\documentclass{gsm-l}

\usepackage{graphicx,amscd,amssymb,array,epsf}
\usepackage{amsthm}
\usepackage{longtable}
\usepackage{color}
\usepackage{mathrsfs}

\input xy
\xyoption{all}
\newtheoremstyle{sec-thm}%
  {}%      Space above (blank means "normal")
  {}%      Space below
  {\itshape}%   Body font
  {}%           Indent amount (empty = no indent, \parindent = para indent)
  {\bfseries}%  Thm head font
  {.}%        Punctuation after thm head
  {.5em}%     Space after thm head
  {}%         Thm head spec (can be left empty, meaning `normal')

\newtheoremstyle{sec-dfn}% name
  {}%      Space above (blank means "normal")
  {}%      Space below
  {}%      Body font
  {}%      Indent amount (empty = no indent, \parindent = para indent)
  {\bfseries}%  Thm head font
  {.}%        Punctuation after thm head
  {.5em}%     Space after thm head
  {}%         Thm head spec (can be left empty, meaning `normal')

\numberwithin{section}{chapter}
\numberwithin{equation}{section}

\numberwithin{figure}{chapter}

\swapnumbers
\theoremstyle{sec-thm}
\newtheorem{theorem}[subsection]{Theorem}
\newtheorem*{xtheorem}{Theorem}
\newtheorem{lemma}[subsection]{Lemma}
\newtheorem*{xlemma}{Lemma}
\newtheorem{proposition}[subsection]{Proposition}
\newtheorem*{xproposition}{Proposition}
\newtheorem{conjecture}[subsection]{Conjecture}
\newtheorem*{xconjecture}{Conjecture}
\newtheorem{corollary}[subsection]{Corollary}
\newtheorem*{xcorollary}{Corollary}

\theoremstyle{sec-dfn}
\newtheorem{definition}[subsection]{Definition}
\newtheorem*{xdefinition}{Definition}
\newtheorem{example}[subsection]{Example}
\newtheorem*{xexample}{Example}
\newtheorem{xca}[subsection]{Exercise}
\newtheorem*{xxca}{Exercise}

\newtheorem{remark}[subsection]{Remark}
\newtheorem*{xremark}{Remark}

\newcommand{\CP}{C} % Conway polynomial
\newcommand{\eps}{\varepsilon}
\newcommand{\MM}{\textrm{MM}}
\newcommand{\isom}{\cong}

\newcommand{\rb}{\raisebox}
\newcommand{\ig}{\includegraphics}

\newcommand{\sm}{\setminus}
\newcommand{\wt}{\widetilde}

\newcommand{\Z}{\mathbb{Z}} % integers
\newcommand{\Q}{\mathbb{Q}} % rationals
\newcommand{\R}{\mathbb{R}} % real
\newcommand{\C}{\mathbb{C}} % complex
\newcommand{\Fi}{\mathbb{F}} % ground field
\def\Ring{\mathcal{R}} % ground ring

\newcommand{\ChD}{\mathbf{A}} % set of chord diagrams
\newcommand{\FD}{\mathbf{C}} % set of closed Jacobi diagrams
\newcommand{\OD}{\mathbf{B}}  % set of open Jacobi diagrams
\newcommand{\GD}{\mathbf{GD}} % set of Gauss diagrams
\newcommand{\A}{\mathcal{A}}  % space of chord diagrams
\newcommand{\F}{\mathcal{C}}  % space of closed Jacobi diagrams
\newcommand{\B}{\mathcal{B}}   % space of open Jacobi diagrams
\newcommand{\Polyak}{\mathcal{P}} % Polyak algebras and such

\newcommand{\Gr}{\mathfrak{G}}  % for the algebra of (all) graphs
\newcommand{\GR}{\mathcal{G}}  % for the set of group-like elements
\newcommand{\AS}{\mathrm{AS}}
\newcommand{\IHX}{\mathrm{IHX}}
\newcommand{\STU}{\mathrm{STU}}
\newcommand{\fT}{\mathrm{4T}}

\newcommand{\I}{\mathcal{I}}
\newcommand{\T}{\mathcal{T}}
\newcommand{\Ha}{\mathcal{H}} % Habiro's modules
\newcommand{\K}{\mathcal{K}} % isotopy classes of knots
\newcommand{\La}{\mathcal{L}} % for Lando's algebra of graphs
\newcommand{\M}{\mathcal{M}} % set of maps
\newcommand{\PR}{\mathcal{P}} % for the primitive subspace
\newcommand{\V}{\mathcal{V}}  % Vassiliev invariants
\newcommand{\X}{\mathcal{X}}  % spaces of fixed diagrams (Vogel)
\newcommand{\g}{\mathfrak{g}}  % for Lie algebra
\newcommand{\gl}{\mathfrak{gl}}
\newcommand{\so}{\mathfrak{so}}
\newcommand{\sP}{\mathfrak{sp}}
\newcommand{\sL}{\mathfrak{sl}} % for the Lie algebra sl
\newcommand{\G}{\Gamma}  % for the algebra of 3-graphs
\renewcommand{\a}{\alpha}
\renewcommand{\b}{\beta}
\newcommand{\e}{\varepsilon}
\newcommand{\f}{\varphi}
\newcommand{\io}{\iota}
\newcommand{\om}{\omega}
\newcommand{\op}{\oplus}
\newcommand{\Op}{\mathop\bigoplus\limits}
\newcommand{\ot}{\otimes}
\renewcommand\t{\theta}

\newcommand{\sz}{\scriptsize}
\newcommand{\ol}{\overline}
\newcommand{\const}{\mathop{\rm const}}
\newcommand{\sign}{\mathop{\rm sign}}
\newcommand{\im}{\mathop{\rm im}}
\newcommand{\id}{\mathrm{id}}
\newcommand{\tr}{\mathop{\rm tr}}
\newcommand{\Hom}{\mathop{\rm Hom}\nolimits}
\newcommand{\End}{\mathop{\rm End}\nolimits}

\newcommand{\Tr}{\mathop{\rm Tr}}

\newcommand{\Di}{\mathcal{D}}
\newcommand\cd[1]{\risS{-10}{#1}{}{24}{20}{15}}
\newcommand{\cdO}{\rb{1pt}{\chd{cd1ch4}}}
\newcommand{\cdWO}{\rb{1pt}{\cd{cd21ch4}}}
\newcommand{\cdWW}{\rb{1pt}{\cd{cd22ch4}}}
\newcommand{\cdTO}{\rb{1pt}{\chd{cd31ch4}}}
\newcommand{\cdTW}{\rb{1pt}{\chd{cd32ch4}}}
\newcommand{\cdTT}{\rb{1pt}{\chd{cd33ch4}}}
\newcommand{\cdTF}{\rb{1pt}{\chd{cd34ch4}}}
\newcommand{\cdTV}{\rb{1pt}{\chd{cd35ch4}}}

\newcommand{\lcdAa}{\rb{-2.5mm}{\ig[width=26mm]{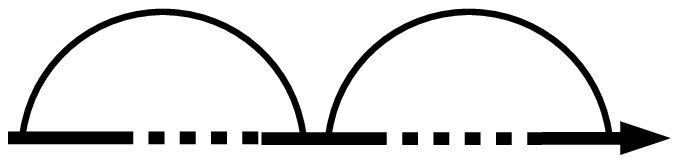}}}
\newcommand{\lcdBa}{\rb{-2.5mm}{\ig[width=26mm]{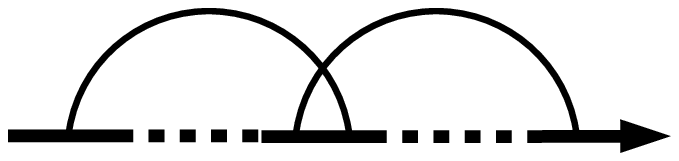}}}
\newcommand{\lcdAb}{\rb{-2.5mm}{\ig[width=26mm]{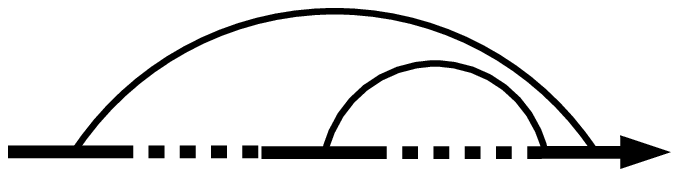}}}
\newcommand{\lcdBb}{\rb{-2.5mm}{\ig[width=26mm]{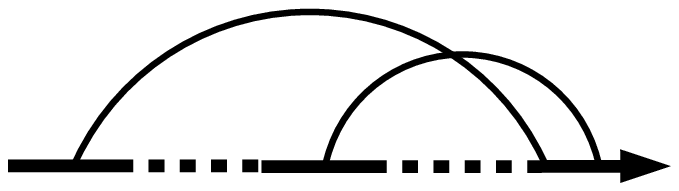}}}
\newcommand{\lcdAc}{\rb{-2.5mm}{\ig[width=26mm]{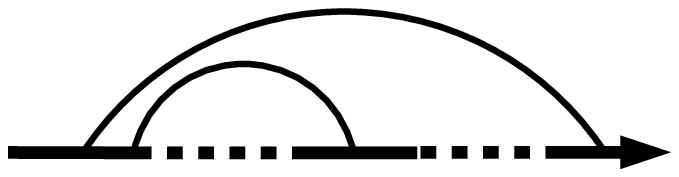}}}
\newcommand{\lcdBc}{\rb{-2.5mm}{\ig[width=26mm]{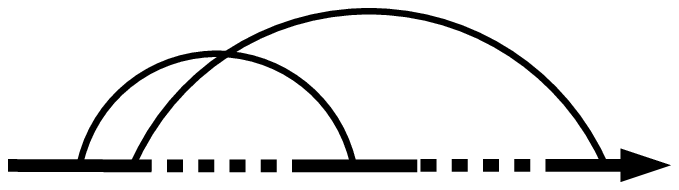}}}

\newcommand{\unkn}{\rb{-4.2mm}{\ig[width=10mm]{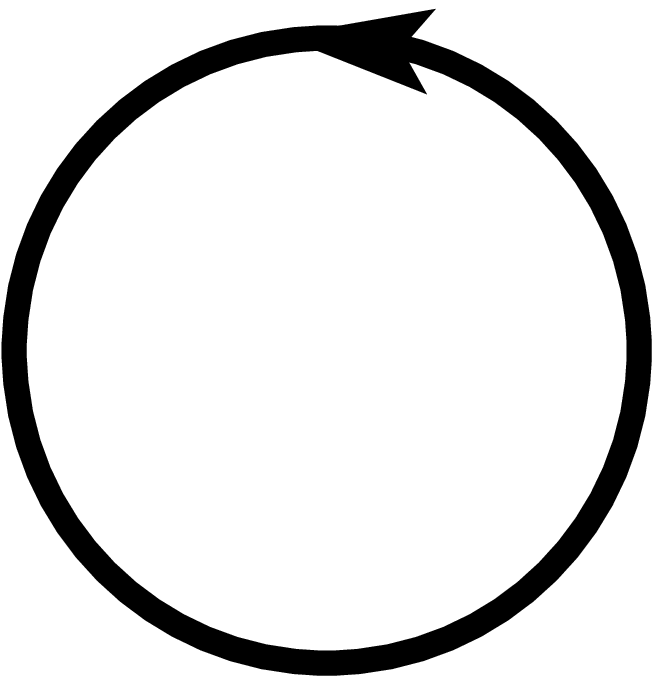}}}
\newcommand{\tallunkn}{\rb{-1.2mm}{\ig[width=10mm]{unkn.eps}}}
\newcommand{\double}{\rb{-4.2mm}{\ig[width=10mm]{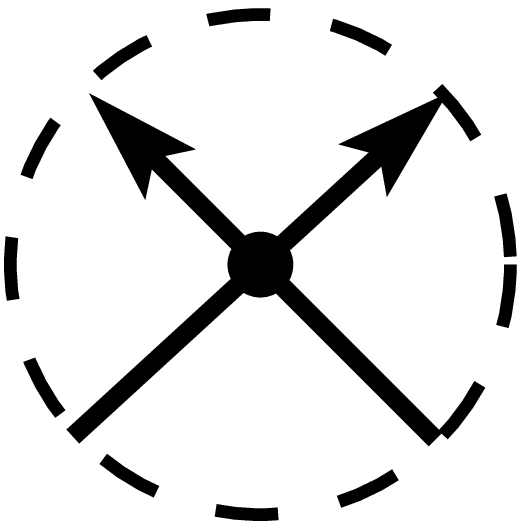}}}
\newcommand{\eight}{\rb{-2.3mm}{\ig[height=6mm]{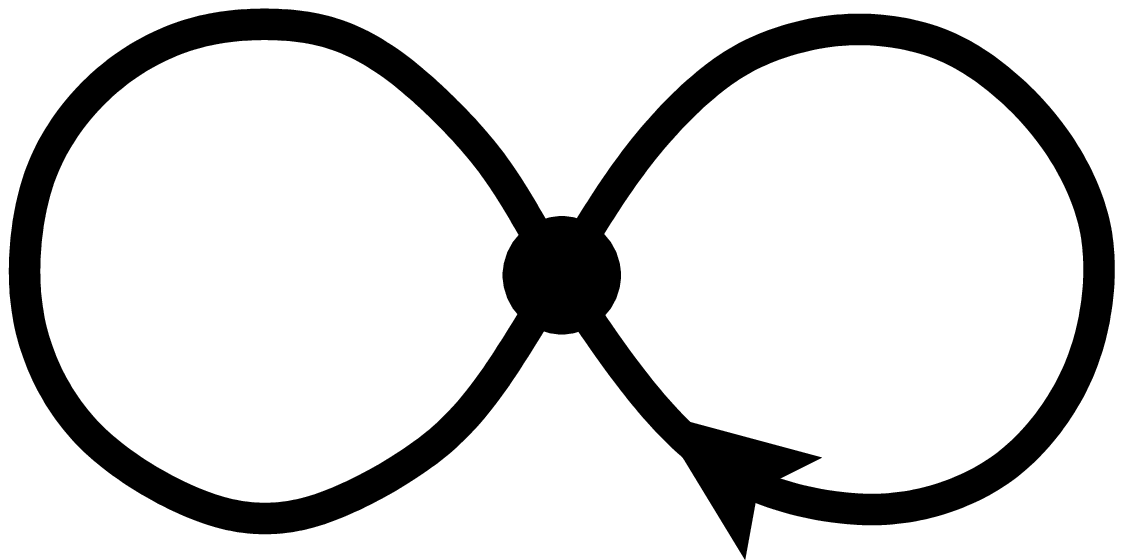}\,}}
\newcommand{\twoup}{\rb{-4.2mm}{\ig[width=10mm]{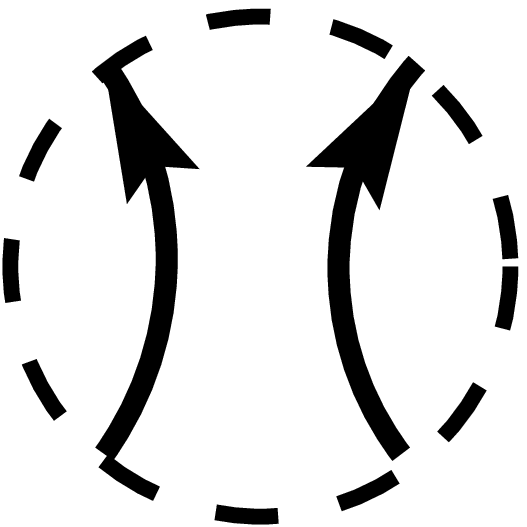}}}
\newcommand{\rlints}{\rb{-4.2mm}{\ig[width=10mm]{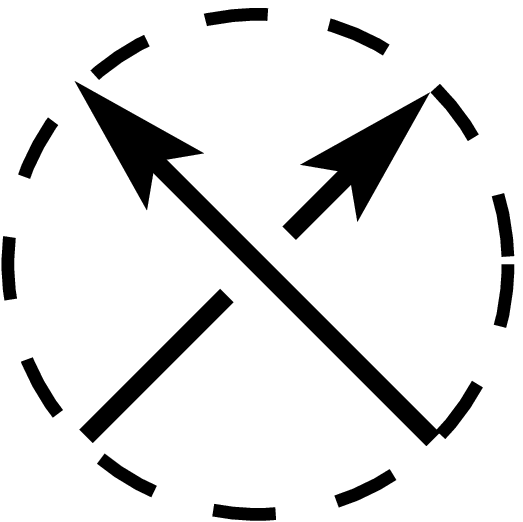}}}
\newcommand{\lrints}{\rb{-4.2mm}{\ig[width=10mm]{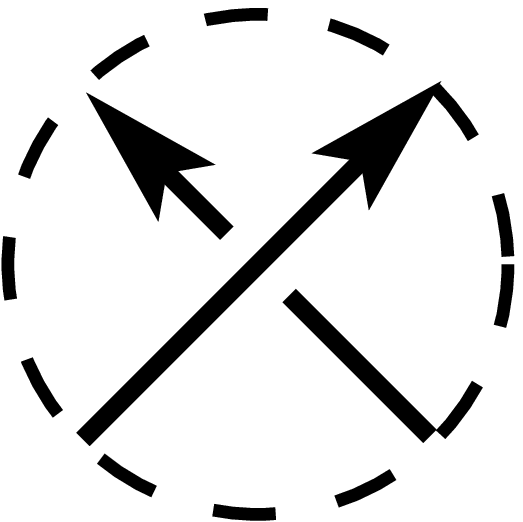}}}
\newcommand{\doublefr}{\rb{-4.2mm}{\ig[width=5mm]{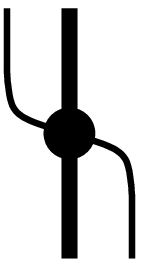}}}
\newcommand{\rlintsfr}{\rb{-4.2mm}{\ig[width=5mm]{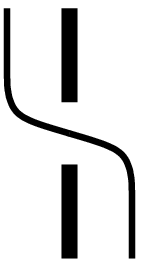}}}
\newcommand{\lrintsfr}{\rb{-4.2mm}{\ig[width=5mm]{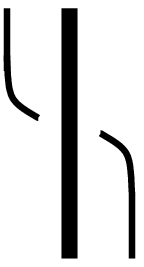}}}
\newcommand{\frameddbl}{\rb{-0.2mm}{\ig[width=15mm]{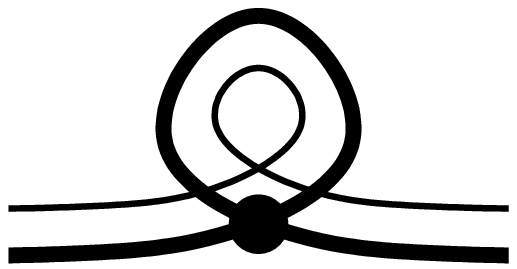}}}
\newcommand{\framedone}{\rb{-1.2mm}{\ig[width=15mm]{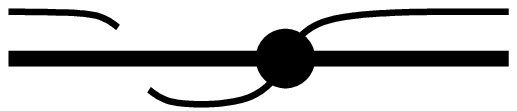}}}
\newcommand{\framedtwo}{\rb{-1.2mm}{\ig[width=15mm]{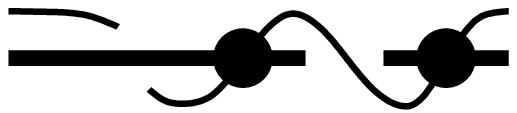}}}

\newcommand{\krestik}{\rb{-1pt}{\ig[width=8pt]{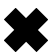}}}

\newcommand{\ccO}[1]{\rb{0.5mm}{\ig[width=10mm]{#1.eps}}}
\newcommand{\rccO}[2]{\rb{#1mm}{\ig[width=10mm]{#2.eps}}}

\newcommand\pil[1]{\rb{-15pt}{\begin{picture}(45,40)(0,0)
      \put(0,0){\ig[width=45pt]{#1.eps}}
     \end{picture}}}
\newcommand{\bub}{
  \begin{picture}(10,10)(-2,2)
      \put(5,5){\circle{10}}
      \put(0,5){\line(1,0){10}}
      \put(0,5){\circle*{2}}
      \put(10,5){\circle*{2}}
  \end{picture}}
\newcommand\risS[6]{\rb{#1pt}[#5pt][#6pt]{\begin{picture}(#4,15)(0,0)
  \put(0,0){\ig[width=#4pt]{#2.eps}} #3
     \end{picture}}}
\newcommand\expic[1]{\risS{-20}{#1}{}{60}{20}{20}}
\newcommand\tanG[4]{\risS{#1}{#2}{#3}{#4}{0}{0}}

\newcommand\ris[7]{\rb{#1pt}[#5pt][#6pt]{
    \begin{picture}(#4,#3)(0,0)
    \put(0,0){\ig[width=#4pt]{#2.eps}} #7
    \end{picture}}}
\newcommand\pmris[1]{\ris{-12}{#1}{30}{27}{20}{16}{}}
\newcommand\bmris[1]{\ris{-12}{#1}{30}{28}{20}{16}{}}
\newcommand\vmris[1]{\ris{-12}{#1}{30}{27}{20}{16}{}}
\newcommand\vbmris[1]{\ris{-11}{#1}{30}{28}{20}{16}{}}
\newcommand\mris[1]{\ris{-12}{#1}{30}{26}{20}{16}{} \hspace{-1pt}}

\def\CAB{\C \langle\!\langle A,B \rangle\!\rangle}
\def\la{\langle\!\langle}
\def\ra{\rangle\!\rangle}
\def\PhiKZ{\Phi_{\mbox{\scriptsize KZ}}}
\def\PhiBN{\Phi_{\mbox{\scriptsize BN}}}
\def\bA{A}    %  was: \def\bA{\overline{A}}
\def\bB{B}    %  was: \def\bB{\overline{B}}

\def\tucpic#1{\risS{-2}{#1}{}{20}{10}{3}}
\def\vstcod{\tucpic{vst4od}}
\def\vstcot{\tucpic{vst4ot}}
\def\vstcdt{\tucpic{vst4dt}}
\def\vstctc{\tucpic{vst4tc}}
\def\vstcdc{\tucpic{vst4dc}}

\def\D{\Delta}
\def\cD{\mathcal{D}}

\def\tupic#1{\risS{-3}{#1}{}{15}{10}{3}}
\def\vstod{\tupic{vstod}}
\def\vstdt{\tupic{vstdt}}
\def\vstemp{\tupic{vstemp}}
\def\vstot{\tupic{vstot}}
\def\vstodod{\tupic{vstod2}}

\def\vstdtdt{\tupic{vstdt2}}

\def\vstdtod{\tupic{vstdtod}}
\def\vstoddt{\tupic{vstoddt}}

\def\o{\omega}

\def\d{\delta}
\def\Ab{\widehat{\A}}
\def\bbp{\mbox{\bf p}}
\def\bbi{\mbox{\bf i}}
\def\bbj{\mbox{\bf j}}
\def\bbq{\mbox{\bf q}}
\def\bbr{\mbox{\bf r}}
\def\bbs{\mbox{\bf s}}
\def\Li{{\rm Li}_2}

\def\taud#1#2{\tau_{#1}\Bigl(#2\Bigr)}

\newcommand{\W}{\mathcal{W}}
\def\bo{\mbox{\bf I}} % delo mastera Bo
\newcommand\figsl[1]{\mbox{\ \begin{picture}(30,20)(0,0)
          \put(0,0){\ig[width=30pt]{#1.eps}}
          \end{picture}\ }}
\newcommand\ad{\mbox{ad}}
\newcommand\Ad{\mbox{Ad}}
\newcommand\symb{\mathrm{symb}}

\long\def\omitfragment#1{}  % to omit portions of text without having to
\newcommand{\chd}[1]{\risS{-10}{#1}{}{25}{20}{15}}

\newcommand{\ccOy}[1]{\rb{0.5mm}{\ig[width=7mm]{#1.eps}}}

\def\clhl{\risS{-.5}{clhl}{}{12}{0}{0}}

\def\sClDx#1{\risS{-8}{#1}{}{8}{5}{4}}
\def\Onechord{\sClDx{onechord}}

\def\sClD#1{\risS{-12}{#1}{}{30}{20}{15}}
\newcommand{\chdh}[1]{\risS{-10}{#1}{}{27}{20}{15}}
\newcommand{\chdn}[2]{\risS{#1}{#2}{}{29}{20}{15}}
\newcommand{\chdr}[2]{\risS{#1}{#2}{}{27}{20}{15}}

\def\boldX{\boldsymbol{X}}
\def\boldY{\boldsymbol{Y}}
\def\xx{\boldsymbol{x}}
\def\yy{\boldsymbol{y}}

\def\SS{\mathcal{S}}
\def\Links{\mathcal{L}}
\def\FLinks{\mathcal{FL}}

\def\maxr{\stackrel{\displaystyle\longrightarrow}{\max}\:}
\def\maxl{\stackrel{\displaystyle\longleftarrow}{\max}\:}
\def\minr{\rb{-6pt}{$\stackrel{\displaystyle\min}{\longrightarrow}$}}
\def\minl{\rb{-6pt}{$\stackrel{\displaystyle\min}{\longleftarrow}$}}
\def\smaxr{\stackrel{\scriptstyle\longrightarrow}{\max}\:}
\def\smaxl{\stackrel{\scriptstyle\longleftarrow}{\max}\:}
\def\sminr{\rb{-5pt}{$\stackrel{\min}{\scriptstyle\longrightarrow}$}}
\def\sminl{\rb{-5pt}{$\stackrel{\min}{\scriptstyle\longleftarrow}$}}
\def\Fxy{\F(\xx\, |\, \yy)}
\def\Bbig{\B^{\circ}}
\def\partbig{\partial^{\circ}}

\def\hl{\risS{-2}{lhl}{}{7}{0}{0}}
\def\dhl{\risS{-2}{ldhl}{}{7}{0}{0}}
\def\Zed{\mathcal{Z}}
\def\Iu#1{\risS{-10}{#1}{}{10.7}{20}{15}}
\def\Id#1{\risS{-10}{#1}{}{18.5}{20}{15}}
\def\eohst{\exp\bigl(\textstyle \frac{1}{2}\,\risS{1.5}{strut}{}{15}{10}{10}\bigr)}
\def\exohy#1#2{\exp\bigl(\textstyle \frac{1}{2}\,\bigl|_{#1}^{#2}\bigr)}
\def\yu#1{\risS{-10}{#1}{}{17}{20}{15}}

\def\neresh{\makebox(0,8){$\hspace*{-7pt}^*$}} % unsolved problem
\def\knt#1{\ig[width=14.5mm,height=14.5mm]{#1.eps}} %knot in table

\def\z{\zeta}
\newcommand{\onept}[1]{{#1}^{\bullet}}
\def\gL_1_1{{\mathfrak{gl}}(1|1)}

\def\kr{{\mathfrak{C}}}
\def\ard#1{\risS{-12}{#1}{}{25}{15}{17}}
\def\plwa{\langle A|S \rangle}
\def\scp#1#2{\langle #1,#2 \rangle}
\newcommand{\gaussAA}{\rb{-3mm}{\ig[width=20mm]{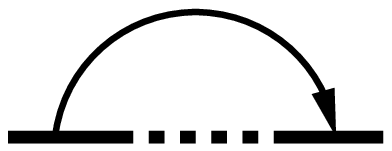}}}
\newcommand{\gaussAC}{\rb{-3mm}{\ig[width=20mm]{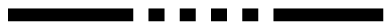}}}

\newcommand{\gaussBA}{\rb{-3mm}{\ig[height=7mm]{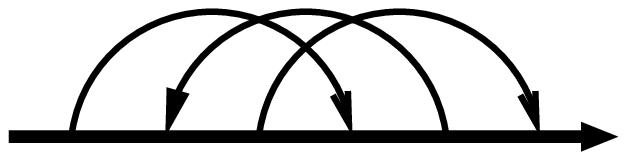}}}
\newcommand{\gaussBC}{\rb{-3mm}{\ig[height=7mm]{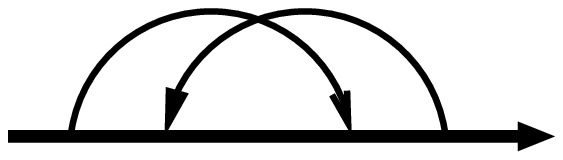}}}
\newcommand{\gaussBD}{\rb{-3mm}{\ig[height=7mm]{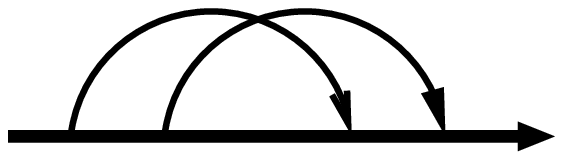}}}
\newcommand{\gaussBE}{\rb{-3mm}{\ig[height=7mm]{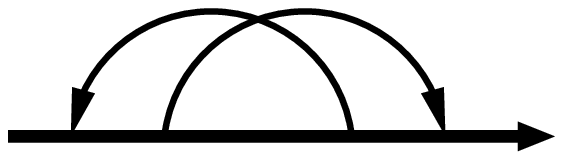}}}
\newcommand{\gaussBF}{\rb{-3mm}{\ig[height=7mm]{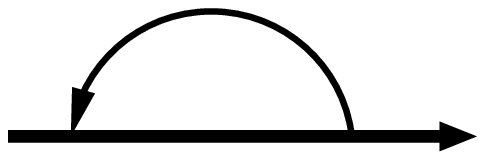}}}
\newcommand{\gaussBG}{\rb{-3mm}{\ig[height=7mm]{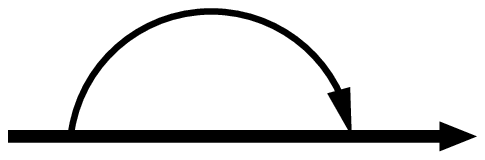}}}
\newcommand{\gaussBH}{\rb{-3mm}{\ig[width=9mm]{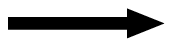}}}

\newcommand{\gaussCA}{\rb{-6mm}{\ig[height=13mm]{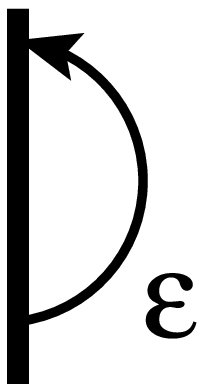}}}
\newcommand{\gaussCB}{\rb{-4mm}{\ig[height=9mm]{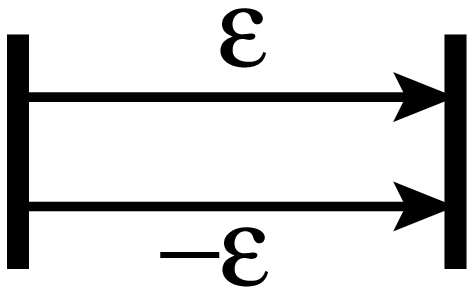}}}
\newcommand{\gaussCC}{\rb{-4mm}{\ig[height=9mm]{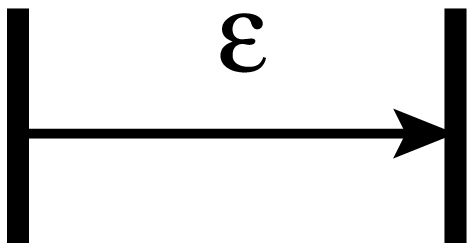}}}
\newcommand{\gaussCD}{\rb{-4mm}{\ig[height=9mm]{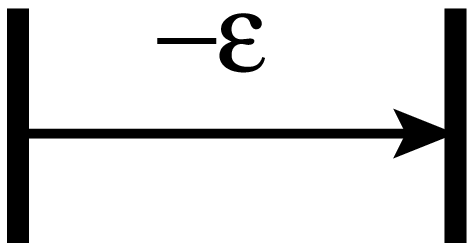}}}
\newcommand{\gaussCE}{\rb{-6mm}{\ig[height=13mm]{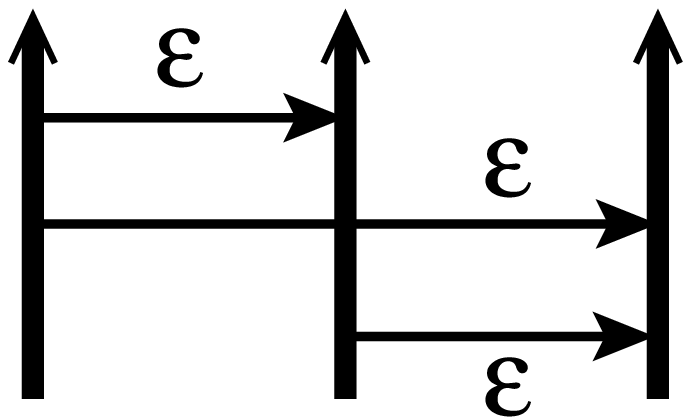}}}
\newcommand{\gaussCF}{\rb{-6mm}{\ig[height=13mm]{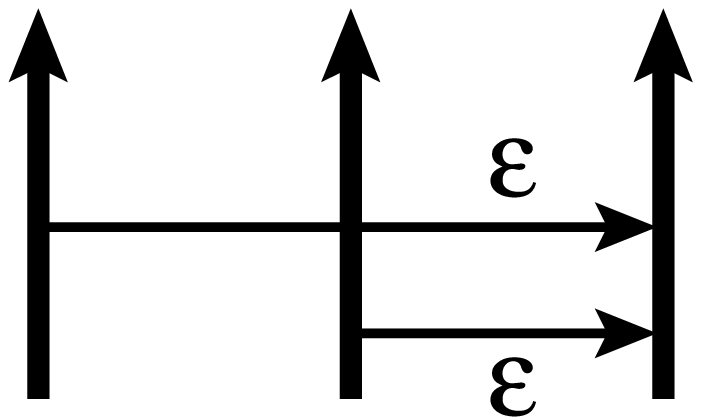}}}
\newcommand{\gaussCG}{\rb{-6mm}{\ig[height=13mm]{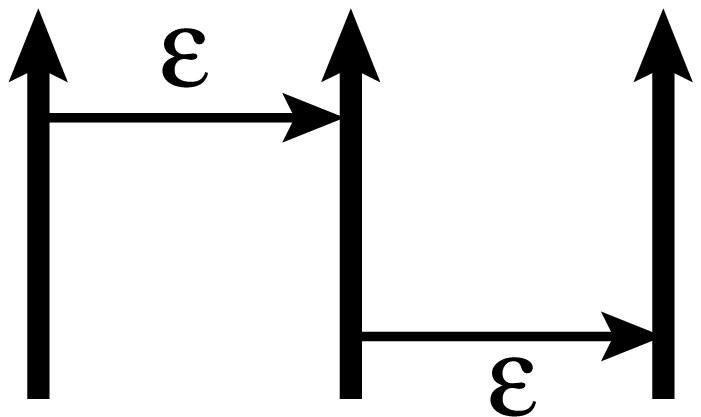}}}
\newcommand{\gaussCH}{\rb{-6mm}{\ig[height=13mm]{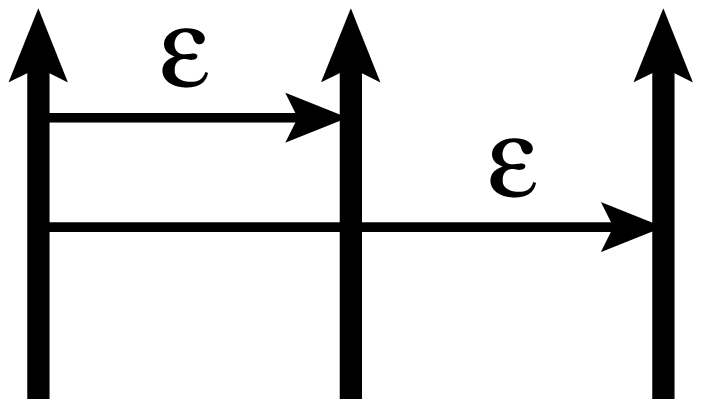}}}
\newcommand{\gaussCI}{\rb{-6mm}{\ig[height=13mm]{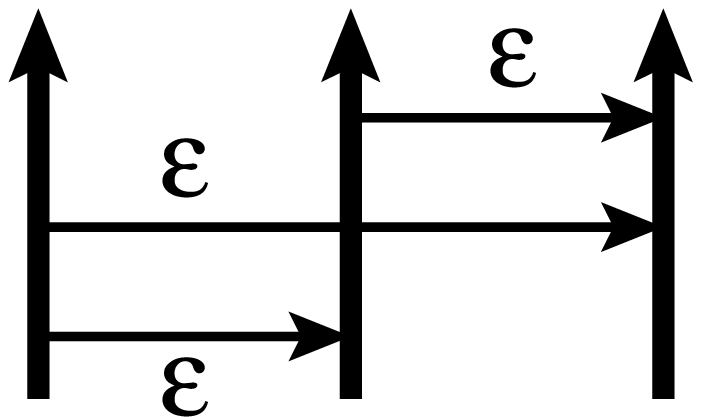}}}
\newcommand{\gaussCJ}{\rb{-6mm}{\ig[height=13mm]{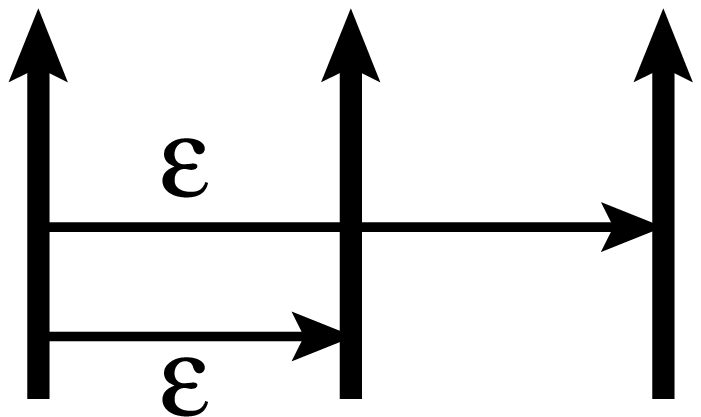}}}
\newcommand{\gaussCK}{\rb{-6mm}{\ig[height=13mm]{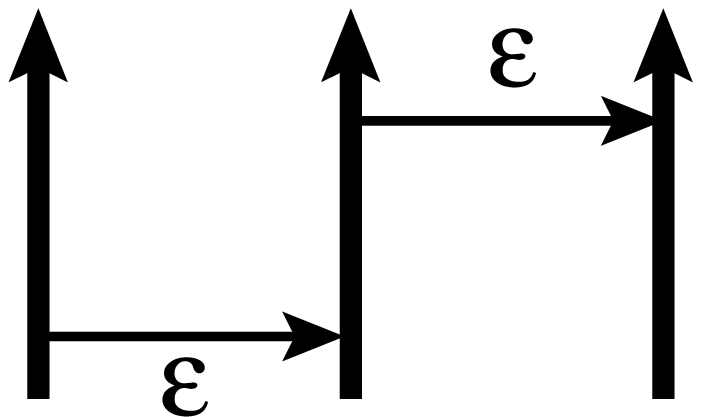}}}
\newcommand{\gaussCL}{\rb{-6mm}{\ig[height=13mm]{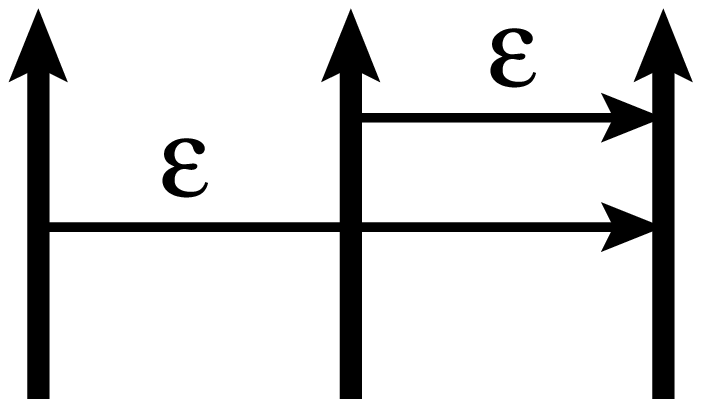}}}

\newcommand{\gaussDA}{\rb{-3mm}{\ig[width=20mm]{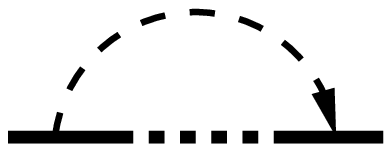}}}
\newcommand{\gaussDB}{\rb{-3mm}{\ig[width=20mm]{gauss-1A.eps}}}
\newcommand{\gaussDC}{\rb{-3mm}{\ig[width=20mm]{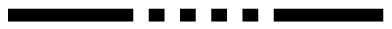}}}

\newcommand{\gaussEA}{\rb{-5mm}{\ig[height=11mm]{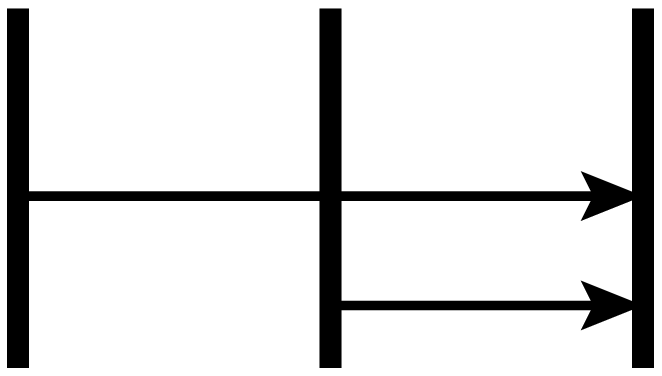}}}
\newcommand{\gaussEB}{\rb{-5mm}{\ig[height=11mm]{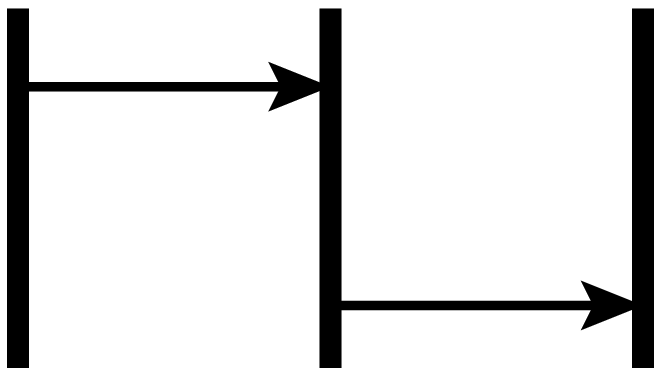}}}
\newcommand{\gaussEC}{\rb{-5mm}{\ig[height=11mm]{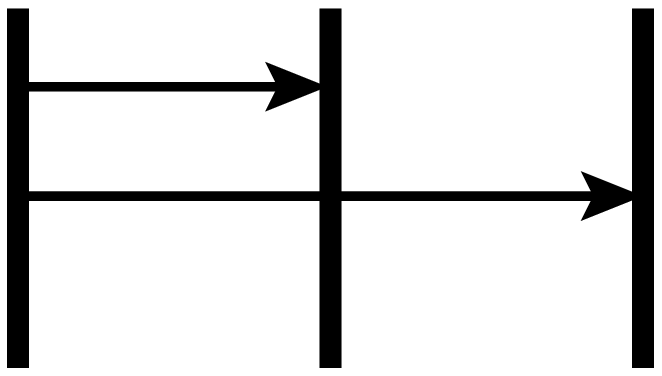}}}
\newcommand{\gaussED}{\rb{-5mm}{\ig[height=11mm]{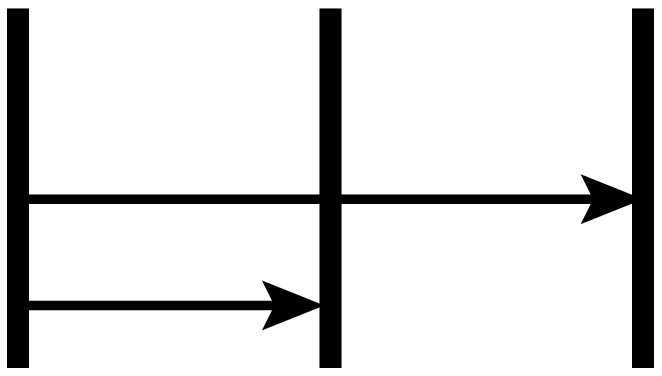}}}
\newcommand{\gaussEE}{\rb{-5mm}{\ig[height=11mm]{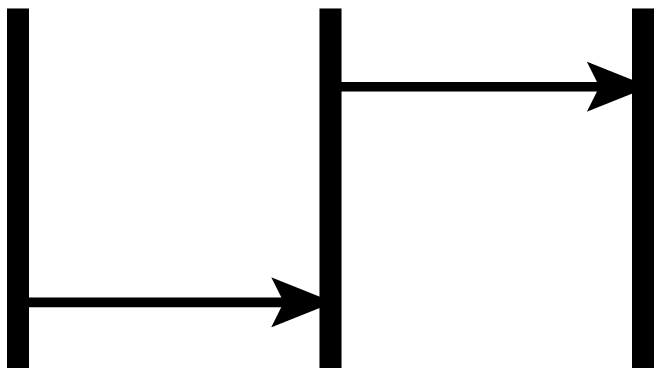}}}
\newcommand{\gaussEF}{\rb{-5mm}{\ig[height=11mm]{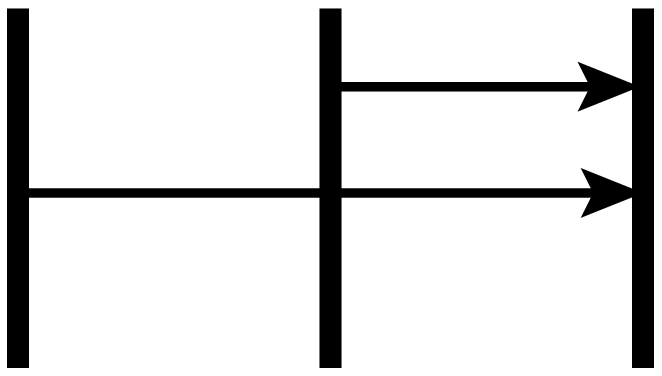}}}

\newcommand{\gaussJA}{\rb{-2mm}{\ig[width=27mm]{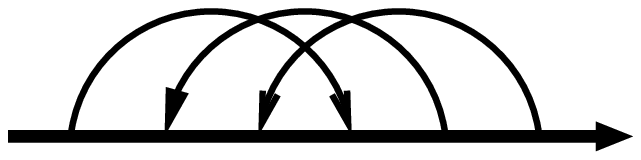}}}
\newcommand{\gaussJB}{\rb{-2mm}{\ig[width=27mm]{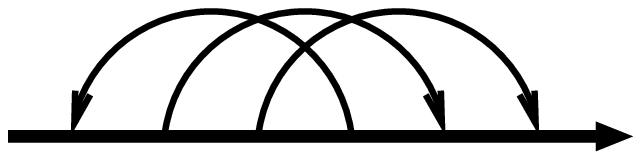}}}
\newcommand{\gaussJC}{\rb{-2mm}{\ig[width=27mm]{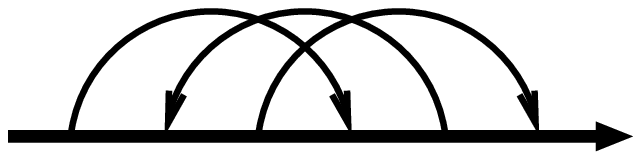}}}
\newcommand{\gaussJD}{\rb{-2mm}{\ig[width=27mm]{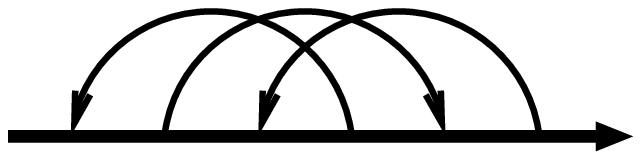}}}
\newcommand{\gaussJE}{\rb{-2mm}{\ig[width=27mm]{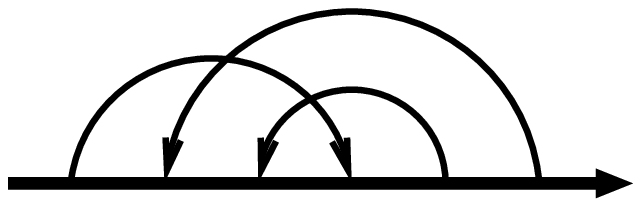}}}
\newcommand{\gaussJF}{\rb{-2mm}{\ig[width=27mm]{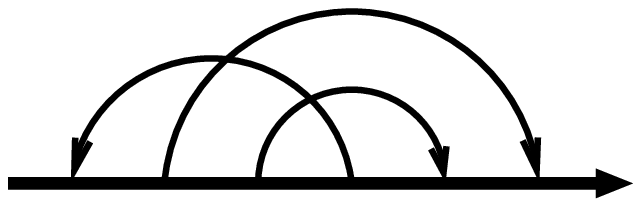}}}
\newcommand{\gaussJG}{\rb{-2mm}{\ig[width=27mm]{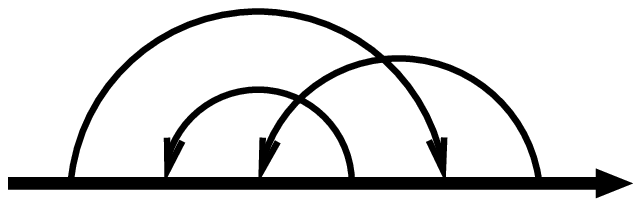}}}
\newcommand{\gaussJH}{\rb{-2mm}{\ig[width=27mm]{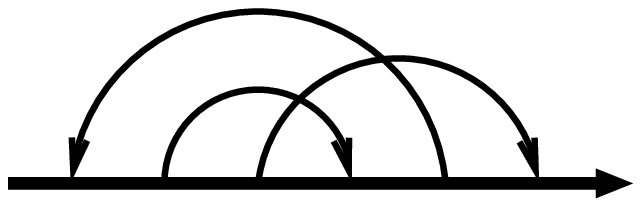}}}
\newcommand{\gaussJI}{\rb{-2.5mm}{\ig[width=27mm]{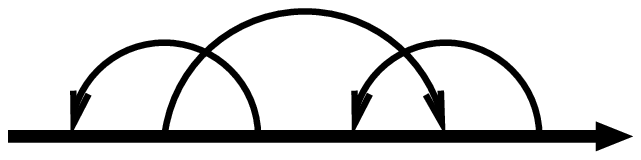}}}
\newcommand{\gaussJJ}{\rb{-2.5mm}{\ig[width=27mm]{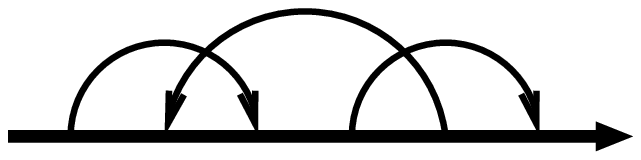}}}
\newcommand{\gaussJK}{\rb{-2.5mm}{\ig[width=21mm]{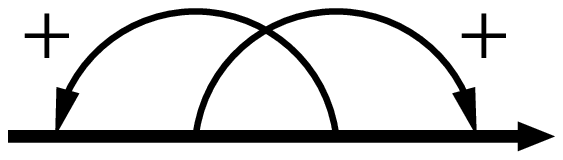}}}
\newcommand{\gaussJL}{\rb{-2.5mm}{\ig[width=21mm]{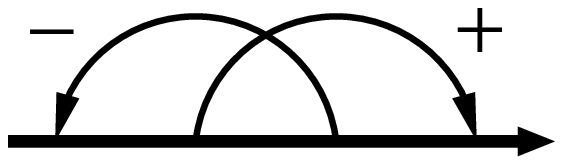}}}

\newcommand{\gaussKA}{\rb{-2mm}{\ig[width=24mm]{gauss-13A.eps}}}
\newcommand{\gaussKB}{\rb{-2mm}{\ig[width=24mm]{gauss-13B.eps}}}
\newcommand{\gaussKC}{\rb{-2mm}{\ig[width=24mm]{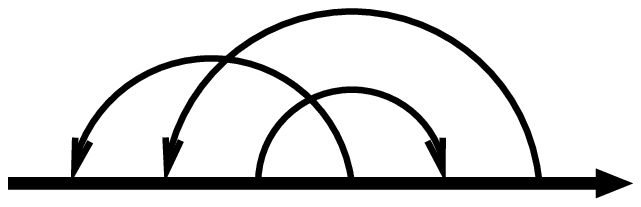}}}
\newcommand{\gaussKD}{\rb{-2mm}{\ig[width=24mm]{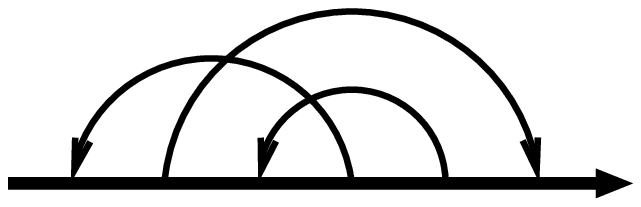}}}
\newcommand{\gaussKE}{\rb{-2mm}{\ig[width=24mm]{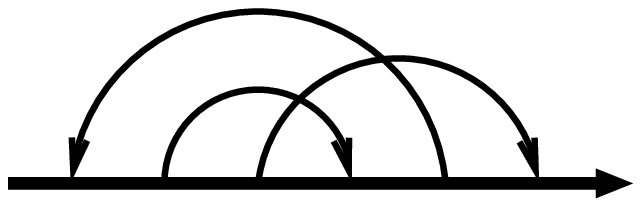}}}
\newcommand{\gaussKF}{\rb{-2mm}{\ig[width=24mm]{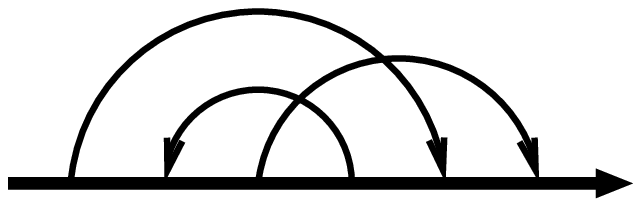}}}
\newcommand{\gaussKG}{\rb{-2mm}{\ig[width=24mm]{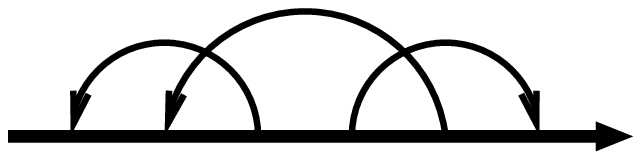}}}
\newcommand{\gaussKH}{\rb{-2mm}{\ig[width=24mm]{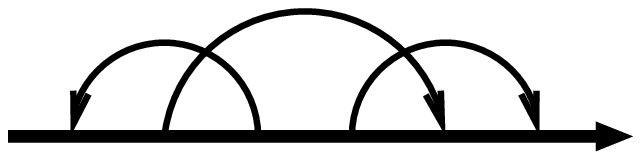}}}
\newcommand{\gaussKI}{\rb{-2mm}{\ig[width=24mm]{gauss-13I.eps}}}
\newcommand{\gaussKJ}{\rb{-2mm}{\ig[width=24mm]{gauss-13J.eps}}}
\newcommand{\gaussKK}{\rb{-2mm}{\ig[width=24mm]{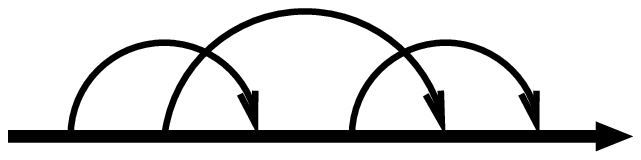}}}

\newenvironment{dedication}
{
   \cleardoublepage
   \thispagestyle{empty}
   \vspace*{\stretch{1}}
   \hfill\begin{minipage}[t]{0.66\textwidth}
   \raggedright
}%
{
   \end{minipage}
   \vspace*{\stretch{3}}
   \clearpage
}

\author{
September 17, 2011 
}

\author{S.~Chmutov}
\address{The Ohio State University, Mansfield Campus,
1680 University Drive, Mansfield, OH 44906, USA}
\email{chmutov@math.ohio-state.edu}
\author{S.~Duzhin}
\address{Steklov Institute of Mathematics, St.~Petersburg Division,
Fontanka 27, St.~Petersburg, 191023, Russia}
\email{duzhin@pdmi.ras.ru}
\author{J.~Mostovoy}
\address{
Departamento de Matem\'aticas, CINVESTAV, Apartado Postal 14-740,
C.P.\ 07000 M\'exico, D.F. Mexico} \email{jacob@math.cinvestav.mx}

\title{Introduction to Vassiliev Knot Invariants}

\makeindex
\begin{document}
\maketitle

\begin{dedication}
{\sl To the memory of V.\ I.\ Arnold}
\end{dedication}

\tableofcontents
\numberwithin{equation}{subsection}
\numberwithin{table}{subsection}
\numberwithin{figure}{subsection}
\renewcommand{\leq}{\leqslant}
\renewcommand{\geq}{\geqslant}
\renewcommand{\le}{\leqslant}
\renewcommand{\ge}{\geqslant}

\section*{Preface}

This book is a detailed introduction to the theory of finite type
(Vassiliev) knot invariants, with a stress on its combinatorial
aspects. It is intended to serve both as a textbook for readers with
no or little background in this area, and as a guide to some of the
more advanced material. Our aim is to lead the reader to
understanding by means of pictures and calculations, and for this
reason we often prefer to convey the idea of the proof on an
instructive example rather than give a complete argument. While we
have made an effort to make the text reasonably self-contained, an
advanced reader is sometimes referred to the original papers for the
technical details of the proofs.

\subsection*{Historical remarks}

The notion of a finite type knot invariant was introduced by Victor
Vassiliev (Moscow) in the end of the 1980's and first appeared in print in
his paper \cite{Va1} (1990).  Vassiliev, at the time, was not specifically
interested in low-dimensional topology.  His main concern was the general
theory of discriminants in the spaces of smooth maps, and his description of
the space of knots was just one, though the most spectacular, application of
a machinery that worked in many seemingly unrelated contexts.  It was
V.~I.~Arnold \cite{Ar2} who understood the importance of finite type
invariants, coined the name ``Vassiliev invariants'' and popularized the
concept; since that time, the term ``Vassiliev invariants'' has become
standard.

A different perspective on the finite type invariants was developed by
Mikhail Goussarov (St.\,Petersburg).  His notion of $n$-equivalence, which
first appeared in print in \cite{G2} (1993), turned out to be useful in
different situations, for example, in the study of the finite type
invariants of 3-manifolds.\footnote{Goussarov cites Vassiliev's works in his
earliest paper \cite{G1}.  Nevertheless, according to O.\ Viro, Goussarov
first mentioned finite type invariants in a talk at the Leningrad
topological seminar as early as in 1987.} Nowadays some people use the
expression ``Vassiliev-Goussarov invariants'' for the finite type
invariants.

Vassiliev's definition of finite type invariants is based on the
observation that knots form a topological space and the knot
invariants can be thought of as the locally constant functions on
this space. Indeed, the space of knots is an open subspace of the
space $M$ of all smooth maps from $S^1$ to $\R^3$; its complement is
the so-called discriminant $\Sigma$ which consists of all maps that
fail to be embeddings.  Two knots are isotopic if and only if they
can be connected in $M$ by a path that does not cross $\Sigma$.

Using simplical resolutions, Vassiliev constructs a spectral
sequence for the homology of $\Sigma$. After applying the Alexander
duality this spectral sequence produces cohomology classes for the
space of knots $M-\Sigma$; in dimension zero these are precisely the
Vassiliev knot invariants.

Vassiliev's approach, which is technically rather demanding, was
simplified by J.\ Birman and X.-S.\ Lin in \cite{BL}.
They explained the relation between the Jones polynomial and finite
type invariants\footnote{independently from Goussarov who was the first to 
discover this relation in \cite{G1}.} and emphasized the role of the algebra of chord
diagrams. M.~Kontsevich showed that the study of real-valued
Vassiliev invariants can, in fact, be reduced entirely to the
combinatorics of chord diagrams \cite{Kon1}. His proof used an
analytic tool ({the Kontsevich integral}) which may be thought of as
a powerful generalization of the Gauss linking formula. Kontsevich also 
defined a coproduct on the algebra of chord diagrams which turns it into 
a Hopf algebra. 

D.~Bar-Natan was the first to give a comprehensive treatment of
Vassiliev knot and link invariants. In his preprint \cite{BN0} and
PhD thesis \cite{BNt} he found the relationship between finite type
invariants and the topological quantum field theory developed by
his thesis advisor E.~Witten \cite{Wit1,Wit2}\index{Witten}.
Bar-Natan's paper \cite{BN1} (whose preprint edition \cite{BN1a}
appeared in 1992) is still the most authoritative source on the
fundamentals of the theory of Vassiliev invariants. About the same
time,  T.~Le and J.~Murakami \cite{LM2}, relying on V.~Drinfeld's
work \cite{Dr1,Dr2}, proved the rationality of the Kontsevich
integral.

Among further developments in the area of finite type knot
invariants let us mention:
\begin{itemize}
\item
 The existence of non-Lie-algebraic weight systems
(P.~Vogel \cite{Vo1}, J.~Lie\-berum \cite{Lieb}) and an interpretation of all weight systems as Lie algebraic weight systems in a suitable category (V.~Hinich, A.~Vaintrob \cite{HV});
\item
J.~Kneissler's analysis \cite{Kn1,Kn2,Kn3} of the structure of the algebra $\Lambda$ introduced by P.~Vogel \cite{Vo1};
\item
The proof by Goussarov \cite{G5} that Vassiliev invariants are polynomials in the gleams for a fixed Turaev shadow;
\item
Gauss diagram formulae of M.~Polyak and O.~Viro \cite{PV1}
and the proof by M.~Goussarov \cite{G3} that all finite type invariants
can be expressed by such formulae;
\item
Habiro's theory of claspers \cite{Ha2} (see also \cite{G4});
\item
V.~Vassiliev's papers \cite{Va4,Va5} where a general technique
for deriving combinatorial formulae for cohomology classes in the
complements to discriminants, and in particular, for finite type
invariants, is proposed;
\item
Explicit formulae for the Kontsevich integral of some knots and
links \cite{BLT,BNL,Roz2,Kri2,Mar,GK};
\item
The interpretation of the Vassiliev spectral sequence in terms of the Hochschild homology
of the Poisson operad by V.~Turchin \cite{Tu1};
\item
The alternative approaches to the topology of the space of knots via configuration spaces and the Goodwillie calculus \cite{Sinha}.
\end{itemize}

One serious omission in this book is the connection between the Vassiliev
invariants and the Chern-Simons theory.\index{Chern-Simons theory} This
connection motivates much of the interest in finite-type invariants and
gives better understanding of the nature of the Kontsevich integral. 
Moreover, it suggests another form of the universal Vassiliev invariant,
namely, the configuration space integral.  There are many texts that explain
this connection with great clarity; the reader may start, for instance, with
\cite{Lab}, \cite{Saw} or \cite{Po2}.  The original paper of Witten
\cite{Wit1} has not lost its relevance and, while it does not deal directly
with the Vassiliev invariants (it dates from 1989), it still is one of the
indispensable references.

An important source of 
information on finite type invariants is the online {\em Bibliography of Vassiliev invariants}
started by D.~Bar-Natan and currently living at
\begin{center}
\verb#http://www.pdmi.ras.ru/~duzhin/VasBib/#
\end{center}
In March, 2011 it contained 644 entries, and this number is increasing. The study of finite type invariants is going on. However,
notwithstanding all efforts, the most important question put forward in
1990:
\begin{quote}
   \textit{Is it true that Vassiliev invariants distinguish knots?}
\end{quote}
is still open. At the moment is is not even known whether the Vassiliev invariants can detect knot orientation. A number of open problems related to finite-type invariants are listed in \cite{Oht3}.

\subsection*{Prerequisites}

We assume that the reader has a basic knowledge of calculus on manifolds
(vector fields, differential forms, Stokes' theorem), general algebra
(groups, rings, modules, Lie algebras, fundamentals of homological algebra),
linear algebra (vector spaces, linear operators, tensor algebra, elementary
facts about representations) and topology (topological spaces, homotopy, homology, 
Euler characteristic). Some of this and more advanced algebraic material (bialgebras, free algebras, universal enveloping algebras etc.) which is of primary importance
in this book, can be found in the Appendix at the end of the book.
No knowledge of knot theory is presupposed, although it may be useful.

\subsection*{Contents}

The book consists of fifteen chapters, which can logically be divided into four parts.

\begin{figure}[ht]
{\bf Chapter dependence}\\
\ \\
$$
\risS{-188}{chapt_depend}{}{270}{150}{188}
$$
\end{figure}

The first part opens with a short introduction into the theory
of knots and their classical polynomial invariants and closes with
the definition of Vassiliev invariants.

In part 2, we systematically study the graded Hopf algebra naturally
associated with the filtered space of Vassiliev invariants,
which appears in three different guises:
as the algebra of chord diagrams $\A^{fr}$,
as the algebra of closed Jacobi diagrams $\F$, and
as the algebra of open Jacobi diagrams $\B$.
After that, we study the auxiliary algebra $\G$ generated by
regular trivalent graphs and closely related to the algebras $\A$, $\B$,
$\mathcal{C}$ as well as to Vogel's algebra $\Lambda$. In the last chapter
we discuss the weight systems defined by Lie algebras, both universal and
depending on a chosen representation.

Part 3 is dedicated to a detailed exposition of the Kontsevich integral;
it contains the proof of the main theorem of the theory of Vassiliev knot
invariants that reduces their study to combinatorics of chord diagrams
and related algebras. Chapters~\ref{chapKI} and \ref{chap:operations} treat 
the Kontsevich integral from the analytic point of view. Chapter~\ref{DA_chap}
is dedicated to the Drinfeld associator and the combinatorial construction of the 
Kontsevich integral. Chapter \ref{advKI} contains some additional
material on the Kontsevich integral: the wheels formula,
the Rozansky rationality conjecture etc.

The last part of the book is devoted to various topics left out in the
previous exposition, such as the Vassiliev invariants for braids, Gauss diagram
formulae, the Melvin--Morton conjecture, the Goussarov--Habiro theory,
the size of the space of Vassiliev invariants etc. The book closes with a description 
of Vassiliev's original construction for the finite type invariants.

The book is intended to be a textbook, so we have included many exercises. Some exercises
are embedded in the text; the others appear in a separate section
at the end of each chapter. Open problems are marked with an asterisk.

\subsection*{Acknowledgements}

The work of the first two authors on this book actually began in August
1992, when our colleague Inna Scherbak returned to Pereslavl-Zalessky from
the First European Mathematical Congress in Paris and brought a photocopy of
Arnold's lecture notes about the newborn theory of Vassiliev knot
invariants. We spent several months filling our waste-paper baskets with
pictures of chord diagrams, before the first joint article \cite{CD1} was
ready.

In the preparation of the present text, we have extensively used our
papers (joint, single-authored and with other coauthors, see
bibliography) and in particular, lecture notes of the course
``Vassiliev invariants and combinatorial structures'' that one of us
(S.~D.) delivered at the Graduate School of Mathematics, University
of Tokyo, in Spring 1999. It is our pleasure to thank V.~I.~Arnold,
D.~Bar-Natan, J.~Birman, C.~De~Concini, O.~Dasbach, A.~Durfee, F.~Duzhin,
V.~Goryunov, O.~Karpenkov, T.~Kerler, T.~Kohno, S.~Lando, M.~Polyak,
I.~Scherbak, A.~Sossinsky, A.~Vaintrob, A.~Varchenko, V.~Vassiliev,
and S.~Willerton
for many useful comments concerning the subjects touched upon in the book.
We are likewise indebted to the anonymous referees whose criticism and
suggestions helped us to improve the text.

Our work was supported by several grants:
INTAS 00-0259, NWO 047.008.005, NSh-709.2008.1, RFFI-05-01-01012
and 08-01-00379 (S.~D.),
Professional Development Grants of OSU, Mansfield (2002, 2004, S.~Ch.),
CONACyT CO2-44100 (J.~M.)
Part of the work was accomplished when the first author was visiting MSRI
(summer 2004), the second author was visiting the
Ohio State University (autumn quarter 2003) and all the three authors,
at various times, visited the Max-Planck-Institut f\"ur Mathematik.
We are grateful to all these institutions for excellent working conditions
and a stimulating atmosphere.
 % preface

%\part{Fundamentals}\label{part_knots}
\chapter{Knots and their relatives} % 01

This book is about knots. It is, however, hardly possible to speak
about knots without mentioning other one-dimensional topological
objects embedded into the three-dimensional space. Therefore, in
this introductory chapter we give basic definitions and
constructions pertaining to knots and their relatives: links, braids
and tangles.

The table of knots at the end of this chapter (page
\pageref{knot_table}) will be used throughout the book as a source
of examples and exercises.

\section{Definitions and examples}

A {\em knot} is a closed non-self-intersecting curve in 3-space. In
this book, we shall mainly study smooth oriented knots. A precise
definition can be given as follows.

\begin{definition} \index{Knot}
A {\em parametrized knot} is an embedding of the circle $S^1$ into the Euclidean
space $\R^3$.
\end{definition}

Recall that an {\em embedding}\index{Embedding} is a smooth map
which is injective and whose differential is nowhere zero. In our
case, the non-vanishing of the differential means that the tangent
vector to the curve is non-zero. In the above definition and
everywhere in the sequel, the word {\em smooth} means {\em
infinitely differentiable}.

A choice of an orientation for the parametrizing circle
$$S^1=\{(\cos t,\sin t)\, \vert\, t\in\R\}\subset\R^2$$ gives an
orientation to all the knots simultaneously. We shall always assume
that $S^1$ is oriented counterclockwise. We shall also fix an
orientation of the 3-space; each time we pick a basis for $\R^3$ we
shall assume that it is consistent with the orientation.

If coordinates $x$, $y$, $z$ are chosen in $\R^3$, a knot can be given by
three smooth periodic functions of one variable $x(t)$, $y(t)$, $z(t)$.

\begin{example}\label{plane_circle}
The simplest knot is represented by a plane circle:
\begin{center}
\parbox{4cm}{
$\begin{array}{ccl}
   x &=& \cos t,\\
   y &=& \sin t,\\
   z &=& 0.
\end{array}$}
\qquad
\parbox{4cm}{\rb{-7mm}{\ig[height=20mm]{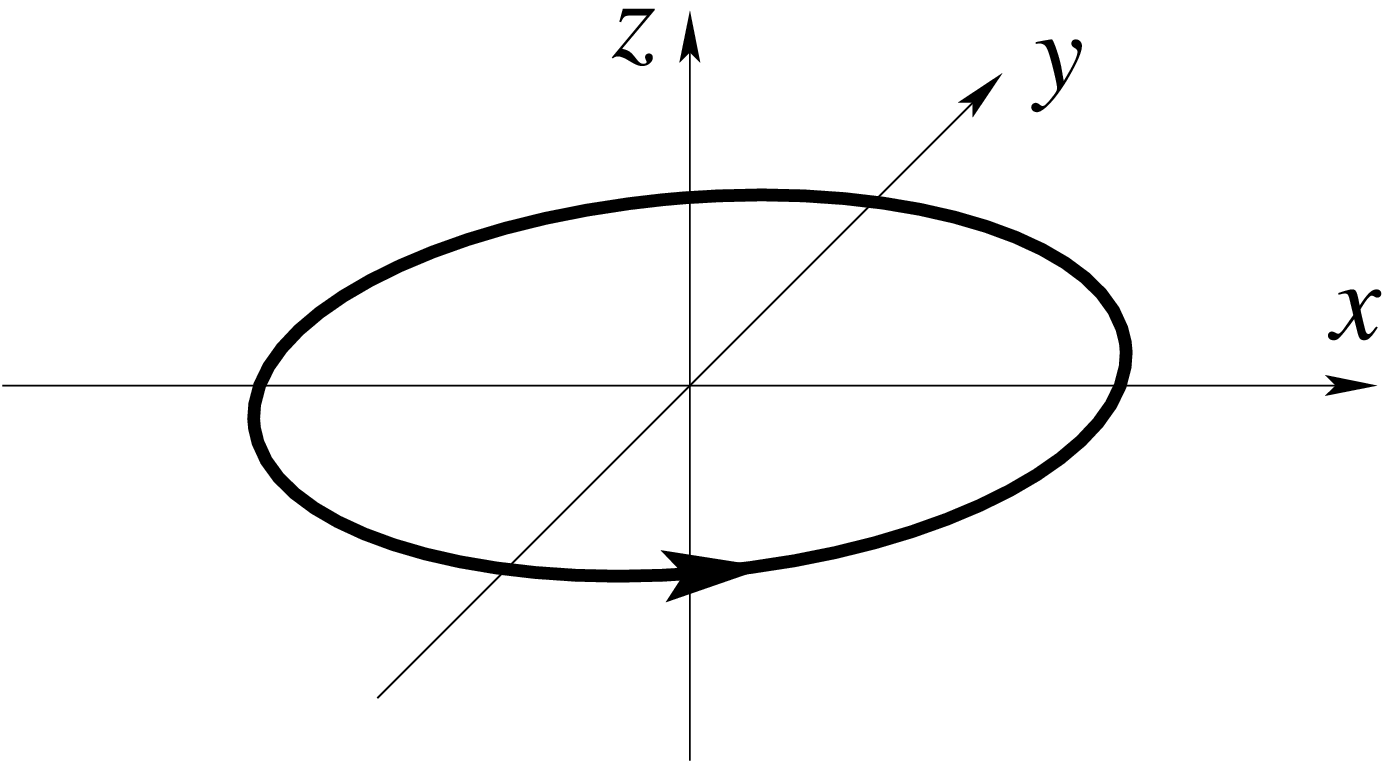}}}
\end{center}
\end{example}

\begin{example}\label{torus23}
The curve that goes 3 times around and 2 times across a standard
torus in $\R^3$ is called the {\em left trefoil knot}, or the
$(2,3)$-{\em torus knot}:
\begin{center}
\parbox{4cm}{
$\begin{array}{ccl}
   x &=& (2+\cos 3t)\cos 2t,\\
   y &=& (2+\cos 3t)\sin 2t,\\
   z &=& \sin 3t.
\end{array}$}
\qquad\qquad
\parbox{4cm}{\rb{-15mm}{\ig[height=18mm]{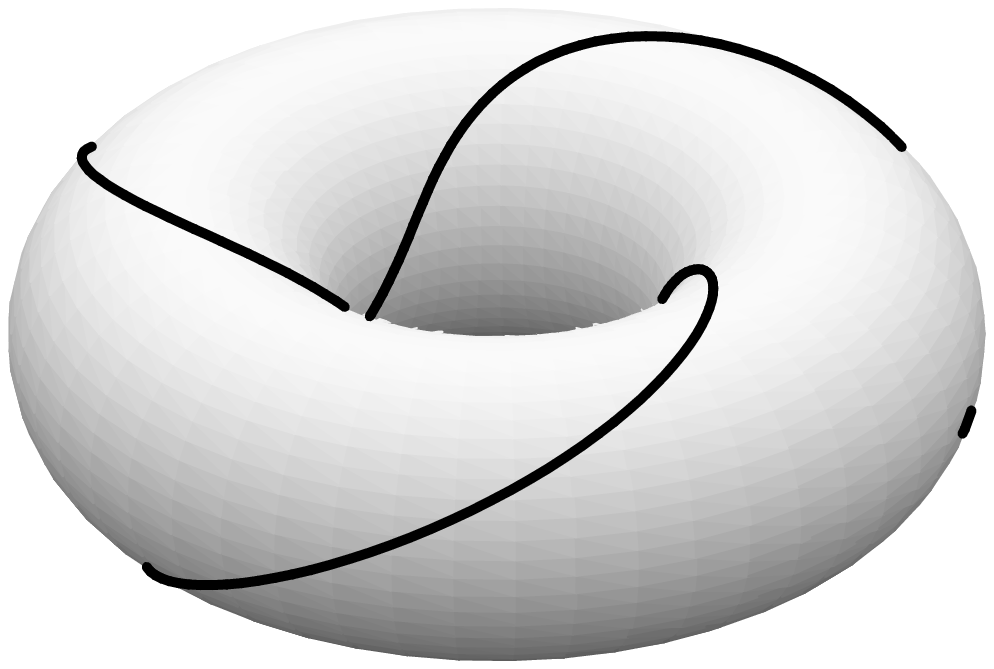}}}
\end{center}
\end{example}

\begin{xca}
Give the definition of a $(p,q)$-torus knot. What are the
appropriate values of $p$ and $q$ for this definition?
\end{xca}

It will be convenient to identify knots that only differ by a change of a parametrization. An {\em oriented knot} is an equivalence class of parametrized knots under orientation-preserving diffeomorphisms of the parametrizing circle.
Allowing {\em all} diffeomorphisms of $S^1$ in this definition, we obtain {\em unoriented knots}.\index{Knot!unoriented}
Alternatively, an unoriented knot can be defined as the {\em image} of
 an embedding of $S^1$ into $\R^3$; an oriented knot is then an image
of such an embedding together with the choice of one of the two
possible directions on it.

We shall distinguish oriented/unoriented knots from parametrized knots in
the notation: oriented and unoriented knots will be usually denoted by
capital letters, while for the individual embeddings lowercase letters will
be used. As a rule, the word ``knot'' will mean ``oriented knot'', unless it
is clear from the context that we deal with unoriented knots, or consider
a specific choice of parametrization.

\section{Isotopy}

The study of parametrized knots falls within the scope of differential
geometry.
The {\em topological} study of knots requires an equivalence relation
which would not only discard the specific choice of parametrization, but also
model the physical transformations of a closed piece of rope in
space.

By a {\em smooth family of maps}, or a {\em map smoothly depending
on a parameter}, we understand a smooth map $F: S^1\times I \to
\R^3$, where $I\subset\R$ is an interval. Assigning a fixed value
$a$ to the second argument of $F$, we get a map  $f_a: S^1\to\R^3$.

\begin{definition}\label{isotopy}\index{Isotopy}
A smooth isotopy of a knot $f: S^1 \to \R^3$, is a smooth family of
knots $f_u$, with $u$ a real parameter, such that for some value $u=a$ we have $f_{a}=f$.
\end{definition}

For example, the formulae
\begin{equation*}
\begin{array}{ccl}
   x &=& (u+\cos{3t})\cos{2t},\\
   y &=& (u+\cos{3t})\sin{2t},\\
   z &=& \sin{3t},
\end{array}
\end{equation*}
where $u\in(1,+\infty)$, represent a smooth isotopy of the trefoil
knot \ref{torus23}, which corresponds to $u=2$. In the pictures below the
space curves are shown by their projection to the $(x,y)$ plane:

\begin{center}
\begin{tabular}{cccc}
{\ig[width=15mm,height=15mm]{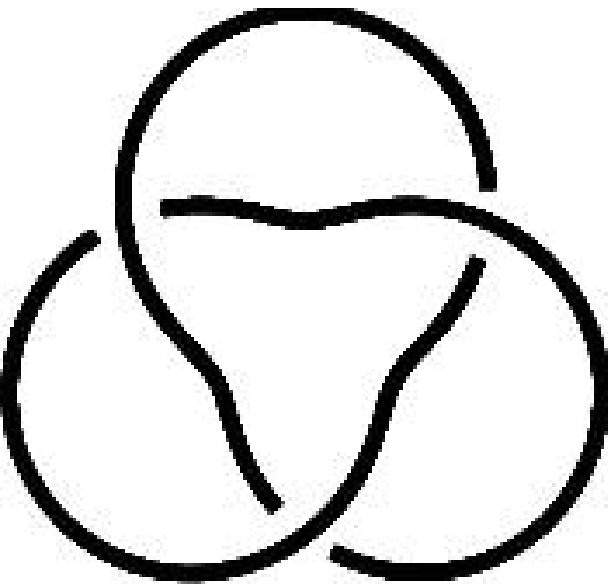}}&
{\ig[width=15mm,height=15mm]{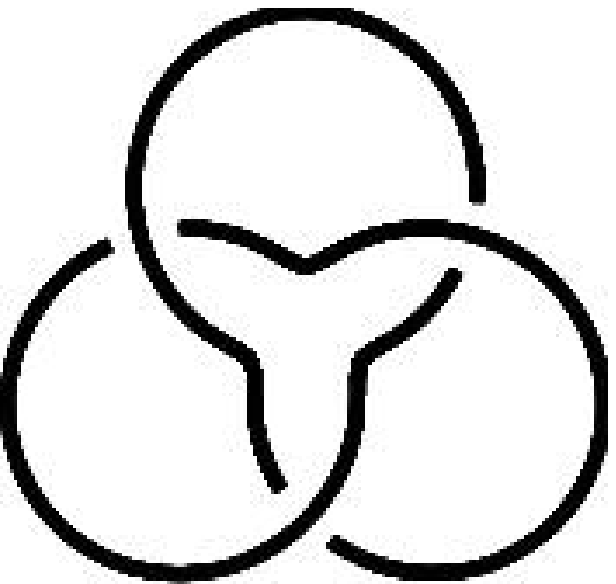}}&
{\ig[width=15mm,height=15mm]{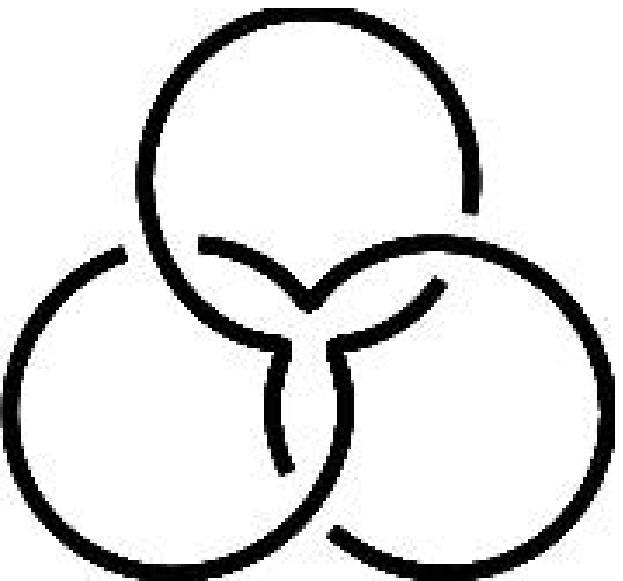}}&
{\ig[width=15mm,height=15mm]{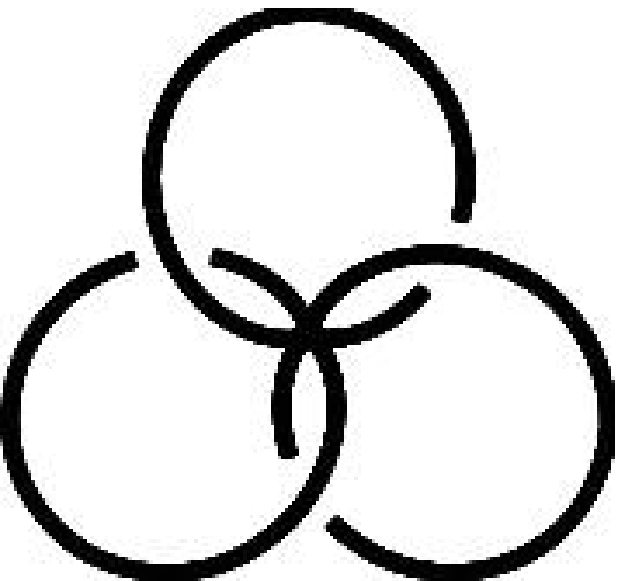}}\\
$u=2$ & $u=1.5$ & $u=1.2$ & $u=1$
\end{tabular}
\end{center}

For any $u>1$ the resulting curve is smooth and has no
self-intersections, but as soon as the value $u=1$ is reached we get
a singular curve with three coinciding cusps\footnote{A {\em cusp}
of a spatial curve is a point where the curve can be represented as
$x=s^2$, $y=s^3$, $z=0$ in some local coordinates.}
corresponding to the values $t=\pi/3$, $t=\pi$ and $t=5\pi/3$. This
curve is not a knot.

\begin{definition}
Two parametrized knots are said to be {\em isotopic} if one can be transformed
into another by means of a smooth isotopy. Two oriented knots
are isotopic if they represent the classes of isotopic parametrized knots; the same definition is valid for unoriented knots.
\end{definition}

\begin{xexample}\label{deform8}
This picture shows an isotopy of the figure eight knot
into its mirror image:
\begin{center}
  \ig[width=\textwidth]{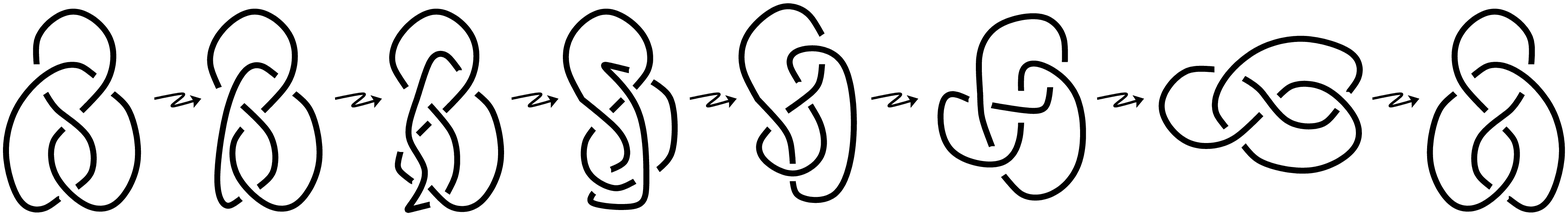}
\end{center}
\end{xexample}

\subsection{}
There are other notions of knot equivalence, namely, {\em ambient
equivalence} and {\em ambient isotopy}, which, for smooth knots, are
the same thing as isotopy. Here are the definitions. A proof that
they are equivalent to our definition of isotopy can be found in
\cite{BZ}.

\begin{xdefinition}\index{Knot!ambient equivalence}
Two parametrized knots, $f$ and $g$, are {\em ambient equivalent} if there is a
commutative diagram
$$
\CD
  S^1    @>f>>  \mathbb R^3   \\
@V\f VV       @VV\psi V \\
  S^1    @>g>>  \mathbb R^3
\endCD
$$
where $\f$ and $\psi$ are orientation preserving diffeomorphisms of
the circle and the 3-space, respectively.
\end{xdefinition}

\begin{xdefinition}\index{Knot!ambient isotopy}
Two parametrized knots, $f$ and $g$, are {\em ambient isotopic} if there is a
smooth family of diffeomorphisms of the 
3-space $\psi_t:\R^3\to\R^3$
with $\psi_0=\id$ and $\psi_1\circ f=g$.
\end{xdefinition}

\subsection{}\label{unknots}
A knot, equivalent to the plane circle of Example \ref{plane_circle}
is referred to as a {\em trivial knot}, or an {\em unknot}.

Sometimes, it is not immediately clear from a diagram of a trivial knot that it is
indeed trivial:

\begin{center}\label{unknots-diag}
  \ig[height=17mm]{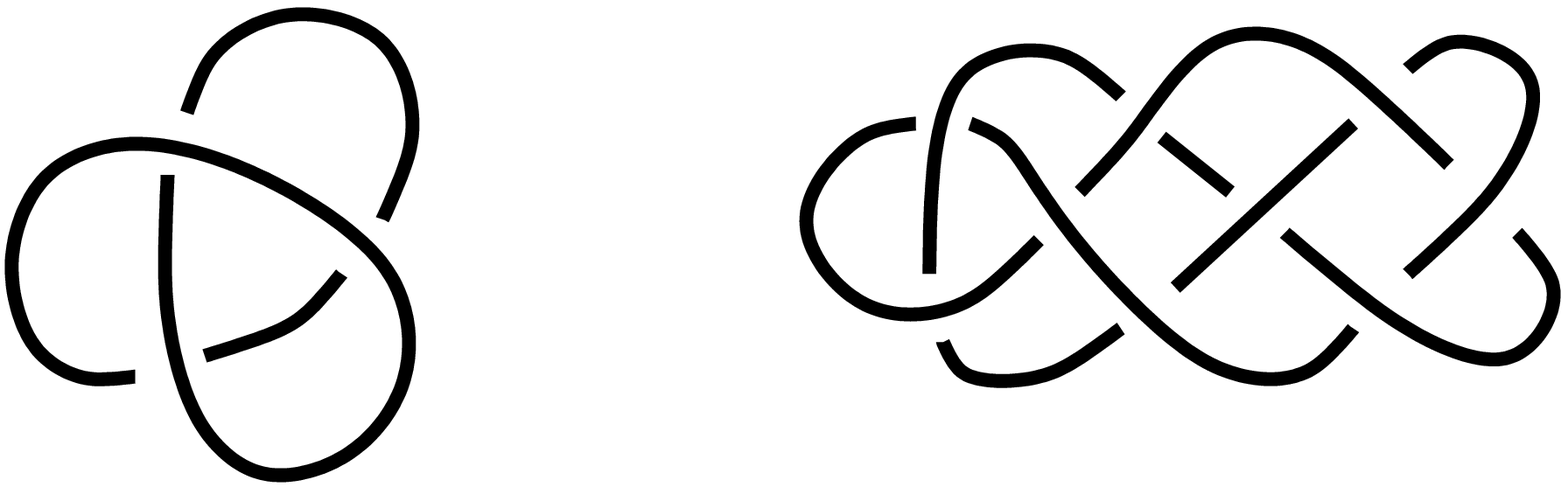}\qquad\qquad
  \risS{7}{unknot-g3}{}{75}{0}{0}\\[-2mm]
Trivial knots
\end{center}

There are algorithmic procedures to detect whether a given knot
diagram represents an unknot. One of them, based on W.~Thurston's
ideas, is implemented in J.~Weeks' computer program
\texttt{SnapPea}, see \cite{Wee}; another algorithm, due to I.~Dynnikov, 
is described in \cite{Dyn}.

Here are several other examples of knots.

\begin{center}\label{named-knots}
\begin{tabular}{c@{\quad}c@{\quad}c@{\quad}c@{\quad}c}
\ig[height=16mm]{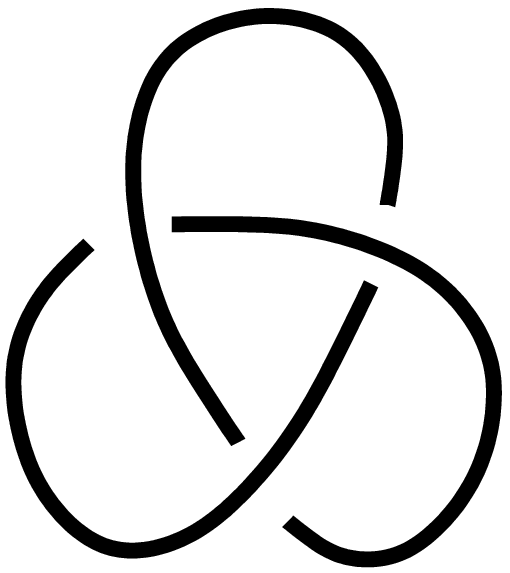} & \ig[height=16mm]{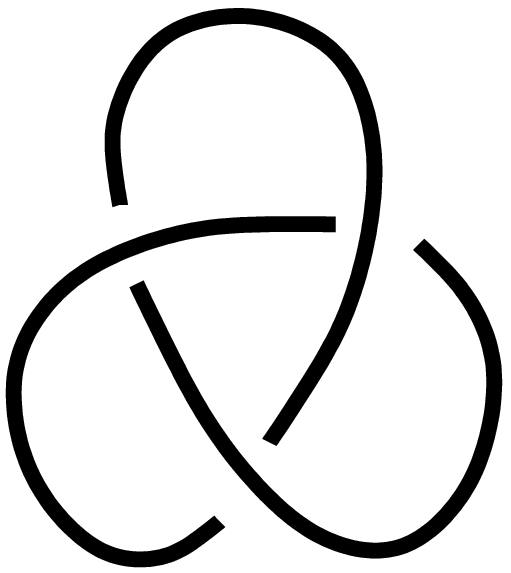} &
\ig[height=16mm]{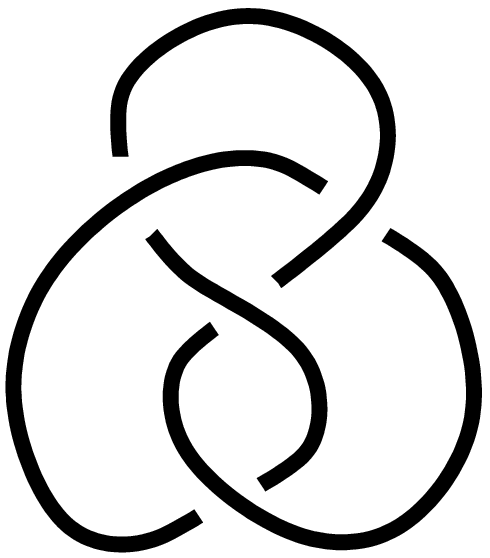} & \ig[height=16mm]{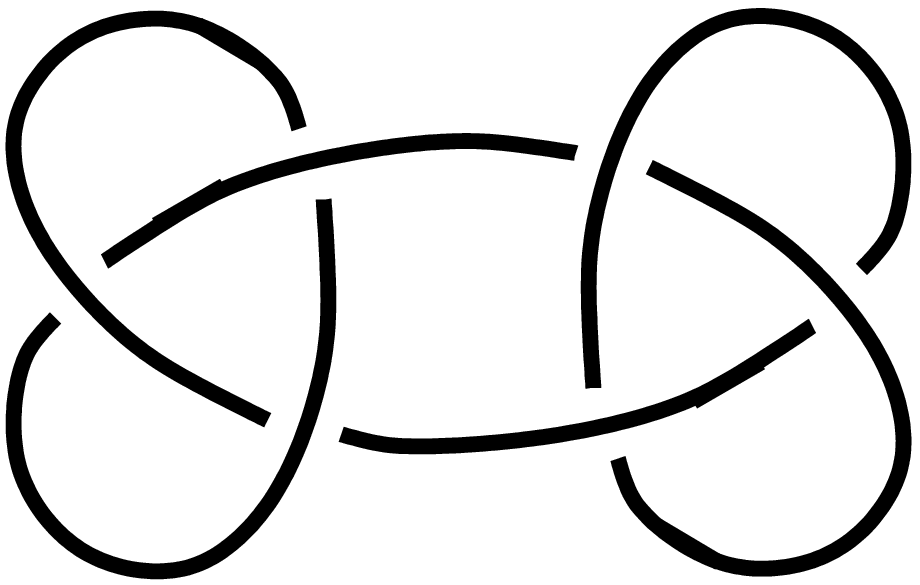} &
\ig[height=16mm]{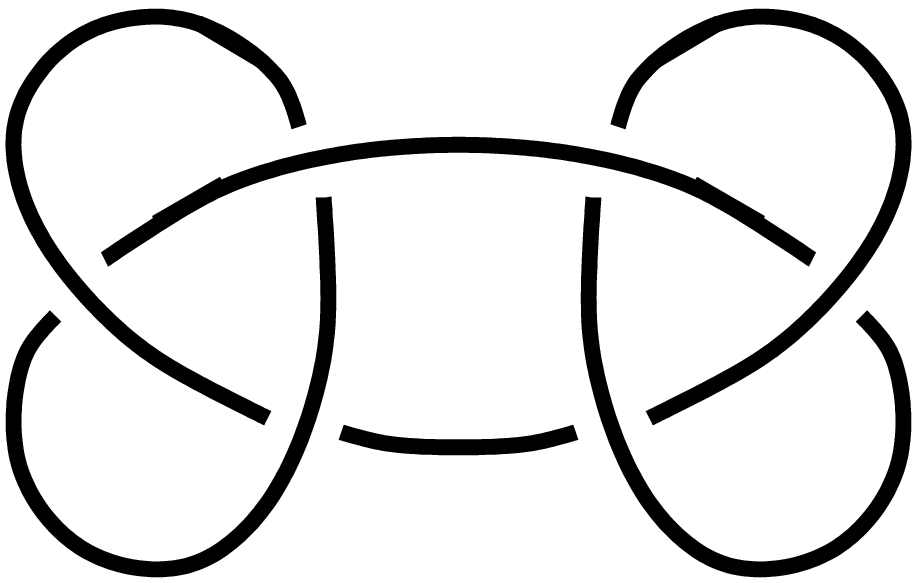} \\
\scriptsize Left trefoil & \scriptsize Right trefoil & \scriptsize Figure 8 knot &
\scriptsize Granny knot & \scriptsize Square knot
\index{Knot!figure eight}\index{Trefoil}\index{Knot!trefoil}
\index{Granny knot}\index{Knot!granny}
\index{Square knot}\index{Knot!square}
\end{tabular}
\end{center}

Knots are a special case of links.

\begin{definition} \label{def:link}\index{Link}
A {\em link} is a smooth embedding $S^1\sqcup\dots\sqcup
S^1\to\R^3$, where $S^1\sqcup\dots\sqcup S^1$ is the disjoint union
of several circles.
\end{definition}

\begin{center}
\index{Borromean rings} \index{Hopf link} \index{Whitehead link}
\begin{tabular}{c@{\quad}c@{\quad}c@{\quad}c}
\ig[height=16mm]{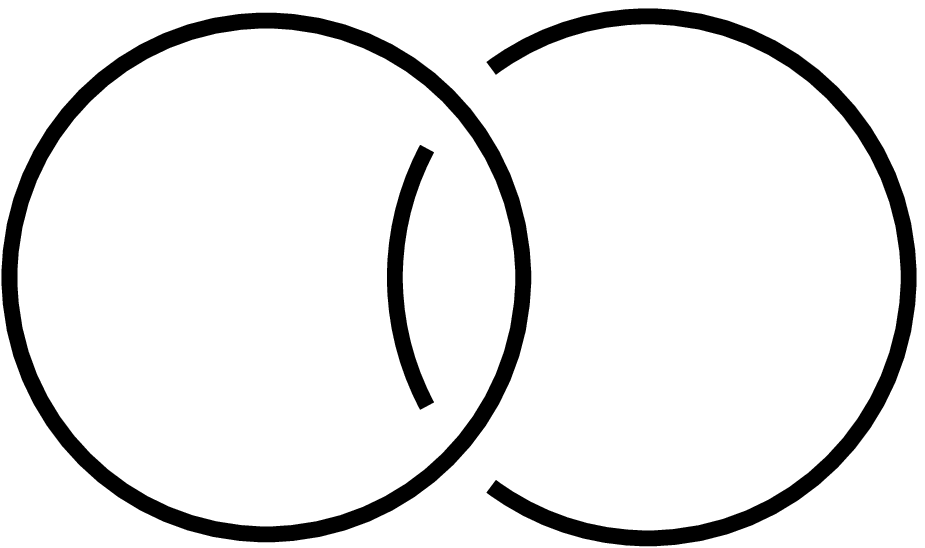} & \ig[height=16mm]{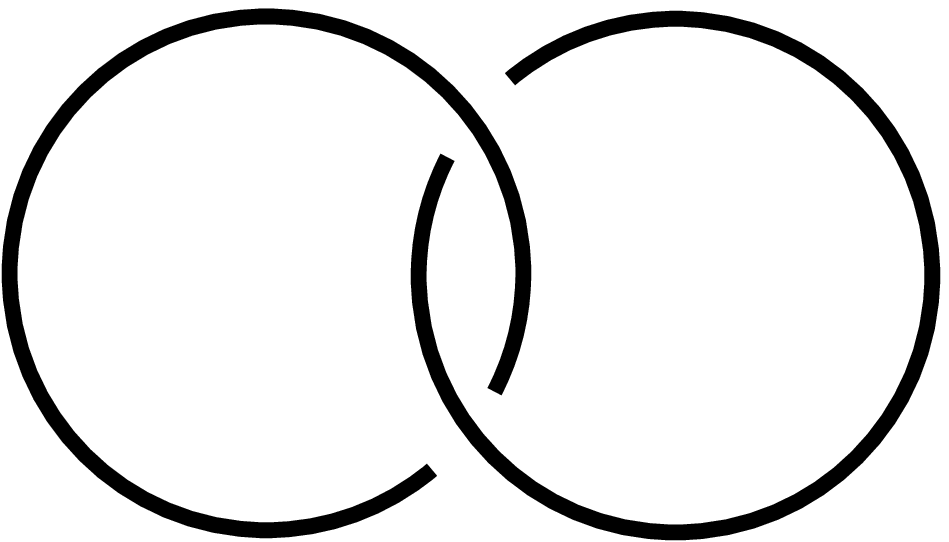} &
\ig[height=16mm]{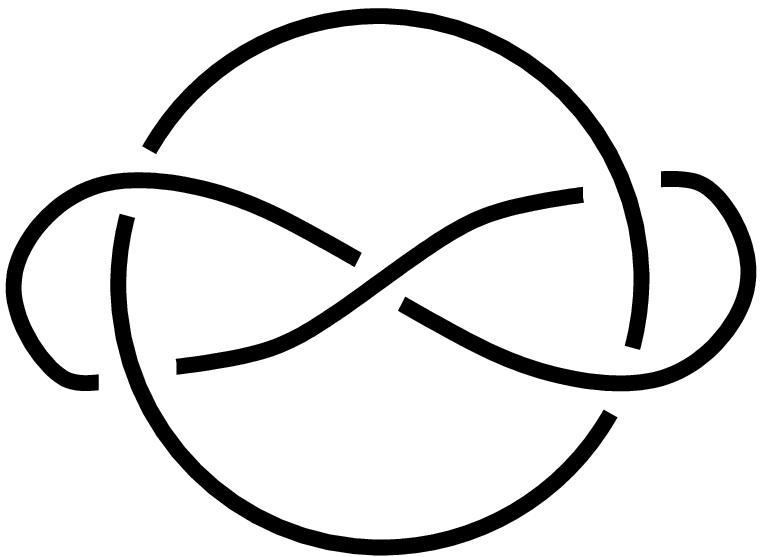} & \ig[height=16mm]{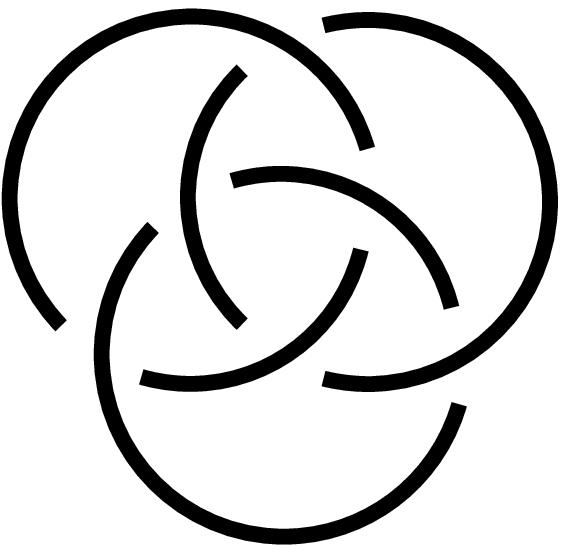}\\
\scriptsize Trivial 2-component link & \scriptsize Hopf link &
\scriptsize Whitehead link & \scriptsize Borromean rings
\index{Link!trivial}\index{Link!Hopf}\index{Link!Whitehead}\index{Link!Borromean
rings} \index{Hopf link}\index{Whitehead link}\index{Borromean
rings} \label{link_ex}
\end{tabular}
\end{center}

Equivalence of links is defined in the same way as for knots ---
with the exception that now one may choose whether to distinguish or
not between the components of a link and thus speak about the
equivalence of links with {\em numbered} or {\em unnumbered}
components.

In the future, we shall often say ``knot (link)'' instead of
``equivalence class'', or ``topological type of knots (links)''.

\section{Plane knot diagrams}
\label{pl_diag}

Knots are best represented graphically by means of {\em knot
diagrams}\index{Knot diagram}.  A knot diagram is a plane curve
whose only singularities are transversal double points (crossings),
together with the choice of one branch of the curve at each
crossing. The chosen branch is called an {\em overcrossing}; the
other branch is referred to as an {\em undercrossing}. A knot
diagram is thought of as a projection of a knot along some
``vertical'' direction; overcrossings and undercrossings indicate
which branch is ``higher'' and which is ``lower''.
To indicate the orientation, an arrow is added to the knot diagram.

\begin{theorem}[Reidemeister \cite{Rei}, proofs can be found in
\cite{PS, BZ,Mura}]
\label{ReiTh}\index{Theorem!Reidemeister} Two unoriented knots
$K_1$ and $K_2$, are equivalent if and only if a diagram of $K_1$
can be transformed into a diagram of $K_2$ by a sequence of ambient
isotopies of the plane and local moves of the following three types:
$$
\index{Reidemeister moves}\index{Moves!Reidemeister}
\begin{array}{ccccccc}
\Omega_1 &&& \Omega_2 &&& \Omega_3 \\[3mm]
\ig[width=90pt,height=40pt]{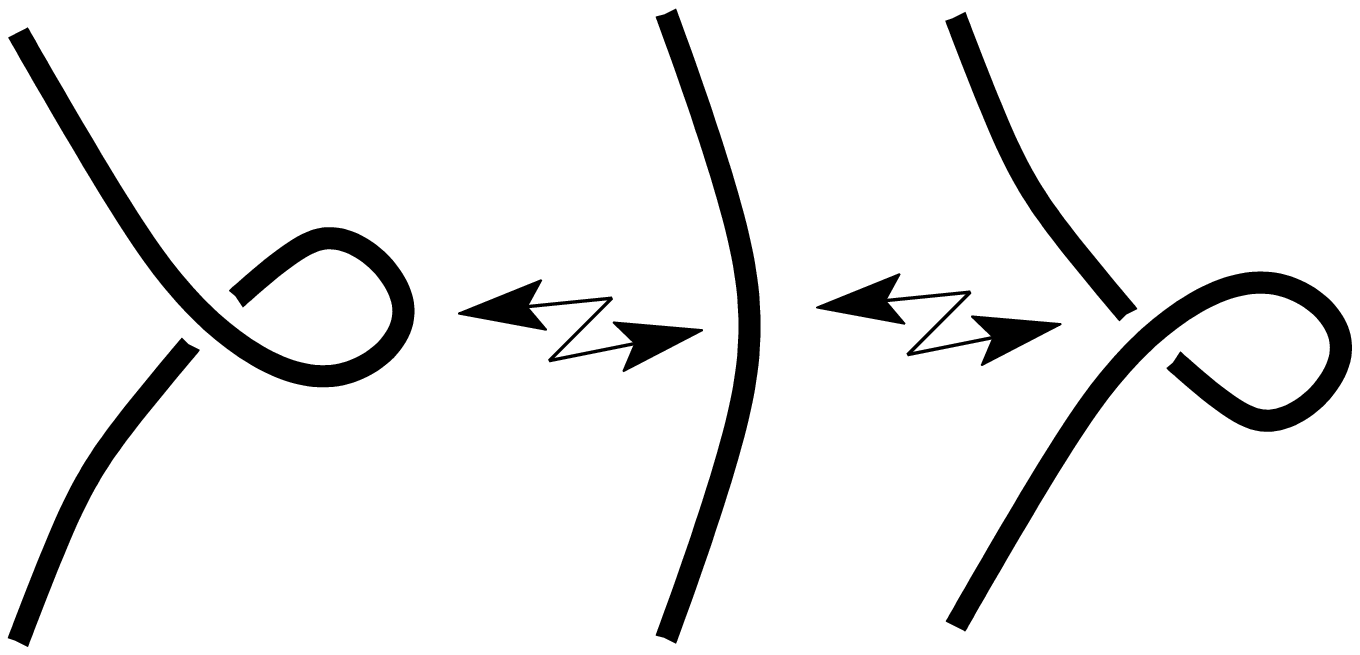} &&&
\ig[width=90pt,height=40pt]{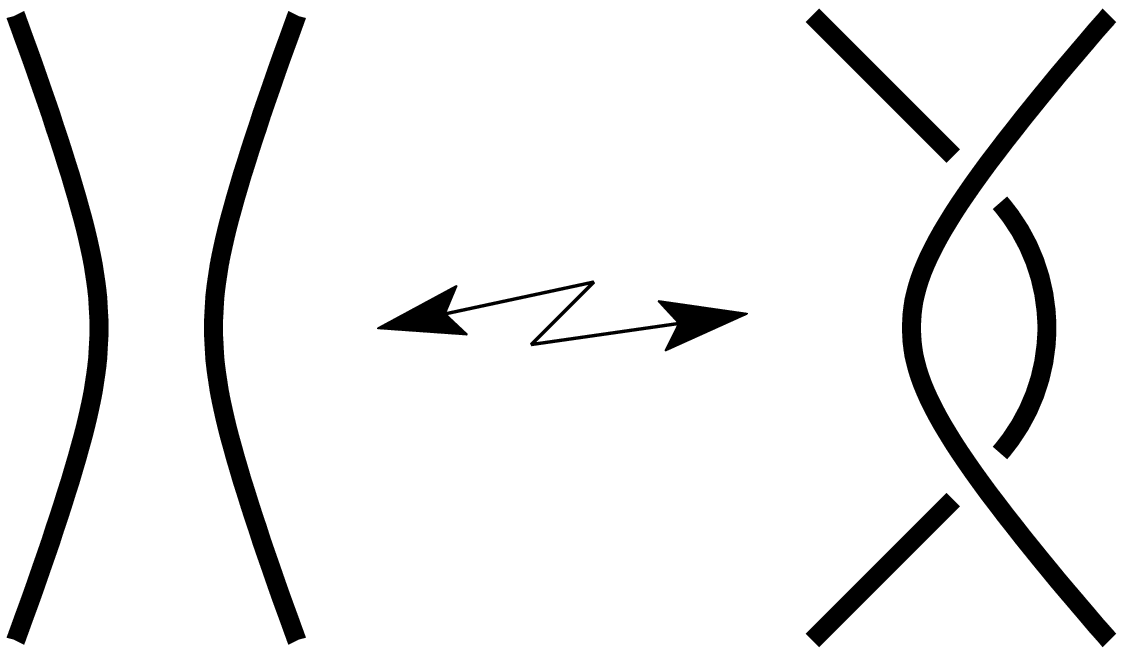} &&&
\ig[width=105pt,height=40pt]{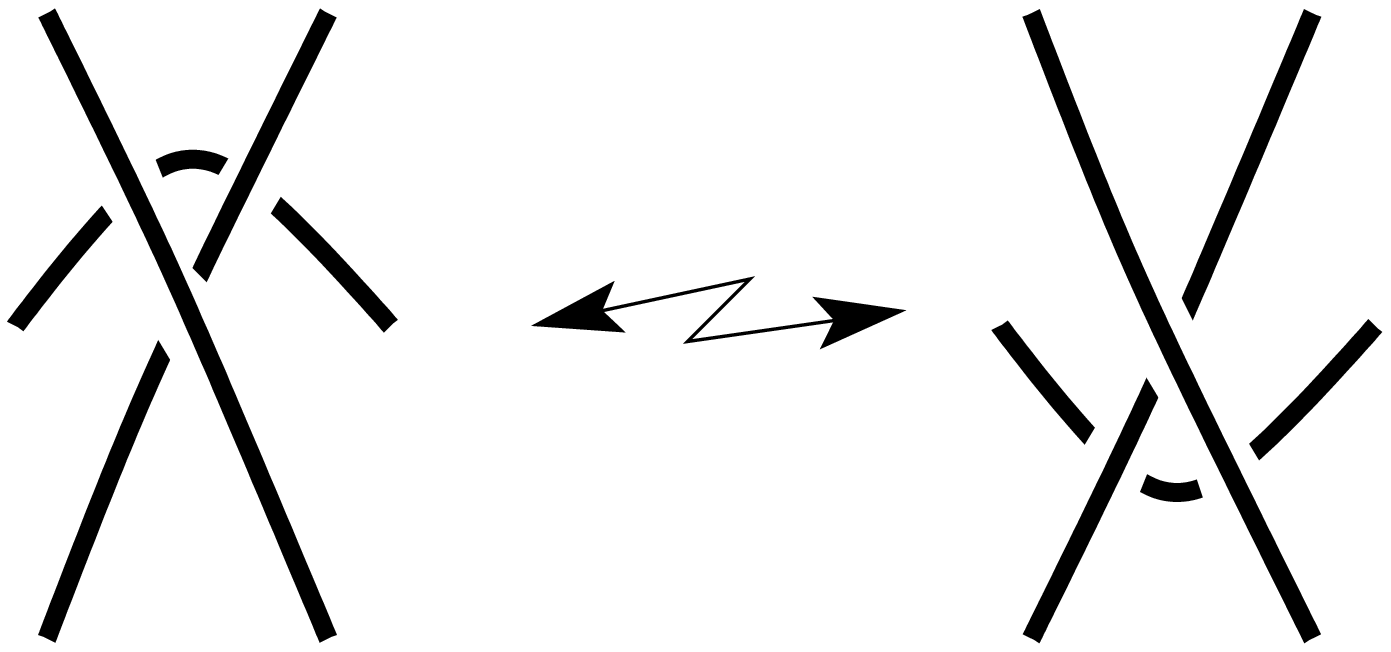} \\[3mm]
\multicolumn{7}{c}{\mbox{Reidemeister moves}}
\end{array}
$$
\end{theorem}

To adjust the assertion of this theorem to the oriented case, each
of the three Reidemeister moves has to be equipped with orientations
in all possible ways. Smaller sufficient sets of oriented moves exist;
one such set will be given later in terms of Gauss diagrams (see p.
\pageref{gd_moves}).
\medskip
\noindent

\textbf{Exercise.}
Determine the sequence of Reidemeister moves that
relates the two diagrams of the trefoil knot below:

\begin{center}
\ig[height=20mm]{k3-1.eps}
\hspace{20mm}
\ig[height=20mm]{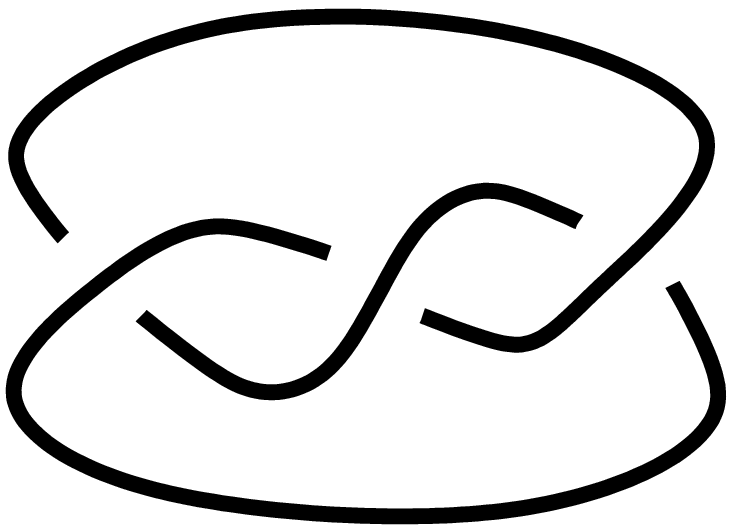}
\end{center}

\subsection{Local writhe}\label{writhe}

Crossing points on a diagram come in two species, positive and
negative:
\begin{center}\label{loc_writhe}\index{Writhe!local}
\begin{tabular}{ccc}
\ig[width=15mm]{lrints.eps} &
\qquad &
\ig[width=15mm]{rlints.eps} \\
\parbox{40mm}{
\begin{center}
Positive crossing
\end{center}} &
\qquad &
\parbox{40mm}{
\begin{center}
Negative crossing
\end{center}}
\end{tabular}
\end{center}

Although this sign is defined in terms of the knot orientation, it
is easy to check that it does not change if the orientation is
reversed. For \textit{links} with more than one component, the
choice of orientation is essential.

The {\em local writhe} of a crossing is defined as
$+1$ or $-1$ for positive or negative points, respectively. The
{\em writhe}\index{Writhe!of a knot diagram} (or {\em total
writhe})\index{Writhe!total} of a diagram is the sum of the writhes
of all crossing points, or, equivalently, the difference between the
number of positive and negative crossings. Of course, the same knot
may be represented by diagrams with different total writhes. In
Chapter \ref{kn_inv} we shall see how the writhe can be used to
produce knot invariants.

\subsection{Alternating knots}\label{altkn}
\index{Knot!alternating} A knot diagram is called {\it
alternating}\index{Knot diagram!alternating} if its overcrossings
and undercrossing alternate as we travel along the knot. A knot
is called {\it alternating} if it has an alternating diagram. A knot
diagram is called \textit{reducible}\index{Knot diagram!reducible}
if it becomes disconnected after the removal of a small
neighbourhood of some crossing.

The number of crossings in a reducible diagram can
be decreased by a move shown in the picture:
$$\ris{-10}{nred_kn_d}{80}{115}{60}{20}{
  \put(10,-10){\mbox{\scriptsize reducible diagram}}}\qquad
\ris{18}{totor}{0}{40}{0}{0}{}\qquad
\ris{-2}{red_kn_d}{0}{100}{0}{0}{
  \put(20,-15){\mbox{\scriptsize reduction}}}\qquad
$$

A diagram which is not reducible is called \textit{reduced}.
As there is no immediate way to simplify a reduced diagram,
the following conjecture naturally arises
(P.~G.~Tait, 1898).

{\bf The Tait conjecture.}\index{Conjecture!Tait}\index{Tait conjecture}
{\it A reduced alternating diagram has the minimal number of crossings
among all diagrams of the given knot.}

This conjecture stood open for almost 100 years. It was proved only
in 1986 (using the newly invented Jones polynomial) simultaneously
and independently by  L.~Kauffman \cite{Ka6}, K.~Murasugi
\cite{Mur}, and M.~Thistlethwaite \cite{Th} (see Exercise
(\ref{ex_KMT_theorem}) in Chapter \ref{kn_inv}).

\section{Inverses and mirror images}
\label{orientation}

Change of orientation (taking the inverse) and taking the mirror
image are two basic operations on knots which are induced by
orientation reversing smooth involutions on $S^1$ and $\R^3$
respectively. Every such involution on $S^1$ is conjugate to the
reversal of the parametrization; on $\R^3$ it is conjugate to a
reflection in a plane mirror.

Let $K$ be a knot. Composing the parametrization reversal of $S^1$
with the map $f:S^1\to\R^3$ representing $K$, we obtain the {\em
inverse} $K^*$ of $K$. The {\em mirror image} of $K$, denoted by
$\ol{K}$, is a composition of the map $f:S^1\to\R^3$ with a
reflection in $\R^3$. Both change of orientation and taking the
mirror image are involutions on the set of (equivalence classes of)
knots. They generate a group isomorphic to $\Z_2\oplus\Z_2$; the
symmetry properties of a knot $K$ depend on the subgroup that leaves
it invariant. The group $\Z_2\oplus\Z_2$ has 5 (not necessarily
proper) subgroups, which give rise to 5 symmetry classes of knots.

\begin{definition}\index{Knot!invertible}
\index{Knot!amphicheiral}\label{amphicheiral} \index{Knot!symmetric}
\index{Knot!asymmetric} A knot is called:
\begin{itemize}
\item \textit{invertible}, if $K^*=K$,
\item \textit{plus-amphicheiral}, if $\ol{K}=K$,
\item \textit{minus-amphicheiral}, if $\ol{K}=K^*$,
\item \textit{fully symmetric}, if $K=K^*=\ol{K}=\ol{K}^*$,
\item \textit{totally asymmetric}, if all knots $K$, $K^*$, $\ol{K}$,
$\ol{K}^*$ are different.
\end{itemize}
\end{definition}

The word {\em amphicheiral} means either plus- or
minus-amphicheiral. For invertible knots, this is the same.
Amphicheiral and non-amphicheiral knots are also referred to as {\em
achiral}\index{Knot!achiral} and {\em chiral} knots,
respectively.\index{Knot!chiral}

The 5 symmetry classes of knots are summarized in the following
table. The word ``minimal'' means ``with the minimal number of
crossings''; $\sigma$ and $\tau$ \label{inv-s-t}
denote the involutions of taking
the mirror image and the inverse respectively.

\begin{center}
\def\lis{\rb{-6pt}{\makebox(0,20){ }}}
\begin{tabular}{|c|c|c|c|}
\hline
{\bf Stabiliser} & {\bf Orbit} & {\bf Symmetry type} $\lis$& {\bf Min example} \\
\hline\hline $\{1\}$ & $\{K,\ol{K},K^*,\ol{K}^*\}$ & totally
asymmetric & $9_{32}$,\quad $9_{33}
               $\lis \\
\hline
$\{1,\sigma\}$ & $\{K,K^*\}$ & $+$amphicheiral, non-inv & $12^a_{427}$\lis \\
\hline
$\{1,\tau\}$ & $\{K,\ol{K}\}$  & invertible, chiral & $3_{1}$\lis \\
\hline
$\{1,\sigma\tau\}$ & $\{K,K^*\}$ & $-$amphicheiral, non-inv & $8_{17}$\lis \\
\hline
$\{1,\sigma,\tau,\sigma\tau\}$ & $\{K\}$ & fully symmetric & $4_{1}$\lis \\
\hline
\end{tabular}
\end{center}

\begin{xexample}
The trefoil knots are invertible, because the rotation through
$180^\circ$ around an axis in $\R^3$ changes the direction of the
arrow on the knot.
\end{xexample}

The existence of non-invertible knots was first proved by H.~Trotter
\cite{Tro} in 1964. The simplest instance of Trotter's theorem is a
{\em pretzel knot\index{Pretzel knot}\index{Knot!pretzel}
with parameters $(3,5,7)$}:
$$
\ig[height=25mm]{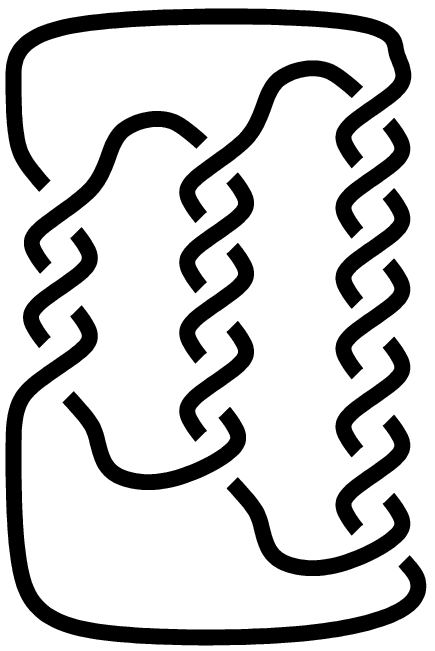}
$$
Among the knots with up to 8 crossings (see
Table \ref{knot_table} on page \pageref{knot_table}) there is only
one non-invertible knot: $8_{17}$, which is, moreover, minus-amphicheiral.
These facts were proved in 1979 by A.~Kawauchi \cite{Ka1}.

\begin{xexample}
The trefoil knots are not amphicheiral, hence the distinction
between the left and the right trefoil. A proof of this fact, based
on the calculation of the Jones polynomial, will be given in Sec.
\ref{Jones}.
\end{xexample}

\begin{xremark}
Knot tables only list knots up to taking inverses and mirror images.
In particular, there is only one entry for the trefoil knots. Either
of them is often referred to as {\em the trefoil}.
\end{xremark}

\begin{xexample}
The figure eight knot is amphicheiral. The isotopy between this knot
and its mirror image is shown on page \pageref{deform8}.
\end{xexample}

Among the 35 knots with up to 8 crossings shown in Table \ref{knot_table},
there are exactly 7 amphicheiral knots: $4_1$, $6_3$, $8_3$, $8_9$, $8_{12}$,
$8_{17}$, $8_{18}$, out of which $8_{17}$ is minus-amphicheiral,
the rest, as they are invertible, are both plus- and minus-amphicheiral.

The simplest totally asymmetric knots appear in 9 crossings, they are
$9_{32}$ and $9_{33}$. The following are all non-equivalent:

$$\begin{array}{ccccccc}
\ig[height=16mm]{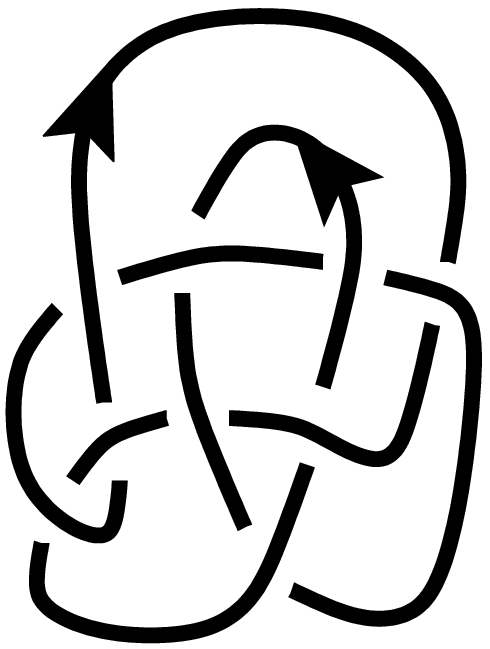}&\qquad&
\ig[height=16mm]{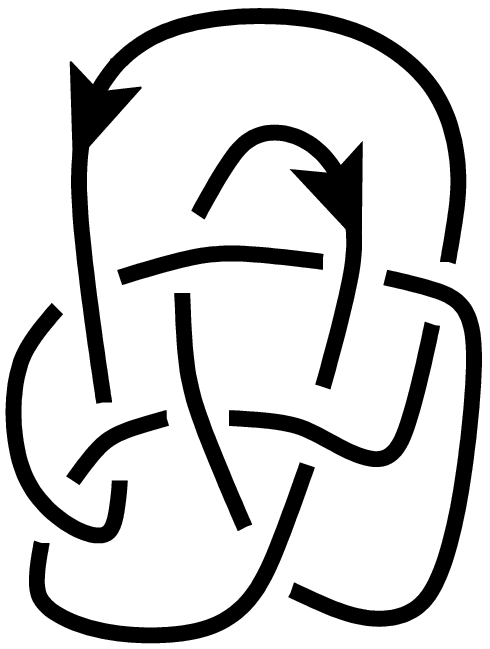}&\qquad&
\ig[height=16mm]{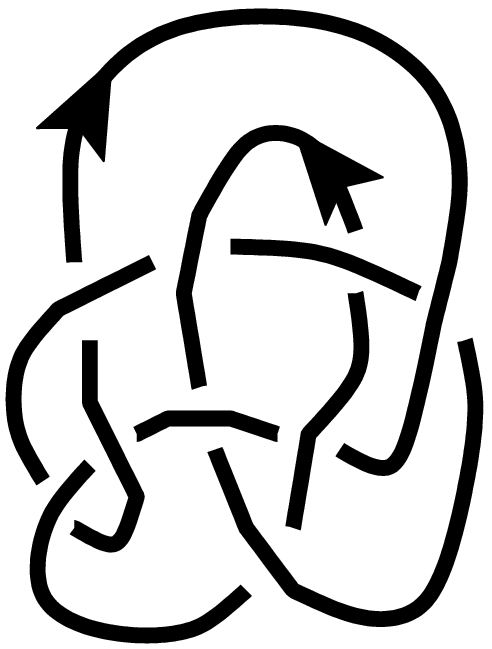}&\qquad&
\ig[height=16mm]{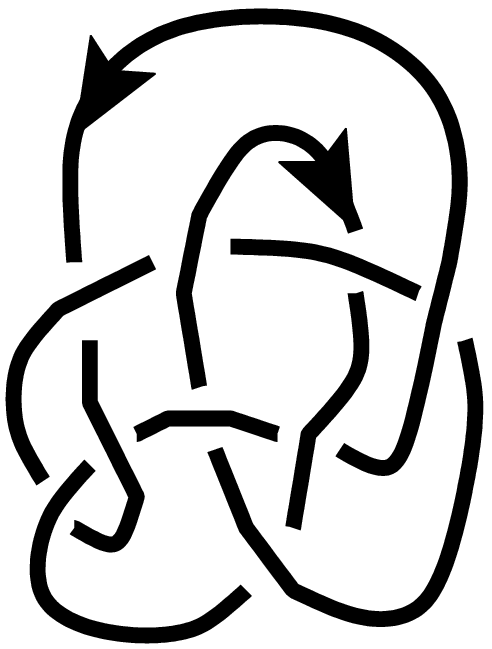}\\
9_{33} && 9_{33}^* && \ol{9}_{33} && \ol{9}_{33}^*
\end{array}$$

Here is the simplest plus-amphicheiral non-invertible knot, together with its
inverse:
$$\begin{array}{ccc}
\ig[height=20mm]{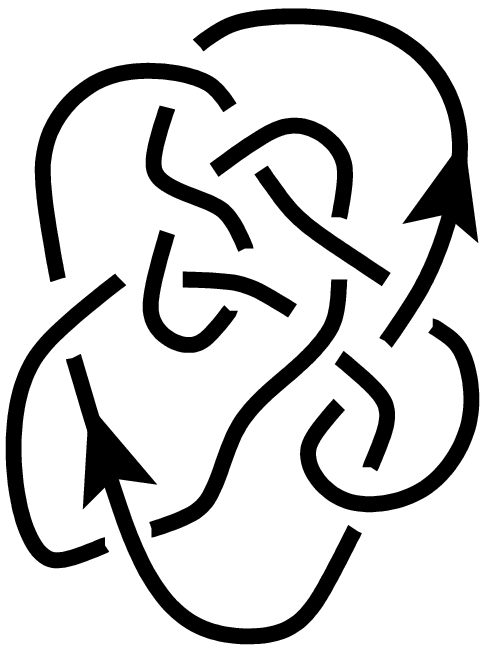} & \qquad & \ig[height=20mm]{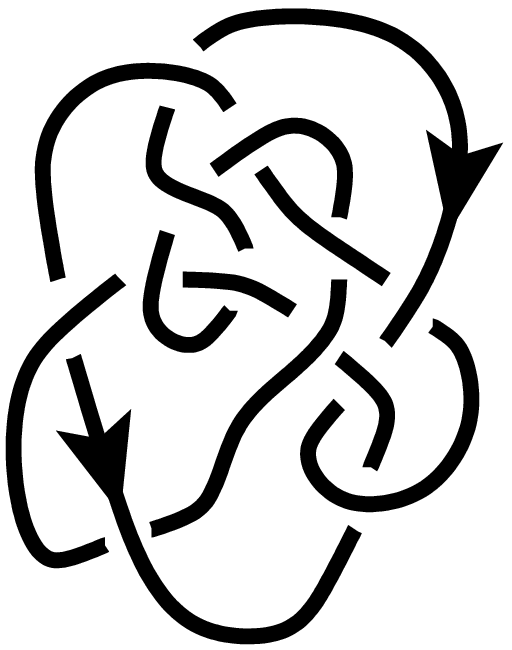}\\
12^a_{427} && {12^{a*}_{427}}
\end{array}$$

In practice, the easiest way to find the symmetry type of a given
knot or link is by using the computer program Knotscape \cite{Knsc},
which can handle link diagrams with up to 49 crossings.

\def\boldX{\boldsymbol{X}}

\section{Knot tables}

\subsection{Connected sum}
There is a natural way to fuse two knots into one: cut each of the
two knots at some point, then connect the two pairs of loose ends.
This must be done with some caution: first, by a smooth isotopy,
both knots should be deformed so that for a certain plane projection
they look as shown in the picture below on the left, then they
should be changed inside the dashed disk as shown on the right:

$$
  {\ig[width=100pt]{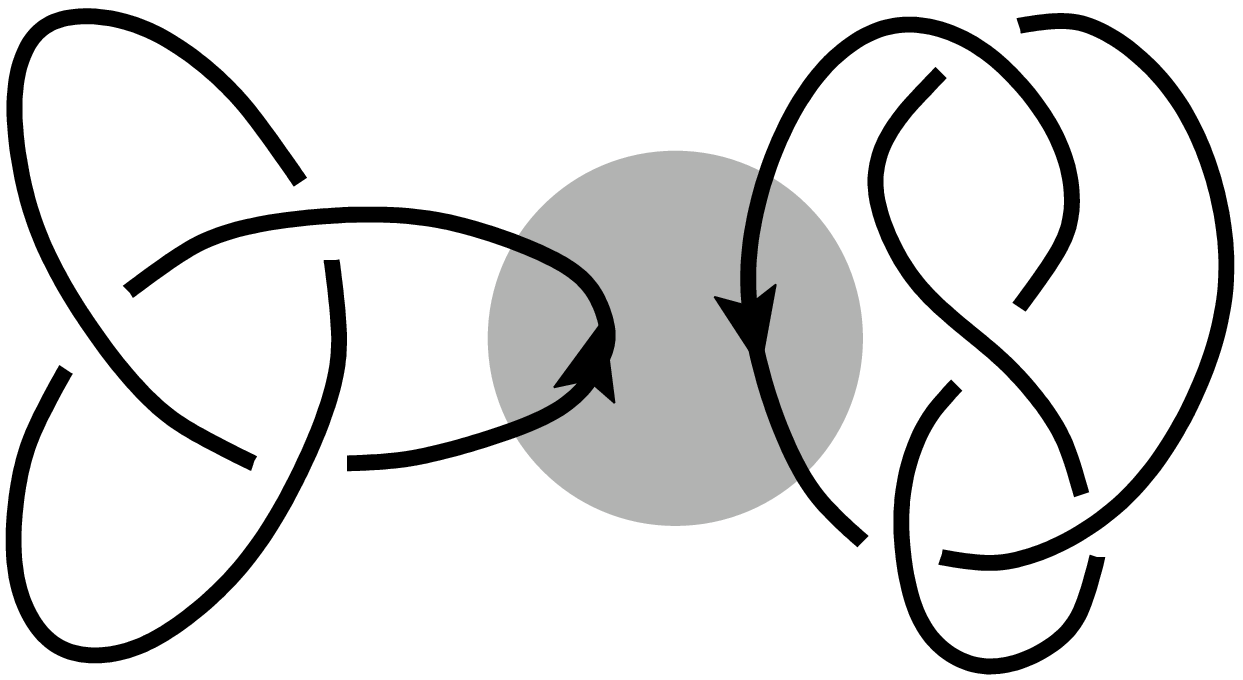}}\qquad
  \rb{20pt}{\ig[width=40pt]{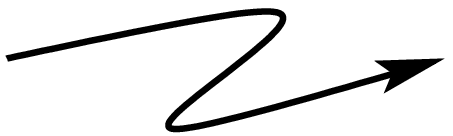}}\qquad
  {\ig[width=100pt]{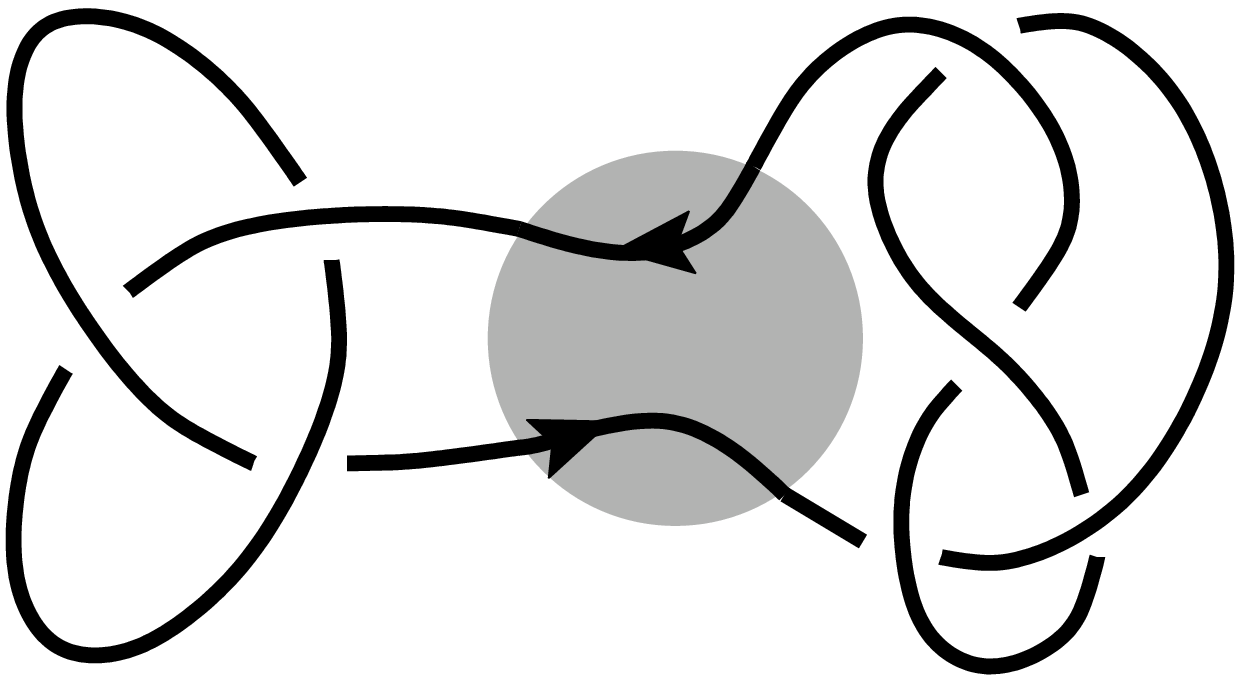}}
$$

The connected sum makes sense only for oriented knots.
It is well-defined and commutative on the equivalence classes of knots.
The connected sum of knots $K_1$ and $K_2$ is denoted by $K_1\#K_2$.\index{Connected sum! of knots}\label{connectedsum}

\begin{definition}
A knot is called {\em prime} if it cannot be represented as the
connected sum of two nontrivial knots.
\end{definition}

Each knot is a connected sum of prime knots, and this decomposition
is unique (see \cite{CrF} for a proof). In particular, this means
that a trivial knot cannot be decomposed into a sum of two
nontrivial knots. Therefore, in order to classify all knots, it is
enough to have a table of all prime knots.

Prime knots are tabulated according to the minimal number of crossings that
their diagrams can have.  Within each group of knots with the same crossing
number, knots are numbered in some, usually rather arbitrary, way. In 
Table~\ref{knot_table}, we use the widely adopted numbering that goes back to the
table compiled by Alexander and Briggs in 1927 \cite{AB}, then repeated
(in an extended and modified way) by D.~Rolfsen in \cite{Rol}. We also
follow Rolfsen's conventions in the choice of the version of
non-amphicheiral knots: for example, our
$3_1$ is the left, not the right, trefoil.

\begin{table}[htb]
$$\begin{array}{ccccccc}\label{kn_table}
\knt{k3-1}&\knt{k4-1}&\knt{k5-1}&\knt{k5-2}& 
   \knt{k6-1}&\knt{k6-2}&\knt{k6-3}\\
3_1&4_1(a)&5_1&5_2&6_1&6_2&6_3(a)\\[1mm]
\hline\\
\knt{k7-1}&\knt{k7-2}&\knt{k7-3}&\knt{k7-4}&
   \knt{k7-5}&\knt{k7-6}&\knt{k7-7}\\
7_1&7_2&7_3&7_4&7_5&7_6&7_7\\[1mm]
\hline\\
\knt{k8-01}&\knt{k8-02}&\knt{k8-03}&\knt{k8-04}&
   \knt{k8-05}&\knt{k8-06}&\knt{k8-07}\\
8_1&8_2&8_3(a)&8_4&8_5&8_6&8_7\\[1mm]
\hline\\
\knt{k8-08}&\knt{k8-09}&\knt{k8-10}&\knt{k8-11}&
   \knt{k8-12}&\knt{k8-13}&\knt{k8-14}\\
8_8&8_9(a)&8_{10}&8_{11}&8_{12}(a)&8_{13}&8_{14}\\[1mm]
\hline\\
\knt{k8-15}&\knt{k8-16}&\knt{k8-17}&\knt{k8-18}&
   \knt{k8-19}&\knt{k8-20}&\knt{k8-21}\\
8_{15}&8_{16}&8_{17}(na-)&8_{18}(a)&8_{19}&8_{20}&8_{21}\\[1mm]
\hline\\
\end{array}$$
\caption{Prime knots, up to orientation and mirror images, with at most 
8 crossings. Amphicheiral knots are marked by `a', the (only)
non-invertible minus-amphicheiral knot by `na-'.}
\label{knot_table}\index{Knot!table}\index{Table of!knots}
\end{table}

Rolfsen's table of knots, authoritative as it is, contained an
error. It is the famous {\em Perko pair} (knots $10_{161}$ and
$10_{162}$ in Rolfsen) --- two equivalent knots that were thought to
be different for 75 years since 1899:

\begin{center}
 \ig[height=2cm]{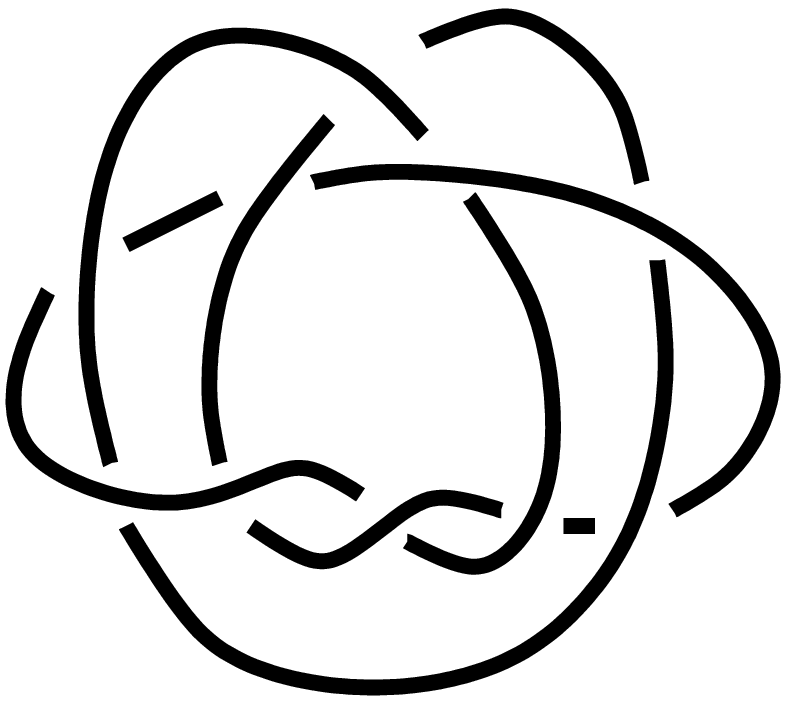}
 \hspace{2cm}
 \ig[height=2cm]{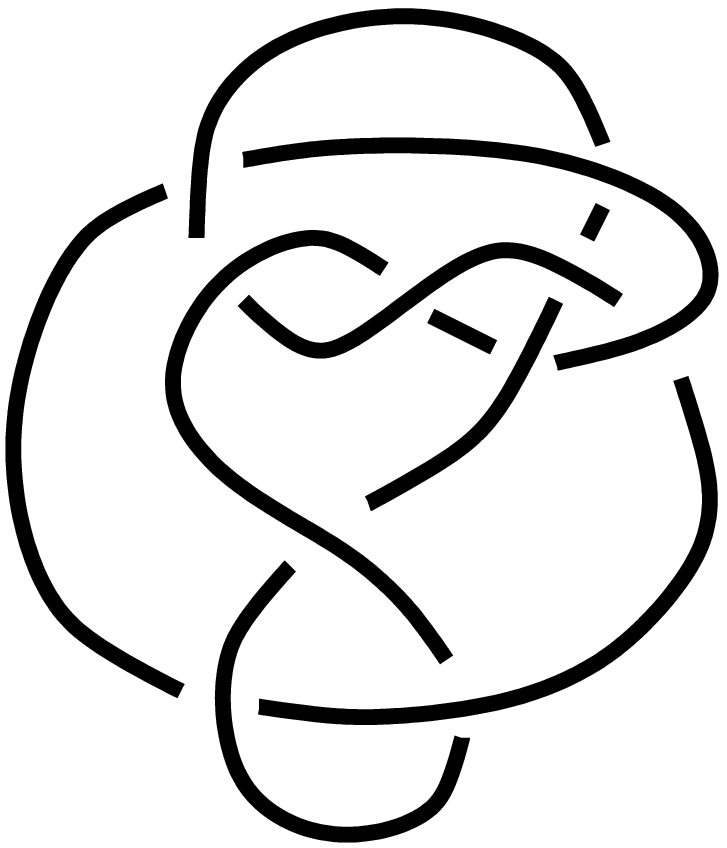}
\end{center}

The equivalence of these two knots was established in 1973 by
K.~A.~Perko \cite{Per1}, a lawyer from New York who
studied mathematics at Princeton in 1960--1964 \cite{Per2} but later
chose jurisprudence to be his profession.\footnote{The combination of a
professional lawyer and an amateur mathematician in one person is
not new in the history of mathematics (think of Pierre Fermat!).}

Complete tables of knots are currently known up to crossing number 16
\cite{HTW}.
For knots with 11 through 16 crossings it is nowadays customary to use
the numbering of Knotscape \cite{Knsc} where the tables are built into the
software.
For each crossing number, Knotscape has a separate list of alternating and
non-alternating knots. For example, the notation $12^a_{427}$ used in
Section \ref{orientation}, refers to the item number 427 in the list of
alternating knots with 12 crossings.

\section{Algebra of knots}
\label{alg_knots}

Denote by $\K$ the set of the equivalence classes of knots.
\label{set-of-knots}
It forms a commutative monoid (semigroup with a unit)
under the connected sum of knots,
and, therefore we can construct the monoid algebra
$\Z\K$ of $\K$.
By definition, elements of
$\Z\K$ are formal finite linear combinations $\sum\lambda_i K_i$,
$\lambda_i\in\Z$, $K_i\in\K$, the
product is defined by $(K_1,K_2) \mapsto K_1\#K_2$ on knots
and then extended by linearity to the entire space $\Z\K$. This algebra
$\Z\K$ will be referred to as the {\em algebra of knots}.
\index{Algebra!of knots}

The algebra of knots provides a convenient language for the study of
{\em knot invariants}\index{Knot!invariant} (see the next chapter):
in these terms, a knot invariant is nothing but a linear functional
on $\Z\K$. Ring homomorphisms from $\Z\K$ to some ring are referred
to as {\em multiplicative invariants};
\index{Knot!invariant!multiplicative} later, in Section
\ref{bialg_knots}, we shall see the importance of this notion.

In the sequel, we shall introduce more operations in this algebra,
as well as in the dual algebra of knot invariants. We shall also
study a filtration on $\Z\K$ that will give us the notion of a
finite type knot invariant.

\section{Tangles, string links and braids}
\label{tangles}

A {\em tangle} is a generalization of a knot which at the same time
is simpler and more complicated than a knot: on one hand, knots are
a particular case of tangles, on the other hand, knots can be
represented as combinations of (simple) tangles.

\begin{definition}\index{Tangle}
A {\em (parametrized) tangle} is a smooth embedding of a
one-dimensional compact oriented manifold, $\boldX$, possibly with
boundary, into a {\em box}
$$\{(x,y,z)\, |\ w_0\leq x\leq w_1\, ,\,  -1\leq y\leq 1\, ,\, h_0\leq z\leq h_1 \} \subset
\R^3,$$ where $w_0,w_1,h_0,h_1\in\R$, such that the boundary of
$\boldX$ is sent into the intersection of the (open) upper and lower
faces of the box with the plane $y=0$. An {\em oriented tangle} is a tangle considered up to an orientation-preserving change of parametrization; an {\em unoriented tangle} is an image of a parametrized tangle.

The boundary points of $\boldX$ are divided into the top and the
bottom part; within each of these groups the points are ordered,
say, from the left to the right. The manifold $\boldX$, with the set
of its boundary points divided into two ordered subsets, is called
the {\em skeleton}\label{skeleton} of the tangle.

The number $w_1-w_0$ is called the {\em width}, and the number
$h_1-h_0$ is the {\em height} of the tangle.
\end{definition}

Speaking of embeddings of manifolds with boundary, we mean that such
embedding send boundaries to boundaries and interiors --- to
interiors. 
Here is an example of a tangle, shown together with its box:
$$\risS{-20}{tang-in-cub}{}{90}{55}{20}
$$
Usually the boxes will be omitted in the pictures.

We shall always identify tangles obtained by translations of boxes.
Further, it will be convenient to have two notions of equivalences
for tangles. Two tangles will be called {\em fixed-end isotopic} if
one can be transformed into the other by a boundary-fixing isotopy
of its box. We shall say that two tangles are simply {\em isotopic},
or {\em equivalent} if they become fixed-end isotopic after a
suitable re-scaling of their boxes of the form
$$(x,y,z) \to (f(x),y,g(z)),$$
where $f$ and $g$ are strictly increasing functions.

\subsection{Operations}
In the case when the bottom of a tangle $T_1$ coincides with the top
of another tangle $T_2$ of the same width (for oriented tangles we
require the consistency of orientations, too), one can define the
{\em product} $T_1\cdot T_2$ by putting $T_1$ on top of $T_2$ (and,
if necessary, smoothing out the corners at the joining points):
$$
  T_1=\tanG{-10}{tang2}{}{60}; \qquad T_2=\ \tanG{-10}{tang1}{}{60};
         \qquad\qquad T_1\cdot T_2 =\ \risS{-20}{tang3}{}{60}{20}{15} \ .
$$
\index{Product!of tangles}\index{Tangle!product}

Another operation, {\em tensor product}, is defined by placing one
tangle next to the other tangle of the same height:

\smallskip

$$
T_1\ot T_2=\tanG{-10}{tang2}{}{60}\tanG{-10}{tang1}{}{60}\ .
$$
\index{Tangle!tensor product}

Both operations give rise to products on equivalence classes of
tangles. The product of two equivalence classes is defined whenever
the bottom of one tangle and the top of the other consist of the
same number of points (with matching orientations in the case of
oriented tangles), the tensor product is defined for any pair of
equivalence classes.

\subsection{Special types of tangles}\label{spec-tangles}

Knots, links and braids are particular cases of tangles. For
example, an $n$-component link is just a tangle whose skeleton is a
union of $n$ circles (and whose box is disregarded).

Let us fix $n$ distinct points $p_i$ on the top boundary of a box of
unit width and let $q_i$ be the projections of the $p_i$ to the
bottom boundary of the box. We choose the points $p_i$ (and, hence,
the $q_i$) to lie in the plane $y=0$.

\begin{xdefinition}
A {\em string link on $n$ strings} \index{String link} is an
(unoriented) tangle whose skeleton consists of $n$ intervals, the
$i$th interval connecting $p_i$ with $q_i$. A string link on one
string is called a {\em long knot} \index{Knot!long}.
\end{xdefinition}

\begin{xdefinition}
A string link on $n$ strings whose tangent vector is never
horizontal is called a {\em pure braid on $n$
strands}.\index{Braid!pure}
\end{xdefinition}

One difference between pure braids and string links is that the
components of a string link can be knotted. However, there are
string links with unknotted strands that are not equivalent to
braids.

Let $\sigma$ be a permutation of the set of $n$ elements.
\begin{xdefinition}
A {\em braid on $n$ strands} \index{Braid} is an (unoriented)
tangle whose skeleton consists of $n$ intervals, the $i$th interval
connecting $p_i$ with $q_{\sigma (i)}$, with the property that the
tangent vector to it is never horizontal.
\end{xdefinition}
Pure braids are a specific case of braids with $\sigma$ the identity
permutation. Note that with our definition of equivalence, an
isotopy between two braids can pass through tangles with points
where the tangent vector is horizontal. Often, in the definition of
the equivalence for braids it is required that that an isotopy
consist entirely of braids; the two approaches are equivalent.

The above definitions are illustrated by the following comparison
chart:
$$
  \ig[height=27mm]{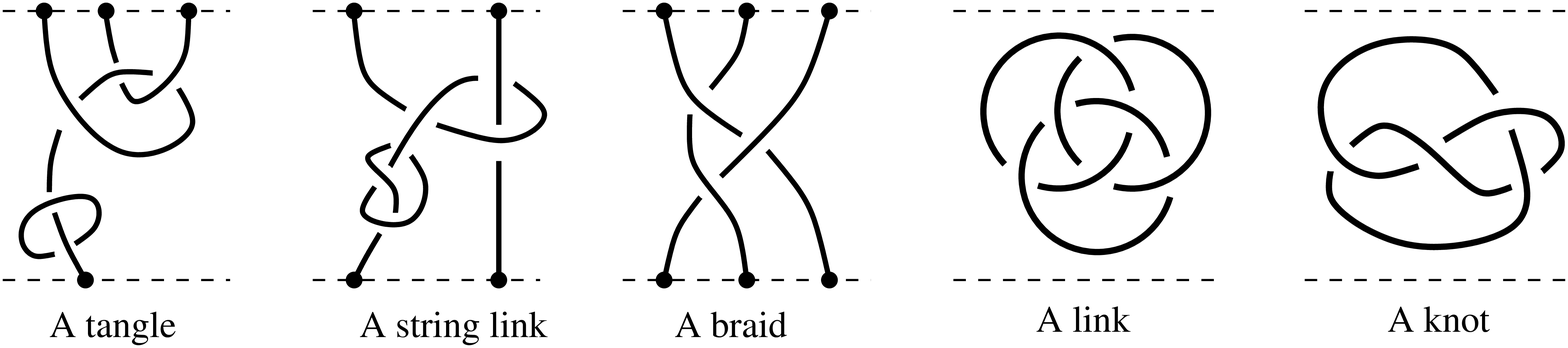}
$$

\subsection{Braids}\label{braid}
Braids are useful in the study of links, because any link can be
represented as a \textit{closure} of a braid (Alexander's theorem
\cite{Al1}\index{Theorem!Alexander}):
\begin{center}
\includegraphics[height=25mm]{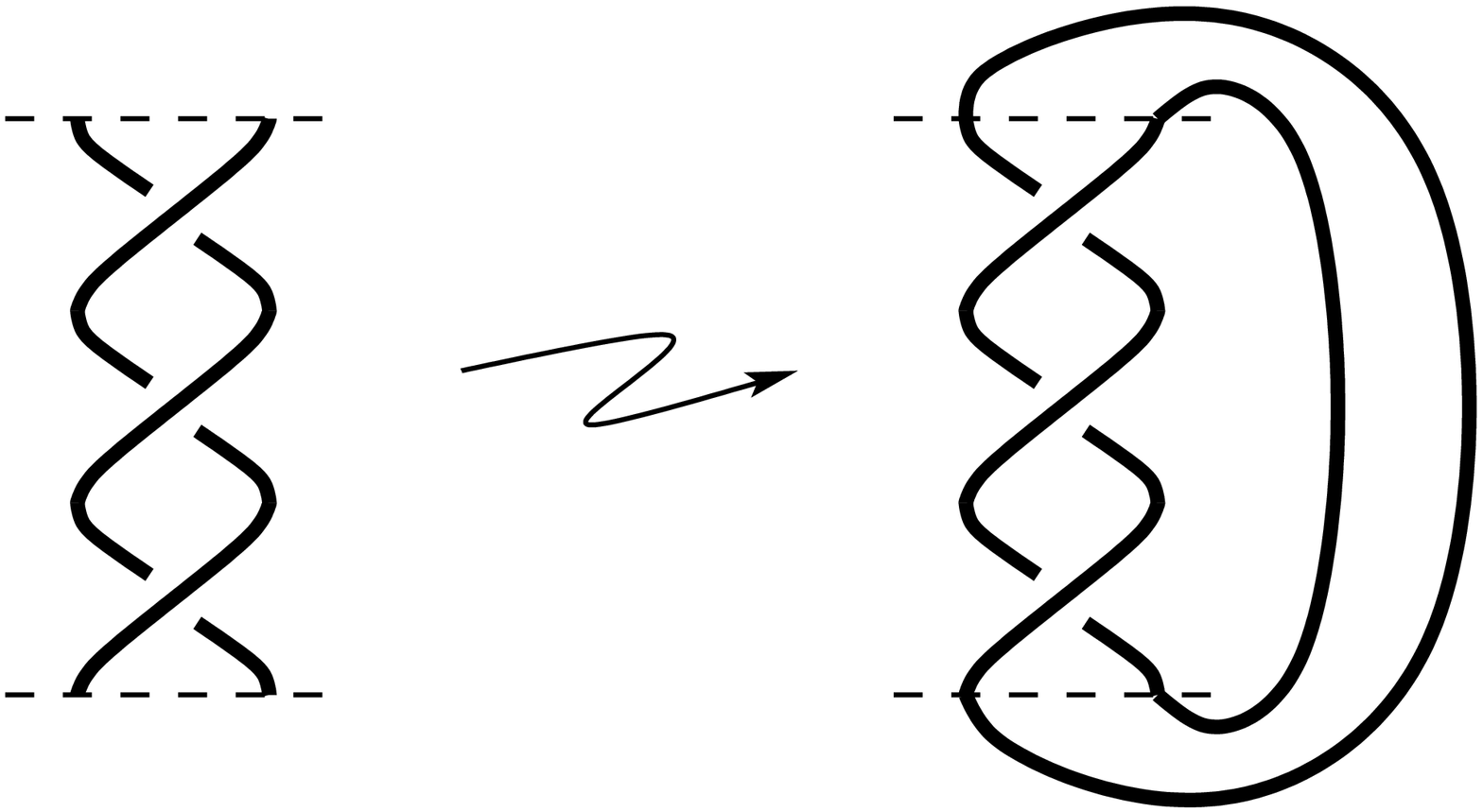}
\end{center}

Braids are in many respects easier to work with, as they form groups
under tangle multiplication: the set of equivalence classes of
braids on $n$ strands is the {\em braid group} denoted by $B_n$. A
convenient set of generators for the
group $B_n$ consists of the elements $\sigma_i$, $i=1,\dots,
n-1$:
\begin{center}\index{Braid!generators}
\ig[height=8mm]{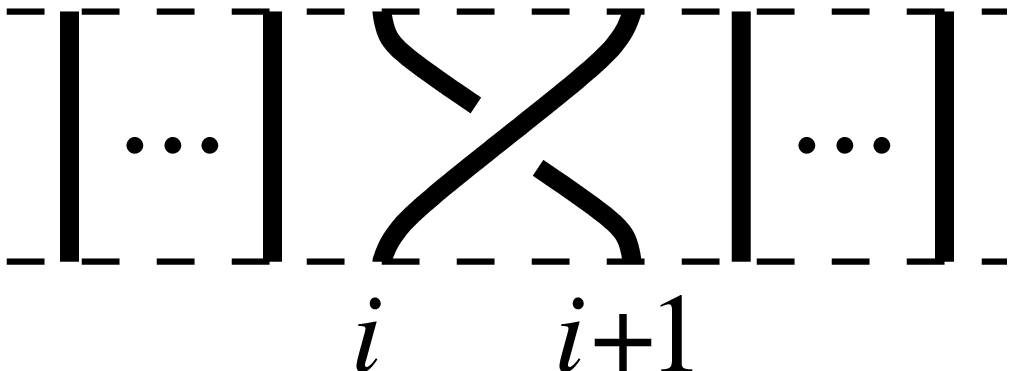}
\end{center}
which satisfy the following complete set of relations.

{\it Far commutativity},\quad\index{Braid!relations}
$\sigma_i\sigma_j=\sigma_j\sigma_i,\quad\mbox{for}\quad |i-j|>1$.
\begin{center}
\ig[height=12mm]{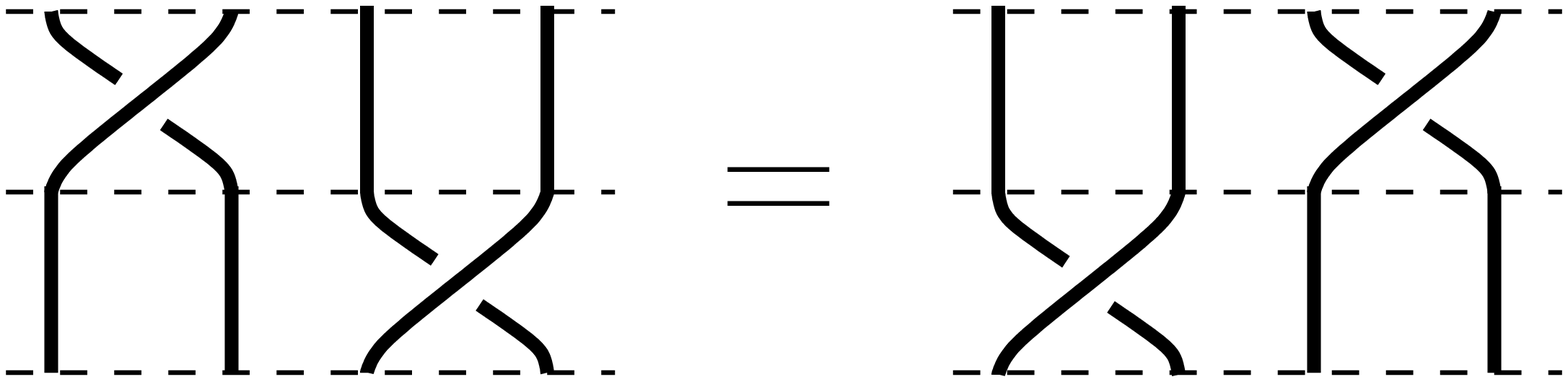}
\end{center}

{\it Braiding relation},\quad\index{Braid!relation}
$\sigma_i\sigma_{i+1}\sigma_i=\sigma_{i+1}\sigma_i\sigma_{i+1},
  \quad\mbox{for}\quad i=1,2,\dots,n-2$.
\begin{center}
\ig[height=18mm]{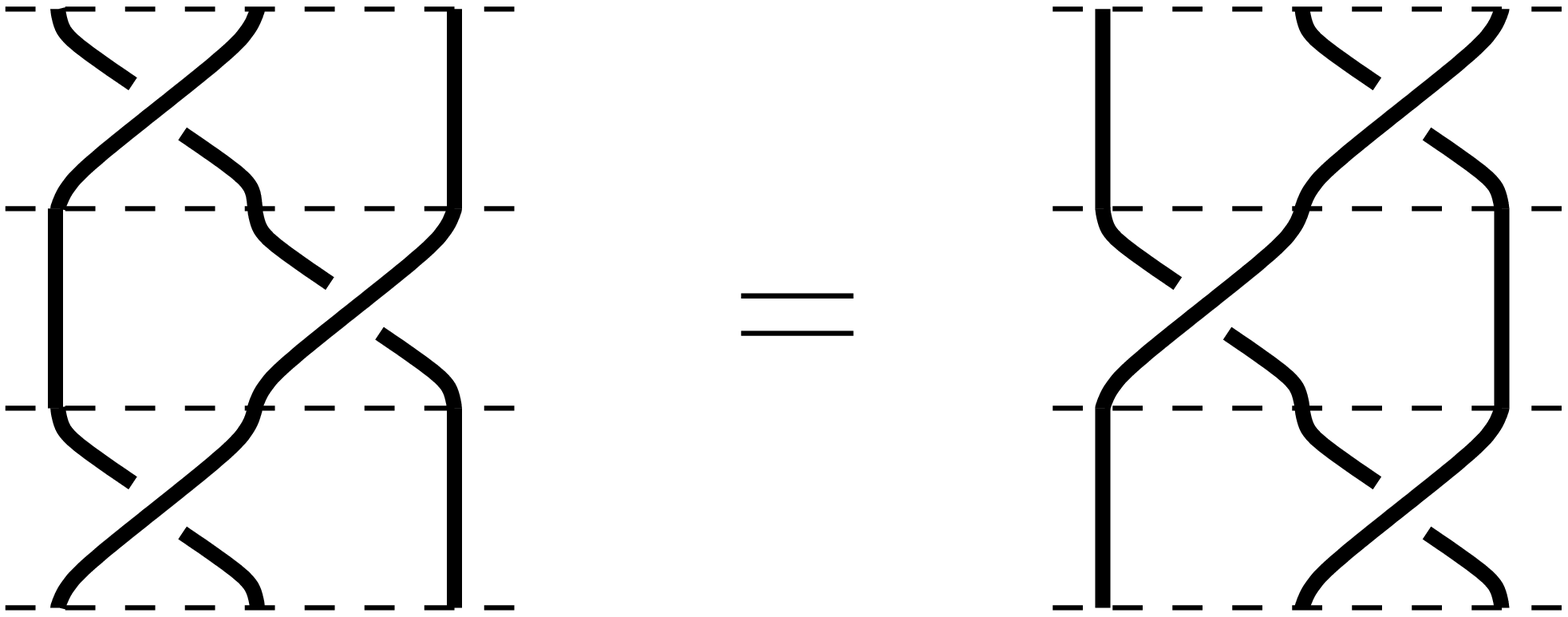}
\end{center}

Assigning to each braid in $B_n$ the corresponding permutation
$\sigma$, we get an epimorphism $\pi: B_n\to S_n$ of the braid group
on $n$ strands onto the symmetric group on $n$ letters. The kernel
of $\pi$ consists of pure braids and is denoted by $P_n$.

\begin{xtheorem}[Markov \cite{Mark,Bir1}]
\label{MarTh}\index{Theorem!Markov} Two closed braids are equivalent
(as links) if and only if the braids are related by a finite
sequence of the following Markov moves:
\index{Markov moves}\index{Moves!Markov}
\begin{itemize}
\item[(M1)] $b\longleftrightarrow aba^{-1}$ for any $a,b\in B_n$;
           \vspace{10pt}
\item[(M2)] $B_n\ni\tanG{-15}{marmo21}{\put(12,19){$b$}}{33.5}
\longleftrightarrow
\tanG{-15}{marmo22}{\put(12,19){$b$}}{49}\!\in B_{n+1}$\ ,\quad
$\tanG{-15}{marmo21}{\put(12,19){$b$}}{33.5}\longleftrightarrow
\tanG{-15}{marmo22m}{\put(12,19){$b$}}{49}$\ .\vspace{10pt}
\end{itemize}
\end{xtheorem}

\subsection{Elementary tangles}
A link 
can be cut into several {\em simple} tangles\index{Tangle!simple} by
a finite set of horizontal planes, and the link is equal to the
product of all such tangles. Every simple tangle is a tensor product
of the following {\em elementary} tangles. \index{Tangle!elementary}

\def\smtan#1{\rb{-5pt}{\ig[height=6mm]{#1.eps}}}
\noindent
Unoriented case:
$$\id:=\smtan{unort_id},\quad
  X_+:=\smtan{unort_xp},\quad
  X_-:=\smtan{unort_xm},\quad
  \max:=\smtan{unort_max},\quad
  \min:=\smtan{unort_min}\ .
$$
Oriented case:
$$\id:=\smtan{ort_id},\quad
  \id^*:=\smtan{ort_ids},\quad
  X_+:=\smtan{ort_xp},\quad
  X_-:=\smtan{ort_xm},
$$
$$\maxr:=\smtan{ort_maxr},\quad
  \maxl:=\smtan{ort_maxl},\quad
  \minr:=\smtan{ort_minr},\quad
  \minl:=\smtan{ort_minl}\ .
$$

For example, the generator $\sigma_i\in B_n$ of the braid group is a simple
tangle represented as the tensor product,
$\sigma_i = \id^{\ot(i-1)}\ot X_+\ot\id^{\ot(n-i-1)}$.

\begin{xca}
Decompose the tangle \smtan{ort_xp_d} into elementary tangles.
\end{xca}

\subsection{The Turaev moves}\label{Turmov}
Having presented a tangle as a product of simple tangles it is
natural to ask for an analogue of Reidemeister's (\ref{ReiTh}) and
Markov's (\ref{MarTh}) theorems, that is, a criterion for two such
presentations to give isotopic tangles. Here is the answer.

\newcommand\risH[6]{\rb{#1pt}[#5pt][#6pt]{\begin{picture}(30,#4)(0,0)
  \put(0,0){\ig[height=#4pt]{#2.eps}} #3
     \end{picture}}}
\begin{xtheorem}[\cite{Tur3}]\label{TurTh}\index{Theorem!Turaev}
Two products of simple tangles are isotopic if and only if they are
related by a finite sequence of the following Turaev moves.
\index{Turaev moves}

\noindent Unoriented case:
\begin{itemize}\label{Tur-mov}\index{Turaev moves!unoriented}
\index{Moves!Turaev!unoriented}
\item[(T0)]
$\risH{-12}{tmuno1-1}{\put(5,19){$\scriptstyle T_1$}
                      \put(18,5){$\scriptstyle T_2$}}{30}{15}{15}
\longleftrightarrow
\risH{-12}{tmuno1-2}{\put(5,5){$\scriptstyle T_1$}
                      \put(18,19){$\scriptstyle T_2$}}{30}{10}{10}$\qquad
\parbox{2.8in}{\scriptsize Note that the number of strands at top or bottom of either tangle $T_1$ or $T_2$, or both might be zero.}

\item[(T1)]
$\risH{-12}{trmuno1-1}{}{30}{15}{15}
\longleftrightarrow\hspace{-5pt}
\risH{-12}{tmuno2-2}{}{30}{10}{10}\hspace{-5pt}
\longleftrightarrow
\risH{-12}{trmuno1-3}{}{30}{10}{10}$\quad
$\begin{array}{l}\scriptstyle
  (\id\ot\max)\cdot(X_+\ot\id)\cdot(\id\ot\min)=\id= \\
\scriptstyle\hspace{45pt}=(\id\ot\max)\cdot(X_-\ot\id)\cdot(\id\ot\min)
\end{array}$

\item[(T2)]
$\risH{-12}{trmuno2-1}{}{30}{15}{15}
\longleftrightarrow\hspace{-5pt}
\risH{-12}{trmuno2-2}{}{30}{10}{10}\hspace{-5pt}
\longleftrightarrow
\risH{-12}{trmuno2-3}{}{30}{10}{10}$\quad
$\scriptstyle X_+\cdot X_-=\id\ot\id=X_-\cdot X_+$

\item[(T3)]
$\risH{-12}{trmuno3-1}{}{30}{15}{15}
\longleftrightarrow
\risH{-12}{trmuno3-2}{}{30}{10}{10}$\quad
$\scriptstyle (X_+\ot\id)\cdot(\id\ot X_+)\cdot(X_+\ot\id)=
  (\id\ot X_+)\cdot(X_+\ot\id)\cdot(\id\ot X_+)$

\item[(T4)]
$\risH{-12}{tmuno2-1}{}{30}{15}{15}
\longleftrightarrow\hspace{-5pt}
\risH{-12}{tmuno2-2}{}{30}{10}{10}\hspace{-5pt}
\longleftrightarrow
\risH{-12}{tmuno2-3}{}{30}{10}{10}$\quad
$\scriptstyle (\max\!\ot\id)\cdot(\id\ot\min)=\id=
 (\id\ot\max)\cdot(\min\!\ot\id)$

\item[(T5)]
$\risH{-12}{tmuno5-1}{}{30}{15}{15}
\longleftrightarrow
\risH{-12}{tmuno5-2}{}{30}{10}{10}$\qquad
$\scriptstyle (\id\ot\max)\cdot(X_+\ot\id)=
  (\max\!\ot\id)\cdot(\id\ot X_-)$

\item[(T5$'$)]
$\risH{-12}{tmuno5-3}{}{30}{15}{15}
\longleftrightarrow
\risH{-12}{tmuno5-4}{}{30}{10}{10}$\qquad
$\scriptstyle (\id\ot\max)\cdot(X_-\ot\id)=
  (\max\!\ot\id)\cdot(\id\ot X_+)$
\end{itemize}
Oriented case:
\begin{itemize}\label{Tur-mov-or}\index{Turaev moves!oriented}
\index{Moves!Turaev!oriented}
\item[(T0)]\hspace*{-5pt}
Same as in the unoriented case with arbitrary orientations of participating
 strings.

\item[(T1\ ]\hspace{-7pt}\mbox{\rm ---\ T3)}
Same as in the unoriented case with orientations of all strings from
bottom to top.\vspace{5pt}

\item[(T4)]
$\risH{-12}{tmor2-1}{}{30}{15}{15}
\longleftrightarrow\hspace{-5pt}
\risH{-12}{tmor2-2}{}{30}{10}{10}\hspace{-5pt}
\longleftrightarrow
\risH{-12}{tmor2-3}{}{30}{10}{10}$\quad
$\scriptstyle (\smaxr\!\ot\id)\cdot(\id\ot\sminr)=\id=
 (\id\ot\smaxl\!)\cdot(\sminl\ot\id)$

\item[(T4$'$)]
$\risH{-12}{tmor2-4}{}{30}{15}{15}
\longleftrightarrow\hspace{-10pt}
\risH{-12}{tmor2-5}{}{30}{10}{10}\hspace{-10pt}
\longleftrightarrow
\risH{-12}{tmor2-6}{}{30}{10}{10}$\quad
$\scriptstyle (\smaxl\!\ot\id^*\!)\cdot(\id^*\!\ot\sminl)=\id^*\!=
 (\id^*\!\ot\smaxr\!)\cdot(\sminr\ot\id^*\!)$

\item[(T5)]
$\risH{-18.5}{tmor5-1}{}{42.5}{22}{20}
\longleftrightarrow\hspace{-5pt}
\risH{-18.5}{tmor5-2}{}{42.5}{10}{10}$\quad
$\begin{array}{l}\scriptstyle
  (\smaxl\!\ot\id\ot\id^*\!)\cdot(\id^*\!\ot X_-\ot\id^*\!)\cdot
  (\id^*\!\ot\id\ot\sminl)\cdot \\
\scriptstyle\hspace{15pt}\cdot (\id^*\!\ot\id\ot\smaxr\!)\cdot
  (\id^*\!\ot X_+\ot\id^*\!)\cdot(\sminr\ot\id\ot\id^*\!)=\id\ot\id^*
\end{array}$

\item[(T5$'$)]
$\risH{-18.5}{tmor5-1pr}{}{42.5}{22}{15}
\longleftrightarrow\hspace{-5pt}
\risH{-18.5}{tmor5-2pr}{}{42.5}{10}{10}$\quad
$\begin{array}{l}\scriptstyle
  (\id^*\!\ot\id\ot\smaxr\!)\cdot
  (\id^*\!\ot X_+\ot\id^*\!)\cdot(\sminr\ot\id\ot\id^*\!)\cdot \\
\scriptstyle\hspace{15pt}\cdot(\smaxl\!\ot\id\ot\id^*\!)\cdot
      (\id^*\!\ot X_-\ot\id^*\!)\cdot(\id^*\!\ot\id\ot\sminl)=\id^*\ot\id
\end{array}$

\item[(T6)]
$\risH{-15}{tmor6-1}{}{35.5}{30}{15}\hspace{20pt}
\longleftrightarrow
\risH{-15}{tmor6-2}{}{35.5}{10}{0}$\hspace{117pt}
\rb{-20pt}{\makebox(0,0){$\begin{array}{l}\scriptstyle
  (\smaxl\!\ot\id^*\!\ot\id^*\!)\cdot
  (\id^*\!\ot\smaxl\!\ot\id\ot\id^*\!\ot\id^*\!)\cdot\\
\scriptstyle\hspace{15pt}\cdot
  (\id^*\!\ot\id^*\!\ot X_{\pm}\ot\id^*\!\ot\id^*)\cdot\\
\scriptstyle\hspace{15pt}\cdot
  (\id^*\!\ot\id^*\!\ot\id\ot\sminl\ot\id^*\!)\cdot
  (\id^*\!\ot\id^*\!\ot\sminl)\ =\vspace{5pt}\\
\scriptstyle\hspace{-10pt}=\
  (\id^*\!\ot\id^*\!\ot\smaxr\!)\cdot
  (\id^*\!\ot\id^*\!\ot\id\ot\smaxr\!\ot\id^*\!)\cdot\\
\scriptstyle\hspace{15pt}\cdot
  (\id^*\!\ot\id^*\!\ot X_{\pm}\ot\id^*\!\ot\id^*)\cdot\\
\scriptstyle\hspace{15pt}\cdot
  (\id^*\!\ot\sminr\ot\id\ot\id^*\!\ot\id^*\!)\cdot
  (\sminr\ot\id^*\!\ot\id^*\!)
\end{array}$}}

\item[(T6$'$)]
$\risH{-15}{tmor6-1pr}{}{35.5}{15}{15}\hspace{20pt}
\longleftrightarrow
\risH{-15}{tmor6-2pr}{}{35.5}{0}{10}\hspace{20pt}$
\end{itemize}

\end{xtheorem}

\section{Variations}

\subsection{Framed knots}\label{subsection:framed}

A {\em framed knot}\index{Framed knot} is a knot equipped with a
{\em framing},\index{Framing} that is, a smooth family of non-zero
vectors perpendicular to the knot. Two framings are considered as
equivalent, if one can be transformed to another by a smooth
deformation. Up to this equivalence relation, a framing is uniquely
determined by one integer: the linking number between the knot
itself and the curve formed by a small shift of the knot in the
direction of the framing. This integer, called the {\em self-linking
number}\index{Self-linking number}, can be arbitrary. The framing
with self-linking number $n$ will be called the {\em $n$-framing}
and a knot with the $n$-framing will be referred to as {\em
$n$-framed}.

One way to choose a framing is to use the {\em blackboard framing},\index{Framing!blackboard}\index{Blackboard framing}
defined by a plane knot projection, with the vector field everywhere
parallel to the projection plane, for example

\begin{center}\index{Knot!framed}
\ig[height=30mm]{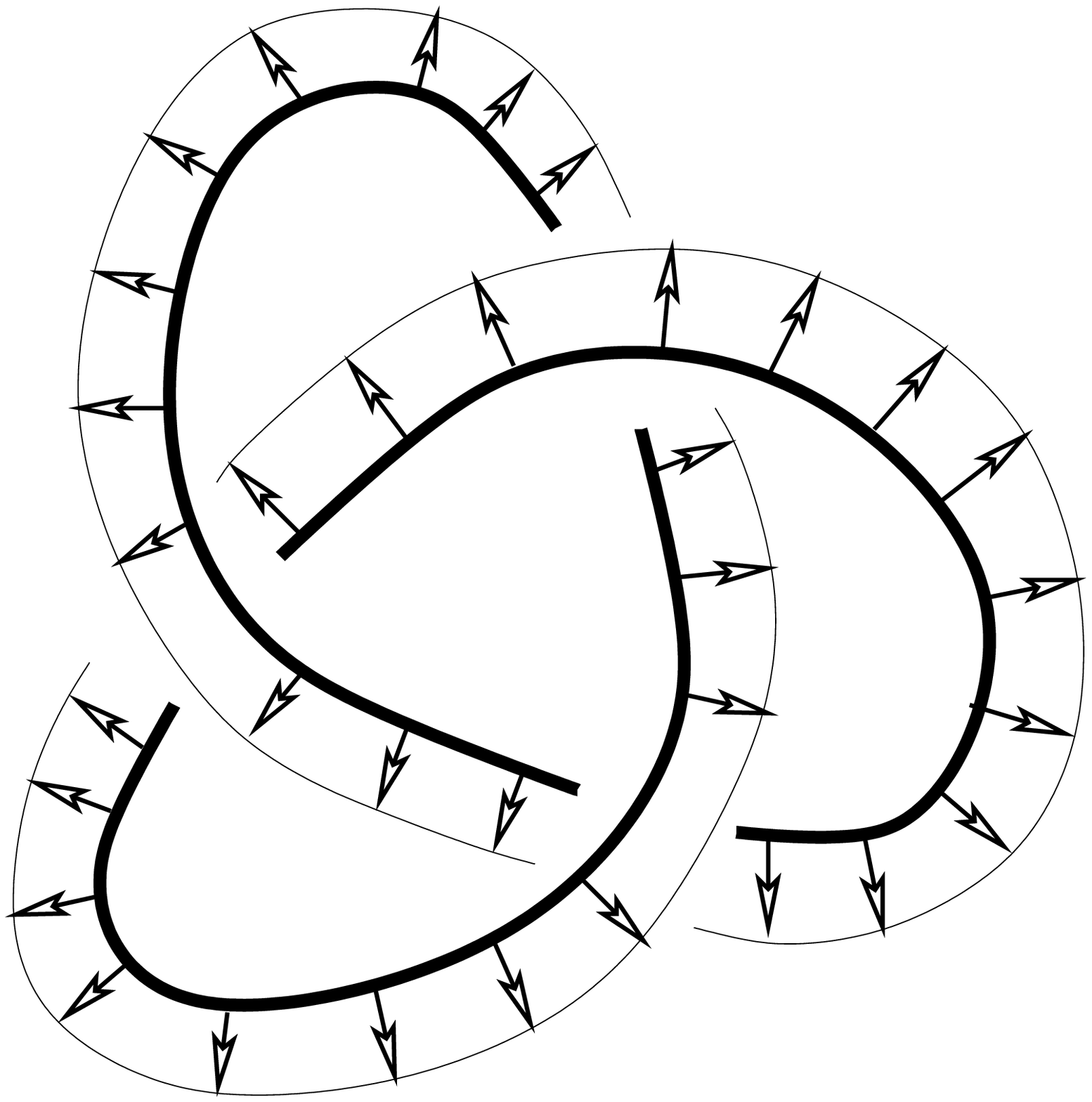}\hspace{2cm}
\ig[height=30mm]{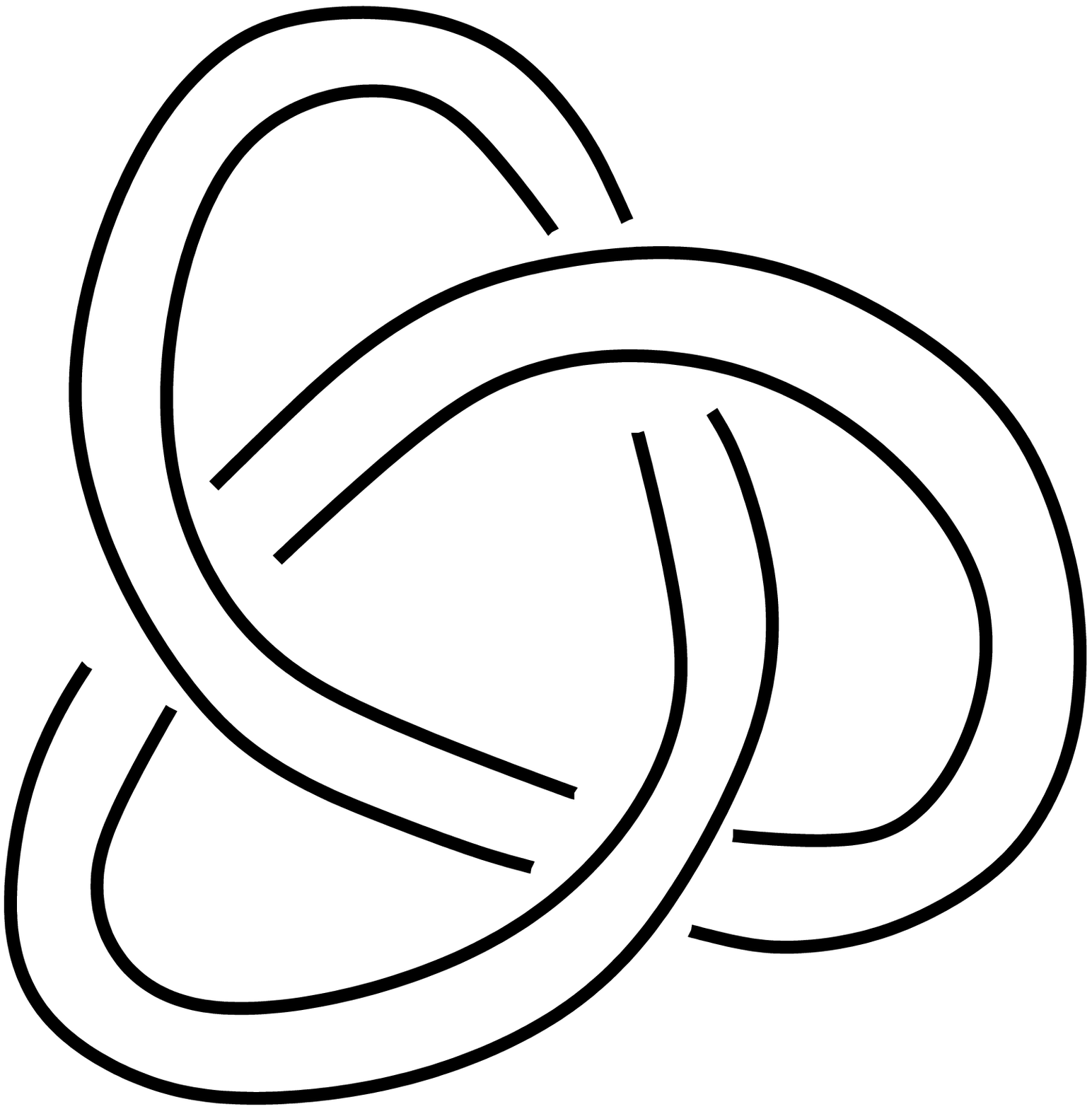}
\end{center}

A framed knot can also be visualized as a {\em ribbon knot},
\index{Knot!ribbon}\index{Ribbon knot} that is, a narrow knotted
strip (see the right picture above).

An arbitrary framed knot can be represented by a plane diagram with
the blackboard framing. This is achieved by choosing an arbitrary
projection and then performing local moves to straighten out the
twisted band:
\begin{center}
\rb{-20pt}{\ig[height=40pt]{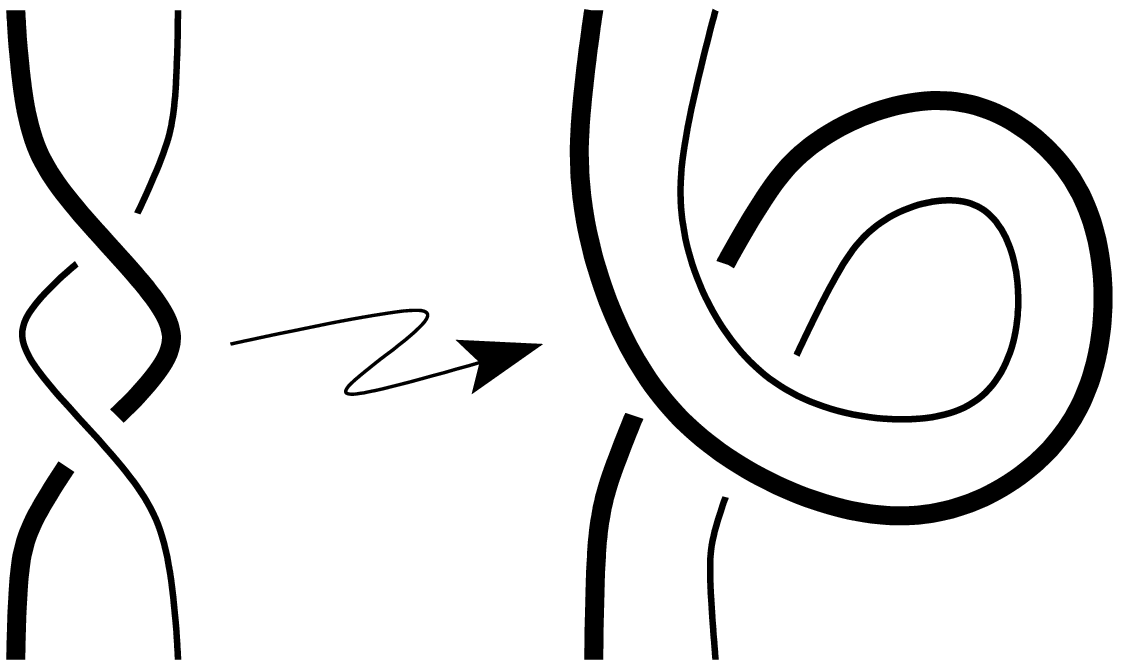}\quad ,} \hspace{1cm}
\rb{-30pt}{\ig[height=60pt]{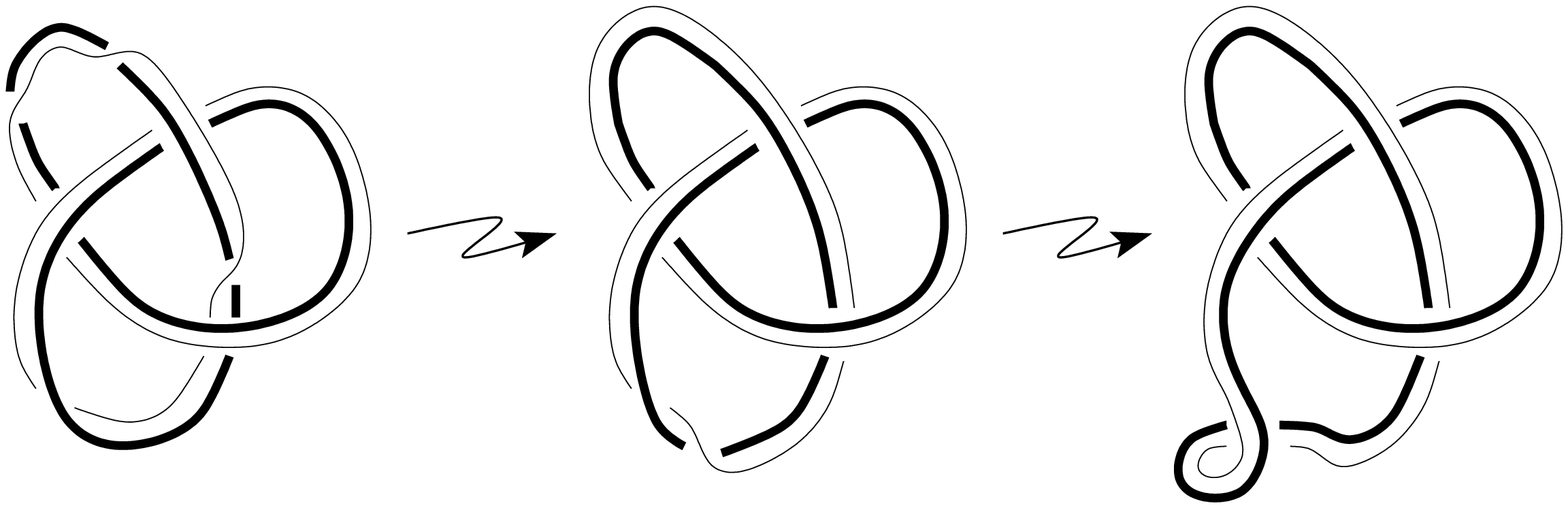}}   
\end{center}
For framed knots (with blackboard framing) the Reidemeister theorem
\ref{ReiTh} does not hold since the first Reidemeister move
$\Omega_1$ changes the blackboard framing. Here is an appropriate
substitute.
\begin{theorem}[{\bf framed Reidemeister theorem}]
\label{frReiTh}\index{Theorem!Reidemeister!framed}
\quad
Two knot diagrams with blackboard framing $D_1$ and $D_2$ are
equivalent if and only if $D_1$ can be transformed into $D_2$ by a
sequence of plane isotopies and local moves of
three types $F\Omega_1$, $\Omega_2$, and $\Omega_3$, where
$$\index{Reidemeister moves!framed}\index{Moves!Reidemeister!framed}
F\Omega_1:\qquad
\rb{-18pt}{\ig[height=40pt]{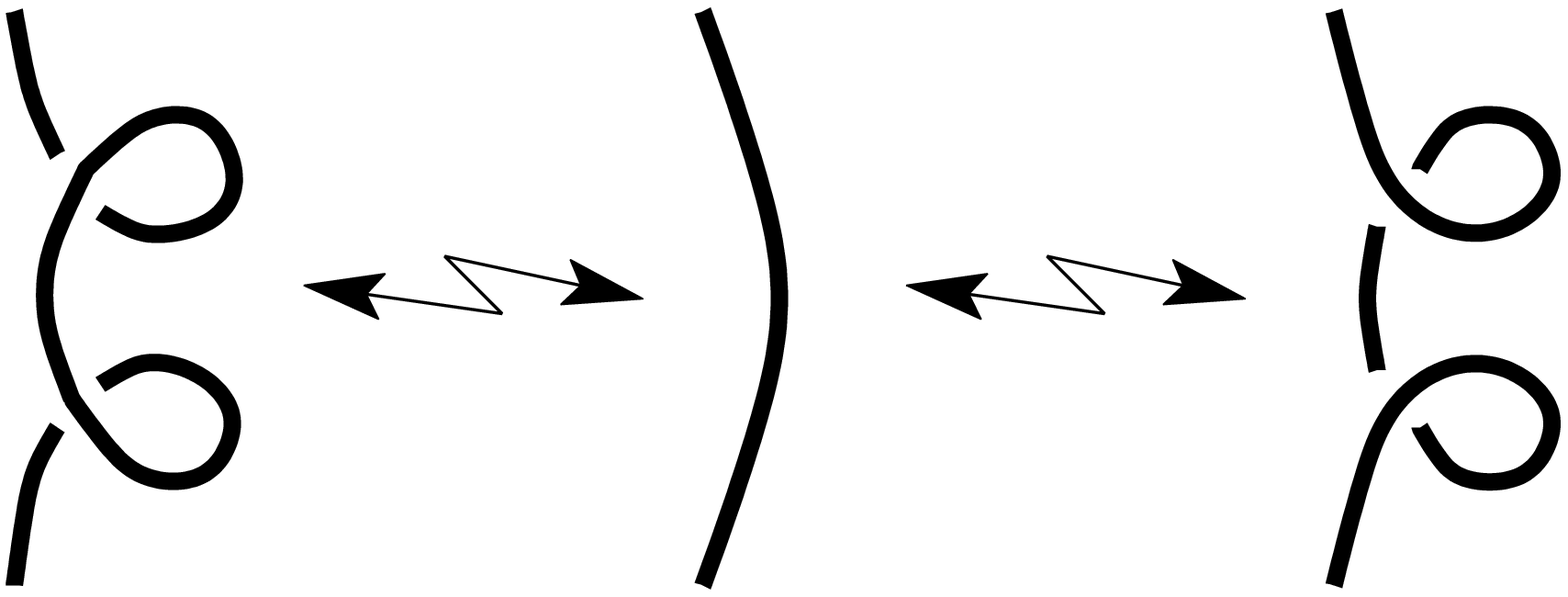}} $$
while $\Omega_2$ and $\Omega_3$ are usual Reidemeister moves defined in \ref{ReiTh}.
\end{theorem}

One may also consider framed tangles. These are defined in the same
manner as framed knots, with the additional requirement that at each
boundary point of the tangle the normal vector is equal to
$(\e,0,0)$ for some $\e>0$. Framed tangles can be represented by
tangle diagrams with blackboard framing. For such tangles there is
an analogue of Theorem \ref{TurTh} --- the Turaev move (T1) should
be replaced by its framed version that mimics the move $F\Omega_1$.
\index{Moves!Turaev!framed}\index{Turaev moves!framed}

\subsection{Long knots}
\label{long}

Recall that a long knot is a string link on one string. A long knot
can be converted into a usual knot by choosing an orientation (say,
upwards) and joining the top and the bottom points by an arc of a
sufficiently big circle. It is easy to prove that this construction
provides a one-to-one correspondence between the sets of equivalence
classes of long knots and knots, and, therefore the two theories are
isomorphic.

Some constructions on knots look more natural in the context of long
knots. For example, the cut and paste procedure for the connected
sum becomes a simple concatenation.

\subsection{Gauss diagrams and virtual knots}

Plane knot diagrams are convenient for presenting knots graphically,
but for other purposes, such as coding knots in a
computer-recognizable form, {\em Gauss diagrams} are suited better.

\begin{xdefinition}\index{Gauss diagram}
A Gauss diagram is an oriented circle with a distinguished set of
distinct points divided into ordered pairs, each pair carrying a
sign $\pm 1$.
\end{xdefinition}

Graphically, an ordered pair of points on a circle can be
represented by a chord with an arrow connecting them and pointing,
say, to the second point. Gauss diagrams are considered up to
orientation-preserving homeomorphisms of the circle. Sometimes, an
additional basepoint is marked on the circle and the diagrams are
considered up to homeomorphisms that keep the basepoint fixed. In
this case, we speak of {\em based Gauss diagrams}.

To a plane knot diagram one can associate a Gauss diagram as
follows. Pairs of points on the circle correspond to the values of
the parameter where the diagram has a self-intersection, each arrow
points from the overcrossing to the undercrossing and its sign is
equal to the local writhe at the crossing.

Here is an example of a plane knot diagram and the corresponding
Gauss diagram:
\begin{center}
\ig[width=4cm]{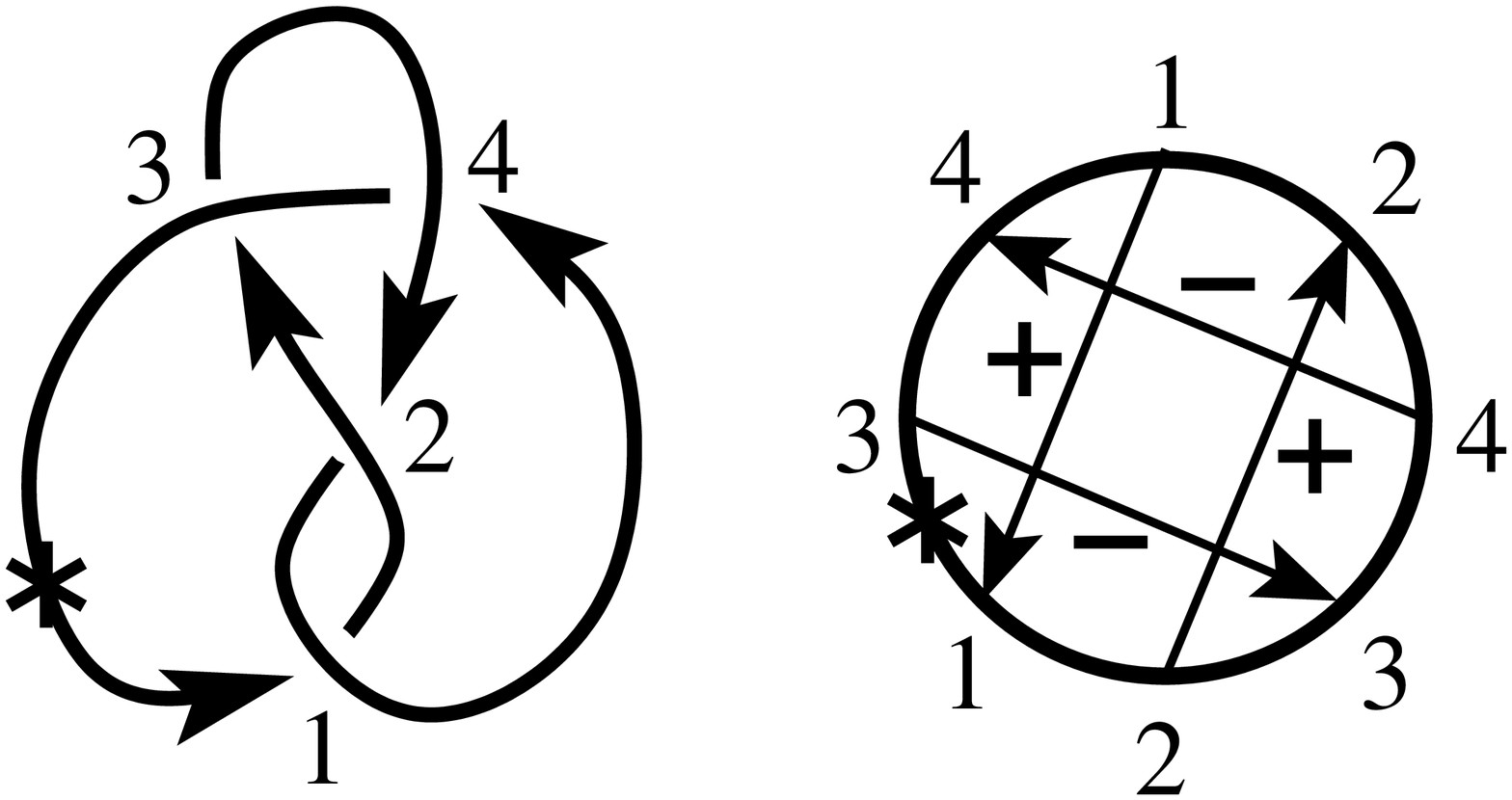}
\end{center}

\begin{xca}
What happens to a Gauss diagram, if (a) the knot is mirrored, (b) the knot
is reversed?
\end{xca}

A knot diagram can be uniquely reconstructed from the corresponding
Gauss diagram. We call a Gauss diagram {\em realizable},\index{Gauss diagram!realizable} if it comes
from a knot. Not every Gauss diagram is realizable, the simplest example being
$$\chd{gauss-nr}.$$

As we know, two oriented knot diagrams give the same knot type if and only
if they are related by a sequence of oriented Reidemeister moves. The
corresponding moves translated into the language of Gauss diagrams look as
follows:

$V\Omega_1:\qquad \risS{-20}{virrI}{
            \put(24,15){\mbox{$\scriptstyle \e$}}
            \put(180,14){\mbox{$\scriptstyle -\e$}}
             }{200}{20}{35}$

$V\Omega_2:\qquad  \risS{-20}{virrII1}{
            \put(8,18){\mbox{$\scriptstyle \e$}}
            \put(22,18){\mbox{$\scriptstyle -\e$}}
             }{120}{20}{35}\qquad\quad
\risS{-20}{virrII2}{\put(4,16){\mbox{$\scriptstyle -\e$}}
            \put(23,14){\mbox{$\scriptstyle \e$}}
              }{120}{20}{35}$

$V\Omega_3:\qquad  \risS{-20}{virrIII1}{}{120}{20}{35}\qquad\quad
\risS{-20}{virrIII2}{}{120}{20}{35}
$\label{gd_moves}
\\
In fact, the two moves $V\Omega_3$ do not exhaust all the
possibilities for representing the third Reidemeister move on
Gauss diagrams. It can be shown, however, that all the other versions of
the third move are combinations of the moves $V\Omega_2$
and $V\Omega_3$, see Exercises \ref{ex_rei_moves} -- \ref{ex_rei2_moves_cr} 
on page \pageref{ex_rei_moves} for examples and \cite{Oll} for a proof.

These moves, of course, have a geometric meaning only for
realizable diagrams. However, they make sense for all Gauss
diagrams, whether realizable or not. In particular, a realizable
diagram may be equivalent to non-realizable one:
$$\chd{cd0}\ \sim\ \risS{-10}{gauss-nr}{
   \put(3,13){\mbox{\scriptsize $-$}}\put(16,12){\mbox{\scriptsize $+$}}
    }{25}{20}{15}\,.
$$

\begin{xdefinition} A {\em virtual knot}\index{Knot!virtual} is a
Gauss diagram considered up to the Reidemeister moves $V\Omega_1$,
$V\Omega_2$, $V\Omega_3$. A {\em long}, or {\em based} virtual knot
is a based Gauss diagram, considered up to Reidemeister moves that
do not involve segments with the basepoint on them. Contrary to the case of usual knots, the theories of circular and long virtual knots differ.

\end{xdefinition}
It can be shown that the isotopy classes of knots form a subset of
the set of virtual knots. In other words, if there is a chain of
Reidemeister moves connecting two realizable Gauss diagrams, we can
always modify it so that it goes only though realizable diagrams.

Virtual knots were introduced by L.~Kauffman \cite{Ka5}. Almost at
the same time, they turned up in the work of M.~Goussarov,
M.~Polyak, O.~Viro \cite{GPV}. There are various geometric
interpretations of virtual knots. Many knot invariants are known to
extend to invariants of virtual knots.

\subsection{Knots in arbitrary manifolds}

We have defined knots as embeddings of the circle into the Euclidean
space $\R^3$. In this definition $\R^3$ can be replaced by the
3-sphere $S^3$, since the one-point compactification $\R^3\to S^3$
establishes a one-to-one correspondence between the equivalence
classes of knots in both manifolds. Going further and replacing
$\R^3$ by an arbitrary 3-manifold $M$, we can arrive to a theory of
knots in $M$ which may well be different from the usual case of
knots in $\R^3$; see, for instance, \cite{Kal,Va6}.

If the dimension of the manifold $M$ is bigger than 3, then all
knots in $M$ that represent the same element of the fundamental
group $\pi_1(M)$, are isotopic. It does not mean, however, that the
theory of knots in $M$ is trivial: the space of all embeddings
$S^1\to M$ may have non-trivial higher homology groups. These
homology groups are certainly of interest in dimension 3 too; see
\cite{Va6}. Another way of doing knot theory in higher-dimensional
manifolds is studying multidimensional knots, like embeddings
$S^2\to\R^4$, see, for example, \cite{Rol}. An analogue of knot
theory for 2-manifolds is Arnold's theory of immersed curves
\cite{Ar3}.

\ \\

\begin{xcb}{Exercises}
\begin{enumerate}

\item Find the following knots in the knot table (page \pageref{knot_table}):
\begin{center}
\begin{tabular}{ccccccccc}
\raisebox{12mm}{(a)} &
\includegraphics[width=14.5mm,height=14.5mm]{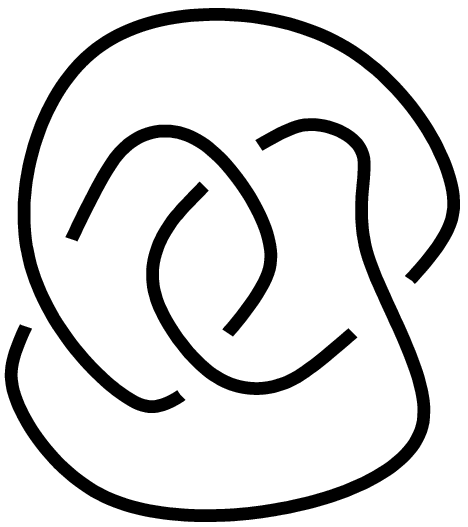}&\qquad&
\raisebox{12mm}{(b)\!} &
\includegraphics[width=14.5mm,height=14.5mm]{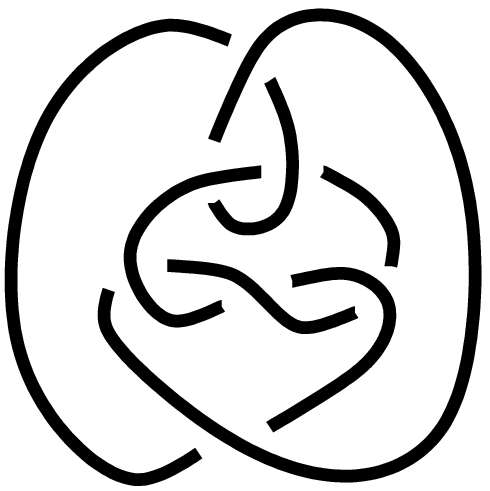}&\qquad&
\raisebox{12mm}{(c)\!} &
\includegraphics[width=14.5mm,height=14.5mm]{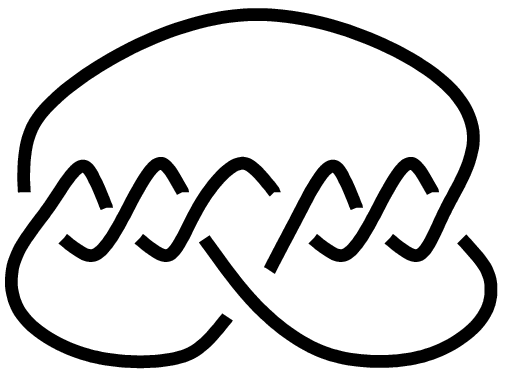}&\qquad
\end{tabular}
\end{center}

\item
Can you find the following links in the picture on page \pageref{link_ex}?
\begin{center}
\ig[height=20mm]{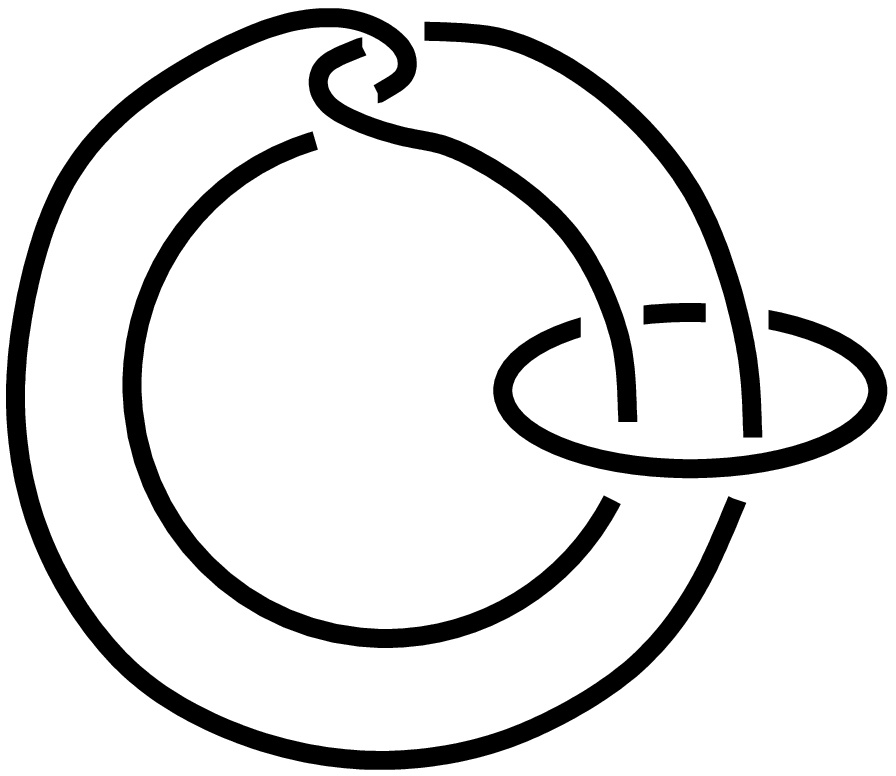}\qquad
\ig[height=20mm]{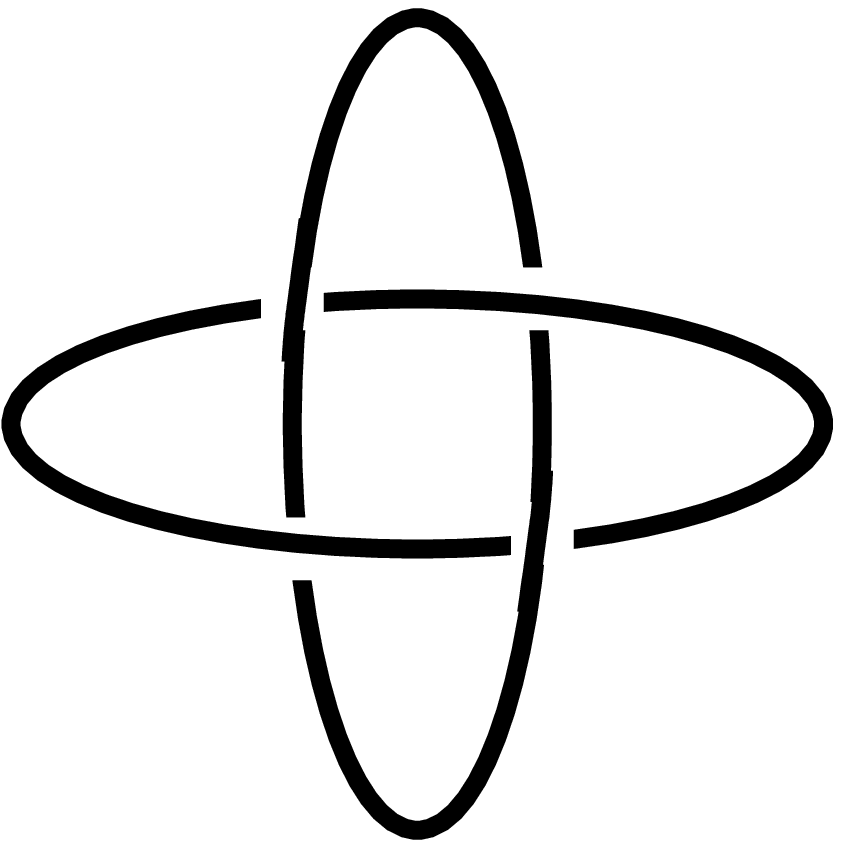}
\end{center}

\item
Borromean rings (see page \pageref{link_ex}) have the property that after deleting any component
the remaining two-component link becomes trivial. Links with such property are
called {\em Brunnian}.\index{Brunnian link}\index{Link!Brunnian} Find a Brunnian link with 4 components.

\item
Table \ref{knot_table} shows 35 topological types of knots up to
change of orientation and taking the mirror images. How many
distinct knots do these represent?

\item
 Find an isotopy that transforms the knot $6_{3}$
   into its mirror image $\overline{6_{3}}$.

\item Repeat Perko's achievement:
find an isotopy that transforms one of the knots of the Perko
pair into another one.

\item
Let $G_n$ be the {\em Goeritz diagram} \cite{Goer} with $2n+5$ crossings, as in the figure below.\index{Goeritz diagram}\index{Unknot!Goeritz}\vspace{-10pt}\\
\begin{enumerate}
\item \parbox[t]{1.8in}{
Show that $G_n$ represents a trivial knot.}\quad
\rb{-10pt}{$G_n=\ \risS{-20}{unknot-gn}{
         \put(8,36){\mbox{\scriptsize $n$ crossings}}
         \put(78,36){\mbox{\scriptsize $n\!+\!1$ crossings}}}{135}{0}{8}$}

\item \parbox[t]{1.8in}{Prove that for $n\geqslant 3$ in any}\\
 sequence of the Reidemeister moves\\ transforming
$G_n$ into the plane circle there is an intermediate knot diagram with more than $2n+5$ crossings.
\item Find a sequence of 23 Reidemeister moves (containing the $\Omega_1$ move 5 times, the $\Omega_2$ move 7 times, and the $\Omega_3$ move 11 times)
transforming $G_3$ into the plane circle. See the picture of $G_3$ in \ref{unknots} on page \pageref{unknots-diag}.
\end{enumerate}

\item
\parbox[t]{3.5in}{Decompose the knot on the right into a connected sum of
prime knots.}\quad
$\risS{-12}{k3_1+4_1}{}{35}{10}{15}$

\item \label{ex_cr_change1}
Show that by changing some crossings from overcrossing to
undercrossing or vice versa, any knot diagram can be transformed
into a diagram of the unknot.

\item ({\it C.~Adams \cite{AdC}})
Show that by changing some crossings from overcrossing to
undercrossing or vice versa, any knot diagram can be made
alternating.

\item
Represent the knots $4_1$, $5_1$, $5_2$ as closed braids.

\item
Analogously to the braid closure, one can define the closure of a
string link. Represent the Whitehead link and the Borromean rings
from Section \ref{def:link} (page \pageref{def:link}) as closures of
string links on 2 and 3 strands respectively.

\item
Find a sequence of Markov moves that transforms the closure of the braid
\def\s{\sigma}
$\s_1^2\s_2^3\s_1^4\s_2$ into the closure of the braid
$\s_1^2\s_2\s_1^4\s_2^3$.

\item
\parbox[t]{3.5in}{Garside's {\it fundamental braid}
\index{Braid!fundamental} $\D\in B_n$
is defined as
$\D:=(\s_1\s_2\dots\s_{n-1})(\s_1\s_2\dots\s_{n-2})\dots
      (\s_1\s_2)(\s_1)\ .
$}
\qquad\qquad
\rb{6pt}{\makebox(0,0){$\D=\risS{-18}{delta5}{}{35}{20}{15}$}}
\begin{enumerate}
\item Prove that $\s_i\D=\D\s_{n-i}$ for every standard generator
   $\s_i\in B_n$.
\item Prove that $\D^2=(\s_1\s_2\dots\s_{n-1})^n$.
\item Check that $\D^2$ belongs to the centre $Z(B_n)$ of the braid group.
\item Show that any braid can be represented as a product
of a certain
  power (possibly negative) of $\D$ and a {\em positive braid}, that is, a braid that
  contains only positive powers of standard generators $\s_i$.
\end{enumerate}
In fact, for $n\ge3$, the centre $Z(B_n)$ is the infinite cyclic
group generated by $\D^2$. The word and conjugacy problems in the
braid group were solved by F.~Garside \cite{Gar}. The structure of
positive braids that occur in the last statement was studied in
\cite{Adya,ElMo}.

\item
\begin{enumerate}
\item
Prove that the sign of the permutation corresponding to a braid $b$
is equal to the parity of the number of crossings of $b$, that is
$(-1)^{\ell(b)}$, where $\ell(b)$ is the length of $b$ as a word in
generators $\s_1,\dots,\s_{n-1}$.
\item \parbox[t]{3.2in}{Prove that
the subgroup $P_n$ of pure braids is generated by the braids
$A_{ij}$ linking the $i$th and $j$th strands with each other
behind all other strands.} \qquad\qquad
\rb{-12pt}{\makebox(0,0){$A_{ij}=\risS{-22}{braid_Aij}{}{30}{20}{15}$}}
\end{enumerate}
\smallskip

\item\label{rep_br_cub}
Let $V$ be a vector space of dimension $n$ with a distinguished
basis $e_1,\dots,e_n$, and let $\Xi_i$ be the counterclockwise
$90^\circ$ rotation in the plane $\langle e_i, e_{i+1}\rangle$:
$\Xi_i(e_i)=e_{i+1}$,  $\Xi_i(e_{i+1})=-e_i$, $\Xi_i(e_j)=e_j$ for
$j\ne i,i+1$. Prove that sending each elementary generator $\s_i\in
B_n$ to $\Xi_i$ we get a representation $B_n\to GL_n(\R)$ of the
braid group.

\item{\bf Burau representation.}\label{rep_burau}\index{Burau representation}
Consider the free module over the ring of Laurent polynomials
$\Z[x^{\pm1}]$ with a basis $e_1,\dots,e_n$. The {\it Burau
representation} $B_n\to GL_n(\Z[x^{\pm1}])$ sends $\s_i\in B_n$ to
the linear operator that transforms $e_i$ into $(1-x)e_i+e_{i+1}$, and
$e_{i+1}$ into $x e_i$.
\begin{enumerate}
\item Prove that it is indeed a representation of the braid group.
\item The Burau representation is reducible. It splits into the trivial
one-dimensional representation and an $(n-1)$-dimensional
irreducible representation which is called the {\em reduced Burau
representation}. \index{Burau representation!reduced}

Find a basis of the reduced Burau representation where the matrices have the form
$$
\s_1\mapsto\left(\begin{smallmatrix}
-x&x&\ldots&0\\
 0&1&\ldots&0\vspace{-4pt}\\
 \vdots&\vdots&\ddots&\vdots\\
  0&0&\ldots&1
  \end{smallmatrix}\right),\quad
\s_i\mapsto\left(\begin{smallmatrix}
1\hspace{-3pt}&      & &  & &      &\vspace{-8pt}\\
 &\ddots\hspace{-3pt}& &  & &      &\\
 &      &1& 0&0&      &\\
 &      &1&-x&x&      &\\
 &      &0& 0&1&      &\vspace{-8pt}\\
 &      & &  & &\hspace{-4pt}\ddots\hspace{-4pt}&\vspace{-2pt}\\
 &      & &  & &      &1
 \end{smallmatrix}\right),\quad
\s_{n-1}\mapsto\left(\begin{smallmatrix}
1&\ldots&0&0\vspace{-4pt}\\
\vdots&\ddots&\vdots&\vdots\\
0&\ldots&1&0\\
0&\ldots&1&-x
  \end{smallmatrix}\right)
$$
{\it Answer.} $\{xe_1-e_2, xe_2-e_3,\dots,xe_{n-1}-e_n\}$
\end{enumerate}
The Burau representation is faithful for $n\leqslant 3$ \cite{Bir1},
and not faithful for $n\geqslant 5$ \cite{Big1}. The case $n=4$ remains open.

\item{\bf Lawrence--Krammer--Bigelow representation.}
\label{rep_Kr_Big}\index{Lawrence--Krammer--Bigelow representation}
Let $V$ be a free $\Z[q^{\pm1},t^{\pm1}]$ module of dimension
$n(n-1)/2$ with a basis $e_{i,j}$ for $1\leqslant i<j\leqslant n$.
The {\it Lawrence--Krammer--Bigelow representation} can be defined
via the action of $\s_k\in B_n$ on $V$: {\scriptsize
$$\s_k(e_{i,j})=
\begin{cases}
e_{i,j} & \text{if $k<i-1$ or $k>j$,}\\
e_{i-1,j}+(1-q)e_{i,j} & \text{if $k=i-1$,}\\
tq(q-1)e_{i,i+1}+qe_{i+1,j} & \text{if $k=i<j-1$,}\\
tq^2e_{i,j} & \text{if $k=i=j-1$,}\\
e_{i,j}+tq^{k-i}(q-1)^2e_{k,k+1} & \text{if $i<k<j-1$,}\\
e_{i,j-1}+tq^{j-i}(q-1)e_{j-1,j} & \text{if $i<k=j-1$,}\\
(1-q)e_{i,j}+qe_{i,j+1} & \text{if $k=j$.}
\end{cases}
$$}
Prove that this assignment determines a representation of the braid group.
It was shown in \cite{Big2,Kram} that this representation is faithful for
any $n\geqslant1$. Therefore the braid group is a linear group.

\item
Represent the knots $4_1$, $5_1$, $5_2$ as
products of simple tangles.

\item
Consider the following two knots given as products of simple tangles:
$$\scriptstyle
(\smaxl\ot\smaxr)\cdot
(\id^*\ot X_+\ot\id^*)\cdot (\id^*\ot X_+\ot\id^*)\cdot
(\id^*\ot X_+\ot\id^*)\cdot(\sminr\ot\sminl) $$
and
$$\scriptstyle
\smaxr\cdot(\id\ot\smaxr\ot\id^*)\cdot
(X_+\ot\id^*\ot\id^*)\cdot
(X_+\ot\id^*\ot\id^*)\cdot
(X_+\ot\id^*\ot\id^*)\cdot
(\id\ot\sminl\ot\id^*)\cdot\sminl$$

\begin{enumerate}
\item Show that these two knots are equivalent.
\item Indicate a sequence of the Turaev moves that transforms one product
into another.
\item Forget about the orientations and consider the corresponding
unoriented tangles. Find a sequence of unoriented Turaev moves that
transforms one product into another.
\end{enumerate}

\item
\parbox[t]{3in}{Represent the oriented tangle move on the right as
a sequence of oriented Turaev moves from page \pageref{Tur-mov-or}.} \qquad
\rb{-12pt}{$\risH{-18.5}{tmor5-3}{}{42.5}{22}{25}
\longleftrightarrow\hspace{-5pt}
\risH{-18.5}{tmor5-2pr}{}{42.5}{10}{10}$}

\item \label{Wh_trick}\index{Whitney trick}
\parbox[t]{3.2in}{{\bf Whitney trick.}
Show that the move $F\Omega_1$ in the framed Reidemeister Theorem
\ref{frReiTh} can be replaced by the move shown on the
right.}\qquad\qquad\quad
\rb{-10pt}{\makebox(0,0){$\ig[height=40pt]{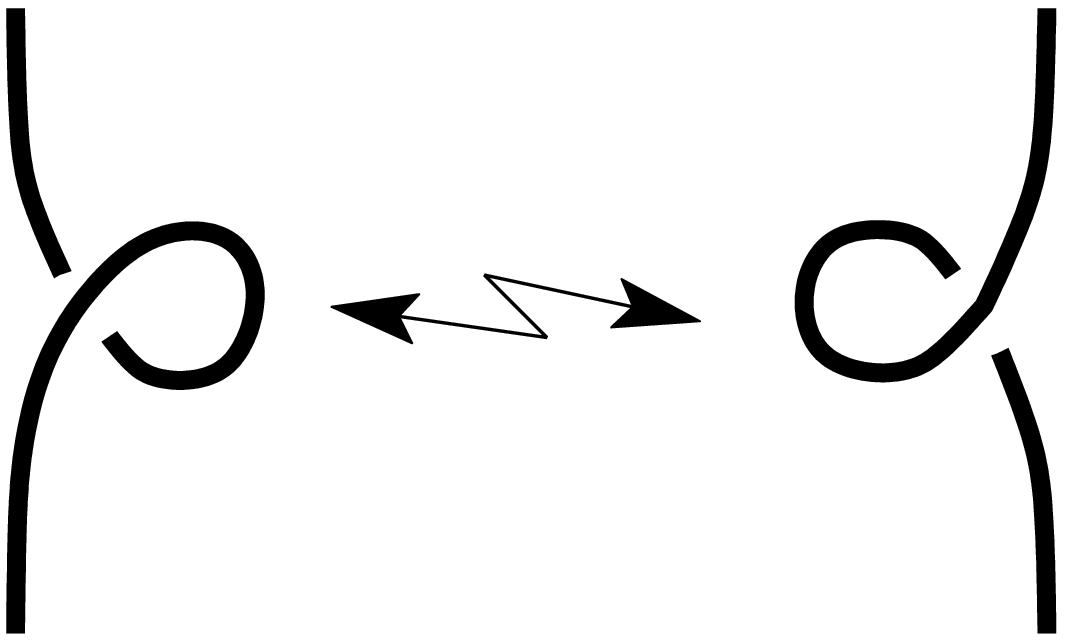}$}}

\item
The group $\Z_2^{k+1}$ acts on oriented $k$-component links,
changing the orientation of each component and taking the mirror
image of the link. How many different links are there in the orbit
of an oriented Whitehead link under this action?

\item \label{ex_rei_moves}
Show that each of the moves $V\Omega_3$ can be obtained as a combination of
the moves $V\Omega_2$ with the moves $V\Omega_3'$ below:
$$V\Omega_3':\quad  \risS{-18}{virrIII4}{}{120}{20}{25}\qquad\quad
\risS{-18}{virrIII3}{}{120}{0}{0}\ .
$$
Conversely, show that the moves $V\Omega_3'$ can be obtained as combinations of
the moves $V\Omega_2$ and $V\Omega_3$.

\item \label{ex_rei3_moves}
Show that the following moves are equivalent modulo $V\Omega_2$.
$$\risS{-18}{virrIII3}{}{120}{20}{25}\qquad\qquad
\risS{-18}{virrIII5}{}{120}{0}{0}\ .
$$
This means that either one can be obtained as a combination of another one
with the $V\Omega_2$ moves.

\item \label{ex_rei2_moves_cr}
({\em O.-P.~\"{O}stlund \cite{Oll}}) Show that the second version of
$V\Omega_2$:
$$\risS{-20}{virrII2}{\put(4,16){\mbox{$\scriptstyle -\e$}}
            \put(23,14){\mbox{$\scriptstyle \e$}}
              }{120}{20}{25}
$$ is redundant. It can be obtained as a combination of the first version,
$$\risS{-20}{virrII1}{
            \put(8,18){\mbox{$\scriptstyle \e$}}
            \put(22,18){\mbox{$\scriptstyle -\e$}}
             }{120}{20}{25}\ ,
$$
with the moves $V\Omega_1$  and $V\Omega_3$.

\item \label{min_rei_moves}
({\em M.Polyak \cite{Po3}}) Show that the following moves
\\
$V\Omega_1:\ \risS{-20}{virrI}{
            \put(24,15){\mbox{$\scriptstyle \e$}}
            \put(180,14){\mbox{$\scriptstyle -\e$}}
             }{200}{20}{35}$
\\
$\hspace*{-10pt} V\Omega_2^{\uparrow\downarrow}:\  \risS{-20}{virrII1}{
            \put(8,18){\mbox{$\scriptstyle \e$}}
            \put(22,18){\mbox{$\scriptstyle -\e$}}
             }{120}{20}{35}\qquad
V\Omega_3^{+++}:\  \risS{-18}{virrIII4}{}{120}{20}{25}$\\
are sufficient to generate all Reidemeister moves $V\Omega_1$,
$V\Omega_2$, $V\Omega_3$.

\end{enumerate}
\end{xcb}
 %1 knots
\chapter{Knot invariants} % 02
\label{kn_inv}

Knot invariants are functions of knots that do not change under
isotopies. The study of knot invariants is at the core of knot
theory; indeed, the isotopy class of a knot is, tautologically, a
knot invariant.

\section{Definition and first examples}

Let $\K$ be the set of all equivalence classes of knots.

\begin{xdefinition} \index{Knot invariant}
A {\em knot invariant with values in a set $S$} is a map from $\K$
to $S$.
\end{xdefinition}
\noindent In the same way one can speak of invariants of links, framed knots, etc.

\subsection{Crossing number}
Any knot can be represented by a plane diagram in infinitely many ways.

\begin{xdefinition}
The {\em crossing number} \index{Crossing number}$c(K)$ of a knot
$K$ is the minimal number of crossing points in a plane diagram of
$K$.
\end{xdefinition}

\noindent
\textbf{Exercise.}
Prove that if $c(K)\le2$, then the knot $K$ is trivial.
\medskip

It follows that the minimal number of crossings required to draw a
diagram of a nontrivial knot is at least 3. A little later we shall
see that the trefoil knot is indeed nontrivial.

Obviously, $c(K)$ is a knot invariant taking values in the set of
non-negative integers.

\subsection{Unknotting number}

Another integer-valued invariant of knots which admits a simple
definition is the unknotting number.

Represent a knot by a plane diagram. The diagram can be transformed
by plane isotopies, Reidemeister moves and crossing changes:
\begin{center}
\ig[height=10mm]{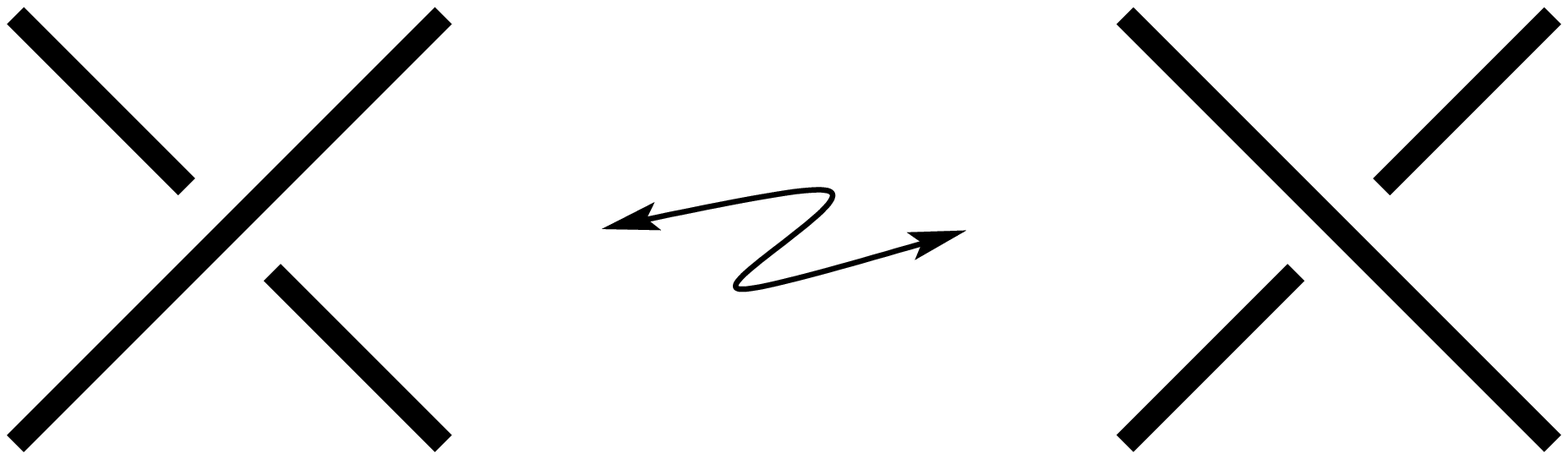}
\end{center}
As we know, modifications of the first two kinds preserve the
topological type of the knot, and only crossing switches can change
it.

\begin{xdefinition}
The {\em unknotting number} \index{Unknotting number} $u(K)$ of a
knot $K$ is the minimal number of crossing changes in a plane
diagram of $K$ that convert it to a trivial knot, provided that any
number of plane isotopies and Reidemeister moves is also allowed.
\end{xdefinition}

\noindent
\textbf{Exercise.}
What is the unknotting number of the knots $3_1$ and $8_3$?
\medskip

Finding the unknotting number, if it is greater than 1, is a
difficult task; for example, the second question of the previous
exercise was answered only in 1986 (by T.~Kanenobu and H.~Murakami).

\subsection{Knot group}
The {\em knot group} is the fundamental group of the complement to
the knot in the ambient space: $\pi(K)=\pi_1(\R^3\sm K)$. The knot
group is a very strong invariant. For example, a knot is trivial if
and only if its group is infinite cyclic. More generally, two prime
knots with isomorphic fundamental groups are isotopic. For a
detailed discussion of knot groups see \cite{Lik}.

\medskip

\noindent
\textbf{Exercise.}
Prove that
\begin{enumerate}
\setlength{\itemsep}{1pt plus 1pt minus 1pt}
\item
the group of the trefoil is generated by two elements $x$, $y$
with one relation $x^2=y^3$;
\item
this group is isomorphic to the braid group $B_3$
(in terms of $x,y$ find another pair of generators $a,b$
that satisfy $aba=bab$).
\end{enumerate}

\section{Linking number}
\index{Linking number}\label{link_num}

The {\em linking number} is an example of a Vassiliev invariant of
two-component links; it has an analog for framed knots, called {\em self-linking number}.

Intuitively, the linking number $lk(A,B)$ of two oriented spatial
curves $A$ and $B$ is the number of times $A$ winds around $B$. To
give a precise definition, choose an oriented disk $D_A$ immersed in
space so that its oriented boundary is the curve $A$ (this means
that the ordered pair consisting of an outward-looking normal vector
to $A$ and the orienting tangent vector to $A$ gives a positive
basis in the tangent space to $D_A$). The linking number $lk(A,B)$
is then defined as the intersection number\index{Intersection
number} of $D_A$ and $B$. To find the intersection number, if
necessary, make a small perturbation of $D_A$ so as to make it meet
the curve $B$ only at finitely many points of transversal
intersection. At each intersection point, define the sign to be
equal to $\pm1$ depending on the orientations of $D_A$ and $B$ at
this point. More specifically, let $(e_1,e_2)$ be a positive pair of
tangent vectors to $D_A$, while $e_3$ a positively directed tangent
vector to $B$ at the intersection point; the sign is set to $+1$ if
and only if  the frame $(e_1,e_2,e_3)$ defines a positive
orientation of $\R^3$. Then the linking number $lk(A,B)$ is the sum
of these signs over all intersection points $p\in D_A\cap B$. One
can prove that the result does not depend on the choice of the
surface $D_A$ and that $lk(A,B)=lk(B,A)$.

\noindent
{\bf Example.} The two curves shown in the picture
\begin{center}
\ig[width=25mm]{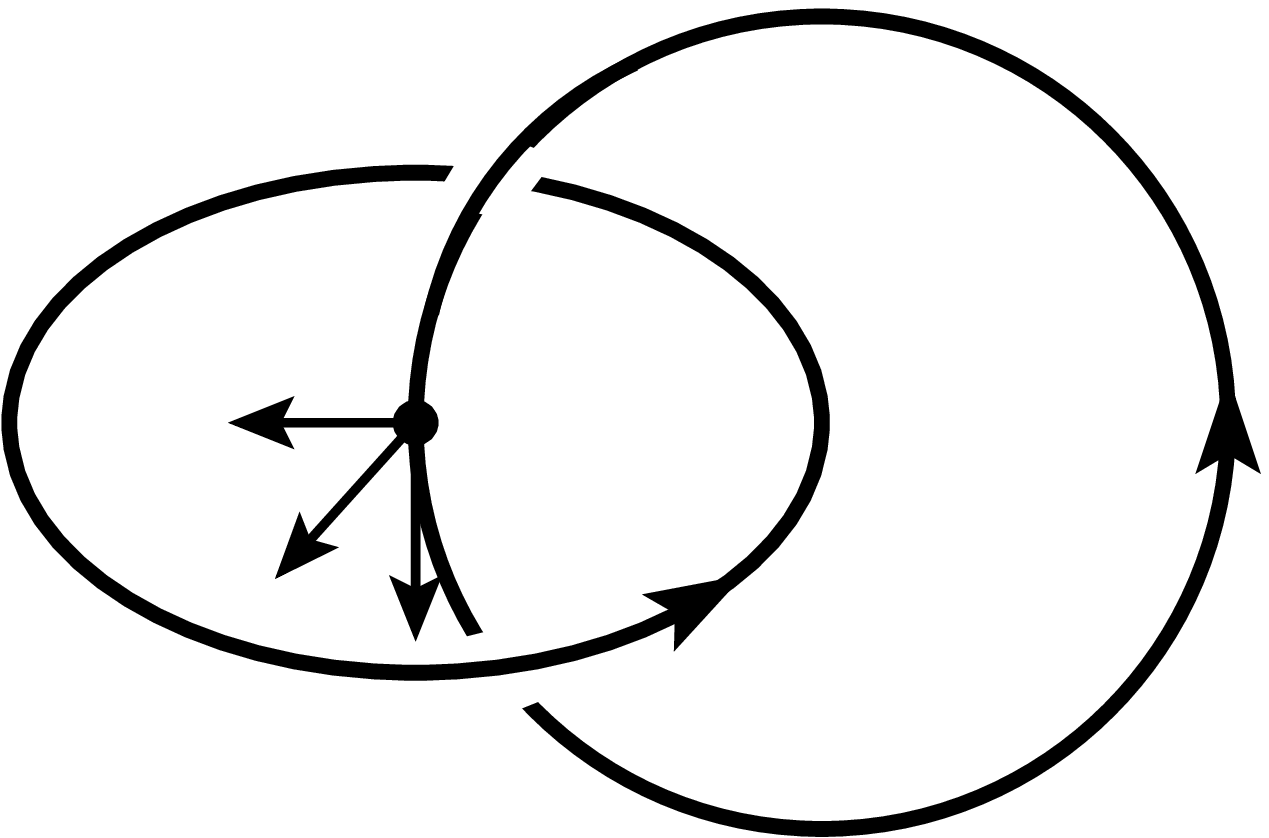}
\end{center}
have their linking number equal to $-1$.

Given a plane diagram of a two-component link, there is a simple
combinatorial formula for the linking number. Let $I$ be the set of
crossing points involving branches of both components $A$ and $B$
(crossing points involving branches of only one component are
irrelevant here). Then $I$ is the disjoint union of two subsets
$I_B^A$ (points where $A$ passes over $B$) and $I_A^B$ (where $B$
passes over $A$).

\begin{proposition}
$$
  lk(A,B)=\sum_{p\in I_B^A} w(p)=\sum_{p\in I_A^B} w(p)
=\frac{1}{2}\sum_{p\in I}w(p)
$$
where $w(p)=\pm 1$ is the local writhe of the crossing point.
\label{comb_lk_num}
\end{proposition}

\begin{proof}
Crossing changes at all points $p\in I_A^B$ make the two components
unlinked. Call the new curves $A'$ and $B'$, then $lk(A',B')=0$. It
is clear from the pictures below that each crossing switch changes
the linking number by $-w$ where $w$ is the local writhe:
\begin{center}
\rb{-12mm}{\ig[height=25mm]{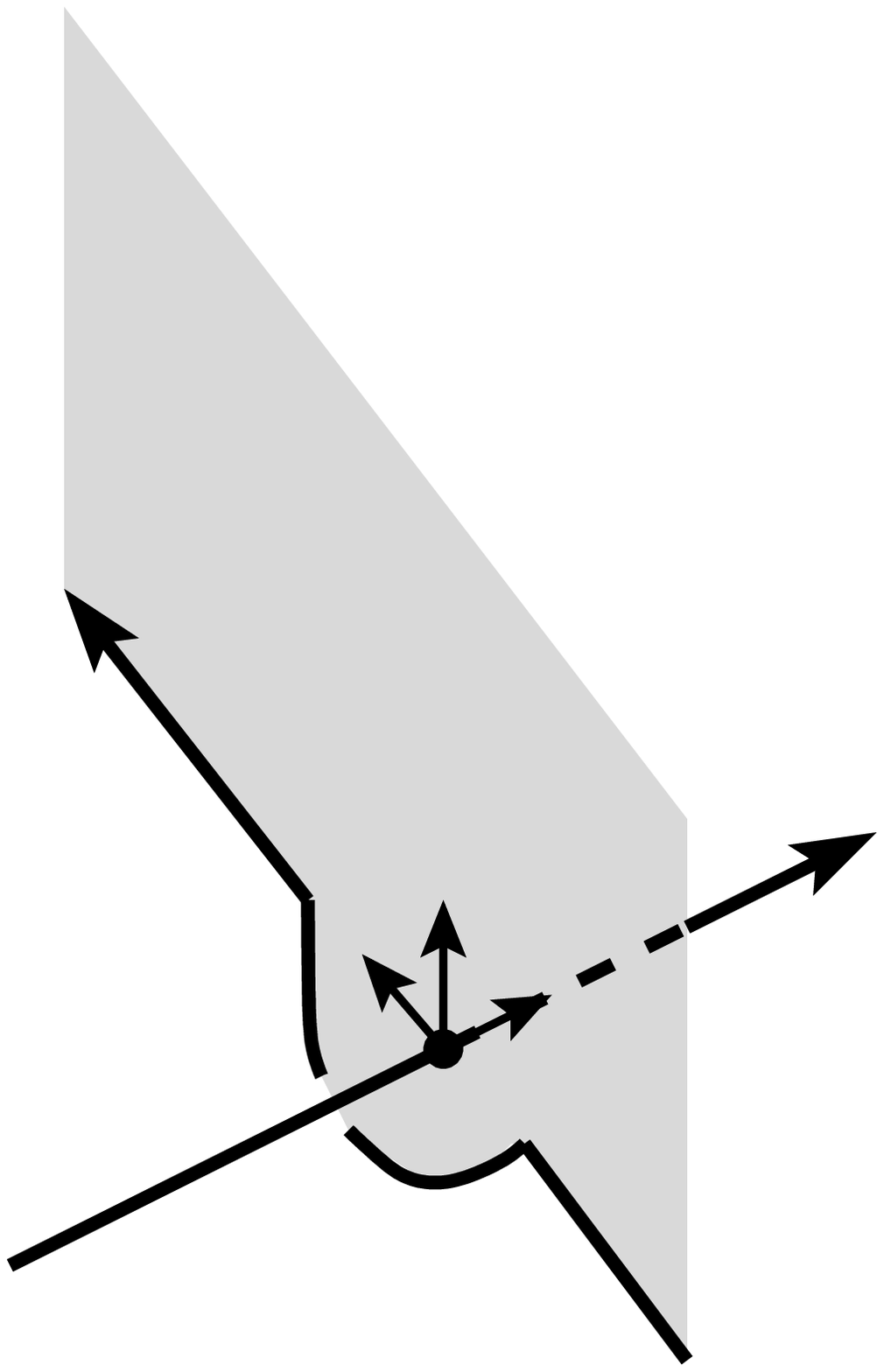}}\qquad
\rb{-2mm}{\ig[height=4mm]{totor.eps}}\qquad
\rb{-12mm}{\ig[height=25mm]{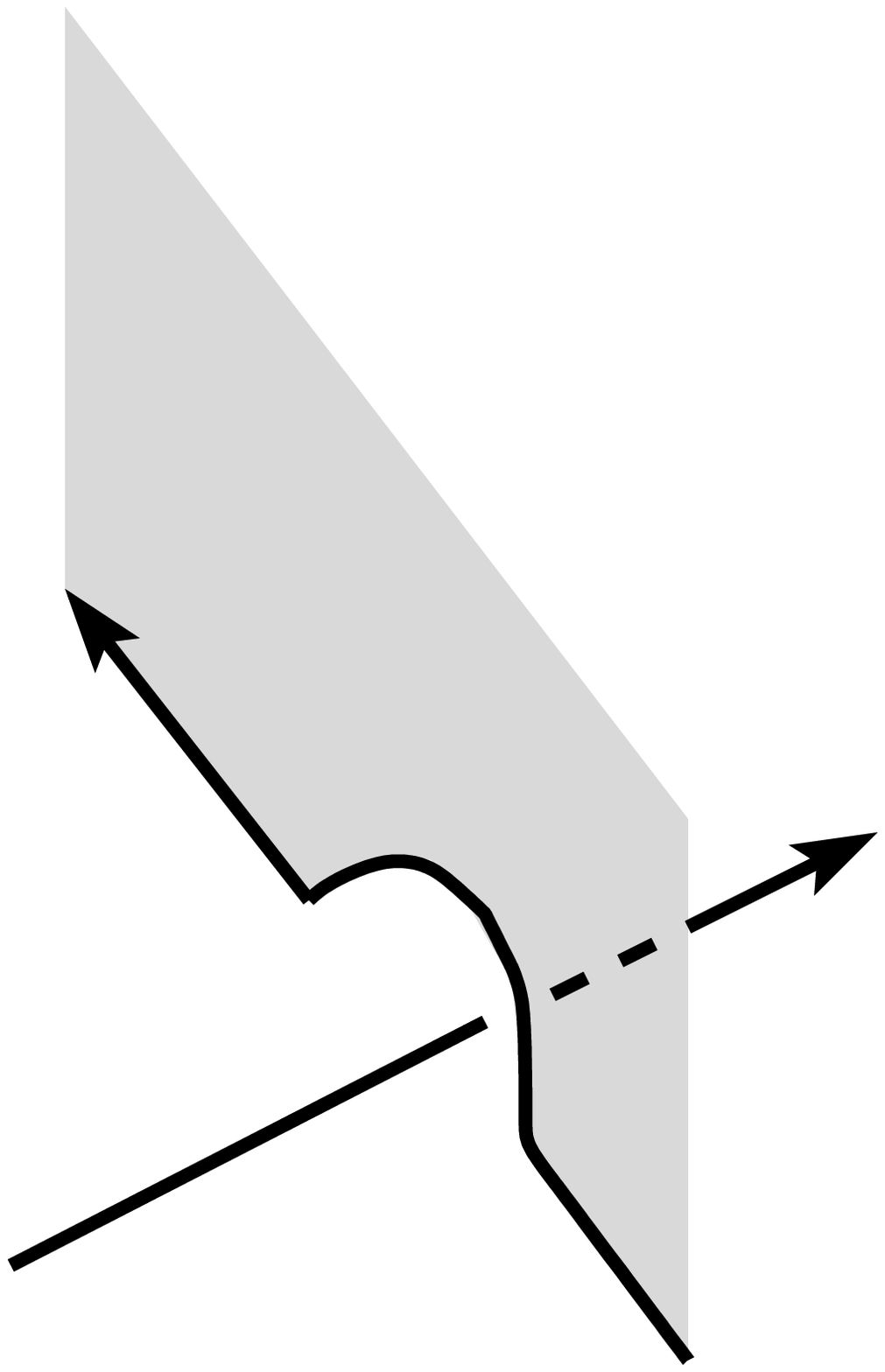}}
\end{center}
\begin{center}
\rb{-12mm}{\ig[height=25mm]{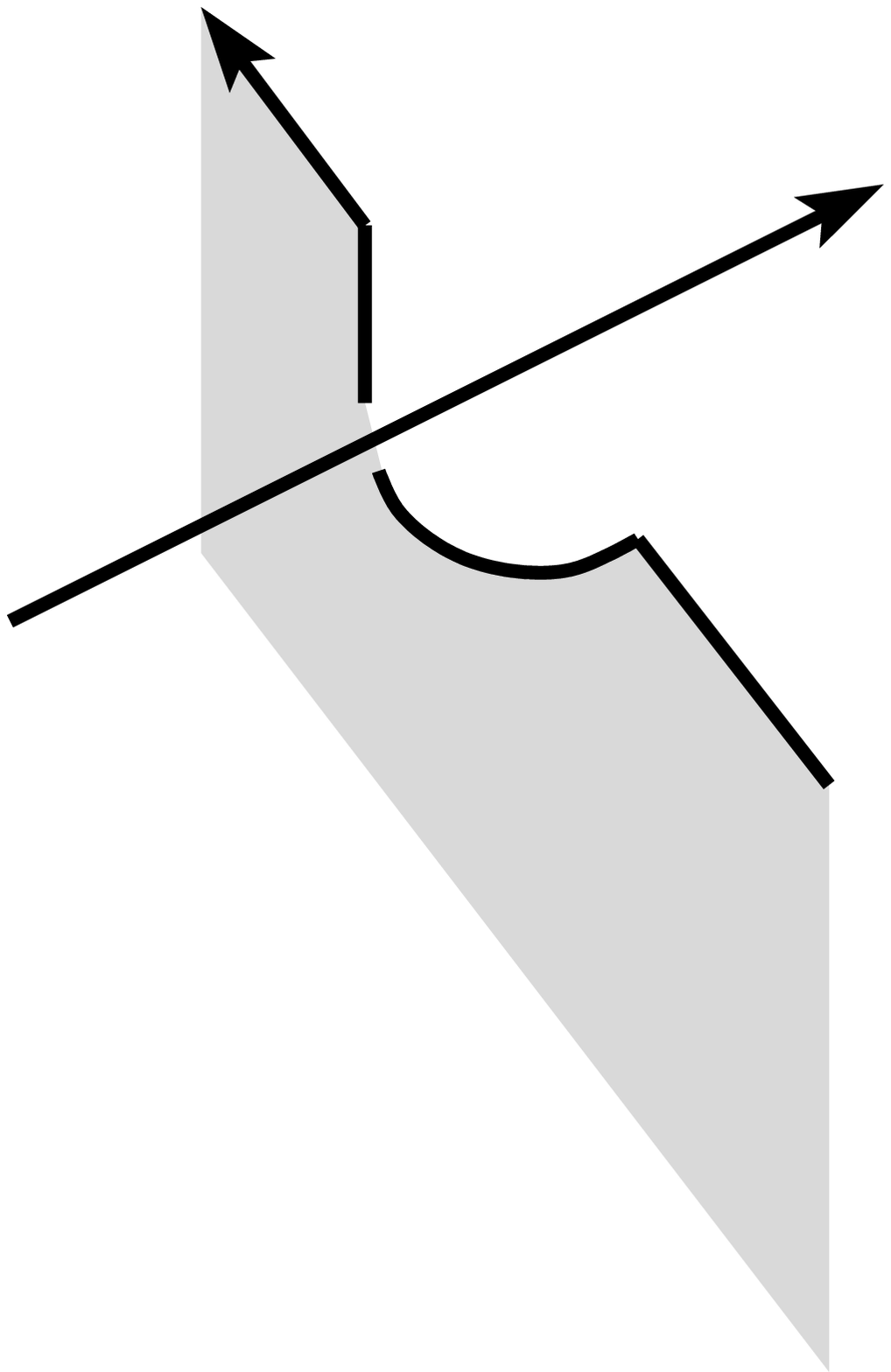}}\qquad
\rb{-2mm}{\ig[height=4mm]{totor.eps}}\qquad
\rb{-12mm}{\ig[height=25mm]{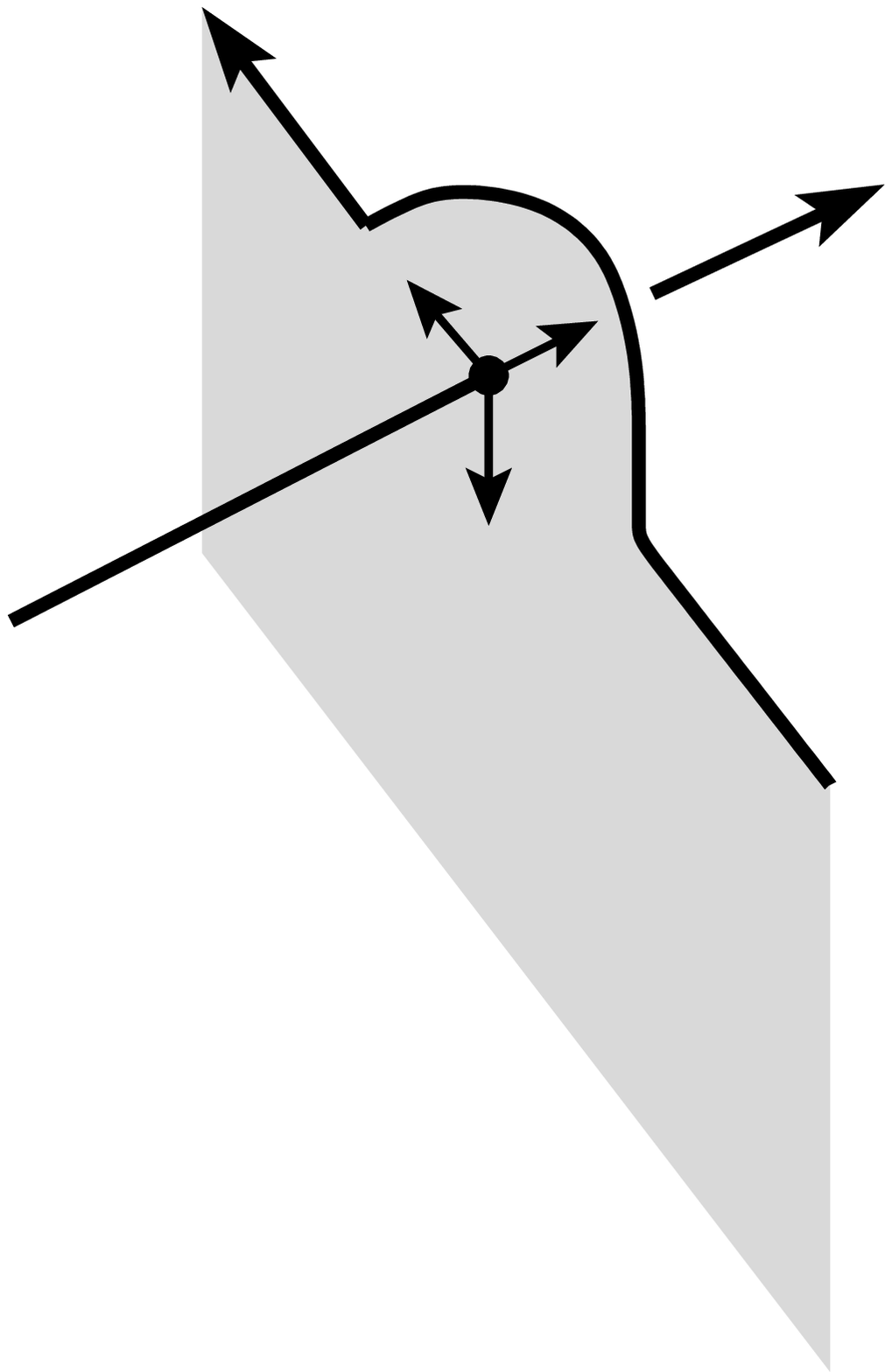}}
\end{center}
Therefore, $lk(A,B)-\sum_{p\in I_A^B}w(p)=0$,
and the assertion follows.
\end{proof}

\noindent {\bf Example}. For the two curves below both ways to
compute the linking number give $+1$:
$$\risS{-30}{v2proof3}{\put(-2,45){$A$}\put(27,-3){$B$}}{60}{18}{35}
$$

\subsection{Integral formulae}

There are various integral formulae for the linking number.
The most famous formula was found by Gauss (see
\cite{Spi} for a proof).

\begin{xtheorem}\label{gauss_f}
Let $A$ and $B$ be two non-intersecting curves in $\R^3$,
parameterized, respectively, by the smooth functions $\a,\b:S^1\to\R^3$. Then
$$
  lk(A,B)=\frac{1}{4\pi}
  \int_{S^1\times S^1} \frac{(\b(v)-\a(u),du,dv)}{|\b(v)-\a(u)|^3},
$$
where the parentheses in the numerator stand for the mixed product
of 3 vectors.
\end{xtheorem}

Geometrically, this formula computes the degree of the Gauss map
from $A\times B = S^1\times S^1$ to the 2-sphere $S^2$, that is, the
number of times the normalized vector connecting a point on $A$ to a
point on $B$ goes around the sphere.

A different integral formula for the linking number will be stated
and proved in Chapter \ref{chapKI}, see page
\pageref{link_num_formula}. It represents the simplest term of the
{\em Kontsevich integral}, which encodes all Vassiliev invariants.

\subsection{Self-linking}
\label{self_link}
\index{Self-linking number}

Let $K$ be a framed knot and let $K'$ be the knot obtained from $K$
by a small shift in the direction of the framing.

\begin{xdefinition}
The {\em self-linking number} of $K$ is the linking number of $K$ and $K'$.
\end{xdefinition}

Note, by the way, that the linking number is the same if $K$ is
shifted in the direction, {\em opposite} to the framing.

\begin{xproposition}
The self-linking number of a framed knot given by a diagram $D$
with blackboard framing is equal to the total writhe of the diagram $D$.
\end{xproposition}

\begin{proof}
Indeed, in the case of blackboard framing, the only crossings of $K$
with $K'$ occur near the crossing points of $K$. The neighbourhood
of each crossing point looks like
\begin{center}
\ig[height=20mm]{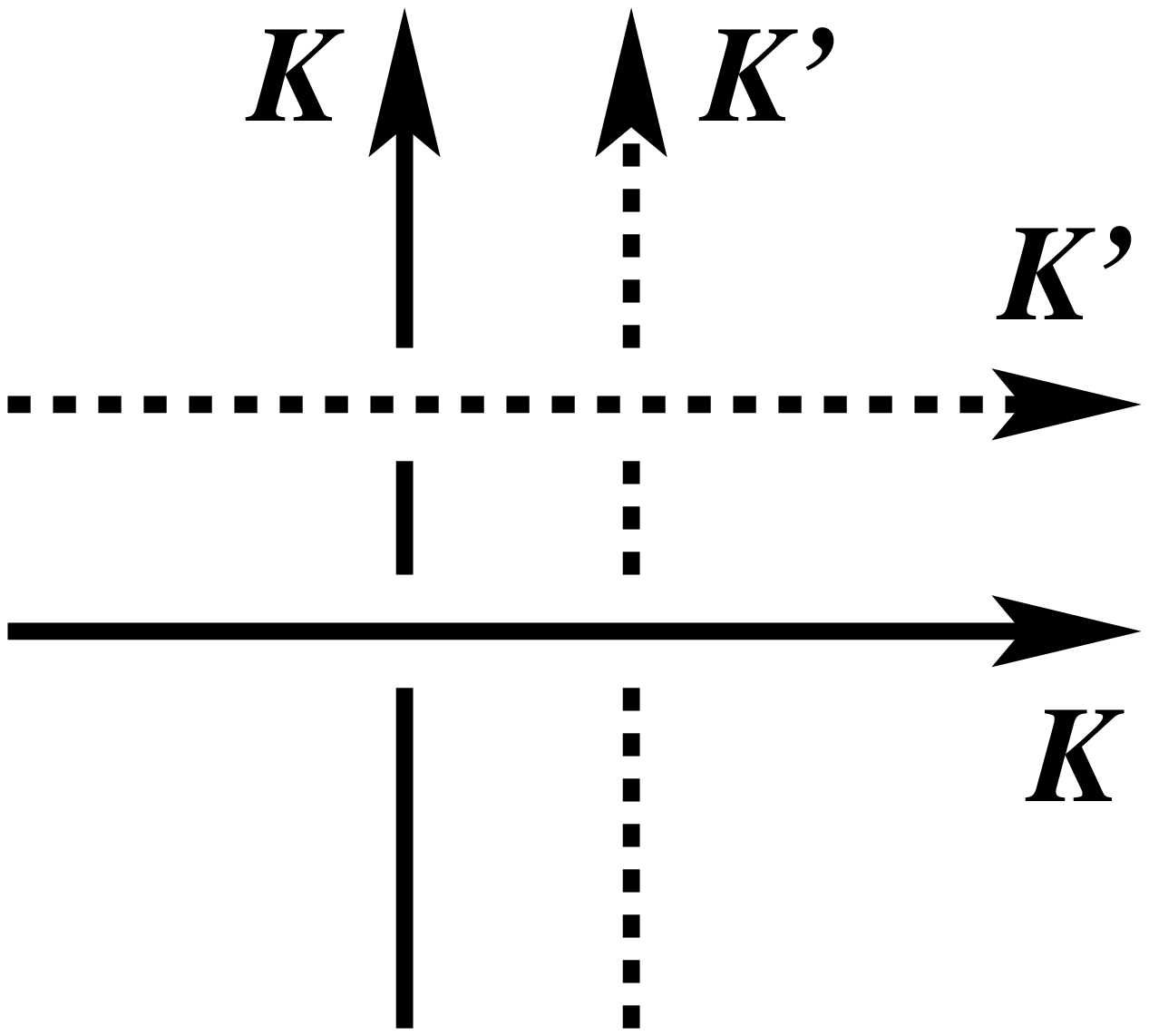}
\end{center}
The local writhe of the crossing where $K$ passes over $K'$ is the
same as the local writhe of the crossing point of the knot $K$ with
itself. Therefore, the claim follows from the combinatorial formula
for the linking number (Proposition \ref{comb_lk_num}).
\end{proof}

\section{The Conway polynomial}
\label{conway}

In what follows we shall usually consider invariants with values in
a commutative ring. Of special importance in knot theory are {\em
polynomial knot invariants} taking values in the rings of
polynomials (or Laurent polynomials\footnote{A {\em Laurent
polynomial} in $x$ is a polynomial in $x$ and $x^{-1}$.}) in one or
several variables, usually with integer coefficients.

Historically, the first polynomial invariant for knots was the {\em
Alexander polynomial}\index{Alexander polynomial} $A(K)$ introduced
in 1928 \cite{Al2}. See \cite{CrF,Lik,Rol} for a discussion of the
beautiful topological theory related to the Alexander polynomial. In
1970 J.~Conway \cite{Con} found a simple recursive construction of a
polynomial invariant $C(K)$ which differs from the Alexander
polynomial only by a change of variable, namely,
$A(K)=C(K)\mid_{t\mapsto x^{1/2}-x^{-1/2}}$. In this book, we only
use Conway's normalization. Conway's definition, given in terms of
plane diagrams, relies on
crossing point resolutions that may take a knot diagram into a link
diagram; therefore, we shall speak of links rather than knots.

\begin{definition} \index{Conway polynomial}\label{axiconw}
The Conway polynomial $\CP$ is an invariant of oriented links (and,
in particular, an invariant of oriented knots) taking values in the
ring $\Z[t]$ and defined by the two properties:

$$\begin{array}{l}
\CP\Bigl(\unkn\Bigr) = 1,\vspace{10pt}\\
\CP\Bigl(\lrints\Bigr)-\CP\Bigl(\rlints\Bigr) = t\CP\Bigl(\twoup\Bigr)\ .
\end{array}$$
\end{definition}

Here $\tallunkn$ stands for the unknot \index{Unknot} (trivial knot)
while the three pictures in the second line stand for three diagrams
that are identical everywhere except for the fragments shown. The
second relation is referred to as {\em Conway's skein
relation}.\index{Skein relation!Conway's} Skein relations are
equations on the values of some functions on knots (links, etc.)
represented by diagrams that differ from each other by local changes
near a crossing point. These relations often give a convenient way
to work with knot invariants.
\medskip

It is not quite trivial to prove the existence of an invariant
satisfying this definition, but as soon as this fact is established,
the computation of the Conway polynomial becomes fairly easy.

\begin{example}\ \\

$\displaystyle\mbox{(i)}\quad
   \CP\Bigl(\rb{-4mm}{\ig[width=22mm]{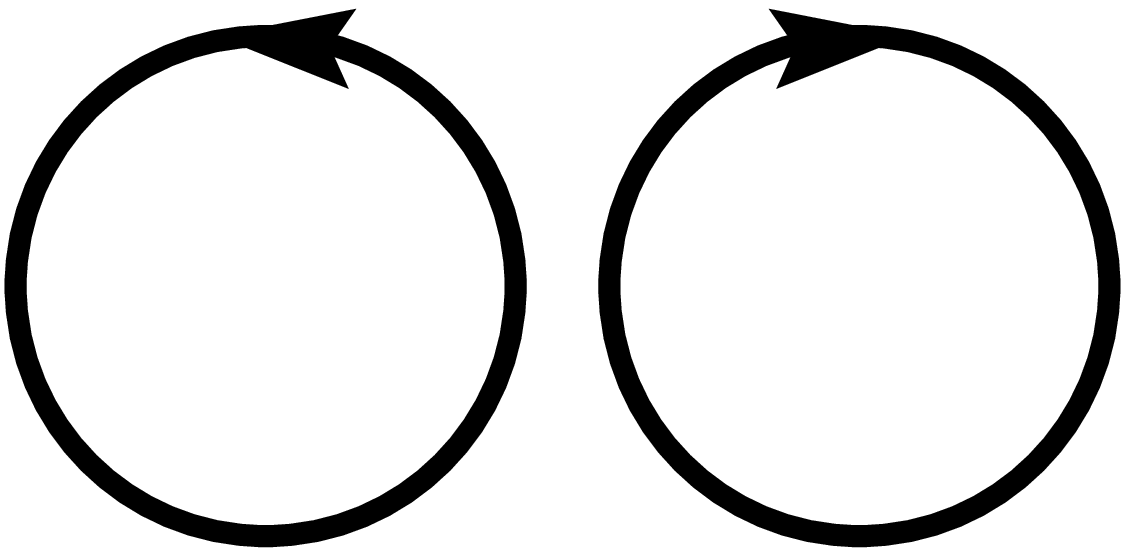}}\Bigr)
  \ =\ \frac{1}{t}
   \CP\Bigl(\rb{-4mm}{\ig[width=20mm]{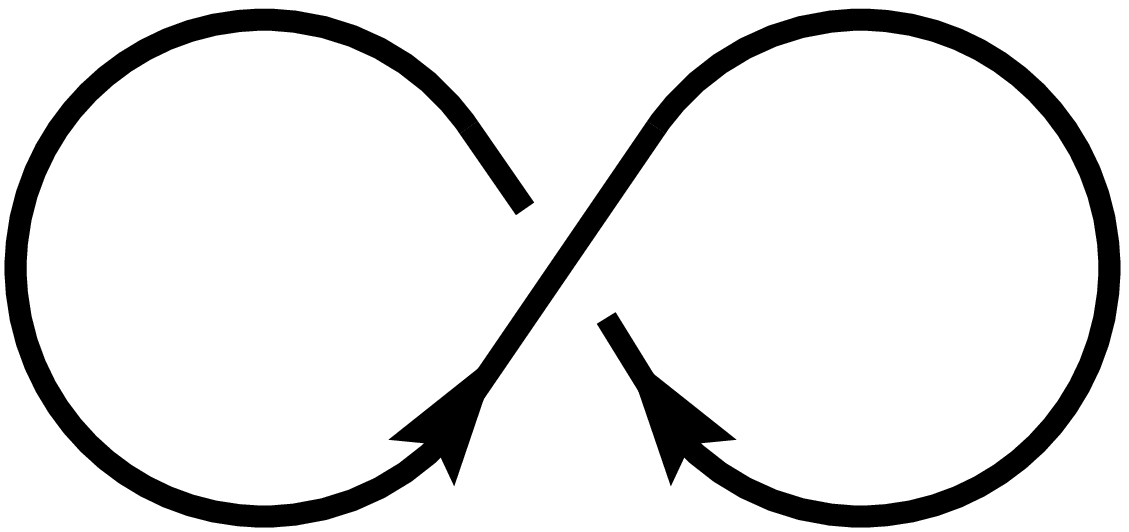}}\Bigr) -
   \frac{1}{t}
   \CP\Bigl(\rb{-4mm}{\ig[width=20mm]{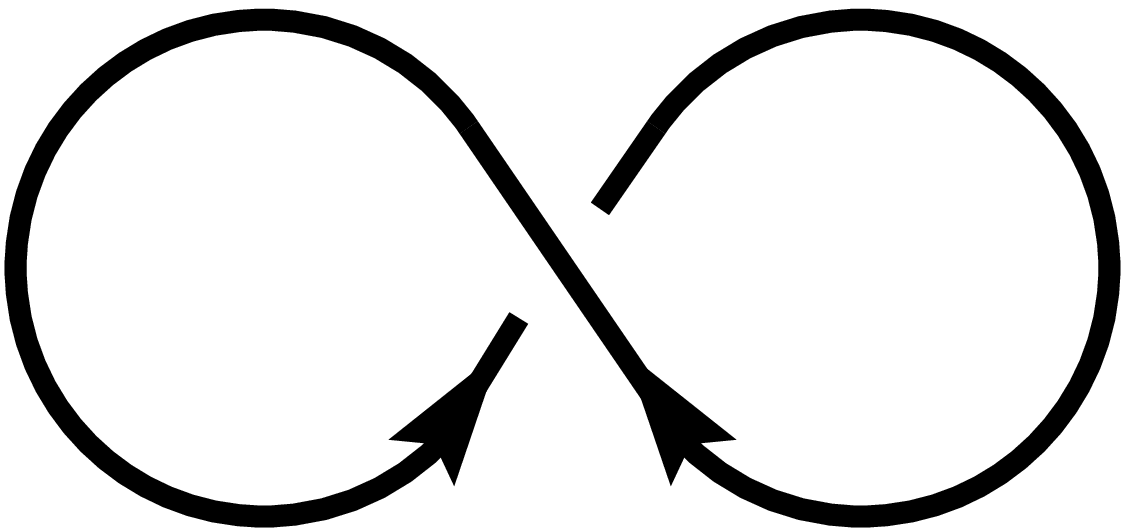}}\Bigr)
  \ =\ 0,
$\vspace{10pt}\\
because the two knots on the right are equivalent (both are trivial).\vspace{10pt}

$\begin{array}{rcl}
\mbox{(ii)}\quad
   \CP\Bigl(\rb{-4mm}{\ig[width=18mm]{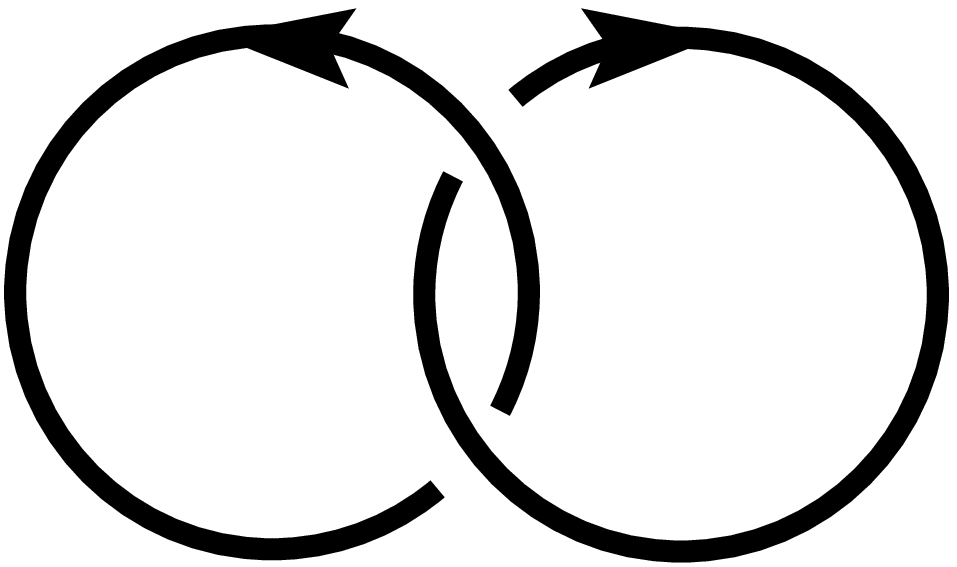}}\Bigr)
  &=&\CP\Bigl(\rb{-4mm}{\ig[width=18mm]{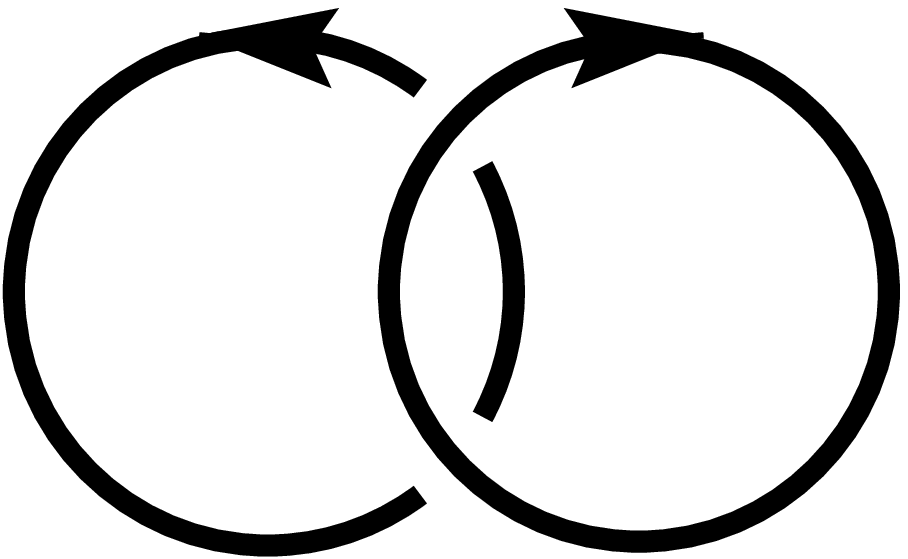}}\Bigr) -
    t \CP\Bigl(\rb{-4mm}{\ig[width=18mm]{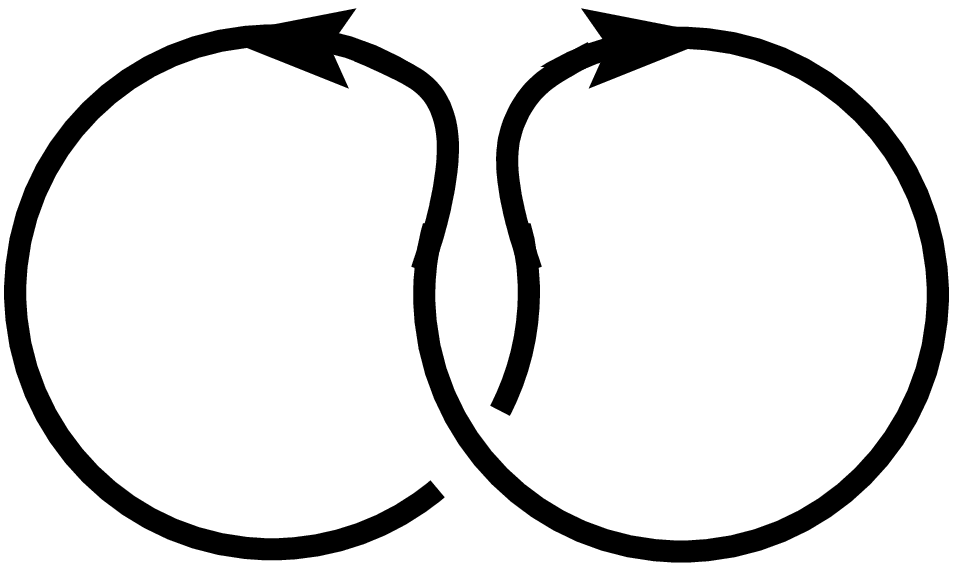}}\Bigr) \vspace{5pt}\\
  &=&\CP\Bigl(\rb{-4mm}{\ig[width=22mm]{2circ.eps}}\Bigr) -
    t \CP\Bigl(\unkn\Bigr)\ =\ -t\ .
\end{array}$\vspace{10pt}

$\begin{array}{rcl}
\mbox{(iii)}\quad
   \CP\Bigl(\rb{-4mm}{\ig[width=13mm]{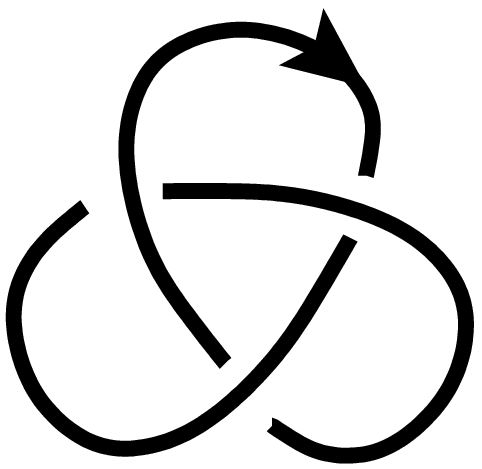}}\Bigr)
  &=&\CP\Bigl(\rb{-4mm}{\ig[width=13mm]{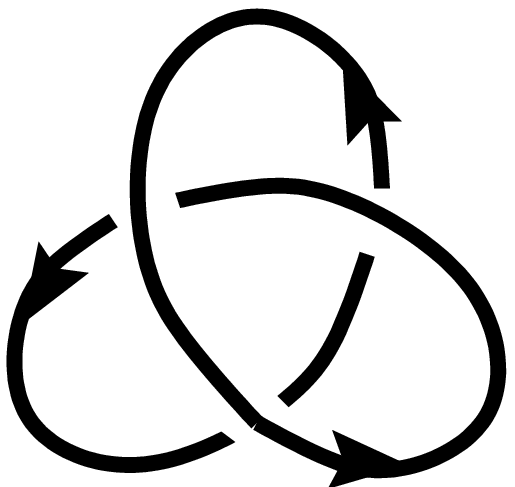}}\Bigr) -
    t \CP\Bigl(\rb{-4mm}{\ig[width=13mm]{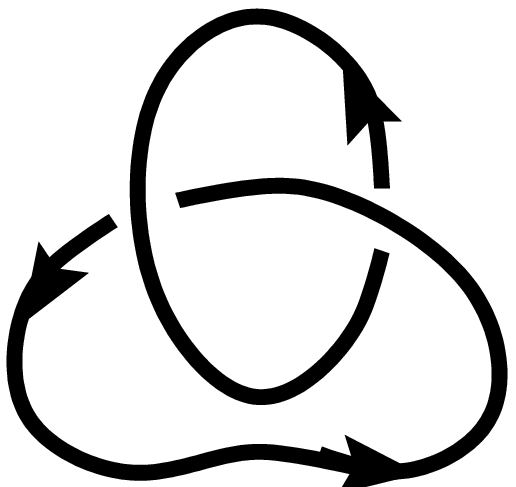}}\Bigr) \vspace{8pt}\\
  &=&\CP\Bigl(\unkn\Bigr) -
    t \CP\Bigl(\rb{-4mm}{\ig[width=18mm]{ornontrivl.eps}}\Bigr)
     \ =\ 1+t^2\ .
\end{array}$\vspace{10pt}
\end{example}

\subsection{}\label{conway_tabl}
The values of the Conway polynomial on knots with up to 8 crossings
are given in Table \ref{conw_table}. Note that the Conway polynomial
of the inverse knot $K^*$ and the mirror knot $\ol{K}$ coincides
with that of knot $K$.

\begin{table}[htb]
$$\begin{array}{c|l||c|l||c|l}
K   & \CP(K)       & K      & \CP(K)       & K      & \CP(K)      \\
\hline
3_1 & 1+ t^2          & 7_6    & 1+ t^2 -t^4     & 8_{11} & 1 -t^2 -2t^4   \\
4_1 & 1 -t^2          & 7_7    & 1 -t^2 +t^4     & 8_{12} & 1 -3t^2 +t^4   \\
5_1 & 1+ 3t^2 +t^4    & 8_1    & 1 -3t^2         & 8_{13} & 1+ t^2+ 2t^4   \\
5_2 & 1+ 2t^2         & 8_2    & 1 -3t^4 -t^6    & 8_{14} & 1 -2t^4        \\
6_1 & 1 -2t^2         & 8_3    & 1 -4t^2         & 8_{15} & 1+ 4t^2 +3t^4  \\
6_2 & 1 -t^2 -t^4     & 8_4    & 1 -3t^2 -2t^4   & 8_{16} & 1+t^2+2t^4+t^6 \\
6_3 & 1+ t^2 +t^4     & 8_5    & 1-t^2-3t^4-t^6  & 8_{17} & 1-t^2-2t^4-t^6 \\
7_1 & 1+6t^2+5t^4+t^6 & 8_6    & 1 -2t^2 -2t^4   & 8_{18} & 1+t^2-t^4-t^6  \\
7_2 & 1+ 3t^2         & 8_7    & 1+2t^2+3t^4+t^6 & 8_{19} & 1+5t^2+5t^4+t^6\\
7_3 & 1+ 5t^2+ 2t^4   & 8_8    & 1+ 2t^2+ 2t^4   & 8_{20} & 1+ 2t^2 +t^4   \\
7_4 & 1+ 4t^2         & 8_9    & 1-2t^2-3t^4-t^6 & 8_{21} & 1 -t^4         \\
7_5 & 1+ 4t^2 +2t^4   & 8_{10} & 1+3t^2+3t^4+t^6 &        &
\end{array}$$
\vspace{3mm}
\caption{Conway polynomials of knots with up to 8 crossings}\label{conw_table}
\index{Conway polynomial!table}
\index{Table of!Conway polynomials}
\end{table}

For every $n$, the coefficient $c_n$ \label{c_n-Conw-coef}
of $t^n$ in $\CP$ is a numerical invariant of the knot.

\subsection{}
The behaviour of the Conway polynomial under the change of orientation of
one component of a link does not follow any known rules. Here is an example.
$$\begin{array}{ccccccc}\label{solomon}
\rb{-10mm}{\ig[height=20mm]{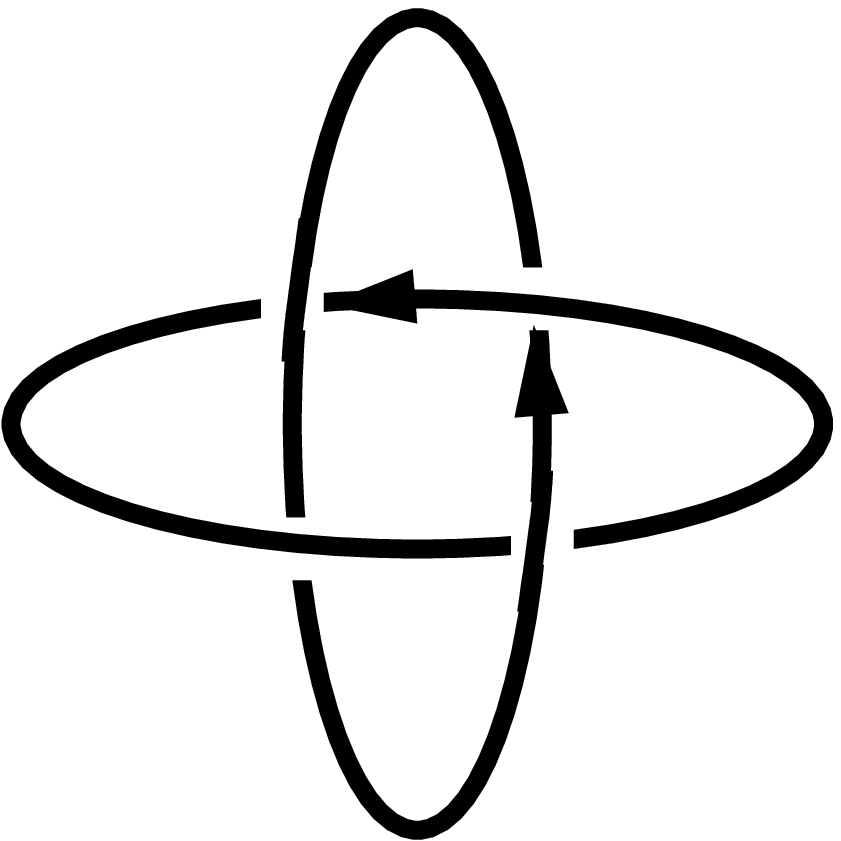}} &
\quad &
\rb{-10mm}{\ig[height=20mm]{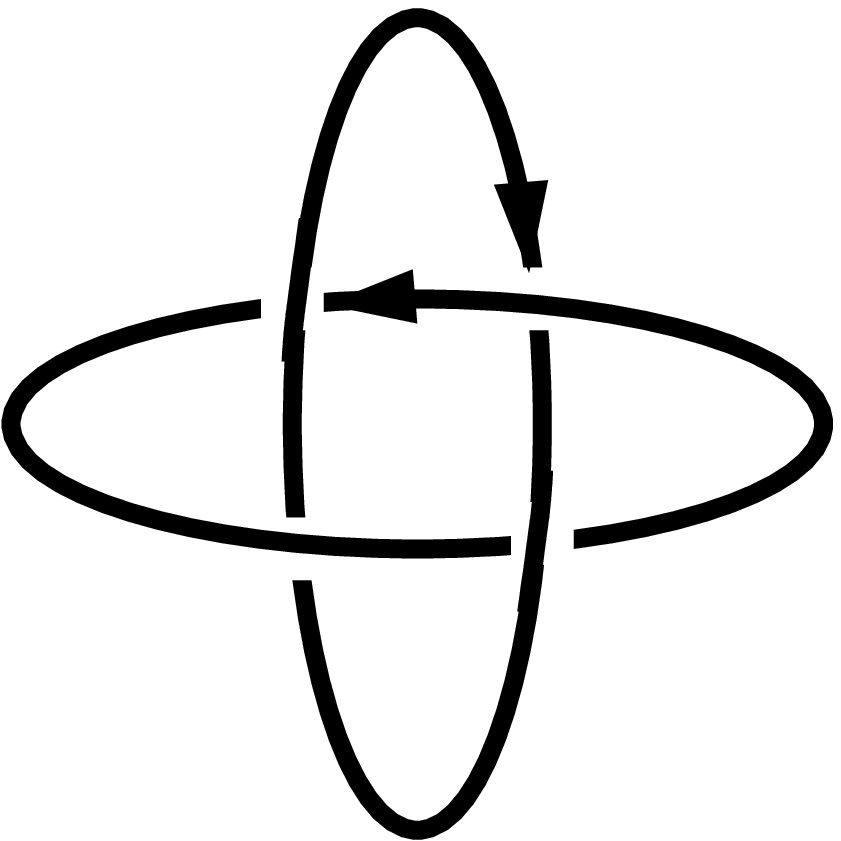}} &
\quad &
\rb{-10mm}{\ig[height=20mm]{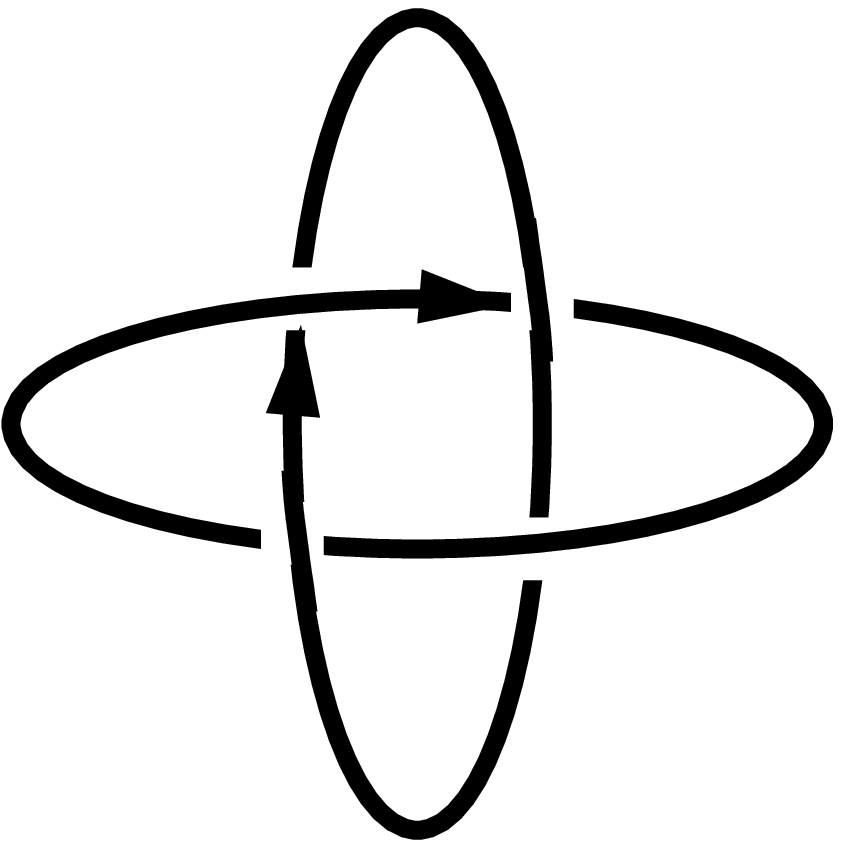}} &
\quad &
\rb{-10mm}{\ig[height=20mm]{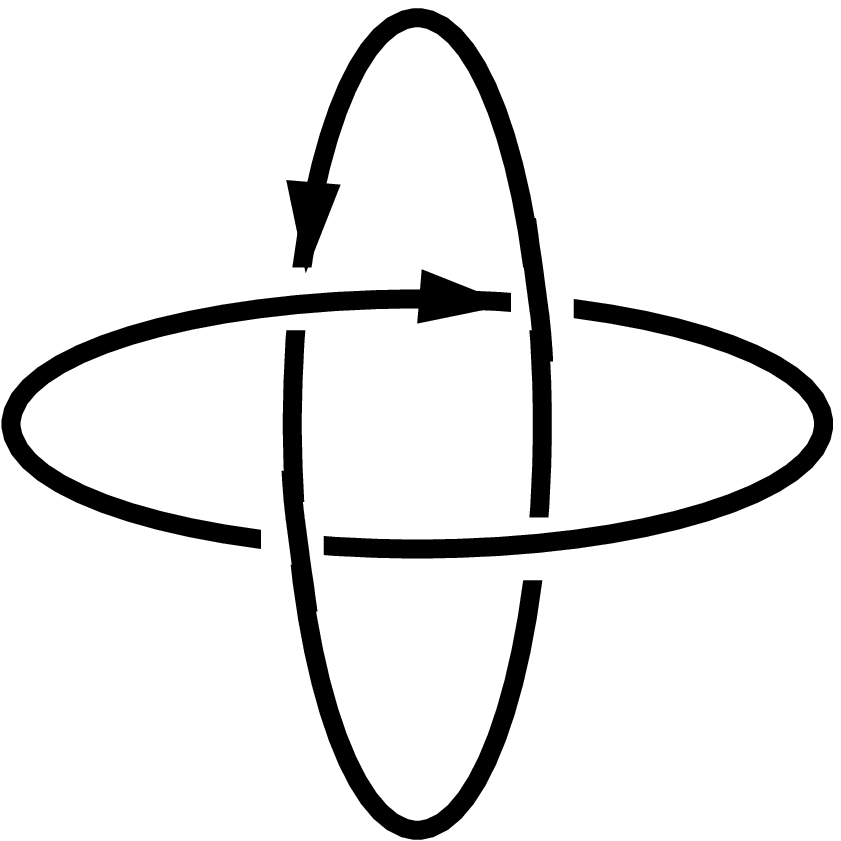}} \\
\rule{0pt}{20pt} -t^3-2t  &  \quad & 2t & \quad &  t^3+2t & \quad &  -2t
\end{array}$$

\section{The Jones polynomial}
\label{Jones}

The invention of the Jones polynomial \cite{Jo1} in 1985 produced a
genuine revolution in knot theory. The original construction of
V.~Jones was given in terms of state sums and von Neumann algebras.
It was soon noted, however, that the Jones polynomial can be defined
by skein relations, in the spirit of Conway's definition
\ref{axiconw}.

Instead of simply giving the corresponding formal equations, we
explain, following L.~Kauffman \cite{Ka6}, {\em how} this definition
could be invented. As with the Conway polynomial, the construction
given below requires that we consider invariants on the totality of
all links, not only knots, because the transformations used may turn
a knot diagram into a link diagram with several components.

Suppose that we are looking for an invariant of unoriented links,
denoted by angular brackets, that has a prescribed behaviour with
respect to the resolution of diagram crossings and the addition of a
disjoint copy of the unknot:
\newcommand{\rcross}{\rb{-3.5mm}{\ig[height=8mm]{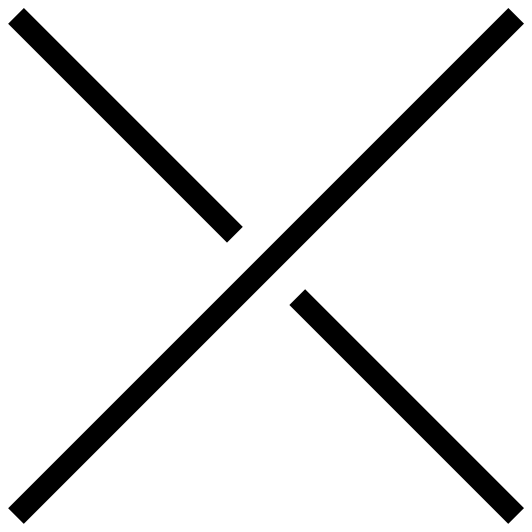}}}
\newcommand{\twovert}{\rb{-3.5mm}{\ig[height=8mm]{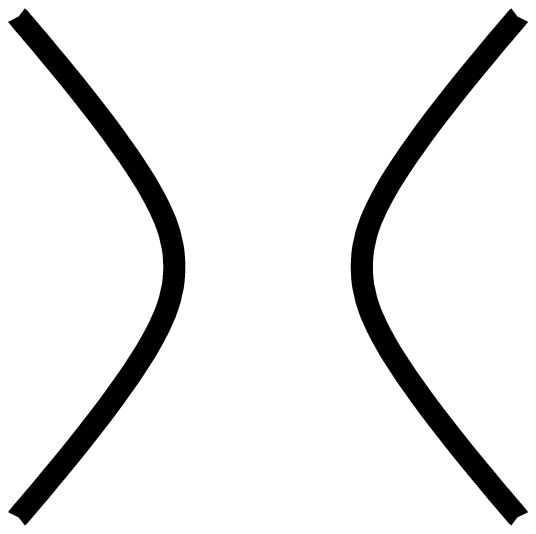}}}
\newcommand{\twohor}{\rb{-3.5mm}{\ig[height=8mm]{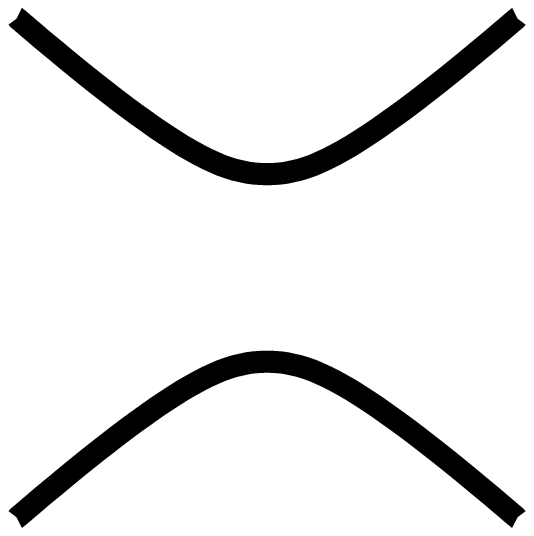}}}
\newcommand{\nunkn}{\rb{-3mm}{\ig[height=8mm]{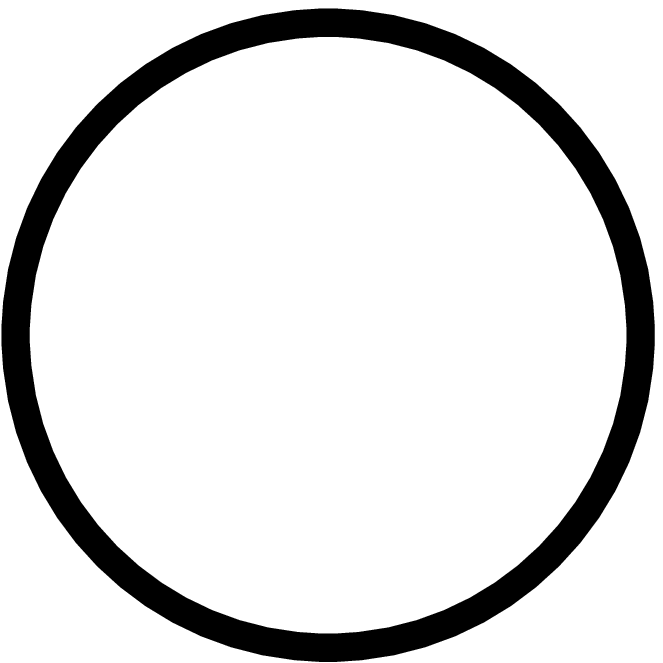}}}
\begin{gather*}
  \langle\;\rcross\;\rangle = a\,\langle\;\twovert\;\rangle
                        + b\,\langle\;\twohor\;\rangle,\\
  \langle\,L\sqcup\nunkn\,\rangle = c\,\langle\,L\,\rangle,
\end{gather*}
where $a$, $b$ and $c$ are certain fixed coefficients.

For the bracket $\langle\,,\,\rangle$ to be a link invariant, it
must be stable under the three Reidemeister moves $\Omega_1$,
$\Omega_2$, $\Omega_3$ (see Section~\ref{pl_diag}).

\begin{xca}
Show that the bracket $\langle\, , \rangle$ is $\Omega_2$-invariant
if and only if $b=a^{-1}$ and $c=-a^2-a^{-2}$. Prove that
$\Omega_2$-invariance in this case implies $\Omega_3$-invariance.
\end{xca}

\begin{xca}
  Suppose that $b=a^{-1}$ and $c=-a^2-a^{-2}$.
Check that the behaviour of the bracket with respect to the first
Reidemeister move is described by the equations
\newcommand{\hill}{\rb{-3.5mm}{\ig[height=8mm]{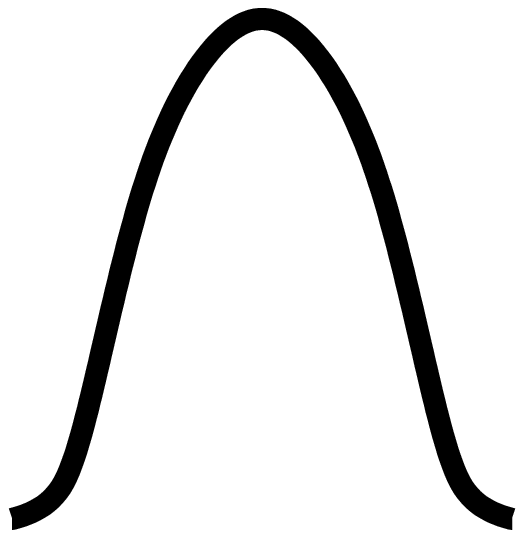}}}
\newcommand{\rhill}{\rb{-3.5mm}{\ig[height=8mm]{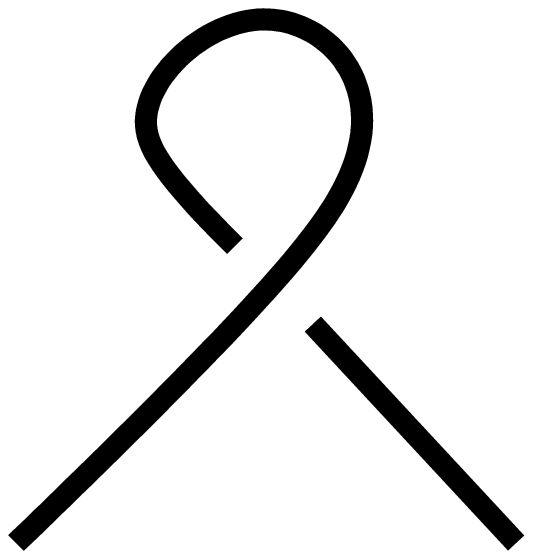}}}
\newcommand{\lhill}{\rb{-3.5mm}{\ig[height=8mm]{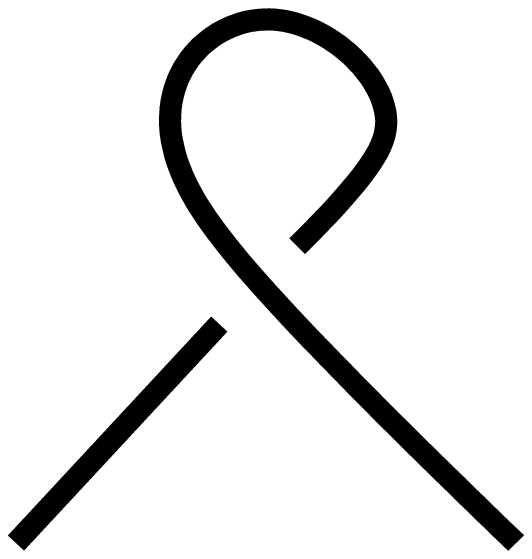}}}
\begin{gather*}
  \langle\;\rhill\;\rangle = -a^{-3}\,\langle\;\hill\;\rangle, \\
  \langle\;\lhill\;\rangle = -a^3\,\langle\;\hill\;\rangle.
\end{gather*}
\end{xca}

In the assumptions $b=a^{-1}$ and $c=-a^2-a^{-2}$, the bracket
polynomial $\langle L\rangle$ normalized by the initial condition
$$
  \langle\,\nunkn\,\rangle = 1
$$
is referred to as the {\em Kauffman
bracket\/}\label{kauf_br}\index{Kauffman bracket} of $L$. We see
that the Kauffman bracket changes only under the addition (or
deletion) of a small loop, and this change depends on the local
writhe of the corresponding crossing. It is easy, therefore, to
write a formula for a quantity that would be invariant under all
three Reidemeister moves:
$$
  J(L) = (-a)^{-3w} \langle L\rangle,
$$
where $w$ is the total writhe of the diagram (the difference between
the number of positive and negative crossings).

The invariant $J(L)$ is a Laurent polynomial called the {\em Jones
polynomial} \index{Jones polynomial}  (in {\em $a$-normalization}).
The more standard {\em $t$-normalization} is obtained by the
substitution $a=t^{-1/4}$. Note that the Jones polynomial is an
invariant of an oriented link, although in its definition we use the
Kauffman bracket which is determined by a diagram without
orientation.

\begin{xca}\label{jones_skein}
Check that the Jones polynomial is uniquely determined by the
skein relation
$$
  t^{-1}J(\lrints) - tJ(\rlints) = (t^{1/2}-t^{-1/2})J(\twoup) \eqno(1)
$$
and the initial condition
$$
  J(\unkn) = 1. \eqno(2)
$$
\end{xca}

\begin{example}
Let us compute the value of the Jones polynomial on the
left trefoil $3_1$. The calculation requires several steps, each consisting
of one application of the rule (1) and some applications of rule (2)
and/or using the results of the previous steps. We leave the details
to the reader.

(i)\quad
$J\left(\rb{-2mm}{\ig[height=6mm]{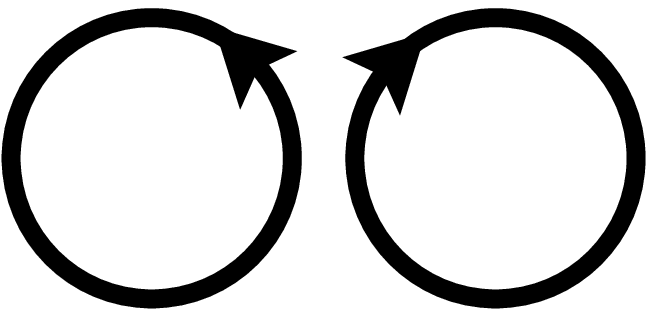}}\right)=-t^{1/2}-t^{-1/2}$.

(ii)\quad
$J\left(\rb{-2mm}{\ig[height=6mm]{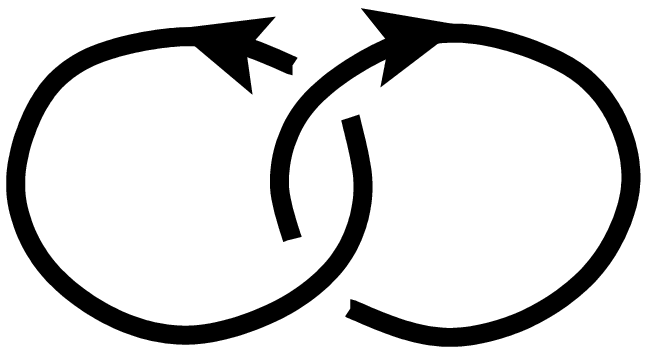}}\right)=-t^{1/2}-t^{5/2}$.

(iii)\quad
$J\left(\rb{-2mm}{\ig[height=6mm]{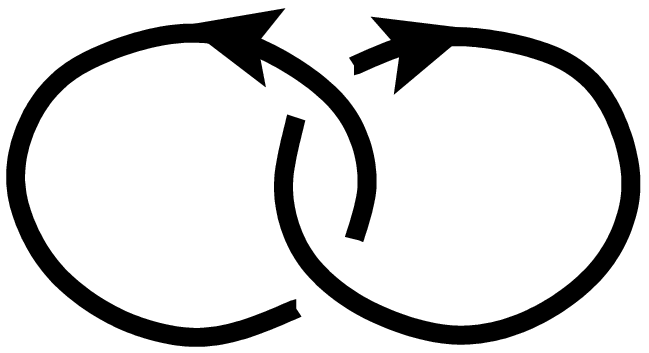}}\right)=-t^{-5/2}-t^{-1/2}$.

(iv)\quad
$J\left(\rb{-4mm}{\ig[height=10mm]{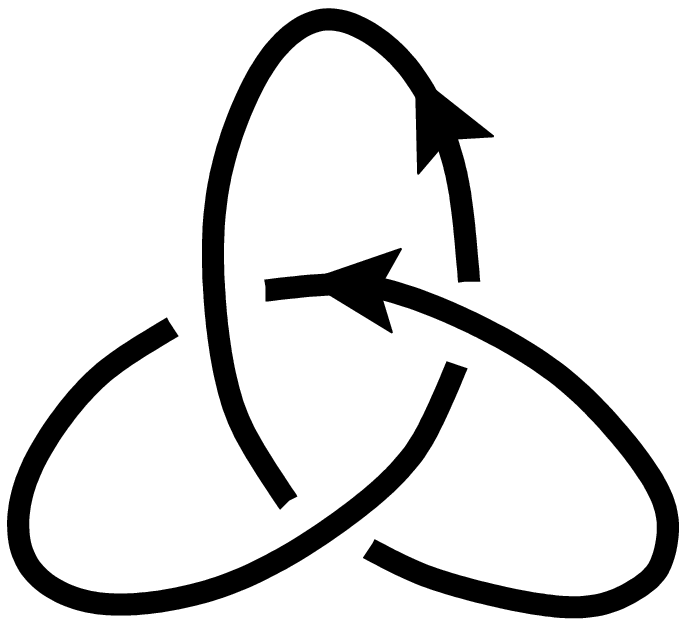}}\right)=-t^{-4}+t^{-3}+ t^{-1}$.
\vspace{10pt}
\end{example}

\begin{xca}
Repeat the previous calculation for the right trefoil and prove that
$J(\ol{3_1})=t+t^3-t^4$.
\end{xca}

We see that the Jones polynomial $J$ can tell apart two knots which
the Conway polynomial $\CP$ cannot. This does not mean, however,
that $J$ is stronger than $\CP$. There are pairs of knots, for
example, $K_1=10_{71}$, $K_2=10_{104}$ such that $J(K_1)=J(K_2)$,
but $\CP(K_1)\ne\CP(K_2)$ (see, for instance, \cite{Sto2}).

\subsection{}
The values of the Jones polynomial on standard knots with up to 8
crossings are given in Table \ref{jones_table}. The Jones polynomial
does not change when the knot is inverted (this is no longer true
for links), see Exercise~\ref{Jones_reversal}. The behaviour of the
Jones polynomial under mirror reflection is described in
Exercise~\ref{Jones_mirror}.

\begin{table}[htb]
$$\begin{array}{c|l}
3_1 & -t^{-4}+t^{-3}+t^{-1}             \\
4_1 & t^{-2}-t^{-1}+1-t+t^2 \\
5_1 & -t^{-7}+t^{-6}-t^{-5}+t^{-4}+t^{-2} \\
5_2 & -t^{-6}+t^{-5}-t^{-4}+2t^{-3}-t^{-2}+t^{-1}\\
6_1 & t^{-4}-t^{-3}+t^{-2}-2t^{-1}+2-t+t^2  \\
6_2 & t^{-5}-2t^{-4}+2t^{-3}-2t^{-2}+2t^{-1}-1+t \\
6_3 & -t^{-3} +2t^{-2} -2t^{-1} +3 -2t +2t^2 -t^3\\
7_1 & -t^{-10}+t^{-9}-t^{-8}+t^{-7}-t^{-6}+t^{-5}+t^{-3} \\
7_2 & -t^{-8}+t^{-7}-t^{-6}+2t^{-5}-2t^{-4}+2t^{-3}-t^{-2}+t^{-1} \\
7_3 & t^2 -t^3 +2t^4 -2t^5 +3t^6 -2t^7 +t^8 -t^9\\
7_4 & t -2t^2 +3t^3 -2t^4 +3t^5 -2t^6 +t^7 -t^8\\
7_5 & -t^{-9}+2t^{-8}-3t^{-7}+3t^{-6}-3t^{-5}+3t^{-4}-t^{-3}+t^{-2}  \\
7_6 & -t^{-6}+2t^{-5}-3t^{-4}+4t^{-3}-3t^{-2}+3t^{-1}-2+t \\
7_7 & -t^{-3}+3t^{-2}-3t^{-1}+4-4t+3t^2-2t^3+t^4 \\
8_1 & t^{-6}-t^{-5}+t^{-4}-2t^{-3}+2t^{-2}-2t^{-1}+2-t+t^2 \\
8_2 & t^{-8}-2t^{-7}+2t^{-6}-3t^{-5}+3t^{-4}-2t^{-3}+2t^{-2}-t^{-1}+1 \\
8_3 & t^{-4} -t^{-3} +2t^{-2} -3t^{-1} +3 -3t +2t^2 -t^3 +t^4 \\
8_4 & t^{-5} -2t^{-4} +3t^{-3} -3t^{-2} +3t^{-1} -3 +2t -t^2 +t^3 \\
8_5 & 1 -t +3t^2 -3t^3 +3t^4 -4t^5 +3t^6 -2t^7 +t^8 \\
8_6 & t^{-7}-2t^{-6}+3t^{-5}-4t^{-4}+4t^{-3}-4t^{-2}+3t^{-1}-1+t \\
8_7 & -t^{-2}+2t^{-1}-2+4t-4t^2+4t^3-3t^4+2t^5-t^6   \\
8_8 & -t^{-3}+2t^{-2}-3t^{-1}+5-4t+4t^2-3t^3+2t^4-t^5 \\
8_9 & t^{-4} -2t^{-3} +3t^{-2} -4t^{-1} +5 -4t +3t^2 -2t^3 +t^4 \\
8_{10} & -t^{-2}+2t^{-1}-3+5t-4t^2+5t^3-4t^4+2t^5-t^6 \\
8_{11} & t^{-7}-2t^{-6}+3t^{-5}-5t^{-4}+5t^{-3}-4t^{-2}+4t^{-1}-2+t \\
8_{12} & t^{-4} -2t^{-3} +4t^{-2} -5t^{-1} +5 -5t +4t^2 -2t^3 +t^4 \\
8_{13} & -t^{-3}+3t^{-2}-4t^{-1}+5-5t+5t^2-3t^3+2t^4-t^5 \\
8_{14} & t^{-7}-3t^{-6}+4t^{-5}-5t^{-4}+6t^{-3}-5t^{-2}+4t^{-1}-2+t  \\
8_{15} &t^{-10}-3t^{-9}+4t^{-8}-6t^{-7}+6t^{-6}-5t^{-5}+5t^{-4}-2t^{-3}+t^{-2} \\
8_{16} & -t^{-6}+3t^{-5} -5t^{-4} +6t^{-3} -6t^{-2} +6t^{-1} -4 +3t -t^2 \\
8_{17} & t^{-4} -3t^{-3} +5t^{-2} -6t^{-1} +7 -6t +5t^2 -3t^3 +t^4 \\
8_{18} & t^{-4} -4t^{-3} +6t^{-2} -7t^{-1} +9 -7t +6t^2 -4t^3 +t^4\\
8_{19} & t^3 +t^5 -t^8 \\
8_{20} & -t^{-5} +t^{-4} -t^{-3} +2t^{-2} -t^{-1} +2 -t \\
8_{21} & t^{-7}-2t^{-6}+2t^{-5}-3t^{-4}+3t^{-3}-2t^{-2}+2t^{-1}
\end{array}$$
\caption{Jones polynomials of knots with up to 8 crossings}\label{jones_table}
\index{Jones polynomial!table}
\index{Table of!Jones polynomials}
\end{table}

\section{Algebra of knot invariants}
\label{alg_inv}\index{Algebra!of knot invariants}

Knot invariants with values in a given commutative ring $\Ring$
form an algebra $\I$ over that ring with respect to usual
pointwise operations on functions
\begin{gather*}
   (f+g)(K)=f(K)+g(K),\\
   (fg)(K)=f(K)g(K).
\end{gather*}

Extending knot invariants by linearity to the whole algebra of knots
we see that
$$
  \I = \Hom_\Z(\Z\K,\Ring).
$$
In particular, as an $\Ring$-module (or a vector space, if $\Ring$
is a field) $\I$ is dual to the algebra $\Ring\K:=
\Z\K\otimes\Ring$, where $\Z\K$ is the algebra of knots introduced
in Section~\ref{alg_knots}. It turns out (see page~\pageref{knotcoproduct}) 
that the product on $\I$ corresponds under this duality to the {\em coproduct} on the algebra $\Ring\K$ of knots.

\section{Quantum invariants}
\label{qi}

The subject of this section is not entirely elementary. However, we
are not going to develop here a full theory of quantum groups and
corresponding invariants, confining ourselves to some basic ideas
which can be understood without going deep into complicated details.
The reader will see that it is possible to use quantum invariants
without even knowing what a quantum group is!

\subsection{}
The discovery of the Jones polynomial inspired many people to search
for other skein relations compatible with Reidemeister moves and
thus defining knot polynomials. This lead to the introduction of the
HOMFLY (\cite{HOM,PT}) and Kauffman's (\cite{Ka3,Ka4}) polynomials.
It soon became clear that all these polynomials are the first
members of a vast family of knot invariants called {\em quantum
invariants}.

The original idea of quantum invariants (in the case of 3-manifolds)
was proposed by E.~Witten in the famous paper \cite{Wit1}. Witten's
approach coming from physics was not completely justified from the
mathematical viewpoint. The first mathematically impeccable
definition of quantum invariants of links and 3-manifolds was given
by Reshetikhin and Turaev \cite{RT1,Tur2}, who used in their
construction the notion of {\em quantum groups} introduced shortly
before that by V.~Drinfeld in \cite{Dr4} (see also \cite{Dr3}) and
M.~Jimbo in \cite{Jimb}. In fact, a {\em quantum group\/} is not a
group at all. Instead, it is a family of algebras, more precisely,
of {\em Hopf algebras} (see Appendix \ref{hopf_alg}), depending on a
complex parameter $q$ and satisfying certain axioms. The {\it
quantum group} \index{Quantum group} $U_q\g$ of a semisimple Lie
algebra $\g$ is a remarkable deformation of the universal enveloping
algebra (see Appendix \ref{uea}) of $\g$ (corresponding to the value
$q=1$) in the class of Hopf algebras.

In this section, we show how the Jones polynomial $J$ can be
obtained by the techniques of quantum groups, following the approach
of Reshetikhin and Turaev. It turns out that $J$ coincides, up to
normalization, with the quantum invariant corresponding to the Lie
algebra $\g=\sL_2$ in its standard two-dimensional representation
(see Appendix \ref{stand_repr}). Later in the book, we shall
sometimes refer to the ideas illustrated in this section. For
detailed expositions of quantum groups, we refer the interested
reader to \cite{Jan,Kas,KRT}.

\subsection{}
Let $\g$ be a semisimple Lie algebra and let $V$ be its 
finite-dimensional representation. One can view $V$ as a representation of
the universal enveloping algebra $U(\g)$ (see Appendix,
page~\pageref{uea}). It is remarkable that this representation can
also be deformed with parameter $q$ to a representation of the
quantum group $U_q\g$. The vector space $V$ remains the same, but
the action now depends on $q$. For a generic value of $q$ all
irreducible representations of $U_q\g$ can be obtained in this way.
However, when $q$ is a root of unity the representation theory is
different and resembles the representation theory of $\g$ in finite
characteristic. It can be used to derive quantum invariants of
3-manifolds. For the purposes of knot theory it is enough to use
generic values of $q$, that is, those which are not roots of unity.

\subsection{} \label{qis2}
An important property of quantum
groups is that every representation gives rise to a solution $R$
of the {\em
quantum Yang--Baxter equation} \index{Yang--Baxter equation}
$$
 (R\ot\id_V)(\id_V\ot R)(R\ot\id_V) = (\id_V\ot R)(R\ot\id_V)(\id_V\ot R)
$$
where $R$ (the {\it $R$-matrix} \index{R-matrix}\label{R-matrix})
is an invertible linear operator $R:V\ot V\to V\ot V$,
and both sides of the equation are understood as linear
transformations $V\ot V\ot V\to V\ot V\ot V$.

\begin{xxca}
Given an $R$-matrix, construct a representation of the braid group $B_n$ in
the space $V^{\otimes n}$.
\end{xxca}

There is a procedure to construct an $R$-matrix associated with a
representation of a Lie algebra. We are not going to describe it in
general, confining ourselves just to one example: the Lie algebra
$\g=\sL_2$ and its standard two dimensional representation $V$
\index{Quantum invariant!$\sL_2$} (for $\sL_N$ case see exercise
(\ref{q-sl_N}) on page \pageref{q-sl_N}). In this case the
associated $R$-matrix has the form
$$R:\left\{\begin{array}{ccl}
e_1\ot e_1& \mapsto&\ \ q^{1/4} e_1\ot e_1\\
e_1\ot e_2& \mapsto& q^{-1/4} e_2\ot e_1\\
e_2\ot e_1& \mapsto& q^{-1/4} e_1\ot e_2 + (q^{1/4}-q^{-3/4}) e_2\ot e_1\\
e_2\ot e_2& \mapsto&\ \ q^{1/4} e_2\ot e_2
\end{array}\right.\label{Rmatix}$$
for an appropriate basis $\{e_1,e_2\}$ of the space $V$.
The inverse of $R$ (we shall need it later) is given by the formulae
$$R^{-1}:\left\{\begin{array}{ccl}
e_1\ot e_1 &\mapsto& q^{-1/4} e_1\ot e_1\\
e_1\ot e_2 &\mapsto&\ \  q^{1/4} e_2\ot e_1 + (-q^{3/4}+q^{-1/4}) e_1\ot e_2\\
e_2\ot e_1 &\mapsto&\ \  q^{1/4} e_1\ot e_2\\
e_2\ot e_2 &\mapsto& q^{-1/4} e_2\ot e_2
\end{array}\right.$$

\subsection{Exercise.}
Check that this operator $R$ satisfies the quantum Yang-Baxter equation.

\subsection{} \label{qis4}

The general procedure of constructing quantum invariants is
organized as follows (see details in \cite{Oht1}). Consider a knot
diagram in the plane and take a generic horizontal line. To each
intersection point of the line with the diagram assign either the
representation space $V$ or its dual $V^*$ depending on whether the
orientation of the knot at this intersection is directed upwards or
downwards. Then take the tensor product of all such spaces over the
whole horizontal line. If the knot diagram does not intersect the
line, then the corresponding vector space is the ground field $\C$.

A portion of a knot diagram between two such horizontal lines
represents a tangle $T$ (see the general definition in
Section~\ref{tangles}). We assume that this tangle is framed by the
blackboard framing. With $T$ we associate a linear transformation
$\t^{fr}(T)$ from the vector space corresponding to the bottom of
$T$ to the vector space corresponding to the top of $T$. The
following three properties hold for the linear transformation
$\t^{fr}(T)$:
\begin{center}\begin{tabular}{cl}
$\bullet$ & $\t^{fr}(T)$ is an invariant of the isotopy class of the framed
tangle $T$;\\
$\bullet$ & $\t^{fr}(T_1\cdot T_2) = \t^{fr}(T_1)\circ\t^{fr}(T_2)$;\\
$\bullet$ & $\t^{fr}(T_1\ot T_2) = \t^{fr}(T_1)\ot\t^{fr}(T_2)$.
\end{tabular}\end{center}
$$
\rb{-50pt}{\ig[width=120pt]{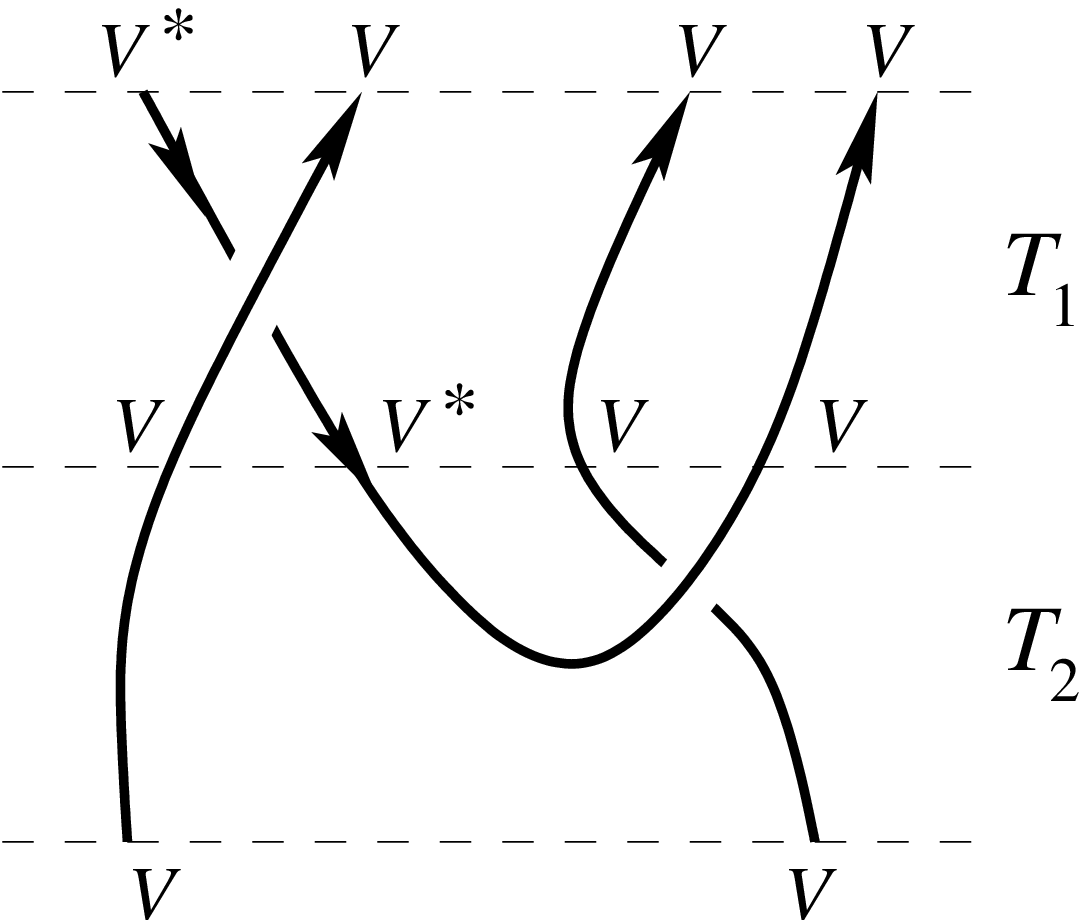}}
\quad\ \ig[width=25pt]{totor.eps}\ \quad
\rb{40pt}{\xymatrix{
    V^*\ot V\ot V\ot V \\
    V\ot V^*\ot V\ot V \ar[u]^{\t^{fr}(T_1)} \\
    V\ot V \ar[u]^{\t^{fr}(T_2)}
           \ar@/_60pt/[uu]_{\t^{fr}(T_1\cdot T_2)} }}
$$
Now we can define a knot invariant $\t^{fr}(K)$ regarding the knot
$K$ as a tangle between the two lines below and above $K$. In this
case $\t^{fr}(K)$ would be a linear transformation from $\C$ to
$\C$, that is, multiplication by a number. Since our linear
transformations depend on the parameter $q$, this number is actually
a function of $q$. \index{Quantum invariant!framed}

\subsection{} \label{qis5}
Because of the multiplicativity property $\t^{fr}(T_1\cdot
T_2)=\t^{fr}(T_1)\circ\t^{fr}(T_2)$ it is enough to define
$\t^{fr}(T)$ only for elementary tangles $T$ such as a crossing, a
minimum or a maximum point. This is precisely where quantum groups
come in. Given a quantum group $U_q\g$ and its finite-dimensional
representation $V$, one can associate certain linear transformations
with elementary tangles in a way consistent with the Turaev oriented
moves from page \pageref{Tur-mov-or}. The $R$-matrix appears here as
the linear transformation corresponding to a positive crossing,
while $R^{-1}$ corresponds to a negative crossing. Of course, for a
trivial tangle consisting of a single string connecting the top and
bottom, the corresponding linear operator should be the identity
transformation. So we have the following correspondence valid for
all quantum groups:
$$\begin{array}{c|@{\qquad}c@{\qquad}} \hline
\rb{-30pt}[36pt][35pt]{\ig[width=65pt]{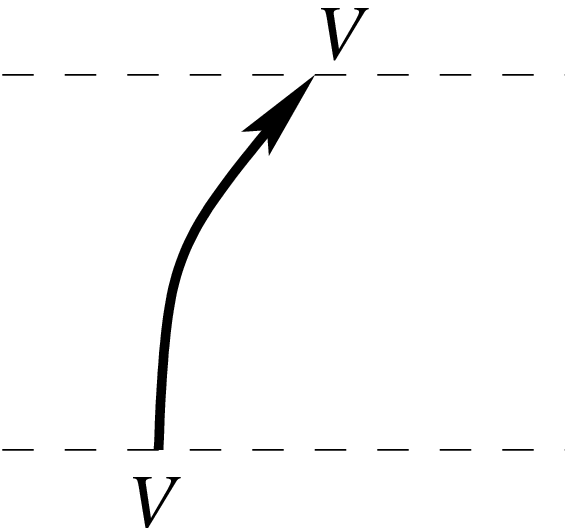}} & \begin{CD}
                              V \\
                              @AA id_V A\\
                              V
                           \end{CD} \\ \hline
\end{array}\qquad\qquad
\begin{array}{c|@{\qquad}c@{\qquad}}  \hline
\rb{-30pt}[36pt][35pt]{\ig[width=65pt]{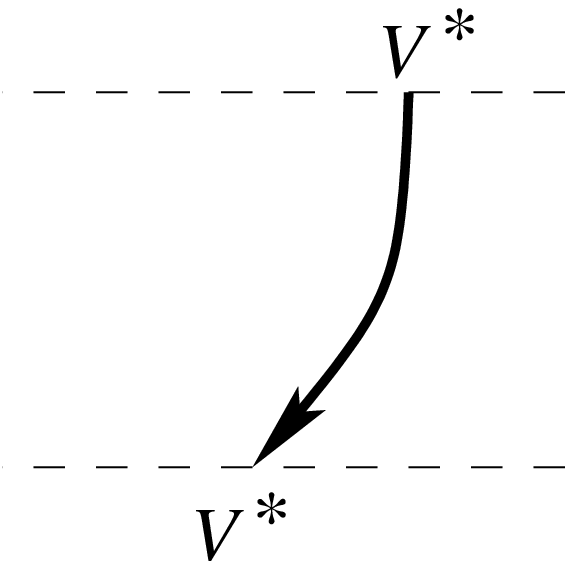}} & \begin{CD}
                              V^* \\
                              @AA id_{V^*} A\\
                              V^*
                           \end{CD} \\  \hline
\end{array}
$$
$$\begin{array}{c|@{\qquad}c@{\qquad}} \hline
\rb{-30pt}[36pt][35pt]{\ig[width=65pt]{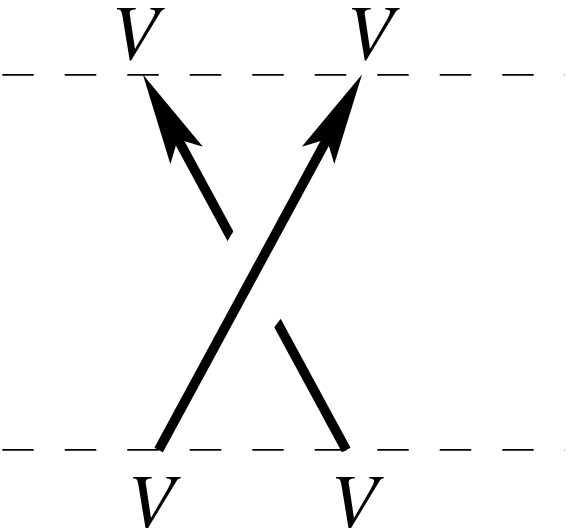}} & \begin{CD}
                              V\ot V \\
                              @AA R A\\
                              V\ot V
                           \end{CD}  \\  \hline
\end{array}\qquad\qquad
\begin{array}{c|@{\qquad}c@{\qquad}} \hline
\rb{-30pt}[36pt][35pt]{\ig[width=65pt]{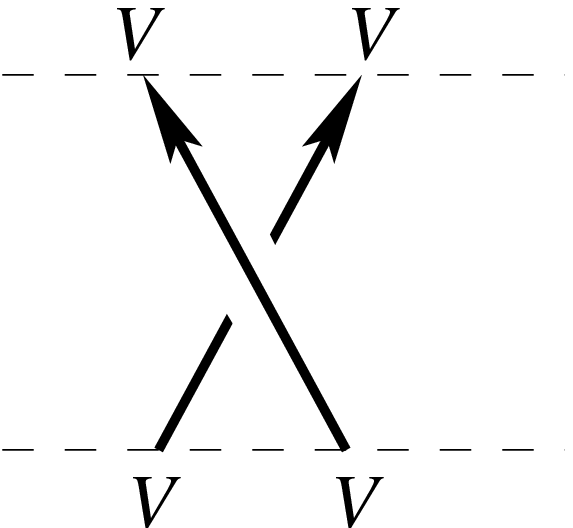}} & \begin{CD}
                              V\ot V \\
                              @AA R^{-1} A\\
                              V\ot V
                           \end{CD}  \\  \hline
\end{array}
$$

Using this we can easily check that the invariance of a quantum
invariant under the third Reidemeister move is nothing else but the
quantum Yang--Baxter equation:
$$\begin{array}{rcl}
\begin{CD}
   V\ot V\ot V \\
   @A R\ot id_V AA\\
   V\ot V\ot V \\
   @A id_V\ot R AA\\
   V\ot V\ot V \\
   @A R\ot id_V AA\\
   V\ot V\ot V
\end{CD}\quad
\rb{-70pt}{\ig[width=85pt]{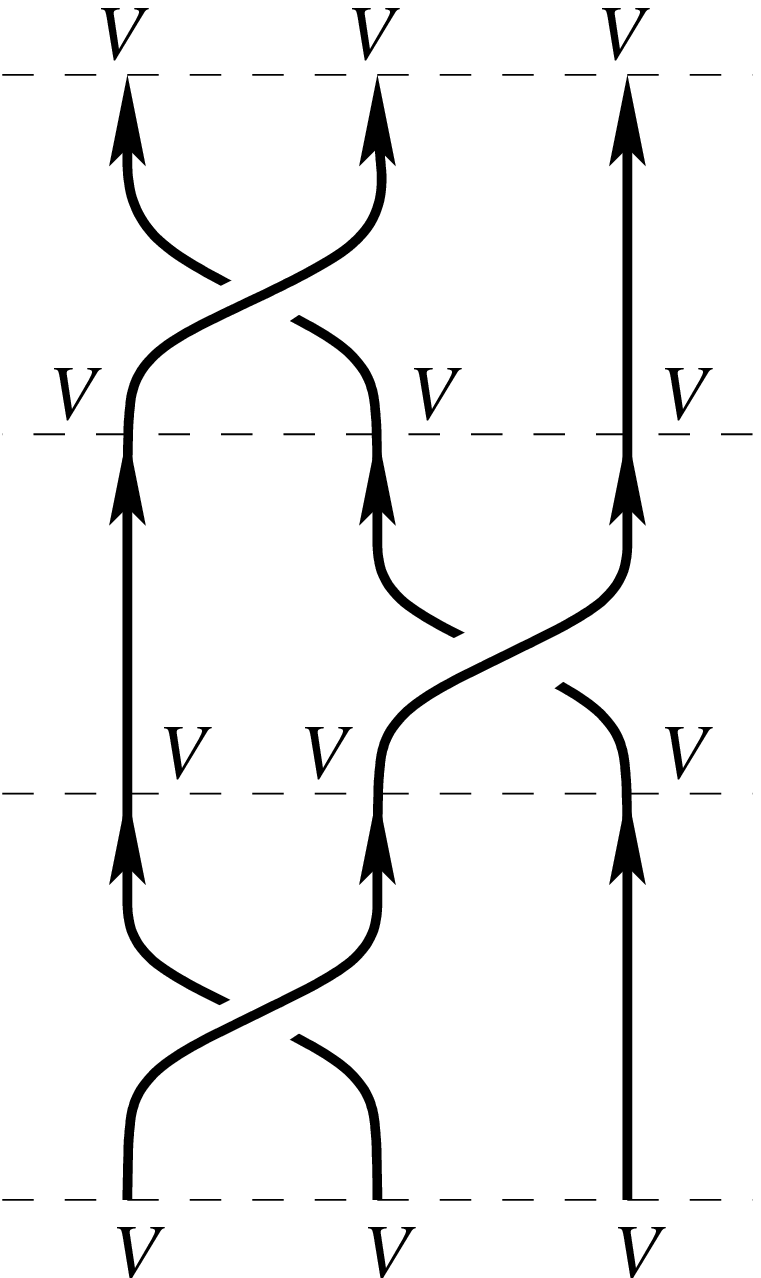}}
&=&
\rb{-70pt}{\ig[width=85pt]{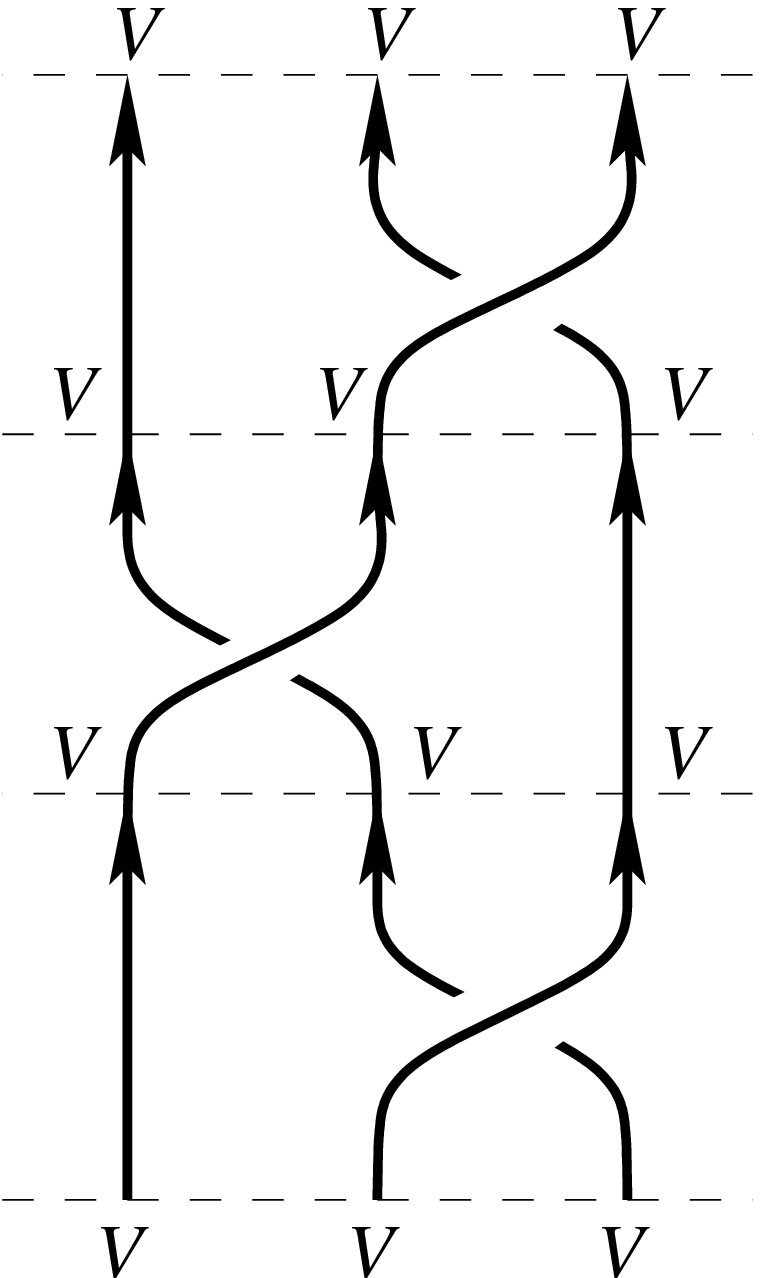}}
\quad
\begin{CD}
   V\ot V\ot V \\
   @AA id_V\ot R A\\
   V\ot V\ot V \\
   @AA R\ot id_V A\\
   V\ot V\ot V \\
   @AA id_V\ot R A\\
   V\ot V\ot V
\end{CD}\vspace{15pt}\\
(R\ot id_V)\!\circ\!(id_V\ot R)\!\circ\!(R\ot id_V)&=&
(id_V\ot R)\!\circ\!(R\ot id_V)\!\circ\!(id_V\ot R)
\end{array}
$$

Similarly, the fact that we assigned mutually inverse operators ($R$
and $R^{-1}$) to positive and negative crossings implies the
invariance under the second Reidemeister move. (The first
Reidemeister move is treated in Exercise \ref{exR1} below.)

To complete the construction of our quantum invariant we should
assign appropriate operators to the minimum and maximum points.
These depend on all the data involved: the quantum group, the
representation and the $R$-matrix. For the quantum group $U_q\sL_2$,
its standard two dimensional representation $V$ and the $R$-matrix
chosen in \ref{qis2} these operators are:
$$\label{min-max-oper}
\begin{array}{c|c|c} \hline
\rb{-22pt}{\makebox(15,10){$\minr=$}}\rb{-35pt}{\ig[width=65pt]{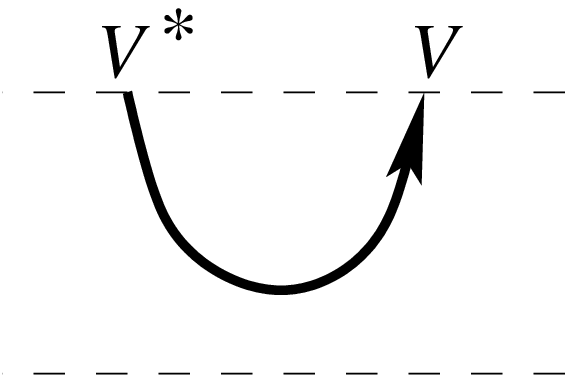}} &
{\xymatrix{V^*\ot V\ \\
          \C\ar@{->}[u]}}&
{\xymatrix{q^{-1/2}e^1\ot e_1 + q^{1/2}e^2\ot e_2\ \ \\
              1\ar@{|->}[u]}}\\ \hline
\end{array}
$$
$$
\begin{array}{c|c|c} \hline
\rb{-22pt}{\makebox(15,10){$\minl=$}}\rb{-35pt}{\ig[width=65pt]{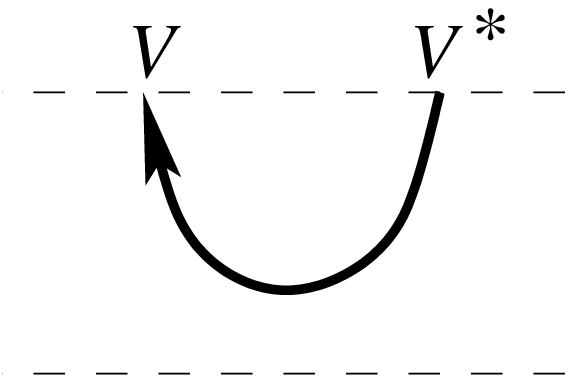}} &
{\xymatrix{\ V\ot V^*\\
          \C\ar@{->}[u]}}&
{\xymatrix{e_1\ot e^1 + e_2\ot e^2 \\
              1\ar@{|->}[u]}}\\ \hline
\end{array}
$$
$$
\begin{array}{c|c|c} \hline
\rb{-22pt}{\makebox(15,10){$\maxr=$}}\rb{-40pt}{\ig[width=65pt]{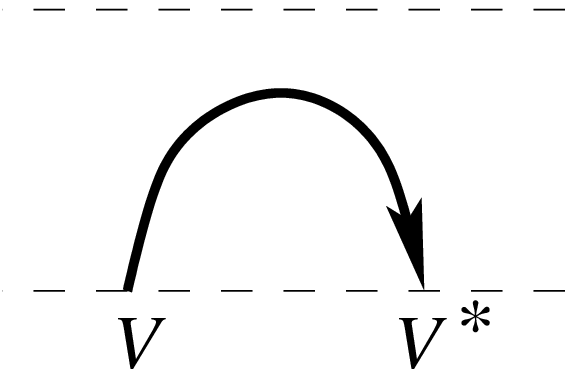}} &
{\xymatrix{\C \\
          \ V\ot V^*\ar@{->}[u]}}&
{\xymatrix{\ \ q^{1/2}\\ e_1\ot e^1\ar@{|->}[u]}}\ \
\rb{-2.5pt}{$\xymatrix{   0    \\ e_1\ot e^2\ar@{|->}[u]}$}\ \
\rb{-2.5pt}{$\xymatrix{   0    \\ e_2\ot e^1\ar@{|->}[u]}$}\ \
{\xymatrix{\quad q^{-1/2}\\ e_2\ot e^2\ar@{|->}[u]}}\\ \hline
\end{array}
$$
$$
\begin{array}{c|c|c} \hline
\rb{-22pt}{\makebox(15,10){$\maxl=$}}\rb{-40pt}{\ig[width=65pt]{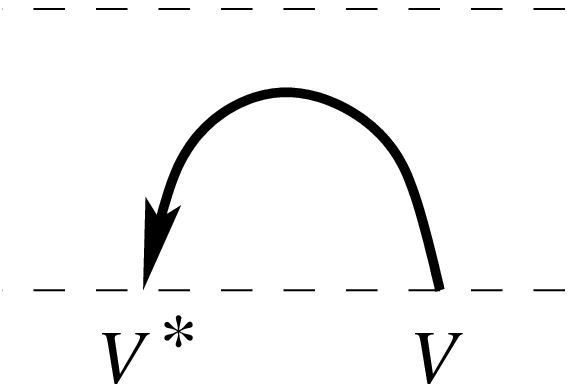}} &
{\xymatrix{\C \\
          V^*\ot V\ \ar@{->}[u]}}&
{\xymatrix{1\\ e^1\ot e_1\ar@{|->}[u]}}\ \
{\xymatrix{0\\ e^1\ot e_2\ar@{|->}[u]}}\ \
{\xymatrix{0\\ e^2\ot e_1\ar@{|->}[u]}}\ \
{\xymatrix{1\\ e^2\ot e_2\ar@{|->}[u]}}\\ \hline
\end{array}
$$
where $\{e^1,e^2\}$ is the basis of $V^*$ dual to the basis
$\{e_1,e_2\}$ of the space $V$.

We leave to the reader the exercise to check that these operators
are consistent with the oriented Turaev moves from page
\pageref{Tur-mov-or}. See  Exercise~\ref{q-sl_N} for their
generalization to $\sL_N$.

\subsection{Example.} \label{qis6}
Let us compute the $\sL_2$-quantum invariant of the unknot.
Represent the unknot as a product of two tangles and compute the
composition of the corresponding transformations
$$
\rb{-45pt}{\ig[width=80pt]{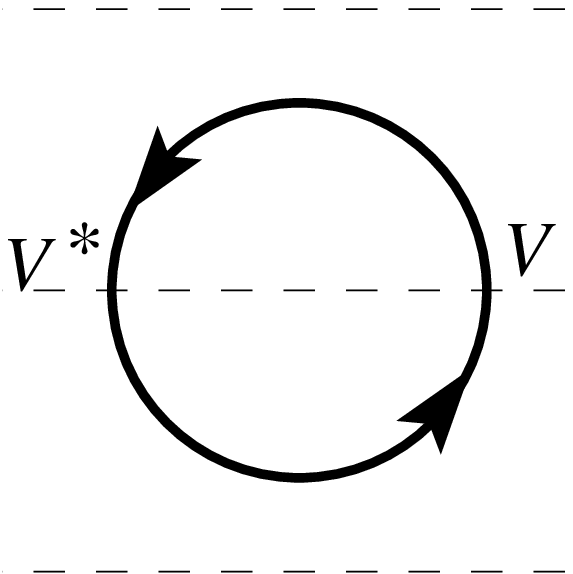}}\qquad\qquad
\rb{35pt}{\xymatrix{\C \\ V^*\ot V\ \ar@{->}[u] \\ \C\ar@{->}[u]}}\qquad
\begin{array}{c}{\underbrace{\begin{array}{ccc}
    {\underbrace{q^{-1/2}}}&+&{\underbrace{q^{1/2}}} \vspace{-6pt}\\
    {\xymatrix{\\ \ar@{|->}[u]}}&&{\xymatrix{\\ \ar@{|->}[u]}} \vspace{-4pt}\\
    {\overbrace{q^{-1/2}e^1\ot e_1}}&+&{\overbrace{q^{1/2}e^2\ot e_2}}
                       \end{array}}} \vspace{-6pt}\\
{\xymatrix{ \\ 1\ar@{|->}[u]}}
\end{array}
$$
So $\t^{fr}(\mbox{unknot})= q^{1/2}+q^{-1/2}$. Therefore, in order
to normalize our invariant so that its value on the unknot is equal
to 1, we must divide $\t^{fr}(\cdot)$ by $q^{1/2}+q^{-1/2}$. We
denote the normalized invariant by
$\wt{\t}^{fr}(\cdot)=\frac{{\t^{fr}}(\cdot)}{q^{1/2}+q^{-1/2}}$.

\subsection{Example.} \label{qis7}
Let us compute the quantum invariant for the left trefoil.
Represent the diagram of the trefoil as follows.
$$\rb{-140pt}{\ig[width=170pt]{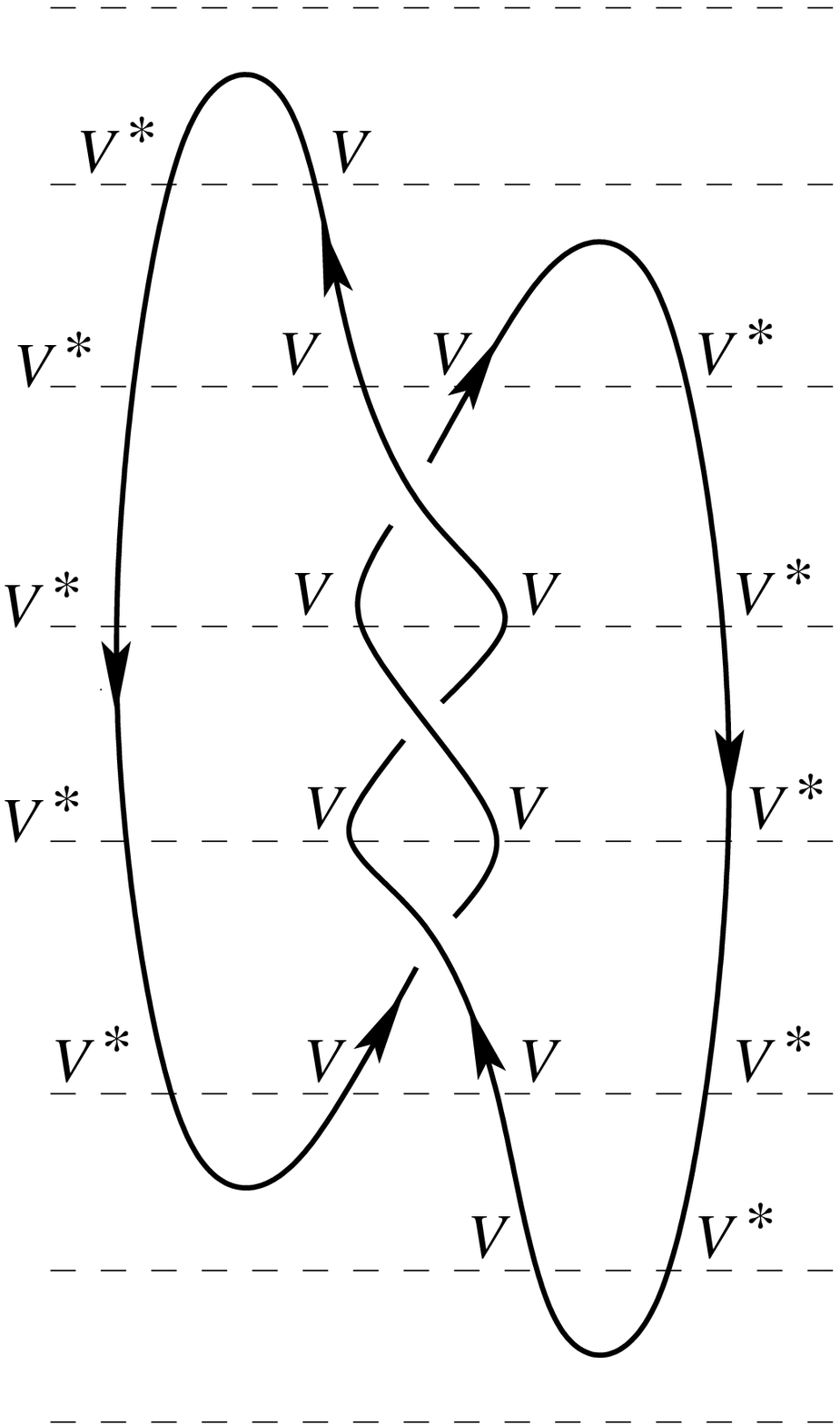}}\qquad\qquad
\begin{CD}
   \C\\@AAA\\ V^*\ot V\\@AAA\\ V^*\ot V\ot V\ot V^*\\
       @AA \id_{V^*}\ot R^{-1}\ot \id_{V^*} A \\ V^*\ot V\ot V\ot V^*\\
       @AA \id_{V^*}\ot R^{-1}\ot \id_{V^*} A \\ V^*\ot V\ot V\ot V^*\\
       @AA \id_{V^*}\ot R^{-1}\ot \id_{V^*} A \\ V^*\ot V\ot V\ot V^*\\
   @AAA\\ V\ot V^*\\ @AAA\\ \C
\end{CD}$$
Two maps at the bottom send $1\in \C$ into the tensor
$$\begin{array}{rcl}
1&\mapsto& q^{-1/2}e^1\ot e_1\ot e_1\ot e^1\quad +\quad
           q^{-1/2}e^1\ot e_1\ot e_2\ot e^2 \vspace{5pt}\\
&&\hspace{2cm}+\quad q^{1/2}e^2\ot e_2\ot e_1\ot e^1\quad
+\quad q^{1/2}e^2\ot e_2\ot e_2\ot e^2\ .
\end{array}$$
Then applying $R^{-3}$ to two tensor factors in the middle we get
$$\begin{array}{l}
 q^{-1/2}e^1\ot \Bigl(q^{-3/4} e_1\ot e_1\Bigr)\ot e^1 \\
+q^{-1/2}e^1\ot \Bigl(
   \bigl(-q^{9/4}+q^{5/4}-q^{1/4}+q^{-3/4}\bigr)e_1\ot e_2
  \\ \hspace{87pt}
  + \bigl(-q^{7/4}-q^{3/4}-q^{-1/4}\bigr)e_2\ot e_1\Bigr)\ot e^2 \\
+q^{1/2}e^2\ot \Bigl(
   \bigl(q^{7/4}-q^{3/4}+q^{-1/4}\bigr)e_1\ot e_2 +
   \bigl(-q^{5/4}+q^{1/4}\bigr)e_2\ot e_1\Bigr)\ot e^1 \\
+q^{1/2}e^2\ot \Bigl(q^{-3/4} e_2\ot e_2\Bigr)\ot e^2\ .
\end{array}$$
Finally, the two maps at the top contract the whole tensor into a number
$$\begin{array}{ccl}
\t^{fr}(3_1)&=&
 q^{-1/2}q^{-3/4}q^{1/2}
+q^{-1/2}\bigl(-q^{9/4}+q^{5/4}-q^{1/4}+q^{-3/4}\bigr)q^{-1/2}\vspace{8pt} \\
&&+q^{1/2}\bigl(-q^{5/4}+q^{1/4}\bigr)q^{1/2}
  +q^{1/2}q^{-3/4}q^{-1/2} \vspace{10pt}\\
&=&2q^{-3/4} -q^{5/4}+q^{1/4}-q^{-3/4}+q^{-7/4}-q^{9/4}+q^{5/4}
         \vspace{10pt}\\
&=&q^{-7/4}+q^{-3/4}+q^{1/4}-q^{9/4}
\end{array}$$
Dividing by the normalizing factor $q^{1/2}+q^{-1/2}$ we get
$$\frac{\t^{fr}(3_1)}{q^{1/2}+q^{-1/2}}
  \quad=\quad q^{-5/4} + q^{3/4} - q^{7/4}\,.$$

The invariant $\t^{fr}(K)$ remains unchanged under the second and
third Reidemeister moves. However, it varies under the first
Reidemeister move and thus depends on the framing. One can {\em
deframe}\index{Deframing!of framed knot invariants} it, that is,
manufacture an invariant of unframed knots out of it, according to
the formula
$$\t(K) = q^{-\frac{c\cdot w(K)}{2}}\t^{fr}(K)\,,
\index{Quantum invariant!unframed}$$ where $w(K)$ is the writhe of
the knot diagram and $c$ is the quadratic Casimir number (see
Appendix \ref{casimir_num}) defined by the Lie algebra $\g$ and its
representation. For $\sL_2$ and the standard 2-dimensional
representation $c=3/2$. The writhe of the left trefoil in our
example equals $-3$. Hence for the unframed normalized quantum
invariant we have
$$\wt{\t}(3_1)=\frac{\t(3_1)}{q^{1/2}+q^{-1/2}}
  \quad=\quad q^{9/4}\Bigl(q^{-5/4} + q^{3/4} - q^{7/4}\Bigr)
  \quad =\quad q + q^3 - q^4\,.$$
The substitution $q=t^{-1}$ gives the Jones polynomial \quad
$t^{-1} + t^{-3} - t^{-4}$.

\section{Two-variable link polynomials}
\subsection{HOMFLY polynomial}\label{HOMFLY-def}

The {\em HOMFLY polynomial}\index{HOMFLY polynomial} $P(L)$ is an
unframed link invariant. It is defined as the Laurent polynomial in
two variables $a$ and $z$ with integer coefficients satisfying the
following skein relation and the initial condition:
$$aP(\lrints)\ -\  a^{-1}P(\rlints)\ = \  zP(\twoup)\ ;\qquad
P(\unkn)\quad = \quad 1\,.$$ The existence of such an invariant is a
difficult theorem. It was established simultaneously and
independently by five groups of authors \cite{HOM,PT} (see also
\cite{Lik}). The HOMFLY polynomial is equivalent to the collection
of quantum invariants associated with the Lie algebra $\sL_N$ and
its standard $N$-dimensional representation for all values of $N$
(see Exercise~\ref{q-sl_N} on page~\pageref{q-sl_N} for details).
\begin{table}
$$\begin{array}{c|l}
3_1 & (2a^2 -a^4) + a^2z^2 \\
4_1 & (a^{-2}-1+a^2) - z^2 \\
5_1 & (3a^4-2a^6) + (4a^4 -a^6)z^2 + a^4z^4 \\
5_2 & (a^2+a^4-a^6) + (a^2+a^4)z^2 \\
6_1 & (a^{-2}-a^2+a^4) + (-1-a^2)z^2 \\
6_2 & (2-2a^2+a^4) + (1-3a^2+a^4)z^2 - a^2z^4 \\
6_3 & (-a^{-2}+3-a^2) + (-a^{-2}+3-a^2)z^2 + z^4 \\
7_1 & (4a^6-3a^8) + (10a^6-4a^8)z^2 + (6a^6-a^8)z^4 + a^6z^6 \\
7_2 & (a^2+a^6-a^8) + (a^2+a^4+a^6)z^2 \\
7_3 & (a^{-4}+2a^{-6}-2a^{-8}) + (3a^{-4}+3a^{-6}-a^{-8})z^2 +
(a^{-4}+a^{-6})z^4 \\
7_4 & (2a^{-4}-a^{-8}) + (a^{-2}+2a^{-4}+a^{-6})z^2 \\
7_5 & (2a^4-a^8) + (3a^4+2a^6-a^8)z^2 + (a^4+a^6)z^4 \\
7_6 & (1-a^2+2a^4-a^6) + (1-2a^2+2a^4)z^2 - a^2z^4 \\
7_7 & (a^{-4}-2a^{-2}+2) + (-2a^{-2}+2-a^2)z^2 + z^4 \\
8_1 & (a^{-2}-a^4+a^6) + (-1-a^2-a^4)z^2 \\
8_2 & (3a^2-3a^4+a^6) + (4a^2-7a^4+3a^6)z^2 + (a^2-5a^4+a^6)z^4 -a^4z^6 \\
8_3 & (a^{-4}-1+a^4) + (-a^{-2}-2-a^2)z^2 \\
8_4 & (a^4-2+2a^{-2}) + (a^4-2a^2-3+a^{-2})z^2 + (-a^2-1)z^4 \\
8_5 & (4a^{-2}-5a^{-4}+2a^{-6}) + (4a^{-2}-8a^{-4}+3a^{-6})z^2 \\
       & + (a^{-2}-5a^{-4}+a^{-6})z^4 -a^{-4}z^6 \\
8_6 & (2-a^2-a^4+a^6) + (1-2a^2-2a^4+a^6)z^2 + (-a^2-a^4)z^4 \\
8_7 & (-2a^{-4}+4a^{-2}-1) + (-3a^{-4}+8a^{-2}-3)z^2 + (-a^{-4}+5a^{-2}-1)z^4 \\
       & + a^{-2}z^6 \\
8_8 & (-a^{-4}+a^{-2}+2-a^2) + (-a^{-4}+2a^{-2}+2-a^2)z^2 + (a^{-2}+1)z^4 \\
8_9 & (2a^{-2}-3+2a^2) + (3a^{-2}-8+3a^2)z^2 + (a^{-2}-5+a^2)z^4 - z^6 \\
8_{10} & (-3a^{-4}+6a^{-2}-2) + (-3a^{-4}+9a^{-2}-3)z^2 + (-a^{-4}+5a^{-2}-1)z^4\\
       & + a^{-2}z^6 \\
8_{11} & (1+a^2-2a^4+a^6) + (1-a^2-2a^4+a^6)z^2 + (-a^2-a^4)z^4 \\
8_{12} & (a^{-4}-a^{-2}+1-a^2+a^4) + (-2a^{-2}+1-2a^2)z^2 + z^4 \\
8_{13} & (-a^{-4}+2a^{-2}) + (-a^{-4}+2a^{-2}+1-a^2)z^2 + (a^{-2}+1)z^4 \\
8_{14} & 1 + (1-a^2-a^4+a^6)z^2 + (-a^2-a^4)z^4 \\
8_{15} & (a^4+3a^6-4a^8+a^{10}) + (2a^4+5a^6-3a^8)z^2 + (a^4+2a^6)z^4 \\
8_{16} & (-a^4+2a^2) + (-2a^4+5a^2-2)z^2 + (-a^4+4a^2-1)z^4 + a^2z^6 \\
8_{17} & (a^{-2}-1+a^2) + (2a^{-2}-5+2a^2)z^2 + (a^{-2}-4+a^2)z^4 - z^6 \\
8_{18} & (-a^{-2}+3-a^2) + (a^{-2}-1+a^2)z^2 + (a^{-2}-3+a^2)z^4 - z^6 \\
8_{19} & (5a^{-6}-5a^{-8}+a^{-10}) + (10a^{-6}-5a^{-8})z^2 +
(6a^{-6}-a^{-8})z^4 + a^{-6}z^6 \\
8_{20} & (-2a^4+4a^2-1) + (-a^4+4a^2-1)z^2 + a^2z^4 \\
8_{21} & (3a^2-3a^4+a^6) + (2a^2-3a^4+a^6)z^2 - a^4z^4\vspace{10pt}
\end{array}$$
\caption{HOMFLY polynomials of knots with up to 8 crossings}\label{HOMFLY_table}
\index{HOMFLY polynomial!table}
\index{Table of!HOMFLY polynomials}
\end{table}

Important properties of the HOMFLY polynomial are contained in the
following exercises.

\begin{xca}\

\begin{enumerate}
\item
Prove the uniqueness of such an invariant.
In other words,
prove that the relation above are sufficient
to compute the HOMFLY polynomial.
\item
\parbox[t]{3in}{Compute the HOMFLY polynomial for the knots
$3_1$, $4_1$ and compare your results with those given in
Table~\ref{HOMFLY_table}.}\\[1mm]
\item \parbox[t]{3in}{
Compare the HOMFLY polynomials of the Conway 
and Kinoshita-Terasaka knots on the right ({see, for instance, \cite{Sos}}).
\index{Conway knot}\index{Knot!Conway}
\index{Kinoshita--Terasaka knot}\index{Knot!Kinoshita--Terasaka}
}
\qquad\qquad\quad
$\makebox(0,0){\rb{50pt}{$\begin{array}{r} C=\risS{-15}{sos_kn2}{}{35}{10}{30}\\
                KT=\risS{-15}{sos_kn1}{}{35}{20}{45}\end{array}$}}
$\label{C-KT1}
\end{enumerate}
\end{xca}

\begin{xca}
Prove that the HOMFLY polynomial of a link is preserved when the
orientation of all components is reversed.
\end{xca}

\begin{xca}
({\em W.~B.~R.~Lickorish \cite{Lik}}) Prove that
\begin{enumerate}
\item 
$P(\ol{L}) = \ol{P(L)}$, where $\ol{L}$ is the mirror reflection
of $L$ and $\ol{P(L)}$ is the polynomial obtained from $P(L)$ by
substituting $-a^{-1}$ for $a$;

\item $P(K_1\#K_2)=P(K_1)\cdot P(K_2)$;
\item
\parbox[t]{3in}{$\displaystyle P(L_1\sqcup L_2)=
       \frac{a-a^{-1}}{z}\cdot P(L_1)\cdot P(L_2)$, where
$L_1\sqcup L_2$ means the split union of links\index{Split union}
(that is, the union of $L_1$ and $L_2$ such that each of these two links is
contained inside its own ball, and the two balls do not have common points);}
\hspace{55pt}
$\makebox(0,0){$\begin{array}{l}
  \rb{-50pt}{$8_8 =\ $} \risS{-70}{k8-08}{}{35}{20}{15}\vspace*{-12pt}\\
\hspace{-10pt}\rb{-50pt}{$\ol{10}_{129} =\ $} \risS{-70}{k10-129}{}{35}{20}{15}
 \end{array}$}
$
\item $P(8_8)=P(\ol{10}_{129})$.\\
\parbox[t]{3in}{These knots can be distinguished by the two-variable Kauffman
polynomial defined below.}
\end{enumerate}
\end{xca}

\subsection{Two-variable Kauffman polynomial}\label{t-var-Kauf}
\index{Kauffman polynomial}

In \cite{Ka4}, L.~Kauffman found another invariant Laurent
polynomial $F(L)$ in two variables $a$ and $z$. Firstly, for a
unoriented link diagram $D$ we define a polynomial $\Lambda(D)$
which is invariant under Reidemeister moves $\Omega_2$ and
$\Omega_3$ and satisfies the relations\vspace{-10pt}
$$\Lambda(\rb{-4.2mm}{\ig[width=10mm]{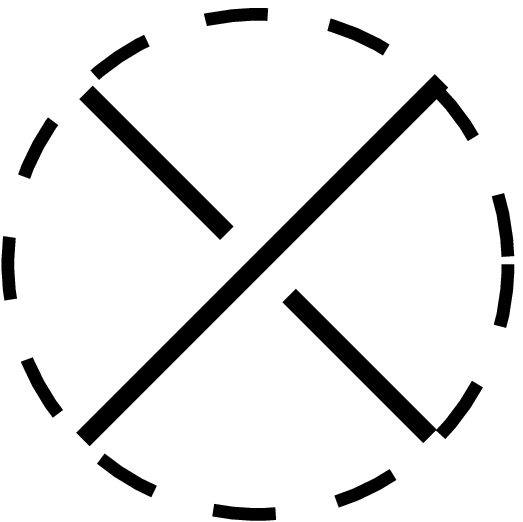}})\ +\
\Lambda(\rb{-12pt}{\rotatebox{90}{\ig[width=10mm]{unor_plcr.eps}}})
    \ = \  z\Bigl(
\Lambda(\rb{-4.2mm}{\ig[width=10mm]{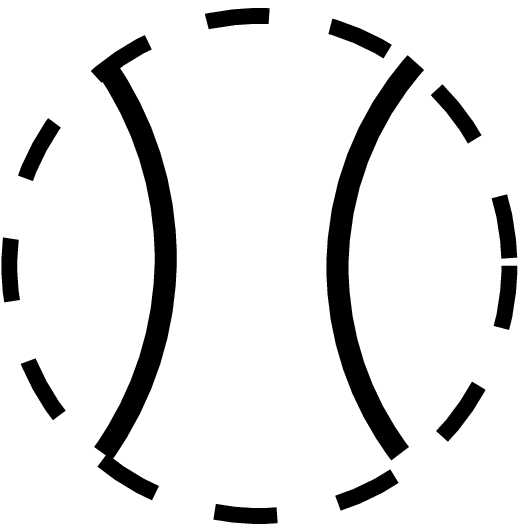}})\ +\
\Lambda(\rb{-12pt}{\rotatebox{90}{\ig[width=10mm]{unor_zero.eps}}})
\Bigr)\,,\index{Kauffman polynomial}
$$\vspace{-5pt}
$$\Lambda(\rb{-4.2mm}{\ig[width=10mm]{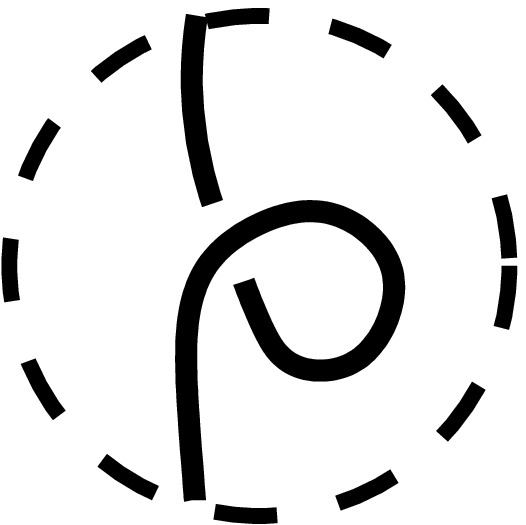}})
    \ = \  a
\Lambda(\rb{-4.2mm}{\ig[width=10mm]{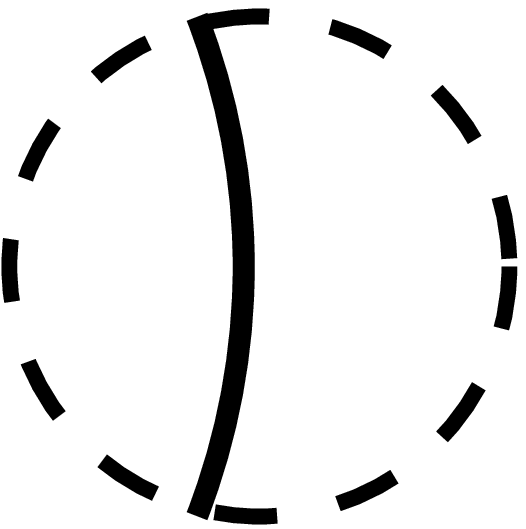}})\,,
\qquad\qquad
\Lambda(\rb{-4.2mm}{\ig[width=10mm]{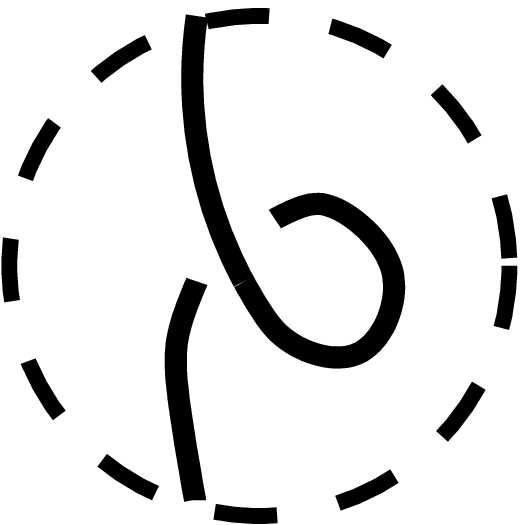}})
    \ = \  a^{-1}
\Lambda(\rb{-4.2mm}{\ig[width=10mm]{unor_oneup.eps}})
\,,
$$\vspace{-5pt}
and the initial condition $\Lambda(\chd{cd0})=1$.

Now, for any diagram $D$ of an oriented link $L$ we put
$$F(L) := a^{-w(D)} \Lambda(D).$$
It turns out
that this polynomial is equivalent to the collection of the quantum
invariants associated with the Lie algebra $\so_N$ and its standard
$N$-dimensional representation for all values of $N$ (see \cite{Tur3}).
\setlength{\intextsep}{2pt}
\setlength{\floatsep}{0pt}
\setlength{\textfloatsep}{0pt}
\begin{table}[h!]
$$\begin{array}{c|l}
3_1 & (-2a^2 - a^4) + (a^3 + a^5)z + (a^2 + a^4)z^2\\

4_1 & (-a^{-2} - 1 - a^2) + (-a^{-1} - a)z + (a^{-2} + 2 + a^2) z^2
+ (a^{-1} + a)z^3\\

5_1 & (3a^4 +2a^6)+ (-2a^5 -a^7 +a^9)z+ (-4a^4 -3a^6 +a^8)z^2\\
    &\quad + (a^5 +a^7)z^3 +(a^4 +a^6)z^4\\

5_2 & (-a^2 + a^4 + a^6) + (-2a^5 - 2a^7)z + (a^2 - a^4 - 2a^6)z^2\\
    &\quad + (a^3 + 2a^5 + a^7)z^3 + (a^4 + a^6)z^4\\

6_1 & (-a^{-2} + a^2 + a^4) + (2a+ 2a^3)z + (a^{-2} - 4a^2 - 3a^4)z^2\\
    &\quad + (a^{-1}- 2a-3a^3)z^3 + (1 + 2a^2 + a^4)z^4 + (a+ a^3)z^5\\

6_2 & (2 + 2a^2 + a^4) + (-a^3 - a^5)z + (-3 - 6a^2 - 2a^4 + a^6)z^2 \\
    &\quad + (-2a+ 2a^5)z^3 + (1 + 3a^2 + 2a^4)z^4 + (a+ a^3)z^5\\

6_3 & (a^{-2} + 3 + a^2)+ (-a^{-3} - 2a^{-1} - 2a - a^3)z
+ (-3a^{-2} - 6 - 3a^2)z^2\\
    &\quad + (a^{-3} + a^{-1} + a+ a^3)z^3 + (2a^{-2} + 4 + 2a^2)z^4 + (a^{-1} + a)z^5\\

7_1 & (-4a^6 - 3a^8) + (3a^7 + a^9 - a^{11} + a^{13})z
+ (10a^6 + 7a^8 - 2a^{10} + a^{12})z^2\\
    &\quad + (-4a^7 - 3a^9 + a^{11})z^3 + (-6a^6 - 5a^8 + a^{10})z^4 +
       (a^7 + a^9)z^5 \\
    &\quad + (a^6 + a^8)z^6\\

7_2 & (-a^2 - a^6 - a^8) + (3a^7 + 3a^9)z + (a^2 + 3a^6 + 4a^8)z^2 \\
    &\quad + (a^3 - a^5 - 6a^7 - 4a^9)z^3 + (a^4 - 3a^6 - 4a^8)z^4
+ (a^5 + 2a^7 + a^9)z^5 \\
    &\quad + (a^6 + a^8)z^6\\

7_3 & (-2a^{-8} - 2a^{-6} + a^{-4}) + (-2a^{-11} + a^{-9} + 3a^{-7})z\\
    &\quad + (-a^{-10} + 6a^{-8} + 4a^{-6} - 3a^{-4})z^2  + (a^{-11} - a^{-9} - 4a^{-7} - 2a^{-5})z^3 \\
    &\quad + (a^{-10} - 3a^{-8} - 3a^{-6} + a^{-4})z^4
+ (a^{-9} + 2a^{-7} + a^{-5})z^5\\
    &\quad + (a^{-8} + a^{-6})z^6\\

7_4 & (-a^{-8} +2a^{-4})+ (4a^{-9} +4a^{-7})z
+ (2a^{-8} -3a^{-6} -4a^{-4} +a^{-2})z^2\\
    &\quad + (-4a^{-9} -8a^{-7} - 2a^{-5} + 2a^{-3})z^3 + (-3a^{-8} + 3a^{-4})z^4 \\
    &\quad + (a^{-9} + 3a^{-7} + 2a^{-5})z^5 + (a^{-8} + a^{-6})z^6\\

7_5 & (2a^4 - a^8) + (-a^5 + a^7 + a^9 - a^{11})z + (-3a^4 + a^8 -2a^{10})z^2 \\
    &\quad + (-a^5 - 4a^7 - 2a^9 + a^{11})z^3 + (a^4 - a^6 + 2a^{10})z^4 \\
    &\quad + (a^5 + 3a^7 + 2a^9)z^5 + (a^6 + a^8)z^6\\

7_6 & (1 + a^2 + 2a^4 + a^6) + (a+ 2a^3 - a^7)z + (-2 - 4a^2 - 4a^4 - 2a^6)z^2\\
    &\quad + (-4a-6a^3 - a^5 + a^7)z^3 + (1 + a^2 + 2a^4 + 2a^6)z^4 \\
    &\quad + (2a+ 4a^3 + 2a^5)z^5 + (a^2 + a^4)z^6\\

7_7 & (a^{-4} + 2a^{-2} + 2) + (2a^{-3} + 3a^{-1} + a)z
+ (-2a^{-4} - 6a^{-2} - 7 - 3a^2)z^2\\
    &\quad + (-4a^{-3} - 8a^{-1} - 3a+ a^3)z^3 + (a^{-4} + 2a^{-2} + 4 + 3a^2)z^4\\
    &\quad + (2a^{-3} + 5a^{-1}+ 3a)z^5 + (a^{-2} + 1)z^6\\
\end{array}$$
\caption{Kauffman polynomials of knots with up to 7 crossings}
\label{Kauffman_table}\index{Kauffman polynomial!table}
\index{Table of!Kauffman polynomials}
\end{table}
\begin{table}[h!]
$$\begin{array}{c|l}
8_1 & (-a^{-2} - a^4 - a^6)+ (-3a^3 - 3a^5)z+ (a^{-2} + 7a^4 + 6a^6)z^2\\
    & + (a^{-1} - a + 5a^3 + 7a^5)z^3 + (1 - 2a^2 - 8a^4 - 5a^6)z^4\\
    & + (a- 4a^3 - 5a^5)z^5
+ (a^2 + 2a^4 + a^6)z^6 + (a^3 + a^5)z^7\\

8_2 & (-3a^2 - 3a^4 - a^6) + (a^3 + a^5 - a^7 - a^9)z \\
    & + (7a^2 + 12a^4 + 3a^6 -a^8 + a^{10})z^2 + (3a^3 - a^5 - 2a^7 + 2a^9)z^3 \\
    & + (-5a^2 - 12a^4 - 5a^6 + 2a^8)z^4 + (-4a^3 - 2a^5 +2a^7)z^5 \\
    & + (a^2 + 3a^4 + 2a^6)z^6 + (a^3 + a^5)z^7\\

8_3 & (a^{-4} - 1 + a^4) + (-4a^{-1} - 4a)z + (-3a^{-4} + a^{-2} + 8 + a^2 - 3a^4)z^2\\
    &+ (-2a^{-3} + 8a^{-1} + 8a- 2a^3)z^3 + (a^{-4} - 2a^{-2} - 6 - 2a^2 + a^4)z^4 \\
    &+ (a^{-3} - 4a^{-1} - 4a+ a^3)z^5 + (a^{-2} + 2 + a^2)z^6 + (a^{-1} + a)z^7\\

8_4 & (-2a^{-2} - 2 + a^4) + (-a^{-1} + a+ 2a^3)z \\
    &+ (7a^{-2} + 10 - a^2 - 3a^4 +a^6)z^2 + (4a^{-1} - 3a- 5a^3 + 2a^5)z^3 \\
    &+ (-5a^{-2} - 11 - 3a^2 + 3a^4)z^4 + (-4a^{-1} - a+ 3a^3)z^5 \\
    &+ (a^{-2} + 3 + 2a^2)z^6 + (a^{-1} + a)z^7\\

8_5 & (-2a^{-6} - 5a^{-4} - 4a^{-2}) + (4a^{-7} + 7a^{-5} + 3a^{-3})z \\
    &+ (a^{-10} - 2a^{-8} + 4a^{-6} + 15a^{-4} + 8a^{-2})z^2
+ (2a^{-9} - 8a^{-7} - 10a^{-5})z^3\\
    & + (3a^{-8} - 7a^{-6} - 15a^{-4} - 5a^{-2})z^4
+ (4a^{-7} + a^{-5} - 3a^{-3})z^5 \\
    &+ (3a^{-6} + 4a^{-4} + a^{-2})z^6 + (a^{-5} + a^{-3})z^7\\

8_6 & (2 + a^2 - a^4 - a^6) + (-a- 3a^3 - a^5 + a^7)z \\
    &+ (-3 - 2a^2 + 6a^4 + 3a^6 - 2a^8)z^2 + (-a+ 5a^3 + 2a^5 - 4a^7)z^3\\
    & + (1 - 6a^4 - 4a^6 + a^8)z^4 + (a- 2a^3 - a^5 + 2a^7)z^5\\
    & + (a^2 + 3a^4 + 2a^6)z^6 + (a^3 + a^5)z^7\\

8_7 & (-2a^{-4} - 4a^{-2} - 1) + (-a^{-7} + 2a^{-3} + 2a^{-1} + a)z \\
    &+ (-2a^{-6} + 4a^{-4} + 12a^{-2} + 6)z^2 +
       (a^{-7} - a^{-5} - 2a^{-3} - 3a^{-1} - 3a)z^3 \\
    & + (2a^{-6} - 3a^{-4} - 12a^{-2} - 7)z^4 + (2a^{-5} - a^{-1} + a)z^5 \\
    &+ (2a^{-4} + 4a^{-2} + 2)z^6 + (a^{-3} + a^{-1})z^7\\

8_8 & (-a^{-4} - a^{-2} + 2 + a^2) + (2a^{-5} + 3a^{-3} + a^{-1} - a- a^3)z \\
    &+ (4a^{-4} + 5a^{-2} - 1 - 2a^2)z^2 + (-3a^{-5} - 5a^{-3} - 3a^{-1} + a^3)z^3\\
    & + (-6a^{-4} - 9a^{-2} - 1 + 2a^2)z^4 + (a^{-5} + a^{-1} + 2a)z^5 \\
    &+ (2a^{-4} + 4a^{-2} + 2)z^6 + (a^{-3} + a^{-1})z^7\\

8_9 & (-2a^{-2} - 3 - 2a^2) + (a^{-3} + a^{-1} + a+ a^3)z \\
    &+ (-2a^{-4} + 4a^{-2} + 12+ 4a^2 - 2a^4)z^2 + (-4a^{-3} - a^{-1} - a- 4a^3)z^3 \\
    &+ (a^{-4} - 4a^{-2} - 10 - 4a^2 + a^4)z^4 + (2a^{-3} + 2a^3)z^5 \\
    &+ (2a^{-2} + 4 + 2a^2)z^6 + (a^{-1} + a)z^7\\

8_{10} & (-3a^{-4} - 6a^{-2} - 2) + (-a^{-7} + 2a^{-5} + 6a^{-3} + 5a^{-1} + 2a)z \\
    &+ (-a^{-6} + 6a^{-4} + 12a^{-2} + 5)z^2 +
               (a^{-7} - 3a^{-5} - 9a^{-3} - 8a^{-1} - 3a)z^3 \\
    &+ (2a^{-6} - 5a^{-4} - 13a^{-2} - 6)z^4 + (3a^{-5} + 3a^{-3} + a^{-1} + a)z^5 \\
    &+ (3a^{-4} + 5a^{-2} + 2)z^6 + (a^{-3} + a^{-1})z^7\vspace{10pt}\\
    \end{array}$$
\setcounter{table}{0}
\caption{Kauffman polynomials of knots with 8 crossings}
\index{Kauffman polynomial!table}\index{Table of!Kauffman polynomials}
\end{table}

\begin{table}[p]
$$\begin{array}{c|l}
8_{11} & (1 - a^2 - 2a^4 - a^6) + (a^3 + 3a^5 + 2a^7)z + (-2 + 6a^4 + 2a^6 - 2a^8)z^2\\
    &+ (-3a- 2a^3 - 3a^5 - 4a^7)z^3 + (1 - 2a^2 - 7a^4 - 3a^6 + a^8)z^4 \\
    &+ (2a+ a^3 + a^5 + 2a^7)z^5 + (2a^2 + 4a^4 + 2a^6)z^6 + (a^3 + a^5)z^7 \\

8_{12} & (a^{-4} + a^{-2} + 1 + a^2 + a^4) + (a^{-3} + a^3)z \\
    &+ (-2a^{-4} -2a^{-2} - 2a^2 - 2a^4)z^2 + (-3a^{-3} - 3a^{-1} - 3a- 3a^3)z^3 \\
    &+ (a^{-4} - a^{-2} - 4 - a^2 + a^4)z^4 + (2a^{-3} + 2a^{-1} + 2a+ 2a^3)z^5 \\
    &+ (2a^{-2} + 4 + 2a^2)z^6 + (a^{-1} + a)z^7\\

8_{13} & (-a^{-4} - 2a^{-2}) + (2a^{-5} + 4a^{-3} + 3a^{-1} + a)z
+ (5a^{-4} +7a^{-2} - 2a^2)z^2 \\
    &+ (-3a^{-5} - 7a^{-3} - 9a^{-1} - 4a+ a^3)z^3
+ (-6a^{-4} - 11a^{-2} - 2 + 3a^2)z^4 \\
    &+ (a^{-5} +a^{-3} + 4a^{-1} + 4a)z^5
+ (2a^{-4} + 5a^{-2} + 3)z^6 + (a^{-3} + a^{-1})z^7\\

8_{14} & 1 + (a+ 3a^3 + 3a^5 + a^7)z + (-2 - a^2 + 3a^4 + a^6 - a^8)z^2 \\
    &+ (-3a- 6a^3 - 8a^5 - 5a^7)z^3 + (1 - a^2 - 7a^4 - 4a^6 + a^8)z^4 \\
    &+ (2a+ 3a^3 + 4a^5 + 3a^7)z^5 + (2a^2 + 5a^4 +3a^6)z^6 + (a^3 + a^5)z^7\\

8_{15} & (a^4 - 3a^6 - 4a^8 - a^{10}) + (6a^7 + 8a^9 + 2a^{11})z \\
    &+ (-2a^4 + 5a^6 + 8a^8 -a^{12})z^2 + (-2a^5 - 11a^7 - 14a^9 - 5a^{11})z^3 \\
    &+ (a^4 - 5a^6 - 10a^8 - 3a^{10} + a^{12})z^4
+ (2a^5 + 5a^7 + 6a^9 + 3a^{11})z^5 \\
    &+ (3a^6 + 6a^8 + 3a^{10})z^6 + (a^7 + a^9)z^7\\

8_{16} & (-2a^2 -a^4)+ (a^{-1} +3a +4a^3 +2a^5)z+ (5+10a^2 +4a^4 -a^6)z^2\\
    &+ (-2a^{-1} - 6a- 10a^3 - 5a^5 + a^7)z^3 + (-8 - 18a^2 - 7a^4 + 3a^6)z^4 \\
    &+ (a^{-1} - a+ 3a^3 + 5a^5)z^5 + (3 + 8a^2 + 5a^4)z^6 + (2a+ 2a^3)z^7\\

8_{17} & (-a^{-2} - 1 - a^2) + (a^{-3} + 2a^{-1} + 2a+ a^3)z \\
    &+ (-a^{-4} + 3a^{-2} + 8 + 3a^2 - a^4)z^2 + (-4a^{-3} - 6a^{-1} - 6a- 4a^3)z^3 \\
    &+ (a^{-4} - 6a^{-2} - 14 - 6a^2 + a^4)z^4 + (3a^{-3} + 2a^{-1} + 2a+ 3a^3)z^5 \\
    &+ (4a^{-2} + 8 + 4a^2)z^6 + (2a^{-1} + 2a)z^7\\

8_{18} & (a^{-2} + 3 + a^2) + (a^{-1} + a)z + (3a^{-2} + 6 + 3a^2)z^2 \\
    &+ (-4a^{-3} - 9a^{-1} - 9a-4a^3)z^3 + (a^{-4} - 9a^{-2} - 20 - 9a^2 + a^4)z^4 \\
    &+ (4a^{-3} + 3a^{-1} + 3a+ 4a^3)z^5
+ (6a^{-2} +12 + 6a^2)z^6 + (3a^{-1} + 3a)z^7\\

8_{19} & (-a^{-10} -5a^{-8} -5a^{-6})+ (5a^{-9} +5a^{-7})z
+ (10a^{-8} +10a^{-6})z^2\\
    &+ (-5a^{-9} -5a^{-7})z^3 + (-6a^{-8} - 6a^{-6})z^4 +
         (a^{-9} + a^{-7})z^5 \\
    &+ (a^{-8} + a^{-6})z^6\\

8_{20} & (-1 - 4a^2 - 2a^4) + (a^{-1} + 3a+ 5a^3 + 3a^5)z
+ (2 + 6a^2 + 4a^4)z^2 \\
    &+ (-3a-7a^3 - 4a^5)z^3 + (-4a^2 - 4a^4)z^4 + (a+ 2a^3 + a^5)z^5 \\
    &+ (a^2 + a^4)z^6\\

8_{21} & (-3a^2 - 3a^4 - a^6) + (2a^3 + 4a^5 + 2a^7)z
+ (3a^2 + 5a^4 - 2a^8)z^2 \\
    &+ (-a^3 -6a^5 - 5a^7)z^3 + (-2a^4 - a^6 + a^8)z^4 + (a^3 + 3a^5 + 2a^7)z^5 \\
    &+ (a^4 + a^6)z^6\vspace{10pt}\\
\end{array}$$
\setcounter{table}{0}
\caption{Kauffman polynomials of knots with 8 crossings (Continuation)}
\index{Kauffman polynomial!table}\index{Table of!Kauffman polynomials}
\end{table}

As in the previous section, we conclude with a series of exercises
with additional information on the Kauffman polynomial.

\begin{xca}
Prove that the defining relations are sufficient to compute the
Kauffman polynomial.
\end{xca}
\begin{xca}
Compute the Kauffman polynomial for the knots
$3_1$, $4_1$ and compare the results with those given in the above table.
\end{xca}
\begin{xca} Prove that the Kauffman polynomial of a knot is preserved when
the knot orientation is reversed.
\end{xca}
\begin{xca}
({\it W.~B.~R.~Lickorish \cite{Lik}}) Prove that
\begin{enumerate}
\item \parbox[t]{3in}{
$F(\ol{L}) = \ol{F(L)}$, where $\ol{L}$ is the mirror reflection
of $L$, and $\ol{F(L)}$ is the polynomial obtained from $F(L)$ by
substituting $a^{-1}$ for $a$;}
\item $F(K_1\#K_2)=F(K_1)\cdot F(K_2)$;
\item
\parbox[t]{3in}{$\displaystyle F(L_1\sqcup L_2)=
 \Bigl((a+a^{-1})z^{-1}-1\Bigr)\cdot F(L_1)\cdot F(L_2)$,
where $L_1\sqcup L_2$ means the split union of links;}
\qquad\qquad\quad
$\rb{-25pt}{\makebox(0,0){$\stackrel{
  \rb{35pt}{$11^a_{30} =\ $} \risS{20}{k11-30}{}{35}{20}{0}}
{\rb{90pt}{$11^a_{189} =\ $} \risS{70}{k11-189}{}{35}{0}{15}}$}}
$
\item $F(11^a_{30})=F(11^a_{189})$;\\
(these knots can be distinguished by the Conway and, hence, by the
HOMFLY polynomial; note that we use the Knotscape numbering of 
knots \cite{Knsc}, while in \cite{Lik} the old Perko's notation is used).

\item $F(L^*) = a^{4 lk(K, L-K)} F(L)$, where the link
$L^*$ is obtained from an oriented link $L$ by
reversing the orientation of a connected component $K$.
\end{enumerate}
\end{xca}

\subsection{Comparative strength of polynomial invariants}

Let us say that an invariant $I_1$ {\em dominates} an invariant
$I_2$, if the equality $I_1(K_1)=I_1(K_2)$ for any two knots $K_1$ and
$K_2$ implies the equality $I_2(K_1)=I_2(K_2)$. Denoting this
relation by arrows, we have the following comparison
chart:
$$\xymatrix{
&&\fbox{\mbox{HOMFLY}}
   \ar[lldd]_(.6){\begin{array}{rr}\scriptstyle a=1\vspace{-3pt}\\
             \scriptstyle z=x^{1/2}-x^{-1/2}\hspace*{15pt}\end{array}\hspace*{-17pt}}
   \ar[dd]^{\hspace{-5pt}\begin{array}{ll}\scriptstyle a=1\vspace{-5pt}\\
                             \scriptstyle z=t\end{array}}
   \ar[rrdd]^(.6){\hspace*{-20pt}\begin{array}{ll}\scriptstyle
                                  z=t^{1/2}-t^{-1/2}\vspace{-3pt}\\
             \scriptstyle \hspace*{15pt}a=t^{-1}\end{array}}
 &&& \fbox{\mbox{Kauffman}}
   \ar[ldd]^(.3){\hspace*{-22pt}\begin{array}{rr}
             \scriptstyle a=-t^{-3/4}\ \vspace{-3pt}\\
             \scriptstyle z=t^{1/4}+t^{-1/4}\end{array}} \\
&&&&&\\
\fbox{\mbox{Alexander}}
   \ar@/^1.5pc/[rr]_{\rb{-18pt}{$\ \ \scriptstyle x^{1/2}-x^{-1/2}=t$}} &&
\fbox{\mbox{Conway}} \ar@/^1.5pc/[ll] && \fbox{\mbox{Jones}} &
}
$$
(the absence of an arrow between the two invariants means that neither
of them dominates the other).

{\bf Exercise.}
Find in this chapter all the facts sufficient to justify this chart.

\begin{xcb}{Exercises}

\begin{enumerate}

\item {\bf Bridge number.} \label{ex:bridge-n}
The {\em bridge number} \index{Bridge number} $b(K)$ of a knot $K$
can be defined as the minimal number of local maxima of the
projection of the knot onto a straight line, where the minimum is
taken over all projections and over all closed curves in $\R^3$
representing the knot. Show that that
$$b(K_1\#K_2)=b(K_1)+b(K_2)-1\ .$$
Knots of bridge number 2 are also called {\em rational knots}.
\index{Rational knots}\index{Knot!rational}

\item Prove that the Conway and the Jones polynomials of a knot are
preserved when the knot orientation is reversed.
\index{Conway polynomial}\index{Jones polynomial}

\item
Compute the Conway and the Jones polynomials for the links from
Section \ref{def:link}, page \pageref{def:link}, with some
orientations of their components.

\item
A link is called {\it split}\index{Link!split} if it is equivalent to a link
which has some components in a ball while the other components are located
outside of the ball. Prove that the Conway polynomial of a split link is
trivial: $\CP(L)=0$.

\item
For a split link $L_1\sqcup L_2$ prove that
$$J(L_1\sqcup L_2)=(-t^{1/2}-t^{-1/2})\cdot J(L_1)\cdot J(L_2)\ .$$

\item\label{sum_conway}
Prove that $\CP(K_1\#K_2)=\CP(K_1)\cdot \CP(K_2)$.

\item\label{sum_jones}
Prove that $J(K_1\#K_2)=J(K_1)\cdot J(K_2)$.

\item (cf. {\it J.~H.~Conway} \cite{Con})\label{ex:C-rels}
Check that the Conway polynomial satisfies the following relations.
\begin{enumerate}
\item $C\Bigl(\chd{conpo-pp}\Bigr) + C\Bigl(\chd{conpo-mm}\Bigr)\
=\  (2+t^2) C\Bigl(\chd{conpo}\Bigr)$\ ;

\item $C\Bigl(\chd{conpr-mm}\Bigr) + C\Bigl(\chd{conpr-pp}\Bigr)\
=\  2 C\Bigl(\chd{conpr}\Bigr)$\ ;

\item $C\Bigl(\chd{tr_p_res1}\Bigr) + C\Bigl(\chd{tr_p_res3}\Bigr)\
=\  C\Bigl(\chd{tr_p_res4}\Bigr) + C\Bigl(\chd{tr_p_res2}\Bigr)$\ .
\end{enumerate}

\item Compute the Conway polynomials of the Conway and the
Kinoshita--Terasaka knots (see page \pageref{C-KT1}).

\item Prove that for any knot $K$ the Conway polynomial $\CP(K)$ is an even
polynomial in $t$ and its constant term is equal to 1:
$$\CP(K)=1+c_2(K)t^2+c_4(K)t^4+\dots$$

\item
Let $L$ be a link with two components $K_1$ and $K_2$. Prove that the Conway polynomial
$\CP(L)$ is an odd polynomial in $t$ and its lowest coefficient is equal to the linking
number $lk(K_1,K_2)$:\index{Linking number}
$$\CP(L) = lk(K_1,K_2)t + c_3(L)t^3 + c_5(L)t^5+\dots$$

\item
Prove that for a link $L$ with $k$ components the Conway polynomial
$\CP(L)$ is divisible by $t^{k-1}$ and is odd or even depending on the parity of $k$:
$$\CP(L) = c_{k-1}(L)t^{k-1} + c_{k+1}(L)t^{k+1} + c_{k+3}(L)t^{k+3}+\dots$$

\item
For a knot $K$, show that $\CP(K)\bigr|_{t=2i}\equiv 1 \mbox{\ or\ }
5\ (\mbox{mod}\ 8)$ depending of the parity of $c_2(K)$. The
reduction of $c_2(K)$ modulo 2 is called the {\it Arf invariant} of
$K$.\index{Arf invariant}

\item
Show that $J(L)\bigr|_{t=-1}=\CP(L)\bigr|_{t=2i}$ for any link $L$.
The absolute value of this number is called the {\em determinant} \index{Determinant of a
link}\index{Link! determinant} of the link $L$. 

{\sl Hint.} Choose $\sqrt{t}$ in such a way that $\sqrt{-1}=-i$.

\item\label{rel:switch-J}
Check the following {\it switching formula} for the Jones polynomial.
\index{Jones polynomial!switching formula}
$$J(\lrints) - tJ(\rlints) = t^{3\lambda_0}(1-t)J(\risS{-12}{k-infty}{}{28}{0}{0})\ ,
$$
where $\lambda_0$ is the linking number of two components of the
link, $\risS{-2}{twoup}{}{29}{0}{0}$, obtained by smoothing the
crossing according to the orientation. Note that the knot in the right
hand side of the formula is unoriented. That is because such a
smoothing destroys the orientation. Since the Jones polynomial does
not distinguish the orientation of a knot, we may choose it arbitrarily.

\item\label{rel:int-cr-J}
{\bf Interlacing crossings formulae.} \index{Jones
polynomial!interlacing crossings formulae} Suppose $K_{++}$ is a
knot diagram with two positive crossings which are {\em interlaced}.
That means when we trace the knot we first past the first crossing,
then the second, then again the first, and after that the second.
Consider the following four knots and one link:
$$\hspace*{30pt}
  \risS{-40}{k-int-pp}{\put(10,70){$K_{++}$}}{42}{35}{45}\qquad
  \risS{-40}{k-int-zz}{\put(12,70){$K_{00}$}}{42}{0}{0}\qquad
  \risS{-40}{k-int-zi}{\put(10,70){$K_{0\infty}$}}{42}{0}{0}\qquad
  \risS{-40}{k-int-im}{\put(10,70){$K_{\infty-}$}}{42}{0}{0}\qquad\quad
  \risS{-40}{k-int-zp}{\put(12,70){$L_{0+}$}}{42}{0}{0}
$$
Check that the Jones polynomial satisfies the relation
$$J(K_{++}) =  t J(K_{00}) + t^{3\lambda_{0+}}
     \Bigl(J(K_{0\infty})- t J(K_{\infty-})\Bigr)\ ,
$$
where $\lambda_{0+}$ is the linking number of two components of the link $L_{0+}$.

Check the similar relations for $K_{+-}$ and $K_{--}$:
$$J(K_{+-}) =  J(K_{00}) + t^{3\lambda_{0-}+1}
   \Bigl(J(K_{0\infty})-J(K_{\infty+})\Bigr)\ ,
$$
$$J(K_{--}) =  t^{-1}J(K_{00}) + t^{3\lambda_{0-}}
   \Bigl(J(K_{0\infty})-t^{-1}J(K_{\infty+})\Bigr)\ .
$$
If a knot diagram does not contain interlacing crossings then it
represents the unknot. Thus the three relations above allow to
compute the Jones polynomial for knots recursively without referring
to links.

\item \label{ex:J-rels}
Show that the Jones polynomial satisfies the following relations.
\begin{enumerate}
\item $t^{-2}J\Bigl(\chd{conpo-pp}\Bigr) + t^2J\Bigl(\chd{conpo-mm}\Bigr)\
   =\  (t+t^{-1}) J\Bigl(\chd{conpo}\Bigr)$\ ;

\item $tJ\Bigl(\chd{conpr-mm}\Bigr) + t^{-1}J\Bigl(\chd{conpr-pp}\Bigr)\
=\  (t+t^{-1}) J\Bigl(\chd{conpr}\Bigr)$\ ;

\item
  $t^2J\Bigl(\chd{tr_p_res1}\Bigr) + t^{-2}J\Bigl(\chd{tr_p_res3}\Bigr)\
=\ t^{-2}J\Bigl(\chd{tr_p_res4}\Bigr) + t^2J\Bigl(\chd{tr_p_res2}\Bigr)$\ .
\end{enumerate}
Compare these relations with those of Exercise~\ref{ex:C-rels} for
the Conway polynomial.

\item \label{pr:J-pol}
Prove that for a link $L$ with an odd number of components, $J(L)$ is
a polynomial in $t$ and $t^{-1}$, and for a link $L$ with an even number of components
$J(L)=t^{1/2}\cdot(\text{a polynomial in }t\text{ and }t^{-1}\text{)}$.

\item Prove that for a link $L$ with $k$ components
$J(L)\bigr|_{t=1} = (-2)^{k-1}$. In particular,
$J(K)\bigr|_{t=1} = 1$ for a knot $K$.

\item
Prove that\quad
$\displaystyle \frac{d(J(K))}{dt}\Biggr|_{t=1} = 0$\quad for any knot $K$.

\item
Evaluate the Kauffman bracket $\langle L\rangle$ at $a=e^{\pi i/3}$,
$b=a^{-1}$, $c=-a^2-a^{-2}$. 
Deduce from here that $J(L)\bigr|_{t=e^{2\pi i/3}} = 1$.

{\sl Hint.} $\sqrt{t}=a^{-2}=e^{4\pi i/3}$.

\item
Let $L$ be a link with $k$ components. For odd (resp. even) $k$ let
$a_j$ ($j=0,1,2,\mbox{\ or\ }3$) be the sum of the coefficients of
$J(L)$ (resp. $J(L)/\sqrt{t}$, see problem \ref{pr:J-pol}) at $t^s$
for all $s\equiv j\ (\mbox{mod}\ 4)$.
\begin{enumerate}
\item For odd $k$, prove that $a_1=a_3$.
\item For even $k$, prove that $a_0+a_1=a_2+a_3$.
\end{enumerate}

\item ({\it W.~B.~R.~Lickorish \cite[Theorem 10.6]{Lik}})
Let $t=i$ with $t^{1/2}=e^{\pi i/4}$. Prove that for a knot $K$,
$J(K)\bigr|_{t=i} = (-1)^{c_2(K)}$.

\item\label{Jones_mirror}
For the mirror reflection $\overline{L}$  of a link $L$ prove that
$J(\overline{L})$ is obtained from $J(L)$ by substituting $t^{-1}$ for
$t$.

\item\label{Jones_reversal} For the link $L^*$ obtained from an oriented link $L$ by
reversing the orientation of one of its components $K$, prove that
$J(L^*) = t^{-3 lk(K, L-K)} J(L)$.

\item\neresh
Find a non-trivial knot $K$ with $J(K)=1$.

\item ({\it L.~Kauffman} \cite{Ka6}, {\it K.~Murasugi} \cite{Mur},
        {\it M.~Thistlethwaite} \cite{Th}).\label{ex_KMT_theorem}
Prove the  for a reduced alternating knot diagram $K$
(Section~\ref{altkn}) the number of crossings is equal to
$span(J(K))$, that is, to the difference beween the maximal and
minimal degrees of $t$ in the Jones polynomial $J(K)$. (This
exercise is not particularly difficult, although it solves a one
hundred years old conjecture of Tait\index{Conjecture!Tait}. Anyway,
the reader can find a rather simple solution in \cite{Tur1}.)

\item
Let $L$ be a link with $k$ components. Show that its HOMFLY
polynomial $P(L)$ is an even function in each of the variables $a$
and $z$ if $k$ is odd, and it is an odd function if $k$ is even.

\item
For a link $L$ with $k$ components, show that the lowest power of $z$
in its HOMFLY polynomial is $z^{-k+1}$. In particular the HOMFLY
polynomial $P(K)$ of a knot $K$ is a genuine polynomial in $z$. This
means that it does not contain terms with $z$ raised to a negative
power.

\item
For a knot $K$ let $p_0(a):=P(K)|_{z=0}$ be the constant term of the
HOMFLY polynomial. Show that its derivative at $a=1$ equals zero.

\item
Let $L$ be a link with two components $K_1$ and $K_2$. Consider
$P(L)$ as a Laurent polynomial in $z$ with coefficients in Laurent
polynomials in $a$. Let $p_{-1}(a)$ and $p_1(a)$ be the coefficients
of $z^{-1}$ and $z$. Check that $p_{-1}\bigr|_{a=1}=0$,\qquad
$p'_{-1}\bigr|_{a=1}=2$,\qquad $p''_{-1}\bigr|_{a=1}=-8 lk(K_1,K_2)
- 2$,\index{Linking number} and $p_1\bigr|_{a=1}= lk(K_1,K_2)$.

\item
Compute the HOMFLY polynomial of the four links shown on page~\pageref{solomon}. Note that,
according to the result,  the behaviour of the HOMFLY
polynomial under the change of orientation of one component is rather unpredictable.
(The same is true for the Conway polynomial, but not true
for the Jones and the Kauffmann polynomials.)

\item ({\it W.~B.~R.~Lickorish \cite{Lik}})
Prove that for an oriented link $L$ with $k$ components,\vspace{-10pt}
$$\Bigl.(J(L))^2\Bigr|_{t=-q^{-2}} =
  (-1)^{k-1}\Biggl.F(L)\Biggr|_{\rb{10pt}{$\begin{array}{ll}
             \scriptstyle a=q^3\ \vspace{-5pt}\\
             \scriptstyle z=q+q^{-1}\end{array}$}}\ ,\vspace{-5pt}
$$
where $J(L)$ is the Jones polynomial and $F(L)$ is the two-variable
Kauffman polynomial defined on page \pageref{t-var-Kauf}.

\item
Let $L$ be a link with $k$ components. Show that its two-variable
Kauffman polynomial $F(L)$ is an even function of both variables $a$
and $z$ (that is, it consists of monomials $a^iz^j$ with $i$ and $j$
of the same parity) if $k$ is odd, and it is an odd function
(different parities of $i$ and $j$) if $k$ is even.

\item
Prove that the Kauffman polynomial $F(K)$ of a knot $K$ is a genuine
polynomial in $z$.

\item
For a knot $K$ let $f_0(a):=F(K)|_{z=0}$ be the constant term of the
Kauffman polynomial. Show that it is related to the constant term of
the HOMFLY polynomial of $K$ as $f_0(a)=p_0(\sqrt{-1}\cdot a)$.

\item
{\bf Quantum $\sL_2$-invariant.} \index{Quantum invariant!$\sL_2$}
Let $\t(\cdot)$ and $\t^{fr}(\cdot)$ be the quantum invariants
 constructed in Sections~\ref{qis2}
and \ref{qis5} for the Lie algebra $\sL_2$ and its standard 2-dimensional
representation.
\begin{enumerate}
\item\label{exR1}
Prove the following dependence of $\t^{fr}(\cdot)$ on the first Reidemeister
move
$$\t^{fr}(\rb{-4.2mm}{\ig[width=10mm]{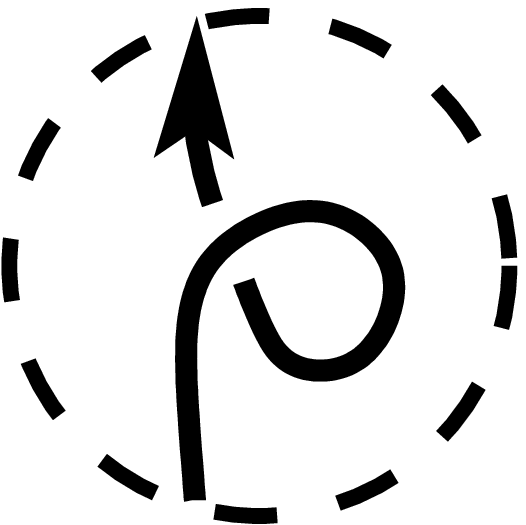}})=
 q^{3/4} \t^{fr}(\rb{-4.2mm}{\ig[width=10mm]{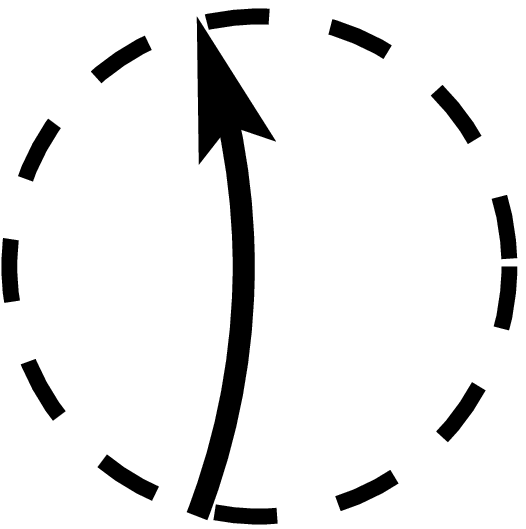}})\,.
$$
\item
Prove that $\t(\cdot)$ remains unchanged under the first Reidemeister
move.
\item Compute the value $\t(4_1)$.
\item Show that the $R$-matrix defined in page \pageref{Rmatix} satisfies the
equation
$$q^{1/4}R - q^{-1/4}R^{-1} \quad = \quad
  (q^{1/2}-q^{-1/2}) id_{_{V\ot V}}\,.$$
\item Prove that $\t^{fr}(\cdot)$ satisfies the skein relation
$$q^{1/4}\t^{fr}(\lrints)\ -\  q^{-1/4}\t^{fr}(\rlints)
  \quad = \quad (q^{1/2}-q^{-1/2})\t^{fr}(\twoup)\,.
$$
\item Prove that $\t(\cdot)$ satisfies the skein relation
$$q\t(\lrints)\ -\  q^{-1}\t(\rlints)
  \quad = \quad (q^{1/2}-q^{-1/2})\t(\twoup)\,.
$$
\item For any link $L$ with $k$ components prove that
$$\t^{fr}(L) = (-1)^k (q^{1/2}+q^{-1/2})\cdot
  \langle L\rangle \Bigr|_{a= -q^{1/4}}\,,$$
where $\langle\cdot\rangle$ is the Kauffman bracket defined on page
\pageref{kauf_br}.
\end{enumerate}

\item\label{q-sl_N}
{\bf Quantum $\sL_N$ invariants.} \index{Quantum invariant!$\sL_N$}
Let $V$ be an $N$ dimensional vector space of the standard
representation of the Lie algebra $\sL_N$ with a basis $e_1,\dots,
e_N$. Consider the operator $R:V\ot V\to V\ot V$ given by the
formulae
$$R(e_i\ot e_j)=\left\{\begin{array}{ll}
q^{\frac{-1}{2N}}e_j\ot e_i& \mbox{if\quad} i<j \vspace{8pt}\\
q^{\frac{N-1}{2N}}e_i\ot e_j& \mbox{if\quad} i=j \vspace{8pt}\\
q^{\frac{-1}{2N}}e_j\ot e_i +
    \Bigl(q^{\frac{N-1}{2N}}-q^{\frac{-N-1}{2N}}\Bigr)e_i\ot e_j& \mbox{if\quad} i>j
\end{array}\right.\label{Rmatrix_for_sl_N}\index{R-matrix!for $\sL_N$}$$
which for $N=2$ coincides with the operator from Section~\ref{qis2},
page \pageref{Rmatix}.
\begin{enumerate}
\item Prove that it satisfies the quantum Yang--Baxter equation
$$R_{12}R_{23}R_{12}=R_{23}R_{12}R_{23}\ ,$$
where $R_{ij}$ is the operator $R$ acting on the $i$th and $j$th
factors of $V\ot V\ot V$, that is, $R_{12}=R\ot\id_V$ and
$R_{23}=\id_V\ot R$.

\item Show that its inverse is given by the formulae
$$R^{-1}(e_i\ot e_j)=\left\{\begin{array}{ll}
q^{\frac1{2N}}e_j\ot e_i +
  \Bigl(-q^{\frac{N+1}{2N}}+q^{\frac{-N+1}{2N}}\Bigr)e_i\ot e_j& \mbox{if\quad} i<j
           \vspace{8pt}\\
q^{\frac{-N+1}{2N}}e_i\ot e_j& \mbox{if\quad} i=j \vspace{8pt}\\
q^{\frac1{2N}}e_j\ot e_i& \mbox{if\quad} i>j
\end{array}\right.$$

\item Check that\quad
$q^{\frac1{2N}}R - q^{\frac{-1}{2N}}R^{-1} = (q^{1/2}-q^{-1/2}) \id_{V\ot V}$.

\item Extending the assignments of operators for maximum/minimum tangles from page
\pageref{min-max-oper} we set:
$$\begin{array}{l@{\qquad}l}
\minr:\C\to V^*\ot V, &
     \minr(1):=\sum\limits_{k=1}^N q^{\frac{-N-1}2+k} e^k\ot e_k\ ;\vspace{5pt}\\
\minl:\C\to V\ot V^*, &
     \minl(1):=\sum\limits_{k=1}^N  e^k\ot e_k\ ;\vspace{5pt}\\
\maxr:V\ot V^*\to \C, &
     \maxr(e_i\ot e^j):=\left\{\begin{array}{ll}0&\mbox{if\ } i\not=j\\
             q^{\frac{N+1}2-i}&\mbox{if\ } i=j\end{array}\right.  \ ;\vspace{5pt}\\
\maxl:V^*\ot V\to \C, &
     \maxl(e^i\ot e_j):=\left\{\begin{array}{ll}0&\mbox{if\ } i\not=j\\
             1&\mbox{if\ } i=j\end{array}\right.  \ .\vspace{5pt}\\
\end{array}
$$
Prove that all these operators are consistent in the sense that
their appropriate combinations are consistent with the oriented
Turaev moves from page \pageref{Tur-mov-or}. Thus we get a link
invariant denoted by $\t^{fr,St}_{\sL_N}$.\label{q-sl_N-fr}

\item Show the $\t^{fr,St}_{\sL_N}$ satisfies the following
skein relation
$$q^{\frac1{2N}}\t^{fr,St}_{\sL_N}(\lrints)\ -\
  q^{-\frac1{2N}}\t^{fr,St}_{\sL_N}(\rlints)\quad = \quad
  (q^{1/2}-q^{-1/2})\t^{fr,St}_{\sL_N}(\twoup)
$$
and the following framing and initial conditions
$$\begin{array}{rcl}\displaystyle
\t^{fr,St}_{\sL_N}(\rb{-4.2mm}{\ig[width=10mm]{plkink.eps}})&=&
 \displaystyle  q^{\frac{N-1/N}{2}}
  \t^{fr,St}_{\sL_N}(\rb{-4.2mm}{\ig[width=10mm]{oneup.eps}})
      \vspace{10pt}\\ \displaystyle
\t^{fr,St}_{\sL_N}(\unkn) &=& \displaystyle
   \frac{q^{N/2} - q^{-N/2}}{q^{1/2}-q^{-1/2}}\,.
\end{array}
$$

\item
The quadratic Casimir number for the standard $\sL_N$ representation
is equal to $N-1/N$. Therefore, the deframing of this invariant
gives $\displaystyle\ \t^{St}_{\sL_N} := q^{-\frac{N-1/N}{2}\cdot w}
\t^{fr,St}_{\sL_N}$\label{q-sl_N-unfr} which satisfies
$$q^{N/2}\t^{St}_{\sL_N}(\lrints)\ -\
  q^{-N/2}\t^{St}_{\sL_N}(\rlints)\quad = \quad
  (q^{1/2}-q^{-1/2})\t^{St}_{\sL_N}(\twoup)\,;
$$
$$\t^{St}_{\sL_N}(\unkn) =
  \frac{q^{N/2} - q^{-N/2}}{q^{1/2}-q^{-1/2}}\,.
$$
Check that this invariant is essentially a specialization of the HOMFLY polynomial,
$$\t^{St}_{\sL_N}(L) = \frac{q^{N/2} - q^{-N/2}}{q^{1/2}-q^{-1/2}}
  \Biggl.P(L)\Biggr|_{\rb{10pt}{$\begin{array}{ll}
             \scriptstyle a=q^{N/2}\ \vspace{-5pt}\\
             \scriptstyle z=q^{1/2}-q^{-1/2}\end{array}$}}\ .$$
Prove that the set of invariants $\{\t^{St}_{\sL_N}\}$ for all
values of $N$ is equivalent to the HOMFLY polynomial. Thus
$\{\t^{fr,St}_{\sL_N}\}$ may be considered as a framed version of
the HOMFLY polynomial.
\end{enumerate}

\item A different framed version of the HOMFLY polynomial is defined in
\cite[page 54]{Ka7}: $P^{fr}(L):= a^{w(L)}P(L)$.
\index{HOMFLY polynomial!framed}\label{fr-HOMFLY}
Show that $P^{fr}$ satisfies the following skein relation
$$P^{fr}(\lrints)\ -\ P^{fr}(\rlints)\quad = \quad zP^{fr}(\twoup)
$$
and the following framing and initial conditions
$$\begin{array}{rcl@{\qquad}rcl}\displaystyle
  P^{fr}(\rb{-4.2mm}{\ig[width=10mm]{plkink.eps}})&=&
 \displaystyle  a P^{fr}(\rb{-4.2mm}{\ig[width=10mm]{oneup.eps}})
      \ , &
  P^{fr}(\rb{-4.2mm}{\ig[width=10mm]{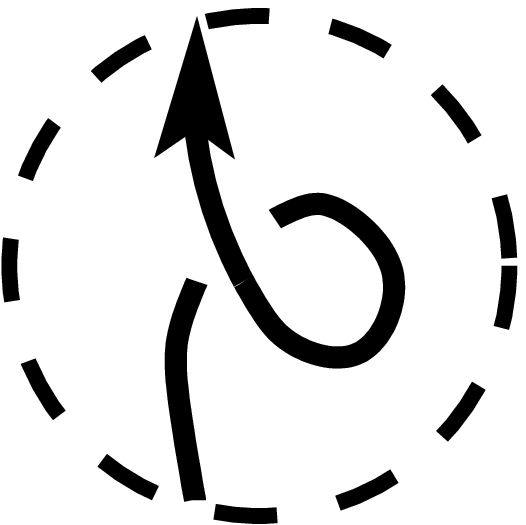}})&=&
 \displaystyle  a^{-1} P^{fr}(\rb{-4.2mm}{\ig[width=10mm]{oneup.eps}})
      \vspace{10pt}\\ \displaystyle
 P^{fr}(\unkn) &=& 1\,. &&&
\end{array}
$$

\end{enumerate}

\end{xcb}
 %2 knot invariants
\chapter{Finite type invariants} %03
\label{FT_inv}
\newcommand\risA[4]{\mbox{\begin{picture}(#2,#3)(0,0)
                   \put(0,#4){\ig[width=#2pt]{#1.eps}}
                \end{picture}}}
\newcommand\kn[2]{\rb{#2mm}{\ig[width=9.5mm]{#1}}}

In this chapter we introduce the main protagonist of this book: the
{\em finite type}, or {\em Vassiliev} knot invariants.

First we define the Vassiliev skein relation and extend, with its
help, arbitrary knot invariants to knots with double points. A
Vassiliev invariant of order at most $n$ is then defined as a knot
invariant which vanishes identically on knots with more than $n$
double points.

After that, we introduce a combinatorial object of great importance:
the chord diagrams. Chord diagrams serve as a means to describe the
symbols (highest parts) of the Vassiliev invariants.

Then we prove that classical invariant polynomials are all, in a
sense, of finite type, explain a simple method of calculating the
values of Vassiliev invariants on any given knot, and give a table
of basis Vassiliev invariants up to degree 5.

Finally, we show how Vassiliev invariants can be defined for framed
knots and for arbitrary tangles.

\section{Definition of Vassiliev invariants}

\subsection{} The original definition of finite type knot invariants was
just an application of the general machinery developed by
V.Vassiliev to study {\em complements of discriminants} in spaces of
maps.

The discriminants in question are subspaces of maps with
singularities of some kind. In particular, consider the space of all
smooth maps of the circle into $\R^3$. Inside this space, define the
discriminant as the subspace formed by maps that fail to be
embeddings, such as curves with self-intersections, cusps etc. Then
the complement of this discriminant can be considered as {\em the
space of knots}. The connected components of the space of knots are
precisely the isotopy classes of knots; knot invariants are locally
constant functions on the space of knots.

Vassiliev's machinery produces a spectral sequence that may (or may
not, nobody knows it yet) converge to the cohomology of the space of
knots. The zero-dimensional classes produced by this spectral
sequence correspond to knot invariants which are now known as
Vassiliev invariants.

This approach is indispensable if one wants to understand the higher
cohomology of the space of knots. However, if we are only after the
zero-dimensional classes, that is, knot invariants, the definitions
can be greatly simplified. In this chapter we follow the easy path
that requires no knowledge of algebraic topology whatsoever. For the
reader who is not intimidated by spectral sequences we outline
Vassiliev's construction in Chapter~\ref{chap_VSS}.

\subsection{Singular knots and the Vassiliev skein relation}

A {\em singular knot} is a smooth map $S^1\to\R^3$ that fails to be
an embedding. We shall only consider singular knots with the
simplest singularities, namely transversal self-intersections, or
{\em double points}.

\begin{xdefinition}
Let $f$ be a map of a one-dimensional manifold to $\R^3$. A point
$p\in\im(f)\subset\R^3$ is a {\em double point} \index{Double point}
of $f$ if $f^{-1}(p)$ consists of two points $t_1$ and $t_2$ and the
two tangent vectors $f'(t_1)$ and $f'(t_2)$ are linearly
independent. Geometrically, this means that in a neighbourhood of
the point $p$ the curve $f$ has two branches with non-collinear
tangents.
\end{xdefinition}

\begin{center}
\begin{tabular}{c}
\ig[width=25mm]{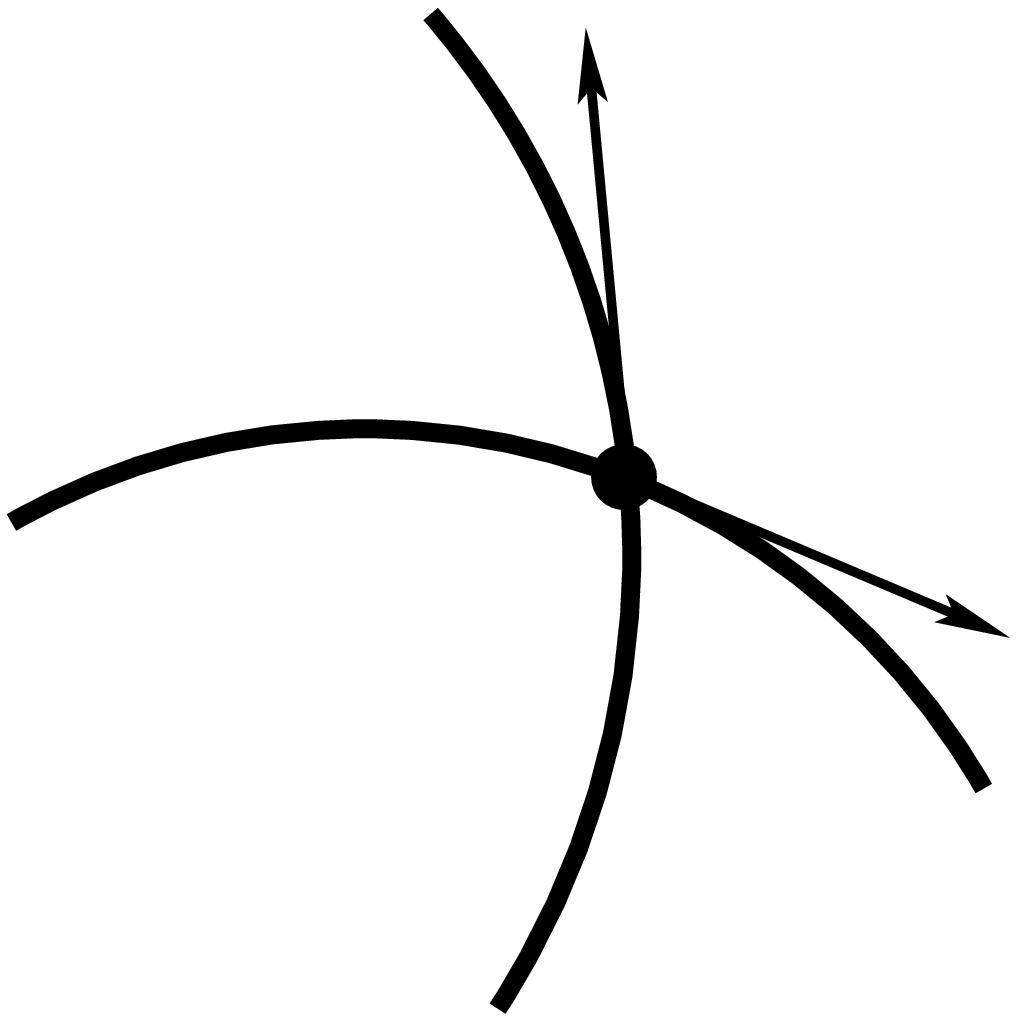} \\
A double point
\end{tabular}
\end{center}

\begin{xremark}
In fact, we gave a definition of a {\em simple double point}. We
omit the word ``simple'' since these are the only double points we
shall see.
\end{xremark}

Any knot invariant can be extended to knots with double points by
means of the {\em Vassiliev skein relation}:
\index{Vassiliev!extension}
\index{Vassiliev!skein relation}
\begin{equation} \label{Vskein}
  v(\double) = v(\lrints)-v(\rlints).
\end{equation}
Here $v$ is the knot invariant with values in some abelian group,
the left-hand side is the value of $v$ on a singular knot $K$ (shown
in a neighbourhood of a double point) and the right-hand side is the
difference of the values of $v$ on (possibly singular) knots
obtained from $K$ by replacing the double point with a positive and
a negative crossing respectively. The process of applying the skein
relation is also referred to as {\em resolving a double point}. It
is clearly independent of the plane projection of the singular knot.

Using the Vassiliev skein relation recursively, we can extend any
knot invariant to knots with an arbitrary number of double points.
There are many ways to do this, since we can choose to resolve
double points in an arbitrary order. However, the result is
independent of any choice. Indeed, the calculation of the value of
$v$ on a singular knot $K$ with $n$ double points is in all cases
reduced to the {\em complete resolution}%
\label{resolution}\index{Resolution!complete} of the knot $K$ which
yields an alternating sum
\begin{equation}
  v(K)=\sum_{\e_1=\pm1, \ldots ,\e_n=\pm1}(-1)^{|\e|}v(K_{\e_1,...,\e_n}),
\label{alt_sum}
\end{equation}
where $|\e|$ is the number of $-1$'s in the sequence $\e_1,\ldots,
\e_n$, and $K_{\e_1,\ldots,\e_n}$ is the knot obtained from $K$ by a
positive or negative resolution of the double points according to
the sign of $\e_i$ for the point number $i$.

\label{def_VI}
\begin{definition} \textrm{(}V.~Vassiliev \cite{Va1}\textrm{)}.
\index{Vassiliev!invariant}\index{Finite type invariant}%
\index{Knot invariant!Vassiliev}\index{Knot invariant!finite type}%
A knot invariant is said to be a {\em Vassiliev invariant} (or a
{\em finite type invariant}) of order (or degree) $\le n$ if its
extension vanishes on all singular knots with more than $n$ double
points. A Vassiliev invariant is said to be of order (degree) $n$ if
it is of order $\le n$ but not of order $\le n-1$.
\end{definition}

In general, a Vassiliev invariant may take values in an arbitrary
abelian group. In practice, however, all our invariants will take
values in commutative rings and it will be convenient to make this
assumption from now on.

\noindent{\bf Notation.} We shall denote by $\V_n$ \label{d:Vn} the
set of Vassiliev invariants of order $\le n$ with values in a ring
$\Ring$ \label{Ring}. Whenever necessary, we shall write
$\V_n^\Ring$ to indicate the range of the invariants explicitly. It
follows from the definition that, for each $n$, the set $\V_n$ is an
$\Ring$-module. Moreover, $\V_n\subseteq\V_{n+1}$, so we have an increasing
filtration
$$\V_0\subseteq\V_1\subseteq\V_2\subseteq \dots \subseteq\V_n
\subseteq\dots\subseteq\V:=\bigcup_{n=0}^\infty\V_n\ .
\label{d:Vspace}
$$

We shall further discuss this definition in the next section. First,
let us see that there are indeed many (in fact, infinitely many)
independent Vassiliev invariants.

\begin{example}\label{conw_as_vass} (\cite{BN0}).
The $n$th coefficient of the Conway polynomial is a
Vassiliev invariant of order $\le n$.
\smallskip

Indeed, the definition of the Conway polynomial, together with the
Vassiliev skein relation, implies that
$$\CP(\double)=t\CP(\twoup).$$ Applying this relation several times,
we get
$$
  \CP(\double\dots\double)=t^k\CP(\twoup\dots\twoup)
$$
for a singular knot with $k$ double points. If $k\ge n+1$, then the
coefficient at $t^n$ in this polynomial is zero.
\end{example}

\section{Algebra of Vassiliev invariants}\label{sec:alg_vi}

\subsection{The singular knot filtration}
\label{sing_knot_filtr}

Consider the ``tautological knot invariant'' $\K\to\Z\K$ which sends
a knot to itself. Applying the Vassiliev skein relation, we extend
it to knots with double points; a knot with $n$ double points is
then sent to an alternating sum of $2^n$ genuine knots.

Recall that we denote by $\Z\K$ the free abelian group spanned by
the equivalence classes of knots with multiplication induced by the
connected sum of knots. Let $\K_n$ be the $\Z$-submodule of the
algebra $\Z\K$ spanned by the images of knots with $n$ double
points.

\noindent
\textbf{Exercise.} Prove that $\K_n$ is an ideal of $\Z\K$.

A knot with $n+1$ double points gives rise to a difference of two
knots with $n$ double points in $\Z\K$; hence, we have the
descending {\em singular knot filtration}
$$\Z\K=\K_0\supseteq\K_1\supseteq\ldots\supseteq\K_n\supseteq\ldots$$
The definition of Vassiliev invariants can now be re-stated in the
following terms:
\begin{xdefinition}
Let $\Ring$ be a commutative ring. A Vassiliev invariant of order
$\le n$ is a linear function $\Z\K\to\Ring$ which vanishes on
$\K_{n+1}$.
\end{xdefinition}

According to this definition, the module of $\Ring$-valued Vassiliev
invariants of order $\le n$ is naturally isomorphic to the space of
linear functions $\Z\K/\K_{n+1}\to\Ring$. So, in a certain sense, the
study of the Vassiliev invariants is equivalent to studying the
filtration $\K_n$. In the next several chapters we shall mostly
speak about invariants, rather than the filtration on the algebra of
knots. Nevertheless, the latter approach, developed by Goussarov
\cite{G2} is important and we cannot skip it here altogether.
\begin{xdefinition}
Two knots $K_1$ and $K_2$ are {\em n-equivalent}
\index{n-equivalence} if they cannot be distinguished by Vassiliev
invariants  of degree $n$
and smaller with values in an arbitrary abelian group. 
A knot that is $n$-equivalent to the trivial knot is
called {\em $n$-trivial}.\index{n-triviality}
\end{xdefinition}
In other words, $K_1$ and $K_2$ are $n$-equivalent if and only if
$K_1-K_2\in\K_{n+1}$.
\begin{xdefinition}
Let $\G_n\K$ be the set of $(n-1)$-trivial knots. The {\em Goussarov
filtration}\index{Goussarov!filtration} on $\K$ is the descending filtration
$$\K=\G_1\K\supseteq\G_2\K\supseteq\ldots\supseteq\G_n\K\supseteq\ldots$$
\end{xdefinition}
The sets $\G_n\K$ are, in fact, abelian monoids under the connected
sum of knots (this follows from the fact that each $\K_n$ is a
subalgebra of $\Z\K$). Goussarov proved that the monoid quotient
$\K/\Gamma_n\K$ is an (abelian) group. We shall consider
$n$-equivalence in greater detail in Chapters~\ref{chapBr} and
\ref{chapMisc}.

\subsection{Vassiliev invariants as polynomials}\label{ss:vas-pol}
A useful way to think of Vassiliev invariants is as follows. Let $v$
be an invariant of singular knots with $n$ double points and
$\nabla(v)$ \label{nabla} be the extension of $v$  to singular knots
with $n+1$ double points using the Vassiliev skein relation. We can
consider $\nabla$ as an operator between the corresponding spaces of
invariants. Now, a function $v:\K\to\Ring$ is a Vassiliev invariant
of degree $\le n$, if it satisfies the difference equation
$\nabla^{n+1}(v)=0$. This can be seen as an analogy between
Vassiliev invariants as a subspace of all knot invariants and
polynomials as a subspace of all smooth functions on a real line:
the role of differentiation is played by the operator $\nabla$. It
is well known that continuous functions on a real line can be
approximated by polynomials. The main open problem of the theory of
finite type invariants is to find an analogue of this statement in
the knot-theoretic context, namely, to understand to what extent an
arbitrary numerical knot invariant can be approximated by Vassiliev
invariants. More on this in Section~\ref{ss:approx}.

\subsection{The filtration on the algebra of Vassiliev invariants}
\label{sec:grading_on_alg_vi}
\index{Algebra!of Vassiliev invariants}
\index{Vassiliev!invariant!algebra of}

The set of all Vassiliev invariants forms a commutative filtered
algebra with respect to the usual (pointwise) multiplication of
functions.\index{Multiplication!of Vassiliev invariants}

\begin{xtheorem}
The product of two Vassiliev invariants of degrees $\le p$ and $\le
q$ is a Vassiliev invariant of degree $\le p+q$.
\end{xtheorem}

\begin{proof}

Let $f$ and $g$ be two invariants with values in a ring $\Ring$, of
degrees $p$ and $q$ respectively. Consider a singular knot $K$ with
$n=p+q+1$ double points. The complete resolution of $K$ via the
Vassiliev skein relation gives
$$
(fg)(K)=\sum_{\e_1=\pm1, \ldots
,\e_n=\pm1}(-1)^{|\e|}f(K_{\e_1,...,\e_n})g(K_{\e_1,...,\e_n})
$$
in the notation of (\ref{alt_sum}). The alternating sum on the
right-hand side is taken over all points of an $n$-dimensional
binary cube
$$Q_n=\{(\e_1,\ldots,\e_n)\, |\, \e_i=\pm 1 \}.$$
In general, given a function $v$ on $Q_n$ and a subset $S\subseteq
Q_n$, the {\em alternating sum} of $v$ over $S$ is defined as
$\sum_{\e\in S} (-1)^{|\e|}v(\e)$.

If we set
$$f({\e_1,...,\e_n})=f(K_{\e_1,...,\e_n})$$ and define
$g({\e_1,...,\e_n})$ similarly, we can think of $f$ and $g$ as
functions on $Q_n$. The fact that $f$ is of degree $p$  means that
the alternating sum of $f$ on each $(p+1)$-face of $Q_n$ is zero.
Similarly, on each $(q+1)$-face of $Q_n$ the alternating sum of $g$
vanishes. Now, the theorem is a consequence of the following lemma.

\begin{xlemma}
Let $f,g$ be functions on $Q_{n}$, where $n=p+q+1$. If the
alternating sums of $f$ over any $(p+1)$-face, and of $g$ over any
$(q+1)$-face of $Q_{n}$ are zero, so is the alternating sum of the
product $fg$ over the entire cube $Q_{n}$.
\end{xlemma}
\noindent\textbf{Proof of the lemma}. Use induction on $n$. For
$n=1$ we have $p=q=0$ and the premises of the lemma read
$f(-1)=f(1)$ and $g(-1)=g(1)$. Therefore, $(fg)(-1)=(fg)(1)$, as
required.

For the general case, denote by $\mathcal{F}_n$ the space of
functions $Q_n\to\Ring$. We have two operators
$$\rho_-,\rho_+:\mathcal{F}_n\to\mathcal{F}_{n-1}$$
which take a function $v$
to its restrictions to the $(n-1)$-dimensional faces $\e_1=-1$ and
$\e_1=1$ of $Q_n$:
$$\rho_-(v)(\e_2,\dots,\e_n)=v(-1,\e_2,\dots,\e_n)$$ and
$$\rho_+(v)(\e_2,\dots,\e_n)=v(1,\e_2,\dots,\e_n).$$ Let
$$\delta=\rho_+-\rho_-.$$
Observe that if the alternating sum of $v$ over any $r$-face of
$Q_n$ is zero, then the alternating sum of $\rho_\pm(v)$
(respectively, $\delta(v)$) over any $r$-face (respectively,
$(r-1)$-face) of $Q_{n-1}$ is zero.

A direct check shows that the operator $\delta$ satisfies the
following Leibniz rule:
$$\delta(fg)=\rho_+(f)\cdot \delta(g) + \delta(f)\cdot \rho_-(g).
$$
Applying the induction assumption to each of the two summands on the
right-hand side, we see that the alternating sum of $\delta(fg)$
over the cube $Q_{n-1}$ vanishes. By the definition of $\delta$,
this sum coincides with the alternating sum of $fg$ over $Q_n$.
\end{proof}

\begin{xremark}
The existence of the filtration on the algebra of Vassiliev invariants
can be thought of as a manifestation of their polynomial character.
Indeed, a polynomial of degree $\le n$ in one variable can be
defined as a function whose $n+1$st derivative is identically zero.
Then the fact that a product of polynomials of degrees $\le p$ and
$\le q$ has degree $\le p+q$ can be proved by induction using the
Leibniz formula. In our argument on Vassiliev invariants we have
used the very same logic. A further discussion of the Leibniz
formula for finite type invariants can be found in \cite{Wil4}.
\end{xremark}

\subsection{Approximation by Vassiliev invariants}\label{ss:approx}

The analogy between finite type invariants and polynomials would be
even more satisfying if there existed a Stone-Weierstra\ss\ type
theorem for knot invariants that would affirm that any invariant can
be approximated by Vassiliev invariants. At the moment no such
statement is known. In fact, understanding the strength of the class
of finite type invariants is the main problem in the theory.

There are various ways of formulating this problem as a precise
question. Let us say that a class $\mathcal{U}$ of knot invariants
is {\em complete} if for any finite set of knots the invariants from
$\mathcal{U}$ span the space of all functions on these knots. We say
that invariants from $\mathcal{U}$ {\em distinguish knots} if for
any two different knots $K_1$ and $K_2$ there exists
$f\in\mathcal{U}$ such that $f(K_1)\neq f(K_2)$. Finally, the class
$\mathcal{U}$ {\em detects the unknot} if any knot can be
distinguished from the trivial knot by an invariant from
$\mathcal{U}$. A priori, completeness is the strongest of these
properties. In this terminology, the main outstanding problem in the
theory of finite type invariants is to determine whether the
Vassiliev invariants distinguish knots. While it is conjectured that
the set of rational-valued Vassiliev knot invariants is complete, it
is not even known if the class of all Vassiliev knot invariants
detects the unknot.

Note that the rational-valued Vassiliev invariants are complete if
and only if the intersection $\cap\K_n$ of all the terms of the
singular knot filtration is zero.
Indeed, a non-zero element of $\cap\K_n$ produces a universal relation
among the values of the invariants on a certain set of knots. On the
other hand, let $\cap\K_n=0$. Then the map $\K\to\Z\K/\K_{n+1}$ is a
Vassiliev invariant of order $n$ whose values on any given set of knots
become linearly independent as $n$ grows.
As for the Goussarov filtration $\G_n\K$, the intersection of all of
its terms consists of the trivial knot if and only if the Vassiliev
invariants detect the unknot.

There are knot invariants, of which we shall see many examples,
which are not of finite type, but, nevertheless, can be approximated
by Vassiliev invariants in a certain sense.
These are the {\em polynomial} and the {\em power series Vassiliev
invariants}.\index{Vassiliev!invariant!power series} A polynomial
Vassiliev invariant is an element of the vector
space\label{pol-Vas-inv}
$$\V_{\bullet}=\bigoplus_{n=0}^{\infty} \V_n.$$
Since the product of two invariants of degrees $m$ and $n$ has
degree at most $m+n$, the space $\V_{\bullet}$ is, in fact, a
commutative graded algebra. The power series Vassiliev invariants
are, by definition, the elements of its {\em graded completion}
$\widehat{\V}_{\bullet}$\label{V-grad_compl}
(see Appendix~\ref{filtration}, page~\pageref{gr-compl}).

The Conway polynomial $C$ is an example of a power series
invariant. Observe that even though for any knot $K$ the value
$C(K)$ is a polynomial, the Conway polynomial $C$ is not a
polynomial invariant according to the definition of this paragraph.

Power series Vassiliev invariants are just one possible approach to
defining approximation by finite type invariants. A wider class of
invariants are those {\em dominated by Vassiliev invariants}. We say
that a knot invariant $u$ is dominated by Vassiliev invariants if
$u(K_1)\neq u(K_2)$ for some knots $K_1$ and $K_2$ implies that
there is a Vassiliev knot invariant $v$ with $v(K_1)\neq v(K_2)$.
Clearly, if Vassiliev invariants distinguish knots, then each knot
invariant is dominated by Vassiliev invariants. At the moment,
however, it is an open question whether, for instance, the signature
of a knot \cite{Rol} is dominated by Vassiliev invariants.

\section{Vassiliev invariants of degrees 0, 1 and 2}

\begin{proposition}\label{prop:v0}
$\V_0 = \{\const\}$, $\dim \V_0 = 1$.
\end{proposition}

\begin{proof}
Let $f\in\V_0$. By definition, the value of (the extension of) $f$
on any singular knot with one double point is 0. Pick an arbitrary
knot $K$. Any diagram of $K$ can be turned into a diagram of the
trivial knot $K_0$ by crossing changes done one at a time. By assumption, the jump of $f$ at every crossing
change is 0, therefore, $f(K)=f(K_0)$. Thus $f$ is constant.
\end{proof}

\begin{proposition}\label{V1}\label{prop:v1}
$\V_1 = \V_0$. 
\end{proposition}

\begin{proof}
A singular knot with one double point is divided by the double point
into two closed curves. An argument similar to the last proof shows
that the value of $v$ on any knot with one double point is equal to
its value on the ``figure infinity'' singular knot and, hence, to 0:
\begin{equation}\label{v1_eq_v0}
   v(\double)=v(\eight)=0
\end{equation}
Therefore, $\V_1=\V_0$.
\end{proof}

The first non-trivial Vassiliev invariant appears in degree 2: it is
the second coefficient $c_2$ of the Conway polynomial, also known as
the {\em Casson invariant}.

\begin{proposition}\label{prop:v2}
$\dim \V_2 = 2$.
\end{proposition}

\begin{proof}
Let us explain why the argument of the proof of
Propositions~\ref{prop:v0} and \ref{prop:v1} does not work in this
case. Take a knot with two double points and try to
transform it into some fixed knot with two double points using
smooth deformations and crossing changes. It is easy to see that any
knot with two double points can be reduced to one of the following
two basic knots:
\begin{center}
\begin{tabular}{ccc}
\risS{-10}{kn21}{}{35}{0}{10} & \qquad & \risS{-10}{kn22}{}{35}{0}{0} \\
Basic knot $K_1$ & \qquad & Basic knot $K_2$
\end{tabular}\label{fig_k1_k2}
\end{center}
--- but these two knots cannot be obtained one from the other!
The essential difference between them is in the
order of the double points on the curve.

Let us label the double points of $K_1$ and $K_2$, say, by 1 and 2.
When travelling along the first knot, $K_1$, the two double points
are encountered in the order 1122 (or 1221, 2211, 2112 if you start
from a different initial point). For the knot $K_2$ the sequence is
1212 (equivalent to 2121). The two sequences 1122 and 1212 are
different even if cyclic permutations are allowed.

Now take an arbitrary singular knot $K$ with two double points. If
the cyclic order of these points is 1122, then we can transform the
knot to $K_1$, passing in the process of deformation through some
singular knots with three double points; if the order is 1212, we
can reduce $K$ in the same way to the second basic knot $K_2$.

The above argument shows that, to any $\Ring$-valued order 2
Vassiliev invariant there corresponds a function on the set of two
elements $\{K_1,K_2\}$ with values in $\Ring$. We thus obtain a
linear map $\V_2\to\Ring^2$. The kernel of this map is equal to $\V_1$:
indeed, the fact that a given invariant $f\in\V_2$ satisfies
$f(K_1)=f(K_2)=0$ means that it vanishes on {\em any} singular knot
with 2 double points, which is by definition equivalent to saying
that $f\in\V_1$.

On the other hand, the image of this linear map is no more than
one-dimensional, since for {\em any} knot invariant $f$ we have
$f(K_1)=0$. This proves that $\dim \V_2 \le 2$. In fact,
$\dim\V_2=2$, since the second coefficient $c_2$ of the Conway
polynomial is not constant (see Table \ref{conway_tabl}).
\end{proof}

\section{Chord diagrams}
\label{ch_diag}

Now let us give a formal definition of the combinatorial structure
which is implicit in the proof of Proposition~\ref{prop:v2}.

\begin{xdefinition}\index{Chord diagram}\index{Order}\index{Degree}
A {\em chord diagram} of order $n$ (or degree $n$) is an oriented
circle with a distinguished set of $n$ disjoint pairs of distinct
points, considered up to orientation preserving diffeomorphisms of
the circle. The set of all chord diagrams of order $n$ will be
denoted by $\ChD_n$.\label{ChD_n}
\end{xdefinition}

We shall usually omit the orientation of the circle in pictures of
chord diagrams, assuming that it is oriented counterclockwise.
\medskip

\noindent
{\bf Examples.}
\begin{eqnarray*}
\ChD_1&=&\{\cdO\},\\
\ChD_2&=&\{\cdWO,\cdWW\},\\
\ChD_3&=&\{\cdTO,\cdTW,\cdTT,\cdTF,\cdTV\}.
\end{eqnarray*}

\begin{xremark}
Chord diagrams that differ by a mirror reflection are, in general,
different:
$$
  \rb{-7mm}{\ig[height=15mm]{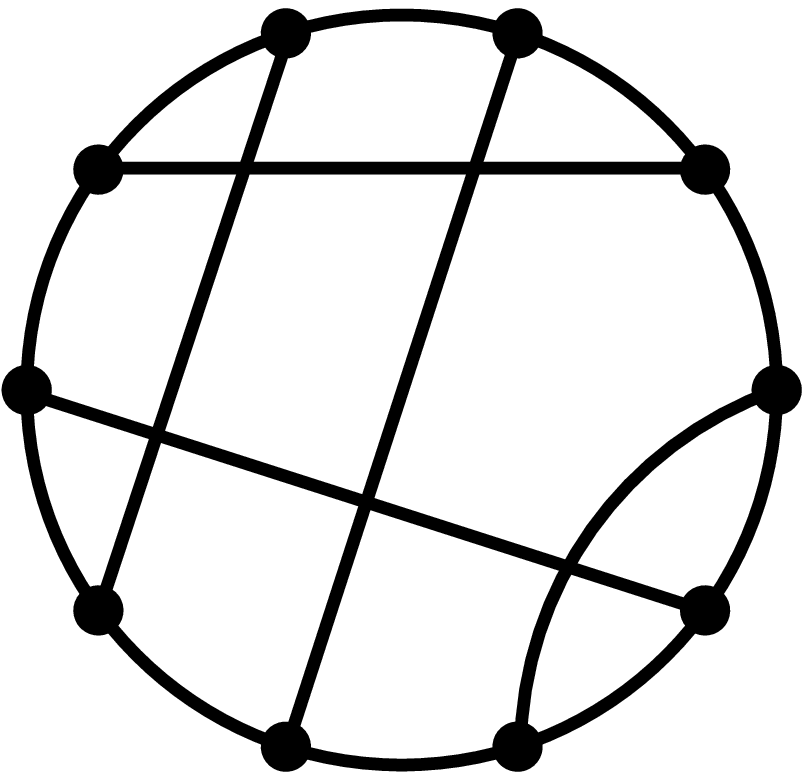}}
  \quad\ne\quad
  \rb{-7mm}{\ig[height=15mm]{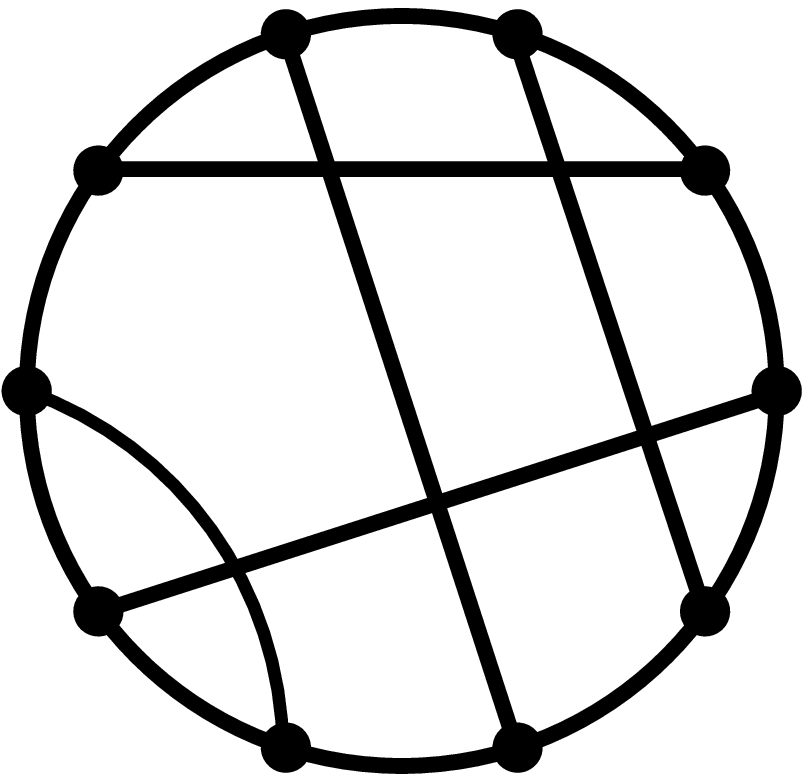}}
$$
\end{xremark}
This observation reflects the fact that we are studying
\textit{oriented} knots.

\subsection{The chord diagram of a singular knot}

Chord diagrams are used to code certain information about singular
knots.

\begin{xdefinition}\index{Chord diagram!of a singular knot}
The  chord diagram $\sigma (K)\in \ChD_n$ of a singular knot with
$n$ double points is obtained by marking on the parametrizing
circle $n$ pairs of points whose images are the $n$ double points of
the knot.
\end{xdefinition}

\noindent
{\bf Examples.}\vspace*{-10pt}
$$\sigma\Bigl(\ \risS{-10}{kn21}{}{30}{0}{20}\ \Bigr)\ =\ \cdWO\ ,\qquad
\sigma\Bigl(\ \risS{-10}{kn22}{}{30}{0}{0}\ \Bigr)\ =\ \cdWW\ .
$$

\begin{proposition} {\rm (V.~Vassiliev \cite{Va1})}.
\label{fact_thru} The value of a Vassiliev invariant $v$ of order
$\le n$ on a knot $K$ with $n$ double points depends only on the
chord diagram of $K$:
$$
  \sigma (K_1)=\sigma (K_2) \Rightarrow v(K_1)=v(K_2).
$$
\end{proposition}

\begin{proof} Suppose that $\sigma (K_1)=\sigma (K_2)$.
Then there is a one-to-one correspondence between the chords of both
chord diagrams, and, hence, between the double points of $K_1$ and
$K_2$. Place $K_1,K_2$ in $\R^3$ so that the corresponding double
points coincide together with both branches of the knot in the
vicinity of each double point.

\begin{center}
\begin{tabular}{c@{\qquad\qquad}c}
\risS{-10}{kn22}{}{35}{18}{15} & \risS{-8}{kn22a}{}{40}{0}{0} \\
{\tt Knot} $K_1$ & {\tt Knot} $K_2$
\end{tabular}
\end{center}

Now we can deform $K_1$ into $K_2$ in such a way that some small
neighbourhoods of the double points do not move. We can assume that
the only new singularities created in the process of this
deformation are a finite number of double points, all at distinct
values of the deformation parameter.
By the Vassiliev skein relation, in each of these events the value
of $v$ does not change, and this implies that
 $v(K_1)=v(K_2)$.
\end{proof}

Proposition~\ref{fact_thru}
shows that there is a well defined map $\alpha_n:\V_n \to
\Ring\ChD_n$ \label{alpha_n}\label{Ring(ChD_n)} (the $\Ring$-module
of $\Ring$-valued functions on the set $\ChD_n$):
$$
  \alpha_n(v)(D)=v(K),
$$
where $K$ is an arbitrary knot with $\sigma (K)=D$.

We want to understand the size and the structure of the space
$\V_n$, so it would be of use to have a description of the kernel
and the image of $\alpha_n$.

The description of the kernel follows immediately from the
definitions: $\ker\alpha_n=\V_{n-1}$. Therefore, we obtain an
injective homomorphism
\begin{equation}\label{alfa_n}
  \ol{\alpha}_n: \V_n/\V_{n-1} \to \Ring \ChD_n.
\end{equation}
The problem of describing the image of $\alpha_n$ is much more
difficult. The answer to it will be given in Theorem~\ref{fund_thm} on
page~\pageref{fund_thm}.

Since there is only a finite number of diagrams of each order,
Proposition~\ref{fact_thru} implies the following
\begin{corollary}\label{fingen}
The module of $\Ring$-valued Vassiliev invariants of degree at most
$n$ is finitely generated over $\Ring$.
\end{corollary}

Since the map $\alpha_n$ discards the order $(n-1)$ part of a
Vassiliev invariant $v$, we can, by analogy with differential
operators, call the function $\alpha_n(v)$ on chord diagrams the
\textit{symbol} of the Vassiliev invariant $v$:
\index{Symbol!of a Vassiliev invariant}
$$
  \symb(v) = \alpha_n(v),
$$\label{Symbol}\index{Vassiliev!invariant!symbol of}
where $n$ is the order of $v$.

\begin{xexample}
The symbol of the Casson invariant is equal to 0 on
the chord diagram with two parallel chords, and to 1 on the chord
diagram with two interecting chords.
\end{xexample}

\begin{remark}\label{remark:cd} It may be instructive to state all the above in the
dual setting of the singular knot filtration. The argument in the
proof of Proposition~\ref{fact_thru} essentially
says that $\ChD_n$ is the set of singular knots with $n$ double
points modulo isotopies and crossing changes. In terms of the
singular knot filtration, we have shown that if two knots with $n$
double points have the same chord diagram, then their difference
lies in $\K_{n+1}\subset\Z\K$. Since $\K_n$ is spanned by the
complete resolutions of knots with $n$ double points, we have a
surjective map
$$
\Z\ChD_n\to\K_n/\K_{n+1}.
$$
The kernel of this map , after tensoring with the
rational numbers, is spanned by the so-called {\em 4T} and {\em 1T
relations}, defined in the next chapter. This is the content of
Theorem~\ref{fund_thm}.
\end{remark}

\section{Invariants of framed knots}

A {\em singular framing} on a closed curve immersed in $\R^3$ is a
smooth normal vector field with a finite number of
simple zeroes on this curve. A {\em singular framed knot} is a knot with simple
double points in $\R^3$ equipped with a singular framing whose set
of zeroes is disjoint from the set of double points.

Invariants of framed knots are extended to singular framed knots by
means of the Vassiliev skein relation; for double points it has the
same form as before, and for the zeroes of the singular framing it
can be drawn as
$$v\Bigl(\doublefr\Bigr) = v\Bigl(\lrintsfr\Bigr)-v\Bigl(\rlintsfr\Bigr).
$$
An invariant of framed knots is of order $\le n$ if
its extension vanishes on knots with more than $n$ singularities
(double points or zeroes of the framing).

Let us denote the space of invariants of order $\le n$ by
$\V_n^{fr}$.
\label{fr-vi}
\index{Vassiliev!invariant!framed}
There is a natural
inclusion $i:\V_n \to \V_n^{fr}$ defined by setting $i(f)(K)=f(K')$
where $K$ is a framed knot, and $K'$ is the same knot without
framing. It turns out that this is a proper inclusion for all
$n\ge1$.

Let us determine the framed Vassiliev invariants of small degree.
Any invariant of degree zero is, in fact, an unframed knot invariant
and, hence, is constant. Indeed, increasing the framing by one can
be thought of as passing a singularity of the framing, and this does
not change the value of a degree zero invariant.

\begin{xca}
(1) Prove that $\dim\V_1^{fr}=2$, and that $\V_1^{fr}$ is spanned
by the constants and the self-linking number.\\
(2) Find the dimension and a basis of the vector space $\V_2^{fr}$.
\end{xca}

\begin{xca}
Let $v$ be a framed Vassiliev invariant degree $n$, and $K$ an
unframed knot. Let $v(K,k)$ be the value of $v$ on $K$ equipped with
a framing with self-linking number $k$. Show that $v(K,k)$ is a
polynomial in $k$ of degree at most $n$.
\end{xca}

\subsection{Chord diagrams for framed knots}

We have seen that chord diagrams on $n$ chords can be thought of as
singular knots with $n$ double points modulo isotopies and crossing changes.
Following the same logic, we should define a chord diagram for framed knots
as an equivalence class of framed singular knots with $n$ singularities
modulo isotopies, crossing changes and additions of zeroes of the framing.
In this way, the value of a degree $n$ Vassiliev invariant on a singular
framed knot with $n$ singularities will only depend on the chord diagram
of the knot.

As a combinatorial object, a framed chord diagram of degree $n$ can be
defined as a usual chord diagram of degree $n-k$ together with $k$ dots
marked on the circle. The chords correspond to the double points of a
singular knot and the dots represent the zeroes of the framing.

In the sequel we shall not make any use of diagrams with dots, for
the following reason. If $\Ring$ is a ring where $2$ is invertible, a zero
of the framing on a knot with $n$ singularities can be replaced, modulo knots
with $n+1$ singularities, by ``half of a double point'':
$$v(\framedone)=\frac{1}{2} v(\frameddbl)-\frac{1}{2} v(\framedtwo)$$
for any invariant $v$.
In particular, if we replace a dot with a chord whose endpoints are
next to each other on some diagram,
the symbol of any  Vassiliev invariant on this diagram
is simply multiplied by 2.

On the other hand, the fact that we can use the same chord diagrams for both
framed and unframed knots does not imply that the corresponding theories of
Vassiliev invariants are the same. In particular, we shall see that the symbol
of any invariant of unframed knots vanishes on a diagram which has a
chord that has no intersections with other chords. This does not hold for an
arbitrary framed invariant.

\begin{xexample} The symbol of the self-linking number is the
function equal to 1 on the chord diagram with one chord.
\end{xexample}

\section{Classical knot polynomials as Vassiliev invariants}
\label{knpol_as_vi}

In Example \ref{conw_as_vass}, we have seen that the coefficients of
the Conway polynomial are Vassiliev invariants. The Conway
polynomial, taken as a whole, is not, of course, a finite type
invariant, but it is an infinite linear combination of such; in
other words, it is a power series Vassiliev invariant. This property
holds for all classical knot polynomials --- but only after a
suitable substitution.

\subsection{}
Modify the Jones polynomial of a knot $K$ substituting $t=e^h$ and
then expanding it into a formal power series in $h$. \index{Jones
polynomial!modified} Let $j_n(K)$ be the coefficient of $h^n$ in
this expansion.

\begin{xtheorem}[\cite{G1,BL,BN1}]\label{jones_vi}
The coefficient $j_n(K)$ is a Vassiliev invariant of order $\leq n$.
\end{xtheorem}

\begin{proof}
Plugging $t=e^h=1+h+\dots$ into the skein relation from
Section~\ref{jones_skein} we get
$$
  (1-h+\dots)\cdot J(\lrints) -
  (1+h+\dots)\cdot J(\rlints) =
  (h+\dots)\cdot J(\twoup)\ .
$$
We see that the difference $$J(\lrints) - J(\rlints)=J(\double)$$ is
congruent to $0$ modulo $h$. Therefore, the Jones polynomial of a
singular knot with $k$ double points is divisible by $h^k$. In
particular, for $k\geq n+1$ the coefficient of $h^n$ equals zero.
\end{proof}

Below we shall give an explicit description of the symbols of the
finite type invariants $j_n$; the similar description for the Conway
polynomial is left as an exercise (no.\ \ref{ex_consymb} at the end
of the chapter, page \pageref{ex_consymb}).

\subsection{Symbol of the Jones invariant $j_n(K)$}
\label{symb_j_n}\index{Weight system!of the Jones coefficients}
\index{Symbol!of the Jones coefficients}
To find the symbol of $j_n(K)$, we must compute the
coefficient of $h^n$ in the Jones polynomial $J(K_n)$ of a
singular knot $K_n$ with $n$ double points in terms of its chord
diagram $\sigma(K_n)$. Since
$$
  J(\double)\! = J(\lrints) -\! J(\rlints)\! =\!
   h\Bigl( j_0(\twoup) + j_0(\lrints) + j_0(\rlints)\Bigr) + \dots
$$
the contribution of a double point of $K_n$ to the coefficient
$j_n(K_n)$ is the sum of the values of $j_0(\cdot)$ on the three
links in the parentheses above. The values of $j_0(\cdot)$ for the
last two links are equal, since, according to Exercise~4 to this
chapter, to $j_0(L) = (-2)^{\mbox{\#(components of\ } L)-1}$. So it
does not depend on the specific way $L$ is knotted and linked and we
can freely change the under/over-crossings of $L$. On the level of
chord diagrams these two terms mean that we just forget about the
chord corresponding to this double point. The first term,
$j_0(\twoup)$, corresponds to the smoothing of the double point
according to the orientation of our knot (link). On the level of
chord diagrams this corresponds to the doubling of a chord:
$$\risS{-12}{symbjn1}{}{30}{20}{15}\quad
  \risS{-2}{totonew}{}{25}{20}{15}\quad
  \risS{-12}{symbjn2}{}{30}{20}{15}\ .
$$
This leads to the following procedure of computing the value of the
symbol of $j_n(D)$ on a chord diagram $D$. Introduce {\it a state}
$s$ for $D$ as an arbitrary function on the set chords of $D$ with
values in the set $\{1,2\}$. With each state $s$ we associate an
immersed plane curve obtained from $D$ by resolving (either doubling or
deleting) all its chords according to $s$:
$$\risS{-12}{symbjn1}{\put(14,17){\mbox{$\scriptstyle c$}}}{30}{20}{15}\
  \risS{-2}{totonew}{}{25}{20}{15}\
  \risS{-12}{symbjn2}{}{30}{20}{15}\ ,\mbox{\ if\ }s(c)=1;\qquad
 \risS{-12}{symbjn1}{\put(14,17){\mbox{$\scriptstyle c$}}}{30}{20}{15}\
  \risS{-2}{totonew}{}{25}{20}{15}\
 \risS{-12}{symbjn3}{}{30}{20}{15}\ ,\mbox{\ if\ }s(c)=2.
$$
Let $|s|$ denote the number of components of the curve obtained in
this way. Then\index{Jones polynomial!symbol}
$$\symb(j_n)(D) = \sum_s\ \Bigl(\prod_c s(c)\Bigr)\ (-2)^{|s|-1}\ ,$$
where the product is taken over all $n$ chords of $D$, and the sum
is taken over all $2^n$ states for $D$.

For example, to compute the value of the symbol of $j_3$ on the
chord diagram $\risS{-12}{symcd4}{}{30}{20}{15}$ we must consider
8 states:
\def\stsjn#1#2#3{\begin{array}{c}
         \risS{-12}{#1}{}{30}{20}{15}\\
         \scriptstyle \prod s(c) = #2 \\
         \scriptstyle |s| = #3
         \end{array}}
$$\stsjn{symst1}{1}{2}\qquad \stsjn{symst2}{2}{1}\qquad
  \stsjn{symst3}{2}{1}\qquad \stsjn{symst4}{2}{3}
$$
$$\stsjn{symst5}{4}{2}\qquad \stsjn{symst6}{4}{2}\qquad
  \stsjn{symst7}{4}{2}\qquad \stsjn{symst8}{8}{1}
$$
Therefore,
$$\symb(j_3)\bigl(\risS{-12}{symcd4}{}{30}{20}{15}\bigr) =
-2+2+2+2(-2)^2+4(-2)+4(-2)+4(-2)+8=-6
$$
Similarly one can compute the values of $\symb(j_3)$ on all chord
diagrams with three chords. Here is the result:
$$\begin{array}{c||c|c|c|c|c}
D& \risS{-12}{symcd1}{}{30}{20}{15} & \risS{-12}{symcd2}{}{30}{20}{15}
 & \risS{-12}{symcd3}{}{30}{20}{15} & \risS{-12}{symcd4}{}{30}{20}{15}
 & \risS{-12}{symcd5}{}{30}{20}{15} \\ \hline
\symb(j_3)(D) & 0 & 0 & 0 & -6 & -12
\end{array}
$$

This function on chord diagrams, as well as the whole Jones
polynomial, is closely related to the Lie algebra $\sL_2$ and its
standard 2-dimensional representation. We shall return to this
subject several times in the sequel (see Sections
\ref{ws_sl_2_on_A}, \ref{ws_sl_n_St}, etc).

\subsection{}
According to Exercise~\ref{Jones_mirror}
(page~\pageref{Jones_mirror}), for the mirror reflection
$\overline{K}$ of a knot $K$ the power series expansion of
$J(\overline{K})$ can be obtained from the series $J(K)$ by
substituting $-h$ for $h$. This means that
$j_{2k}(\overline{K})=j_{2k}(K)$ and
$j_{2k+1}(\overline{K})=-j_{2k+1}(K)$.

\subsection{}
Table \ref{taylor_jones} displays the first five terms of the
power series expansion of the Jones polynomial after the
substitution $t=e^h$.

\begin{longtable}{c|lllll@{\quad+\dots\makebox(0,13){}}}
$3_1 $&$ 1 $&$ -3h^2 $&$ +6h^3 $&$ -(29/4)h^4 $&$ +(13/2)h^5 $\\
$4_1 $&$ 1 $&$ +3h^2 $&$       $&$ +(5/4)h^4 $&$ $\\
$5_1 $&$ 1 $&$ -9h^2 $&$ +30h^3 $&$ -(243/4)h^4 $&$ +(185/2)h^5 $\\
$5_2 $&$ 1 $&$ -6h^2 $&$ +18h^3 $&$ -(65/2)h^4 $&$ +(87/2)h^5 $\\
$6_1 $&$ 1 $&$ +6h^2 $&$ -6h^3 $&$ +(17/2)h^4 $&$ -(13/2)h^5 $\\
$6_2 $&$ 1 $&$ +3h^2 $&$ -6h^3 $&$ +(41/4)h^4 $&$ -(25/2)h^5 $\\
$6_3 $&$ 1 $&$ -3h^2 $&$       $&$ -(17/4)h^4 $&$ $\\
$7_1 $&$ 1 $&$ -18h^2 $&$ +84h^3 $&$ -(477/2)h^4 $&$ +511h^5 $\\
$7_2 $&$ 1 $&$ -9h^2  $&$ +36h^3 $&$ -(351/4)h^4 $&$ +159h^5 $\\
$7_3 $&$ 1 $&$ -15h^2 $&$ -66h^3 $&$ -(697/4)h^4 $&$ -(683/2)h^5 $\\
$7_4 $&$ 1 $&$ -12h^2 $&$ -48h^3 $&$ -113h^4 $&$ -196h^5 $\\
$7_5 $&$ 1 $&$ -12h^2 $&$ +48h^3 $&$ -119h^4 $&$ +226h^5 $\\
$7_6 $&$ 1 $&$ -3h^2  $&$ +12h^3 $&$ -(89/4)h^4 $&$ +31h^5 $\\
$7_7 $&$ 1 $&$ +3h^2  $&$ +6h^3 $&$ +(17/4)h^4 $&$ +(13/2)h^5 $\\
$8_1 $&$ 1 $&$ +9h^2  $&$ -18h^3 $&$ +(135/4)h^4 $&$ -(87/2)h^5 $\\
$8_2 $&$ 1 $&$        $&$ -6h^3 $&$ +27h^4 $&$ -(133/2)h^5 $\\
$8_3 $&$ 1 $&$ +12h^2 $&$       $&$ +17h^4 $&$ $\\
$8_4 $&$ 1 $&$ +9h^2 $&$ -6h^3 $&$ +(63/4)h^4 $&$ -(25/2)h^5 $\\
$8_5 $&$ 1 $&$ +3h^2 $&$ +18h^3 $&$ +(209/4)h^4 $&$ +(207/2)h^5 $\\
$8_6 $&$ 1 $&$ +6h^2 $&$ -18h^3 $&$ +(77/2)h^4 $&$ -(123/2)h^5 $\\
$8_7 $&$ 1 $&$ -6h^2 $&$ -12h^3 $&$ -(47/2)h^4 $&$ -31h^5 $\\
$8_8 $&$ 1 $&$ -6h^2 $&$ -6h^3 $&$ -(29/2)h^4 $&$ -(25/2)h^5 $\\
$8_9 $&$ 1 $&$ +6h^2 $&$       $&$ +(23/2)h^4 $&$ $\\
$8_{10} $&$ 1 $&$ -9h^2 $&$ -18h^3 $&$ -(123/4)h^4 $&$ -(75/2)h^5 $\\
$8_{11} $&$ 1 $&$ +3h^2 $&$ -12h^3 $&$ +(125/4)h^4 $&$ -55h^5 $\\
$8_{12} $&$ 1 $&$ +9h^2 $&$       $&$ +(51/4)h^4 $&$ $\\
$8_{13} $&$ 1 $&$ -3h^2 $&$ -6h^3 $&$ -(53/4)h^4 $&$ -(25/2)h^5 $\\
$8_{14} $&$ 1 $&$       $&$       $&$ +6h^4 $&$ -18h^5 $\\
$8_{15} $&$ 1 $&$ -12h^2 $&$ +42h^3 $&$ -80h^4 $&$ +(187/2)h^5 $\\
$8_{16} $&$ 1 $&$ -3h^2 $&$ +6h^3 $&$ -(53/4)h^4 $&$ +(37/2)h^5 $\\
$8_{17} $&$ 1 $&$ +3h^2 $&$       $&$ +(29/4)h^4 $&$ $\\
$8_{18} $&$ 1 $&$ -3h^2 $&$       $&$ +(7/4)h^4 $&$ $\\
$8_{19} $&$ 1 $&$ -15h^2 $&$ -60h^3 $&$ -(565/4)h^4 $&$ -245h^5 $\\
$8_{20} $&$ 1 $&$ -6h^2 $&$ +12h^3 $&$ -(35/2)h^4 $&$ +19h^5 $\\
$8_{21} $&$ 1 $&$       $&$ -6h^3 $&$ +21h^4 $&$ -(85/2)h^5 $\vspace{10pt}\\
\caption{Taylor expansion of the modified Jones polynomial}
\label{taylor_jones}\index{Jones polynomial!modified!table}
\index{Table of!Jones polynomials!modified}
\end{longtable}

\subsection{Example}
\label{ex_J_sing}
In the following examples the $h$-expansion of the Jones polynomial
starts with a power of $h$ equal to the number of double points in a
singular knot, in compliance with Theorem \ref{jones_vi}.
\def\risn#1#2#3{\rb{#3pt}{\ig[width=#2pt]{#1.eps}}}
\def\jn#1{J\!\left(\risn{#1}{30}{-10}\right)}
\def\jne#1{J\!\left(\risn{#1}{30}{-15}\right)}
$$\begin{array}{ccl}
\jn{rtfssp}&=&\underbrace{\jn{rtfsmr}}_{\begin{array}{c}
                     {\scriptstyle ||}\\0\end{array}}-
  \jn{rtfsl}=-\underbrace{\jn{rtfm}}_{\begin{array}{c}
                     {\scriptstyle ||}\\1\end{array}}+\jn{31}\\
&=& -1+J(3_1)=-3h^2 +6h^3 -\frac{29}{4}h^4 +\frac{13}{2}h^5+\dots
\end{array}$$

Similarly,
$$\jn{rtfss1}\ =\ J(\overline{3_1}) -1 =
 -3h^2 -6h^3 -\frac{29}{4}h^4 -\frac{13}{2}h^5+\dots
$$

Thus we have
$$\jn{rtfsss}\ =\ \jn{rtfss1}-\jn{rtfssp} =
  -12h^3 -13h^5+\dots
$$

\subsection{}
\label{qift} J.~Birman and X.-S.~Lin proved in \cite{Bir2,BL} that
all quantum invariants produce Vassiliev invariants in the same way
as the Jones polynomial. More precisely, let $\theta(K)$ be the
quantum invariant constructed as in Section~\ref{qi}. It is a
polynomial in $q$ and $q^{-1}$. Now let us make substitution
$q=e^h$ and consider the coefficient $\theta_n(K)$ of $h^n$ in the
Taylor expansion of $\theta(K)$.

\begin{xtheorem}[\cite{BL,BN1}]
\index{Theorem!Birman--Lin}
The coefficient $\theta_n(K)$ is a Vassiliev invariant of order $\leq n$.
\end{xtheorem}

\begin{proof}
The argument is similar to that used in Theorem \ref{jones_vi}: it
is based on the fact that an $R$-matrix $R$ and its inverse $R^{-1}$
are congruent modulo $h$.
\end{proof}

\subsection{The Casson invariant}\label{subsection:Casson}

The second coefficient of the Conway polynomial, or the {\em Casson
invariant}\index{Casson invariant}, can be computed directly from
any knot diagram by counting (with signs) pairs of crossings of
certain type\footnote{The Casson invariant was defined in 1985 by Casson as an invariant of homology 3-spheres. The Casson invariant of a knot can be interpreted as the difference between the Casson invariants of the homology spheres obtained by surgeries on the knot with different framings, see \cite{AM}.}.

Namely, fix a based Gauss diagram $G$ of a knot $K$, with an
arbitrary basepoint, and consider all pairs of arrows of $G$ that
form a subdiagram of the following form:\vspace{-10pt}
\begin{equation}\label{fig:PV}
\risS{-12}{ad2}{\put(-2,2){$\scriptstyle \e_1$}
                \put(21,2){$\scriptstyle \e_2$}}{25}{15}{15}
\end{equation}
The Casson invariant $a_2(K)$ is defined as the number of such pairs
of arrows with $\e_1\e_2=1$ minus the number of
pairs of this form with $\e_1\e_2=-1$.

\begin{xtheorem}\label{v2_th}
The Casson invariant coincides with the second coefficient of the
Conway polynomial $c_2$.
\end{xtheorem}

\begin{proof}
We shall prove that the Casson invariant as defined above, is a
Vassiliev invariant of degree 2. It can be checked directly that it
vanishes on the unknot and is equal to 1 on the left trefoil. Since
the same holds for the invariant $c_2$ and $\dim\V_2=2$, the
assertion of the theorem will follow.

First, let us verify that $a_2$ does not depend on the location of
the basepoint on the Gauss diagram. It is enough to prove that
whenever the basepoint is moved over the endpoint of one arrow, the
value of $a_2$ remains the same.

Let $c$ be an arrow of some Gauss diagram. For another arrow $c'$ of
the same Gauss diagram with the sign $\varepsilon(c')$, the {\em
flow} of $c'$ through $c$ is equal to $\varepsilon(c')$ if $c'$
intersects $c$, and is equal to 0 otherwise. The {\em flow to the
right} through $c$ is the sum of the flows through $c$ of all arrows
$c'$ such that $c'$ and $c$, in this order, form a positive basis of
$\R^2$. The {\em flow to the left} is defined as the sum of the
flows of all $c'$ such that $c',c$ form a negative basis. The {\em
total flow} through the arrow $c$ is the difference of the right and
the left flows through $c$.

Now, let us observe that if a Gauss diagram is realizable, then the
total flow through each of its arrows is equal to zero. Indeed, let
us cut and re-connect the branches of the knot represented by the
Gauss diagram in the vicinity of the crossing point that corresponds
to the arrow $c$. What we get is a two-component link:
$$\risS{-30}{v2proof1}{}{60}{18}{35}\hspace{100pt}
\risS{-30}{v2proof3}{\put(-2,45){$A$}\put(27,-3){$B$}}{60}{0}{0}
$$
It is easy to see that the two ways of computing the linking number
of the two components $A$ and $B$ (see Section \ref{link_num}) are
equal to the right and the left flow through $c$ respectively. Since
the linking number is an invariant, the difference of the flows is
0.

Now, let us see what happens when the basepoint is moved over an
endpoint of an arrow $c$. If this endpoint corresponds to an
overcrossing, this means that the arrow $c$ does not appear in any
subdiagram of the form (\ref{fig:PV}) and, hence, the value of $a_2$
remains unchanged. If the basepoint of the diagram is moved over an
undercrossing, the value of $a_2$ changes by the amount that is
equal to the number of all subdiagrams of $G$ involving $c$, counted
with signs. Taking the signs into the account, we see that this
amount is equal to the total  flow through the chord $c$ in $G$,
that is, zero.

Let us now verify that $a_2$ is invariant under the Reidemeister
moves. This is clear for the move $V\Omega_1$, since an arrow with
adjacent endpoints cannot participate in a subdiagram of the form
(\ref{fig:PV}).

The move $V\Omega_2$ involves two arrows; denote them by $c_1$ and
$c_2$. Choose the basepoint ``far'' from the endpoints of $c_1$ and
$c_2$, namely, in such a way that it belongs neither to the interval
between the sources of $c_1$ and $c_2$, nor to the interval between
the targets of these arrows. (Since $a_2$ does not depend on the
location of the basepoint, there is no loss of generality in this
choice.) Then the contribution to $a_2$ of any pair that contains
the arrow $c_1$ cancels with the corresponding contribution for
$c_2$.

The moves of type 3 involve three arrows. If we choose a basepoint
far from all of these endpoints, only one of the three distinguished
arrows can participate in a subdiagram of the from (\ref{fig:PV}).
It is then clear that exchanging the endpoints of the three arrows
as in the move $V\Omega_3$ does not affect the value of $a_2$.

It remains to show that $a_2$ has degree $2$. Consider a knot with 3
double points. Resolving the double point, we obtain an alternating
sum of eight knots whose Gauss diagrams are the same except for the
directions and signs of 3 arrows. Any subdiagram of the form
(\ref{fig:PV}) fails to contain at least one of these three arrows.
It is, therefore clear that for each instance that the Gauss diagram
of one of the eight knots contains the diagram (\ref{fig:PV}) as a
subdiagram, there is another occurrence of (\ref{fig:PV}) in another
of the eight knots, counted in $a_2$ with the opposite sign.
\end{proof}

{\bf Remark.} This method of calculating $c_2$ (invented by Polyak
and Viro \cite{PV1, PV2}) is an example of a {\em Gauss diagram
formula}.  See Chapter \ref{chapGD} for details and for more
examples.

\section{Actuality tables}\label{section:actuality_tables}

In general, the amount of information needed to describe a knot
invariant $v$ is infinite, since $v$ is a function on an infinite
domain: the set of isotopy classes of knots. However, Vassiliev
invariants require only a finite amount of information for their
description. We already mentioned the analogy between Vassiliev
invariants and polynomials. A polynomial of degree $n$ can be
described, for example, using the Lagrange interpolation formula, by
its values in $n+1$ particular points. In a similar way a given
Vassiliev invariant can be described by its values on a finitely
many knots. These values are organized in the {\em actuality
table}\index{Actuality table} (see \cite{Va1,BL,Bir2}).

\subsection{Basic knots and actuality tables}
To construct the actuality table we must choose a representative
({\em basic}) singular knot for every chord diagram. A possible
choice of basic knots up to degree 3 is shown in the table.
$$
\begin{array}{c||c||c|c||c|c|c|c|c}
\ChD_0&\ChD_1&\multicolumn{2}{c||}{\ChD_2}&\multicolumn{5}{c}{\ChD_3} \\[2mm]
\hline \cd{cd0}&\chd{cd1ch4}&\cd{cd21ch4}&\cd{cd22ch4}&
\chd{cd31ch4}&\chd{cd32ch4}&\chd{cd33ch4}&
   \chd{cd34ch4}&\chd{cd35ch4} \\[5mm]
\hline
\kn{kn0}{-8}&\kn{eight}{-6}&\kn{kn21}{-8.5}&\kn{kn22}{-8}&
\kn{kn31}{-9}&\kn{kn32}{-8}&\kn{kn33}{-8.5}&\kn{kn34}{-8}&
\kn{kn35}{-8}
\end{array}
\label{basic_knots}
$$

\bigskip
The actuality table for a particular invariant $v$ of order $\leq n$
consists of the set of its values on the set of all basic knots with
at most $n$ double points. The knowledge of this set is sufficient
for calculating $v$ for any knot.

Indeed, any knot $K$ can be transformed into any other knot, in
particular, into the basic knot with no singularities (in the table
above this is the unknot), by means of crossing changes and
isotopies. The difference of two knots that participate in a
crossing change is a knot with a double point, hence in $\Z\K$ the
knot $K$ can be written as a sum of the basic non-singular knot and
several knots with one double point. In turn, each knot with one
double point can be transformed, by crossing changes and isotopies,
into the basic singular knot {\em with the same chord diagram}, and
can be written, as a result, as a sum of a basic knot with one
double point and several knots with two double points. This process
can be iterated until we obtain a representation of the knot $K$ as
a sum of basic knots with at most $n$ double points and several
knots with $n+1$ double points. Now, since $v$ is of order $\leq n$,
it vanishes on the knots with $n+1$ double points, so $v(K)$ can be
written as a sum of the values of $v$ on the basic knots with at
most $n$ singularities.

By Proposition~\ref{fact_thru}, the values of $v$ on the knots with
precisely $n$ double points depend only on their chord diagrams. For
a smaller number of double points, the values of $v$ in the actuality
table depend not only on chord diagrams, but also on the basic
knots. Of course, the values in the actuality table cannot be
arbitrary. They satisfy certain relations which we shall discuss
later (see Section~\ref{4Tsec}). The simplest of these relations,
however, is easy to spot from the examples: the value of any
invariant on a diagram with a chord that has no intersections with
other chords is zero.

\begin{example}\label{Con2} The second coefficient $c_2$ of the Conway
polynomial (Section~\ref{def_VI}) is a Vassiliev invariant of order
$\leq 2$. Here is an actuality table for it.
$$
  c_2:\qquad\qquad\begin{array}{c||c||c|c}0&0&0&1\end{array}
$$
The order of the values in this table corresponds to the order of
basic knots in the table on page~\pageref{basic_knots}.
\end{example}

\begin{example}
A Vassiliev invariant of order 3 is given by the third coefficient
$j_3$ of the Taylor expansion of Jones polynomial
(Section~\ref{Jones}). The actuality table for $j_3$ looks as
follows.
$$j_3:\qquad\qquad\begin{array}{c||c||c|c||c|c|c|c|c}
0&0&0&6&0&0&0&-6&-12 \end{array}$$
\end{example}

\subsection{}
To illustrate the general procedure of computing the value of a
Vassiliev invariant on a particular knot by means of actuality
tables let us compute the value of $j_3$ on the right-hand trefoil.
The right-hand trefoil is an ordinary knot, without singular points,
so we have to deform it (using crossing changes) to our basic knot
without double points, that is, the unknot. This can be done by one
crossing change, and by the Vassiliev skein relation
we have
$$j_3\left(\risA{rtf}{30}{15}{-10}\right) =
  j_3\left(\risA{rtfm}{30}{15}{-10}\right) +
  j_3\left(\risA{rtfs}{30}{15}{-10}\right) =
  j_3\left(\risA{rtfs}{30}{15}{-10}\right) $$
because $j_3(\mathit{unknot})=0$ in the actuality table. Now the
knot with one double point we got is not quite the one from our
basic knots. We can deform it to a basic knot changing the upper
right crossing.
$$j_3\left(\risA{rtfs}{30}{15}{-10}\right) =
  j_3\left(\risA{rtfsm}{30}{15}{-10}\right) +
  j_3\left(\risA{rtfss}{30}{15}{-10}\right) =
  j_3\left(\risA{rtfss1}{30}{15}{-10}\right)$$
Here we used the fact that {\em any} invariant vanishes on the basic
knot with a single double point. The knot with two double points on
the right-hand side of the equation still differs by one crossing
from the basic knot with two double points. This means that we have
to do one more crossing change. Combining these equations together
and using the values from the actuality table we get the final
answer
$$j_3\left(\risA{rtf}{30}{15}{-10}\right) =
  j_3\left(\risA{rtfss1}{30}{15}{-10}\right) =
  j_3\left(\risA{rtfssp}{30}{15}{-10}\right) +
  j_3\left(\risA{rtfsss}{30}{15}{-10}\right) = 6-12=-6$$

\subsection{The first ten Vassiliev invariants}

Using actuality tables, one can find the values of the Vassiliev
invariants of low degree. Table \ref{vi_table} 
uses a certain basis in the space of Vassiliev invariants up to
degree 5. It represents an abridged version of the table compiled by
T.~Stanford \cite{Sta1}, where the values of invariants up to degree
6 are given on all knots with at most 10 crossings.

Some of the entries in Table \ref{vi_table} are different from
\cite{Sta1}, this is due to the fact that, for some non-amphicheiral
knots, Stanford uses mirror reflections of the Rolfsen's
knots shown in Table \ref{knot_table}.
\begin{longtable}{|r|c|||r||r||r||r|r|r||r|r|r|r|}
\hline
&&$v_0$&$v_2$&$v_3$&$v_{41}$&$v_{42}$&$v_2^2$\makebox(0,12){}
                &$v_{51}$&$v_{52}$&$v_{53}$&$v_2v_3$\\
\hline
\hline%          v2    v3    v41   v42  v2^2   v51  v52   v53  v2v3
$0_1$&$++$&$ 1 $&$ 0 $&$ 0 $&$ 0 $&$ 0 $&$ 0 $&$ 0 $&$ 0$&$ 0$&$ 0$\\
\hline\hline
$3_1$&$-+$&      $1$ &  $1$ &  $-1$ &    $1$ &  $-3$ &  $1$ &  $-3$ &   $1$&   $-2$& $-1$\\
\hline\hline
$4_1$&$++$&$      1 $&$ -1 $&$   0 $&$   -2 $&$   3 $&$  1 $&$   0 $&$   0$&$    0$&$  0$\\
\hline\hline
$5_1$&$-+$&$      1 $&$  3 $&$  -5 $&$    1 $&$  -6 $&$  9 $&$ -12 $&$   4$&$   -8$&$-15$\\
\hline
$5_2$&$-+$&$      1 $&$  2 $&$  -3 $&$    1 $&$  -5 $&$  4 $&$  -7 $&$   3$&$   -5$&$ -6$\\
\hline\hline
$6_1$&$-+$&$      1 $&$ -2 $&$   1 $&$   -5 $&$   5 $&$  4 $&$   4 $&$  -1$&$    2$&$ -2$\\
\hline
$6_2$&$-+$&$      1 $&$ -1 $&$   1 $&$   -3 $&$   1 $&$  1 $&$   3 $&$  -1$&$    1$&$ -1$\\
\hline
$6_3$&$++$&$      1 $&$  1 $&$   0 $&$    2 $&$  -2 $&$  1 $&$   0 $&$   0$&$    0$&$  0$\\
\hline\hline
$7_1$&$-+$&$      1 $&$  6 $&$ -14 $&$   -4 $&$  -3 $&$ 36 $&$ -21 $&$   7$&$  -14$&$-84$\\
\hline
$7_2$&$-+$&$      1 $&$  3 $&$  -6 $&$    0 $&$  -5 $&$  9 $&$  -9 $&$   6$&$   -7$&$-18$\\
\hline
$7_3$&$-+$&$      1 $&$  5 $&$  11 $&$   -3 $&$  -6 $&$ 25 $&$  16 $&$  -8$&$   13$&$ 55$\\
\hline
$7_4$&$-+$&$      1 $&$  4 $&$   8 $&$   -2 $&$  -8 $&$ 16 $&$  10 $&$  -8$&$   10$&$ 32$\\
\hline
$7_5$&$-+$&$      1 $&$  4 $&$  -8 $&$    0 $&$  -5 $&$ 16 $&$ -14 $&$   6$&$   -9$&$-32$\\
\hline
$7_6$&$-+$&$      1 $&$  1 $&$  -2 $&$    0 $&$  -3 $&$  1 $&$  -2 $&$   3$&$   -2$&$ -2$\\
\hline
$7_7$&$-+$&$      1 $&$ -1 $&$  -1 $&$   -1 $&$   4 $&$  1 $&$   0 $&$   2$&$    0$&$  1$\\
\hline\hline
$8_1$&$-+ $&$     1 $&$ -3 $&$   3 $&$   -9 $&$   5 $&$  9 $&$  12 $&$  -3 $&$   5$&$ -9$ \\
\hline
$8_2$&$-+ $&$     1 $&$  0 $&$   1 $&$   -3 $&$  -6 $&$  0 $&$   2 $&$   0 $&$  -3$&$  0$ \\
\hline
$8_3$&$++ $&$     1 $&$ -4 $&$   0 $&$  -14 $&$   8 $&$ 16 $&$   0 $&$   0 $&$   0$&$  0$ \\
\hline
$8_4$&$-+ $&$     1 $&$ -3 $&$   1 $&$  -11 $&$   4 $&$  9 $&$   0 $&$  -2 $&$  -1$&$ -3$ \\
\hline
$8_5$&$-+ $&$     1 $&$ -1 $&$  -3 $&$   -5 $&$  -5 $&$  1 $&$  -5 $&$   3 $&$   2$&$  3$ \\
\hline
$8_6$&$-+ $&$     1 $&$ -2 $&$   3 $&$   -7 $&$   0 $&$  4 $&$   9 $&$  -3 $&$   2$&$ -6$ \\
\hline
$8_7$&$-+ $&$     1 $&$  2 $&$   2 $&$    4 $&$  -2 $&$  4 $&$   7 $&$  -1 $&$   3$&$  4$ \\
\hline
$8_8$&$-+ $&$     1 $&$  2 $&$   1 $&$    3 $&$  -4 $&$  4 $&$   2 $&$  -1 $&$   1$&$  2$ \\
\hline
$8_9$&$++ $&$     1 $&$ -2 $&$   0 $&$   -8 $&$   1 $&$  4 $&$   0 $&$   0 $&$   0$&$  0$ \\
\hline
$8_{10}$&$-+ $&$  1 $&$  3 $&$   3 $&$    3 $&$  -6 $&$  9 $&$   5 $&$  -3 $&$   3$&$  9$ \\
\hline
$8_{11}$&$-+ $&$  1 $&$ -1 $&$   2 $&$   -4 $&$  -2 $&$  1 $&$   8 $&$  -1 $&$   2$&$ -2$ \\
\hline
$8_{12}$&$++ $&$  1 $&$ -3 $&$   0 $&$   -8 $&$   8 $&$  9 $&$   0 $&$   0 $&$   0$&$  0$ \\
\hline
$8_{13}$&$-+ $&$  1 $&$  1 $&$   1 $&$    3 $&$   0 $&$  1 $&$   6 $&$   0 $&$   3$&$  1$ \\
\hline
$8_{14}$&$-+ $&$  1 $&$  0 $&$   0 $&$   -2 $&$  -3 $&$  0 $&$  -2 $&$   0 $&$  -3$&$  0$ \\
\hline
$8_{15}$&$-+ $&$  1 $&$  4 $&$  -7 $&$    1 $&$  -7 $&$ 16 $&$ -16 $&$   5 $&$ -10$&$-28$ \\
\hline
$8_{16}$&$-+ $&$  1 $&$  1 $&$  -1 $&$    3 $&$   0 $&$  1 $&$   2 $&$   2 $&$   2$&$ -1$ \\
\hline
$8_{17}$&$+- $&$  1 $&$ -1 $&$   0 $&$   -4 $&$   0 $&$  1 $&$   0 $&$   0 $&$   0$&$  0$ \\
\hline
$8_{18}$&$++ $&$  1 $&$  1 $&$   0 $&$    0 $&$  -5 $&$  1 $&$   0 $&$   0 $&$   0$&$  0$ \\
\hline
$8_{19}$&$-+ $&$  1 $&$  5 $&$  10 $&$    0 $&$  -5 $&$ 25 $&$  18 $&$  -6 $&$  10$&$ 50$ \\
\hline
$8_{20}$&$-+ $&$  1 $&$  2 $&$  -2 $&$    2 $&$  -5 $&$  4 $&$  -1 $&$   3 $&$  -1$&$ -4$ \\
\hline
$8_{21}$&$-+ $&$  1 $&$  0 $&$   1 $&$   -1 $&$  -3 $&$  0 $&$   1 $&$  -1 $&$  -1$&$  0$ \\
\hline
\caption{Vassiliev invariants of order $\leq5$} \label{vi_table}
\end{longtable}

\noindent
The two signs after the knot number refer to their symmetry
properties: a plus in the first position means that the knot is
amphicheiral, a plus in the second position means that the knot is
invertible.

\section{Vassiliev invariants of tangles}

Knots are tangles whose skeleton is a circle. A theory of Vassiliev
invariants, similar to the theory for knots, can be constructed for
isotopy classes of tangles with any given skeleton $\boldX$.

Indeed, similarly to the case of knots, one can introduce {\em
tangles with double points}, with the only extra assumption that the
double points lie in the interior of the tangle box. Then, any
invariant of tangles can be extended to tangles with double points
with the help of the Vassiliev skein relation. An invariant of
tangles is a Vassiliev invariant of degree $\le n$ if it vanishes on
all tangles with more that $n$ double points.

We stress that we define Vassiliev invariants separately for each
skeleton $\boldX$. Nevertheless, there are relations among
invariants of tangles with different skeleta.

\begin{xexample} Assume that the isotopy classes of tangles with
the skeleta $\boldX_1$ and $\boldX_2$ can be multiplied. Given a
tangle $T$ with skeleton $\boldX_1$ and a Vassiliev invariant $v$ of
tangles with skeleton $\boldX_1\boldX_2$, we can define an invariant
of tangles on $\boldX_2$ of the same order as $v$ by composing a
tangle with $T$ and applying $v$.
\end{xexample}

\begin{xexample} In the above example the product of tangles can
be replaced by their tensor product. (Of course, the condition that
$\boldX_1$ and $\boldX_2$ can be multiplied is no longer necessary
here.)
\end{xexample}

In particular, the Vassiliev invariants of tangles whose skeleton
has one component, can be identified with the Vassiliev  invariants
of knots.

\begin{xexample} Assume that $\boldX'$ is obtained from $\boldX$ by
dropping one or several components. Then any Vassiliev invariant
$v'$ of tangles with skeleton $\boldX'$ gives rise to an invariant
$v$ of tangles on $\boldX$ of the same order; to compute $v$ drop
the components of the tangle that are not in $\boldX'$ and apply
$v'$.
\end{xexample}

This example immediately produces a lot of tangle invariants of
finite type: namely, those coming from knots. The simplest example
of a Vassiliev invariant that does not come from knots is the
linking number of two components of a tangle. So far, we have
defined the linking number only for pairs of closed curves. If one
or both of the components are not closed, we can use the
constructions above to close them up in some fixed way.

\begin{xlemma}
The linking number of two components of a tangle is a Vassiliev
invariant of order 1.
\end{xlemma}

\begin{proof}
Consider a two-component link with one double point. This double
point can be of two types: either it is a self-intersection point of a
single component, or it is an intersection of two different
components. Using the Vassiliev skein relation and the
formula~\ref{comb_lk_num}, we see that in the first case the linking
number vanishes, while in the second case it is equal to 1. It
follows that for a two-component link with two double points the
linking number is always zero.
\end{proof}

Among the invariants for all classes of tangles, the string link
invariants have attracted most attention. Two particular classes of
string link invariants are the knot invariants (recall that string
links on one strand are in one-to-one correspondence with knots) and
the invariants of pure braids. We shall treat the Vassiliev
invariants of pure braids in detail in Chapter~\ref{chapBr}.

\begin{xcb}{Exercises}

\begin{enumerate}
\item
Using the actuality tables, compute the value of $j_3$ on the
left-hand trefoil.

\item
Choose the basic knots with four double points and construct the
actuality tables for the fourth coefficients $c_4$ and $j_4$ of the
Conway and Jones polynomials.

\item Prove that $j_0(K)=1$ and $j_1(K)=0$ for any knot $K$.

\item Show that the value of $j_0$ on a link with $k$ components
is equal to $(-2)^{k-1}$.

\item
For a link $L$ with two components $K_1$ and $K_2$ prove that\\
$j_1(L) = -3\cdot lk(K_1,K_2)$.\index{Linking number} In other
words,
$$J(L) = -2 -3\cdot lk(K_1,K_2)\cdot h + j_2(L)\cdot h^2 +
    j_3(L)\cdot h^3+\dots\ .
$$

\item\label{div_by_6}
Prove that for any knot $K$ the integer $j_3(K)$ is divisible by 6.

\item For a knot $K$, find the relation between the second coefficients
$c_2(K)$ and $j_2(K)$ of the Conway and Jones polynomials.

\item
Prove that $v(3_1\#3_1)=2v(3_1)-v(0)$, where $0$ is the trivial knot,
for any Vassiliev invariant $v\in\V_3$.

\item Prove that for a knot $K$ the $n$th derivative at $1$ of the Jones
polynomial
$$\left.\frac{d^n(J(K))}{dt^n}\right|_{t=1}$$
is a Vassiliev invariant of order $\leq n$. Find the relation
between these invariants and $j_1,\dots,j_n$ for small values of
$n$.

\item
Express the coefficients $c_2$, $c_4$, $j_2$, $j_3$, $j_4$, $j_5$ of the
Conway and Jones polynomials in terms of the basis
Vassiliev invariants from Table \ref{vi_table}.

\item
Find the symbols of the Vassiliev invariants from Table
\ref{vi_table}.

\item
Express the invariants of Table \ref{vi_table} through the coefficients
of the Conway and the Jones polynomials.

\item
Find the actuality tables for some of the Vassiliev invariants appearing in
Table \ref{vi_table}.

\item
Explain the correlation between the first sign and the zeroes in the
last four columns of Table \ref{vi_table}.

\item
Check that Vassiliev invariants up to order 4 are enough to distinguish,
up to orientation, all knots with at most 8 crossings from
Table \ref{knot_table} on page \pageref{knot_table}.

\item \label{ex_consymb}
\index{Conway polynomial!symbol}\index{Symbol!of the Conway coefficients}
Prove that the symbol of the
coefficient $c_n$ of the Conway polynomial can be calculated as
follows.  Double every chord of a given chord diagram $D$ as in
Section~\ref{symb_j_n}, and let $|D|$ be equal to the number of
components of the obtained curve. Then
$$\symb(c_n)(D)=\left\{\begin{array}{cl}
  1, & \mbox{if\ } |D|=1\\
  0, & \mbox{otherwise}\ .
  \end{array}\right.
$$

\item
Prove that $c_n$ is a Vassiliev invariant of degree exactly $n$.

\item
Prove that there is a well-defined extension of knot invariants to
singular knots with a non-degenerate \textit{triple} point according
to the rule
$$f(\chd{tr_point})=f(\chd{tr_p_res1})-f(\chd{tr_p_res2})\ .
$$
Is it true that, according to this extension, a Vassiliev invariant
of degree 2 is equal to 0 on any knot with a triple point?

Is it possible to use the same method to define an extension of knot
invariants to knots with self-intersections of multiplicity higher
than 3?

\item
Following Example \ref{ex_J_sing}, find the power series expansion
of the modified Jones polynomial of the singular knot
\rb{-8mm}{\ig[height=11mm]{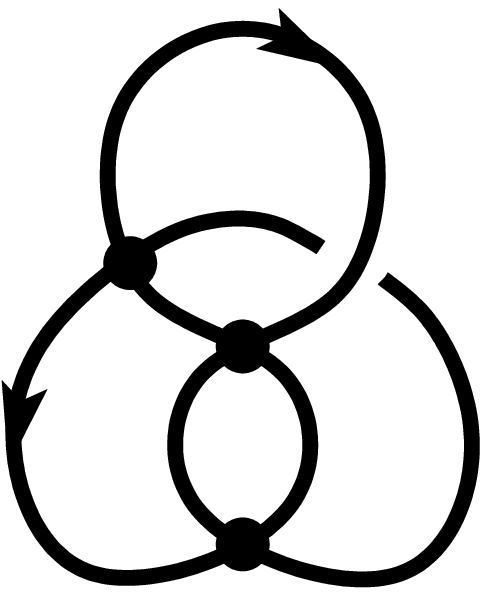}}\ .

\item
Prove the following relation between the Casson knot invariant $c_2$,
extended to singular knots, and the linking number of two curves.
Let $K$ be a knot with one double point. Smoothing the double point by
the rule $\double \mapsto \twoup$, one obtains a 2-component link $L$.
Then $lk(L)=c_2(K)$.

\item
Is there a prime knot $K$ such that $j_4(K)=0$?

\item {\bf Vassiliev invariants from the HOMFLY polynomial.}
\label{pr:vas-homfly} \index{HOMFLY polynomial!Vassiliev invariants}
\index{Vassiliev!invariant!from HOMFLY} For a link $L$ make a
substitution $a=e^h$ in the HOMFLY polynomial $P(L)$ and take the
Taylor expansion in $h$. The result will be a Laurent polynomial in
$z$ and a power series in $h$. Let $p_{k,l}(L)$ be its coefficient
at $h^kz^l$.
\begin{itemize}
\item[(a)] Show that  for any link $L$ the total degree $k+l$ is not negative.
\item[(b)] If $l$ is odd, then $p_{k,l}=0$.
\item[(c)] Prove that $p_{k,l}(L)$ is a Vassiliev invariant of order $\leqslant k+l$.
\item[(d)] Describe the symbol of $p_{k,l}(L)$.
\end{itemize}

\end{enumerate}
\end{xcb}
 %3 finite type inv
\chapter{Chord diagrams} % 04
\label{ChD}
\newcommand{\cld}[1]{\risS{-9.5}{#1}{}{23}{20}{15}}

A chord diagram encodes the order of double points along a singular
knot. We saw in the last chapter that a Vassiliev invariant of order
$n$ gives rise to a function on chord diagrams with $n$ chords. Here
we shall describe the conditions, called one-term and four-term
relations, that a function on chord diagrams should satisfy in order
to come from a Vassiliev invariant. We shall see that the vector
space spanned by chord diagrams modulo these relations has the
structure of a Hopf algebra. This Hopf algebra turns out to be dual
to the algebra of the Vassiliev invariants.

\section{Four- and one-term relations}
\label{4Tsec}

Recall that $\Ring$ denotes a commutative ring and $\V_n$ is the
space of $\Ring$-valued Vassiliev invariants of order $\le n$. Some
of our results will only hold when $\Ring$ is a field of
characteristic 0; sometimes we shall take $\Ring=\C$. On
page~\pageref{alfa_n} in Section~\ref{def_VI} we constructed a
linear inclusion (the symbol of an invariant)
$$
  \ol{\alpha}_n: \V_n/\V_{n-1} \to \Ring\ChD_n,
$$
where  $\Ring\ChD_n$ is the space of $\Ring$-valued functions on the
set $\ChD_n$ of chord diagrams of order $n$.

To describe the image of $\ol{\alpha}_n$, we need the following
definition.

\begin{definition}\index{Weight system}
\index{Relation!four-term}
\index{Four-term relation!for chord diagrams}
A function $f\in\Ring\ChD_n$ is said to satisfy the {\em 4-term (or
4T) relations} if the alternating sum of the values of $f$ is zero
on the following quadruples of diagrams:
\begin{equation}\label{4Teq}
  f(\chd{4T1})-f(\chd{4T2})+f(\chd{4T4})-f(\chd{4T3})=0.
\end{equation}
In this case $f$ is also called a {\em (framed) weight system} of
order $n$.
\end{definition}

Here it is assumed that the diagrams in the pictures may have other
chords with endpoints on the dotted arcs, while all the endpoints of
the chords on the solid portions of the circle are explicitly shown.
For example, this means that in the first and second diagrams the
two bottom points are adjacent. The chords omitted from the pictures
should be the same in all the four cases.

\begin{xexample}
Let us find all 4-term relations for chord diagrams of order 3. We
must add one chord in one and the same way to all the four terms of
Equation~(\ref{4Teq}). Since there are 3 dotted arcs, there are 6
different ways to do that, in particular,\vspace{-7pt}
$$f(\chd{4T1-3a})-f(\chd{4T2-3a})+f(\chd{4T4-3a})-f(\chd{4T3-3a})=0
\vspace{-10pt}$$
and\vspace{-7pt}
$$f(\chd{4T1-3b})-f(\chd{4T2-3b})+f(\chd{4T3-3b})-f(\chd{4T4-3b})=0
$$
Some of the diagrams in these equations are equal, and the relations can be
simplified as
$f(\chd{cd31ch4})=f(\chd{cd32ch4})$\ and\
$f(\chd{cd35ch4})-2f(\chd{cd34ch4})+f(\chd{cd33ch4})=0$.\label{cd-rel-3}
The reader is invited to check that the remaining four 4-term relations
(we wrote only 2 out of 6) are either trivial or coincide with one of these
two.
\end{xexample}

It is often useful to look at a 4T relation from the following point
of view. We can think that one of the two chords that participate in
equation~(\ref{4Teq}) is fixed, and the other is moving. One of the
ends of the moving chord is also fixed, while the other end travels
around the fixed chord stopping at the four locations adjacent to
its endpoints. The resulting four diagrams are then summed up with
alternating signs. Graphically,\vspace{-6pt}
\begin{equation} \label{leg_motion}
  f(\chd{4Tm1})-f(\chd{4Tm2})+f(\chd{4Tm3})-f(\chd{4Tm4})=0.\vspace{-8pt}
\end{equation}
where the fixed end of the moving chord is marked by \krestik.

Another way of writing the 4T relation, which will be useful in
Section~\ref{algFD}, is to split the four terms into two pairs:
$$f(\chd{4T1})-f(\chd{4T2})= f(\chd{4T3})-f(\chd{4T4})\ .
$$

Because of the obvious symmetry, this can be completed as follows:
\begin{equation}\label{3pairs}
f(\chd{4T1})-f(\chd{4T2})= f(\chd{4T5})-f(\chd{4T6})\ .
\end{equation}

Note that for each order $n$
the choice of a specific
4-term relation depends on the following data:
\begin{itemize}
\setlength{\itemsep}{1pt plus 1pt minus 1pt}
\item
   a diagram of order $n-1$,
\item
   a distinguished chord of this diagram (``fixed chord''), and
\item
   a distinguished arc on the circle of this diagram (where the fixed
endpoint of the ``moving chord'' is placed).
\end{itemize}

There are 3 fragments of the circle that participate in a 4-term
relation, namely, those that are shown by solid lines in the
equations above. If these 3 fragments are drawn as 3 vertical line
segments, then the 4-term relation can be restated as follows:
\vspace{3mm}
\begin{multline}\label{horizontal4T}
     (-1)^\downarrow\
     f\Bigl(\ \rb{-6mm}{\ig[height=12mm]{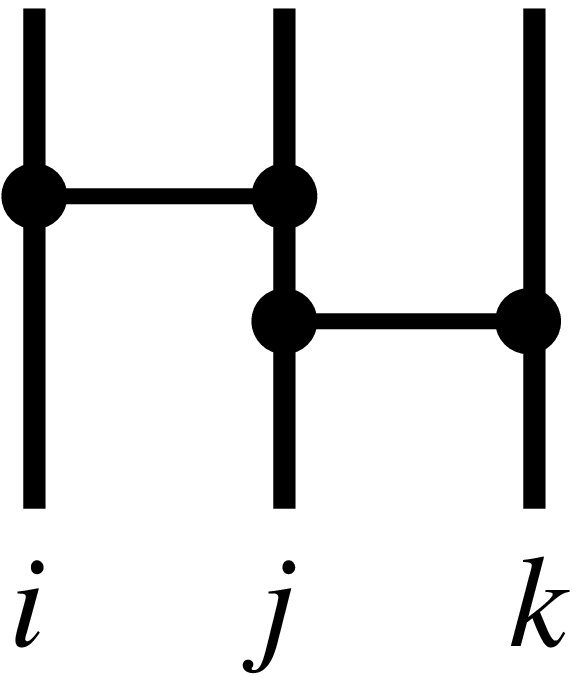}}\ \Bigr)
\ -\
(-1)^\downarrow\
 f\Bigl(\ \rb{-6mm}{\ig[height=12mm]{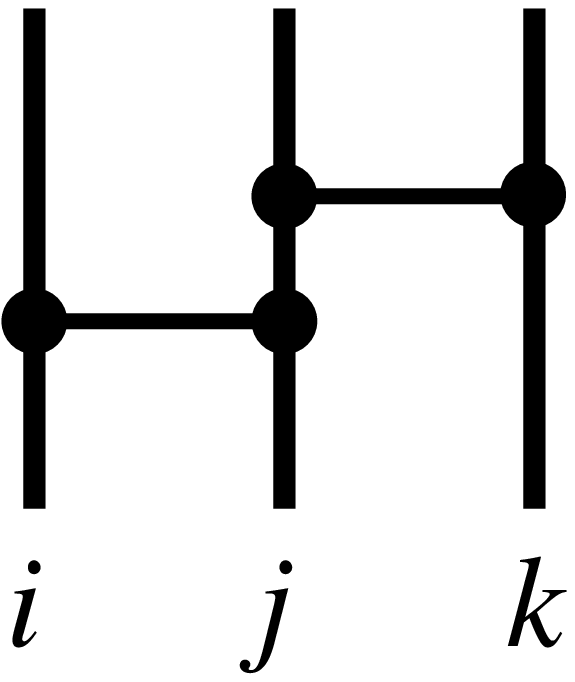}}\ \Bigr)
\\ +\
(-1)^\downarrow\
f\Bigl(\ \rb{-6mm}{\ig[height=12mm]{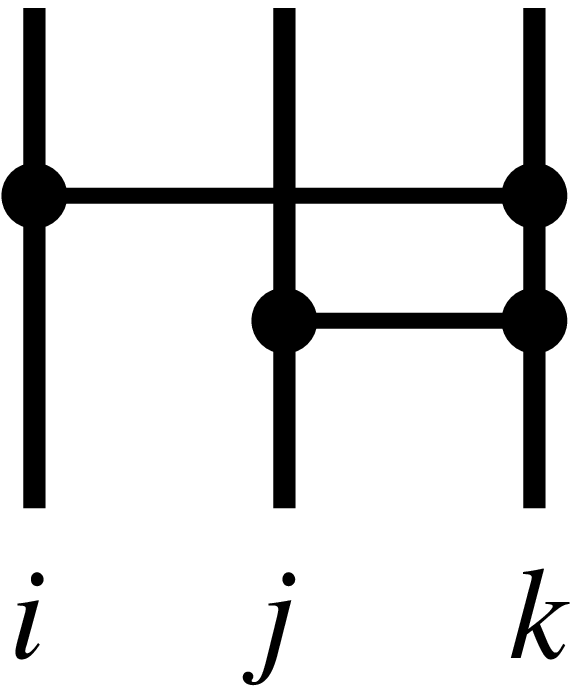}}\ \Bigr) \ -\
(-1)^\downarrow\
f\Bigl(\ \rb{-6mm}{\ig[height=12mm]{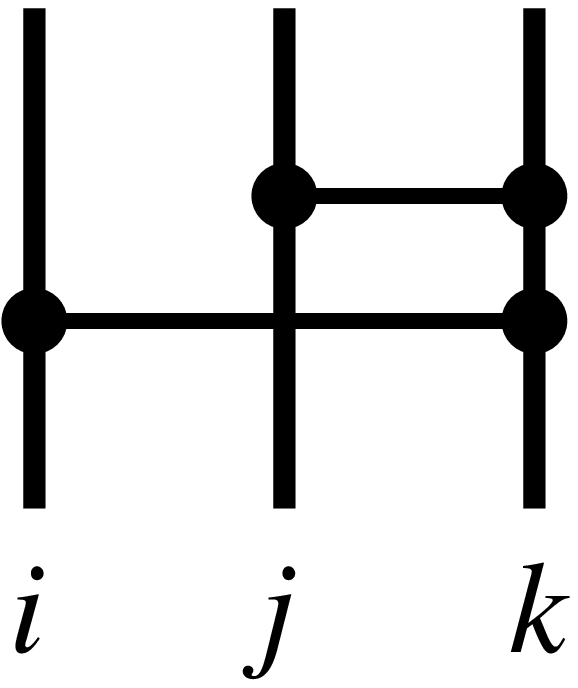}}\ \ \Bigr)\ =\ 0\,.
\end{multline}
where $\downarrow$ stands for the number of endpoints of the chords in which
the orientation of the strands is directed downwards.
This form of a 4T relation is called a {\em horizontal 4T relation}
\index{Four-term relation!horizontal}. (See also Section \ref{line_CD}).
It first appeared, in a different context, in the work by T.~Kohno
\cite{Koh2}.

\begin{xca} Choose some orientations of the three fragments of the
circle, add the portions necessary to close it up and check that the
last form of the 4-term relation carries over into the ordinary
four-term relation.

Here is an example:
$$f\Bigl(\ \rb{-4mm}{\ig[height=11mm]{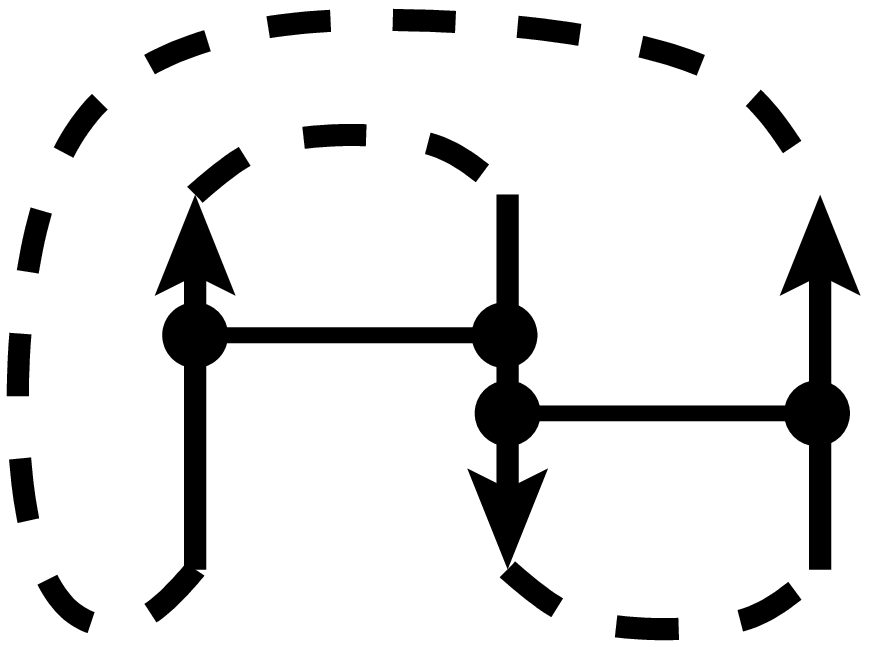}}\ \Bigr)
\ -\ f\Bigl(\ \rb{-4mm}{\ig[height=11mm]{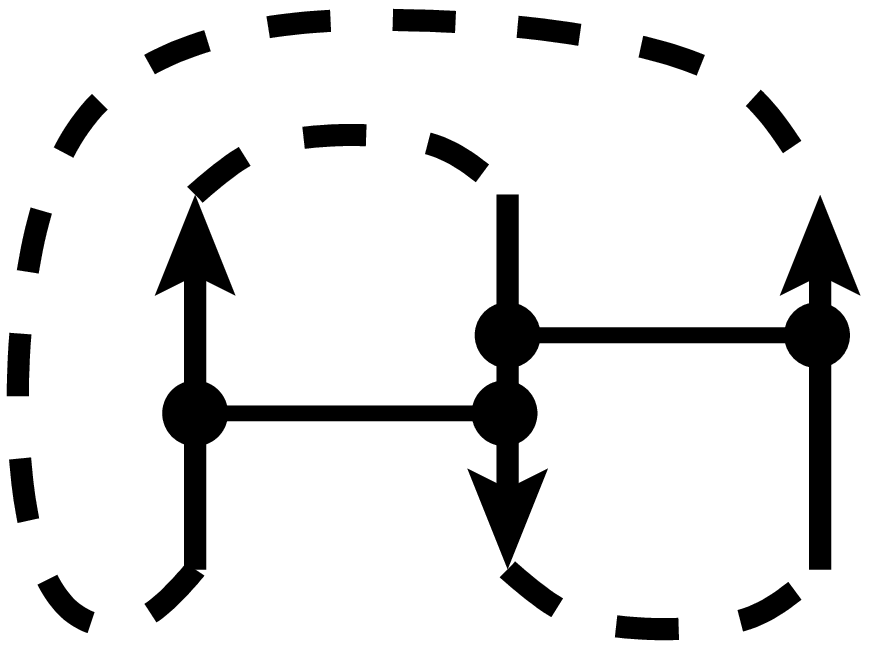}}\ \Bigr) \ -\
f\Bigl(\ \rb{-4mm}{\ig[height=11mm]{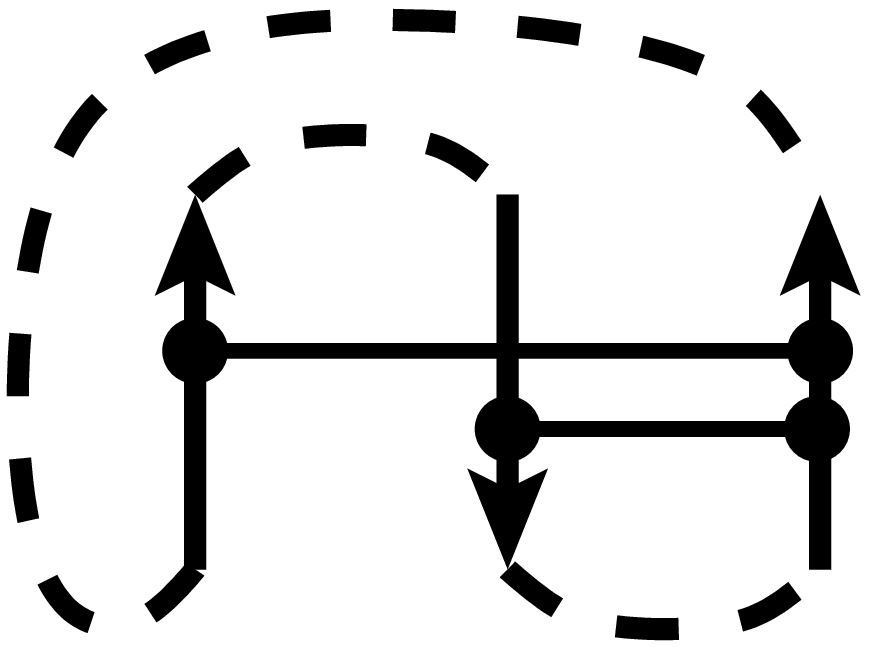}}\ \Bigr) \ +\
f\Bigl(\ \rb{-4mm}{\ig[height=11mm]{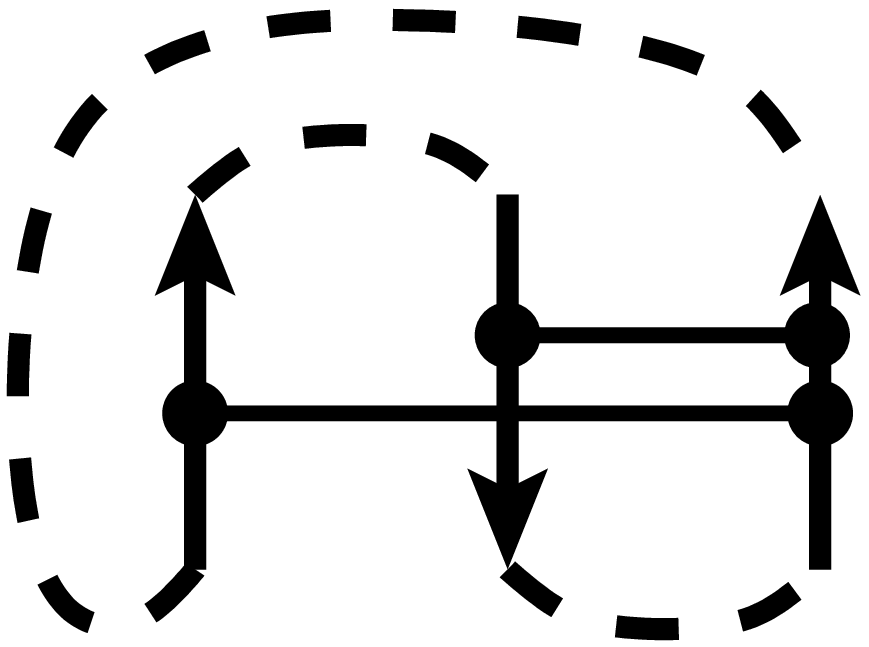}}\ \Bigr)\ =\ 0\,.
$$
\end{xca}

We shall see in the next section that the four-term relations are
always satisfied by the symbols of Vassiliev invariants, both in the
usual and in the framed case. For the framed knots, there are no
other relations; in the unframed case, there is another set of
relations, called {\em one-term}, or {\em framing independence}
relations.

\begin{definition}\index{Weight system!unframed}
An {\em isolated chord\/} \index{Chord!isolated} is a chord that
does not intersect any other chord of the diagram. A function
$f\in\Ring\ChD_n$ is said to satisfy the {\em 1-term relations} if
it vanishes on every chord diagram with an isolated chord. An {\em
unframed weight system\/} of order $n$ is a weight system that
satisfies the 1-term relations.
\index{Relation!one-term}\index{One-term relation}
\end{definition}

Here is an example of a 1T relation:\ $f(\chd{isol}) = 0$\ .

\subsection{Notation.}\label{W_n}\label{W_n^{fr}}
We denote by $\W_n^{fr}$ the subspace of $\Ring\ChD_n$ consisting of
all (framed) weight systems of order $n$ and by
$\W_n\subset\W_n^{fr}$ the space of all unframed weight systems of
order $n$. 

\section{The Fundamental Theorem}
\label{MT_sec}

In Section \ref{ch_diag} we showed that the symbol of an invariant
gives an injective map $\ol{\alpha}_n: \V_n/\V_{n-1}
\to\Ring\ChD_n$. The Fundamental Theorem on Vassiliev invariants
describes its image.

\begin{theorem}[Vassiliev--Kontsevich]
\label{fund_thm}\index{Theorem!Vassiliev--Kontsevich} For $\Ring=\C$
the map $\ol{\alpha}_n$ identifies $\V_n/\V_{n-1}$ with the subspace
of unframed weight systems $\W_n\subset\Ring\ChD_n$. In other words,
the space of unframed weight systems
is isomorphic to the graded space
associated with the filtered space
of Vassiliev invariants,
$$\W=\bigoplus\limits_{n=0}^\infty\W_n\ \isom\
\bigoplus\limits_{n=0}^\infty\V_n/\V_{n+1}\ .$$
\end{theorem}

The theorem consists of two parts:
\begin{itemize}
\setlength{\itemsep}{1pt plus 1pt minus 1pt}
\item
  (V.~Vassiliev) The symbol of every Vassiliev invariant is an unframed weight
system.
\item
  (M.~Kontsevich) Every unframed weight system is the symbol of a certain
Vassiliev invariant.
\end{itemize}

We shall now prove the first (easy) part of the theorem. The second
(difficult) part will be proved later (in Section~\ref{pkt}) using
the Kontsevich integral.
\medskip

The first part of the theorem consists of two assertions, and we prove them
one by one.

\subsection{First assertion:} \label{first_assertion}
{\it any function $f\in\Ring\ChD_n$ coming from an invariant
$v\in\V_n$ satisfies the 1-term relations.}

\begin{proof}
Let $K$ be a singular knot whose chord diagram contains an isolated
chord. The double point $p$ that corresponds to the isolated chord
divides the knot into two parts: $A$ and $B$.

\begin{center}
\ig[height=20mm]{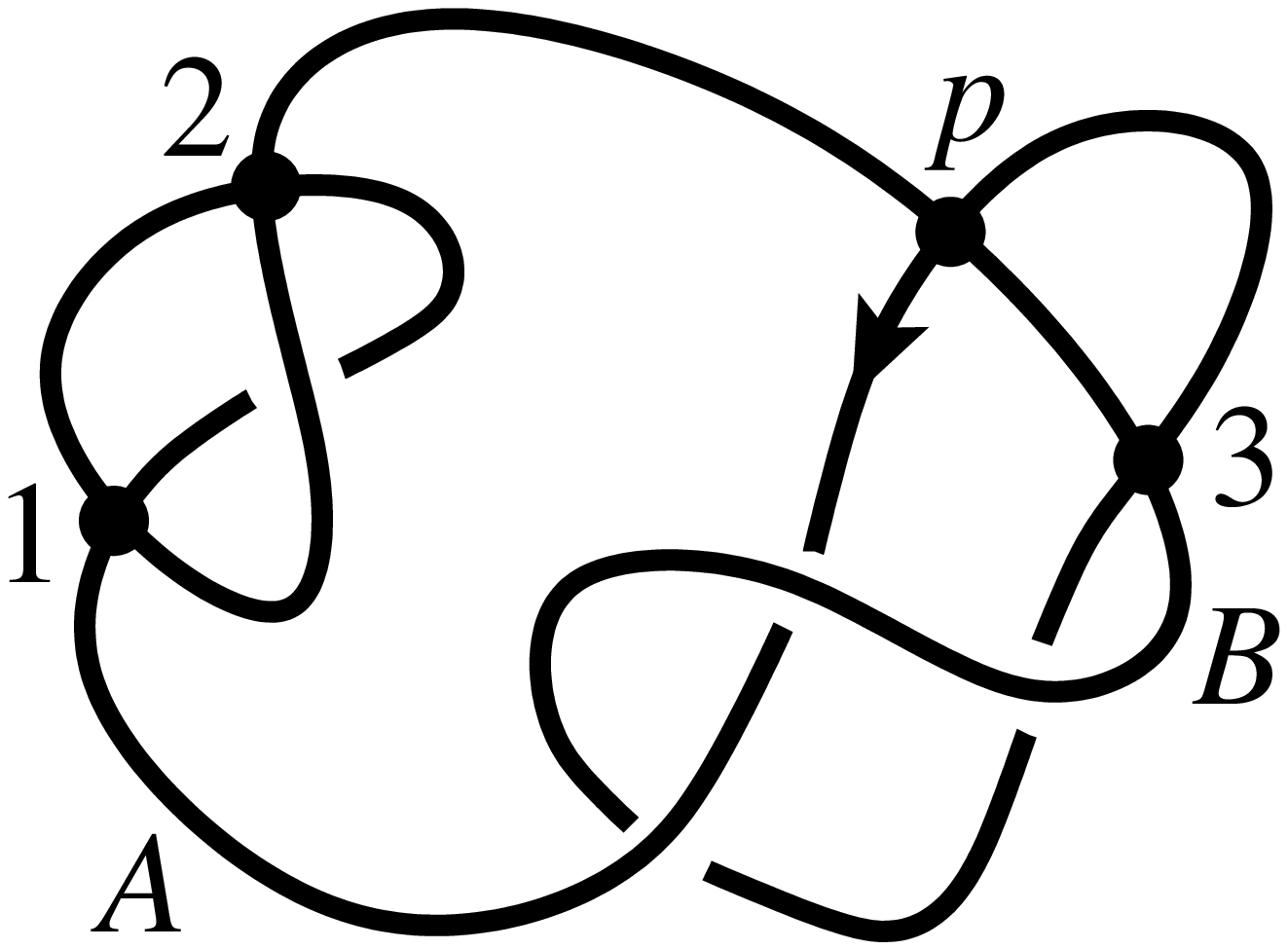}
\hspace{2cm}
\ig[height=20mm]{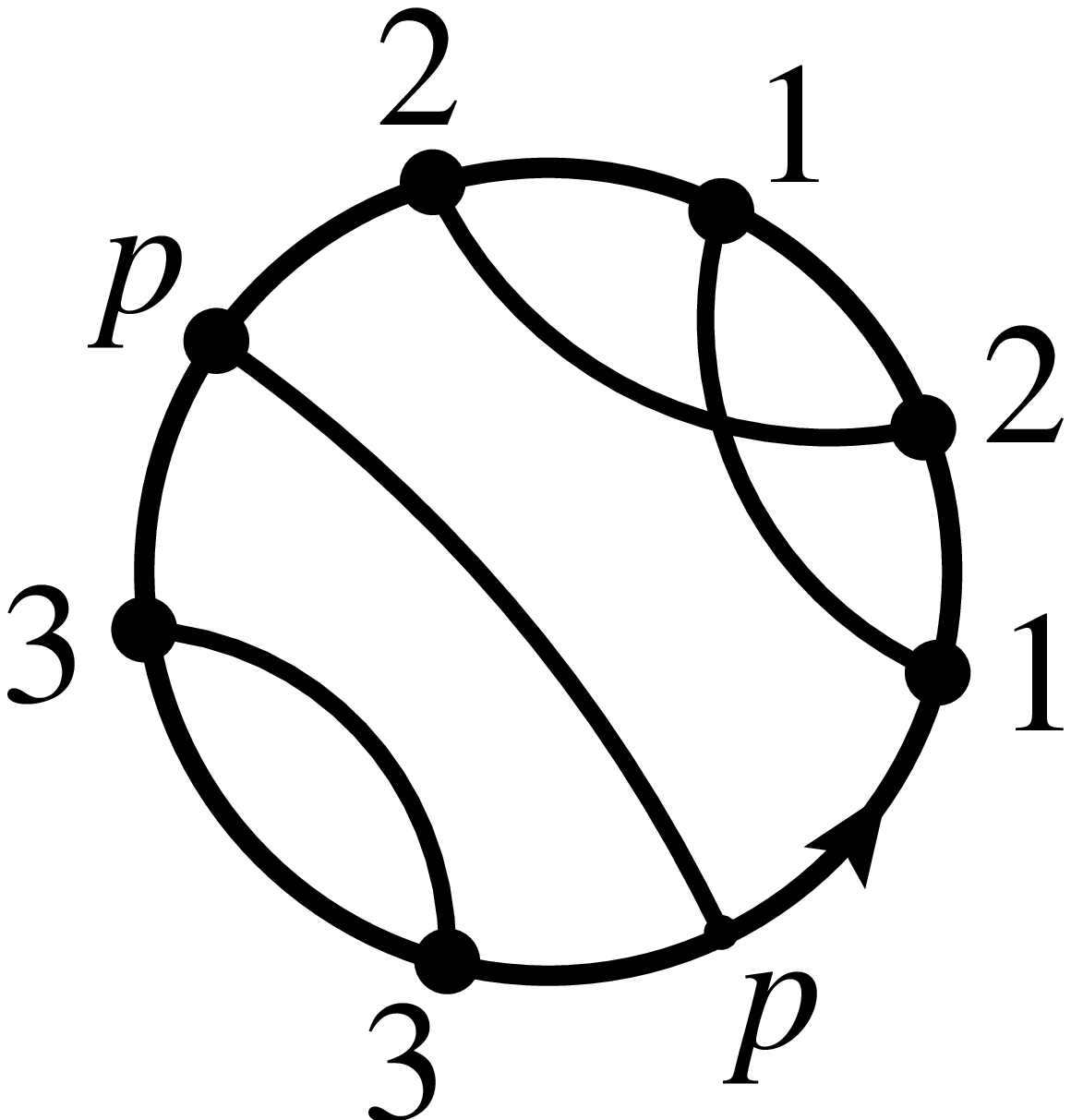}
\end{center}

The fact that the chord is isolated means that $A$ and $B$ do not
have common double points. There may, however, be crossings
involving branches from both parts. 
By crossing changes, we can untangle part $A$ from part $B$ thus
obtaining a singular knot $K'$ with the same chord diagram as $K$
and with the property that the two parts lie on either side of some
plane in $\R^3$ that passes through the double point $p$\ :
\begin{center}
\ig[width=25mm]{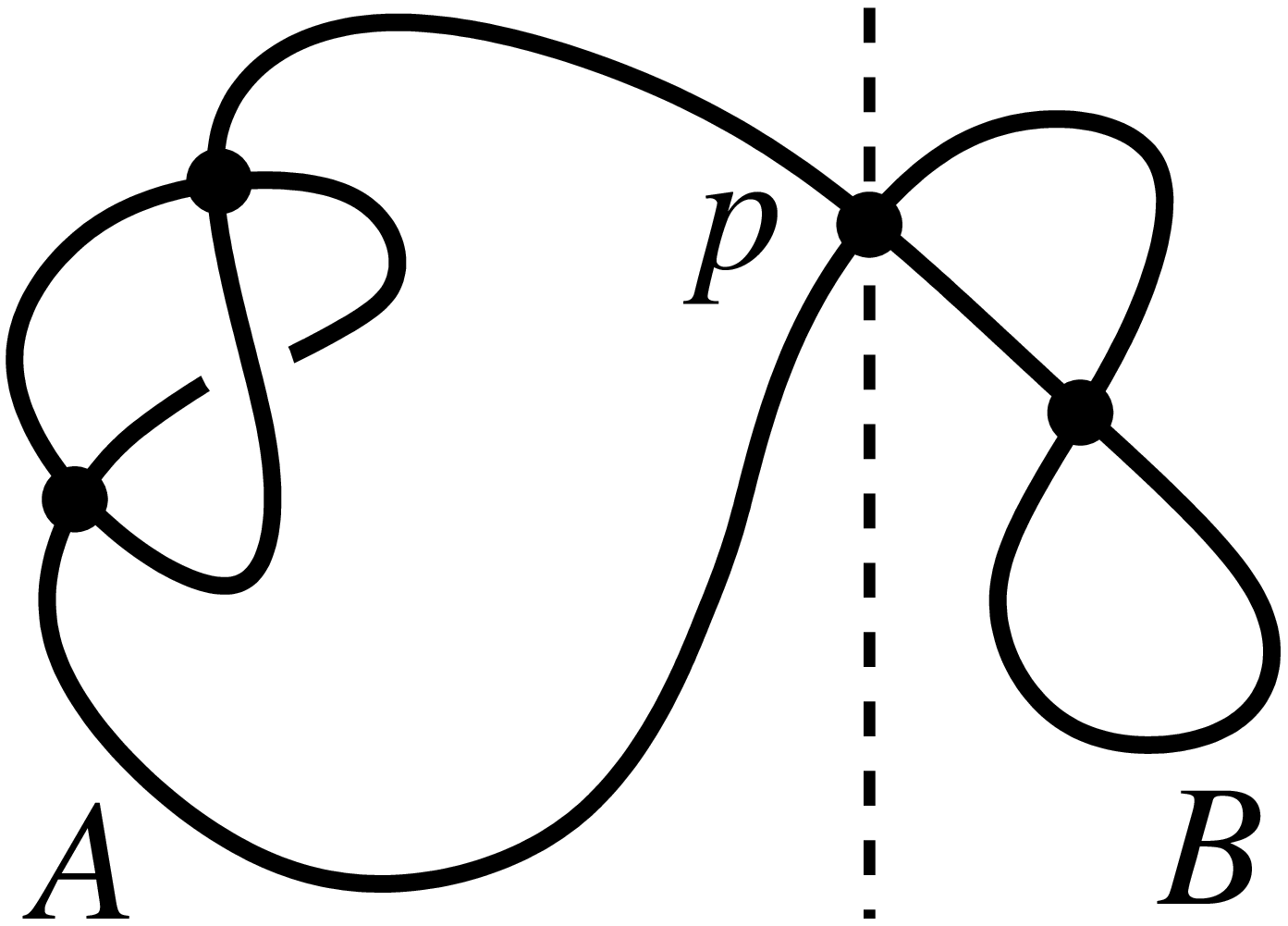}
\end{center}
Here it is obvious that the two resolutions of the double point $p$
give equivalent singular knots, therefore
$v(K)=v(K')=v(K'_{+})-v(K'_{-})=0$.
\end{proof}

\subsection{Second assertion:}
{\it any function $f\in\Ring\ChD_n$ coming from an invariant
$v\in\V_n$ satisfies the 4-term relations.}

\begin{proof}
We shall use the following lemma.

\begin{xlemma}[4-term relation for knots]
\index{Four-term relation!for knots}
Any Vassiliev invariant satisfies
$$
   f\bigl(\rb{-13pt}{\ig[width=30pt]{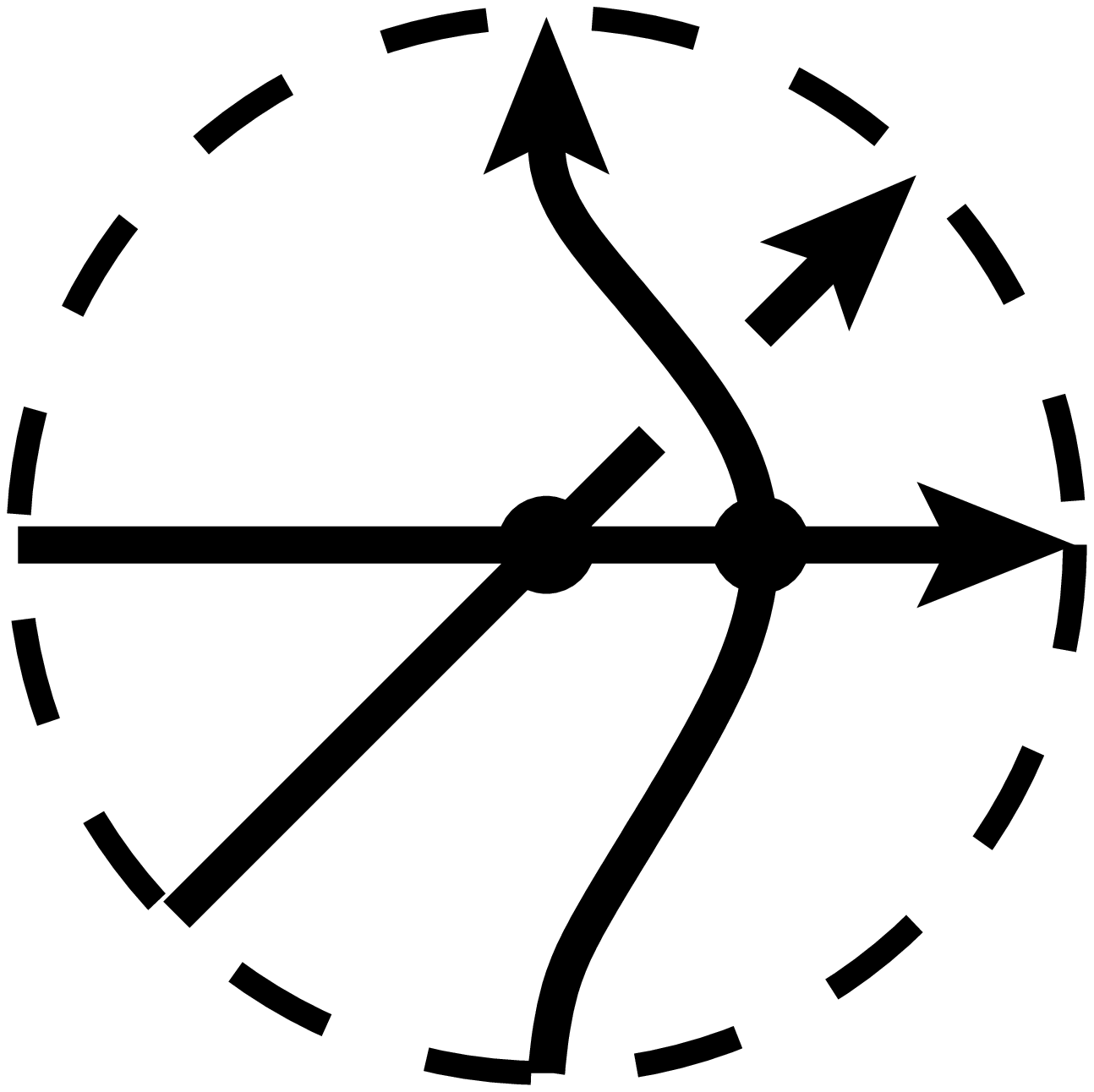}}\bigr)
 - f\bigl(\rb{-13pt}{\ig[width=30pt]{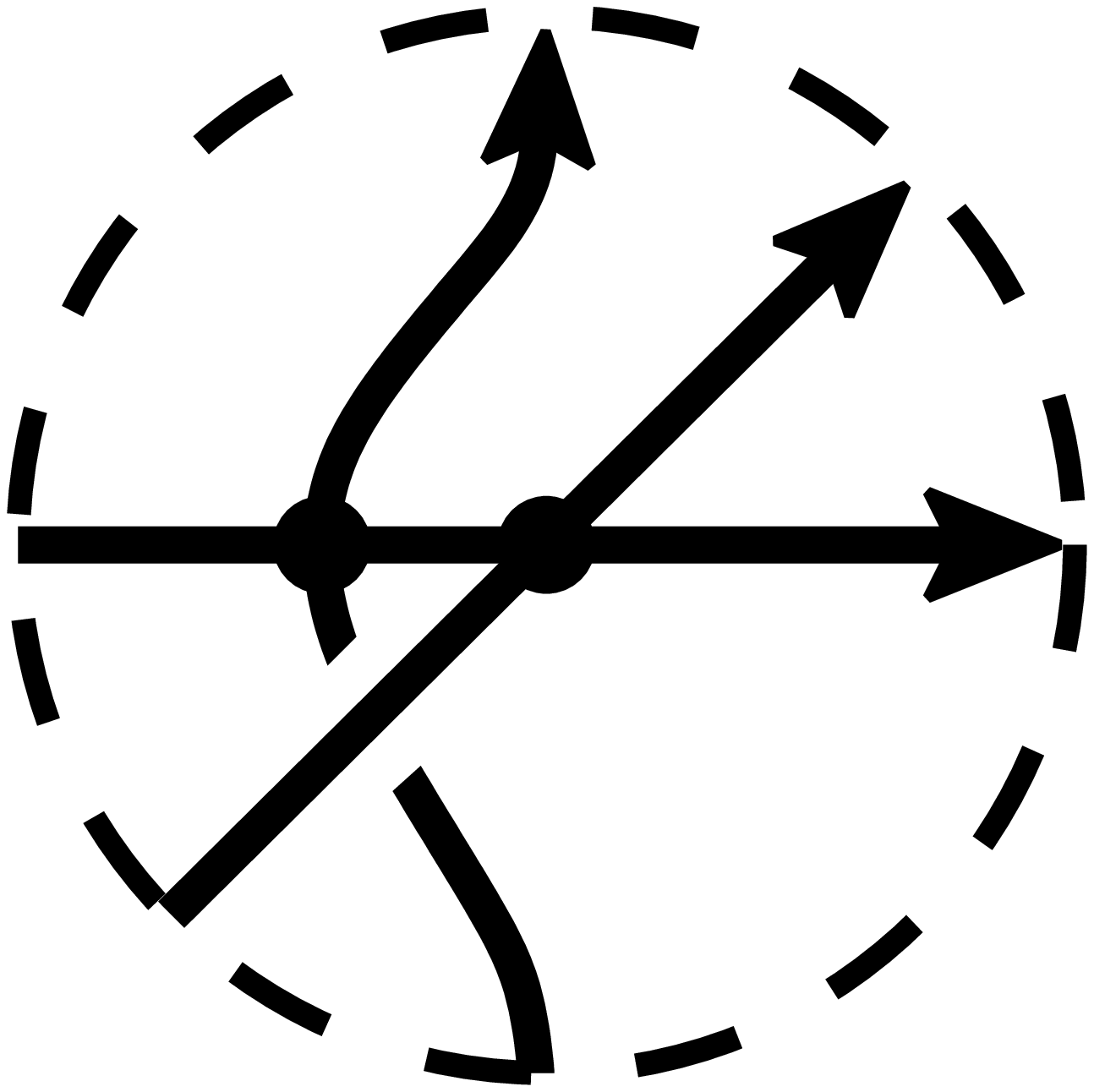}}\bigr)
 + f\bigl(\rb{-13pt}{\ig[width=30pt]{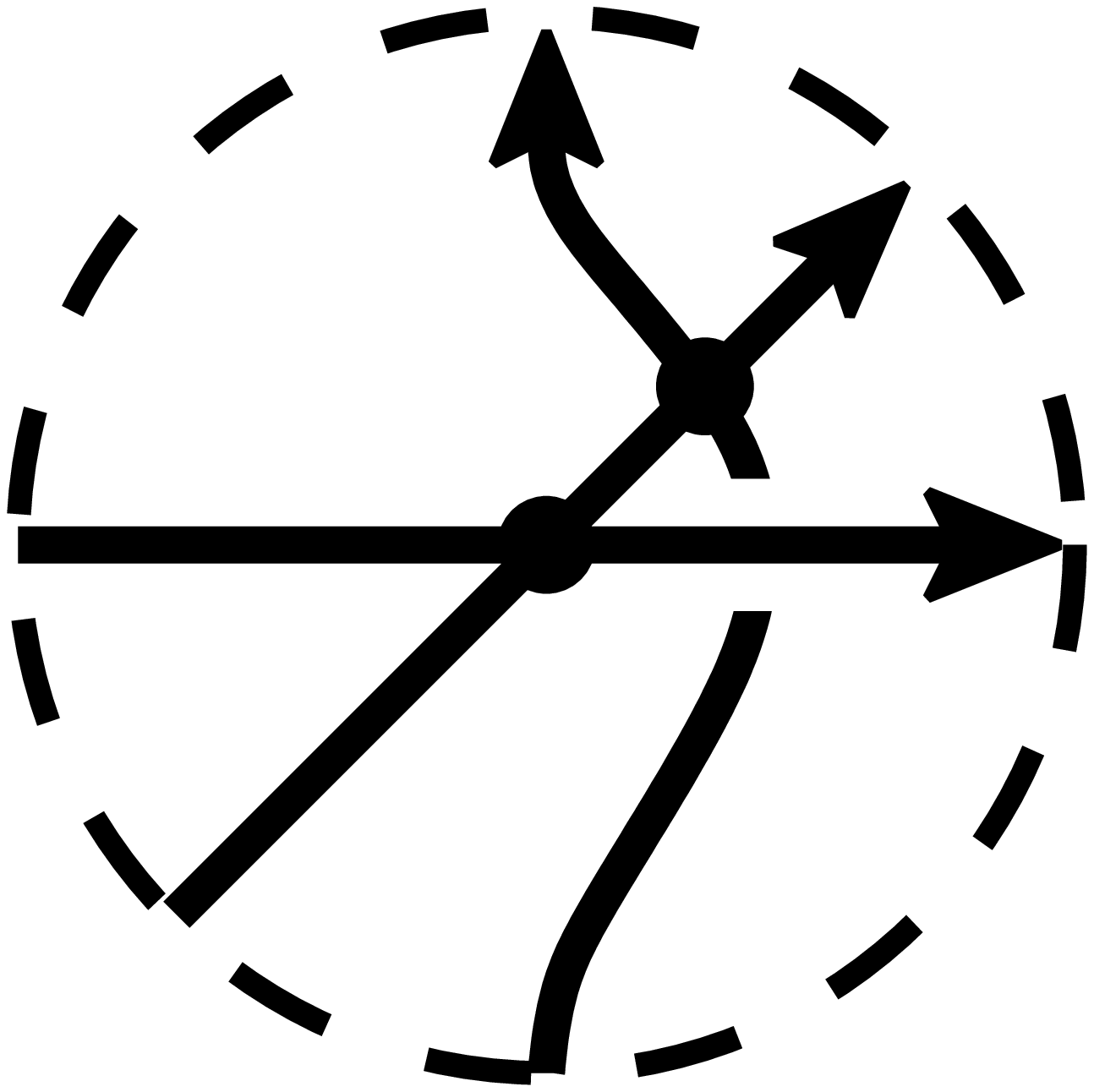}}\bigr)
 - f\bigl(\rb{-13pt}{\ig[width=30pt]{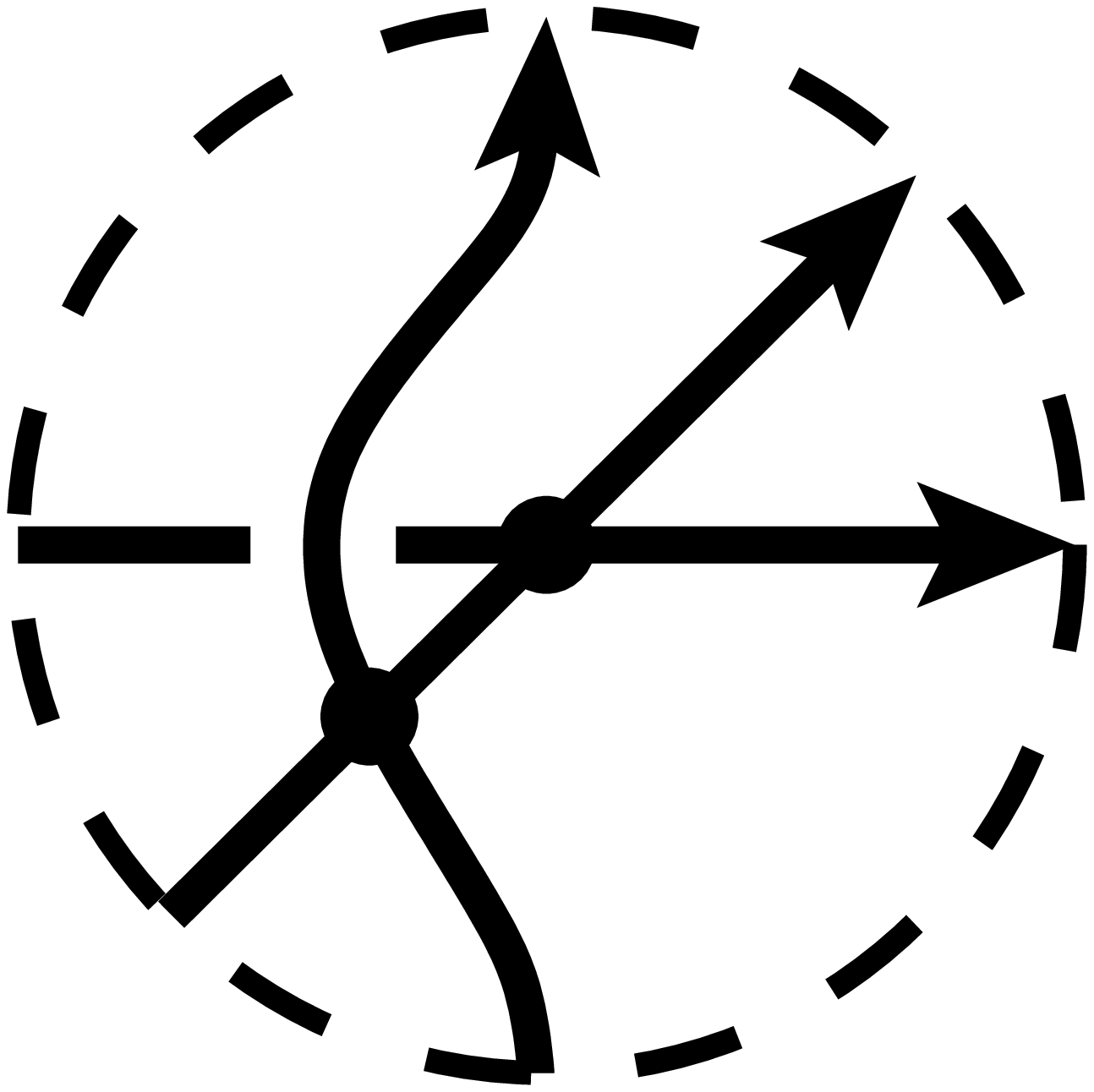}}\bigr)\ =\ 0,
$$
\end{xlemma}

\begin{proof}
By the Vassiliev skein relation,
\begin{eqnarray*}
f\bigl(\rb{-13pt}{\ig[width=30pt]{4T1k.eps}}\bigr) &=&
  f\bigl(\rb{-13pt}{\ig[width=30pt]{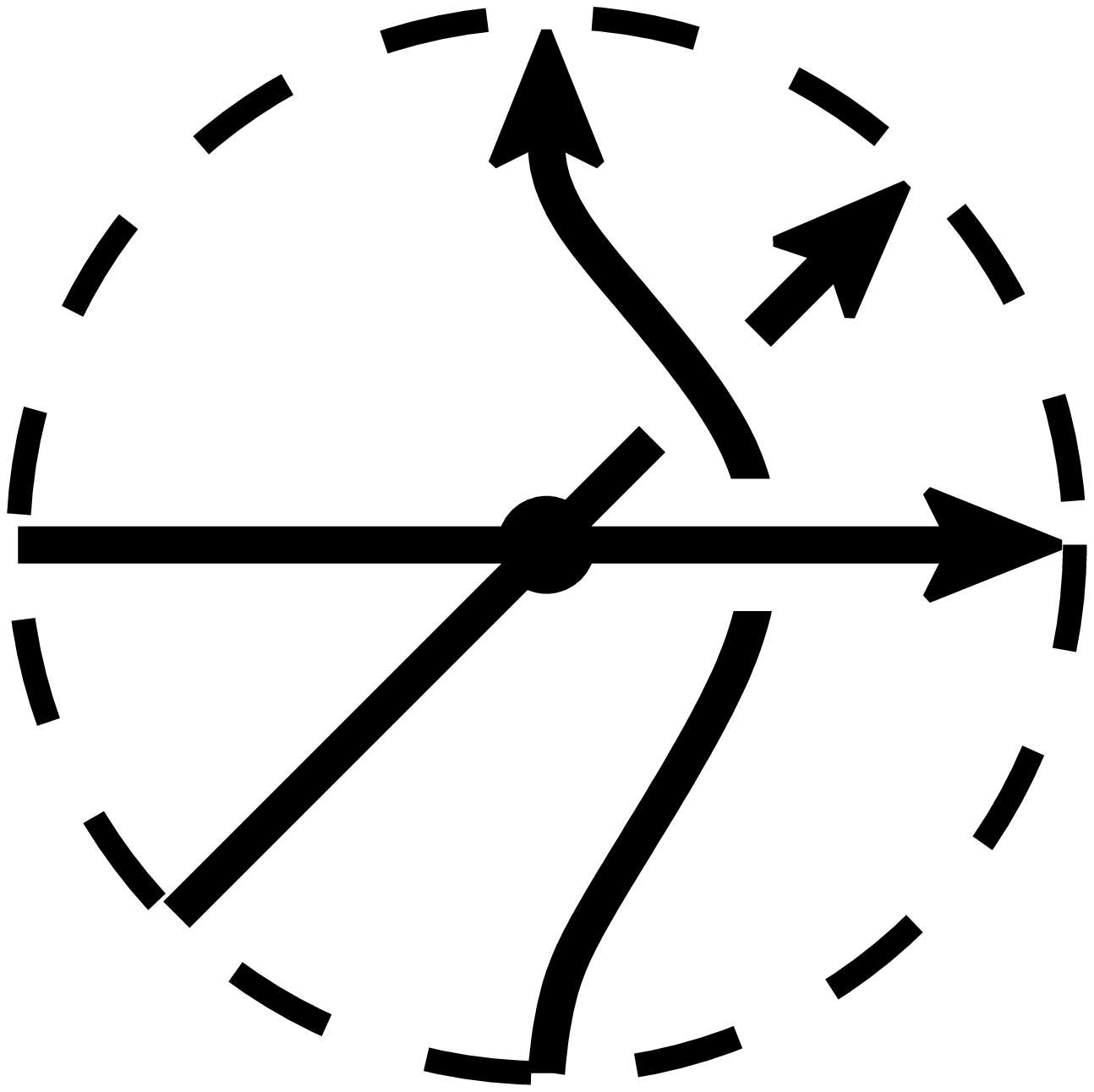}}\bigr)
- f\bigl(\rb{-13pt}{\ig[width=30pt]{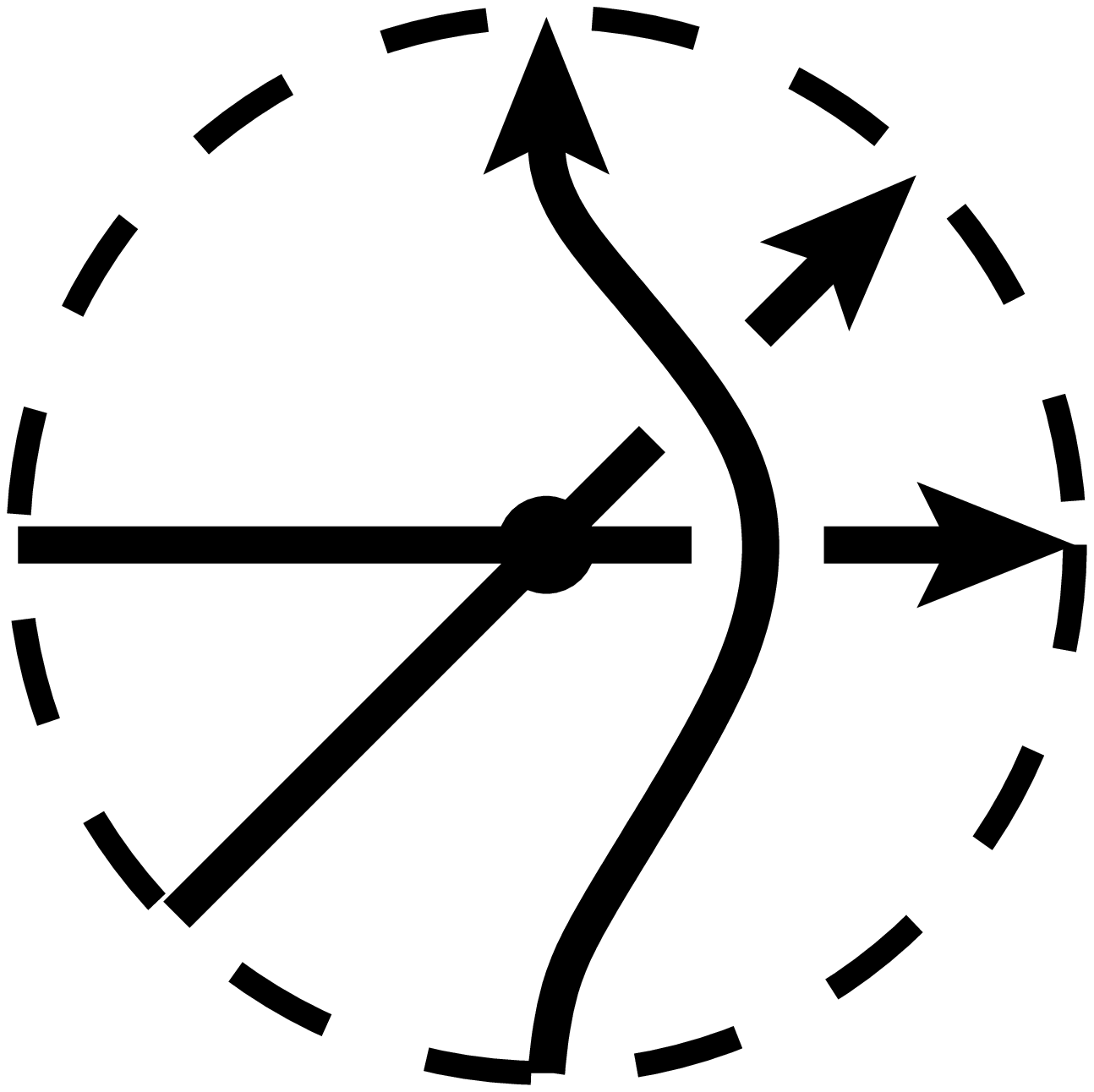}}\bigr)\ =\ a-b, \\
f\bigl(\rb{-13pt}{\ig[width=30pt]{4T2k.eps}}\bigr) &=&
  f\bigl(\rb{-13pt}{\ig[width=30pt]{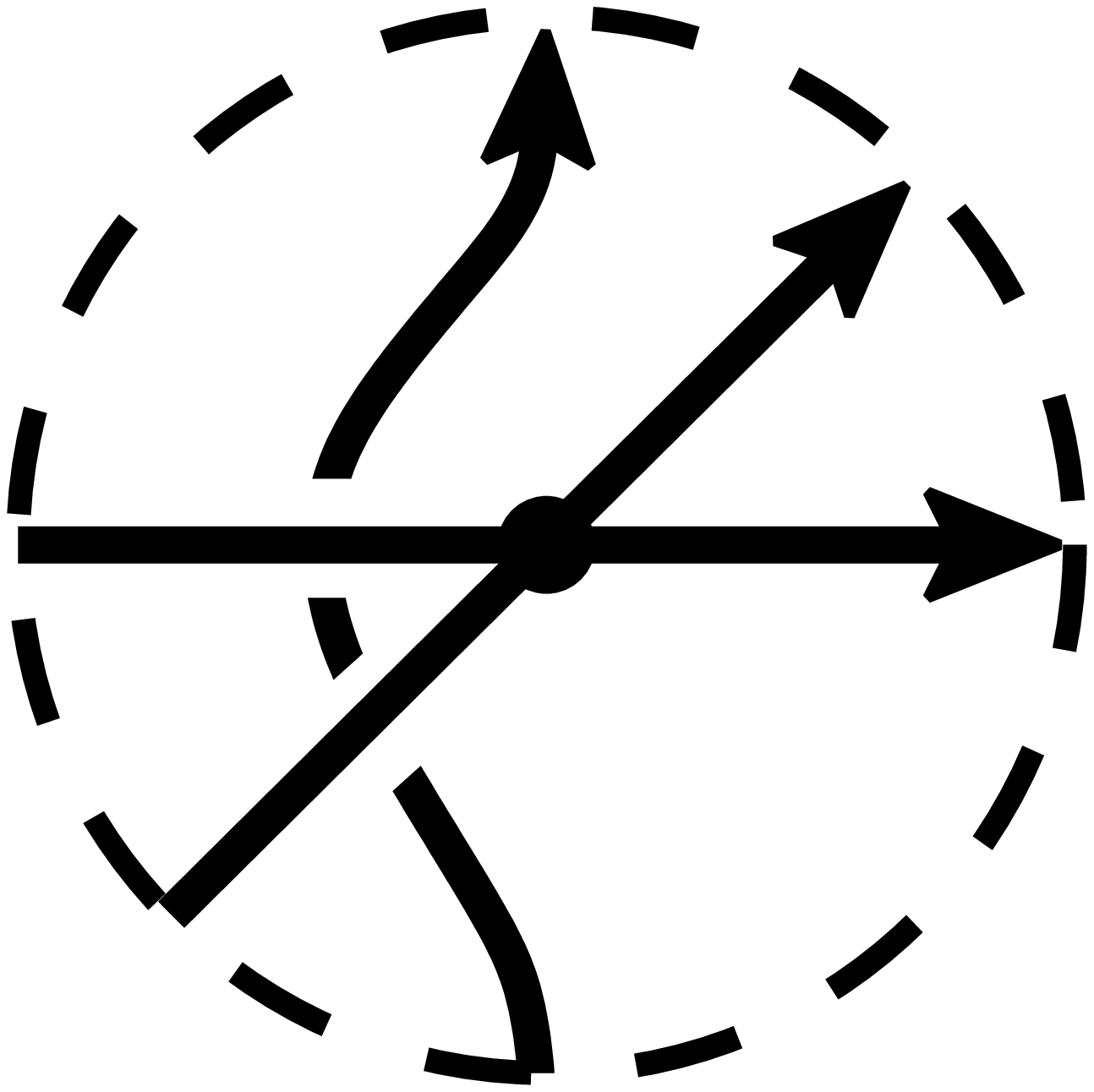}}\bigr)
- f\bigl(\rb{-13pt}{\ig[width=30pt]{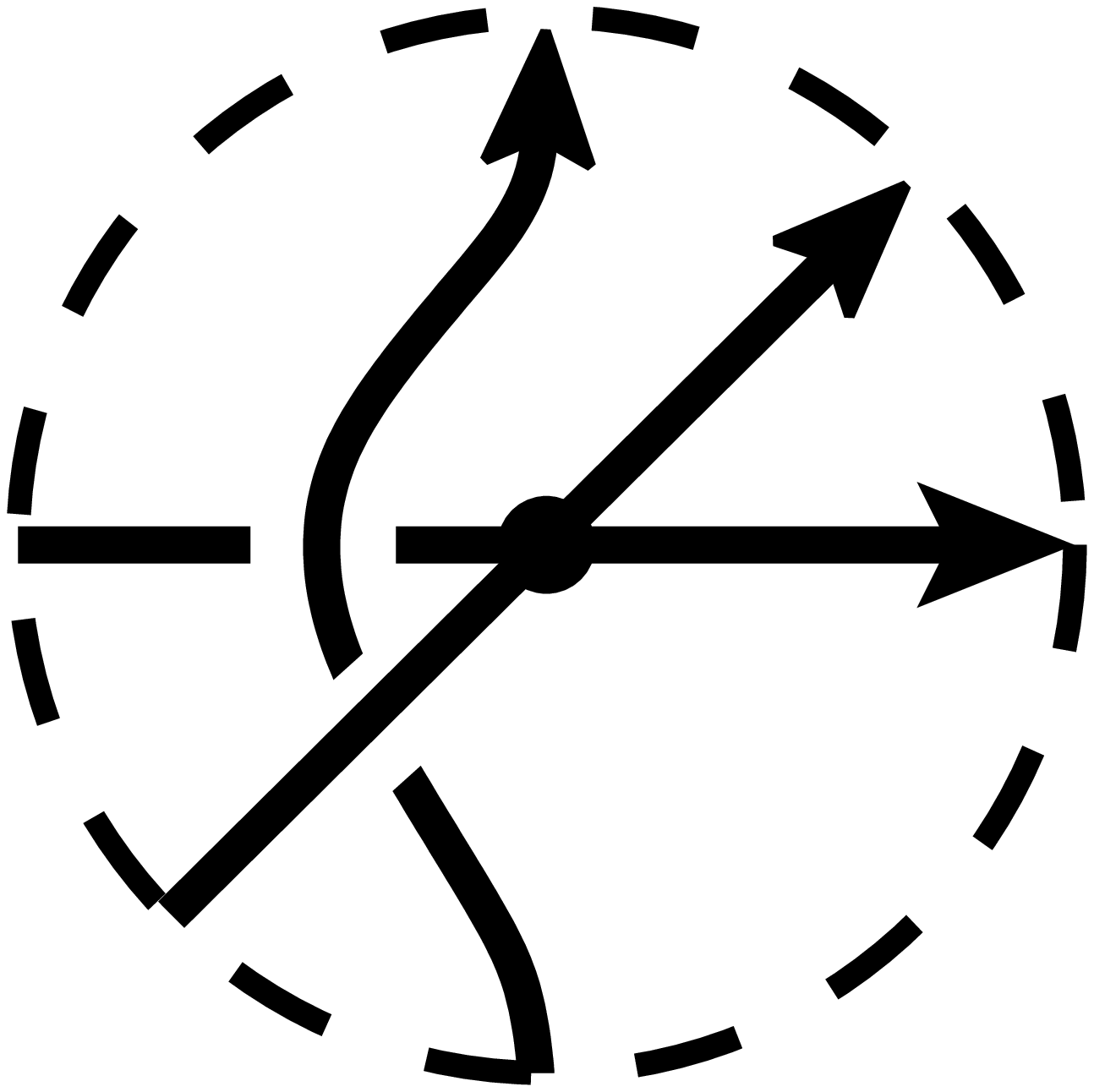}}\bigr)\ =\ c-d, \\
f\bigl(\rb{-13pt}{\ig[width=30pt]{4T3k.eps}}\bigr) &=&
  f\bigl(\rb{-13pt}{\ig[width=30pt]{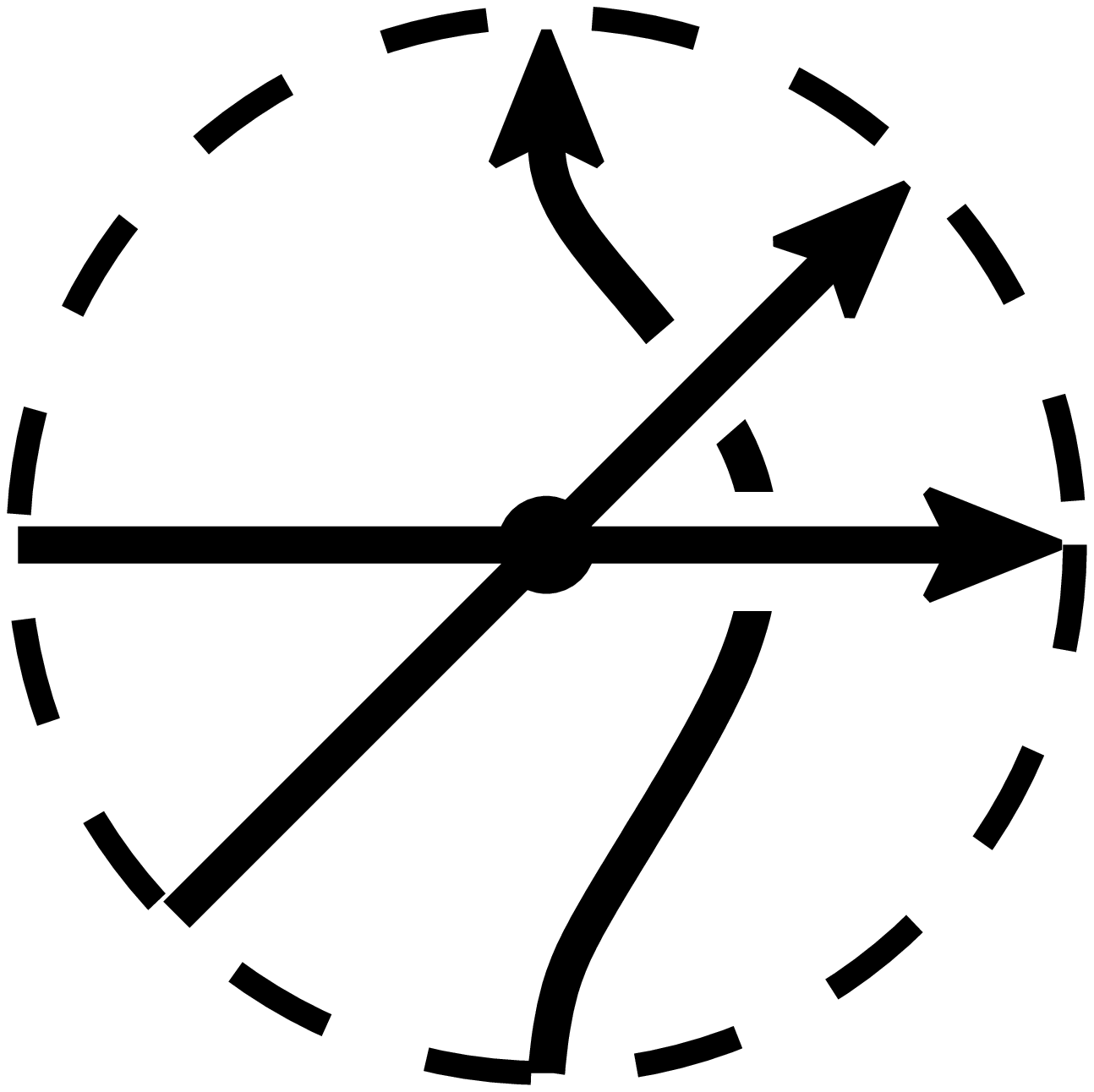}}\bigr)
- f\bigl(\rb{-13pt}{\ig[width=30pt]{4T1a.eps}}\bigr)\ =\ c-a, \\
f\bigl(\rb{-13pt}{\ig[width=30pt]{4T4k.eps}}\bigr) &=&
  f\bigl(\rb{-13pt}{\ig[width=30pt]{4T2d.eps}}\bigr)
- f\bigl(\rb{-13pt}{\ig[width=30pt]{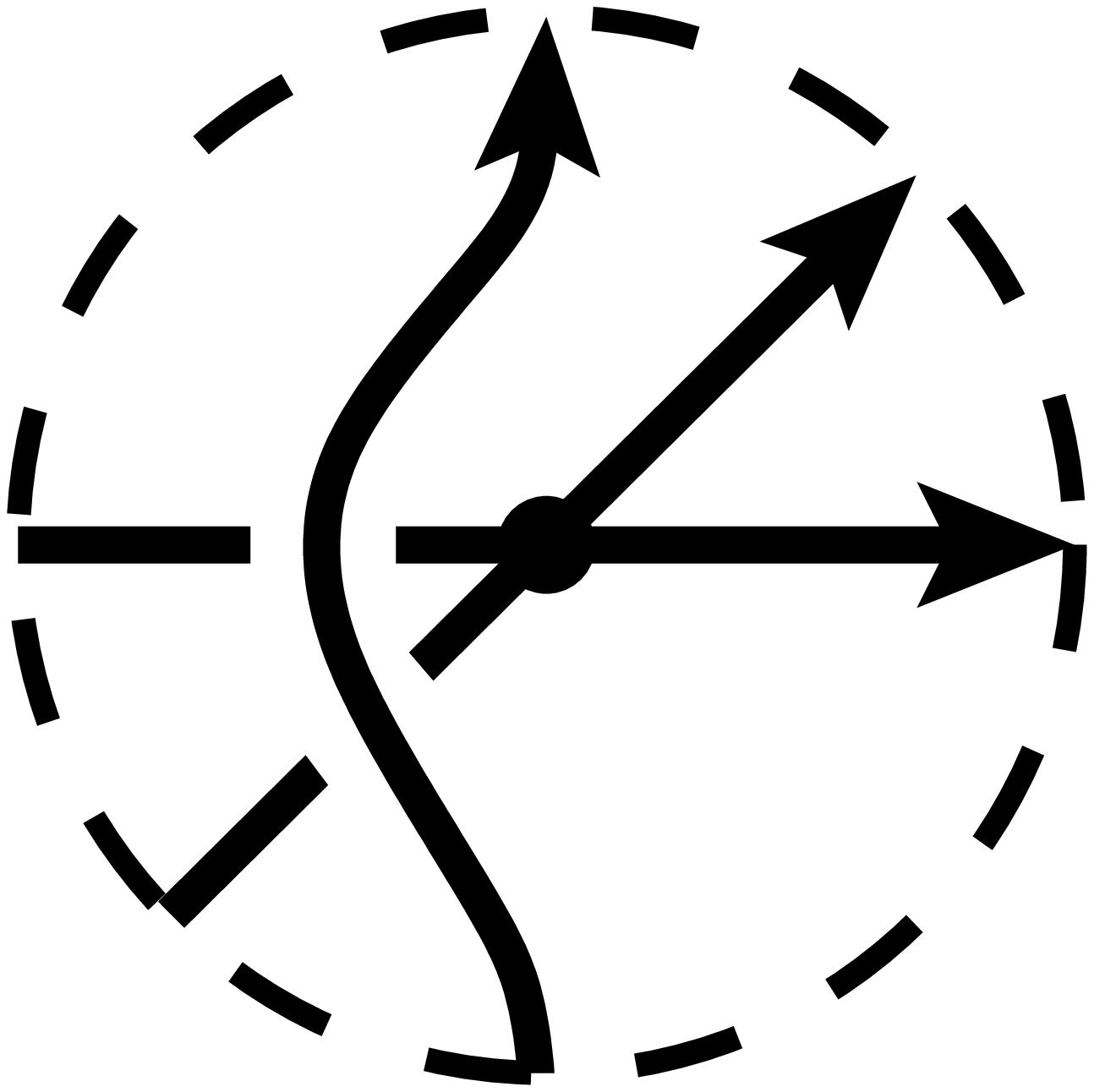}}\bigr)\ =\ d-b.
\end{eqnarray*}
\par The alternating sum of these expressions is
$(a-b)-(c-d)+(c-a)-(d-b)=0$, and the lemma is proved.
\end{proof}

Now, denote by $D_1, \ldots, D_4$ the four diagrams in a 4T
relation. In order to prove the 4-term relation for the symbols of
Vassiliev invariants, let us choose for the first diagram $D_1$ an
arbitrary singular knot $K_1$ such that $\sigma (K_1)=D_1$:
$$
\sigma \Bigl(\ \risS{-22}{4T1kn}{
     \put(28,38){\mbox{$K_1$}}
     \put(95,38){\mbox{$D_1$}}}{45}{20}{20}\ \Bigr)\ =\ \chd{4T2}\ \ .
$$
Then the three remaining knots $K_2$, $K_3$, $K_4$ that participate
in the 4-term relation for knots, correspond to the three
remaining chord diagrams of the 4-term relation for chord diagrams,
and the claim follows from the lemma.
$$
\sigma \Bigl(\risS{-22}{4T2kn}{
     \put(28,38){\mbox{$K_2$}}
     \put(60,38){\mbox{$D_2$}}}{45}{25}{5}\,\Bigr) = \chd{4T1}\ ,\quad
\sigma \Bigl(\risS{-22}{4T3kn}{
     \put(28,38){\mbox{$K_3$}}
     \put(60,38){\mbox{$D_3$}}}{45}{20}{5}\,\Bigr) = \chd{4T3}\ ,\quad
\sigma \Bigl(\risS{-22}{4T4kn}{
     \put(28,38){\mbox{$K_4$}}
     \put(60,38){\mbox{$D_4$}}}{45}{20}{5}\,\Bigr) = \chd{4T4}\ .
$$
\end{proof}

\subsection{The case of framed knots}

As in the case of usual knots, for the invariants of framed knots we
can define a linear map  $\V_n^{fr}/\V^{fr}_{n-1}\to\Ring\ChD_n$.
This map satisfies the 4T relations, but does {\em not} satisfy the
1T relation, since the two knots differing by a crossing change (see
the proof of the first assertion in \ref{first_assertion}), are not
equivalent as framed knots (the two framings differ by 2). 
The Fundamental Theorem also holds, in fact, for framed knots: we
have the equality $$\V^{fr}_n/\V^{fr}_{n-1}=\W^{fr}_n;$$ it can be
proved using the Kontsevich integral for framed knots (see
Section~\ref{frKI}).

This explains why the 1-term relation for the Vassiliev invariants
of (unframed) knots is also called the {\em framing independence
relation}
\index{Framing!independence}
\index{Relation!framing independence}.

\subsection{}
We see that, in a sense, the 4T relations are more fundamental than
the 1T relations. Therefore, in the sequel we shall mainly study
combinatorial structures involving the 4T relations only. In any
case, 1T relations can be added at all times, either by considering
an appropriate subspace or an appropriate quotient space (see
Section~\ref{defram_cd}). This is especially easy to do in terms of
the primitive elements (see page~\pageref{primA}): the problem
reduces to simply leaving out one primitive generator.

\section{Bialgebras of knots and Vassiliev knot invariants}
\label{bialg_knots}

Prerequisites on bialgebras can be found in the Appendix (see page
\pageref{bialg}). In this section it will be assumed that
$\Ring=\Fi$, a field of characteristic zero.

In Section~\ref{alg_inv} we noted that the algebra of knot
invariants $\I$, as a vector space, is dual to the algebra of knots
$\Fi\K=\Z\K\otimes\Fi$. This duality provides the algebras of knot and knot invariants with additional structure. Indeed, the dual map to a product $V\ot V\to V$ on a vector space $V$ is a map $V^*\to (V\ot V)^*$; when $V$ is finite-dimensional it is a {\em coproduct} $V^*\to V^*\ot V^*$.
This observation does not apply to the algebras of knots and knot invariants directly, since they are not finite-dimensional. Nevertheless, the coproduct on the algebra of knots exists and is given by  an explicit formula
$$
  \d(K)=K\ot K
$$
for any knot $K$; by linearity this map extends to the entire space $\Fi\K$. Note that its dual 
is precisely the product in $\I$.\label{knotcoproduct}

\begin{xca}
Show that with this coproduct  $\Fi\K$ is a bialgebra. (For this, define the counit and check the compatibility conditions for the product and the coproduct.)
\end{xca}

The {\em singular knot filtration} $\K_n$ on $\Fi\K$ is obtained
from the singular knot filtration on $\Z\K$ (page~\pageref{sing_knot_filtr}) simply by tensoring it
with the field $\Fi$.

\begin{theorem}\label{filt_alg_knots}
The bialgebra of knots $\Fi\K$ considered with the singular knot
filtration  is a bialgebra with a decreasing filtration \textrm{(}page
\pageref{d-filt-bialg}\textrm{)}.
\end{theorem}

\begin{proof}
There are two assertions to prove:
\begin{enumerate}
\setlength{\itemsep}{1pt plus 1pt minus 1pt}
\item
   If $x\in\K_m$ and $y\in\K_n$, then $xy\in\K_{m+n}$,
\item
   If $x\in\K_n$, then $\d(x)\in\sum\limits_{p+q=n}\K_p\ot\K_q$.
\end{enumerate}

The first assertion was proved in Chapter 3.

In order to prove (2), first let us introduce some additional notation.

Let $K$ be a knot given by a plane diagram with $\ge n$ crossings
out of which exactly $n$ are distinguished and numbered. Consider
the set $\hat{K}$ of $2^n$ knots that may differ from $K$ by
crossing changes at the distinguished points and the vector space
$X_K\subset\Fi\K$ spanned by $\hat{K}$. The group $\Z_2^n$ acts on
the set $\hat{K}$; the action of $i$th generator $s_i$ consists in
the flip of under/overcrossing at the distinguished point number
$i$. We thus obtain a set of $n$ commuting linear operators
$s_i:X_K\to X_K$. Set $\sigma_i=1-s_i$. In these terms, a typical
generator $x$ of $\K_n$ can be written as
$x=(\sigma_1\circ\dots\circ\sigma_n)(K)$. To evaluate $\d(x)$, we
must find the commutator relations between the operators $\delta$
and $\sigma_i$.

\begin{lemma}
$$
  \d\circ\sigma_i = (\sigma_i\ot\id+s_i\ot\sigma_i)\circ\d,
$$
where both the left-hand side and the right-hand side are understood as
linear operators from $X_K$ into $X_K\ot X_K$.
\end{lemma}

\begin{proof}
Just check that the values of both operators on an arbitrary
element of the set $\hat{K}$ are equal.
\end{proof}

A successive application of the lemma yields:
\begin{align*}
  \d\circ\sigma_1\circ\dots\circ\sigma_n
  &=\left(\prod_{i=1}^n (\sigma_i\ot\id+s_i\ot\sigma_i)\right)\circ\d \\
  &=(\sum_{I\subset\{1,\dots,n\}}\prod_{i\in I}\sigma_i\prod_{i\not\in I}s_i
   \ot\prod_{i\not\in I}\sigma_i)\circ\d\,.
\end{align*}
Therefore, an element $x=(\sigma_1\circ\dots\circ\sigma_n)(K)$
satisfies
$$
  \d(x)=\sum_{I\subset\{1,\dots,n\}}
        (\prod_{i\in I}\sigma_i\prod_{i\not\in I}s_i)(K)
        \ot(\prod_{i\not\in I}\sigma_i)(K)
$$
which obviously belongs to $\sum_{p+q=n}\Z\K_p\ot\Z\K_q$.
\end{proof}

\subsection{} In contrast with the knot algebra, the algebra of invariants does not have a 
natural coproduct. The map dual to the product in $\Fi\K$ is given by
$$
  \d(f)(K_1\ot K_2)=f(K_1\#K_2)
$$
for an invariant $f$ and any pair of knots $K_1$ and $K_2$. It sends $\I=(\Fi\K)^*$ to 
$(\Fi\K\ot \Fi\K)^*$ but its image is not contained in $\I\ot \I$.
\begin{xca} 
Find a knot invariant whose image under $\d$ is not  in $\I\ot \I$.
\end{xca}

Even though the map $\d$ is not a coproduct, it becomes one if we restrict our attention to the subalgebra $\V^{\Fi}\subset\I$ consisting of all $\Fi$-valued Vassiliev invariants. 
\begin{proposition}\label{fa_vi}
The algebra of $\Fi$-valued Vassiliev knot invariants $\V^{\Fi}$  is a bialgebra with an
increasing filtration \textrm{(}page
\pageref{i-filt-bialg}\textrm{)}.
\end{proposition}
Indeed, the algebra of $\V^{\Fi}$  is a bialgebra is dual {\em as a filtered bialgebra} to the bialgebra of knots with the singular knot filtration. The filtrations on  $\V^{\Fi}$ and $\Fi\K$ are of finite type by Corollary~\ref{fingen} and, hence, the Proposition follows from 
Theorem~\ref{filt_alg_knots} and Proposition~\ref{dual-bialgebra} on page~\pageref{dual-bialgebra}.

\subsection{} Let us now find all the {\em primitive} and the {\em group-like} elements in
the algebras $\Fi\K$ and $\V^{\Fi}$ (see definitions in Appendix
\ref{ap:pr-gr-el} on page \pageref{ap:pr-gr-el}). 
As for the algebra of knots $\Fi\K$, both structures are quite poor: it follows
from the definitions that $\PR(\Fi\K)=0$, while $\GR(\Fi\K)$
consists of only one element: the trivial knot. (Non-trivial
knots are semigroup-like, but not group-like!)

The case of the algebra of Vassiliev invariants is more interesting. As a
consequence of Proposition~\ref{prim_grlike_in_dual} we obtain a
description of primitive and group-like Vassiliev knot invariants:
\index{Vassiliev!invariant!primitive}
\index{Vassiliev!invariant!group-like}
these
are nothing but the {\em additive} and the {\em multiplicative}
invariants, respectively, that is, the invariants satisfying the
relations
\begin{align*}
  & f(K_1\#K_2)=f(K_1)+f(K_2),\\
  & f(K_1\#K_2)=f(K_1)f(K_2),
\end{align*}\label{add-mult-prop}
respectively, for any two knots $K_1$ and $K_2$.

As in the case of the knot algebra, the group-like elements of $\V^{\Fi}$ are scarce:
\begin{xca}\label{grlike-vi}
Show that the only group-like Vassiliev invariant is the constant 1.
\end{xca}
In contrast, we shall see that primitive Vassiliev invariants abound. 

\subsection{}
The bialgebra structure of the Vassiliev invariants extends naturally to the 
power series Vassiliev invariants term by term.  In this framework, there are many more group-like invariants.
\begin{xexample}
According to Exercise~\ref{sum_conway}
to Chapter~\ref{kn_inv}, the Conway polynomial is a group-like power series Vassiliev 
invariant. Taking its logarithm one obtains a primitive power series Vassiliev
invariant. For example, the coefficient $c_2$ (the Casson
invariant) is primitive.
\end{xexample}
\begin{xca}
Find a finite linear combination of coefficients $j_n$ of the Jones
polynomial that gives a primitive Vassiliev invariant.
\end{xca}

\section{Bialgebra of chord diagrams}
\label{bialgCD}
\index{Bialgebra!of chord diagrams}

\subsection{The vector space of chord diagrams}
\label{vscd}

A dual way to define the weight systems is to introduce the 1- and 4-term
relations directly in the vector space spanned by chord diagrams.
\medskip

\begin{definition}\index{Vector space!of chord diagrams}\label{A_n^{fr}}
The space  $\A_n^{fr}$ of chord diagrams of order $n$ is the vector
space generated by the set $\ChD_n$ (all diagrams of order $n$)
modulo the subspace spanned by all 4-term linear combinations\vspace{-5pt}
$$\chd{4T1}\ -\ \chd{4T2}\ +\ \chd{4T3}\ -\ \chd{4T4}\ .\vspace{-5pt}
$$
The space $\A_n$ of \textit{unframed} chord diagrams
\index{Vector space!of unframed chord diagrams}\label{A_n}
\index{Unframed chord diagrams}
of order $n$ is the quotient of $\A_n^{fr}$ by the subspace spanned
by all diagrams with an isolated chord.
\end{definition}

In these terms, the space of framed weight systems $\W_n^{fr}$ is
dual to the space of framed chord diagrams $\A_n^{fr}$, and the
space of unframed weight systems $\W_n$ --- to that of unframed
chord diagrams $\A_n$:
\begin{eqnarray*}
  \W_n &=& \mathop{\rm Hom}(\A_n,\Ring),\\
  \W_n^{fr} &=& \mathop{\rm Hom}(\A_n^{fr},\Ring).
\end{eqnarray*}

\medskip

Below, we list the dimensions and some bases of the spaces
$\A_n^{fr}$ for $n=1$, 2 and 3:

$\A_1^{fr}=\bigl<\chd{cd1ch4}\bigr>$, $\dim\A^{fr}_1=1$.

$\A_2^{fr}=\bigl<\cld{cd21ch4},\cld{cd22ch4}\bigr>$,
$\dim\A^{fr}_2=2$, since the only 4-term relation involving chord
diagrams of order 2 is trivial.

$\A_3^{fr}=\bigl<\chd{cd32ch4},\chd{cd33ch4},\chd{cd34ch4}\bigr>$,
$\dim\A^{fr}_3=3$, since $\ChD_3$ consists of 5 elements, and there
are 2 independent 4-term relations (see
page~\pageref{cd-rel-3}):\vspace{-10pt}
$$\chd{cd31ch4}\ =\ \chd{cd32ch4}\qquad\mbox{and}\qquad
\chd{cd35ch4}\ -\ 2\,\chd{cd34ch4}\ +\ \chd{cd33ch4}\ =\ 0.\vspace{-8pt}
$$

Taking into account the 1-term relations, we get the following
result for the spaces of unframed chord diagrams of small orders:

$\A_1=0$, $\dim\A_1=0$.

$\A_2=\bigl<\cld{cd22ch4}\bigr>$, $\dim\A_2=1$.

$\A_3=\bigl<\chd{cd34ch4}\bigr>$, $\dim\A_3=1$.
\bigskip

The result of similar calculations for order 4 diagrams is presented
in Table \ref{cd4tab}. 
In this case $\dim\A^{fr}_4=6$; the set
$\{d^4_3,d^4_6,d^4_7,d^4_{15},d^4_{17},d^4_{18}\}$ is used in the
table as a basis. The table is obtained by running Bar-Natan's
computer program available at \cite{BN5}. The numerical notation for
chord diagrams like $[12314324]$ is easy to understand: one writes
the numbers on the circle in the positive direction and connects
equal numbers by chords. Of all possible codes we choose the
lexicographically minimal one.

\begin{table}[ht]
\begin{tabular}{|c||l||c||l|}
\hline
CD&Code and expansion&CD&Code and expansion \\ \hline\hline
\chd{cd4-01} &
\parbox{2in}{
$d^4_1=[12341234]$\\ $\phantom{d^4_1}=d^4_3+2d^4_6-d^4_7-2d^4_{15}+d^4_{17}$}
     &%
\chd{cd4-02} &
\parbox{2in}{
$d^4_2=[12314324]$\\ $\phantom{d^4_2}=d^4_3-d^4_6+d^4_7$}
     \\ \hline
\chd{cd4-03} &
\parbox{2in}{
$d^4_3=[12314234]$\\ $\phantom{d^4_3}=d^4_3$}
     &%
\chd{cd4-04} &
\parbox{2in}{
$d^4_4=[12134243]$\\ $\phantom{d^4_4}=d^4_6-d^4_7+d^4_{15}$}
     \\ \hline
\chd{cd4-05} &
\parbox{2in}{
$d^4_5=[12134234]$\\ $\phantom{d^4_5}=2d^4_6-d^4_7$}
     &%
\chd{cd4-06} &
\parbox{2in}{
$d^4_6=[12132434]$\\ $\phantom{d^4_6}=d^4_6$}
     \\ \hline
\chd{cd4-07} &
\parbox{2in}{
$d^4_7=[12123434]$\\ $\phantom{d^4_7}=d^4_7$}
     &%
\chd{cd4-08} &
\parbox{2in}{
$d^4_8=[11234432]$\\ $\phantom{d^4_8}=d^4_{18}$}
     \\ \hline
\chd{cd4-09} &
\parbox{2in}{
$d^4_9=[11234342]$\\ $\phantom{d^4_9}=d^4_{17}$}
     &%
\chd{cd4-10}  &
\parbox{2in}{
$d^4_{10}=[11234423]$\\ $\phantom{d^4_{10}}=d^4_{17}$}
     \\ \hline
\chd{cd4-11} &
\parbox{2in}{
$d^4_{11}=[11234324]$\\ $\phantom{d^4_{11}}=d^4_{15}$}
     &%
\chd{cd4-12}  &
\parbox{2in}{
$d^4_{12}=[11234243]$\\ $\phantom{d^4_{12}}=d^4_{15}$}
     \\ \hline
\chd{cd4-13}  &
\parbox{2in}{
$d^4_{13}=[11234234]$\\ $\phantom{d^4_{13}}=2d^4_{15}-d^4_{17}$}
     &%
\chd{cd4-14}  &
\parbox{2in}{
$d^4_{14}=[11232443]$\\ $\phantom{d^4_{14}}=d^4_{17}$}
     \\ \hline
\chd{cd4-15} &
\parbox{2in}{
$d^4_{15}=[11232434]$\\ $\phantom{d^4_{15}}=d^4_{15}$}
     &%
\chd{cd4-16}  &
\parbox{2in}{
$d^4_{16}=[11223443]$\\ $\phantom{d^4_{16}}=d^4_{18}$}
     \\ \hline
\chd{cd4-17} &
\parbox{2in}{
$d^4_{17}=[11223434]$\\ $\phantom{d^4_{17}}=d^4_{17}$}
     &%
\chd{cd4-18} &
\parbox{2in}{
$d^4_{18}=[11223344]$\\ $\phantom{d^4_{18}}=d^4_{18}$}
\\
\hline
\end{tabular}\vspace{5pt}
\label{cd4tab}\caption{Chord diagrams of order 4} \index{Table of!
chord diagrams}
\end{table}

\subsection{Multiplication of chord diagrams}
\label{algCD}

Now we are ready to define the structure of an algebra in the vector
space $\A^{fr} = \Op_{k\ge0} \A^{fr}_k$\label{A^{fr}}\vspace{-5pt}
of chord diagrams.

\begin{xdefinition}\label{mult-in-A}
\index{Chord diagram!product}\index{Product!in $\A$} The {\em
product} of two chord diagrams $D_1$ and $D_2$ is defined by cutting
and glueing the two circles as shown:\vspace{-10pt}
$$\chd{cd34ch4}\ \cdot\ \chd{cd34ch4}\ =\
\rb{-3.9mm}{\ig[width=22mm]{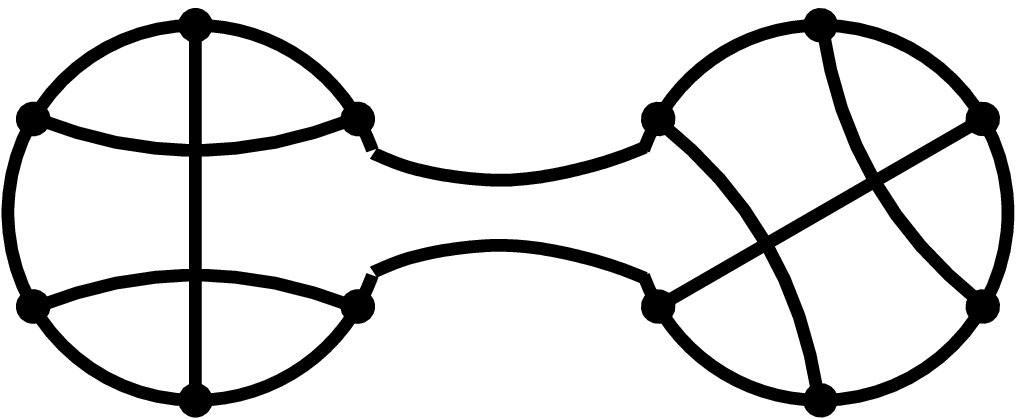}}\ =\ \chd{cd_prod3}\vspace{-8pt}
$$
\end{xdefinition}
This map is then extended by linearity to
$$
  \mu:\A^{fr}_m\otimes \A^{fr}_n \to \A^{fr}_{m+n}.
$$
Note that the product of diagrams depends on the choice of the
points where the diagrams are cut: in the example above we could
equally well cut the circles in other places and get a different
result: $\chd{prod_other}$\ .

\begin{xlemma}\label{mult_cd}
The product is well-defined modulo 4T relations.
\end{xlemma}

\begin{proof}
We shall show that the product of two diagrams is well-defined; it
follows immediately that this is also true for linear combinations
of diagrams. It is enough to prove that if one of the two diagrams,
say $D_2$, is turned inside the product diagram by one ``click''
with respect to $D_1$, then the result is the same modulo 4T
relations.

Note that such rotation is equivalent to the following
transformation. Pick a chord in $D_2$ with endpoints $a$ and $b$
such that $a$ is adjacent to $D_1$. Then, fixing the endpoint $b$,
move $a$ through the diagram $D_1$. In this process we obtain $2n+1$
diagrams $P_0$, $P_1$, ..., $P_{2n}$, where $n$ is the order of
$D_1$, and we must prove that $P_0\equiv P_{2n}\mod 4T$. Now, it is
not hard to see that the difference $P_0-P_{2n}$ is, in fact, equal
to the sum of all $n$ four-term relations which are obtained by
fixing the endpoint $b$ and all chords of $D_1$, one by one. For
example, if we consider the two products shown above and use the
following notation:
$$
\begin{array}{c@{\qquad}c@{\qquad}c@{\qquad}c@{\qquad}c@{\qquad}c@{\qquad}c}
\risS{-10}{prod0}{\put(15,26){\mbox{$\scriptstyle a$}}
                  \put(23,-2){\mbox{$\scriptstyle b$}}}{25}{18}{15} &
\risS{-10}{prod1}{\put(3,24){\mbox{$\scriptstyle a$}}
                  \put(23,-2){\mbox{$\scriptstyle b$}}}{25}{10}{15} &
\risS{-10}{prod2}{\put(-4,19){\mbox{$\scriptstyle a$}}
                  \put(23,-2){\mbox{$\scriptstyle b$}}}{25}{10}{15} &
\risS{-10}{prod3}{\put(-5,10){\mbox{$\scriptstyle a$}}
                  \put(23,-2){\mbox{$\scriptstyle b$}}}{25}{10}{15} &
\risS{-10}{prod4}{\put(-3,3){\mbox{$\scriptstyle a$}}
                  \put(23,-2){\mbox{$\scriptstyle b$}}}{25}{10}{15} &
\risS{-10}{prod5}{\put(2,-3){\mbox{$\scriptstyle a$}}
                  \put(23,-2){\mbox{$\scriptstyle b$}}}{25}{10}{15} &
\risS{-10}{prod6}{\put(10,-5){\mbox{$\scriptstyle a$}}
                  \put(23,-2){\mbox{$\scriptstyle b$}}}{25}{10}{15} \\
P_0 & P_1 & P_2 & P_3 & P_4 & P_5 & P_6
\end{array}
$$
then we must take the sum of the three linear combinations
\begin{eqnarray*}
   && P_0-P_1+P_2-P_3, \\
   && P_1-P_2+P_4-P_5, \\
   && P_3-P_4+P_5-P_6,
\end{eqnarray*}
and the result is exactly $P_0-P_6$.

\end{proof}

\begin{xxca}
Show that the multiplication of chord diagrams corresponds to the connected sum operation on knots in the following sense: if $K_1$ and $K_2$ are two singular knots and $D_1$ and $D_2$ are their chord diagrams, there exists a singular knot, equal to $K_1\# K_2$ as an element of $\Z\K$, whose diagram is $D_1\cdot D_2$.    
\end{xxca}
In view of this exercise, the product of chord diagrams $D_1$ and $D_2$ is sometimes referred to as their {\em connected sum} and denoted by $D_1\# D_2$.\label{connectedsumdiags}\index{Connected sum!of diagrams}

\subsection{Comultiplication of chord diagrams}
\label{comult_cd}
\index{Chord diagram!coproduct}%
\index{Coproduct!in $\A$} The {\em coproduct}\label{cmult-in-A} in
the algebra $\A^{fr}$
$$
  \delta:\A^{fr}_n \to \bigoplus_{k+l=n}\A^{fr}_k\otimes \A^{fr}_l
$$
is defined as follows. For a diagram $D\in\A^{fr}_n$ we put
$$
  \delta(D) :=\sum_{J\subseteq [D]} D_J\ot D_{\overline J},
$$
the summation taken over all subsets $J$ of the set of chords of
$D$. Here $D_J$ is the diagram consisting of the chords that belong
to $J$ and $\overline J=[D]\setminus J$ is the complementary subset
of chords. To the entire space $\A^{fr}$ the operator $\delta$ is
extended by linearity.

If $D$ is a diagram of order $n$, the total number of summands in
the right-hand side of the definition is $2^n$.
\begin{xexample}
$$\begin{array}{l}
  \d\bigl(\chd{cd34ch4}\bigr) \ =\ \cld{cd0}\ot\chd{cd34ch4} + \chd{cd1ch4}\ot\cld{cd21ch4}
        + \cld{cd3y2}\ot\chd{cd3y13} + \cld{cd3y3}\ot\chd{cd3y12} \\
\hspace{54pt}+ \cld{cd21ch4}\ot\chd{cd1ch4} + \chd{cd3y13}\ot\cld{cd3y2}
          + \chd{cd3y12}\ot\cld{cd3y3} + \chd{cd34ch4}\ot\cld{cd0} \vspace{10pt}\\
\hspace{47pt}
  =\ \cld{cd0}\ot\chd{cd34ch4} + 2\chd{cd1ch4}\ot\cld{cd22ch4} +
           \chd{cd1ch4}\ot\cld{cd21ch4} \\
\hspace{54pt}+ \chd{cd34ch4}\ot\cld{cd0} + 2\cld{cd22ch4}\ot\chd{cd1ch4} +
           \cld{cd21ch4}\ot\chd{cd1ch4}
\end{array}
$$
\end{xexample}

\begin{xlemma}
The coproduct $\delta$ is well-defined modulo 4T relations.
\end{xlemma}

\begin{proof}
Let $D_1-D_2+D_3-D_4=0$ be a 4T relation. We must show that the sum
$\delta(D_1)-\delta(D_2)+\delta(D_3)-\delta(D_4)$ can be written as
a combination of 4T relations. Recall that a specific four-term
relation is determined by the choice of a moving chord $m$ and a
fixed chord $a$. Now, take one and the same splitting $A\cup B$ of
the set of chords in the diagrams $D_i$, the same for each $i$, and
denote by $A_i$, $B_i$ the resulting chord diagrams giving the
contributions $A_i\otimes B_i$ to $\delta(D_i)$, $i=1,2,3,4$.
Suppose that the moving chord $m$ belongs to the subset $A$. Then
$B_1=B_2=B_3=B_4$ and $A_1\otimes B_1-A_2\otimes B_2+A_3\otimes
B_3-A_4\otimes B_4 =(A_1-A_2+A_3-A_4)\otimes B_1$. If the fixed
chord $a$ belongs to $A$, then the $A_1-A_2+A_3-A_4$ is a four-term
combination; otherwise it is easy to see that $A_1=A_2$ and
$A_3=A_4$ for an appropriate numbering. The case when $m\in B$ is
treated similarly.
\end{proof}

The {\em unit} and the {\em counit} in $\A^{fr}$ are defined as
follows:\vspace{-5pt}\index{Unit!in $\A$}\index{Counit!in $\A$}
$$\begin{array}{lcl} \label{u-cu-in-A}
\index{Unit!in $\A$}\index{Counit!in $\A$}
  \io:\Ring\to\A^{fr} &,& \io(x)=x\ \cld{cd0}\ ,\vspace{-10pt}\\
  \e:\A^{fr}\to\Ring  &,& \e\bigl(x\ \cld{cd0}\ +\ ...\bigr)=x\ .\vspace{-10pt}
\end{array}
$$

\begin{xxca} Check the axioms of a bialgebra for
$\A^{fr}$ and verify that it is commutative, cocommutative and
connected.
\end{xxca}

\subsection{Deframing the chord diagrams}
\label{defram_cd} \index{Deframing!of chord diagrams} The space of
unframed chord diagrams $\A$ was defined as the quotient of the
space $\A^{fr}$ by the subspace spanned by all diagrams with an
isolated chord. In terms of the multiplication in $\A^{fr}$, this
subspace can be described as the ideal of $\A^{fr}$ generated by
$\Theta$\label{Theta}, the chord diagram with one chord, so that we
can write:
$$\label{A}
  \A = \A^{fr}/(\Theta).
$$
It turns out that there is a simple explicit formula for a linear
operator $p:\A^{fr}\to\A^{fr}$ whose kernel is the ideal $(\Theta)$.
Namely, define $p_n:\A^{fr}_n\to\A^{fr}_n$ by
$$
  p_n(D):=\sum_{J\subseteq [D]} (-\Theta)^{n-|J|}\cdot D_J\ ,
$$
where, as earlier, $[D]$ stands for the set of chords in the diagram $D$
and $D_J$ means the subdiagram of $D$ with only the chords from $J$ left.
The sum of $p_n$ over all $n$ is the operator $p:\A^{fr}\to\A^{fr}$.

\begin{xca}
Check that
\begin{enumerate}
\setlength{\itemsep}{1pt plus 1pt minus 1pt}
\item
  $p$ is a homomorphism of algebras,
\item
  $p(\Theta)=0$ and hence $p$ takes the entire ideal $(\Theta)$ into 0.
\item
  $p$ is a projector, that is, $p^2=p$.
\item
  the kernel of $p$ is exactly $(\Theta)$.
\end{enumerate}
\end{xca}

We see, therefore, that the quotient map
$\bar{p}:\A^{fr}/(\Theta)\to\A^{fr}$\label{map:bar-p} is the
isomorphism of $\A$ onto its image and we have a direct
decomposition $\A^{fr}=\bar{p}(\A) \op (\Theta)$. Note that the
first summand here is different from the subspace spanned merely by
all diagrams without isolated chords!

For example, $p(\A^{fr}_3)$ is spanned by the two
vectors\vspace{-8pt}
$$\begin{array}{lcl}
  p\bigl(\chd{cd34ch4}\bigr)&=&\chd{cd34ch4} -2\, \chd{cd33ch4} +\chd{cd31ch4}\ ,
                   \vspace{-5pt}\\
  p\bigl(\chd{cd35ch4}\bigr)&=&\chd{cd35ch4} -3\, \chd{cd33ch4} +2\,\chd{cd31ch4}\ = 2p\bigl(\chd{cd34ch4}\bigr).
\end{array}\vspace{-8pt}
$$
while the subspace generated by the
elements $\chd{cd34ch4}$ and $\chd{cd35ch4}$ is 2-dimensional and has a nonzero
intersection with the ideal $(\Theta)$.

\section{Bialgebra of weight systems}
\label{bialgWS}
\index{Bialgebra!of weight systems}

According to \ref{A_n^{fr}} the vector space $\W^{fr}$ is dual to
the space $\A^{fr}$. Since  now $\A^{fr}$ is equipped with the
structure of a Hopf algebra, the general construction of Section
\ref{dual_Ha} supplies the space $\W^{fr}$ with the same structure.
In particular, weight systems can be multiplied: $(w_1\cdot w_2)(D)
:= (w_1\ot w_2)(\d(D))$ \index{Product!in $\W$} and comultiplied:
$(\d(w))(D_1\ot D_2) := w(D_1\cdot D_2)$. \index{Coproduct!in $\W$}
The {\em unit} of $\W^{fr}$ is the weight system $\bo_0$ which takes
\index{Unit!in $\W$} value 1 on the chord diagram without chords and
vanishes elsewhere. The {\em counit} sends a weight system to its
value on on the chord diagram without chords.\index{Counit!in $\W$}

For example, if $w_1$ is a weight system which takes
value $a$ on the chord diagram $\cld{cd22ch4}$\ \vspace{-8pt}, and
zero value on all other chord diagrams, and $w_2$ takes value $b$ on
$$\chd{cd1ch4}$$ and vanishes elsewhere, then\vspace{-10pt}
$$(w_1\cdot w_2)(\chd{cd34ch4}) = (w_1\ot w_2)(\delta(\chd{cd34ch4}))
      = 2w_1(\chd{cd22ch4})\cdot w_2(\chd{cd1ch4}) = 2ab\ .
$$

\begin{proposition}\label{sym-hom}
The symbol $\symb:\V^{fr}\to \W^{fr}$ commutes with multiplication
and comultiplication.
\end{proposition}

\begin{proof}[Proof of the proposition.]
Analyzing the proof of Theorem \ref{sec:grading_on_alg_vi} one can
conclude that for any two Vassiliev invariants of orders $\leq p$
and $\leq q$ the symbol of their product is equal to the product of
their symbols. This implies that the map $\symb$ respects the
multiplication. Now we prove that $\symb(\d(v))=\d(\symb(v))$ for a
Vassiliev invariant $v$ of order $\leq n$. Let us apply both parts
of this equality to the tensor product of two chord diagrams $D_1$
and $D_2$ with the number of chords $p$ and $q$ respectively where
$p+q=n$. We have
$$\symb(\d(v))\bigl(D_1\ot D_2\bigr) =
\d(v)\bigl(K^{D_1}\ot K^{D_2}\bigr) = v\bigl(K^{D_1}\# K^{D_2}\bigr)\ ,
$$
where the singular knots $K^{D_1}$ and $K^{D_2}$ represent chord
diagrams $D_1$ and $D_2$. But the singular knot $K^{D_1}\# K^{D_2}$
represents the chord diagram $D_1\cdot D_2$. Since the total number
of chords in $D_1\cdot D_2$ is equal to $n$, the value of $v$ on the
corresponding singular knot would be equal to the value of its
symbol on the chord diagram:
$$v\bigl(K^{D_1}\# K^{D_2}\bigr) = \symb(v)\bigl(D_1\cdot D_2\bigr)
= \d(\symb(v))\bigl(D_1\ot D_2\bigr)\ .$$
\end{proof}

\begin{xremark} The map $\symb:\V^{fr}\to \W^{fr}$ is not a bialgebra
homomorphism because it does not respect the addition. Indeed, the
sum of two invariants $v_1+v_2$ of different orders $p$ and $q$
with, say $p>q$ has the order $p$. That means
$\symb(v_1+v_2)=\symb(v_1)\ne \symb(v_1)+\symb(v_2)$.

However, we can extend the map $\symb$ to power series Vassiliev
invariants by sending the invariant
$\prod v_i\in \widehat{\V}_{\bullet}^{fr}$
to the element $\sum \symb (v_i)$ of the graded
completion $\widehat{\W}^{fr}$. Then the above Proposition implies
that the map $\symb:\widehat{\V}^{fr}_{\bullet}\to\widehat{\W}^{fr}$ is a
graded bialgebra homomorphism.
\end{xremark}

\subsection{}
We call a weight system $w$ {\em multiplicative} \index{Weight
system! multiplicative} if for any two chord diagrams $D_1$ and
$D_2$ we have $$w(D_1\cdot D_2)=w(D_1) w(D_2).$$ This is the same as
to say that $w$ is a semigroup-like \index{Group-like element} element
in the bialgebra of weight systems (see Appendix \ref{gr_like}).
Note that a multiplicative weight system always takes value 1 on the
chord diagram with no chords. The unit $\bo_0$ is the only
group-like element of the bialgebra $\W^{fr}$ (compare with
Exercise~\ref{grlike-vi} on page \pageref{grlike-vi}). However, the
graded completion $\widehat{\W}^{fr}$\label{grad_compl}  contains many interesting
group-like elements.

 \noindent{\bf Corollary of Proposition
\ref{sym-hom}.} {\it Suppose that $$v=\prod_{n=0}^\infty
v_n\in\widehat{\V}^{fr}_{\bullet}$$ is multiplicative.
Then its symbol is also multiplicative.}

Indeed any homomorphism of bialgebras sends group-like elements to
group-like elements.

\subsection{}
A weight system that belongs to a homogeneous component $\W_n^{fr}$
of the space $\W^{fr}$ is said to be {\em homogeneous of degree
$n$}. \index{Weight system! homogeneous} Let $w\in\widehat{\W^{fr}}$
be an element with homogeneous components $w_i\in\W_i^{fr}$ such
that $w_0=0$. Then the exponential of $w$ can be defined as the
Taylor series
$$
  \exp(w) = \sum_{k=0}^\infty \frac{w^k}{k!}.
$$
This formula makes sense because only a finite number of operations
is required for the evaluation of each homogeneous component of this
sum. One can easily check that the weight systems $\exp(w)$ and
$\exp(-w)$ are inverse to each other:
$$\exp(w)\cdot \exp(-w) = \bo_0.$$

By definition, a {\em primitive} weight system $w$ satisfies
$$w(D_1\cdot D_2) = \bo_0(D_1)\cdot w(D_2) + w(D_1)\cdot\bo_0(D_2).$$
(In particular, a primitive weight system is always zero on a
product of two nontrivial diagrams $D_1\cdot D_2$.) The exponential
$\exp(w)$ of a primitive weight system $w$ is multiplicative
(group-like). Note that it always belongs to the completion
$\widehat{\W}^{fr}$, even if $w$ belongs to $\W^{fr}$.

A simple example of a homogeneous weight system of degree $n$ is
provided by the function on the set of chord diagrams which is equal
to 1 on any diagram of degree $n$ and to 0 on chord diagrams of all
other degrees. This function clearly satisfies the four-term
relations. Let us denote this weight system by $\bo_n$.\label{bo_n}

\begin{lemma} $\bo_n\cdot\bo_m = \binom{m+n}{n}\bo_{n+m}$.
\end{lemma}

This directly follows from the definition of the multiplication for
weight systems.

\begin{corollary}

{\rm(i)} \quad $\frac{\bo_1^n}{n!}=\bo_n$;

{\rm(ii)}\quad If we set\
$\bo=\sum_{n=0}^\infty\bo_n$ (that is, $\bo$ is the weight system
that is equal to 1 on every chord diagram), then
$$\exp({\bo_1})=\bo.$$

\end{corollary}

Strictly speaking, $\bo$ is not an element of $\W^{fr}=\oplus_n
\W_n^{fr}$ but of the graded completion $\widehat{\W}^{fr}$. Note
that $\bo$ is not the unit of $\widehat{\W}^{fr}$. Its unit, as well
as the unit of $\W$ itself, is represented by the element $\bo_0$.

\subsection{Deframing the weight systems}
\label{defram_ws} \index{Deframing!of weight systems} Since
$\A=\A^{fr}/(\Theta)$ is a quotient of $\A^{fr}$, the corresponding
dual spaces are embedded one into another, $\W\subset \W^{fr}$. The
elements of $\W$ take zero values on all chord diagrams with an
isolated chord. In Section~\ref{4Tsec} they were called {\em
unframed weight systems}.  The deframing procedure for chord
diagrams (Section \ref{defram_cd}) leads to a deframing procedure
for weight systems. By duality, the projector $p:\A^{fr}\to\A^{fr}$
gives rise to a projector $p^*:\W^{fr}\to\W^{fr}$ whose value on an
element $w\in\W_n^{fr}$ is defined by
$$
  w'(D)= p^*(w)(D):=w(p(D))
=\sum_{J\subseteq [D]} w\left((-\Theta)^{n-|J|}\cdot D_J\right)\ .
$$
Obviously, $w'(D)=0$ for any $w$ and any chord diagram $D$ with an
isolated chord. Hence the operator $p^*:w\mapsto w'$ is a projection
of the space $\widehat{\W}^{fr}$ onto its subspace $\widehat{\W}$
consisting of unframed weight systems.

The deframing operator looks especially nice for multiplicative
weight systems.

\begin{xca}
Prove that for any number $\theta\in\Fi$ the exponent
$e^{\theta\bo_1}\in \widehat{W}$ is a multiplicative weight system.
\end{xca}

\begin{lemma} \label{defr-mul-ws}
Let $\theta=w(\Theta)$ for a multiplicative weight
system $w$. Then its deframing is $w'=e^{-\theta\bo_1}\cdot w$.
\end{lemma}

We leave the proof of this lemma to the reader as an exercise.
The lemma, together with the previous exercise, implies that the deframing
of a multiplicative weight system is again multiplicative.

\section{Primitive elements in $\A^{fr}$}
\index{Primitive space!in $\A$}
\label{primA}

The algebra of chord diagrams $\A^{fr}$ is commutative,
cocommutative and connected. Therefore, by the Milnor-Moore Theorem
\ref{MMthm}, any element of $\A^{fr}$ is uniquely represented as a
polynomial in basis primitive elements. Let us denote the $n$th
homogeneous component of the primitive subspace by
$\PR_n=\A^{fr}_n\cap \PR(\A^{fr})$ \label{PR_n} and find an explicit
description of $\PR_n$ for small $n$.

\underline{$\dim=1$.} $\PR_1=\A^{fr}_1$ is one-dimensional and
spanned by $\chd{cd1ch4}$\ .

\bigskip
\underline{$\dim=2$.} Since
\begin{eqnarray*}
\d\bigl(\cld{cd22ch4}\bigr)&=&\cld{cd0}\ot\cld{cd22ch4} +
             2\chd{cd1ch4}\ot\chd{cd1ch4} +\cld{cd22ch4}\ot\cld{cd0}\ ,\\
\d\bigl(\cld{cd21ch4}\bigr)&=&\cld{cd0}\ot\cld{cd21ch4} +
             2\chd{cd1ch4}\ot\chd{cd1ch4} +\cld{cd21ch4}\ot\cld{cd0}\ ,
\end{eqnarray*}
the element $\cld{cd22ch4}-\cld{cd21ch4}$ is primitive. It constitutes
a basis of $\PR_2$.

\bigskip
\underline{$\dim=3$.} The coproducts of the 3 basis elements of $\A^{fr}_3$ are
\begin{align*}
\d\bigl(\chd{cd34ch4}\bigr)&=\cld{cd0}\ot\chd{cd34ch4} +
        2\chd{cd1ch4}\ot\cld{cd22ch4} +\chd{cd1ch4}\ot\cld{cd21ch4} +\dots,\\
\d\bigl(\chd{cd33ch4}\bigr)&=\cld{cd0}\ot\chd{cd33ch4} +
         \chd{cd1ch4}\ot\cld{cd22ch4}
       +2\chd{cd1ch4}\ot\cld{cd21ch4} +\dots,\\
\d\bigl(\chd{cd31ch4}\bigr)&=\cld{cd0}\ot\chd{cd31ch4} +
        3\chd{cd1ch4}\ot\cld{cd21ch4} + \dots
\end{align*}
(Here the dots stand for the terms symmetric to the terms that are
shown explicitly.) Looking at these expressions, it is easy to check
that the element $$\chd{cd34ch4}-2\chd{cd33ch4}+\chd{cd31ch4}$$ is
the only, up to multiplication by a scalar, primitive element of
$\A^{fr}_3$.

\bigskip
The exact dimensions of $\PR_n$ are currently (2011) known up to $n=12$
(the last three values, corresponding to
$n=10,11,12$, were found by J.~Kneissler \cite{Kn0}):
\begin{center}
\begin{tabular}{c||c|c|c|c|c|c|c|c|c|c|c|c}
$n$         & 1 & 2 & 3 & 4 & 5 & 6 & 7 & 8  & 9  & 10 & 11 & 12\\
\hline
$\dim \PR_n$ & 1 & 1 & 1 & 2 & 3 & 5 & 8 & 12 & 18 & 27 & 39 & 55
\end{tabular}
\end{center}
We shall discuss the sizes of the spaces $\PR_n$, $\A_n$ and $\V_n$
in more detail later (see Sections~\ref{FDprim} and \ref{bounds}).

If the dimensions of $\PR_n$ were known for all $n$, then the dimensions of
$\A_n$ would also be known.
\medskip

{\bf Example.} Let us find the dimensions of $\A^{fr}_n$, $n\le5$,
assuming that we know the values of $\dim\PR_n$ for $n=1,2,3,4,5$,
which are equal to $1,1,1,2,3$, respectively. Let $p_i$ be the basis
element of $\PR_i$, $i=1,2,3$ and denote the bases of $\PR_4$ and
$\PR_5$ as $p_{41},p_{42}$ and $p_{51},p_{52},p_{53}$, respectively.
Nontrivial monomials up to degree 5 that can be made out of these
basis elements are:

Degree 2 monomials (1): \ $p_1^2$.

Degree 3 monomials (2): \ $p_1^3$, $p_1p_2$.

Degree 4 monomials (4): \ $p_1^4$, $p_1^2p_2$, $p_1p_3$, $p_2^2$.

Degree 5 monomials (7): \ $p_1^5$, $p_1^3p_2$, $p_1^2p_3$, $p_1p_2^2$,
$p_1p_{41}$, $p_1p_{42}$, $p_2p_3$.

A basis of each $\A^{fr}_n$ can be made up of the primitive elements
and their products of the corresponding degree.
For $n=0,1,2,3,4,5$ we get:
$\dim\A^{fr}_0=1$,
$\dim\A^{fr}_1=1$,
$\dim\A^{fr}_2=1+1=2$,
$\dim\A^{fr}_3=1+2=3$,
$\dim\A^{fr}_4=2+4=6$,
$\dim\A^{fr}_5=3+7=10$.

The partial sums of this sequence give the dimensions of the spaces
of framed Vassiliev invariants:
$\dim\V^{fr}_0=1$,
$\dim\V^{fr}_1=2$,
$\dim\V^{fr}_2=4$,
$\dim\V^{fr}_3=7$,
$\dim\V^{fr}_4=13$,
$\dim\V^{fr}_5=23$.
\medskip

\begin{xca}
Let $p_n$ be the sequence of dimensions
of primitive spaces in a Hopf algebra
and $a_n$ the sequence of dimensions
of the entire algebra.
Prove the relation between the generating functions
$$
1+a_1t+a_2t^2+\dots=
\frac{1}{(1-t)^{p_1}(1-t^2)^{p_2}(1-t^3)^{p_3}\dots}.
$$
\end{xca}

Note that primitive elements of $\A^{fr}$ are represented by rather
complicated linear combinations of chord diagrams. A more concise
and clear representation can be obtained via connected closed
diagrams, to be introduced in the next chapter (Section
\ref{FDprim}).

\section{Linear chord diagrams}
\label{line_CD}

The arguments of this chapter, applied to {\em long} knots (see
\ref{long}), lead us naturally to considering the space of {\em
linear} chord diagrams, that is, diagrams on an oriented line:

\begin{center}
\ig[width=50mm]{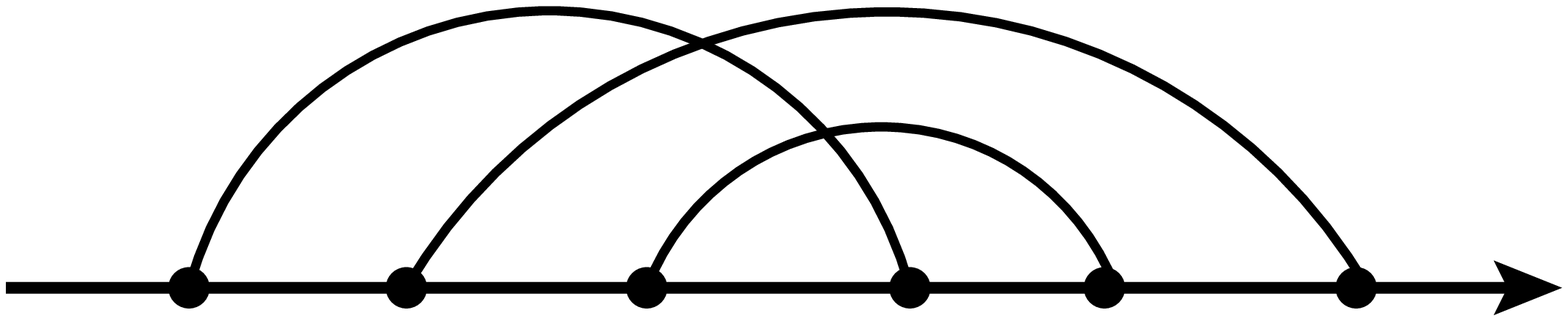}
\end{center}
subject to the 4-term relations:
\begin{align*}
\lcdAa-\lcdBa &=\lcdAb-\lcdBb\\
&=\lcdAc-\lcdBc.
\end{align*}

Let us temporarily denote the space of linear chord diagrams with
$n$ chords modulo the 4-term relations by $(\A^{fr}_n)^{long}$. The
space $(\A^{fr})^{long}$ of such chord diagrams of all degrees
modulo the 4T relations is a bialgebra; the product in
$(\A^{fr})^{long}$ can be defined simply by concatenating the
oriented lines.

If the line is closed into a circle, linear 4-term relations become
circular (that is, usual) 4-term relations; thus, we have a linear
map $(\A^{fr}_n)^{long}\to\A^{fr}_n$. This map is evidently onto, as
one can find a preimage of any circular chord diagram by cutting the
circle at an arbitrary point. This preimage, in general, depends on
the place where the circle is cut, so it may appear that this map
has a non-trivial kernel. For example, the linear diagram shown
above closes up to the same diagram as the one drawn below:
\begin{center}
\ig[width=50mm]{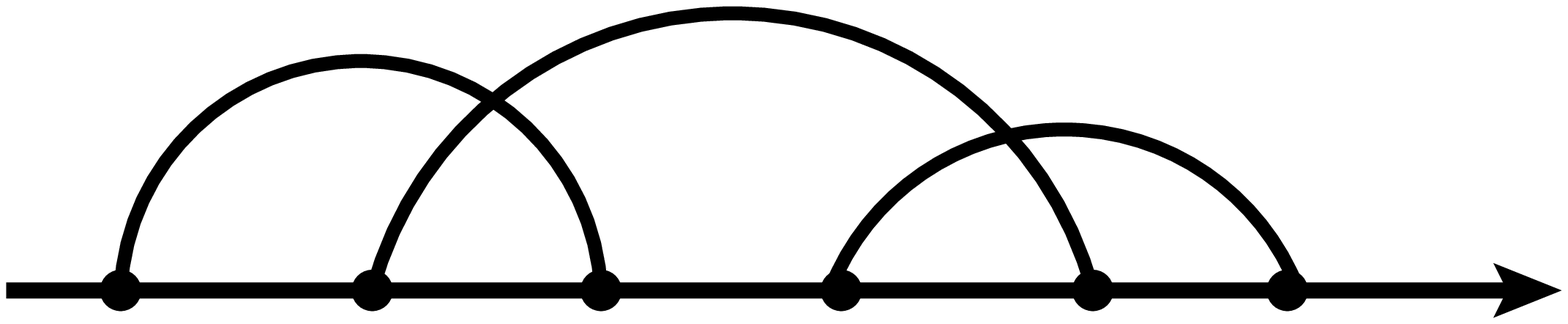}
\end{center}

Remarkably, {\em modulo 4-term relations}, all the preimages of any
circular chord diagram are equal in $(\A^{fr}_3)^{long}$ (in
particular, the two diagrams in the above pictures give the same
element of $(\A^{fr}_3)^{long}$). This fact is proved by exactly the
same argument as the statement that the product of chord diagrams is
well-defined (Lemma \ref{mult_cd}); we leave it to the reader as an
exercise.

Summarizing, we have:

\begin{xproposition}
Closing up the line into the circle gives rise to a vector space
isomorphism $(\A^{fr})^{long}\to\A^{fr}$. This isomorphism is
compatible with the multiplication and comultiplication and thus
defines an isomorphism of bialgebras.
\end{xproposition}

A similar statement holds for diagrams modulo 4T and 1T relations.
Further, one can consider chord diagrams (and 4T relations) with
chords attached to an arbitrary one-dimensional oriented manifold
--- see Section \ref{tangl_Jac}.

\section{Intersection graphs}
\index{Intersection graph}
\label{intgraphs}

\begin{definition} (\cite{CD1}) The {\em intersection graph} $\G(D)$
\label{G(D)}\index{Graph!intersection} of a chord diagram $D$ is the
graph whose vertices correspond to the chords of $D$ and whose edges
are determined by the following rule: two vertices are connected by
an edge if and only if the corresponding chords intersect, and
multiple edges are not allowed. (Two chords, $a$ and $b$, are said
to intersect if their endpoints $a_1$, $a_2$ and $b_1$, $b_2$ appear
in the interlacing order $a_1$, $b_1$, $a_2$, $b_2$ along the
circle.)
\end{definition}

For example,\vspace{-10pt}
\begin{center}\tt
\begin{tabular}{ccc}
\rb{-12mm}{\ig[width=25mm]{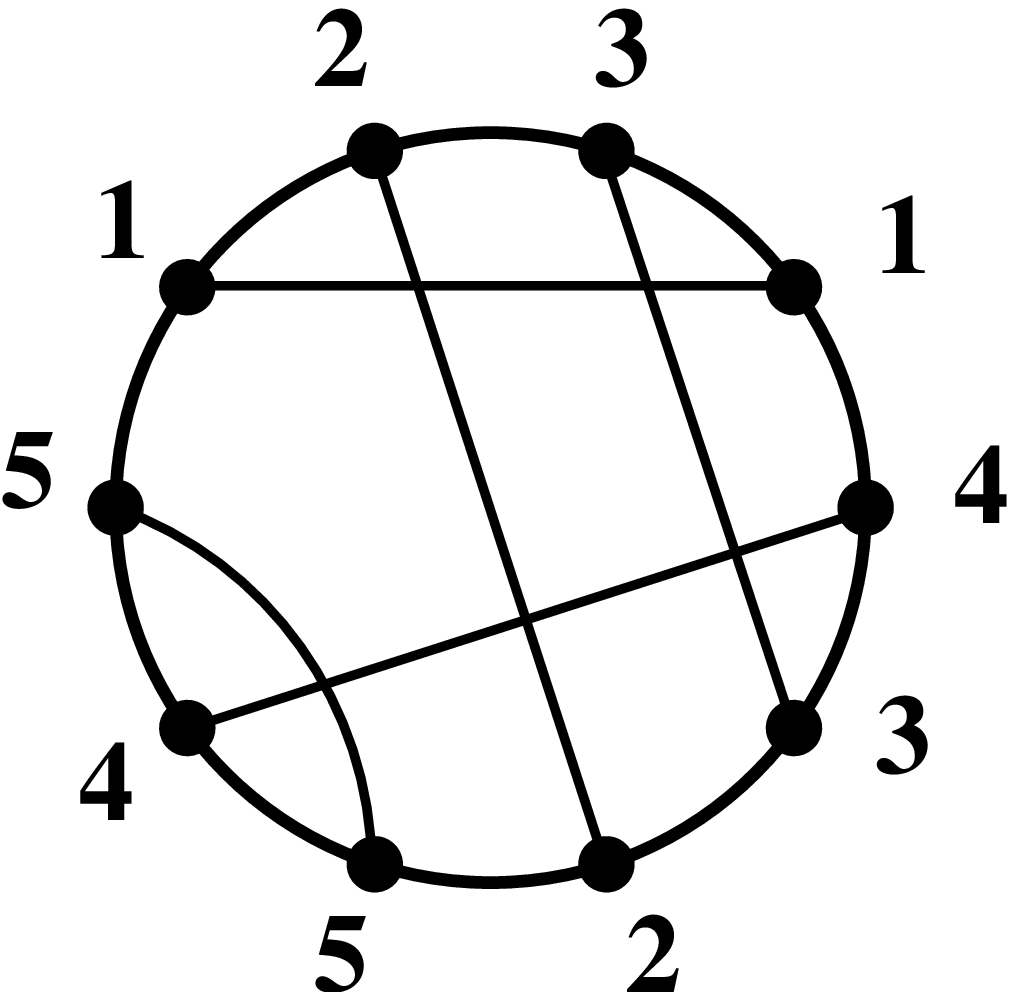}}& $\longrightarrow$&
       \rb{-7.5mm}{\ig[width=25mm]{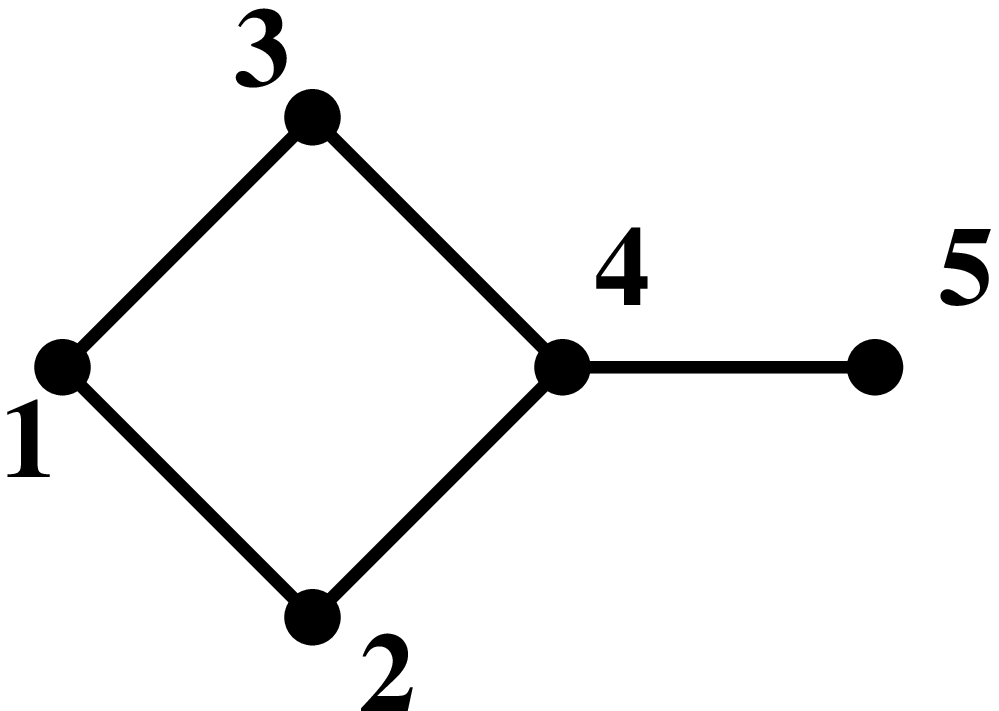}} 
\end{tabular}
\end{center}

The intersection graphs of chord diagrams are also called
{\em circle graphs}\index{Graph!circle}\index{Circle graph} or
{\em alternance graphs}. (See \cite{Bou1}).
\index{Graph!alternance}\index{Alternance graph}

Note that not every graph can be represented as the intersection graph of a
chord diagram. For example, the following graphs are not intersection
graphs: \qquad \chd{pent_c}\ ,\qquad \chd{hex_c}\ ,\qquad \chd{seveng_c}\ .

\subsection{Exercise.}
Prove that all graphs with no more than 5 vertices are intersection
graphs.

\medskip
On the other hand, distinct diagrams may have coinciding intersection graphs.
For example, there are three different diagrams\vspace{-5pt}
$$\chd{cd_A5a}\qquad\qquad \chd{cd_A5b} \qquad\qquad \chd{cd_A5c}\vspace{-5pt}
$$
with the same intersection
graph\qquad $\rb{2pt}{\ig[width=25mm]{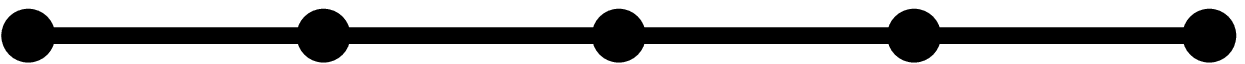}}$\ .

A complete characterization of those graphs that can be realized as
intersection graphs was given by A.~Bouchet \cite{Bou2}.

With each chord diagram $D$ we can associate an oriented surface
$\Sigma_D$ by attaching a disc to the circle of $D$ and thickening
the chords of $D$. Then the chords determine a basis in
$H_1(\Sigma_D,\Z_2)$ as in the picture below. The intersection
matrix for this basis coincides with the adjacency matrix of $\G_D$.
Using the terminology of  singularity theory we may say that the
intersection graph $\G_D$ is the {\em Dynkin diagram} \index{Dynkin
diagram}\index{Diagram!Dynkin} of the intersection form in
$H_1(\Sigma_D,\Z_2)$ constructed for the basis of
$H_1(\Sigma_D,\Z_2)$.
$$D=\risS{-10}{cd34ch4}{}{25}{30}{0}\quad \risS{-2}{totor}{}{20}{0}{0}\quad
\Sigma_D=\risS{-25}{basket}{}{50}{0}{0}\quad \risS{-2}{totor}{}{20}{0}{0}\quad
\G_D=\ \risS{0}{gr_A3}{}{30}{0}{20}\label{basket}
$$

\subsection{Some weight systems}

Intersections graphs are useful, for one thing, because they provide a
simple way to define some weight systems.
We will describe two framed weight systems which depend only on the
intersection graph. The reader is invited to find their deframings, using
the formulae of section \ref{defram_ws}.

1. Let $\nu$ be number of intersections of chords in a chord diagram (or, if
you like, the number of edges in its intersection graph).
This number satisfies the four-term relations and thus descends to a
well-defined mapping $\nu:\A^{fr}\mapsto\Z$.

2. Let $\chi(G)$ be a function of a natural number $n$ which is equal to the
number of ways to colour the vertices of a graph $G$ by $t$ colours (not
necessarily using all the colours) so that the endpoints of any edge are
coloured differently. It is easy to see \cite{Har} that $\chi(G)$ is a  
polynomial in $t$ called the {\it chromatic
polynomial}\label{chrom_pol}\index{Chromatic polynomial} of $G$.
If $D$ is a chord diagram, then keeping the notation $\chi(D)$ for the
chromatic
polynomial of $\Gamma(D)$, one can prove that this function satisfies the 4T   
relations and therefore produces a weight system $\chi:\A^{fr}\to\Z[t]$
(this follows from the deletion--contraction relation
for the chromatic polynomial and relation \ref{fourtlg} in Section
\ref{algLando}). 

\begin{xxca}
Prove that the primitivization (see Exercise \ref{primitiviz_cd} at the end
of the chapter)
of the chord diagram with complete intersection graph provides one non-zero
primitive element of $\A^{fr}$ in each degree thus giving the first
non-trivial lower estimate on the dimensions of the spaces $\A^{fr}$.
\end{xxca}

\bigskip
Intersection graphs contain a good deal of information about chord diagrams.
In \cite{CDL1} the following conjecture was stated.

\subsection{Intersection graph conjecture.}\label{IGC}
\index{Intersection graph!conjecture}
{\it If $D_1$ and $D_2$ are two chord diagrams whose intersection
graphs are equal, $\G(D_1) = \G(D_2)$, then $D_1 = D_2$ as elements
of $\A^{fr}$ {\rm(}that is, modulo four-term relations{\rm)}.}
\medskip

Although wrong in general (see Section~\ref{IGCwrong}), this assertion
is true in some particular situations:

(1) for all diagrams $D_1$, $D_2$ with  up to 10 chords (a direct
computer check \cite{CDL1} up to 8 chords and \cite{Mu} for 9 and 10
chords);

(2) when $\G(D_1) = \G(D_2)$ is a tree (see \cite{CDL2}) or, more
generally, $D_1$, $D_2$ belong to the forest subalgebra (see
\cite{CDL3});

(3) when $\G(D_1)=\G(D_2)$ is a graph with a single loop (see
\cite{Mel1});

(4) for weight systems $w$ coming from standard representations of
Lie algebras $\gl_N$ or $\so_N$. This means that $\G(D_1) = \G(D_2)$
implies $w(D_1) = w(D_2)$; see Chapter \ref{LAWS}, proposition on
page \pageref{IGglNst} and exercise \ref{ex_so_N_intgr} on page
\pageref{ex_so_N_intgr} of the same chapter;

(5) for the universal $\sL_2$ weight system and the weight system
coming from the standard representation of the Lie superalgebra
$\gl(1|1)$ (see \cite{CL}).

In fact, the intersection graph conjecture can be refined to the
following theorem which covers items (4) and (5) above.

\begin{xtheorem}[\cite{CL}] \label{IGTM}  The symbol of a Vassiliev
invariant that does not distinguish mutant knots depends on the
intersection graph only.
\end{xtheorem}

We postpone the discussion of mutant knots, the proof of this
theorem and its converse to Section~\ref{mutation}.

\subsection{Chord diagrams representing a given graph}\label{cd-w-ig}

To describe all chord diagrams representing a given intersection
graph we need the notion of a {\em share} \cite{CDL1,CL}. 
Informally, a {\em share} of a chord diagram is a subset of chords
whose endpoints are separated into at most two parts by the
endpoints of the complementary chords. More formally,

\medskip
\begin{xdefinition}\index{Share} A {\em share} is a part of a
chord diagram consisting of two arcs of the outer circle with the
following property: each chord one of whose ends belongs to these
arcs has both ends on these arcs.
\end{xdefinition}

\medskip
Here are some examples:
$$\risS{-20}{share1}{
       \put(5,-10){\mbox{\scriptsize\tt A share}}}{40}{25}{30}\hspace{2cm}
  \risS{-20}{share2}{
       \put(-3,-10){\mbox{\scriptsize\tt Not a share}}}{40}{20}{20}\hspace{2cm}
  \risS{-20}{share3}{
       \put(2,-10){\mbox{\scriptsize\tt Two shares}}}{48}{20}{20}
$$
The complement of a share also is a share.
The whole chord diagram is its own share whose complement contains no chords.

\medskip
{\bf Definition.}\index{Mutation!of a chord diagram}\index{Chord
diagram!mutation} A {\em mutation of a chord diagram} is another
chord diagram obtained by a flip of a share.

\medskip
For example, three mutations of the share in the first chord diagram above
produce the following  chord diagrams:\label{st-example}
$$\risS{-20}{sh-mut1}{}{40}{25}{20}\hspace{2cm}
  \risS{-20}{sh-mut2}{}{40}{20}{20}\hspace{2cm}
  \risS{-20}{sh-mut3}{}{40}{20}{20}
$$

Obviously, mutations preserve the intersection graphs of chord diagrams.

\begin{xtheorem}\label{cd-mut-theorem} Two chord diagrams have
the same intersection graph if and only if they are related by a
sequence of mutations.
\end{xtheorem}

This theorem is contained implicitly in papers \cite{Bou1,GSH} where
chord diagrams are written as {\em double occurrence words}.

\begin{proof}[Proof of the theorem.]
The proof uses Cunningham's theory of graph decompositions
\cite{Cu}.

A {\em split}\index{Graph!split}\index{Split of a graph} of a
(simple) graph $\Gamma$ is a disjoint bipartition $\{V_1,V_2\}$ of
its set of vertices $V(\Gamma)$ such that each part contains at
least 2 vertices, and with the property that there are subsets
$W_1\subseteq V_1$, $W_2\subseteq V_2$ such that all the edges of
$\Gamma$ connecting $V_1$ with $V_2$ form the complete bipartite
graph $K(W_1,W_2)$ with the parts $W_1$ and $W_2$. Thus for a split
$\{V_1,V_2\}$ the whole graph $\Gamma$ can be represented as a union
of the induced subgraphs $\Gamma(V_1)$ and $\Gamma(V_2)$ linked by a
complete bipartite graph.

Another way to think about splits, which is sometimes more
convenient and which we shall use in the pictures below, is as
follows. Consider two graphs $\Gamma_1$ and $\Gamma_2$ each with a
distinguished vertex $v_1\in V(\Gamma_1)$ and $v_2\in V(\Gamma_2)$,
respectively, called {\em markers}. Construct the new graph $$\Gamma
= \Gamma_1 \boxtimes_{(v_1,v_2)} \Gamma_2$$ whose set of vertices is
$V(\Gamma) = \{V(\Gamma_1)-v_1\}\cup \{V(\Gamma_2)-v_2\}$, and whose
set of edges is
$$\begin{array}{l}
E(\Gamma) =
    \{(v'_1,v''_1)\!\in\!E(\Gamma_1)\!:\!v'_1\not=v_1\not=v''_1 \} \cup
  \{(v'_2,v''_2)\!\in\!E(\Gamma_2)\!:\!v'_2\not=v_2\not=v''_2 \} \vspace{8pt}\\
\hspace{40pt} \cup\
  \{(v'_1,v'_2): (v'_1,v_1)\in E(\Gamma_1)\ \mbox{and}\ (v_2,v'_2)\in E(\Gamma_2)\}\ .
\end{array}$$
Representation of $\Gamma$ as $\Gamma_1 \boxtimes_{(v_1,v_2)}
\Gamma_2$ is called a {\em decomposition} of $\Gamma$, the graphs
$\Gamma_1$ \index{Graph!decomposition}\index{Decomposition of
graphs} and $\Gamma_2$ are called the {\em components} of the
decomposition. The partition $\{V(\Gamma_1)-v_1, V(\Gamma_2)-v_2\}$
is a split of $\Gamma$. Graphs $\Gamma_1$ and $\Gamma_2$ might be
decomposed further giving a finer decomposition of the initial graph
$\Gamma$. Graphically, we represent a decomposition by pictures of
its components where the corresponding markers are connected by a
dashed edge.

A {\em prime} \index{Graph!prime}\index{Prime graph} graph is a
graph with at least three vertices admitting no splits. A
decomposition of a graph is said to be {\em canonical}
\index{Graph!canonical decomposition}\index{Canonical decomposition
of graphs} if the following conditions are satisfied:
\begin{itemize}
\item[(i)] each component is either a prime graph, or a complete
graph $K_n$, or a star $S_n$, which is the tree with a vertex, the
{\em centre}, adjacent to $n$ other vertices; \item[(ii)] no two
components that are complete graphs are neighbours, that is, their
markers are not connected by a dashed edge; \item[(iii)] the markers
of two components that are star graphs connected by a dashed edge
are either both centres or both not centres of their components.
\end{itemize}

W.~H.~Cunningham proved \cite[Theorem 3]{Cu} that each graph with
at least six vertices possesses a unique canonical decomposition.

Let us illustrate the notions introduced above by an example of
canonical decomposition of an intersection graph.
We number the chords and the corresponding vertices in
our graphs, so that the unnumbered vertices are the markers of the
components.
$$\risS{-20}{cd-st-ex}{
    \put(0,-10){\mbox{\scriptsize\tt A chord diagram}}}{55}{35}{40}\hspace{1.5cm}
  \risS{-15}{ig-st-ex}{
    \put(-3,-15){\mbox{\scriptsize\tt The intersection graph}}}{60}{20}{20}\hspace{1.5cm}
  \risS{-20}{candec-st-ex}{
    \put(20,-10){\mbox{\scriptsize\tt The canonical decomposition}}}{140}{20}{20}
$$

The key observation in the proof of the theorem
is that components of the canonical decomposition of any
intersection graph admit a unique representation by chord
diagrams. For a complete graph and star components, this is
obvious. For a prime component, this was proved by A.~Bouchet
\cite[Statement 4.4]{Bou1} (see also \cite[Section 6]{GSH} for an
algorithm finding such a representation for a prime graph).

Now, in order to describe all chord diagrams with a given
intersection graph, we start with a component of its canonical
decomposition. There is only one way to realize the component by a
chord diagram. We draw the chord corresponding to the marker as a
dashed chord and call it the {\it marked chord}. This chord
indicates the places where we must cut the circle removing the
marked chord together with small arcs containing its endpoints. As a
result we obtain a chord diagram on two arcs. Repeating the same
procedure with the next component of the canonical decomposition, we
get another chord diagram on two arcs. We have to glue the arcs of
these two diagrams together in the alternating order. There are four
possibilities to do this, and they
differ by mutations of the share corresponding to one of the two components. This completes the proof of the theorem.
\end{proof}

\medskip
To illustrate the last stage of the proof consider our standard example and take the star
2-3-4 component first and then the triangle component. We get
$$\risS{-12}{re-st-com}{}{160}{25}{8} \ \mbox{and}\
  \risS{-12}{re-tr-com}{}{160}{25}{8}\ .
$$
Because of the symmetry, the four ways of glueing these diagrams
produce only two distinct chord diagrams with a marked chord:
$$\risS{-12}{re-st-tr}{}{100}{15}{8} \qquad\mbox{and}\qquad
  \risS{-12}{re-st-rt}{}{100}{15}{8}\ .
$$
Repeating the same procedure with the marked chord for the last
1-6 component of the canonical decomposition, we get
$$\risS{-12}{re-os-com}{}{160}{20}{10}\ .
$$
Glueing this diagram into the previous two in all possible ways we
get the four mutant chord diagrams from page \pageref{st-example}.

\subsection{2-term relations and the genus of a diagram}
\label{2term}\index{Two-term relations}\index{Relation!two-term}
A 2-term (or \textit{endpoint sliding}) relation for chord diagrams
has the form
\begin{align*}
  \ig[height=10mm]{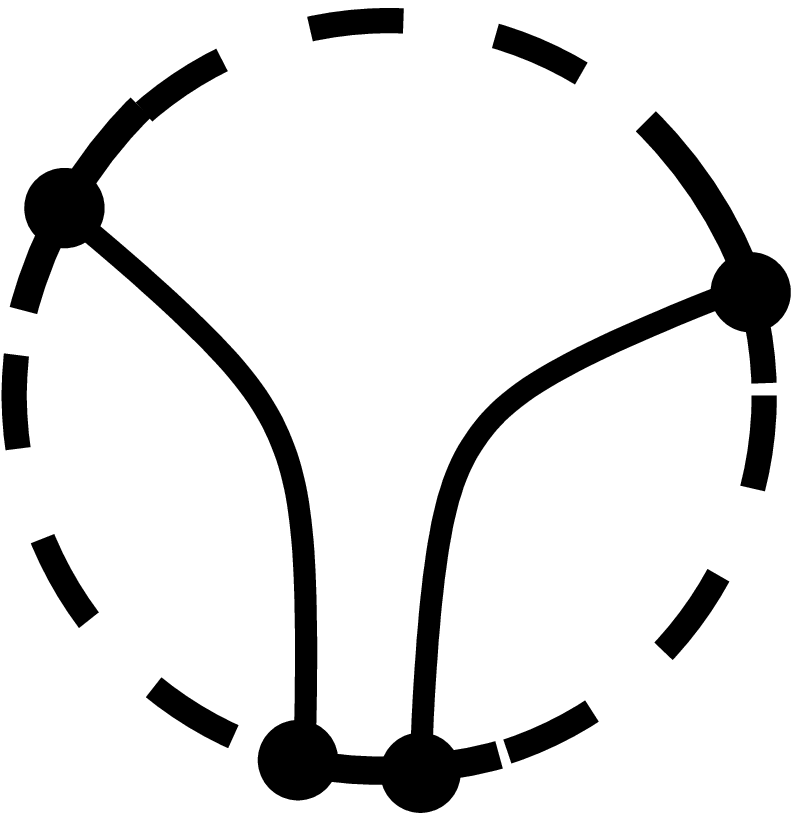}\quad &\rb{5mm}{=}
             \quad \ig[height=10mm]{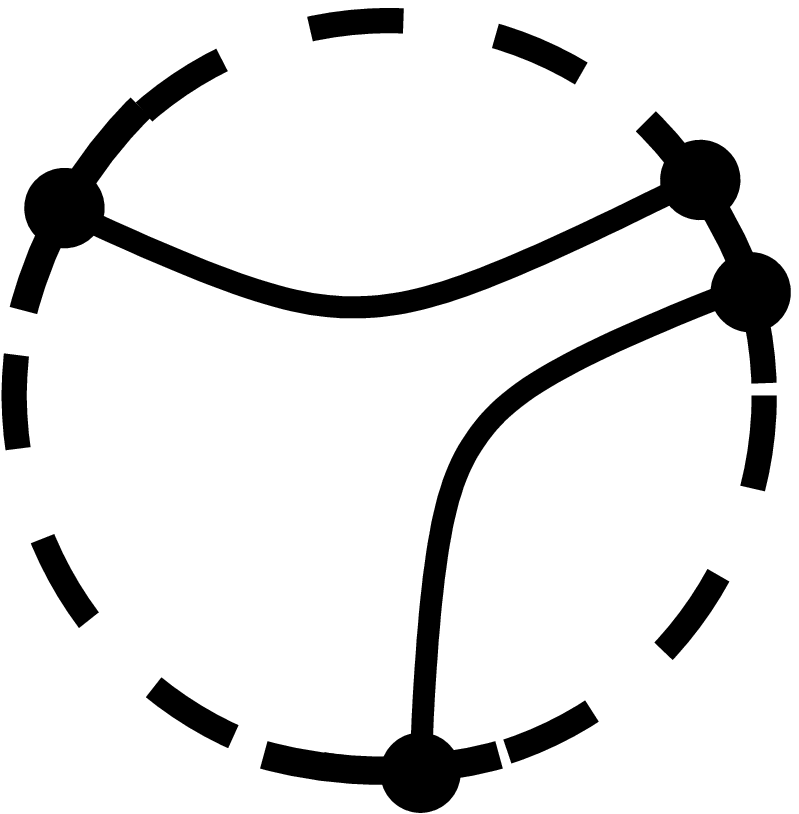}\ \rb{4mm}{,}\\
  \ig[height=10mm]{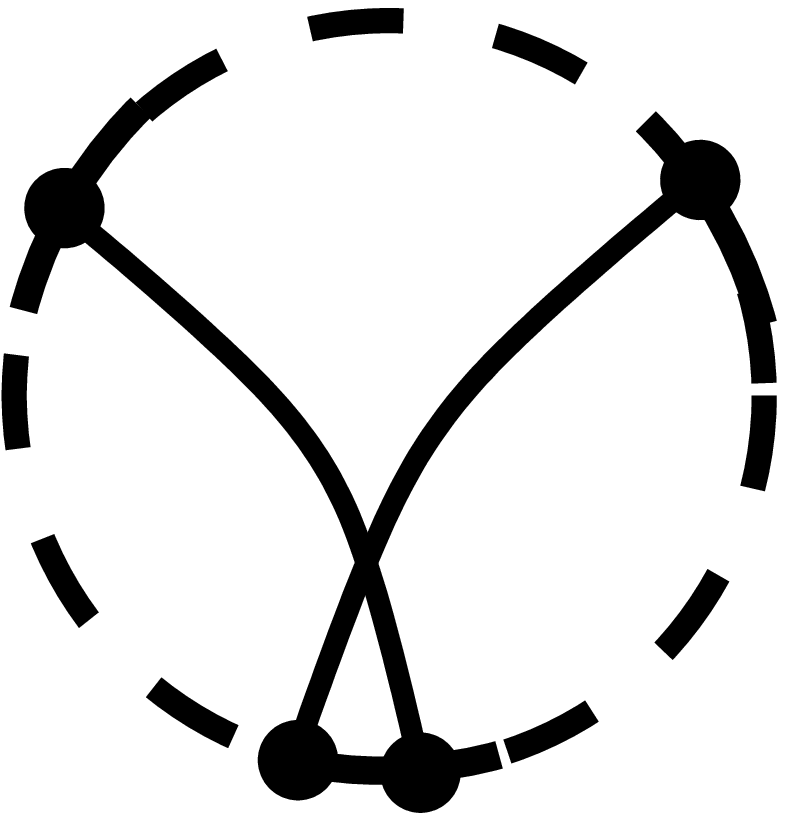}\quad &\rb{5mm}{=}
             \quad \ig[height=10mm]{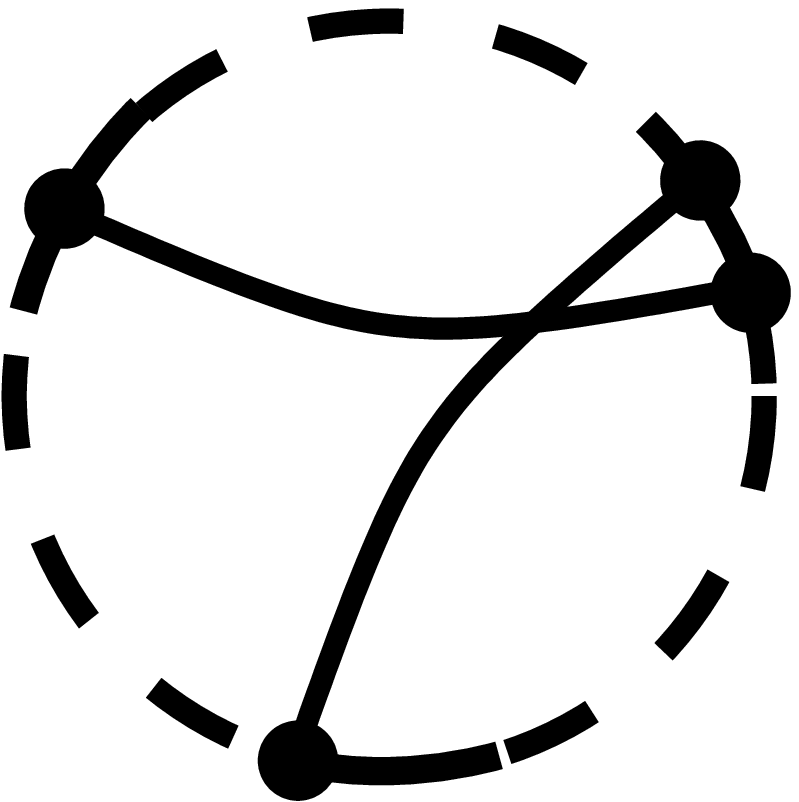}\ \rb{4mm}{.}
\end{align*}

The 4-term relations are evidently a consequence of the 2-term relations;
therefore, any function on chord diagrams that satisfies 2-term
relations is a weight system. An example of such a weight system is
the {\em genus of a chord diagram} defined as follows.

Replace the outer circle of the chord diagram and all its chords by
narrow untwisted bands --- this yields an orientable surface with
boundary. Attaching a disk to each boundary component gives a closed
orientable surface. This is the same as attaching disks to the
boundary components of the surface $\Sigma_D$ from page
\pageref{basket}. Its genus is by definition the genus of the chord
diagram. The genus can be calculated from the number of boundary
components using Euler characteristic. Indeed, the Euler
characteristic of the surface with boundary obtained by above
described procedure from a chord diagram of degree $n$ is equal to
$-n$. If this surface has $c$ boundary components and genus $g$,
then we have $-n = 2-2g-c$ while $g=1+(n-c)/2$. For example, the two
chord diagrams of degree 2 have genera 0 and 1, because the number
of connected components of the boundary is 4 and 2, respectively, as
one can see in the following picture:
$$\chd{cd21ch4}\quad \risS{-2}{totor}{}{20}{0}{0}\quad
\risS{-12}{genera1}{}{30}{0}{0}\qquad\qquad
\chd{cd22ch4}\quad \risS{-2}{totor}{}{20}{0}{0}\quad
\risS{-12}{genera2}{}{30}{0}{0}$$

The genus of a chord diagram satisfies 2-term relations, since
sliding an endpoint of a chord along another adjacent chord does not
change the topological type of the corresponding surface with
boundary.

An interesting way to compute the genus from the intersection graph
of the chord diagram was found by Moran (see \cite{Mor}). Moran's
theorem states that the genus of a chord diagram is half the rank of
the adjacency matrix over $\Z_2$ of the intersection graph. This
theorem can be proved transforming a given chord diagram to the
canonical form using the following two exercises.

\begin{xxca} Let $D_1$ and $D_2$ be two chord diagrams differing by a
2-term relation. Check that the corresponding adjacency matrices
over $\Z_2$ are conjugate (one is obtained from the other by adding
the $i$th column to the $j$th column and the $i$th row to the $j$th
row).
\end{xxca}
\begin{xxca} A {\em caravan of $m_1$ ``one-humped camels'' and
$m_2$ ``two-humped camels} is the product of $m_1$ diagrams with one
chord and $m_2$ diagrams with 2 crossing chords:
$$\risS{-12}{caravan}{\put(18,-8){$m_1$}\put(75,-8){$m_2$}}{114}{0}{20}
$$
Show that any chord diagram is equivalent, modulo 2-term relations,
to a caravan. Show that the caravans form a basis in the vector
space of chord diagrams modulo 2-term relations.
\end{xxca}

The algebra generated by caravans is thus a quotient algebra of the
algebra of chord diagrams. 
\begin{xremark}
The last exercise is, essentially, equivalent to the classical topological classification of  compact oriented surfaces with boundary by the genus, $m_2$, and the number of boundary components,
$m_1+1$.
\end{xremark}

\begin{xcb}{Exercises}

\begin{enumerate}

\item
A {\em short} chord is a chord whose endpoints are adjacent, that
is, one of the arcs that it bounds contains no endpoints of other
chords. In particular, short chords are isolated. Prove that the
linear span of all diagrams with a short chord and all four-term
relation contains all diagrams with an isolated chord. This means
that the restricted one-term relations (only for diagrams with a
short chord) imply general one-term relations provided that the
four-term relations hold.

\item
   Find the number of different chord diagrams of order $n$ with $n$
isolated chords. Prove that all of them are equal to each other modulo
the four-term relations.

\item
Using Table \ref{cd4tab} on page \pageref{cd4tab}, find the space of
unframed weight systems $\W_4$.

{\sl Answer.} The basis weight systems are:
\begin{center}
\begin{tabular}{ccccccc}
\chd{cd4-07} &
\chd{cd4-06} &
\chd{cd4-04} &
\chd{cd4-05} &
\chd{cd4-02} &
\chd{cd4-03} &
\chd{cd4-01} \\
\hline
1 & 0 & -1 & -1 & 0 & -1 & -2 \\
\hline
0 & 1 & 1 & 2 & 0 & 1 & 3 \\
\hline
0 & 0 & 0 & 0 & 1 & 1 & 1 \\
\hline
\end{tabular}
\end{center}
The table shows that the three diagrams \ \chd{cd4-02}, \chd{cd4-06}
and \chd{cd4-07} \
form a basis in the space $\A_4$.

\item\neresh
Is it true that any chord diagram of order 13 is equivalent to its mirror
image modulo 4-term relations?

\item
Prove that the deframing operator $'$ (Section~\ref{defram_ws}) is a
homomorphism of algebras: $(w_1\cdot w_2)'=w'_1\cdot w'_2$.

\item
Give a proof of Lemma \ref{defr-mul-ws} on page \pageref{defr-mul-ws}.

\item
Find a basis in the primitive space $\PR_4$.

{\sl Answer.}
A possible basis consists of the elements $d^4_6-d^4_7$ and $d^4_2-2d^4_7$
from the table on page \pageref{cd4tab}.

\item\label{ex_defr_prim}
Prove that for any primitive element $p$ of degree $>1$,
$w(p)=w'(p)$ where $w'$ is the deframing of a weight system $w$.

\item
Prove that the symbol of a primitive Vassiliev invariant is a primitive
weight system.

\item\label{primitiviz_cd}
Prove that the projection onto the space of the primitive 
elements (see Exercise~\ref{primitivization} on page~\pageref{primitivization})
in the algebra $\A^{fr}$ can be given by the following explicit formula:
$$
\pi(D)=D-1!\sum D_1D_2+2!\sum D_1D_2D_3-...,
$$
where the sum is taken over all unordered splittings of the set of chords
of $D$ into 2, 3, etc nonempty subsets.

\item
   Let $\Theta$ be the chord diagram with a single chord. By a direct
computation, check that
$\exp(\alpha\Theta) :=
  \sum_{n=0}^\infty \frac{\alpha^n\Theta^n}{n!}\in \widehat{\A}^{fr}$
is a group-like element in the completed Hopf algebra of chord diagrams.

\item
(a) Prove that no chord diagram is equal to 0 modulo 4-term
relations.

\noindent
(b) Let $D$ be a chord diagram without isolated chords. Prove that $D\ne0$
modulo 1- and 4-term relations.

\item
Let $c(D)$ be the number of chord intersections in a chord diagram $D$.
Check that $c$ is a weight system. Find its deframing $c'$.

\item {\bf The generalized 4-term relations.}\\
\noindent(a) Prove the following relation:
$$\index{Four-term relation!generalized}\index{Generalized 4-term relation}
  \rb{-7mm}{\ig[height=15mm]{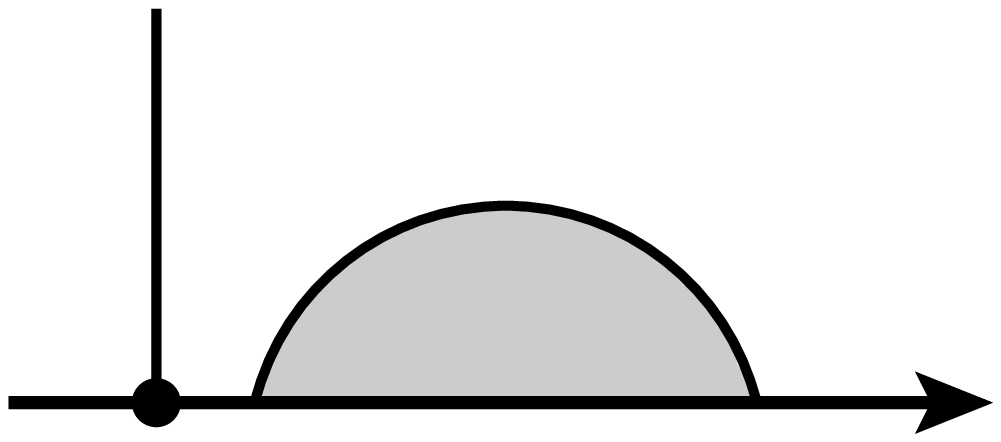}}
  \quad=\quad
  \rb{-7mm}{\ig[height=15mm]{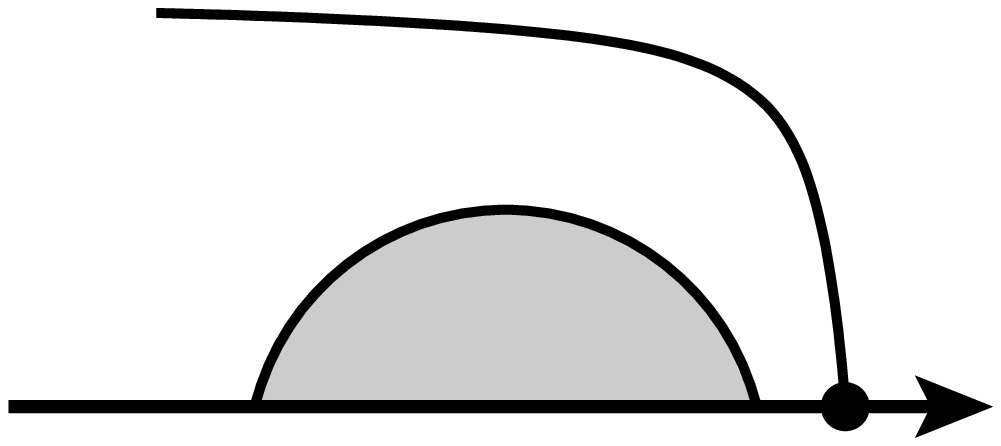}}
$$
Here the horizontal line is a fragment of the circle of the diagram,
while the grey region denotes an arbitrary conglomeration of chords.

\noindent (b) Prove the following relation:

$$\rb{-7mm}{\ig[height=12mm]{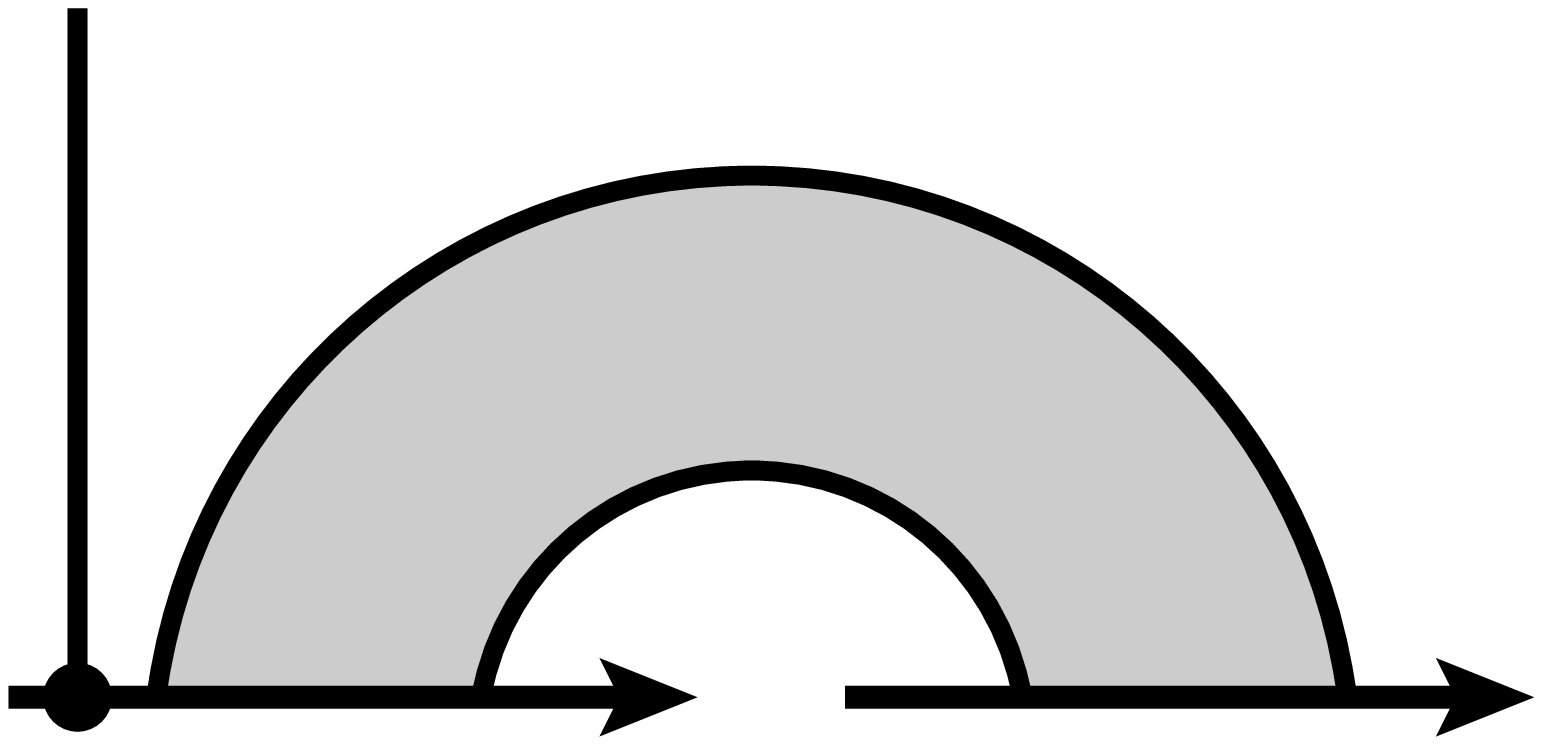}}
  \ +\ \rb{-7mm}{\ig[height=12mm]{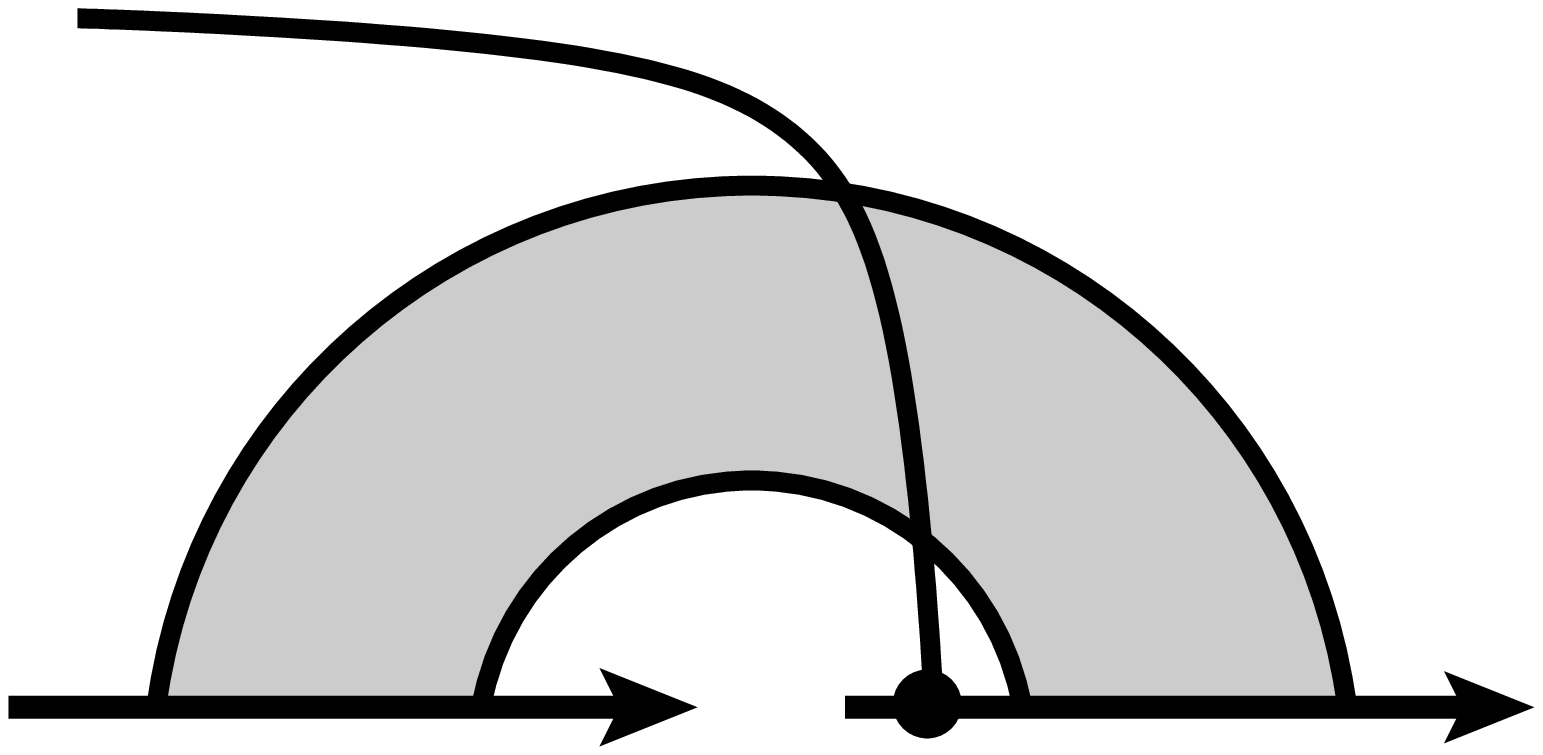}}
  \ =\quad
  \rb{-7mm}{\ig[height=12mm]{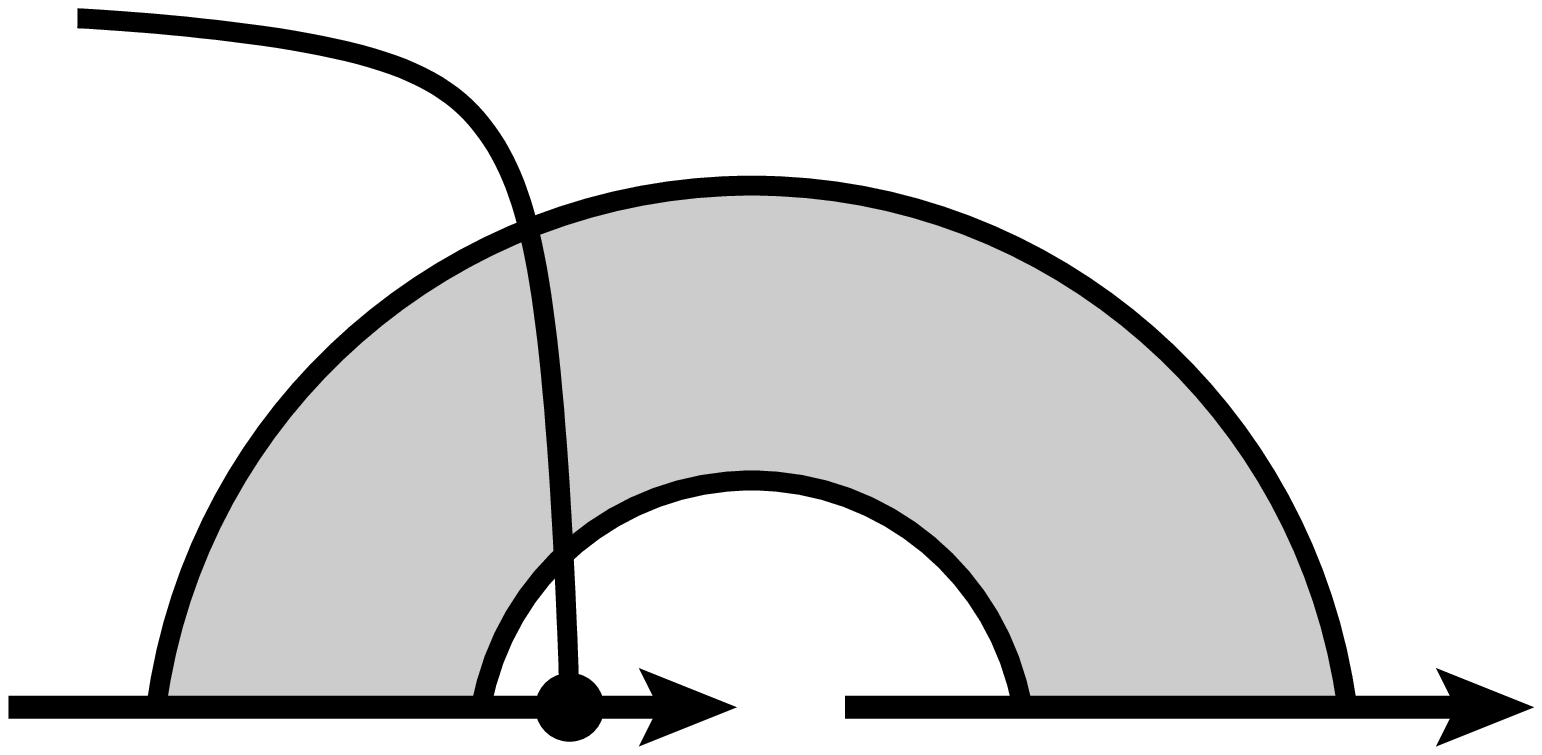}}
  \ +\ \rb{-7mm}{\ig[height=12mm]{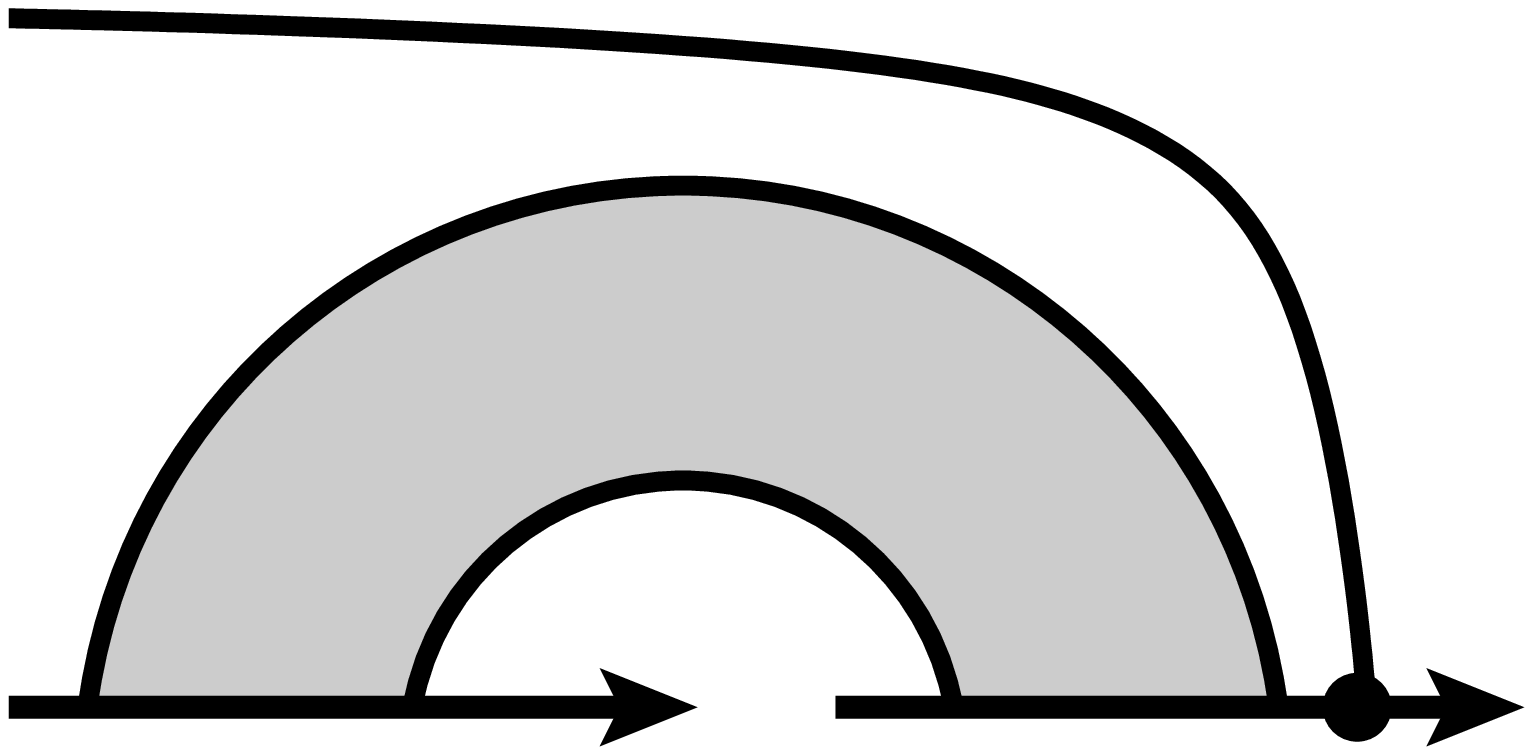}}
$$
\noindent or, in circular form:\qquad $\index{Four-term
relation!generalized}\index{Generalized 4-term relation}
\label{gen4T} \chd{gen4Tc1}\ +\ \chd{gen4Tc4} \quad=\quad
\chd{gen4Tc2}\ +\ \chd{gen4Tc3}\ . $

\item Using the generalized 4-term relation 
prove the following identity:\vspace{-5pt}
$$\rb{-5mm}{\ig[height=12mm]{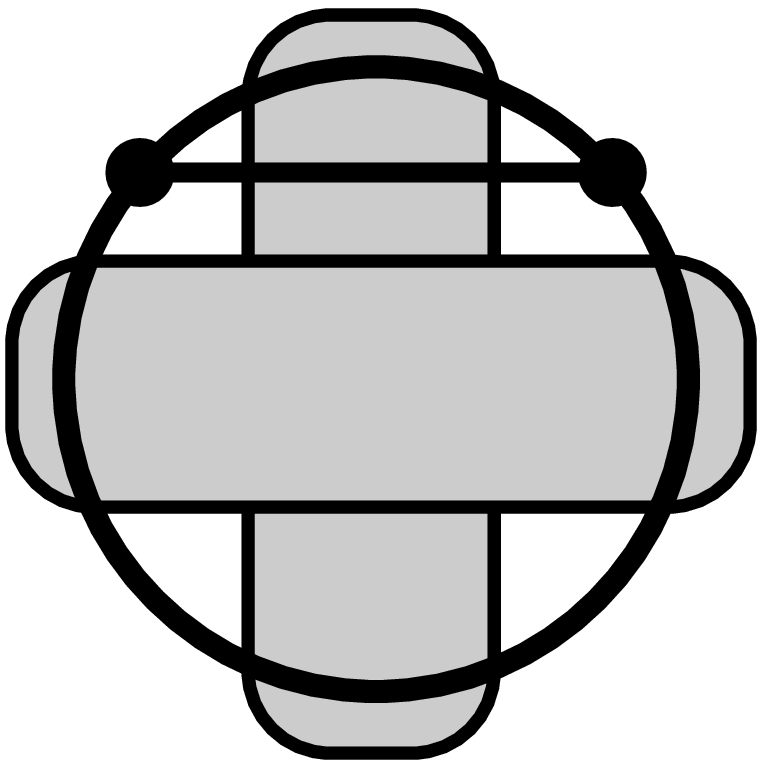}}
  \quad=\quad
  \rb{-5mm}{\ig[height=12mm]{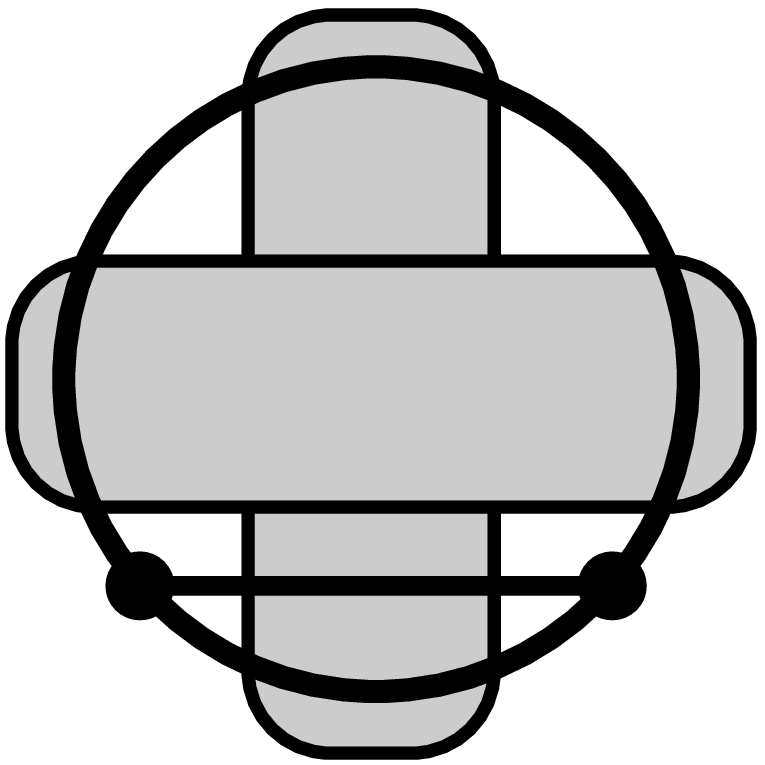}}
$$\vspace{5pt}

\item Prove the proposition of Section~\ref{line_CD}.

\item
Check that for the chord diagram below, the intersection graph and
its canonical decomposition are as shown:
$$\risS{-25}{ccd}{
    \put(10,-10){\mbox{\scriptsize\tt Chord diagram}}}{70}{50}{45}\hspace{4.5cm}
  \risS{-25}{ig-c-ex}{
    \put(20,-10){\mbox{\scriptsize\tt Intersection graph}}}{120}{20}{20}
$$
$$\hspace{-20pt}\risS{-23}{candec-c-ex}{
    \put(45,-12){\mbox{\scriptsize\tt Canonical decomposition}}}{180}{20}{50}
$$

\item
 (\cite[example 6.4.11]{LZ})
Prove that $e^{\symb(c_2)}(D)$ is equal to the number of perfect matchings
of the intersection graph $\Gamma(D)$. (A {\em perfect
matching}\index{Perfect matching} in a graph is a set of disjoint
edges covering all the vertices of the graph.)

\end{enumerate}
\end{xcb}
 %4 Ch.Diag

%\part{Diagrammatic Algebras}
\chapter{Jacobi diagrams} %11
\label{algFDchap}\label{algOD}

In the previous chapter we saw that the study of Vassiliev knot
invariants, at least complex-valued, is largely reduced to the study
of the algebra of chord diagrams. Here we introduce two different
types of diagrams representing elements of this algebra, namely {\em
closed Jacobi diagrams} and {\em open Jacobi diagrams}. These
diagrams provide better understanding of the primitive space $\PR\A$
and bridge the way to the applications of the Lie algebras in the
theory of Vassiliev invariants, see Chapter \ref{LAWS} and
Section~\ref{wheels}.

The name {\em Jacobi diagrams} is justified by a close resemblance
of the basic relations imposed on Jacobi diagrams (STU and IHX) to
the Jacobi identity for Lie algebras.

\section{Closed Jacobi diagrams}
\label{algFD}

\begin{definition}\label{def_FD}
\index{Closed diagram}\index{Jacobi diagram}
\index{Diagram!closed}\index{Order}\index{Degree} A {\em closed
Jacobi diagram} (or, simply, a {\em closed diagram}) is a connected
trivalent graph with a distinguished embedded oriented cycle, called
\textit{Wilson loop}\footnote{this terminology, introduced by Bar-Natan, makes an allusion to field theory where a Wilson loop is an observable that assigns to a connection (field potential) its holonomy along a fixed closed curve.},\index{Wilson loop} and a fixed
cyclic order of half-edges at each vertex not on the Wilson loop.
Half the number of the vertices of a closed diagram is called the
\textit{degree}, or \textit{order}, of the diagram. This number is
always an integer.
\end{definition}

\begin{xremark} Some authors (see, for instance, \cite{HM}) also include the
cyclic order of half-edges at the vertices on the Wilson loop into
the structure of a closed Jacobi diagram; this leads to the same
theory.
\end{xremark}

\begin{xremark} A Jacobi diagram is allowed to have multiple edges and hanging
loops, that is, edges with both ends at the same vertex. It is the
possible presence of hanging loops that requires introducing the
cyclic order on half-edges rather than edges.
\end{xremark}

\begin{xexample} Here is a closed diagram of degree 4:
\begin{center}
  \ig[width=20mm]{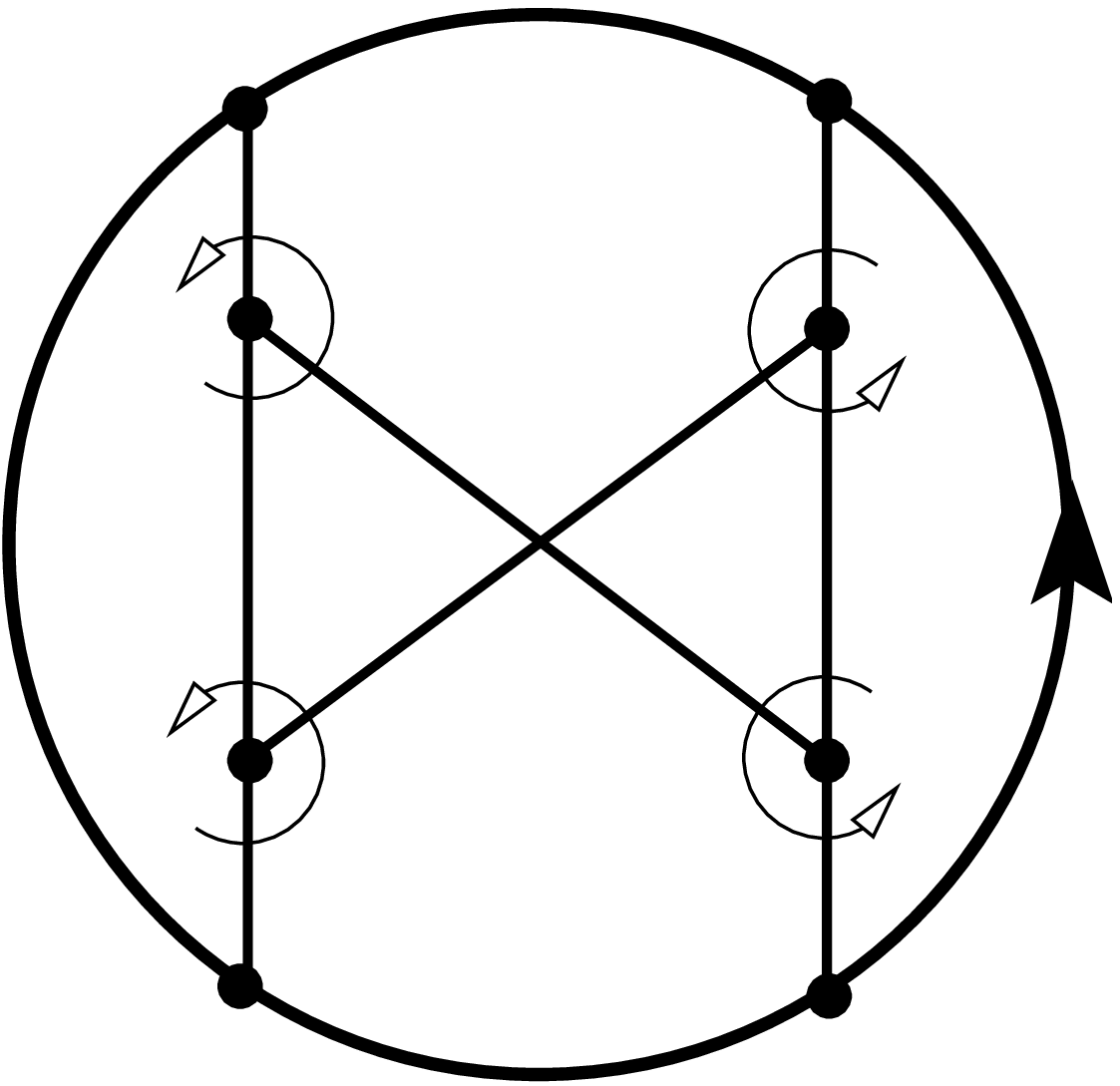}
\end{center}
\end{xexample}
The orientation of the Wilson loop and the cyclic orders of
half-edges at the internal vertices are indicated by arrows. In the
pictures below, we shall always draw the diagram inside its Wilson
loop, which will be assumed to be oriented counterclockwise unless
explicitly specified otherwise. Inner vertices will also be assumed
to be oriented counterclockwise. (This convention is referred to as
the \textit{blackboard orientation}.) Note that the intersection of
two edges in the centre of the diagram above is not actually a
vertex.

Chord diagrams are closed Jacobi diagrams all of whose vertices lie
on the Wilson loop.

Other terms used for closed Jacobi diagrams in the literature
include \textit{Chinese character diagrams} \cite{BN1},
\textit{circle diagrams} \cite{Kn0}, \textit{round diagrams}
\cite{Wil1} and \textit{Feynman diagrams} \cite{KSA}.

\def\dSTU#1#2{\risS{-12}{#1}{#2}{40}{20}{20}}
\def\lClD#1#2{\risS{-17}{#1}{#2}{40}{20}{20}}
\def\ClD#1{\risS{-12}{#1}{}{30}{10}{15}}

\begin{definition}
\index{Vector space!of closed diagrams}
\label{F_n}

The vector space of closed diagrams
$\F_n$ is the space spanned by all closed diagrams of degree $n$
modulo the \textit{STU relations}:
$$\index{Relation!STU}\index{STU}\label{STUrel}
\dSTU{sstuA}{\put(35,-5){\mbox{\scriptsize S}}}\quad =\quad
\dSTU{tstuA}{\put(35,-5){\mbox{\scriptsize T}}}\ -\
\dSTU{ustuA}{\put(35,-5){\mbox{\scriptsize U}}}\ .
$$
The three diagrams S, T and U must be identical outside the shown
fragment. We write $\F$ for the direct sum of the spaces $\F_n$ for
all $n\ge 0$.
\end{definition}

The two diagrams $T$ and $U$ are referred to as the {\em
resolutions} of the diagram $S$. The choice of the plus and minus
signs in front of the two resolutions in the right-hand side of the
STU relation, depends on the orientation for the Wilson loop and on
the cyclic order of the three edges meeting at the internal vertex
of the S-term. Should we reverse one of them, say the orientation of
the Wilson loop, the signs of the T- and U-terms change. Indeed,
$$\dSTU{sstu}{}\ =\ \rb{-5pt}{\dSTU{stuAs}{}}\ \stackrel{\STU}{=}
\ \rb{-5pt}{\dSTU{stuAt}{}}\ - \rb{-5pt}{\dSTU{stuAu}{}}\ =\
\dSTU{ustu}{}\ -\ \dSTU{tstu}{}\ .
$$
This remark will be important in Section \ref{detect_orient} where
we discuss the problem of detecting knot orientation. One may think
of the choice of the direction for the Wilson loop in an $\STU$
relation as a choice of the cyclic order
``forward-sideways-backwards'' at the vertex lying on the Wilson
loop. In these terms, the signs in the $\STU$ relation depend on the
cyclic orders at both vertices of the S-term, the relation above may
be thought of as a consequence of the antisymmetry relation $\AS$
(see \ref{AS_rel}) for the vertex on the Wilson loop, and the $\STU$
relation itself can be regarded as a particular case of the $\IHX$ relation
(see \ref{IHX_rel}).

\subsection{Examples}
There exist two different closed diagrams of order 1:
$\sClD{wssoN1}$\ ,\  $\sClD{fd1b}$\ , one of which vanishes due to
the STU relation:
$$\ClD{fd1b}\quad =\quad \ClD{fd1b+}\ -\ \ClD{fd1b-}\quad =\ 0\ .
$$
There are ten closed diagrams of degree 2:

$$\sClD{fd2a}\ ,\quad\sClD{fd2b}\ ,\quad\risS{-13}{fd2c}{}{33}{20}{15}
\ ,\quad\sClD{fd2d}\ ,
$$
$$\sClD{fd2e}\ ,\quad\risS{-13}{fd2f}{}{33}{20}{15}\ ,\quad\sClD{fd2g}
\ ,\quad\risS{-13}{fd2h}{}{33}{20}{15}\ ,\quad\sClD{fd2i} \
,\quad\sClD{fd2j}\ .
$$

The last six diagrams are zero. This is easy to deduce from the
$\STU$ relations,  but the most convenient way of seeing it is by
using the AS relations which follow from the $\STU$ relations (see
Lemma \ref{STU_AS} below).

Furthermore, there are at least two relations among the first four
diagrams:

\begin{eqnarray*}
  \sClD{fd2d} &=& \sClD{fd2b}\ -\ \sClD{fd2a}\ ; \\
  \risS{-13}{fd2c}{}{33}{20}{15} &=& \sClD{fd2c+}\ -\ \sClD{fd2c-}\
            =\ 2\ \sClD{fd2d}\ .
\end{eqnarray*}

It follows that $\dim\F_2\le2$. Note that the first of the above
equalities gives a concise representation,
\ $\sClD{fd2d}$\ , for the basis primitive element of degree 2.

\medskip

\begin{xca}
Using the $\STU$ relations, rewrite the basis primitive element of order 3
in a concise way.
\end{xca}

{\sl Answer.}
$$\sClD{symcd4}\ -\ 2\ \sClD{symcd3}\ +\ \sClD{symcd1}\quad
   =\quad \sClD{fd34a}\ .
$$

We have already mentioned that chord diagrams are a particular case
of closed diagrams. Using the $\STU$ relations, one can rewrite any
closed diagram as a linear combination of chord diagrams. (Examples
were given just above.)
\medskip

A vertex of a closed diagram that lies on the Wilson loop is called
{\em external};\index{Vertex!external} otherwise it is called {\em
internal}.\index{Vertex!internal} External vertices are also called
{\em legs}.\index{Leg!of a closed diagram} There is an increasing
filtration on the space $\F_n$ by subspaces $\F_n^m$ spanned by
diagrams with at most $m$ external vertices:
$$
  \F_n^1 \subset \F_n^2 \subset ... \subset \F_n^{2n}.
$$

\begin{xca}\label{one_leg_A}
Prove that $\F_n^1=0$.
\end{xca}

{\sl Hint.}
In a diagram with only two legs one of the legs
can go all around the circle and change places with the second.

\section{IHX and AS relations}

\begin{lemma}\label{STU_4T}
The $\STU$ relations imply the $\fT$ relations for chord diagrams.
\end{lemma}

\begin{proof}
Indeed, writing the four-term relation in the form
$$
  \sClD{d4t1} \ -\ \sClD{d4t2} \ =\ \sClD{d4t3} \ -\ \sClD{d4t4}
$$
and applying the $\STU$ relations to both parts of this equation, we
get the same closed diagrams.
\end{proof}

\def\relD#1{\risS{-18}{#1}{}{30}{20}{15}}
\begin{definition}\label{AS_rel}\index{Relation!AS}\index{Antisymmetry relation}
An $\AS$ (={\it antisymmetry}) {\it relation} is:
$$\relD{as1} \quad =\quad -\ \relD{as2}\ .
$$
In other words, a diagram changes sign when the cyclic order of
three edges at a trivalent vertex is reversed.
\end{definition}

\begin{definition}
\label{IHX_rel}\index{IHX relation}\index{Relation!IHX}
An {\it $\IHX$ relation} is:
$$\sClD{ihxI} \quad =\quad \sClD{ihxH}\ -\ \sClD{ihxX}\ .
$$
\end{definition}

As usual, the unfinished fragments of the pictures denote graphs
that are identical (and arbitrary) everywhere but in this explicitly
shown fragment.

\begin{xca}
Check that the three terms of the IHX relation ``have equal
rights''. For example, an H turned 90 degrees looks like an I; write an
IHX relation starting from that I and check that it is the same as
the initial one. Also, a portion of an X looks like an H; write down
an IHX relation with that H and check that it is again the same.
The IHX relation is in a sense unique; this is discussed in
Exercise~\ref{strangeIHX} on page~\pageref{strangeIHX}.
\end{xca}

\begin{lemma}\label{STU_AS}
The $\STU$ relations imply the $\AS$ relations for the internal vertices of a
closed Jacobi diagram.
\end{lemma}

\begin{proof}
Induction on the distance (in edges) of the vertex in question from the
Wilson loop.

{\sl Induction base.}
If the vertex is adjacent to an external vertex, then the assertion follows by
one application of the $\STU$ relation:\vspace{-3pt}
$$\begin{array}{ccl}
\dSTU{sstuA}{} &=& \dSTU{tstuA}{}\ -\ \dSTU{ustuA}{} \vspace{-3pt}\\
\dSTU{sstuB}{} &=& \dSTU{tstuB}{}\ -\ \dSTU{ustuB}{}\ .
\end{array}$$

{\sl Induction step.} Take two closed diagrams $f_1$ and $f_2$ that
differ only by a cyclic order of half-edges at one internal vertex
$v$. Apply STU relations to both diagrams {\it in the same way\/} so
that $v$ gets closer to the Wilson loop.

\end{proof}

\begin{lemma}\label{STU_IHX}
The $\STU$ relations imply the $\IHX$ relations for the internal edges
of a closed diagram.
\end{lemma}

\begin{proof}
The argument is similar to the one used in the previous proof.
We take an $\IHX$ relation somewhere inside a closed diagram and,
applying the same sequence of STU moves to each of the three diagrams,
move the $\IHX$ fragment closer to the Wilson loop. The proof of the induction
base is shown in these pictures:
$$\begin{array}{ccccl}
\sClD{fIHXa} &=& \sClD{fIHXaa} - \sClD{fIHXab}
   &=& \sClD{fIHXi} - \sClD{fIHXii} - \sClD{fIHXiii} + \sClD{fIHXiv}\ , \\
\sClD{fIHXb} &=& \sClD{fIHXba} - \sClD{fIHXbb}
   &=& \sClD{fIHXi} - \sClD{fIHXv} - \sClD{fIHXvi} + \sClD{fIHXiv}\ , \\
\sClD{fIHXc} &=& \sClD{fIHXca} - \sClD{fIHXcb}
   &=& \sClD{fIHXv} - \sClD{fIHXiii} - \sClD{fIHXii} + \sClD{fIHXvi}\ .
\end{array}$$
Therefore,
$$
  \ClD{fIHXa}\quad =\quad \ClD{fIHXb}\ +\ \ClD{fIHXc}\quad .
$$
\end{proof}

\subsection{Other forms of the $\IHX$ relation}

\def\ddD#1#2{\risS{-15}{#1}{#2}{40}{20}{20}}

The $\IHX$ relation can be drawn in several forms, for example:
\begin{itemize}
\setlength{\itemsep}{1pt plus 1pt minus 1pt}
\item (rotationally symmetric form)\vspace{-8pt}
$$\dSTU{ihxRS1}{}\quad +\quad \dSTU{ihxRS2}{}\quad +\quad \dSTU{ihxRS3}{}
   \qquad=\qquad 0\ .
$$
\item (Jacobi form)
$$\ddD{ihxJ2}{}\quad =\quad \ddD{ihxJ1}{}\quad +\quad \ddD{ihxJ3}{}\ .
$$
\item (Kirchhoff form)\vspace{-5pt}
$$\ddD{ihx1}{}\qquad =\qquad \ddD{ihx3}{}\quad +\quad \ddD{ihx2}{}\ .
$$
\end{itemize}

\begin{xca}
By turning your head and pulling the strings of the diagrams, check
that  all these forms are equivalent.
\end{xca}

The Jacobi form of the $\IHX$ relation can be interpreted as
follows. Suppose that to the upper 3 endpoints of each diagram we
assign 3 elements of a Lie algebra, $x$, $y$ and $z$, while every
trivalent vertex, traversed downwards, takes the pair of
``incoming'' elements into their commutator:\vspace{10pt}
$$\ddD{commut}{\put(-8,42){\mbox{$x$}}
               \put(42,42){\mbox{$y$}}
               \put(24,0){\mbox{$[x,y]$}}}\qquad .
$$
Then the IHX relation means that
$$
  [x,[y,z]] = [[x,y],z] + [y,[x,z]],
$$
which is the classical Jacobi identity.
This observation, properly developed, leads to the construction of Lie
algebra weight systems --- see Chapter \ref{LAWS}.
\medskip

The Kirchhoff presentation is reminiscent of the Kirchhoff's law in
electrotechnics. Let us view the portion
\rb{-3mm}{\ig[height=8mm]{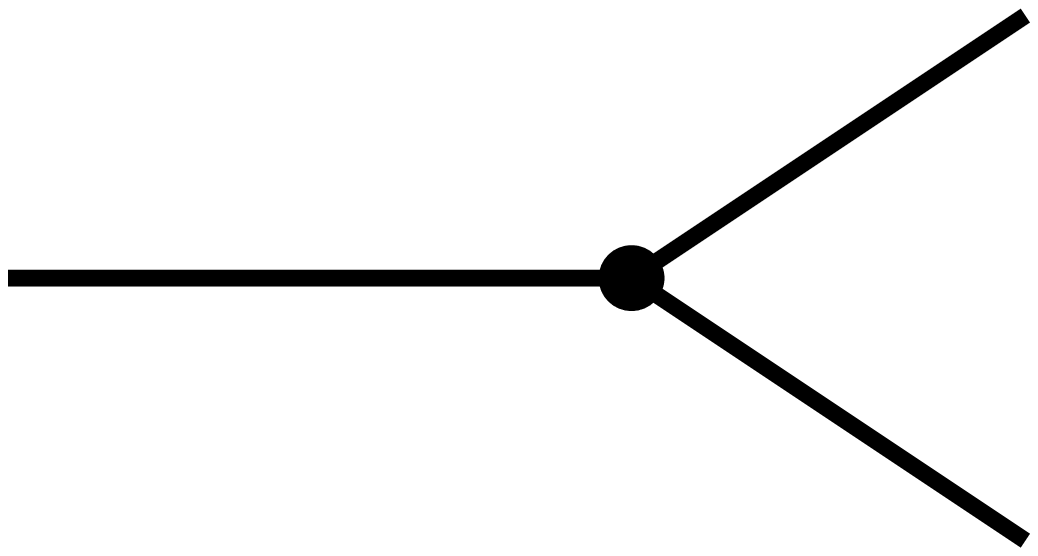}} of the given graph as a piece of
electrical circuit, and the variable vertex as an ``electron'' $e$
with a ``tail'' whose endpoint is fixed. Suppose that the electron
moves towards a node of the circuit:
$$\risS{-15}{kl2}{\put(19,35){\mbox{$e$}}
                  \put(-20,-5){\mbox{\small ``tail"}}}{75}{33}{28}
$$
Then the $\IHX$ relation expresses the well-known Kirchhoff rule:
{\it the sum of currents entering a node is equal to the sum of
currents going out of it}. This electrotechnical analogy is very
useful, for instance, in the proof of the generalized IHX relation
below.
\medskip

The $\IHX$ relation can be generalized as follows:

\begin{lemma} (Kirchhoff law, or generalized IHX relation).
\label{Kirch}\index{Kirchhoff law}
\index{Relation!Kirchhoff}\index{IHX relation!generalized}
The following identity holds:
$$\risS{-35}{gKirch1}{
     \put(82,52){\mbox{$\scriptstyle 1$}}
     \put(82,44){\mbox{$\scriptstyle 2$}}
     \put(77,39){\mbox{$\scriptstyle\cdot$}}
     \put(77,35){\mbox{$\scriptstyle\cdot$}}
     \put(77,31){\mbox{$\scriptstyle\cdot$}}
     \put(82,24){\mbox{$\scriptstyle k$}}}{80}{28}{30}
\qquad =\qquad \sum_{i=1}^k\ \
\risS{-35}{gKirch2}{
     \put(82,52){\mbox{$\scriptstyle 1$}}
     \put(77,48){\mbox{$\scriptstyle\cdot$}}
     \put(77,46){\mbox{$\scriptstyle\cdot$}}
     \put(77,44){\mbox{$\scriptstyle\cdot$}}
     \put(82,38){\mbox{$\scriptstyle i$}}
     \put(77,34){\mbox{$\scriptstyle\cdot$}}
     \put(77,32){\mbox{$\scriptstyle\cdot$}}
     \put(77,30){\mbox{$\scriptstyle\cdot$}}
     \put(82,24){\mbox{$\scriptstyle k$}}}{80}{28}{30}\qquad ,
$$
where the grey box is an arbitrary subgraph which has only 3-valent
vertices.
\end{lemma}

\begin{proof} Fix a horizontal line in the plane
and consider an immersion of the given graph into the plane with
smooth edges, generic with respect to the projection onto this line.
More precisely, we assume that (1) the projections of all vertices
onto the horizontal line are distinct, (2) when restricted to an
arbitrary edge, the projection has only non-degenerate critical
points, and (3) the images of all critical points are distinct and
different from the images of vertices.

{\it Bifurcation points\/} are the images of vertices and critical
points of the projection. Imagine a vertical line that moves from left to
right; for every position of this line take the sum of all diagrams obtained
by attaching the loose end to one of the intersection points.
This sum does not depend on the position of the vertical line, because
it does not change when the line crosses one bifurcation point.

Indeed, bifurcation points fall into six categories:
$$1)\ \ \ClD{bt0}\quad 2)\ \ \ClD{bt1}\quad
  3)\ \ \risS{-15}{bt2}{}{30}{0}{0}\quad
  4)\ \ \risS{-15}{bt3}{}{30}{0}{0}\quad
  5)\ \ \risS{-25}{bt4}{}{30}{0}{0}\quad
  6)\ \ \risS{-25}{bt5}{}{30}{0}{20}\quad .
$$
In the first two cases the assertion follows from the IHX relation, in
cases 3 and 4 --- from the AS relation. Cases 5 and 6 by a
deformation of the immersion are reduced to a combination of the
previous cases (also, they can be dealt with by one application of
the IHX relation in the symmetric form).

\end{proof}

\begin{xexample}\vspace*{-5pt}
\begin{eqnarray*}
\expic{ex1}\ &=& \ \expic{ex2}\ +\ \expic{ex3} \\
&=&\ \expic{ex4}\ +\ \expic{ex6}\ +\ \expic{ex7}
                                        \vspace{15pt}\\
&=&\ \expic{ex5}\ =\ \expic{ex9}\ +\ \expic{ex8}
\end{eqnarray*}
\end{xexample}
\begin{xremark} The difference between inputs and outputs in the
equation of Lemma \ref{Kirch} is purely notational. We may bend the
left-hand leg to the right and move the corresponding term to the
right-hand side of the equation, changing its sign because of the
antisymmetry relation, and thus obtain:
$$
\sum_{i=1}^{k+1}\quad\rb{-12mm}{\ig[height=20mm]{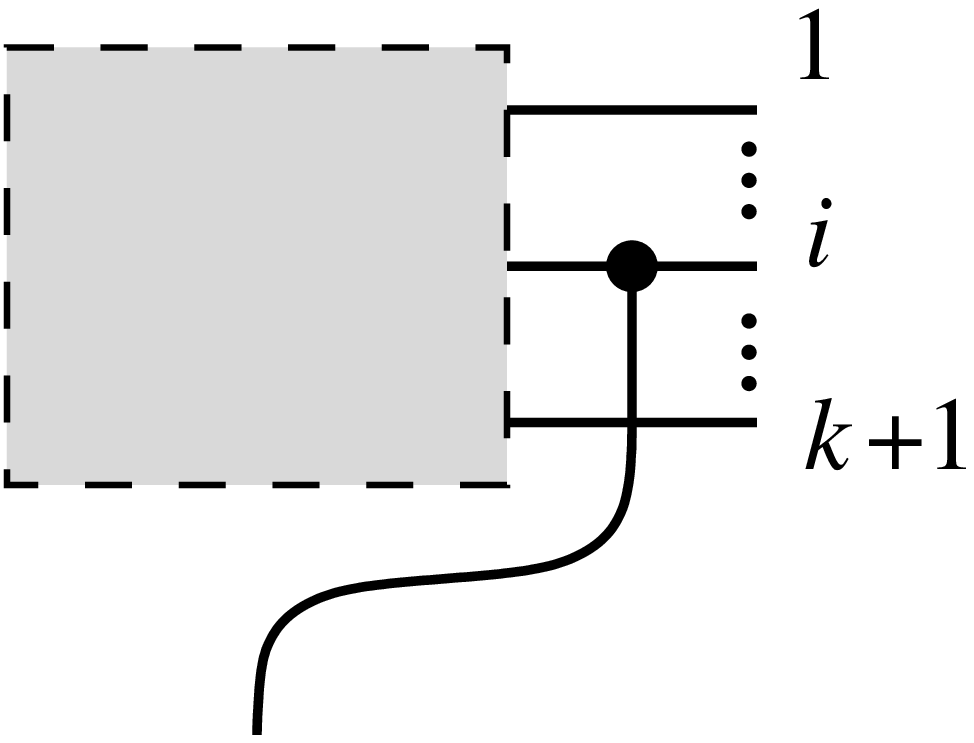}}
\quad=\quad0\ .
$$

Or we may prefer to split the legs into two arbitrary subsets, putting one
part on the left and another on the right. Then:
$$
\sum_{i=1}^{k}\quad\rb{-12mm}{\ig[height=20mm]{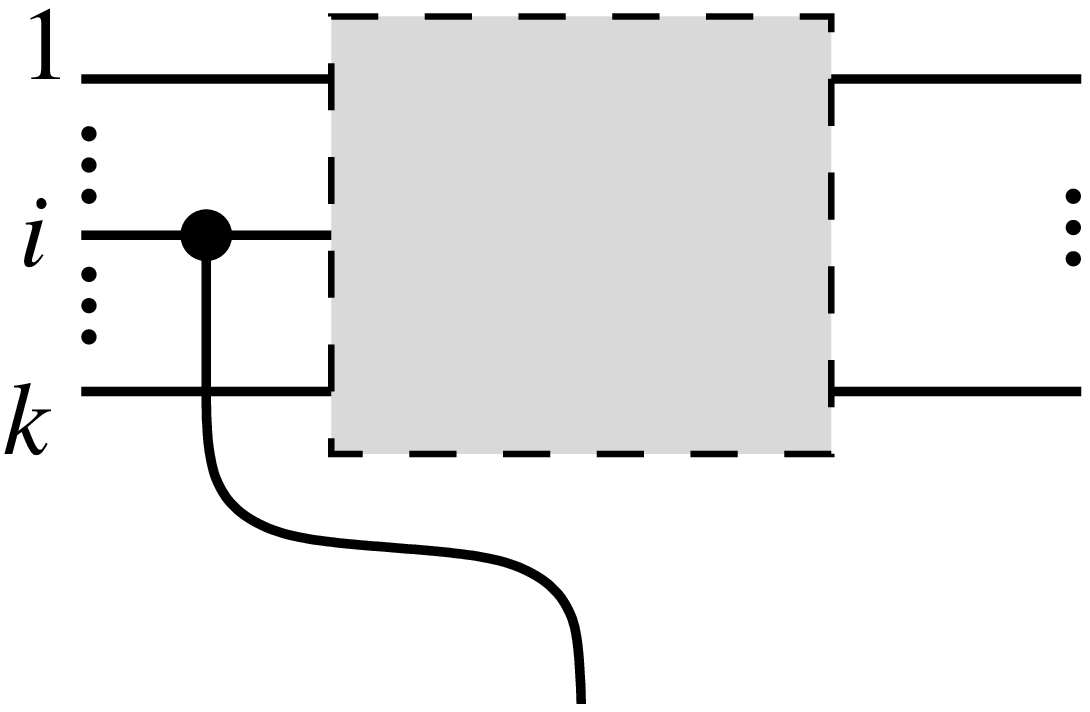}}
\quad=\quad
\sum_{i=1}^{l}\quad\rb{-12mm}{\ig[height=20mm]{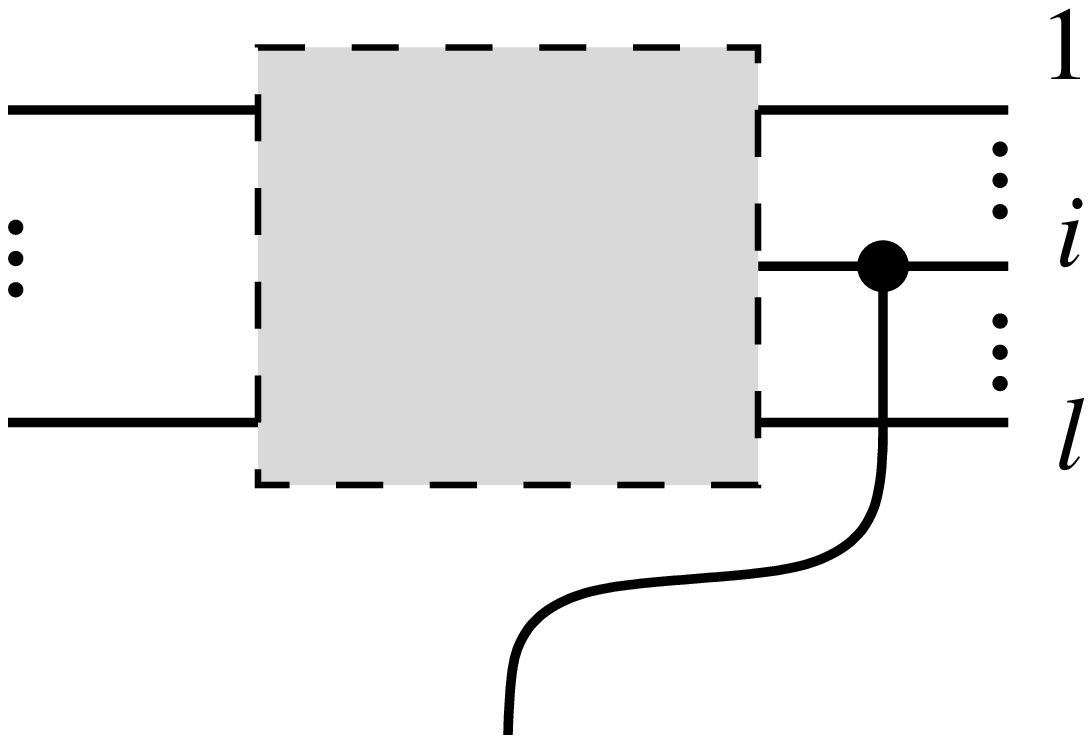}}
\quad.
$$
\end{xremark}

\subsection{A corollary of the AS relation}
\label{AS_corr}

A simple corollary of the antisymmetry relation in the space $\F$ is
that any diagram $D$ containing a hanging loop
$\rb{-3pt}{\ig[height=12pt]{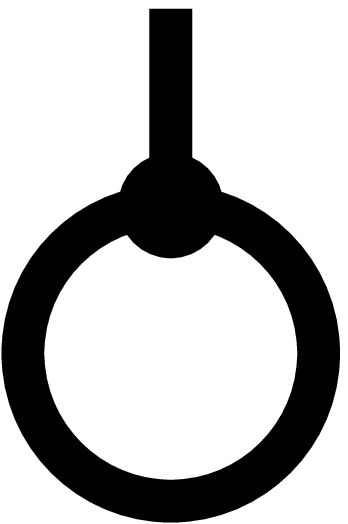}}$ is equal to zero. Indeed,
there is an automorphism of the diagram that changes the two
half-edges of the small circle and thus takes $D$ to $-D$, which
implies that $D=-D$ and $D=0$. This observation also applies to the
case when the small circle has other vertices on it and contains a
subdiagram, symmetric with respect to the vertical axis. In fact,
the assertion is true even if the diagram inside the circle is not
symmetric at all. This is a generalization of
Exercise~\ref{one_leg_A}, but cannot be proved by the same argument.
In Section~\ref{openJD} we shall prove a similar statement
(Lemma~\ref{odd_legs_B}) about {\em open} Jacobi diagrams; that
proof also applies here.

\section{Isomorphism $\A^{fr}\simeq\F$}

Let $\ChD_n$
be the set of chord diagrams of order $n$ and
$\FD_n$ \label{FD_n}
the set
of closed diagrams of the same order. We have a natural inclusion
$\lambda:\ChD_n\to\FD_n$.

\begin{theorem}\label{IsomAF}
\index{Isomorphism!$\A\simeq\F$}
The inclusion $\lambda$ gives rise to an isomorphism of vector
spaces $\lambda:\A^{fr}_n\to\F_n$.
\end{theorem}

\begin{proof}
We must check:
\begin{enumerate}
\item[(A)] that $\lambda$ leads to a well-defined linear map from
$\A_n^{fr}$ to $\F_n$;
\item[(B)] that this map is a linear isomorphism.
\end{enumerate}

\medskip

Part (A) is easy. Indeed,
$\A^{fr}_n=\langle\ChD_n\rangle/\langle\fT\rangle$,
$\F_n=\langle\FD_n\rangle/\langle\STU\rangle$, where angular
brackets denote linear span. Lemma \ref{STU_4T} implies that
$\lambda(\langle\fT\rangle)\subseteq\langle\STU\rangle$, therefore
the map of the quotient spaces is well-defined.
\medskip

(B) We shall construct a linear map $\rho: \F_n\to\A^{fr}_n$ and
prove that it is inverse to $\lambda$.

As we mentioned before, any closed diagram by the iterative use of
$\STU$ relations can be transformed into a combination of chord
diagrams. This gives rise to a map
$\rho:\FD_n\to\langle\ChD_n\rangle$ which is, however, multivalued,
since the result may depend on the specific sequence of relations
used. Here is an example of such a situation (the place where the
$\STU$ relation is applied is marked by an asterisk):
\begin{align}
\ \ \risS{-13}{fd33}{\put(-6,10){*}}{32}{20}{15}\ &\mapsto\
    \sClD{isAC34a} - \sClD{fd34b} = 2 \sClD{isAC34a}\ \mapsto\
    2\Bigl(\sClD{isAC32} - 2\ \sClD{isAC33} + \sClD{isAC34}\Bigr),\notag\\
\risS{-13}{fd33}{\put(31,0){*}}{32}{20}{15}\ &\mapsto\
    \sClD{fd34c} - \sClD{fd34d} =
2\ \ \risS{-13}{fd35a}{\put(-5,8){*}}{32}{20}{15} -
2\ \risS{-13}{fd35b}{\put(33,9){*}}{32}{20}{15} \notag\\
  &\mapsto\ 2\Bigl(\sClD{isAC32} - \sClD{isAC33} - \sClD{isAC34}
                 + \sClD{isAC35}\Bigr).\notag
\end{align}

However, the combination $\rho(C)$ is well-defined as an element of
$\A_n^{fr}$, that is, modulo the $\fT$ relations. The proof of this
fact proceeds by induction on the number $k$ of internal vertices in
the diagram $C$.

If $k=1$, then the diagram $C$ consists of one tripod and several chords
and may look something like this:
$$
  \rb{-10mm}{\ig[width=20mm]{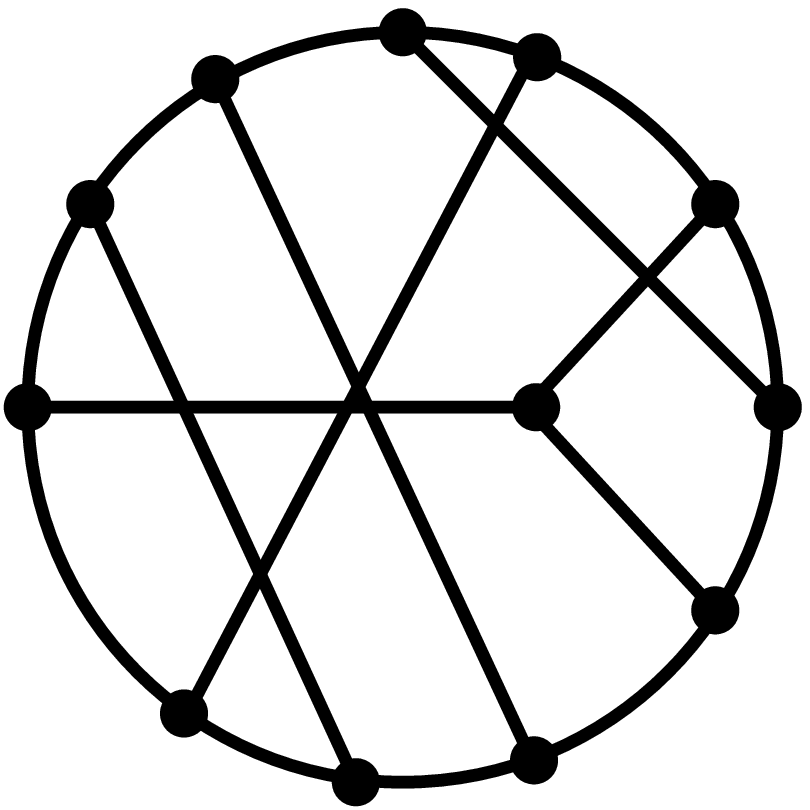}}
$$
There are 3 ways to resolve the internal triple point by an $\STU$
relation, and the fact that the results are the same in $\A^{fr}_n$
is exactly the definition of the 4T relation.

Suppose that $\rho$ is well-defined on closed diagrams with $<k$
internal vertices. Pick a diagram in $\FD_n^{2n-k}$. The process of
eliminating the triple points starts with a pair of neighbouring external
vertices. Let us prove, modulo the inductive hypothesis, that if we
change the order of these two points, the final result will remain
the same.

There are 3 cases to consider: the two chosen points on the Wilson
loop are (1) adjacent to a common internal vertex, (2) adjacent to
neighbouring internal vertices, (3) adjacent to non-neighbouring internal
vertices. The proof for the cases (1) and (2) is shown in the
pictures that follow.

(1)\vspace{-10pt}
$$\dSTU{fdRIa}{\put(7,-6){\mbox{$*$}}}\quad \longmapsto\quad
      \dSTU{fdRIaa}{}\ -\ \dSTU{fdRIab}{}\ ,\vspace{-5pt}
$$
$$\dSTU{fdRIa}{\put(27,-6){\mbox{$*$}}}\quad \longmapsto\quad
      \dSTU{fdRIba}{}\ -\ \dSTU{fdRIab}{}\ ,
$$
The position of an isolated chord does not matter, because, as we
know, the multiplication in $\A^{fr}$ is well-defined.

(2)\vspace{-15pt}
$$\hspace{40pt}\begin{array}{ccl}
\dSTU{fdRIIa}{\put(7,-6){\mbox{$*$}}} &\longmapsto&
          \dSTU{fdRIIaa}{\put(17,-6){\mbox{$*$}}} \ -\
          \dSTU{fdRIIab}{\put(29,-6){\mbox{$*$}}}  \\ &\longmapsto&
      \dSTU{fdRIIaaa}{}\ -\ \dSTU{fdRIIaab}{}\ -\
      \dSTU{fdRIIaba}{}\ +\ \dSTU{fdRIIabb}{}\  ,
\end{array}\vspace{-5pt}
$$
$$\hspace{40pt}\begin{array}{ccl}
\dSTU{fdRIIa}{\put(29,-6){\mbox{$*$}}} &\longmapsto&
          \dSTU{fdRIIba}{\put(7,-6){\mbox{$*$}}} \ -\
          \dSTU{fdRIIbb}{\put(7,-6){\mbox{$*$}}}  \\ &\longmapsto&
      \dSTU{fdRIIbaa}{}\ -\ \dSTU{fdRIIaba}{}\ -\
      \dSTU{fdRIIaab}{}\ +\ \dSTU{fdRIIabb}{}\  .
\end{array}
$$
After the first resolution, we can choose the sequence of further
resolutions arbitrarily, by the inductive hypothesis.

\noindent{\bf Exercise.} Give a similar proof for the case (3).
\medskip

We thus have a well-defined linear map $\rho:\F_n\to\A_n^{fr}$. The fact
that it is two-sided inverse to $\lambda$ is clear.
\end{proof}

\section{Product and coproduct in $\F$}

Now we shall define a bialgebra structure in the space $\F$.

\begin{definition}
The product \index{Product!in $\F$}\index{Closed diagram!product} of
two closed diagrams is defined in the same way as for chord
diagrams: the two Wilson loops are cut at arbitrary places and then
glued together into one loop, in agreement with the orientations:
$$\sClD{fd33}\ \cdot\ \sClD{fd2d}\ =\ \sClD{fd_prod}\ .
$$
\end{definition}

\begin{proposition}
This multiplication is well-defined, that is, it does not depend on
the place of cuts.
\end{proposition}

\begin{proof}
The isomorphism $\A^{fr}\isom\F$ constructed in Theorem
\ref{IsomAF} identifies the product in $\A^{fr}$ with the above
product in $\F$.

Since the multiplication is well-defined in $\A^{fr}$, it is also
well-defined in $\F$.
\end{proof}

To define the {\em coproduct} in the space $\F$, we need the
following definition:

\begin{definition}\label{con_co_ClD}\index{Closed diagram!connected}
The {\em internal graph} \index{Internal graph of a closed diagram}
of a closed diagram is the graph obtained by stripping off the
Wilson loop. A closed diagram is said to be {\em connected} if its
internal graph is connected. The {\it connected components} of a
closed diagram are defined as the connected components of its
internal graph.
\end{definition}

In the sense of this definition, any chord diagram of order $n$ consists
of $n$ connected components --- the maximal possible number.

Now, the construction of the coproduct proceeds in the same way as
for chord diagrams.

\begin{definition}
Let $D$ be a closed diagram and $[D]$ the set of its connected
components. For any subset $J\subseteq[D]$ denote by $D_J$ the
closed diagram with only those components that belong to $J$ and by
$D_{\ol J}$ the ``complementary''
diagram (${\ol J}:=[D]\setminus J$). We set
$$\index{Coproduct!in $\F$}\index{Closed diagram!coproduct}
  \delta(D) := \sum_{J\subseteq[D]} D_J \ot D_{\ol J}.
$$
\end{definition}

\begin{xexample}
$$ \delta\Bigl(\sClD{fd_prod}\Bigr)\ =\
1 \ot \sClD{fd_prod}\ +\ \sClD{fd33}\ot\sClD{fd2d}\ +\
  \sClD{fd2d}\ot\sClD{fd33}\ +\ \sClD{fd_prod}\ot 1.
$$
\end{xexample}

We know that the algebra $\F$, as a vector space, is spanned by
chord diagrams. For chord diagrams, algebraic operations defined in
$\A^{fr}$ and $\F$, tautologically coincide. It follows that the
coproduct in $\F$ is compatible with its product and that the
isomorphisms $\lambda$, $\rho$ are, in fact, isomorphisms of
bialgebras.

\section{Primitive subspace of $\F$}
\label{FDprim}

By definition, connected closed diagrams are primitive with respect to the
coproduct $\delta$. 
It may sound surprising that the converse is also true:

\begin{theorem} \cite{BN1} \label{primF}
\index{Primitive space!in $\F$}
The primitive space $\PR$ of the bialgebra $\F$ coincides with the linear
span of connected closed diagrams.
\end{theorem}

Note the contrast of this straightforward characterization of the
primitive space in $\F$ with the case of chord diagrams.

\begin{proof}
If the primitive space $\PR$ were bigger than the span of connected
closed diagrams, then, according to Theorem \ref{MMthm}, it would
contain an element that cannot be represented as a polynomial in
connected closed diagrams. Therefore, to prove the theorem it is
enough to show that every closed diagram is a polynomial in
connected diagrams. This can be done by induction on the number of
legs  of a closed diagram $C$. Suppose that the diagram $C$ consists
of several connected components (see \ref{con_co_ClD}). The $\STU$
relation tells us that we can freely interchange the legs of $C$
modulo closed diagrams with fewer legs. Using such permutations we
can separate the connected components of $C$. This means that modulo
closed diagrams with fewer legs $C$ is equal to the product of its
connected components.
\end{proof}

\subsection{Filtration of $\PR_n$}
\label{filtr_pr}

The primitive space $\PR_n$ cannot be graded by the number of legs
$k$, because the STU relation is not homogeneous with respect to
$k$. However, it can be \textit{filtered}: \index{Primitive space!
filtration}
$$0=\PR_n^1\subseteq  \PR_n^2\subseteq \PR_n^3\subseteq \dots \subseteq
    \PR_n^{n+1}=\PR_n\ .
$$
where $\PR_n^k$ is the subspace of $\PR_n$ generated by connected
closed diagrams with \textit{at most} $k$ legs.

The connectedness of a closed diagram with $2n$ vertices implies
that the number of its legs cannot be bigger than $n+1$. That is why
the filtration ends at the term $\PR_n^{n+1}$.

The following facts about the filtration are known.
\begin{itemize}
\item \cite{CV} The filtration stabilizes even sooner. Namely, $\PR_n^n=\PR_n$
for even $n$, and $\PR_n^{n-1}=\PR_n$ for odd $n$. Moreover, for
even $n$ the quotient space $\PR_n^n / \PR_n^{n-1}$ has dimension
one and is generated by the {\it wheel} $\ol{w}_n$ with $n$
spokes\index{Wheel}:
$$\ol{w}_n = \quad \risS{-13}{wssl2wn}{
         \put(0,-8){\mbox{\scriptsize $n$ spokes}}}{30}{20}{25}
$$
This fact is related to the Melvin-Morton conjecture (see Section
\ref{melmor} and Exercise~\ref{ex_filt_pr}).

\item \cite{Da1} The quotient space $\PR_n^{n-1} / \PR_n^{n-2}$ has dimension
$[n/6]+1$ for odd $n$, and $0$ for even $n$.

\item \cite{Da2} For even $n$
$$
  \dim(\PR_n^{n-2} / \PR_n^{n-3})= \left[\frac{(n-2)^2+12(n-2)}{48}\right] +1\ .
$$

\item For small degrees the dimensions of the quotient spaces
$\PR_n^k / \PR_n^{k-1}$ were calculated by J.~Kneissler \cite{Kn0}
(empty entries in the table are zeroes):

\begin{center}\index{Primitive space! dimensions}
\index{Table of!dimensions of!the primitive spaces}
\label{dim-of-prim-sp}
\begin{tabular}{c||c|c|c|c|c|c|c|c|c|c|c|c||c}
\makebox(10,10){\begin{picture}(10,10)(0,0)
   \put(-2,12){\line(5,-4){18}}
   \put(0,-2){\mbox{$n$}} \put(9,5){\mbox{$k$}} \end{picture} }
  & 1 & 2 & 3 & 4 & 5 & 6 & 7 & 8 & 9 & 10 & 11 & 12 & $\dim \PR_n$\\
\hline \hline
1 &   & 1 &   &   &   &   &   &   &   &    &    &    &  1 \\ \hline
2 &   & 1 &   &   &   &   &   &   &   &    &    &    &  1 \\ \hline
3 &   & 1 &   &   &   &   &   &   &   &    &    &    &  1 \\ \hline
4 &   & 1 &   & 1 &   &   &   &   &   &    &    &    &  2 \\ \hline
5 &   & 2 &   & 1 &   &   &   &   &   &    &    &    &  3 \\ \hline
6 &   & 2 &   & 2 &   & 1 &   &   &   &    &    &    &  5 \\ \hline
7 &   & 3 &   & 3 &   & 2 &   &   &   &    &    &    &  8 \\ \hline
8 &   & 4 &   & 4 &   & 3 &   & 1 &   &    &    &    & 12 \\ \hline
9 &   & 5 &   & 6 &   & 5 &   & 2 &   &    &    &    & 18 \\ \hline
10&   & 6 &   & 8 &   & 8 &   & 4 &   & 1  &    &    & 27 \\ \hline
11&   & 8 &   &10 &   &11 &   & 8 &   & 2  &    &    & 39 \\ \hline
12&   & 9 &   &13 &   &15 &   &12 &   & 5  &    & 1  & 55
\end{tabular}\end{center}
\end{itemize}

\subsection{Detecting the knot orientation}
\label{detect_orient} One may notice that in the table above  all
entries with odd $k$ vanish. This means that any connected closed
diagram with an odd number of legs is equal to a suitable linear
combination of diagrams with fewer legs. This observation is closely
related to the problem of distinguishing knot orientation by
Vassiliev invariants. The existence of the universal Vassiliev
invariant given by the Kontsevich integral reduces the problem of
detecting the knot orientation to a purely combinatorial problem.
Denote by $\tau$ \label{def:tau-rever} the operation of reversing
the orientation of the Wilson loop of a chord diagram; its action is
equivalent to a mirror reflection of the diagram as a planar
picture. This operation descends to $\A$; we call an element of $\A$
{\em symmetric},\index{Chord diagram!symmetric} if $\tau$ acts on it
as identity. Then,  Vassiliev invariants do not distinguish the
orientation of knots if and only if all chord diagrams are
symmetric: $D=\tau(D)$ for all $D\in\A$. The following theorem
translates this fact into the language of primitive subspaces.

\begin{xtheorem} \label{or_detect_in_C}\index{Orientation!detecting}
Vassiliev invariants do not distinguish the orientation
of knots if and only if $\PR_n^k = \PR_n^{k-1}$ for any odd $k$
and arbitrary $n$.
\end{xtheorem}

To prove the Theorem we need to reformulate the question whether
$D=\tau(D)$ in terms of closed diagrams. Reversing the orientation
of the Wilson loop on closed diagrams should be done with some
caution, see the discussion in \ref{STUrel} on
page~\pageref{STUrel}). The correct way of doing it is carrying the
operation $\tau$ from chord diagrams to closed diagrams by the
isomorphism $\lambda:\A^{fr}\to\F$; then we have the following
assertion:

\begin{xlemma} Let
$P=\lClD{prinor1}{\put(17,18){\mbox{$\scriptstyle P'$}}}$\
be a closed diagram with $k$ external vertices. Then
$$\tau(P)\ =\ (-1)^k\ \lClD{prinor2}{\put(17,18){\mbox{$\scriptstyle P'$}}}\ .
$$
\end{xlemma}

\begin{proof} Represent $P$ as a linear combination of chord
diagrams using $\STU$ relations, and then reverse the orientation of
the Wilson loop of all chord diagrams obtained. After that, convert
the resulting linear combination back to a closed diagram. Each
application of the $\STU$ relation multiplies the result by $-1$
because of the reversed Wilson loop (see page~\pageref{STUrel}). In
total, we have to perform the $\STU$ relation $2n-k$ times, where
$n$ is the degree of $P$. Therefore, the result gets multiplied by
$(-1)^{2n-k}=(-1)^k$.
\end{proof}

In the particular case $k=1$ the Lemma asserts that $\PR_n^1=0$ for
all $n$ --- this fact appeared earlier as Exercise~\ref{one_leg_A}.

The operation $\tau:\F\to\F$ is, in fact, an algebra automorphism,
$\tau(C_1\cdot C_2) = \tau(C_1)\cdot \tau(C_2)$.  Therefore, to
check the equality $\tau=\id_{\F}$ it is enough to check it on the
primitive subspace, that is, determine whether $P=\tau(P)$ for every
connected closed diagram $P$.

\medskip
\noindent{\bf Corollary of the Lemma.}\label{B_coroll} {\it Let
$P\in \PR^k =\bigoplus\limits_{n=1}^\infty \PR_n^k$ be a connected
closed diagram with $k$ legs. Then\quad $\tau(P)\equiv (-1)^kP \mod
\PR^{k-1}$.}

\medskip
\noindent{\bf Proof of the Corollary.} Rotating the Wilson loop in
3-space by $180^\circ$ about the vertical axis, we get:
$$\tau(P)\ =\ (-1)^k\ \lClD{prinor2}{\put(17,18){\mbox{$\scriptstyle P'$}}}
 \  =\ (-1)^k\ \lClD{prinor3}{\put(17,22){\mbox{$\scriptstyle P'$}}}\ .
$$
The $STU$ relations allow us to permute the legs modulo diagrams
with fewer number of legs. Applying this procedure to the last
diagram we can straighten out all legs and get $(-1)^kP$.
\hfill$\square$

\medskip
\noindent{\bf Proof of the Theorem.} Suppose that the Vassiliev
invariants do not distinguish the orientation of knots. Then
$\tau(P)=P$ for every connected closed diagram $P$. In particular,
for a diagram $P$ with an odd number of legs $k$ we have $P\equiv -P
\mod \PR^{k-1}$. Hence, $2P\equiv 0 \mod \PR^{k-1}$, which means
that $P$ is equal to a linear combination of diagrams with fewer
legs, and therefore $\dim(\PR_n^k / \PR_n^{k-1})=0$.

Conversely, suppose that Vassiliev invariants do distinguish the
orientation. Then there is a connected closed diagram $P$ such that
$\tau(P)\not=P$. Choose such $P$ with the smallest possible number of
legs $k$. Let us show that $k$ cannot be even. Consider
$X=P-\tau(P)\not=0$. Since $\tau$ is an involution $\tau(X)=-X$. But, in the
case of even $k$, the non-zero element $X$ has fewer legs than $k$,
and $\tau(X)=-X\not=X$, so $k$ cannot be minimal. Therefore, the minimal
such $k$ is odd, and $\dim(\PR_n^k / \PR_n^{k-1})\not=0$.
\hfill$\square$

\begin{xca}
Check that, for invariants of fixed degree, the theorem can be
specialized as follows. {\it Vassiliev invariants of degree $\le n$
do not distinguish the orientation of knots if and only if $\PR_m^k
= \PR_m^{k-1}$ for any odd $k$ and arbitrary $m\le n$.}
\end{xca}

\begin{xca}
Similarly to the filtration in the primitive space $\PR$, one can
introduce the leg filtration in the whole space $\F$. Prove the
following version of the above theorem: {\it Vassiliev invariants of
degree $n$ do not distinguish the orientation of knots if and only if
$\F_n^k = \F_n^{k-1}$ for any odd $k$ and arbitrary $n$.}
\end{xca}

\section{Open Jacobi diagrams}
\label{openJD}
\label{prodB}

The subject of this section is the combinatorial bialgebra $\B$
which is isomorphic to the bialgebras $\A^{fr}$ and $\F$ as a vector
space and as a coalgebra, but has a different natural
multiplication. This leads to the remarkable fact that in the vector
space $\A^{fr}\simeq\F\simeq\B$ there are two multiplications both
compatible with one and the same coproduct.

\begin{definition}\label{ODdef}\index{Open diagram}\index{Diagram!open}
An {\em open Jacobi diagram} is a graph with 1- and 3-valent
vertices, cyclic order of (half-)edges at every 3-valent vertex and
with at least one 1-valent vertex in every connected component.
\end{definition}

An open diagram is not required to be connected. It may have loops
and multiple edges. We shall see later that, modulo the natural
relations any diagram with a loop vanishes. However, it is important
to include the diagrams with loops in the definition, because the
loops may appear during natural operations on open diagrams, and it
is exactly because of this fact that we introduce the cyclic order
on half-edges, not on whole edges.

\index{Order}\index{Degree} The total number of vertices of an open
diagram is even. Half of this number is called the {\em degree} (or
{\em order}) of an open diagram. We denote the set of all open
diagrams of degree $n$ by $\OD_n$\label{OD_n}. \index{Leg!of an
open diagram} The univalent vertices will sometimes be referred to
as {\em legs}.

In the literature, open diagrams are also referred to as {\em
1-3-valent diagrams}, \index{Diagram!1-3-valent} {\em Jacobi
diagrams}, \index{Diagram!Jacobi} {\em web diagrams}
\index{Diagram!web} and {\em Chinese characters}.\index{Chinese
characters}
\medskip

\begin{xdefinition}
An isomorphism between two open diagrams is a one-to-one correspondence
between their respective sets of vertices and half-edges that preserves the
vertex-edge adjacency and the cyclic order of half-edges at every vertex.
\end{xdefinition}

\begin{xexample} Below is the complete list of open diagrams of
degree 1 and 2, up to isomorphism just introduced.

\begin{eqnarray*}
\OD_1&=&\{\ \ccO{ccOa}\ ,\ \rccO{-1.8}{ccOb}\ \} \\
\OD_2&=&\Bigl\{\ \rccO{-1}{ccWa}\ ,\ \rccO{-3.5}{tripod}\ ,\
\rccO{-2}{ccWTb}\ ,\ \rccO{-.8}{ccWb}\ ,\
\rccO{-2.9}{ccWWa}\ ,\ \rccO{-2.5}{ccWWb}\ ,\ \\
&&\hspace{12pt}\rccO{-1.8}{ccWOa}\ ,\ \rccO{-.4}{ccWOb}\ ,\
\rccO{-.2}{ccWOc}\ \Bigr\}
\end{eqnarray*}
\end{xexample}
Most of the elements listed above will be of no importance to us, as they are
killed by the following definition.

\begin{definition}\index{Vector space!of open diagrams}
The space of open diagrams of degree $n$ is the quotient space
$$\B_n:= \langle \OD_n\rangle /\langle \AS,\ \IHX\rangle,
$$
where $\langle \OD_n\rangle$ is the vector space formally generated
by all open diagrams of degree $n$ and $\langle \AS,\ \IHX\rangle$
stands for the subspace spanned by all AS and IHX relations (see
\ref{AS_rel}, \ref{IHX_rel}). By definition, $\B_0$ is
one-dimensional, spanned by the empty diagram, and
$\B:=\Op_{n=0}^\infty \B_n$\label{alg-B}.
\end{definition}

Just as in the case of closed diagrams (Section~\ref{AS_corr}), the
AS relation immediately implies that any open diagram with a loop
(\rb{-3pt}{\ig[height=12pt]{mhang_loop.eps}}) vanishes in $\B$. Let
us give a most general statement of this observation --- valid, in
fact, both for open and for closed Jacobi diagrams.

\begin{xdefinition}
An {\em anti-automorphism} of a Jacobi diagram $b\in\OD_n$ is a
graph automorphism of $b$ such that the cyclic order of half-edges
is reversed in an odd number of vertices.
\end{xdefinition}

\begin{lemma}\label{anti-auto}
If a diagram $b\in\OD_n$ admits an anti-automorphism, then $b=0$ in
the vector space $\B$.
\end{lemma}

\begin{proof}
Indeed, it follows from the definitions that in this case $b=-b$.
\end{proof}

\begin{xexample}
$$\risS{-19}{zero_diag}{}{92}{10}{30}\quad=\quad 0\ .
$$
\end{xexample}
\noindent{\bf Exercise.} Show that $\dim\B_1=1$, $\dim\B_2=2$.
\medskip

The relations AS and IHX imply the generalized IHX relation, or
Kirchhoff law (Lemma \ref{Kirch}) and many other interesting
identities among the elements of the space $\B$.  Some of them are
proved in the next chapter (Section \ref{some_id}) in the context of
the algebra $\G$. Here is one more assertion that makes sense only
in $\B$, as its formulation refers to univalent vertices (legs).

\begin{xlemma}\label{odd_legs_B}
If $b\in\B$ is a diagram with an odd number of legs, all of which are
attached to one and the same edge, then $b=0$ modulo AS and IHX relations.
\end{xlemma}

\begin{xexample}
$$\risS{-19}{3legs}{}{76}{10}{30}\quad=\quad 0\ .
$$
\end{xexample}
Note that in this example the diagram does not have an
anti-automorphism, so the previous lemma does not apply.

\begin{proof}
Any diagram satisfying the premises of the lemma can be put into the
form on the left of the next picture. Then by the generalized IHX
relation it is equal to the diagram on the right which obviously
possesses an anti-automorphism and therefore is equal to zero:
$$\risS{-19}{odd_legs1}{}{76}{30}{20}\quad=\quad
\risS{-19}{odd_legs2}{}{76}{0}{0}
$$
where the grey region is an arbitrary subdiagram.
\end{proof}

In particular, any diagram with exactly one leg vanishes in $\B$.
This is an exact counterpart of the corresponding property of closed
diagrams (see Exercise~\ref{one_leg_A}); both facts are, furthermore,
equivalent to each other in view of the isomorphism $\F\isom\B$
that we shall speak about later (in Section \ref{A=B}).

\noindent{\bf Conjecture.} \label{conj:oddlegs} Any diagram with an
odd number of legs is 0 in $\B$.

This important conjecture is equivalent to the conjecture that
Vassiliev invariants do not distinguish the orientation of knots
(see Section~\ref{grad_B}).

Relations AS and IHX, unlike STU, preserve the separation of
vertices into 1- and 3-valent. Therefore, the space $\B$ has a much
finer grading than $\A^{fr}$. Apart from the main grading by half
the number of vertices, indicated by the subscript in $\B$, it also
has a grading by the number of univalent vertices
$$\index{Grading!by legs in $\B$}\label{leg_grad_B}
   \displaystyle \B = \Op_n \Op_k \B_n^k,
$$
indicated by the superscript in $\B$, so that $\B_n^k$ is the
subspace spanned by all diagrams with $k$ legs and $2n$ vertices in
total.

For disconnected diagrams the second grading can, in turn, be
refined to a multigrading by the number of legs in each connected
component of the diagram:
$$\displaystyle \B = \Op_n \Op_{k_1\leq\ldots\leq k_m} \B_n^{k_1,...,k_m}.
$$

Yet another important grading in the space $\B$ is the grading by
the number of \textit{loops} in a diagram, that is, by its first
Betti number. In fact, we have a
decomposition:
$$\index{Grading!by loops in $\B$}\label{loop_grad_B}
   \displaystyle \B = \Op_n \Op_k \Op_l\ {^l\!\B_n^{k}},
$$
where $l$ can also be replaced by $m$ (the number of connected components)
because of the relation $l+k = n+m$, which can be proved by a simple
argument involving the Euler characteristic.

The abundance of gradings makes the work with the space $\B$ more
convenient than with $\F$, although both are isomorphic, as we shall
soon see.

\subsection{The bialgebra  structure on $\B$}

Both the product and the coproduct in the vector space $\B$ are
defined in a rather straightforward way. We first define the product
and coproduct on diagrams, then extend the operations by linearity
to the free vector space spanned by the diagrams, and then note that
they are compatible with the AS and IHX relations and thus descend
to the quotient space $\B$.

\begin{definition}\index{Product!in $\B$}
The product of two open diagrams is their disjoint union.
\end{definition}

\begin{xexample}\qquad $\rb{0mm}{\ig[width=12mm]{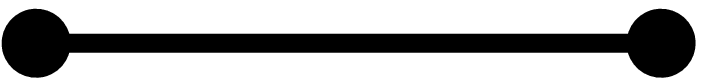}}\ \cdot\
 \rb{-1.5mm}{\ig[width=12mm]{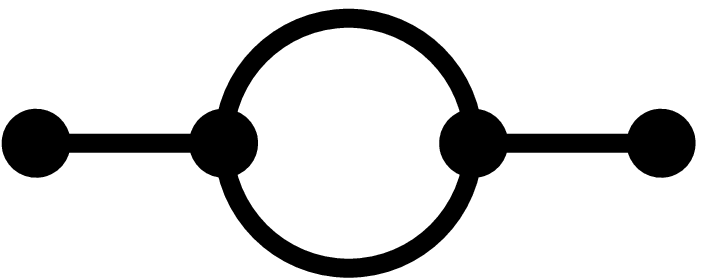}}\ =\
 \rb{-3mm}{\ig[width=12mm]{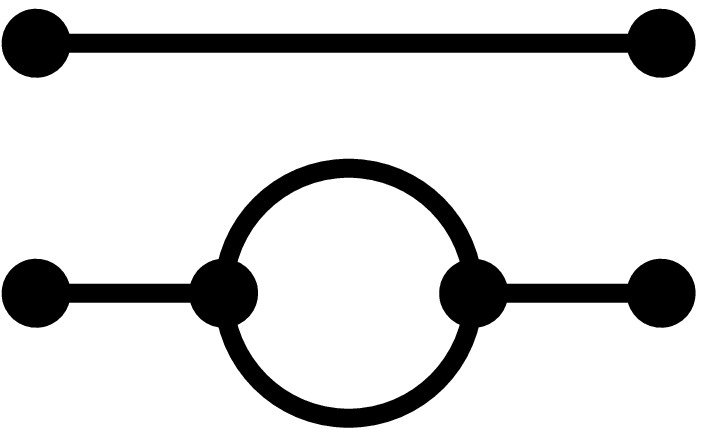}}\ .$
\end{xexample}
\begin{definition}\index{Coproduct!in $\B$}
Let $D$ be an open diagram and $[D]$ --- the set of its connected
components. For a subset $J\subseteq[D]$, denote by $D_J$ the union
of the components that belong to $J$ and by $D_{\ol J}$ --- the
union of the components that do not belong to $J$. We set
$$ \delta(D) := \sum\limits_{J\subseteq [D]}D_J \ot D_{\ol J}.$$
\end{definition}

\begin{xexample}
$$\delta\Bigl(\rccO{-2.6}{ccTb}\Bigr) = 1 \ot \rccO{-2.6}{ccTb} \ +\
 \ccO{ccOa} \ot \rccO{-.8}{ccWb}\ +\ \rccO{-.8}{ccWb} \ot \ccO{ccOa}\
            +\ \rccO{-2.6}{ccTb} \ot 1,
$$
\end{xexample}

As the relations in $\B$ do not intermingle different connected
components of a diagram,  the product of an AS or IHX combination of
diagrams by an arbitrary open diagram belongs to the linear span of
the relations of the same types. Also, the coproduct of any AS or
IHX relation vanishes modulo these relations. Therefore, we have
well-defined algebraic operations in the space $\B$, and they are
evidently compatible with each other. The space $\B$ thus becomes a
graded
bialgebra. \index{Bialgebra!of open diagrams}

\section{Linear isomorphism $\B\simeq\F$}
\label{A=B}

In this section we construct a linear isomorphism between vector
spaces $\B_n$ and $\F_n$. The question whether it preserves
multiplication will be discussed later (Section \ref{twoprodB}). Our
exposition follows \cite{BN1}, with some details omitted, but some
examples added.

To convert an open diagram into a closed diagram, we join all of its
univalent vertices by a Wilson loop. Fix $k$ distinct points on the
circle. For an open diagram with $k$ legs $D\in\OD_n^k$ there are
$k!$ ways of glueing its legs to the Wilson loop at these $k$
points, and we set $\chi(D)$ to be equal to the arithmetic mean of
all the resulting closed diagrams. Thus we get the
\textit{symmetrization map}\index{Symmetrization map}
$$\chi: \OD \to \F.
$$

For example,
\newcommand\kol[1]{\rccO{-4.2}{#1}}
$$
\begin{array}{ccl}
\chi\Bigl(\ \kol{wh4}\ \Bigr)\hspace*{-7pt}&=\hspace*{-8pt}&
\displaystyle\frac{1}{24}\left(
  \kol{kol1}+\kol{kol2}+\kol{kol3}+\kol{kol4}+\kol{kol5}+\kol{kol6}
  \right.\\[5mm]
&&\hspace{11pt}+\
\kol{kol7}+\kol{kol8}+\kol{kol9}+\kol{kol10}+\kol{kol11}+\kol{kol12}\\[5mm]
&&\hspace{11pt}+\
\kol{kol13}+\kol{kol14}+\kol{kol15}+\kol{kol16}+\kol{kol17}+\kol{kol18}\\[5mm]
&&\left. \hspace{11pt}+\
\kol{kol19}+\kol{kol20}+\kol{kol21}+\kol{kol22}+\kol{kol23}+\kol{kol24}
\right).
\end{array}
$$
Scrutinizing these pictures, one can see that 16 out of 24 summands are
equivalent to the first diagram, while the remaining 8 are equivalent to the
second one. Therefore,
$$
\chi(\kol{wh4})\ =\
\displaystyle\frac{1}{3}\ \kol{kol1} + \displaystyle\frac{2}{3}\
\kol{kol2}\ .
$$

\bigskip
\noindent{\bf Exercise.} Express this element via chord diagrams,
using the isomorphism $\F\simeq\A^{fr}$.

{\sl Answer:}\qquad \rb{-5pt}{$\displaystyle
 \kol{cd4-02}\ -\ \frac{10}{3}\ \kol{cd4-06}\ +\ \frac{4}{3}\ \kol{cd4-07}\ .$}

\bigskip
\begin{theorem}\label{BFisom}
\index{Isomorphism!$\B\simeq\F$}
The symmetrization map $\chi:\OD\to\F$ descends to a
linear map $\chi:\B\to\F$, which is a graded isomorphism between the
vector spaces $\B$ and $\F$.
\end{theorem}

The theorem consists of two parts:
\begin{itemize}
\setlength{\itemsep}{1pt plus 1pt minus 1pt}
\item Easy part: $\chi$ is well-defined.
\item Difficult part: $\chi$ is bijective.
\end{itemize}

The proof of bijectivity of $\chi$ is difficult because not every
closed diagram can be obtained by a symmetrization of an open
diagram. For example, the diagram \ \kol{fd2d} \ is not a
symmetrization of any open diagram, even though it looks very much
symmetric. Notice that symmetrizing the internal graph of this
diagram we get 0.

{\bf Easy part of the theorem.} To prove the easy part, we must show
that the $\AS$ and $\IHX$ combinations of open diagrams go to 0 in
the space $\F$. This follows from lemmas \ref{STU_AS} and
\ref{STU_IHX}.

{\bf Difficult part of the theorem.} To prove the difficult part, we
construct a linear map $\tau$ from $\F$ to $\B$, inverse to $\chi$.
\label{inv-chi}
This will be done inductively by the number of legs of the diagrams.
We shall write $\tau_k$ for the restriction of $\tau$ to the
subspace spanned by diagrams with at most $k$ legs.

There is only one way to attach the only leg of an open diagram to
the Wilson loop. Therefore, we can define $\tau_1$ on a closed
diagram $C$ with one leg as the internal graph of $C$. (In fact,
both open and closed diagrams with one leg are all zero in $\B$ and
$\F$ respectively, see Exercise \ref{one_leg_A} and
Lemma \ref{odd_legs_B}). For diagrams with two
legs the situation is similar. Every closed diagram with two legs is
a symmetrization of an open diagram, since there is only one cyclic
order on the set of two elements. For example, \kol{fdWb} is the
symmetrization of the diagram \ \rccO{-1.2}{ccWb}\ ~. Therefore, for
a closed diagram $C$ with two legs we can define $\tau_2(C)$ to be
the internal graph of $C$.

In what follows, we shall often speak of the action of the symmetric
group $S_k$ on closed diagrams with $k$ legs. This action preserves
the internal graph of a closed diagram and permutes the points where
the legs of the internal graph are attached to the Wilson loop.
Strictly speaking, to define this action we need the legs of the
diagrams to be numbered. We shall always assume that such numbering
is chosen; the particular form of this numbering will be irrelevant.

The difference of a closed diagram $D$ and the same diagram whose
legs are permuted by some permutation $\sigma$, is equivalent,
modulo STU relations, to a combination of diagrams with a smaller
number of external vertices. For every given $D$ and $\sigma$ we
fix such a linear combination. 

Assuming that the map $\tau$ is defined for closed diagrams having
less than $k$ legs, we define it for a diagram $D$ with exactly $k$
legs by the formula:
\begin{equation}\label{def_tau_k}
\tau_k(D) = \wt{D} + \frac{1}{k!}\sum_{\sigma\in
S_k}\tau_{k-1}(D-\sigma(D))\ ,
\end{equation}
where $\wt{D}$ is the internal graph of $D$, 
and $D-\sigma(D)$ is represented as a combination of diagrams with
less than $k$ legs according to the choice above.

For example, we know that\quad
$\taud{}{\rccO{-3.8}{fdWb}}\ =\ \rccO{-1.2}{ccWb}$\ ,
and we want to find $\taud{}{\kol{fd2d}}$. By the above formula, we have:
$$\begin{array}{ccl}
  \taud{3}{\kol{fd2d}} &=& \rccO{-3.5}{tripod}
   +\ \frac{1}{6}\ \Biggl(\taud{2}{\kol{fd2d}-\kol{fd2d}}+
                       \taud{2}{\kol{fd2d}-\kol{fd2da}} \vspace{10pt}\\
 &&\hspace{57pt}      +\taud{2}{\kol{fd2d}-\kol{fd2db}}+
                       \taud{2}{\kol{fd2d}-\kol{fd2dc}} \vspace{10pt}\\
 &&\hspace{57pt}      +\taud{2}{\kol{fd2d}-\kol{fd2dd}} +
                       \taud{2}{\kol{fd2d}-\kol{fd2de}}\Biggr)\!\!\vspace{15pt}\\
       &=&\frac{1}{2}\taud{2}{\rccO{-3.8}{fdWb}} = \frac{1}{2}\ \rccO{-1.2}{ccWb}\ .
\end{array}
$$

We have to prove the following assertions:
\begin{itemize}
\item[(i)] 
The value $\tau_{k-1}(D-\sigma(D))$ in the formula (\ref{def_tau_k})
does not depend on the presentation of $D-\sigma(D)$ as a
combination of diagrams with a smaller number of external vertices.
\item[(ii)] The map $\tau$ respects STU relations. 
\item[(iii)] $\chi\circ\tau=\mbox{id}_\F$ and $\tau$ is surjective.
\end{itemize}
The first two assertions imply that $\tau$ is well-defined and the
third means that $\tau$ is an isomorphism. The rest of the section
is dedicated to the proof of these statements.

In the vector space spanned by all closed diagrams (with no
relations imposed) let $\cD^k$ be the subspace spanned by all
diagrams with at most $k$ external vertices. We have a chain of
inclusions
$$\cD^0\subset\cD^1\subset\cD^2\subset\ldots .    $$
We denote by $\I^k$ be the subspace in $\cD^k$ spanned by all STU,
IHX and antisymmetry relations that do not involve diagrams with
more than $k$ external vertices.

\subsection{Action of permutations on closed diagrams}
\label{perm_act_on_C}

The action of the symmetric group $S_k$ on closed diagrams with $k$
legs can be represented graphically as the ``composition'' of a
closed diagram with the diagram of the permutation:
$$k=4\ ;\qquad\qquad \sigma=(4132)= \rccO{-1.3}{perpi}\ ;\qquad\qquad
  D=\kol{kol7}\ ;
$$
$$  \sigma D\ =\ \risS{-18}{actperd1}{}{40}{20}{20}\quad=\quad
     \kol{kol6}\ .
$$

\subsection{}\label{Lemma_5_5_BN} {\bf Lemma.}
Let $D\in \cD^k$. {\it \begin{itemize}
\item Modulo $\I^k$, the difference $D-\sigma D$
belongs to $\cD^{k-1}$.
\item Any choice $U_\sigma$ of a presentation of $\sigma$ as a product
of transpositions determines in a natural way an element
$\G_D(U_\sigma)\in\cD^{k-1}$ such that $$\G_D(U_\sigma)\equiv
D-\sigma D \mod{\I^k}.$$
\item Furthermore, if $U_\sigma$ and $U'_\sigma$ are two such presentations, then
$\G_D(U_\sigma)$ is equal to $\G_D(U'_\sigma)$ modulo $\I^{k-1}$.
\end{itemize}}

This is Lemma 5.5 from \cite{BN1}. Rather than giving the details of
the proof (which can be found in \cite{BN1}) we illustrate it on a
concrete example.

Take the permutation $\sigma=(4132)$ and let $D$ be the diagram
considered above. Choose two presentations of $\sigma$ as a product
of transpositions:
$$
U_\sigma=(34)(23)(34)(12)=\risS{-13}{perpi1}{}{35}{20}{20}\ ;
   \qquad\quad
U'_\sigma=(23)(34)(23)(12)=\risS{-13}{perpi2}{}{35}{20}{20}\ .
$$
For each of these products we represent $D-\sigma D$ as a sum:
\def\per#1{\rb{1.5pt}{$\scriptstyle #1$}}
$$\begin{array}{ccl}
D-\sigma D\! &=&\! (D-\per{(12)}D) + (\per{(12)}D - \per{(34)(12)}D)
  + (\per{(34)(12)}D-\per{(23)(34)(12)}D) \\
&&+ (\per{(23)(34)(12)}D-\per{(34)(23)(34)(12)}D)
\end{array}
$$
and
$$\begin{array}{ccl}
D-\sigma D\! &=&\! (D-\per{(12)}D) + (\per{(12)}D - \per{(23)(12)}D)
          + (\per{(23)(12)}D-\per{(34)(23)(12)}D) \\
&&+ (\per{(34)(23)(12)}D-\per{(23)(34)(23)(12)}D)\ .
\end{array}
$$
Here, the two terms in every pair of parentheses differ only by a
transposition of two neighbouring legs, so their difference is the
right-hand side of an STU relation. Modulo the subspace $\I^4$ each
difference can be replaced by the corresponding left-hand side of
the STU relation, which is a diagram in $\cD^{3}$. We get
$$\G_D(U_\sigma) = \risS{-16}{actperd3}{}{35}{20}{20}
              + \risS{-16}{actperd4}{}{35}{20}{20}
              + \risS{-16}{actperd5}{}{35}{20}{20}
              + \risS{-16}{actperd6}{}{35}{20}{20}
$$
$$\G_D(U'_\sigma) = \risS{-16}{actperd3}{}{35}{20}{20}
              + \risS{-16}{actperd7}{}{35}{20}{20}
              + \risS{-16}{actperd8}{}{35}{20}{20}
              + \risS{-16}{actperd9}{}{35}{20}{20}
$$
Now the difference $\G_D(U_\sigma)-\G_D(U'_\sigma)$ equals
$$
\Bigl(\risS{-16}{actperd4}{}{35}{20}{20}
       -\risS{-16}{actperd9}{}{35}{20}{20}\Bigr) +
\Bigl(\risS{-16}{actperd5}{}{35}{20}{20}
       -\risS{-16}{actperd8}{}{35}{20}{20}\Bigr) +
\Bigl(\risS{-16}{actperd6}{}{35}{20}{20}
       -\risS{-16}{actperd7}{}{35}{20}{20}\Bigr)
$$
Using the STU relation in $\I^3$ we can represent it in the form
$$\G_D(U_\sigma)-\G_D(U'_\sigma) = \risS{-16}{actperd10}{}{35}{20}{20}
            + \risS{-16}{actperd11}{}{35}{20}{20}
            - \risS{-16}{actperd12}{}{35}{20}{20}
 = 0
$$
which is zero because of the IHX relation.

\subsection{Proof of assertions (i) and (ii)}

Let us assume that the map $\tau$, defined by the formula
(\ref{def_tau_k}), is (i) well-defined on $\cD^{k-1}$ and (ii)
vanishes on $\I^{k-1}$.

Define $\tau'(D)$ to be equal to $\tau(D)$ if $D\in\cD^{k-1}$, and
if $D\in \cD^{k}-\cD^{k-1}$ set
$$\tau'(D) = \wt{D}
            + \frac{1}{k!}\sum_{\sigma\in S_k} \tau(\G_D
            (U_\sigma))\ .
$$

Lemma \ref{Lemma_5_5_BN} means that for any given $D\in \cD^k$ with
exactly $k$ external vertices $\tau(\G_D (U_\sigma)))$ does not
depend on a specific presentation $U_\sigma$ of the permutation
$\sigma$ as a product of transpositions. Therefore, $\tau'$ gives a
well-defined map $\cD^k\to\B$.

Let us now show that $\tau'$ vanishes on $\I^k$. It is obvious that
$\tau'$  vanishes on the IHX and antisymmetry relations since these
relations hold in ${\B}$. So we only need to check the STU relation
which relates a diagram $D^{k-1}$ with $k-1$ external vertices and
the corresponding two diagrams $D^k$ and $U_iD^k$ with $k$ external
vertices, where $U_i$ is a transposition $U_i=(i,i+1)$.  Let us
apply $\tau'$ to the right-hand side of the STU relation:
$$\begin{array}{ccl}
\tau'(D^k-U_iD^k) &=& \wt{D^k}
   + \frac{1}{k!}\sum\limits_{\sigma\in S_k} \tau (\G_{D^k} (U_\sigma)) \\
                   &&\hspace{-10pt}  - \wt{U_iD^k}
   - \frac{1}{k!}\sum\limits_{\sigma'\in S_k} \tau
                         (\G_{U_iD^k} (U_{\sigma'}))\ .
\end{array}
$$
Note that $\wt{D^k} = \wt{U_iD^k}$. Reparametrizing the first sum, we get
$$\tau'(D^k-U_iD^k) =
    \frac{1}{k!}\sum_{\sigma\in S_k} \tau (\G_{D^k} (U_\sigma U_i)
                   - \G_{U_iD^k} (U_\sigma))\ .
$$
Using the obvious identity $\G_D(U_\sigma U_i) = \G_D(U_i) +
\G_{U_iD^k} (U_\sigma)$ and the fact that $D^{k-1}=\G_D(U_i)$, we
now obtain
$$\tau'(D^k-U_iD^k) =
    \frac{1}{k!}\sum_{\sigma\in S_k} \tau (D^{k-1})
    = \tau (D^{k-1}) = \tau'(D^{k-1})\ ,
$$
which means that $\tau'$ vanishes on the STU relation, and, hence,
on the whole of $\I^k$.

Now, it follows from the second part of Lemma \ref{Lemma_5_5_BN}
that $\tau'=\tau$ on $\cD^k$. In particular, this means that $\tau$
is well-defined on $\cD^k$ and vanishes on $\I^k$. By induction,
this implies the assertions (i) and (ii).

\subsection{Proof of assertion (iii)}

Assume that $\chi\circ\tau$ is the identity for diagrams with at
most $k-1$ legs. Take $D\in \cD^k$ representing an element of $\F$.
Then
$$\begin{array}{ccl}
(\chi\circ\tau)(D)&=&
   \chi\Bigl(\wt{D}
   + \frac{1}{k!}\sum\limits_{\sigma\in S_k} \tau (\G_{D} (U_\sigma))
  \Bigr) \vspace{8pt}\\
&=& \frac{1}{k!}\sum\limits_{\sigma\in S_k}
   \bigl(\sigma D + (\chi\circ\tau) (\G_{D} (U_\sigma))\bigr)\ .
\end{array}
$$
Since $\G_{D} (U_\sigma)$ is a combination of diagrams with at most
$k-1$ legs, by the induction hypothesis $\chi\circ\tau (\G_{D}
(U_\sigma))=\G_{D} (U_\sigma)$ and, hence,
$$
(\chi\circ\tau)(D) =
  \frac{1}{k!}\sum\limits_{\sigma\in S_k}
     \bigl(\sigma D + \G_{D} (U_\sigma)\bigr) =
  \frac{1}{k!}\sum\limits_{\sigma\in S_k}
     \bigl(\sigma D + D-\sigma D\bigr) = D\ .
$$

\bigskip

The surjectivity of $\tau$ is clear from the definition, so we have
established that $\chi$ is a linear isomorphism between $\B$ and
$\F$. \hfill$\square$

\section{Relation between $\B$ and $\F$}
\label{twoprodB}

It is easy to check that the isomorphism $\chi$ is compatible with
the coproduct in the algebras $\B$ and $\F$. (\textbf{Exercise}:
pick a decomposable diagram $b\in\B$ and check that
$\delta_A(\chi(b))$ and $\chi(\delta_B(b))$ coincide.) However,
$\chi$ is \textit{not} compatible with the product. For example,
$$
  \chi(\ccO{ccOa})\ =\ \sClD{cd1}\ .
$$
The square of the element\quad \ccO{ccOa}\quad in $\B$ is
\quad\rccO{-1}{ccWa}\ . However, the
corresponding element of $\F$
$$\chi(\ \rccO{-1}{ccWa}\ )\ =\ \frac{1}{3}\ \sClD{fd2a}\ +\
   \frac{2}{3}\ \sClD{fd2b}
$$
is not equal to the square of\quad \sClD{cd1}\ .

We can, of course, carry the natural multiplication of the algebra
$\B$ to the algebra $\F$ with the help of the isomorphism $\chi$,
thus obtaining a bialgebra  with two different products, both
compatible with one and the same coproduct.

\medskip

By definition, any connected diagram $p\in\B$ is primitive.
Similarly to Theorem \ref{primF} we have:

\begin{theorem}
The primitive space of the bialgebra $\B$ is spanned by connected
open diagrams.
\end{theorem}

\begin{proof}
The same argument as in the proof of Theorem \ref{primF}, with a
simplification that in the present case we do not have to prove that
every element of $\B$ has a polynomial expression in terms of
connected diagrams: this holds by definition.
\end{proof}

Although the isomorphism $\chi$ does not respect the multiplication,
the two algebras $\B$ and $\F$ are isomorphic. This is clear from
what we know about their structure: by Milnor--Moore theorem both
algebras are commutative polynomial algebras over the corresponding
primitive subspaces. But the primitive subspaces coincide, since
$\chi$ preserves the coproduct! An explicit algebra isomorphism
between $\B$ and $\F$ will be the subject of Section \ref{wheeling}.

Situations of this kind appear in the theory of Lie algebras.
Namely, the bialgebra of invariants in the symmetric algebra of a
Lie algebra $L$ has a natural map into the centre of the universal
enveloping algebra of $L$. This map, which is very similar in spirit
to the symmetrization map $\chi$, is an isomorphism of coalgebras,
but does not respect the multiplication. In fact, this analogy is
anything but superficial. It turns out that the algebra $\F$ is
isomorphic to the centre of the universal enveloping algebra for a
certain Casimir Lie algebra in a certain tensor category. For
further details see \cite{HV}.

\subsection{Unframed version of $\B$}
\index{Jacobi diagram!unframed} The unframed version of the algebras
$\A^{fr}$ and $\F$ are obtained by taking the quotient by the ideal
generated by the diagram with 1 chord $\Theta$. Although the product
in $\B$ is different, it is easy to see that multiplication in $\F$
by $\Theta$ corresponds to multiplication in $\B$ by the
\textit{strut}\index{Strut} $s$: the diagram of degree 1 consisting
of 2 univalent vertices and one edge. Therefore, the unframed
version of the algebra $\B$ is its quotient by the ideal generated
by $s$ and we have: $\B':=\B/(s)\isom\F/(\Theta)=:\F'$.

\subsection{Grading in $\PR\B$ and filtration in $\PR\F$}
\label{grad_B}
The space of primitive elements $\PR\B$ is carried by
 $\chi$ isomorphically onto $\PR\F$. The space $\PR\F=\PR$ is filtered
(see Section~\ref{filtr_pr}), the space $\PR\B$ is  graded
(page~\pageref{leg_grad_B}). It turns out that $\chi$ intertwines
the grading on $\B$ with the filtration on $\F$; as a corollary, the
filtration on $\PR\F$ comes from a grading. Indeed, the definition
of $\chi$ and the construction of the inverse mapping $\tau$ imply
two facts:
\begin{gather*}
  \chi(\PR\B^i) \subset \PR^i \subset \PR^k, \text{ if }i<k,\\
  \tau(\PR^k) \subset \Op_{i=1}^k \PR\B^k.
\end{gather*}
Therefore, we have an isomorphism
$$\tau: \PR_n^k \longrightarrow
        \PR\B_n^1\op \PR\B_n^2\op \dots\op\PR\B_n^{k-1}\op \PR\B_n^k\ .
$$
and, hence, an isomorphism $\PR_n^k/\PR_n^{k-1}\isom \PR\B_n^k$.

Using this fact, we can give an elegant reformulation of the theorem
about detecting the orientation of knots (Section
\ref{detect_orient}, page~\pageref{or_detect_in_C}):

\begin{xcorollary} \index{Orientation!detecting} Vassiliev
invariants do distinguish the orientation of knots if and only if
$\PR\B_n^k\not=0$ for an odd $k$ and some $n$.
\end{xcorollary}

Let us clarify that by saying that Vassiliev invariants do
distinguish the orientation of knots we mean that there exists a
knot $K$ non-equivalent to its inverse $K^*$ and a Vassiliev
invariant $f$ such that $f(K)\ne f(K^*)$.

\begin{xxca}
Check that in the previous statement the letter $\PR$ can be
dropped:  Vassiliev invariants do distinguish the orientation of
knots if and only if $\B_n^k\not=0$ for an odd $k$ and some $n$.
\end{xxca}

The relation between $\F$ and $\B$ in this respect can also be stated in the
form of a commutative diagram:
$$\begin{CD}
\B  @>\chi>>   \F \\
@V{\tau_B}VV     @VV{\tau_C}V \\
\B    @>>\chi>       \F
\end{CD}$$
where $\chi$ is the symmetrization isomorphism, $\tau_C$ is the
orientation reversing map in $\F$ defined by the lemma  in Section
\ref{detect_orient}, while $\tau_B$ on an individual diagram from
$\B$ acts as multiplication by $(-1)^k$ where $k$ is the number of
legs. The commutativity of this diagram is a consequence of the
corollary to the above mentioned lemma (see
page~\pageref{B_coroll}).

\section{The three algebras in small degrees}
\label{ABC_sm_deg} Here is a comparative table which displays some
linear bases of the algebras $\A^{fr}$, $\F$ and $\B$ in small
degrees.

\newcommand\mc{\multicolumn}
\begin{center}
\renewcommand{\arraystretch}{1.5}
\begin{tabular}{|c||cc|c||cc|c||cc|c|}
\hline
$n$ & \mc{3}{c||}{$\A^{fr}$} & \mc{3}{c||}{$\F$} & \mc{3}{c|}{$\B$} \\
\hline
 0 & \mc{3}{c||}{\chd{c25_0}} & \mc{3}{c||}{\chd{c25_0}} &
        \mc{3}{c|}{$\emptyset$} \\
\hline
 1 & \mc{3}{c||}{\chdh{c25_1}} & \mc{3}{c||}{\chdh{c25_1}} &
        \mc{3}{c|}{\ccO{ccOa}}\\
\hline
 2 & \mc{3}{c||}{\begin{tabular}{c|c} \chd{c25_21}\hspace{10pt}
                      &\hspace{7pt}  \chd{c25_22} \end{tabular}}
   & \mc{3}{c||}{\begin{tabular}{c|c} \chd{c25_21}\hspace{10pt}
                      &\hspace{7pt}  \chdh{f25_Wb} \end{tabular}}
   & \mc{3}{c|}{\begin{tabular}{c|c} \chdn{-1}{ccWa}\hspace{10pt}
                      &\hspace{7pt}  \chdn{-.8}{ccWb} \end{tabular}} \\
\hline
 3 & \chdh{c25_31}\! & \chdh{c25_33} & \chdh{c25_34}
   & \chdh{c25_31}\! & \chdh{f25_Tb} & \chdh{f25_Tc}
   & \chdr{-5}{c25_Ta}\! & \chdr{-7}{c25_Tb} & \chdr{-2}{c25_Tc} \\
\hline
 4 & \chd{cd4-08}\! & \chd{cd4-09} & \chd{cd4-06}
   & \chd{cd4-08}\! & \chdh{f25_Fb} & \chdh{f25_Fe}
   & \chdr{-5}{c25_Fa}\! & \chdr{-5}{c25_Fb} & \chdr{-2}{c25_Fe} \\[4mm]
   & \chd{cd4-07}\! & \chd{cd4-15} & \chd{cd4-02}
   & \chdh{f25_Fc}\! & \chdh{f25_Fd} & \chd{f25_Ff}
   & \chdr{-4}{c25_Fc}\! & \chdr{-6}{c25_Fd} & \chdr{-10}{c25_Ff} \\[4mm]
\hline
\end{tabular}
\end{center}

In every order up to 4, for each of the three algebras, this table
displays a basis of the corresponding homogeneous component.
Starting from order 2, decomposable elements (products of elements
of smaller degree) appear on the left, while the new indecomposable
elements appear on the right. The bases of $\F$ and $\B$ are chosen
to consist of primitive elements and their products. We remind that
the difference between the $\A^{fr}$ and $\F$ columns is notational
rather than anything else, since chord diagrams are a special case
of closed Jacobi diagrams, the latter can be considered as linear
combinations of the former, and the two algebras are in any case
isomorphic.

\section{Jacobi diagrams for tangles}
\label{tangl_Jac}
In order to define chord diagrams and, more generally, closed Jacobi
diagrams, for arbitrary tangles it suffices to make only minor
adjustments to the definitions. Namely, one simply replaces the
Wilson loop with an arbitrary oriented one-dimensional manifold (the
{\em skeleton} of the Jacobi diagram). In the 4-term relations the points of
attachment of chords are allowed to belong to different components of the
skeleton, while the STU relations remain the same.

The Vassiliev invariants for tangles with a given skeleton can be
described with the help of chord diagrams or closed diagrams with
the same skeleton; in fact the Vassiliev-Kontsevich Theorem is valid
for tangles and not only for knots.

Open Jacobi diagrams can also be defined for arbitrary tangles. If
we consider tangles whose skeleton is not connected, the legs of
corresponding open diagrams have to be labeled by the connected
components of the skeleton. Moreover, for such tangles there are
mixed spaces of diagrams, some of whose legs are attached to the
skeleton, while others are ``hanging free''. Defining spaces of open
and mixed diagrams for tangles is a more delicate matter than
generalizing chord diagrams: here new relations, called {\em link
relations} may appear in addition to the STU,  IHX and AS relations.

\subsection{Jacobi diagrams for tangles}
\label{section:tang_Jac}

\begin{xdefinition}\index{Jacobi diagram!for tangles}
\index{Tangle!Jacobi diagram} Let $\boldX$ be a tangle skeleton (see
page~\pageref{skeleton}). A {\em tangle closed Jacobi diagram $D$
with skeleton $\boldX$} is a unitrivalent graph with a distinguished
oriented subgraph identified with $\boldX$, a fixed cyclic order of
half-edges at each vertex not on $\boldX$, and such that:
\begin{itemize}
\item it has no univalent vertices other than the boundary points of $\boldX$;
\item each connected component of $D$ contains at least one connected
component of $\boldX$.
\end{itemize}
A tangle Jacobi diagram whose all vertices belong to the skeleton,
is called a {\em tangle chord diagram}.\index{Chord diagram!for
tangles} As with usual closed Jacobi diagrams, half the number of
the vertices of a closed diagram is called the {\em degree}, or {\em
order}, of the diagram.
\end{xdefinition}

\begin{xexample} A tangle diagram whose skeleton consists of a line
segment and a circle:
$$\risS{-22}{tangdiag}{}{60}{20}{25}
$$
\end{xexample}

The vector space of tangle closed Jacobi diagrams with skeleton
$\boldX$ modulo the STU relations is denoted by $\F(\boldX)$, or by
$\F(\xx_1,\ldots,\xx_n)$ where the $\xx_i$ are the connected
components of $\boldX$. The space $\F_n(\boldX)$ is the subspace of
$\F(\boldX)$ spanned by diagrams of degree $n$. It is clear that for
any $\boldX$ the space $\F_n(\boldX)$ is spanned by chord diagrams
with $n$ chords.

Two tangle diagrams are considered to be equivalent if there is a
graph isomorphism between them which preserves the skeleton and the
cyclic order of half-edges at the trivalent vertices outside the
skeleton.

Weight systems of degree $n$ for tangles with skeleton $\boldX$ can
now be defined as linear functions on $\F_n(\boldX)$. The
Fundamental Theorem \ref{fund_thm} extends to the present case:
\begin{xtheorem}\index{Fundamental theorem!for tangles}
Each tangle weight system of degree $n$ is a symbol of some degree
$n$ Vassiliev invariant of framed tangles.
\end{xtheorem}
In fact, we shall prove this, more general version of the
Fundamental Theorem in Chapter~\ref{chapKI} and deduce the
corresponding statement for knots as a corollary.

Now, assume that $\boldX$ is a union of connected components $\xx_i$
and $\yy_j$ and suppose that the $\yy_j$ have no boundary.
\begin{xdefinition}\index{Jacobi diagram!mixed, for tangles}
A {\em mixed tangle Jacobi diagram} is a unitrivalent graph with a
distinguished oriented subgraph (the {\em skeleton}) identified with
$\cup \xx_i$, with all univalent vertices, except those on the
skeleton, labeled by elements of the set $\{\yy_j\}$ and a fixed
cyclic order of edges at each vertex not on the skeleton, and such
that each connected component either contains at least one of the
$\xx_i$, or at least one univalent vertex. A {\em leg} of a mixed
diagram is a univalent vertex that does not belong to the skeleton.
\end{xdefinition}

Here is an example of a mixed Jacobi diagram:

$$\ig{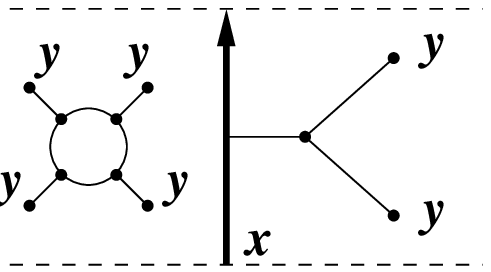}$$

Mixed Jacobi diagrams, apart from the usual STU, IHX and
antisymmetry relations, are subject to a new kind of relations,
called {\em link relations}.
\index{Link relation}\index{Relation!link}
To obtain a link relation, take a mixed
diagram, choose one of its legs and one label $\yy$. For each
$\yy$-labeled vertex, attach the chosen leg to the edge, adjacent
to this vertex, sum all the results and set this sum to be equal to
0. The attachment is done according to the cyclic order as illustrated by the following picture:
$$\rb{-30pt}{\includegraphics[height=1in]{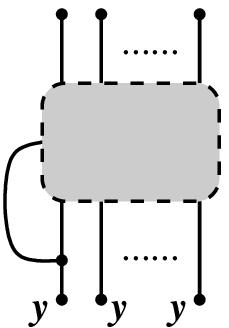}}\ + \
\rb{-30pt}{\includegraphics[height=1in]{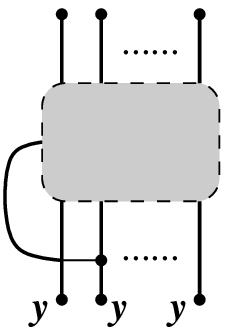}}\ +\dots+\
\rb{-30pt}{\includegraphics[height=1in]{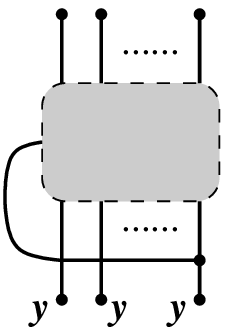}}\ =\ 0.$$
Here the shaded parts of all diagrams coincide, the skeleton is
omitted from the pictures and the unlabeled legs are assumed to
have labels distinct from $\yy$.

 Note that when the skeleton is empty and $\yy$ is the
only label (that is, we are speaking about the usual open Jacobi
diagrams), the link relations are an immediate consequence from the
Kirchhoff law.

Now, define the vector space $\F(\xx_1,\ldots,\xx_n\, |\,
\yy_1,\ldots,\yy_m)$ to be spanned by all mixed diagrams with the
skeleton $\cup\xx_i$ and label in the $\yy_j$, modulo the STU, IHX,
antisymmetry and link relations.\label{sp-mix-J-d}

Both closed and open diagrams are particular cases of this
construction. In particular, $\F(\xx_1,\ldots,\xx_n\, |\,
\emptyset)=\F(\boldX)$ and $\F(\emptyset\, |\, \yy)=\B$. The latter
equality justifies the notation $\B(\yy_1,\ldots,\yy_m)$ or just
$\B(m)$\label{B-of-m} for the space of $m$-coloured open Jacobi diagrams
$\F(\emptyset\, |\, \yy_1,\ldots,\yy_m)$.

Given a diagram $D$ in $\F(\xx_1,\ldots,\xx_n\, |\,
\yy_1,\ldots,\yy_m)$ we can perform ``symmetrization of $D$ with
respect to the label $\yy_m$'' by taking the average of all possible
ways of attaching the $\yy_m$-legs of $D$ to a circle with the label
$\yy_m$. This way we get the map
$$\chi_{\yy_m}:\F(\xx_1,\ldots,\xx_n\, |\, \yy_1,\ldots,\yy_m)\to
\F(\xx_1,\ldots,\xx_n, \yy_m\, |\, \yy_1,\ldots,\yy_{m-1}).
\label{sym-map-C}$$
\begin{xtheorem}
The symmetrization map $\chi_{\yy_m}$ is an isomorphism of vector
spaces.
\end{xtheorem}
In particular, iterating $\chi_{\yy_m}$ we get the isomorphism between
the spaces $\F(\xx_1,\ldots,\xx_n\, |\, \yy_1,\ldots,\yy_m)$ and
$\F(\boldX\cup\boldY)$, where $\boldX=\cup\xx_i$ and $\boldY=\cup\yy_j$.

Let us indicate the idea of the proof; this will also clarify the
origin of the link relations.

Consider the vector space $\F(\xx_1,\ldots,\xx_n\, |\,
\yy_1,\ldots,\yy^*_m)$ defined just like $\F(\xx_1,\ldots,\xx_n\,
|\, \yy_1,\ldots,\yy_m)$ but without the link relations on the
$\yy_m$-legs. Also, define the space $\F(\xx_1,\ldots,\xx_n,
\yy^*_m\, |\, \yy_1,\ldots,\yy_{m-1})$ in the same way as
$\F(\xx_1,\ldots,\xx_n, \yy_m\, |\, \yy_1,\ldots,\yy_{m-1})$ but
with an additional feature that all diagrams have a marked point on
the component $\yy_m$.

Then we have the symmetrization map
$$\chi_{\yy_m^*}: \F(\xx_1,\ldots,\xx_n\, |\,
\yy_1,\ldots,\yy^*_m)\to\F(\xx_1,\ldots,\xx_n, \yy_m^*\, |\,
\yy_1,\ldots,\yy_{m-1})$$ which consists in attaching, in all
possible ways, the $\yy_m$-legs to a pointed circle labeled
$\yy_m$, and taking the average of all the results.

\noindent{\bf Exercise.} Prove
that $\chi_{\yy_m^*}$ is an isomorphism.

Now, consider the map
$$\F(\xx_1,\ldots,\xx_n, \yy_m^*\, |\, \yy_1,\ldots,\yy_{m-1})
\to\F(\xx_1,\ldots,\xx_n, \yy_m\, |\, \yy_1,\ldots,\yy_{m-1})$$ that
simply forgets the marked point on the circle $\yy_m$. The kernel of
this map is spanned by differences of diagrams of the form
$$\rb{-25pt}{\ig[height=2cm]{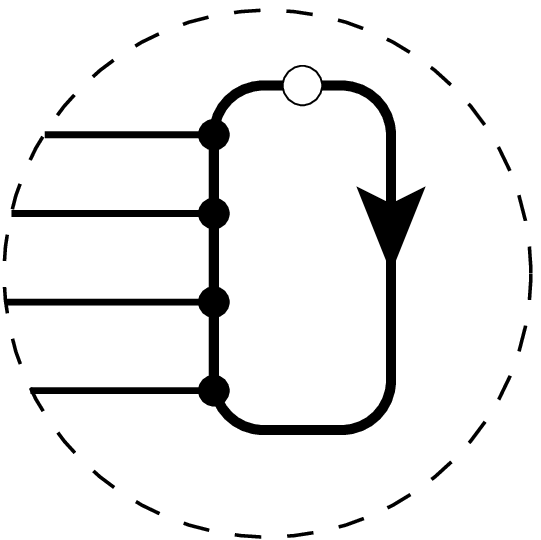}}\quad -\quad
\rb{-25pt}{\ig[height=2cm]{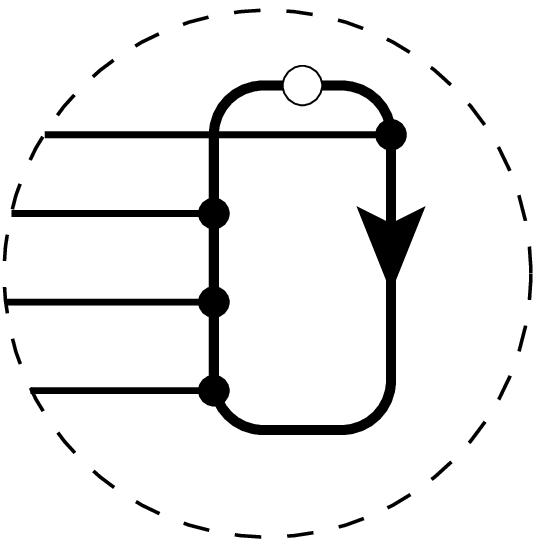}}\ .$$
(The diagrams above illustrate the particular case of 4 legs
attached to the component $\yy_m$.) By the STU relations the above is
equal to the following ``attached link relation'':
$$\rb{-25pt}{\ig[height=2cm]{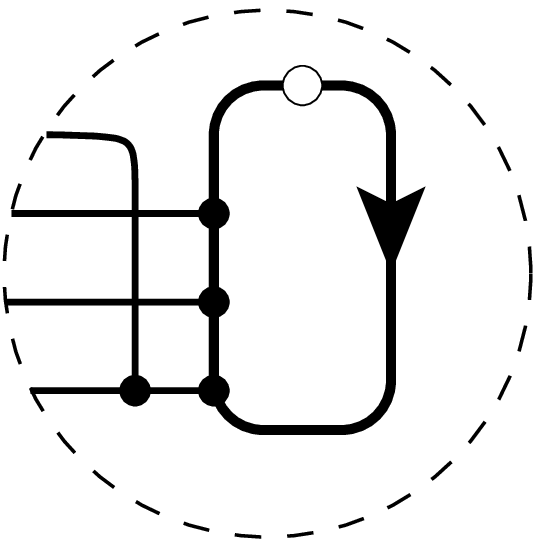}}\quad +\quad
\rb{-25pt}{\ig[height=2cm]{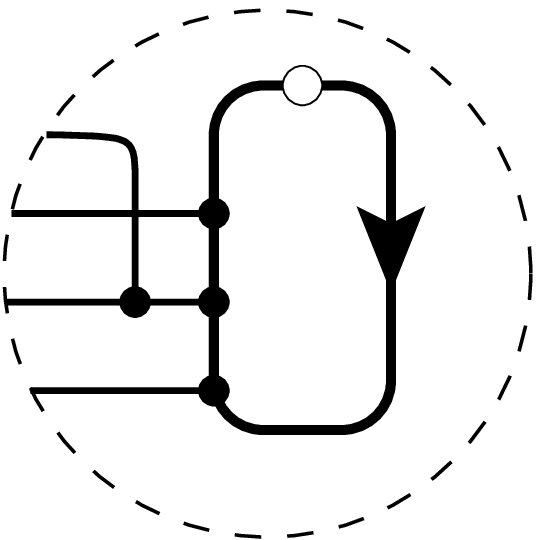}}\quad
+\quad \rb{-25pt}{\ig[height=2cm]{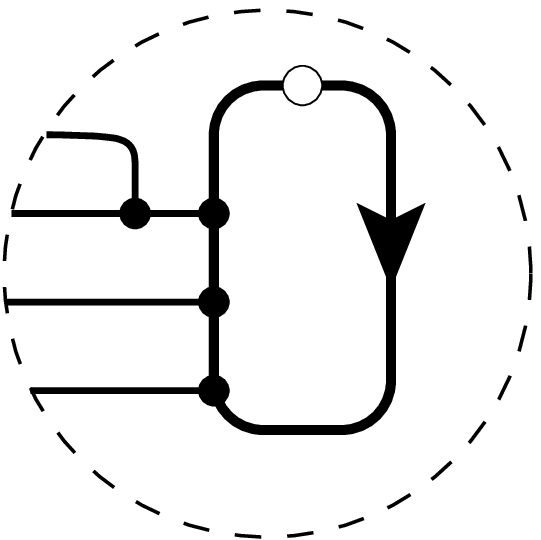}}
.$$

\noindent{\bf Exercise.} Show that the symmetrization map
$\chi_{\yy_m^*}$ identifies the subspace of link relations in
$\F(\xx_1,\ldots,\xx_n\, |\, \yy_1,\ldots,\yy^*_m)$ with the
subspace of $\F(\xx_1,\ldots,\xx_n, \yy_m^*\, |\,
\yy_1,\ldots,\yy_{m-1})$ spanned by all ``attached link relations''.

\subsection{Pairings on diagram spaces}
\label{subsec:pairing}

There are several kinds of pairings on diagram spaces. The first
pairing is induced by the product on tangles; it generalizes the
multiplication in the algebra $\F$. This pairing exists between the
vector spaces $\F(\boldX_1)$ and $\F(\boldX_2)$ such that the bottom
part of $\boldX_1$ coincides with the top part of $\boldX_2$ and
these manifold can be concatenated into an oriented 1-manifold
$\boldX_1\circ\boldX_2$. In this case we have the bilinear map
$$\F(\boldX_1)\otimes\F(\boldX_2)\to\F(\boldX_1\circ\boldX_2),$$
obtained by putting one diagram on top of another. \index{Tangle
diagrams!product}

If $\boldX$ is a collection of $n$ intervals, with one top and one
bottom point on each of them, $\boldX\circ\boldX$ is the same thing
as $\boldX$ and in this case we have an algebra structure on
$\F(\boldX)$. This is the algebra of closed Jacobi diagrams for
string links on $n$ strands. When $n=1$, we, of course, come back to
the algebra $\F$.
\begin{xremark}
While $\F(\boldX)$ is not necessarily an algebra, it is always a
coalgebra with the coproduct defined in the same way as for the
usual closed Jacobi diagrams:
$$\delta(D) := \sum_{J\subseteq[D]} D_J \ot D_{\ol J},$$
where $[D]$ is the set of connected components of the internal graph
of $D$.
\end{xremark}

The second multiplication is the {\em tensor product} of tangle
diagrams. It is induced the tensor product of tangles, and consists
in placing the diagrams side by side. \index{Tangle diagrams!tensor
product}

There is yet another pairing on diagram spaces, which is sometimes
called ``inner product''. For diagrams $C\in\Fxy$ and $D\in\B(\yy)$
define the diagram $\langle C, D\rangle_{\yy}\in \F(\xx)$ as the sum
of all ways of glueing all the $\yy$-legs of $C$ to the $\yy$-legs
of $D$. If the numbers of $\yy$-legs of $C$ and $D$ are not equal,
we set $\langle C, D\rangle_{\yy}$ to be zero. It may happen that in
the process of glueing we get closed circles with no vertices on
them (this happens if $C$ and $D$ contain intervals with both ends
labeled by $\yy$). We set such diagrams containing circles to be
equal to zero.

\begin{lemma}\label{pair-C-B}
The inner product
$$\langle\ ,\, \rangle_{\yy}:\Fxy\otimes\B(\yy)\to\F(\xx)$$
is well-defined.
\end{lemma}
\begin{proof}
We need to show that the class of the resulting diagram in $\F(\xx)$
does not change if we modify the second argument of $\langle \ ,\,
\rangle_{\yy}$ by IHX or antisymmetry relations, and the first
argument --- by STU or link relations. This is clear for the first
three kinds of relations. For link relations it follows from the
Kirchhoff rule and the antisymmetry relation. For example, we have
\begin{multline*}
\rb{-20pt}{\ig[height=1.6cm]{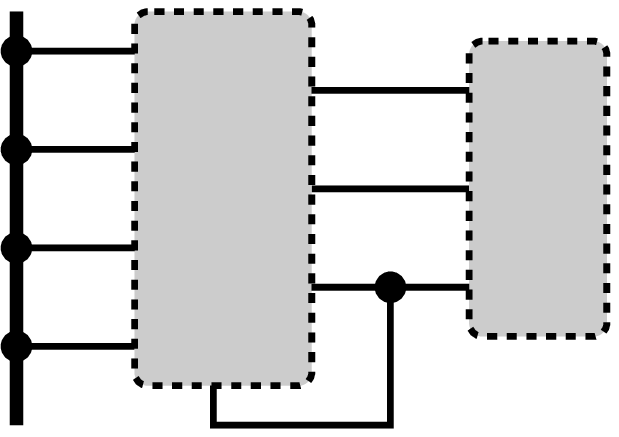}}  \ +\
\rb{-20pt}{\ig[height=1.6cm]{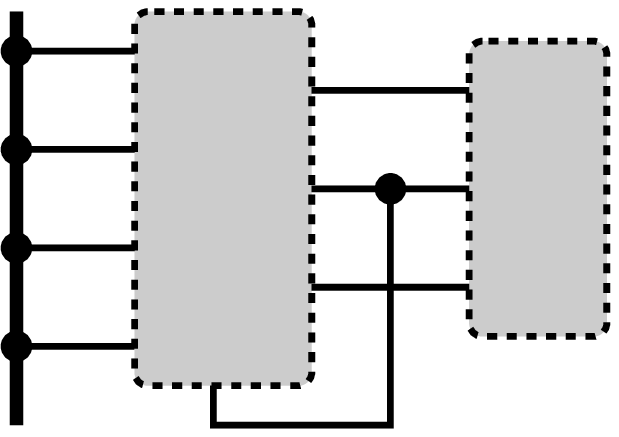}}\ +\
\rb{-20pt}{\ig[height=1.6cm]{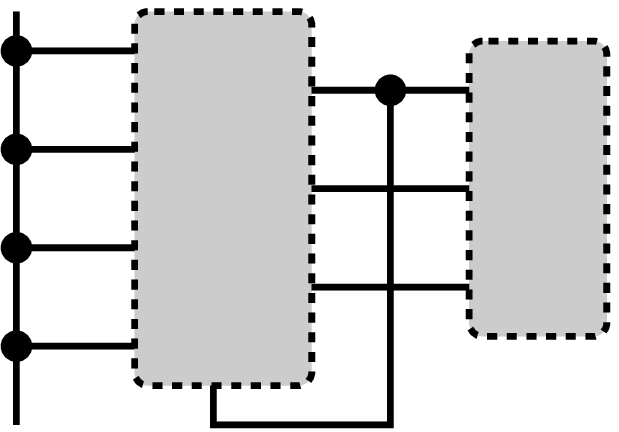}}\quad = \\
=\quad -\ \rb{-20pt}{\ig[height=1.6cm]{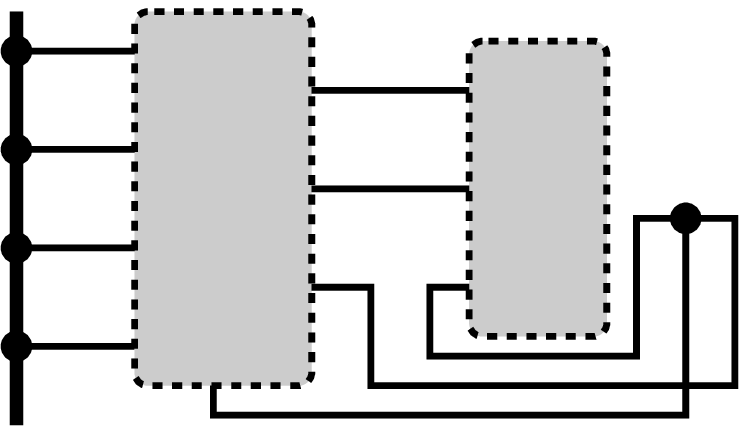}}\ +\
\rb{-20pt}{\ig[height=1.6cm]{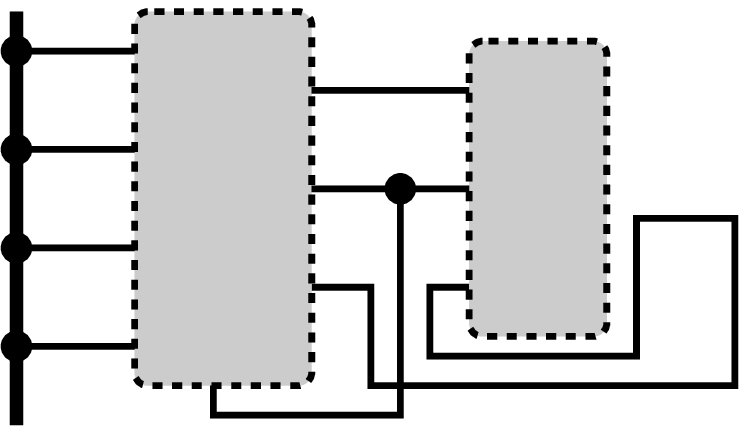}}\ +\
\rb{-20pt}{\ig[height=1.6cm]{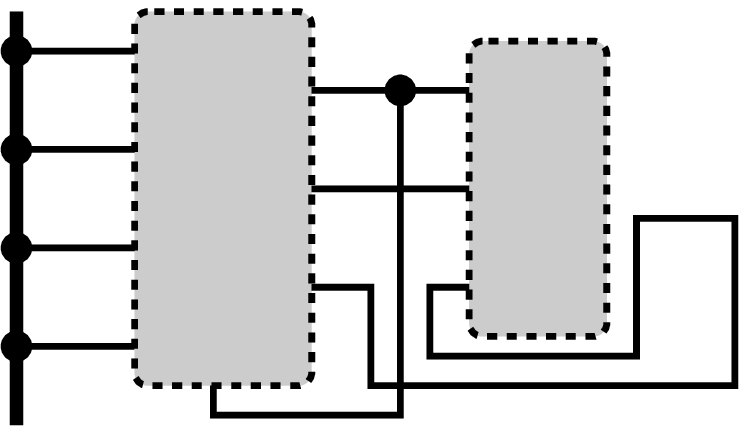}}\quad =
\quad 0.
\end{multline*}
\end{proof}

The definition of the inner product can be extended. For example, if
two diagrams $C,D$ have the same number of $\yy_1$-legs and the same
number of $\yy_2$-legs, they can be glued together along the
$\yy_1$-legs and then along the $\yy_2$-legs. The sum of the results
of all such glueings is denoted by $\langle C ,
D\rangle_{\yy_1,\yy_2}$. This construction, clearly, can be
generalized further.

\subsection{Actions of $\F$ and $\B$ on tangle
diagrams}\label{subsec:modules-tangle} While the coalgebra
$\F(\boldX)$, in general, does not have a product, it carries an
algebraic structure that generalizes the product in $\F$. Namely,
for each component $\xx$ of $\boldX$, there is an action of
$\F(\xx)$ on $\F(\boldX)$, defined as the connected sum along the
component $\xx$. We denote this action by $\#$, as if it were the
usual connected sum. More generally, the spaces of mixed tangle
diagrams $\F(\xx_1, \ldots, \xx_n\, |\, \yy_1, \ldots, \yy_m )$ are
two-sided modules over $\F(\xx_i)$ and $\B(\yy_j)$. The algebra
$\F(\xx_i)$ acts, as before, by the connected sum on the component
$\xx_i$, while the action of $\B(\yy_j)$ consists in taking the
disjoint union with diagrams in $\B(\yy_j)$. We shall denote the
action of $\B(\yy_j)$ by $\cup$.

We cannot expect the relation of the module structures on the space of
mixed diagrams with the symmetrization map to be straightforward,
since the symmetrization map from $\B$ to $\F$ fails to be
multiplicative. We shall clarify this remark in
\ref{subsec:generalwheel}.

\noindent{\bf Exercise.}  Prove that the above actions are
well-defined. In particular, prove that the action of $\F(\xx_i)$
does not depend on the location where the diagram is inserted into
the corresponding component of the tangle diagram, and show that the
action of $\B(\yy_j)$ respects the link relations.

\subsection{Sliding property}
\label{sec:sliding}

There is one important corollary of the IHX relation (Kirchhoff
law), called {\em sliding property} (\cite{BLT}), which holds in the
general context of tangle Jacobi diagrams. To formulate it, we need
to define the operation
$\D_{\xx}^{(n)}:\F(\xx\cup\boldY)\to\F(\xx_1\cup\dots\cup\xx_n\cup\boldY)$.
By definition, $\D_{\xx}^{(n)}(D)$ is the lift of $D$ to the $n$th
disconnected cover of the line $\xx$, that is, for each $\xx$-leg of the
diagram $D$ we take the sum over all ways to attach it to $\xx_i$
for any $i=1,\dots,n$ (the sum consists of $n^k$ terms, if $k$ is
the number of vertices of $D$ belonging to $\xx$). Example:
$$\D_{\xx}^{(2)}\Big(\,\risS{-12}{sl_ex0}{
      \put(4,30){$\scriptstyle \xx$}}{27}{20}{15}\ \Big)
  \quad=\quad \risS{-12}{sl_ex1}{
      \put(-7,32){$\scriptstyle \xx_1$}
      \put(8,32){$\scriptstyle \xx_2$}}{27}{0}{0}\ +\
  \risS{-12}{sl_ex2}{
      \put(-7,32){$\scriptstyle \xx_1$}
      \put(8,32){$\scriptstyle \xx_2$}}{27}{0}{0}\ +\
  \risS{-12}{sl_ex3}{
      \put(-7,32){$\scriptstyle \xx_1$}
      \put(8,32){$\scriptstyle \xx_2$}}{27}{0}{0}\ +\
  \risS{-12}{sl_ex4}{
      \put(-7,32){$\scriptstyle \xx_1$}
      \put(8,32){$\scriptstyle \xx_2$}}{27}{0}{0}\ .
$$

\begin{xproposition} (Sliding relation)
\label{prop:sliding}\index{Sliding relation}\index{Relation!sliding}
Suppose that $D\in\F(\xx\cup\boldY)$; let $D_1=\D_{\xx}^{(n)}(D)$.
Then for any diagram $D_2\in\F(\xx_1\cup\dots\cup\xx_n)$ we have
$D_1D_2=D_2D_1$. In pictures:
$$\risS{-35}{sliding1}{
      \put(0,-4){$\scriptstyle \xx_1\cdots\xx_2$}
      \put(20,12){$D_1$}\put(5,50){$D_2$}
      \put(35,34){$\scriptstyle \boldY$}}{52}{30}{40}
  \quad=\quad
\risS{-35}{sliding2}{
      \put(0,-4){$\scriptstyle \xx_1\cdots\xx_2$}
      \put(20,50){$D_1$}\put(5,13){$D_2$}
      \put(35,30){$\scriptstyle \boldY$}}{52}{0}{0}
$$
\end{xproposition}

\begin{proof}
Indeed, take the leg in $D_1$ which is closest to $D_2$ and consider
the sum of all diagrams on $\xx_1\cup\dots\cup\xx_n\cup\boldY$ where
this leg is attached to $\xx_i$, $i=1,\dots,n$, while all the other
legs are fixed. By Kirchhoff law, this sum is equal to the similar
sum where the chosen leg has jumped over $D_2$. In this way, all the
legs jump over $D_2$ one by one, and the commutativity follows.
\end{proof}

\subsection{Closing a component of a Jacobi diagram} \label{sec:Jac_links}

Recall that long knots can be closed up to produce usual knots. This
closure induces a bijection of the corresponding isotopy classes and
an isomorphism of the corresponding diagram spaces.

This fact can be generalized to tangles whose skeleton consists of
one interval and several circles.

\begin{xtheorem}
Let $\boldX$ be a tangle skeleton with only one interval component,
and $\boldX'$ be a skeleton obtained by closing this component into
a circle. The induced map $$\F(\boldX)\to\F(\boldX')$$ is an
isomorphism of vector spaces.
\end{xtheorem}

The proof of this theorem consists in applying the
Kirchhoff's law and we leave it to the reader.

We should point out that closing one component of a skeleton with
more that one interval component does not produce an isomorphism of
the corresponding diagram spaces. Indeed, let us denote by $\A(2)$
the space of closed diagrams for string links on 2 strands. A direct
calculation shows that the two diagrams of order 2 below on the left
are different in $\A(2)$, while their images under closing one
strand of the skeleton are obviously equal:
$$\risS{-12}{linkd3}{}{16}{15}{15}\ \ne\
\risS{-12}{linkd4}{}{16}{0}{0}\hspace{3cm}
\risS{-12}{linkd2}{}{27}{0}{0}\ =\ \risS{-12}{linkd1}{}{27}{0}{0}
$$

The above statements about tangle diagrams, of course, are not
arbitrary, but reflect the following topological fact that we state
as an exercise:

\noindent{\bf Exercise.} Define the map of closing one component on
isotopy classes of tangles with a given skeleton and show that it is
bijective if and only if it is applied to tangles whose skeleton has
only one interval component.

\section{Horizontal chord diagrams}

There is yet another diagram algebra which will be of great
importance in what follows, namely, the algebra $\A^h(n)$ of {\em
horizontal \index{Algebra!horizontal chord diagrams}\label{A^h(n)}
chord diagrams} on $n$ strands.

A horizontal chord diagram on $n$ strands is a tangle diagram whose
skeleton consists of $n$ vertical intervals (all oriented, say,
upwards) and all of whose chords are horizontal. Two such diagrams
are considered to be equivalent if one can be deformed into the
other through horizontal diagrams.

A product of two horizontal diagrams is clearly a horizontal
diagram; by definition, the algebra $\A^h(n)$ is generated by the
equivalence classes of all such diagrams modulo the horizonal 4T
relations \ref{horizontal4T} (see Section \ref{4Tsec}). We denote by
${\bf 1}_n$ the the empty diagram in $\A^h(n)$ which is the
multiplicative unit.

Each horizontal chord diagram is equivalent to a diagram whose
chords are all situated on different levels, that is, to a product
of diagrams of degree 1. Set
$$
  u_{jk} = \ \rb{-8mm}{\ig[height=15mm]{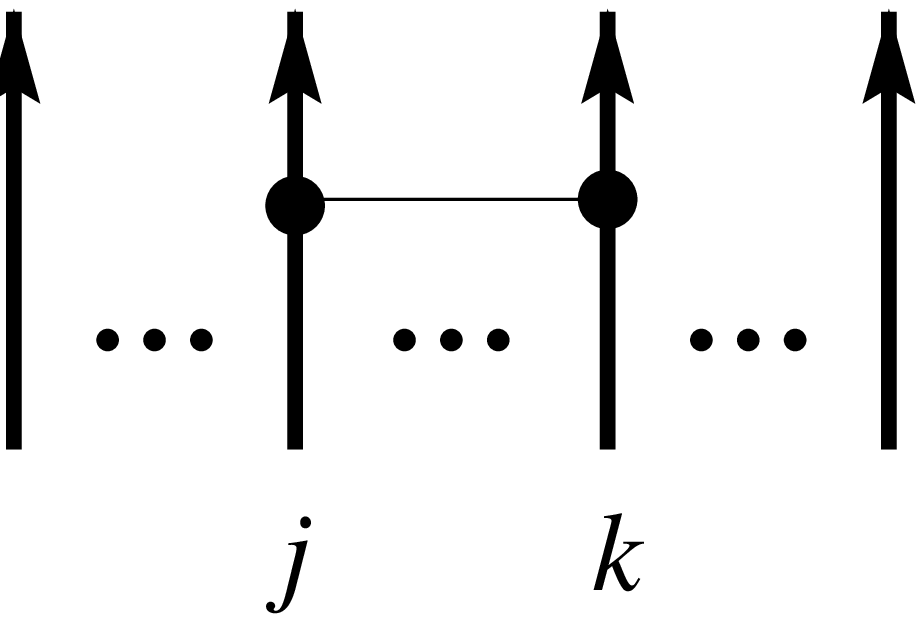}}\ ,\qquad
         1\leqslant j<k\leqslant n,
$$
and for $1\leqslant k<j\leqslant n$ set $u_{jk}=u_{kj}$. Then
$\A^h(n)$ is generated by the $u_{jk}$ subject to the following
relations ({\em infinitesimal pure braid relations}, first appeared
in \cite{Koh2})
\begin{align*}
   &[u_{jk},u_{jl}+u_{kl}] =0,\quad \text{ if } j,k,l \text{ are different},\\
   &[u_{jk},u_{lm}] =0, \quad \text{ if } j,k,l,m \text{ are
   different}.
\end{align*}
Indeed, the first relation is just the horizontal 4T relation. The
second relation is similar to the far commutativity relation in
braids. The products of the $u_{jk}$ up to this relation are
precisely the equivalence classes of horizontal diagrams.

The algebra $\A^h(2)$ is simply the free commutative algebra on one
generator $u_{12}$.
\begin{proposition}\label{prop:a3}
$\A^h(3)$ is a direct product of the free algebra on two generators
$u_{12}$ and $u_{23}$, and the free commutative algebra on one
generator $$u=u_{12}+u_{23}+u_{13}.$$
\end{proposition}
In particular, $\A^h(3)$ is highly non-commutative.
\begin{proof}
Choose $u_{12}$, $u_{23}$ and $u$ as the set of generators for
$\A^h(3)$. In terms of these generators all the relations in
$\A^h(3)$ can be written as
$$[u_{12}, u]=0, \quad \text{and} \quad
[u_{23}, u]=0\ .$$
\end{proof}

For $n>3$ the multiplicative structure of the algebra $\A^h(n)$ is
rather more involved, even though it admits a simple description as
a vector space. We shall treat this subject in more detail in
Chapter~\ref{chapBr}, as the algebra $\A^h(n)$ plays the same role
in the theory of finite type invariants for pure braids as the
algebra $\A$ in the theory of the Vassiliev knot invariants.

We end this section with one property of $\A^h(n)$ which will be
useful in Chapter~\ref{DA_chap}.

\begin{lemma}\label{lemma:hordiag_center}
Let $J,K\subseteq \{1,\ldots, n\}$ be two non-empty subsets with
$J\cap K=\emptyset$. Then the element $\sum_{j\in J, k\in K} u_{jk}$
commutes in $\A^h(n)$ with any generator $u_{pq}$ with $p$ and $q$
either both in $J$ or both in $K$.
\end{lemma}

\begin{proof}
It is clearly sufficient to prove the lemma for the case when $K$
consists of one element, say $k$, and both $p$ and $q$ are in $J$.
Now, any $u_{jk}$ commutes with $u_{pq}$ if $j$ is different from
both $p$ and $q$. But $u_{pk}+u_{qk}$ commutes with $u_{pq}$ by the
horizontal 4T relation, and this proves the lemma.
\end{proof}

\subsection{Horizontal diagrams and string link diagrams}
Denote by $\A(n)$\label{A-ot-p} be the algebra of closed diagrams for string
links. Horizontal diagrams are examples of string link diagrams and
the horizontal 4T relations are a particular case of the usual 4T
relations, and, hence, there is an algebra homomorphism
$$\A^{h}(n)\to \A(n).$$
This homomorphism is injective, but this is a surprisingly
non-trivial fact; see \cite{BN8, HM}. We shall give a proof of this
in Chapter~\ref{chapBr}, see page~\pageref{theorem:stringlinx}.

\noindent{\bf Exercise.}
\begin{itemize}
\item[(a)] Prove that the chord diagram consisting of one chord
connecting the two components of the skeleton belongs to the centre
of the algebra $\A(2)$.

\item[(b)] Prove that any chord diagram consisting of two intersecting
chords belongs to the centre of the algebra $\A(2)$.

\item[(c)] Prove that Lemma~\ref{lemma:hordiag_center} is also valid
for $\A(n)$. Namely, show that the element  $\sum_{j\in J, k\in K}
u_{jk}$ commutes in $\A(n)$ with any chord diagram whose chords have
either both ends on the strands in $J$ or on the strands in $K$.
\end{itemize}

\begin{xcb}{Exercises}
\begin{enumerate}

\item
Prove that $\kol{circdiag6}=\kol{circdiag1}$.

\item
Let\quad $a_1=\kol{circdiag1}\ ,\ a_2=\kol{circdiag2}\ ,\
    a_3=\kol{circdiag3}\ ,\ a_4=\kol{circdiag4},\ a_5=\kol{circdiag5}$.

\begin{itemize}
\item[(a)] Find a relation between $a_1$ and $a_2$.

\item[(b)] Represent the sum $a_3+a_4-2a_5$ as a connected closed diagram.

\item[(c)] Prove the linear independence of $a_3$ and $a_4$ in $\F$.
\end{itemize}

\item
Express the primitive elements $\chdh{f25_Tc}$,
$\chdh{f25_Fe}$ and $\chd{f25_Ff}$ of degrees 3 and 4 as linear
combinations of chord diagrams.

\item\label{A2C}
Prove the following identities in the algebra $\F$:
$$\begin{array}{ccl}
\chd{c25_22} &=&\displaystyle \chd{c25_21}-\frac12\ \chdh{f25_Wb}\ ;\\
\chdh{c25_34} &=&\displaystyle \chdh{c25_31} - \chdh{f25_Tb}
                  +\frac14\ \risS{-10}{f25_Tc}{}{28}{20}{15}\ ;\\
\chd{cd4-04} &=&\displaystyle \chd{cd4-08}-\frac32\ \chdh{f25_Fb}
               +\frac34\ \risS{-10}{f25_Fc}{}{28}{20}{15}
               -\frac18\ \risS{-10}{f25_Fe}{}{28}{20}{15}\ ;\\
\chd{cd4-06} &=&\displaystyle \chd{cd4-08}-\frac32\ \chdh{f25_Fb}
   +\frac12\ \risS{-10}{f25_Fc}{}{28}{20}{15}+\frac14\ \chdh{f25_Fd}
   -\frac18\ \risS{-10}{f25_Fe}{}{28}{20}{15}\ ;\\
\chd{cd4-02} &=&\displaystyle \chd{cd4-08}-2\ \chdh{f25_Fb}
   +\risS{-10}{f25_Fc}{}{28}{20}{15}+\frac12\ \chdh{f25_Fd}
   -\frac12\ \risS{-10}{f25_Fe}{}{28}{20}{15}+\chd{f25_Ff}\ .\\
\end{array}$$

\item
Show that the symbols of the coefficients of the Conway polynomial
\index{Conway polynomial!symbol} (Section~\ref{conway}) take the
following values on the basis primitive diagrams of degree 3 and 4.
$$\begin{array}{c} \symb(c_3)\Bigl(\chdh{f25_Tc}\Bigr) = 0,\\
  \symb(c_4)\Bigl(\chdh{f25_Fe}\Bigr) = 0,\qquad\qquad
  \symb(c_4)\Bigl(\chdh{f25_Ff}\Bigr) = -2.\end{array}
$$

\item
Show that the symbols of the coefficients of the Jones polynomial
\index{Jones polynomial!symbol} (Section~\ref{symb_j_n}) take the
following values on the basis primitive diagrams of degrees 3 and 4.
$$\begin{array}{c} \symb(j_3)\Bigl(\chdh{f25_Tc}\Bigr) = -24,\\
  \symb(j_4)\Bigl(\chdh{f25_Fe}\Bigr) = 96,\qquad\qquad
  \symb(j_4)\Bigl(\chdh{f25_Ff}\Bigr) = 18.\end{array}
$$

\item \label{forest_el}
\parbox[t]{3in}{(\cite{CV}) Let $\ol{t_n}\in \PR_{n+1}$ be the closed diagram
shown on the right. Prove the following identity
$$\ol{t_n}\quad =\quad \frac{1}{2^n}\ \risS{-20}{forest_bn}{
         \put(5,-8){\mbox{\scriptsize $n$ bubbles}}}{40}{20}{25}
$$}\qquad
\parbox[t]{1in}{\rb{-10pt}{$\ol{t_n} =\quad \risS{-15}{wssl2tn}{
         \put(5,-8){\mbox{\scriptsize $n$ legs}}}{30}{20}{15}$}}

\noindent
Deduce that $\ol{t_n}\in \PR_{n+1}^2$.

\item
Express $\ol{t_n}$ as a linear combination of chord diagrams. In
particular, show that the intersection graph of every chord diagram
that occurs in this expression is a forest.

\item (\cite{CV})
Prove the following identity in the space $\F$ of closed diagrams:
$$\risS{-17}{cvpent1}{}{37}{20}{20}\quad =\quad
  \frac{3}{4}\ \risS{-17}{cvpent2}{}{37}{20}{20}\ \ -\
 \frac{1}{12}\ \lClD{cvpent3}{}\ -\ \frac{1}{48}\ \lClD{cvpent4}{}\ .
$$

{\sl Hint.} Turn the internal pentagon of the left-hand side of the
identity %as a rigid object
in the 3-space by $180^\circ$ about the vertical axis. The result
will represent the same graph with the cyclic orders at all five
vertices of the pentagon changed to the opposite:\bigskip
$$\risS{-22}{cvpent5}{}{37}{20}{20}\ =\
  (-1)^5\ \risS{-17}{cvpent6}{}{37}{20}{20}\ =\
  -\ \risS{-17}{cvpent1}{}{37}{20}{20}\ +\
  (\mbox{\scriptsize terms\ with\ at\ most\ 4\ legs})\ .
$$
The last equality follows from the $\STU$ relations which allow us
to rearrange the legs modulo diagrams with a smaller number of legs.
To finish the solution, the reader must figure out the terms in the
parentheses.

\item
Prove the linear independence of the three elements in the
right-hand side of the last equality, using Lie algebra invariants
defined in Chapter \ref{LAWS}.

\item (\cite{CV})
Prove that the primitive space in the algebra $\F$ is generated by
the closed diagrams whose internal graph is a tree.

\item (\cite{CV})
With each permutation $\sigma$ of $n$ objects associate a closed
diagram $P_\sigma$ acting as in Section~\ref{perm_act_on_C} by the
permutation on the lower legs of a closed diagram $P_{(12\dots
n)}=\ol{t_n}$ from problem \ref{forest_el}. Here are some examples:
$$P_{(2143)}\ =\ \sClD{p_2143}\ ;\qquad
P_{(4123)}\ =\ \sClD{p_4123}\ ;\qquad
P_{(4132)}\ =\ \sClD{p_4132}\ .
$$
Prove that the diagrams $P_\sigma$ span the vector space
$\PR_{n+1}$.

\item(\cite{CV}) Prove that
\begin{itemize}
\item\label{ex_filt_pr}   $\PR_n^n=\PR_n$ for
even $n$,
        and $\PR_n^{n-1}=\PR_n$ for odd $n$;
\item
for even $n$ the quotient space $\PR_n^n / \PR_n^{n-1}$
has dimension one and generated by the wheel $\ol{w}_n$.
\end{itemize}

\item
\newcommand{\kolo}[1]{\rb{-4mm}{\ig[height=10mm]{#1.eps}}}
Let\qquad $b_1=\kolo{freediag0}\ ,\quad b_2=\kolo{freediag1}\ ,\quad
     b_3=\kolo{freediag2}\ ,\quad b_4=\kolo{c25_Ff}\ .$\vspace{2mm}
Which of these diagrams are zero in $\B$, that is, vanish modulo AS
and IHX relations?

\item\label{strangeIHX}
Prove that the algebra generated by all open diagrams modulo the AS and
the modified IHX equation
$I=a H- b X$, where $a$ and $b$ are arbitrary complex numbers,
is isomorphic (equal) to $\B$ if and only if $a=b=1$, in all other cases
it is a free polynomial algebra on one generator.

\item
\begin{itemize}
  \item
Indicate an explicit form of the isomorphisms
$\A^{fr}\isom\F\isom\B$ in the bases given in
Section~\ref{ABC_sm_deg}.
  \item
Compile the multiplication table for $\B_m\times\B_n\rightarrow\B_{m+n}$,
$m+n\le4$, for the second product in $\B$ (the one pulled back from
$\F$ along the isomorphism $\F\isom\B$).
  \item
Find some bases of the spaces $\A^{fr}_n$, $\F_n$, $\B_n$ for $n=5$.
\end{itemize}

\item\label{caterpil}
({\it J. Kneissler}). Let $\B_{n}^{u}$ be the space of open diagrams
of degree $n$ with $u$ univalent vertices. Denote by
$\om_{i_1i_2\dots i_k}$ the element of
$\B_{i_1+\dots+i_k+k-1}^{i_1+\dots+i_k}$ represented by a
{\em caterpillar} diagram consisting of $k$ body segments with $i_1$,
\dots, $i_k$ ``legs", respectively.\index{Diagram!caterpillar} Using
the AS and IHX relations, prove that $\om_{i_1i_2\dots i_k}$ is
well-defined, that is, for inner segments it makes no difference on
which side of the body they are drawn. For example,
$$\omega_{0321}\quad=\quad\risS{-20}{om0321-1}{}{65}{20}{20}
               \quad=\quad\risS{-20}{om0321-2}{}{65}{20}{12}
$$

\item\neresh
({\it J. Kneissler})
Is it true that any caterpillar diagram in the algebra $\B$ can be
expressed through caterpillar diagrams with even indices $i_1$, \dots,
$i_k$?
Is it true that the primitive space $\PR(\B)$ (that is, the space spanned
by connected open diagrams) is generated by caterpillar diagrams?

\item
Prove the equivalence of the two claims:
\begin{itemize}
\setlength{\itemsep}{1pt plus 1pt minus 1pt}
\item
all chord diagrams are symmetric modulo one- and four-term
relations.
\item
all chord diagrams are symmetric modulo only four-term relations.
\end{itemize}

\item
Similarly to symmetric chord diagrams
(page~\pageref{def:tau-rever}), we can speak of {\em anti-symmetric}
\index{Chord diagram!anti-symmetric} diagrams: an element $D$ of
$\A$ or $\A^{fr}$ is anti-symmetric if $\tau(D)=-D$. Prove that
under the isomorphism $\chi^{-1}:\A^{fr}\to\B$:
\begin{itemize}
\item
 the image of a symmetric chord diagram is a linear
combination of open diagrams with an even number of legs,
\item
 the image of an anti-symmetric chord diagram in is a linear
combination of open diagrams with an odd number of legs.
\end{itemize}

\item\neresh\label{threelegs}
(The simplest unsolved case of Conjecture \ref{conj:oddlegs}).
Is it true that an open diagram with 3 univalent vertices is always
equal to 0 as an element of the algebra $\B$?

\item
Prove that the diagram \rb{-7mm}{\ig[height=15mm]{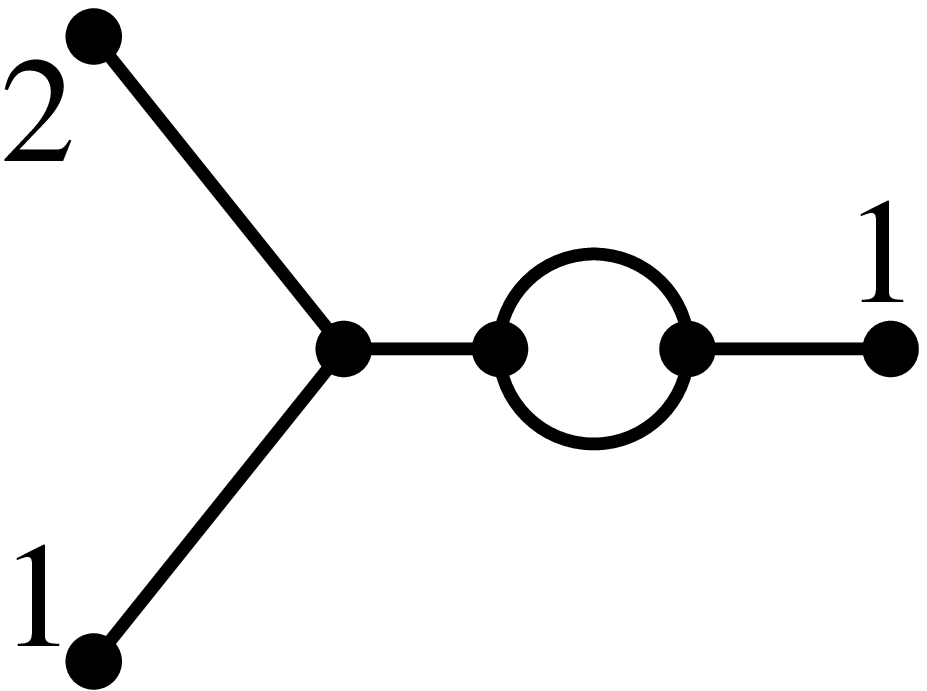}} is
equal to 0 in the space $\B(2)$.

\item
Let $u_{ij}$ be the diagram in $\A(3)$ with one chord connecting the
$i$th and the $j$th component of the skeleton. Prove that for any
$k$ the combination $u_{12}^k+u_{23}^k+u_{13}^k$ belongs to the
centre of $\A(3)$.

\item\label{ex_C2B_34leg}
Let $D^3 = \risS{-7}{mixd3l1}{
         \put(30,0){\mbox{$\yy$}}}{28}{20}{15}\ \ \in \F(\boldX,\yy)$ and
$D^4 = \risS{-7}{mixd4l1}{
         \put(30,0){\mbox{$\yy$}}}{28}{20}{15}\ \ \in \F(\boldX,\yy)$
be tangle diagrams with exactly three and four $\yy$-legs respectively.
Show that
\def\skry{\mbox{$\scriptstyle \yy$}}
\def\mdto{\risS{-7}{mixd3l2}{\put(4,-5){\skry}\put(12,-5){\skry}\put(20,-5){\skry}
                    }{28}{20}{15}\!\!}
\def\mddos#1{\risS{-7}{#1}{\put(8,-5){\skry}\put(18,-5){\skry}}{28}{20}{15}\!\!}
\def\mdft#1{\risS{-7}{#1}{\put(4,-5){\skry}\put(18,-5){\skry}}{28}{20}{15}\!\!}
$$\begin{array}{l}
\F(\boldX|\yy)\ni\chi_{\yy}^{-1}(D^3) = \mdto + \frac12\! \mddos{mixd3l3}\ =
   \mdto + \frac12\! \mddos{mixd3l4}\ = \mdto - \frac12\! \mddos{mixd3l5}\ ; \\
\F(\boldX|\yy)\ni\chi_{\yy}^{-1}(D^4) =
  \risS{-7}{mixd4l2}{\put(2,-5){\skry}\put(8,-5){\skry}\put(14,-5){\skry}
                     \put(20,-5){\skry}}{28}{20}{15}
  + \frac12\! \risS{-7}{mixd4l3}{\put(7,-5){\skry}\put(13,-5){\skry}
                     \put(19,-5){\skry}}{28}{20}{15}\!\!
  + \frac12\! \risS{-7}{mixd4l4}{\put(4,-5){\skry}\put(10,-5){\skry}
                     \put(16,-5){\skry}}{28}{20}{15}\!\!
  + \frac18\! \mdft{mixd4l5} + \frac5{24}\! \mdft{mixd4l6}\ .
\end{array}
$$

{\sl Hint.} Follow the proof of Theorem \ref{BFisom} on page \pageref{def_tau_k} and
then use link relations.
\end{enumerate}
\end{xcb}
 %5 Cl.Op.Diag
\chapter{Lie algebra weight systems} %6
\label{LAWS}
\index{Lie algebra!weight systems}

Given a Lie algebra $\g$ equipped with a non-degenerate invariant
bilinear form, one can construct a weight system with values in the
centre of the universal enveloping algebra $U(\g)$. In a similar
fashion one can define a map from the space $\B$ into the
ad-invariant part of the symmetric algebra $S(\g)$. These
constructions are due to M.~Kontsevich \cite{Kon1}, with basic ideas
already appearing in \cite{Pen}. If, in addition, we have a finite
dimensional representation of the Lie algebra then taking the trace
of the corresponding operator we get a numeric weight system. It
turns out that these weight systems are the symbols of the quantum
group invariants (Section~\ref{qift}). The construction of weight
systems based on representations first appeared in D.~Bar-Natan's
paper \cite{BN0}. The reader is invited to consult the Appendix for
basics on Lie algebras and their universal envelopes.

A useful tool to compute Lie algebra weight systems is Bar-Natan's
computer program called {\tt main.c} and available online at
\cite{BN5}. The tables in this chapter were partially obtained using
that program.

There is another construction of weight systems, also invented by
Kontsevich: the weight systems coming from {\em marked surfaces}. As
proved in \cite{BN1}, this construction gives the same set of weight
systems as the classical Lie algebras, and we shall not speak about
it here.

\section{Lie algebra weight systems for the algebra $\A^{fr}$}
\label{LAWS_A}
\index{Weight system! Lie algebra}

\subsection{Universal Lie algebra weight systems}
\label{univ_lie_ws}

Kontsevich's construction proceeds as follows. Let $\g$ be a
metrized Lie algebra over $\R$ or $\C$, that is, a Lie algebra with
an ad-invariant non-degenerate bilinear form
$\langle\cdot,\cdot\rangle$ (see \ref{metrized}). 
Choose a basis
$e_1,\dots,e_m$ of $\g$ and let $e^*_1,\dots,e^*_m$ be the dual
basis with respect to the form $\langle\cdot,\cdot\rangle$.

Given a chord diagram $D$ with $n$ chords, we first choose a base
point on its Wilson loop, away from the chords of $D$. This gives a
linear order on the endpoints of the chords, increasing in the
positive direction of the Wilson loop. Assign to each chord $a$ an
{\em index}, that is, an integer-valued variable, $i_a$. The values
of $i_a$ will range from 1 to $m$, the dimension of the Lie algebra.
Mark the first endpoint of the chord with the symbol $e_{i_a}$ and
the second endpoint with $e^*_{i_a}$.

Now, write the product of all the $e_{i_a}$ and all the $e^*_{i_a}$,
in the order in which they appear on the Wilson loop of $D$, and
take the sum of the $m^n$ elements of the universal enveloping
algebra $U(\g)$ obtained by substituting all possible values of the
indices $i_a$ into this product. Denote by $\f_\g(D)$ the resulting
element of $U(\g)$.\label{def:phiA}

For example,
$$\f_\g\bigl(\cdO\bigr) = \sum_{i=1}^m e_i e^*_i =: c
$$
is the {\em quadratic Casimir element}\label{QuaCa} \index{Casimir
element} associated with the chosen invariant form. The next theorem
shows, in particular, that the Casimir element does not depend on
the choice of the basis in $\g$. Another example: if
$$
   D\ =\ \rb{-6.5mm}{\ig[width=15mm]{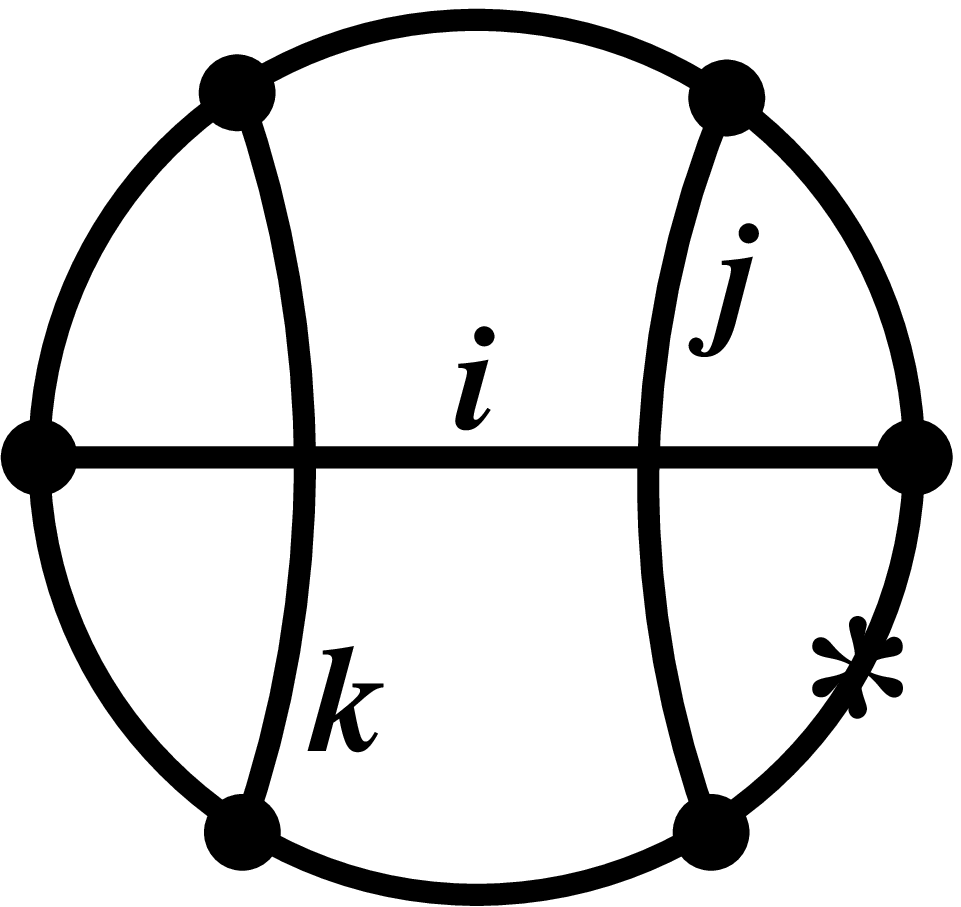}}\ ,
$$
then
$$
   \f_\g(D)=\sum_{i=1}^m\sum_{j=1}^m\sum_{k=1}^m
    e^{}_ie^{}_je^{}_ke^*_ie^*_ke^*_j.
$$

\begin{theorem}\label{phi_prop}
The above construction has the following properties:
\begin{enumerate}
\setlength{\itemsep}{1pt plus 1pt minus 1pt}
\item the element $\f_\g(D)$ does not depend on the choice of the base point
on the diagram;
\item it does not depend on the choice of the basis $\{e_i\}$ of the Lie
algebra;
\item it belongs to the ad-invariant subspace
$$U(\g)^\g=\{x\in U(\g)\mid  xy=yx\ \text{for\ all}\ y\in\g\}$$
of the universal enveloping algebra $U(\g)$ (that is, to the centre
$ZU(\g)$);
\item the function $D\mapsto\f_\g(D)$ satisfies 4-term relations;
\item the resulting map $\f_\g:\A^{fr}\to ZU(\g)$ is a homomorphism of
algebras.
\end{enumerate}
\end{theorem}

\begin{proof}
(1) Introducing a base point means that a circular chord diagram is
replaced by a linear chord diagram (see Section \ref{line_CD}).
Modulo 4-term relations, this map is an isomorphism, and, hence, the
assertion follows from (4).

(2) An exercise in linear algebra: take two different bases
$\{e_i\}$ and $\{f_j\}$ of $\g$ and reduce the expression for
$\f_\g(D)$ in one basis to the expression in another using the
transition matrix between the two bases. Technically, it is enough
to do this exercise only for $m=\dim\g=2$, since the group of
transition matrices $GL(m)$ is generated by linear transformations
in the 2-dimensional coordinate planes.

This also follows from the invariant construction of this weight
system in Section \ref{lie_alg_ws_without_bases} which does not use
any basis.

(3) It is enough to prove that $\f_\g(D)$ commutes with any basis
element $e_r$. By property (2), we can choose the basis to be
orthonormal with respect to the ad-invariant form
$\langle\cdot,\cdot\rangle$, so that $e_i^*=e_i$ for all $i$. Now,
the commutator of $e_r$ and $\f_\g(D)$ can be expanded into a sum of
$2n$ expressions, similar to $\f_\g(D)$, only with one of the $e_i$
replaced by its commutator with $e_r$. Due to the antisymmetry of
the structure constants $c_{ijk}$ (Lemma \ref{skewsymm} on page
\pageref{skewsymm}), these
expressions cancel in pairs that correspond to the ends of each
chord.

To take a concrete example,
\begin{align*}
    &[e_r,\sum_{ij}e_ie_je_ie_j] \\
  = &\sum_{ij}[e_r,e_i]e_je_ie_j
  + \sum_{ij}e_i[e_r,e_j]e_ie_j
  + \sum_{ij}e_ie_j[e_r,e_i]e_j
  + \sum_{ij}e_ie_je_i[e_r,e_j] \\
  = &\sum_{ijk}c_{rik}e_ke_je_ie_j
  + \sum_{ijk}c_{rjk}e_ie_ke_ie_j
  + \sum_{ijk}c_{rik}e_ie_je_ke_j
  + \sum_{ijk}c_{rjk}e_ie_je_ie_k \\
  = &\sum_{ijk}c_{rik}e_ke_je_ie_j
  + \sum_{ijk}c_{rjk}e_ie_ke_ie_j
  + \sum_{ijk}c_{rki}e_ke_je_ie_j
  + \sum_{ijk}c_{rkj}e_ie_ke_ie_j.
\end{align*}
Here the first and the second sums cancel with the third and the fourth
sums, respectively.

(4) We still assume that the basis $\{e_i\}$ is
$\langle\cdot,\cdot\rangle$--orthonormal. Then one of the pairwise
differences of the chord diagrams that constitute the 4 term
relation in equation (\ref{3pairs}) (page \pageref{3pairs}) is sent
by $\f_\g$ to
$$
  \sum c_{ijk} \dots e_{i} \dots e_{j} \dots e_{k} \dots\ ,
$$
while the other goes to
$$
  \sum c_{ijk} \dots e_{j} \dots e_{k} \dots e_{i} \dots
  = \sum c_{kij} \dots e_{i} \dots e_{j} \dots e_{k} \dots
$$
By the cyclic symmetry of the structure constants $c_{ijk}$ in
an orthonormal basis (again see Lemma~\ref{skewsymm}), these two
expressions are equal.

(5) Using property (1), we can place
the base point in the product diagram $D_1\cdot D_2$
between $D_1$ and $D_2$. Then the identity
$\f_\g(D_1\cdot D_2)=\f_\g(D_1)\f_\g(D_2)$ becomes evident.
\end{proof}

\begin{xremark}
 If $D$ is a chord diagram with $n$ chords, then
$$\f_\g(D) = c^n + \{ {\rm terms\ of\ degree\ less\ than\ } 2n {\rm\ in\ }
U(\g) \},
$$
where $c$ is the quadratic Casimir element as on
page~\pageref{QuaCa}. Indeed, we can permute the endpoints of chords
on the circle without changing the highest term of $\f_\g(D)$ since
all the additional summands arising as commutators have degrees
smaller than $2n$. Therefore, the highest degree term of $\f_\g(D)$
does not depend on $D$. Finally, if $D$ is a diagram with $n$
isolated chords, that is, the $n$th power of the diagram with one
chord, then $\f_\g(D) = c^n$.
\end{xremark}

The centre $ZU(\g)$ of the universal enveloping algebra is precisely
the $\g$-invariant subspace $U(\g)^{\g}\subset U(\g)$, where the action of
$\g$ on $U(\g)$
consists in taking the commutator. According to
the Harish-Chandra theorem (see \cite{Hum}), for a semi-simple Lie
algebra $\g$, the centre $ZU(\g)$ is isomorphic to the algebra of
polynomials in certain variables $c_1=c, c_2, \ldots , c_r$, where
$r = rank(\g)$.

\subsection{}\label{lie_alg_ws_without_bases}
The construction of Lie algebra weight systems can be described
without referring to any particular basis.

A based chord diagram $D$ with $n$ chords gives a permutation
$\sigma_D$ of the set $\{1,2,\ldots,2n\}$ as follows. As we have
noted before, the endpoints of chords of a based chord diagram are
ordered, so we can order the chords of $D$ by their first endpoint.
Let us number the chords from 1 to $n$, and their endpoints from 1
to $2n$, in the increasing order. Then, for $1\leq i\leq n$ the
permutation $\sigma_D$ sends $2i-1$ to the (number of the) first
endpoint of the $i$th chord, and $2i$ to the second endpoint of the
same chord. In the terminology of Section~\ref{perm_act_on_C},
page~\pageref{perm_act_on_C}, the permutation $\sigma_D$  sends the
diagram with $n$ consecutive isolated chords into $D$. For instance:
$$\risS{-20}{ya-kart0}{\put(35,52){$D$}\put(14,22){$\sigma_D$}
           }{160}{35}{20} \hspace{1.5cm} \sigma_D=(132546)
$$

The bilinear form $\langle\cdot,\cdot\rangle$ on $\g$ is a tensor in
$\g^*\ot\g^*$. The algebra $\g$ is metrized, so we can identify
$\g^*$ with $\g$ and think of $\langle\cdot,\cdot\rangle$ as an
element of $\g\ot\g$. The permutation $\sigma_D$ acts on $\g^{\ot
2n}$ by interchanging the factors. The value of the universal Lie
algebra weight system $\f_\g(D)$ is then the image of the $n$th
tensor power $\langle\cdot,\cdot\rangle^{\ot n}$ under the map
$$\g^{\ot 2n}\stackrel{\sigma_D}{\longrightarrow}\g^{\ot 2n}\to U(\g),$$
where the second map is the natural projection of the tensor algebra
on $\g$ to its universal enveloping algebra.

\subsection{The universal $\sL_2$ weight system} \label{ws_sl_2_on_A}
\index{Weight system! $\sL_2$}
Consider the Lie algebra $\sL_2$ of $2\times 2$ matrices with zero
trace. It is a three-dimensional Lie algebra spanned by the matrices
$$H=\left(\begin{array}{cc} 1&0\\ 0&-1\end{array}\right)\ ,\qquad
  E=\left(\begin{array}{cc} 0&1\\ 0&0\end{array}\right)\ ,\qquad
  F=\left(\begin{array}{cc} 0&0\\ 1&0\end{array}\right)
$$
with the commutators
$$[H,E]=2E,\qquad [H,F]=-2F,\qquad [E,F]=H\ .$$
\def\bilinf#1#2{\langle #1, #2 \rangle}
We shall use the symmetric bilinear form $\bilinf{x}{y} =
\mbox{Tr}(xy)$:
$$\bilinf{H}{H}=2,\ \bilinf{H}{E}=0,\ \bilinf{H}{F}=0,\
\bilinf{E}{E}=0,\ \bilinf{E}{F}=1,\ \bilinf{F}{F}=0.
$$
One can easily check that it is ad-invariant and non-degenerate.
The corresponding dual basis is
$$H^*=\frac{1}{2}H,\qquad E^*=F,\qquad F^*=E,$$
and, hence, the Casimir element is $c=\frac{1}{2}HH +EF + FE$.

The centre $ZU(\sL_2)$ is
isomorphic to the algebra of polynomials in a single variable $c$.
The value $\f_{\sL_2}(D)$ is thus a polynomial in $c$.
In this section, following \cite{CV}, we explain a combinatorial
procedure to compute this polynomial for a given chord diagram $D$.

The algebra $\sL_2$ is simple, hence, any invariant form is equal to
$\lambda\bilinf{\cdot}{\cdot}$ for some constant $\lambda$. The
corresponding Casimir element $c_\lambda$, as an element of the
universal enveloping algebra, is related to $c=c_1$ by the formula
$c_\lambda = \frac{c}{\lambda}$. Therefore, the weight system
$$\f_{\sL_2}(D) = c^n + a_{n-1}c^{n-1} + a_{n-2}c^{n-2} + \dots
                      + a_2c^2 + a_1c
$$
and the weight system corresponding to
$\lambda\bilinf{\cdot}{\cdot}$
$$\f_{\sL_2,\lambda}(D) = c_\lambda^n + a_{n-1,\lambda}c_\lambda^{n-1}
        + a_{n-2,\lambda}c_\lambda^{n-2} + \dots
        + a_{2,\lambda}c_\lambda^2 + a_{1,\lambda}c_\lambda
$$
are related by the formula
$\f_{\sL_2,\lambda}(D) = \frac{1}{\lambda^n}\cdot
    \f_{\sL_2}(D)\vert_{_{c=\lambda\cdot c_\lambda}}$, or
$$a_{n-1} = \lambda a_{n-1,\lambda},\quad
  a_{n-2} = \lambda^2 a_{n-2,\lambda},\ \dots\
  a_{2} = \lambda^{n-2} a_{2,\lambda},\quad
  a_{1} = \lambda^{n-1} a_{1,\lambda}.
$$

\newsavebox{\pp}
\sbox{\pp}{
   \begin{picture}(1,4)
   \put(-2,0){\line(0,1){4}}
   \put(1,0){\line(0,1){4}}
   \end{picture}
          }
\def\pr{{\usebox{\pp}}}
\def\kr{\times}

\def\chvadth#1#2{\ \risS{-15}{#1}{
             \put(-7,28){\mbox{$\scriptstyle p_i$}}
             \put(-7,8){\mbox{$\scriptstyle p_j$}}
             \put(36,28){\mbox{$\scriptstyle p^*_i$}}
             \put(36,8){\mbox{$\scriptstyle p^*_j$}}
             #2}{35}{30}{20}\quad }

\begin{xtheorem} Let $\f_{\sL_2}$ be the weight system associated
with $\sL_2$, with the invariant form $\bilinf\cdot\cdot$. Take a chord
diagram $D$ and choose a chord $a$ of $D$. Then
$$\f_{\sL_2}(D) = (c-2k) \f_{\sL_2}(D_a) +
   2 \sum_{1 \leq i<j \leq k} \left(\f_{\sL_2}(D^\pr_{i,j}) -
                                    \f_{\sL_2}(D^\kr_{i,j}) \right),$$
where:\\
$\bullet$ $k$ is the number of chords that intersect the  chord $a$;\\
$\bullet$ $D_a$ is the chord diagram obtained from $D$ by
deleting the chord $a$;\\
$\bullet$ $D^\pr_{i,j}$ and $D^\kr_{i,j}$ are the chord diagrams
obtained from $D_a$ in the following way. Draw the diagram $D$ so
that the chord $a$ is vertical. Consider an arbitrary pair of chords
$a_i$ and $a_j$ different from $a$ and such that each of them
intersects $a$. 
Denote by $p_i$ and $p_j$ the endpoints of $a_i$ and $a_j$ that lie
to the left of $a$ and by $p^*_i, p^*_j$ the endpoints of $a_i$ and
$a_j$ that lie to the right. There are three ways to connect the
four points $p_i, p^*_i, p_j, p^*_j$ by two chords. $D_a$ is the
diagram where these two chords are $(p_i,p^*_i), (p_j,p^*_j)$, the
diagram $D^\pr_{i,j}$ has the chords $(p_i,p_j), (p^*_i, p^*_j)$ and
$D^\kr_{i,j}$ has the chords $(p_i,p^*_j),(p^*_i,p_j)$. All other
chords are the same in all the diagrams:
$$D=\chvadth{chvard01}{
             \put(8,28){\mbox{$\scriptstyle a_i$}}
             \put(8,7){\mbox{$\scriptstyle a_j$}}
             \put(19,17){\mbox{$\scriptstyle a$}} };\quad\!
  D_a=\chvadth{chvard02}{};\quad\!
  D^\pr_{i,j}=\chvadth{chvard03}{};\quad\!
  D^\kr_{i,j}=\chvadth{chvard04}{}.
$$
\end{xtheorem}

The theorem allows one to compute $\f_{\sL_2}(D)$ recursively, as
each of the three diagrams $D_a, D^\pr_{i,j}$ and $D^\kr_{i,j}$ has
one chord less than $D$.

\bigskip
\noindent{\bf Examples.}

\def\chvad#1{\risS{-10}{#1}{}{25}{10}{20}}
\def\chvadws#1{\f_{\sL_2}\bigl(\chvad{#1}\bigr)}

\begin{enumerate}
\item $\chvadws{chvard05} = (c-2)c.$ In this case,  $k=1$ and the sum
in the right hand side is zero, since there are no pairs $(i,j)$.

\item
$$\begin{array}{rl}
\ \chvadws{chvard06} &=\ (c-4) \chvadws{chvard07}
   + 2 \chvadws{chvard08} - 2 \chvadws{chvard05} \\
&=\ (c-4)c^2 + 2c^2 - 2(c-2)c\ \ =\ \ (c-2)^2c.
\end{array}$$

\item
$$\begin{array}{rl}
\ \chvadws{chvard09} &=\ (c-4) \chvadws{chvard05}
   + 2 \chvadws{chvard08} - 2 \chvadws{chvard07} \\
&=\ (c-4)(c-2)c + 2c^2 - 2c^2 \ \ =\ \ (c-4)(c-2)c.
\end{array}$$
\end{enumerate}

\medskip
\begin{xremark}
Choosing the invariant form $\lambda\bilinf{\cdot}{\cdot}$, we
obtain a modified relation
$$\f_{\sL_2,\lambda}(D) = \Bigl(c_\lambda - \frac{2k}{\lambda}\Bigr)
                    \f_{\sL_2,\lambda}(D_a) +
\frac{2}{\lambda} \sum_{1 \leq i<j \leq k}
              \left(\f_{\sL_2,\lambda}(D^\pr_{i,j}) -
                    \f_{\sL_2,\lambda}(D^\kr_{i,j})\right).
$$
If $k=1$, the second summand vanishes. In particular, for the
Killing form ($\lambda=4$) and $k=1$ we have
$$\f_\g(D) = (c-1/2) \f_\g(D_a).$$
It is interesting that the last formula is valid for any simple Lie
algebra $\g$ with the Killing form and any chord $a$ which intersects precisely one other chord. 
See Exercise~\ref{ex6_sl2} for a generalization of this fact in the case $\g=\sL_2$.
\end{xremark}

\noindent{\bf Exercise.} Deduce the theorem of this section from the
following lemma by induction (in case of difficulty see the proof in
\cite{CV}).

\begin{xlemma} [6-term relations for the universal $\sL_2$ weight
system] Let $\f_{\sL_2}$ be the weight system associated with
$\sL_2$ and the invariant form $\bilinf{\cdot}{\cdot}$. Then
$$\f_{\sL_2}\Bigl( \chvad{chvard10}\ -\ \chvad{chvard11}\ -\
     \chvad{chvard12}\ +\  \chvad{chvard13} \Bigr)\ =\
  2\f_{\sL_2}\Bigl( \chvad{chvard14}\ -\ \chvad{chvard15} \Bigr)\ ;
$$
$$\f_{\sL_2}\Bigl( \chvad{chvard16}\ -\ \chvad{chvard17}\ -\
     \chvad{chvard18}\ +\  \chvad{chvard19} \Bigr)\ =\
  2\f_{\sL_2}\Bigl( \chvad{chvard14}\ -\ \chvad{chvard20} \Bigr)\ ;
$$
$$\f_{\sL_2}\Bigl( \chvad{chvard21}\ -\ \chvad{chvard22}\ -\
     \chvad{chvard23}\ +\  \chvad{chvard24} \Bigr)\ =\
  2\f_{\sL_2}\Bigl( \chvad{chvard25}\ -\ \chvad{chvard26} \Bigr)\ ;
$$
$$\f_{\sL_2}\Bigl( \chvad{chvard27}\ -\ \chvad{chvard28}\ -\
     \chvad{chvard29}\ +\  \chvad{chvard30} \Bigr)\ =\
  2\f_{\sL_2}\Bigl( \chvad{chvard31}\ -\ \chvad{chvard32} \Bigr)\ .
$$
\end{xlemma}

These relations also provide a recursive way to compute
$\f_{\sL_2}(D)$ as the two chord diagrams on the right-hand side
have one chord less than the diagrams on the left-hand side, and the
last three diagrams on the left-hand side are simpler than the first
one since they have less intersections between their chords. See
Section~\ref{ws_sl_2_on_C} for a proof of this lemma.

\subsection{Weight systems associated with representations}
\label{ws_repres}

The construction of Bar-Natan, in comparison with that of
Kontsevich, uses one additional ingredient: a representation of a
Lie algebra (see \ref{metrized}).

A linear representation $T:\g\to\End(V)$ extends to a homomorphism
of associative algebras $U(T):U(\g)\to\End(V)$. The composition of
following three maps (with the last map being the trace)
$$
 \A \stackrel{\f_\g}{\to} U(\g) \stackrel{U(T)}{\to} \End(V)
    \stackrel{\Tr}{\to} \C
$$
by definition gives the {\em weight system associated with the
representation}\label{def:pr-phiA}
$$
  \f_\g^T = \Tr \circ U(T) \circ \f_\g
$$
(by abuse of notation, we shall sometimes write $\f_\g^V$ instead of
$\f_\g^T$).

The map $\f_\g^T$ is not in general multiplicative (the reader may check
this for the diagram $\t$ and its square in the standard representation of
the algebra $\gl_N$, see Section \ref{ws_glN_on_A}).
However, if the representation $T$ is irreducible, then, according to the
Schur Lemma \cite{Hum}, every element of the centre $ZU(\g)$ is
represented (via $U(T)$) by a scalar operator $\mu\cdot\mbox{id}_V$.
Therefore, its trace equals $\f^T_\g(D) = \mu\dim V$. The number
$\mu=\frac{\f^T_\g(D)}{\dim V}$, as a function of the chord diagram
$D$, is a weight system which is clearly multiplicative.

\subsection{Algebra $\sL_2$ with the standard
representation} Consider the standard 2-dimensional representation
$St$ of $\sL_2$. Then the Casimir element is represented by the
matrix
$$c=\frac{1}{2}HH + EF + FE =
    \left(\begin{array}{cc} 3/2&0\\ 0&3/2 \end{array}\right) =
    \frac{3}{2}\cdot\id_2 .
$$
In degree 3 we have the following weight systems
$$\begin{array}{c||c|c|c|c|c}
D& \risS{-12}{symcd1}{}{30}{20}{15} & \risS{-12}{symcd2}{}{30}{20}{15}
 & \risS{-12}{symcd3}{}{30}{20}{15} & \risS{-12}{symcd4}{}{30}{20}{15}
 & \risS{-12}{symcd5}{}{30}{20}{15} \\ \hline
\f_{\sL_2}(D) & c^3 & c^3 & c^2(c-2) & c(c-2)^2 & c(c-2)(c-4)
         \rb{-6pt}{\makebox(0,17){}}\\ \hline
\f^{St}_{\sL_2}(D) & 27/4 & 27/4 & -9/4 & 3/4 & 15/4
         \rb{-7pt}{\makebox(0,20){}} \\ \hline
\f'^{St}_{\sL_2}(D) & 0 & 0 & 0 & 12 & 24 \makebox(0,13){}
\end{array}
$$
Here the last row represents the unframed weight system obtained
from $\f^{St}_{\sL_2}$ by the deframing procedure from
Section~\ref{defram_ws}. A comparison of this computation with the
one from Section~\ref{symb_j_n} shows that on these elements $\symb(j_3) =
-\frac{1}{2} \f'^{St}_{\sL_2}$. \index{Weight system!of the Jones
coefficients} See Exercises~\ref{ex_defr_sl_N} and
\ref{ex_ws_j_N_and_sl_2} at the end of the chapter for more
information about these weight systems.

\subsection{Algebra $\gl_N$ with the standard representation}
\label{ws_glN_on_A}
\index{Weight system! $\gl_N$}

Consider the Lie algebra $\g=\gl_N$ of all $N\times N$ matrices
and its standard representation $St$. Fix the trace of the product
of matrices as the preferred ad-invariant form: $\langle x,y
\rangle = \mbox{Tr}(xy)$.

The algebra $\gl_N$ is linearly spanned by matrices $e_{ij}$ with 1
on the intersection of $i$th row with $j$th column and zero
elsewhere. We have $\langle e_{ij},e_{kl} \rangle = \delta_i^l
\delta_j^k$, where $\delta$ is the Kronecker delta. Therefore, the
duality between $\gl_N$ and $(\gl_N)^*$ defined by
$\langle\cdot,\cdot\rangle$ is given  by the formula $e_{ij}^* =
e_{ji}$.

\begin{xxca}
Prove that the form $\langle \cdot\, ,\,\cdot
\rangle$ is equal $2(N-1)$ times the Killing form.
({\sl Hint}: It is enough to compute the trace of just one operator
$(\ad e_{11})^2$.)
\end{xxca}

One can verify that $[e_{ij},e_{kl}] \not= 0$ only in the following cases:
\begin{itemize}
  \item $[e_{ij},e_{jk}] = e_{ik}$, if $i\ne k$,
  \item $[e_{ij},e_{ki}] = -e_{kj}$, if $j\ne k$,
  \item $[e_{ij},e_{ji}] = e_{ii}-e_{jj}$, if $i\ne j$,
\end{itemize}
This gives the following formula for the Lie bracket as a tensor in
$\gl_N^*\ot \gl_N^*\ot \gl_N$:
$$[\cdot\, ,\cdot] = \sum_{i,j,k=1}^N (e_{ij}^*\ot e_{jk}^*\ot e_{ik} -
                                    e_{ij}^*\ot e_{ki}^*\ot e_{kj}).
$$
When transferred to $\gl_N\ot \gl_N\ot \gl_N$ via the above mentioned
duality, this tensor takes the form
$$J = \sum_{i,j,k=1}^N (e_{ji}\ot e_{kj}\ot e_{ik} -
                                    e_{ji}\ot e_{ik}\ot e_{kj}).
$$
This formula will be used later in Section \ref{LAWS_C}.

D.\,Bar-Natan found the following elegant way of computing the
weight system $\f^{St}_{\gl_N}$.

\begin{xtheorem}[\cite{BN0}] Denote by $s(D)$ the number of
connected components of the curve obtained by doubling all chords of
a chord diagram $D$.
$$\risS{-12}{symbjn1}{}{30}{20}{15}\quad
  \risS{-2}{totonew}{}{25}{20}{15}\quad
  \risS{-12}{symbjn2}{}{30}{20}{15}\ .
$$
Then $\f^{St}_{\gl_N}(D) = N^{s(D)}\ .$
\end{xtheorem}

\begin{xremark}
By definition, the number $s(D)$ equals $c-1$, where $c$ is the number of
boundary components of the surface described in Section \ref{2term}.
\end{xremark}

\begin{xexample} For $D = \risS{-12}{symcd4}{}{30}{20}{15}$ we obtain
the picture $\risS{-12}{symst1}{}{30}{20}{15}$. Here $s(D)=2$, hence
$\f^{St}_{\gl_N}(D) = N^2$.
\end{xexample}

\begin{proof}
We take the matrices $e_{ij}$ as the chosen basis of $\gl_N$. The
values of the index variables associated with the chords are pairs
$(ij)$; each chord has one end labeled by a matrix $e_{ij}$ and the
other end by $e_{ji}=e_{ij}^*$.

Now, consider the curve $\gamma$ obtained by doubling the chords.
Given a chord whose ends are labeled by $e_{ij}$ and $e_{ji}$, we
can label the two copies of this chord in $\gamma$, as well as the
four pieces of the Wilson loop adjacent to its endpoints, by the
indices $i$ and $j$ as follows:
$$\risS{-12}{symbjn1}{\put(-12,12){\mbox{$\scriptstyle e_{ij}$}}
                      \put(32,12){\mbox{$\scriptstyle e_{ji}$}}
                     }{30}{20}{15}\qquad
  \risS{-2}{totonew}{}{25}{20}{15}\qquad
  \risS{-12}{symbjn2}{\put(-4,18){\mbox{$\scriptstyle i$}}
                      \put(-5,7){\mbox{$\scriptstyle j$}}
                      \put(31,18){\mbox{$\scriptstyle i$}}
                      \put(30,7){\mbox{$\scriptstyle j$}}
                      \put(14,19){\mbox{$\scriptstyle i$}}
                      \put(14,7){\mbox{$\scriptstyle j$}}
                      }{30}{20}{15}\ \ .
$$
To compute the value of the weight system $\f^{St}_{\gl_N}(D)$, we
must sum up the products $\dots e_{ij}e_{kl}\dots$. Since we are
dealing with the standard representation of $\gl_N$, the product
should be understood as genuine matrix multiplication, rather than the
formal product in the universal enveloping algebra. Since $e_{ij}\cdot
e_{kl} = \delta_{jk}\cdot e_{il}$, we get a non-zero summand only if
$j=k$. This means that the labels of the chords must follow the
pattern:\label{ws_glN_proof}
$$\risS{-12}{chvard38}{\put(12,-6){\mbox{$\scriptstyle e_{ij}$}}
                      \put(38,-6){\mbox{$\scriptstyle e_{jl}$}}
                     }{60}{20}{20}\qquad
  \risS{-2}{totonew}{}{25}{20}{15}\qquad
  \risS{-12}{chvard39}{\put(8,-4){\mbox{$\scriptstyle i$}}
                      \put(22,-6){\mbox{$\scriptstyle j$}}
                      \put(32,-6){\mbox{$\scriptstyle j$}}
                      \put(50,-4){\mbox{$\scriptstyle l$}}
                      \put(16,33){\mbox{$\scriptstyle i$}}
                      \put(22,33){\mbox{$\scriptstyle j$}}
                      \put(32,33){\mbox{$\scriptstyle j$}}
                      \put(40,33){\mbox{$\scriptstyle l$}}
                     }{60}{20}{20}\ \ .
$$
Therefore, all the labels on one and the same connected component of
the curve $\gamma$ are equal. If we take the whole product of
matrices along the circle, we get the operator $e_{ii}$ whose trace
is 1. Now, we must sum up the traces of all such operators over all
possible labelings. The number of labelings is equal to the number
of values the indices $i$, $j$, $l,\ldots$ take on the connected
components of the curve $\gamma$. Each component gives exactly $N$
possibilities, so the total number is $N^{s(D)}$.
\end{proof}

\begin{xproposition} The weight system $\f^{St}_{\gl_N}(D)$
depends only on the intersection graph of $D$. \label{IGglNst}
\end{xproposition}

\begin{proof}
The value $\f^{St}_{\gl_N}(D)$ is defined by the number $s(D)=c-1$
(where $c$ has the meaning given on page \pageref{2term}),
therefore it is a function of the genus of the diagram $D$.
In Section \ref{2term} we proved that the genus depends only on the
intersection graph.
\end{proof}

\subsection{Algebra $\sL_N$ with the standard representation}
\label{ws_sl_n_St}\index{Weight system! $\sL_N$}

Here we describe the weight system $\f^{St}_{\sL_N}(D)$ associated with the
Lie algebra $\sL_N$, its standard representation by $N\times N$ matrices
with zero trace and the invariant form $\langle x,y \rangle = \mbox{Tr}(xy)$,
Following Section~\ref{symb_j_n}, introduce {\it a state} $\sigma$
for a chord diagram $D$ as an arbitrary function on the set $[D]$ of
chords of $D$ with values in the set $\{1,-\frac{1}{N}\}$.
With each state $\sigma$ we associate an
immersed plane curve obtained from $D$ by resolutions of all its
chords according to $s$:
$$\risS{-12}{symbjn1}{\put(14,17){\mbox{$\scriptstyle c$}}}{30}{20}{15}\
  \risS{-2}{totonew}{}{25}{20}{15}\
  \risS{-12}{symbjn2}{}{30}{20}{15}\ ,\mbox{\ if\ }\sigma(c)=1;\qquad
 \risS{-12}{symbjn1}{\put(14,17){\mbox{$\scriptstyle c$}}}{30}{20}{15}\
  \risS{-2}{totonew}{}{25}{20}{15}\
 \risS{-12}{symbjn3}{}{30}{20}{15}\ ,\mbox{\ if\ }\sigma(c)=-\frac{1}{N}.
$$
Let $|\sigma|$ denote the number of components of the curve obtained in
this way.

\begin{xtheorem}
$\displaystyle
   \f^{St}_{\sL_N}(D) = \sum_\sigma\ \Bigl(\prod_c \sigma(c)\Bigr)\
N^{|\sigma|}\ ,$ where the product is taken over all $n$ chords of
$D$, and the sum is taken over all $2^n$ states for $D$.
\end{xtheorem}

\medskip
One can prove this theorem in the same way as we did for $\gl_N$ by
picking an appropriate basis for the vector space $\sL_N$ and then
working with the product of matrices (see
Exercise~\ref{ex_ws_sl_N_basis}). However, we prefer to prove it in
a different way, via several reformulations, using the algebra
structure of weight systems which is dual to the coalgebra structure
of chord diagrams (Section~\ref{bialgWS}).

\noindent{\bf Reformulation 1.} {\em For a subset $J\subseteq [D]$
(the empty set and the whole $[D]$ are allowed) of chords of $D$,
denote by $|J|$ the cardinality of $J$, and let
$n-|J|=|\overline{J}|$ stand for the number of chords in
$\ol{J}=[D]\setminus J$. Write $D_J$ for the chord diagram formed by
the chords from $J$, and denote by $s(D_J)$ the number of connected
components of the curve obtained by doubling all the chords of
$D_J$. Then
$$\f^{St}_{\sL_N}(D) = \sum_{J\subseteq [D]}
(-1)^{n-|J|} N^{s(D_J)-n+|J|}.
$$}

This assertion is obviously equivalent to the Theorem: for every
state $s$, the subset $J$ consists of all chords $c$ with value
$s(c)=1$.

Consider the weight system $e^{-\frac{\bo_1}{N}}$ from
Section~\ref{defram_ws}, which is equal to the constant
$\frac{1}{(-N)^n}$ on any chord diagram with $n$ chords.

\medskip

\noindent{\bf Reformulation 2.} {\it
$$\f^{St}_{\sL_N} = e^{-\frac{\bo_1}{N}} \cdot \f^{St}_{\gl_N}\ .$$}

Indeed, by the definition of the product of weight systems
(Section~\ref{bialgWS}),
$$\left(e^{-\frac{\bo_1}{N}} \cdot \f^{St}_{\gl_N}\right)(D) =
\left(e^{-\frac{\bo_1}{N}} \ot \f^{St}_{\gl_N}\right)(\d(D))\ ,
$$
where $\d(D)$ is the coproduct (Section~\ref{bialgCD}) of the chord
diagram $D$. It splits $D$ into two complementary parts
$D_{\overline{J}}$ and $D_J$: $\d(D)=\sum\limits_{J\subseteq [D]}
D_{\overline{J}}\ot D_J$. The weight system $\f^{St}_{\gl_N}(D_J)$
gives $N^{s(D_J)}$. The remaining part is given by
$e^{-\frac{\bo_1}{N}}(D_{\overline{J}})$.

\bigskip
\noindent{\bf Reformulation 3.} {\it
$$\f^{St}_{\gl_N} = e^{\frac{\bo_1}{N}} \cdot \f^{St}_{\sL_N}\ .$$}

The equivalence of this and the foregoing formulae follows from the
fact that the weight systems $e^{{\bo_1}/{N}}$ and $e^{-{\bo_1}/{N}}$
are inverse to each other as elements of the
completed algebra of weight systems.

\begin{proof}
We shall prove the theorem in Reformulation~3. The Lie algebra
$\gl_N$ is a direct sum of $\sL_N$ and the trivial one-dimensional
Lie algebra generated by the identity matrix $\mbox{id}_N$. Its dual
is $\mbox{id}^*_N=\frac{1}{N} \mbox{id}_N$. We can choose a basis
for the vector space $\gl_N$ consisting of the basis for $\sL_N$ and
the unit matrix $\mbox{id}_N$. To every chord we must assign either
a pair of dual basis elements of $\sL_N$, or the pair
$(\mbox{id}_N,\frac{1}{N} \mbox{id}_N)$, which is equivalent to
forgetting the chord and multiplying the obtained diagram by
$\frac{1}{N}$. This means precisely that we are applying the weight
system
$e^{{\bo_1}/{N}}$
to the chord subdiagram
$D_{\overline{J}}$ formed by the forgotten chords, and the weight
system $\f^{St}_{\sL_N}$ to the chord subdiagram $D_J$ formed by the
remaining chords.
\end{proof}

\subsection{Algebra $\so_N$ with the standard representation}
\label{ws_so_n_St}\index{Weight system! $\so_N$} In this case {\it a
state} $\sigma$ for $D$ is a function on the set $[D]$ of chords of
$D$ with values in the set $\{1/2,-1/2\}$. The rule for the
resolution of a chord according to its state is
$$\risS{-12}{symbjn1}{\put(14,17){\mbox{$\scriptstyle c$}}}{30}{20}{15}\
  \risS{-2}{totor}{}{25}{20}{15}\
  \risS{-12}{symbjn2}{}{30}{20}{15}\ ,\mbox{\ if\ }\sigma(c)=\frac{1}{2};\qquad
 \risS{-12}{symbjn1}{\put(14,17){\mbox{$\scriptstyle c$}}}{30}{20}{15}\
  \risS{-2}{totor}{}{25}{20}{15}\
 \risS{-12}{wssoNres}{}{30}{20}{15}\ ,\mbox{\ if\ }\sigma(c)=-\frac{1}{2}.
$$
As before, $|\sigma|$ denotes the number of components of the obtained curve.

\begin{xtheorem}[\cite{BN0,BN1}]
For the invariant form $\bilinf{x}{y} = \mbox{Tr}(xy)$,
$$\f^{St}_{\so_N}(D) = \sum_\sigma\ \Bigl(\prod_c\sigma(c)\Bigr)\
N^{|\sigma|}\ ,
$$
where the product is taken over all $n$ chords of
$D$, and the sum is taken over all $2^n$ states for $D$.
\end{xtheorem}

We leave the proof of this theorem to the reader as an exercise
(number \ref{ex_so_N_basis} at the end of the chapter, to be
precise). Alternatively, one can view a chord diagram as a closed
Jacobi diagram and use the theorem on page~\pageref{soN_stand_C}.

Here is the table of values of $\f^{St}_{\so_N}(D)$ for some basis elements
of $\A^{fr}$ of small degree:
$$\begin{array}{c||c|c|c|c}
D& \risS{-12}{wssoN1}{}{30}{20}{15} & \risS{-12}{chvard05a}{}{30}{20}{15}
 & \risS{-12}{symcd4}{}{30}{20}{15} & \risS{-12}{symcd5}{}{30}{20}{15}
       \\ \hline
\f^{St}_{\so_N}(D)
&\scriptstyle\frac{1}{2}(N^2-N)
&\scriptstyle\frac{1}{4}(N^2-N)
& \scriptstyle \frac{1}{8}(N^2-N)
& \scriptstyle \frac{1}{8}N(-N^2+4N-3)
         \rb{-8pt}{\makebox(0,20){}}\\ \hline\hline
\multicolumn{5}{l}{\hspace{-5pt}\begin{array}{c||c|c|c}
D&   \risS{-12}{cd4-02}{}{30}{20}{15} & \risS{-12}{cd4-03}{}{30}{20}{15}
        & \risS{-12}{cd4-01}{}{30}{20}{15} \\ \hline
\f^{St}_{\so_N}(D) & \scriptstyle \frac{1}{16}N(3N^2-8N+5)
              & \scriptstyle \frac{1}{16}(2N^3-5N^2+3N)
              & \scriptstyle \frac{1}{16}N(N^3-4N^2+6N-3)
         \rb{-8pt}{\makebox(0,20){}}\end{array}}\\ \hline
\end{array}
$$

Exercises \ref{ex_so_N_intgr}--\ref{ex_so_N_circle} contain
additional information about this weight system.

\subsection{Algebra $\sP_{2N}$ with the standard representation}
\index{Weight system! $\sP_{2N}$}

It turns out that $$\f^{St}_{\sP_{2N}}(D) =
(-1)^{n+1}\f_{\so_{-2N}}(D),$$ where the last notation means the
formal substitution of $-2N$ instead of the variable $N$ in the
polynomial $\f_{\so_N}(D)$, and $n$, as usual, is the degree of $D$.
This implies that the weight system $\f^{St}_{\sP_{2N}}$ does not
provide any new knot invariant. Some details about it can be found
in \cite{BN0,BN1}.

\medskip
It would be interesting to find a combinatorial description of the
weight systems for the exceptional simple Lie algebras $E_6$, $E_7$,
$E_8$, $F_4$, $G_2$.

\section{Lie algebra weight systems for the algebra $\F$}
\label{LAWS_C}

Since every closed diagram is a linear combination of chord
diagrams, the weight system $\f_\g$ can be treated as a function on
$\F$ with values in $U(\g)$.
It turns out that $\f_\g$  can be evaluated on any closed
diagram directly, often in a more convenient way.

The $\STU$ relation (Section~\ref{STUrel}), which defines the
algebra $\F$, gives us a hint how to do it. Namely, if we assign
elements $e_i$, $e_j$ to the endpoints of chords of the T- and U-
diagrams from the $\STU$ relations,
$$\risS{-12}{dT_iz_STU}{\put(12,-6){\mbox{$\scriptstyle e_i$}}
                      \put(40,-6){\mbox{$\scriptstyle e_j$}}
                      \put(18,34){\mbox{$\scriptstyle e^*_i$}}
                      \put(38,34){\mbox{$\scriptstyle e^*_j$}}
                      \put(-10,0){\mbox{T}}
                     }{60}{30}{20}\qquad - \qquad
\risS{-12}{dU_iz_STU}{\put(12,-6){\mbox{$\scriptstyle e_j$}}
                      \put(40,-6){\mbox{$\scriptstyle e_i$}}
                      \put(20,34){\mbox{$\scriptstyle e^*_i$}}
                      \put(38,34){\mbox{$\scriptstyle e^*_j$}}
                      \put(-10,0){\mbox{U}}
                     }{60}{20}{20}\qquad = \qquad
\risS{-12}{dS_iz_STU}{\put(20,-6){\mbox{$\scriptstyle [e_i,e_j]$}}
                      \put(20,34){\mbox{$\scriptstyle e^*_i$}}
                      \put(38,34){\mbox{$\scriptstyle e^*_j$}}
                      \put(-10,0){\mbox{S}}
                     }{60}{20}{20}\ ,
$$
then it is natural to assign the commutator $[e_i,e_j]$ to the
trivalent vertex on the Wilson loop of the S-diagram.

Strictly speaking, $[e_i,e_j]$ may not be a basis vector. A diagram
with an endpoint marked by a linear combination of the basis vectors
should be understood as a corresponding linear combination of
diagrams marked by basis vectors. For this reason it will be more
convenient to use the description of $\f_\g$ given in
\ref{lie_alg_ws_without_bases}, which does not depend of the choice
of a basis. The formal construction goes as follows.

Let $C\in \FD_n$ be a closed Jacobi diagram with a base point and
$V=\{v_1,\dots,v_m\}$ be the set of its external vertices ordered
according to the orientation of the Wilson loop. We shall construct
a tensor $T_\g(C)\in\g^{\ot m}$ whose $i$th tensor factor $\g$
corresponds to the element $v_i$ of the set $V$. The weight system
$\f_\g$ evaluated on $C$ is the image of $T_\g(C)$ in $U(\g)$ under
the natural projection.

In order to construct the tensor $T_\g(C)$, consider the internal
graph of $C$ and cut all the edges connecting the trivalent vertices
of $C$. This splits the internal graph of $C$ into a union of
elementary pieces of two types: chords and tripods, the latter
consisting of one trivalent vertex and three legs with a fixed
cyclic order. Here is an example: 
$$\risS{-8}{ya-kart1a}{}{25}{18}{12}\quad\risS{0}{totor}{}{25}{0}{0}
\quad\risS{-8}{ya-kart1b}{}{31}{0}{0}.$$%\ .

To each leg of a chord or of a tripod we associate a copy of $\g$,
marked by this leg. Just as in \ref{lie_alg_ws_without_bases}, to
each chord we can assign the tensor $\langle\cdot\,,\cdot\rangle$
considered as an element of $\g\ot\g$, where the copies of $\g$ in
the tensor product are labeled by the ends of the chord. Similarly,
to a tripod we associate the tensor $-J\in\g\ot\g\ot\g$ defined as
follows. The Lie bracket $[\cdot\,,\cdot]$ is an element of
$\g^*\ot\g^*\ot\g$. Identifying $\g^*$ and $\g$ by means of
$\langle\cdot\,,\cdot\rangle$ we see that it corresponds to a tensor
in $\g\ot\g\ot\g$ which we denote by $J$. The order on the three
copies of $\g$ should be consistent with the cyclic order of legs in
the tripod.

Now, take the tensor product $\widetilde{T}_\g(C)$ of all the
tensors assigned to the elementary pieces of the internal graph of
$C$, with an arbitrary order of the factors.  It is an element of
the vector space $\g^{\ot (m+2k)}$ which has one copy of $\g$ for
each external vertex $v_i$ of $C$ and two copies of $\g$ for each of
the $k$ edges where the internal graph of $C$ has been cut. The form
$\langle\cdot\,,\cdot\rangle$, considered now as a bilinear map of
$\g\ot\g$ to the ground field, induces a map
$$\g^{\ot (m+2k)}\to \g^{\ot m}$$
by contracting a tensor over all pairs of coinciding labels. Apply
this contraction to $\widetilde{T}_\g(C)$; the result is a tensor in
$\g^{\ot m}$ where the factors are indexed by the $v_i$, but
possibly in a wrong order. Finally, re-arranging the factors in
$\g^{\ot m}$ according to the cyclic order of vertices on the Wilson loop,
we obtain the tensor $T_\g(C)$ we were looking for.

\begin{xremark}
Note that we associate the tensor $-J$, not $J$, to each tripod. This
is not a matter of choice, but a reflection of our convention for the
default cyclic order at the 3-valent vertices and the signs in the STU
relation.
\end{xremark}

\begin{xremark}
The construction of the tensor $T_\g(C)$ consists of two steps:
taking the product of the tensors that correspond to the elementary
pieces of the internal graph of $C$, and contracting these tensors
on the coinciding labels. These two steps can be performed
``locally''. For instance, let
$$C = \risS{-8}{ya-kart2}{}{25}{18}{12}.$$
The internal graph of $C$ consists of 3 tripods. To obtain $T_\g(C)$
we first take the tensor product of two copies of $-J$ and contract
the resulting tensor on the coinciding labels, thus obtaining a
tensor in $\g^{\ot 4}$. Graphically, this could be illustrated by
glueing together the two tripods into a graph with four univalent
vertices. Next, this graph is glued to the remaining tripod; this
means taking product of the corresponding tensors and contracting it
on a pair of labels:
$$\risS{-25}{ya-kart3a}{\put(53,37){$\g$}\put(56,27){$\g$}\put(35,7){$\g$}
          \put(8,-1){$\g$}\put(62,-6){$\g$}\put(-6,24){$\g$}
          \put(23,60){$\g$}\put(62,60){$\g$}\put(27,8){$\g$}
          }{65}{40}{35}\qquad\risS{0}{totor}{}{25}{0}{0}\quad
\quad\risS{-25}{ya-kart3b}{\put(28,37){$\g$}\put(31,28){$\g$}
          \put(8,-6){$\g$}\put(37,-6){$\g$}\put(-6,20){$\g$}
          \put(-2,60){$\g$}\put(37,60){$\g$}
          }{40}{0}{0}\qquad\risS{0}{totor}{}{25}{0}{0}\quad
\quad\risS{-12}{ya-kart3c}{\put(8,-7){$\g$}\put(37,-6){$\g$}
          \put(6,40){$\g$}\put(37,39){$\g$}\put(-6,20){$\g$}
          }{40}{0}{0}$$
The result is, of course, the same as if we took first the tensor
product of all 3 copies of $-J$ and then performed all the
contractions.
\end{xremark}

The only choice involved in the construction of $T_\g(C)$ is the
order of the factors in the tensor product $\g^{\ot 3}$ that corresponds
to a tripod. The following exercise shows that this order does not
matter as long as it is consistent with the cyclic order of legs:
\begin{xca}
Use the properties of $[\cdot\,,\cdot]$ and
$\langle\cdot\,,\cdot\rangle$ to prove that the tensor $J$ is
skew-symmetric under the permutations of the three tensor factors
(for the solution see Lemma \ref{skewsymm} on page
\pageref{skewsymm}; note that we have already used this fact earlier
in the proof of Theorem \ref{phi_prop}).
\end{xca}
This shows that $T_\g(C)$ is well-defined. Moreover, it produces a
weight system: the definition of the commutator in the universal
enveloping algebra implies that the element $\f_\g(C)$, which is the
image of $T_\g(C)$ in $U(\g)$, satisfies the STU relation. If $C$ is
a chord diagram, this definition of $\f_\g(C)$ coincides with the
definition given in \ref{lie_alg_ws_without_bases}.

Since the STU relation implies both the AS and the IHX relations,
$\f_\g$ satisfies these relations too. Moreover, it is easy to see
that the AS and the IHX relations are already satisfied for the
function $C\mapsto T_\g(C)$:
\begin{itemize}
\item the AS relation follows from the fact that the tensor $J$ changes
sign under odd permutations of the three factors in $\g\ot\g\ot\g$.
\item
the IHX relation is a corollary of the Jacobi identity in $\g$.
\end{itemize}

\subsection{}
\label{univ_ws_bubble}
Let us show how the construction of $T_\g$ works on an example and
prove the following lemma that relates the tensor corresponding to a
``bubble'' with the quadratic Casimir tensor.

\begin{xlemma}  For the Killing form $\langle \cdot\,,\cdot \rangle^K$ as
the preferred invariant form, the tensor $T_\g$ does not change if
a bubble is inserted into an internal edge of a diagram:
$$T_\g(\risS{-12}{laws_bub1}{}{30}{15}{15})\ =\
 T_\g(\risS{-12}{laws_bub2}{}{30}{0}{0})\ .
$$
\end{xlemma}

\begin{proof}
The fragment of a closed diagram on the right hand side is
obtained from two tripods by contracting the corresponding two
copies of the tensor $-J$. This gives the following tensor written
in an orthonormal basis $\{e_i\}$:
$$\begin{array}{c@{\qquad}c@{\quad}l}
\risS{-11}{laws_bub3}{
      \put(0,19){\mbox{${\scriptstyle e_i}$}}
      \put(16,30){\mbox{\frame{\makebox(25,11){
                                $\!{\scriptstyle e_k\ e_{k'}}$}}}}
      \put(16,-12){\mbox{\frame{\makebox(25,11){
                                $\!{\scriptstyle e_j\ e_{j'}}$}}}}
      \put(56,19){\mbox{${\scriptstyle e_{l}}$}}
}{60}{20}{15} & \risS{0}{totor}{}{25}{30}{25} & \displaystyle
\sum_{i,l}\ \sum_{k,j,k',j'} c_{ijk}c_{lk'j'}
          \langle e_k,e_{k'}\rangle^K  \langle e_j,e_{j'}\rangle^K
                          e_i\ot e_{l} \\
&& \displaystyle =\ \sum_{i,l} \Bigl(\sum_{j,k} c_{ijk}c_{lkj}
\Bigr)
                    e_i\ot e_{l}\ ,
\end{array}
$$
where $c_{ijk}$ are the structure constants:
$J= \sum\limits_{i,j,k=1}^d c_{ijk} e_i\ot e_j\ot e_k$.

To compute the coefficient $\Bigl(\sum\limits_{j,k} c_{ijk}c_{lkj}
\Bigr)$ let us find the value of the Killing form
$$\langle e_i,e_l\rangle^K = \mbox{Tr}( \mbox{ad}_{e_i} \mbox{ad}_{e_l})\ .
$$
Since
$$\mbox{ad}_{e_i}(e_s) = \sum_k c_{isk} e_k\quad
\mbox{and} \quad \mbox{ad}_{e_l}(e_t) = \sum_k c_{ltk} e_k\ ,
$$
the $(j,r)$-entry of the matrix of the product
$\mbox{ad}_{e_i} \mbox{ad}_{e_l}$ will be $\sum\limits_k c_{ikj}c_{lrk}$. Therefore,
$$\langle e_i,e_l\rangle^K = \sum_{k,j} c_{ikj}c_{ljk} =
\sum_{j,k} c_{ijk}c_{lkj}\ .
$$
Orthonormality of the basis $\{e_i\}$ implies that
$$\sum_{j,k} c_{ijk}c_{lkj} = \delta_{i,l}.
$$
This means that the tensor on the left-hand side in the statement of the lemma
equals
$$\sum_i e_i\ot e_i\ ,
$$
which is the quadratic Casimir tensor from the right-hand side.
\end{proof}
\begin{xremark}
If in the above lemma we use the bilinear form $\mu\langle
\cdot,\cdot \rangle^K$ instead of the Killing form, the rule changes
as follows:
$$T_\g(\risS{-12}{laws_bub1}{}{30}{10}{20})\ =\ \frac{1}{\mu}
 T_\g(\risS{-12}{laws_bub2}{}{30}{0}{0}).$$
\end{xremark}

\subsection{The universal $\sL_2$ weight system for $\F$}
\label{ws_sl_2_on_C}
\index{Weight system! $\sL_2$}
\begin{xtheorem}[\cite{CV}] For the invariant form
$\bilinf{x}{y} =\mbox{Tr}(xy)$ the tensor $T_{\sL_2}$ satisfies
the following skein relation:
$$T_{\sL_2}(\risS{-12}{chvar_euler1}{}{30}{15}{15})\ =\
 2T_{\sL_2}(\risS{-12}{chvar_euler2}{}{30}{15}{15}) -
 2T_{\sL_2}(\risS{-12}{chvar_euler3}{}{30}{15}{15})\ .$$
\end{xtheorem}

If the chosen invariant form is $\lambda\bilinf{\cdot}{\cdot}$,
then the coefficient 2 in this equation is replaced by
$\frac{2}{\lambda}$.

\begin{proof}
For the algebra $\sL_2$ the Casimir tensor and the Lie bracket tensor are
$$ C = \frac{1}{2} H\ot H + E\ot F + F\ot E\ ;$$
$$-J = -H\ot F\ot E +F\ot H\ot E +H\ot E\ot F -E\ot H\ot F
       -F\ot E\ot H  +E\ot F\ot H.$$
Then the tensor corresponding to the elementary pieces on the
right-hand side is equal to (we enumerate the vertices according
to the tensor factors)
$$\begin{array}{ccl}
  T_{\sL_2}\Bigl(\ \ \risS{-6}{chvar_euler4}{
      \put(-3,18){\mbox{$\scriptstyle 1$}}
      \put(-3,-6){\mbox{$\scriptstyle 2$}}
      \put(25,18){\mbox{$\scriptstyle 4$}}
      \put(25,-6){\mbox{$\scriptstyle 3$}}
           }{25}{15}{0}\ \ \Bigr) &=& \scriptstyle
  -H\ot F\ot H\ot E +H\ot F\ot E\ot H +F\ot H\ot H\ot E -F\ot H\ot E\ot H  \\
&&\scriptstyle -H\ot E\ot H\ot F
           +H\ot E\ot F\ot H +E\ot H\ot H\ot F -E\ot H\ot F\ot H \\
&&\scriptstyle  +2F\ot E\ot F\ot E -2F\ot E\ot E\ot F
           -2E\ot F\ot F\ot E +2E\ot F\ot E\ot F \vspace{10pt} \\
&&\hspace*{-70pt} =\ \scriptstyle
      2\bigl(\frac{1}{4}H\ot H\ot H\ot H +\frac{1}{2}H\ot E\ot F\ot H
   + \frac{1}{2}H\ot F\ot E\ot H +\frac{1}{2}E\ot H\ot H\ot F \\
&&\hspace*{-60pt}\scriptstyle
   +E\ot E\ot F\ot F + E\ot F\ot E\ot F +\frac{1}{2}F\ot H\ot H\ot E
   +F\ot E\ot F\ot E + F\ot F\ot E\ot E \bigr) \\
&&\hspace*{-68pt} \ \scriptstyle
     - 2\bigl(\frac{1}{4}H\ot H\ot H\ot H +\frac{1}{2}H\ot E\ot H\ot F
   + \frac{1}{2}H\ot F\ot H\ot E +\frac{1}{2}E\ot H\ot F\ot H \\
&&\hspace*{-60pt}\scriptstyle
   +E\ot E\ot F\ot F + E\ot F\ot F\ot E +\frac{1}{2}F\ot H\ot E\ot H
   +F\ot E\ot E\ot F + F\ot F\ot E\ot E \bigr) \vspace{10pt} \\
&&\hspace*{-70pt} =\
  2T_{\sL_2}\Bigl(\ \ \risS{-6}{chvar_euler5}{
      \put(-3,18){\mbox{$\scriptstyle 1$}}
      \put(-3,-6){\mbox{$\scriptstyle 2$}}
      \put(25,18){\mbox{$\scriptstyle 4$}}
      \put(25,-6){\mbox{$\scriptstyle 3$}}
           }{25}{15}{0}\ \ \Bigr) -
  2T_{\sL_2}\Bigl(\ \ \risS{-6}{chvar_euler6}{
      \put(-3,18){\mbox{$\scriptstyle 1$}}
      \put(-3,-6){\mbox{$\scriptstyle 2$}}
      \put(25,18){\mbox{$\scriptstyle 4$}}
      \put(25,-6){\mbox{$\scriptstyle 3$}}
           }{25}{15}{0}\ \ \Bigr)\ .
\end{array}$$
\end{proof}

\begin{xremark} While transforming a closed diagram according to this
theorem a closed circle different from the Wilson loop may occur
(see the example below). In this situation the circle should be
replaced by the numeric factor $3=\dim\sL_2$, which is the trace
of the identity operator in the adjoint representation of $\sL_2$.
\end{xremark}
\begin{xremark}
In the context of weight systems this relation was first noted in
\cite{CV}; afterwards, it was rediscovered several times. In a more
general context of graphical notation for tensors  it appeared
already in R.~Penrose's paper \cite{Pen}. In a certain sense, this
relation goes back to Euler and Lagrange because it is an exact
counterpart of the classical ``$bac-cab$'' rule, $$\mbox{\bf
a}\times (\mbox{\bf b}\times \mbox{\bf c}) =
 \mbox{\bf b}(\mbox{\bf a}\cdot \mbox{\bf c}) -
 \mbox{\bf c}(\mbox{\bf a}\cdot \mbox{\bf b}),$$
for the ordinary cross product of vectors in 3-space.
\end{xremark}
\begin{xexample}
$$\begin{array}{ccl}
  \f_{\sL_2}\Bigl( \risS{-12}{chvar_ex1}{}{30}{20}{15} \Bigr) &=&
 2\f_{\sL_2}\Bigl( \risS{-12}{chvar_ex2}{}{30}{20}{15} \Bigr) -
 2\f_{\sL_2}\Bigl( \risS{-12}{chvar_ex3}{}{30}{20}{15} \Bigr)\ =\
 4\f_{\sL_2}\Bigl( \risS{-12}{chvar_ex4}{}{30}{20}{15} \Bigr)\\
&&
   -4\f_{\sL_2}\Bigl( \risS{-12}{chvar_ex5}{}{30}{20}{15} \Bigr) -
    4\f_{\sL_2}\Bigl( \risS{-12}{chvar_ex6}{}{30}{20}{15} \Bigr) +
    4\f_{\sL_2}\Bigl( \risS{-12}{chvar_ex7}{}{30}{20}{15} \Bigr)\\
&=& 12c^2 -4c^2 -4c^2 +4c^2\ =\ 8c^2\ .
\end{array}$$
\end{xexample}
The next corollary implies the 6-term relation from
Section~\ref{ws_sl_2_on_A}.

\begin{xcorollary}
$$\f_{\sL_2}\Bigl( \chvad{chvar_lem1} \Bigr) =
  2\f_{\sL_2}\Bigl( \chvad{chvard14}\ -\ \chvad{chvard15} \Bigr);\quad
\f_{\sL_2}\Bigl( \chvad{chvar_lem2} \Bigr)=
  2\f_{\sL_2}\Bigl( \chvad{chvard14}\ -\ \chvad{chvard20} \Bigr);
$$
$$\f_{\sL_2}\Bigl( \chvad{chvar_lem3} \Bigr)=
  2\f_{\sL_2}\Bigl( \chvad{chvard25}\ -\ \chvad{chvard26} \Bigr);\quad
\f_{\sL_2}\Bigl( \chvad{chvar_lem4} \Bigr)=
  2\f_{\sL_2}\Bigl( \chvad{chvard31}\ -\ \chvad{chvard32} \Bigr).
$$
\end{xcorollary}

\subsection{The universal  $\gl_N$ weight system for $\F$}
\label{univ_ws_glN}
Let us apply the general procedure of the beginning of this section 
to the Lie algebra $\gl_N$ equipped with the bilinear form
$\langle e_{ij},e_{kl}\rangle=\d_{il}\d_{jk}$ so that $e_{ij}^*=e_{ji}$.
The corresponding universal weight system $\f_{\gl_N}$ can be calculated with the help of a graphical
calculus similar to that invented by R.~Penrose in \cite{Pen}. (A modification of this calculus
is used in \cite{BN1} to treat  the standard representation of $\gl_N$, see Section~\ref{ws_glN_on_C} below).

According to the general procedure, in order to construct $T_{\gl_N}$ we first erase the Wilson loop of the diagram, then place a copy of the tensor
$$
 -J = \sum_{i,j,k=1}^N (e_{ij}\ot e_{jk}\ot e_{ki}
                       - e_{ij}\ot e_{ki}\ot e_{jk})
$$
into each trivalent vertex and, finally, make contractions along all edges.  Any interval component (that is, chord) of the internal graph of the diagram is replaced simply by a copy of the bilinear form understood as the element $\sum e_{ij}\ot e_{ji}$.  The cyclic order of the endpoints is remembered. The universal weight system $\f_{\gl_N}$ is the image of $T_{\gl_N}$ in the universal enveloping algebra $U(\gl_N)$; in order to obtain it we simply omit the symbol of the tensor product in the above expressions:
$$
 -J = \sum_{i,j,k=1}^N (e_{ij}e_{jk}e_{ki} - e_{ij}e_{ki}e_{jk})\ .
$$

Now, the formula for $-J$ can be represented graphically as\label{T-diagram}
$$
-J = \rb{-6.5mm}{\ig[height=15mm]{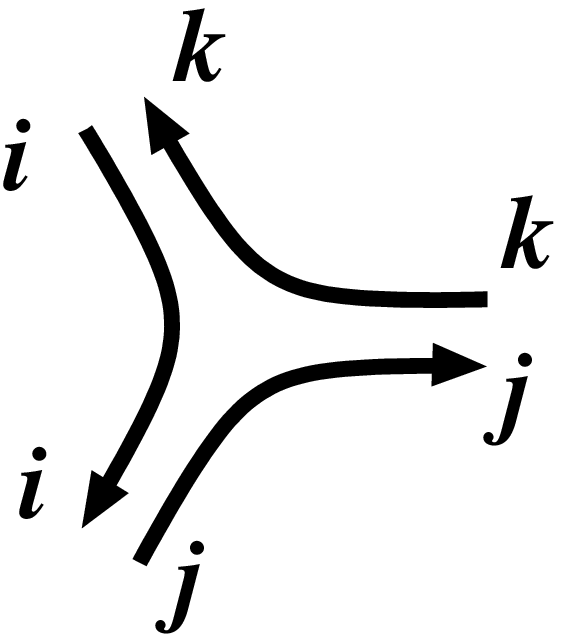}}
   - \rb{-6.5mm}{\ig[height=15mm]{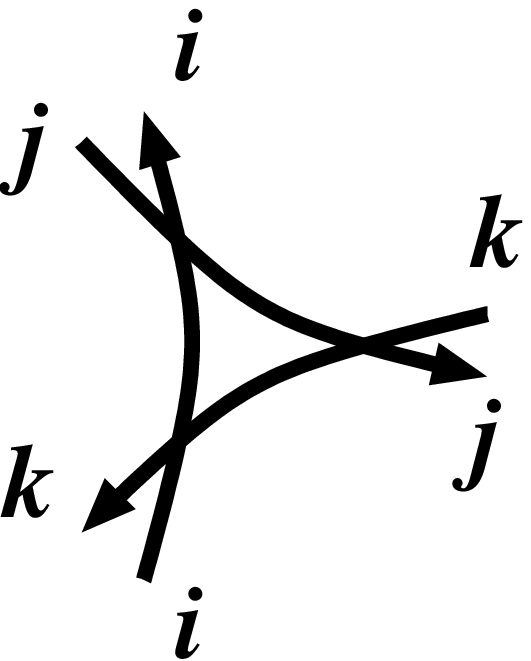}}
$$
One should imagine a basis element $e_{ij}$ attached to any pair of
adjacent endpoints marked $i$ and $j$, with $i$ on the incoming line and $j$
on the outgoing line. More generally, one may encode tensors by pictures as follows: 
specify $k$ pairs of points, each point connected to some other point with an arrow. Each of the $k$ arrows carries an index and each of the $k$ pairs carries the generator $e_{ij}$, where $i$ is the index of the incoming arrow and $j$ is the index of the outgoing arrow. The tensor that corresponds to such a picture is obtained by fixing some order on the set of pairs, taking the product of the $n$ elements $e_{ij}$ that correspond to the pairs, in the corresponding order, and then taking the sum over all the possible values of all the indices.

Choose one of the two pictures as above for each trivalent vertex
(this may be thought of as ``resolving'' the trivalent vertex in a positive or negative way).
The contraction along the edges means that we must glue together the small
pictures. This is done in the following manner. For any edge connecting two trivalent 
vertices, the contraction along it always gives zero except for the case when we have
$\langle e_{ij},e_{ji}\rangle=1$. Graphically, this means that
we must connect the endpoints of the tripods and write one and the same
letter on each connected component of the resulting curve.
Note that the orientations on the small pieces of curves (that come from the
cyclic order of the edges at every vertex) always agree for any set of
resolutions, so that we get a set of oriented curves.
We shall, further, add small intervals at each univalent vertex (now
doubled) thus obtaining one connected oriented curve for every connected
component of the initial diagram. To convert this curve into an element
of the universal enveloping algebra, we write, at every univalent vertex,
the element $e_{ij}$ where the subscripts
$i$ and $j$ are written in the order induced by the orientation on the
curve. Then we take the product of all such elements in the order coming
from the cyclic order of univalent vertices on the Wilson loop. As we know,
the result is invariant under cyclic permutations
of the factors. Finally, we sum up these results over all resolutions of the
triple points.
\smallskip

\begin{xexample}
Let us compute the value of $\f$ on the diagram
$$
  C = \rb{-4.5mm}{\ig[height=11mm]{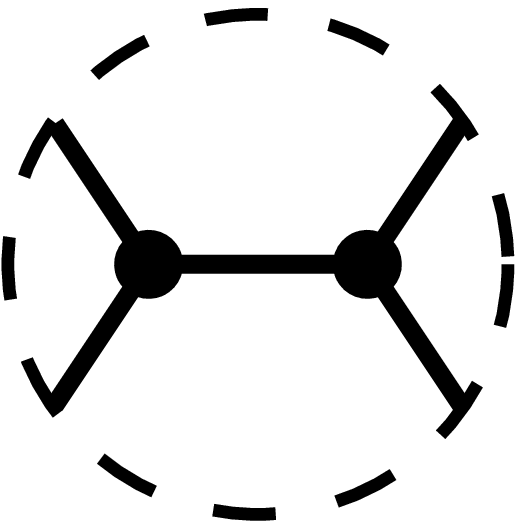}}
$$
We have:
\label{ex_ws_glN}
\begin{eqnarray*}
  C & \longmapsto &
               \rb{-3mm}{\ig[height=8mm]{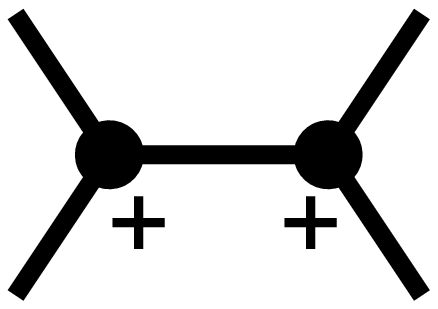}}
          \ -\ \rb{-3mm}{\ig[height=8mm]{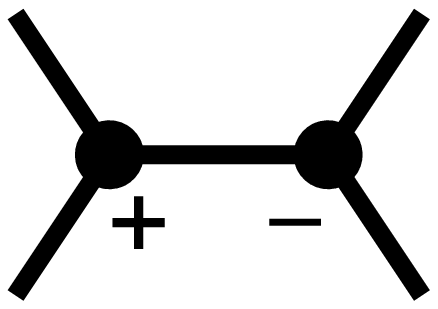}}
          \ -\ \rb{-3mm}{\ig[height=8mm]{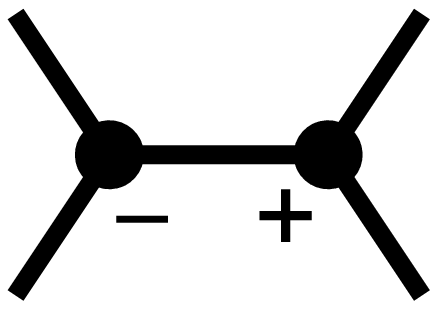}}
          \ +\ \rb{-3mm}{\ig[height=8mm]{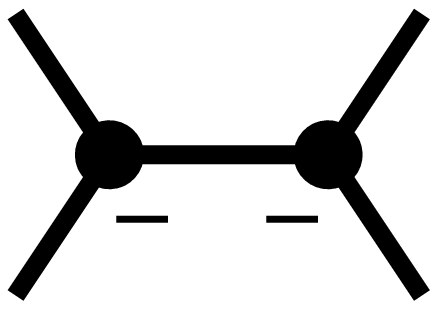}} \\[2mm]
 & \longmapsto &
               \rb{-8mm}{\ig[height=18mm]{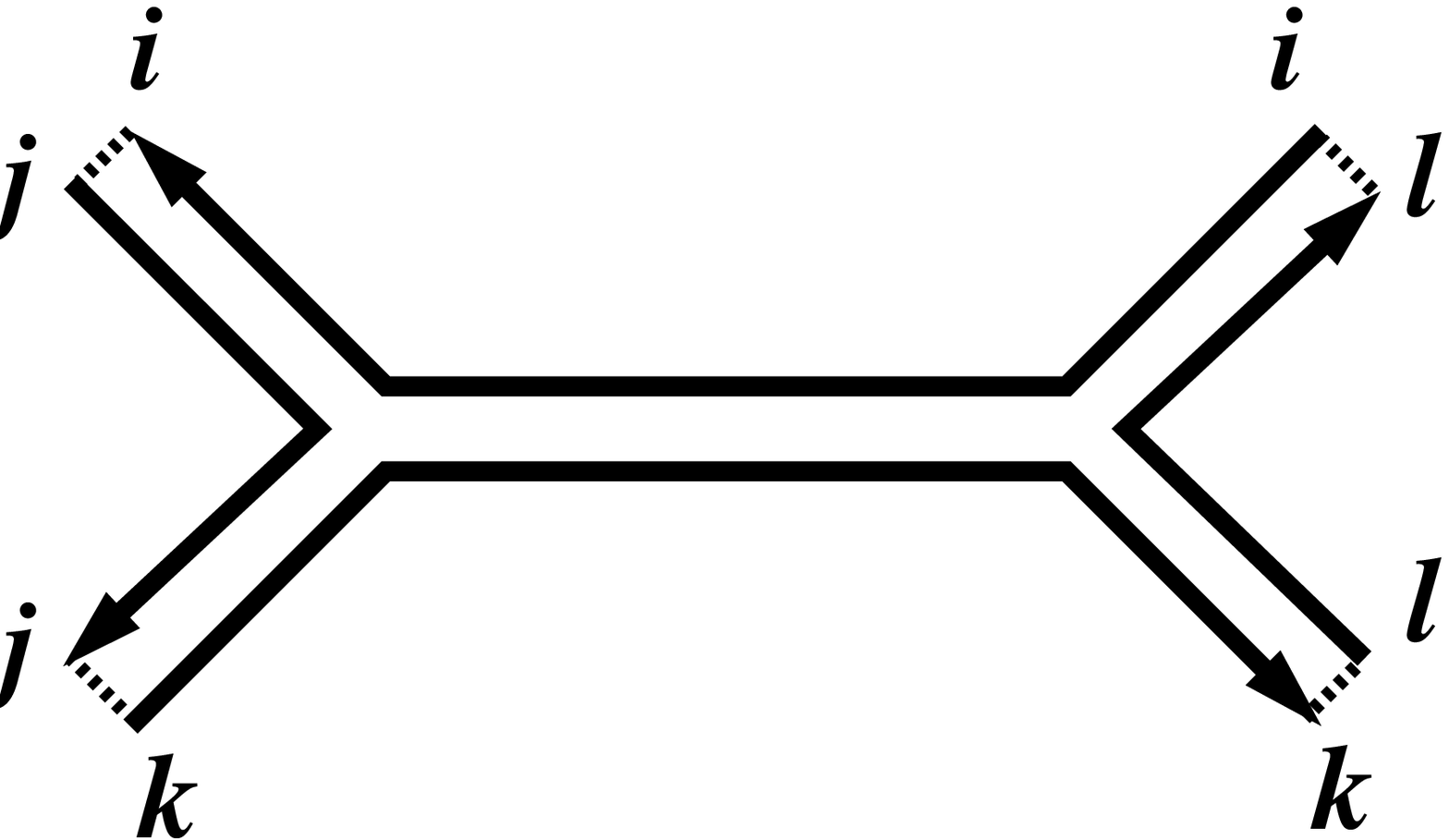}}
          \ -\ \rb{-8mm}{\ig[height=18mm]{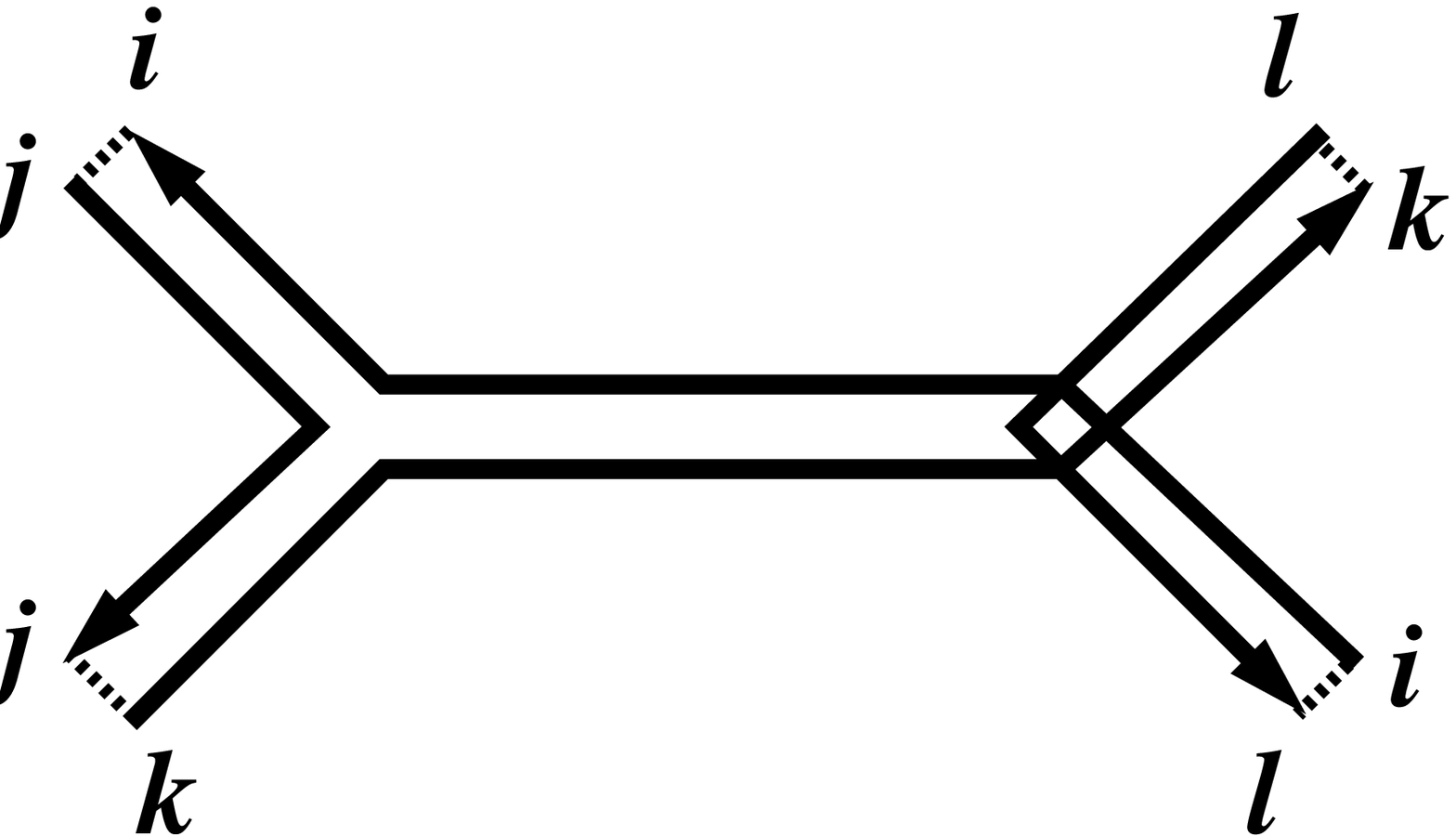}}\\[2mm]
          &&\ -\ \rb{-8.5mm}{\ig[height=18mm]{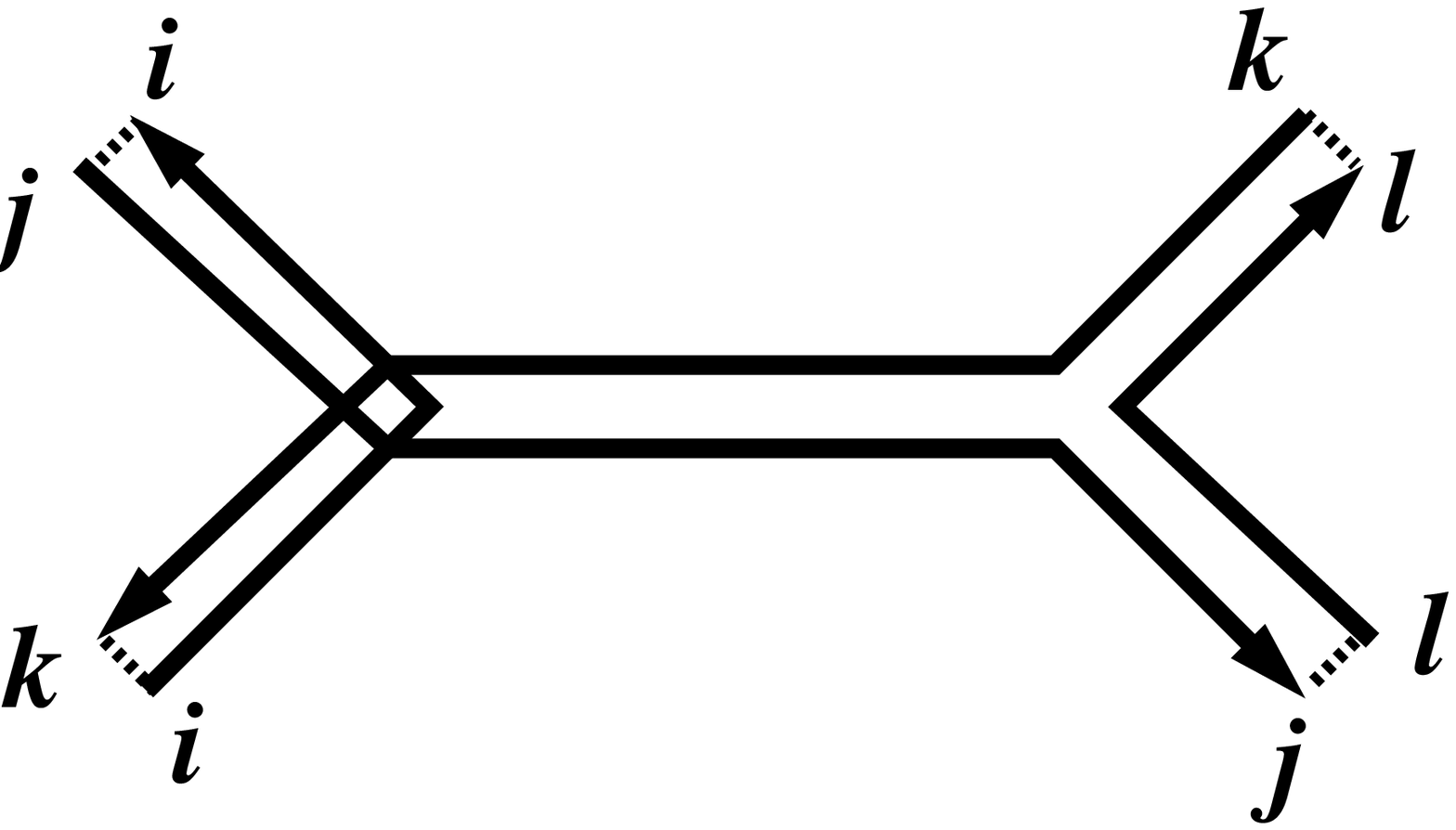}}
          \ +\ \rb{-8mm}{\ig[height=18mm]{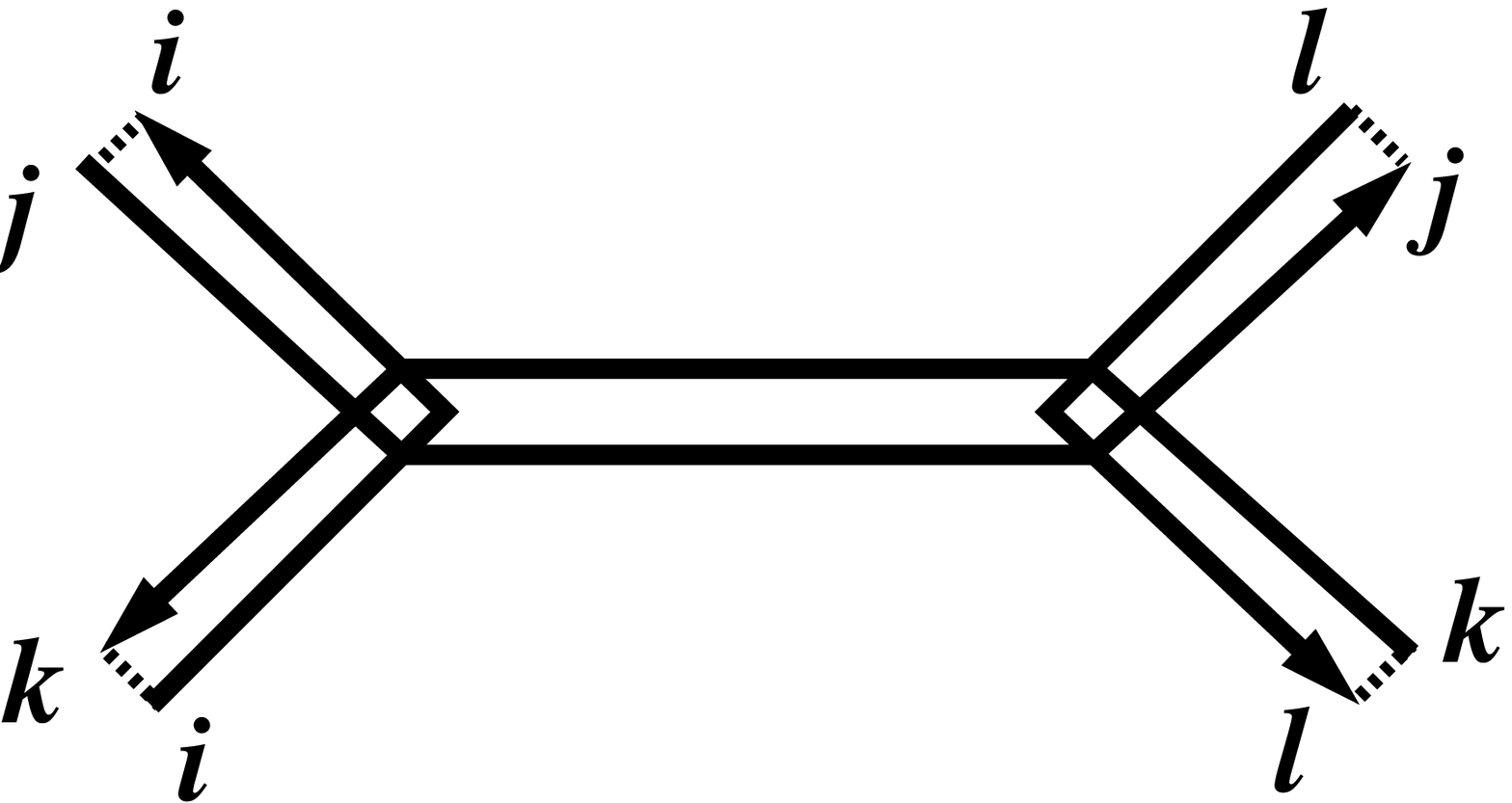}}\\[2mm]
 & \longmapsto &
   \sum_{i,j,k,l=1}^N
(e_{ij}e_{jk}e_{kl}e_{li}
- e_{ij}e_{jk}e_{li}e_{kl}
- e_{ij}e_{ki}e_{jl}e_{lk}
+ e_{ij}e_{ki}e_{lk}e_{jl})\ .
\end{eqnarray*}
\end{xexample}

As we know, $\f_{\gl_N}$ of any diagram always belongs to the centre of $U(\gl_N)$,
so it can be written as a polynomial in $N$ commuting
variables $c_1,\dots,c_N$ (the generalized Casimir elements\index{Casimir
element!generalized}, see, for instance, \cite{Zh}):
$$
  c_s = \sum_{i_1,\dots,i_j=1}^N
        e_{i_1i_2}e_{i_2i_3}\dots e_{i_{s-1}i_s}e_{i_si_1}.
$$
In the graphical notation
$$c_s = \underbrace{\risS{-15}{xm}{}{200}{25}{18}}_{s \mbox{\ {\scriptsize pairs}}}\ .
$$
\newcommand\xone{\risS{-4}{xone}{}{5}{5}{8}}
\newcommand\xzero{\risS{-2}{xzer}{}{15}{5}{8}}
\newcommand\xtwo{\risS{-4}{xtwo}{}{20}{5}{8}}
In particular, $$c_1 = \xone = \sum\limits_{i=1}^N e_{ii}$$ is the
unit matrix (note that it is {\it not\/} the unit of the algebra $U(\g)$),
$$
c_2 = \xtwo = \sum\limits_{i,j=1}^N e_{ij}e_{ji}
$$
is the quadratic Casimir element.
It is convenient to extend the list $c_1,\dots,c_N$ of our variables
by setting $c_0=N$; the graphical notation for $c_0$ will be a
circle. This is especially useful when speaking about the direct limit
$\gl_\infty=\lim_{N\to\infty}\gl_N$.

For instance, the first term in the expansion of $\f(C)$ in the
previous example is nothing but $c_4$; the whole alternating sum,
after some transformations, turns out to be equal to
$c_0^2(c_2-c_1^2)$. Expressing  the values of $\f$ on closed
Jacobi diagrams via the generators $c_i$ is, in general, a
non-trivial operation; a much clearer description exists for the
analog of the map $\f$ defined for the algebra of open diagrams,
see Section~\ref{glN_B}.

\begin{xremark} If the resulting picture contains curves which have
no univalent vertices, then, in the corresponding element of $U(\gl_N)$ 
every such curve is replaced by the numerical factor $N$.
This happens because every such curve leads to a sum where one of the indices does
not appear among the subscripts of the product $e_{i_1i_1}\dots e_{i_sj_s}$,
but the summation over this index still must be done. The proof is
similar to that of the general lemma in Section~\ref{univ_ws_bubble}, where
a different bilinear form is used.
For the diagram $C$ given in that section as an example, we obtain
$$
  \f(C)=\sum_{i,j,k=1}^N e_{ij}e_{ji} = N \sum_{i,j=1}^N e_{ij}e_{ji}.
$$
\end{xremark}

\subsection{Algebra $\gl_N$ with the standard representation}
\label{ws_glN_on_C}
\index{Weight system! $\gl_N$}

The procedure for the closed diagrams repeats what we did with chord diagrams in Section~\ref{ws_glN_on_A}. 
For a closed diagram $C\in \F_n$ with the set $IV$ of $t$ internal
trivalent vertices we double each internal edge and count the number
of components of the resulting curve as before. The only problem
here is how to connect the lines near an internal vertex. This can
be decided by means of a state function $s: IV\to
\{-1,1\}$.\label{gl_N-state}

\begin{xtheorem}[\cite{BN1}] Let $\f^{St}_{\gl_N}$ be the
weight system associated with the standard representation of the Lie
algebra $\gl_N$ with the invariant form $\bilinf{x}{y} =
\mbox{Tr}(xy)$.

For a closed diagram $C$ and a state $s: IV\to \{-1,1\}$ double every
internal edge and connect the lines together in a neighbourhood of a vertex
$v\in IV$ according to the state $s$:
$$\risS{-12}{ws_glN_trv1}{\put(12,19){\mbox{$\scriptstyle v$}}}{29}{20}{10}\ \
  \risS{-2}{totor}{}{25}{0}{0}\
  \risS{-12}{ws_glN_trv2}{}{29}{0}{0}\ ,\mbox{\ if\ }s(v)=1;\qquad
 \risS{-12}{ws_glN_trv1}{\put(12,19){\mbox{$\scriptstyle v$}}}{29}{0}{0}\ \
  \risS{-2}{totor}{}{25}{0}{0}\
 \risS{-12}{ws_glN_trv3}{}{29}{0}{0}\ ,\mbox{\ if\ }s(v)=-1;
$$
and replace each external vertex as follows\quad
$\risS{-5}{ws_glN_ev1}{}{18}{15}{15}\ \
  \risS{0}{totor}{}{25}{0}{0}\
  \risS{-5}{ws_glN_ev2}{}{18}{0}{0}\ .$

Let $|s|$ denote the number of components of the curve obtained in
this way. Then
$$\f^{St}_{\gl_N}(C) = \sum_s\ \Bigl(\prod_v s(v)\Bigr)\ N^{|s|}\ ,
$$
where the product is taken over all $t$ internal vertices of $C$,
and the sum is taken over all the $2^t$ states for $C$.
\end{xtheorem}

A straightforward way to prove this theorem is to use the STU
relation and the theorem of Section~\ref{ws_glN_on_A}. We leave the
details to the reader.

\begin{xexample} Let us compute the value
$\f^{St}_{\gl_N}\Bigl(\chvad{chvar_ex8}\Bigr)$.
\def\stglcd#1#2#3{\begin{array}{c}
         \risS{-12}{#1}{}{30}{20}{15}\\
         \scriptstyle \prod s(v) = #2 \\
         \scriptstyle |s| = #3
         \end{array}}
There are four resolutions of the triple points:
$$\stglcd{chvar_ex9}{1}{4}\qquad \stglcd{chvard_ex10}{-1}{2}\qquad
  \stglcd{chvard_ex11}{-1}{2}\qquad \stglcd{chvard_ex12}{1}{2}
$$
Therefore, $\f^{St}_{\gl_N}\Bigl(\chvad{chvar_ex8}\Bigr) = N^4-N^2$.
\end{xexample}

Other properties of the weight system $\f^{St}_{\gl_N}$ are
formulated in exercises \ref{ex_gl_N_on_C_sl_N} -- \ref{ex_gl_N_on_C_wheel}.

\subsection{Algebra $\so_N$ with standard representation}
\label{ws_soN_on_C} \index{Weight system! $\so_N$} Here, a {\em
state} for $C\in \F_n$ will be a function $s: IE\to \{-1,1\}$ on the
set $IE$ of internal edges (those which are not on the Wilson loop).
The value of a state indicates the way of doubling the corresponding
edge:
$$\risS{1}{ws_soN_edge}{\put(12,4){\mbox{$\scriptstyle e$}}}{29}{10}{8}\ \
  \risS{-2}{totor}{}{25}{0}{0}\ \
  \risS{-2}{pso}{}{29}{0}{0}\ ,\mbox{\ if\ }s(e)=1;\quad
 \risS{1}{ws_soN_edge}{\put(12,4){\mbox{$\scriptstyle e$}}}{29}{0}{0}\ \
  \risS{-2}{totor}{}{25}{0}{0}\ \
 \risS{-2}{mso}{}{29}{0}{0}\ ,\mbox{\ if\ }s(e)=-1.
$$
In the neighbourhoods of trivalent and external vertices we connect
the lines in the standard fashion as before. For example, if the values
of the state on three edges $e_1$, $e_2$, $e_3$ meeting at a vertex $v$
are $s(e_1)=-1$, $s(e_2)=1$, and $s(e_3)=-1$, then we resolve it as follows:
$$\risS{-12}{ws_glN_trv1}{\put(13,19){\mbox{$\scriptstyle v$}}
  \put(-7,28){\mbox{$\scriptstyle e_1$}}
  \put(-8,-2){\mbox{$\scriptstyle e_2$}}
  \put(30,12){\mbox{$\scriptstyle e_3$}}}{29}{20}{15}\qquad
  \risS{-2}{totor}{}{25}{20}{15}\quad
  \risS{-12}{ws_soN_mpm1}{}{29}{20}{15} \ .
$$
As before, $|s|$ denotes the number of components of the curve obtained in
this way.

\begin{xtheorem}[\cite{BN1}]\label{soN_stand_C}
Let $\f^{St}_{\so_N}$ be the weight system associated with the
standard representation of the Lie algebra $\so_N$ with the
invariant form $\bilinf{x}{y} = \mbox{Tr}(xy)$. Then
$$\f^{St}_{\so_N}(C) = 2^{-\deg{C}}
   \sum_s\ \Bigl(\prod_e s(e)\Bigr)\ N^{|s|}\ ,$$
where the product is taken over all internal edges of $C$ and the
sum is taken over all the states $s:IE(C)\to\{1,-1\}$.
\end{xtheorem}

\begin{proof}
First let us note that $\deg{C}={\#(IE) - \#(IV)}$, where $\#(IV)$
and $\#(IE)$ denote the numbers of internal vertices and edges
respectively. We prove the Theorem by induction on $\#(IV)$.

If $\#(IV)=0$ then $C$ is a chord diagram. In this case the Theorem
coincides with the Theorem of Section~\ref{ws_so_n_St}, page
\pageref{ws_so_n_St}.

If $\#(IV)\not=0$ we can use the STU relation to decrease the number
of internal vertices. Thus it remains to prove that the formula for
$\f^{St}_{\so_N}$ satisfies the STU relation. For this we split the
8 resolutions of the $S$ diagram corresponding to the various values
of $s$ on the three edges of $S$ into two groups which can be
deformed to the corresponding resolutions of the $T$ and the $U$
diagrams:
$$\begin{array}{c@{\qquad}c@{\qquad}l}
\risS{-12}{ws_soN_S}{\put(-10,10){$S$}}{29}{20}{20} &
\risS{-2}{totor}{}{25}{0}{0} &
  \Biggl(\
  \risS{-10}{ws_soN_ppp}{}{29}{0}{0}\ -\ \risS{-10}{ws_soN_mpp}{}{29}{0}{0}\
  -\ \risS{-10}{ws_soN_pmp}{}{29}{0}{0}\ +\ \risS{-10}{ws_soN_mmp}{}{29}{20}{15}\
  \Biggr) \\
&&\hspace{-12pt}-\ \Biggl(\
  \risS{-10}{ws_soN_mmm}{}{29}{0}{0}\ -\ \risS{-10}{ws_soN_pmm}{}{29}{0}{0}\
  -\ \risS{-10}{ws_soN_mpm}{}{29}{0}{0}\ +\ \risS{-10}{ws_soN_ppm}{}{29}{0}{0}\
   \Biggr)\ ;
\end{array}
$$\label{8_states_so_N}
$$\begin{array}{c@{\qquad}c@{\qquad}l}
\risS{-12}{ws_soN_T}{\put(-10,10){$T$}}{29}{20}{20} &
\risS{-2}{totor}{}{25}{0}{0} &
  \risS{-10}{ws_soN_Tpp}{}{29}{0}{0}\ -\ \risS{-10}{ws_soN_Tpm}{}{29}{0}{0}\
  -\ \risS{-10}{ws_soN_Tmp}{}{29}{0}{0}\ +\ \risS{-10}{ws_soN_Tmm}{}{29}{20}{15}\ ;
\end{array}
$$
$$\begin{array}{c@{\qquad}c@{\qquad}l}
\risS{-12}{ws_soN_U}{\put(-10,10){$U$}}{29}{20}{20} &
\risS{-2}{totor}{}{25}{0}{0} &
  \risS{-10}{ws_soN_Upp}{}{29}{0}{0}\ -\ \risS{-10}{ws_soN_Upm}{}{29}{0}{0}\
  -\ \risS{-10}{ws_soN_Ump}{}{29}{0}{0}\ +\ \risS{-10}{ws_soN_Umm}{}{29}{0}{0}\ .
\end{array}
$$
\end{proof}

\begin{xexample}
$$\f^{St}_{\so_N}\Bigl( \risS{-12}{ws_soN_ex}{}{30}{20}{15} \Bigr)
  = \frac{1}{4}(N^3 -3N^2 + 2N)\ =\ \frac{1}{4}N(N-1)(N-2)\ .
$$
\end{xexample}

\subsection{A small table of values of $\f$}

The following table shows the values of $\f$ on the generators
of the algebra $\F$ of degrees $\leq4$:
$$
t_1=\rb{-4mm}{\ig[height=10mm]{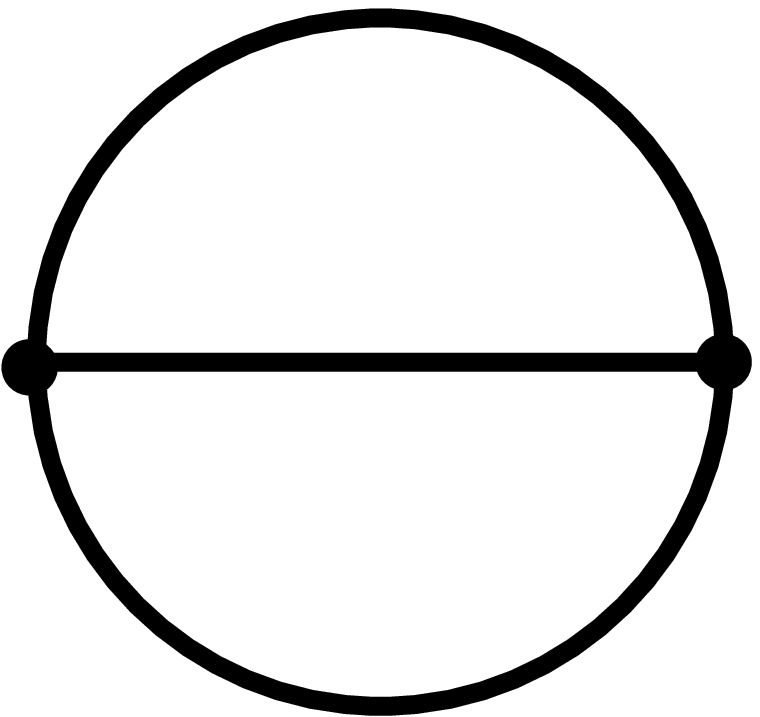}},\
t_2=\rb{-4mm}{\ig[height=10mm]{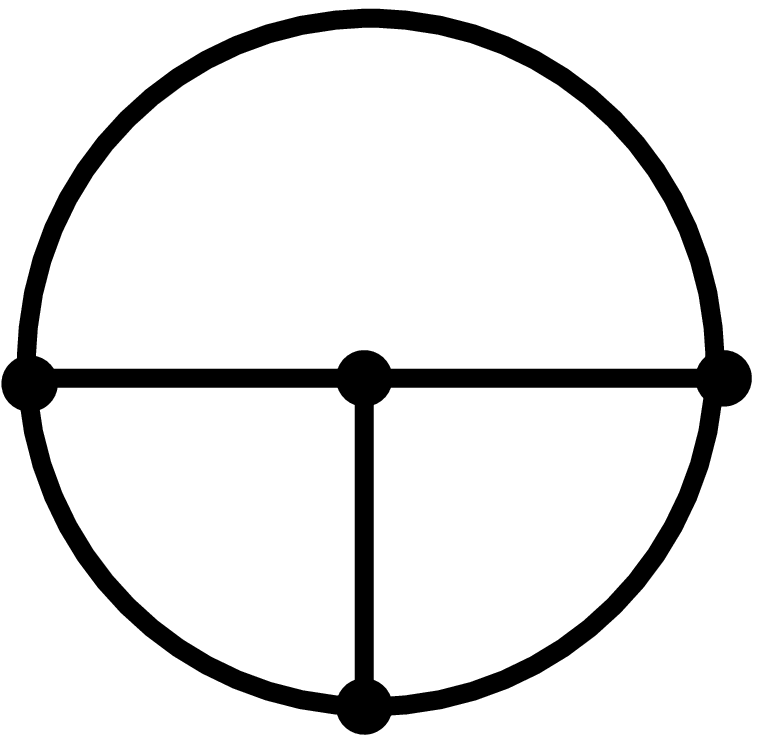}},\
t_3=\rb{-4mm}{\ig[height=10mm]{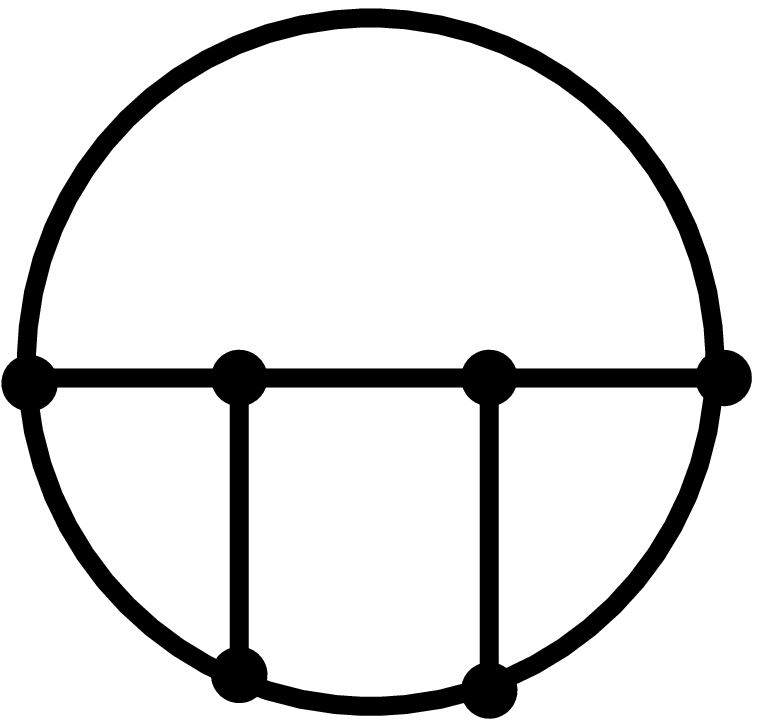}},\
t_4=\rb{-4mm}{\ig[height=10mm]{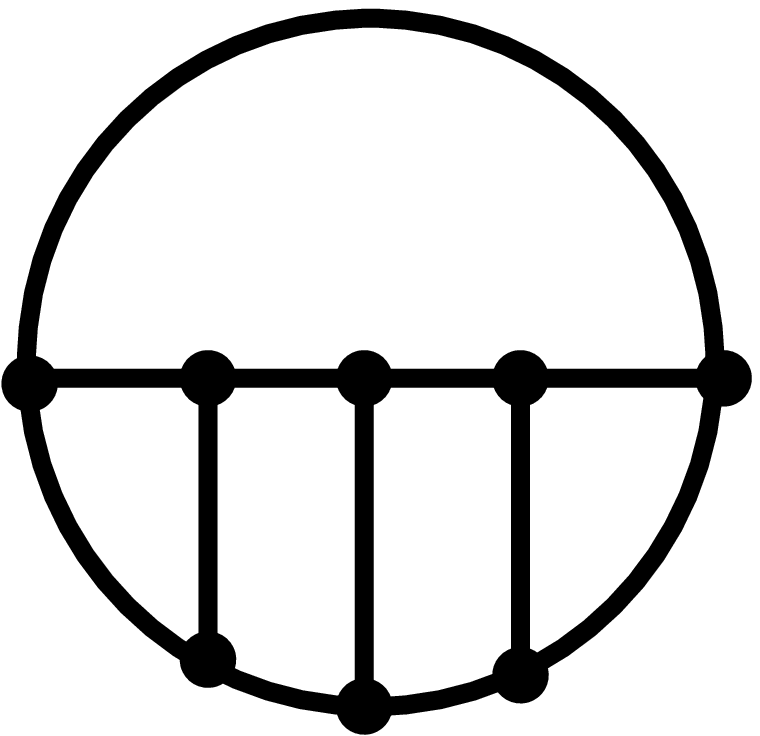}},\
w_4=\rb{-4mm}{\ig[height=10mm]{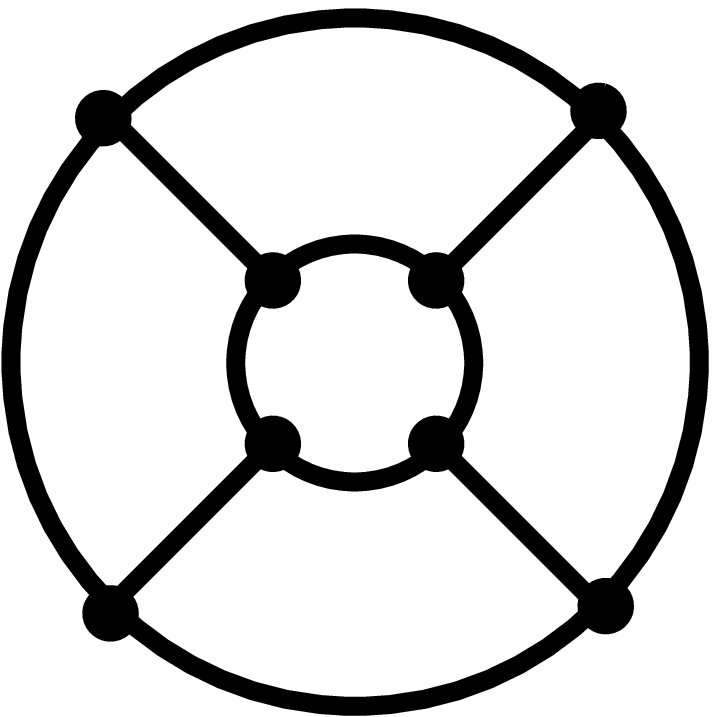}}
$$
for the simple Lie algebras
$A_1$, $A_2$, $A_3$, $A_4$, $B_2$, $B_3$, $C_3$, $D_4$, $G_2$,
computed by A.\,Kaishev \cite{Kai}.

$$\begin{array}{l|c|c|c|c|c}
     & t_1  & t_2    & t_3      & t_4      & w_4 \\ \hline
 A_1 & c & -2 c  & 4 c   & -8 c    & 8 c^2 \\ \hline
 A_2 & c & -3 c  & 9 c   & -27 c   & 9 c^2 + 9 c \\ \hline
 A_3 & c & -4 c  & 16 c  & -64 c   & e\\ \hline
 A_4 & c & -5 c  & 25 c  & -125 c  & e \\ \hline
 B_2 & c & -3/2 c& 9/4 c & -27/8 c & d \\ \hline
 B_3 & c & -5/2 c& 25/4 c& -125/8 c & d \\ \hline
 C_3 & c & -4 c  & 16 c  & -64 c   & d \\ \hline
 D_4 & c & -3 c  & 9 c   & -27 c   & 3 c^2 + 15 c \\ \hline
 G_2 & c & -2 c  & 4 c   & -8 c    & 5/2 c^2 + 11/3 c
\end{array}$$

Here $c$ is the quadratic Casimir element of the corresponding enveloping
algebra $U(\g)$, while $d$ and $e$ are the following (by degree)
independent generators of $ZU(\g)$. Note that in this table all $d$'s and
$e$'s have degree 4 and are defined modulo elements of smaller degrees.
The exact expressions for $d$ and $e$ can be found in \cite{Kai}.

A look at the table shows that the mapping $\f$ for almost
all simple Lie algebras has a non-trivial kernel. In fact,
$\f_\g(t_1t_3-t_2^2)=0$.

\begin{xca}
Find a metrized Lie algebra $\g$ such that the mapping $\f_\g$ has a
non-trivial cokernel.
\end{xca}

\section{Lie algebra weight systems for the algebra $\B$}
\label{LAWS_B}

The construction of the Lie algebra weight systems for open Jacobi
diagrams is very similar to the procedure for closed diagrams. For a
metrized Lie algebra $\g$ we construct a weight system
$\rho_\g:\B\to S(\g)$, defined on the space of open diagrams $\B$
and taking values in the symmetric algebra of the vector space $\g$
(in fact, even in its $\g$-invariant subspace $S(\g)^\g$).

Let $O\in\OD$ be an open diagram. Choose an order on the set of its
univalent vertices; then $O$ can be treated as the internal graph of
some closed diagram $C_O$. Following the recipe of
Section~\ref{LAWS_C}, construct a tensor $T_\g(C_O)\in\g^{\ot m}$,
where $m$ is the number of legs of the diagram $O$. Now we define
$\rho_\g(O)$ as the image of the tensor $T_\g(C_O)$ in $S^m(\g)$
under the natural projection of the tensor algebra on $\g$ onto
$S(\g)$.\label{def:psiB}

The choice of an order on the legs of $O$ is of no importance.
Indeed, it amounts to choosing an order on the tensor factors in the
space $\g^{\ot m}$ to which the tensor $T_\g(C_O)$ belongs. Since
the algebra $S(\g)$ is commutative, the image of $T_\g(C_O)$ is
always the same.

\subsection{The formal PBW theorem}

The relation between the Lie algebra weight systems for the open
diagrams and for the closed diagrams is
expressed by the following theorem.

\begin{xtheorem}
For any metrized Lie algebra $\g$ the diagram
\begin{equation*}\begin{CD}
\B  @>\rho_{\g} >>  S(\g) \\
@V{\chi}VV          @VV{\beta_{\g}}V \\
\F     @>>\f_{\g}>   U(\g)
\end{CD}\end{equation*}
commutes.
\end{xtheorem}

\begin{proof}
The assertion becomes evident as soon as one recalls the definitions
of all the ingredients of the diagram: the isomorphism $\chi$
between the algebras $\F$ and $\B$ described in section \ref{A=B},
the weight systems $\f_{\g}$ and $\rho_{\g}$, defined in sections
\ref{LAWS_A} and \ref{LAWS_B}, and $\beta_{\g}$, the
Poincar\'{e}--Birkhoff--Witt isomorphism taking an element
$x_1x_2...x_n$ into the arithmetic mean of
$x_{i_1}x_{i_2}...x_{i_n}$ over all permutations $(i_1,i_2,...,i_n)$
of the set $\{1,2,\dots,n\}$. Its restriction to the invariant
subspace $S(\g)^\g$ is a vector space isomorphism with the centre of
$U(\g)$.
\end{proof}

\begin{example}
Let $\g$ be the Lie algebra $\so_3$. It has a basis $\{a,b,c\}$
which is orthonormal with respect to the Killing form
$\langle\cdot,\cdot\rangle^K$ and with the commutators $[a,b]=c$,
$[b,c]=a$, $[c,a]=b$. As a metrized Lie algebra $\so_3$ is
isomorphic to the Euclidean 3-space with the cross product as a Lie
bracket. The tensor that we put in every trivalent vertex in this
case is
\begin{eqnarray*}
  -J &=& -a\wedge b\wedge c \\  &\hspace{-40pt}=&\hspace{-20pt}
      - a\ot b \ot c - b\ot c \ot a - c\ot a \ot b
      + b\ot a \ot c + c\ot b \ot a + a\ot c \ot b.
\end{eqnarray*}
Since the basis is orthonormal, the only way to get a non-zero
element in the process of contraction along the edges is to choose
the same basis element on either end of each edge. On the other
hand, the formula for $J$ shows that in every vertex we must choose
a summand with different basis elements corresponding to the 3
edges. This leads to the following algorithm for computing the
tensor $T_{\so_3}(O)$ for a given diagram $O$: one must list all
3-colourings of the edges of the graph by 3 colours $a$, $b$, $c$ such
that the 3 colours at every vertex are always different, then sum up
the tensor products of the elements written on the legs, each taken
with the sign $(-1)^s$, where $s$ is the number of {\it negative\/}
vertices (that is, vertices where the colours, read counterclockwise,
come in the negative order $a$, $c$, $b$).

For example, consider the diagram (the {\it Pont-Neuf diagram with
parameters\index{Diagram!Pont-Neuf} $(1,3)$} in the terminology of
O.~Dasbach \cite{Da3}, see also page~\pageref{pont_neuf} below):
$$\index{Pont-Neuf diagram}
  O\ =\ \rb{-18mm}{\ig[width=35mm]{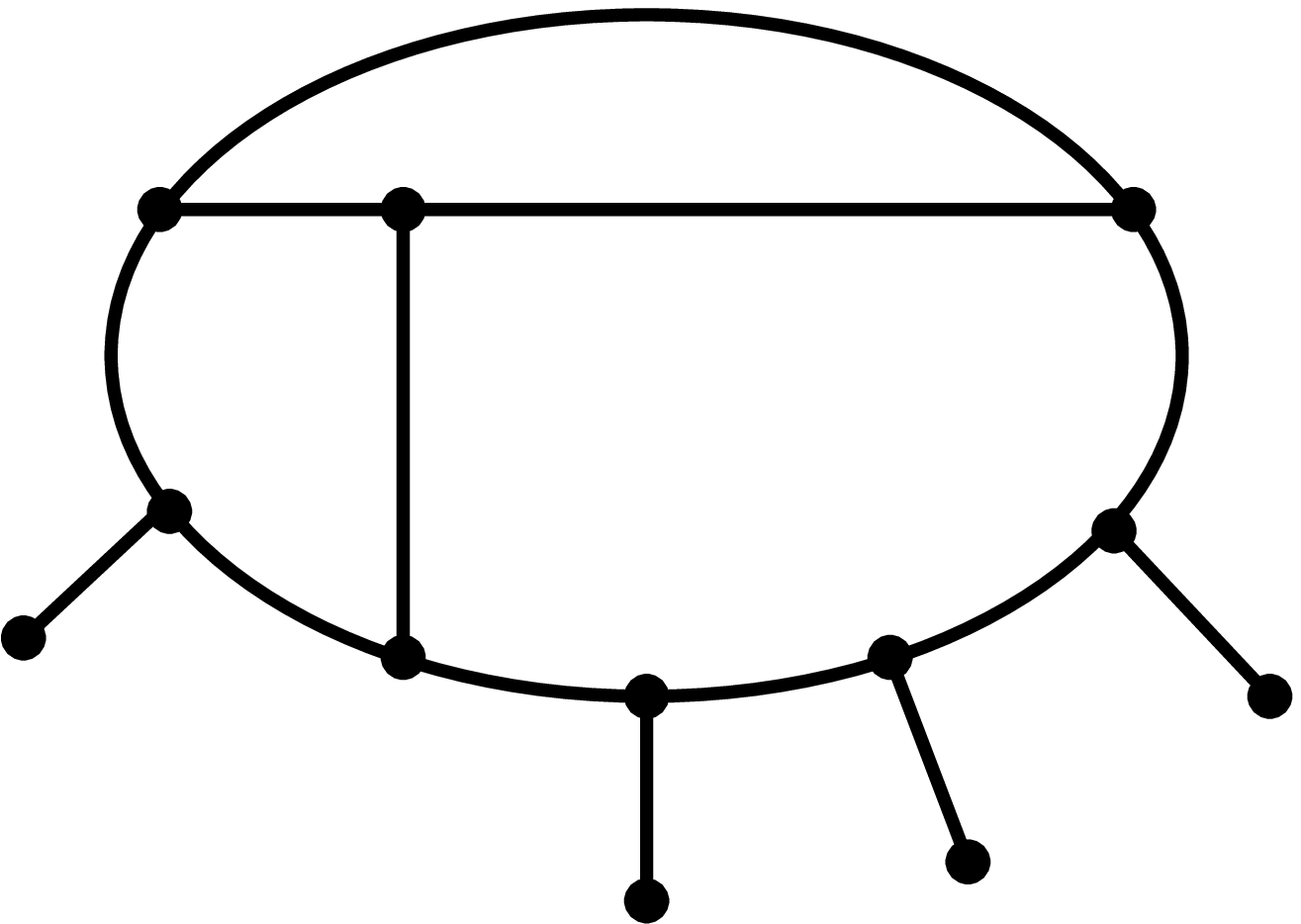}}
$$
It has 18 edge 3-colourings, which can be obtained from the following three
by permutations of $(a,b,c)$:

$$
  \ig[width=35mm]{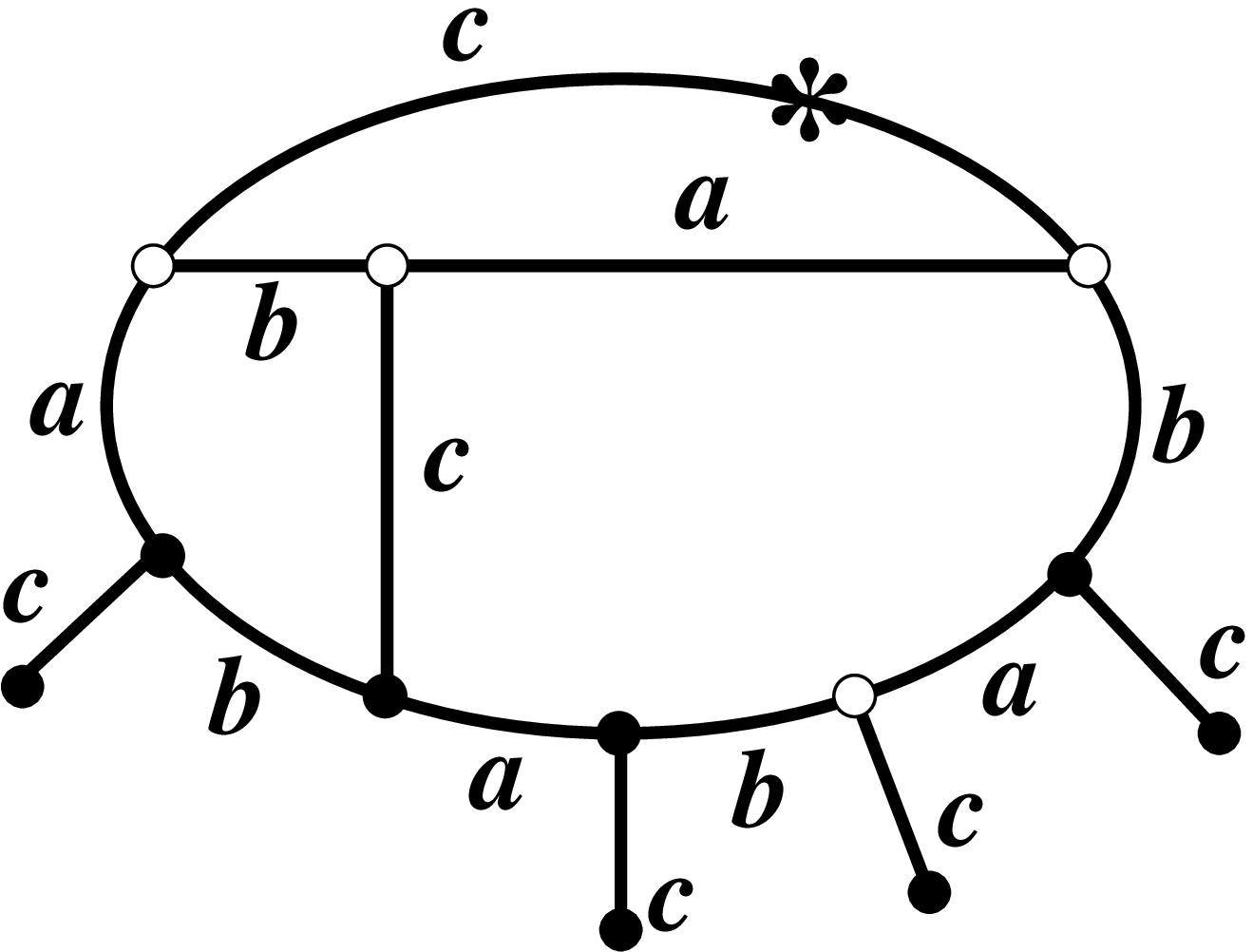}\qquad
  \ig[width=35mm]{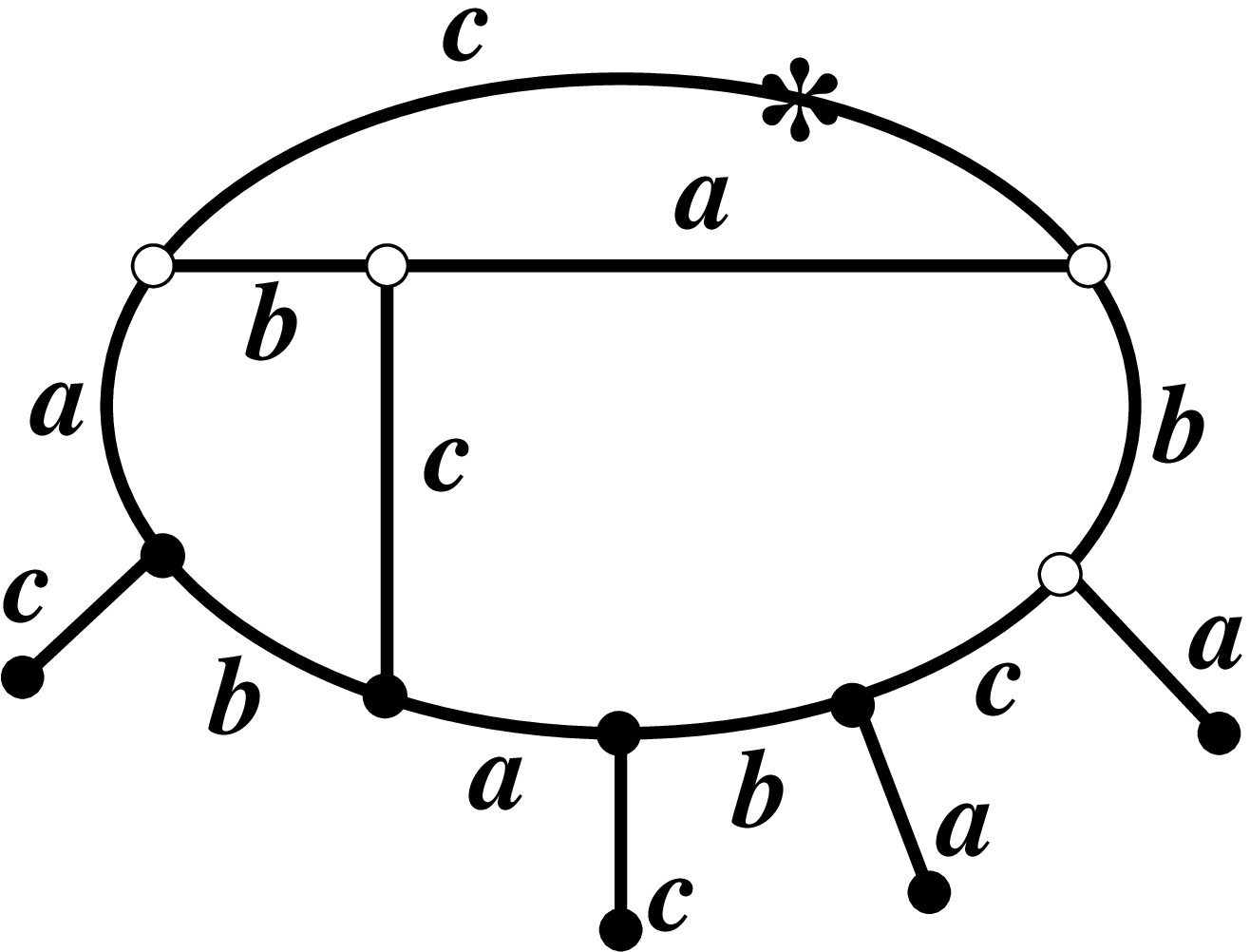}\qquad
  \ig[width=35mm]{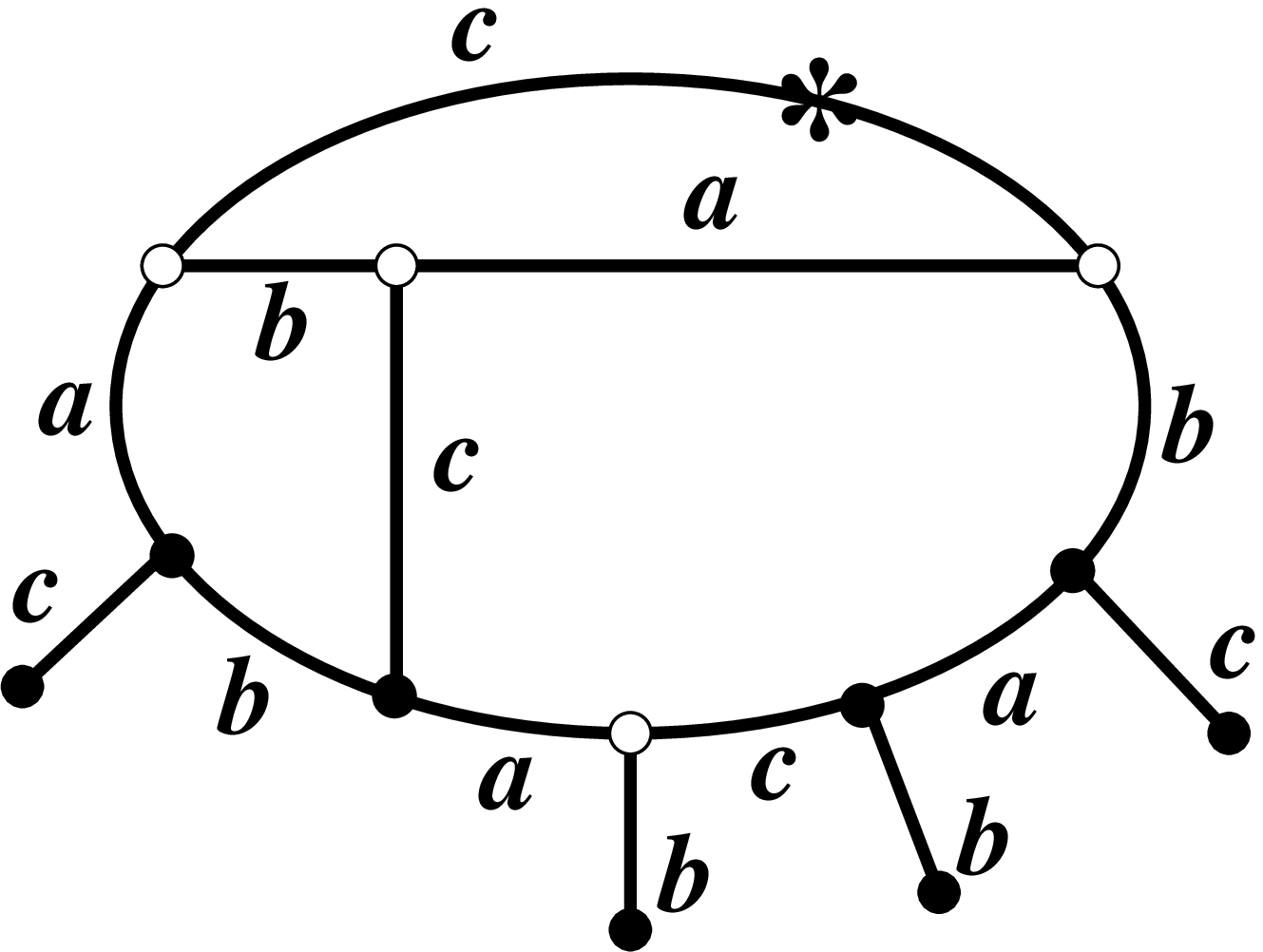}
$$
In these pictures, negative vertices are marked by small empty circles.
Writing the tensors in the counterclockwise order starting from the marked
point, we get:
\begin{align*}
  2 (&a\ot a\ot a\ot a + b\ot b\ot b\ot b+ c\ot c\ot c\ot c) \\
+ &a\ot b\ot b\ot a + a\ot c\ot c\ot a + b\ot a\ot a\ot b \\
+ &b\ot c\ot c\ot b + c\ot a\ot a\ot c + c\ot b\ot b\ot c \\
+ &a\ot a\ot b\ot b+ a\ot a\ot c\ot c+ b\ot b\ot a\ot a \\
+ &b\ot b\ot c\ot c+ c\ot c\ot a\ot a+ c\ot c\ot b\ot b.
\end{align*}
Projecting onto the symmetric algebra, we get:
$$
  \rho_{\so_3}(O)=2(a^2+b^2+c^2)^2.
$$

This example shows that the weight system defined by the Lie algebra
$\so_3$, is closely related to the 4-colour theorem, see \cite{BN3}
for details.
\end{example}

\begin{example}\label{rhootkolesa}
For an arbitrary metrized Lie algebra $\g$ let us calculate
$\rho_{\g}(w_{n})$ where $w_{n}\in\B$ is the \index{Wheel!in $\B$}
{\em wheel with $n$ spokes}:
$$w_{n} :=\quad \risS{-3}{wssl2wn-b}{\put(-3,-8){\mbox{\scriptsize $n$
spokes}}}{30}{20}{15}$$
Note that $n$ must be even; otherwise by Lemma~\ref{anti-auto}
$w_n=0$.

Dividing the wheel into $n$ tripods, contracting the resulting
tensors of rank 3 and projecting the result to $S(\g)$ we get
$$
c_{j_1 i_1 j_2}\ldots c_{j_n i_n j_{1}}\cdot e_{i_1}\ldots e_{i_n}
=\Tr{(\ad{\, e_{i_1}}\ldots \ad{\, e_{i_n}})}\cdot e_{i_1}\ldots
e_{i_n},$$ where $\{e_i \}$ is an orthonormal basis for $\g$, and
the summation by repeating indices is implied.
\end{example}
\smallskip

\subsection{The universal $\gl_N$ weight system for the algebra $\B$}
\label{glN_B}

The $\gl_N$ weight system for the algebra $\B$ of open Jacobi diagrams
is computed in exactly the same way as for the closed diagrams (see Section
\ref{univ_ws_glN}), only now we treat the variables $e_{ij}$ as commuting
elements of $S(\gl_N)$.
For instance, the diagram $B=\rb{-3mm}{\ig[height=7mm]{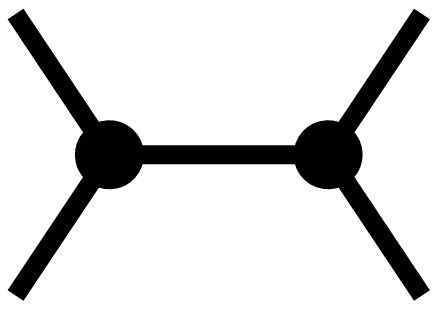}}$
obtained by stripping the Wilson loop off the diagram $C$ of
Section \ref{univ_ws_glN}, goes to 0 under the mapping $\rho_{\gl_N}$,
because all the four summands in the corresponding alternating sum
become now equal.

In general, as we know that the invariant part  $S(\gl_N)^{\gl_N}$ of the 
algebra $S(\gl_N)$ is isomorphic to the centre of $U(\gl_N)$, it is also freely 
generated by the Casimir elements $c_1,\dots,c_N$. Here is an example, 
where we, as above, write $c_0$ instead of $N$:

\begin{xexample}
\newcommand\figb[1]{\risS{-10}{#1}{}{30}{10}{15}}
$$\label{B2_ex}
\begin{array}{rcl}
\rho_{gl_N}\Bigl(\ \figb{bt}\ \Bigr) &=& \figb{btpp}\ -\
\figb{btmp}\ -\
       \figb{btpm}\ +\ \figb{btmm} \\ 
&=& \xzero\ \xtwo\ -\ \xone\ \xone\ -\ \xone\ \xone\ +\ \xtwo\ \xzero
\ =\ 2(c_0c_2 - c_1^2).
\end{array}
$$
\end{xexample}

\subsection{Invariants of string links and the algebra of necklaces}

Recall that the algebra $\A(n)$ of closed diagrams for string links on $n$ strands (see \ref{A-ot-p}).
has a $\B$-analog, denoted by $\B(n)$ and called the {\em algebra of coloured open Jacobi diagrams}, see page~\pageref{B-of-m}.
In this section we shall describe the weight system generalizing $\rho_{\gl_N}:\B\to S(\gl_N)$ to a mapping $$\rho^{(n)}_{\gl_N}:\B(n)\to S(\gl_N)^{\otimes n}.$$

A diagram in $\B(n)$ is an open Jacobi diagram with univalent
vertices marked by numbers between 1 and $n$ (or coloured by $n$ colours).
The vector space spanned by these elements modulo AS and IHX relations
is what we call $\B(n)$. The colour-respecting averaging map
$\chi_n:\B(n)\to\A(n)$, defined similarly to the simplest case
$\chi:\B\to\A$ (see Section \ref{A=B}), is a linear isomorphism (see \cite{BN4}).

Given a coloured open Jacobi diagram, we consider positive and negative
resolutions of all its $t$ trivalent vertices and get the alternating sum of
$2^t$ pictures as on page \pageref{ex_ws_glN} with the univalent legs
marked additionally by the colours. For each resolution, mark the connected
components by different variables $i$, $j$ etc, then add small arcs near the univalent vertices and obtain a set of oriented closed curves. To each small arc 
(which was a univalent vertex before) there corresponds
a pair of indices, say $i$ and $j$. Write $e_{ij}$ in the tensor factor of
$S(\gl_N)^{\otimes n}$ whose number is the number of that univalent vertex,
and where $i$ and $j$ go in the order consistent with the orientation on the
curve. Then take the sum over all subscripts from 1 to $N$.

To make this explanation clearer, let us illustrate it on a concrete example.
Take the coloured diagram
$$
  D\quad = \quad\rb{-10mm}{\ig[height=21mm]{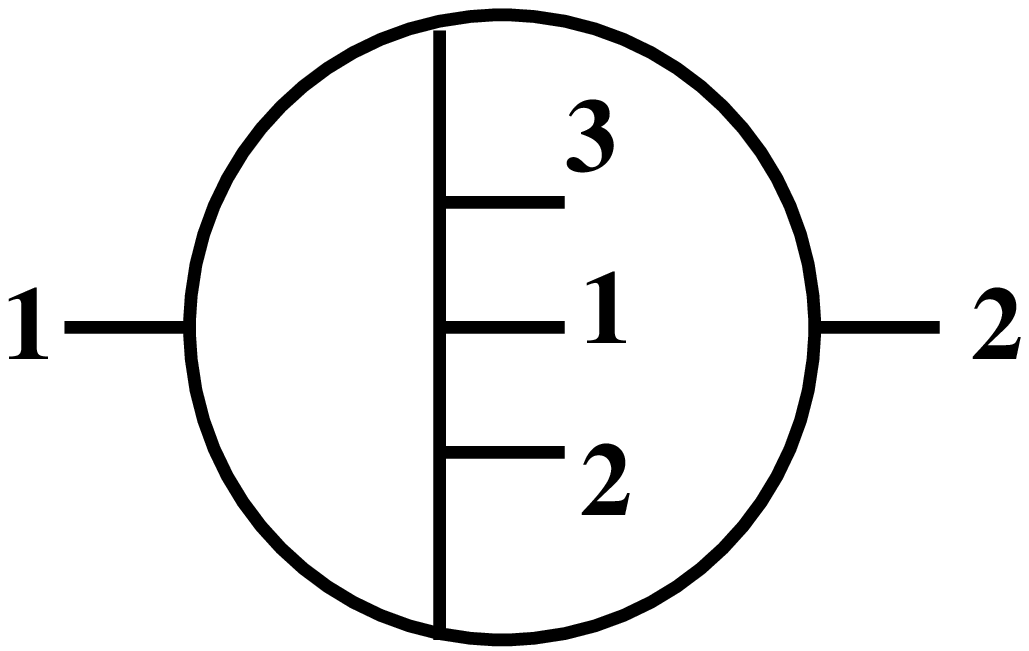}}
$$
with the blackboard (counterclockwise) cyclic order of edges meeting at
trivalent vertices. Resolving all the trivalent vertices positively, we get the following
collection of directed curves:
$$
  \rb{-11mm}{\ig[height=23mm]{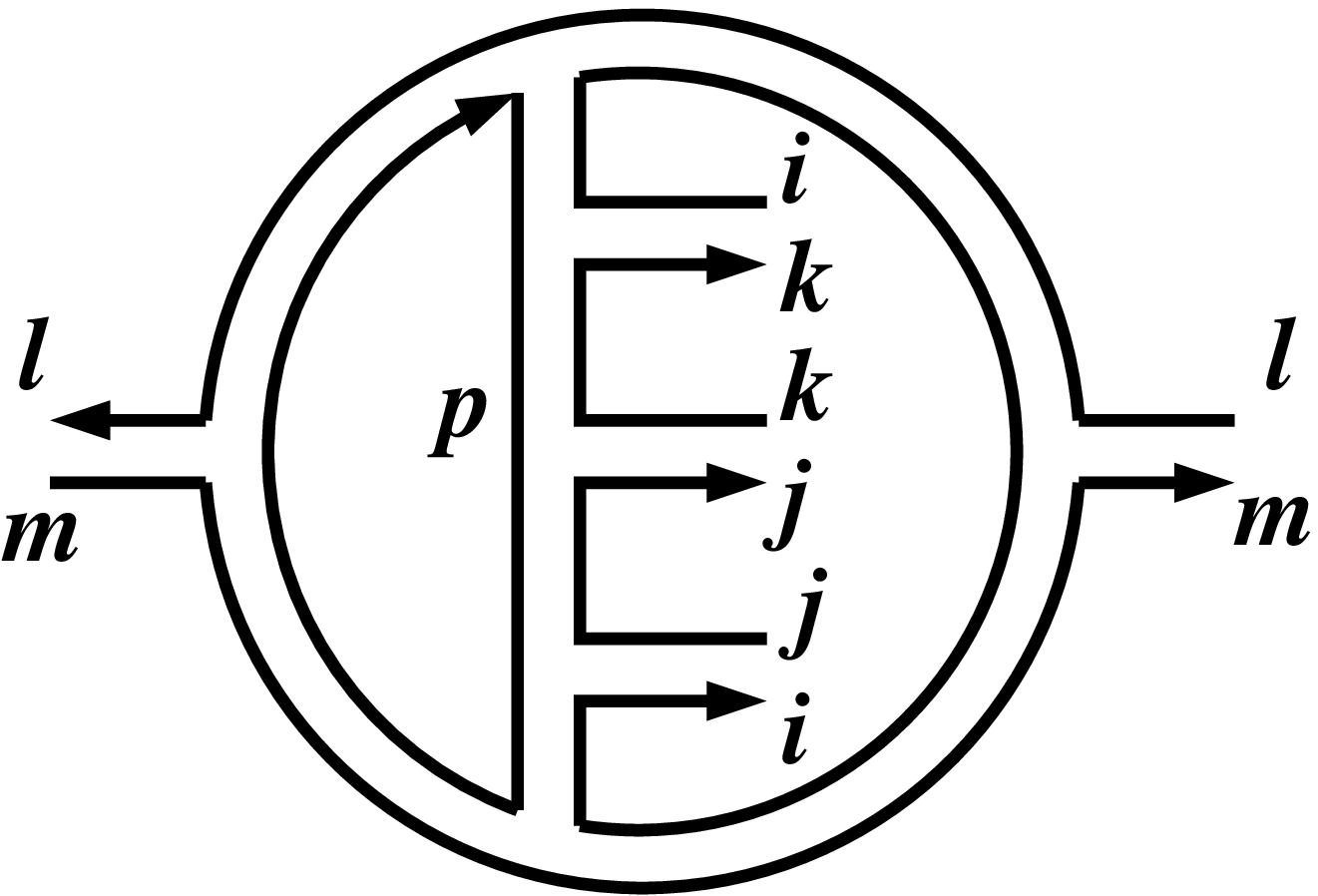}}
$$
which, according to the above procedure, after filling in the gaps at
univalent vertices, transcribes as the following
element of $S(\gl_N)^{\otimes3}$:
$$
  \sum_{i,j,k,l,m,p=1}^N e_{lm}e_{jk}\ot e_{ml}e_{ij}\ot e_{ki}
=N\cdot\sum_{i,j,k=1}^N e_{jk}\ot e_{ij}\ot e_{ki}
  \cdot\sum_{l,m=1}^N e_{lm}\ot e_{ml}\ot 1.
$$

We see that the whole expression is the product of three elements
corresponding to the three connected components of the closed curve.
In particular, the factor $N$ corresponds to the circle without univalent
vertices and can be represented alternatively as multiplication
by $\sum_{n=1}^N 1\ot1\ot1$.

As the choice of notations for the summation indices does not matter,
we can write the obtained formula schematically as the product of three
{\em necklaces}:
$$
       \rb{-5mm}{\ig[height=11mm]{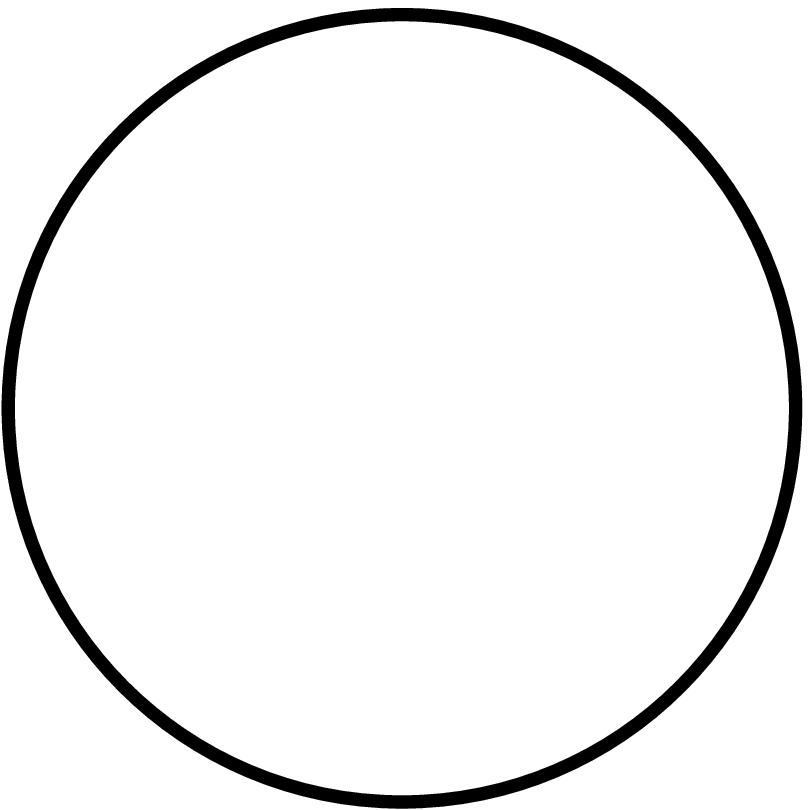}}
  \ \cdot\ \rb{-5mm}{\ig[height=11mm]{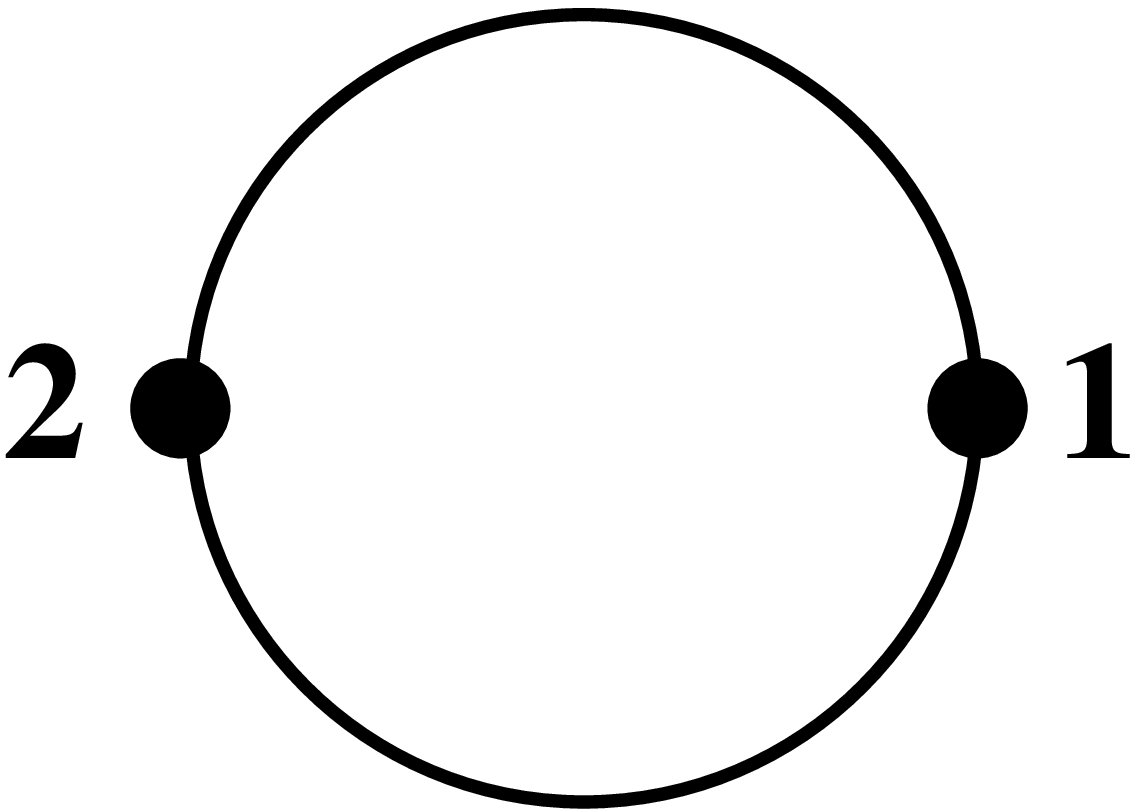}}
  \ \cdot\ \rb{-5mm}{\ig[height=11mm]{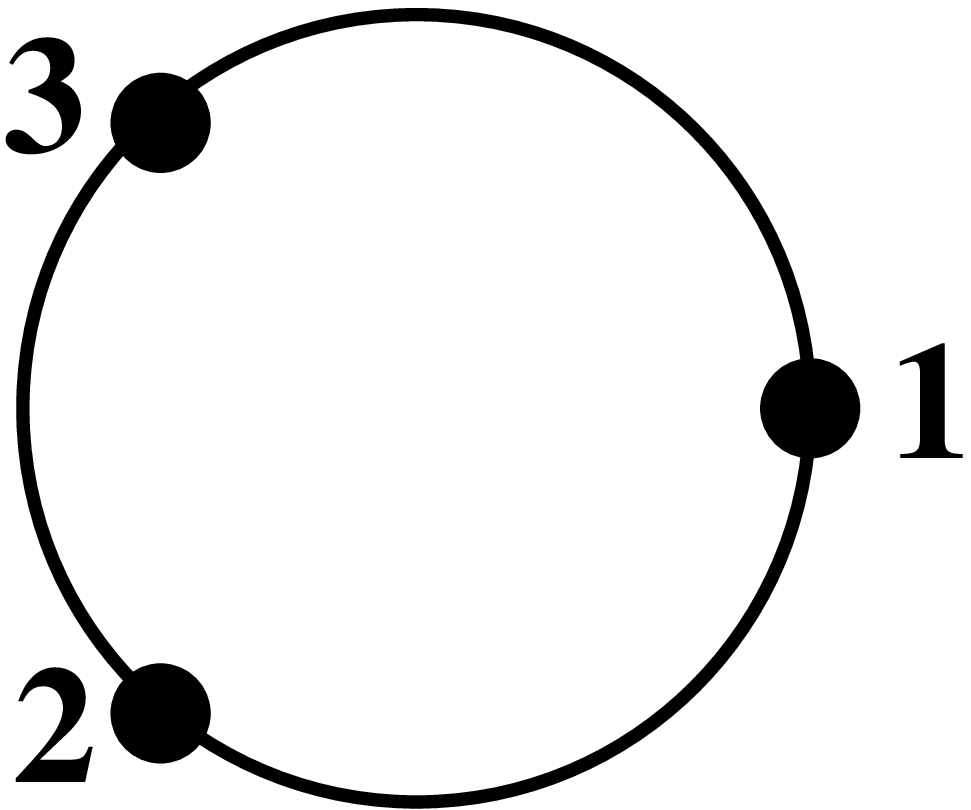}}
$$

An $n$-coloured {necklace}\label{necklace}\index{Necklace} is an
arrangement of several {\em beads}, numbered between 1 to $n$, along
an oriented circle (the default orientation is counterclockwise).
A necklace can be uniquely denoted by a letter, say $x$,
with a subscript consisting of the sequence of bead numbers chosen to be
lexicographically smallest among all its cyclic shifts.
Any $n$-coloured necklace corresponds to an element of
the tensor power of $S(\gl_N)$ according to the following rule.
Mark each arc of the circle between two beads by a different integer variable
$i$, $j$, etc. To each bead we assign the element $e_{ij}$, where $i$ is
the variable written on the incoming arc and $j$, on the outgoing arc.
Then compose the tensor product of all these $e_{ij}$'s putting each into
the tensor factor of $S(\gl_N)^{\ot n}$ whose number is the number of the
bead under consideration, and take the sum of these expressions where each
integer variable runs from 1 to $N$.

Examples (for $n=3$):
\begin{align*}
  &x_{123}:=\ \rb{-4.5mm}{\ig[height=11mm]{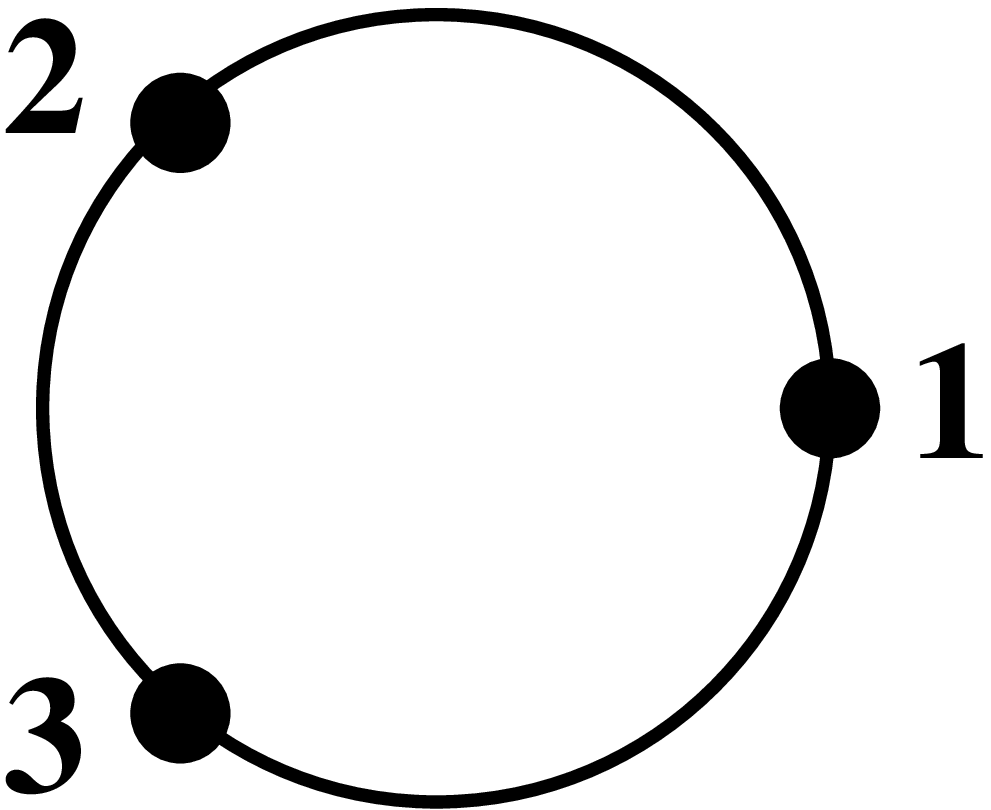}}\ \mapsto
          \ \sum_{i,j,k=1}^N e_{ij}\ot e_{jk}\ot e_{ki} \\
  &x_{132}:=\ \rb{-4.5mm}{\ig[height=11mm]{nl132.eps}}\ \mapsto
          \ \sum_{i,j,k=1}^N e_{jk}\ot e_{ij}\ot e_{ki} \\
  &x_{12123}:=\ \rb{-5.5mm}{\ig[height=13mm]{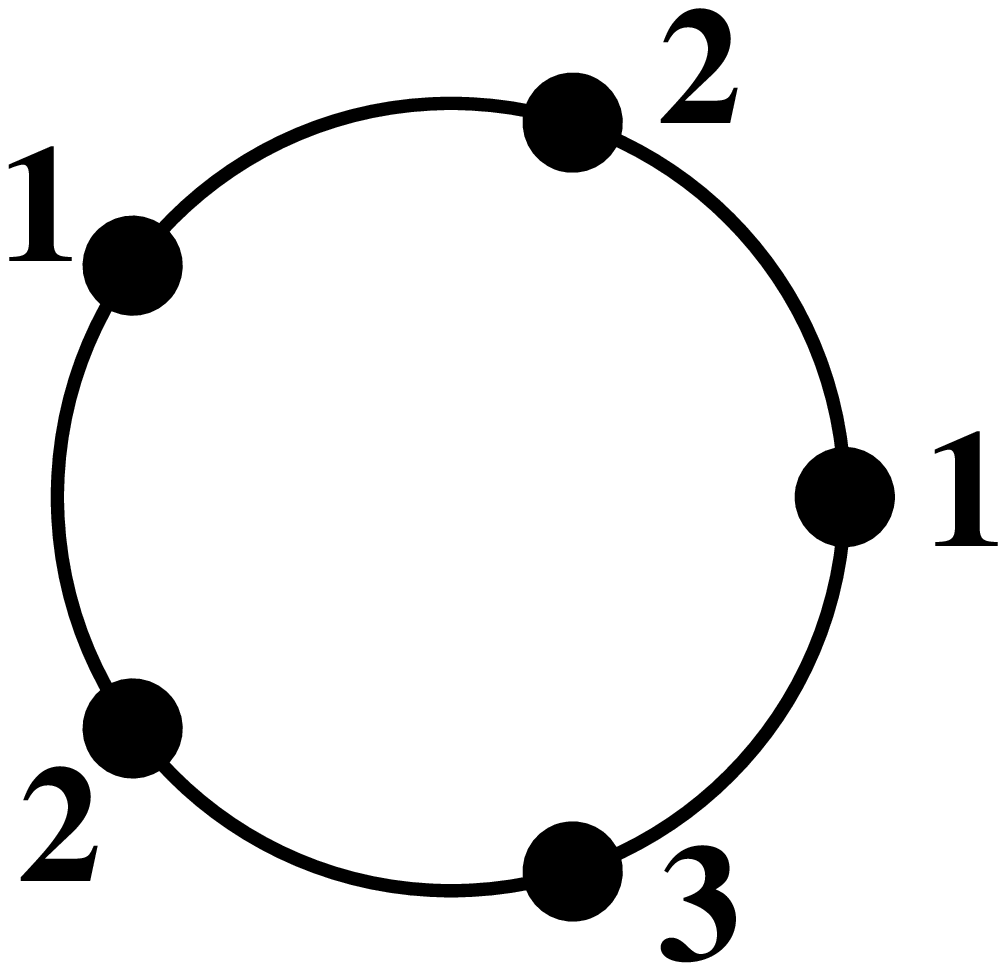}}\ \mapsto
          \ \sum_{i,j,k,l.m=1}^N e_{ij}e_{kl}\ot e_{jk}e_{lm}\ot e_{mi}
\end{align*}
(All the circles are oriented counterclockwise.)

We will call such elements of $S(\gl_N)^{\ot n}$ the \textit{necklace
elements}\index{Necklace!element}.
By a theorem of S.\,Donkin \cite{Don}, the $\gl_N$-invariant subspace of
the algebra $S(\gl_N)^{\otimes n}$ is generated by the necklace elements,
and the algebraic relations between them may exist for small values of
$N$, but disappear as $N\to\infty$, so that the invariant subspace of the
direct limit $S(\gl_\infty)^{\otimes n}$ is isomorphic to the free polynomial
algebra generated by $n$-coloured necklaces.

Summing up, we can formulate the algorithm of finding the image of any given
diagram in $S(\gl_\infty)^{\otimes n}$ immediately in terms of necklaces.
For a given coloured $\B$-diagram, take the alternating sum over all
resolutions of the triple points. For each resolution convert the obtained
picture into a collection of oriented closed curves, put the numbers
(1,\dots,$n$) of the univalent vertices on the places where they were
before closing and thus get a product of necklaces. For instance:
$$
  \rb{-9mm}{\ig[height=19mm]{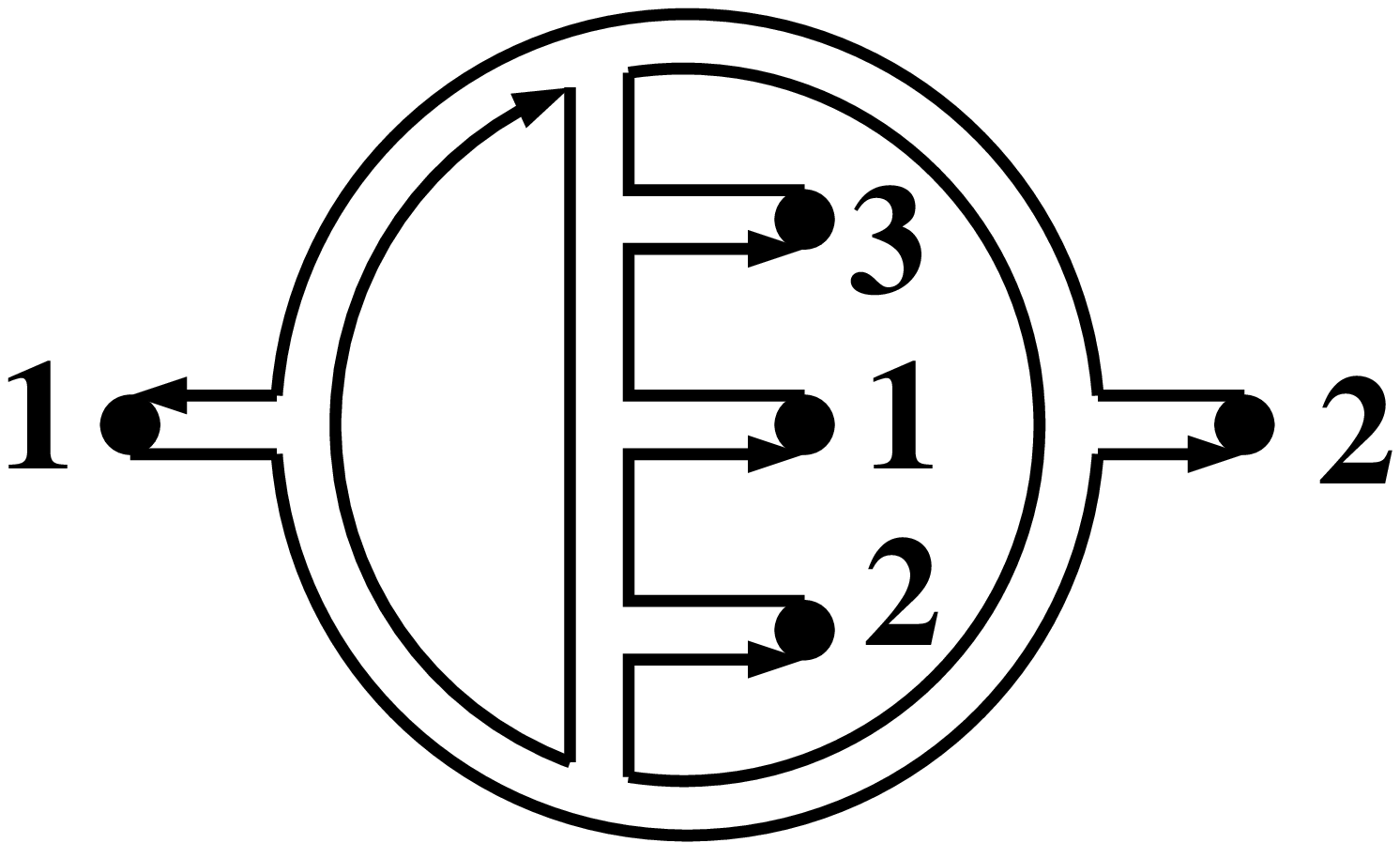}}
  \ \longmapsto
  \ x x_{12} x_{132}.
$$

\smallskip

\begin{xxca}
Find a direct proof, without appealing to Lie algebras,
that the described mapping $\rho$ into the necklace algebra provides
a weight system, that is, satisfies the AS and IHX relations.
{\sl Hint}: see \cite{Da3}, where this is done for the case $n=1$.
\end{xxca}

One application of (unicoloured) necklace weight system is the lower bound on
the dimensions of the spaces $\V_n$ for knots, see Section
\ref{lower_bound}.

Another application --- of the 2-coloured necklace weight system ---
is the proof that there exists a degree 7 Vassiliev invariant
that is capable to detect the change of orientation in
two-component string links, see \cite{DK}. This fact follows from the
computation
$$
\rho\Bigl(\,\rb{-9mm}{\ig[height=20mm]{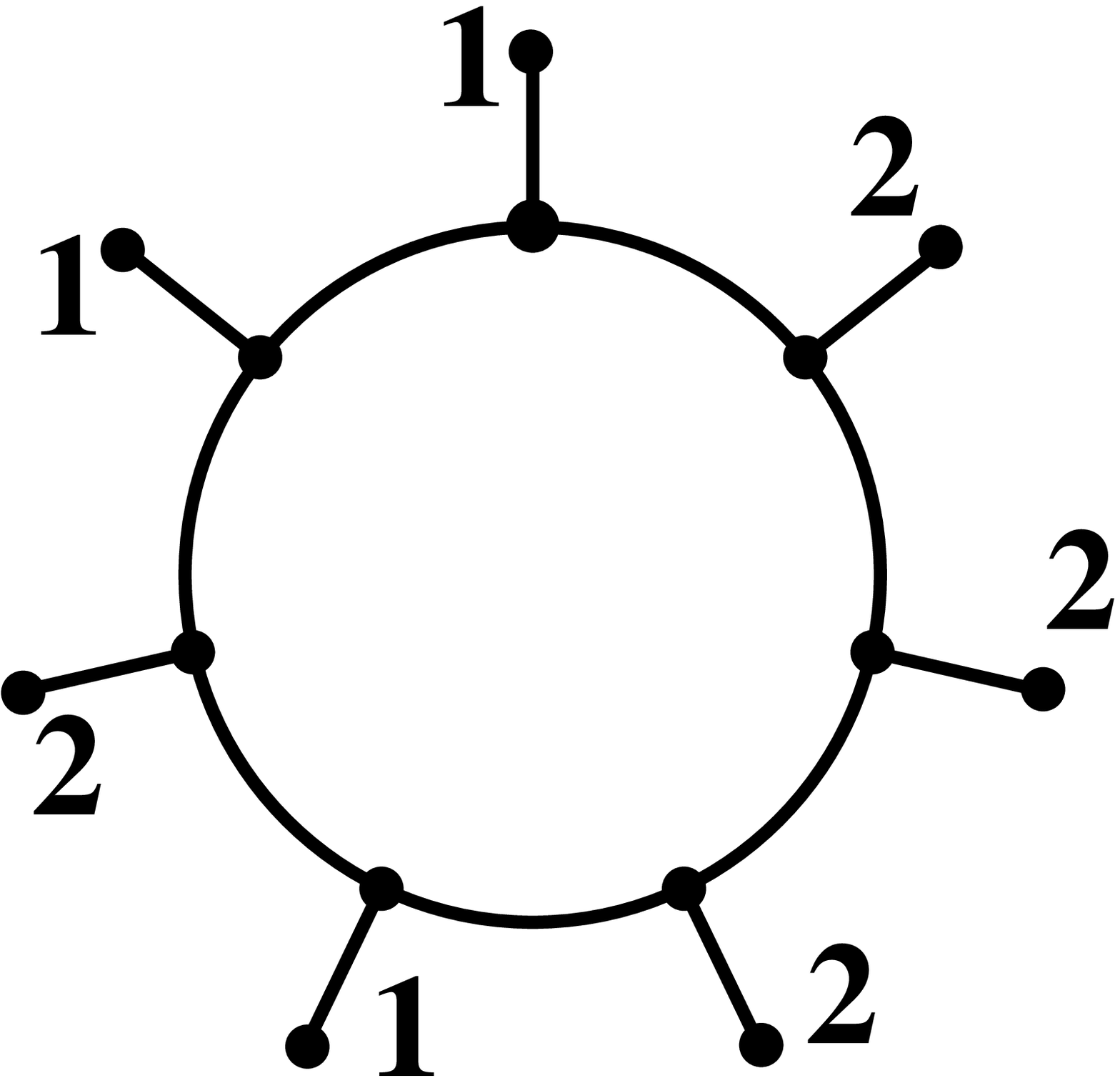}}\Bigr)
\ =
\ x(x_{1121222}-x_{1122212})+3x_2(x_{112212}-x_{112122})
$$
which implies that the depicted diagram is non-zero in $\B(2)$.

\section{Lie superalgebra weight systems}

The construction of Lie algebra weight systems works for algebraic
structures more general than Lie algebras \cite{Vai3, FKV, HV},
namely for the analogs of metrized Lie algebras in categories more
general than the category of vector spaces.  An example of such a
category is that of super vector spaces; Lie algebras in this
category are called {\em Lie superalgebras}.\index{Lie superalgebra}
The definition and basic properties of Lie superalgebras are
discussed in Appendix~\ref{lie_superalgebras}; we refer the reader
to \cite{Kac1, Kac2} for more details.

\subsection{Weight systems for Lie superalgebras}
Recall the construction of the Lie algebra weight systems for the
closed diagrams as described in
Sections~\ref{lie_alg_ws_without_bases} and \ref{LAWS_C}. It
consists of several steps. First, the internal graph of the diagram
is cut into tripods and chords. Then to each tripod we assign a
tensor in $\g^{\ot 3}$ coming from the Lie bracket, and to each
chord -- a tensor in $\g^{\ot 2}$ coming from the invariant form.
Next, we take the tensor product  $\widetilde{T}_\g$ of all these
tensors and perform contractions on the pairs of indices
corresponding to the points where the diagram was cut. Finally, we
re-arrange the factors in the tensor product and this gives the
tensor $T_\g$ whose image in $U(\g)$ is the weight system we were
after.

If $\g$ is a metrized Lie superalgebra, the very same construction
works with only one modification: re-arranging the factors in the
final step should be done with certain care. Instead of simply
permuting the factors in the tensor product one should use a certain
representation of the symmetric group $S_m$ on $m$ letters that
exists on the $m$th tensor power of any super vector space.

This representation is defined as follows. Let
$$S:\g\ot\g\to\g\ot\g$$ be the linear map that sends $u\ot v$ to
$(-1)^{p(u)p(v)} v\ot u$, where  $u,v$ are homogeneous (that is,
purely even or purely odd) elements of $\g$ and $p(x)$ stands for
the parity of $x$. The map $S$ is an involution; in other words,
it defines a representation of the symmetric group $S_2$ on the
vector space $\g^{\ot 2}$. More generally, the representation of
$S_m$ on $\g^{\ot m}$ is defined by sending the transposition
$(i,i+1)$ to ${\rm id}^{\ot i-1}\ot S \ot {\rm id}^{m-i-1}$. If
the odd part of $\g$ is zero, this representation simply permutes
the factors in the tensor product.

We shall use the same notation $\f_\g$ for the resulting weight
system.

\noindent{\bf Example.} Let $\g$ be a metrized Lie superalgebra with
the orthonormal bases $e_1, \ldots, e_m$ and $f_1,\ldots, f_r$ for
the even and the odd parts, respectively. Denote by $D$ the
diagram\quad $\risS{-12}{cd22Indexed}{}{25}{17}{12}$\ . Then
$$
   \f_\g(D)=\sum_{i=1}^m\sum_{j=1}^m  e_i e_j e_i e_j -
    \sum_{i=1}^m\sum_{j=1}^r (e_i f_j e_i f_j + f_j e_i f_j e_i) +
    \sum_{i=1}^r\sum_{j=1}^r  f_i f_j f_i f_j.
$$

\noindent{\bf Exercise.} Write down the expression for
$\f_\g\Bigl(\risS{-10}{cdIndexed1}{}{27}{0}{12}\Bigr)$.
This exercise is resolved in \cite{FKV} (Example~2 in Section~1.3)
though with a different base point and a not necessarily
orthonormal basis.

\noindent{\bf Exercise.} Show that $\f_\g$ is a well-defined weight
system with values in the (super) centre of $U(\g)$. In particular,
prove that $\f_\g$ satisfies the 4T relation.

\subsection{The $\gL_1_1$ weight system}
The simplest non-trivial example of a Lie superalgebra is the space
$\gL_1_1$ of endomorphisms of the super vector space of dimension
$1+1$. The universal weight system for $\gL_1_1$ can be calculated
with the help of a recursive formula similar to the formula for
$\sL_2$ (see Section~\ref{ws_sl_2_on_A}).

The (super) centre of $U(\gL_1_1)$ is a polynomial algebra in two
generators $c$ and $h$, where $c$ is the quadratic Casimir element
and $h\in\gL_1_1$ is the identity matrix.
\begin{xtheorem}[\cite{FKV}]
Let $\f_{\gL_1_1}$ be the weight system associated with $\gL_1_1$
with the invariant form $\bilinf{x}{y} = {\rm sTr}(xy)$.

Take a chord diagram $D$ and choose a chord $a$ of $D$. Then
\begin{multline*}
\f_{\gL_1_1}(D) = c\f_{\gL_1_1}(D_a) +h^2 \sum_{1\leq i\leq k}
\f_{\gL_1_1}(D_{i})\\
   -h^2 \sum_{1 \leq i<j \leq k} \left(\f_{\gL_1_1}(D^{+-}_{i,j}) +
\f_{\gL_1_1}(D^{-+}_{i,j})
   - \f_{\gL_1_1}(D^{l}_{i,j}) - \f_{\gL_1_1}(D^{r}_{i,j}) \right),
\end{multline*}
where:\\
$\bullet$ $k$ is the number of chords that intersect the  chord $a$;\\
$\bullet$ $D_a$ is the chord diagram obtained from $D$ by
deleting the chord $a$;\\
$\bullet$ for each chord $a_i$ that intersects $a$, the diagram
$D_i$ is obtained from $D$ by deleting the chords $a$ and $a_i$;\\
$\bullet$ $D^{+-}_{i,j}$, $D^{-+}_{i,j}$, $D^{l}_{i,j}$ and
$D^{r}_{i,j}$ are the chord diagrams obtained from $D_a$ in the
following way. Draw the diagram $D$ so that the chord $a$ is
vertical. Consider an arbitrary pair of chords $a_i$ and $a_j$
different from $a$ and such that each of them intersects $a$. Denote
by $p_i$ and $p_j$ the endpoints of $a_i$ and $a_j$ that lie to the
left of $a$ and by $p^*_i, p^*_j$ the endpoints of $a_i$ and $a_j$
that lie to the right. Delete from $D$ the chords $a$, $a_i$ and
$a_j$ and insert one new chord: $(p_i,p^*_j)$ for $D^{+-}_{i,j}$,
$(p_j,p^*_i)$ for $D^{-+}_{i,j}$, $(p_i,p_j)$ for $D^{l}_{i,j}$ and
$(p^*_i,p^*_j)$ for $D^{r}_{i,j}$:

$$D=\chvadth{chvard01}{
             \put(8,28){\mbox{$\scriptstyle a_i$}}
             \put(8,7){\mbox{$\scriptstyle a_j$}}
             \put(19,17){\mbox{$\scriptstyle a$}} };\quad\!
  D_a=\chvadth{chvard02}{};\quad\!
  D_i=\chvadth{cdvai1}{};\quad\!
$$
$$D^{+-}_{i,j}=\risS{-15}{cdvai2}{}{35}{30}{20}\ ;\quad\!
  D^{-+}_{i,j}=\risS{-15}{cdvai3}{}{35}{0}{0}\ ;\quad\!
  D^{l}_{i,j}=\risS{-15}{cdvai4}{}{35}{0}{0}\ ;\quad\!
  D^{r}_{i,j}=\risS{-15}{cdvai5}{}{35}{0}{0}\ .
$$

In particular, $\f_{\gL_1_1}(D)$ is a polynomial in $c$ and $h^2$.
\end{xtheorem}

We refer to \cite{FKV} for the proof.

\subsection{Invariants not coming from Lie algebras}
Lie algebra weight systems produce infinite series of examples of
Vassiliev invariants. J.\ Kneissler has shown in \cite{Kn0} that all
invariants up to order 12 come from Lie algebras. However,  in
general, this is not the case. P.\ Vogel \cite{Vo1} has used the
family of Lie superalgebras $D(1,2,\alpha)$ depending on the
parameter $\alpha$; he showed that these algebras produce invariants
which cannot be expressed as combinations of invariants coming from
Lie algebras. (J.\ Lieberum \cite{Lieb} gave an example of an order
17 closed diagram detected by $D(1,2,\alpha)$ but not by semisimple
Lie algebra weight systems.) Moreover, there exist Vassiliev
invariants that do no come from Lie (super) algebras \cite{Vo1,
Lieb}.

The main technical tool for proving these results is the algebra
$\Lambda$ constructed by Vogel. In the next chapter we shall
consider the {\em algebra of 3-graphs} closely related to Vogel's
algebra $\Lambda$.

\begin{xcb}{Exercises}

\begin{enumerate}
\item Let $(\g_1,\bilinf{\cdot}{\cdot}_1)$ and $(\g_2,\bilinf{\cdot}{\cdot}_2)$
be two metrized Lie algebras. Then their direct sum $\g_1\op \g_2$ is also a metrized
Lie algebra with respect to  the form
$\bilinf{\cdot}{\cdot}_1 \op \bilinf{\cdot}{\cdot}_2$. Prove that
$\f_{\g_1\op \g_2} = \f_{\g_1}\cdot \f_{\g_2}$.

\medskip
{\setlength{\baselineskip}{12pt} 
\scriptsize The general aim of exercises (\ref{ex1_sl2})-(\ref{ex6_sl2})
is  to compare the behaviour of
$\f_{\sL_2}(D)$ with that of
the chromatic polynomial of a graph. In these exercises we use the form
$2\bilinf{\cdot}{\cdot}$ as the invariant form.}

\item\label{ex1_sl2} ({\it S.~Chmutov, S.Lando \cite{CL}}).
 Prove that $\f_{\sL_2}(D)$ depends only on the intersection graph
$\Gamma(D)$ of the chord diagram $D$.

\item
Prove that the polynomial $\f_{\sL_2}(D)$ has alternating coefficients.

\item Show that for any chord diagram $D$ the polynomial
$\f_{\sL_2}(D)$ is divisible by $c$.

\item\neresh
Prove that the sequence of coefficients of the polynomial
$\f_{\sL_2}(D)$ is unimodal (that is, its absolute values form a
sequence with only one maximum).

\item
Let $D$ be a chord diagram with $n$ chords for which $\Gamma(D)$ is a tree.
Prove that $\f_{\sL_2}(D) = c(c-2)^{n-1}$.

\item
Prove that the highest three terms of the polynomial $\f_{\sL_2}(D)$
are
$$\f_{\sL_2}(D) = c^n - e\cdot c^{n-1} +
 ( e(e-1)/2 - t + 2q )\cdot c^{n-2} - \dots\ ,
$$
where $e$ is the number of intersections of chords of $D$; $t$ and
$q$ are the numbers of {\em triangles} and {\em quadrangles} of $D$
respectively. A triangle is a subset of three chords of $D$ with all
pairwise intersections. A quadrangle of $D$ is an unordered subset
of four chords
 $a_1, a_2, a_3, a_4$
which form a cycle of length four. This means that, after a suitable
relabeling, $a_1$ intersects $a_2$ and $a_4$, $a_2$ intersects $a_3$
and $a_1$, $a_3$ intersects $a_4$ and $a_2$, $a_4$ intersects $a_1$
and $a_3$ and any other intersections are allowed. For example,
$$e\Bigl( \chvad{chvard33} \Bigr) = 6,\qquad
  t\Bigl( \chvad{chvard33} \Bigr) = 4,\qquad
  q\Bigl( \chvad{chvard33} \Bigr) = 1\ .
$$

\item\label{ex6_sl2} ({\it A.~Vaintrob \cite{Vai2}}).
Define {\it vertex multiplication} of chord diagrams as follows:
\vspace{-10pt}
$$ \chvad{chvard34} \vee \chvad{chvard35}\quad :=\quad
   \risS{-11}{chvard36}{}{55}{20}{20}\quad =\quad \chvad{chvard37}\ .
$$
Of course, the result depends of the choice of vertices where multiplication is
performed. Prove that for any choice
$$\f_{\sL_2}(D_1\vee D_2) =
    \frac{\f_{\sL_2}(D_1)\cdot \f_{\sL_2}(D_2)}{c}\,.$$

\item
({\it D.~Bar-Natan, S.~Garoufalidis \cite{BNG}}) \label{Weight
system!of the Conway coefficients} Let $c_n$ be the coefficient of
$t^n$ in the Conway polynomial and $D$ a chord diagram of degree
$n$. Prove that $\symb(c_n)(D)$ is equal, modulo 2, to the
determinant of the adjacency matrix for the intersection graph
$\Gamma(D)$.

\item\label{ex_gl_N_circle}
\parbox[t]{3in}{Let $D_n$ be the chord diagram with $n$ chords whose
intersection graph is a circle, $n\geqslant3$. Prove that
$\f^{St}_{\gl_N}(D_n) = \f^{St}_{\gl_N}(D_{n-2})$. Deduce that
$\f^{St}_{\gl_N}(D_n)= N^2$ for odd $n$ and
$\f^{St}_{\gl_N}(D_n)= N^3$ for even $n$.
}\qquad
$\rb{-20pt}{$D_n =$}\quad \risS{-38}{wssoNDn}{
         \put(8,-10){\mbox{\scriptsize $n$ chords}}}{50}{20}{15}$

\item\label{ex_ws_sl_N_basis}
Work out a proof of the theorem from Section~\ref{ws_sl_n_St} about
the $\sL_N$ weight system with standard representation, similar to the
one given in Section~\ref{ws_glN_on_A}. Use the basis of the vector
space $\sL_N$ consisting of the matrices $e_{ij}$ for $i\not= j$ and
the matrices $e_{ii} - e_{i+1,i+1}$.

\item\label{ex_defr_sl_N}
Prove that $\f'^{St}_{\sL_N}\equiv \f'^{St}_{\gl_N}$.

{\sl Hint.}\hspace{2cm} $\f'^{St}_{\sL_N} =
e^{-\frac{N^2-1}{N}\bo_1}\cdot \f^{St}_{\sL_N}
     = e^{-N\bo_1}\cdot \f^{St}_{\gl_N} =\f'^{St}_{\gl_N}$.

\item\label{ex_ws_j_N_and_sl_2}
\index{Weight system!of the Jones coefficients} Compare the symbol
of the coefficient $j_n$ of the Jones polynomial
(Section~\ref{symb_j_n}) with the weight system coming form $\sL_2$,
and prove that
$$\symb(j_n)= \frac{(-1)^n}{2} \f'^{St}_{\sL_2}\ .$$

{\sl Hint.} Compare the formula for $\f'^{St}_{\sL_2}$ from the
previous problem and the formula for $\symb(j_n)$ from
Section~\ref{symb_j_n}, and prove that
$$(|s|-1) \equiv \mbox{\#\{chords $c$ such that $s(c)=1$\}} \mod{2}\ .$$

\item\label{ex_so_N_basis}
Work out a proof of the theorem from Section~\ref{ws_so_n_St} about
the $\so_N$ weight system in standard representation. Use the basis of
$\so_N$ formed by matrices $e_{ij}-e_{ji}$ for $i<j$. (In case of
difficulty consult \cite{BN0,BN1}.)

\item\label{ex_so_N_intgr}
Work out a proof, similar to the proof of the Proposition from
Section~\ref{ws_glN_on_A}, that $\f^{St}_{\so_N}(D)$ depends only on
the intersection graph of $D$.

\item\label{problem_Mellor}
({\it B.~Mellor \cite{Mel2}}). For any subset $J\subseteq [D]$, let
$M_J$ denote the {\em marked adjacency matrix} of the intersection
graph of $D$ over the filed $\Fi_2$ , that is the adjacency matrix
$M$ with each diagonal element corresponding to an element of $J$
replaced by 1. Prove that
$$\f^{St}_{\so_N}(D) = \frac{N^{n+1}}{2^n} \sum\limits_{J\subseteq [D]}
                (-1)^{|J|}N^{-\mbox{rank}(M_J)}\ ,$$
where the rank is computed as the rank of a matrix over $\Fi_2$.
This gives another proof of the fact that $\f^{St}_{\so_N}(D)$
depends only on the intersection graph $\Gamma(D)$.

\item\label{ex_so_N_roots}
Show that $N=0$ and $N=1$ are roots of the polynomial
$\f^{St}_{\so_N}(D)$ for any chord diagram $D$.

\item\label{ex_so_N_tree}
Let $D$ be a chord diagram with $n$ chords, such that the intersection graph
$\Gamma(D)$ is a tree. Show that $\f^{St}_{\so_N}(D) = \frac{1}{2^n}N(N-1)$.

\item\label{ex_so_N_circle}
Let $D_n$ be the chord diagram from Exercise~\ref{ex_gl_N_circle}.
Prove that
\begin{enumerate}
\item[(a)] $\f^{St}_{\so_N}(D_n)=\frac{1}{2}\left(
        \f^{St}_{\so_N}(D_{n-2}) - \f^{St}_{\so_N}(D_{n-1})\right)$;\vspace{8pt}

\item[(b)] $ \f^{St}_{\so_N}(D_n)=
              \frac{1}{(-2)^n} N(N-1)(a_{n-1}N-a_n),$
where the recurrent sequence $a_n$ is defined by $a_0=0, a_1=1,
a_n=a_{n-1}+2a_{n-2}$.
\end{enumerate}

\item\label{ex_sl2_on_cincub_cubinb}
Compute the values of $\f_{\sL_2}$ on the closed diagrams
$\risS{-12}{fdbw4_ex_sl2}{}{29}{20}{15}$ and
$\risS{-12}{fdw4b_ex_sl2}{}{30}{20}{15}$, and show that these two
diagrams are linearly independent.

{\sl Answer}: $16c^2$, $64c$.

\item\label{ex_sl_2_on_C_tree}
\parbox[t]{3in}{Let $\ol{t_n}\in\F_{n+1}$ be a closed diagram with $n$ legs as
shown in the figure.\\ Show that $\f_{\sL_2}(\ol{t_n}) = 2^n c$.
}\qquad $\rb{-3pt}{$\ol{t_n} :=$}\quad \risS{-17}{wssl2tn}{
         \put(5,-8){\mbox{\scriptsize $n$ legs}}}{30}{20}{15}$

\item\label{ex_sl_2_on_C_wheel}\index{Wheel!in $\F$}
\parbox[t]{3in}{Let $\ol{w_n}\in\F_n$ be a wheel with $n$ spokes.\\
Show that\\
$\f_{\sL_2}(\ol{w_2}) = 4c$,\quad $\f_{\sL_2}(\ol{w_3}) = 4c$,\ and
}\qquad $\rb{-3pt}{$\ol{w_n} :=$}\quad \risS{-17}{wssl2wn}{
         \put(0,-8){\mbox{\scriptsize $n$
spokes}}}{30}{20}{15}$\vspace{5pt}\\
$\f_{\sL_2}(\ol{w_n}) = 2c\cdot\f_{\sL_2}(\ol{w_{n-2}}) +
   2\f_{\sL_2}(\ol{w_{n-1}}) - 2^{n-1}c$.

\item\label{ex_sl_2_on_B_wheel}
\parbox[t]{3in}{Let $w_{2n}\in\B_{2n}$ be a wheel with $2n$ spokes and
$(\risS{1.5}{strut}{}{15}{10}{0})^n \in\B_n$ be the $n$th power of
the element $\risS{1.5}{strut}{}{15}{0}{0}$ in the algebra $\B$.\\
Show that for the tensor $T_{\sL_2}$ as in \ref{ws_sl_2_on_C} the }
\parbox[t]{2in}{\qquad $\rb{7pt}{$w_{2n} :=$}\quad \risS{-3}{wssl2wn-b}{
         \put(-3,-8){\mbox{\scriptsize $2n$ spokes}}}{30}{20}{15}$

\noindent
$\rb{-3pt}{$(\risS{1.5}{strut}{}{15}{0}{0})^n :=$}
             \rb{-2pt}{$\left.\begin{picture}(15,12)(0,0)
                      \put(0,-5){\ig[width=15pt]{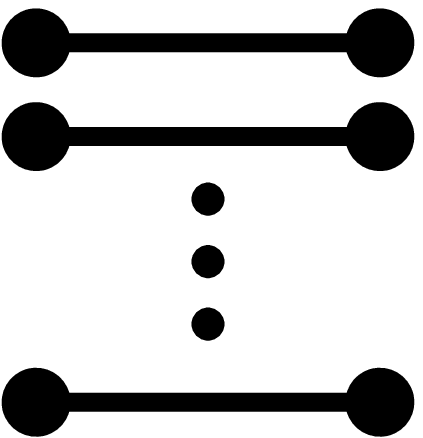}}
                  \end{picture}\right\}$}
 \rb{-1pt}{\scriptsize $n$ segments}$
               }\\
following equality holds: $T_{\sL_2}(w_{2n}) = 2^{n+1}\
T_{\sL_2}((\risS{1.5}{strut}{}{15}{0}{0})^n)$. Therefore,
$\rho_{\sL_2}(w_{2n}) = 2^{n+1}\
\rho_{\sL_2}((\risS{1.5}{strut}{}{15}{0}{10})^n).$

\item\label{ex_sl_2_on_Pkn}
Let $p\in \PR^k_n\subset\F_n$ be a primitive element of degree $n>1$
with at most $k$ external vertices. Show that $\f_{\sL_2}(p)$ is a
polynomial in $c$ of degree $\leqslant k/2$.

{\sl Hint.} Use the theorem from \ref{ws_sl_2_on_C} and the
calculation of $\f_{\sL_2}(\ol{t_3})$ from Exercise
(\ref{ex_sl_2_on_C_tree}).

\item\label{ex_sl_2unfr}
Let $\f'_{\sL_2}$ be the deframing of the weight system $\f_{\sL_2}$
according to the procedure of Section \ref{defram_ws}. Show that for
any element $D\in \A_n$,  the value $\f'_{\sL_2}(D)$ is a polynomial
in $c$ of degree $\leqslant [n/2]$.

{\sl Hint.} Use the previous exercise, Exercise \ref{ex_defr_prim}
of Chapter~\ref{ChD}, and Section~\ref{filtr_pr}.

\item\label{ex_sl_2_col_Jones}
Denote by $V_k$ the $(k+1)$-dimensional irreducible
representation of $\sL_2$ (see Appendix \ref{irr_sl2}).
Let $\f'^{V_k}_{\sL_2}$ be the
corresponding weight system. Show that for any
element $D\in \A_n$ of degree $n$, $\f'^{V_k}_{\sL_2}(D)/k$ is a
polynomial in $k$ of degree at most $n$.

{\sl Hint.} The Casimir number (see
page~\pageref{casimir_num}) in this case is
$\frac{k^2-1}{2}$.

\item\label{ex_gl_N_on_C_sl_N}
Let $D\in \F_n$ ($n>1$) be a connected closed diagram. Prove that
$\f^{St}_{\gl_N}(D) = \f^{St}_{\sL_N}(D)$.

{\sl Hint.} For the Lie algebra $\gl_N$ the tensor $J\in \gl_N^{\ot
3}$ lies in the subspace $\sL_N^{\ot 3}$.

\item\label{ex_gl_N_on_C_surface}
Consider a closed diagram $D\in \F_n$ and a $\gl_N$-state $s$ for it
(see page~\pageref{gl_N-state}). Construct a surface $\Sigma_s(D)$
by attaching a disk to the Wilson loop, replacing each edge by a
narrow band and glueing the bands together at the trivalent vertices
with a twist if $s=-1$, and without a twist if $s=1$. Here is an
example:
$$\qquad D =\ \risS{-24}{trtrvert}{
         \put(9,15){\mbox{\scriptsize $+$}}
         \put(9,36){\mbox{\scriptsize $-$}}
         \put(27,27){\mbox{\scriptsize $-$}}
         \put(36,20){\mbox{\scriptsize $+$}}
       }{50}{27}{22} \quad
\risS{-2}{totonew}{}{25}{10}{10}\quad
\risS{-24}{cdtr1}{}{50}{27}{20} \quad
\risS{-2}{totonew}{}{25}{10}{10}\quad
\risS{-27}{sdtr}{}{60}{30}{25}\ =: \Sigma_s(D)\ .
$$
\begin{enumerate}
\item[(a)] Show that the surface $\Sigma_s(D)$ is orientable.
\item[(b)] Compute the Euler characteristic of $\Sigma_s(D)$ in terms of $D$,
and show that it depends only on the degree $n$ of $D$.
\item[(c)] Prove that $\f^{St}_{\gl_N}(D)$ is an odd polynomial for even $n$,
   and it is an even polynomial for odd $n$.
\end{enumerate}

\item\label{ex_gl_N_on_C_root}
Show that $N=0$, $N=-1$, and $N=1$ are roots of the polynomial
$\f^{St}_{\gl_N}(D)$ for any closed diagram $D\in \F_n$ ($n>1$).

\item\label{ex_gl_N_on_C_tree}
Compute $\f^{St}_{\gl_N}(\ol{t_n})$, where $\ol{t_n}$ is the closed
diagram from Exercise~\ref{ex_sl_2_on_C_tree}.

{\sl Answer.} For $n\geq 1$, $\f^{St}_{\gl_N}(\ol{t_n})=
N^n(N^2-1)$.

\item\label{ex_gl_N_on_C_wheel}
For the closed diagram $\ol{w_n}$ as in
Exercise~\ref{ex_sl_2_on_C_wheel}, prove that
$\f^{St}_{\gl_N}(\ol{w_n}) = N^2(N^{n-1}-1)$ for odd $n$, and
$\f^{St}_{\gl_N}(\ol{w_n}) = N(N^n +N^2-2)$ for even $n$.

{\sl Hint.} Prove the recurrent formula $\f^{St}_{\gl_N}(\ol{w_n})=
N^{n-1}(N^2-1) +
 \f^{St}_{\gl_N}(\ol{w_{n-2}})$ for $n\geq 3$.

\item\label{ex_conw_on_C}
\index{Weight system!of the Conway coefficients} Extend the
definition of the weight system $\symb(c_n)$ of the coefficient
$c_n$ of the Conway polynomial to $\F_n$, and prove that
$$\symb(c_n)(D) = \sum_s\ \Bigl(\prod_v s(v)\Bigr) \delta_{1,|s|}\ ,$$
where the states $s$ are precisely the same as in the theorem of
Section~\ref{ws_glN_on_C} for the weight system $\f^{St}_{\gl_N}$.
In other words, prove that $\symb(c_n)(D)$ is equal to the coefficient
of $N$ in the polynomial $\f^{St}_{\gl_N}(D)$. In particular, show
that $\symb(c_n)(\ol{w_n}) = -2$ for even $n$, and
$\symb(c_n)(\ol{w_n}) = 0$ for odd $n$.

\item\label{ex_so_N_on_C_root2}
\begin{enumerate}
\item[(a)] Let $D\in \F$ be a closed diagram with at least one internal trivalent vertex.
  Prove that $N=2$ is a root of the polynomial $\f^{St}_{\so_N}(D)$.
\item[(b)] Deduce that $\f^{St}_{\so_2}(D) = 0$ for any primitive
  closed diagram $D$.
\end{enumerate}

{\sl Hint.} Consider the eight states that differ only on three
edges meeting at an internal vertex (see
page~\pageref{8_states_so_N}). Show that the sum over these eight
states, $\sum\sign(s)2^{|s|}$, equals zero.

\item\label{ex_so_N_on_C_tree}
Prove that $\f^{St}_{\so_N}(\ol{t_n}) =
\frac{N-2}{2}\f^{St}_{\so_N}(\ol{t_{n-1}})$
for $n>1$, where $\ol{t_n}$ is as in Exercise~\ref{ex_sl_2_on_C_tree}. \\
In particular, $\f^{St}_{\so_N}(\ol{t_n}) =
\frac{(N-2)^n}{2^{n+1}}N(N-1)$.

\item
Using some bases in $\F_2$ and $\B_2$, find the matrix of the
isomorphism $\chi$, then calculate (express as polynomials in the
standard generators) the values on the basis elements of the weight
systems $\f_\g$ and $\rho_\g$ for the Lie algebras $\g=\so_3$ and
$\g=\gl_N$ and check the validity of the relation
$\beta\circ\rho=\f\circ\chi$ in this particular case.

\item
Prove that the map $\rho:\B\to S(\g)$ is well-defined.

\end{enumerate}
\end{xcb}
 %6 3-graphs
\chapter{Algebra of 3-graphs} %13
\label{alg3g}

The {\em algebra of 3-graphs} $\G$, introduced in \cite{DKC}, is
related to the diagram algebras $\F$ and $\B$. The difference
between 3-graphs and closed diagrams is that 3-graphs do not have a
distinguished cycle (Wilson loop); neither they have univalent
vertices, which distinguishes them from open diagrams. Strictly
speaking, there are two different algebra structures on the space of
3-graphs, given by the edge (Section~\ref{ed_prod_in_G}) and the
vertex (Section~\ref{prod_in_G}) products. The space $\G$ is closely
related to the Vassiliev invariants in several ways:

\begin{itemize}
\item
The vector space $\G$ is isomorphic to the subspace $\PR^2$ of the
primitive space $\PR\subset\F$  spanned by the connected diagrams
with 2 legs (Section~\ref{iso_P2_and_G}).
\item
The algebra $\G$ acts on the primitive space $\PR$ in two ways, via
the edge, and via the vertex products (see Sections~\ref{G_ed_act_on
P} and \ref{G_v_act_on_P}). These actions behave nicely with respect
to Lie algebra weight systems (see Chapter~\ref{LAWS}); as a
consequence, the algebra $\G$ is as good a tool for the proof of
existence of non-Lie-algebraic weight systems as the algebra
$\Lambda$ in Vogel's original approach (Section~\ref{ws_not_Lie}).
\item
The vector space
$\G$ describes the combinatorics of finite type
invariants of integral homology 3-spheres in much the same way as
the space of chord diagrams describes the combinatorics of Vassiliev
knot invariants. This topic, however, lies outside of the scope
of our book and we refer an interested reader to \cite{Oht1}.
\end{itemize}

Unlike $\F$ and $\B$, the algebra $\G$ does not have any natural
coproduct.

\section{The space of 3-graphs}

A {\em 3-graph} \index{Three-graph} is a connected 3-valent graph
with a fixed cyclic order of half-edges at each vertex. Two 3-graphs
are {\em isomorphic} if there exists a graph isomorphism between
them that preserves the cyclic order of half-edges at every vertex.
The {\em degree}, or {\em order} \index{Order}\index{Degree} of a
3-graph is defined as half the number of its vertices. It will be
convenient to consider a circle with no vertices on it as a 3-graph
of degree 0 (even though, strictly speaking, it is not a graph).

\begin{xexample}
Up to an isomorphism, there are three different 3-graphs
of degree 1:
$$\risS{-12}{wssoN1}{}{30}{10}{15}\qquad\qquad
  \risS{-8}{g3-11}{}{30}{10}{15}\qquad\qquad
  \risS{-3}{g3-12}{}{30}{10}{15}
$$
\end{xexample}

\begin{xremark} Graphs with a cyclic order of half-edges at each vertex
are often called {\em ribbon graphs} (see \cite{LZ}),\index{Ribbon
graph} as every such graph can be represented as an orientable
surface with boundary obtained by ``thickening'' the graph:
$$
  \rb{-4mm}{\ig[height=10mm]{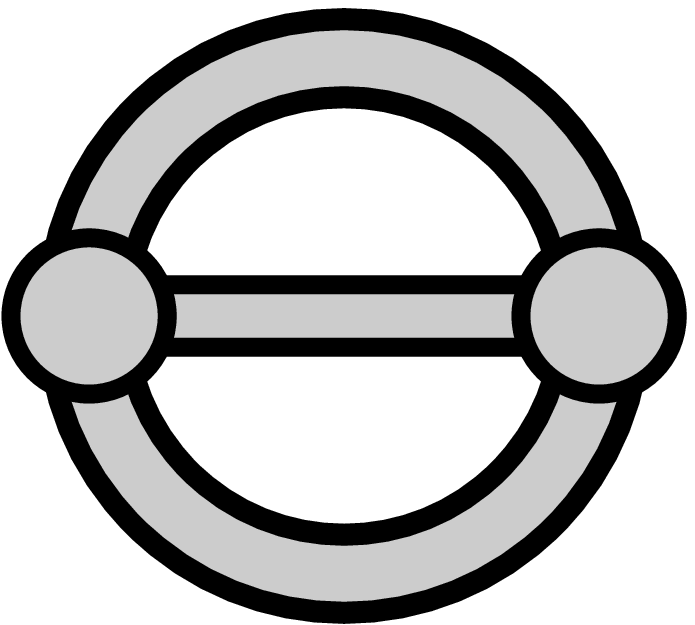}}\qquad\qquad
  \rb{-4mm}{\ig[height=10mm]{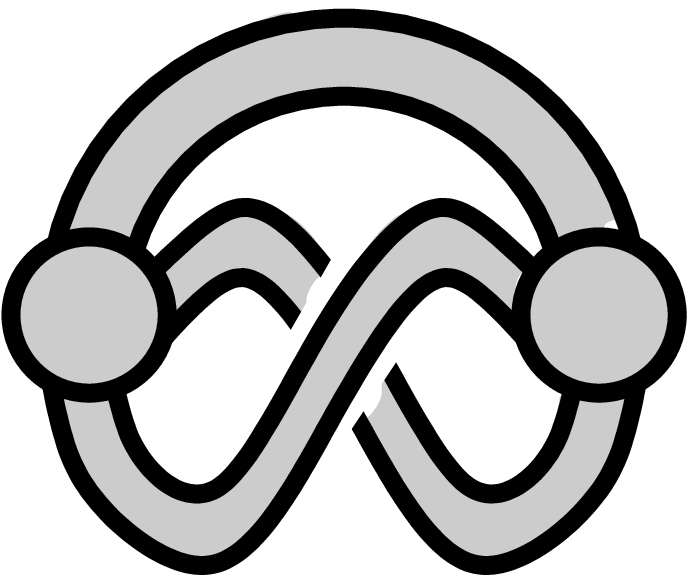}}\qquad\qquad
  \rb{-3mm}{\ig[height=8mm]{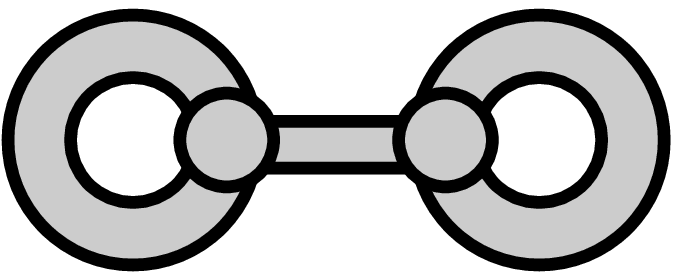}}
$$
To be more precise, given a graph, we replace each of its vertices
and each of its edges by an oriented disk (imagine that the disks
for the vertices are ``round'' while the disks for the edges are
``oblong''). The disks are glued together along segments of their
boundary in agreement with the orientation and with the prescribed
cyclic order at each vertex; the cyclic order at a vertex is taken
in the positive direction of the vertex-disk boundary.
\end{xremark}

\begin{xdefinition}\label{alg-3-graphs}
The {\em space of 3-graphs} $\G_n$ is the $\Q$-vector space spanned
by all 3-graphs of degree $n$ modulo the AS and IHX
relations (see page~\pageref{AS_rel}). 
\end{xdefinition}
In particular, the space $\G_0$ is one-dimensional and spanned by
the circle.

\noindent{\bf Exercise.} \mbox{ }\\
\parbox[t]{4in}{
Check that the 3-graph on the right is equal to zero as an element of the space
$\G_3$.}\qquad
\risS{-8}{g3-6gon}{}{30}{10}{15}

\section{Edge multiplication}
\label{ed_prod_in_G}

In the graded space
$$
  \G = \G_0 \op \G_1 \op \G_2 \op \G_3 \op \dots
$$\index{Algebra!of 3-graphs}
there is a natural structure of a commutative algebra.

Let $G_1$ and $G_2$ be two 3-graphs. Choose arbitrarily an edge in
$G_1$ and an edge in $G_2$. Cut each of these two edges in the
middle and re-connect them in any other way so as to get a 3-graph.
$$\risS{-30}{g1prod}{\put(20,51){\mbox{$G_1$}} \put(20,8){\mbox{$G_2$}}
           }{50}{35}{35} \qquad
 \risS{-4}{totonew}{}{25}{20}{15}\qquad
 \risS{-30}{g2prod}{\put(20,51){\mbox{$G_1$}} \put(20,8){\mbox{$G_2$}}
           }{50}{10}{15} \quad = G_1\cdot G_2\ .$$
The resulting 3-graph is called the {\em edge product} \index{Edge
product} \index{Product!in $\G$} of $G_1$ and $G_2$.

The edge product of 3-graphs can be thought of as the connected sum
of $G_1$ and $G_2$ along the chosen edges, or as the result of
insertion of one graph, say $G_1$, into an edge of $G_2$.

\begin{xremark} The product of two connected graphs may yield a
disconnected graph, for example:
$$ \risS{-17}{g3-13}{}{15}{20}{20}\ \times\
   \risS{-17}{g3-13}{}{15}{10}{15}\quad =\quad
   \risS{-17}{g3-14}{}{33}{10}{15}
$$
This happens, however, only in the case when each of the
two graphs becomes disconnected after cutting the chosen edge,
and in this case both graphs are 0 modulo AS and IHX relations
(see Lemma \ref{lemma2b}(b) below).
\end{xremark}
\begin{theorem}
The edge product of 3-graphs, viewed as an element of the space
$\G$, is well-defined.
\end{theorem}

Note that, as soon as this assertion is proved, one immediately sees
that the edge product is commutative.

The claim that the product is well-defined consists of two parts.
Firstly, we need to prove that modulo the AS and the IHX relations
the product does not depend on the {\it choice\/} of the two edges
of $G_1$ and $G_2$ which are cut and re-connected. Secondly, we must
show that the product does not depend on the {\it way\/} they are
re-connected (clearly, the two loose ends of $G_1$ can be glued to
the two loose ends of $G_2$ in two different ways). These two facts
are established in the following two lemmas.

\begin{lemma}\label{lemma1a} Modulo the AS and the IHX relations,
a subgraph with two legs can be carried through a vertex: 
$$\risS{-20}{l13}{\put(12,28){\mbox{$G$}}}{50}{30}{15}\quad = \quad
  \risS{-20}{l14}{\put(31,28){\mbox{$G$}}}{50}{0}{0}.$$
\end{lemma}

\begin{proof}
This lemma is a particular case ($k=1$) of the Kirchhoff law (see
page \pageref{Kirch}).
\end{proof}

Lemma~\ref{lemma1a} shows that given an insertion of a 3-graph $G_1$
into an edge of $G_2$, there exists an equivalent insertion of $x$
into any adjacent edge. Since $G_2$ is connected, it only remains to
show that the two possible insertions of $G_1$ into an edge of $G_2$
give the same result.

\begin{lemma}\label{lemma2a}
The two different ways to re-connect two 3-graphs produce the same
element of the space $\G$:
$$\risS{-20}{l21}{\put(20,40){\mbox{$G_1$}} \put(20,7){\mbox{$G_2$}}
           }{50}{30}{20}\quad = \quad
  \risS{-20}{l22}{\put(20,40){\mbox{$G_1$}} \put(20,7){\mbox{$G_2$}}
           }{50}{0}{10}\ .
$$
\end{lemma}

\begin{proof}
At a vertex of $G_1$ which lies next to the subgraph $G_2$ in the
product, one can, by Lemma \ref{lemma1a}, perform the following
manoeuvres:
$$\risS{-20}{man1}{\put(10,20){\mbox{$\scriptstyle G_1$}}
                   \put(42,8){\rotatebox{180}{$\scriptstyle G_2$}}
           }{60}{30}{20}\ = \
  \risS{-20}{man2}{\put(10,20){\mbox{$\scriptstyle G_1$}}
                   \put(38,25){\rotatebox{90}{$\scriptstyle G_2$}}
           }{60}{30}{20}\ = \
   \risS{-25}{man3}{\put(10,25){\mbox{$\scriptstyle G_1$}}
                   \put(15,11){\mbox{$\scriptstyle G_2$}}
           }{60}{30}{20}\ = \
  \risS{-20}{man4}{\put(10,20){\mbox{$\scriptstyle G_1$}}
                   \put(43,4){\mbox{$\scriptstyle G_2$}}
           }{60}{30}{20}\ .
$$
Therefore,
$$\risS{-20}{l21}{\put(20,42){\mbox{$G_1$}} \put(20,8){\mbox{$G_2$}}
           }{50}{20}{20}\quad = \quad
  \risS{-20}{l21}{\put(20,42){\mbox{$G_1$}} \put(20,12){\rotatebox{180}{$G_2$}}
           }{50}{20}{20}\quad = \quad
  \risS{-20}{l22}{\put(20,42){\mbox{$G_1$}} \put(20,8){\mbox{$G_2$}}
           }{50}{0}{10}\ .
$$
The lemma is proved, and the edge multiplication of 3-graphs is thus
well-defined.
\end{proof}

The edge product of 3-graphs extends by linearity to the whole space
$\G$.

\begin{xcorollary}The edge product in $\G$ is well-defined and
associative.
\end{xcorollary}

This follows from the fact that a linear combination of either AS or
IHX type relations, when multiplied by an arbitrary graph, is a
combination of the same type. The associativity is obvious.

\subsection{Some identities}
\label{some_id}

There are two very natural operations defined on the space $\G$: the
{\em insertion of a bubble into an edge}:
$$\risS{1}{l42-b1}{}{45}{5}{10}\quad\risS{-2}{totor}{}{35}{0}{0}\quad
  \risS{-1}{l42-b2}{}{45}{0}{0}\ .
$$
and the {\em insertion of a triangle into a vertex}:
\index{Three-graph!bubble}
$$\risS{-12}{l42-tr}{}{45}{25}{20}\quad\risS{-2}{totor}{}{35}{0}{0}\quad
  \risS{-12}{l41}{}{45}{0}{0}\ .
$$
Inserting a bubble into an edge of a 3-graph is the same thing as
multiplying this graph by $\beta=\bub\ \in\G_1$. In particular, this
operation is well-defined and does not depend on the edge where the
bubble is created. Inserting a triangle into a vertex can be
expressed in a similar fashion via the {\em vertex multiplication}
discussed below in \ref{prod_in_G}. The following lemma implies that
inserting a triangle into a vertex is a well-defined operation:
\begin{lemma}\label{lemma1b}
A triangle is equal to one half of a bubble:
$${\displaystyle \pil{l41} \quad=\quad \frac12 \ \pil{l42}
\quad=\quad \frac12\ \pil{l43} \quad=\quad \frac12\
\rb{-3pt}{\pil{l44}}}.
$$
\end{lemma}

\begin{proof}
\begin{align*}
\pil{l41} \quad &=\quad \pil{l45}\ +\ \pil{l46}
\quad=\quad \pil{l42}\ -\ \pil{l47}\\
&= \quad \pil{l42}\ -\ \pil{l41}
\end{align*}\vspace*{-10pt}
\end{proof}

\begin{xremark}
It was proved by Pierre Vogel \cite{Vo2} that the operator of bubble
insertion has non-trivial kernel. He exhibited an element of degree
15 which is killed by inserting a bubble.
\end{xremark}

The second lemma describes two classes of 3-graphs which are equal
to 0 in the algebra $\G$, that is, modulo the AS and IHX relations.

\begin{lemma}\label{lemma2b}\mbox{ }\\
\parbox[c]{2.8in}{{\rm ({\bf a})} A graph with a loop is $0$ in $\G$.}\hspace{50pt}
$\risS{-5}{g3-loop}{}{50}{10}{10}\quad = \quad 0$ \\
\parbox{2.8in}{{\rm ({\bf b})} More generally, if the edge connectivity of the
graph $\gamma$ is $1$, that is, if it becomes disconnected after
removal of an edge, then $\gamma=0$ in $\G$.}\qquad $\gamma\ =\
\risS{-10}{l51}{}{70}{30}{15}\ = \ 0$
\end{lemma}

\begin{proof}
(a) A graph with a loop is zero because of the antisymmetry
relation. Indeed, changing the cyclic order at the vertex of the
loop produces a graph which is, on one hand, isomorphic to the
initial graph, and on the other hand, differs from it by a sign.

(b) Such a graph can be represented as a product of two graphs, one
of which is a graph with a loop that vanishes according to (a):
$$\gamma\ =\
    \risS{-10}{l51}{\put(4,10){\mbox{$\scriptstyle G_1$}} \put(60,10){\mbox{$\scriptstyle G_2$}}}{70}{15}{0}\ = \
             \risS{-10}{l53}{\put(4,10){\mbox{$\scriptstyle G_1$}}}{62}{15}{0}\times
             \risS{-10}{l52}{\put(11,10){\mbox{$\scriptstyle G_2$}}}{22}{15}{0} \ = \ 0\ .\vspace*{-10pt}
$$
\end{proof}

\subsection{The Zoo}

Table~\ref{gen-Gamma} shows the dimensions $d_n$ and displays the
bases of the vector spaces $\G_n$ for $n\leq 11$, obtained by
computer calculations.

\begin{table}[htb]
$$\begin{array}{c|c|l}
n&d_n&\mbox{additive generators} \\ \hline\hline
1&1&\hspace{-5pt}\pmris{b1}\\ \hline
2&1&\hspace{-5pt}\pmris{b2}\\ \hline
3&1&\hspace{-5pt}\pmris{b3}\\ \hline
4&2&\hspace{-5pt}\pmris{b4}\pmris{w4}\\ \hline
5&2&\hspace{-5pt}\pmris{b5}\pmris{w4b1}\\ \hline
6&3&\hspace{-5pt}\pmris{b6}\pmris{w4b2}\vmris{w6} \\ \hline
7&4&\hspace{-5pt}\pmris{b7}\vbmris{w4b3}\vmris{w6b1}\vmris{w7}\\ \hline
8&5&\hspace{-5pt}\pmris{b8}\bmris{w4b4}\pmris{w6b2}\pmris{w7b1}\vmris{w8}
              \\ \hline
9&6&\hspace{-5pt}\pmris{b9}\bmris{w4b5}\vmris{w6b3}\vmris{w7b2}
                 \vmris{w8b1}\vmris{w9}\\ \hline
10&8&\hspace{-5pt}\pmris{b10}\bmris{w4b6}\pmris{w6b4}\pmris{w7b3}
                  \pmris{w8b2}\pmris{w9b1}\vmris{w10}\pmris{dodek}\\ \hline
11&9&\hspace{-5pt}\pmris{b11}\pmris{w4b7}\vmris{w6b5}\vbmris{w7b4}\vmris{w8b3}
                  \vmris{w9b2}\vmris{w10b1}\pmris{dodb1}\vmris{w11}
\end{array}$$
\caption{Additive generators of the algebra of 3-graphs $\G$}
\label{gen-Gamma}
\index{Table of!generators of $\G$}
\end{table}

Note that the column for $d_n$ coincides with the column for $k=2$
in the table of primitive spaces on page \pageref{dim-of-prim-sp}.
This will be proved in Proposition~\ref{iso_P2_and_G}.

One can see from the table that the multiplicative generators of the
algebra $\G$ up to degree 11 can be chosen as follows (here $\beta$
stands for ``bubble'', $\om_i$
--- for ``wheels'', $\delta$ --- for ``dodecahedron''):
\index{Three-graph!wheel}\index{Three-graph!dodecahedron}
$$\begin{array}{c|c|c|c|c|c|c|c}
1&4&6&7&8&9&10&11 \\ \hline
\hspace{-4pt}\mris{b1}
&\hspace{-4pt}\mris{w4}
&\hspace{-4pt}\mris{w6}
&\hspace{-4pt}\mris{w7}
&\hspace{-4pt}\mris{w8}
&\hspace{-4pt}\mris{w9}
&\hspace{-4pt}\mris{w10}\hspace{4pt}\mris{dodek}
&\hspace{-4pt}\mris{w11}\\ \hline
\beta&\om_4&\om_6&\om_7&\om_8&\om_9&\om_{10}\quad\quad \delta&\om_{11}
\end{array}$$

The reader may have noticed that the table of additive generators
does not contain the elements $\om_4^2$ of degree 8 and $\om_4\om_6$
of degree 10. This is due to the following relations (found by
A.~Kaishev \cite{Kai}) in the algebra $\G$:
\begin{eqnarray*}
\om_4^2 &=& \frac{5}{384}\beta^8 - \frac{5}{12}\beta^4\om_4 +
       \frac{5}{2}\beta^2\om_6 - \frac32 \beta\om_7, \\
\om_4\om_6 &=& \frac{305}{27648}\beta^{10} - \frac{293}{864}\beta^6\om_4
+ \frac{145}{72}\beta^4\om_6 - \frac{31}{12}\beta^3\om_7 + 2\beta^2\om_8
- \frac34 \beta\om_9.
\end{eqnarray*}

In fact, as we shall see in Section~\ref{prod_in_G}, it is true in
general that the product of an arbitrary pair of homogeneous
elements of $\G$ of positive degree belongs to the ideal generated
by $\beta$.

Since there are non-trivial relations between the generators, the
algebra of 3-graphs, in contrast to the algebras $\A$, $\B$ and
$\F$, is commutative but not free commutative and, hence, does not
possess the structure of a Hopf algebra.

\section{Vertex multiplication}
\label{prod_in_G}

Apart from the edge product, the space
$$
  \G_{\geq 1} = \G_1 \op \G_2 \op \G_3 \op \G_4 \op \dots\
$$
spanned by all the 3-graphs of non-zero degree has another
commutative and associative product.

Let $G_1$ and $G_2$ be two 3-graphs of positive degree. Choose
arbitrarily a vertex in $G_1$ and a vertex in $G_2$. Cut out each of
these two vertices and attach the three loose ends that appear on
$G_1$ to the three loose ends on $G_2$. There are six possible ways
of doing this. Take the alternating average of all of them,
assigning the negative sign to those three cases where the cyclic
order on the loose ends of $G_1$ agrees with that for $G_2$, and the
positive sign to the other three cases. This alternating average is
called the {\em vertex product} \index{Vertex
product}\index{Product!in $\G$} of $G_1$ and $G_2$.

Pictorially, if the graphs $G_1$, $G_2$ are drawn as
$G_1=\risS{-9}{g1}{\put(2,16){\mbox{$\scriptstyle
G_1$}}}{15}{20}{10}$,\
$G_2=\risS{-11}{g2}{\put(2,5){\mbox{$\scriptstyle
G_2$}}}{15}{20}{12}$, then, in order to draw their vertex product we
have to merge them, inserting a permutation of the three strands in
the middle. Then we take the result with the sign of the permutation
and average it over all six permutations:

\def\gonegtwo#1{\risS{-18}{G1G2#1}{\put(2,5){\mbox{$\scriptstyle G_2$}}
              \put(2,33){\mbox{$\scriptstyle G_1$}}
                }{15}{20}{15}}

$$\risS{-9}{g1}{\put(2,16){\mbox{$\scriptstyle G_1$}}}{15}{20}{10}\ \vee\
  \risS{-11}{g2}{\put(2,5){\mbox{$\scriptstyle G_2$}}}{15}{20}{12}\quad =\quad
 \frac{1}{6}\Biggl[\ \ \gonegtwo{a}\ -\ \gonegtwo{b}\ -\ \gonegtwo{c}\ -\ \gonegtwo{d}\
                +\ \gonegtwo{e}\ +\ \gonegtwo{f}\ \ \Biggr]\ .
$$
As an example, let us compute the vertex product with the theta
graph:

\def\bubg#1{\risS{-13}{bubG2#1}{\put(4,5){\mbox{$\scriptstyle G$}}
                }{15}{15}{15}}

$$\beta\vee \risS{-11}{g2}{\put(4,5){\mbox{$\scriptstyle G$}}}{15}{15}{12} =
 \risS{-5}{g3_vpr_bub}{}{9}{20}{12}\vee
       \risS{-11}{g2}{\put(4,5){\mbox{$\scriptstyle G$}}}{15}{15}{12} =
 \frac{1}{6}\Biggl[\ \bubg{a} - \bubg{b} - \bubg{c} - \bubg{d} +
     \bubg{e} + \bubg{f} \ \Biggr] =
    \risS{-11}{g2}{\put(4,5){\mbox{$\scriptstyle G$}}}{15}{15}{12}\ ,
$$
since all the summands in the brackets (taken with their signs) are
equal to each other due to the AS relation. Therefore, $\beta$ will
be the unit for the vertex product on $\G_{\geq 1}$.

In order to simplify the notation, we shall use diagrams with
shaded disks, understanding them as alternating linear
combinations of six graphs as above. For example:
$$\risS{-9}{g1}{\put(2,16){\mbox{$\scriptstyle G_1$}}}{15}{20}{10}\ \vee\
  \risS{-11}{g2}{\put(2,5){\mbox{$\scriptstyle G_2$}}}{15}{20}{12}\quad =\quad
  \risS{-18}{G1sG2}{\put(2,5){\mbox{$\scriptstyle G_2$}}
              \put(2,33){\mbox{$\scriptstyle G_1$}}
                }{15}{20}{15} \quad =\quad
  \risS{-18}{G1G2s}{\put(2,5){\mbox{$\scriptstyle G_2$}}
              \put(2,33){\mbox{$\scriptstyle G_1$}}
                }{15}{20}{15}\ .
$$

\begin{theorem} \label{vertex_well_def}
The vertex product in $\G_{\geq 1}$ is well-defined, commutative
and associative.
\end{theorem}

\begin{proof}
It is sufficient to prove that the the AS and the IHX relations
imply the following equality:
$$
 X_1\ =\ \rb{-9mm}{\ig[height=20mm]{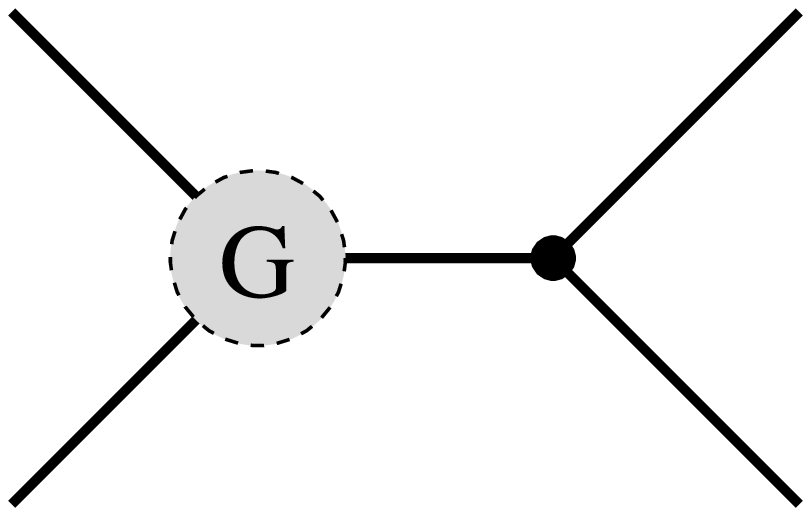}}
 \ =\ \rb{-9mm}{\ig[height=20mm]{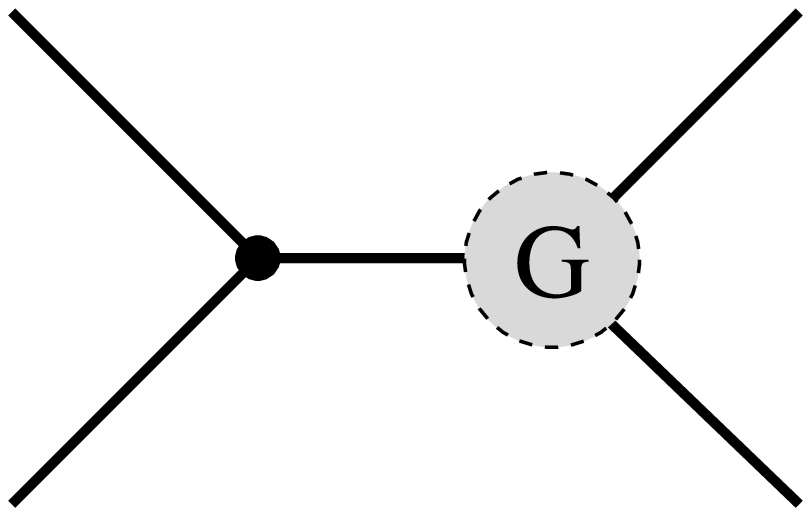}}\ =\ X_2.
$$
where $G$ denotes an arbitrary subgraph with three legs (and each picture
is the alternating sum of six diagrams).

By the Kirchhoff law we have:
\begin{align*}
 & \rb{-7mm}{\ig[height=15mm]{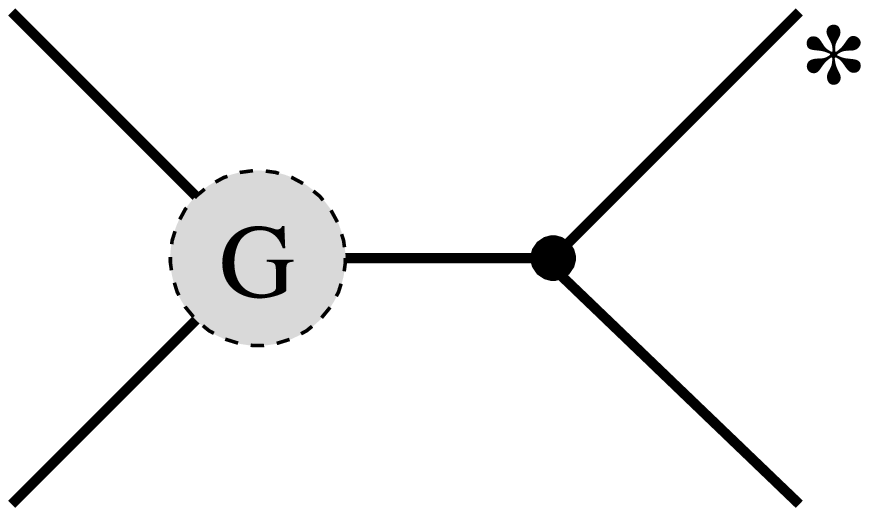}}
 \ =\ \rb{-7mm}{\ig[height=15mm]{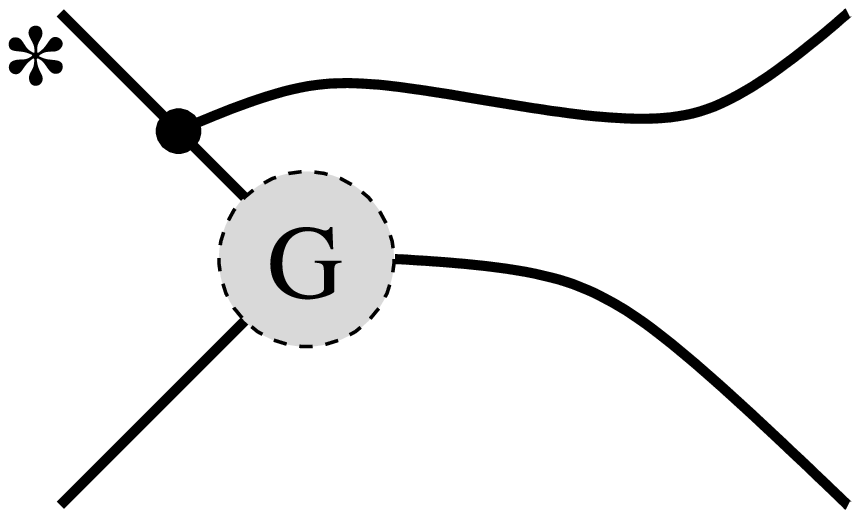}}
 \ +\ \rb{-7mm}{\ig[height=15mm]{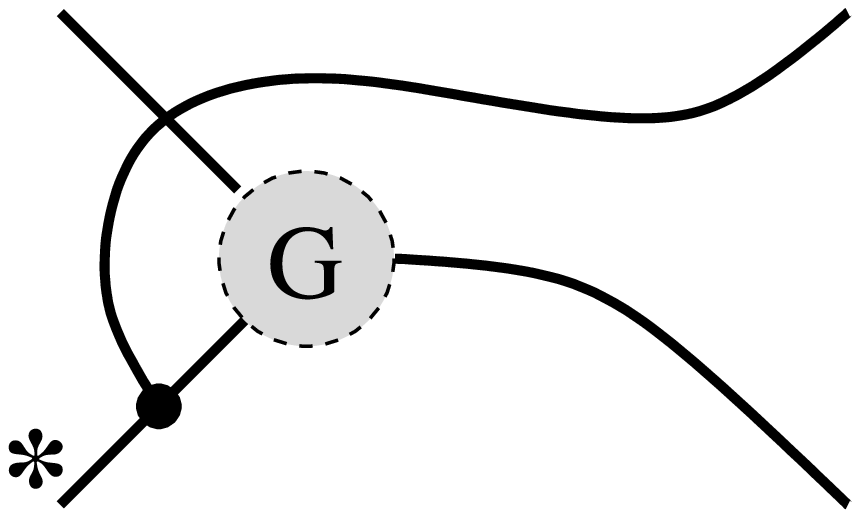}} \\[3mm]
 & =\ \rb{-7mm}{\ig[height=15mm]{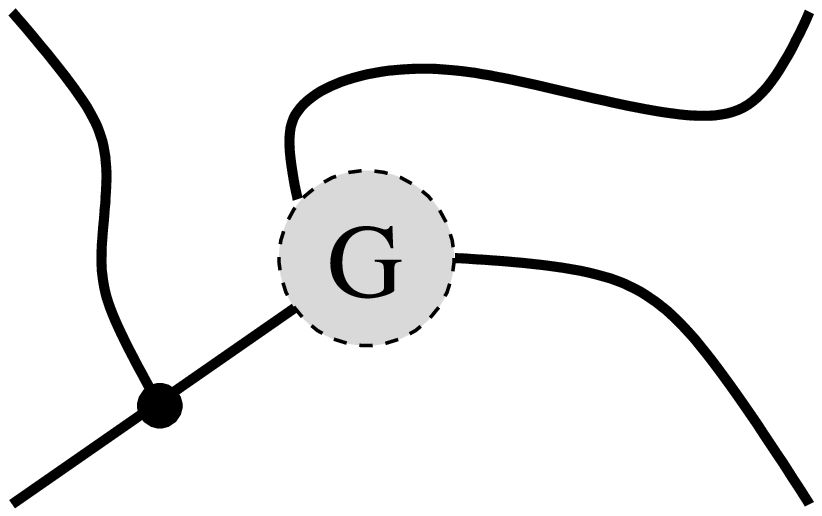}}
 \ +\ \rb{-7mm}{\ig[height=15mm]{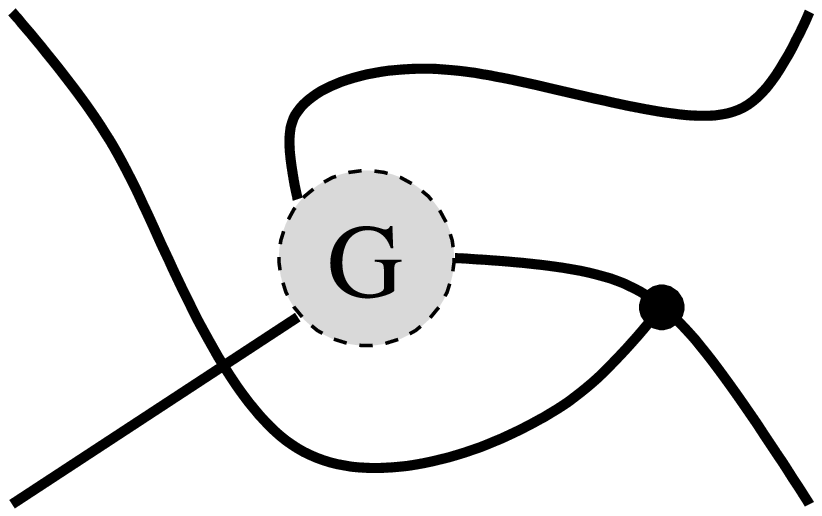}}
 \ +\ \rb{-7mm}{\ig[height=15mm]{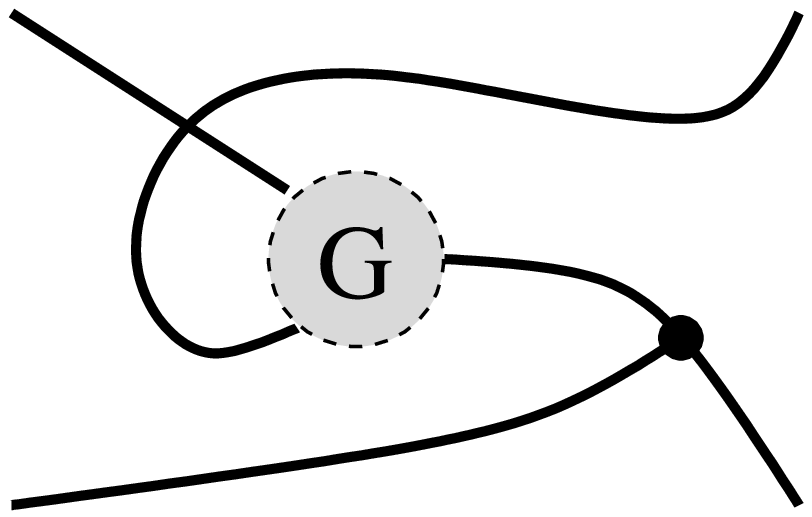}}
 \ +\ \rb{-7mm}{\ig[height=15mm]{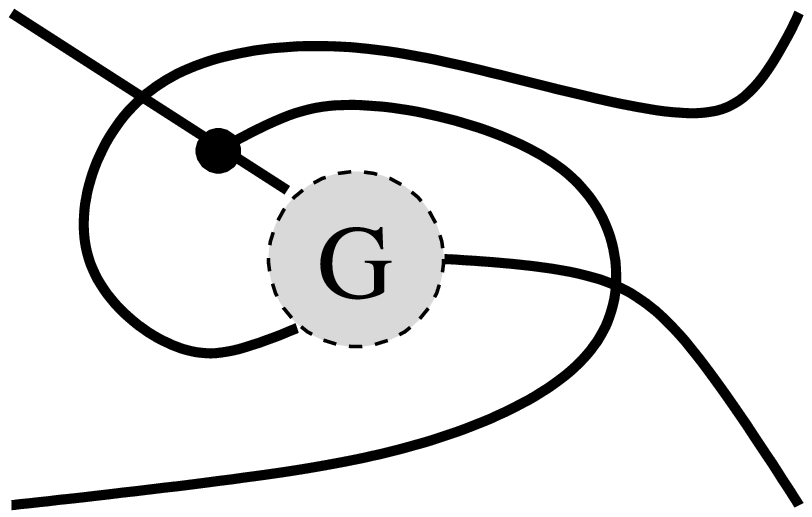}}
\end{align*}
(the stars indicate the place where the tail of the ``moving
electron'' is fixed in Kirchhoff's relation). Now, in the last line
the first and the fourth diagrams are equal to $X_2$, while the sum
of the second and the third diagrams is equal to $-X_1$ (again, by
an application of Kirchhoff's rule). We thus have $2X_1=2X_2$ and
therefore $X_1=X_2$.

Commutativity and associativity are obvious.
\end{proof}

\begin{remark}
Unlike the edge product, which respects the grading on $\G$, the
vertex multiplication is an operation of degree $-1$:
$$\G_n\vee \G_m \subset \G_{n+m-1}.$$
\end{remark}

\subsection{Relation between the two products in $\G$}
\begin{xproposition}\label{v_prod_in_G}
The edge product $\cdot$ in the algebra of 3-graphs $\G$ is
related to the vertex product $\vee$ on $\G_{\geq 1}$  as follows:
$$
  G_1\cdot G_2 = \beta\cdot(G_1 \vee G_2).
$$
\end{xproposition}

\begin{proof}
Choose a vertex in each of the given graphs $G_1$ and $G_2$ and call
its complement  $G_1'$ and $G_2'$, respectively:
$$G_1 = \risS{-11}{Gvpr1}{\put(29,23){\mbox{$\scriptstyle G_1'$}}}{50}{20}{12} =
        \risS{-11}{Gvpr1s}{\put(29,23){\mbox{$\scriptstyle G_1'$}}}{50}{20}{12}\ ,\qquad\qquad
G_2 = \risS{-11}{Gvpr1}{\put(29,23){\mbox{$\scriptstyle
G_2'$}}}{50}{20}{12}\ ,
$$
where, as explained above, the shaded region indicates the
alternating average over the six permutations of the three legs.

Then, by Theorem~\ref{vertex_well_def} we have:
$$\begin{array}{ccl}
G_1\cdot G_2&=& \risS{-20}{Gprod1}{\put(37,30){\mbox{$\scriptstyle
G'_1$}}
               \put(0,45){\mbox{$\scriptstyle G_1$}}
               \put(91,30){\mbox{$\scriptstyle G'_2$}}
               \put(110,45){\mbox{$\scriptstyle G_2$}}
               }{120}{35}{12}\ =\
    \risS{-20}{Gprod2}{\put(43,31){\mbox{$\scriptstyle G'_1$}}
              \put(70,31){\mbox{$\scriptstyle G'_2$}}}{95}{20}{12} \\
&=&\beta\cdot
    \risS{-20}{Gprod3}{\put(12,28){\mbox{$\scriptstyle G'_1$}}
              \put(41,28){\mbox{$\scriptstyle G'_2$}}}{65}{45}{12}
     \ =\ \beta\cdot (G_1\vee G_2)\ .
\end{array}
$$
\end{proof}

\section{Action of $\G$ on the primitive space $\PR$}

\subsection{Edge action of $\G$ on $\PR$}
\label{G_ed_act_on P}

As we know (Section~\ref{FDprim}) the space $\PR$ of the primitive
elements in the algebra $\F$ is spanned by connected diagrams, that
is, closed Jacobi diagrams which remain connected after the Wilson
loop is stripped off. It is natural to define the {\em edge action}
of $\G$ on a primitive diagram $D\in\PR$ simply by taking the edge
product of a graph $G\in \G$ with $D$ as if $D$ were a 3-graph,
using an internal edge of $D$.
The resulting graph $G\cdot D$ is again a closed diagram; moreover,
it lies in $\PR$. Since $D$ is connected, Lemmas~\ref{lemma1a} and
\ref{lemma2a} imply that $G\cdot D$ does not depend on the choice of
the edge in $G_1$ and of the internal edge in $D$.  Therefore, we
get a well-defined action of $\G$ on $\PR$, which is clearly
compatible with the gradings:
$$\G_n\cdot \PR_m \subset \PR_{n+m}.$$

\begin{proposition} \label{iso_P2_and_G}
The vector space $\G$ is isomorphic as a graded vector space (with
the grading shifted by one) to the subspace $\PR^2\subset \PR$ of
primitive closed diagrams spanned by connected diagrams with 2
legs: $\G_n \cong \PR^2_{n+1}$ for all $n\geq 0$.
\end{proposition}

\begin{proof}
The isomorphism $\G\to \PR^2$ is given by the edge action of
3-graphs on the element $\Theta\in\PR^2$ represented by the chord
diagram with a single chord, $G\mapsto G\cdot \Theta$. The inverse
map is equally simple. For a connected closed diagram $D$ with two
legs strip off the Wilson loop and glue together the two loose ends
of the resulting diagram, obtaining a 3-graph of degree one less
than $D$. Obviously, this map is well-defined and inverse to the
edge action on $\Theta$.
\end{proof}

\subsection{Vertex action of $\G$ on $\PR$}
\label{G_v_act_on_P}

In order to perform the vertex multiplication, we need at least one
vertex in each of the factors. Therefore, we shall define an action
of the algebra $\G_{\geq 1}$ (with the vertex product) on the space
$\PR_{> 1}$ of primitive elements of degree strictly greater than 1.

The action $G\vee D$ of $G$ on $D$ is the alternated average over
all six ways of inserting $G$, with one vertex removed, into $D$
with one internal vertex taken out. Again, since $D$ is connected,
the proof of Theorem~\ref{vertex_well_def} works to show that this
action is well-defined. Note that the vertex action decreases the
total grading by 1 and preserves the number of legs:  $$\G_n\vee
\PR^k_m \subset \PR^k_{n+m-1}.$$

The simplest element of $\PR$ on which $\G_{\geq 1}$  acts in this
way is the ``Mercedes-Benz diagram'' $$\ol{t}_1 =
\risS{-12}{mercedes}{}{30}{20}{15}.$$

\begin{lemma}\label{G_v_mercedes}
\ \\
\begin{itemize}
\item[(a)]
The map $\G_{\geq 1}\to \PR$ defined as $G \mapsto G\vee \ol{t}_1$
is injective.
\item[(b)]
For all $G\in \G_{\geq q}$ we have $$ G\vee \ol{t}_1=\frac{1}{2}
G\cdot \Theta.$$
\end{itemize}
\end{lemma}

\begin{proof}
Indeed, $\ol{t}_1=\frac{1}{2}\beta\cdot \Theta$. Therefore,
$$G\vee\ol{t}_1=\frac{1}{2}(G\vee\beta)\cdot\Theta =
\frac{1}{2}G\cdot\Theta.$$ Since the map $G\mapsto G\cdot \Theta$ is
an isomorphism (Proposition~\ref{iso_P2_and_G}), the map $G\mapsto
G\vee \ol{t}_1$ is also an isomorphism $\G_{\geq 1}\isom \PR^2_{>
1}$.
\end{proof}

\subsection{A product on the primitive space $\PR$}
In principle, the space of primitive elements $\PR$ of the algebra
$\F$ does not possess any a priori defined multiplicative structure.
Primitive elements only generate the algebra $\F$ much in the same
way as the variables $x_1,\dots,x_n$ generate the polynomial algebra
$\R[x_1,\dots,x_n]$. However, the link between the space $\PR$ and
the algebra of 3-graphs $\G$ allows to introduce  a
(non-commutative) multiplication in $\PR$.

There is a projection $\pi: \PR_n \to \G_n$, which consists in
introducing a cyclic order on the half-edges at the vertices of the
Wilson loop according to the rule ``forward--sideways--backwards''
and then forgetting the fact that the Wilson loop was distinguished.
The edge action $\G \times \PR \to \PR$ then gives rise to an
operation $*: \PR \times \PR \to \PR$ defined by the rule
$$
  p*q = \pi(p)\cdot q.
$$
where $\pi:\PR\to\G$ is the homomorphism of forgetting the Wilson loop
defined above.

The operation $*$ is associative, but, in general, non-commutative:
$$
\rb{-4mm}{\ig[height=10mm]{fdb_ex_sl2}}
\ *\ \rb{-4mm}{\ig[height=10mm]{fdw4_ex_sl2}}\ =\
  \rb{-4mm}{\ig[height=10mm]{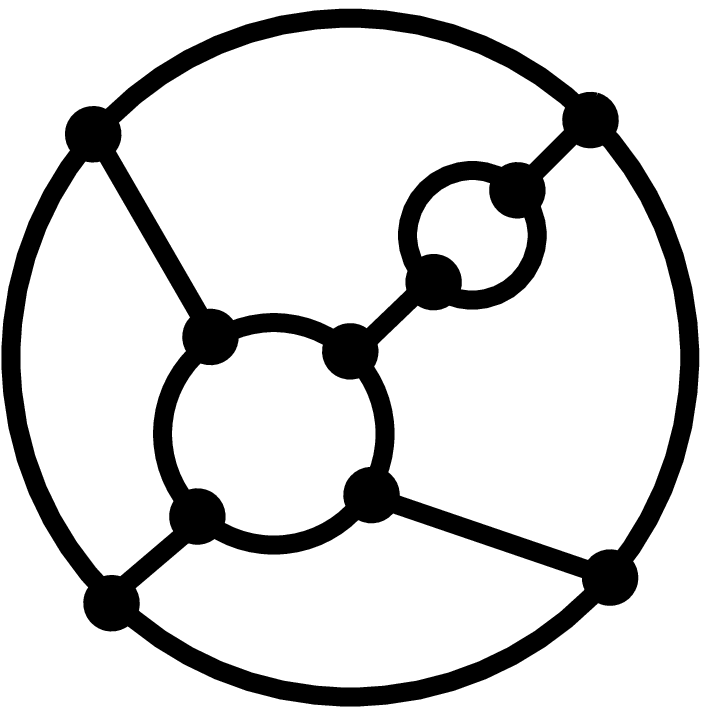}}\ ,\quad\mbox{but}\quad
\rb{-4mm}{\ig[height=10mm]{fdw4_ex_sl2}}
\ *\ \rb{-4mm}{\ig[height=10mm]{fdb_ex_sl2}}\ =\
  \rb{-4mm}{\ig[height=10mm]{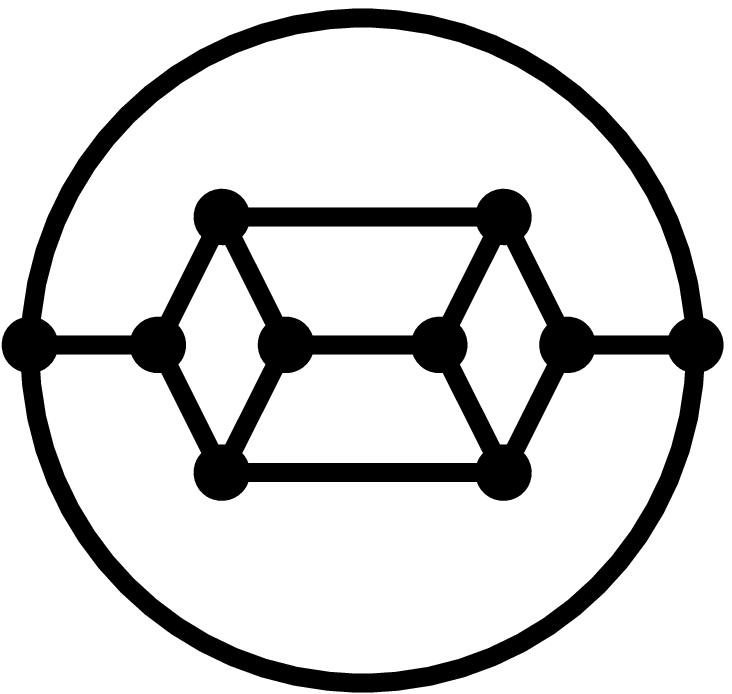}}\ .
$$

These two elements of the space $\PR$ are different; they can be
distinguished, for instance, by the $\sL_2$-invariant (see
Exercise~\ref{ex_sl2_on_cincub_cubinb} at the end of
Chapter~\ref{LAWS}). However, $\pi$ projects these two elements into
the same element $\beta\cdot\om_4\in\G_5$.

\section{Lie algebra weight systems for the algebra $\G$}
\label{LAWS_3G}

A {\em weight system for 3-graphs} is a function on 3-graphs that
satisfies the antisymmetry and the IHX relations. Lie algebras (and
super Lie algebras) give rise to weight systems for 3-graphs in the
very same fashion as for the algebra $\F$ (Section~\ref{LAWS_C}),
using the structure tensor $J$. Since 3-graphs have no univalent
vertices, these weight systems take values in the ground field (here
assumed to be $\C$). For a graph $G\in\G$ we put $$\f_\g(G) :=
T_\g(G)\in \g^0\cong\C.$$

In particular, the weight system $\f_\g$ evaluated on the circle
(the 3-graph without vertices, which is the unit in $\G$) gives
the dimension of the Lie algebra.
\begin{xremark}
This circle should not be confused with the circles appearing in
the state sum formulae for $\f_{\sL_N}$ and $\f_{\so_N}$ from
Sections~\ref{ws_glN_on_C} and \ref{ws_soN_on_C}. The contribution
of these circles to the values of $\f_{\sL_N}$ and $\f_{\so_N}$ is
equal to $N$ while the dimensions of the corresponding Lie
algebras are $N^2-1$ and $\frac{1}{2}N(N-1)$, respectively.
\end{xremark}

\subsection{Changing the bilinear form}
From the construction of $\f_\g$ it is easy to see that the
function $\f_{\g, \lambda}$ corresponding to the form
$\lambda\langle\cdot,\cdot\rangle$ is a multiple of $\f_\g$:
$$\f_{\g,{\lambda}}(G) = \lambda^{-n}\f_\g(G)
$$
for $G\in\G_n$.

\subsection{Multiplicativity with respect to
the edge product in $\G$}
\begin{xproposition} For a simple Lie algebra $\g$ and any choice of
an $\ad$-invariant non-degenerate symmetric bilinear form
$\langle\cdot,\cdot\rangle$ the function $\frac{1}{\dim \g}
\f_\g:\G\to\C$ is multiplicative with respect to the edge product in
$\G$.
\end{xproposition}

\begin{proof}
This is a consequence of the fact that, up to a constant, the
quadratic Casimir tensor of a simple Lie algebra is the only
$\ad$-invariant, symmetric, non-degenerate tensor in $\g\ot\g$.

Consider two graphs $G_1,G_2\in\G$ and chose an orthonormal basis
$e_i$ for the Lie algebra $\g$. Cut an arbitrary edge of the graph
$G$ and consider the tensor that corresponds to the resulting graph
$G_1'$ with two univalent vertices. This tensor is a scalar multiple
of the quadratic Casimir tensor $c\in \g\ot\g$:
$$a\cdot c = a \sum_{i=1}^{\dim \g} e_i\ot e_i\ .
$$
Now, $\f_{\g,K}(G_1)$ is obtained by contracting these two tensor
factors.  This gives $\f_{\g,K}(G_1) = a\dim \g$, and
$a=\frac{1}{\dim \g}\f_\g(G_1)$. Similarly, for the graph $G_2$ we
get the tensor $\frac{1}{\dim \g}\f_\g(G_2)\cdot c$. Now, if we join
together one pair of univalent vertices of the graphs $G_1'$ and
$G_2'$ (where $G_2'$ is obtained from $G_2$ by cutting an edge), the
partial contraction of the element $c^{\ot 2}\in\g^{\ot 4}$ will
give $$\frac{1}{(\dim\g)^2}\, \f_\g(G_1)\f_\g(G_2) \cdot
c\in\g\ot\g.$$ But, on the other hand, this tensor equals
$$\frac{1}{\dim \g}\, \f_\g(G_1\cdot G_2)\cdot  c\in\g\ot\g.$$ This
shows that $\frac{1}{\dim\g}\f_\g$ is multiplicative.
\end{proof}

\subsection{Compatibility with the edge action of $\G$ on $\F$}
\label{comp_e_act}

Recall the definition of the edge action of 3-graphs on closed
diagrams (see Section~\ref{G_ed_act_on P}). We choose an edge in
$G\in\G$ and an internal edge in $D\in\F$, and then take the
connected sum of $G$ and $D$ along the chosen edges. In fact, this
action depends on the choice of the connected component of the
internal graph of $D$ containing the chosen edge. It is well defined
only on the primitive subspace $\PR\subset\F$. In spite of this
indeterminacy we have the following lemma.

\begin{xlemma}
For any choice of the glueing edges, $\f_\g(G\cdot D) =
\frac{\f_\g(G)}{\dim \g}\, \f_\g(D)$.
\end{xlemma}
\begin{proof}
Indeed, in order to compute $\f_\g(D)$ we assemble the tensor
$T_\g(D)$ from tensors that correspond to tripods and chords. The
legs of these elementary pieces are glued together by contraction
with the quadratic Casimir tensor $c$, which corresponds to the
metric on the Lie algebra. By the previous argument, to compute the
tensor $T_\g(G\cdot D)$ one must use for the chosen edge the tensor
$\frac{1}{\dim \g}\f_\g(G)\cdot c$ instead of $c$. This gives the
coefficient $\frac{1}{\dim \g}\f_\g(G)$ in the expression for
$\f_\g(G\cdot D)$ as compared with $\f_\g(D)$.
\end{proof}

One particular case of the edge action of $\G$ is especially
interesting: when the graph $G$ varies, while $D$ is fixed and equal
to $\Theta$, the chord diagram with only one chord. In this case the
action is an isomorphism of the vector space $\G$ with the subspace
$\PR^2$ of the primitive space $\PR$ generated by connected closed
diagrams with 2 legs (section \ref{iso_P2_and_G}).

\begin{xcorollary}  For the weight systems associated with a
simple Lie algebra $\g$ and the Killing form
$\langle\cdot,\cdot\rangle^K$:
$$\f_{\g,K}(G) = \f^{ad}_{\g,K}(G\cdot \Theta)\ ,$$
where $\f^{ad}_{\g,K}$ is the weight system corresponding to the
adjoint representation of $\g$.
\end{xcorollary}

\begin{proof}
Indeed, according to the Lemma, for the universal enveloping algebra invariants we have
$$\f_{\g,K}(G\cdot \Theta) = \frac{1}{\dim \g}\f_\g(G) \f_{\g,K}(\Theta)
  =\frac{1}{\dim \g}\f_\g(G) \sum_{i=1}^{\dim \g} e_i e_i\ ,$$
where $\{e_i\}$ is a basis orthonormal with respect to the Killing
form. Now to compute $\f^{ad}_{\g,K}(G\cdot \Theta)$ we take the
trace of the product of operators in the adjoint representation:
$$\f^{ad}_{\g,K}(G\cdot \Theta) = \frac{1}{\dim \g}\f_\g(G)
       \sum_{i=1}^{\dim \g} \mbox{Tr}( \mbox{ad}_{e_i} \mbox{ad}_{e_i} ) = \f_\g(G)
$$
by the definition of the Killing form.
\end{proof}

\subsection{Multiplicativity with respect to the vertex product
in $\G$}
\begin{xproposition} Let $w:\G\to\C$ be an edge-multiplicative weight system,
and $w(\beta)\not= 0$. Then $\frac{1}{w(\beta)} w: \G\to\C$ is
multiplicative with respect to the vertex product. In particular,
for a simple Lie algebra $\g$, $\frac{1}{\f_\g(\beta)}\f_\g$ is
vertex-multiplicative.
\end{xproposition}

\begin{proof}
According to \ref{v_prod_in_G} the edge product is related to the
vertex product as $G_1\cdot G_2 = \beta\cdot(G_1\vee G_2)$.
Therefore,
$$w(G_1)\cdot w(G_2) = w(\beta\cdot(G_1\vee G_2)) = w(\beta)\cdot w(G_1\vee G_2)\ .
$$
This means that the weight system $\frac{1}{w(\beta)} w: \G\to\C$ is multiplicative
with respect to the vertex product.
\end{proof}

\begin{xcorollary} The weight systems $\frac{1}{2N(N^2-1)} \f_{\sL_N},\
                \frac{2}{N(N-1)(N-2)} \f_{\so_N}:\G\to\C$
associated with the ad-invariant form $\langle x,y \rangle =
\mbox{\rm Tr}(xy)$ are multiplicative with respect to the vertex
product in $\G$.
\end{xcorollary}

This follows from a direct computation for the ``bubble'':
$$\f_{\sL_N}(\beta)= 2N(N^2-1),\mathrm{\ and\ } \f_{\so_N}(\beta)=
\frac{1}{2}N(N-1)(N-2).$$

\subsection{Compatibility with the vertex action of $\G$ on $\F$}
\label{comp_v_act} The vertex action $G\vee D$ of a 3-graph $G\in
\G$ on a closed diagram $D\in\F$ with at least one vertex (see
Section~\ref{G_v_act_on_P}) is defined as the alternating sum of 6
ways to glue the graph $G$ with the closed diagram $D$ along chosen
internal vertices in $D$ and $G$. Again, this action is well-defined
only on the primitive space $\PR_{>1}$.

\begin{xlemma}  Let $\g$ be a simple Lie algebra. Then for any choice
of the glueing vertices in $G$ and $D$:
$$\f_\g(G\vee D) = \frac{\f_\g(G)}{\f_\g(\beta)} \f_\g(D)\ .$$
\end{xlemma}

\begin{proof}
Using the edge action (Section~\ref{comp_e_act}) and its relation to
the vertex action we can write
$$\frac{\f_\g(G)}{\dim \g} \f_\g(D) =  \f_\g(G\cdot D)
   = \f_\g(\beta\cdot(G_1\vee D))
   = \frac{\f_\g(\beta)}{\dim \g} \f_\g(G_1\vee D)\ ,
$$
which is what we need.
\end{proof}

\subsection{$\sL_N$- and $\so_N$-polynomials}
The $\sL_N$- and $\so_N$-{\em polynomials} are the weight systems
$\f_{\sL_N}$ (with respect to the bilinear form $\bilinf{x}{y} =
\mbox{Tr}(xy)$), and $\f_{\so_N}$  (with the bilinear form
$\bilinf{x}{y} = \frac{1}{2}\mbox{Tr}(xy)$). In the case of
$\so_N$ the choice of {\em half} the trace as the bilinear form is
more convenient since it gives polynomials in $N$ with integral
coefficients. In particular, for this form $\f_{\so_N}(\beta)=
N(N-1)(N-2)$, and in the state sum formula from the Theorem of
Section~\ref{ws_soN_on_C} the coefficient in front of the sum
equals 1.

The polynomial $\f_{\sL_N}(G)\ (= \f_{\gl_N}(G))$ is divisible by
$2N(N^2-1)$ (Exercise~\ref{ex_gl_N_on_G_root} in the end of this
chapter) and the quotient is a multiplicative function with respect
to the vertex product. We call this quotient the {\em reduced
$\sL$-polynomial\/} and denote it by $\wt{\sL}(G)$. Dividing the
$\so$-polynomial $\f_{\so_N}(G)$ by $N(N-1)(N-2)$ (see
Exercise~\ref{ex_so_N_on_G_root}), we obtain the {\em reduced
$\so$-polynomial\/} $\wt{\so}(G)$, which is also multiplicative with
respect to the vertex product.

A.~Kaishev \cite{Kai} computed the values of $\wt{\sL}$-, and
$\wt{\so}$-polynomials on the generators of $\G$ of small degrees
(for $\wt{\so}$-polynomial the substitution $M=N-2$ is used in the
table):

\def\slpol#1{\parbox[t]{87pt}{$\scriptstyle #1$}}
\def\sopol#1{\parbox[t]{192pt}{$\scriptstyle #1$\vspace{6pt}}}

$$\begin{array}{|c|c|l|l|}
\hline
\mbox{deg} &   &  \wt{\sL}\mbox{-polynomial}\makebox(0,14){}
               &  \wt{\so}\mbox{-polynomial}\\
\hline
 1 & \beta & 1 & 1 \\
\hline
 4 & \om_4 & \slpol{N^3 +12N}
    & \sopol{M^3 - 3M^2 + 30M - 24} \\
\hline

 6 & \om_6 & \slpol{N^5+32N^3+48N}
    & \sopol{M^5 - 5M^4 + 80M^3 - 184M^2 + 408M - 288}    \\
\hline

 7 & \om_7 & \slpol{N^6+64N^4+64N^2}
    & \sopol{M^6 - 6M^5 + 154M^4 - 408M^3 + 664M^2 - 384} \\
\hline

 8 & \om_8 & \slpol{N^7+128N^5+128N^3 \\ +192N}
    & \sopol{M^7 - 7M^6 + 294M^5 - 844M^4 + 1608M^3 - 2128M^2 \\
    + 4576M - 3456} \\
\hline

 9 & \om_9 & \slpol{N^8+256N^6+256N^4 \\ +256N^2}
    & \sopol{M^8 - 8M^7 + 564M^6 - 1688M^5 + 3552M^4 - 5600M^3 \\
    - 5600M^3 + 6336M^2 + 6144M - 9216} \\
\hline

10 & \om_{10} & \slpol{N^9+512N^7+512N^5 \\ +512N^3+768N}
    & \sopol{M^9 - 9M^8 + 1092M^7 - 3328M^6 + 7440M^5 - 13216M^4 \\
    + 18048M^3 - 17920M^2 + 55680M - 47616} \\
\hline

10 & \delta & \slpol{N^9+11N^7+114N^5 \\ -116N^3}
    & \sopol{M^9 - 9M^8 + 44M^7 - 94M^6 + 627M^5 + 519M^4 \\
    - 2474M^3 - 10916M^2 + 30072M - 17760} \\
\hline

11 & \om_{11} & \slpol{ N^{10}+1024N^8+1024N^6 \\ +1024N^4+1024N^2}
    & \sopol{M^{10} - 10 M^9 + 2134 M^8 - 6536 M^7 + 15120 M^6 \\
    - 29120 M^5 + 45504 M^4 - 55040 M^3 + 48768 M^2 \\+ 145408 M - 165888} \\
\hline
\end{array}\index{Table of!Lie algebra weight systems on $\G$}$$

There are recognizable patterns in this table. For example, we see
that
\begin{center}
\begin{tabular}{l}
$\wt{\sL}(\om_n) = N^{n-1} + 2^{n-1}(N^{n-3} + \dots + N^2)$, for odd $n>5$; \\
$\wt{\sL}(\om_n) = N^{n-1} + 2^{n-1}(N^{n-3} + \dots + N^3) +
                   2^{n-2}3N$, for even $n\geq 4$.
\end{tabular}\end{center}
It would be interesting to know if these observations are particular
cases of some general statements.

\section{Vogel's algebra $\Lambda$} %15
\label{lambda}

Diagrams with 1- and 3-valent vertices can be considered with
different additional structures on the set of univalent vertices. If
there is no structure, then we get the notion of an open Jacobi
diagram; open diagrams are considered modulo AS and IHX relations.
If the legs are attached to a circle or a line, then we obtain
closed Jacobi diagrams; for the closed diagrams, AS, IHX and STU
relations are used. Connected diagrams with a linear order
(numbering) on the set of legs and a cyclic order on the half-edges
at each 3-valent vertex, considered modulo AS and IHX, but without
STU relations, will be referred to as {\em fixed diagrams}.
\index{Diagram!fixed} The set ${\bf X}$ of all fixed diagrams has
two gradings: by the number of legs (denoted by a superscript) and
by half the total number of vertices (denoted by a subscript).

\begin{xdefinition}
The $\Q$-vector space spanned by fixed diagrams with $k$ legs modulo
the usual AS and IHX relations
$$
   \X^k_n = \langle {\bf X}^k_n\rangle/\langle AS,IHX\rangle,
$$
is called the {\em space of fixed diagrams} of degree $n$ with $k$
legs.
\end{xdefinition}
We shall write $\X^k$ for the direct sum $\oplus_n \X^k_n$.

\begin{xremark} The spaces
$\X^k$ for different values of $k$ are related by various
operations. For example, one may think about the diagram
$\risS{-2}{mum_yum3b}{}{15}{0}{0}$ as of a linear operator from
$\X^4$ to $\X^3$. Namely, it acts on an element $G$ of $\X^4$ as
follows:

$$\risS{-2}{mum_yum3b}{}{15}{0}{0}:\
    \risS{-7}{el1_of_x4}{\put(12,10){\mbox{$G$}}
               \put(1,-6){\mbox{$\scriptstyle 1$}}
               \put(6,-6){\mbox{$\scriptstyle 2$}}
               \put(11,-6){\mbox{$\scriptstyle 3$}}
               \put(25,-6){\mbox{$\scriptstyle 4$}}}{30}{0}{0}\ \ \mapsto\ \
    \risS{-15}{el1_of_x3}{\put(12,25){\mbox{$G$}}
               \put(12,-6){\mbox{$\scriptstyle 1$}}
               \put(18,-6){\mbox{$\scriptstyle 2$}}
               \put(23,-6){\mbox{$\scriptstyle 3$}}}{30}{20}{12}\ .
$$
\end{xremark}

\begin{xca}\label{ex_x3_x4}
\ \\
\begin{itemize}
\item[(a)] Prove the following relation
$$
  \rb{-5mm}{\ig[width=15mm]{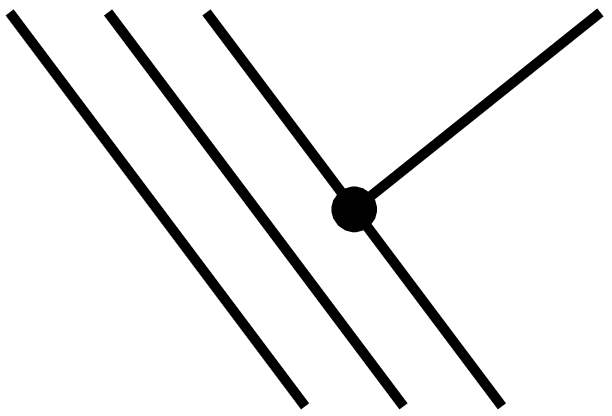}}
  +\rb{-5mm}{\ig[width=15mm]{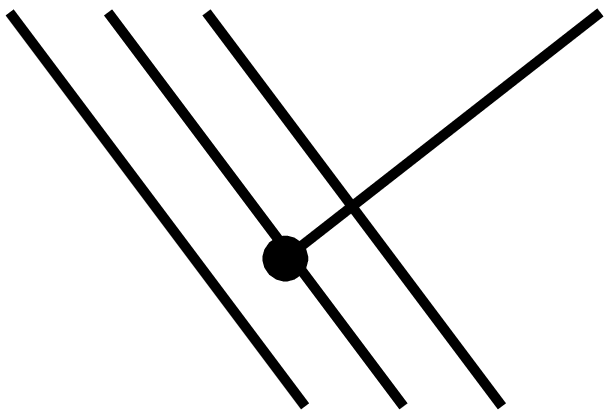}}
  +\rb{-5mm}{\ig[width=15mm]{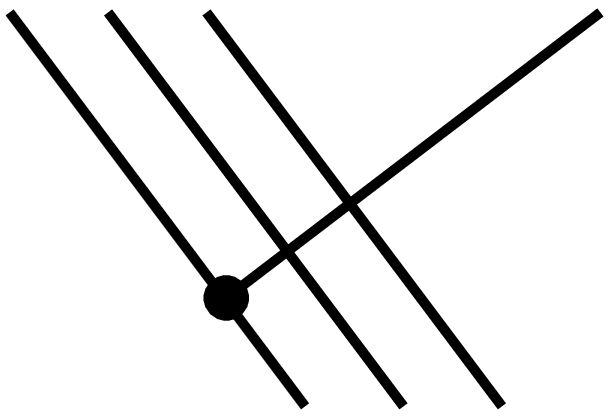}} = 0,
$$
among the three linear operators from $\X^4$ to $\X^3$.
\item[(b)] Prove that
$$
  \rb{-5mm}{\ig[width=15mm]{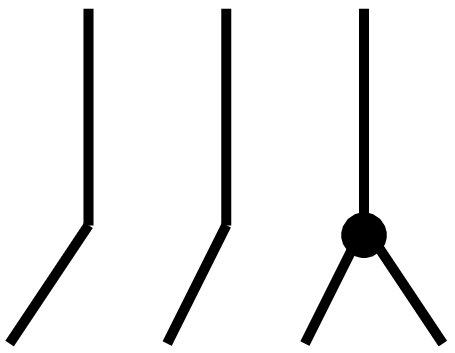}}
  =\rb{-5mm}{\ig[width=15mm]{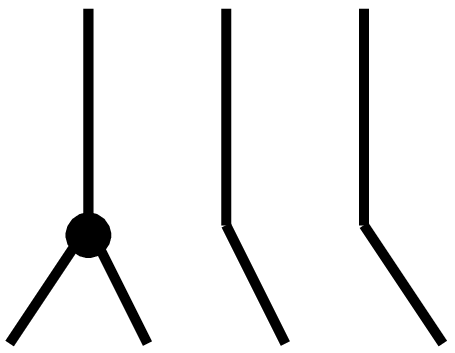}}
$$
as linear maps from $\X^3$ to $\X^4$.
\end{itemize}
\end{xca}

The space of open diagrams $\B^k$ studied in Chapter \ref{algOD} is
the quotient of $\X^k$ by the action of the symmetric group $S_k$
which permutes the legs of a fixed diagram. The quotient map
$\X\to\B$ has a nontrivial kernel; for example, a tripod, which is
nonzero in $\X^3$, becomes zero in $\B$:

$$
  0\ \ne\ \rb{-4mm}{\ig[height=10mm]{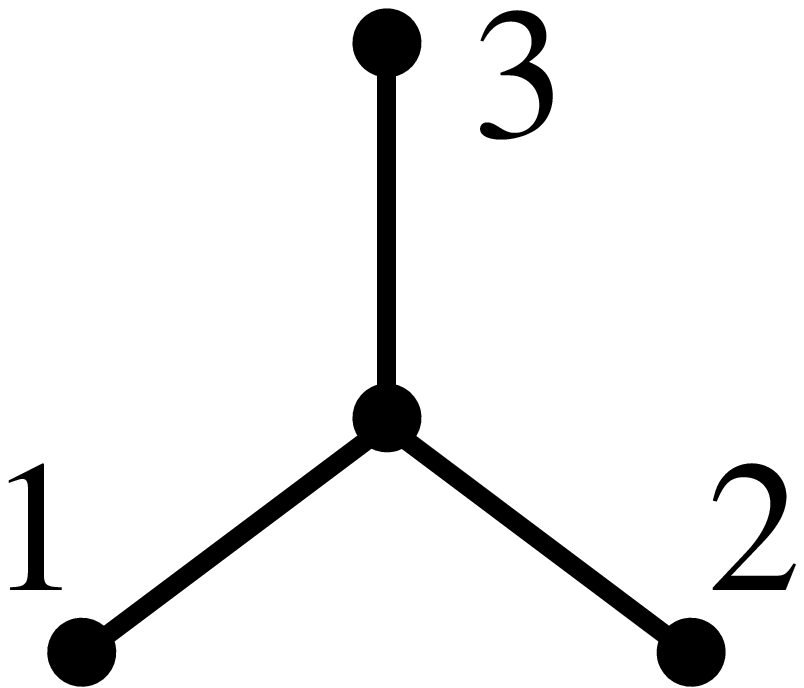}}
\longmapsto \rb{-4mm}{\ig[height=10mm]{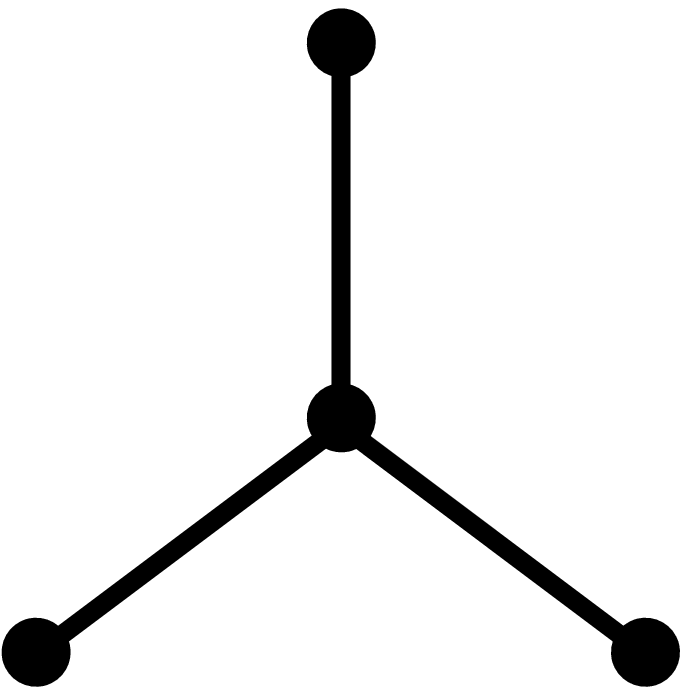}} =\ 0.
$$

\subsection{The algebra $\Lambda$}

\label{Vogel-algebra}
\index{Vogel!algebra $\Lambda$}
\index{Algebra!Vogel's $\Lambda$}
The algebra $\Lambda$ is the subspace of $\X^3$ that consists of all
elements antisymmetric with respect to permutations of their legs.
The product in $\Lambda$ is similar to the vertex product in $\G$.
Given a connected fixed diagram and an element of $\Lambda$, we
remove an arbitrary vertex in the diagram and insert the element of
$\Lambda$ instead --- in compliance with the cyclic order at the
vertex. This operations extends to a well-defined product on
$\Lambda$, and this fact is proved in the same way as for the vertex
multiplication in $\Gamma$. Since antisymmetry is presupposed, we do
not need to take the alternated average over the six ways of
insertion, as in $\Gamma$,
--- all the six summands will be equal to each other.

\def\vogt#1{\risS{-15}{#1}{\put(-2,-6){\mbox{$\scriptstyle 1$}}
          \put(28,-6){\mbox{$\scriptstyle 2$}}
          \put(10,38){\mbox{$\scriptstyle 3$}}}{30}{20}{30}}
\begin{xexample}

$$\vogt{vog_t} \quad\vee\quad \vogt{vog_t} \quad=\quad \vogt{vog_tt}\ .
$$
\end{xexample}

Conjecturally, the antisymmetry requirement in this definition is
superfluous:
\begin{conjecture}\label{nechet_X3}
$\Lambda=\X^3$, that is, any fixed diagram with 3 legs is
antisymmetric with respect to leg permutations.
\end{conjecture}

\begin{xremark}
Note the sign difference in the definitions of the product in
$\Lambda$ and the vertex product in $\Gamma$: in $\Gamma$ when two
graphs are glued together in compliance with the cyclic order of
half-edges, the corresponding term is counted with a {\em negative}
sign.
\end{xremark}

\begin{xremark}
Vogel in \cite{Vo1} defines the spaces $\X^k$ and the algebra
$\Lambda$ over the integers, rather than over $\Q$. In this
approach, the equality (b) of Exercise~\ref{ex_x3_x4} no longer
follows from the AS and IHX relations. It has to be postulated
separately as one of the equations defining $\Lambda$ in $\X^3$ in
order to make the product in $\Lambda$ well-defined.
\end{xremark}

The product in the algebra $\Lambda$ naturally generalizes to the
action of $\Lambda$ on different spaces generated by 1- and 3-valent
diagrams, such as the space of connected open diagrams $\PR\B$ and
the space of 3-graphs $\G$. The same argument as above shows that
these actions are well-defined.

\subsection{Relation between $\Lambda$ and $\G$}

Recall that the space $\G_{\ge 1}$ of all 3-graphs of degree at
least one is an algebra with respect to the vertex product.

\begin{xproposition}
The algebra $\Lambda$ is isomorphic to $(\G_{\ge 1},\vee)$.
\end{xproposition}

\begin{proof}
There two mutually inverse maps between $\Lambda$ and $\G_{\ge 1}$.
The map from $\Lambda$ to $\G_{\ge 1}$ glues the three legs of an
element of $\Lambda$ to a new vertex so that the cyclic order of
edges at this vertex is opposite to the order of legs:
$$\risS{-27}{LtoG1}{\put(-2,6){\mbox{$\scriptstyle 1$}}
          \put(48,-6){\mbox{$\scriptstyle 2$}}
          \put(48,54){\mbox{$\scriptstyle 3$}}}{50}{30}{30}\quad\mapsto\quad
 (-1)\cdot\risS{-17}{LtoG2}{}{50}{30}{30}\quad .
$$
In order to define a map from $\G_{\ge 1}$ to $\Lambda$, we choose
an arbitrary vertex of a 3-graph, delete it and antisymmetrize:

\def\GtoL#1#2#3{\risS{-12}{GtoL1}{\put(10,32){\mbox{$\scriptstyle #1$}}
          \put(14,20){\mbox{$\scriptstyle #2$}}
          \put(27,26){\mbox{$\scriptstyle #3$}}}{30}{0}{0}}

$$\risS{-12}{GtoL0}{}{30}{25}{10}\quad\mapsto\quad \frac{1}{6}\Biggl[\
     -\GtoL123 + \GtoL213 + \GtoL132 + \GtoL321 - \GtoL231 - \GtoL312\ \Biggr]\ .
$$
It is easy to see that this is indeed a well-defined map (Hint: use
part (b) of Exercise~\ref{ex_x3_x4}).

It is evident from the definitions that both maps are inverse to
each other and send products to products.
\end{proof}

\begin{xremark} If Conjecture~\ref{nechet_X3} is true for $k=3$
then all the six terms (together with their signs) in the definition
of the map $\G_{\ge 1}\to\Lambda$ are equal to each other. This
means that there is no need to antisymmetrize. What we do is remove
one vertex (with a small neighbourhood) and number the three legs
obtained according to their cyclic ordering at the deleted vertex.
This would also simplify the definition of the vertex product in
Section~\ref{prod_in_G} as in this case
$$\risS{-8}{Gshed}{\put(4,13){\mbox{$\scriptstyle G$}}}{15}{10}{15} \quad =\quad
  \risS{-8}{Gclean}{\put(4,13){\mbox{$\scriptstyle G$}}}{15}{0}{0}\ ,
$$
and we simply insert one graph in a vertex of another.
\end{xremark}

\begin{xconjecture}[\cite{Vo1}]
\label{vog_conj}% conj 1
The algebra $\Lambda$ is generated by the elements $t$ and $x_k$
with odd $k=3,5,...$:\vspace{10pt}
$$\begin{array}{c@{\qquad}c@{\qquad}c@{\qquad}c@{\qquad}c}
\vogt{vog_t} & \vogt{vog_x3} & \vogt{vog_x4} & \vogt{vog_x5} & \cdots \\
t & x_3 & x_4 & x_5 & \cdots
\end{array}
$$
\end{xconjecture}

\subsection{Weight systems not coming from Lie algebras}
\label{ws_not_Lie}

In order to construct a weight system which would not be a
combination of Lie algebra weight systems, it is sufficient to find
a non-zero element in $\F$ on which all the Lie algebra weight
systems vanish. The same is true, of course, for super Lie algebras.

In \cite{Vo1} Vogel produces diagrams that cannot be detected by
(super) Lie algebra weight systems. Vogel's work involves heavy
calculations of which we shall give no details here. His
construction can be (very briefly) described as follows.

First, he gives a list of super Lie algebras with the property that
whenever all the weight systems for the algebras from this list
vanish on an element of $\Lambda$, this element cannot be detected
by any (super) Lie algebra weight system. This list includes a
certain super Lie algebra $D(2,1,\alpha)$; this algebra detects an
element of $\Lambda$ which all other algebras from the list do not
detect. Making this element act by the vertex action on the
``Mercedes-Benz'' closed diagram $\ol{t}_1$ he obtains a closed
diagram which is non-zero because of Lemma~\ref{G_v_mercedes} but
which cannot be detected by any super Lie algebra weight system. We
refer to \cite{Vo1} and \cite{Lieb} for the details.

\begin{xcb}{Exercises}

\begin{enumerate}

\item
Find an explicit chain of IHX and AS relations that proves the
following equality in the algebra $\G$ of 3-graphs:
$$\risS{-14}{g3_ex_g1}{}{65}{20}{20} \quad=\quad\risS{-22}{g3_ex_g2}{}{65}{20}{22}
$$

\item
Let $\tau_2:\X^2\to\X^2$ be the transposition of legs in a fixed
diagram. Prove that $\tau_2$ is the identity. {\sl Hint}: (1) prove
that a ``hole'' can be dragged through a trivalent vertex
(2) to change the numbering of the two legs, use manoeuvres like in
Lemma~\ref{lemma2a} with $G_2=\emptyset$).

\item${}^*$
Let $\G$ be the algebra of 3-graphs.
   \begin{itemize}
\item
  Is it true that $\G$ is generated by plane graphs?
\item
  Find generators and relations of the algebra $\G$.
\item
  Suppose that a graph $G\in\G$ consists of two parts $G_1$ and $G_2$
  connected by three edges. Is the following equality:
$$\risS{-10}{g3_ex_g3}{\put(7,11){\mbox{$\scriptstyle G_1$}}
                   \put(50,11){\mbox{$\scriptstyle G_2$}}}{65}{18}{20} \quad=\quad
  \risS{-10}{g3_ex_g4}{\put(7,11){\mbox{$\scriptstyle G_1$}}
                   \put(50,11){\mbox{$\scriptstyle G_2$}}}{65}{0}{22}
$$
true?
\end{itemize}

\item${}^*$
Is it true that the algebra of primitive elements $\PR$ has no
divisors of zero with respect to the product $*$?

\item
Let $\X^{k}$ be the space of 1- and 3-valent graphs with $k$
numbered legs. Consider the transposition of two legs of an element of $\X^{k}$.
\begin{itemize}
  \item Give a example of a non-zero element of $\X^{k}$ with even $k$ which
is changed under such a transposition.
  \item${}^*$ Is it true that any such transposition
changes the sign of the element if $k$ is odd?
(The first nontrivial case is when $k=3$ --- this is Conjecture \ref{nechet_X3}.)
\end{itemize}

\item${}^*$
Let $\Lambda$ be Vogel's algebra, that is, the subspace of $\X^{3}$
consisting of all antisymmetric elements.
   \begin{itemize}
\item
   Is it true that $\Lambda=\X^{3}$ (this is again Conjecture \ref{nechet_X3})?
\item
   Is it true that $\Lambda$ is generated by the elements $t$ and $x_k$
(this is the Conjecture \ref{vog_conj}; see also Exercises \ref{t_x_etc}
and \ref{dodec_etc})?
   \end{itemize}

\item\label{t_x_etc}
Let $t$, $x_3$, $x_4$, $x_5$, ...  be the elements of the space
$\X^{3}$ defined above.
   \begin{itemize}
\item
   Prove that $x_i$'s belong to Vogel's algebra $\Lambda$, that is,
   that they are antisymmetric with respect to permutations of legs.
\item
   Prove the relation $x_4=-\frac{4}{3}t\vee x_3-\frac{1}{3}t^{\vee 4}$.
\item
   Prove that $x_k$ with an arbitrary even $k$ can be expressed through
   $t$, $x_3$, $x_5$, ...
   \end{itemize}

\item\label{dodec_etc}
Prove that the dodecahedron
$$d\quad=\quad\risS{-20}{vog-dodec}{\put(-3,-6){\mbox{$\scriptstyle 1$}}
          \put(61,-6){\mbox{$\scriptstyle 2$}}
          \put(23,50){\mbox{$\scriptstyle 3$}}}{60}{10}{20}$$
belongs to $\Lambda$, and express it as a vertex polynomial in $t$, $x_3$, $x_5$, $x_7$, $x_9$.

\item${}^*$
The group $S_3$ acts in the space of fixed diagrams with 3 legs
$\X^{3}$, splitting it into 3 subspaces:
   \begin{itemize}
\item
  symmetric, which is isomorphic to $\B^{3}$
(open diagrams with 3 legs),
\item
  totally antisymmetric, which is Vogel's $\Lambda$ by definition, and
\item
  some subspace $Q$, corresponding to a 2-dimensional irreducible
     representation of $S_3$.
   \end{itemize}
Is it true that $Q=0$?

\item\label{ex_gl_N_on_G_root}
Show that $N=0$, $N=-1$, and $N=1$ are roots of the polynomial
$\f_{\gl_N}(G)$ for any 3-graph $G\in \G_n$ ($n>1$).

\item\label{ex_so_N_on_G_root}
Show that $N=0$, $N=1$ and $N=2$ are roots of polynomial
$\f_{\so_N}(G)$ for any 3-graph $G\in \G_n$ ($n>0$).

\end{enumerate}
\end{xcb}
 %7 LieAlg

%\part{The Kontsevich Integral}\label{part_ki}
\chapter{The Kontsevich integral} % 08
\label{chapKI}

The Kontsevich integral appeared in the paper \cite{Kon1} by
M.~Kontsevich as a tool to prove the Fundamental Theorem of the
theory of Vassiliev invariants (that is, Theorem~\ref{fund_thm}).
Any Vassiliev knot invariant with coefficients in a field of
characteristic 0 can be factored through the universal invariant
defined by the Kontsevich integral.

Detailed (and different) expositions of the construction and
properties of the Kontsevich integral can be found in \cite{BN1, CD3, Les}.
Other important references are \cite{Car1}, \cite{LM1},
\cite{LM2}.

About the notation: in this chapter we shall think of $\R^3$ as the
product of a (horizontal) complex plane $\C$ with the complex
coordinate $z$ and a (vertical) real line $\R$ with the coordinate
$t$. All Vassiliev invariants are always thought of having values in
the complex numbers.

\section{First examples}
\label{link_num_formula}

\parbox[t]{3.4in}{We start with two examples where the Kontsevich integral
appears in a simplified form and with a clear geometric meaning.}
\vspace*{-10pt}
\subsection{The braiding number of a 2-braid}\ \\
\parbox[t]{3.4in}{\hspace*{20pt}A braid on two strands has a complete
invariant: the number of full twists that one strand makes around
the other.

\hspace*{20pt}Let us consider the horizontal coordinates of
}\hspace{10pt}
\parbox{1in}{$\risS{-30}{2braid}{\put(15,1){$\scriptstyle z(t)$}
    \put(48,28){$\scriptstyle w(t)$}\put(80,22){$\C$}
    \put(43,110){$t$}}{100}{0}{0}$}\vspace{3pt}\\
points on the strands, $z(t)$ and $w(t)$, as functions of the vertical coordinate $t$, $0\le t\le 1$, then the number of
full twists can be computed by the integral formula
$$
  \frac{1}{2\pi i}\int_0^1\frac{dz-dw}{z-w}.
$$
Note that the number of full twists is not necessarily an integer;
however, the number of {\em half-twists} always is.

\subsection{Kontsevich type formula for the linking number}
\label{subsec:linking_number}

The Gauss integral formula for the linking number of two spatial curves
$lk(K,L)$ (discussed in Section~\ref{gauss_f}) involves integration
over a torus (namely, the product of the two curves).
Here we shall give a different integral formula for the same invariant,
with the integration over an interval, rather than a torus. This formula
generalizes the expression for the braiding number of a braid on two strands
and, as we shall later see, gives the first term of the Kontsevich integral of a
two-component link.

\begin{xdefinition} A link in $\R^3$ is a {\em Morse} link \index{Link!Morse}
if the function $t$ (the vertical coordinate) on it has only non-degenerate critical
points. A Morse link is a {\em strict} Morse link \index{Link!strict Morse} if the
critical values of the vertical coordinate are all distinct. Similarly one speaks of
{\em Morse tangles} and {\em strict Morse tangles}.
\end{xdefinition}

\begin{xtheorem}\label{th:ki-lk}
Suppose that two disjoint connected curves $K$, $L$ are embedded into $\R^3$ as a strict Morse link.
$$\risS{-30}{lk_num}{\put(-24,40){$z_j(t)$}
    \put(68,40){$w_j(t)$}}{65}{60}{30}
$$
Then
$$\index{Linking number}
  lk(K,L)=\frac{1}{2\pi i}
          \int \sum_j (-1)^{\downarrow_j}
          \frac{d(z_j(t)-w_j(t))}{z_j(t)-w_j(t)},
$$
where the index  $j$ enumerates all possible choices
of a pair of strands on the link as functions $z_j(t)$, $w_j(t)$
corresponding
to $K$ and $L$, respectively, and the integer $\downarrow_j$ is the
number of strands in the pair which are oriented downwards.
\end{xtheorem}

\begin{xremark} In fact, the condition that the link in question is a
strict Morse link can be relaxed. One may consider piecewise linear
links with no horizontal segments, or smooth links whose vertical
coordinate function has no flattening points (those where all the
derivatives vanish).
\end{xremark}

\begin{proof}
The proof consists of three steps which --- in a more elaborate setting ---
will also appear in the full construction of the Kontsevich integral.

{\bf Step 1.} {\em The value of the sum in the right hand side is an integer.}
Note that for a strict Morse link with two components $K$ and $L$, the
{\em configuration space} of all horizontal chords joining $K$ and $L$ is a
closed one-dimensional manifold, that is, a disjoint union of several circles.

For example, assume that two adjacent critical values $m$ and $M$ (with $m<M$) of the vertical
coordinate correspond to a minimum on the component $K$ and a maximum on the component $L$
respectively:
\begin{center}
\includegraphics[height=35mm]{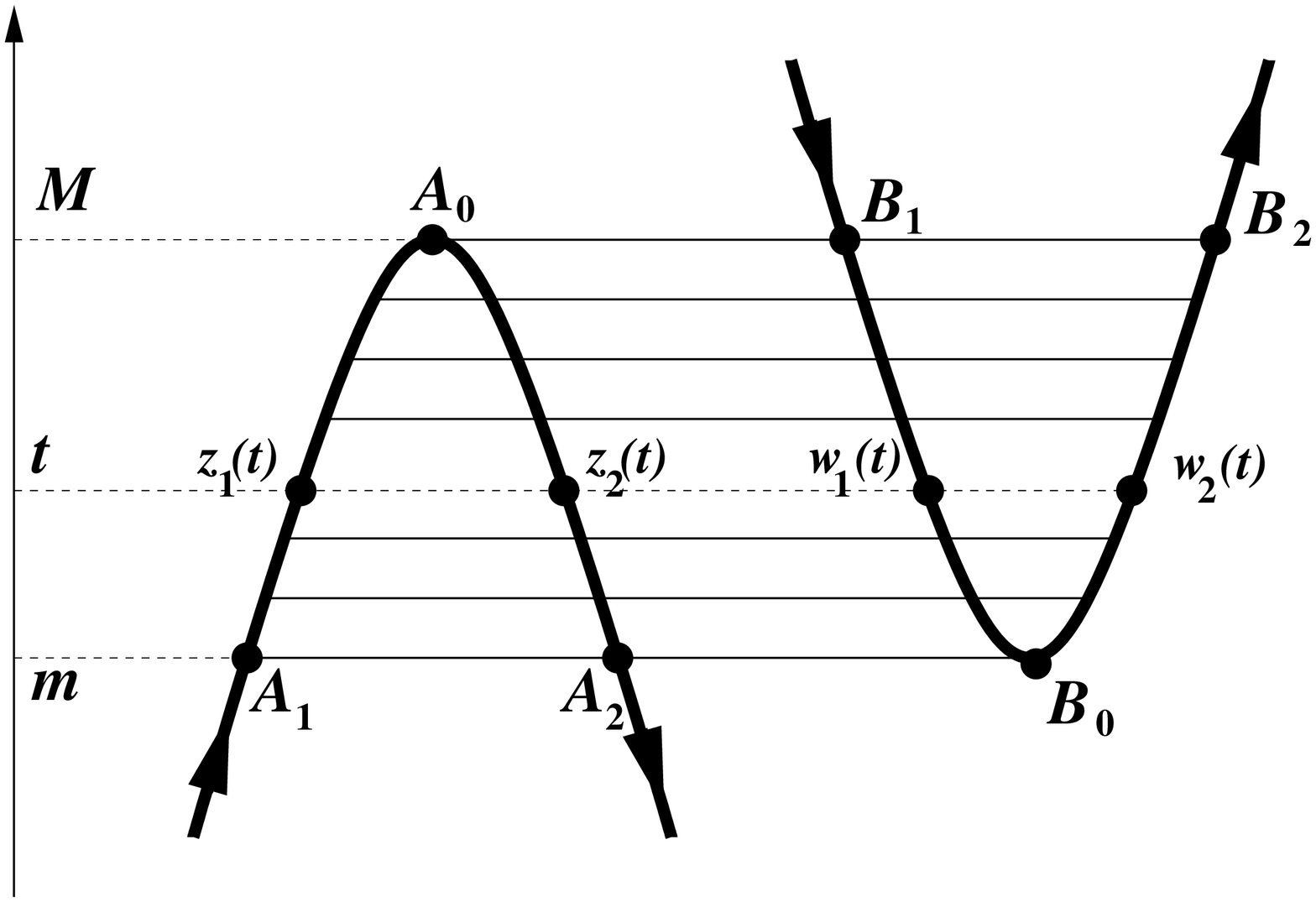}
\end{center}
The space of all horizontal chords that join the shown parts of $K$
and $L$ consists of four intervals which join together to form a
circle. The motion along this circle starts, say, at a chord
$A_1B_0$  and proceeds as $$A_1B_0 \to A_0B_1 \to A_2B_0 \to A_0B_2
\to A_1B_0.$$ Note that when the moving chord passes a critical
level (either $m$ or $M$), the direction of its motion changes, and
so does the sign $(-1)^{\downarrow_j}$. (Exercise
(\ref{pr:hor-hor-conf-sp}) on page~\pageref{pr:hor-hor-conf-sp}
deals with a more complicated example of the  configuration space of
horizontal chords.)

It is now clear that our integral formula counts the number of
complete turns made by the horizontal chord while running through
the whole configuration space of chords with one end $(z_j(t),t)$ on
$K$ and the other end $(w_j(t),t)$ on $L$. This is, clearly, an
integer.

{\bf Step 2.} {\em The value of the right hand side remains
unchanged under a continuous horizontal deformation of the link.}
(By a horizontal deformation we mean a deformation of a link which
moves every point in a horizontal plane $t=\const$.) The assertion
is evident, since the integral changes continuously while always
remaining an integer. Note that this is true even if we allow
self-intersections within each of the components; this does not
influence the integral because $z_j(t)$ and $w_j(t)$ lie on the
different components.

{\bf Step 3.} {\em Reduction to the combinatorial formula for the
linking number (Section \ref{link_num})}. Choose a vertical plane in
$\R^3$ and represent the link by a generic projection to that plane.
By a horizontal deformation, we can flatten the link so that it lies
in the plane completely, save for the small fragments around the
diagram crossings between $K$ and $L$ (as we noted above, {\em
self-intersections} of each component are allowed). Now, the
rotation of the horizontal chord for each crossing is by $\pm\pi$,
and the signs are in agreement with the number of strands oriented
downwards. The reader is invited to draw the two different possible
crossings, then, for each picture, consider the four possibilities
for the orientations of the strands and make sure that the sign of
the half-turn of the moving horizontal chord always agrees with the
factor $(-1)^{\downarrow_j}$. (Note that the integral in the theorem
is computed over $t$, so that each specific term computes the angle
of rotation of the chord as it moves {from bottom to top}.)
\end{proof}

The Kontsevich integral can be regarded as a generalization of
this formula. Here we kept track of one horizontal chord
moving along the two curves. The full Kontsevich integral
keeps track of how {\em finite sets} of horizontal chords on the knot
(or a tangle) rotate when moved in the vertical direction.
This is the somewhat na{\"\i}ve approach that we use in the next section.
Later, in Section~\ref{KZ}, we shall adopt a more sophisticated point of view,
interpreting the Kontsevich integral as the monodromy of the
Knizhnik--Zamolodchikov connection in the complement
of the union of diagonals in $\C^n$.

\section{The construction}
\label{defki}

Let us recall some notation and terminology of the preceding section.
For points of $\R^3$ we use coordinates $(z,t)$ with $z$ complex and $t$ real;
the planes $t=\mathrm{const}$ are thought of being horizontal. Having chosen the
coordinates, we can speak of {\em strict Morse knots}, \index{Knot!strict Morse}
namely, knots with the property that the coordinate $t$ restricted to the knot
has only non-degenerate critical points with distinct critical values.

We define the Kontsevich integral for strict Morse knots.
Its values belong to the graded completion $\widehat{\A}$ \label{hat-A}
of the algebra of chord diagrams with 1-term relations $\A=\A^{fr}/(\Theta)$.
(By definition, the elements of a graded algebra are finite linear
combinations of homogeneous elements.
The {\em graded completion}\index{Graded completion}
consists of all infinite combinations of such elements.)

\begin{definition}
\label{formulaKI} \index{Kontsevich integral} %
The  {\em Kontsevich integral} $Z(K)$ of a strict Morse knot $K$ is
given by the following formula:
$$\label{ki-Z(K)}
  Z(K) = \sum_{m=0}^\infty \frac{1}{(2\pi i)^m}
         \int\limits_{\substack{
    t_{\mbox{\tiny min}}<t_m<\dots<t_1<t_{\mbox{\tiny max}}\\
    t_j\mbox{\tiny\ are noncritical}}}
         \sum_{P=\{(z_j,z'_j)\}} (-1)^{\downarrow_P} D_P
         \bigwedge_{j=1}^m \frac{dz_j-dz'_j}{z_j-z'_j} \ .
$$
\end{definition}

The ingredients of this formula have the following meaning.
\medskip

The real numbers $t_{\mbox{\scriptsize min}}$ and
$t_{\mbox{\scriptsize max}}$ are the minimum and the maximum of
the function $t$ on $K$.

\medskip
The integration domain is the set of all points of the $m$-dimensional simplex
$t_{\mbox{\scriptsize min}}<t_m<\dots<t_1<t_{\mbox{\scriptsize max}}$
none of whose coordinates $t_i$ is a critical value of $t$. The $m$-simplex is
divided by the critical values into several {\em connected components}. \label{conncomp}
For example, for the following embedding of the unknot and $m=2$ the corresponding
integration domain has six connected components and looks like
 \begin{center}
     \begin{picture}(280,95)(0,0)
      \put(0,0){\ig[width=280pt]{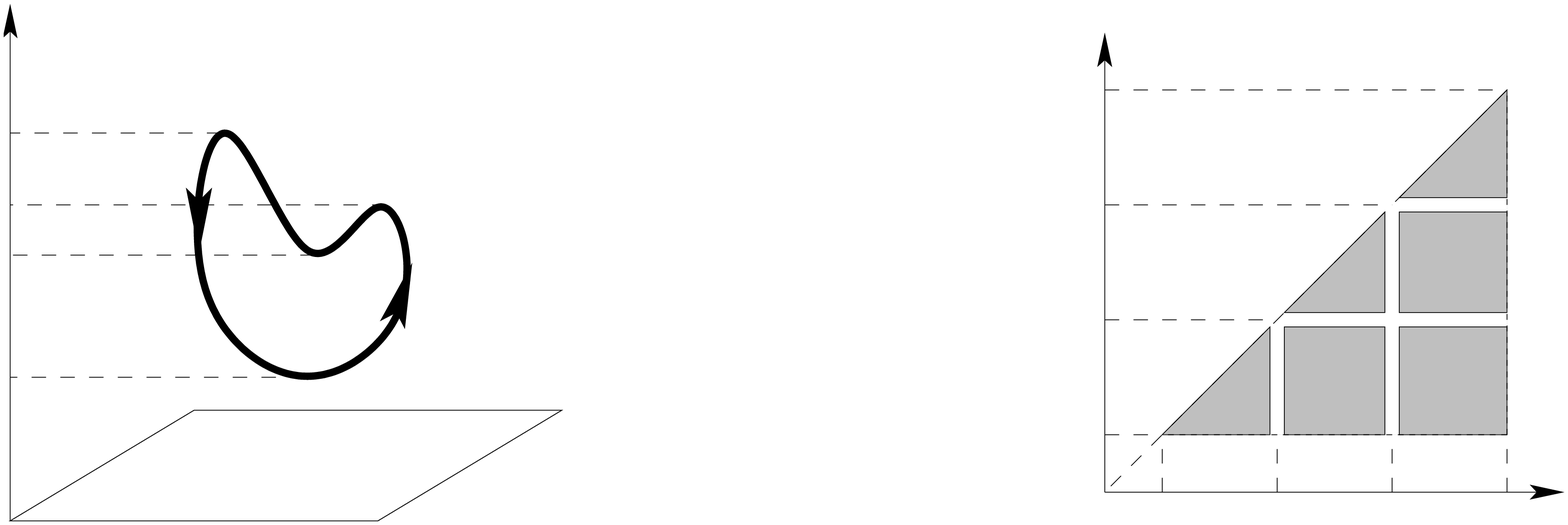}}
      \put(-7,90){\mbox{$t$}}
      \put(-20,71){\mbox{$t_{\mbox{\tiny max}}$}}
      \put(-14,57){\mbox{$t_{c_1}$}}
      \put(-14,42){\mbox{$t_{c_2}$}}
      \put(-20,22){\mbox{$t_{\mbox{\tiny min}}$}}
      \put(77,0){\mbox{$z$}}
      \put(187,92){\mbox{$t_2$}}
      \put(177,74){\mbox{$t_{\mbox{\tiny max}}$}}
      \put(183,57){\mbox{$t_{c_1}$}}
      \put(183,36){\mbox{$t_{c_2}$}}
      \put(177,15){\mbox{$t_{\mbox{\tiny min}}$}}
      \put(285,2){\mbox{$t_1$}}
      \put(266,-5){\mbox{$t_{\mbox{\tiny max}}$}}
      \put(246,-5){\mbox{$t_{c_1}$}}
      \put(225,-5){\mbox{$t_{c_2}$}}
      \put(203,-5){\mbox{$t_{\mbox{\tiny min}}$}}
     \end{picture}
  \end{center}

\medskip
The number of summands in the integrand is constant in each connected
component of the integration domain, but can be different for different
components.
In each plane $\{t=t_j\}\subset\R^3$ choose an unordered pair of distinct
points $(z_j,t_j)$ and $(z'_j,t_j)$ on $K$, so that $z_j(t_j)$
and $z'_j(t_j)$ are continuous functions.
We denote by $P=\{(z_j,z'_j)\}$ the set of such pairs for $j=1,\dots,m$ and call it a
{\em pairing}.\index{Pairing}

\begin{center}\label{summands}
  \ig[width=12cm]{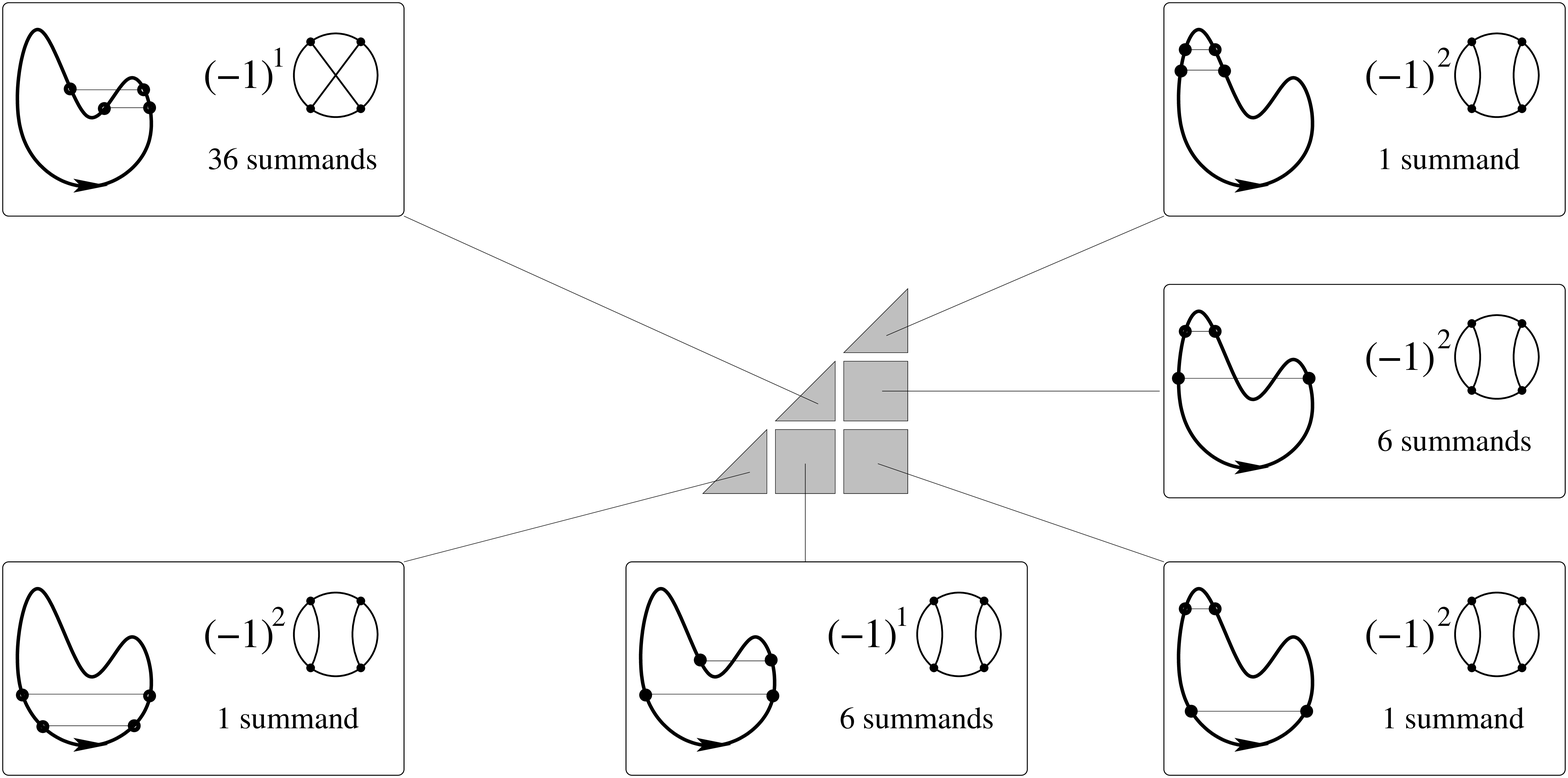}
\end{center}

The integrand is the sum over all choices of the pairing $P$.
In the example above for the component
$\{t_{c_1}<t_1<t_{\mbox{\scriptsize max}}$,
$t_{\mbox{\scriptsize min}}<t_2<t_{c_2}\}$, in the bottom right corner,
we have only one possible pair of points on the levels $\{t=t_1\}$ and
$\{t=t_2\}$. Therefore, the sum over $P$ for this component consists of
only one summand. In contrast, in the component next to it,
$\{t_{c_2}<t_1<t_{c_1}$, $t_{\mbox{\scriptsize min}}<t_2<t_{c_2}\}$,
we still have only one possibility for the chord $(z_2,z'_2)$ on the level $\{t=t_2\}$, but the plane
$\{t=t_1\}$ intersects our knot $K$ in four points. So we have
$\binom{4}{2}=6$ possible pairs $(z_1,z'_1)$ and the total number
of summands here is six (see the picture above).

\medskip
For a pairing $P$ the symbol `$\downarrow_P$' denotes the number of points
$(z_j,t_j)$ or $(z'_j,t_j)$ in $P$ where the coordinate $t$ decreases
as one goes along $K$.

\medskip
Fix a pairing $P$. Consider the knot $K$ as an oriented circle and connect
the points $(z_j,t_j)$ and $(z'_j,t_j)$ by a chord. We obtain a chord diagram
with $m$ chords. (Thus, intuitively, one can think of a pairing as a way of
inscribing a chord diagram into a knot in such a way that all chords are
horizontal and are placed on different levels.)
The corresponding element of the algebra $\A$ is denoted by $D_P$.
In the picture below, for each connected component in our example, we show
one of the possible pairings, the corresponding chord diagram with the sign
$(-1)^{\downarrow_P}$ and the number of summands of the integrand
(some of which are equal to zero in $\A$ due to the one-term relation).

Over each connected component, $z_j$ and $z'_j$ are smooth functions
in $t_j$. By
$$
{\displaystyle \bigwedge_{j=1}^m \frac{dz_j-dz'_j}{z_j-z'_j}}
$$ 
we mean the
{pullback} of this form to the integration domain of the variables
$t_1,\dots,t_m$. The integration domain is considered with the
{orientation} of the space $\R^m$ defined by the natural order of the
coordinates $t_1, \ldots, t_m$.

By convention, the term in the Kontsevich integral corresponding to $m=0$
is the (only) chord diagram of order 0 taken with coefficient one. It is
the unit of the algebra $\widehat{\A}$.

\subsection{Basic properties}

We shall see later in this chapter that the Kontsevich integral has the
following basic properties:
\begin{itemize}
\item
   $Z(K)$ converges for any strict Morse knot $K$.
\item
It is invariant under the deformations of the knot in the class
of (not necessarily strict) Morse knots.
\item
It behaves in a predictable way under the deformations that add a pair
of new critical points to a Morse knot.
\end{itemize}

Let us explain the last item in more detail. While the Kontsevich
integral is indeed an invariant of Morse knots, it is not preserved by
deformations that change the number of critical points of $t$.
However, the following formula shows how the integral changes
when a new pair of critical points is added to the knot:
\begin{equation}\label{addcritF}
  Z\biggl(\,\rb{-12pt}{\ig[width=30pt]{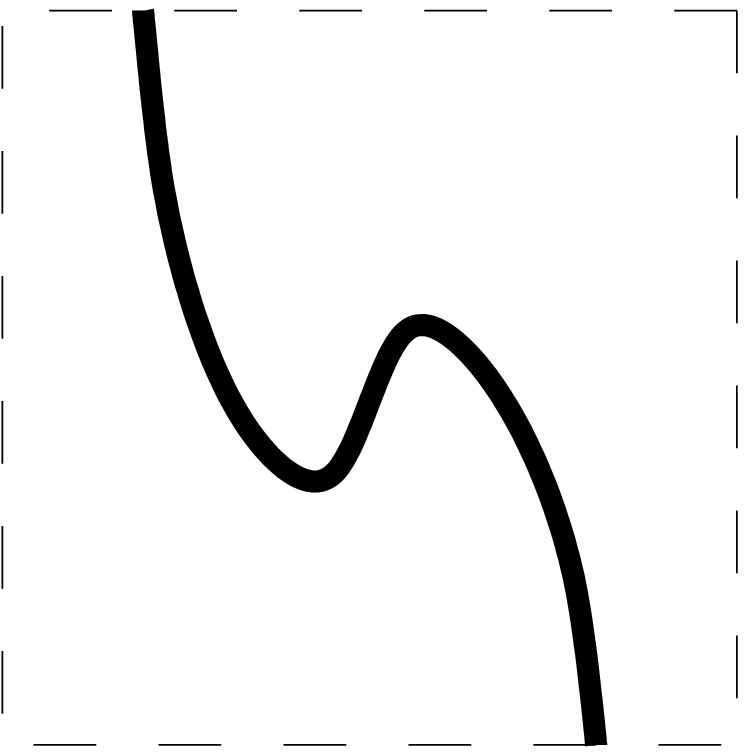}}\,\biggr)
  = Z(H) \cdot
  Z\biggl(\,\rb{-12pt}{\ig[width=30pt]{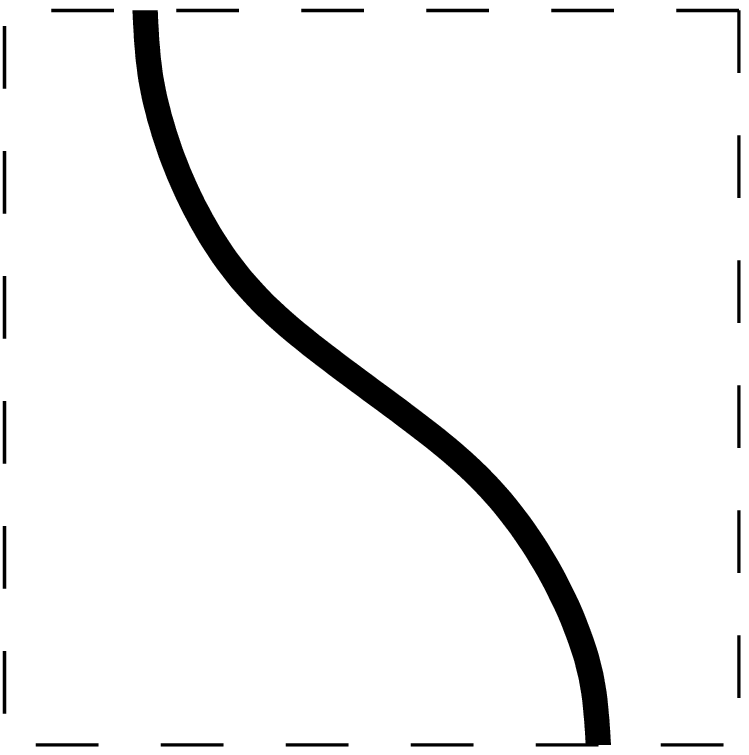}}\,\biggr).
\end{equation}
Here the first and the third pictures represent two embeddings of an
arbitrary knot that coincide outside the fragment shown,
$$
  H\ :=\ \risS{-12}{humpk}{}{25}{15}{10}
$$
is the {\em hump}\index{Hump}\label{hump}
(an unknot with two maxima),
and the product is the product in the completed algebra $\widehat\A$
of chord diagrams. The equality (\ref{addcritF}) allows to define a
genuine knot invariant by the formula
$$\label{ki-I(K)}
  I(K) = \frac{Z(K)}{Z(H)^{c/2}},
$$
where $c$ denotes the number of critical points of $K$ and the ratio
means the division in the algebra $\widehat\A$ according to the rule
$(1+a)^{-1} = 1-a+a^2-a^3+\dots$
The knot invariant $I(K)$ is sometimes referred to as the {\em final} Kontsevich
\index{Kontsevich integral!final}%
integral as opposed to the {\em preliminary} Kontsevich integral
\index{Kontsevich integral!preliminary}%
$Z(K)$.

The central importance of the final Kontsevich integral in the theory of finite type
invariants is that it is a universal Vassiliev invariant in the following sense.

Consider an unframed weight system $w$ of degree $n$ (that is, a function on the set
of chord diagrams with $m$ chords satisfying one- and four-term relations).
Applying $w$ to the $m$-homogeneous part of the series $I(K)$,
we get a numerical knot invariant $w(I(K))$. This invariant is a Vassiliev
invariant of order $m$ and such invariants span the space of all finite type
invariants. This argument
will be used to prove the Fundamental Theorem on Vassiliev Invariants, see
Section~\ref{pkt}.

The Kontsevich integral has many interesting properties that we
shall describe in this and in the subsequent chapters. Among these
are its behaviour with respect to the connected sum of knots
(Section~\ref{ki_tangles} and \ref{long2closed}) to the coproduct in
the Hopf algebra of chord diagrams (Section~\ref{glki}), cablings
(Chapter~\ref{chap:operations}), mutation
(Section~\ref{ki_mutation}). We shall see that it can be computed
combinatorially (Section~\ref{comb_ki}) and has rational
coefficients (Section~\ref{kifieqki}).

\section{Example of calculation}
\label{excalc}

Here we shall calculate the coefficient of the chord diagram
$\risS{-3}{cd22-ki}{}{17}{0}{0}$ in $Z(H)$, where $H$ is the hump (plane curve with 4 critical points,
as in the previous section) directly from the definition of the Kontsevich integral.
The following computation is valid for an arbitrary shape of the curve,
provided that the length of the segments $a_1a_2$ and $a_3a_4$
(see picture below) decreases with $t_1$, while that of the segment
$a_2a_3$ increases.

First of all, note that out of the total number of 51 pairings shown
in the picture on page~\pageref{summands}, the following 16 contribute
to the coefficient of $\risS{-6}{cd22-ki}{}{17}{0}{0}$:
\newcommand\tw[1]{\rb{-15pt}{\ig[width=30pt]{#1.eps}}}
\begin{gather*}
\tw{t44}\quad \tw{t43}\quad \tw{t42}\quad \tw{t41}
\quad\tw{t24}\quad \tw{t23}\quad \tw{t22}\quad \tw{t21} \\
\tw{t34}\quad \tw{t33}\quad \tw{t32}\quad \tw{t31}
\quad\tw{t14}\quad \tw{t13}\quad \tw{t12}\quad \tw{t11}
\end{gather*}

We are, therefore, interested only in the band between the critical values
$c_1$ and $c_2$. Denote by $a_1$, $a_2$, $a_3$, $a_4$ (resp.
$b_1$, $b_2$, $b_3$, $b_4$) the four points of intersection of the knot
with the level $\{t=t_1\}$ (respectively, $\{t=t_2\}$):
$$
  \ig[height=4cm]{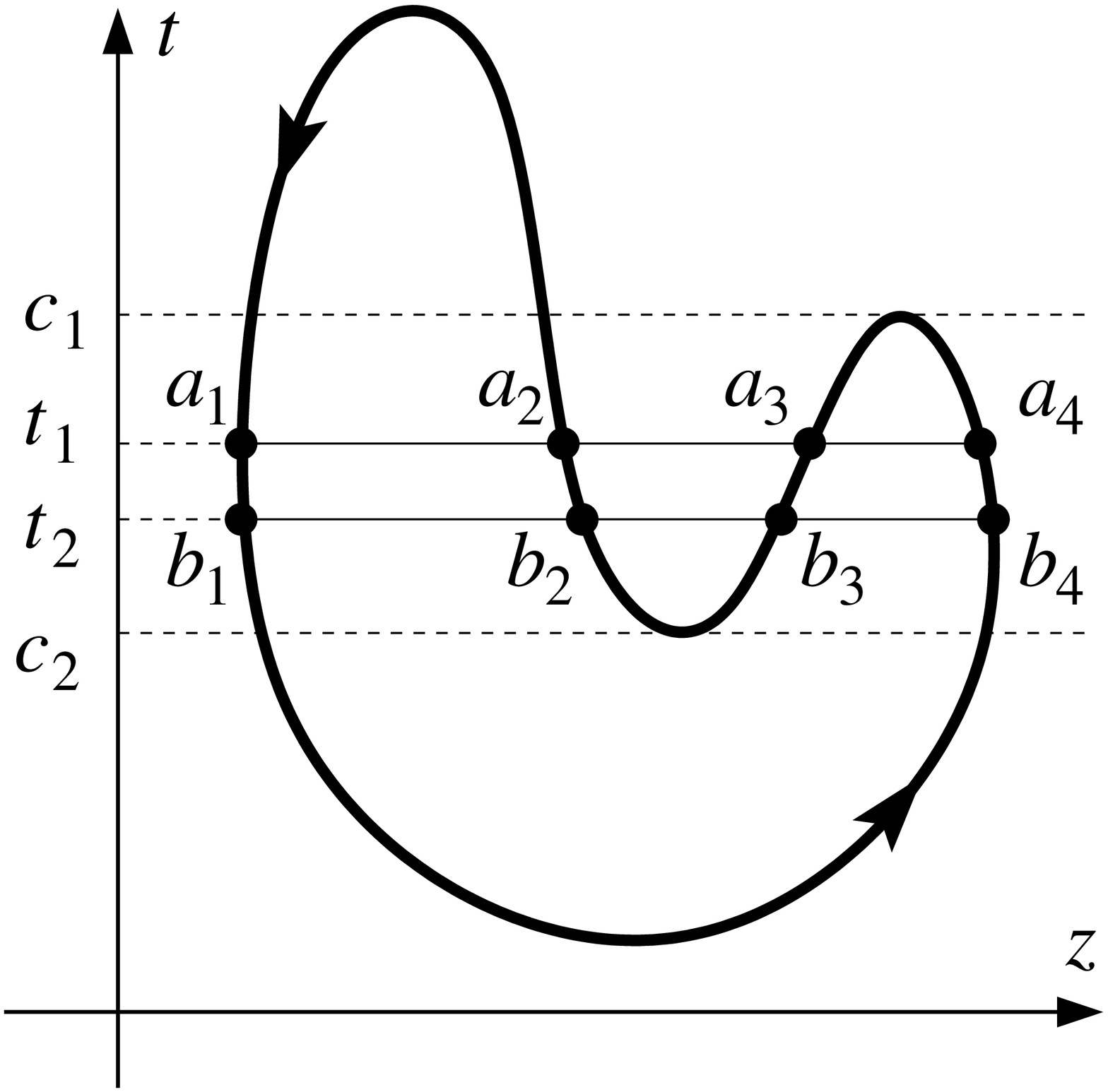}
$$
The sixteen pairings shown in the picture above correspond to
the differential forms
$$
  (-1)^{j+k+l+m}d\ln a_{jk} \wedge d\ln b_{lm},
$$
where $a_{jk}=a_k-a_j$, $b_{lm}=b_m-b_l$, and the pairs $(jk)$ and $(lm)$
can take 4 different values each: $(jk)\in\{(13),(23),(14),(24)\}=:A$,
$(lm)\in\{(12),(13),(24),(34)\}=:B$. The sign $(-1)^{j+k+l+m}$ is equal
to $(-1)^{\downarrow_P}$, because in our case the upward oriented strings have
even numbers, while the downward oriented strings have odd numbers.

The coefficient of $\risS{-5}{cd22-ki}{}{17}{0}{0}$ is, therefore, equal to
\begin{align*}
  &\frac{1}{(2\pi i)^2} \int\limits_\Delta
  \sum_{(jk)\in A}\sum_{(lm)\in B}(-1)^{j+k+l+m} d\ln a_{jk} \wedge
     d\ln b_{lm}\\
= &-\frac{1}{4\pi^2} \int\limits_\Delta
  \sum_{(jk)\in A}(-1)^{j+k+1} d\ln a_{jk}
  \wedge\sum_{(lm)\in B}(-1)^{l+m-1}  d\ln b_{lm}\\
= &-\frac{1}{4\pi^2} \int\limits_\Delta
   d\ln\frac {a_{14}a_{23}}{a_{13}a_{24}}
    \wedge d\ln \frac {b_{12}b_{34}}{b_{13}b_{24}},
\end{align*}
where the integration domain $\Delta$ is the triangle described by the
inequalities $c_2<t_1<c_1$, $c_2<t_2<t_1$. Assume the following notation:
$$
  u=\frac {a_{14}a_{23}}{a_{13}a_{24}},\qquad
  v=\frac {b_{12}b_{34}}{b_{13}b_{24}}.
$$
It is easy to see that $u$ is an increasing function of $t_1$ ranging from 0
to 1, while $v$ is an decreasing function of $t_2$ ranging from 1 to 0.
Therefore, the mapping $(t_1,t_2)\mapsto(u,v)$ is a diffeomorphism with
a negative Jacobian, and after the change of variables the integral we are
computing becomes
$$
 \frac{1}{4\pi^2} \int\limits_{\Delta'}
   d\ln u \wedge d\ln v
$$
where $\Delta'$ is the image of $\Delta$. It is obvious that the boundary
of $\Delta'$ contains the segments $u=1, 0\le v\le 1$ and $v=1, 0\le u\le 1$
that correspond to $t_1=c_1$ and $t_2=c_2$. What is not immediately evident
is that the third side of the triangle $\Delta$ also goes into a straight
line, namely, $u+v=1$. Indeed, if $t_1=t_2$, then all $b$'s are equal to the
corresponding $a$'s and the required fact follows from the identity
$a_{12}a_{34}+a_{14}a_{23}=a_{13}a_{24}$.
$$
  \ig[height=2cm]{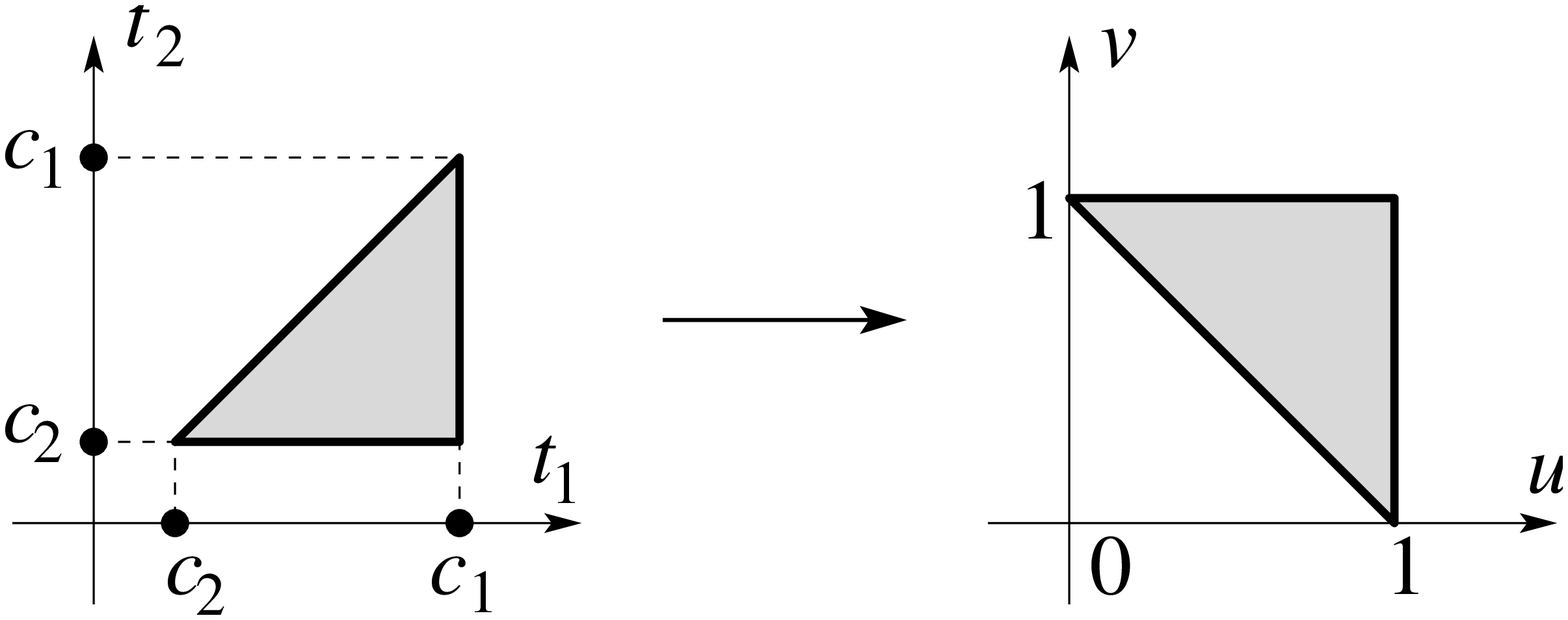}
$$

Therefore,
\begin{eqnarray*}
 \frac{1}{4\pi^2} \int\limits_{\Delta'}
   d\ln u \wedge d\ln v &=& \frac{1}{4\pi^2}
\int\limits_0^1\left(\ \int\limits_{1-u}^1 d\ln v\right)
                          \frac{du}{u}\\
 &=&
-\frac{1}{4\pi^2} \int\limits_0^1\ln(1-u)\frac{du}{u}\ .
\end{eqnarray*}

Taking the Taylor expansion of the logarithm we get
$$
\frac{1}{4\pi^2}\sum_{k=1}^\infty\int\limits_0^1\frac{u^k}{k}\frac{du}{u}
\ =\  \frac{1}{4\pi^2} \sum_{k=1}^\infty \frac{1}{k^2}
\ =\  \frac{1}{4\pi^2} \zeta(2)\ =\  \frac{1}{24}.
$$

Two things are quite remarkable in this answer: (1) that it is
expressed via a value of the zeta function, and (2) that the answer
is rational. In fact, for any knot $K$ the coefficient of any chord
diagram in $Z(K)$ is rational and can be computed through the values
of multivariate $\zeta$-functions:
$$\zeta(a_1,\dots,a_n)\ =\
     \sum_{0<k_1<k_2<\dots<k_n} k_1^{-a_1}\dots k_n^{-a_n}.
$$
We shall speak about that in more detail in Section \ref{zetanumb}.

For a complete formula for $Z(H)$ see Section~\ref{wheels}.

\section{The Kontsevich integral for tangles}
\label{ki_tangles}\index{Kontsevich integral!for tangles}

The definition of the preliminary Kontsevich integral for knots (see
Section~\ref{defki}) makes sense for an arbitrary strict Morse
tangle $T$. One only needs to replace the completed algebra
$\widehat\A$ of chord diagrams by the graded completion of the
vector space of tangle chord diagrams on the skeleton of $T$, and
take $t_{\mbox{\sz min}}$ and $t_{\mbox{\sz max}}$ to correspond to
the bottom and the top of $T$, respectively. In the section
\ref{s:conver-ki} we shall show that the coefficients of the chord
diagrams in the Kontsevich integral of any (strict Morse) tangle
actually converge.

In particular, one can speak of the Kontsevich integral of links or
braids.

\begin{xca}
For a two-component link, what is the coefficient in the Kontsevich integral
of the chord diagram of degree 1 whose chord has ends on both components?
\end{xca}
{\sl Hint:} see Section~\ref{subsec:linking_number}.

\begin{xca}\label{comp_ki_R}
Compute the integrals
$$Z\left(\ \risS{-12}{crt}{}{30}{20}{15}\ \right)\qquad\mbox{and}\qquad
Z\left(\ \risS{-12}{crtm}{}{30}{20}{15}\ \right).$$ \label{ki-R}
\end{xca}

{\sl Answer:}
$$\rb{-3mm}{\ig[height=9mm]{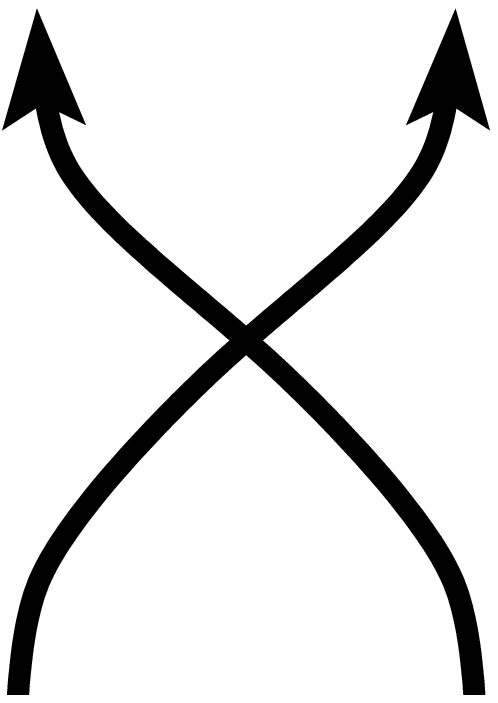}}
\cdot\exp\Bigl(\frac{\risS{-3}{doc}{}{17}{0}{0}}{2}\Bigr)
             \quad\mbox{and}\quad
\rb{-3mm}{\ig[height=9mm]{dzcvirt.eps}} \cdot\exp\Bigl(-
\frac{\risS{-3}{doc}{}{17}{0}{0}}{2}\Bigr),\quad\mbox{respectively,}
$$
where $\exp a$ is the series $1+a+\frac{a^2}{2!}+\frac{a^3}{3!}+\ldots$.

Strictly speaking, before describing the properties of the Kontsevich integral
we need to show that it is always well-defined. This will be done in the following
section. Meanwhile, we shall assume that this is indeed the case for all the tangles
in question.

\begin{proposition}\label{multTKI}
The Kontsevich integral for tangles is multiplicative:
$$
  Z(T_1\cdot T_2) = Z(T_1)\cdot Z(T_2)
$$
whenever the product $T_1\cdot T_2$ is defined.
\end{proposition}

\begin{proof}
Let $t_{\mbox{\sz min}}$ and $t_{\mbox{\sz max}}$ correspond to the
bottom and the top of $T_1\cdot T_2$, respectively, and let
$t_{\mbox{\sz mid}}$ be the level of the top of $T_2$ (or the bottom
of $T_1$, which is the same). In the expression for the Kontsevich
integral of the tangle $T_1\cdot T_2$ let us remove from the domain
of integration all points with at least one coordinate $t$ equal to
$t_{\mbox{\sz mid}}$. This set is of codimension one, so the value
of integral remains unchanged. On the other hand, the connected
components of the new domain of integration are precisely all
products of the connected components for $T_1$ and $T_2$, and the
integrand for $T_1\cdot T_2$ is the exterior product of the
integrands for $T_1$ and $T_2$. The Fubini theorem on multiple
integrals implies that $Z(T_1\cdot T_2) = Z(T_1)\cdot Z(T_2)$.
\end{proof}

The behaviour of the Kontsevich integral under the tensor product of
tangles is more complicated. In the expression for $Z(T_1\ot T_2)$
indeed there are terms that add up to the tensor product $Z(T_1)\ot
Z(T_2)$: they involve pairings without chords that connect $T_1$
with $T_2$. However, the terms with pairings that do have such
chords are not necessarily zero and we have no effective way of
describing them. Still, there is something we can say but we need a
new definition for this.

\subsection{Parametrized tensor products}
By a (horizontal) {\em $\e$-rescaling} of $\R^3$ we mean the map
sending $(z,t)$ to $(\e z, t)$. For $\e>0$ it induces an operation
on tangles; we denote by $\e T$ the result of an $\e$-rescaling
applied to $T$. Note that $\e$-rescaling of a tangle does not change
its Kontsevich integral.

Let $T_1$ and $T_2$ be two tangles such that $T_1\ot T_2$ is
defined. For $0<\e\leq 1$ we define the {\em
$\e$-parametrized tensor product}
$T_1\ot_{\e} T_2$ as the result of placing $\e
T_1$ next to $\e T_2$ on the left, with the distance of
$1-\e$ between the two tangles:

$$
T_1=\tanG{-15}{tang2sub}{}{59}; \quad T_2=\
\tanG{-15}{tang1sub}{}{49};\quad T_1\ot_{\e}
T_2=\tanG{-15}{tang4sub}{}{80}\ .
$$
\index{Tangle!tensor product!parametrized} \label{eparam}
\medskip

More precisely, let $\mathbf{0}_{1-\e}$ be the empty tangle
of width $1-\e$ and the same height and depth as
$\e T_1$ and $\e T_2$. Then
$$T_1\ot_{\e} T_2= \e T_1\ot \mathbf{0}_{1-\e}\ot\e T_2.$$
When $\e=1$ we get the usual tensor product. Note that when
$\e<1$, the parametrized tensor product is, in general,
not associative.

\begin{proposition}\label{hormult}
The Kontsevich integral for tangles is asymptotically multiplicative with respect
to the parametrized tensor product:
$$
  \lim_{\e\to 0} Z(T_1\ot_{\e} T_2) = Z(T_1)\ot Z(T_2)
$$
whenever the product $T_1\ot T_2$ is defined. Moreover, the
difference
$$Z(T_1\ot_{\e} T_2) - Z(T_1)\otimes Z(T_2)$$
as $\e$ tends to 0 is of the same or smaller order of magnitude as
$\e$.
\end{proposition}

\begin{proof}
As we have already noted before, $Z(T_1\ot_{\e} T_2)$ consists of two parts: the terms that do not involve
chords that connect $\e T_1$ with $\e T_2$, and the terms that do. The first part does not depend on
$\e$ and is equal to $Z(T_1)\otimes Z(T_2)$, and the second part tends to 0 as $\e\to 0$.

Indeed, each pairing $P=\{(z_j,z'_j)\}$ for $T_1\ot T_2$ give rise to a continuous family of pairings
$P_{\e}=\{(z_j(\e),z'_j(\e))\}$ for $T_1\otimes_{\e} T_2$.
Consider one such  family $P_{\e}$. For all $k$
$$ dz_{k}(\e)-dz'_{k}(\e)=\e (dz_{k}-dz'_{k}).$$
If the $k$th chord has has both ends on $\e T_1$ or on $\e T_2$, we have
$$z_{k}(\e)-z'_{k}(\e)=\e (z_{k}-z'_{k})$$
for all $\e$. Therefore the limit of the first part is equal to
$Z(T_1)\ot Z(T_2)$.

On the other hand, if $P_{\e}$ has at least one chord
connecting the two factors, we have
$|z_{k}(\e)-z'_{k}(\e)|\to 1$
as $\e\to 0$. Thus the integral corresponding to the
pairing $P_{\e}$ tends to zero as $\e$ gets smaller, and we
see that the whole second part of the Kontsevich integral of
$T_1\otimes_{\e} T_2$ vanishes in the limit at least as
fast as $\e$:
$$Z(T_1\ot_{\e} T_2)=Z(T_1)\ot Z(T_2)+O(\e)\ .$$
\end{proof}

\section{Convergence of the integral}\label{s:conver-ki}

\begin{proposition}\index{Kontsevich integral!convergence}
For any strict Morse tangle $T$, the Kontsevich integral $Z(T)$ converges.
\end{proposition}

\begin{proof}
The integrand of the Kontsevich integral may have singularities near the boundaries
of the connected components. This happens near a critical point of a tangle
when the pairing includes a ``short'' chord whose ends are on the branches of the
tangle that come together at a critical point.

Let us assume that the tangle $T$ has exactly one critical point. This is
sufficient since any strict Morse tangle can be decomposed as a product of
such tangles (and the case when there are no critical points at all, is
trivial).  The argument in the proof of Proposition~\ref{multTKI} shows
that the Kontsevich integral of a product converges whenever the integral of
the factors does.

Suppose, without loss of generality, that $T$ has a critical point which is a maximum %$c$
with the value $t_c$. Then we only need to consider pairings with no chords above
$t_c$. Indeed, for any pairing its coefficient in the Kontsevich integral of $T$ is
a product of two integrals: one corresponding to the chords above $t_c$, and the
other - to the chords below $t_c$. The first integral obviously converges since the integrand
has no singularities, so it is sufficient to consider the factor with chords below $t_c$.

Essentially, there are two cases.

1) An isolated chord $(z_1,z'_1)$ tends to zero:
$$\begin{picture}(50,22)(0,5)
  \put(0,0){\ig[width=40pt]{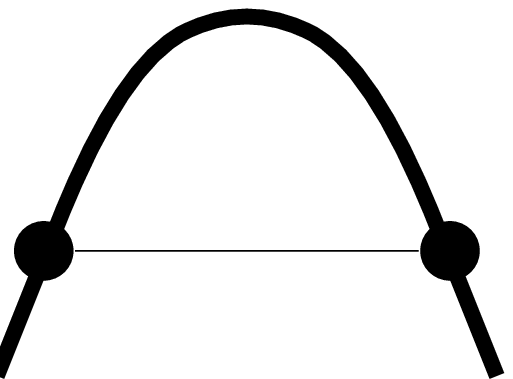}}
  \put(-8,10){\mbox{${\scriptstyle z_1}$}}
  \put(40,10){\mbox{${\scriptstyle z'_1}$}}
\end{picture}$$
In this case the corresponding chord diagram $D_P$
is equal to zero in $\A$ by the one-term relation.

2) A chord $(z_j,z'_j)$ tends to zero near a critical point but
is separated from that point by one or more other chords:

$$
\begin{picture}(200,35)(0,5)
  \put(0,0){\ig[width=200pt]{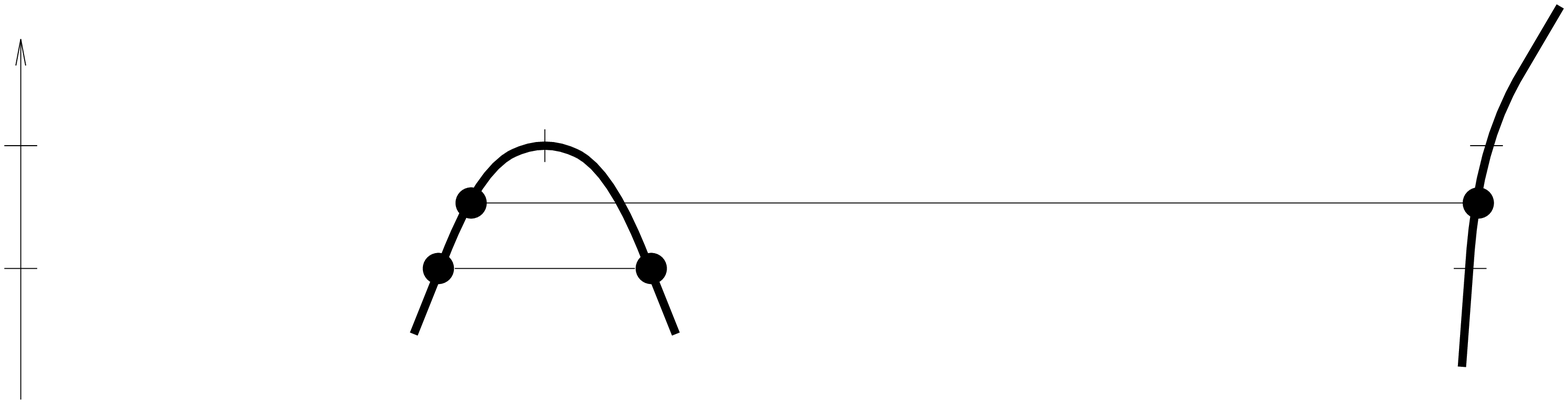}}
  \put(-8,32){\mbox{${\scriptstyle t_c}$}}
  \put(-8,14){\mbox{${\scriptstyle t_2}$}}
  \put(45,14){\mbox{${\scriptstyle z_2}$}}
  \put(88,14){\mbox{${\scriptstyle z'_2}$}}
  \put(47,25){\mbox{${\scriptstyle z_1}$}}
  \put(64,38){\mbox{${\scriptstyle z_c}$}}
  \put(198,35){\mbox{${\scriptstyle z'_c}$}}
  \put(195,25){\mbox{${\scriptstyle z'_{1}}$}}
  \put(191,11){\mbox{${\scriptstyle z''_2}$}}
\end{picture}
$$

Consider, for example, the case shown on the figure, where the ``short'' chord $(z_2,z'_2)$
is separated from the critical point by another, ``long'' chord
$(z_{1},z'_{1})$. We have:
\begin{eqnarray*}
  \left|\ \int\limits_{t_2}^{t_c}
         \frac{dz_{1}-dz'_{1}} {z_{1}-z'_{1}}\ \right|
&\leq& C\left|\ \int\limits_{t_2}^{t_c} d(z_{1}-z'_{1})\ \right|\\
&=& C\left|(z_c-z_2)-(z'_c-z''_2)\right|\ \leq\ C' |z_2-z'_2|
\end{eqnarray*}
for some positive constants $C$ and $C'$. This integral is of the
same order as $z_2-z'_2$ and this compensates the denominator
corresponding to the second chord.

More generally, one shows by induction that if a
``short'' chord $(z_j,z'_j)$ is separated from the maximum by
$j-1$ chords, the first of which (that is, the nearest to the maximum)
is ``long'', the integral
$$       \int\limits_{t_{j}<t_{j-1}<\dots<t_1<t_{c}}
         \bigwedge_{i=1}^{j-1} \frac{dz_i-dz'_i}{z_i-z'_i}
$$
is of the same order as $z_j-z'_j$. This implies the convergence of
the Kontsevich integral.
\end{proof}

\section{Invariance of the integral}

\begin{theorem}\label{horinv}
\index{Kontsevich integral!invariance}%
The Kontsevich integral is invariant under the deformations in the
class of (not necessarily strict) Morse knots.
\end{theorem}

The proof of this theorem spans the whole of this section.
\medskip

Any deformation of a knot within the class of Morse knots can be
approximated by a sequence of deformations of three types:
orientation-preserving reparametrizations, {\em horizontal
deformations} and {\em movements of critical points}.

The invariance of the Kontsevich integral under
orientation-preserving reparametrizations is immediate since the
parameter plays no role in the definition of the integral apart from
determining the orientation of the knot.

A horizontal deformation is an isotopy of a knot in $\R^3$ which
preserves all horizontal planes $\{t=\const\}$ and leaves all the
critical points (together with some small neighbourhoods) fixed. The
invariance under horizontal deformations is the most essential point
of the theory. We prove it in the next subsection.

A  movement of a critical point $C$ is an isotopy which is identical
everywhere outside a small neighbourhood of $C$ and does not
introduce new critical points on the knot.
The invariance of the Kontsevich integral under the movements of critical points will be considered in \ref{ssmcrit}.

As we mentioned before, the Kontsevich integral is not invariant
under isotopies that change the number of critical points. Its
behaviour under such deformations will be discussed in Section
\ref{change_no_crit_pts}.

\subsection{Invariance under horizontal deformations}
\label{hor_inv}

Let us decompose the given knot into a product of tangles without
critical points of the function $t$ and very thin tangles containing
the critical levels. A horizontal deformation keeps fixed the
neighbourhoods of the critical points, so, due to multiplicativity,
it is enough to prove that the Kontsevich integral for a tangle
without critical points is invariant under horizontal deformations
that preserve the boundary pointwise.

\begin{xproposition}
Let $T_0$ be a tangle without critical points
and $T_\lambda$, a horizontal deformation
of $T_0$ to $T_1$ (preserving the top and the bottom of the tangle).
Then $Z(T_0) = Z(T_1)$.
\end{xproposition}

\begin{proof}
Denote by $\om$ the integrand form in the $m$th term of the
Kontsevich integral:
$$
  \om = \sum_{P=\{(z_j,z'_j)\}}\  (-1)^\downarrow D_P \quad
         \bigwedge_{j=1}^m \frac{dz_j-dz'_j} {z_j-z'_j}\ .
$$
Here the functions $z_j$, $z'_j$ depend not only on $t_1$, ..., $t_m$, but
also on $\lambda$, and {\it all differentials are understood as complete
differentials with respect to all these variables}. This means that the
form $\om$ is not exactly the form which appears in the Kontsevich's
integral (it has some additional $d\lambda$'s), but this does not change the
integrals over the simplices
$$
  \Delta_\lambda=\{t_{\mbox{\sz min}}<t_m<\dots<t_1<t_{\mbox{\sz max}}\}
   \times\{\lambda\},
$$
because the value of $\lambda$ on such a simplex is fixed.

We must prove that the integral of $\om$ over $\Delta_0$ is equal to its
integral over $\Delta_1$.

Consider the product polytope
$$
  \Delta = \Delta_0\times[0,1] = \quad
\risS{-15}{del}{
   \put(-5,-7){\mbox{${\scriptstyle \Delta_0}$}}
   \put(67,-7){\mbox{${\scriptstyle \Delta_1}$}}
      }{90}{20}{20}\ .
$$
By Stokes' theorem, we have \quad
${\displaystyle \int\limits_{\partial\Delta} \om\quad =
                      \quad\int\limits_\Delta d\om}\ $.

The form $\om$ is closed: $d\om =0$.
The boundary of the integration domain is
$\partial\Delta = \Delta_0 - \Delta_1 +\sum\mbox{\{\it faces\}}$.
The theorem will follow from the fact that
$\om\vert_{\mbox{\small\{{\it face}\}}}=0$.
To show this, consider two types of faces.

\medskip
The first type corresponds to $t_m = t_{\mbox{\sz min}}$ or $t_1 =
t_{\mbox{\sz max}}$ . In this situation, $dz_j=dz'_j=0$ for $j=1$ or
$m$, since $z_j$ and $z'_j$  do not depend on $\lambda$.

\medskip
The faces of the second type are those where we have $t_k = t_{k+1}$
for some $k$. In this case we have to choose the $k$th and $(k+1)$st
chords on the same level $\{t=t_k\}$. In general, the endpoints of
these chords may coincide and we do not get a chord diagram at all.
Strictly speaking, $\om$ and $D_P$  do not extend to such a face so
we have to be careful.
Extend $D_P$ to this face in the following manner: if 
some endpoints of $k$th and $(k+1)$st chords belong to the same
string (and therefore coincide) we place $k$th chord a little higher
than $(k+1)$st chord, so that its endpoint differs from the endpoint
of $(k+1)$st chord. This trick yields a well-defined prolongation of
$D_P$ and $\om$ to the face, and we use it here.

All summands of $\om$ are divided into three parts:

 (1) $k$th and $(k+1)$st chords connect the same two strings;

 (2) $k$th and $(k+1)$st chords
are chosen in such a way that their endpoints belong to four different
strings;

 (3) $k$th and $(k+1)$st chords
are chosen in such a way that there exist exactly three different strings
containing their endpoints.

Consider all these cases one by one.

1) We have $z_k=z_{k+1}$ and $z'_k=z'_{k+1}$ or vice versa. So 
$$
d(z_k-z'_k)\wedge d(z_{k+1}-z'_{k+1}) = 0
$$ 
and, therefore, the restriction of $\om$ to the face is zero.

2) All choices of chords in this part of $\om$ appear in mutually
canceling pairs. Fix four strings and number them by 1, 2, 3, 4.
Suppose that for a certain choice of the pairing, the $k$th chord
connects the first two strings and $(k+1)$st chord connects the last
two strings. Then there exists another choice for which, on the contrary, the $k$th chord connects the last two strings and $(k+1)$st
chord connects the first two strings. These two choices give two
summands of $\om$ differing by a sign:
$$
  \cdots d(z_k-z'_k)\wedge d(z_{k+1}-z'_{k+1}) \cdots +
  \cdots d(z_{k+1}-z'_{k+1})\wedge d(z_k-z'_k) \cdots = 0.
$$

3) This is the most difficult case. The endpoints of $k$th and
$(k+1)$st chords have exactly one string in common. Call the three
relevant strings 1, 2, 3 and denote by $\om_{ij}$ the 1-form
${\displaystyle \frac{dz_i-dz_j} {z_i-z_j}}$. Then $\om$ is the
product of a certain $(m-2)$-form and the sum of the following six
\label{sixpossum} 2-forms:
\newcommand\ttc[1]{(-1)^\downarrow\ \rb{-4.5mm}{\ig[height=11mm]{tijkN#1.eps}}}
\newcommand\ftc[2]{\om_{#1}\wedge\om_{#2}}
$$\begin{array}{rcl}
 \ttc{2}\ \ftc{12}{23} &+&  \ttc{6}\ \ftc{12}{13} \\[7mm]
+\ttc{5}\ \ftc{13}{12} &+&  \ttc{4}\ \ftc{13}{23} \\[7mm]
+\ttc{1}\ \ftc{23}{12} &+&  \ttc{3}\ \ftc{23}{13}\ .
\end{array}$$

Using the fact that $\om_{ij}=\om_{ji}$, we can rewrite this as follows:
$$
 \begin{array}{rc}
   \biggl(\ttc{2}\ -\ \ \ttc{1} \biggr)\ \ftc{12}{23} & \\[7mm]
  +\biggl(\ttc{3}\ -\ \ \ttc{4} \biggr)\ \ftc{23}{31} & \\[7mm]
  +\biggl(\ttc{5}\ -\ \ \ttc{6} \biggr)\ \ftc{31}{12} & .
 \end{array}
$$
The four-term relations in horizontal form (page
\pageref{horizontal4T}) say that the expressions in parentheses are
one and the same element of $\A$, hence, the whole sum is equal to
$$
  \biggl(\ttc{2}\ -\ \ \ttc{1} \biggr)\
  (\ftc{12}{23}+\ftc{23}{31}+\ftc{31}{12}).
$$
The 2-form that appears here is actually zero! This simple, but
remarkable fact, known as {\it Arnold's identity} (see \cite{Ar1})
can be put into the following form:
$$
  f+g+h=0\ \Longrightarrow\ \frac{df}{f}\wedge\frac{dg}{g}
  +\frac{dg}{g}\wedge\frac{dh}{h}+\frac{dh}{h}\wedge\frac{df}{f}=0
$$
(in our case $f=z_1-z_2$, $g=z_2-z_3$, $h=z_3-z_1$)
and verified by a direct computation.

This finishes the proof.
\end{proof}

\begin{xremark} The Kontsevich integral of a tangle may change, if the
boundary points are moved. Examples may be found below in
Exercises~\ref{kiasst}---\ref{kimint}.
\end{xremark}

\subsection{Moving the critical points}
\label{ssmcrit}

Let $T_0$ and $T_1$ be two tangles which are identical except a
sharp ``tail'' of width $\e$, which may be twisted:
$$
  \rb{-13mm}{\ig[height=28mm]{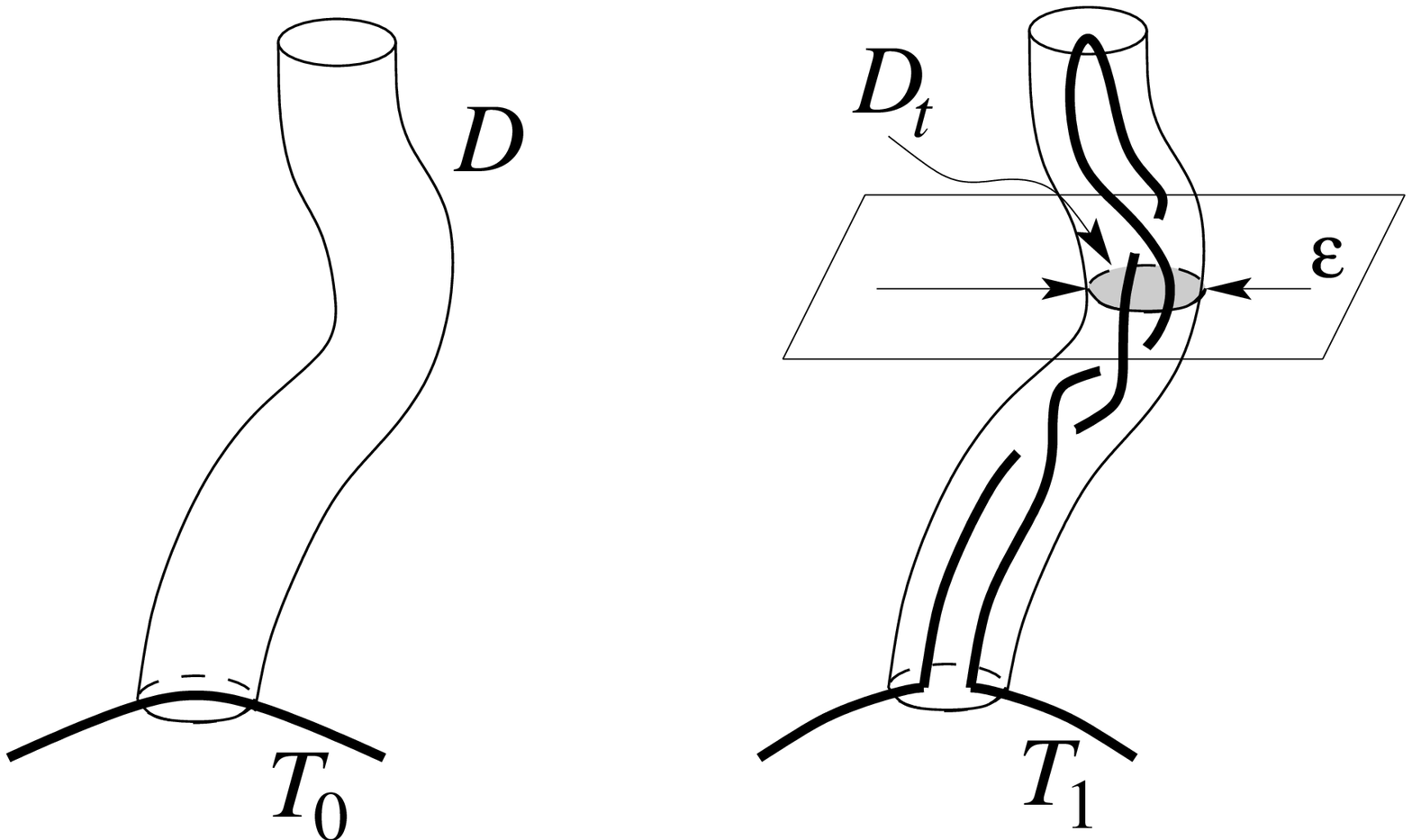}}
$$
More exactly, we assume that (1) $T_1$ is different from $T_0$ only
inside a region $D$ which is the union of disks $D_t$ of diameter
$\e$ lying in horizontal planes with fixed $t\in[t_1,t_2]$, (2) each
tangle $T_0$ and $T_1$ has exactly one critical point in $D$, and
(3) each tangle $T_0$ and $T_1$ intersects every disk $D_t$ at most
in two points. We call the passage from $T_0$ to $T_1$ a {\em
special movement\/}\label{spec_mov} of the critical point. To prove
Theorem~\ref{horinv} it is sufficient to show the invariance of the
Kontsevich integral under such movements. Note that special
movements of critical points may take a Morse knot out of the class
of strict Morse knots.

\begin{xproposition}
The Kontsevich integral remains unchanged under a special movement
of the critical point: $Z(T_0) = Z(T_1)$.
\end{xproposition}

\begin{proof}
The difference between $Z(T_0)$ and $Z(T_1)$ can come only from the
terms with a chord ending on the tail.

Consider the tangle $T_1$ ($T_0$ is considered similarly.)
If the highest of such chords connects the two sides of the tail, then the
corresponding tangle chord diagram is zero by a one-term relation.

So we can assume that the highest, say, the $k$th,  chord is a
``long'' chord, which means that it connects the tail with another
part of $T_1$. Suppose the endpoint of the chord belonging to the
tail is $(z'_k,t_k)$. Then there exists another choice for $k$th
chord which is almost the same but ends at another point of the tail
$(z''_k,t_k)$ on the same horizontal level:
$$
  \rb{-13mm}{\ig[height=28mm]{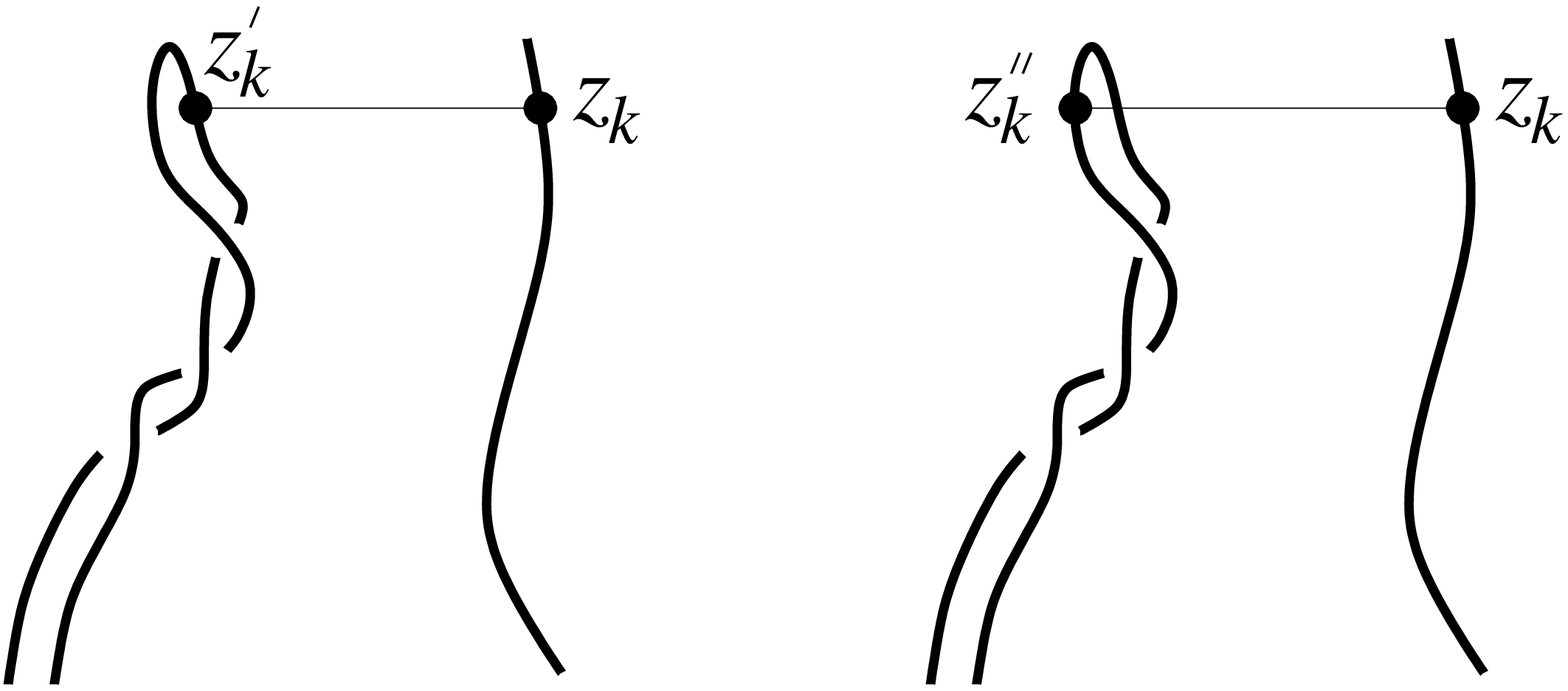}}
$$
The corresponding two terms appear in $Z(T_1)$ with the opposite
signs due to the sign $(-1)^\downarrow$.

Let us estimate the difference of the integrals corresponding to
such $k$th chords:
\begin{gather*}
\left|\  \int\limits_{t_{k+1}}^{t_c} d(\ln(z'_k-z_k))\ -
        \ \int\limits_{t_{k+1}}^{t_c} d(\ln(z''_k-z_k)) \right| =
\left|\ \ln\biggl(
     \frac{z''_{k+1}-z_{k+1}}{z'_{k+1}-z_{k+1}}\biggr) \right|
\\
=\left|\ \ln\biggl(
     1+\frac{z''_{k+1}-z'_{k+1}}{z'_{k+1}-z_{k+1}}\biggr) \right|
    \sim \left| z''_{k+1}-z'_{k+1}\right| \leq\e
\end{gather*}
(here $t_c$ is the value of $t$ at the uppermost point of the tail).

Now, if the next $(k+1)$st chord is also long, then, similarly, it
can be paired with another long chord so that they give a
contribution to the integral proportional to $\left|
z''_{k+2}-z'_{k+2}\right| \leq\e$.

In the case the $(k+1)$st chord is short (that is, it connects two
points $z''_{k+1}$, $z'_{k+1}$ of the tail) we have the following
estimate for the double integral corresponding to $k$th and
$(k+1)$st chords:
\begin{align*}
&\left|\ \displaystyle \int\limits_{t_{k+2}}^{t_c} \Biggl(
        \ \int\limits_{t_{k+1}}^{t_c} d(\ln(z'_k-z_k))\ -
        \ \int\limits_{t_{k+1}}^{t_c} d(\ln(z''_k-z_k))\ \Biggr)\
     \displaystyle \frac{dz''_{k+1}-dz'_{k+1}}{z''_{k+1}-z'_{k+1}} \right| \\
&\leq\
\const\cdot\left|\ \displaystyle \int\limits_{t_{k+2}}^{t_c}
        \left| z''_{k+1}-z'_{k+1}\right|
   \displaystyle \frac{dz''_{k+1}-dz'_{k+1}}{\left| z''_{k+1}-z'_{k+1}\right|}
\right|
\end{align*}

$\displaystyle\ \
=\ \const\cdot\left|\ \displaystyle \int\limits_{t_{k+2}}^{t_c}
d(z''_{k+1}-z'_{k+1}) \right|
           \sim \left| z''_{k+2}-z'_{k+2}\right| \leq\e\ .$

Continuing this argument, we see that the difference between
$Z(T_0)$ and $Z(T_1)$ is $O(\e)$. Now, by horizontal deformations we
can make $\e$ tend to zero. This proves the theorem and completes
the proof of the Kontsevich integral's invariance in the class of
knots with nondegenerate critical points.
\end{proof}

\section{Changing the number of critical points}\label{change_no_crit_pts}

The multiplicativity of the Kontsevich integral for tangles (Propositions~\ref{multTKI} and
\ref{hormult}) have several immediate consequences for knots.

\subsection{From long knots to usual knots}
\label{long2closed}

A long (Morse) knot can be closed up so as to produce a usual (Morse) knot:
$$\ig[width=5cm]{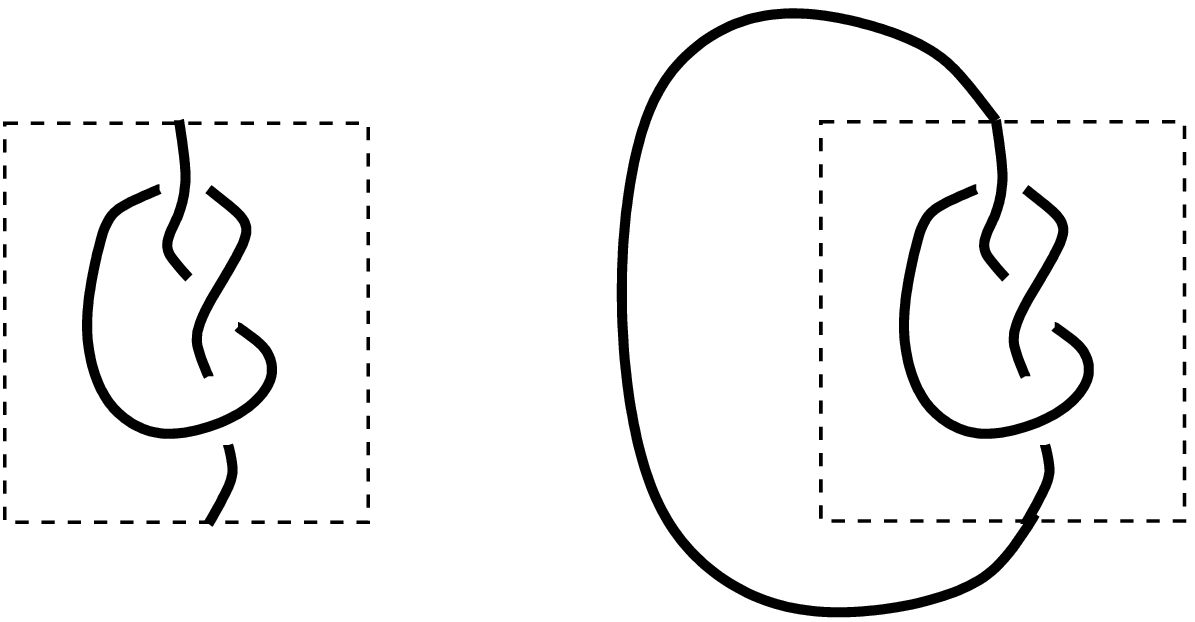}$$

Recall that the algebras of chord diagrams for long knots and for usual
knots are essentially the same; the isomorphism is given by closing up a linear chord diagram.
\begin{xproposition}
The Kontsevich integral of a long knot $T$ coincides with that of its closure $K_T$.
\end{xproposition}
\begin{proof}
Denote by ${\id}$ the tangle consisting of one vertical strand. Then
$K_T$ can be written as $T_{max}\cdot (T\otimes_{\e} {\id})\cdot
T_{min}$ where $T_{max}$ and $T_{min}$ are a maximum and a minimum
respectively, and $0<\e\leq 1$.

Since the Kontsevich integral of $K_T$ does not depend on $\e$, we can
take $\e\to 0$. Therefore,
$$Z(K_T)=Z(T_{max})\cdot (Z(T)\otimes Z({\id})) \cdot Z(T_{min}).$$
However, the Kontsevich integrals of $T_{max}$, $T_{min}$ and
${\id}$ consist of one diagram with no chords, and the Proposition
follows.
\end{proof}

A corollary of this is the formula (\ref{addcritF}) (page~\pageref{addcritF})
which describes the behaviour of the Kontsevich integral under the addition
of a pair of critical points. Indeed, adding a pair of critical points to a
long knot $T$ is the same as multiplying it by
$$\ig[width=30pt]{humpt.eps},$$
and (\ref{addcritF}) then follows from the multiplicativity of the Kontsevich
integral for tangles.

\subsection{The universal Vassiliev invariant}
\label{univ_vi} \index{Vassiliev!invariant!universal for knots}
\index{Universal Vassiliev invariant} The formula~(\ref{addcritF})
allows one to define the {\em universal Vassiliev invariant} by
either
$$
  I(K) = \frac{Z(K)}{Z(H)^{c/2}}
$$
or
$$
  I'(K) = \frac{Z(K)}{Z(H)^{c/2-1}},
$$ \label{Iprime}
where $c$ denotes the number of critical points of $K$ in an
arbitrary Morse representation, and the quotient means division in
the algebra $\widehat\A$: $(1+a)^{-1} = 1-a+a^2-a^3+\dots$.

Any isotopy of a knot in $\R^3$ can be approximated by a sequence
consisting of isotopies within the class of (not necessarily strict)
Morse knots and insertions/deletions of ``humps'', that is, pairs of
adjacent maxima and minima. Hence, the invariance of $Z(K)$ in the
class of Morse knots and the formula~(\ref{addcritF}) imply that
both $I(K)$ and $I'(K)$ are invariant under an arbitrary deformation
of $K$. (The meaning of the ``universality'' will be explained in
Section \ref{univki}.)

The version $I'(K)$ has the advantage of being multiplicative with
respect to the connected sum of knots; in particular, it vanishes
(more precisely, takes the value 1) on the unknot. However, the
version $I(K)$ is also used as it has a direct relationship with the
quantum invariants (see \cite{Oht1}). In particular, we shall use the term
``Kontsevich integral of the unknot''; this, of course, refers to
$I$, and not $I'$.

\section{Proof of the Kontsevich theorem}
\label{pkt}

First of all we reformulate the Kontsevich theorem (or, more exactly,
Kontsevich's part of the Vassiliev--Kontsevich theorem \ref{fund_thm})
as follows.

\begin{theorem}\label{Ko_part_fund_th}
Let $w$ be an unframed weight system of order $n$. Then there exists
a Vassiliev invariant of order $n$ whose symbol is $w$.
\end{theorem}

\begin{proof} The desired knot invariant is given by the formula
$$K \longmapsto w(I(K)).$$
Let $D$ be a chord diagram of order $n$ and let $K_D$ be a singular
knot with chord diagram $D$. The theorem follows from the fact that
$I(K_D) = D +(\mbox{terms of order $> n$})$. Since the denominator
of $I(K)$ starts with the unit of the algebra $\A$, it is sufficient
to prove that
\begin{equation}\label{ki_sing_knot}
  Z(K_D) = D + (\mbox{terms of order $>n$}).
\end{equation}
In fact, we shall establish (\ref{ki_sing_knot}) for $D$ an
arbitrary {\em tangle} chord diagram and $K_D=T_D$ a singular
tangle with the diagram $D$.

If $n=0$, the diagram $D$ has no chords and $T_D$ is non-singular.
For a non-singular tangle the Kontsevich integral starts with a
tangle chord diagram with no chords, and (\ref{ki_sing_knot})
clearly holds. Note that the Kontsevich integral of any singular
tangle (with at least one double point) necessarily starts with
terms of degree at least 1.

Consider now the case $n=1$. If $T_D$ is a singular 2-braid, there
is only one possible term of degree 1, namely the chord diagram with
the chord connecting the two strands. The coefficients of this
diagram in $Z(T_{+})$ and $Z(T_{-})$, where $T_{+}-T_{-}$ is a
resolution of the double point of $T_D$, simply measure the number
of full twists in $T_{+}$ and $T_{-}$ respectively. The difference
of these numbers is 1, so in this case (\ref{ki_sing_knot}) is also
true.

Now, let $T_D$ be an arbitrary singular tangle with exactly one
double point, and $V_{\e}$ be the
$\e$-neighbourhood of the singularity. We can assume that
the intersection of  $T_D$ with $V_{\e}$ is a singular
2-braid, and that the double point of $T_D$ is resolved as
$T_D=T^{\e}_{+}-T^{\e}_{-}$ where
$T^{\e}_{+}$ and $T^{\e}_{-}$ coincide with $T$
outside $V_{\e}$.

Let us write the degree 1 part of $Z(T^{\e}_{\pm})$ as a
sum $Z'_{\pm}+Z''_{\pm}$ where $Z'_{\pm}$ is the integral over all
chords whose both ends are contained in $V_{\e}$ and
$Z''_{\pm}$ is the rest, that is, the integral over the chords with
at least one end outside $V_{\e}$. As $\e$ tends
to 0, $Z''_{+}-Z''_{-}$ vanishes. On the other hand, for all
$\e$ we have that $Z''_{+}-Z''_{-}$ is equal to the diagram
$D$ with the coefficient 1. This settles the case $n=1$.

Finally, if $n>1$, using a suitable deformation, if necessary, we
can always achieve that $T_D$ is a product of $n$ singular tangles
with one double point each. Now (\ref{ki_sing_knot}) follows from
the multiplicativity of the Kontsevich integral for tangles.
\end{proof}

\bigskip

\subsection{Universality of $I(K)$}\label{univki}

In the proof of the Kontsevich Theorem we have seen that for a
singular knot $K$ with $n$ double points, $I(K)$ starts with terms
of degree $n$. This means that if $I_n(K)$ denotes the $n$th graded
component of the series $I(K)$, then the function $K \mapsto I_n(K)$
is a Vassiliev invariant of order $n$.

In some sense, all Vassiliev invariants are of this type:

\begin{proposition}
Any Vassiliev invariant can be factored through $I$:
for any $v\in\V$ there exists a linear function $f$ on $\widehat\A$
such that $v=f\circ I$.
\end{proposition}

\begin{proof}
Let $v\in\V_n$. By the Kontsevich theorem we know that there is a
function $f_0$ such that $v$ and $f_0\circ I_n$ have the same
symbol. Therefore, the highest part of the difference $v-f_0\circ I_n$ 
belongs to
$\V_{n-1}$ and is thus representable as $f_1\circ I_{n-1}$.
Proceeding in this way, we shall finally obtain:
$$
   v = \sum_{i=1}^n f_i \circ I_{n-i}.
$$
\end{proof}

\begin{xremark} The construction of the foregoing proof shows that the
universal Vassiliev invariant induces a splitting of the filtered
space $\V$ into a direct sum with summands isomorphic to the factors
$\V_n/\V_{n-1}$. Elements of these subspaces are referred to as {\em
canonical Vassiliev
invariants}.\index{Vassiliev!invariant!canonical} We shall speak
about them in more detail later in Section \ref{canvi}.
\end{xremark}

As a corollary, we get the following statement:
\begin{theorem}
The universal Vassiliev invariant $I$ is exactly as strong as the set of
all Vassiliev invariants: for any two knots $K_1$ and $K_2$ we have
$$
  I(K_1)=I(K_2) \quad \Longleftrightarrow \quad
  \forall v\in\V \quad v(K_1)=v(K_2).
$$
\end{theorem}

\section{Symmetries and the group-like property of $Z(K)$}
\subsection{Reality}
\label{reality} \index{Kontsevich integral!reality} Choose a basis
in the vector space $\A$ consisting of chord diagrams. A priori, the
coefficients of the Kontsevich integral of a knot $K$ with respect
to this basis are complex numbers.
\begin{xtheorem} All coefficients of the Kontsevich integral with
respect to a basis of chord diagrams are real.
\end{xtheorem}
\begin{xremark} Of course, this fact is a consequence of the
Le--Murakami Theorem stating that these coefficients are rational
(Section~\ref{kifieqki}). However, the rationality of the Kontsevich
integral is a very non-trivial fact, while the proof that its
coefficients are real is very simple.
\end{xremark}
\begin{proof}
Rotate the coordinate frame in $\R^3$ around the real axis $x$ by
$180^\circ$; denote the new coordinates by $t^{\star}=-t$,
$z^{\star}=\ol{z}$. If $K$ is a Morse knot, it will still be a Morse
knot, with the same number of maxima, with respect to the new
coordinates, and its Kontsevich integral, both preliminary and
final, will be the same in both coordinate systems. Let us denote
the preliminary Kontsevich integral with respect to the starred
coordinates by $Z^{\star}(K)$.

For each pairing $P=\{(z_j,z'_j)\}$ with $m$ chords that appears in
the formula for $Z(K)$, there is a pairing
$P^{\star}=\{(z^{\star}_j,{z_j^{\star}}')\}$ that appears in the
formula for $Z^{\star}(K)$. It consists of the very same chords as
$P$ but taken in the starred coordinate system: $z^{\star}_j =
\ol{z}_{m-j+1}$ and ${z_j^{\star}}' = \ol{z'}_{m-j+1}$. The
corresponding chord diagrams are, obviously, equal:
$D_P=D_{P^{\star}}$. Moreover, $\downarrow^{\star} =
2m-\!\!\downarrow$ and, hence,
$(-1)^{\downarrow^{\star}}=(-1)^{\downarrow}$. The simplex
$\Delta=t_{\mbox{\sz min}}<t_1<\dots<t_m<t_{\mbox{\sz max}}$ for the
variables $t_i$ corresponds to the simplex $\Delta^{\star}
=-t_{\mbox{\sz max}}<t^{\star}_m<\dots<t^{\star}_1<-t_{\mbox{\sz
min}}$ for the variables $t^{\star}_i$. The coefficient of
$D_{P^{\star}}$ in $Z^{\star}(K)$ is
$$
  c(D_{P^{\star}})\ =\
  \frac{(-1)^{\downarrow}}{(2\pi i)^m}
  \int \bigwedge_{j=1}^m d\ln(z^{\star}_j-{z_j^{\star}}') \ ,
$$
where $z^{\star}_j$ and ${z_j^{\star}}'$ are understood as functions
in $t^{\star}_1$, \dots, $t^{\star}_m$ and the integral is taken
over a connected component in the simplex $\Delta^{\star}$. In the
last integral we make the change of variables according to the
formula $t^{\star}_j=-t_{m-j+1}$. The Jacobian of this
transformation is equal to $(-1)^{m(m+1)/2}$. Therefore,
$$
  c(D_{P^{\star}})\ =\
  \frac{(-1)^{\downarrow}}{(2\pi i)^m}\int (-1)^{m(m+1)/2}
    \bigwedge_{j=1}^m d\ln(\ol{z}_{m-j+1}-\ol{z'}_{m-j+1})
$$
(integral over the corresponding connected component in the simplex
$\Delta$). Now permute the differentials to arrange the subscripts
in the increasing order. The sign of this permutation is
$(-1)^{m(m-1)/2}$. Note that
$(-1)^{m(m+1)/2}\cdot(-1)^{m(m-1)/2}=(-1)^m$. Hence,
$$\begin{array}{rcl}
  c(D_{P^{\star}})\ &=&\ \displaystyle
\frac{(-1)^{\downarrow}}{(2\pi i)^m}(-1)^m\ \int\
    \bigwedge_{j=1}^m d\ln(\ol{z}_j-\ol{z'}_j) \vspace{10pt}\\
  &=&\ \displaystyle
\frac{(-1)^{\downarrow}}{(2\pi i)^m}\ \int\
    \bigwedge_{j=1}^m \ol{d\ln(z_j-z'_j)}\quad =\quad \ol{c(D_P)} .
\end{array}$$
Since any chord diagram $D_P$ can be expressed as a combination of
the basis diagrams with real coefficients, this proves the theorem.
\end{proof}

\subsection{The group-like property}
\label{glki}
\begin{xtheorem}\index{Kontsevich integral!group-like}
For any Morse tangle $T$ with skeleton $\boldX$ the Kontsevich
integral $Z(T)$ is a group-like element \index{Group-like element}
in the graded completion of the coalgebra $\F(\boldX)$:
$$
\d(Z(T)) = Z(T)\ot Z(T)\,.
$$
In particular, if $K$ is a knot, $Z(K)$ is a group-like element in
$\Ab$.
\end{xtheorem}
\begin{proof}
The proof consists of a direct comparison of the coefficients on
both sides. Write the Kontsevich integral as
$$Z(K)=\sum_{P} c_P D_P,$$
where the sum is over all possible pairings $P$ for $K$. Then, for
given pairings $P_1$ and $P_2$ we have
$$\sum_{P\backslash J=P_1,\, J=P_2} c_P= c_{P_1}c_{P_2}.$$
Indeed, the chords in $P_1$ and $P_2$ vary over the simplices
$\Delta_1=\{t_{min}<t_{|P_1|}< \ldots <t_1<t_{max}\}$ and
$\Delta_2=\{t_{min}<\tilde{t}_{|P_2|}< \ldots
<\tilde{t}_1<t_{max}\}$ respectively. The product $\Delta_1\times
\Delta_2$ is subdivided into simplices $\Delta_{P,J}$, each of
dimension $|P_1|+|P_2|$, by hyperplanes $t_i=\tilde{t}_k$; each of
these simplices can be thought of a way of fusing $P_1$ and $P_2$
into one pairing $P$, with $P\backslash J=P_1$ and $ J=P_2$. Then
each term on the left-hand side is the integral
$$
\frac{(-1)^{\downarrow}}{(2\pi i)^{|P_1|+|P_2|}} \int_{\Delta_{P,J}}
\bigwedge_{j=1}^{|P_1|+|P_2|} d\ln(z_j-z_j') \ ,
$$
while $c_{P_1}c_{P_2}$ is the integral of the same expression over
the product $\Delta_1\times\Delta_2=\cup \Delta_{P,J}$.

It follows now that
$$\d(Z(K))=\sum_{P,J\subseteq P} c_P D_{P\backslash J}\ot D_J,$$
is equal to
$$Z(K)\ot Z(K)= \sum_{P_1,P_2} c_{P_1}c_{P_2} D_{P_1}\ot D_{P_2}.$$
\end{proof}

Group-like elements in a bialgebra form a group and this implies
that the final Kontsevich integral is also group-like.

\subsection{Change of orientation}
\label{invertib}
\begin{xtheorem}  The Kontsevich integral commutes
with the operation $\tau$ of orientation reversal:
$$
  Z(\tau(K))=\tau(Z(K)).
$$
\end{xtheorem}
\begin{proof} Changing the orientation of $K$ has the following effect
on the formula for the Kontsevich integral on
page~\pageref{formulaKI}: each diagram $D$ is replaced by $\tau(D)$
and the factor $(-1)^\downarrow$ is replaced by $(-1)^\uparrow$,
where by $\uparrow$ we mean, of course, the number of points
$(z_j,t_j)$ or $(z'_j,t_j)$ in a pairing $P$ where the coordinate
$t$ grows along the parameter of $K$. Since the number of points in
a pairing is always even, $(-1)^\downarrow=(-1)^\uparrow$, so that
$\tau(D)$ appears in  $Z(\tau(K))$ with the same coefficient as $D$
in $Z(K)$. The theorem is proved.
\end{proof}
\begin{xcorollary} The following two assertions are equivalent:
\begin{itemize}
\setlength{\itemsep}{1pt plus 1pt minus 1pt}
\item
Vassiliev invariants do not distinguish the orientation of knots,
\item
all chord diagrams are symmetric: $D=\tau(D)$ modulo one- and
four-term relations.
\end{itemize}
\end{xcorollary}
The calculations of \cite{Kn0} show that up to order 12 all chord
diagrams are symmetric. For bigger orders the problem is still open.

\subsection{Mirror images}
\label{mirsym}

Recall that $\sigma$ is the operation sending a knot to its mirror
image (see \ref{orientation}). Define the corresponding operation
$\sigma:\A\to \A$ by sending a chord diagram $D$ to $(-1)^{\deg{D}}
D$. It extends to a map $\Ab\to \Ab$ which we also denote by
$\sigma$.
\begin{xtheorem} \label{mir_refl}
The Kontsevich integral commutes with $\sigma$:
$$Z(\sigma(K))=\sigma(Z(K))\,.$$
\end{xtheorem}
\begin{proof} Let us realize the operation $\sigma$ on knots by the reflection
of $\R^3$ coming from the complex conjugation in $\C$:
$(z,t)\mapsto(\bar{z},t)$. Then the Kontsevich integral for
$\sigma(K)$ can be written as
$$\begin{array}{rcl}
Z(\sigma(K)) &=& \displaystyle\sum_{m=0}^\infty \ \frac1{(2\pi i)^m}
            \int \sum_P\  (-1)^\downarrow D_P
            \bigwedge_{j=1}^m d\ln(\ol{z_j-z'_j}) \vspace{20pt}\\
        &=& \displaystyle
            \sum_{m=0}^\infty \ (-1)^m \ol{\frac1{(2\pi i)^m}
            \int \sum_P\  (-1)^\downarrow D_P
            \bigwedge_{j=1}^m d\ln(z_j-z'_j)} \ .
\end{array}$$
Comparing this with the formula for $Z(K)$ we see that the terms of
$Z(\sigma(K))$ with an even number of chords coincide with those of
$\ol{Z(K)}$ and terms of $Z(\sigma(K))$ with an odd number of chords
differ from the corresponding terms of $\ol{Z(K)}$ by a sign. Since
$Z(K)$ is real, this implies the theorem.
\end{proof}

Since the Kontsevich integral is equivalent to the totality of all
finite type invariants, Theorem \ref{mir_refl} implies that if  $v$
is a Vassiliev invariant of degree $n$,  $K$ is a singular knot with
$n$ double points and $\ol{K}=\sigma(K)$ its mirror image, then
$v(K)=v(\ol{K})$ for even $n$ and $v(K)=-v(\ol{K})$ for odd $n$.

\noindent{\bf Exercise.} Prove this statement without using the
Kontsevich integral.
\medskip

Recall (page~\pageref{amphicheiral}) that a knot $K$ is called {\em
plus-amphicheiral}, \index{Knot!plus-amphicheiral} if it is
equivalent to its mirror image as an oriented knot: $K=\sigma(K)$,
and {\em minus-amphicheiral} if it is equivalent to the inverse of
the mirror image: $K=\tau{\sigma{K}}$. Write $\tau$ for the mirror
reflection on chord diagrams (see \ref{detect_orient}), and recall
that an element of $\A$ is called {\em symmetric},\index{Chord
diagram!symmetric} ({\em antisymmetric}), \index{Chord
diagram!anti-symmetric} if $\tau$ acts on it as identity, (as
multiplication by $-1$, respectively).
\begin{xcorollary}
The Kontsevich integral $Z(K)$ of a plus-amphicheiral knot $K$
consist only of even order terms. For a minus-amphicheiral knot $K$
the Kontsevich integral $Z(K)$ has the following property: its
even-degree part consists only of symmetric chord diagrams, while
the odd-degree part consists only of anti-symmetric diagrams. The
same is true for the universal Vassiliev invariant $I(K)$.
\end{xcorollary}
\begin{proof}
For a plus-amphicheiral knot, the theorem implies that
$Z(K)=\sigma(Z(K))$, hence all the odd order terms in the series
$Z(K)$ vanish. The quotient of two even series in the graded
completion $\widehat\A$ is obviously even, therefore the same
property holds for $I(K)=Z(K)/Z(H)^{c/2}$.

For a minus-amphicheiral knot $K$, we have
$Z(K)=\tau(\sigma(Z(K)))$, which implies the second assertion.
\end{proof}

Note that it is an open question whether non-symmetric chord
diagrams exist. If they do not, then, of course, both assertions of
the theorem, for plus- and minus-amphicheiral knots, coincide.

\section{Towards the combinatorial Kontsevich integral}
\label{comb_ki_intro}

Since the Kontsevich integral comprises all Vassiliev invariants,
calculating it explicitly is a very important problem. Knots are,
essentially, combinatorial objects so it is not surprising that the
Kontsevich integral, which we have defined analytically, can be
calculated combinatorially from the knot diagram. Different versions
of such combinatorial definition were proposed in several papers
(\cite{BN2,Car1,LM1,LM2,Piu}) and treated in several books
(\cite{Kas,Oht1}). Such a definition will be given in
Chapter~\ref{DA_chap}; here we shall explain the idea behind it.

The multiplicativity of the Kontsevich integral hints at the
following method of computing it: present a knot as a product of
several standard tangles whose Kontsevich integral is known and then
multiply the corresponding values of the integral. This method works
well for the quantum invariants, see Sections \ref{qis4} and
\ref{qis5}; however, for the Kontsevich integral it turns out to be
too na\"\i ve to be of direct use.

Indeed, in the case of quantum invariants we decompose the knot into
elementary tangles, that is, crossings, max/min events and pieces of
vertical strands using both the usual product and the tensor product
of tangles. While the Kontsevich integral behaves well with respect
to the usual product of tangles, there is no simple expression for
the integral of the tensor product of two tangles, even if one of
the factors is a trivial tangle. As a consequence, the Kontsevich
integral is hard to calculate even for the generators of the braid
group, not to mention other possible candidates for ``standard''
tangles.

Still, we know that the Kontsevich integral is {\em asymptotically}
multiplicative with respect to the parametrized tensor product.
This suggests the following procedure.

Write a knot $K$ as a product of tangles $K=T_1\cdot\ldots\cdot T_n$
where each $T_i$ is simple, that is, a tensor product of several
elementary tangles. Let us think of each $T_i$ as of an
$\e$-parametrized tensor product of elementary tangles with $\e=1$.
We want to vary this $\e$ to make it very small. There are two
issues here that should be taken care of.

Firstly, the $\e$-parametrized tensor product is not
associative for $\e\neq 1$, so we need a {\em
parenthesizing} on the factors in $T_i$. We choose the
parenthesizing arbitrarily on each $T_i$ and denote by
$T_i^{\e}$ the tangle obtained from $T_i$ by replacing
$\e=1$ by an arbitrary positive $\e\leq 1$.

Secondly, even though the tangles $T_i$ and $T_{i+1}$ are
composable, the tangles $T_i^{\e}$ and
$T_{i+1}^{\e}$ may fail to be composable for
$\e<1$. Therefore, for each $i$ we have to choose a family
of {\em associating} \index{Associating
tangle}\index{Tangle!associating} tangles without crossings
$Q_i^{\e}$ which connect the bottom endpoints of
$T_i^{\e}$ with the corresponding top endpoints of
$T_{i+1}^{\e}$.

Now we can define a family of knots $K^{\e}$ as
$$K^{\e}=T_1^{\e}\cdot
Q_1^{\e}\cdot T_2^{\e}\cdot \ldots\cdot
Q_{n-1}^{\e}\cdot T_n^{\e}.$$ The following
picture illustrates this construction on the example of a trefoil
knot:

\begin{figure}[h]
 \begin{center}
     \begin{picture}(310,200)(0,0)
      \put(0,0){\includegraphics[width=150pt]{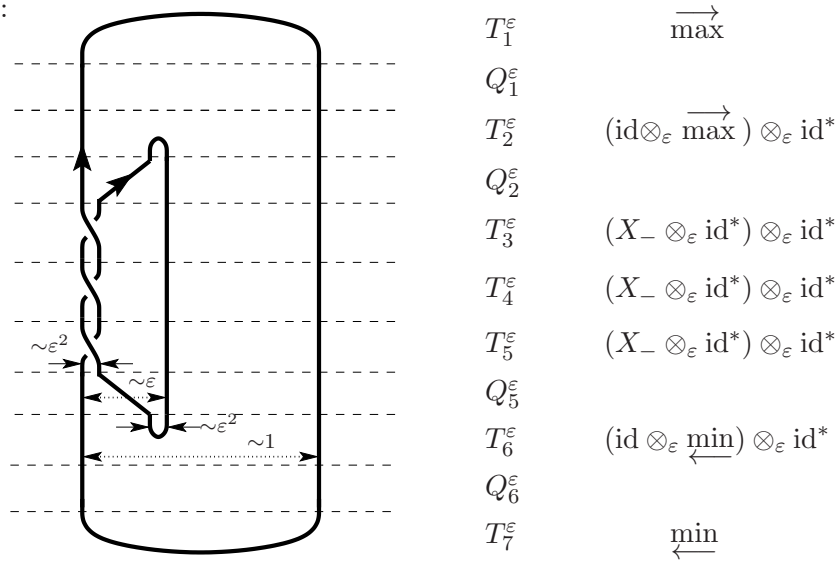}}

      \put(250,5){\mbox{$\minl$}}
      \put(225,41){\mbox{$(\id\ot_{\e}\minl)\ot_{\e}\id^*$}}
      \put(225,77){\mbox{($X_-\ot_{\e}\id^*)\ot_{\e}\id^*$}}
      \put(225,98){\mbox{($X_-\ot_{\e}\id^*)\ot_{\e}\id^*$}}
      \put(225,120){\mbox{($X_-\ot_{\e}\id^*)\ot_{\e}\id^*$}}
      \put(225,158){\mbox{$(\id\ot_{\e}\maxr)\ot_{\e}\id^*$}}
      \put(250,196){\mbox{$\maxr$}}

      \put(180,5){\mbox{$T_7^{\e}$}}
      \put(180,23){\mbox{$Q_6^{\e}$}}
      \put(180,41){\mbox{$T_6^{\e}$}}
      \put(180,59){\mbox{$Q_5^{\e}$}}
      \put(180,77){\mbox{$T_5^{\e}$}}
      \put(180,98){\mbox{$T_4^{\e}$}}
      \put(180,120){\mbox{$T_3^{\e}$}}
      \put(180,139){\mbox{$Q_2^{\e}$}}
      \put(180,158){\mbox{$T_2^{\e}$}}
      \put(180,177){\mbox{$Q_1^{\e}$}}
      \put(180,196){\mbox{$T_1^{\e}$}}

      \put(8,77){\mbox{$\scriptstyle\sim\e^2$}}
      \put(45,64){\mbox{$\scriptstyle\sim\e$}}
      \put(72,47){\mbox{$\scriptstyle\sim\e^2$}}
      \put(90,41){\mbox{$\scriptstyle\sim 1$}}
     \end{picture}
 \end{center}
\caption{A decomposition of the trefoil into associating tangles and
$\e$-parametrized tensor products of elementary tangles, with the
notations from Section~\ref{Turmov}. The associating tangles between
$T_3^{\e}$, $T_4^{\e}$ and $T_5^{\e}$ are omitted since these
tangles are composable for all $\e$.} \label{parenth_tref_ch8}
\label{parenth_tref} 
\end{figure}

Since for each $\e$ the knot $K^{\e}$ is isotopic to $K$ it is
tempting to take $\e\to 0$, calculate the limits of the Kontsevich
integrals of the factors and then take their product. The Kontsevich
integral of any elementary tangle, and, hence, of the limit
$$\lim_{\e\to 0}Z(T_i^{\e})$$ is easily evaluated, so it only
remains to calculate the limit of $Z(Q_i^{\e})$ as $\e$ tends to
zero.

\def\ZsmT#1{Z\bigl(\tupic{#1}\bigr)}

Calculating this last limit is not a straightforward task, to say
the least. In particular, if $Q^{\e}$ is the simplest
associating tangle

$$ \risS{-40}{assatel}{
         \put(7,1){\mbox{$\scriptstyle 0$}}
         \put(7,21){\mbox{$\scriptstyle \e$}}
         \put(22,2){\mbox{$\scriptstyle \e$}}
         \put(45,1){\mbox{$\scriptstyle 1-\e$}}
         \put(65,1){\mbox{$\scriptstyle 1$}}
         \put(-3,52){\mbox{$\scriptstyle 1-\e$}}
         \put(8,67){\mbox{$\scriptstyle 1$}}
         \put(18,80){\mbox{$t$}}
         \put(93,5){\mbox{$z$}}
                 }{90}{45}{50}
$$
we shall see in Chapter \ref{DA_chap} that asymptotically, as $\e\to 0$
we have
$$\ZsmT{asst0}\simeq \e^{\frac{1}{2\pi i}\vstdt}\cdot\PhiKZ\cdot\e^{-\frac{1}{2\pi i}\vstod},$$
where $\e^x$ is defined as the formal power series $\exp(x\log{\e})$
and $\PhiKZ$ is the power series known as the {\em
Knizhnik-Zamolodchikov associator}. Similar formulae can be written
for other associating tangles.

There are two difficulties here. One is that the integral
$Z(Q^{\e})$ does not converge as $\e$ tends to 0. However, all the
divergence is hidden in the terms $\e^{\frac{1}{2\pi i}\vstdt}$ and
$\e^{-\frac{1}{2\pi i}\vstod}$ and careful analysis shows that all
such terms from all associating tangles cancel each other out in the
limit, and can be omitted. The second problem is to calculate the
associator. This a highly non-trivial task, and is the main subject
of Chapter~\ref{DA_chap}.

\begin{xcb}{Exercises}
\begin{enumerate}

\item\label{pr:hor-hor-conf-sp}
\parbox[t]{2.9in}{For the link with two components $K$ and $L$ shown on the right draw
the configuration space of horizontal chords joining $K$ and $L$ as
in the proof of the linking number theorem from Section
\ref{subsec:linking_number} (see page \pageref{th:ki-lk}). Compute
the linking number of $K$ and $L$ using this theorem. } \qquad
\parbox{1.4in}{$\risS{-80}{l-3-1-t}{
         \put(5,65){\mbox{$\scriptstyle K$}}
         \put(70,5){\mbox{$\scriptstyle L$}}
                 }{75}{0}{10}$}\\

\item
Is it true that $Z(\ol{H})=Z(H)$, where $H$ is the hump as shown in page
\pageref{hump} and $\ol{H}$ is the same hump reflected in a horizontal line?

\item\label{ki_kont} 
M.~Kontsevich in his pioneering paper \cite{Kon1} and some of his followers
(for example, \cite{BN1,CD3}) defined the Kontsevich integral slightly
differently, numbering the chords upwards. Namely,\qquad $Z_{Kont}(K) = $
$$
  = \sum_{m=0}^\infty \frac{1}{(2\pi i)^m}
         \int\limits_{\substack{
    t_{\mbox{\tiny min}}<t_1<\dots<t_m<t_{\mbox{\tiny max}}\\
    t_j\mbox{\tiny\ are noncritical}}}
         \sum_{P=\{(z_j,z'_j)\}} (-1)^{\downarrow_P} D_P
         \bigwedge_{j=1}^m \frac{dz_j-dz'_j}{z_j-z'_j} \ .
$$
Prove that for any tangle $T$, $Z_{Kont}(T)=Z(T)$, as series of
tangle chord diagrams.

{\sl Hint.} Change of variables in multiple integrals.

\item
Express the integral over the cube
$$Z_{\square}(K) :=
  \sum_{m=0}^\infty \frac{1}{(2\pi i)^m}
         \int\limits_{\substack{
    t_{\mbox{\tiny min}}<t_1,\dots,t_m<t_{\mbox{\tiny max}}\\
    t_j\mbox{\tiny\ are noncritical}}}
         \sum_{P=\{(z_j,z'_j)\}} (-1)^{\downarrow_P} D_P
         \bigwedge_{j=1}^m \frac{dz_j-dz'_j}{z_j-z'_j}
$$
in terms of $Z(K)$.

\item
Compute the Kontsevich integral of the tangles \smtan{ort_xp_d} and
\smtan{ort_xp_l}\,.

\item\label{kislant}
\parbox[t]{2.8in}{Prove that for the tangle $\ \tupic{kis-slant}\ $ shown on the right\quad
$Z(\ \tupic{kis-slant}\ )=\exp\Bigl(
\frac{\risS{-3}{doc}{}{17}{18}{3}}{2\pi i}\cdot\ln\e\Bigr)\ .$
}
\qquad\qquad
\parbox[t]{1.4in}{$\risS{-50}{ki-slant}{
         \put(5,3){\mbox{$\scriptstyle 0$}}
         \put(20,3){\mbox{$\scriptstyle \e$}}
         \put(58,3){\mbox{$\scriptstyle 1$}}
         \put(80,4){\mbox{$\scriptstyle z$}}
         \put(-9,52){\mbox{$\scriptstyle 1-\e$}}
         \put(4,65){\mbox{$\scriptstyle t$}}
         \put(50,55){\mbox{$\scriptstyle t=1-z$}}
                 }{85}{25}{50}$}\\

\item\label{EuLi}
The {\it Euler dilogarithm}\index{Dilogarithm}\label{E-Li}
\label{dilogarithm}
is defined by the power series
$$
\Li(z)=\sum_{k=1}^\infty \frac{z^k}{k^2}
$$ 
for $|z|\leq 1$. Prove the following identities
$$\begin{array}{c}
\Li(0)=0;\qquad \Li(1)=\frac{\pi^2}{6};\qquad
    \Li'(z)=-\frac{\ln(1-z)}{z}\ ;\\
\frac{d}{dz}\Bigl(\Li(1-z)+\Li(z)+\ln z\ln(1-z)\Bigr)=0\ ;\\
\Li(1-z)+\Li(z)+\ln z\ln(1-z) = \frac{\pi^2}{6}\ .
\end{array}$$
About these and other remarkable properties of $\Li(z)$ see \cite{Lew,Kir,Zag2}.

\item\label{kiasst}
\parbox[t]{2.9in}{Consider the associating tangle
\index{Associating tangle}\index{Tangle!associating}
$\ \tupic{asst0}\ $ shown on the right. Compute $\ZsmT{asst0}$ up to
the second order.

{\sl Answer.}
$\tupic{s2td0}\ -\
\frac{1}{2\pi i} \ln\left(\frac{1-\e}{\e}\right)
     \Bigl(\tupic{s2td12} - \tupic{s2td23} \Bigr)$
\noindent
$
-\ \frac{1}{8\pi^2} \ln^2\left(\frac{1-\e}{\e}\right)
     \Bigl(\tupic{s2td1212} + \tupic{s2td2323} \Bigr)$
}
\qquad
\parbox{1.4in}{\makebox(100,30){\rb{-230pt}{\quad$\risS{-25}{asst}{
         \put(5,5){\mbox{$\scriptstyle 0$}}
         \put(20,5){\mbox{$\scriptstyle \e$}}
         \put(48,5){\mbox{$\scriptstyle 1-\e$}}
         \put(69,5){\mbox{$\scriptstyle 1$}}
         \put(80,5){\mbox{$\scriptstyle z$}}
         \put(4,20){\mbox{$\scriptstyle \e$}}
         \put(-7,64){\mbox{$\scriptstyle 1-\e$}}
         \put(4,83){\mbox{$\scriptstyle t$}}
         \put(40,77){\mbox{$\scriptstyle z=t$}}
                 }{85}{90}{100}$}}}\\
$+ \frac{1}{4\pi^2} \Bigl(\ln(1-\e)\ln\left(\frac{1-\e}{\e}\right)
   + \Li(1-\e) - \Li(\e) \Bigr) \tupic{s2td1223}$\\
$- \frac{1}{4\pi^2} \Bigl(\ln(\e)\ln\left(\frac{1-\e}{\e}\right)
   + \Li(1-\e) - \Li(\e) \Bigr) \tupic{s2td2312} $

 The calculation here uses the
 dilogarithm function defined in problem (\ref{EuLi}). Note that the Kontsevich
 integral diverges as $\e\to0$.

\item
Make the similar computation $\ZsmT{asst0rl}$ for the reflected tangle.
Describe the difference with the answer to the previous problem.

\item\label{kimaxt}
\parbox[t]{3.1in}{Compute the Kontsevich integral
$\ZsmT{maxt0}$ of the maximum tangle shown on the right.\\
{\sl Answer.}
$\tupic{maxt0}\ +\
\frac{1}{2\pi i} \ln(1-\e) \tupic{maxt12}$\\
$
+\ \frac{1}{4\pi^2} \Bigl(\Li\bigl(\frac{\e}{2-\e}\bigr) -
        \Li\bigl(\frac{-\e}{2-\e}\bigr)\Bigr) \tupic{maxt1223}$
}
\qquad
\parbox{1.2in}{\makebox(85,0){\rb{-135pt}{$\risS{-25}{maxt}{
         \put(12,3){\mbox{$\scriptstyle 0$}}
         \put(32,3){\mbox{$\scriptstyle 1-\e$}}
         \put(56,3){\mbox{$\scriptstyle 1$}}
         \put(75,3){\mbox{$\scriptstyle z$}}
         \put(-7,35){\mbox{$\scriptstyle 1-\e$}}
         \put(4,57){\mbox{$\scriptstyle 1$}}
         \put(4,68){\mbox{$\scriptstyle t$}}
         \put(25,67){\mbox{$\scriptstyle t=-z^2+(2-\e)z$}}
                 }{90}{70}{55}$}}}\\
$+\ \frac{1}{8\pi^2} \Bigl(\ln^2 2 -\ln^2\left(\frac{1-\e}{2-\e}\right)
      +2\Li\bigl(\frac{1}{2}\bigr) -
       2\Li\bigl(\frac{1-\e}{2-\e}\bigr)\Bigr) \tupic{maxt1212}$\\
$+\ \frac{1}{8\pi^2} \Bigl(\ln^2 2 -\ln^2(2-\e)
      +2\Li\bigl(\frac{1}{2}\bigr) -
       2\Li\bigl(\frac{1}{2-\e}\bigr)\Bigr) \tupic{maxt1213}$

\item\label{kimint}
\parbox[t]{3.1in}{Compute the Kontsevich integral
$\ZsmT{mint0}$ of the minimum tangle shown on the right.\\
{\sl Answer.}
$\tupic{mint0}\ -\
\frac{1}{2\pi i} \ln(1-\e) \tupic{mint23}$\\
$+\ \frac{1}{4\pi^2} \Bigl(\Li\bigl(\frac{\e}{2-\e}\bigr) -
        \Li\bigl(\frac{-\e}{2-\e}\bigr)\Bigr) \tupic{mint1223}$
}
\qquad
\parbox{1.2in}{\makebox(85,0){\rb{-110pt}{$\risS{-25}{mint}{
         \put(12,4){\mbox{$\scriptstyle 0$}}
         \put(30,4){\mbox{$\scriptstyle \e$}}
         \put(62,4){\mbox{$\scriptstyle 1$}}
         \put(77,4){\mbox{$\scriptstyle z$}}
         \put(4,30){\mbox{$\scriptstyle \e$}}
         \put(4,53){\mbox{$\scriptstyle t$}}
         \put(34,50){\mbox{$\scriptstyle t=z^2-\e z+\e$}}
                 }{90}{70}{40}$}}}\\
$+\ \frac{1}{8\pi^2} \Bigl(\ln^2 2 -\ln^2\left(\frac{1-\e}{2-\e}\right)
      +2\Li\bigl(\frac{1}{2}\bigr) -
       2\Li\bigl(\frac{1-\e}{2-\e}\bigr)\Bigr) \tupic{mint2323}$\\
$+\ \frac{1}{8\pi^2} \Bigl(\ln^2 2 -\ln^2(2-\e)
      +2\Li\bigl(\frac{1}{2}\bigr) -
       2\Li\bigl(\frac{1}{2-\e}\bigr)\Bigr) \tupic{mint1323}$

Note that all nontrivial terms in the last two problems tend to zero
as $\e\to 0$.

\item\label{kihump}
Express the Kontsevich integral of the hump as the product of tangle
chord diagrams from problems \ref{kiasst}, \ref{kimaxt},
\ref{kimint}:
$$Z\Bigl(\ \risS{-15}{humppr0}{}{15}{10}{3}\ \Bigr)  =
\ZsmT{maxt0}\cdot \ZsmT{asst0}\cdot \ZsmT{mint0}\ .$$
To do this introduce shorthand notation for the coefficients:\\
$\ZsmT{maxt0}=\tupic{maxt0}\ +\ A\tupic{maxt12}\ +\ B\tupic{maxt1223}\
    +\ C\,\tupic{maxt1212}\ +\ D\tupic{maxt1213} \vspace{10pt}$\\
$\ZsmT{asst0}=\tupic{asst0}\
        +\ E\left(\tupic{s2td12} - \tupic{s2td23}\right)
      \ +\ F\left(\tupic{s2td1212} + \tupic{s2td2323}\right)
  +\ G\,\tupic{s2td1223}\ +\ H\,\tupic{s2td2312}
      \vspace{10pt}$\\
$\ZsmT{mint0}=\tupic{mint0}\ +\ I\,\tupic{mint23}\ +\ J\,\tupic{mint1223}\
    +\ K\,\tupic{mint2323}\ +\ L\tupic{mint1323}\ .$

Show that the order 1 terms of the product vanish.

The only nonzero chord diagram of order 2 on the hump is the cross
(diagram without isolated chords). The coefficient of this diagram
is $B+D+G+J+L-AE+AI+EI$. Show that it is equal to
$$\textstyle
\frac{\Li\bigl(\frac{\e}{2-\e}\bigr)
 - \Li\bigl(\frac{-\e}{2-\e}\bigr) + \Li\bigl(\frac{1}{2}\bigr)
 - \Li\bigl(\frac{1}{2-\e}\bigr) - \Li(\e)}{2\pi^2}
 + \frac{\ln^2 2 - \ln^2(2-\e)}{4\pi^2} + \frac{1}{24}\ .
$$
Using the properties of the dilogarithm mentioned in problem \ref{EuLi}
prove that the last expression equals $\frac{1}{24}$. This is also a
consequence of the remarkable {\it Rogers five-term relation}
\index{Rogers five-term relation}
(see, for example, \cite{Kir})
$$\textstyle
\Li x+ \Li y -\Li xy= \Li\frac{x(1-y)}{1-xy}
+  \Li\frac{y(1-x)}{1-xy}
+ \ln\frac{(1-x)}{1-xy} \ln\frac{(1-y)}{1-xy}
$$
and the {\it Landen connection formula}
\index{Landen connection formula}
(see, for example, \cite{Roos})
$$\textstyle
\Li z + \Li\frac{-z}{1-z} = -\frac{1}{2}\ln^2(1-z)\ .$$

\item\label{op_S_i}
Let $S_i$ be the operation of reversing the orientation of the $i$th
component of a tangle $T$. Denote by the same symbol $S_i$ the
operation on tangle chord diagrams which (a) reverses the $i$th
component of the skeleton of a diagram; (b) multiplies the diagram
by $-1$ raised to the power equal to the number of chord endpoints
lying on the $i$th component. \label{op-S_i}\index{Operation!$S_i$
on tangle chord diagrams} Prove that
$$Z(S_i(T))=S_i(Z(T))$$
We shall use this operation in Chapter \ref{DA_chap}.

\item\label{ki_Phi}
\parbox[t]{2.2in}{Compute the Kontsevich integral
$Z(AT^t_{b,w})$ up to the order 2. Here $\e$ is a small parameter, and
$w$, $t$, $b$ are natural numbers subject to inequalities $w<b$ and $w<t$.
\\
{\sl Answer.}\qquad $Z(AT^t_{b,w})\ =\ \vstemp\ +$
}
\qquad
\parbox{2.1in}{\makebox(145,0){\rb{-120pt}{$
 AT^t_{b,w} = \risS{-50}{assate}{
         \put(20,18){\mbox{$\scriptstyle \e^b$}}
         \put(43,11){\mbox{$\scriptstyle \e^w$}}
         \put(69,84){\mbox{$\scriptstyle \e^t$}}
                 }{90}{40}{50}$}}}\\
$\hspace{25pt}+\
\frac{1}{2\pi i} \ln\left(\frac{\e^w-\e^t}{\e^b}\right) \vstod\ -\
\frac{1}{2\pi i} \ln\left(\frac{\e^w-\e^b}{\e^t}\right) \vstdt$\\
$\hspace{25pt}-\
\frac{1}{8\pi^2}\ln^2\left(\frac{\e^w-\e^t}{\e^b}\right) \vstodod\ -\
\frac{1}{8\pi^2}\ln^2\left(\frac{\e^w-\e^b}{\e^t}\right) \vstdtdt$\\
$\hspace{25pt} - \frac{1}{4\pi^2}\left(
   \ln(\e^{b-w})\ln\left(\frac{\e^w-\e^b}{\e^t}\right)
   + \Li(1-\e^{b-w}) - \Li(\e^{t-w})
   \right)\vstoddt$\\
$\hspace{25pt} + \frac{1}{4\pi^2}\left(
   \ln(1-\e^{t-w})\ln\left(\frac{\e^w-\e^b}{\e^t}\right)
   + \Li(1-\e^{b-w}) - \Li(\e^{t-w})
   \right)\vstdtod\ .$

\medskip
\item\label{parprfo}
The set of elementary tangles can be expanded by adding crossings
with arbitrary orientations of strands.  Express the figure eight
knot $4_1$ in terms of associating tangles and $\e$-parametrized
tensor products of elementary (in this wider sense) tangles in the
same manner as the trefoil $3_1$ is described in Figure \ref{parenth_tref_ch8}.\\
{\sl Answer:} A possible answer is shown in
Figure~\ref{parenth_eight}.
\begin{figure}[ht]
\begin{center}
\begin{picture}(110,195)(0,5)
{\put(0,0){\ig[width=150pt]{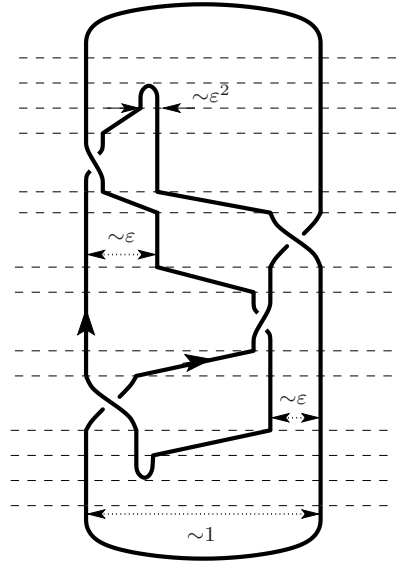}}
      \put(67,8){\mbox{$\scriptstyle \sim 1$}}
      \put(102,60){\mbox{$\scriptstyle \sim\e$}}
      \put(37,121){\mbox{$\scriptstyle \sim\e$}}
      \put(69,174){\mbox{$\scriptstyle \sim\e^2$}}
}
\end{picture}
\end{center}
\caption{The figure eight knot in terms of elementary tangles and
associating tangles.} \label{parenth_eight}
\end{figure}

\end{enumerate}
\end{xcb}
 %8 KI
\chapter{Framed knots and cabling operations} %9
\label{chap:operations}

In this chapter we show how to associate to a framed knot $K$ an
infinite set of framed knots and links, called the {\em
$(p,q)$-cables} of $K$. The operations of taking the $(p,q)$-cable
respect the Vassiliev filtration, and give rise to operations on
Vassiliev invariants and on chord diagrams. We shall give explicit
formulae that describe how the Kontsevich integral of a framed knot
changes under the cabling operations. As a corollary, this will give
an expression for the Kontsevich integral of all torus knots.

\section{Framed version of the Kontsevich integral}
\label{frKI} \index{Kontsevich integral!framed}

In order to describe a framed knot one only needs to specify the
corresponding unframed knot and the self-linking number. This
suggests that there should be a simple formula to define the
universal Vassiliev invariant for a framed knot via the Kontsevich
integral of the corresponding unframed knot. This is, indeed, the
case, as we shall see in Section~\ref{subsec:fra-unfra}. However,
for our purposes it will be more convenient to use a definition of
the framed Kontsevich integral given by  V.~Goryunov in \cite{Gor1},
which is in the spirit of the original formula of Kontsevich
described in Section~\ref{defki}.

\begin{xremark} For framed knots and links, the universal Vassiliev
invariant was first defined by Le and Murakami \cite{LM2} who gave a
combinatorial construction of it using the Drinfeld associator (see
Chapter~\ref{DA_chap}). Goryunov used his framed Kontsevich integral
in \cite{Gor2} to study Arnold's $J^+$-theory of plane curves (or,
equivalently, Legendrian knots in a solid torus).
\end{xremark}

\subsection{}\label{Gappr}
Let $K_\e$ be a copy of $K$ shifted by a small distance $\e$ in the
direction of the framing. We assume that both $K$ and $K_\e$ are in
general position with respect to the height function $t$ as in
Section~\ref{defki}. Then we construct the (preliminary) integral
$Z(K,K_\e)$ defined by the formula
$$Z(K,K_\e)=
\sum_{m=0}^\infty \frac{1}{(4\pi i)^m}  \hspace{-8pt}
         \int\limits_{\substack{\vspace{5pt} \\
    t_{\mbox{\tiny min}}<t_m<\dots<t_1<t_{\mbox{\tiny max}}\\
    t_j\mbox{\tiny\ are noncritical}}}
         \sum_{P=\{(z_j,z'_j)\}}\hspace{-8pt}  (-1)^\downarrow D_P
         \bigwedge_{j=1}^m \frac{dz_j-dz'_j}{z_j-z'_j} \ ,
$$
whose only difference with the formula for the unframed Kontsevich
integral is the numerical factor in front of the integral. However,
the notation here has a different meaning. The class of the diagram
$D_P$ is taken in $\A^{fr}$ rather than in $\A$. We consider only
those pairings $P=\{(z_j,z'_j)\}$ where $z_j$ lies on $K$ while
$z'_j$ lies on $K_\e$. In order to obtain the chord diagram $D_P$ we
project the chord endpoints that lie on $K_{\e}$ back onto $K$ along
the framing. If an endpoint $z'_j$ projects exactly to the point
$z_j$ on $K$, we place a ``small'' isolated chord in a neighbourhood
of $z_j$. The following picture illustrates this definition:
$$\risS{-110}{frtrdi}{
      \put(-4,125){\mbox{$\scriptstyle t$}}
      \put(68,15){\mbox{$\scriptstyle K$}}
      \put(13,30){\mbox{$\scriptstyle K_\e$}}
      \put(-6,45){\mbox{$\scriptstyle t_3$}}
      \put(85,51){\mbox{$\scriptstyle z_3$}}
      \put(102,45){\mbox{$\scriptstyle z'_3$}}
      \put(-6,55){\mbox{$\scriptstyle t_2$}}
      \put(15,59){\mbox{$\scriptstyle z'_2$}}
      \put(48,55){\mbox{$\scriptstyle z_2$}}
      \put(-6,85){\mbox{$\scriptstyle t_1$}}
      \put(107,85){\mbox{$\scriptstyle z'_1$}}
      \put(29,88){\mbox{$\scriptstyle z_1$}}
      \put(170,25){\mbox{$\scriptstyle D_P$}}
      \put(195,55){\mbox{$\scriptstyle z_3$}}
      \put(187,45){\mbox{$\scriptstyle z'_3$}}
      \put(155,97){\mbox{$\scriptstyle z'_2$}}
      \put(155,42){\mbox{$\scriptstyle z_2$}}
      \put(194,85){\mbox{$\scriptstyle z'_1$}}
      \put(138,65){\mbox{$\scriptstyle z_1$}}
             }{200}{20}{110}$$
Now the framed Kontsevich integral can be defined as
$$Z^{fr}(K) = \lim\limits_{\e\to 0} Z(K,K_\e)\,.$$

In \cite{Gor1} V.~Goryunov proved that the limit does exist and is
invariant under the deformations of the framed knot $K$ in the class
of framed Morse knots. The proof is, essentially, the same as that
of the unframed case, with only one difference: since a diagram with
an isolated chord is not necessarily zero in $\A^{fr}$ we have to be
careful about the case when an isolated chord approaches a maximum
or minimum. The key observation here is that the diverging pairs
$(z,z')$ and $(z,z'')$ cancel each other in the limit:
$$\risS{-15}{frcan}{
      \put(-3,7){\mbox{$\scriptstyle z$}}
      \put(49,7){\mbox{$\scriptstyle z'$}}
      \put(17,10){\mbox{$\scriptstyle z''$}}
             }{50}{0}{15}$$
We refer to \cite{Gor1} for details.
\begin{xexample}
Let $O^{+m}$ be the $m$-framed unknot:
$$\ig[height=1.5in]{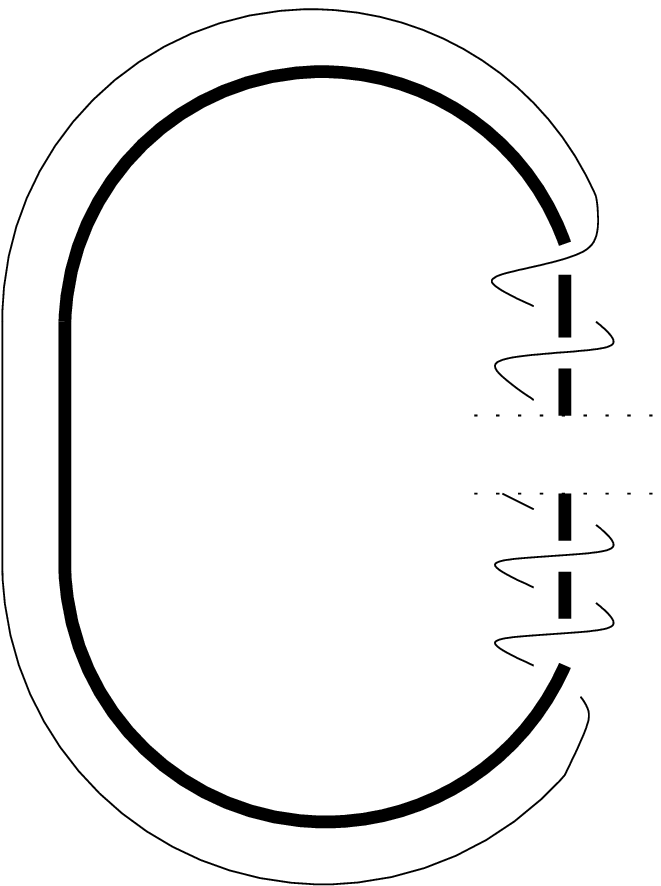}$$
Then $$Z^{fr}(O^{+m})=\exp{\frac{m\Theta}{2}}.$$
\end{xexample}
\begin{xexample}
The integral formula for the linking number in
\ref{subsec:linking_number} shows that the coefficient of the
diagram $\Theta$ in $Z^{fr}(K)$ is equal to $w(K)/2$ where $w(K)$ is
the self-linking number of $K$.
\end{xexample}

Define the final framed Kontsevich integral as
$$I^{fr}(K) = \frac{Z^{fr}(K)}{Z^{fr}(H)^{c/2}}\ ,$$
where $H$ is the zero-framed hump unknot (see page~\pageref{hump}).
With its help one proves the framed version of
Theorem~\ref{Ko_part_fund_th}:
\begin{xtheorem}\label{frKo_part_fund_th}
Let $w$ be a framed $\C$-valued weight system of order $n$. Then
there exists a framed Vassiliev invariant of order $\leq n$ whose
symbol is $w$.
\end{xtheorem}

\subsection{The relation with the unframed integral}\label{subsec:fra-unfra}
\begin{xproposition}
The image of the framed Kontsevich integral $Z^{fr}(K)$ under the
quotient map $\widehat{\A}^{fr}\to \widehat{\A}$ is the unframed
Kontsevich integral $Z(K)$.
\end{xproposition}
\begin{proof}
Each horizontal chord with endpoints on $K$ can be lifted to a chord
with one end on $K$ and the other on $K_\e$ in two possible ways.
Therefore, each pairing $P$ with $m$ chords for the unframed
Kontsevich integral comes from $2^m$ different pairings for the
framed integral. As $\e$ tends to zero, each of these pairings gives
the same contribution to the integral as $P$ and its coefficient is
precisely $(2\pi i)^{-m}/2^m=(4\pi i)^{-m}$.
\end{proof}

In fact, we can prove a much more precise statement. As we have seen
in Section~\ref{defram_cd}, the algebra of chord diagrams $\A$ can
be considered as a direct summand of $\A^{fr}$. This allows us to
compare the framed and the unframed Kontsevich integrals directly.
\begin{xtheorem}\label{th:framed} Let $K$ be a framed knot with self-linking number $w(K)$.
Then $$Z^{fr}(K)=Z(K)\cdot\exp{\frac{w(K)\Theta}{2}}$$ where $Z(K)$
is considered as an element of $\widehat{\A}^{fr}$.
\end{xtheorem}
This statement can be taken as a definition of the framed Kontsevich
integral.

\begin{proof}
Recall that $\A$ is identified with a direct summand of $\A^{fr}$ by
means of the algebra homomorphism $p:\A^{fr}\to\A^{fr}$ whose kernel
is the ideal generated by the diagram $\Theta$, and which is defined
on a diagram $D$ as
$$
  p(D)=\sum_{J\subseteq [D]} (-\Theta)^{\deg{D}-|J|}\cdot D_J\ ,
$$
see Section~\ref{defram_cd}. We shall prove that
\begin{equation}\label{formula:pK}
p(Z^{fr}(K))=Z^{fr}(K)\cdot\exp{-\frac{w(K)\Theta}{2}},
\end{equation}
which will imply the statement of the theorem.

Write $p(D)$ as a sum $\sum_k (-1)^k \Theta^k\cdot p_{(k)}(D)$ where
the action of $p_{(k)}$ consists in omitting $k$ chords from a
diagram in all possible ways: $$p_{(k)}(D)=\sum_{J\subseteq [D],\,
\deg{D}-|J|=k} D_J.$$ We have $p_{(k)}(Z^{fr}(K))= \sum c_P D_P$
where the sum is taken over all possible pairings $P$. The
coefficient $c_P$ is equal to the sum of all the coefficients in
$Z^{fr}(K)$ that correspond to pairings $P'$ obtained from $P$ by
adding $k$ chords. These chords can be taken arbitrarily, so,
writing $m$ for the degree of $P$ we have
$$c_P=\frac{1}{(4\pi i)^{m+k}} \hspace{-25pt}
         \int\limits_{\substack{\vspace{9pt} \\
    t_{\mbox{\tiny min}}<t_{m}<\dots<t_1<t_{\mbox{\tiny max}}\\
    t_{\mbox{\tiny min}}<\tau_{k}<\dots<\tau_1<t_{\mbox{\tiny max}}}}
         \hspace{-25pt} \sum_{P=\{(z_j,z'_j)\}}  \hspace{-8pt} (-1)^\downarrow
         \bigwedge_{j=1}^{m} d\ln{(z_j-z'_j)} \wedge \bigwedge_{i=1}^{k}
         d\ln{(\zeta_i-\zeta'_i)}\,
         ,
$$
where all $t_j$ and $\tau_i$ are non-critical and distinct, $z_j$
and $z'_j$ depend on $t_j$ and $\zeta_i$ and $\zeta'_i$ --- on
$\tau_i$. This expression is readily seen to be a product of two
factors: the coefficient at $D_P$ in $Z^{fr}(K)$ and
$$\frac{1}{(4\pi i)^{k}} \hspace{-15pt}
         \int\limits_{\substack{\vspace{5pt} \\t_{\mbox{\tiny min}}<\tau_{k}<\dots<\tau_1<t_{\mbox{\tiny
         max}}}}
         \hspace{-10pt} \sum_{P'=\{(\zeta_i,\zeta'_i)\}}  \hspace{-15pt} (-1)^\downarrow \bigwedge_{i=1}^{k}
         d\ln{(\zeta_i-\zeta'_i)}\, .
$$
The latter expression is equal to
$$
         \frac{1}{k!\cdot(4\pi i)^{k}} \hspace{-15pt}
         \int\limits_{\substack{\vspace{5pt} \\t_{\mbox{\tiny min}}<\tau_{1},\dots,\tau_k <t_{\mbox{\tiny
         max}}}}
         \hspace{-10pt} \sum_{P'=\{(\zeta_i,\zeta'_i)\}}  \hspace{-15pt} (-1)^\downarrow \bigwedge_{i=1}^{k}
         d\ln{(\zeta_i-\zeta'_i)}= \frac{1}{k!}\Bigl(\frac{w(K)}{2}\Bigr)^k,
$$
so that
$$p_{(k)}(Z^{fr}(K))=\frac{1}{k!}\Bigl(\frac{w(K)}{2}\Bigr)^k\cdot Z^{fr}(K),$$
and (\ref{formula:pK}) follows.
\end{proof}

\subsection{The case of framed tangles}\label{subsec:kifrt}
The above methods produce the Kontsevich integral not just for
framed knots, but, more generally, for framed tangles. The
preliminary integral $Z^{fr}(T)$ of a tangle $T$ can be constructed
just as in the case of knots, and the final integral $I^{fr}(T)$ is
defined as
$$I^{fr}(T)=Z^{fr}(H)^{-m_1}\#\ldots\# Z^{fr}(H)^{-m_k}\#
Z^{fr}(T),$$ where  $m_i$ is the number of maxima on the $i$th
component of $T$ and $Z^{fr}(H)^{-m_i}$ acts on the $i$th component
of $Z^{fr}(T)$ as defined in \ref{subsec:modules-tangle}. Here $k$ is the
number of components of $T$.

Note that the final integral $I^{fr}(T)$ is multiplicative with
respect to the tangle product, but not the connected sum of knots.

\section{Cabling operations}
\subsection{Cabling operations on framed knots}
\label{def:cabling} Let $p,q$ be two coprime integers with $p\neq
0$, and $K$ be a framed knot given by an embedding $f:S^1\to\R^3$
with the framing vector $v(\theta)$ for $\theta\in S^1$. Denote by
$r_\alpha v(\theta)$ the rotation of the vector $v(\theta)$ by the
angle $\alpha$ in the plane orthogonal to the knot. Then, for all
sufficiently small values of $\e$, the map
$$\theta \to f(p\theta)+\e\cdot r_{q\theta} v(p\theta)$$ is actually a knot. This knot is called
the {\em $(p,q)$-cable of $K$}; we denote it by $K^{(p,q)}$. Note
that $q$ is allowed to be zero: $K^{(1,0)}$ is $K$ itself, and
$K^{(-1,0)}$ is the inverse $K^*$.

\begin{xexample} Here is the left trefoil with
the blackboard framing and its $(3,1)$-cable:
$$\ig[height=20mm]{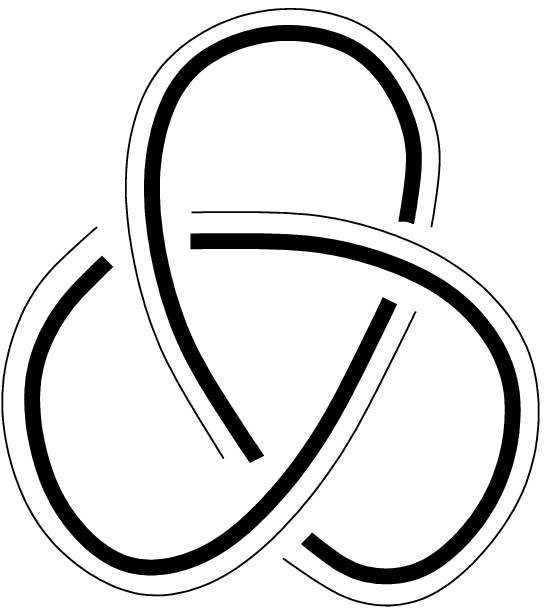}\qquad\quad\ig[height=20mm]{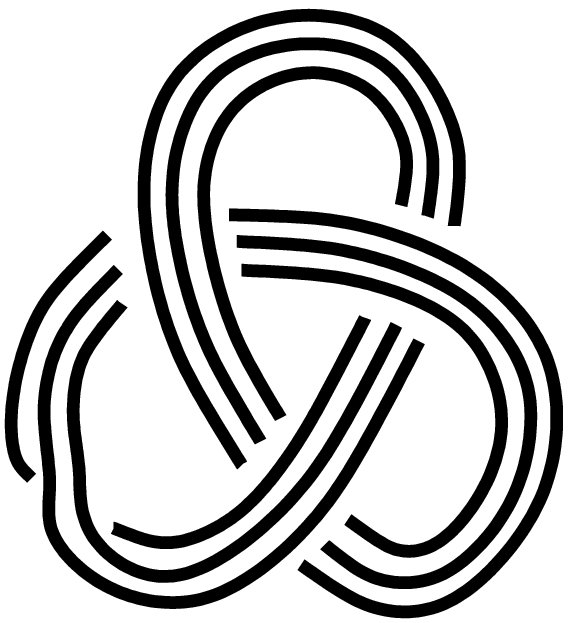}$$
\end{xexample}

The $(p,q)$-cables can, in fact, be defined for {\em arbitrary}
integers $p,q$ with $p\neq 0$, as follows. Take a small tubular
neighbourhood $N$ of $K$. On its boundary there are two
distinguished simple closed oriented curves: the {\em meridian},
which bounds a small disk perpendicular to the knot\footnote{This
defines the meridian up to isotopy.} and is oriented so as to have
linking number one with $K$, and the {\em longitude}, which is
obtained by shifting $K$ to $\partial N$ along the framing. A closed
oriented curve without self-intersections on $\partial N$ which
represents $p$ times the class of the longitude plus $q$ times the
class of the meridian in $H_1(\partial N)$ is the $(p,q)$-cable of
$K$. In general, $K^{(p,q)}$ is a knot if and only if $p$ and $q$
are relatively prime; otherwise, it is a link with more than one
component. The number of components of the resulting links is
precisely the greatest common divisor of $p$ and $q$. Sometimes, the
$(k,0)$-cable of $K$ is called {\em the $k$th disconnected cabling
of $K$} and the $(k,1)$-cable {\em the $k$th connected cabling
of $K$}.

We shall consider $K^{(p,q)}$ as a framed link with the framing
normal to $\partial N$ and pointing outwards:
$$\ig[height=35mm]{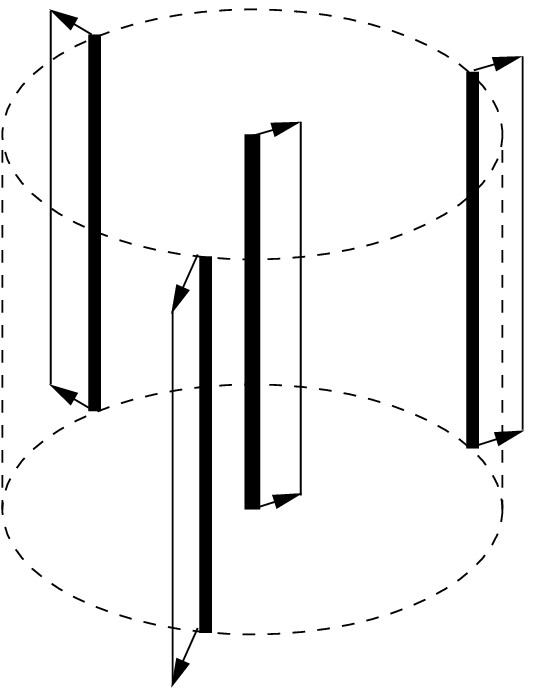}$$

\begin{xexample} The $(p,q)$-torus knot
(link) \index{Knot!torus} can be defined as the $(p,q)$-cable of the
zero-framed unknot.
\end{xexample}

\subsection{Cables and Vassiliev invariants}
\label{vassiliev_cables} Composing a link invariant with a cabling
operation on knots we obtain a new invariant of (framed) knots.
\begin{xproposition}
Let $p,q$ be a pair of integers and $r$ be their greatest common
divisor. If $v$ is a Vassiliev invariant of framed $r$-component
links whose degree is at most $n$, the function $v^{(p,q)}$ sending
a framed knot $K$ to $v(K^{(p,q)})$ is an invariant of degree $\leq
n$.
\end{xproposition}
\begin{proof}
Indeed, the operation of taking the $(p,q)$-cable sends the singular
knot filtration on $\Z\K^{fr}$ into the filtration by singular links
on the free abelian group generated by the $r$-component framed
links, since the difference
$$\rb{-9mm}{\ig[height=20mm]{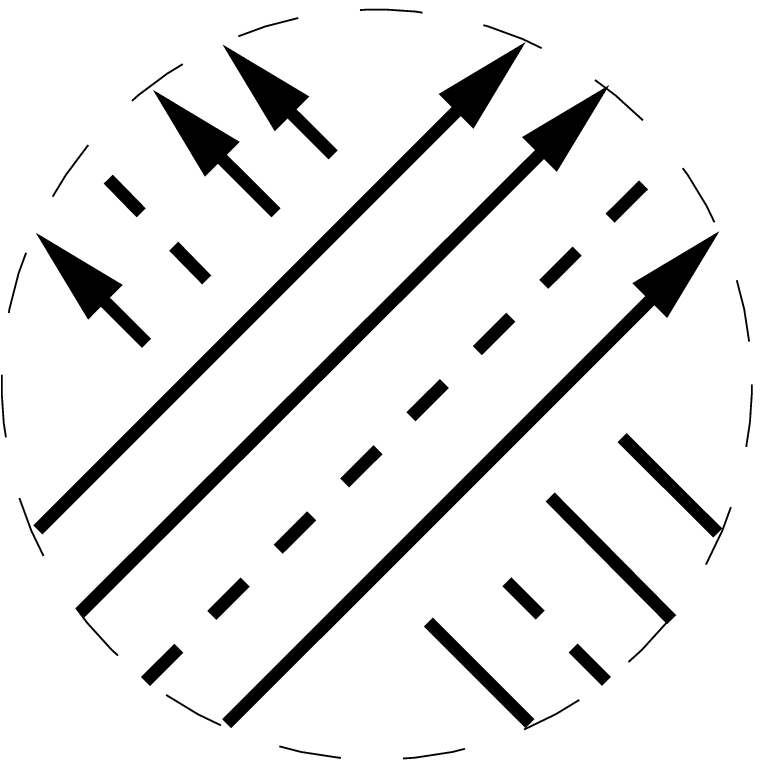}}\ -\ \rb{-9mm}{\ig[height=20mm]{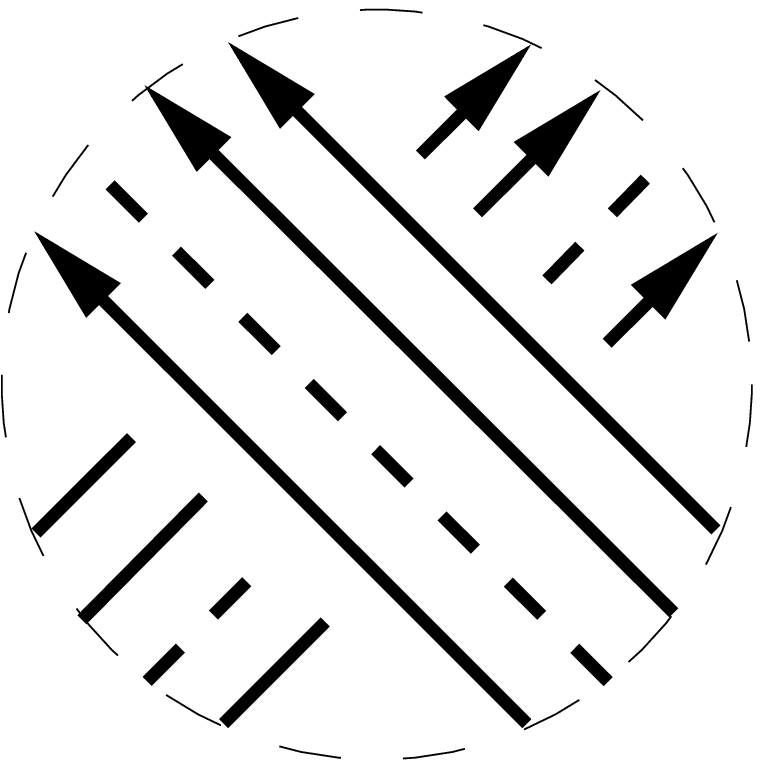}}$$
can be written as a sum of several double points. For instance,
$$\rb{-5mm}{\ig[height=12mm]{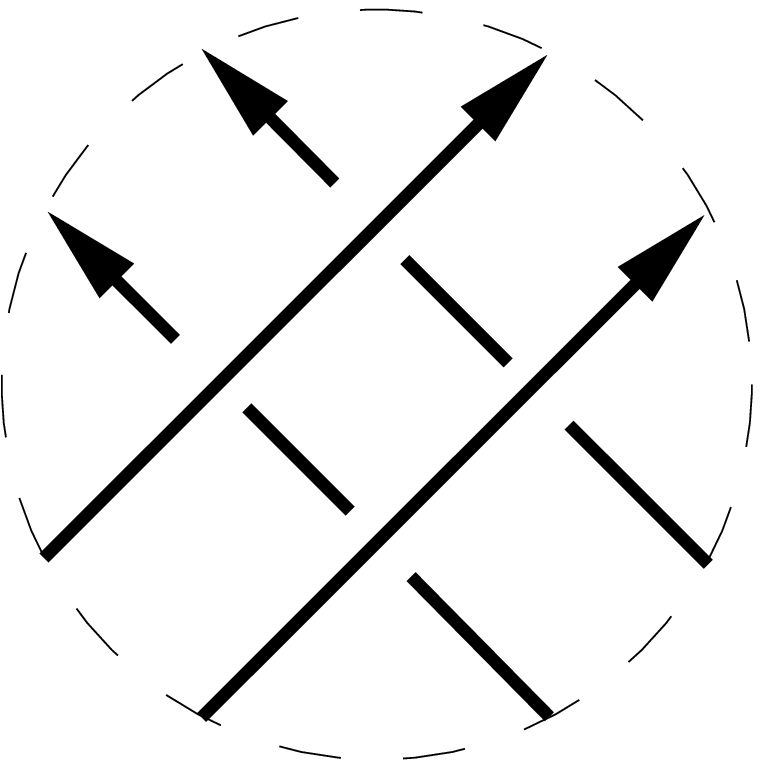}}\ -\ \rb{-5mm}{\ig[height=12mm]{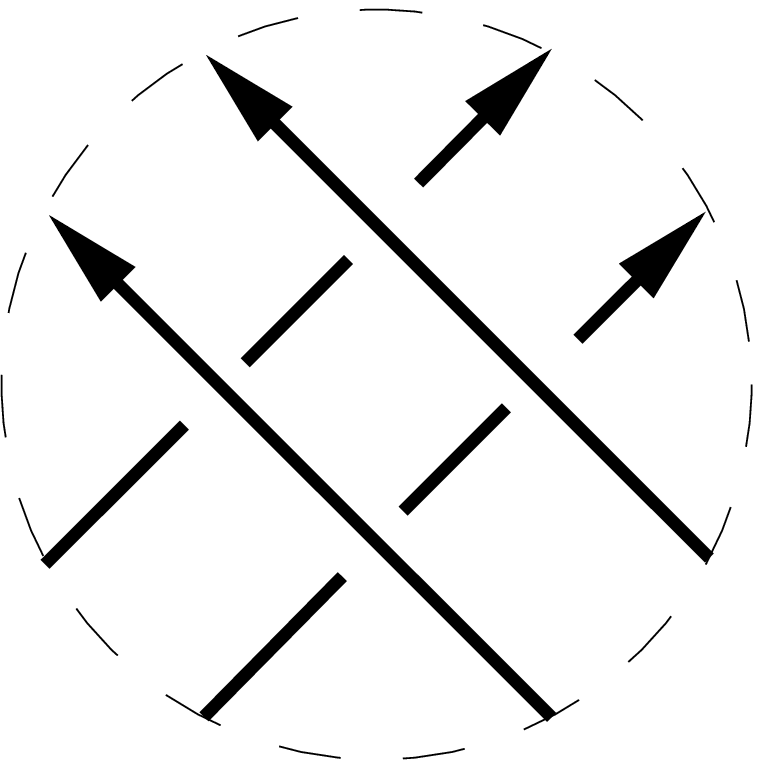}}\ =
\ \rb{-5mm}{\ig[height=12mm]{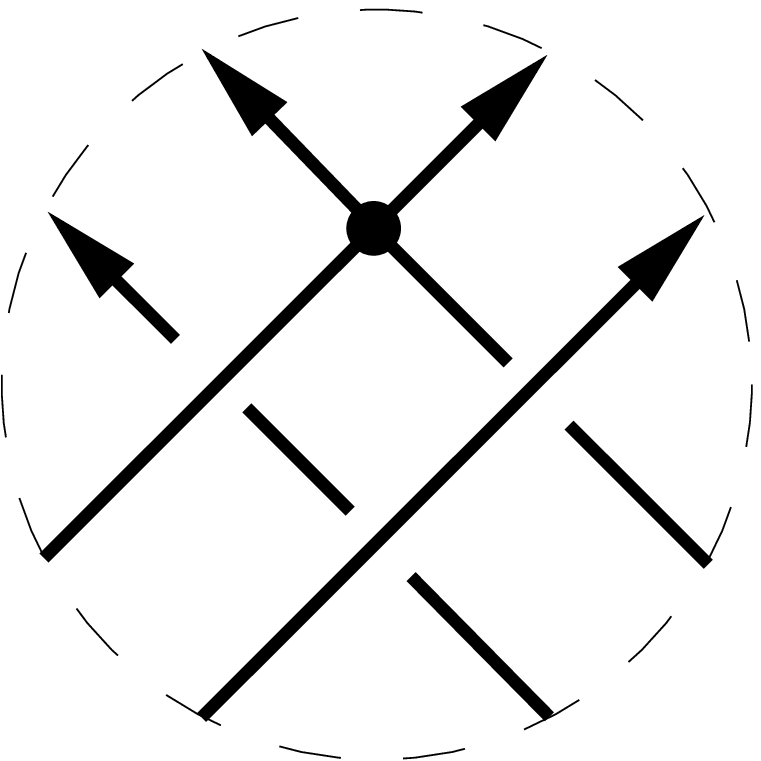}}\ +\
\rb{-5mm}{\ig[height=12mm]{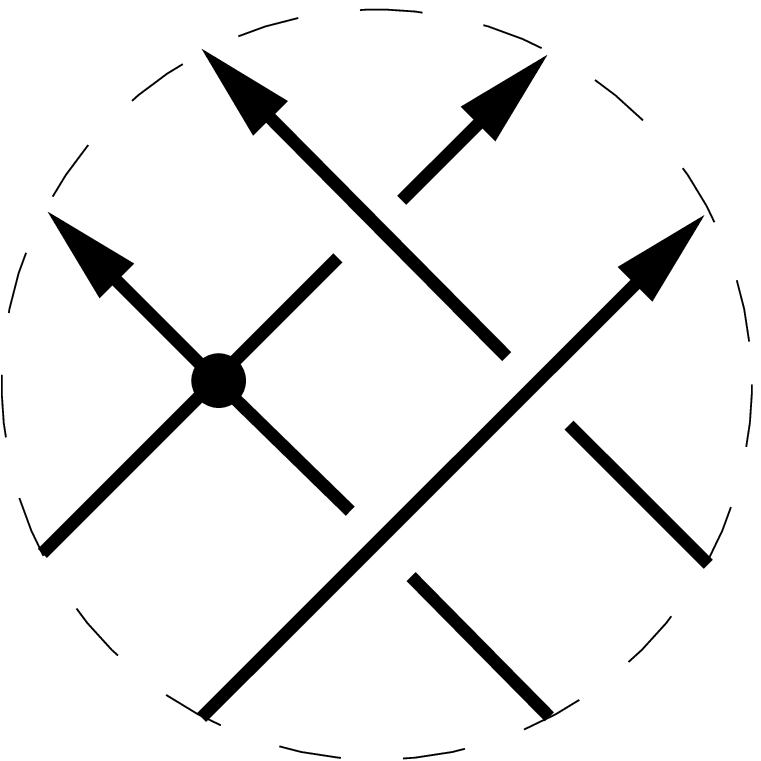}}\ + \
\rb{-5mm}{\ig[height=12mm]{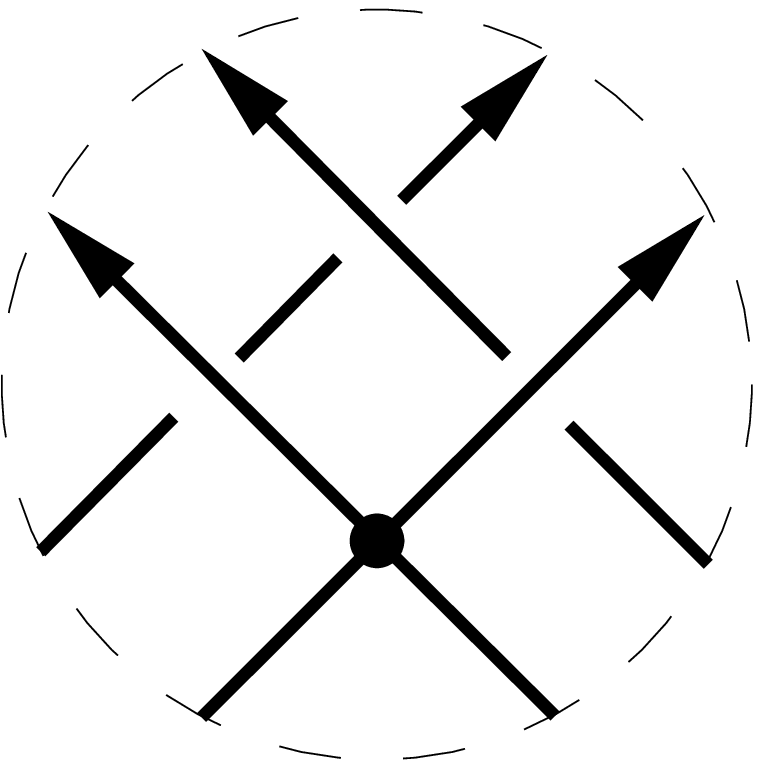}}\ +\
\rb{-5mm}{\ig[height=12mm]{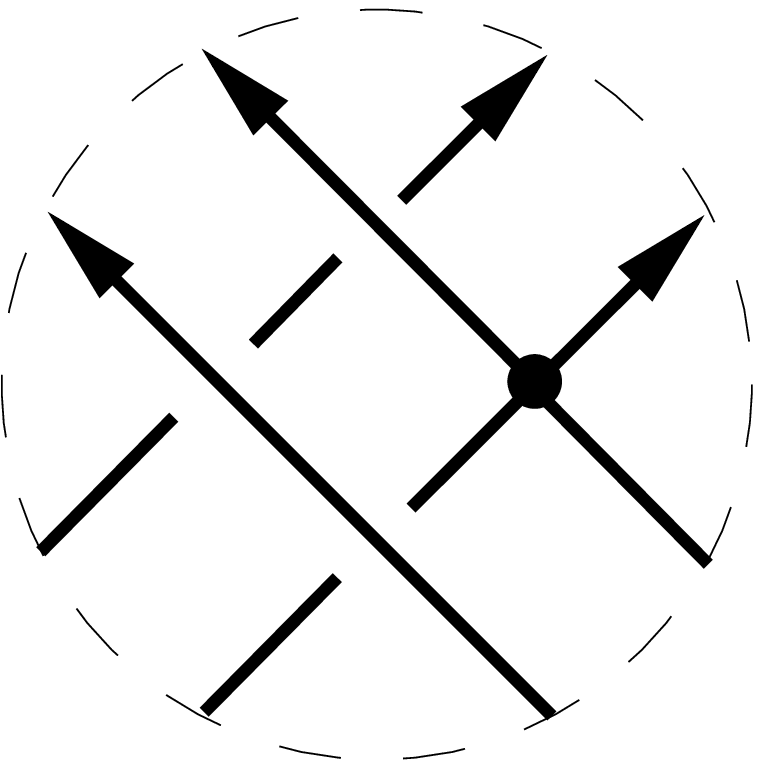}}\, .
$$
\end{proof}
It is clear from the above argument what effect the cabling
operation has on chord diagrams. Consider first the case of $p$ and
$q$ coprime, when the $(p,q)$-cabling gives an operation on knot
invariants. For a chord diagram $D$ define $\psi^p(D)$ to be the sum
of chord diagrams obtained by all possible ways of lifting the ends
of the chords to the $p$-sheeted connected covering of the Wilson
loop of $D$.\label{psiop}

\newcommand{\cabldp}[1]{\rb{-4.2mm}{\ig[width=10mm]{#1.eps}}}
\newcommand{\cabldps}[1]{\rb{-1mm}{\ig[width=10mm]{#1.eps}}}
\newcommand{\cabldpn}[1]{\rb{-4.2mm}{\ig[width=7mm]{#1.eps}}}
\newcommand{\wtcabl}[1]{\rb{-4.2mm}{\ig[width=10mm]{#1.eps}\hspace{-11.5mm}
                                    \ig[width=10mm]{wt_cab_psi2.eps}}}

\noindent{\bf Examples.}
$$\begin{array}{rcl}
\psi^2(\cabldp{fd2b})
&=&\cabldp{cabcd1}+\cabldp{cabcd2}+\cabldp{cabcd3}+
  \cabldp{cabcd4}+\cabldp{cabcd5}+\cabldp{cabcd6}  \vspace{8pt}\\
&&+\cabldp{cabcd7}+\cabldp{cabcd8}+\cabldp{cabcd9}+\cabldp{cabcd10}
  +\cabldp{cabcd11}+\cabldp{cabcd12}   \vspace{8pt}\\
&&+\cabldp{cabcd13}+\cabldp{cabcd14}+\cabldp{cabcd15}+\cabldp{cabcd16}
         \vspace{10pt}\\
&=& 12\,\cabldp{fd2b}\ +\ 4\,\cabldp{fd2a}\,. \vspace{10pt}\\
\psi^2(\cabldp{fd2a})&=&\ 8\ \cabldp{fd2b}\ +\ 8\,\cabldp{fd2a}\,.
\end{array}
$$ It is a simple exercise to see that $\psi^p$
respects the 4T relations; hence, it gives a linear map
$\psi^p:\A^{fr}\to\A^{fr}.$ We have the following
\begin{xproposition}  $\symb(v^{p,q})= \symb (v) \circ \psi_p.$
\end{xproposition}
Note that the symbol of $v^{(p,q)}$ does not depend on $q$.

The case when $p$ and $q$ are not coprime and the $(p,q)$-cable is a
link with at least two components, is very similar. We shall treat
this case in a slightly more general setting in
Section~\ref{diagram_cabling_general}.

\subsection{Cabling operations in $\F$ and $\B$}
\label{diagram_cabling_general} \index{Cabling!of
closed diagrams} \index{Cabling!of open diagrams} The map $\psi^p$
is defined on general closed diagrams in the very same way as on
chord diagrams: it is the sum of all possible liftings of the legs
of a diagram to the $p$-fold connected cover of its Wilson loop. It
is not hard to see that $\psi^p$ defined in this manner respects the
STU relation. For instance,
\def\cabClD#1{\risS{-12}{#1}{}{40}{20}{20}}
$$\begin{array}{rcl}
\psi^2\Bigl(\cabClD{sstuA}\Bigr) &=&
\psi^2\Bigl(\cabClD{tstuA}\Bigr)
        \ -\ \psi^2\Bigl(\cabClD{ustuA}\Bigr) \\
&=& \cabClD{cabstu1}\ +\ \cabClD{cabstu2}\ +\ \cabClD{cabstu3}\ +\
    \cabClD{cabstu4} \\
&&-\ \cabClD{cabstu5}\ -\ \cabClD{cabstu6}\ -\ \cabClD{cabstu7}\ -\
    \cabClD{cabstu8} \\
&=& \cabClD{cabstu1}\ -\ \cabClD{cabstu5}\ +\ \cabClD{cabstu4}\ -\
    \cabClD{cabstu8} \\
&=& \cabClD{cabstu9}\ +\ \cabClD{cabstu10}\ .
\end{array}$$
Therefore, $\psi^p$ is a well-defined map of $\F$ to itself. Note
that $\psi^p$ is a coalgebra map; however, it does not respect the
product in $\F$. This is hardly surprising since the cabling maps in
general do not respect the connected sum of knots.

The algebra $\B$ is better suited for working with the cabling
operations than $\F$: the map $\psi^p$ applied to an open diagram
with $k$ legs simply multiplies this diagram by $p^k$. Indeed, the
isomorphism $\chi:\B\isom\F$ takes an open diagram $B$ with $k$ legs
into the average of the $k!$ closed diagrams obtained by all
possible ways of attaching the legs of $B$ to a Wilson loop. Lifting
this average to the $p$-fold covering of the Wilson loop we get the
same thing as $p^k \chi (B)$. We arrive to the following
\begin{xproposition}
The operation $\psi^p:\B\to\B$ is an Hopf algebra map. In particular, the subspace
$\B^k$ of diagrams with $k$ legs is the eigenspace for $\psi^p$ with
the corresponding eigenvalue $p^k$.
\end{xproposition}
The fact that $\psi^p$ respects the product on $\B$ follows from the
second part of the Proposition.

\begin{remark}\label{rem:cablehump}
Even though $\psi^p$ does not respect the product in $\F$, there is
one instance where it is exhibits simple behaviour with respect to
the connected sum. Namely, if $H$ is the zero-framed hump unknot, we
have $$\psi^p (D \# Z^{fr}(H))=\psi^p(D) \# Z^{fr}(H)^{\# p}$$ as an
element of $\widehat{\F}$. This is due to the fact that if $K$ is a
framed Morse knot, $(K\# H)^{(p,q)}$ is obtained from $K^{(p,q)}$ by
taking the connected sum with $p$ copies of $H$.
\end{remark}

\subsection{Cablings on tangle diagrams}
So far we have only considered the effect of the $(p,q)$-cables on
chord diagrams for coprime $p$ and $q$. However, there is no
difficulty in extending our results to the case of arbitrary $p$ and
$q$.

Given a framed tangle $T$ with a closed component $\yy$, we can
define its {\em $(p,q)$-cable along $\yy$}, denoted by
$T^{(p,q)}_{\yy}$ in the same manner as for knots. If $p,q$ are
coprime the result will have the same skeleton as the original
tangle, otherwise the component $\yy$ will be replaced by several
components whose number is the greatest common divisor of $p$ and
$q$.

If $p'=rp$ and $q'=rq$ with $p$ and $q$ coprime, the map
$\psi^{r\cdot p}_{\yy}$ corresponding to the $(p',q')$-cable on the
space of closed Jacobi diagrams with the skeleton $\boldX\cup\yy$
can be described as follows. Consider the map
$$\boldX\cup \yy_1\cup\ldots\cup \yy_r\to \boldX\cup \yy$$
where $\yy_i$ are circles, which sends $\boldX$ to $\boldX$ by the
identity map and maps each $\yy_i$ to $\yy$ as a $p$-fold covering. Then
$\psi^{r\cdot p}_{\yy}$ of a closed diagram $D$ is the sum of all
the different ways of lifting the legs of $D$ to $\boldX\cup_i
\yy_i$. For example,
\newcommand{\cabldph}[1]{\rb{-4.2mm}{\ig[height=10mm]{#1.eps}}}
$$
\begin{array}{c}
{
\begin{array}{rcl}
\psi^{2\cdot 1}\Bigl(\cabldp{fd2b}\Bigr)
&=&\cabldp{cabdcd1}+\cabldp{cabdcd2}+\cabldp{cabdcd3}+
  \cabldp{cabdcd4}+\cabldp{cabdcd5}+\cabldp{cabdcd6}  \vspace{8pt}\\
&&+\cabldp{cabdcd7}+\cabldp{cabdcd8}+\cabldp{cabdcd9}+\cabldp{cabdcd10}
  +\cabldp{cabdcd11}+\cabldp{cabdcd12}   \vspace{8pt}\\
&&+\cabldp{cabdcd13}+\cabldp{cabdcd14}+\cabldp{cabdcd15}+\cabldp{cabdcd16}
         \vspace{10pt}
\end{array}
}\\
{=2\ \cabldph{licd2-20a} + 8\ \cabldph{licd2-3-1a} + 2\
\cabldph{licd2-22a}
     + 4\ \cabldph{licd2-2-2a}\ .}
\end{array}
$$
and
\begin{multline*}
\psi^{2\cdot 1}\Bigl(\cabldp{fd2a}\Bigr)\ =\ 2\ \cabldph{licd2-xoa}
+ 8\,\cabldph{licd2-3-1a}
 \\ + 2\ \cabldph{licd2-22a} +
 4\ \cabldph{licd2-2-2a}\ .
\end{multline*}
Here we have omitted the subscript indicating the component $\yy$,
since the original diagram had only one component. In what follows,
we shall write $\psi^{p}_{\yy}$ instead of $\psi^{1\cdot p}_{\yy}$.

As in Section~\ref{vassiliev_cables}, the $(p,q)$-cable along $\yy$
composed with a Vassiliev invariant $v$ of degree $n$ is again a
Vasiliev invariant of the same degree, whose symbol is obtained by
composing $\psi^{r\cdot p}_{\yy}$ with the symbol of $v$. The map
$\psi^{r\cdot p}_{\yy}$ satisfies the 4T relations and gives rise to
a coalgebra map on the spaces of closed diagrams.

\subsection{Disconnected cabling in $\B$}\label{disccableB}

Just as with connected cabling, disconnected cabling looks very
simple in the algebra $\B$. Composing $\psi^{r\cdot 1}_{\yy}$ with
$\chi$ we immediately get the following
\begin{xproposition}
The disconnected cabling operation $\psi^{r\cdot 1}_{\yy}$ sends an
open diagram in $\B(\yy)$ with $k$ legs to the sum of all $r^k$ ways
of replacing one label $\yy$ by $r$ labels $\yy_1,\ldots,\yy_r$.
\end{xproposition}
A similar statement holds, of course, for diagrams with more than
one skeleton component.

The twofold disconnected cabling $\psi^{2\cdot 1}=\psi^{2\cdot
1}_{\yy}$ can be thought of as a coproduct: it is a map
$$\B\simeq\B(\yy)\to\B(\yy_1,\yy_2)\simeq\B\ot\B.$$
This coproduct is dual, in a certain sense, to the product in $\B$;
we shall make this statement precise in Section~\ref{pairandcable}.

\section{The Kontsevich integral of a $(p,q)$-cable}
\label{cabelKI}\ \index{Cabling!of the Kontsevich integral} The
Kontsevich integral is well-behaved with respect to taking
$(p,q)$-cables for all values of $p$ and $q$.

\begin{xtheorem}[\cite{LM5, BLT}] %\label{thm:disc-cable}
For a framed tangle $T$ with a closed component $\yy$ and $p,q,r$
integers such that $r$ is the greatest common divisor of $p,q$
$$I^{fr}(T^{(p,q)}_{\yy}) =
\psi^{r\cdot p/r}_{\yy}\Bigl(
I^{fr}(T)\#_{\yy}\exp(\frac{q}{2p}\Theta)\Bigr)\,,$$ where
$\#_{\yy}$ denotes the action of $\F$ on the tangle chord diagrams
by taking the connected sum along the component $\yy$.
\end{xtheorem}
\begin{xremark}
At the first sight the formula of \cite{BLT} for the Kontsevich
integral of a $(p,1)$-cable may seem to disagree with the above
theorem. This is due to a different choice of framing on the
$(p,1)$-cable of a knot in \cite{BLT}.
\end{xremark}
\begin{proof}
Given a framing on $T$, the component $\yy$ can be drawn in such a
manner that the framing curve on $\yy$ is produced by shifting $\yy$
along some small {\em constant} vector $(z,0)$ with no vertical
component. (Here we use the coordinates $(z,t)$ in $\R^3$ with
$z\in\C$ and $t\in\R$.) For instance, if $T$ is represented by a
plane diagram, shifting $\yy$ in the direction perpendicular to the
projection plane we get the blackboard framing. Using this
observation, it will be convenient to modify the definition of the
$(p,q)$-cable along $\yy$ as follows.

For $q=0$ we simply replace the component $\yy$ by the $p$ curves
obtained from it by shifting along the vector $(z\exp{\frac{2\pi k
i}{p}},0)$, for $0\leq k< p$. If these curves have non-empty
intersection, they can be made disjoint by a small rotation of the
vector $(z,0)$. Moreover, we can assume these curves remain disjoint
if the vector $(z,0)$ is replaced by $(\e z, 0)$ for any $0<\e<1$.
The framing on each of these curves is chosen so that the
translation by $(z\exp{\frac{2\pi k i}{p}},0)$ gives the framing
curve. The $(p,q)$-cable with $q\neq 0$ is then obtained from the
$(p,0)$-cable by inserting a ``twist'' by $\exp{\frac{2\pi q
i}{p}}$:
$$\ig[height=32mm]{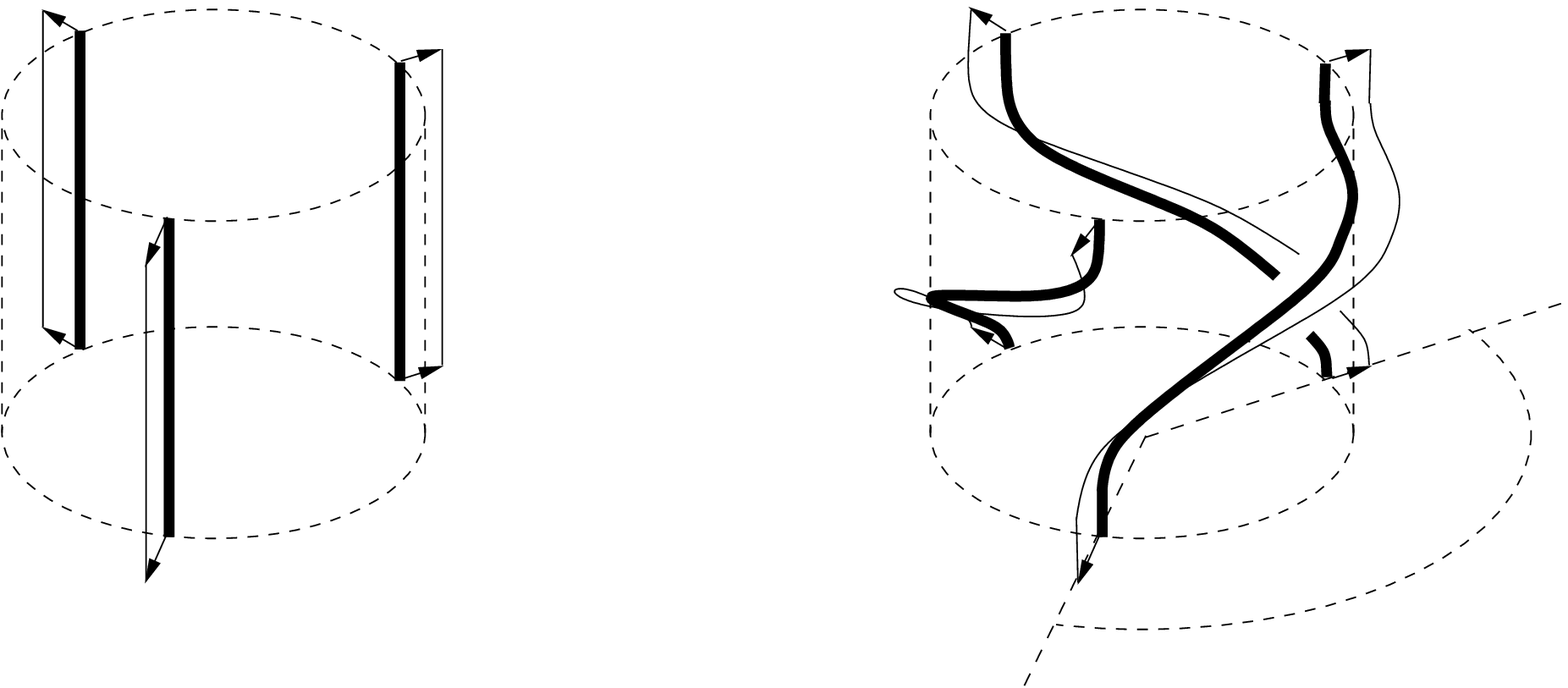}\put(0,30){$\scriptstyle \frac{2\pi q}{p}$}
$$
The main difference of this definition with the one given on
page~\pageref{def:cabling} is that the $p$ branches of the cable are
obtained from $\yy$ by shifting in the horizontal planes rather than
in the planes perpendicular to $\yy$. This difference becomes
significant near the critical points of the height on $\yy$:
$$\ig[height=30mm]{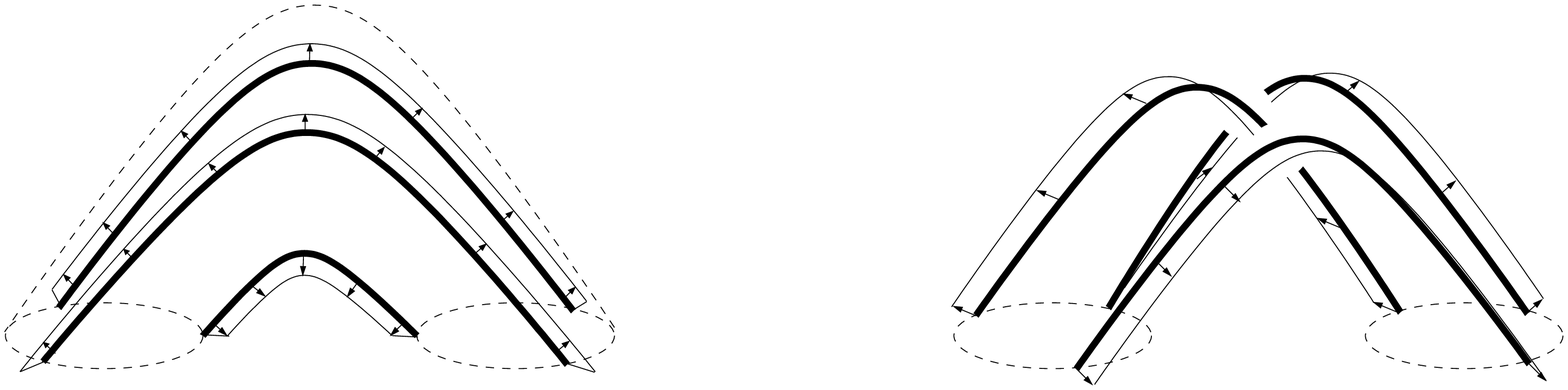}$$
Nevertheless, it is easy to see that this definition of the
$(p,q)$-cable coincides with the old definition up to isotopy.

In view of the Remark~\ref{rem:cablehump} we can replace $I^{fr}$ in
the statement of the theorem by $Z^{fr}$. Let us first consider the
easiest case when $q=0$.

Let us assume that the framing on $\yy$ is given, as above, using a
shift by a constant vector $(z,0)$, and, moreover, that the same
framing is produced by the vectors $(\e z, 0)$ for all $0<\e<1$.
Take an arbitrary diagram $D$ that participates in $Z^{fr}(T)$, and
let $\widetilde{D}$ be a lifting of $D$ from $\F(T)$ to
$\F(T^{(p,0)}_{\yy})$. Then, as $\e$ tends to 0, the coefficient of
$\widetilde{D}$ in the integral $Z^{fr}\bigl(T^{(p,0)}_{\yy}\bigr)$
tends to the coefficient of $D$ in $Z^{fr}(T)$. Since the Kontsevich
integral is an invariant, this gives the statement of the theorem
for $q=0$.

Now, consider the case of coprime $p$ and $q$. The proof of this
case is a twist on the proof of the case $q=0$. Write $\tilde{\yy}$
for the cabling of $\yy$. Let $T'$ be the part of $T^{(p,q)}_{\yy}$
consisting of $1/p$th part of  $\tilde{\yy}$ together with the rest
of $T$:
$$\ig[height=40mm]{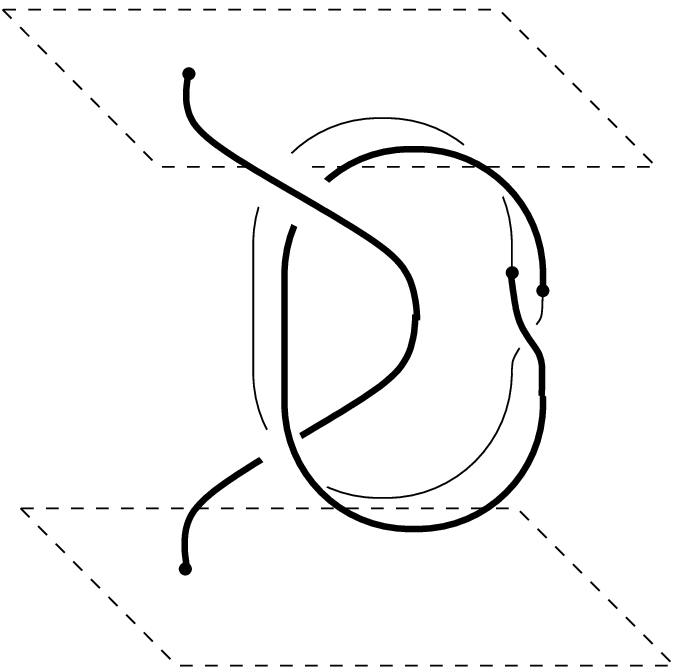}$$
Here $T'$ is indicated by a thicker line. The tangle
$T^{(p,q)}_{\yy}$ projects onto $T'$, this projection
$\pi$ 
sends a point of $T$
to the nearest point of $T'$ in
the same horizontal plane. (Note that there are points where this
projection is not continuous, but this is of no importance).

Define $Z'$ in the same way as the framed Kontsevich integral of
$T^{(p,q)}_{\yy}$ but only integrating over those pairings that have
all their chord endpoints on $T'$. For every such pairing with $s$
endpoints on $T'\cap\tilde{\yy}$ there are $p^s$ lifts, with respect
to $\pi$, to a pairing on $T^{(p,q)}_{\yy}$. As $\e$ tends to $0$,
the coefficient at each of these lifts in $Z^{fr}(T^{(p,q)}_{\yy})$
tends to the coefficient of the original pairing in $Z'$. Moreover,
the sum of all chord diagrams that correspond to all the lifts of a
pairing $P$ is nothing but $\psi^{p}_{\yy}(D_P)$, and, hence,
$$Z^{fr}(T^{(p,q)}_{\yy}) = \psi^{p}_{\yy}(Z').$$
It remains to show that
$$Z'=Z^{fr}(T)\#_{\yy}\exp(\frac{q}{2p}\Theta).$$

Let us say, within this proof, that a chord in a pairing is {\em
short} if it connects a point on a framed tangle to its shift on the
framing curve.
Observe that the chord diagram of a pairing modulo 4T relations does
not depend on the position of the isolated chords, and, in
particular, on the position of short chords. This means that we can
write $Z'=Z'_1 \#_{\yy} Z'_2$ where $Z'_1$ is the part of $Z'$ which
involves only the pairings without short chords, and $Z'_2$ is the
part involving pairings that consist just of short chords.

Note that there is no contribution of short chords in $Z^{fr}(T)$
and, as $\e\to 0$, the limit of $Z'_1$ is precisely $Z^{fr}(T)$. The
only non-trivial contribution in $Z'_2$ comes from the pairings that
have all their chords in the twisted part of $T'$, and this is
readily seen to be $\exp(\frac{q}{2p}\Theta)$. This settles the case
of $p$ and $q$ coprime.

Finally, if $p$ and $q$ are not coprime and $r>1$ is their greatest
common divisor, the argument above applies almost word for word.
Apart from replacing $\psi^{p}_{\yy}$ by $\psi^{r\cdot p/r}_{\yy}$,
the only difference is that now $\tilde{\yy}$ has $r$ components,
and the definition of $Z'$ requires one more step: one has to delete
in all the diagrams participating in $Z'$ the $r-1$ components of
$\tilde\yy$ that do not intersect $T'$.
\end{proof}

\subsection{Torus knots}
The $(p,q)$-torus knot is the $(p,q)$-cable of the unknot, and,
therefore, the formula for the cables of the Kontsevich integral as
a particular case gives an expression for the Kontsevich integral of
torus knots. An essential ingredient of this expression is the
Kontsevich integral $I^{fr}$ of the unknot, which will be treated
later in Chapter~\ref{advKI}.

As we noted in the proof of the cabling formula, it remains true
with $Z^{fr}$ in place of $I^{fr}$. Thus we have an expression for
the preliminary Kontsevich integral of the $(p,q)$-torus knot with
$2p$ critical points of the height function:
$$Z^{fr}(O^{(p,q)})=\psi^{p}\left(\exp\left(\frac{q}{2p}\Theta\right)\right).$$
This formula is not quite explicit, since it involves the operations $\psi^{p}$. We refer the reader to the paper \cite{Mar} by J.~March\'e where the Kontsevich integral of torus knots is computed 
without the help of the operations $\psi^{p}$.

\section{Cablings of the Lie algebra weight systems}\label{cabelWS}
\index{Cabling!of weight systems}

In Chapter~\ref{LAWS_A} we have seen how a semi-simple Lie algebra
$\g$ gives rise to the universal Lie algebra weight system
$\f_\g:\A^{fr}\to U(\g)$, and how a representation $V$ of $\g$
determines a numeric weight system $\f_\g^V:\A^{fr}\to\C$. The
interaction of $\f_\g$ with the operation $\psi^p$ is rather
straightforward.

Define $\mu^p\colon U(\g)^{\ot p}\to U(\g)$ and $\d^p\colon U(\g)\to
U(\g)^{\ot p}$ by
$$\mu^p(x_1\ot x_2\ot\ldots\ot x_p)= x_1 x_2\ldots x_p$$
for $x_i\in U(\g)$ and
$$\d^p(g) = g\ot1\ot\dots\ot 1 + 1\ot g\ot\dots\ot 1 + \dots +
            1\ot1\ot\dots\ot g\,,$$
where $g\in \g$.
\begin{proposition} For $D\in \A^{fr}$ we have
$$(\f_\g\circ\psi^{p})(D) = (\mu^p\circ\d^p)(\f_\g(D))\,.$$
\end{proposition}

\begin{proof}

The construction of the universal Lie algebra weight system
(Section~\ref{univ_lie_ws}) consists in assigning the basis vectors
$e_{i_a}\in\g$ to the endpoints of each chord $a$, then taking their
product along the Wilson loop and summing up over each index $i_a$.
For the weight system $\f_\g\circ\psi^p$, to each endpoint of a
chord we assign not only a basis vector, but also the sheet of the
covering to which that particular point is lifted. (Since the
construction of Lie algebra weight systems uses based diagrams, the
sheets of the covering can actually be enumerated.) To form an
element of the universal enveloping algebra we must read the letters
$e_{i_a}$ along the circle $n$ times. On the first pass we read only
those letters which are related to the first sheet of the covering,
omitting all the others. Then read the circle for the second time
and now collect only the letters from the second sheet, etc. up to
the $p$th reading.  The products of $e_{i_a}$'s thus formed are
summed up over all the $i_a$ and over all the ways of lifting the
endpoints to the covering.

On the other hand, the operation $\mu^p\circ\d^p:U(\g)\to U(\g)$ can
be described as follows. If $A$ is an ordered set of elements of
$\g$, let us write $\prod_A\in U(\g)$ for the product of all the
elements of $A$, according to the order on $A$. Let $x=\prod_A$ for
some $A$. To obtain $\mu^p\circ\d^p(x)$ we take all possible
decompositions of $A$ into an ordered set of $n$ disjoint subsets
$A_i$, with $1\leq i\leq n$, and take the sum of
$\prod_{A_1}\prod_{A_2}\ldots\prod_{A_p}$ over all these
decompositions.

When applied to $\f_\g$, the sets $A_k$ are the sets of $e_{i_a}$
corresponding to the endpoints that are lifted to the $k$th sheet of
the $p$-fold covering. This establishes a bijection between the 
summands on the two sides of the formula.
\end{proof}

\begin{xcb}{Exercises}
\begin{enumerate}

\item

Define the connected sum of two framed knots as their usual
connected sum with the framing whose self-linking number is the sum
of the self-linking numbers of the summands. Show that the framed
Kontsevich integral is multiplicative with respect to the connected
sum.

\item
Prove that the framed Kontsevich integral $Z^{fr}(K)$ is a
group-like element \index{Group-like element} of the Hopf algebra
$\widehat\A^{fr}$.

\item\label{strutsinifr}
Let $K$ be a framed knot. Consider the Kontsevich integral
$I^{fr}(K)$ as an element of $\widehat{\B}$, and show that if at
least one of the diagrams participating in it contains a strut (an
interval component) then $K$ has non-zero framing.

{\sl Hint.} Use the group-like property of $Z^{fr}$.

\item\label{ex_cabl_4t}
Check that the maps $\psi^p$, and, more generally, $\psi^{r\cdot
p}_{\yy}$ are compatible with the four-term relations.

\item
Compute $\psi^3(\cabldp{fd2b})$ and $\psi^3(\cabldp{fd2a})$.

\item
Compute the eigenvalues and eigenvectors of $\psi^3|_{_{\A_2}}$.

\item
Compute $\psi^22(\cabldp{symcd1})$, $\psi^2(\cabldp{symcd2})$,
$\psi^2(\cabldp{symcd3})$, $\psi^2(\cabldp{symcd4})$, and
$\psi^2(\cabldp{symcd5})$.

\item
Compute the eigenvalues and eigenvectors of $\psi^2|_{_{\A_3}}$.

\item\label{ex_cabl_isoch}
\parbox[t]{3in}{Compute $\psi^2(\Theta^m)$, where $\Theta^m$
is a chord diagram with $m$ isolated chords, such as the one shown
on the right. }\qquad $\rb{0pt}{$\Theta^m =$}\quad
\risS{-12}{cabl_isoch}{
         \put(-5,-10){\mbox{\scriptsize $m$ chords}}}{30}{20}{15}$

\item\label{ex_cabl_comul}
Prove that $\psi^p$ commutes with the comultiplication of chord
diagrams. In other words, show that in the notation of Section~\ref{comult_cd}, page~\pageref{comult_cd}, the identity
$$
  \delta(\psi^p(D)) = \sum_{J\subseteq [D]}
          \psi^p(D_J)\ot \psi^p(D_{\overline J})
$$
holds for any chord diagram $D$.

\item
({\em D.~Bar-Natan \cite{BN1}}). Prove that $\psi^p\circ \psi^q =
\psi^{pq}$.

\item Prove the Proposition from Section~\ref{vassiliev_cables}:
$$\symb((v)^{(p,q)})(D) = \symb(v)(\psi^p (D))\,.$$

\item\label{op_Delta}
Let $T$ be a tangle with $k$ numbered components, all of them
intervals, and assume that the $i$th component connects the $i$th
point on the upper boundary with the $i$th point on the lower
boundary. (String links and pure braids are examples of such
tangles.) Let $\Delta^{n_1,\ldots,n_k}_{\e}$ be the operation of
replacing, for each $i$, the $i$th component of $T$ by $n_i$
parallel copies of itself with the distance $\e$ between each copy,
as on Figure~\ref{fig:Delta_2_3}.
\begin{figure}[ht]
\ig[height=3cm]{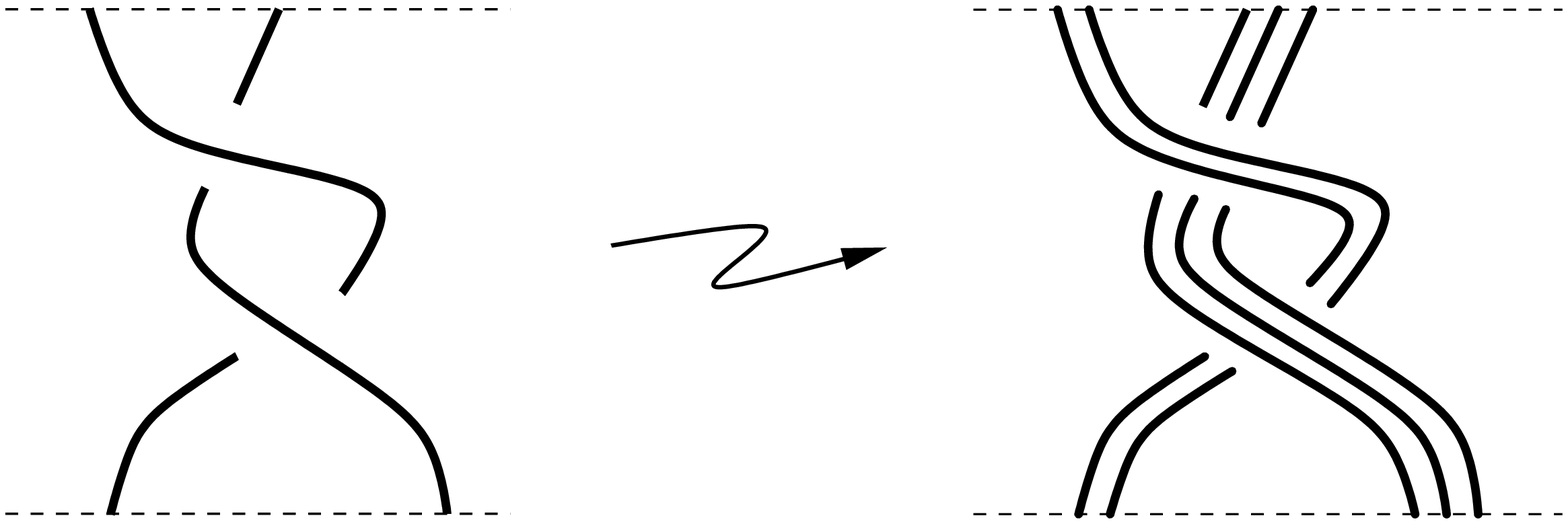}\caption{The effect of
$\Delta^{2,3}_{\e}$}\label{fig:Delta_2_3}
\end{figure}
Denote by $\Delta^{n_1,\ldots,n_k}$ \label{Delta_n1nk} the following
operation on the corresponding tangle chord diagrams: for each $i$
the $i$th strand is replaced by $n_i$ copies of itself and a chord
diagram is sent to the sum of all of its liftings to the resulting
skeleton. Prove that
$$\lim_{\e\to 0}Z(\Delta^{n_1,\ldots,n_k}_{\e}(T))=\Delta^{n_1,\ldots,n_k}Z(T).$$

\item\label{op_Delta_2}
Let $T_{\e}$ be the following family of tangles depending on a
parameter $\e$:
$$T_{\e}=\rb{-14mm}{\ig[height=30mm]{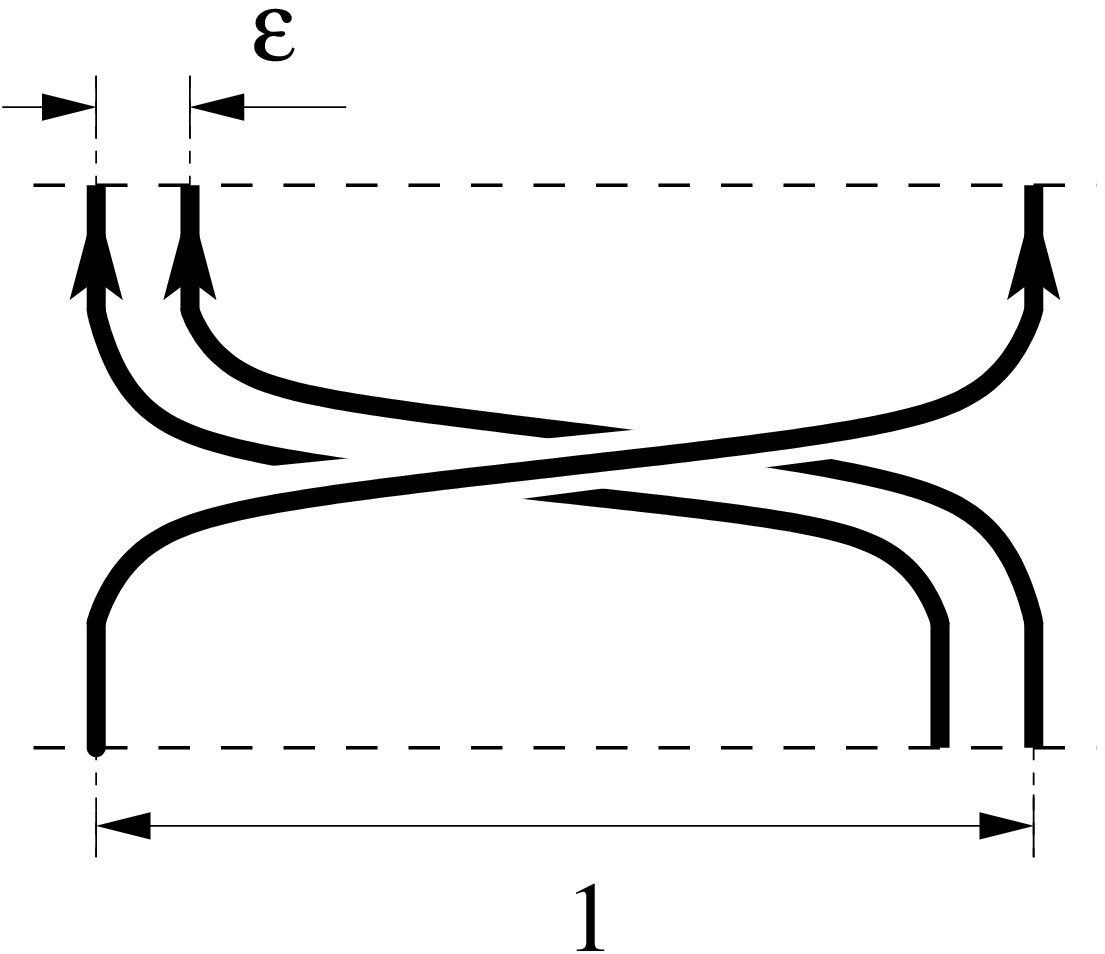}}$$
Show that $$\lim_{\e\to
0}Z(T_{\e})=\rb{-3mm}{\ig[height=9mm]{dzcvirt.eps}}
\cdot\Delta^{1,2}\exp\Bigl(\frac{\risS{-3}{doc}{}{17}{0}{0}}{2}\Bigr).$$

\end{enumerate}
\end{xcb}
 %9 Operations and KI
\chapter{The Drinfeld associator} % 10
\label{DA_chap}
\def\Alg{\widehat{\A}}

\def\smtan#1{\rb{-5pt}{\ig[height=6mm]{#1.eps}}}

In this chapter we give the details of the combinatorial
construction for the Kontsevich integral. The main ingredient of
this construction is the power series known as the {\em Drinfeld
associator} $\PhiKZ$. Here the subscript ``KZ'' indicates that the
associator comes from the solutions to the Knizhnik-Zamolodchikov
equation. The Drinfeld associator enters the theory as a
(normalized) Kontsevich integral for a special tangle without
crossings, which is the simplest {\em associating tangle}.

The associator $\PhiKZ$ is an infinite series in two non-commuting
variables whose coefficients are combinations of multiple zeta
values. In the construction of the Kontsevich integral only some
properties of $\PhiKZ$ are used; adopting them as axioms, we arrive
at the general notion of an associator that appeared in Drinfeld's
papers \cite{Dr1,Dr2} in his study of quasi-Hopf algebras. These
axioms actually describe a large collection of associators belonging
to the completed algebra of chord diagrams on three strands. Some of
these associators have rational coefficients, and this implies the
rationality of the Kontsevich integral.

\section[The Knizhnik--Zamolodchikov equation]{The KZ
equation and iterated integrals}
\label{KZ}

In this section, we give the original Drinfeld's definition of the
associator in terms of the solutions of the simplest
Knizhnik--Zamolodchikov equation.

The Knizhnik--Zamolodchikov (KZ) equation appears in the
Wess--Zumino--Witten model of conformal field theory \cite{KnZa}.
The theory of KZ type equations has been developed in the contexts
of mathematical physics, representation theory and topology
\cite{EFK, Var, Kas, Koh4, Oht1}. Our exposition follows the
topological approach and is close to that of the last three books.

\subsection{General theory}
\label{KZgen}

Let $X$ be a smooth manifold and $\Alg$ a completed graded algebra
over the complex numbers. Choose a set $\omega_1,\ldots,\omega_p$ of
$\C$-valued closed differential 1-forms on $X$ and a set
$c_1,\ldots,c_p$ of homogeneous elements of $\Alg$ of degree 1.
Consider the closed 1-form
$$\omega=\sum_{j=1}^p \omega_j c_j$$
with values in $\Alg$. The Knizhnik-Zamolodchikov equation is a
particular case of the following very general equation:
\begin{equation}\label{eq:Chen_connection}
dI= \omega\cdot I,
\end{equation}
where $I:X\to\A$ is the unknown function.

\noindent{\bf Exercise.} One may be tempted to solve the above
equation as follows: $d\log(I)=\omega$, therefore
$I=\exp\int\omega$. Explain why this is wrong.

The form $\omega$ must satisfy certain conditions so that
Equation~\ref{eq:Chen_connection} may have non-zero solutions.
Indeed, taking the differential of both sides of
(\ref{eq:Chen_connection}), we get that $0=d(\omega I)$. Applying
the Leibniz rule, using the fact that $d\omega=0$ and substituting
$dI=\omega I$, we see that a necessary condition for integrability
can be written as
\begin{equation} \label{intKZ}
\omega\wedge\omega=0
\end{equation}

It turns out that this condition is not only necessary, but also
sufficient for {\em local} integrability: if it holds near a point
$x_0\in X$, then (\ref{eq:Chen_connection}) has the unique solution
$I_0$ in a small neighbourhood of $x_0$, satisfying the initial
condition $I_0(x_0)=a_0$ for any  $a_0\in\Alg$. This fact is
standard in differential geometry where it is called the {\em
integrability of flat connections} (see, for instance, \cite{KN}).
A direct {\em ad hoc} proof can be found in \cite{Oht1},
Proposition~5.2.

\subsection{Monodromy}

Assume that the integrability condition \ref{intKZ} is satisfied for
all points of $X$. Given a (local) solution $I$ of
Equation~\ref{eq:Chen_connection} and $a\in\Alg$, the product $Ia$
is also a (local) solution. Therefore, germs of local solutions at a
point $x_0$ form an $\Alg$-module. This module is free of rank one;
it is generated by the germ of a local solution taking value
$1\in\Alg$ at $x_0$.

The reason to consider germs rather that global solutions is that
the global solutions of (\ref{eq:Chen_connection}) are generally
multivalued, unless $X$ is simply-connected. Indeed, one can extend
a local solution at $x_0$ along any given path which starts at $x_0$
by patching together local solutions at the points on the path. (One
can think of this extension as something like an analytic
continuation of a holomorphic function.) Extending in this way a
local solution $I_0$ at the point $x_0$ along a closed loop
$\gamma:[0,1]\to X$ we arrive to another local solution $I_1$, also
defined in a neighbourhood of $x_0$.

Let $I_1(x_0)=a_\gamma$. Suppose that $a_0$ is an invertible element
of $\Alg$. The fact that the local solutions form a free
one-dimensional $\Alg$-module implies that the two solutions $I_0$
and $I_{1}$ are proportional to each other: $I_1=I_0 a_0^{-1}
a_\gamma$. The coefficient $a_0^{-1} a_\gamma$ does not depend on a
particular choice of the invertible element $a_0\in\Alg$ and the
loop $\gamma$ within a fixed homotopy class. Therefore, we get a
homomorphism $\pi_1(X)\to\Alg^*$ from the fundamental group of $X$
into the multiplicative group of the units of $\Alg$, called the
{\em monodromy representation}.

\subsection{Iterated integrals}
\label{iter_int}

Both the continuation of the solutions and the monodromy
representation can be expressed in terms of the 1-form $\omega$.
Choose a path $\gamma:[0,1]\to X$, not necessarily closed, and
consider the composition $I\circ\gamma$. This is a function
$[0,1]\to\Alg$ which we denote by the same letter $I$; it satisfies
the ordinary differential equation
\begin{equation}\label{eq_on_I}
  \frac{d}{dt}I(t)= \omega(\dot{\gamma}(t))\cdot I(t),
  \qquad\qquad I(0)=1\ .
\end{equation}
The function $I$ takes values in the completed graded algebra
$\Alg$, and it can be expanded in an infinite series according to
the grading:
$$I(t)=I_0(t)+I_1(t)+I_2(t)+\dots\ ,$$
where each term $I_m(t)$ is the homogeneous degree $m$ part of $I(t)$.

The form $\omega$ is homogeneous of degree 1 (recall that
$\omega=\sum c_j\omega_j$, where $\omega_j$ are $\C$-valued 1-forms
and $c_j$'s are elements of $\A_1$). Therefore
Equation~\ref{eq_on_I} is equivalent to an infinite system of
ordinary differential equations
$$\begin{array}{lll}
  I_0'(t)=0,                &\quad& I_0(0)=1, \\
  I_1'(t)=\omega(t) I_0(t), &\quad& I_1(0)=0, \\
  I_2'(t)=\omega(t) I_1(t), &\quad& I_2(0)=0, \\
  \dots\dots&&\dots\dots
\end{array}$$
where $\omega(t)=\gamma^*\omega$ is the pull-back of the 1-form to the
interval $[0,1]$.

These equations can be solved iteratively, one by one. The first one
gives $I_0=\const$, and the initial condition implies
$I_0(t)=1$. Then,\\
$\displaystyle I_1(t)=\int_0^t \omega(t_1)$.
Here $t_1$ is an auxiliary variable that ranges from 0 to $t$.
Coming to the next equation, we now get:
$$
I_2(t)=\int_0^t \omega(t_2)\cdot I_1(t_2) =
    \int_0^t \omega(t_2) \left(\int_0^{t_2} \omega(t_1)\right)
   = \int\limits_{0<t_1<t_2<t}
           \omega(t_2)\wedge\omega(t_1),
$$
Proceeding in the same way, for an arbitrary $m$ we obtain
$$
I_m(t)=\int\limits_{0<t_1<t_2<\dots<t_m<t}
       \omega(t_m)\wedge\omega(t_{m-1})\wedge\dots\wedge\omega(t_1)
$$
In what follows, it will be more convenient to use this formula
with variables renumbered:
\begin{equation} \label{iterint}
I(t)=1+\sum\limits_{m=1}^\infty \
               \int\limits_{0<t_m<t_{m-1}<\dots<t_1<t}
               \omega(t_1)\wedge\omega(t_2)\wedge\dots\wedge\omega(t_m)
\end{equation}
The value $I(1)$ represents the monodromy of the solution over the
loop $\gamma$. Each iterated integral $I_m(1)$ is a homotopy
invariant (of ``order $m$'') of $\gamma$. Note the resemblance of
these expressions to the Kontsevich integral; we shall come back to
that again later.

\begin{xremark}
One may think of the closed 1-form $\omega$ as of an $\Alg$-valued
connection on $X$.  Then the condition $\omega\wedge\omega=0$ means
that this connection is flat. The monodromy $I(t)$ represents the
parallel transport. In this setting the presentation of the parallel
transport as a series of iterated integrals was described by
K.-T.~Chen \cite{Chen2}.
\end{xremark}

\subsection{The Knizhnik-Zamolodchikov equation}

\label{kzcd} Let ${\mathcal H}=\bigcup_{j=1}^p H_j$ be a collection
of affine hyperplanes in $\C^n$. Each hyperplane $H_j$ is defined by
a (not necessarily homogeneous) linear equation $L_j=0$. A {\em
Knizhnik--Zamolodchikov}, or simply {\em KZ}, equation
\index{Knizhnik--Zamolodchikov!equation} \index{KZ equation} is an
equation of the form (\ref{eq:Chen_connection}) with
$$X=\C^n-{\mathcal H}$$
and
$$\omega_j=d\log{L_j}$$
for all $j$.

Many of the KZ equations are related to Lie algebras and their
representations. This class of equations has attracted the most
attention in the literature; see, for example,
\cite{KnZa,Oht1,Koh4}. We are specifically interested in the
following situation.

Suppose that $X=\C^n\setminus\mathcal{H}$ where $\mathcal{H}$ is the
union of the diagonal hyperplanes $\{z_j=z_k\}$, $1\le j<k\le n$,
and the algebra $\Alg=\Ab^h(n)$ is the completed algebra of
horizontal chord diagrams, see page~\pageref{A^h(n)}.
\index{Algebra!horizontal chord diagrams} Recall that $\Ab^h(n)$ is
spanned by the diagrams on $n$ vertical strands (which we assume to
be oriented upwards) whose all chords are horizontal.
Multiplicatively, $\Ab^h(n)$ is generated by the degree-one elements
$u_{jk}$ for all $1\leq j<k\leq n$, (which are simply the horizontal
chords joining the $j$th and the $k$th strands) subject to the
infinitesimal pure braid relations
\begin{align*}
   &[u_{jk},u_{jl}+u_{kl}] =0,\quad \text{ if } j,k,l \text{ are different},\\
   &[u_{jk},u_{lm}] =0, \quad \text{ if } j,k,l,m \text{ are different},
\end{align*}
where, by definition, $u_{jk}=u_{kj}$.

Consider the $\A^h(n)$-valued 1-form $\displaystyle
\omega=\frac{1}{2\pi i}
   \sum\limits_{1\leqslant j<k\leqslant n}
   u_{jk}\frac{dz_j-dz_k}{z_j-z_k}$
and the corresponding KZ equation
\begin{equation}\label{KZeqCD}
dI=\frac{1}{2\pi i}\Bigl(\sum\limits_{1\leqslant j<k\leqslant n}
   u_{jk}\frac{dz_j-dz_k}{z_j-z_k}\Bigr)\cdot I\ .
\end{equation}
This case of the Knizhnik-Zamolodchikov equation is referred to as
the {\em formal KZ equation}\index{KZ equation!formal}.

The integrability condition \ref{intKZ} for the formal KZ equation
is the following identity on the 1-form $\omega$ on $X$ with values
in the algebra $\A^h(n)$:
$$\omega\wedge\omega
=\sum\limits_{\substack{1\leqslant j<k\leqslant n\\
                        1\leqslant l<m\leqslant n}}
 u_{jk}u_{lm}\frac{dz_j-dz_k}{z_j-z_k}\wedge\frac{dz_l-dz_m}{z_l-z_m}
\quad=\quad 0\ .$$ This identity, in a slightly different notation,
was actually proved in Section~\ref{hor_inv} when we checked the
horizontal invariance of the Kontsevich integral.

The space $X=\C^n\setminus{\mathcal H}$ is the configuration space
of $n$ distinct (and distinguishable) points in $\C$. A loop
$\gamma$ in this space may be identified with a pure
braid\index{Braid!pure} (that is, a braid that does not permute the
endpoints of the strands), and the iterated integral formula
\ref{iterint} gives
$$
  I(1)=\sum_{m=0}^\infty \frac{1}{(2\pi i)^m}
         \int\limits_{0<t_m<\dots<t_1<1}\
         \sum_{P=\{(z_j,z'_j)\}}\  D_P
         \bigwedge_{j=1}^m \frac{dz_j-dz'_j}{z_j-z'_j} \ ,
$$
where $P$ (a {\em pairing}) is a choice of $m$ pairs of points on
the braid, with $j$th pair lying on the level $t=t_j$, and $D_P$ is
the product of $m$ $T$-chord diagrams of type $u_{jj'}$
corresponding to the pairing $P$. We can see that the monodromy of
the KZ equation over $\gamma$ coincides with the Kontsevich integral
of the corresponding braid (see Section \ref{ki_tangles}).

\subsection{The case $n=2$}
For small values of $n$ Equation \ref{KZeqCD} is easier to handle.
In the case $n=2$ the algebra $\Ab^h(2)$ is free commutative on one
generator and everything is very simple, as the following exercise
shows.

\noindent{\bf Exercise.} Solve explicitly Equation \ref{KZeqCD} and
find the monodromy representation in the case $n=2$.

\subsection{The case $n=3$}
\label{kzdrass}
The formal KZ equation for $n=3$ has the form
$$dI=\frac{1}{2\pi i} \Bigl(u_{12}d\log(z_2-z_1) +
     u_{13}d\log(z_3-z_1) + u_{23}d\log(z_3-z_2)\Bigr)\cdot I\ ,
$$
which is a partial differential equation in 3 variables. It turns
out that it can be reduced to an ordinary differential equation.

Indeed, make the substitution
$$I=(z_3-z_1)^{\frac{u}{2\pi i}}\cdot
   { G}\ ,$$
where $u:=u_{12}+u_{13}+u_{23}$ and we understand the factor in
front of $G$ as a power series in the algebra $\Ab^h(3)$:
\begin{multline*}
(z_3-z_1)^{\frac{u}{2\pi i}}\ =\
 \exp\left(\frac{\log(z_3-z_1)}{2\pi i}u\right) \\
 = 1+ \frac{\log(z_3-z_1)}{2\pi i}u
+\frac{1}{2!}\frac{\log^2(z_3-z_1)}{(2\pi i)^2}u^2
+ \frac{1}{3!}\frac{\log^3(z_3-z_1)}{(2\pi i)^3}u^3 + \dots
\end{multline*}

By Proposition \ref{prop:a3} on page \pageref{prop:a3}, the algebra
$\Ab^h(3)$ is a direct product of the free algebra on $u_{12}$ and
$u_{23}$, and the free commutative algebra generated by $u$. In
particular, $u$ commutes with all elements of $\Ab^h(3)$. Taking
this into the account we see that the differential equation for $G$
can be simplified so as to become
$$d{ G}=\frac{{1}}{2\pi i}\left(
    u_{12}d\log\Bigl(\frac{z_2-z_1}{z_3-z_1}\Bigr) +
    u_{23}d\log\Bigl(1 - \frac{z_2-z_1}{z_3-z_1} \Bigr) \right) G\ .$$

Denoting $\frac{z_2-z_1}{z_3-z_1}$ simply by $z$, we see that the
function $G$ depends only on $z$ and  satisfies the following
ordinary differential equation (the {\em reduced KZ equation})
\begin{equation}\label{redKZ}\index{KZ equation!reduced}
\frac{dG}{dz} =\left(\frac{A}{z} + \frac{B}{z-1}\right) G
\end{equation}
where $A:=\frac{u_{12}}{2\pi i}$,  $B:=\frac{u_{23}}{2\pi i}$. As
defined, $G$ takes values in the algebra $\Ab^h(3)$ with three
generators $A$, $B$, $u$. However, the space of local solutions of
this equation is a free module over $\Ab^h(3)$ of rank 1, so the
knowledge of just one solution is enough. Since the coefficients of
Equation \ref{redKZ} do not involve $u$, the equation does have a
solution with values in the ring of formal power series $\CAB$ in
two non-commuting variables $A$ and $B$.\index{$\CAB$}

\subsection{The reduced KZ equation}\label{defDrass}
The reduced KZ equation~\ref{redKZ} is a particular case of the
general KZ equation defined by the data $n=1$,
$X=\C\setminus\{0,1\}$, $L_1=z$, $L_2=z-1$, $\A=\CAB$, $c_1=A$,
$c_2=B$.

Although (\ref{redKZ}) is a first order ordinary differential
equation, it is hardly easier to solve than the general KZ equation.
In the following exercises we invite the reader to try out two
natural approaches to the reduced KZ equation.

\begin{xxca}
Try to find the general solution of Equation~\ref{redKZ} by
representing it as a series
$$G=G_0+G_1A+G_2B+G_{11}A^2+G_{12}AB+G_{21}BA+\ldots,$$ where the $G$'s
with subscripts are complex-valued functions of $z$.
\end{xxca}

\begin{xxca}
Try to find the general solution of Equation~\ref{redKZ} in the form
of a Taylor series $G=\sum_k G_k(z-\frac{1}{2})^k$, where the
$G_k$'s are elements of the algebra $\CAB$. (Note that it is not
possible to expand the solutions at $z=0$ or $z=1$, because they
have essential singularities at these points.)
\end{xxca}

These exercises show that direct approaches do not give much insight
into the nature of the solutions of (\ref{redKZ}). Luckily, we know
that at least one solution exists (see Section~\ref{KZgen}) and that
any solution can be obtained from one basic solution via
multiplication by an element of the algebra $\CAB$. The Drinfeld
associator appears as a coefficient between two remarkable
solutions.

\begin{xdefinition}
The (Knizhnik-Zamolodchikov) {\em Drinfeld associator}
\index{Drinfeld associator}
\index{Knizhnik-Zamolodchikov!associator} $\PhiKZ$ \index{$\PhiKZ$}
is the ratio $\PhiKZ= G_1^{-1}(z)\cdot G_0(z)$ of two special
solutions $G_0(z)$ and $G_1(z)$ of this equation described in the
following Lemma.
\end{xdefinition}

\begin{lemma}[\cite{Dr1,Dr2}]\label{dr_lemma}
There exist unique solutions $G_0(z)$ and $G_1(z)$ of
Equation~\ref{redKZ}, analytic in the domain $\{z\in\C \mid |z|<1,
|z-1|<1\}$ and with the following asymptotic behaviour:
$$
  G_0(z)\sim z^{\bA} \mbox{\ as\ } z\to0\qquad\mbox{and}\qquad
  G_1(z)\sim (1-z)^{\bB} \mbox{\ as\ } z\to1\ ,
$$
which means that
$$G_0(z)=f(z)\cdot z^{\bA} \qquad\mbox{and}\qquad
  G_1(z)=g(1-z)\cdot (1-z)^{\bB} \ ,
$$
where $f(z)$ and $g(z)$ are analytic functions in a neighbourhood of
$0\in\C$ with values in $\CAB$ such that $f(0)=g(0)=1$, and the
(multivalued) exponential functions are understood as formal power
series, that is, $z^A=\exp(A\log z)=\sum_{k\ge0}(A\log z)^k/k!$
\end{lemma}

\begin{xremark}
It is sometimes said that the element $\PhiKZ$ represents the
monodromy of the KZ equation over the horizontal interval from 0 to
1. This phrase has the following meaning. In general, the monodromy
along a path $\gamma$ connecting two points $p$ and $q$, is the
value at $q$ of the solution, analytical over $\gamma$ and taking
value 1 at $p$. If $f_p$ and $f_q$ are two solutions analytical over
$\gamma$ with initial values $f_p(p)=f_q(q)=1$, then the monodromy
is the element $f_q^{-1}f_p$. The reduced KZ equation has no
analytic solutions at the points $p=0$ and $q=1$, and the usual
definition of the monodromy cannot be applied directly in this case.
What we do is we choose some natural solutions with reasonable
asymptotics at these points and define the monodromy as their ratio
in the appropriate order.
\end{xremark}

\begin{proof}
Plugging the expression $G_0(z)=f(z)\cdot z^{\bA}$ into
Equation~\ref{redKZ} we get
$$f'(z)\cdot z^{\bA} + f\cdot \frac{\bA}{z}\cdot z^{\bA} =
  \left(\frac{\bA}{z} + \frac{\bB}{z-1}\right)\cdot f\cdot z^{\bA}\ ,
$$
hence $f(z)$ satisfies the differential equation
$$f' -\frac{1}{z}[\bA,f]=\frac{-\bB}{1-z}\cdot f \ .
$$
Let us look for a formal power series solution
$f=1+\sum_{k=1}^\infty f_kz^k$ with coefficients
 $f_k\in \CAB$.
  We have the following
recurrence equation for the coefficient of $z^{k-1}$:
$$
kf_k-[\bA,f_k]=(k-\ad_{\bA})(f_k)=-\bB(1+f_1+f_2+\dots+f_{k-1})\ ,
$$
where $\ad_{\bA}$ denotes the operator $x\mapsto [A,x]$.
The operator $k-\ad_{\bA}$ is invertible:
$$(k-\ad_{\bA})^{-1} =
 \sum_{s=0}^\infty \frac{\ad_{\bA}^s}{k^{s+1}}
$$
(the sum is well-defined because the operator $\ad_{\bA}$ increases the
grading), so the recurrence can be solved
$$f_k =  \sum_{s=0}^\infty \frac{\mbox{ad}_{\bA}^s}{k^{s+1}}
       \bigl(-\bB(1+f_1+f_2+\dots+f_{k-1})\bigr)\ .$$
Therefore the desired solution does exist among formal power series.
Since the point 0 is a regular singular point of
Equation~\ref{redKZ}, it follows from the general theory (see \cite{Wal})
that this power series converges for $|z|<1$. We thus get an analytic solution
$f(z)$.

To prove the existence of the second solution, $G_1(z)$, it is best
to make the change of undependent variable $z\mapsto 1-z$ which
transforms Equation~\ref{redKZ} into a similar equation with $A$ and
$B$ swapped.
\end{proof}

\begin{xremark}
If the variables $A$ and $B$ were commutative, then the function
explicitly given as the product $z^A(1-z)^B$ would be a solution of
Equation \ref{redKZ} satisfying both asymptotic conditions of Lemma
\ref{dr_lemma} at once, so that the analogs of $G_0$ and $G_1$ would
coincide. Therefore, the image of $\PhiKZ$ under the
abelianization map $\CAB\to\C[[A,B]]$ is equal to 1.
\end{xremark}

The next lemma gives another expression for the associator in terms
of the solutions of Equation~\ref{redKZ}.
\begin{lemma}[\cite{LM2}] \label{lim_lemma}
Suppose that $\e\in\R$, $0<\e<1$. Let $G_\e(z)$ be the
unique solution of Equation~\ref{redKZ} satisfying the initial
condition $G_\e(\e)=1$. Then
$$\PhiKZ = \lim\limits_{\e\to0}\
    \e^{-\bB}\cdot G_\e(1-\e)\cdot \e^{\bA}\ .$$
\end{lemma}

\begin{proof}
We rely on, and use the notation of, Lemma \ref{dr_lemma}.
The solution $G_\e$ is proportional to the distinguished solution $G_0$:
$$
G_\e(z)=G_0(z)G_0(\e)^{-1}=G_0(z)\cdot \e^{-\bA}f(\e)^{-1}
=G_1(z)\cdot \PhiKZ\cdot \e^{-\bA}f(\e)^{-1}\
$$
(the function $f$, as well as $g$ mentioned below, was defined in Lemma
\ref{dr_lemma}).
In particular,
$$G_\e(1-\e) =
G_1(1-\e)\cdot \PhiKZ\cdot \e^{-\bA}f(\e)^{-1} = g(\e)\e^{\bB}\cdot
\PhiKZ\cdot
   \e^{-\bA}f(\e)^{-1}\ .$$
We must compute the limit
$$\lim\limits_{\e\to0}\
  \e^{-\bB}g(\e)\e^{\bB}\cdot \PhiKZ\cdot
   \e^{-\bA}f(\e)^{-1}\e^{\bA}\ ,$$
which obviously equals
$\PhiKZ$
because $f(0)=g(0)=1$ and $f(z)$ and $g(z)$ are analytic functions in
a neighbourhood of zero. The lemma is proved.
\end{proof}

\subsection{The Drinfeld associator and the Kontsevich integral}
\label{kieqkz}
Consider the reduced KZ equation \ref{redKZ} on the real interval
$[0,1]$ and apply the techniques of iterated integrals from
Section \ref{iter_int}. Let the path $\gamma$ be the identity inclusion
$[\e,1]\to\C$. Then the solution $G_\e$ can be written as
$$G_\e(t) = 1 + \sum_{m=1}^\infty \ \
          \int\limits_{\e<t_m<\dots<t_2<t_1<t}
          \omega(t_1)\wedge\omega(t_2)\wedge\dots\wedge\omega(t_m)\ .
$$
The lower limit in the integrals is $\e$ because the parameter on the
path $\gamma$ starts from this value.

We are interested in the value of this solution at $t=1-\e$:
$$G_\e(1-\e) = 1 + \sum_{m=1}^\infty \ \
          \int\limits_{\e<t_m<\dots<t_2<t_1<1-\e}
          \omega(t_1)\wedge\omega(t_2)\wedge\dots\omega(t_m)\ .
$$
We claim that this series literally coincides with the Kontsevich
integral of the following tangle
$$ AT_{\e}\ =\quad \risS{-40}{assatel}{
         \put(7,1){\mbox{$\scriptstyle 0$}}
         \put(7,21){\mbox{$\scriptstyle \e$}}
         \put(22,2){\mbox{$\scriptstyle \e$}}
         \put(45,1){\mbox{$\scriptstyle 1-\e$}}
         \put(65,1){\mbox{$\scriptstyle 1$}}
         \put(-3,52){\mbox{$\scriptstyle 1-\e$}}
         \put(8,67){\mbox{$\scriptstyle 1$}}
         \put(18,80){\mbox{$t$}}
         \put(93,5){\mbox{$z$}}
                 }{90}{45}{50}
$$
under the identification $A=\frac{1}{2\pi i}\vstod$,
$B=\frac{1}{2\pi i}\vstdt$. Indeed, on every level $t_j$ the
differential form $\omega(t_j)$ consists of two summands. The first
summand $\bA\frac{dt_j}{t_j}$ corresponds to the choice of a pair
$P=(0,t_j)$ on the first and the second strings and is related to
the chord diagram $A=\vstod$. The second summand
$\bB\frac{d(1-t_j)}{1-t_j}$ corresponds to the choice of a pair
$P=(t_j,1)$ on the second and third strings and is related to the
chord diagram $B=\vstdt$. The pairing of the first and the third
strings does not contribute to the Kontsevich integral, because
these strings are parallel and the correspoding differential
vanishes. We have thus proved the following proposition.

\begin{xproposition} \label{reg_KI}
The value of the solution $G_\e$ at $1-\e$ is equal to the
Kontsevich integral $G_\e(1-\e)=Z(AT_{\e})$. Consequently,
by Lemma \ref{lim_lemma}, the KZ
associator coincides with the regularization of the Kontsevich
integral of the tangle $AT_\e$:
$$
  \PhiKZ = \lim\limits_{\e\to0}\
    \e^{-\bB}\cdot Z(AT_{\e})\cdot \e^{\bA},
$$
where $A=\frac{1}{2\pi i}\vstod$ and $B=\frac{1}{2\pi i}\vstdt$.
\end{xproposition}

\section{Calculation of the KZ Drinfeld associator}
\label{zetanumb}

In this section, following \cite{LM4}, we deduce an explicit formula
for the Drinfeld associator $\PhiKZ$. It turns out that all the
coefficients in the expansion of $\PhiKZ$ as a power series in $A$
and $B$ are values of multiple zeta functions \index{MZV} (see
Section \ref{secetazeta}) divided by powers of $2\pi i$.

\subsection{}

Put $\o_0(z)=\frac{dz}{z}$ and $\o_1(z)=\frac{d(1-z)}{1-z}$. Then
the 1-form $\omega$ studied in \ref{kieqkz} is the linear
combination $\o(z)=\bA\o_0(z)+\bB\o_1(z)$, where
$\bA=\frac{\vstod}{2\pi i}$ and $\bB=\frac{\vstdt}{2\pi i}$. By
definition the terms of the Kontsevich integral $Z(AT_{\e})$
represent the monomials corresponding to all choices of one of the
two summands of $\o(t_j)$ for every level $t_j$. The coefficients of
these monomials are integrals over the simplex
$\e<t_m<\dots<t_2<t_1<1-\e$ of all possible products of the forms
$\o_0$ and $\o_1$. The coefficient of the monomial
$\bB^{i_1}\bA^{j_1}\dots\bB^{i_l}\bA^{j_l}$ ($i_1\ge0$, $j_1>0$,
\dots, $i_l>0$, $j_l\ge0$)  is
$$\begin{array}{l}\displaystyle
\int\limits_{\e<t_m<\dots<t_2<t_1<1-\e}
\underbrace{\o_1(t_1)\wedge\dots\wedge\o_1(t_{i_1})}_{i_1}
\wedge
\underbrace{\o_0(t_{i_1+1})\wedge\dots\wedge\o_0(t_{i_1+j_1})}_{j_1}
\wedge\dots \vspace{10pt}\\
\hspace{164pt}\wedge
\underbrace{\o_0(t_{i_1+\dots+i_l+1})\wedge\dots\wedge
         \o_0(t_{i_1+\dots+j_l}) }_{j_l}\ ,
\end{array}
$$
where $m=i_1+j_1+\dots+i_l+j_l$.
For example, the coefficient of $\bA\bB^2\bA$ equals
$$
\displaystyle
\int\limits_{\e<t_4<t_3<t_2<t_1<1-\e}
\o_0(t_1)\wedge\o_1(t_2)\wedge\o_1(t_3)\wedge\o_0(t_4)\,.
$$

We are going to divide the sum of all monomials into two parts,
``convergent'' $Z^{conv}$ and ``divergent'' $Z^{div}$, depending on
the behaviour of the coefficients as $\e\to0$. We shall have
$Z(AT_{\e})=Z^{conv}+Z^{div}$ and
\begin{equation}\label{condiv}
\Phi = \lim\limits_{\e\to0}\
    \e^{-\bB}\cdot Z^{conv} \cdot \e^{-\bA} +
  \lim\limits_{\e\to0}\
    \e^{-\bB}\cdot Z^{div} \cdot \e^{-\bA}\ .
\end{equation}
Then we shall prove that the second limit equals zero and find an
explicit expression for the first one in terms of multiple zeta
values. We shall see that although the sum $Z^{conv}$ does not
contain any divergent monomials, the first limit in (\ref{condiv})
does.

We pass to exact definitions.

\begin{definition}
A non-unit monomial in letters $\bA$ and $\bB$ with positive powers
is said to be {\it convergent}
\index{Convergent monomial} if it
starts with an $\bA$ and ends with a $\bB$. Otherwise the monomial is
said to be {\it divergent}
\index{Divergent monomial}. We regard the unit monomial $1$ as convergent.
\end{definition}

\begin{example}
The integral
$$\int\limits_{a<t_p<\dots<t_2<t_1<b}
  \o_1(t_1)\wedge\dots\wedge\o_1(t_p) =
    \frac{1}{p!}\log^p\Bigl(\frac{1-b}{1-a}\Bigr)$$
diverges as $b\to1$. It is the coefficient of the monomial $\bB^q$ in
$G_\e(1-\e)$ when $a=\e$, $b=1-\e$, and this is the reason to
call monomials that start with a $\bB$ divergent.

Similarly, the integral
$$\int\limits_{a<t_q<\dots<t_2<t_1<b}
  \o_0(t_1)\wedge\dots\wedge\o_0(t_q) =
    \frac{1}{q!}\log^q\Bigl(\frac{b}{a}\Bigr)$$
diverges as $a\to0$. It is the coefficient of the monomial $\bA^p$ in
$G_\e(1-\e)$ when $a=\e$, $b=1-\e$, and this is the reason to
call monomials that end with an $\bA$ divergent.
\end{example}

Now consider the general case: integral of a product that contains both
$\o_0$ and $\o_1$. For $\d_j=0$ or $1$ and $0<a<b<1$, introduce the
notation
$$I_{\d_1\dots\d_m}^{a,b} =
\int\limits_{a<t_m<\dots<t_2<t_1<b}
  \o_{\d_1}(t_1)\wedge\dots\wedge\o_{\d_m}(t_m)\ .$$

\begin{lemma}
\begin{itemize}
\item[(i)]  If $\d_1=0$, then the integral $I_{\d_1\dots\d_m}^{a,b}$ converges
to a non-zero constant as $b\to1$, and it grows as a power of $\log(1-b)$ if
$\d_1=1$.

\item[(ii)] If $\d_m=1$, then the integral $I_{\d_1\dots\d_m}^{a,b}$ converges
to a non-zero constant as $a\to0$, and it grows as a power of
$\log a$ if $\d_m=0$.
\end{itemize}
\end{lemma}

\begin{proof} Induction on the number of chords $m$. If $m=1$ then
the integral can be calculated explicitly like in the previous example,
and the lemma follows from the result. Now suppose that the lemma is
proved for $m-1$ chords.
By the Fubini theorem the integral can be represented as
$$I_{1\d_2\dots\d_m}^{a,b} =
  \int\limits_{a<t<b} I_{\d_2\dots\d_m}^{a,t}\cdot \frac{dt}{t-1}\ ,
\qquad
I_{0\d_2\dots\d_m}^{a,b} =
  \int\limits_{a<t<b} I_{\d_2\dots\d_m}^{a,t}\cdot \frac{dt}{t}\ ,
$$
for the cases $\d_1=1$ and $\d_1=0$ respectively. By the
induction assumption
$0<c<\left|I_{\d_2\dots\d_m}^{a,t}\right|<\left|\log^k(1-t)\right|$
for some constants $c$ and $k$. The comparison test implies that
the integral $I_{0\d_2\dots\d_m}^{a,b}$ converges as $b\to1$ because
$I_{\d_2\dots\d_m}^{a,t}$ grows slower than any power of $(1-t)$. Moreover,
$\left|I_{0\d_2\dots\d_m}^{a,b}\right|>c\int_a^1 \frac{dt}{t}
 = -c\log(a)>0$ because $0<a<b<1$.

In the case $\d_1=1$ we have $$c\log(1-b)=c\!\!\int_0^b\!
\frac{dt}{t-1}<
  \left|I_{1\d_2\dots\d_m}^{a,b}\right|<
  \left|\int_0^b\!\log^k(1-t) d(\log(1-t))\right| =
       \left|{\textstyle \frac{\log^{k+1}(1-b)}{k+1}}\right|\ ,$$
which proves assertion (i). The proof of assertion (ii) is similar.
\end{proof}

\subsection{}
Here is the plan of our subsequent actions.

Let $\Ab^{conv}$ ($\Ab^{div}$) be the subspace of $\Ab=\CAB$ spanned
by all convergent (respectively, divergent) monomials. We are going
to define a certain linear map $f:\Ab\to\Ab$
which kills divergent monomials and
preserves the associator $\Phi$. Applying $f$ to both parts of
Equation~\ref{condiv} we shall have
\begin{equation}\label{fphi}
\Phi=f(\Phi)=f\left(\lim\limits_{\e\to0}\
    \e^{-\bB}\cdot Z^{conv} \cdot \e^{\bA}\right) =
             f\left(\lim\limits_{\e\to0} Z^{conv}\right)\ .
\end{equation}
The last equality here follows from the fact that only the unit terms of
$\e^{-\bB}$ and $\e^{\bA}$ are convergent and therefore survive under the
action of $f$.

The convergent improper integral
\begin{equation}\label{limZconv}
\lim\limits_{\e\to0} Z^{conv}=1+\sum_{m=2}^\infty\ \ \sum_{\d_2,\dots,\d_{m-1}=0,1}
              I^{0,1}_{0\d_2\dots\d_{m-1}1}\cdot
     \bA C_{\d_2}\dots C_{\d_{m-1}} \bB
\end{equation}
can be computed explicitly (here $C_j=\bA$ if $\d_j=0$ and $C_j=\bB$
if $\d_j=1$). Combining Equations~\ref{fphi} and \ref{limZconv} we
get
\begin{equation}\label{phif}
\Phi=1+\sum_{m=2}^\infty\ \ \sum_{\d_2,\dots,\d_{m-1}=0,1}
              I^{0,1}_{1\d_2\dots\d_{m-1}0}\cdot
     f(\bA C_{\d_2}\dots C_{\d_{m-1}} \bB)
\end{equation}
The knowledge of how $f$ acts on the monomials from $\Ab$ leads to
the desired formula for the associator.

\subsection{Definition of the linear map $f:\Ab\to\Ab$}

Consider the algebra $\Ab[\alpha,\beta]$ of polynomials in two commuting variables
$\alpha$ and $\beta$ with coefficients in $\Ab$. Every monomial in $\Ab[\alpha,\beta]$
can be written uniquely as $\beta^p M\alpha^q$, where $M$ is a monomial in $\Ab$.
Define a $\C$-linear map $j:\Ab[\alpha,\beta]\to\Ab$ by
$j(\beta^p M\alpha^q)=\bB^p M \bA^q$. Now for any element
$\Gamma(\bA,\bB)\in\Ab$ let
$$f(\Gamma(\bA,\bB))=j(\Gamma(\bA-\alpha,\bB-\beta))\ .$$

\begin{lemma}
If $M$ is a divergent monomial in $\Ab$, then $f(M)=0$.
\end{lemma}

\begin{proof} Consider the case where $M$ starts with $\bB$, say
$M=\bB C_2\dots C_m$, where each $C_j$ is either $\bA$ or $\bB$. Then
$$f(M) = j( (\bB-\beta)M_2 ) =j(\bB M_2) - j(\beta M_2)\ ,
$$
where $M_2=(C_2-\gamma_2)\dots(C_m-\gamma_m)$ with
$\gamma_j=\alpha$ or $\gamma_j=\beta$ depending on $C_j$.
But $j(\bB M_2)$ equals $j(\beta M_2)$ by the definition of $j$ above.
The case where $M$ ends with an $\bA$ can be done similarly.
\end{proof}

\subsection{}
One may notice that for any monomial $M\in\Ab$ we have $f(M)=M +
(\mbox{sum of divergent monomials})$. Therefore, by the lemma, $f$
is an idempotent map, $f^2=f$, that is, $f$ is a projection along
$\Ab^{div}$ (but not onto $\Ab^{conv}$).

\begin{proposition} $f(\Phi)=\Phi$.
\end{proposition}

\begin{proof}
We use the definition of the associator $\Phi$ as the KZ Drinfeld
associator from Section~\ref{defDrass} (see Proposition in
Section~\ref{kieqkz}).

It is the ratio $\Phi(\bA,\bB)=G_1^{-1}\cdot G_0$ of two solutions
of the differential equation (\ref{redKZ}) from
Section~\ref{defDrass}
$$G'=\left(\frac{\bA}{z} + \frac{\bB}{z-1}\right)\cdot G
$$
with the asymptotics
$$G_0(z)\sim z^{\bA} \mbox{\ as\ } z\to0\qquad\mbox{and}\qquad
  G_1(z)\sim (1-z)^{\bB} \mbox{\ as\ } z\to1\ .
$$

Consider the differential equation
$$H'=\left(\frac{\bA-\alpha}{z} + \frac{\bB-\beta}{z-1}\right)\cdot H\ .
$$
A direct substitution shows that the functions
$$H_0(z)=z^{-\alpha}(1-z)^{-\beta}\cdot G_0(z) \qquad\mbox{and}\qquad
  H_1(z)=z^{-\alpha}(1-z)^{-\beta}\cdot G_1(z)
$$
are its solutions with the asymptotics
$$H_0(z)\sim z^{\bA-\alpha} \mbox{\ as\ } z\to0\qquad\mbox{and}\qquad
  H_1(z)\sim (1-z)^{\bB-\beta} \mbox{\ as\ } z\to1\ .
$$
Hence we have
$$\Phi(\bA-\alpha, \bB-\beta)=H_1^{-1}\cdot H_0=G_1^{-1}\cdot G_0
   = \Phi(\bA,\bB)\ .
$$
Therefore,
$$f(\Phi(\bA,\bB))=j(\Phi(\bA-\alpha, \bB-\beta))=j(\Phi(\bA,\bB))=\Phi(\bA,\bB)
$$
because $j$ acts as the identity map on the subspace
$\Ab\subset\Ab[\alpha,\beta]$. The proposition is proved.
\end{proof}

\subsection{}
In order to compute $\Phi$ according to (\ref{phif}) we must find
the integrals $I^{0,1}_{0\d_2\dots\d_{m-1}1}$ and the action of $f$
on the monomials. Let us compute $f(\bA C_{\d_2}\dots C_{\d_{m-1}}
\bB)$ first.

Represent the monomial $M=\bA C_{\d_2}\dots C_{\d_{m-1}} \bB$ in the form
$$M=\bA^{p_1}\bB^{q_1}\dots\bA^{p_l}\bB^{q_l}$$
for some positive integers
$p_1, q_1, \dots, p_l,q_l$. Then
$$f(M)=j((\bA-\alpha)^{p_1}(\bB-\beta)^{q_1}\dots
         (\bA-\alpha)^{p_l}(\bB-\beta)^{q_l})\ .
$$
We are going to expand the product, collect all $\beta$'s on the left and
all $\alpha$'s on the right, and then replace $\beta$ by $\bB$ and
$\alpha$ by $\bA$.
To this end let us introduce the following  multi-index notations:
$$\begin{array}{c@{\qquad}c}
\bbr=(r_1,\dots,r_l);\quad \bbi=(i_1,\dots,i_l);&
   \bbs=(s_1,\dots,s_l);\quad \bbj=(j_1,\dots,j_l);
       \vspace{10pt}\\
\bbp=\bbr+\bbi=(r_1+i_1,\dots,r_l+i_l);&
\bbq=\bbs+\bbj=(s_1+j_1,\dots,s_l+j_l); \vspace{10pt}\\
|\bbr|=r_1+\dots+r_l;& |\bbs|=s_1+\dots+s_l;  \vspace{10pt}\\
\dbinom{\bbp}{\bbr} = \dbinom{p_1}{r_1} \dbinom{p_2}{r_2}
                   \dots \dbinom{p_l}{r_l};&
\dbinom{\bbq}{\bbs} = \dbinom{q_1}{s_1} \dbinom{q_2}{s_2}
                   \dots \dbinom{q_l}{s_l}; \vspace{10pt}\\
\multicolumn{2}{c}{
(\bA,\bB)^{(\bbi,\bbj)} = \bA^{i_1}\cdot \bB^{j_1}\cdot\dots\cdot
                      \bA^{i_l}\cdot \bB^{j_l}}
\end{array}
$$
We have
$$\begin{array}{l}
(\bA-\alpha)^{p_1}(\bB-\beta)^{q_1}\dots
(\bA-\alpha)^{p_l}(\bB-\beta)^{q_l}
        = \vspace{10pt}\\ \displaystyle\hspace{90pt}
\sum_{\substack{0\leqslant\bbr\leqslant\bbp \\
                0\leqslant\bbs\leqslant\bbq}}
 (-1)^{|\bbr|+|\bbs|} \dbinom{\bbp}{\bbr} \dbinom{\bbq}{\bbs}
    \cdot \beta^{|\bbs|} (\bA,\bB)^{(\bbi,\bbj)} \alpha^{|\bbr|}\ ,
\end{array}
$$
where the inequalities
$0\leqslant\bbr\leqslant\bbp$ and $0\leqslant\bbs\leqslant\bbq$ mean
$0\leqslant r_1\leqslant p_1$, \dots, $0\leqslant r_l\leqslant p_l$, and
$0\leqslant s_1\leqslant q_1$,\dots, $0\leqslant s_l\leqslant q_l$.
Therefore
\begin{equation}\label{fofM}
f(M)=
 \sum_{\substack{0\leqslant\bbr\leqslant\bbp \\
                0\leqslant\bbs\leqslant\bbq}}
 (-1)^{|\bbr|+|\bbs|} \dbinom{\bbp}{\bbr} \dbinom{\bbq}{\bbs}
    \cdot \bB^{|\bbs|} (\bA,\bB)^{(\bbi,\bbj)} \bA^{|\bbr|}\ .
\end{equation}

\subsection{}\label{secetazeta}

To complete the formula for the associator we need to compute the
coefficient $I^{0,1}_{1\d_2\dots\d_{m-1}0}$ of $f(M)$.
It turns out that, up to a sign, they are equal to some values
of the multivariate $\z$-function
\index{Multivariate $\z$-function}
\index{Multiple zeta values}
\index{MZV}
$$\z(a_1,\dots,a_n)\ =\
     \sum_{0<k_1<k_2<\dots<k_n} k_1^{-a_1}\dots k_n^{-a_n}$$
where $a_1$, ..., $a_n$ are positive integers (see \cite{LM1}).
Namely, the coefficients in question
are equal, up to a sign, to the values of $\z$ at integer points
$(a_1,\dots,a_n)\in\Z^n$, which are called ({\it multiple zeta values},
or MZV for short). Multiple zeta values for $n=2$ were first studied by
L.~Euler in 1775. His
paper \cite{Eu} contains several dozens interesting relations between
MZVs and values of the univariate Riemann's zeta function.
Later, this subject was almost forgotten for more than 200 years
until M.~Hoffman and D.~Zagier revived a general interest to MZVs by their
papers \cite{Hoff}, \cite{Zag3}.

\noindent\textbf{Exercise.} The sum in the definition of the
multivariate $\z$-function converges if and only if
$a_n\geqslant2$.

\begin{remark}
Two different conventions about the order of arguments in $\z$ are
in use: we follow that of D.\,Zagier \cite{Zag3}, also used by
P.\,Deligne, A.\,Goncharov and Le--Murakami \cite{LM1,LM2,LM3,LM4}.
The opposite school that goes back to L.\,Euler \cite{Eu} and
includes J.~Borwein, M.~Hoffman, M.~Petitot, P.~Cartier, writes
$\z(2,1)$ for what we would write as $\z(1,2)$.
\end{remark}

\begin{proposition} \label{propetazeta}
{\it For $\bbp>0$ and $\bbq>0$ let
\begin{equation}\label{etazeta}
\eta(\bbp,\bbq):=
 \z(\underbrace{1,\dots,1}_{q_l-1},\ p_l+1,\
    \underbrace{1,\dots,1}_{q_{l-1}-1},\ p_{l-1}+1,\ \dots\
    \underbrace{1,\dots,1}_{q_1-1},\ p_1+1)\ .
\end{equation}\index{$\eta(\bbp,\bbq)$}
Then
\begin{equation}\label{intzetaval}
I^{0,1}_{\underbrace{\scriptstyle 0\dots0}_{p_1}
           \underbrace{\scriptstyle 1\dots1}_{q_1}\dots\dots
           \underbrace{\scriptstyle 0\dots0}_{p_l}
           \underbrace{\scriptstyle 1\dots1}_{q_l}}
  = (-1)^{|\bbq|} \eta(\bbp,\bbq)\ .
\end{equation}}
\end{proposition}

The calculations needed to prove the proposition, are best organised
in terms of the (univariate)
\textit{polylogarithm}\index{Polylogarithm}\footnote{It is a
generalization of Euler's dilogarithm ${\rm Li}_2(z)$ we used on
page \pageref{dilogarithm}, and a specialization of the multivariate
polylogarithm
$${\rm Li}_{a_1,\dots,a_n}(z_1,\dots,z_n) = \sum_{0<k_1<k_2<\dots<k_n}
    \frac{z_1^{k_1}\dots z_n^{k_n}}{k_1^{a_1}\dots k_n^{a_n}}$$
introduced by A.~Goncharov in \cite{Gon1}. }
\index{Polylogarithm}\index{Multiple polylogarithm}
function defined by the series
\begin{equation}\label{mplog}
{\rm Li}_{a_1,\dots,a_n}(z) =
    \sum_{0<k_1<k_2<\dots<k_n} \frac{z^{k_n}}{k_1^{a_1}\dots k_n^{a_n}}\ ,
\end{equation}
which obviously converges for $|z|<1$.

\begin{lemma} For $|z|<1$
$${\rm Li}_{a_1,\dots,a_n}(z) = \left\{ \begin{array}{ll}
\displaystyle \int_0^z {\rm Li}_{a_1,\dots,a_n-1}(t)
         \frac{dt}{t}\ ,&\quad \mbox{if\ } a_n>1 \vspace{10pt}\\
\displaystyle -\int_0^z {\rm Li}_{a_1,\dots,a_{n-1}}(t)
         \frac{d(1-t)}{1-t}\ ,&\quad \mbox{if\ } a_n=1\ .
\end{array}\right. $$
\end{lemma}

\begin{proof} The lemma follows from the identities below, whose proofs we
leave to the reader as an exercise.
$$\begin{array}{l}
\frac{d}{dz} {\rm Li}_{a_1,\dots,a_n}(z) = \left\{
    \begin{array}{ll}
       \frac{1}{z}\cdot {\rm Li}_{a_1,\dots,a_n-1}(z)
                  \ ,&\quad \mbox{if\ } a_n>1 \vspace{8pt}\\
       \frac{1}{1-z}\cdot {\rm Li}_{a_1,\dots,a_{n-1}}(z)
                  \ ,&\quad \mbox{if\ } a_n=1\ ;
   \end{array}\right. \vspace{10pt} \\
\frac{d}{dz} {\rm Li}_1(z)= \frac{1}{1-z}\ .
\end{array}$$
\end{proof}

\subsection{Proof of Proposition \ref{propetazeta}} From the previous
lemma we have
$$\begin{array}{l}
{\rm Li}_{\underbrace{\scriptstyle 1,1,\dots,1}_{q_l-1},\ p_l+1,\
    \underbrace{\scriptstyle 1,1,\dots,1}_{q_{l-1}-1},\  p_{l-1}+1,\
         \dots,\
    \underbrace{\scriptstyle 1,1,\dots,1}_{q_1-1},\ p_1+1}(z) =
             \vspace{12pt}\\
\displaystyle
= (-1)^{q_1+\dots+q_l} \int\limits_{0<t_m<\dots<t_2<t_1<z}
\underbrace{\o_0(t_1)\wedge\dots\wedge\o_0(t_{p_1})}_{p_1}
     \wedge\vspace{8pt}\\
\hspace{5pt} \wedge
\underbrace{\o_1(t_{p_1+1})\wedge\dots\wedge\o_1(t_{p_1+q_1})}_{q_1}
         \wedge\dots\wedge
\underbrace{\o_1(t_{p_1+\dots+p_l+1})\wedge\dots\wedge
         \o_1(t_{p_1+\dots+q_l}) }_{q_l} =
            \vspace{12pt}\\
\displaystyle  \hspace{190pt}
= (-1)^{|\bbq|} I^{0,z}_{\underbrace{\scriptstyle 0\dots0}_{p_1}
           \underbrace{\scriptstyle 1\dots1}_{q_1}\dots\dots
           \underbrace{\scriptstyle 0\dots0}_{p_l}
           \underbrace{\scriptstyle 1\dots1}_{q_l}}\ .
\end{array}
$$
Note that the multiple polylogarithm series (\ref{mplog})
converges for $z=1$ in the case $a_n>1$. This implies that if
$p_1\geqslant1$ (which holds for a convergent monomial), then we have
$$\begin{array}{ccl}
\eta(\bbp,\bbq) &=&
 \z(\underbrace{1,\dots,1}_{q_l-1},\ p_l+1,\
    \underbrace{1,\dots,1}_{q_{l-1}-1},\ p_{l-1}+1,\ \dots\
    \underbrace{1,\dots,1}_{q_1-1},\ p_1+1)  \vspace{12pt}\\
&=&
 {\rm Li}_{\underbrace{\scriptstyle 1,1,\dots,1}_{q_l-1},\ p_l+1,\
    \underbrace{\scriptstyle 1,1,\dots,1}_{q_{l-1}-1},\  p_{l-1}+1,\ \dots,\
    \underbrace{\scriptstyle 1,1,\dots,1}_{q_1-1},\ p_1+1}(1)
             \vspace{12pt} \\
&=& (-1)^{|\bbq|} I^{0,1}_{\underbrace{\scriptstyle 0\dots0}_{p_1}
           \underbrace{\scriptstyle 1\dots1}_{q_1}\dots\dots
           \underbrace{\scriptstyle 0\dots0}_{p_l}
           \underbrace{\scriptstyle 1\dots1}_{q_l}}\ .
\end{array}
$$
The Proposition is proved.

\subsection{Explicit formula for the associator}

\label{LMformula} Combining equations (\ref{phif}), (\ref{fofM}),
and (\ref{intzetaval}) we get the following formula for the
associator:

\bigskip

\noindent \fbox{$\displaystyle
\PhiKZ=1\!+\!\!\sum_{m=2}^\infty\  
         \sum_{\substack{ 0<\bbp,0<\bbq\\
                          |\bbp|+|\bbq|=m}}\hspace{-5pt}
\eta(\bbp,\bbq)\cdot\hspace{-8pt}
     \sum_{\substack{0\leqslant\bbr\leqslant\bbp \\
                0\leqslant\bbs\leqslant\bbq}}
 (-1)^{|\bbr|+|\bbj|} \dbinom{\bbp}{\bbr} \dbinom{\bbq}{\bbs}
    \cdot \bB^{|\bbs|} (\bA,\bB)^{(\bbi,\bbj)} \bA^{|\bbr|}
$}

\bigskip

\noindent where $\bbi$ and $\bbj$ are multi-indices of the same
length, $\bbp=\bbr+\bbi$, $\bbq=\bbs+\bbj$, and $\eta(\bbp,\bbq)$ is
the multiple zeta value given by (\ref{etazeta}).

This formula was obtained by Le and Murakami in \cite{LM4}.

\subsection{Example. Degree 2 terms of the associator.}
There is only one possibility to represent $m=2$ as the sum of two
positive integers: $2=1+1$. So we have only one possibility for
$\bbp$ and $\bbq$: $\bbp=(1)$, $\bbq=(1)$. In this case
$\eta(\bbp,\bbq)=\z(2)=\pi^2/6$ according to (\ref{etazeta}).
The multi-indices $\bbr$ and $\bbs$ have length 1 and thus consist
of a single number
$\bbr=(r_1)$ and $\bbs=(s_1)$. There are two possibilities for
each of them: $r_1=0$ or $r_1=1$, and $s_1=0$ or $s_1=1$. In all
these cases the binomial coefficients $\dbinom{\bbp}{\bbr}$ and
$\dbinom{\bbq}{\bbs}$ are equal to 1. We arrange all the possibilities in
the following table.
$$\begin{array}{c|c|c|c||c}
r_1 & s_1 & i_1 & j_1 &
    (-1)^{|\bbr|+|\bbj|}\cdot \bB^{|\bbs|}(\bA,\bB)^{(\bbi,\bbj)} \bA^{|\bbr|} \\
\hline\hline
0 & 0 & 1 & 1 & -\bA \bB \makebox(0,12){}\\  \hline
0 & 1 & 1 & 0 &  \bB \bA \makebox(0,12){}\\  \hline
1 & 0 & 0 & 1 &  \bB \bA \makebox(0,12){}\\  \hline
1 & 1 & 0 & 0 & -\bB \bA \makebox(0,12){}
\end{array}
$$
Hence, for the degree 2 terms of the associator we get the formula:
$$-\z(2)[\bA,\bB]=-\frac{\z(2)}{(2\pi i)^2}[a,b]=
           \frac{1}{24}[a,b]\ ,$$
where $a=(2\pi i)A=\vstod$, and $b=(2\pi i)B=\vstdt$.

\subsection{Example. Degree 3 terms of the associator.}
There are two ways to represent $m=3$ as the sum of two positive
integers: $3=2+1$ and $3=1+2$. In each case either $\bbp=(1)$ or
$\bbq=(1)$. Hence $l=1$ and both multi-indices consist of just
one number $\bbp=(p_1)$, $\bbq=(q_1)$. Therefore all other
multi-indices $\bbr$, $\bbs$, $\bbi$, $\bbj$ also consist of one
number.

Here is the corresponding table.
$$\begin{array}{@{\!}c|c|c||c|c|c|c||c}
p_1 & q_1 & \eta(\bbp,\bbq) & r_1 & s_1 & i_1 & j_1 &
    (-1)^{|\bbr|+|\bbj|} \dbinom{\bbp}{\bbr} \dbinom{\bbq}{\bbs}
    \cdot \bB^{|\bbs|} (\bA,\bB)^{(\bbi,\bbj)} \bA^{|\bbr|}
            \rb{-15pt}{\makebox(0,0){}}\\
\hline\hline
&&& 0 & 0 & 2 & 1 & -\bA \bA \bB \makebox(0,12){}\\
&&& 0 & 1 & 2 & 0 & \bB \bA \bA \makebox(0,12){}\\
&&& 1 & 0 & 1 & 1 & 2\bA \bB \bA \makebox(0,12){}\vspace{-10pt}\\
\makebox(4,0){2} & \makebox(4,0){1} &\makebox(4,0){$\z(3)$} &&&&
    \vspace{-4pt}\\
&&& 1 & 1 & 1 & 0 & -2\bB \bA \bA \makebox(0,12){}\\
&&& 2 & 0 & 0 & 1 & -\bB \bA \bA \makebox(0,12){}\\
&&& 2 & 1 & 0 & 0 & \bB \bA \bA \makebox(0,12){}\\ \hline
&&& 0 & 0 & 1 & 2 &  \bA \bB \bB \makebox(0,12){}\\
&&& 1 & 0 & 0 & 2 &  -\bB \bB \bA \makebox(0,12){}\\
&&& 0 & 1 & 1 & 1 &  -2\bB \bA \bB \makebox(0,12){}\vspace{-10pt}\\
\makebox(4,0){1} & \makebox(4,0){2} &\makebox(4,0){$\z(1,2)$} &&&& \vspace{-4pt}\\
&&& 1 & 1 & 0 & 1 & 2\bB \bB \bA \makebox(0,12){}\\
&&& 0 & 2 & 1 & 0 & \bB \bB \bA \makebox(0,12){}\\
&&& 1 & 2 & 0 & 0 & -\bB \bB \bA \makebox(0,12){}\\
\end{array}
$$
Using the Euler identity $\z(1,2)=\z(3)$ (see
section \ref{relMZV}) we can sum up the degree 3 part of $\Phi$ into the
formula
$$\begin{array}{l}
\z(3)\Bigl( -\bA \bA \bB +2\bA \bB \bA - \bB \bA \bA
       +\bA \bB \bB -2\bB \bA \bB +\bB \bB \bA \Bigr) \vspace{8pt}\\
\hspace{50pt} \displaystyle
    =\z(3)\left(-\Bigl[\bA,\bigl[\bA,\bB\bigr]\Bigr] -
                   \Bigl[\bB,\bigl[\bA,\bB\bigr]\Bigr] \right)
    =-\frac{\z(3)}{(2\pi i)^3}\bigl[a+b,[a,b]\bigr]\ .
\end{array}
$$

\subsection{Example. Degree 4 terms of the associator.}
There are three ways to represent $m=4$ as the sum of two positive
integers: $4=3+1$, $4=1+3$, and $4=2+2$. So we have the following
four possibilities for $\bbp$ and $\bbq$:
$$\begin{array}{c|cccc}
 \bbp & (1) & (3) & (2) & (1,1) \\ \hline
 \bbq & (3) & (1) & (2) & (1,1)
\end{array}$$

The table for the multi-indices $\bbr$, $\bbs$, $\bbp$, $\bbq$ and
the corresponding term
$$ T = (-1)^{|\bbr|+|\bbj|} \dbinom{\bbp}{\bbr} \dbinom{\bbq}{\bbs}
    \cdot \bB^{|\bbs|} (\bA,\bB)^{(\bbi,\bbj)} \bA^{|\bbr|}\ $$
is shown on page \pageref{deg4}.

\begin{figure}\label{deg4}
$$\begin{array}{c|c|c||c|c|c|c||c}
\bbp & \bbq & \eta(\bbp,\bbq) & \bbr & \bbs & \bbi & \bbj & T\\
\hline\hline
&&& (0) & (0) & (1) & (3) & - \bA \bB \bB \bB\\
&&& (0) & (1) & (1) & (2) & +3 \bB \bA \bB \bB\\
&&& (0) & (2) & (1) & (1) & -3 \bB \bB \bA \bB\\
&&& (0) & (3) & (1) & (0) & + \bB \bB \bB \bA\vspace{-10pt}\\
\makebox(10,0){(1)} & \makebox(10,0){(3)} & \makebox(4,0){$\z(1,1,2)$} &&& \vspace{-4pt}\\
&&& (1) & (0) & (0) & (3) & + \bB \bB \bB \bA\\
&&& (1) & (1) & (0) & (2) & -3 \bB \bB \bB \bA\\
&&& (1) & (2) & (0) & (1) & +3 \bB \bB \bB \bA\\
&&& (1) & (3) & (0) & (0) & - \bB \bB \bB \bA\\ \hline
&&& (0) & (0) & (3) & (1) & - \bA \bA \bA \bB\\
&&& (0) & (1) & (3) & (0) & + \bB \bA \bA \bA\\
&&& (1) & (0) & (2) & (1) & +3 \bA \bA \bB \bA\\
&&& (1) & (1) & (2) & (0) & -3 \bB \bA \bA \bA\vspace{-10pt}\\
\makebox(10,0){(3)} & \makebox(10,0){(1)} & \makebox(45,0){$\z(4)$} &&& \vspace{-4pt}\\
&&& (2) & (0) & (1) & (1) & -3 \bA \bB \bA \bA\\
&&& (2) & (1) & (1) & (0) & +3 \bB \bA \bA \bA\\
&&& (3) & (0) & (0) & (1) & + \bB \bA \bA \bA\\
&&& (3) & (1) & (0) & (0) & - \bB \bA \bA \bA\\ \hline
&&& (0) & (0) & (2) & (2) & + \bA \bA \bB \bB\\
&&& (0) & (1) & (2) & (1) & -2 \bB \bA \bA \bB\\
&&& (0) & (2) & (2) & (0) & + \bB \bB \bA \bA\\
&&& (1) & (0) & (1) & (2) & -2 \bA \bB \bB \bA\vspace{-10pt}\\
\makebox(10,0){(2)} & \makebox(10,0){(2)} & \makebox(4,0){$\z(1,3)$} &&& \vspace{-4pt}\\
&&& (1) & (1) & (1) & (1) & +4 \bB \bA \bB \bA\\
&&& (1) & (2) & (1) & (0) & -2 \bB \bB \bA \bA\\
&&& (2) & (0) & (0) & (2) & + \bB \bB \bA \bA\\
&&& (2) & (1) & (0) & (1) & -2 \bB \bB \bA \bA\\
&&& (2) & (2) & (0) & (0) & + \bB \bB \bA \bA\\ \hline
&&& (0,0) & (0,0) & (1,1) & (1,1) & + \bA \bB \bA \bB\\
&&& (0,0) & (0,1) & (1,1) & (1,0) & - \bB \bA \bB \bA\\
&&& (0,0) & (1,0) & (1,1) & (0,1) & - \bB \bA \bA \bB\\
&&& (0,0) & (1,1) & (1,1) & (0,0) & + \bB \bB \bA \bA\\
&&& (0,1) & (0,0) & (1,0) & (1,1) & - \bA \bB \bB \bA\\
&&& (0,1) & (0,1) & (1,0) & (1,0) & + \bB \bA \bB \bA\\
&&& (0,1) & (1,0) & (1,0) & (0,1) & + \bB \bA \bB \bA\\
&&& (0,1) & (1,1) & (1,0) & (0,0) & - \bB \bB \bA \bA\vspace{-10pt}\\
 \makebox(20,0){(1,1)} & \makebox(20,0){(1,1)} &
 \makebox(4,0){$\z(2,2)$} &&& \vspace{-4pt}\\
&&& (1,0) & (0,0) & (0,1) & (1,1) & - \bB \bA \bB \bA\\
&&& (1,0) & (0,1) & (0,1) & (1,0) & + \bB \bB \bA \bA\\
&&& (1,0) & (1,0) & (0,1) & (0,1) & + \bB \bA \bB \bA\\
&&& (1,0) & (1,1) & (0,1) & (0,0) & - \bB \bB \bA \bA\\
&&& (1,1) & (0,0) & (0,0) & (1,1) & + \bB \bB \bA \bA\\
&&& (1,1) & (0,1) & (0,0) & (1,0) & - \bB \bB \bA \bA\\
&&& (1,1) & (1,0) & (0,0) & (0,1) & - \bB \bB \bA \bA\\
&&& (1,1) & (1,1) & (0,0) & (0,0) & + \bB \bB \bA \bA\\ \hline
\end{array}
$$
\caption{Degree 4 terms of the associator}
\end{figure}

Combining the terms into the commutators we get the degree 4 part of
the associator $\Phi$:
$$\begin{array}{l} \displaystyle
 \z(1,1,2)\left[\bB,\Bigl[\bB,\bigl[\bB,\bA\bigr]\Bigr]\right]
+\z(4)\left[\bA,\Bigl[\bA,\bigl[\bB,\bA\bigr]\Bigr]\right]
  \vspace{10pt}\\  \hspace{80pt}
\displaystyle
+\z(1,3)\left[\bB,\Bigl[\bA,\bigl[\bB,\bA\bigr]\Bigr]\right]
  +(2\z(1,3)+\z(2,2))\bigl[\bB,\bA\bigr]^2
\end{array}
$$

Recalling that $A=\frac{1}{2\pi i}a$ and $B=\frac{1}{2\pi i}b$,
where $a$ and $b$ are the basic chord diagrams with one chord,
and using the identities from Section~\ref{relMZV}:
$$\z(1,1,2)=\z(4)=\pi^4/90,\qquad \z(1,3)=\pi^4/360,\qquad
  \z(2,2)=\pi^4/120\ ,$$
we can write out the associator $\Phi$ up to degree 4:
$$\fbox{\hspace{-3pt}$\begin{array}{ccl}
\PhiKZ&=&\displaystyle 1\ +\ \frac{1}{24}[a,b]\ -\
  \frac{\z(3)}{(2\pi i)^3}\bigl[a+b,[a,b]\bigr]\ -\
  \frac{1}{1440}\left[a,\bigl[a,[a,b]\bigr]\right]
          \vspace{10pt}\\ &&\displaystyle -\
  \frac{1}{5760}\left[a,\bigl[b,[a,b]\bigr]\right]\ -\
  \frac{1}{1440}\left[b,\bigl[b,[a,b]\bigr]\right]\ +\
  \frac{1}{1152}[a,b]^2  \vspace{10pt}\\
&& +\ (\mbox{terms of order $> 4$})\ .
\end{array}$\hspace{-3pt}}
$$

\subsection{Multiple zeta values}
\label{relMZV}

There are many  relations among MZV's and powers of $\pi$. Some of
them, like $\z(2)=\frac{\pi^2}{6}$ or $\z(1,2)=\z(3)$, were
already known to Euler. The last one can be obtained in the
following way. According to (\ref{etazeta}) and (\ref{intzetaval})
we have
$$\begin{array}{ccl}
\z(1,2)&=&\displaystyle \eta((1),(2))=I^{0,1}_{011} =
  \int\limits_{0<t_3<t_2<t_1<1}
     \o_0(t_1)\wedge\o_1(t_2)\wedge\o_1(t_3) \vspace{5pt}
     \\
 &=&\displaystyle \int\limits_{0<t_3<t_2<t_1<1}
     \frac{dt_1}{t_1}\wedge \frac{d(1-t_2)}{1-t_2}\wedge
           \frac{d(1-t_3)}{1-t_3}\ .
\end{array}
$$
The change of variables $(t_1,t_2,t_3)\mapsto (1-t_3,1-t_2,1-t_1)$
transforms the last integral to
$$\begin{array}{l} \displaystyle
\int\limits_{0<t_3<t_2<t_1<1}
     \frac{d(1-t_3)}{1-t_3}\wedge\frac{dt_2}{t_2}\wedge\frac{dt_1}{t_1}
      \\
 \displaystyle \hspace{20pt}= - \int\limits_{0<t_3<t_2<t_1<1}
     \o_0(t_1)\wedge\o_0(t_2)\wedge\o_1(t_3)
 = -I^{0,1}_{001} = \eta((2),(1)) = \z(3)\ .
\end{array}
$$
In the general case a similar change of variables
$$(t_1,t_2,\dots,t_m)\mapsto (1-t_m,\dots,1-t_2,1-t_1)$$
gives the identity
$$I^{0,1}_{\underbrace{\scriptstyle 0\dots0}_{p_1}
           \underbrace{\scriptstyle 1\dots1}_{q_1}\dots\dots
           \underbrace{\scriptstyle 0\dots0}_{p_l}
           \underbrace{\scriptstyle 1\dots1}_{q_l}}
= (-1)^m I^{0,1}_{\underbrace{\scriptstyle 0\dots0}_{q_l}
           \underbrace{\scriptstyle 1\dots1}_{p_l}\dots\dots
           \underbrace{\scriptstyle 0\dots0}_{q_l}
           \underbrace{\scriptstyle 1\dots1}_{p_l}}\ .
$$

By (\ref{intzetaval}), we have
\begin{eqnarray*}
I^{0,1}_{\underbrace{\scriptstyle 0\dots0}_{p_1}
           \underbrace{\scriptstyle 1\dots1}_{q_1}\dots\dots
           \underbrace{\scriptstyle 0\dots0}_{p_l}
           \underbrace{\scriptstyle 1\dots1}_{q_l}}
&=& (-1)^{|q|}\eta(\bbp,\bbq),\\
I^{0,1}_{\underbrace{\scriptstyle 0\dots0}_{q_l}
           \underbrace{\scriptstyle 1\dots1}_{p_l}\dots\dots
           \underbrace{\scriptstyle 0\dots0}_{q_l}
           \underbrace{\scriptstyle 1\dots1}_{p_l}}
&=& (-1)^{|p|}\eta(\ol\bbq,\ol\bbp),
\end{eqnarray*}
where the bar denotes the inversion of a sequence:
$\overline{\bbp}=(p_l,p_{l-1},\dots, p_1)$,
$\overline{\bbq}=(q_l,q_{l-1},\dots, q_1)$.

Since $|p|+|q|=m$, we deduce that
$$
   \eta(\bbp,\bbq)=\eta(\ol\bbq,\ol\bbp),
$$
This relation is called the {\it duality relation} \index{Duality
relation} between MZV's. After the conversion from $\eta$ to $\z$
according to Equation~\ref{etazeta}, the duality relations become
picturesque and unexpected.

\begin{figure}[hb]
 \ig[width=200pt]{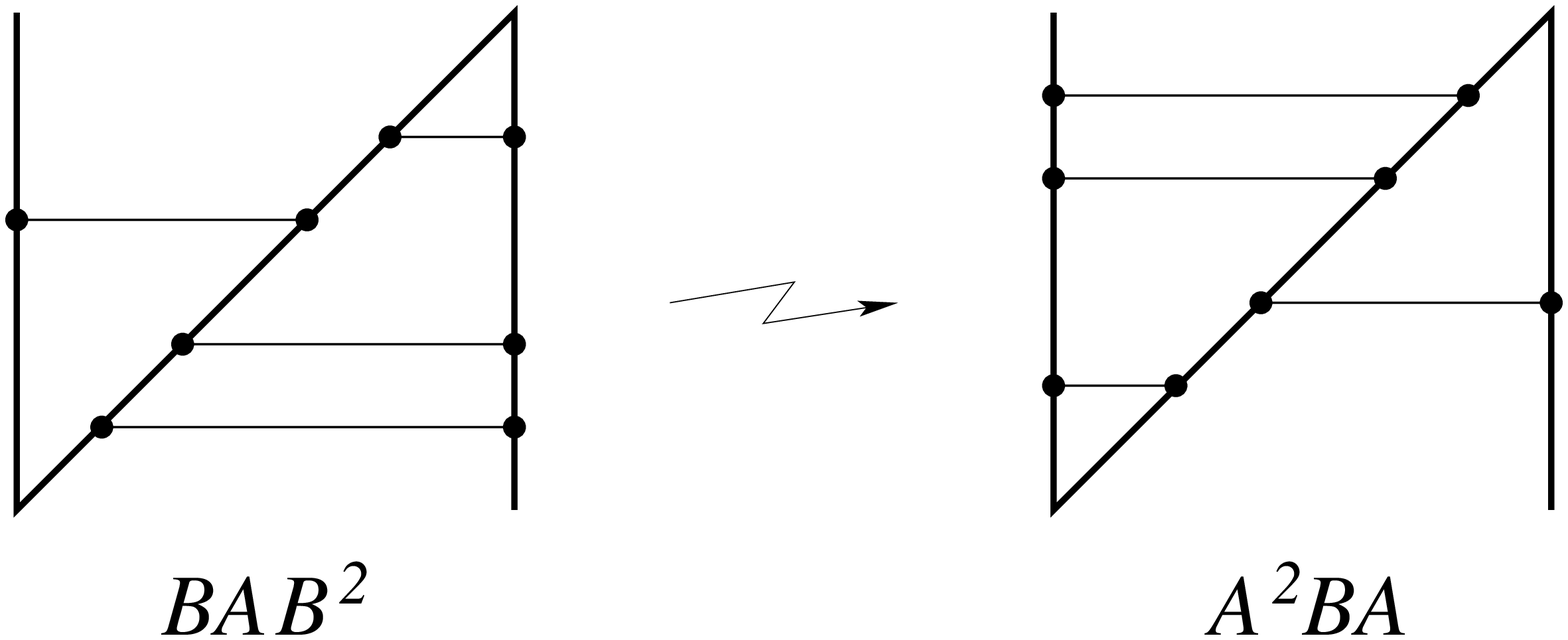}
\caption{} \label{rot_tangle}
\end{figure}

\begin{xxca} Relate the duality to the rotation of a chord diagram by
$180^\circ$ as in Figure~\ref{rot_tangle}.
\end{xxca}

As an example, we give a table of all nontrivial duality relations
of weight $m\le5$:

\def\strut{\rule[-5pt]{0pt}{18pt}}
$$\begin{array}{|c|c|c|c|c|}
\hline
\strut\bbp   & \bbq   &  \bar{\bbq}  &  \bar{\bbp}  & \mathrm{relation} \\
\hline
\strut(1)   & (2)   &  (2)  &  (1)  & \z(1,2)=\z(3) \\
\hline
\strut(1)   & (3)   &  (3)  &  (1)  & \z(1,1,2)=\z(4)\\
\hline
\strut(1)   & (4)   &  (4)  &  (1)  & \z(1,1,1,2)=\z(5)\\
\hline
\strut(2)   & (3)   &  (3)  &  (2)  & \z(1,1,3)=\z(1,4)\\
\hline
\strut(1,1) & (1,2) & (2,1) & (1,1) & \z(1,2,2)=\z(2,3)\\
\hline
\strut(1,1) & (2,1) & (1,2) & (1,1) & \z(2,1,2)=\z(3,2)\\
\hline
\end{array}$$

The reader may want to check this table by way of exercise.

There are other relations between the multiple zeta values
that do not follow from the duality law. Let us quote just a few.

1. Euler's relations:
\begin{gather}
\z(1,n-1) + \z(2,n-2) + \dots + \z(n-2,2) = \z(n),\label{rel1}\\
\z(m)\cdot\z(n)=\z(m,n)+\z(n,m)+\z(m+n)\label{rel2}\ .
\end{gather}

2. Relations obtained by Le and Murakami \cite{LM1}
computing the Kontsevich integral of the unknot by the combinatorial
procedure explained below in Section \ref{comb_ki}
(the first one was earlier proved by M.~Hoffman \cite{Hoff}):
\begin{gather}
\z(\underbrace{2,2,\dots,2}_{m})=\frac{\pi^{2m}}{(2m+1)!}\label{rel3}\\
\Bigl(\frac{1}{2^{2n-2}}-1\Bigr)\z(2n) - \z(1,2n-1)
   + \z(1,1,2n-1) - \dots \label{rel4} \\
   + \z(\underbrace{1,1,\dots,1}_{2n-2},2)= 0\ .\notag
\end{gather}

These relations are sufficient to express all multiple zeta values
with the sum of arguments equal to 4 via powers of $\pi$. Indeed, we
have:
\begin{gather*}
  \z(1,3)+\z(2,2)=\z(4),\\
  \z(2)^2=2\z(2,2)+\z(4),\\
  \z(2,2)=\frac{\pi^4}{120},\\
  -\frac{3}{4}\z(4)-\z(1,3)+\z(1,1,2)=0.
\end{gather*}
Solving these equations one by one and using the identity
$\z(2)=\pi^2/6$, we find the values of all MZVs of weight 4:
$\z(2,2)=\pi^4/120$,
$\z(1,3)=\pi^4/360$,
$\z(1,1,2)=\z(4)=\pi^4/90$.

There exists an extensive literature about the relations between MZV's,
for instance \cite{BBBL,Car2,Hoff,HoOh,OU}, and the interested
reader is invited to consult it.

\medskip

An attempt to overview the whole variety of relations between MZV's was
undertaken by D.~Zagier \cite{Zag3}.
Call the \textit{weight}\index{Weight!of a MZV} of a multiple zeta value
$\z(n_1,\dots,n_k)$ the sum of all its arguments $w=n_1+\dots+n_k$.
Let ${\mathcal Z}_w$ be the vector subspace of the reals $\R$ over the
rationals $\Q$ spanned by all MZV's of a fixed weight $w$.
For completeness we put ${\mathcal Z}_0=\Q$ and ${\mathcal Z}_1=0$.
No inhomogeneous relations between the MZV's of different weight
are known, so that conjecturally the sum of all ${\mathcal Z}_i$'s is
direct. In any case, we can consider the \textit{formal}
direct sum of all ${\mathcal Z}_w$
$$
{\mathcal Z}_\bullet:=\bigoplus\limits_{w\geq0} {\mathcal Z}_w.
$$

\begin{xproposition} The vector space
${\mathcal Z}_\bullet$ forms a graded algebra over $\Q$, i.e.
${\mathcal Z}_u\cdot {\mathcal Z}_v \subseteq {\mathcal Z}_{u+v}$.
\end{xproposition}

Euler's product formula (\ref{rel2}) illustrates this statement.
A proof can be found in \cite{Gon1}. D.~Zagier made a conjecture about the
Poincar\'e series of this algebra.

\begin{xconjecture}[\cite{Zag3}]
$$\sum_{w=0}^\infty \dim_\Q({\mathcal Z}_w)\cdot t^w
      =\frac{1}{1-t^2-t^3}\ ,
$$
which is equivalent to say that
$\dim{\mathcal Z}_0=\dim{\mathcal Z}_2=1$, $\dim{\mathcal Z}_1=0$ and
$\dim{\mathcal Z}_w=\dim{\mathcal Z}_{w-2}+\dim{\mathcal Z}_{w-3}$
for all $w\geq3$.
\end{xconjecture}

This conjecture turns out to be related  to the
dimensions of various subspaces in the primitive space of the chord
diagram algebra $\A$ (see \cite{Br,Kre}) and also to the Drinfeld's
conjecture about the structure of the Lie algebra of the
Grothendieck--Teichm\"{u}ller group \cite{ES}.

It is known (\cite{Gon2, Ter}) that Zagier's sequence gives an upper
bound on the dimension of ${\mathcal Z}_w$; in fact, up to weight 12
any zeta-number can be written as a rational polynomial in
$$
\z(2),\z(3),\z(5),\z(7),\z(2,6),\z(9),\z(2,8),
\z(11),\z(1,2,8),\z(2,10),\z(1,1,2,8).
$$
More information about the generators of the algebra ${\mathcal Z}$
is available on the web pages of M.\,Petitot \cite{Pet} and J.\,Vermaseren
\cite{Ver}.

\subsection{Logarithm of the KZ associator}
\label{ln_ass}

The associator $\PhiKZ$ is group-like \index{Group-like element}
(see Exercise~\ref{ass-gr-l} at the end of the chapter). Therefore
its logarithm can be expressed as a Lie series in variables $A$ and
$B$. Let $L$\index{$L$} be the completion of a free Lie algebra
\index{Free Lie algebra $L$}\index{Lie algebra!free $L$} generated
by $A$ and $B$.

An explicit expression for $\log\PhiKZ$ up to degree 6 was first
written out in \cite{MPH}, and up to degree 12, in \cite{Du3}. The
last formula up to degree 7 in the variables $A$ and $B$ is shown on
page \pageref{ass7}. We use the shorthand notations
$$
\zeta_n=\zeta(n)\ ,\ C_{kl}=\mbox{ad}_B^{k-1}\ \mbox{ad}_A^{l-1}\
[A,B].
$$

\begin{figure}\label{ass7}
\begin{align*}   \log(\PhiKZ)
 &= -\z_2\, C_{11} -\z_3\,\left(C_{12}+C_{21}\right)\\
 & -\frac{2}{5}\z_2^2\,\left(C_{13}+C_{31}\right)-\frac{1}{10}\z_2^2\, C_{22}\\
 &-\z_5 \,\left(C_{14}+C_{41}\right)+\left(\z_2\z_3-2\z_5\right)\,\left(C_{23}+C_{32}\right)\\
 &+\frac{\z_2\z_3-\z_5}{2} \,[C_{11},C_{12}]
 +\frac{\z_2\z_3-3\z_5}{2} \,[C_{11},C_{21}]\\
&-\frac{8}{35}\z_2^3 \,\left(C_{15}+C_{51}\right)
+\left(\frac{1}{2}\z_3^2-\frac{6}{35}\z_2^3\right)
\,\left(C_{24}+C_{42}\right)
+\left(\z_3^2-\frac{23}{70}\z_2^3\right)\, C_{33}\\
&+\left(-\frac{19}{105}\z_2^3+\z_3^2\right) \,[C_{11},C_{13}]
+\left(-\frac{69}{140}\z_2^3+\frac{3}{2}\z_3^2\right) \,[C_{11},C_{22}]\\
&+\left(-\frac{17}{105}\z_2^3\right) \,[C_{11},C_{31}]
+\left(\frac{2}{105}\z_2^3-\frac{1}{2}\z_3^2\right) \,[C_{12},C_{21}]\\
& -\z_7 \,\left(C_{16}+C_{61}\right)
+\left(\frac{2}{5}\z_3\z_2^2+\z_2\z_5-3\z_7\right) \,\left(C_{25}+C_{52}\right)\\
&
+\left(\frac{1}{2}\z_3\z_2^2+2\z_2\z_5-5\z_7\right) \,\left(C_{34}+C_{43}\right)\\
&
+\left(\frac{6}{5}\z_3\z_2^2+\frac{1}{2}\z_2\z_5-4\z_7\right) \,[C_{11},C_{14}]\\
&
+\left(\frac{11}{5}\z_3\z_2^2+\frac{7}{2}\z_2\z_5-13\z_7\right) \,[C_{11},C_{23}]\\
&+\left(\frac{3}{10}\z_3\z_2^2+\frac{13}{2}\z_2\z_5-12\z_7\right)
\,[C_{11},C_{32}]
+\left(\frac{5}{2}\z_2\z_5-5\z_7\right) \,[C_{11},C_{41}]\\
& +\left(\z_3\z_2^2-3\z_7\right) \,[C_{12},C_{13}]
+\left(\frac{23}{20}\z_3\z_2^2-\frac{61}{16}\z_7\right) \,[C_{12},C_{22}]\\
&
+\left(-\frac{3}{10}\z_3\z_2^2-\frac{1}{2}\z_2\z_5+\frac{19}{16}\z_7\right)
\,[C_{12},C_{31}]\\
&
+\left(\frac{4}{5}\z_3\z_2^2+\frac{5}{2}\z_2\z_5-\frac{99}{16}\z_7\right)
\,[C_{21},C_{13}]\\
&
+\left(\frac{7}{20}\z_3\z_2^2+6\z_2\z_5-\frac{179}{16}\z_7\right) \,[C_{21},C_{22}]\\
&
+\left(-\frac{1}{5}\z_3\z_2^2+2\z_2\z_5-3\z_7\right) \,[C_{21},C_{31}]\\
&
+\left(\frac{67}{60}\z_3\z_2^2+\frac{1}{4}\z_2\z_5-\frac{65}{16}\z_7\right)
\,[C_{11},[C_{11},C_{12}]]\\
&
+\left(-\frac{1}{12}\z_3\z_2^2+\frac{3}{4}\z_2\z_5-\frac{17}{16}\z_7\right)
\,[C_{11},[C_{11},C_{21}]]+\dots
\end{align*}
\caption{Logarithm of the Drinfeld associator up to degree 7}
\end{figure}

\begin{xremark} We have expanded the associator with respect to the {\em Lyndon
basis} of the free Lie algebra (see \cite{Reu}). There is a
remarkable one-to-one correspondence between the Lyndon words and
the irreducible polynomials over the field of two elements $\Fi_2$,
so that the associator may be thought of as a mapping from the set
of irreducible polynomials over $\Fi_2$ into the algebra of multiple
zeta values.
\end{xremark}

Now let $L'':=[[L,L],[L,L]]$ be the second commutant of the algebra $L$.
We can consider $L$ as a subspace of $\CAB$. V.~Drinfeld \cite{Dr2} proved
the following formula
$$\log\PhiKZ =
\sum_{k,l\geqslant 0} c_{kl}\ \mbox{ad}_B^l\ \mbox{ad}_A^k\ [A,B]\qquad
  (\mbox{mod}\ L'')\ ,
$$
where the coefficients $c_{kl}$ are defined by the generating function
$$1+\sum_{k,l\geqslant 0} c_{kl} u^{k+1}v^{l+1} =
\exp\left(\sum_{n=2}^\infty \frac{\z(n)}{(2\pi i)^n n}
 \bigl(u^n+v^n-(u+v)^n\bigr)\right)$$
expressed in terms of the univariate zeta function\index{Univariate zeta function}
$\z(n):=\sum_{k=1}^\infty k^{-n}$. In particular,
$c_{kl}=c_{lk}$ and $c_{k0}=c_{0k}=-\frac{\z(k+2)}{(2\pi i)^{k+2}}$.

\section{Combinatorial construction of the Kontsevich integral}
\label{comb_ki}

In this section we fulfil the promise of Section~\ref{comb_ki_intro}
and describe in detail a combinatorial construction for the
Kontsevich integral of knots and links. The associator $\PhiKZ$ is
an essential part of this construction. In 
Section~\ref{zetanumb} we shall give formulae for $\PhiKZ$; using
these expressions one can perform explicit calculations, at least in
low degrees.

\subsection{Non-associative monomials} A
{\em non-associative monomial}\, in one variable
\index{Non-associative monomial} is simply a choice of an order
(that is, a choice of parentheses) of multiplying $n$ factors; the
number $n$ is referred to as the {\em degree} of a non-associative
monomial. The only such monomial in $x$ of degree 1 is $x$ itself.
In degree 2 there is also only one monomial, namely $xx$, in degree
3 there are two monomials $(xx)x$ and $x(xx)$, in degree 4 we have
$((xx)x)x$, $(x(xx))x$, $(xx)(xx)$, $x((xx)x)$ and $x(x(xx))$, et
cetera. Define the product $u\cdot v$ of two non-associative monomials
$u$ and $v$ as their concatenation with each factor of length more than one
surrounded by an extra pair of parentheses, for instance $x\cdot x = xx$,
$xx\cdot xx =(xx)(xx)$.

For each pair $u, v$ of non-associative monomials of the same degree
$n$ one can define the element $\Phi(u,v)\in \Ab^h(n)$ as follows. If $n<3$ we
set $\Phi(u,v)={\bf 1}_n$, the unit in $\Ab^h(n)$. Assume $n\geq 3$.
Then $\Phi(u,v)$ is determined by the following properties:
\begin{enumerate}
\item
If $u=w_1\cdot(w_2\cdot w_3)$ and $v=(w_1\cdot w_2)\cdot w_3$ where
$w_1,w_2,w_3$ are
monomials of degrees $n_1,n_2$ and $n_3$ respectively, then
$$\Phi(u,v)=\Delta^{n_1,n_2,n_3}\PhiKZ,$$
where $\Delta^{n_1,n_2,n_3}$ is the cabling-type operation defined
in Exercise~\ref{op_Delta} to Chapter~\ref{chap:operations}.
\item If $w$ is monomial of degree $m$,
$$\Phi(w\cdot u,w\cdot v)={\bf 1}_m\otimes \Phi(u,v)$$
and
$$\Phi(u\cdot w,v\cdot w)=\Phi(u,v)\otimes {\bf 1}_m;$$
\item If $u,v,w$ are monomials of the same degree,
then
$$\Phi(u,v)=\Phi(u,w)\Phi(w,v).$$
\end{enumerate}
These properties are sufficient to determine $\Phi(u,v)$ since each
non-associative monomial in one variable can be obtained from any
other such monomial of the same degree by moving the parentheses in
triple products. It is not immediate that $\Phi(u,v)$ is
well-defined, however. Indeed, according to (3), we can define
$\Phi(u,v)$ by choosing a sequence of moves that shift one pair of
parentheses at a time, and have the effect of changing $u$ into $v$.
A potential problem is that there may be more than one such
sequence; however, let us postpone this matter for the moment and
work under the assumption that $\Phi(u,v)$ may be multivalued (which
it is not, see page~\pageref{rem:phiwelldef}).

Recall from Section~\ref{tangles} the notion of an elementary
tangle: basically, these are maxima, minima, crossings and vertical
segments.

Take a tensor product of several elementary tangles and choose the
brackets in it, enclosing each elementary tangle other than a
vertical segment in its own pair of parentheses. This choice of
parentheses is encoded by a non-associative monomial $w$, where each
vertical segment is represented by an $x$ and each crossing or a
critical point --- by the product $xx$. Further, we have two more
non-associative monomials, $\ol{w}$ and $\underline{w}$: $\ol{w}$ is
formed by the top boundary points of the tangle, and $\underline{w}$
is formed by the bottom endpoints. For example, consider the
following tensor product, parenthesized as indicated, of three
elementary tangles:
$$
\ig[height=2.4cm]{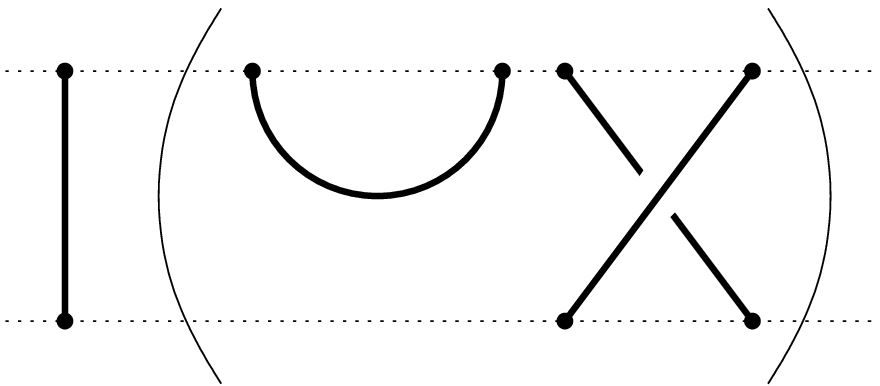}
$$
The parentheses in the product are coded by $w=x((xx)(xx))$. The top
part of the boundary gives $\ol{w}=w$, and the bottom part produces
$\underline{w}=x(xx)$.

Note that here it is important that the factors in the product are
not arbitrary, but elementary tangles, since each elementary tangle
has at most two upper and at most two lower boundary points.

\subsection{The construction}
\label{subsec:constrcombKI} First, recall that in
Exercise~\ref{op_S_i} on page~\pageref{op_S_i} we defined the
operations $S_k$ which describe how the Kontsevich integral changes
when one of the components of a tangle is reversed. Assume that the
components of the diagram skeleton are numbered. Then $S_k$ changes
the direction of the $k$th component and multiplies the diagram by
$-1$ if the number of chord endpoints lying on the $k$th component
is odd.

Represent a given knot $K$ as a product of tangles
$$K=T_1T_2\ldots T_n$$
so that each $T_i$ as a tensor product of elementary tangles:
$$T_i=T_{i,1}\otimes\dots\otimes T_{i,k_i}.$$
Write $Z_i$ for the tensor product of the Kontsevich integrals of
the elementary tangles $T_{i,j}$:
$$Z_i=Z(T_{i,1})\otimes\dots\otimes Z(T_{i,k_{i}}).$$
Note that the only elementary tangles for which the Kontsevich
integral is non-trivial are the crossings $X_-$ and $X_+$, and for
them
$$Z\left( X_+ \right) = \rb{-3mm}{\ig[height=9mm]{dzcvirt.eps}}
\cdot\exp\Bigl(\frac{\rb{-2mm}{\ig[height=4mm]{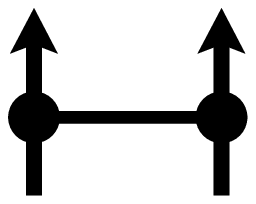}}}{2}\Bigr)
              \ ,\qquad
Z\left( X_- \right)= \rb{-3mm}{\ig[height=9mm]{dzcvirt.eps}}
\cdot\exp\Bigl(-
\frac{\rb{-2mm}{\ig[height=4mm]{doc.eps}}}{2}\Bigr).
$$
For all other elementary tangles the Kontsevich integral consists of
a diagram with no chords:
$$Z(\smaxr)=\smtan{ort_maxr},\qquad Z(\id)= \smtan{ort_id},$$
and so on. We remind that $Z_i$ in general does not coincide with
$Z(T_i)$.

For each simple tangle $T_i$ choose the parentheses in the tensor
product, and represent this choice by a non-associative monomial
$w_i$. Then the {\em combinatorial Kontsevich integral}
\index{Kontsevich integral!combinatorial} $Z_{comb}(K)$ is defined
as
$$Z_{comb}(K)= Z_1\cdot
\Phi(\underline{w}_1, \ol{w}_2)_{\downarrow^1}\cdot Z_2\cdot\ldots
\cdot Z_{n-1} \cdot \Phi(\underline{w}_{n-1},
\ol{w}_n)_{{\downarrow}^{n-1}}\cdot Z_n,$$ where
$\Phi(\underline{w}_i, \ol{w}_{i+1})_{\downarrow^i}$ is the result
of applying to $\Phi(\underline{w}_i, \ol{w}_{i+1})$ all the
operations $S_{k}$ such that at the $k$th point on the bottom of
$T_i$ (or on the top of $T_{i+1}$) the corresponding strand is
oriented downwards.

The combinatorial Kontsevich integral at the first glance may seem
to be a complicated expression. However, it is built of only two
types of elements: the exponential of $\ig[height=3mm]{doc.eps}/2$
and the Drinfeld associator $\PhiKZ$ which produces all the
$\Phi(\underline{w}_i, \ol{w}_{i+1})$.

\begin{xremark} The definition of an elementary tangle in
Section~\ref{tangles} is somewhat restrictive. In particular, of all
types of crossings only $X_+$ and $X_-$ are considered to be
elementary tangles. Note that rotating $X_+$ and $X_-$ by $\pm\pi/2$
and by $\pi$ we get tangles whose Kontsevich integral is an
exponential of the same kind as for $X_+$ and $X_-$. It will be
clear from our argument that we can count these tangles as
elementary for the definition of the combinatorial Kontsevich
integral.
\end{xremark}

\subsection{Example of computation}\label{cki-3_1}

Let us see how the combinatorial Kontsevich integral can be
computed, up to order 2, on the example of the left trefoil $3_1$.
Explicit formulae for the associator will be proved in 
Section~\ref{zetanumb}. In particular, we shall see that
$$\PhiKZ=1+\frac{1}{24}\bigl(\ \vstdtod\ -\ \vstoddt\ \bigr)+\dots.$$
Decompose the left trefoil into elementary tangles as shown below
and choose the parentheses in the tensor product as shown in the
second column:

$$\label{leftrefdecomp}
\rb{-22mm}{\ig[height=46mm]{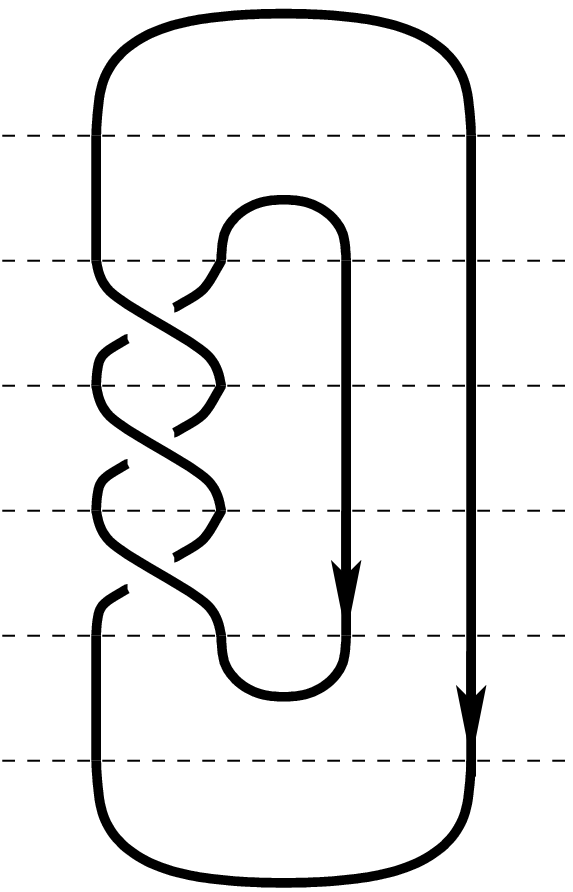}}\quad\quad
\begin{array}{l}
\medskip
\smaxr\\
\medskip
(\id\ot\smaxr)\ot\id^*\\
\medskip
(X_-\ot\id^*)\ot\id^*\\
\medskip
(X_-\ot\id^*)\ot\id^*\\
\medskip
(X_-\ot\id^*)\ot\id^*\\
\medskip
(\id\ot\minl)\ot\id^*\\
\medskip
\minl
\end{array}
$$

The combinatorial Kontsevich integral may then be represented as

$$\label{z_comb_3_1}
Z_{comb}(3_1) \quad=\quad \rb{-20mm}{\ig[height=42mm]{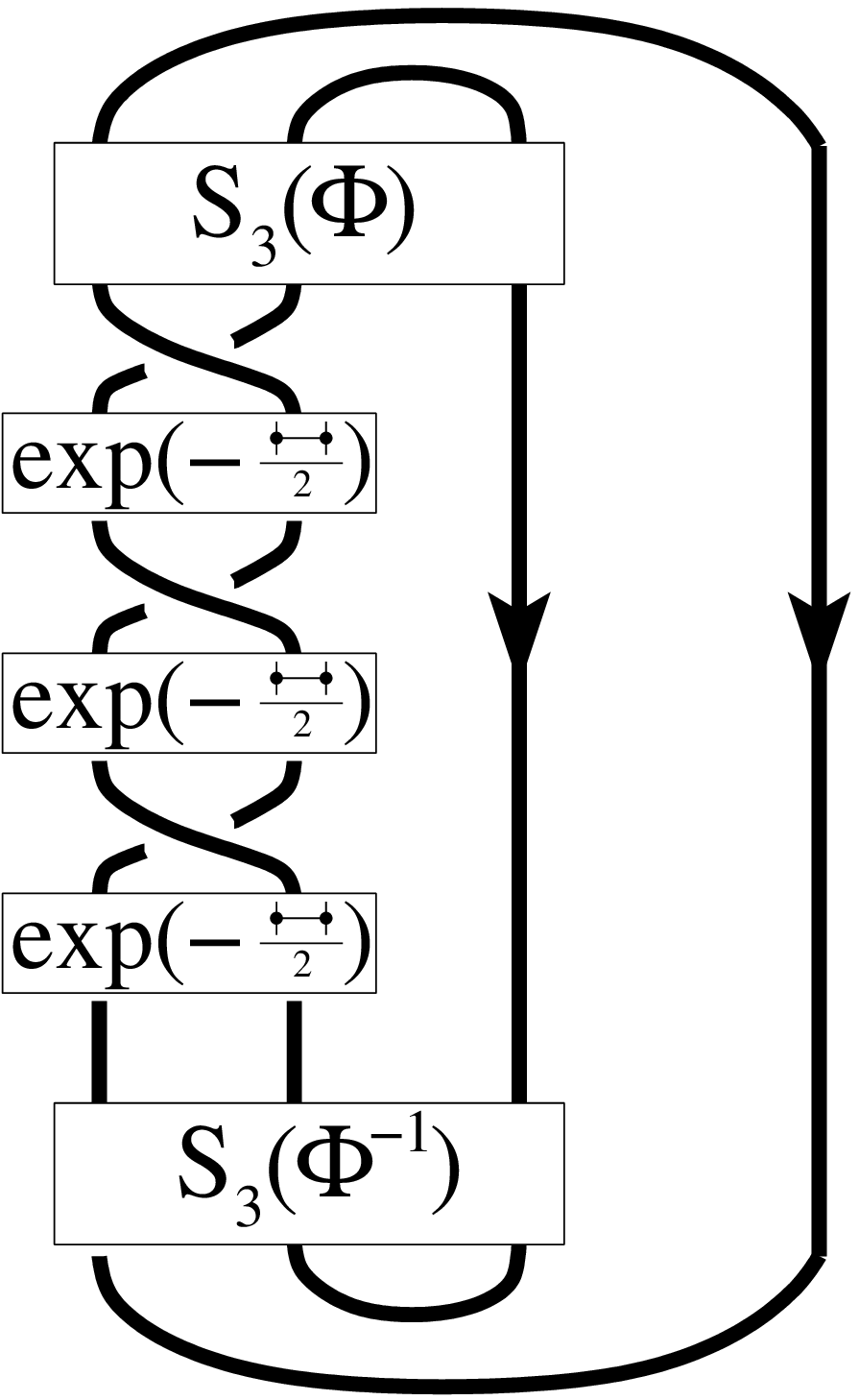}}
$$
where 
$S_3(\vstod)=\tupic{vstod-1}$ and $S_3(\vstdt)=-\tupic{vstdt-1}$.
The crossings in the above picture are, of course, irrelevant since
it shows chord diagrams and not knot diagrams.

We have that
$$S_3(\PhiKZ^{\pm 1})=1\pm \frac{1}{24}\bigl(\ \tupic{vstdtod-1}\ -\ \tupic{vstoddt-1}\ \bigr)+\dots$$
and
$$ \exp\left(\pm \frac{\risS{-3}{doc}{}{10}{0}{0}}{2}\right)=1\pm \frac{\risS{-3}{doc}{}{10}{0}{0}}{2}
+\frac{\ \risS{-3}{doc}{}{10}{0}{0}^2}{8}+\dots$$ Plugging these
expressions into the diagram above we see that, up to degree 2, the
combinatorial Kontsevich integral of the left trefoil is
$$Z_{comb}(3_1)=1+\frac{25}{24}\cdWW+\dots$$
Representing the hump $H$ as
$$\label{humpcomp}
\rb{-12mm}{\ig[height=24mm]{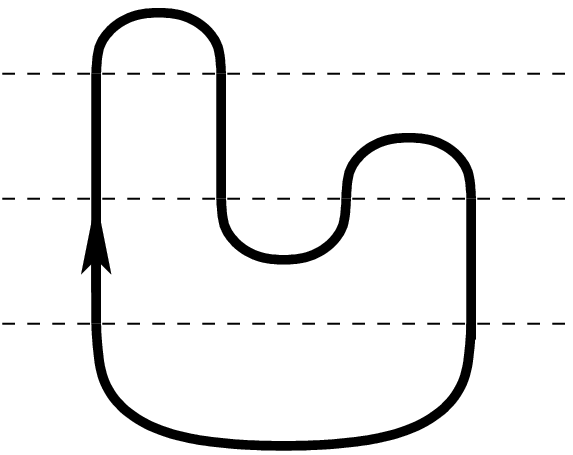}}\quad\quad
\begin{array}{l}
\medskip
\smaxr\\
\medskip
\id\ot(\id^*\ot\smaxr)\\
\medskip
\id\ot(\minr\ot\id^*)\\
\medskip
\minl
\end{array}
$$
we have
$$Z_{comb}(H)=1+\frac{1}{24}\cdWW+\dots$$
and we can speak of the {\em final} combinatorial Kontsevich
integral of the trefoil (for instance, in the multiplicative
normalization as on page~\pageref{Iprime})
\begin{multline*}
I'_{comb}(3_1)=Z_{comb}(3_1)/Z_{comb}(H)\\
=\Bigl(1+\frac{25}{24}\cdWW+\dots\Bigr)
\Bigl(1+\frac{1}{24}\cdWW+\dots\Bigr)^{-1}=1+\cdWW+\dots
\end{multline*}

\subsection{Equivalence of the
combinatorial and analytic definitions} \label{combkieqki} The main
result about the combinatorial Kontsevich integral is the following
theorem:
\begin{xtheorem}[\cite{LM3}] 
The combinatorial Kontsevich integral of a knot or a link is equal
to the usual Kontsevich integral:
$$Z_{comb}(K)=Z(K)\ .$$
\end{xtheorem}

The idea of the proof was sketched the idea in
Section~\ref{comb_ki_intro}.
The most important part of the argument consists in expressing the
Kontsevich integral of an associating tangle via $\PhiKZ$. As often
happens in our setting, this argument will give results about
objects more general than links, and we shall prove it in this
greater generality.

\subsection{The combinatorial integral for parenthesized tangles}
\label{subsec:parenthtangles} The combinatorial construction for the
Kontsevich integral on page~\pageref{subsec:constrcombKI} can also
be performed for arbitrary oriented tangles, in the very same manner
as for knots or links. However, the result of this construction can
be manifestly non-invariant.

\begin{xexample}
Take the trivial tangle ${\id}^{\ot 3}$ on 3 strands and write it as
${\id}^{\ot 3}=T_1 T_2$ where $T_1={\id}\ot (\id\ot\id)$ and
$T_2=(\id\ot\id)\ot\id$. With this choice of the parentheses, the
combinatorial Kontsevich integral of ${\id}^{\ot 3}$ is equal to the
Drinfeld associator $\PhiKZ$. On the other hand, the calculation for
${\id}^{\ot 3}=T_1$ simply gives ${\bf 1}_3$.
\end{xexample}

It turns out that the combinatorial Kontsevich integral is an
invariant of {\em parenthesized tangles}.\index{Tangle!parenthesized} A parenthesized tangle
$(T, u, v)$ is an oriented tangle $T$ together with two
non-associative monomials $u$ and $v$ in one variable, such that the
degrees of $u$ and $v$ are equal to the number of points in the
upper, and, respectively, lower, parts of the boundary of $T$. One
can think of these monomials as sets of parentheses on the boundary
of $T$.

The combinatorial Kontsevich integral $Z_{comb}$ of a parenthesized
tangle $(T, u, v)$ is defined in the same way as the Kontsevich
integral of knots or links, by decomposing $T$ into a product of
simple tangles $T_1\ldots T_n$ and choosing parentheses on the
$T_i$. The only difference is that now we require that the
bracketing chosen on $T_1$ give rise to the monomial $u$ on the top
part of $T_1$ and that the parentheses of $T_n$ produce $v$ on the
bottom of $T_n$. As usual (for instance, in
Section~\ref{subsec:kifrt}), we can define
$$I_{comb}(T)=Z(H)^{-m_1}\#\ldots\# Z(H)^{-m_k}\# Z_{comb}(T),$$
where $m_i$ is the number of maxima on the $i$th component of $T$.

It turns out that the combinatorial Kontsevich integral
$I_{comb}(T,u,v)$ depends only on the isotopy class of $T$ and on
the monomials $u,v$. This can be proved by relating
$Z_{comb}(T,u,v)$ to the Kontsevich integral of a certain family of
tangles. In order to write down the exact formula expressing this
relation, we need to define its ingredients first.

\noindent{\bf Remarks.} For tangles whose upper and lower boundary
consist of at most two points there is only one way to choose the
parentheses on the top and on the bottom, and therefore, in this
case $I_{comb}$ is an invariant of usual (not parenthesized)
tangles.

Note also that $I_{comb}$ of parenthesized tangles is preserved by
{\em all} isotopies, while the analytic Kontsevich integral is only
constant under fixed-end isotopies.

\subsection{Deformations associated with monomials and regularizing factors}

Let $t$ be set of $n$ distinct points in an interval $[a,b]$. To
each non-associative monomial $w$ of degree $n$ we can associate a
deformation $t_{\e}^w$, with $0<\e\leq 1$, as follows.

If $t_1$ and $t_2$ are two configurations of distinct points, in the
intervals $[a_1,b_1]$ and $[a_2,b_2]$ respectively, we can speak of
their $\e$-parametrized tensor product: it is obtained by rescaling
both $t_1$ and $t_2$ by $\e$ and placing the resulting intervals at
the distance $1-\e$ from each other:
\begin{multline*}\ig[height=2mm]{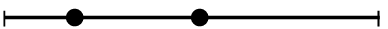}\ot_{\e}\ig[height=2mm]{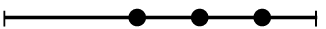}
\\
=\ig[height=2mm]{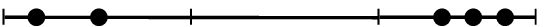}
\end{multline*} This is completely
analogous to the $\e$-parametrized tensor product of tangles on
page~\pageref{eparam}. Just as for tangles, the $\e$-parametrized
tensor product of configurations of distinct points is not
associative, and defined only up to a translation.

Now let us consider our configuration of points $t$. Divide the
interval $[a,b]$ into $n$ smaller intervals so that there is exactly
one point of $t$ in each of them, and take their $\e$-parametrized
tensor product in the order prescribed by the monomial $w$. Call the
result $t_{\e}^w$.

\noindent{\bf Exercise.} Show that $t_{\e}^w$ only depends on $t$,
$w$ and $\e$. In particular, it does not depend on the choice of the
decomposition of $t$ into $n$ intervals.
\medskip

Now, let $(T,u,v)$ be a parenthesized tangle. Denote by $s$ and $t$
the sets of top and bottom boundary points, respectively. A
continuous deformation of the boundary of a tangle can always be
extended to an horizontal deformation of the whole tangle. We shall
denote by $T_{\e,u,v}$ the family of (non-parenthesized) tangles
obtained by deforming $s$ by means of $s_{\e}^u$ and, at the same
time, deforming $t$ by means of $t_{\e}^v$.

The second ingredient we shall need is a certain function from
non-associative monomials of degree $n$  in one variable to
$\Ab^h(n)$.

First, we define for each integer $i\geq 0$ and each non-associative
monomial $w$ in $x$ the element $c_i(w)\in\Ab^h(n)$, where $n$ is
the degree of $w$, by setting
\begin{itemize}
\item $c_i(x)=0$ for all $i$;
\item $c_0(w_1w_2)=\Delta^{n_1,n_2}\bigl(\frac{\risS{-1}{doc}{}{10}{0}{0}}{2\pi i}\bigr)$ if $w_1,w_2\neq
1$, where $n_1$ and $n_2$ are the degrees of $w_1$ and $w_2$
respectively;
\item $c_i(w_1w_2)=c_{i-1}(w_1)\otimes {\bf 1}_{n_2}+
{\bf 1}_{n_1}\otimes c_{i-1}(w_2)$ if $w_1,w_2\neq 1$ with
$\deg{w_{1}}=n_{1}$, $\deg{w_{2}}=n_{2}$  and $i>0$.
\end{itemize}
It is easy to see that for each $w$ all the $c_i(w)$ commute with
each other (this follows directly from
Lemma~\ref{lemma:hordiag_center}) and that only a finite number of
the $c_i$ is non-zero. Now, we set
$$\rho^{\e}(w)=\prod_{k=1}^{\infty}\e^{k c_k(w)}.$$
This product is, of course, finite since almost all terms in it are
equal to the unit in $\Ab^h(n)$. The element $\rho^{\e}(w)$ is
called the {\em regularizing factor}\index{Regularizing factor} of $w$.

\begin{theorem}\label{thm:combki} For a parenthesized tangle $(T,u,v)$
$$\lim_{\e\to 0} \rho^{\e}(u)_{\downarrow^{t}}^{-1}\cdot Z(T_{\e,u,v})
\cdot\rho^{\e}(v)_{\downarrow^{b}}= Z_{comb}(T,u,v),$$ where
$\rho^{\e}(u)_{\downarrow^{t}}$ is the result of applying to
$\rho^{\e}(u)$ all the operations $S_{k}$ such that at the $k$th
point on the top of $T$ the corresponding strand is oriented
downwards, and $_{\downarrow^{b}}$ denotes the same operation at the
bottom of the tangle $T$.
\end{theorem}

In particular, it follows that for a knot or a link both definitions
of the Kontsevich integral coincide, since the boundary of a link is
empty.

\begin{xexample}
Let $(T,u,v)$ be the parenthesized tangle with $T$ being the trivial braid on 4 strands, 
all oriented upwards, $u=(x(xx))x$ and $v=((xx)x)x$. We have 
$$c_1(u)=\frac{1}{2\pi i}\bigl(\ \vstcod\ +\ \vstcot\ \bigr),\qquad 
c_1(v)=\frac{1}{2\pi i}\bigl(\ \vstcdt +\ \vstcot\ \bigr), $$
$$c_2(u)= \frac{1}{2\pi i}\vstcdt, \qquad c_2(v)=\frac{1}{2\pi i}\vstcod.$$
The combinatorial Kontsevich integral of $(T,u,v)$ equals to $\PhiKZ\ot \id$, and we have
\begin{multline*}
\lim_{\e\to 0} 
\e^{-\frac{1}{2\pi i}\bigl(\, \vstcod\, +\, \vstcot\, \bigr)} \e^{- \frac{1}{\pi i}\vstcdt}
\cdot Z(T_{\e,u,v})\cdot
\e^{\frac{1}{\pi i}\vstcod}\e^{\frac{1}{2\pi i}\bigl(\, \vstcdt +\, \vstcot\, \bigr)}\\
= \PhiKZ\ot \id.
\end{multline*}
This is a particular case of Exercise~\ref{regass-3-1-st} on page~\pageref{regass-3-1-st}.

\end{xexample}

\subsection{Proof of Theorem~\ref{thm:combki}} \label{sec:combkieqki}

Let $(T_1, u_1,v_1)$ and $(T_2, u_2, v_2)$ be two parenthesized
tangles with $v_1=u_2$ and such that the orientations of the strands
on the bottom of $T_1$ agree with those on top of $T_2$. Then we can
define their {\em product} to be the parenthesized tangle $(T_1T_2,
u_1, v_2)$. Every parenthesized tangle is a product of tangles
$(T,u,v)$ of three types:
\begin{enumerate}
\item {\em associating} tangles with $T$ trivial (all strands vertical, though
with arbitrary orientations) and $u\neq v$;
\item tangles where $T$ is a tensor product, in some order, of one crossing
and several vertical strands;
\item tangles where $T$ is a tensor product, in some order, of one
critical point and several vertical strands.
\end{enumerate}
In the latter two cases we require that $u$ and $v$ come from the
same choice of brackets on the elementary factors of $T$. For the
tangles of type (2) this implies that $u=v$; in the case of type (3)
tangles one monomial is obtained from the other by deleting one
factor of the form $(xx)$.

The two expressions on both sides of the equality in
Theorem~\ref{thm:combki} are multiplicative with respect to this
product, so it is sufficient to consider the three cases separately.

Let us introduce, for this proof only, the following notation. If
$x$ and $y$ are two elements of $\Ab^h(n)$ that depend on a
parameter $\e$, by saying that $x\sim y$ as $\e\to 0$ we shall mean
that in some fixed basis of $\Ab^h(n)$ (and, hence, in any basis of
this algebra) the coefficient of each diagram in $x-y$ is of the
same or smaller order of magnitude than $\e\ln^{N}{\e}$ for some
non-negative integer $N$ that may depend on the diagram. Note that
for any non-negative $N$ the limit of $\e\ln^{N}{\e}$ as $\e\to 0$
is equal to 0.

First let us consider the associating tangles. Without loss of
generality we can assume that all the strands of the tangle are
oriented upwards. We need to show that if $I =\id^{\ot n}$ is a
trivial tangle,
\begin{equation}\label{eq_zat}
\rho^{\e}(u)^{-1}\cdot Z(I_{\e,u,v})\cdot \rho^{\e}(v)\sim \Phi(u,v)
\end{equation}
 as $\e\to 0$.

\begin{xremark}\label{rem:phiwelldef}
An important corollary of the above formula is that $\Phi(u,v)$ is
well-defined, since the left-hand side is.
\end{xremark}

Let $w$ be a non-associative word and $t$ - a configuration of
distinct points in an interval. We denote by $\e t$ a configuration
of the same cardinality and in the same interval as $t$ but whose
distances between points are equal to the corresponding distances in
$t$, multiplied by $\e$. (There are many such configurations, of
course, but this is of no importance in what follows.)

Write $N_{\e}(w)$\label{ne} for a tangle with no crossings which has
$\e t_{\e}^w$ and $t_{\e}^w$ as its top and bottom configurations of
boundary points respectively, and each of whose strands connects one
point on the top to one on the bottom:
$$\ig[height=3cm]{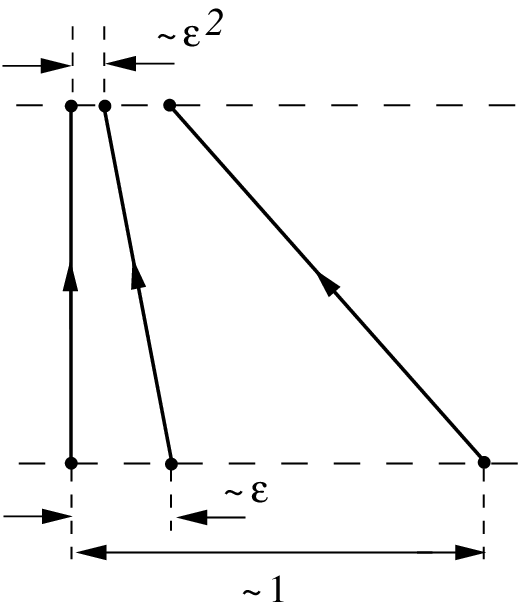}$$
As $\e$ tends to 0, the Kontsevich integral of $N_{\e}(w)$ diverges.
We have the following asymptotic formula:
\begin{equation}\label{eq_N}
Z(N_{\e}(w))\sim \prod_{k=0}^{\infty}\e^{c_k(w)}.
\end{equation}
If $t$ is a two-point configuration this formula is exact, and
amounts to a straightforward computation (see Exercise~\ref{kislant}
to Chapter~\ref{chapKI}). In general, if $w=w_1w_2$ and
$n_i=\deg{w_i}$ we can write $N_{\e}(w)$ as a product in the
following way:
$$\ig[height=3cm]{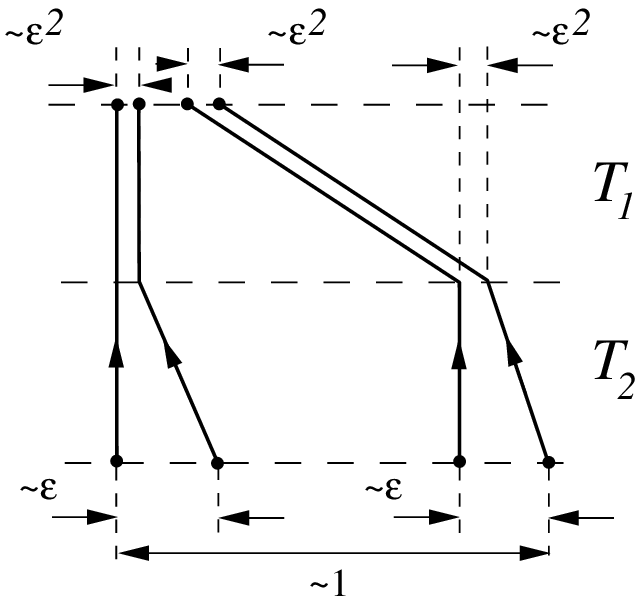}$$
As $\e$ tends to 0, we have
$$Z(T_1)\sim\Delta^{n_1,n_2}\e^{\risS{-1}{doc}{}{10}{0}{0}/2\pi i}=\e^{c_0(w_1w_2)},$$
see Exercise~\ref{op_Delta} on page~\pageref{op_Delta}, and
$$Z(T_2)\sim Z(N_{\e}(w_1))\otimes Z(N_{\e}(w_2)).$$ Using induction
and the definition of the $c_i$ we arrive to the formula
(\ref{eq_N}).

Now, notice that it is sufficient to prove (\ref{eq_zat}) in the
case when $u=w_1(w_2w_3)$ and $v=(w_1w_2)w_3$. Let us draw ${\bf
1}_{\e,u,v}$ as a product $T_1 \cdot T_2 \cdot T_3$ as in the
picture:
$$\ig[height=4cm]{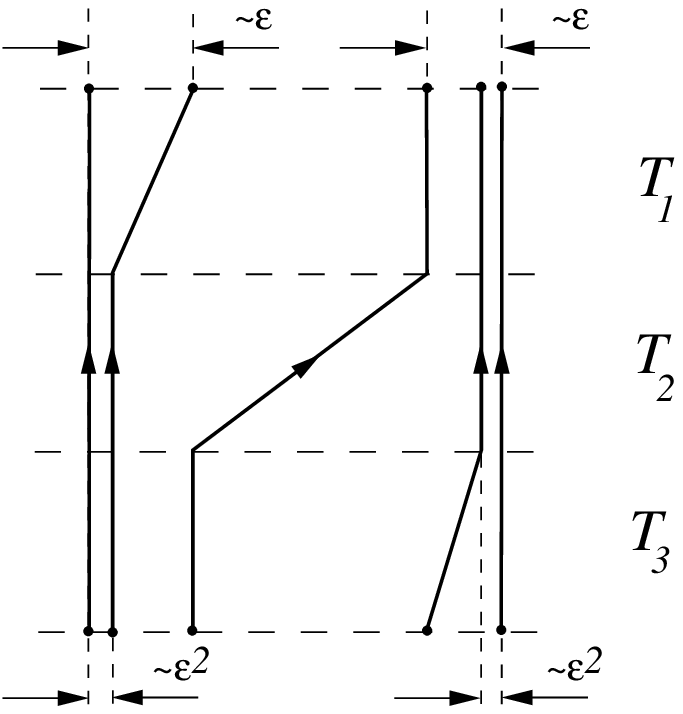}$$
As $\e \to 0$ we have:
\begin{itemize}
\item $Z(T_1)\sim Z(N_{\e}(w_1))^{-1}\otimes {\bf 1}_{n_2+n_3}$;
\item $Z(T_2) \sim ({\bf 1}_{n_1}\otimes c_0 (w_2w_3))\cdot \Delta^{n_1,n_2,n_3}\PhiKZ\cdot(c_0(w_1w_2)\otimes {\bf
1}_{n_3})^{-1}$;
\item $Z(T_3)\sim {\bf 1}_{n_1+n_2}\otimes Z(N_{\e}(w_3))$,
\end{itemize}
where $n_i=\deg{w_i}$. Notice that these asymptotic expressions for
$Z(T_1)$, $Z(T_2)$ and $Z(T_3)$ all commute with each other. Now
(\ref{eq_zat}) follows from (\ref{eq_N}) and the definition of
$\rho^{\e}(w)$.

Let us now consider the case when $T$ is a tensor product of one
crossing and several vertical strands. In this case $T_{\e,u,u}$ is
an iterated $\e$-parametrized tensor product, so
Proposition~\ref{hormult} gives $Z(T_{\e,u,u})\sim Z_{comb}(T,u,u)$.
From the definition of $c_k(u)$ we see that
$$\rho^{\e}(u)_{{\downarrow}^t}\cdot Z_{comb}(T,u,u)=
Z_{comb}(T,u,u) \cdot\rho^{\e}(u)_{{\downarrow}^b},$$ and we are
done.

Finally, let $T$ be a tensor product of one critical point (say,
minimum) and several vertical strands. For example, assume that
$T=\id^{\ot n} \ot \minr$, and that $u=v\cdot(xx)$. Now
$T_{\e,u,v}=T'_{\e}\cdot N_{\e}(v)$ where $T'_{\e}$ is the iterated
$\e$-parametrized tensor product corresponding to the monomial $u$,
so that, in particular, $T'_1=T$, and $N_{\e}(v)$ is as on
page~\pageref{ne}. As $\e\to 0$ we have
$$T'_{\e}\sim Z_{comb}(T).$$
As for the regularizing factors,
$$\rho^{\e}(u)=\rho^{\e}(v(xx))=\Bigl(\rho^{\e}(v)\ot {\id}^{\ot 2}\Bigr)\cdot
\Bigl(\prod_k
\e^{c_k(v)}\ot\e^{\frac{\risS{-1}{doc}{}{10}{0}{0}}{2\pi i}}\Bigr)$$
and we see that they cancel out together with $Z(N_{\e}(u))$.

The general case for tangle of type (3) is entirely similar.

\section{General associators}\label{axass}

We have seen that the Kontsevich integral of a knot is assembled
from Knizhnik-Zamolodchikov associators and exponentials of a
one-chord diagram on two strands. Given that the coefficients of
$\PhiKZ$ are multiple $\z$-values, the following theorem may come
as a surprise:
\begin{theorem}\label{thm:rat}
For any knot (or link) $K$ the coefficients of the Kontsevich
integral $Z(K)$, in an arbitrary basis of $\A$ consisting of chord
diagrams, are rational.
\end{theorem}
The proof of this important theorem is rather involved, and we shall
not give it here. Nevertheless, in this section we sketch very
briefly some ideas central to the argument of the proof.

\subsection{Axioms for associators}
One may ask what properties of $\PhiKZ$ imply that the combinatorial
construction indeed produces the Kontsevich integral for links. Here
we shall give a list of such properties.

Consider the algebra $\A(n)$ of tangle chord diagrams on $n$
vertical strands. (Recall that, unlike in $\A^h(n)$, the chords of
diagrams in $\A(n)$ need not be horizontal.) There are various
homomorphisms between the algebras $\A(n)$, some of which we have
already seen. Let us introduce some notation.

\begin{xdefinition}
\label{eps_i}  The operation $\e_i:\A(n)\to\A(n-1)$
\index{Operation!$\e_i$ on tangle chord diagrams} sends a tangle
chord diagram $D$ to $0$ if at least one chord of $D$ has an
endpoint on the $i$th strand; otherwise $\e_i(D)$ is obtained from
$D$ by removing the $i$th strand.
\end{xdefinition}
\noindent{\bf Examples}. $\e_i(\vstemp) =
\risS{-2}{vstwoemp}{}{13}{10}{3}$ for any $i$,
$\e_1(\vstod)=\e_2(\vstod)=0$,
$\e_3(\vstod)=\risS{-2}{vstwo}{}{13}{10}{3}$.

\medskip

The following notation is simply shorthand for
$\Delta^{1,\ldots,1,2,1,\ldots,1}$:
\begin{xdefinition}  The operation $\Delta_i:A(n)\to\A(n+1)$
\index{Operation!$\Delta_i$ on tangle chord diagrams} consists in
doubling the $i$th strand of a tangle chord diagram $D$ and taking
the sum over all possible lifts of the chord endpoints of $D$ from
the $i$th strand to one of the two new strands.
\end{xdefinition}

The symmetric group on 3 letters acts on $\A(3)$ by permuting the
strands. The action of $\sigma$ can be thought of as conjugation
$$D\to \sigma D \sigma^{-1}$$ by a strand-permuting diagram with no chords
whose $i$th point on the bottom is connected with the $\sigma(i)$th
point on top. For $D\in\A(3)$ and $\{i,j,k\}=\{1,2,3\}$ we shall
write $D^{ijk}$ for $D$ conjugated by the permutation that sends
$(\,1\,2\,3\, )$ to $(\,i\,j\,k\, )$.

All the above operations can be extended to  the graded completion
$\Ab(n)$ of the algebra $\A(n)$ with respect to the number of
chords.

Finally, we write $R$ for
$\exp\left(\frac{\risS{-1}{doc}{}{9}{0}{0}}{2}\right)\in \Ab(2)$ and
$R^{ij}$ for $\exp(u_{ij}/2)\in\Ab(3)$ where $u_{ij}$ has only one
chord that connects the strands $i$ and $j$.

\begin{xdefinition}
\index{Associator!axiomatic definition} An {\em associator} $\Phi$
is an element of the algebra $\Ab(3)$ satisfying the following
axioms:
\begin{itemize}
\item ({\em strong invertibility})
$$\e_1(\Phi)=\e_2(\Phi)=\e_3(\Phi)=1$$ (this property, in particular,
implies that the series $\Phi$ starts with 1 and thus represents an
invertible element of the algebra $\Ab(3)$);
\item ({\em skew symmetry})
$$\Phi^{-1}=\Phi^{321};$$
\item({\em pentagon relation})
\index{Pentagon relation}
$$(\id\ot\Phi)\cdot(\Delta_2\Phi)\cdot(\Phi\ot\id) =
(\Delta_3\Phi)\cdot(\Delta_1\Phi);$$
\item({\em hexagon relation})\label{hex_rel}
\index{Hexagon relation} 
$$\Phi^{231}\cdot(\Delta_2 R)\cdot\Phi =
 R^{13}\cdot\Phi^{213}\cdot R^{12}.$$
\end{itemize}
A version of the last two relations appears in abstract category
theory where they form part of the definition of a monoidal category
(see \cite[Sec.XI.1]{ML}).

\end{xdefinition}

\begin{theorem}
The Knizhnik--Zamolodchikov Drinfeld associator $\PhiKZ$ satisfies
the axioms above.
\end{theorem}

\begin{proof}
The main observation is that the pentagon and the hexagon relations
hold for $\PhiKZ$ as it can be expressed via the Kontsevich
integral. The details of the proof are as follows.

\noindent{\bf Property 1} immediately follows from the explicit
formula \ref{LMformula} for the associator $\PhiKZ$, which shows
that the series starts with 1 and every term appearing with non-zero
coefficient has endpoints of chords on each of the three strands.
\medskip

\noindent{\bf Property 2.} Notice that $\Phi^{321}$ is obtained from
$\Phi$ simply by flipping $\Phi$ about a vertical axis.
Now, consider the following tangle:
$$\ig[width=3cm]{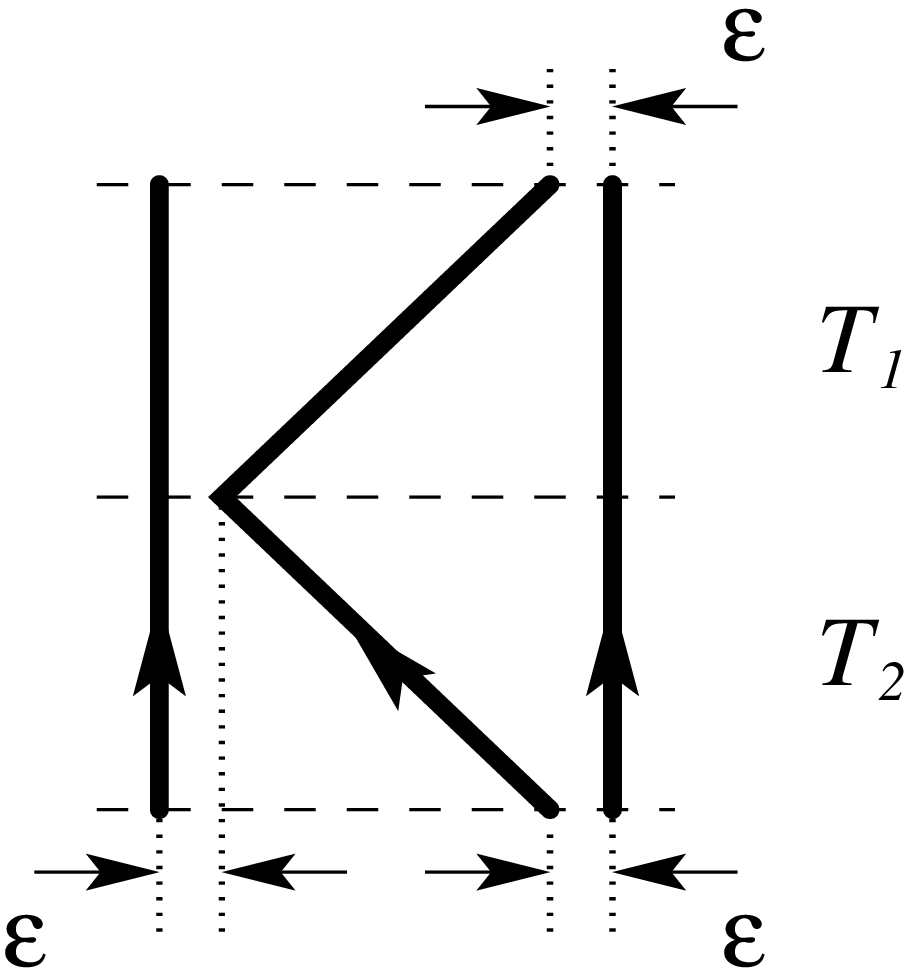}$$
It is isotopic to a tangle whose all strands are vertical, so its
Kontsevich integral is equal to 1. As we know from
Section~\ref{kieqkz}, the Kontsevich integral of the two halves of
this tangle can be expressed as
$$Z(T_1)=
   \lim\limits_{\e\to0}\
    \e^{\frac{1}{2\pi i}\vstdt}\cdot \PhiKZ\cdot \e^{-\frac{1}{2\pi i}\vstod},$$
and
$$Z(T_2)=
   \lim\limits_{\e\to0}\
    \e^{\frac{1}{2\pi i}\vstod}\cdot \PhiKZ^{321}\cdot \e^{-\frac{1}{2\pi i}\vstdt},$$
since $T_2$ is obtained from $T_1$ by flipping it about a vertical
axis. We see that the regularizing factors cancel out and
$$\PhiKZ\cdot\PhiKZ^{321}=1.$$

\noindent{\bf Property 3}. The pentagon relation for $\PhiKZ$ can be
represented by the following diagram:
$$\risS{-65}{pentagon}{
      \put(81,85){\mbox{$\scriptstyle a$}}
      \put(112,85){\mbox{$\scriptstyle b$}}
      \put(100,50){\mbox{$\scriptstyle c$}}
      \put(138,50){\mbox{$\scriptstyle d$}}
      \put(135,7){\mbox{$\scriptstyle ((ab)c)d$}}
      \put(3,57){\mbox{$\scriptstyle (a(bc))d$}}
      \put(50,115){\mbox{$\scriptstyle a((bc)d)$}}
      \put(220,115){\mbox{$\scriptstyle a(b(cd))$}}
      \put(266,57){\mbox{$\scriptstyle (ab)(cd)$}}
      \put(40,23){\mbox{$\Phi\ot\id$}}
      \put(14,88){\mbox{$\Delta_2\Phi$}}
      \put(134,120){\mbox{$\id\ot\Phi$}}
      \put(225,23){\mbox{$\Delta_1\Phi$}}
      \put(265,88){\mbox{$\Delta_3\Phi$}}
             }{300}{70}{80}
$$
Both sides of this relation are, actually, two expressions for  the
combinatorial Kontsevich integral of the trivial tangle
parenthesized by $x(x(xx))$ at the top and $((xx)x)x$ at the bottom.
On the left-hand side it is written as a product of three trivial
tangles with the monomials $x((xx)x)$ and $(x(xx))x$ in the middle.
On the right-hand side it is a product of two trivial tangles, the
monomial in the middle being $(xx)(xx)$.

\noindent{\bf Property 4}, the hexagon relation, is illustrated by
the following diagram:
$$\risS{-65}{hexagon}{
      \put(73,45){\mbox{$\scriptstyle a$}}
      \put(139,45){\mbox{$\scriptstyle c$}}
      \put(100,52){\mbox{$\scriptstyle b$}}
      \put(140,7){\mbox{$\scriptstyle (ab)c$}}
      \put(8,38){\mbox{$\scriptstyle a(bc)$}}
      \put(8,94){\mbox{$\scriptstyle (bc)a$}}
      \put(140,125){\mbox{$\scriptstyle b(ca)$}}
      \put(271,94){\mbox{$\scriptstyle b(ac))$}}
      \put(271,38){\mbox{$\scriptstyle (ba)c$}}
      \put(70,15){\mbox{$\Phi$}}
      \put(-6,67){\mbox{$\Delta_2 R$}}
      \put(70,115){\mbox{$\Phi$}}
      \put(225,15){\mbox{$R\ot\id$}}
      \put(283,67){\mbox{$\Phi$}}
      \put(225,115){\mbox{$\id\ot R$}}
             }{300}{75}{80}
$$
On both sides we have the combinatorial Kontsevich integral of the
tangle
$$T= \rb{-4mm}{\ig[height=1cm]{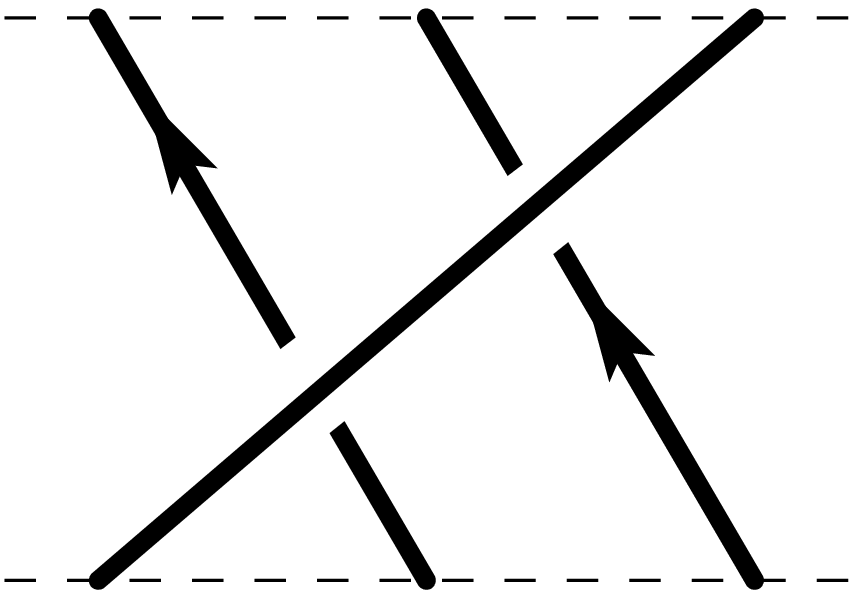}}$$ parenthesized at the top as
$x(xx)$ and at the bottom as $(xx)x$. On the right-hand side this
integral is calculated by decomposing $T$ into a product of two
crossings. On the left-hand side we use Theorem~\ref{thm:combki} and
the expression for the Kontsevich integral of $T_{\e}$ from
Exercise~\ref{op_Delta_2} on page~\pageref{op_Delta_2}. We have
$$
\begin{picture}(240, 132)(0,0)
      \put(0,0){\ig[width=240pt]{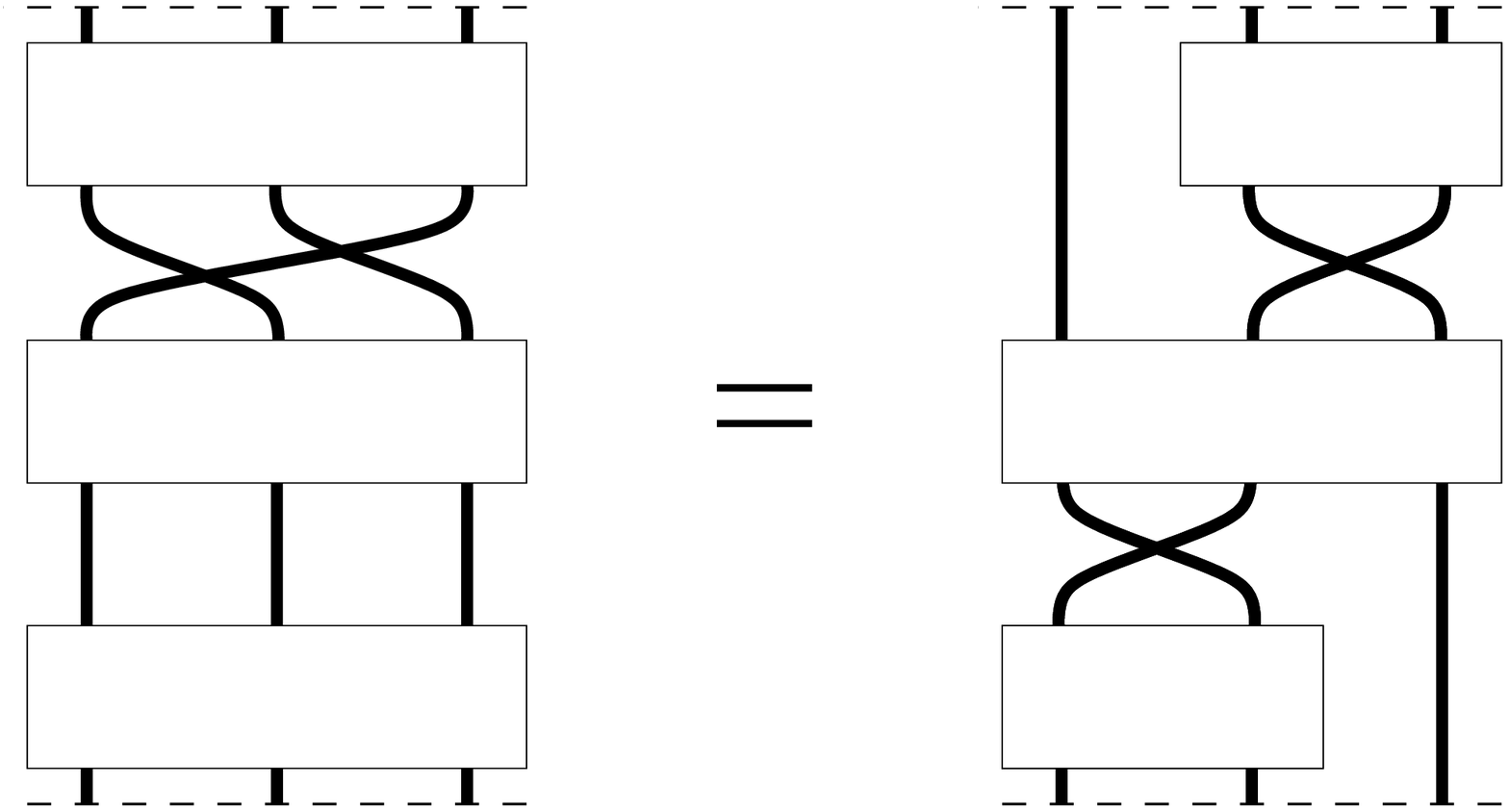}}
      \put(40,15){\mbox{$\Phi$}}
      \put(40,107){\mbox{$\Phi$}}
      \put(30,59){\mbox{$\Delta_2(R)$}}
      \put(178,15){\mbox{$R$}}
      \put(208,107){\mbox{$R$}}
      \put(194,60){\mbox{$\Phi$}}
\end{picture}
$$
which gives the hexagon relation.
\end{proof}

\subsection{The set of all associators.}\label{allass}
Interestingly, the axioms do not define the associator uniquely. The
following theorem describes the totality of all associators.
\begin{xtheorem}[\cite{Dr1,LM2}]
Let $\Phi$ be an associator and $F\in \Ab(2)$ an invertible element.
Then
$$\widetilde{\Phi}= (\id\ot F^{-1})\cdot\Delta_2(F^{-1})
   \cdot\Phi\cdot \Delta_1(F)\cdot(F\ot\id)
$$
is also an associator. Conversely, for any two associators $\Phi$
and $\widetilde{\Phi}$ there exists $F\in \Ab(2)$ invertible so that
$\Phi, \widetilde{\Phi}$ and $F$ are related as above. Moreover, $F$
can be assumed to be symmetric, that is, invariant under conjugation
by the permutation of the two strands.
\end{xtheorem}
The operation $\Phi\mapsto\widetilde{\Phi}$ is called {\em twisting}
by $F$. Diagrammatically, it looks as follows:
$$\risS{-15}{assabox}{
      \put(21,14){\mbox{$\widetilde{\Phi}$}}
             }{47}{20}{20}
\qquad=\qquad
\risS{-60}{assatwi}{
      \put(23,104){\mbox{$F^{-1}$}}
      \put(2,82){\mbox{$\Delta_2(F^{-1})$}}
      \put(20,60){\mbox{$\Phi$}}
      \put(10,38){\mbox{$\Delta_1(F)$}}
      \put(8,16){\mbox{$F$}}
             }{47}{70}{60}
$$
Twisting and the above theorem were discovered by V.~Drinfeld
\cite{Dr1} in the context of quasi-triangular quasi-Hopf algebras,
and adapted for chord diagrams in \cite{LM2}. We refer to \cite{LM2}
for the proof.

\noindent{\bf Exercise.}
 Prove that the twist by an element
$F=\exp(\a\,\risS{-3}{doc}{}{15}{0}{0}^{m})$ is the identity on any
associator for any $m$.

\medskip

\noindent{\bf Exercise.} Prove that
\begin{enumerate}
\item
twisting  by $1+\rb{-5mm}{\ig[height=12mm]{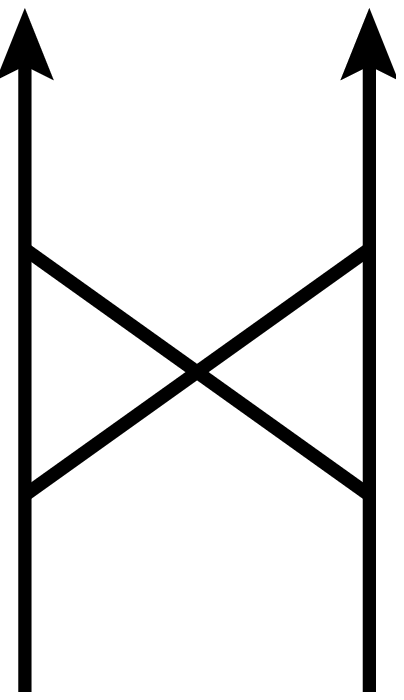}}$ adds
$2([a,b]-ac+bc)$ to the degree 2 term of an associator, where
$a=\vstod$, $b=\vstdt$ and $c=\vstot$.
\item
twisting  by $1+\rb{-5mm}{\ig[height=12mm]{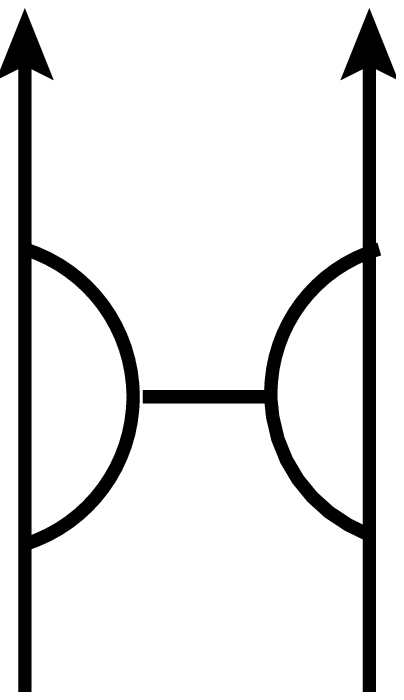}}$ does
not change the degree 3 term of an associator.
\end{enumerate}

\begin{xexample}
Let  $\Phi_{\mbox{\scriptsize BN}}$ be the rational associator
described in the next section. It is remarkable that both
$\Phi_{\mbox{\scriptsize BN}}$ and $\PhiKZ$ are horizontal, that is,
they belong to the subalgebra $\A^h(3)$ of horizontal diagrams, but
can be converted into one another only by a non-horizontal twist.
For example, twisting $\PhiBN$ by the element
$$
  F\  =\
1\ +\ \a\ \rb{-5mm}{\ig[height=12mm]{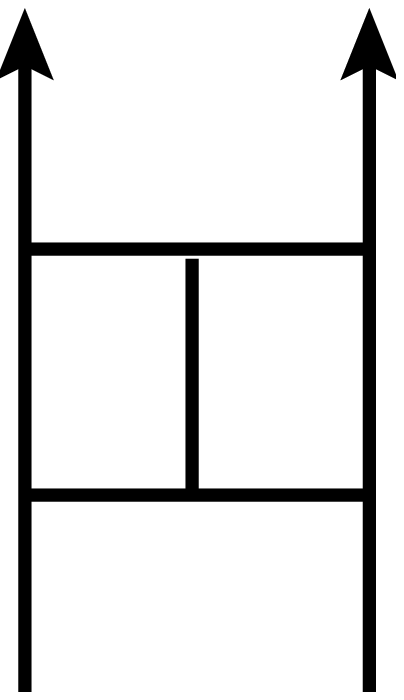}}
$$
with an appropriate constant $\a$ ensures the coincidence with
$\PhiKZ$ up to degree 4.
\end{xexample}

On the other hand, the set of all {\em horizontal} associators can also be
described in terms of the action of the so-called {\em
Grothendieck-Techm\"uller group(s)}, see \cite{Dr2, BN6}.  

V.~Kurlin \cite{Kur} described all {\em group-like} associators modulo the
second commutant.

\subsection{Rationality of the Kontsevich integral} \label{kifieqki}

Let us replace $\PhiKZ$ in the combinatorial construction of the
Kontsevich integral for parenthesized tangles by an arbitrary
associator $\Phi$; denote the result by $Z_\Phi(K)$.

\begin{xtheorem}[\cite{LM2}] For any two associators  $\Phi$
and $\widetilde{\Phi}$ the corresponding combinatorial integrals
coincide for any link $K$: $Z_\Phi(K) = Z_{\widetilde{\Phi}}(K)$.
\end{xtheorem}

A more precise statement is that for any parenthesized tangle
$(T,u,v)$ the integrals $Z_\Phi(T)$ and $Z_{\widetilde{\Phi}}(T)$
are conjugate in the sense that $Z_{\widetilde{\Phi}}(T) = {\mathcal
F}_u\cdot Z_\Phi(T)\cdot {\mathcal F}^{-1}_v$, where the elements
${\mathcal F}_u$ and ${\mathcal F}_v$ depend only on $u$ and $v$
respectively. This can be proved in the same spirit as
Theorem~\ref{thm:combki} by decomposing a parenthesized tangle into
building blocks for which the statement is easy to verify. Then,
since a link has empty boundary, the corresponding combinatorial
integrals are equal.

The fact that the Kontsevich integral does not depend on the
associator used to compute it is the key step to the proof of
Theorem~\ref{thm:rat}. Indeed, V.~Drinfeld \cite{Dr2} (see also
\cite{BN6}) showed that there exists an associator $\Phi_\Q$ with
rational coefficients. Therefore, $Z(K)=Z_{comb}(K)=Z_{\Phi_\Q}(K)$.
The last combinatorial integral has rational coefficients.

We should stress here that the existence of a rational associator
$\Phi_\Q$ is a highly non-trivial fact, and that computing it is a
difficult task. In \cite{BN2} D.~Bar-Natan, following \cite{Dr2},
gave a construction of $\Phi_\Q$  by induction on the degree. He
implemented the inductive procedure in {\tt Mathematica}
(\cite{BN5}) and computed the logarithm of the associator up to
degree 7. With the notation $a=\vstod$, $b=\vstdt$ his answer, which
we denote by $\Phi_{\mbox{\scriptsize BN}}$ is as follows:
$$\begin{array}{ccl}
\log\Phi_{\mbox{\scriptsize BN}}&=& \frac{1}{48}[ab] -
\frac{1}{1440}[a[a[ab]]]
           - \frac{1}{11520}[a[b[ab]]]  \vspace{10pt}\\  &&
+ \frac{1}{60480}[a[a[a[a[ab]]]]] +
\frac{1}{1451520}[a[a[a[b[ab]]]]]
     \vspace{10pt}\\  &&
+ \frac{13}{1161216}[a[a[b[b[ab]]]]] +
\frac{17}{1451520}[a[b[a[a[ab]]]]]
     \vspace{10pt}\\  &&
+ \frac{1}{1451520}[a[b[a[b[ab]]]]]   \vspace{10pt}\\ &&
-(\mbox{similar terms with $a$ and $b$ interchanged})+\dots
\end{array}
$$

\begin{xremark}
This expression is obtained from $\PhiKZ$ expanded to degree 7 (see
the formula on page~\pageref{ass7}) by substitutions $\z(3)\to0$,
$\z(5)\to0$, $\z(7)\to0$.
\end{xremark}

\begin{xcb}{Exercises}
\begin{enumerate}

\item
Find the monodromy of the reduced KZ equation (page~\pageref{redKZ})
around the points 0, 1 and $\infty$.

\item
Using the action of the symmetric group $S_n$ on the configuration
space $X=\C^n\setminus{\mathcal H}$ determine the algebra of values
and the KZ equation on the quotient space $X/S_n$ in such a way that
the monodromy gives the Kontsevich integral of (not necessarily
pure) braids. Compute the result for $n=2$ and compare it with
Exercise~\ref{comp_ki_R} on page~\pageref{comp_ki_R}.

\item\label{ass-gr-l}
Prove that the associator $\PhiKZ$ is group-like.

{\sl Hint.} Use the fact the Kontsevich integral is group-like.

\item Find $Z_{comb}(3_1)$ up to degree 4 using the parenthesized
presentation of the trefoil knot given in Figure~\ref{parenth_tref}
(page~\pageref{parenth_tref}).

\item Compute the Kontsevich integral of the knot $4_1$
up to degree 4, using the parenthesized presentation of the knot
$4_1$ from Exercise~\ref{parprfo} to Chapter~\ref{chapKI}
(page~\pageref{parprfo}).

\item Prove that the condition $\e_2(\Phi)=1$ and the pentagon relation imply
the other two equalities for strong invertibility:
$\e_1(\Phi)=1$ and $\e_3(\Phi)=1$.

\item
Prove the second hexagon relation
$$(\Phi^{312})^{-1}\cdot(\Delta_1 R)\cdot\Phi^{-1} =
 R^{13}\cdot(\Phi^{132})^{-1}\cdot R^{23}$$
for an arbitrary associator $\Phi$.

\item
Any associator $\Phi$ in the algebra of horizontal diagrams
$\A^h(3)$ can be written as a power series in non-commuting
variables $a=\vstod$, $b=\vstdt$, $c=\vstot$:  $\Phi=\Phi(a,b,c)$.
\begin{enumerate}
\item[(a)] Check that the skew-symmetry axiom is equivalent to the identity
$\Phi^{-1}(a,b,c)=\Phi(b,a,c)$. In particular, for an associator
$\Phi(A,B)$ with values in $\CAB$ (like $\Phi_\text{BN}$, or
$\PhiKZ$), we have $\Phi^{-1}(A,B)=\Phi(B,A)$.

\item[(b)] Prove that the hexagon relation from page~\pageref{hex_rel}
can be written in the form
$$\Phi(a,b,c)\exp\bigl(\frac{b+c}{2}\bigr)\Phi(c,a,b) =
  \exp\bigl(\frac{b}{2}\bigr)\Phi(c,b,a)\exp\bigl(\frac{c}{2}\bigr).$$

\item[(c)] ({\em V.~Kurlin \cite{Kur}})
Prove that for a horizontal associator the hexagon relation is equivalent to the relation
$$\Phi(a,b,c)e^{\frac{-a}{2}}\Phi(c,a,b)e^{\frac{-c}{2}}
  \Phi(b,c,a)e^{\frac{-b}{2}} = e^{\frac{-a-b-c}{2}}\ .
$$

\item[(d)] Show that for a horizontal associator $\Phi$,
$$\Phi\cdot\Delta_2(R)\cdot \Phi\cdot\Delta_2(R)\cdot \Phi\cdot\Delta_2(R)
=\exp(a+b+c)\ .
$$
\end{enumerate}

\item
   Express $Z(H)$ via $\PhiKZ$.
Is it true that the resulting power series contains only even degree terms?

\item\label{LM-3-18}
Prove that
$$\qquad\Phi=\lim\limits_{\e\to0}\
 \e^{-\frac{w}{2\pi i}\cdot\bigl(\vstod\ +\ \vstot\bigr)}
 \e^{-\frac{t}{2\pi i}\cdot\vstdt}\cdot
Z(AT^t_{b,w})\cdot
 \e^{\frac{b}{2\pi i}\cdot\vstod}\cdot
 \e^{\frac{w}{2\pi i}\cdot\bigl(\vstot\ +\ \vstdt\bigr)}\ ,
$$
where the tangle $AT^t_{b,w}$ is as in Exercise~\ref{ki_Phi} on
page~\pageref{ki_Phi}.

\item\label{regass-3-1-st}
\parbox[t]{2.5in}{Prove that for the tangle $T^t_{m,b,w}$ in the picture on the right \vspace{10pt}\\
$\lim\limits_{\e\to0}\
 \e^{-\frac{b-w}{2\pi i}\bigl(\vstcod\ +\ \vstcot\bigr)}\cdot
 \e^{-\frac{t-w}{2\pi i}\vstcdt}\cdot
$\vspace{-7pt}\\
$\hspace*{70pt}\cdot Z(T^t_{m,b,w})\cdot$ \vspace{4pt}
}
\qquad
\parbox{1.8in}{\makebox(115,0){\rb{-100pt}{$\qquad\quad
 \risS{-35}{assa31st}{
         \put(26,12){\mbox{$\scriptstyle \e^b$}}
         \put(60,7){\mbox{$\scriptstyle \e^w$}}
         \put(-2,34){\mbox{$\scriptstyle \e^m$}}
         \put(38,63){\mbox{$\scriptstyle \e^t$}}
                 }{110}{30}{40}$}}}\\
$\hspace*{18pt}\cdot\e^{\frac{m-w}{2\pi i}\vstcod}\cdot
 \e^{\frac{b-w}{2\pi i}\bigl(\vstcot\ +\ \vstcdt\bigr)}
\ =\ \Phi\ot\id$.

\item\label{regass-1-3-st}
\parbox[t]{2.5in}{Prove that for the tangle $T^{t,m}_{b,w}$ in the picture on the right \vspace{10pt}\\
$\lim\limits_{\e\to0}\
 \e^{-\frac{t-w}{2\pi i}\bigl(\vstcdt\ +\ \vstcdc\bigr)}\cdot
 \e^{-\frac{m-w}{2\pi i}\vstctc}\cdot
$\vspace{-7pt}\\
$\hspace*{70pt}\cdot Z(T^{t,m}_{b,w})\cdot$ \vspace{4pt}
}
\qquad
\parbox{1.8in}{\makebox(115,0){\rb{-100pt}{$\quad
 \risS{-35}{assa13}{
         \put(69,15){\mbox{$\scriptstyle \e^b$}}
         \put(40,7){\mbox{$\scriptstyle \e^w$}}
         \put(110,52){\mbox{$\scriptstyle \e^m$}}
         \put(80,64){\mbox{$\scriptstyle \e^t$}}
                 }{110}{20}{40}$}}}\\
$\hspace*{18pt}\cdot\e^{\frac{b-w}{2\pi i}\vstcdt}\cdot
 \e^{\frac{t-w}{2\pi i}\bigl(\vstcdc\ +\ \vstctc\bigr)}
\ =\ \id\ot\Phi$.

\end{enumerate}
\end{xcb}
 %12 Drinfeld

%\part{Other Topics}
\chapter{The Kontsevich integral: advanced features} % 11
\label{advKI}

\section{Mutation}
\label{mutation} The purpose of this section is to prove that the
Kontsevich integral commutes with the operation of mutation (this
fact was first noticed by T.~Le). As an application, we construct a
counterexample to the original intersection graph conjecture
(page~\pageref{IGC}) and describe, following \cite{CL}, all
Vassiliev invariants which do not distinguish mutants.

\subsection{Mutation of knots}
\label{mutkn}

Suppose we have a knot $K$ with a distinguished tangle $T$ whose
boundary consists of two points at the bottom and two points at the
top. If the orientations of the strands of $T$ agree both at the top
and the bottom of $T$, we can cut out the tangle, rotate it by
$180^\circ$ around a vertical axis and insert it back. This
operation $M_T$ is called {\em mutation}\index{Mutation} and the
knot $M_T(K)$ thus obtained is called a {\em mutant} of $K$.

Here is a widely known pair of mutant knots, $11^n_{34}$ and $11^n_{42}$,
which are mirrors of the {\em Conway} and {\em Kinoshita--Terasaka}
knots respectively:\label{C_KT_kn}
\index{Knot!Conway}\index{Knot!Kinoshita--Terasaka} \index{Conway
knot}\index{Kinoshita--Terasaka knot}
$$11^n_{34}=\ol{C} =  \rb{-13mm}{\ig[height=28mm]{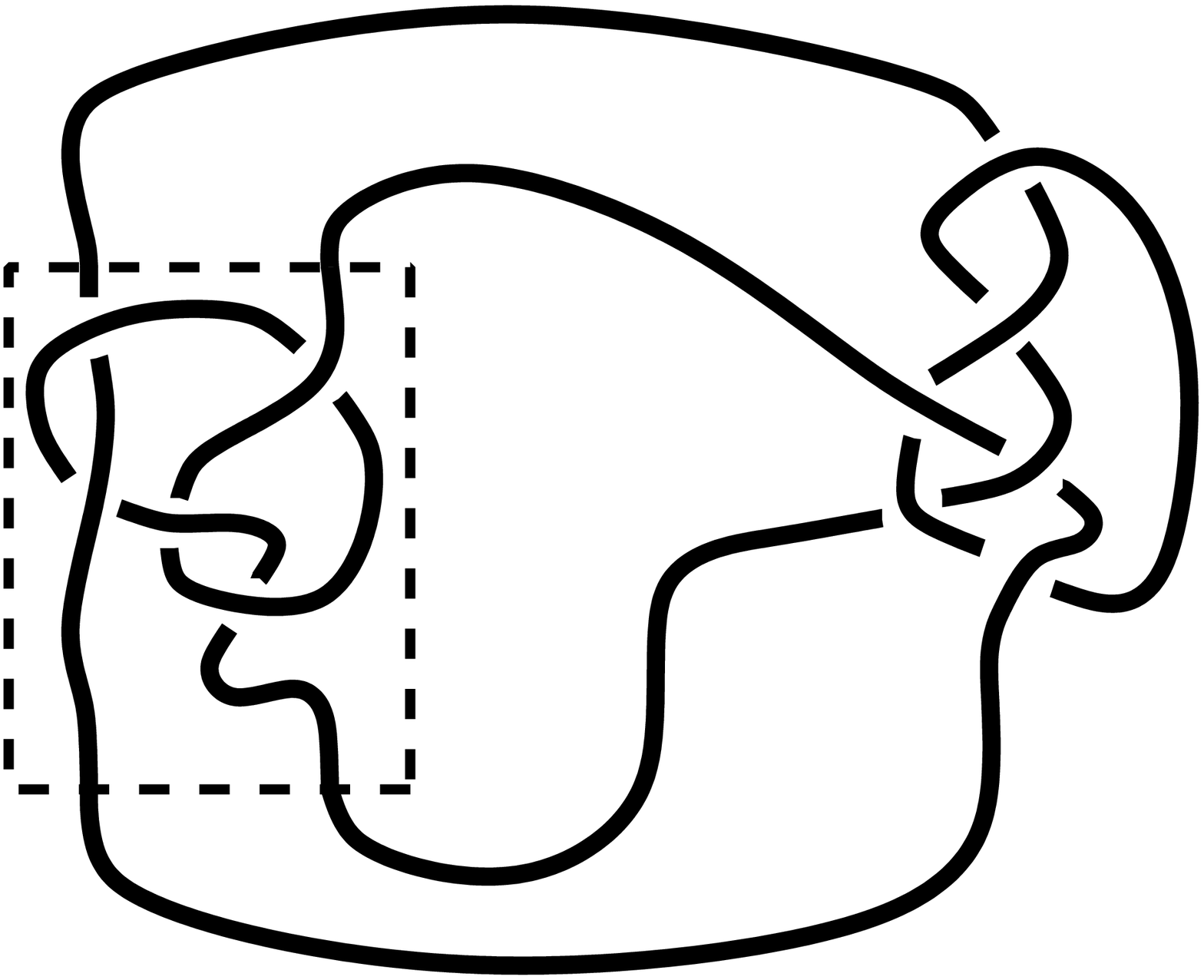}} \qquad
  11^n_{42}=\ol{KT} = \rb{-13mm}{\ig[height=28mm]{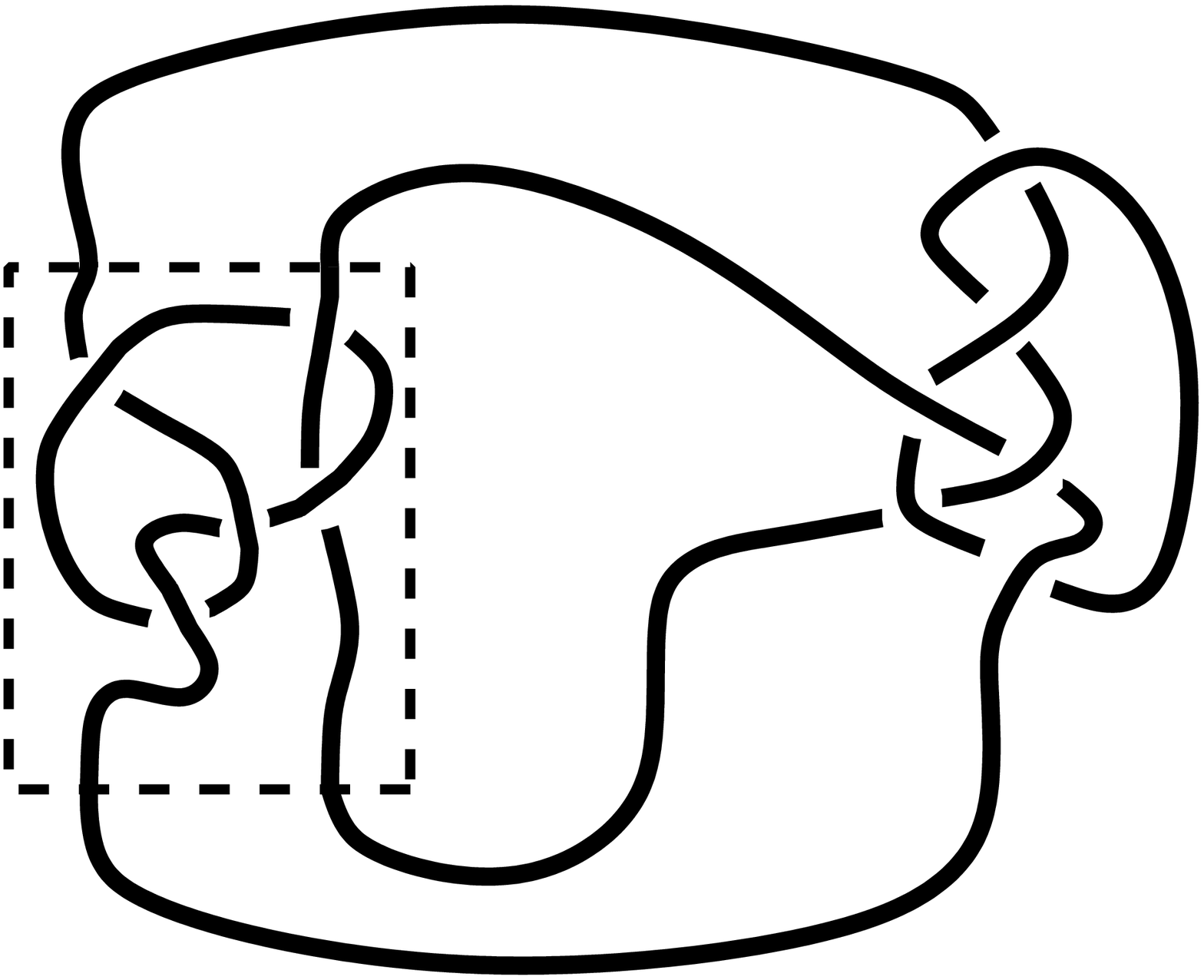}}
$$

\begin{theorem}[\cite{MC}] \label{MCthm}
There exists a Vassiliev invariant $v$ of order 11 such that
$v(C)\not=v(KT)$.
\end{theorem}

Morton and Cromwell manufactured the invariant $v$ using the Lie
algebra $\gl_N$ with a nonstandard representation (or, in other
words, the HOMFLY polynomial of certain cablings of the knots).

J.~Murakami \cite{Mu} showed that any invariant or order at most 10
does not distinguish mutants. So order 11 is the smallest where
Vassiliev invariants detect mutants.

\subsection{Mutation of the Kontsevich integral}
\label{ki_mutation}

Let us describe the behaviour of the Kontsevich integral with
respect to knot mutation.

First, recall the definition of a {\em share} from
Section~\ref{cd-w-ig}:  it is a part of the Wilson loop of a chord
diagram, consisting of two arcs, such that each chord has either
both or no endpoints on it. A mutation of a chord diagram is an
operation of flipping the share with all the chords on it.

In the construction of the Kontsevich integral of a knot $K$ the
Wilson loop of the chord diagrams is parametrized by the same
circle as $K$. For each chord diagram participating in $Z(K)$, the
mutation of $K$ with respect to a subtangle $T$ gives rise to a flip
of two arcs on the Wilson loop.

\begin{xtheorem}[\cite{Le}]\label{th:Le}
Let $M_T(K)$ be the mutant of a knot $K$ with respect to a subtangle
$T$.  Then $Z(K)$ consists only of diagrams for which the part of
the Wilson loop that corresponds to $T$ is a share. Moreover, if
$M_T(Z(K))$ is obtained from $Z(K)$ by flipping the $T$-share of
each diagram, we have
$$
  Z(M_T(K))=M_T(Z(K)).
$$
\end{xtheorem}

\begin{proof} The proof is a straightforward application of the
combinatorial construction of the Kontsevich integral. Write $K$ as
a product $K=A\cdot (T\ot B) \cdot C$ where $A,B,C$ are some
tangles. Then the mutation operation consists in replacing $T$ in
this expression by its flip $T'$ about a vertical axis.

First, observe that rotating a parenthesized tangle with two points
at the top and two points at the bottom by $180^{\circ}$ about a
vertical axis results in the same operation on its combinatorial
Kontsevich integral. Moreover, since there is only one choice of
parentheses for a product of two factors, the non-associative
monomials on the boundary of $T$ are the same as those of $T'$ (all
are equal to $(xx)$). Choose the non-associative monomials for $B$
to be $u$ at the top and $v$ at the bottom. Then
$$Z(K)=Z(A,1,(xx)u)\cdot \Bigl( Z(T,(xx),(xx))
\ot Z(B,u,v)\Bigr) \cdot Z(C, (xx)v, 1),$$ where we write simply $Z$
for $Z_{comb}$, and
$$Z(M_T(K))=Z(A,1,(xx)u)\cdot \Bigl( Z(T',(xx),(xx))
\ot Z(B,u,v)\Bigr) \cdot Z(C, (xx)v, 1).$$ Both expressions only
involve diagrams for which the part of the Wilson loop that
corresponds to $T$ is a share; they differ exactly by the mutation
of all the $T$-shares of the diagrams.
\end{proof}

\subsection{Counterexample to the Intersection Graph Conjecture}
\label{IGCwrong} \index{Intersection graph!conjecture}

It is easy to see that the mutation of chord diagrams does not
change the intersection graph. Thus, if the intersection graph
conjecture (see \ref{IGC}) were true, the Kontsevich integrals of
mutant knots would coincide, and all Vassiliev invariants would take
the same value on mutant knots. But this contradicts Theorem
\ref{MCthm}.

\subsection{}
Now we can prove the theorem announced on page~\pageref{IGTM}:
\begin{xtheorem}[\cite{CL}] The symbol of a Vassiliev
invariant that does not distinguish mutant knots depends on the
intersection graph only.
\end{xtheorem}

\begin{proof}
The idea of the proof can be summarized in one sentence: a mutation
of a chord diagram is always induced by a mutation of a singular
knot.

Let $D_1$ and $D_2$ be chord diagrams of degree $n$ with the same
intersection graph. We must prove that if a Vassiliev knot invariant
$v$, of order at most $n$, does not distinguish mutants, then the
symbol of $v$ takes the same value on $D_1$ and $D_2$.

According to the theorem of Section~\ref{cd-w-ig}
(page~\pageref{cd-mut-theorem}), $D_2$ can be obtained from $D_1$ by
a sequence of mutations. It is sufficient to consider the case when
$D_1$ and $D_2$ differ by a single mutation in a share $S$.

Let $K_1$ be a singular knot with $n$ double points whose chord
diagram is $D_1$. The share $S$ corresponds to two arcs on $K_1$;
the double points on these two arcs correspond to the chords with
endpoints on $S$. Now, shrinking and deforming the two arcs, if
necessary, we can find a ball in $\R^3$ whose intersection with
$K_1$ consists of these two arcs and a finite number of other arcs.
These other arcs can be pushed out of the ball, though not
necessarily by an isotopy, that is, passing through
self-intersections. The result is a new singular knot $K_1'$ with
the same chord diagram $D_1$, whose double points corresponding to
$S$ are collected in a tangle $T_S$. Performing an appropriate
rotation of $T_S$ we obtain a singular knot $K_2$ with the chord
diagram $D_2$. Since $v$ does not distinguish mutants, its values on
$K_1$ and $K_2$ are equal. The theorem is proved.
\end{proof}

\medskip

To illustrate the proof, consider the chord diagram $D_1$ below.
Pick a singular knot $K_1$ representing $D_1$.
$$D_1 =\ \risS{-25}{esche-cd-st-ex}{}{55}{30}{30}\hspace{3cm}
  K_1 =\ \risS{-25}{esche-kn-sin1}{}{80}{15}{8}
$$
By deforming $K_1$ we achieve that the two arcs of the share form a
tangle (placed on its side in the pictures below), and then push all
other segments of the knot out of this subtangle:
$$\risS{-25}{esche-kn-sin2}{
    \put(-5,-10){\mbox{\scriptsize deforming the knot to form the subtangle}}
    }{100}{40}{50}\hspace{2.4cm}
  \risS{-25}{esche-kn-sin3}{
    \put(5,-10){\mbox{\scriptsize pushing out other segments}}
    }{100}{40}{50}
$$

Combining the last theorem with \ref{th:Le} we get the following
corollary.

\begin{xcorollary} Let $w$ be a weight system on chord diagrams with $n$ chords.
Consider a Vassiliev invariant 
$v(K):=w\circ I(K)$. Then
$v$ does not distinguish mutants if and only if the weight system
$w$ depends only on the intersection graph.
\end{xcorollary}

\section{Canonical Vassiliev invariants}\label{canvi}

Theorem~\ref{Ko_part_fund_th}) on the universality of the Kontsevich
integral and its framed version in Section~\ref{frKo_part_fund_th}
provide a means to recover a Vassiliev invariant of order $\leqslant
n$ from its symbol, up to invariants of smaller order. It is natural
to consider those remarkable Vassiliev invariants whose recovery
gives precisely the original invariant.

\begin{definition}\label{def_cani} (\cite{BNG})
\index{Vassiliev!invariant!canonical}
\index{Canonical!invariant}
A (framed) Vassiliev invariant $v$ of order $\leqslant n$ is called
{\it canonical} if for every (framed) knot $K$,
$$v(K) = \symb(v)(I(K))\,.$$
In the case of framed invariants one should write $I^{fr}(K)$
instead of $I(K)$.

A power series invariant $f=\sum_{n=0}^\infty f_n h^n$, with $f_n$
of order $\leq n$ is called {\em canonical} \index{Canonical!series}
if
$$f(K) = \sum_{n=0}^\infty \bigl(w_n(I(K))\bigr) h^n\,$$
for every knot $K$, where $w=\sum_{n=0}^\infty w_n$ is the symbol of
$f$. And, again, in the framed case one should use $I^{fr}(K)$
instead of $I(K)$.
\end{definition}

Recall that the power series invariants were defined on
page~\pageref{pol-Vas-inv} and their symbols
--- in the remark after Proposition~\ref{sym-hom}.

Canonical invariants define a grading in the filtered space of
Vassiliev invariants which is consistent with the filtration.

\begin{xexample} The trivial invariant of order $0$ which is
identically equal to 1 on all knots is a canonical invariant. Its
weight system is equal to $\bo_0$ in the notation of Section
\ref{bialgWS}.
\end{xexample}

\begin{xexample} The Casson invariant $c_2$ is canonical. This follows
from the explicit formula \ref{v2_th} that defines it in terms of
the knot diagram.
\end{xexample}

\noindent{\bf Exercise.} Prove that the invariant $j_3$ (see
\ref{jones_vi}) is canonical.

\medskip

Surprisingly many of the classical knot invariants discussed in Chapters
\ref{kn_inv} and \ref{FT_inv} turn out to be canonical.

The notion of a canonical invariant allows one to reduce various
relations between Vassiliev knot invariants to some combinatorial
relations between their symbols, which gives a powerful tool to
study knot invariants. This approach will be used in
Section~\ref{melmor} to prove the Melvin--Morton conjecture. Now we
shall give examples of canonical invariants following \cite{BNG}.

\subsection{Quantum invariants}\label{can_qi}

Building on the work of Drinfeld \cite{Dr1,Dr2} and Kohno \cite{Koh2},
T.~Le and J.~Murakami \cite[Th 10]{LM3}, and C.~Kassel
\cite[Th XX.8.3]{Kas} \index{Theorem!Le--Murakami--Kassel}
(see also \cite[Th 6.14]{Oht1})
proved that the quantum knot invariants $\theta^{fr}(K)$ and
$\theta(K)$ introduced in Section \ref{qi} become canonical series after
substitution $q=e^h$ and expansion into a power series in $h$.

The initial data for these invariants is a semi-simple Lie algebra $\g$
and its finite dimensional irreducible representation $V_\lambda$,
where $\lambda$ is its highest weight. %!!! explain in the Appendix
To emphasize this data, we shall
write $\theta^{V_\lambda}_{\g}(K)$ for $\theta(K)$ and
$\theta^{fr,V_\lambda}_{\g}(K)$ for $\theta^{fr}(K)$ .
\index{Quantum invariant}

The quadratic Casimir element $c$ (see Section~\ref{LAWS_A}) acts on
$V_\lambda$ as multiplication by a constant, call it $c_\lambda$.
The relation between the framed and unframed quantum invariants is
$$\theta^{fr,V_\lambda}_{\g}(K) = q^{\frac{c_\lambda\cdot w(K)}{2}}
  \theta^{V_\lambda}_{\g}(K)\,,$$
where $w(K)$ is the writhe of $K$.

Set $q=e^h$. Write $\theta^{fr,V_\lambda}_{\g}$ and
$\theta^{V_\lambda}_{\g}$  as power series in $h$:
$$\theta^{fr,V_\lambda}_{\g} =
        \sum_{n=0}^\infty \theta^{fr,\lambda}_{\g,n}h^n \qquad
 \theta^{V_\lambda}_{\g} = \sum_{n=0}^\infty \theta^{\lambda}_{\g,n}h^n.
$$

According to the Birman--Lin theorem (\ref{qift}),
the coefficients $\theta^{fr,\lambda}_{\g,n}$ and
$\theta^{\lambda}_{\g,n}$ are Vassiliev
invariants of order $n$. The Le--Murakami---Kassel Theorem states
that they both are canonical series.

It is important that the symbol of $\theta^{fr,V_\lambda}_{\g}$ is
precisely the weight system $\f_\g^{V_\lambda}$ described in Chapter
\ref{LAWS}. The symbol of $\theta^{V_\lambda}_{\g}$ equals
$\f'^{V_\lambda}_\g$. In other words, it is obtained from
$\f_\g^{V_\lambda}$ by the deframing procedure of
Section~\ref{defram_ws}. Hence,  knowing the Kontsevich integral
allows us to restore the quantum invariants
$\theta^{fr,V_\lambda}_{\g}$ and $\theta^{V_\lambda}_{\g}$ from
these weight systems without the quantum procedure of
Section~\ref{qi}.

\subsection{Coloured Jones polynomial}\label{colJones}
\index{Jones polynomial!coloured}
\index{Coloured Jones polynomial}

    The {\it coloured Jones polynomials}
$J^k := \theta^{V_\lambda}_{\sL_2}$ and
$J^{fr,k}:=\theta^{fr,V_\lambda}_{\sL_2}$ are
particular cases of quantum invariants for $\g=\sL_2$. For this
Lie algebra, the highest weight is an integer $\lambda=k-1$, where $k$
is the dimension of the representation, so in our notation we may
use $k$ instead of $\lambda$. The quadratic Casimir number
in this case
is $c_\lambda=\frac{k^2-1}{2}$, and the relation between the framed
and unframed coloured Jones polynomials is
$$J^{fr,k}(K)= q^{\frac{k^2-1}{4}\cdot w(K)} J^k(K)\,.$$
The ordinary Jones polynomial of Section \ref{Jones} corresponds to
the case $k=2$, that is, to  the standard $2$-dimensional
representation of the Lie algebra $\sL_2$.

Set $q=e^h$. Write $J^{fr,k}$ and $J^k$ as power series in $h$:
$$J^{fr,k} = \sum_{n=0}^\infty J^{fr,k}_{n}h^n \qquad
J^k = \sum_{n=0}^\infty J^k_n h^n.$$
Both series are canonical with the symbols
$$\symb(J^{fr,k}) = \f_{\sL_2}^{V_k}\,,\qquad
  \symb(J^k) = \f'^{V_k}_{\sL_2}
$$
defined in Sections \ref{ws_sl_2_on_A} and \ref{ws_sl_2_on_C}.

\subsection{Alexander--Conway polynomial}\label{conconway}
\index{Conway polynomial}
\index{Alexander--Conway polynomial}

Consider the unframed quantum invariant $\theta^{St}_{\sL_N}$ as a
function of the parameter $N$. Let us think of $N$ not as a discrete
parameter but rather as a continuous variable, where for non integer
$N$ the invariant $\theta^{St}_{\sL_N}$ is defined by the skein and
initial relations above. Its symbol
$\f'^{St}_{\sL_N}=\f'^{St}_{\gl_N}$ (see Exercise~\ref{ex_defr_sl_N}
to Chapter~\ref{LAWS}) also makes sense for all real values of $N$,
because for every chord diagram $D$, $\f'^{St}_{\gl_N}(D)$ is a
polynomial of $N$. Even more, since this polynomial is divisible by
$N$, we may consider the limit
$$\lim_{N\to 0} \frac{\f'^{St}_{\sL_N}}{N}\,.$$

\noindent{\bf Exercise.} Prove that the weight system defined by
this limit coincides with the symbol of the Conway polynomial,
$\symb(\CP)=\sum_{n=0}^\infty \symb(c_n).$
\index{Conway polynomial!symbol}

{\sl Hint.} Use Exercise~\ref{ex_consymb} to Chapter \ref{FT_inv}.

Make the substitution $\theta^{St}_{\sL_N}\bigr|_{q=e^h}$. The skein
and initial relations for $\theta^{St}_{\sL_N}$ allow us to show
(see Exercise~\ref{ex_ex_con_lim} to this chapter) that the
limit\label{Alex-Conv-pol}
$$A := \lim_{N\to 0} \frac{\theta^{St}_{\sL_N}\bigr|_{q=e^h}}{N}$$
does exist and satisfies the relations
\begin{equation}\label{eq_skein_AC}
  A\Bigl(\lrints\Bigr)\ -\  A\Bigl(\rlints\Bigr)\quad = \quad
  (e^{h/2}-e^{-h/2})A\Bigl(\twoup\Bigr)\,;
\end{equation}
\begin{equation}\label{eq_init_AC}
A\Bigl(\unkn\Bigr) = \frac{h}{e^{h/2}-e^{-h/2}}\,.
\end{equation}

A comparison of these relations with the defining relation for the
Conway polynomial \ref{axiconw} shows that
$$A=\frac{h}{e^{h/2}-e^{-h/2}} \CP\bigr|_{t=e^{h/2}-e^{-h/2}}\,.$$
Despite of the fact that the Conway polynomial $\CP$ itself is not a
canonical series, it becomes canonical after the substitution
$t=e^{h/2}-e^{-h/2}$ and multiplication by
$\frac{h}{e^{h/2}-e^{-h/2}}$. The weight system of this canonical
series is the same as for the Conway polynomial. Or, in other words,
$$\frac{h}{e^{h/2}-e^{-h/2}} \CP\bigr|_{t=e^{h/2}-e^{-h/2}}(K) =
  \sum_{n=0}^\infty \bigl(\symb(c_n)\circ I(K)\bigr) h^n\,.$$

\begin{xremark} We cannot do the same for framed invariants because none of
the limits
$$\lim_{N\to 0} \frac{\theta^{fr,St}_{\sL_N}\bigr|_{q=e^h}}{N},\qquad
  \quad \lim_{N\to 0} \frac{\f^{St}_{\sL_N}}{N}$$
exists.
\end{xremark}

\section{Wheeling}\label{wheeling}
\index{Wheeling}

We mentioned in Section~\ref{twoprodB} of Chapter~\ref{algFDchap}
that the relation between the algebras $\B$ and $\F$ is similar to
the relation between the invariants in the symmetric algebra of a
Lie algebra and the centre of its universal enveloping algebra. One
may then expect that there exists an algebra isomorphism between
$\B$ and $\F$ similar to the {Duflo isomorphism} \index{Duflo
isomorphism} for Lie algebras (see page~\pageref{duflo}).

This isomorphism indeed exists. It is called {\em wheeling} and we
describe it in this section. It will be used in the next section to
deduce an explicit formula for the Kontsevich integral of the
unknot.

\subsection{Diagrammatic differential operators and the wheeling map}
\label{subsec:ddo}

For an open diagram $C$ with $n$ legs, let us define the {\em
diagrammatic differential operator} \index{Diagrammatic differential
operator!on $\B$}\label{d-d-o}
$$\partial_{C}:\B\to\B.$$
Take an open diagram $D$. If $D$ has at most $n$ legs, set
$\partial_{C}D=0$. If $D$ has more than $n$ legs, we define
$\partial_{C}(D)\in \B$ as the sum of all those ways of glueing all
the legs of $C$ to some legs of $D$ that produce diagrams having at
least one leg on each connected component. For example, if $w_2$
stands for the diagram $\rccO{-1.2}{ccWb}$, we have
$$\partial_{w_2}(\rccO{-4}{ccw4}) =
  \ 8 \ \rb{-4mm}{\ig[height=10mm]{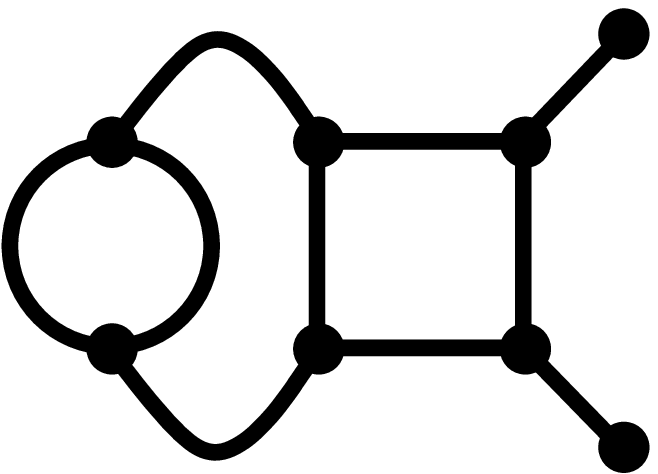}}
  \ +\ 4 \ \rb{-4mm}{\ig[height=10mm]{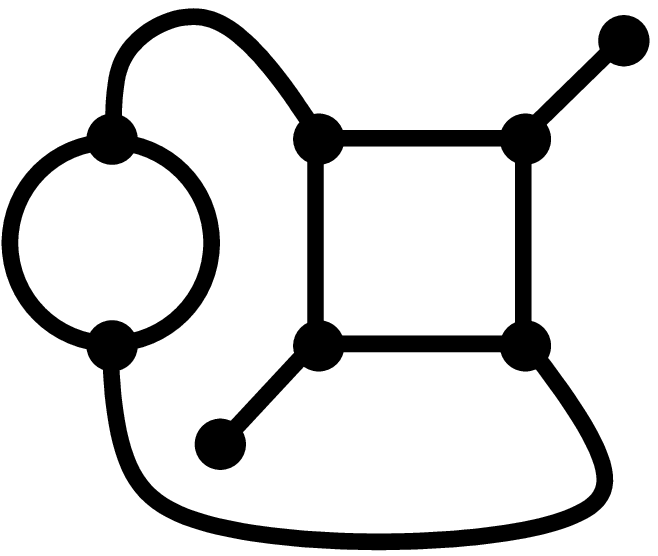}}
  \ =\ 10\ \rb{-2.5mm}{\ig[height=6mm]{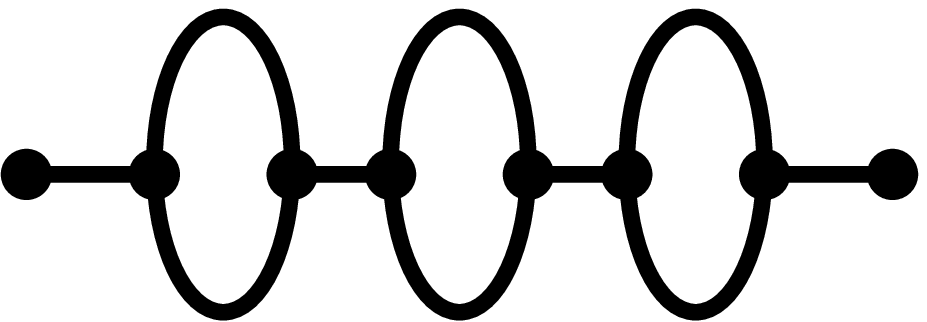}}\ .
$$
Also,
$$
\partial_{w_2}(\ \rccO{-1}{ccWa}\ )\ =\ 8\ \rccO{-1.2}{ccWb}\ ,
$$
since the other four ways of glueing $w_2$ into $\ \rccO{-1}{ccWa}\
$ produce diagrams one of whose components has no legs (see page
\pageref{eg-do}).

Extending the definition by linearity, we can replace the diagram
$C$ in the definition of $\partial_{C}$  by any linear combination
of diagrams. Moreover, $C$ can be taken to be a formal power series
in diagrams, {\em with respect to the grading by the number of
legs}. Indeed, for any given diagram $D$ almost all terms in such
formal power series would have at least as many legs as $D$.

Recall that the {\em wheel} \index{Wheel!in $\B$} $w_n$ in the
algebra $\B$ is the diagram
$$w_n = \quad \risS{-13}{wssl2wn-b}{
         \put(0,-8){\mbox{\scriptsize $n$ spokes}}}{30}{20}{25}
$$
The wheels $w_n$ with $n$ odd are equal to zero; this follows
directly from Lemma~\ref{anti-auto}.

\begin{xdefinition}
The {\em wheels element} $\Omega$ is the formal power series
$$
\Omega=\exp{\sum_{n=1}^{\infty} b_{2n} w_{2n}}
$$
where $b_{2n}$ are the {\em modified Bernoulli numbers}, and the
products are understood to be in the algebra $\B$.
\end{xdefinition}

The modified Bernoulli numbers $b_{2n}$ are the coefficients at
$x^{2n}$ in the Taylor expansion of the function
$$
\frac{1}{2}\ln{\frac{\sinh{x/2}}{x/2}}.
$$
We have $b_2=1/48$, $b_4=-1/5760$ and $b_6=1/362880$. In general,
$$
b_{2n}=\frac{B_{2n}}{4n\cdot(2n)!},
$$
where $B_{2n}$ are the usual Bernoulli numbers. 

\begin{xdefinition} The {\em wheeling map} is the map
$$\label{wheel-map}\index{Wheeling map}
\partial_{\Omega}=\exp{\sum_{n=1}^{\infty} b_{2n} \partial_{w_{2n}}}
$$
\end{xdefinition}
The wheeling map is a degree-preserving linear map $\B\to\B$. It is,
clearly, a vector space isomorphism since $\partial_{\Omega^{-1}}$
is an inverse for it.

\begin{theorem}[Wheeling Theorem]\label{thm:wheeling}
\index{Theorem!wheeling}\index{Wheeling Theorem} The map
$\chi\circ\partial_{\Omega}:\B\to\F$ is an algebra isomorphism.
\end{theorem}

There are several approaches to the proof of the above theorem. It
has been noted by Kontsevich \cite{Kon3} that the Duflo-Kirillov
isomorphism holds for a Lie algebra in any rigid tensor category;
Hinich and Vaintrob showed in \cite{HV} that the wheeling map can be
interpreted as a particular case of such a situation. Here, we shall
follow the proof of Bar-Natan, Le and Thurston \cite{BLT}.

\begin{xexample}
At the beginning of Section~\ref{twoprodB} (page~\pageref{twoprodB})
we saw that $\chi$ is not compatible with the multiplication. Let us
check the multiplicativity of $\chi\circ\partial_{\Omega}$ on the
same example:
$$\begin{array}{ccl}
  \chi\circ\partial_{\Omega}(\ \rccO{-1}{ccWa}\ ) &=&
  \chi\circ(1+b_2\partial_{w_2})(\ \rccO{-1}{ccWa}\ ) \vspace{5pt}\\
&=& \chi\bigl(\rccO{-1}{ccWa} + \frac{1}{48}\cdot 8\cdot
\rccO{-1.2}{ccWb}\bigr) \ =\
    \chi\bigl(\rccO{-1}{ccWa} + \frac{1}{6} \rccO{-1.2}{ccWb}\bigr) \vspace{5pt} \\
&=&\frac{1}{3}\ \rccO{-3.8}{fd2a}\ +\ \frac{2}{3}\
\rccO{-3.8}{fd2b}\
        +\ \frac{1}{6}\ \rccO{-3.8}{fdWb} \ =\ \rccO{-3.8}{fd2b}\ ,
\end{array}
$$
which is the square of the element
$\chi\circ\partial_{\Omega}(\ccO{ccOa}) = \sClD{cd1}$ in the algebra
$\F$.
\end{xexample}

\subsection{The algebra $\Bbig$}\label{subsection:beprime}
\index{Algebra!$\Bbig$}

\def\Fbig{\F^{\circ}}

For what follows it will be convenient to enlarge the algebras $\B$
and $\F$ by allowing diagrams with components that have no legs.

A diagram in the enlarged algebra $\Bbig$ \label{B-big} is a union
of a unitrivalent graph with a finite number of circles with no
vertices on them; a cyclic order of half-edges at every trivalent
vertex is given. The algebra  $\Bbig$ is spanned by all such
diagrams modulo IHX and antisymmetry relations. The multiplication
in $\Bbig$ is the disjoint union. The algebra $\B$ is the subalgebra
of $\Bbig$ spanned by graphs which have at least one univalent
vertex in each connected component. Killing all diagrams which have
components with no legs, we get a homomorphism $\Bbig\to\B$, which
restricts to the identity map on $\B\subset\Bbig$.

The algebra of 3-graphs $\G$ from Chapter \ref{alg3g} is also a
subspace of $\Bbig$. In fact, the algebra $\Bbig$ is the tensor
product of $\B$, the symmetric algebra $S(\G)$ of the vector space
$\G$ and the polynomial algebra in one variable (which counts the
circles with no vertices on them).

The reason to consider $\Bbig$ instead of $\B$ can be roughly
explained as follows. One of our main tools for the study of $\B$ is
the universal weight system
$$\rho_{\g}:\B\to S(\g)$$
with the values in the symmetric algebra of a Lie algebra $\g$. In
fact, much of our intuition about $\B$ comes from Lie algebras,
since $\rho_{\g}$ respects some basic constructions. For instance,
glueing together the legs of two diagrams  corresponds to a
contraction of the corresponding tensors. However, there is a very
simple operation in $S(\g)$ that cannot be lifted to $\B$ via
$\rho_{\g}$. Namely, there is a pairing $S^n(\g)\ot S^n(\g)\to \C$
which extends the invariant form $\g\ot\g\to\C$. Roughly, if the
elements of $S(\g)$ are thought of as symmetric tensors (see
page~\pageref{symtens}), this pairing consists in taking the sum of
all possible contractions of two tensors of the same rank. This
operation is essential if one works with differential operators on
$S(\g)$, and it cannot be lifted to $\B$ since glueing two diagrams
with the same number of legs produces a diagram with no univalent
vertices. The introduction of $\Bbig$ remedies this problem.

Indeed, the map $\rho_{g}$ naturally extends to $\Bbig$. On a
connected diagram with no legs it coincides with the $\C$-valued
weight system $\f_{\g}$ for 3-graphs described in
Section~\ref{LAWS_3G}; in particular, it takes value $\dim{\g}$ on a
circle with no vertices. Finally, $\rho_g$ is multiplicative with
respect to the disjoint union of diagrams.

There is a bilinear symmetric pairing \index{Pairing
$\Bbig\otimes\Bbig\to\G$} on $\Bbig$ whose image lies in the
subspace spanned by legless diagrams.
\begin{xdefinition}
For two diagrams $C, D\in \Bbig$ with the same number of legs we
define $\langle C, D\rangle$ to be the sum of all ways of glueing
all legs of $C$ to those of $D$. If the numbers of legs of $C$ and
$D$ do not coincide we set $\langle C, D\rangle =0$.
\end{xdefinition}
Now, if $C$ and $D$ are two diagrams with the same number of legs,
$\rho_{\g}(\langle C, D\rangle)$ is the sum of all possible
contractions of $\rho_{\g}(C)$ and $\rho_{\g}(D)$ considered as
symmetric tensors.

\begin{xdefinition}
Let $C$ be an open diagram. The {\em diagrammatic differential
operator} \index{Diagrammatic differential operator!on
$\Bbig$}\label{d-d-o-O}
$$\partbig_C:\Bbig\to\Bbig$$
sends $D\in\Bbig$ to  the sum of all ways of glueing the legs of $C$
to those of $D$, if $D$ has at least as many legs as $C$; if $D$ has
less legs than $C$, then $\partbig_C(D)=0$.
\end{xdefinition}
For example,
$$\label{eg-do}
\partbig_{w_2}(\ \risS{0}{strut-2}{}{15}{10}{10}\ )\ =\ 8\ \rccO{-1.2}{ccWb}\ +\ 4\
        \risS{-5}{strut-theta}{}{15}{10}{10}\ .
$$
This definition of diagrammatic operators is consistent with the
definition of diagrammatic operators in $\B$. Namely, if $C$ is a
diagram in $\B$ and $p:\Bbig\to\B$ is the projection, we have
$$\partial_C\circ p = p\circ\partbig_C.$$
Note that while $\partial_C$ and $\partbig_C$ are compatible with
the projection $p$, they are not compatible with the inclusion
$\B\to\Bbig$.

\medskip

Similarly to the algebra $\Bbig$ one defines the algebra $\Fbig$ by
considering not necessarily connected trivalent graphs in the
definition of $\F$. The vector space isomorphism $\chi:\F\to\B$
extends to an isomorphism $\Fbig\simeq\Bbig$ whose definition
literally coincides with that of $\chi$ (and which we also denote by
$\chi$). In particular, for a legless diagram in $\Fbig$ the map
$\chi$ consists in simply erasing the Wilson loop.

Our method of proving the Wheeling Theorem will be to prove it for
the algebras $\Bbig$ and $\Fbig$, with the diagrammatic operator
$\partbig_{\Omega}:\Bbig\to\Bbig$. Then the version for $\B$ and
$\F$ will follow immediately by applying the projection map. First,
however, let us explain the connection of the Wheeling Theorem with
the Duflo isomorphism for Lie algebras.

\subsection{The Duflo isomorphism}\label{duflo}
The wheeling map is a diagrammatic analogue of the
\index{Duflo-Kirillov map}\index{Map!Duflo-Kirillov} {\em
Duflo-Kirillov map} for metrized Lie algebras.

Recall that for a Lie algebra $\g$ the Poincar\'e-Birkhoff-Witt
isomorphism
$$S(\g)\simeq U(\g)$$
sends a commutative monomial in $n$ variables to the average of all
possible {\em non-commutative} monomials in the same variables, see
page~\pageref{PBW}. It is not an algebra isomorphism, of course,
since $S(\g)$ is commutative and $U(\g)$ is not (unless $\g$ is
abelian); however, it is an isomorphism of $\g$-modules. In
particular, we have an isomorphism of vector spaces
$$S(\g)^{\g}\simeq U(\g)^{\g}=Z(U(\g))$$
between the subalgebra of invariants in the symmetric algebra and
the centre of the universal enveloping algebra. This map does not
respect the product either, but it turns out that $S(\g)^{\g}$ and
$Z(U(\g))$ are actually isomorphic as commutative algebras. The
isomorphism between them, known as the {\em Duflo isomorphism},
\index{Duflo isomorphism} is given by the {\em Duflo-Kirillov map},
which is described in the Appendix, see page~\pageref{duflok}.

\begin{xlemma} The wheeling map $\partbig_{\Omega}: \Bbig\to\Bbig$ is
taken by the universal Lie algebra weight system $\rho_{\g}$ to the
Duflo-Kirillov map:
$${\sqrt{j}}\circ\rho_{\g}=\rho_{\g}\circ\partbig_{\Omega}.$$
\end{xlemma}
\begin{proof} Observe that a diagrammatic operator $\partbig_C:\Bbig\to\Bbig$
is taken by $\rho_{\g}$ to the corresponding differential operator
$\partbig_{\rho_{\g}(C)}:S(\g)\to S(\g)$ in the sense that
$$\rho_{\g}\circ\partbig_C=\partial_{\rho_{\g}(C)}\circ\rho_{\g}.$$
This simply reflects the fact that glueing the legs of two diagrams
corresponds to a contraction of the corresponding tensors.

On page~\pageref{rhootkolesa} we have calculated the value of
$\rho_{\g}$ on the wheel $w_{k}$:
$$\rho_{\g}(w_{k})=
\sum_{i_1, \ldots , i_k}\Tr{(\ad{\,e_{i_1}}\ldots
\ad{\,e_{i_k}})}\cdot e_{i_1}\ldots e_{i_k},$$ where $\{e_i\}$ is a
basis for $\g$. In order to interpret this expression as an element
of $S^k(\g^*)$, we must contract it with $k$ copies of $x\in \g$.
The resulting homogeneous polynomial of degree $k$ on $\g$ sends
$x\in\g$ to $\Tr{(\ad{\, x})^{k}}$.

Since $\rho_{\g}$ is multiplicative, it carries the wheeling map
$\partbig_{\Omega}$ to
\begin{multline*}
\exp \left(\sum_n b_{2n} \Tr{(\ad{\, x})^{2n}}\right)= \exp
\Tr{\left(\sum_n b_{2n}(\ad{\, x})^{2n}\right)}\\ = \det \exp
\left(\frac{1}{2}\ln
\frac{\sinh{\frac{1}{2}{\ad{\,x}}}}{{\frac{1}{2}\ad{\,x}}}\right)=\sqrt{j},
\end{multline*}
that is, to the Duflo-Kirillov map.
\end{proof}

The Duflo isomorphism is a rather mysterious fact. Remarkably, more
than one of its proofs involve diagrammatic techniques: apart from
being a consequence of the Wheeling Theorem, it follows from
Kontsevich's work on deformation quantization \cite{Kon3}. In fact,
in \cite{Kon3} the Duflo isomorphism is generalized to a sequence of
isomorphisms between $H^i(\g, S(\g))$ and $H^i(\g, U(\g))$ with the
usual Duflo isomorphism being the case $i=0$.

\subsection{Pairings on diagram spaces and cabling
operations}\label{pairandcable} Everything we said about the
algebras $\Bbig$ and $\Fbig$ can be generalized for the case of
tangles with several components. In particular, there is a bilinear
pairing
$$\Fbig(\xx\,|\,\yy)\otimes\Bbig(\yy)\to\Fbig(\xx).$$
For diagrams $C\in\Fbig(\xx\,|\,\yy)$ and $D\in\Bbig(\yy)$ define
the diagram $\langle C, D\rangle_{\yy}\in \Fbig(\xx)$ as the sum of
all ways of glueing all the $\yy$-legs of $C$ to the $\yy$-legs of
$D$. If the numbers of $\yy$-legs of $C$ and $D$ are not equal, we
set $\langle C, D\rangle_{\yy}$ to be zero. This is a version of the
inner product for diagram spaces defined in
Section~\ref{subsec:pairing}; Lemma~\ref{pair-C-B} shows that it is
well-defined. In what follows we shall indicate by a subscript the
component along which the inner product is taken.

The inner product can be used to express the diagrammatic
differential operators in $\Bbig$ via the disconnected cabling
operations. These are defined for $\Bbig$ in the same way as for
$\B$; for instance, $\psi^{2\cdot 1}(D)$ is the sum of all diagrams
obtained from $D$ by replacing the label (say, $\yy$) on its
univalent vertices by one of the two labels $\yy_1$ or $\yy_2$.  If
$D\in \B(\yy)$, the labels $\yy_1,\yy_2$ are obtained by doubling
$\yy$ and the diagram $C$ is considered as an element of
$\B(\yy_1)$, we have
$$\partbig_{C}(D)=\left\langle C, \psi^{2\cdot 1}(D)\right\rangle_{\yy_1}.$$
The proof consists in simply comparing the diagrams on both sides.

The cabling operation $\psi^{2\cdot 1}$ can be thought of as a
coproduct in $\Bbig$, dual to the disjoint union with respect to the
inner product:
\begin{equation}\label{cableprodduality}
\left\langle C, D_1\cup D_2\right\rangle= \left\langle \psi^{2\cdot
1}(C),D_1\ot D_2\right\rangle_{\yy_1,\yy_2},
\end{equation}
where on the right-hand side $\yy_1$ and $\yy_2$ are the two labels
for the legs of $\psi^{2\cdot 1}(C)$, and $D_1\ot D_2$ belongs to
$\Bbig(\yy_1)\ot\Bbig(\yy_2)$. The proof of this last formula is
also by inspection of both sides.

\subsection{The Hopf link and the map $\Phi_0$}
\label{subsec:KI-and-Phi}

In what follows we shall often write $\#$ for the connected sum and
$\cup$ --- for the disjoint union product, in order to avoid
confusion.

Consider the framed Hopf link $\hl$ \label{Hopf-hl} with one
interval component labelled $\xx$, one closed component labelled
$\yy$, zero framing, and orientations as indicated: \index{Hopf
link!$\hl$}
$$\quad \risS{-13}{temp3}{
             \put(42,57){\mbox{$\yy$}}
             \put(28,2){\mbox{$\xx$}}}{50}{60}{15}$$
The framed Kontsevich integral $I^{fr}(\hl)$ lives in $\F(\xx,\yy)$
or, via the isomorphism
$$\chi^{-1}_{\yy}:\F(\xx,\yy)\to \Fxy,$$
in $\Fxy$.

Let us write $\Zed(\hl)$ \label{Zed-of-hl} for
the image of $I^{fr}(\hl)$ in $\Fxy$. For any diagram
$D\in\Bbig(\yy)$, the pairing $\langle \Zed(\hl), D \rangle_{\yy}$
is well-defined and lives in $\Fbig(\xx)$. Identifying
$\Bbig({\yy})$ with $\Bbig$ and $\Fbig(\xx)$ with $\Fbig$, we obtain
a map
$$\Phi:\Bbig\to\Fbig$$ \index{Map!$\Phi:\B\to\F$}\label{Phi:B-to-F}
defined by
$$D\to \langle \Zed(\hl), D \rangle_{\yy}.$$

\begin{xlemma}
The map $\Phi:\Bbig\to\Fbig$ is a homomorphism of algebras.
\end{xlemma}

\begin{proof}
Taking the disconnected cabling of the Hopf link $\hl$ along the
component ${\yy}$, we obtain a link $\dhl$\label{Hopf-dhl} with one
interval component labelled $\xx$ and two closed parallel components
labelled ${\yy}_1$ and ${\yy}_2$:
$$\quad \risS{-13}{temp4}{
             \put(49,57){\mbox{$\yy_1$}}
             \put(49,26){\mbox{$\yy_2$}}
             \put(28,2){\mbox{$\xx$}}}{50}{65}{15}$$
In the same spirit as $\Phi$, we define the map
$$\Phi_2:\Bbig\otimes\Bbig\to\Fbig$$
using \index{Map!$\Phi_2:\B\otimes\B\to\F$}\label{Phi-2}  the link
$\dhl$ instead of $\hl$. Namely, given two diagrams,
$D_1\in\Bbig({\yy}_1)$ and $D_2\in\Bbig({\yy}_2)$ we have
$$D_1\otimes D_2\in \Bbig({\yy}_1,{\yy}_2).$$
Write $\Zed(\dhl)$ \label{Zed-of-dhl} for the image of the
Kontsevich integral $I^{fr}(\dhl)$ under the map
$$\chi_{{\yy}_1,{\yy}_2}^{-1}: \F(\xx,{\yy}_1,{\yy}_2)\to \F(\xx\,|\,{\yy}_1,{\yy}_2).$$
Identify $\Bbig({\yy}_1,{\yy}_2)$ with $\Bbig\otimes\Bbig$ and
$\Fbig(\xx)$ with $\Fbig$, and define $\Phi_2$ as
$$D_1\otimes D_2\stackrel{\Phi_2}{\longrightarrow} \langle \Zed(\dhl), D_1\otimes D_2\rangle_{\yy_1,\yy_2}.$$
The map $\Phi_2$ glues the legs of the diagram $D_1$ to the
${\yy}_1$-legs of $\Zed(\dhl)$, and the legs of $D_2$
--- to the $\yy_2$-legs of $\Zed(\dhl)$.

There are two ways of expressing $\Phi_2(D_1\otimes D_2)$ in terms
of $\Phi(D_i)$. First, we can use the fact that $\dhl$ is a product
(as tangles) of two copies of the Hopf link $\hl$. Since the legs
of $D_1$ and $D_2$ are glued independently to the legs corresponding
to $\yy_1$ and $\yy_2$, it follows that
$$
\Phi_2(D_1\otimes D_2)=\Phi(D_1)\# \Phi(D_2).
$$

On the other hand, we can apply the formula (\ref{cableprodduality})
that relates the disjoint union multiplication with disconnected
cabling. We have
\begin{align*}
\Phi_2(D_1\otimes D_2)&=\langle \psi^{2\cdot 1}_{\yy} (\Zed(\hl)),
D_1\otimes D_2\rangle_{\yy_1,\yy_2}\\
&=\langle  \Zed(\hl),
D_1\cup D_2\rangle_{\yy}\\
&=\Phi(D_1\cup D_2),
\end{align*}
and, therefore,
$$ \Phi(D_1)\# \Phi(D_2)=\Phi(D_1\cup D_2).$$
\end{proof}

Given a diagram $D\in \Fxy$, the map $\Bbig\to\Fbig$ given by
sending $C\in\Bbig(\yy)$ to $\langle D, C\rangle_{\yy} \in
\Fbig(\xx)$ shifts the degree of $C$ by the amount equal to the
degree of $D$ minus the number of $\yy$-legs of $D$. If $D$ appears
in $\Zed(\hl)$ with a non-zero coefficient, this difference is
non-negative. Indeed, the diagrams participating in $\Zed(\hl)$
contain no struts (interval components) both of whose ends are
labelled with $\yy$, since the $\yy$ component of $\hl$ comes with
zero framing (see Exercise~\ref{strutsinifr} on
page~\pageref{strutsinifr}). Also, if two $\yy$-legs are attached to
the same internal vertex, the diagram is zero, because of the
antisymmetry relation, and therefore, the number of inner vertices
of $D$ is at least as big as the number of $\yy$-legs.

It follows that the Kontsevich integral $\Zed(\hl)$ can be written
as $\Zed_0(\hl)+\Zed_1(\hl)+\ldots$, where $\Zed_i(\hl)$
\label{Zed-i-of-hl} is the
part consisting of diagrams whose degree exceeds the number of $\yy$
legs by $i$. We shall be interested in the term $\Zed_0(\hl)$ of
this sum.

Each diagram that appears in this term is a union of a {\em comb}
\index{Comb with $n$ teeth} with some wheels:
$$\includegraphics{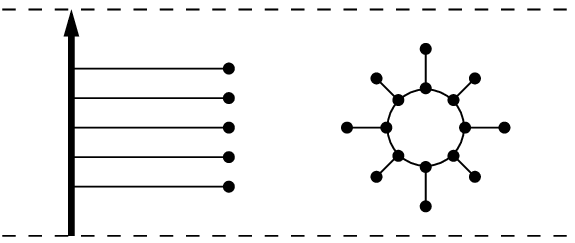}$$

Indeed, each vertex of such diagram is either a $\yy$-leg, or is
adjacent to exactly one $\yy$-leg.

Denote a comb with $n$ teeth by $u^n$. Strictly speaking, $u^n$ is
not really a product of $n$ copies of $u$ since $\Fxy$ is not an
algebra. However, we can introduce a Hopf algebra structure in the
space of all diagrams in $\Fxy$ that consist of combs and wheels.
The product of two diagrams is the disjoint union of all components
followed by the concatenation of the combs; in particular
$u^ku^{m}=u^{k+m}$. The coproduct is the same as in $\Fxy$. This
Hopf algebra is nothing else but the free commutative Hopf algebra
on a countable number of generators.

The Kontsevich integral is group-like, and this implies that
$$\delta(\Zed_0(\hl))=\Zed_0(\hl)\ot \Zed_0(\hl).$$
Group-like elements in the completion of the free commutative Hopf
algebra are the exponentials of linear combinations of generators
and, therefore
$$\Zed_0(\hl)=\exp(cu \cup \sum_n a_{2n}w_{2n}),$$
where $c$ and $a_{2n}$ are some constants.

In fact, the constant $c$ is precisely the linking number of the
components $\xx$ and $\yy$, and, hence is equal to $1$. We can write
$$\Zed_0(\hl)=\sum_n\frac{u^n}{n!}\cup \Omega',$$
where $\Omega'$ the part of $\Zed_0(\hl)$ containing wheels:
$$\Omega'=\exp{\sum_n a_{2n}w_{2n}}.
  \label{omega-pr-def}$$

Define the map $\Phi_0:\Bbig\to\Fbig$
\index{Map!$\Phi_0:\B\to\F$}\label{Phi-0} by taking the pairing of a
diagram in $\Bbig(\yy)$ with $\Zed_0(\hl)$:
$$D\to \langle \Zed_0(\hl), D \rangle_{\yy}.$$
The map $\Phi_0$ can be thought of as the part of $\Phi$ that shifts
the degree by zero. Since $\Phi$ is multiplicative $\Phi_0$ also is.
In fact, we shall see later that  $\Phi_0=\Phi$.

\begin{lemma}\label{lem:partbigomega}
$\Phi_0=\chi\circ\partbig_{\Omega'}.$
\end{lemma}

\begin{proof}
Let us notice first that if $C\in \Fxy$ and $D\in\Bbig(\yy)$, we
have
$$\langle C\cup w_{2n}, D \rangle_{\yy} = \langle C, \partbig_{w_{2n}}(D)
\rangle_{\yy}.$$ Also, for any $D\in \Bbig$ the expression
$$\left\langle \sum_n \frac{u^n}{n!}, D\right\rangle_{\yy}$$
is precisely the average of all possible ways of attaching the legs
of $D$ to the line $\xx$.

Therefore, for $D\in\Bbig(\yy)$

$$\Phi_0(D)=\left\langle\, \sum_n\frac{u^n}{n!}\cup\Omega',\, D\, \right\rangle_{\yy}=
\left\langle\, \sum_n\frac{u^n}{n!},\, \partbig_{\Omega'} D\,
\right\rangle_{\yy}=\chi\circ\partbig_{\Omega'}D.$$
\end{proof}

\subsection{The coefficients of the wheels in $\Phi_0$}
If $D\in\B$ is a diagram, we shall denote by $D_{\yy}$ the result of
decorating all the legs of $D$ with the label $\yy$; the same
notation will be used for linear combinations of diagrams.

First, let us observe that $\Omega'$ is group-like with respect to
the coproduct $\psi^{2\cdot 1}$:

\begin{xlemma}\label{lemma:pseudo-lin}
$$\psi^{2\cdot 1}_{\yy} \Omega'_{\yy}=\Omega'_{\yy_1}\otimes\Omega'_{\yy_2},$$
where $\yy_1$ and $\yy_2$ are obtained by doubling $\yy$.
\end{xlemma}

\begin{proof}
We again use the fact that the sum of the Hopf link $\hl$ with
itself coincides with its two-fold disconnected cabling along the
closed component $\yy$. Since the Kontsevich integral is
multiplicative, we see that
$$\psi^{2\cdot 1}_{\yy} \Zed_0(\hl)_{\xx,\yy}=\Zed_0(\hl)_{\xx,\yy_1}\Zed_0(\hl)_{\xx,\yy_2},$$
where the subscripts indicate the labels of the components and
product on the right-hand side lives in the graded completion of
$\F(\xx\, |\, \yy_1, \yy_2)$. Now, if we factor out on both sides
the diagrams that have at least one vertex on the $\xx$ component,
we obtain the statement of the lemma.
\end{proof}

\begin{xlemma}
For any $D\in\Bbig$
$$\partbig_D(\Omega')
=\langle D, \Omega' \rangle\, \Omega'.$$
\end{xlemma}

\begin{proof}
It is clear from the definitions and from the preceding lemma that
$$\partbig_D(\Omega')=
\langle D_{{\yy}_1}, \psi^{2\cdot 1}_{\yy} \Omega' \rangle= \langle
D_{\yy_1}, \Omega'_{\yy_1}\Omega'_{\yy_2} \rangle= \langle D,
\Omega' \rangle\,\Omega'.$$
\end{proof}

\begin{xlemma}
The following holds in $\Bbig$:
$$\langle\, \Omega', (\ccOy{ccOa})^n\ \rangle=\left(\frac{1}{24} \sClD{cd1}\right)^n.$$
\end{xlemma}

\begin{proof}
According to Exercise~\ref{ex:zed-hl-in-B} on page
\pageref{ex:zed-hl-in-B}, the Kontsevich integral of the Hopf link
$\hl$ up to degree two is equal to
$$\Zed_0(\hl) =  \Iu{I}\hspace{-5pt} + \Iu{I-u}+\frac12\ \Iu{I-uu}
+\frac1{48}\ \Id{I-o}\ .$$
It follows that the coefficient $a_2$ in
$\Omega'$ is equal to $1/48$ and that
$$\langle\, \Omega', \ccOy{ccOa}\, \rangle=\frac{1}{24}\ \sClD{cd1}.$$
This establishes the lemma for $n=1$.  Now, use induction:
\begin{align*}
\langle\, \Omega', (\ccOy{ccOa})^n\ \rangle&=
             \langle\,\partbig_{\risS{4}{strut}{}{15}{10}{10}}\Omega',
                (\ccOy{ccOa})^{n-1}\ \rangle\\
&=\frac{1}{24}\ \sClD{cd1}\cdot \langle\, \Omega', (\ccOy{ccOa})^{n-1} \rangle\\
&=\left(\frac{1}{24}\ \sClD{cd1}\right)^n.
\end{align*}

The first equality follows from the obvious identity
valid for arbitrary
$A,B,C\in\Bbig$:
$$\langle C, A\cup B \rangle = \langle \partbig_B(C), A \rangle.$$
The second equality follows from the preceding lemma, and the third
is the induction step.
\end{proof}

In order to establish that $\Omega'=\Omega$ we have to show that the
coefficients $a_{2n}$ in the expression $\Omega'=\exp{\sum_n a_{2n}w_{2n}}$ 
are equal to the modified Bernoulli numbers $b_{2n}$. In other
words, we have to prove that
\begin{equation}\label{bernoulli-eq}
\sum_n a_{2n}x^{2n} = \frac{1}{2}\ln{\frac{\sinh{x/2}}{x/2}},\mbox{\quad or\quad}
\exp\Bigl(2\sum_n a_{2n}x^{2n}\Bigr) = \frac{\sinh{x/2}}{x/2}.
\end{equation}

To do this we compute the value of the $\sL_2$-weight system
$\f_{\sL_2}$ from Section~\ref{LAWS_3G} (page~\pageref{LAWS_3G}) on
the 3-graph $\langle \Omega',
(\risS{1.5}{strut}{}{15}{10}{10})^n\rangle\in\G$ in two different
ways.

Using the last lemma and Theorem \ref{ws_sl_2_on_C} on page
\pageref{ws_sl_2_on_C} we have
$$\f_{\sL_2}\Bigl(\langle \Omega', (\risS{1.5}{strut}{}{15}{10}{10})^n\rangle\Bigr)
= \f_{\sL_2}\Biggl(\Bigl(\frac{1}{24} \sClD{cd1}\Bigr)^n\Biggr) =
\frac{1}{2^n}\ .
$$

On the other hand, according to Exercise~\ref{ex_sl_2_on_B_wheel} on
page~\pageref{ex_sl_2_on_B_wheel}, the value of the $\sL_2$ weight
system on the wheel $w_{2n}$ is equal to $2^{n+1}$ times its value
on $\ (\risS{1.5}{strut}{}{15}{10}{10})^n$. Therefore,
$$\f_{\sL_2}\Bigl(\langle\, \Omega', (\risS{1.5}{strut}{}{15}{10}{10})^n\rangle\Bigr)
= \f_{\sL_2}\Bigl(\langle\,\exp{\sum_m a_{2m} 2^{m+1}\
(\risS{1.5}{strut}{}{15}{10}{10})^m},
  (\risS{1.5}{strut}{}{15}{10}{10})^n\rangle\Bigr)\ .
$$
Denote by $f_n$ the coefficient at $z^n$ of the power series expansion of the function
$\exp\bigl(2\sum_n a_{2n}z^n\bigr) = \sum_n f_n z^n$. We get
$$\f_{\sL_2}\Bigl(\langle\, \Omega', (\risS{1.5}{strut}{}{15}{10}{10})^n\rangle\Bigr)
= 2^n f_n\ \f_{\sL_2}\Bigl(\langle
(\risS{1.5}{strut}{}{15}{10}{10})^n,
    (\risS{1.5}{strut}{}{15}{10}{10})^n\rangle\Bigr)
$$
Now, using Exercise~\ref{ex-pair-In-In} on
page~\pageref{ex-pair-In-In} and the fact that for the circle
without vertices $$\f_{\sL_2}(\risS{-2}{o}{}{10}{10}{10})=3,$$ (see
page~\pageref{LAWS_3G}) we obtain
$$\begin{array}{ccl}
\f_{\sL_2}\Bigl(\langle (\risS{1.5}{strut}{}{15}{10}{10})^n,
    (\risS{1.5}{strut}{}{15}{10}{10})^n\rangle\Bigr) &=&
  (2n+1)(2n)\ \f_{\sL_2}\Bigl(\langle (\risS{1.5}{strut}{}{15}{10}{10})^{n-1},
    (\risS{1.5}{strut}{}{15}{10}{10})^{n=1}\rangle\Bigr) = \dots \\
&=& (2n+1)!\ .
\end{array}
$$
Comparing these two calculations we find that
$$f_n=\frac{1}{4^n(2n+1)!},$$ which is the coefficient at $z^n$ of the
power series expansion of $$\frac{\sinh{\sqrt{z}/2}}{\sqrt{z}/2}.$$
Hence,
$$\exp\Bigl(2\sum_n a_{2n}z^n\Bigr) = \frac{\sinh{\sqrt{z}/2}}{\sqrt{z}/2}\ .$$
Substituting $z=x^2$ we get the equality (\ref{bernoulli-eq}) which
establishes that $\Omega'=\Omega$ and completes the proof of the
Wheeling Theorem \ref{thm:wheeling}.

\subsection{Wheeling for tangle diagrams}\label{subsec:generalwheel}

A version of the Wheeling Theorem exists for more general spaces of
tangle diagrams. For our purposes it is sufficient to consider the
spaces of diagrams for links with two closed components $\xx$ and
$\yy$.

For $D\in\B$ define the operator
$$(\partial_D)_{\xx}:\B(\xx,\yy)\to\B(\xx,\yy)$$
as the sum of all possible ways of glueing all the legs of $D$ to
some of the $\xx$-legs of a diagram in $\B(\xx,\yy)$ that do not
produce components without legs.

\noindent{\bf Exercise.} Show that $(\partial_D)_{\xx}$ respects the
link relations, and, therefore, is well-defined.

Define the wheeling map $\Phi_{\xx}$ as $\chi_{\xx} \circ
(\partial_\Omega)_{\xx}$. (Strictly speaking, we should call it
$(\Phi_0)_{\xx}$ since we have not yet proved that $\Phi=\Phi_0$.)
The Wheeling Theorem can now be generalized as follows:

\begin{xtheorem}
The map
$$\Phi_{\xx}:\B(\xx,\yy)\to\Fxy$$
identifies the $\B(\xx)$-module $\B(\xx,\yy)$ with the
$\F(\xx)$-module $\Fxy.$
\end{xtheorem}

The proof is, essentially, the same as the proof of the Wheeling
Theorem, and we leave it to the reader.

\section{The unknot and the Hopf link}
\label{wheels}

The arguments similar to those used in the proof of the Wheeling
Theorem allow us to write down an explicit expression for the framed
Kontsevich integral of the zero-framed unknot $O$. Let us denote by
$\Zed(O)$ the image $\chi^{-1}I^{fr}(O)$ of the Konstevich integral
of $O$ in the graded completion of $\B$. (Note that we use the
notation $\Zed(\hl)$ in a similar, but not exactly the same
context.)

\begin{theorem}\label{th:ki-unkn} 
\begin{equation}\label{ki-unkn}
\index{Kontsevich integral!of the unknot} \index{Unknot!Kontsevich
integral} \index{Wheels formula} \Zed(O)=\Omega= \exp
\sum_{n=1}^\infty b_{2n} w_{2n}.
\end{equation}
\end{theorem}
\noindent A very similar formula holds for the Kontsevich integral
of the Hopf link $\hl$:
\begin{theorem}\label{thm:hopflink}
$$\Zed(\hl)=\sum_n\frac{u^n}{n!}\cup\Omega.$$
\end{theorem}

This formula implies that the maps $\Phi$ and $\Phi_0$ of the previous
section, in fact, coincide.

We start the proof with a lemma.

\begin{lemma}\label{lemma:nogi}
If $C_1,\ldots,C_n$ are non-trivial elements of the algebra $\F$,
then $\chi^{-1}(C_1\#\ldots\#C_n)$ is a combination of diagrams in
$\B$ with at least $n$ legs.
\end{lemma}
\begin{proof}
We shall use the same notation as before. If $D\in\B$ is a diagram,
we denote by $D_{\yy}$ the result of decorating all the legs of $D$
with the label $\yy$. Applying the operation $\psi^{2\cdot 1}_{\yy}$
the components obtained from $\yy$ will be called $\yy_1$ and
$\yy_2$

Let $D_i=(\chi\circ\partial_{\Omega})^{-1}(C_i)\in\B$ be the inverse
of $C_i$ under the wheeling map. By the Wheeling Theorem we have
that
$$\chi^{-1}(C_1\#\ldots\#C_n)=
\partial_{\Omega}(D_1\cup\ldots\cup D_n)= \langle \Omega_{\yy_1},
\psi^{2\cdot 1}_{\yy} (D_1\cup\ldots\cup D_n)\rangle.$$ Decompose
$\psi^{2\cdot 1}_{\yy} (D_i)$ as a sum $(D_i)_{\yy_1}+D'_i$ where
$D'_i$ contains only diagrams with at least one leg labelled by
$\yy_2$.

Recall that in the completion of the algebra $\Bbig$ we have
$\partbig_D(\Omega)=\langle D, \Omega \rangle\, \Omega.$ By
projecting this equality to $\B$ we see that $\partial_D(\Omega)$
vanishes unless $D$ is empty. Hence,
$$\langle \Omega_{\yy_1}, (D_1)_{\yy_1}\cup\psi^{2\cdot 1}_{\yy} (D_2\cup\ldots\cup D_n)
\rangle= \langle(\partial_{D_1}\Omega) _{\yy_1}, \psi^{2\cdot
1}_{\yy} (D_2\cup\ldots\cup D_n) \rangle=0.$$ As a result we have
\begin{align*}
\partial_{\Omega}(D_1\cup\ldots\cup D_n)&= \langle \Omega_{\yy_1},
 (D_1)_{\yy_1}\cup\psi^{2\cdot 1}_{\yy} (D_2\cup\ldots\cup D_n)\rangle\\
&\qquad\qquad{+\langle\Omega_{\yy_1},
D'_1\cup\psi^{2\cdot 1}_{\yy} (D_2\cup\ldots\cup D_n)\rangle}\\
&= \langle \Omega_{\yy_1},
D'_1\cup\psi^{2\cdot 1}_{\yy} (D_2\cup\ldots\cup D_n)\rangle\\
&= \langle \Omega_{\yy_1}, D'_1\cup\ldots\cup D'_n\rangle.
\end{align*}
Each of the $D'_i$ has at least one leg labelled $\yy_2$, and these
legs are preserved by taking the pairing with respect to the label
$\yy_1$.
\end{proof}

\subsection{The Kontsevich integral of the unknot}

The calculation of the Kontsevich integral for the unknot is based
on the following geometric fact: the $n$th connected cabling of the
unknot is again an unknot.

The cabling formula on page~\pageref{cabelKI} in this case reads
\begin{equation}\label{eq:cabling-the-unknot}
\psi^{n}\Bigl((I^{fr}(O)\#
\exp{\Bigl(\textstyle\frac{1}{2n}\,\Onechord}\Bigr)\Bigr)
=I^{fr}(O^{(p,1)})=I^{fr}(O)\#\exp{\Bigl(\textstyle\frac{n}{2}\,\Onechord}\Bigr).
\end{equation}

In each degree, the right-hand side of this formula depends on $n$
polynomially. The  term of degree 0 in $n$ is precisely the
Kontsevich integral of the unknot $I^{fr}(O)$.

As a consequence, the left-hand side also contains only non-negative
powers of $n$. We shall be specifically interested in the terms that
are of degree 0 in $n$.

The operator $\psi^{n}$ has a particularly simple form in the
algebra $\B$ (see Section~\ref{diagram_cabling_general}): it
multiplies a diagram with $k$ legs by $n^k$. Let us expand the
argument of $\psi^{n}$ into a power series and convert it to $\B$
term by term.

It follows from Lemma~\ref{lemma:nogi} that if a diagram $D$ is contained in
$$\chi^{-1}\Bigl(I^{fr}(O)\# {\Bigl(\textstyle\frac{1}{2n}\,\Onechord\Bigr)^k}\Bigr),$$
then it has $k'\geq k$ legs. Moreover, by the same lemma it can have
precisely $k$ legs only if it is contained in
$$\chi^{-1}{\Bigl(\textstyle\frac{1}{2n}\,\Onechord\Bigr)^k}.$$
Applying $\psi^{n}$, we multiply $D$ by $n^{k'}$, hence the
coefficient of $D$ on the left-hand side of
(\ref{eq:cabling-the-unknot}) depends on $n$ as $n^{k'-k}$. We see
that if  the coefficient of $D$ is of degree 0 as a function of $n$,
then the number of legs of $D$ must be equal to the degree of $D$.

Thus we have proved that $\Zed(O)$ contains only diagrams whose
number of legs is equal to their degree. We have seen in
\ref{subsec:KI-and-Phi} that the part of the Kontsevich integral of
the Hopf link that consists of such diagrams has the form
$\sum\frac{u^n}{n!}\cup\Omega$. Deleting from this expression the
diagrams with legs attached to the interval component, we obtain
$\Omega$. On the other hand, this is the Kontsevich integral of the
unknot $\Zed(O)$.

\subsection{The Kontsevich integral of the Hopf link}

The Kontsevich integral of the Hopf link both of whose components
are closed with zero framing is computed in \cite{BLT}. Such Hopf
link can be obtained from the zero-framed unknot in three steps:
first, change the framing of the unknot from $0$ to $+1$, then take
the disconnected twofold cabling, and, finally, change the framings
of the resulting components from $+1$ to $0$. We know how the
Kontsevich integral behaves under all these operations and this
gives us the following theorem (see \ref{subsec:generalwheel} for
notation):

\begin{xtheorem}
Let $\clhl$ \label{Hopf-clhl}\index{Link!Hopf} be the Hopf link both
of whose components are closed with zero framing and oriented
counterclockwise. Then
$$I^{fr}(\clhl)=(\Phi_{\xx}\circ\Phi_{\yy}) (\exp{|^{\xx}_{\yy}}),$$
where $|^{\xx}_{\yy}\in\B(\xx,\yy)$ is an interval with one $\xx$
leg and one $\yy$-leg.
\end{xtheorem}

We shall obtain Theorem~\ref{thm:hopflink} from the above statement.

\begin{proof}
Let $O^{+1}$ be the unknot with $+1$ framing. Its Kontsevich
integral is related to that of the zero-framed unknot as in the
theorem on page~\pageref{th:framed}:
$$I^{fr}(O^{+1})=I^{fr}(O)\#\exp{\Bigl(\textstyle\frac{1}{2}\,\Onechord}\Bigr).$$
Applying the inverse of the wheeling map we get
\begin{align*}
\partial_{\Omega}^{-1}\Zed(O^{+1})&=\partial^{-1}_{\Omega}\Zed(O)\cup
\exp\left(\frac{1}{2}\partial^{-1}_{\Omega}\chi^{-1}\Bigl(\Onechord\Bigr)\right)\\
&=\Omega\cup\exp\left(\frac{1}{2}\partial^{-1}_{\Omega}(\risS{1.5}{strut}{}{15}{10}{10})\right).
\end{align*}
Recall that in the proof of Lemma~\ref{lemma:nogi} we have seen that
$\partial_D(\Omega)=0$ unless $D$ is empty. In particular,
$\partial^{-1}_{\Omega}(\Omega)=\Omega$. We see that
\begin{equation}\label{eq:ifroplus}
\Zed(O^{+1})=\partial_{\Omega}\bigl(\Omega \cup \eohst \bigr),
\end{equation}
since $\partial^{-1}_{\Omega}(\risS{1.5}{strut}{}{15}{10}{10})=
\risS{1.5}{strut}{}{15}{10}{10}$.

Our next goal is the following formula:
\begin{equation}\label{eq:exphalfinterval}
\partial^{-2}_{\Omega}(\Zed(O^{+1}))= \eohst\ .
\end{equation}

Applying $\partial_{\Omega}$ to both sides of this equation and
using (\ref{eq:ifroplus}), we obtain an equivalent form of
(\ref{eq:exphalfinterval}):
$$\partial_{\Omega}\bigl(\eohst\bigr)=
\Omega\cup \eohst\ .$$ To prove it, we observe that
\begin{align*}
\partial_{\Omega}\bigl(\eohst\bigr)&=
\langle\Omega_{\yy_1}, \psi^{2\cdot 1}_{\yy}  \eohst  \rangle_{\yy_1}   \\
&=\langle\Omega_{\yy_1}, \exohy{\yy_1}{\yy_1}
\exp\bigl(\bigl|_{\yy_1}^{\yy_2}\bigr)
    \exohy{\yy_2}{\yy_2} \rangle_{\yy_1}.
\end{align*}
The pairing $\B(\yy_1,\yy_2)\otimes\B(\yy_1)\to\B(\yy_2)$ satisfies
$$\langle C, A\cup B \rangle_{\yy_1} = \langle \partial_B(C), A  \rangle_{\yy_1}.$$
for all $A,B\in\B(\yy_1)$, $C\in\B(\yy_1,\yy_2)$.
Therefore, the last expression can be re-written as
$$\langle\partial_{\exp{(\frac{1}{2}|_{\yy_1}^{\yy_1})}}\Omega_{\yy_1}, \exp(\bigl|_{\yy_1}^{\yy_2})\rangle_{\yy_1}
\cup\exohy{\yy_2}{\yy_2}\ .$$ Taking into the account the fact that
$\partial_D(\Omega)=0$ unless $D$ is empty, we see that this is the
same thing as
$$\langle\Omega_{\yy_1}, \exp(\bigl|_{\yy_1}^{\yy_2})\rangle_{\yy_1}\cup\exohy{\yy_2}{\yy_2}
=\Omega\cup \eohst,$$ and this proves (\ref{eq:exphalfinterval}).

To proceed, we need the following simple observation:
\begin{xlemma}
$$\psi^{2\cdot 1}_{\yy} \partial_C(D)=(\partial_C)_{\yy_1}(\psi^{2\cdot 1}_{\yy}  (D))=
(\partial_C)_{\yy_2}(\psi^{2\cdot 1}_{\yy}  (D)).$$
\end{xlemma}

Now, let $\clhl^{+1}$ be the Hopf link both of whose components are
closed with $+1$ framing. Since
$\partial_{\Omega}^{-1}=\partial_{\Omega^{-1}}$, the above lemma and
the cabling formula on page~\pageref{cabelKI} imply that
$$\psi^{2\cdot 1}\partial^{-2}_{\Omega}(\Zed(O^{+1}))=
(\partial_{\Omega})_{\yy_1}^{-1}(\partial_{\Omega})_{\yy_2}^{-1}\bigl(\chi^{-1}_{\yy_1,\yy_2}
I^{fr} (\clhl^{+1})\bigr).$$ On the other hand, this, by
(\ref{eq:exphalfinterval}) is equal to
$$\psi^{2\cdot 1}\eohst =
\exp(\bigl|_{\yy_1}^{\yy_2})\cdot\exohy{\yy_1}{\yy_1}\cdot\exohy{\yy_2}{\yy_2}\
.$$ Applying $\Phi_{\yy_1}\circ \Phi_{\yy_2}$ to the first
expression, we get exactly $I^{fr}(\clhl^{+1})$. On the second
expression, this evaluates to
$$\Phi_{\yy_1}(\Phi_{\yy_2}(\exp(\bigl|_{\yy_1}^{\yy_2}))\#
\exp_{\#}\Bigl(\textstyle\frac{1}{2}\,\Onechord_{\yy_1}\Bigr)\#
\exp_{\#}\Bigl(\frac{1}{2}\,\Onechord_{\yy_2}\Bigr).$$
Changing the framing, we see that
$$I^{fr}(\clhl)=(\Phi_{\yy_1}\circ \Phi_{\yy_2}) (\exp{\bigl|^{\yy_1}_{\yy_2}}).$$
The statement of the theorem follows by a simple change of notation.
\end{proof}

\begin{proof}[Proof of Theorem~\ref{thm:hopflink}]

First, let us observe that
for any diagram $D\in\B$ we have
$$(\partial_D)_{\xx} \exp{\bigl|_{\xx}^{\yy}}=D_{\yy}\cup
\exp{\bigl|_{\xx}^{\yy}}.$$
Now, we have
\begin{align*}
I^{fr}(\hl) \# \chi_{\xx}(\Omega_{\xx}) &= I^{fr}(\clhl)\\
 &= \Phi_{\xx}(\Phi_{\yy}(\exp{\bigl|_{\xx}^{\yy}}))\\
 &= \Phi_{\xx}(\exp{\bigl|_{\xx}^{\yy}} \cup\Omega_{\xx})
 && \text{by the observation above}.
\end{align*}
Since $\partial_{\Omega}=\Omega$, it follows that
$\chi_{\xx}(\Omega_{\xx}=\Phi_{\xx}(\Omega_{\xx}))$ and
\begin{align*}
I^{fr}(\hl) &= \Phi_{\xx}(\Omega_{\xx}^{-1})\#\Phi_{\xx}(\exp{\bigl|_{\xx}^{\yy}} \cup\Omega_{\xx})\\
    &=\Phi_{\xx}(\exp{\bigl|_{\xx}^{\yy}} \cup\Omega_{\xx}\cup\Omega_{\xx}^{-1})
    && \text{by the Wheeling Theorem}\\
    &=\Phi_{\xx}(\exp{\bigl|_{\xx}^{\yy}})\\
    &=\chi_{\xx}(\Omega\cup\exp{\bigl|_{\xx}^{\yy}}).
\end{align*}

\end{proof}

\section{Rozansky's rationality conjecture}\label{Roz_r_conj}

This section concerns a generalization of the wheels formula for the
Kontsevich integral of the unknot to arbitrary knots. The
generalization is, however, not complete -- the Rozansky--Kricker
theorem does not give an explicit formula, it only suggests that
$I^{fr}(K)$ can be written in a certain form.

It turns out that the terms of the Kontsevich integral $I^{fr}(K)$
with values in $\B$ can be rearranged into lines corresponding to
the number of loops in open diagrams from $\B$. Namely, for any term
of $I^{fr}(K)$, shaving off all legs of the corresponding diagram
$G\in\B$, we get a trivalent graph $\gamma$.
Infinitely many terms of $I^{fr}(K)$ give rise to the same $\gamma$.
It turns out that these terms behave in a regular fashion, so that
it is possible to recover all of them from $\gamma$ and some finite
information.

To make this statement precise we introduce {\em marked open
diagrams} which are represented by a trivalent graph whose edges are
marked by power series (it does matter on which side of the edge the
mark is located, and we shall indicate the side in question by a
small leg near the mark). We use such marked open diagrams to
represent infinite series of open diagrams which differ by the
number of legs. More specifically, an edge marked by a power series
$f(x)=c_0+c_1x+c_2x^2+c_3x^3+\dots$ stands for the following series
of open diagrams:
\newcommand\vmark[2]{\rb{-18pt}[20pt][20pt]{
    \begin{picture}(20,40)(0,0)
    \put(0,0){\ig[height=40pt]{#1.eps}} #2
    \end{picture}}}
$$\vmark{maropdi-f}{\put(7,22){\mbox{$\scriptstyle f(x)$}}}\quad :=
\quad c_0 \vmark{maropdi-0}{} +
      c_1 \vmark{maropdi-1}{} + c_2 \vmark{maropdi-2}{} +
      c_3 \vmark{maropdi-3}{} + \dots
$$
In this notation the wheels formula (Theorem \ref{th:ki-unkn}) can
be written as

$$\ln I^{fr}(O) = \rb{-12pt}[10pt][20pt]{
    \begin{picture}(20,40)(0,0)
    \put(0,0){\ig[width=20pt]{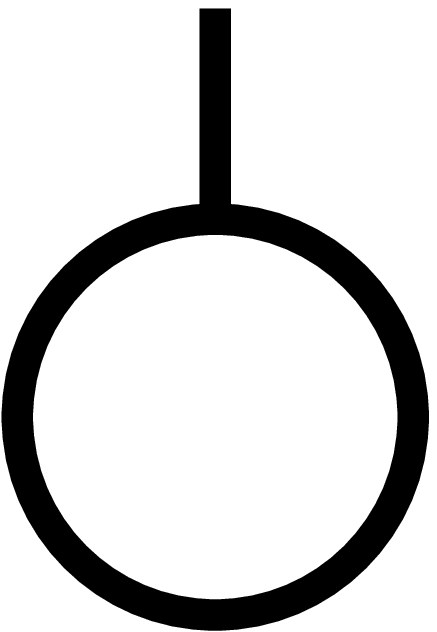}}
    \put(12,25){\mbox{$\scriptstyle \frac{1}{2}
                \ln\frac{e^{x/2}-e^{-x/2}}{x}$}}
    \end{picture}}
$$
Now we can state the

\noindent{\bf Rozansky's rationality conjecture.} \cite{Roz2} {\it
$$\begin{array}{rcl}
\ln I^{fr}(K) &=&\displaystyle \rb{-20pt}[40pt][30pt]{
    \begin{picture}(120,60)(0,0)
    \put(0,0){\ig[width=120pt]{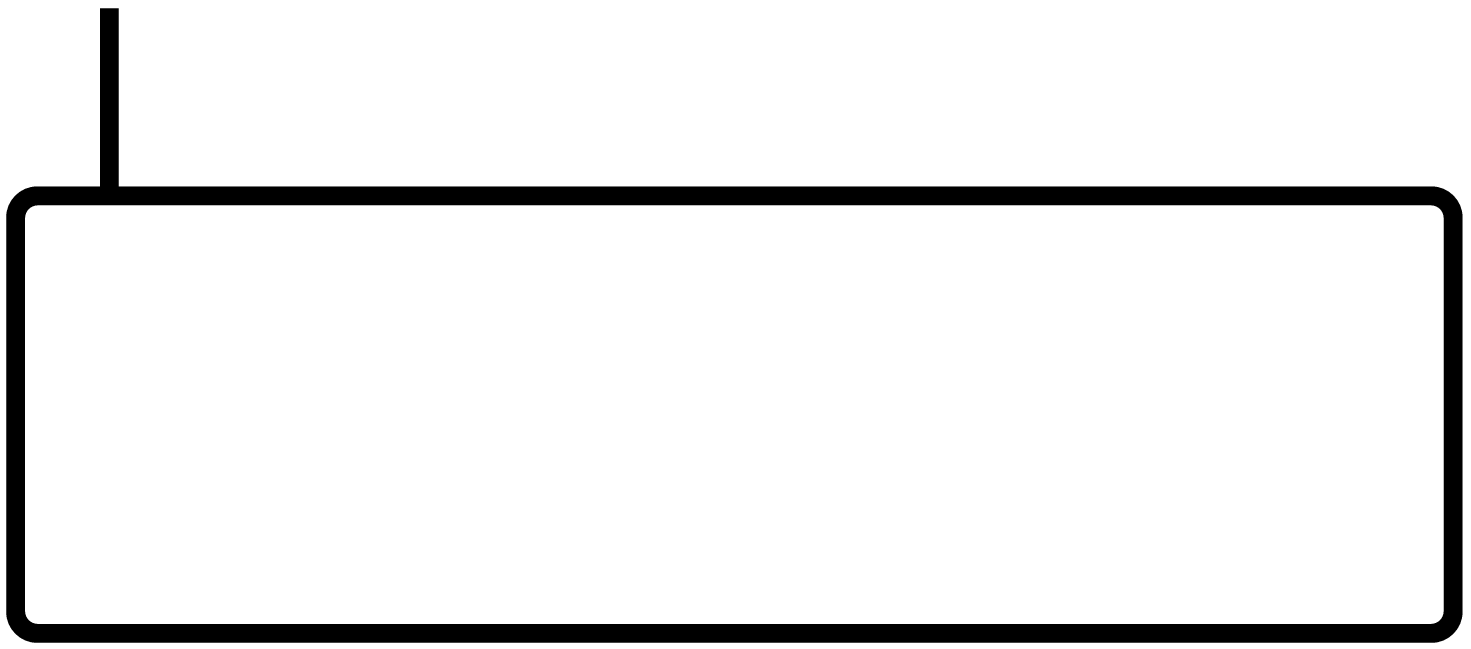}}
    \put(12,45){\mbox{$\scriptstyle
       \frac{1}{2} \ln\frac{e^{x/2}-e^{-x/2}}{x} -
       \frac{1}{2} \ln A_K(e^x)$}}
    \end{picture}} + \sum_i^{\mbox{\scriptsize\rm finite}}
    \rb{-20pt}[40pt][20pt]{
    \begin{picture}(120,80)(0,0)
    \put(0,0){\ig[width=120pt]{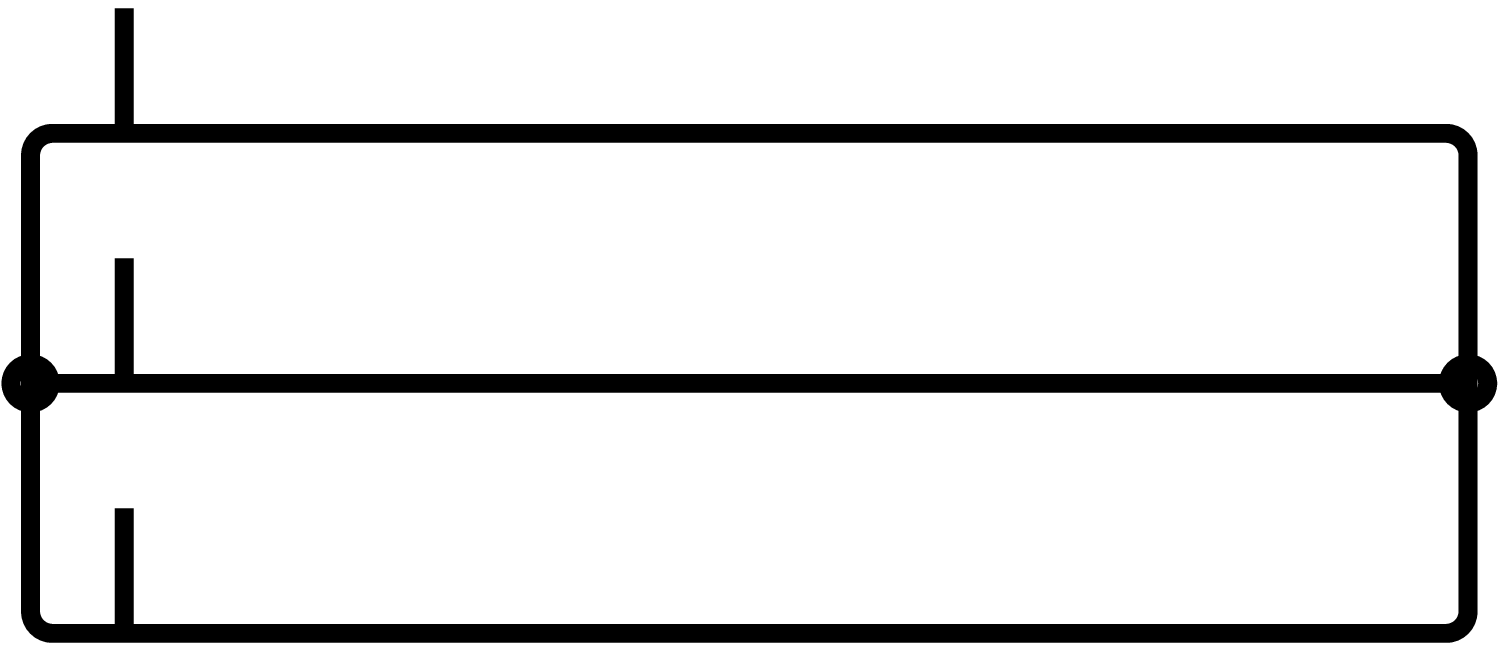}}
    \put(14,48){\mbox{$\scriptstyle p_{i,1}(e^x)/A_K(e^x)$}}
    \put(14,27){\mbox{$\scriptstyle p_{i,2}(e^x)/A_K(e^x)$}}
    \put(14,7){\mbox{$\scriptstyle p_{i,3}(e^x)/A_K(e^x)$}}
    \end{picture}}\\
&& +\ (\mbox{\rm terms with $\geqslant 3$ loops})\ ,
\end{array}
$$
where $A_K(t)$ is the Alexander polynomial
of $K$ normalized so that $A_K(t)=A_K(t^{-1})$ and $A_K(1)=1$,
$p_{i,j}(t)$ are polynomials, and the higher loop terms mean the sum
over marked trivalent graphs (with finitely many copies of each
graph) whose edges are marked by a polynomial in $e^x$ divided by
$A_K(t)$.}

The word ``rationality'' refers to the fact that the labels on all
3-graphs of degree $\ge1$ are rational functions of $e^x$. The
conjecture was proved by A.~Kricker in \cite{Kri2}. Due to $\AS$ and
$\IHX$ relations the specified presentation of the Kontsevich
integral is not unique. Hence the polynomials $p_{i,j}(t)$
themselves cannot be knot invariants. However, there are certain
combinations of these polynomials
 that are genuine knot invariants.
For example, consider the polynomial
$$\Theta'_K(t_1,t_2,t_3)=\sum_i p_{i,1}(t_1)p_{i,2}(t_2)p_{i,3}(t_3)\ .$$
Its symmetrization,
$$\Theta_K(t_1,t_2,t_3)=\sum_{\substack{\e=\pm1\\ \{i,j,k\}=\{1,2,3\}}}
\Theta'_K(t^\e_i,t^\e_j,t^\e_k) \in
 \Q[t^{\pm1}_1,t^{\pm1}_2,t^{\pm1}_3]/(t_1t_2t_3=1)\ ,
$$
over the order 12 group of symmetries of the theta graph, is a knot
invariant. It is called the {\it 2-loop polynomial} of $K$.
\index{Knot invariant!2-loop polynomial}\index{Two-loop polynomial}
Its values on knots with  few crossings are tabulated in
\cite{Roz2}. T.~Ohtsuki \cite{Oht2} found a cabling formula for the
2-loop polynomial and its values on torus knots $T(p,q)$.

\begin{xcb}{Exercises}
\begin{enumerate}

\item $^*$
Find two chord diagrams with 11 chords which have the same
intersection graph but unequal modulo four- and one-term relations.
According to Section \ref{mutkn}, eleven is the least number of
chords for such diagrams. Their existence is known, but no explicit
examples were found yet.

\item $^*$
In the algebra $\A$ consider the subspace $\A^M$ generated by those
chord diagrams whose class in $\A$ is determined by their
intersection graph only. It is natural to regard the quotient space
$\A/\A^M$ as the space of chord diagram {\em distinguishing
mutants}. \index{Chord diagram!distinguishing mutants} Find the
dimension of $\A_n/\A^M_n$. It is known that it is zero for $n\leq
10$ and greater than zero for $n=11$. Is it true that
$\dim(\A_{11}/\A^M_{11})=1$?

\item Find a basis in the space of canonical invariants of degree 4.

{\sl Answer:} $j_4, c_4+c_2/6, c2^2$.

\item
Show that the self-linking number defined in Section~\ref{self_link}
is a canonical framed Vassiliev invariant of order 1.

\item\label{ex_ex_con_lim}
Show the existence of the limit from Section~\ref{conconway}
$$A = \lim_{N\to 0} \frac{\theta_{{\sL_N},V}\bigr|_{q=e^h}}{N}\,.$$

{\sl Hint.}  Choose a complexity function
on link diagrams in such a way that two of the diagrams participating in
the skein relation for $\theta_{{\sL_N},V}$ are strictly simpler
then the third one. Then use induction on complexity.

\item\label{ex_ws_subst}
Let $\displaystyle f(h)=\sum_{n=0}^\infty f_n h^n$ and
$\displaystyle g(h)=\sum_{n=0}^\infty g_n h^n$ be two power series
Vassiliev invariants. 
\begin{enumerate}
\item[(a)] Prove that their product $f(h)\cdot g(h)$ as formal power
series in $h$ is a Vassiliev power series invariant, and
$$\symb(f\cdot g) = \symb(f)\cdot \symb(g)\,.$$
\item[(b)]Suppose that $f$ and $g$ are
related to each other via substitution and multiplication:
$$f(h)=\b(h)\cdot g\bigl(\a(h)\bigr)\,,$$
where $\a(h)$ and $\b(h)$ are formal power series in $h$, and
$$\a(h) = ah + (\mbox{\scriptsize terms of degree $\geqslant 2$})\,,
    \qquad
  \b(h) = 1 + (\mbox{\scriptsize terms of degree $\geqslant 1$})\,.
$$
Prove that $\symb(f_n)=a^n \symb(g_n)$.
\end{enumerate}

\item
Prove that a canonical Vassiliev invariant is primitive if its
symbol is primitive.

\item\label{ex_prod_can}
Prove that the product of any two canonical Vassiliev power series
is a canonical Vassiliev power series.

\item
If $v$ is a canonical Vassiliev invariant of odd order and $K$ an
amphicheiral knot, then $v(K)=0$.

\item
Let $\kappa\in\W_n$ be a weight system of degree $n$. Construct
another weight system $(\kappa\circ\psi^p)'\in\W_n$, where $\psi^p$
is the $p$th connected cabling operator, and $(\cdot)'$ is the
deframing operator from Section~\ref{defram_ws}. We get a function
$f_\kappa\colon p\mapsto (\kappa\circ\psi^p )'$ with values in
$\W_n$. Prove that

\begin{enumerate}
\item[(a)] $f_\kappa(p)$ is a polynomial in $p$ of degree
             $\leqslant n$ if $n$ is even, and of degree
       $\leqslant n-1$ if $n$ is odd.

\item[(b)] \parbox[t]{3in}{The $n$th degree term of the polynomial
$f_\kappa(p)$ is equal to $\displaystyle -\frac{\kappa(\ol{w}_n)}{2}
\symb(c_n) p^n\,,$ where  $\ol{w}_n$ is }\qquad $\rb{10pt}{$\ol{w}_n
=$}\quad \risS{-4}{wssl2wn}{
\put(0,-8){\mbox{\scriptsize $n$ spokes}}}{30}{20}{15}$\vspace{2pt}\\
 the wheel with $n$ spokes, and $c_n$ is the
$n$th coefficient of the Conway polynomial.

\end{enumerate}

\item\label{ex-fr-ki-hump}
Find the framed Kontsevich integrals
$Z^{fr}(H)$ and $I^{fr}(h)$
for the hump unknot with zero framing up to order 4.

\medskip
{\sl Answer.}
$$\begin{array}{ccl}
Z^{fr}\Bigl(\rb{-6pt}{\ig[width=15pt]{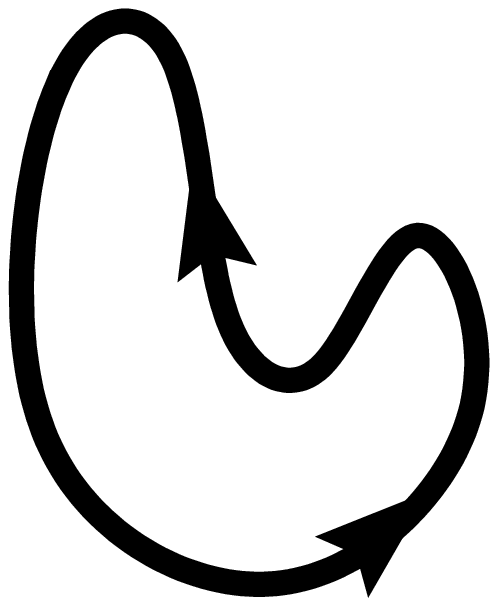}}\Bigr)
   &=& 1 - \frac{1}{24} \cdWO + \frac{1}{24} \cdWW
       + \frac{7}{5760} \rb{-4mm}{\ig[height=9mm]{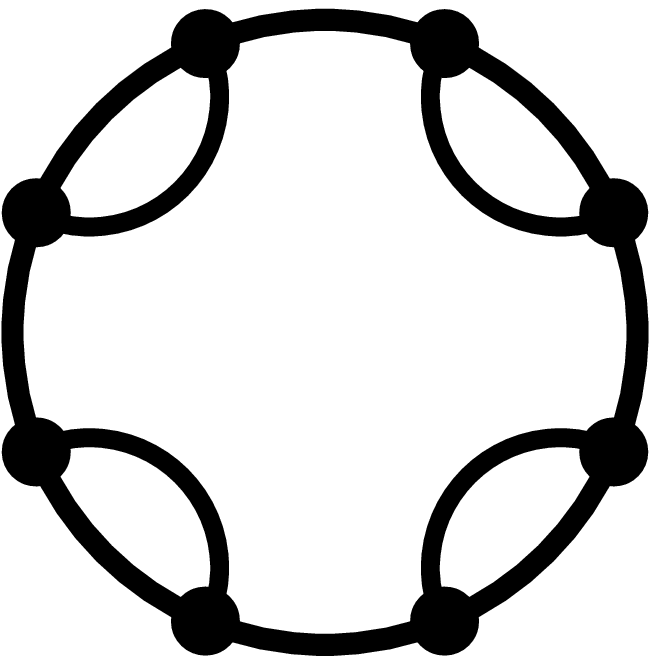}}
       - \frac{17}{5760} \rb{-4mm}{\ig[height=9mm]{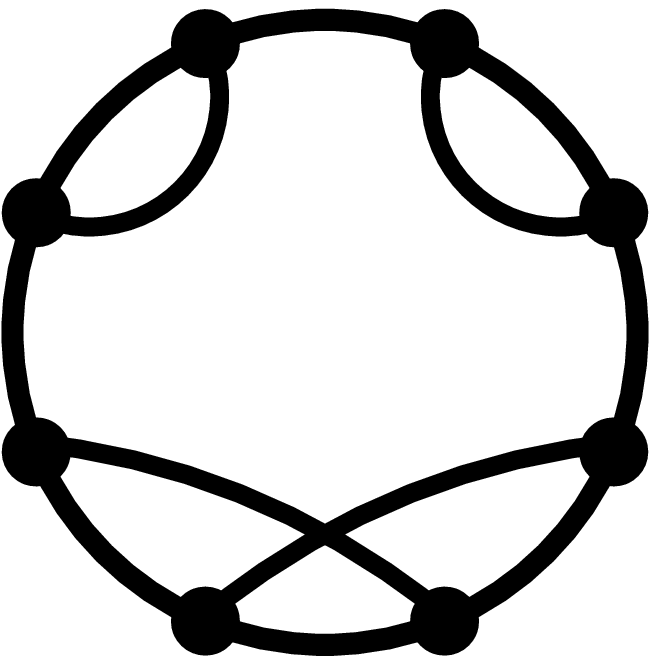}}
   \vspace{10pt} \\ &&
       + \frac{7}{2880} \rb{-4mm}{\ig[height=9mm]{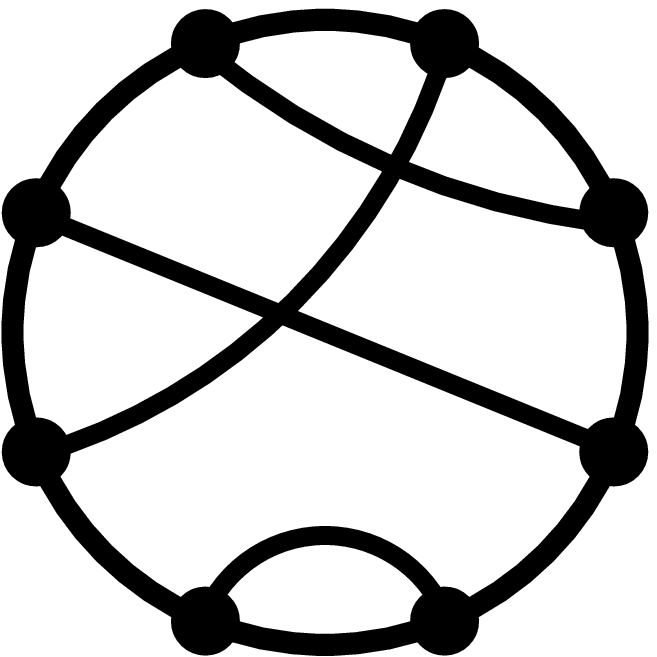}}
       - \frac{1}{720} \rb{-4mm}{\ig[height=9mm]{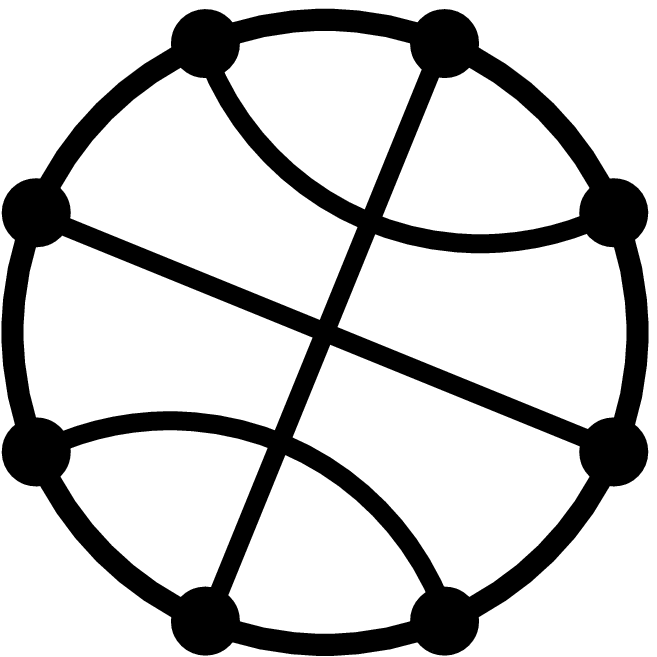}}
       + \frac{1}{1920} \chd{cd4-06} + \frac{1}{5760} \chd{cd4-02}\ .
   \vspace{15pt} \\
I^{fr}\Bigl(\rb{-6pt}{\ig[width=15pt]{humpm.eps}}\Bigr)&=&
 1\Big/Z^{fr}\Bigl(\rb{-6pt}{\ig[width=15pt]{humpm.eps}}\Bigr)\ .
\end{array}
$$

\item
Using Exercise~\ref{A2C} to Chapter \ref{algFDchap}
(page~\pageref{A2C}) show that up to degree 4
$$\begin{array}{ccl}
Z^{fr}\Bigl(\rb{-6pt}{\ig[width=15pt]{humpm.eps}}\Bigr)
   &=& 1-\frac1{48}\ \chdh{f25_Wb}+\frac1{4608}\ \chdh{f25_Fd}
               +\frac1{46080}\ \risS{-10}{f25_Fe}{}{28}{20}{15}
               +\frac1{5760}\ \chd{f25_Ff}\ , \\
I^{fr}\Bigl(\rb{-6pt}{\ig[width=15pt]{humpm.eps}}\Bigr)
   &=& 1+\frac1{48}\ \chdh{f25_Wb}+\frac1{4608}\ \chdh{f25_Fd}
               -\frac1{46080}\ \risS{-10}{f25_Fe}{}{28}{20}{15}
               -\frac1{5760}\ \chd{f25_Ff}\ .
\end{array}
$$

\item
Using the previous problem and Exercise~\ref{ex_C2B_34leg} to
Chapter \ref{algFDchap} (page~\pageref{ex_C2B_34leg}) prove that up
to degree 4
$$\Zed(O) = \chi^{-1}I^{fr}(O) = 1+\frac1{48}\ \chdn{-.8}{ccWb}
   +\frac1{4608}\ \chdr{-6}{c25_Fd} -\frac1{5760}\ \chdr{-10}{c25_Ff}\ .
$$
This result confirms Theorem \ref{th:ki-unkn} from page
\pageref{th:ki-unkn} up to degree 4.

\item
\parbox[t]{3in}{Compute the framed Kontsevich integral $Z^{fr}(\hl)$
\index{Kontsevich integral! of the Hopf link} up to degree 4 for the
Hopf link $\hl$ with one vertical interval component $\xx$ and one
closed component $\yy$ represented by the tangle on the picture.
Write the result as an element of $\F(\xx,\yy)$.}\qquad
\parbox[t]{1in}{$\rb{-30pt}{$\hl =$}\quad \risS{-60}{phl}{
      \put(3,-5){\mbox{$\xx$}}
      \put(32,40){\mbox{$\yy$}}}{30}{20}{15}$}
\medskip

{\sl Answer.}
$$\textstyle
Z^{fr}(\hl) = \yu{yu-ed} + \yu{yu-u} +\frac{1}{2}\ \yu{yu-uu}
     +\frac{1}{6}\ \yu{yu-uuu} -\frac{1}{24}\ \chd{yu-bu}
     +\frac{1}{24}\ \yu{yu-uuuu} -\frac{1}{48}\ \chd{yu-buu}\ .
$$

\item
Compute the final framed Kontsevich integral $I^{fr}(\hl)$ up to
degree 4:
$$\begin{array}{ccl}
I^{fr}(\hl) &=& \yu{yu-ed} + \yu{yu-u} +\frac{1}{2}\ \yu{yu-uu}
        +\frac{1}{48}\ \yu{yu-edo} +\frac{1}{6}\ \yu{yu-uuu}
        -\frac{1}{24}\ \chd{yu-bu}+\frac{1}{48}\ \yu{yu-uo} \\
&&\hspace{-20pt}+\frac{1}{24}\ \yu{yu-uuuu} -\frac{1}{48}\
\chd{yu-buu}
  +\frac1{4608}\ \yu{yu-edoo} -\frac1{46080}\ \yu{yu-edooo}
  -\frac1{5760}\ \Id{yu-edq} +\frac{1}{96}\ \risS{-10}{yu-uuo}{}{18}{20}{15} \ .
\end{array}
$$

\item\label{ex:zed-hl-in-B}
Using the previous problem and Exercise~\ref{ex_C2B_34leg} to
Chapter~\ref{algFDchap} (page~\pageref{A2C}) prove that up to degree
4
$$\begin{array}{ccl}
\Zed(\hl) = \chi^{-1}_{\yy} I^{fr}(\hl) &=& \Iu{I}\hspace{-5pt} +
\Iu{I-u}+\frac12\ \Iu{I-uu} +\frac16\ \Iu{I-uuu} +\frac1{24}\
\Iu{I-uuuu} +\frac1{48}\ \Id{I-o}
+\frac1{48}\ \Id{I-uo} \\
&&+\frac1{96}\ \Id{I-uuo} -\frac1{5760}\ \Id{I-q} +\frac1{4608}\
\Id{I-oo} +\frac1{384}\ \Id{I-bo}\ .
\end{array}
$$
Indicate the parts of this expression forming $\Zed_0(\hl)$,
$\Zed_1(\hl)$,
 $\Zed_2(\hl)$ up to degree 4.
This result confirms Theorem \ref{thm:hopflink} from page
\pageref{thm:hopflink} up to degree 4.

\item
Prove that $\chi\circ\partial_{\Omega}:\B\to\F$ is a bialgebra
isomorphism.

\item
Compute
$\chi\circ\partial_{\Omega}(\risS{-2}{strut-3}{}{15}{10}{10})$,\quad
$\chi\circ\partial_{\Omega}(\risS{-1}{strut-w2}{}{15}{10}{10})$,\quad
$\chi\circ\partial_{\Omega}(w_6)$.

\item\label{ex-pair-In-In}
Show that the pairing $\langle
(\risS{1.5}{strut}{}{15}{10}{10})^n,(\risS{1.5}{strut}{}{15}{10}{10})^n\rangle$
satisfies the recurrence relation
$$\langle (\risS{1.5}{strut}{}{15}{10}{10})^n,(\risS{1.5}{strut}{}{15}{10}{10})^n\rangle
=  2n\cdot \bigl( \risS{-2}{o}{}{10}{10}{10} + 2n-2 \bigr)\cdot
    \langle (\risS{1.5}{strut}{}{15}{10}{10})^{n-1},
            (\risS{1.5}{strut}{}{15}{10}{10})^{n-1}\rangle\ ,
$$
where $\risS{-2}{o}{}{10}{10}{10}$ is a 3-graph in
$\G_0\subset\G\subset\Bbig$ of degree 0 represented by a circle
without vertices and multiplication is understood in algebra $\Bbig$
(disjoint union).

\item
Prove that, after being carried over from $\B$ to $\A^{fr}$, the
right-hand side of Equation~\ref{ki-unkn} (page~\pageref{ki-unkn})
belongs in fact to the subalgebra $\A\subset\A^{fr}$.
Find an explicit expression of the series through some basis of $\A$
up to degree 4.
\medskip

{\sl Answer.} The first terms of the infinite series giving the Kontsevich integral of
the unknot, are:
$$
  I(O) = 1 - \frac{1}{24} \cdWW
   - \frac{1}{5760} \chd{cd4-02}
   + \frac{1}{1152} \chd{cd4-06}
   + \frac{1}{2880} \chd{cd4-07} +\dots
$$
Note that this agrees with the answer to
Exercise~\ref{ex-fr-ki-hump}.

\end{enumerate}

\end{xcb}
 %10 FineKI
\chapter{Braids and string links} % 12
\label{chapBr}

Essentially, the theory of Vassiliev invariants of braids is a
particular case of the Vassiliev theory for tangles, and the main
constructions are very similar to the case of knots. There is,
however, one big difference: many of the questions that are still
open for knots are rather easy to answer in the case of braids.
This, in part, can be explained by the fact that braids form as
group, and it turns out that the whole Vassiliev theory for braids
can be described in group-theoretic terms. In this chapter we shall
see that the Vassiliev filtration on the pure braid groups coincides
with the filtrations coming from the nilpotency theory of groups. In
fact, for any given group the nilpotency theory could be thought of
as a theory of finite type invariants.

The group-theoretic techniques of this chapter can be used to study
knots and links. One such application is the theorem of Goussarov
which says that $n$-equivalence classes of string links on $m$
strands form a group. Another application of the same methods is a
proof that Vassiliev invariants of pure braids extend to invariants
of string links of the same order. In order to make these
connections we shall describe a certain braid closure that produces
string links out of pure braids.

The theory of the finite type invariants for braids was first
developed by T.~Kohno \cite{Koh1, Koh2} several years before
Vassiliev knot invariants were introduced. The connection between
the theory of commutators in braid groups and the Vassiliev knot
invariants was first made by T.~Stanford \cite{Sta4}.

\section{Basics of the theory of nilpotent groups}
\label{section:nilpotent_groups} We shall start by reviewing some
basic notions related to nilpotency in groups. These will not only
serve us a technical tool: we shall see that the theory of Vassiliev
knot invariants has the logic and structure similar to those of the
theory of nilpotent groups. The classical reference for nilpotent
groups are the lecture notes of P.\ Hall \cite{Hall}. For modern
aspects of the theory see \cite{MPbook}. Iterated integrals are
described in the papers of K.T.Chen \cite{Chen1, Chen2}.

\subsection{The dimension series}
Let $G$ be an arbitrary group. The {\em group algebra}\index{Group algebra}
$\Z G$ consists of finite linear combinations $\sum a_ig_i$ where
$g_i$ are elements of $G$ and $a_i$ are integers. The product in $\Z
G$ is the linear extension of the product in $G$. The group itself
can be considered as a subset of the group algebra if we identify
$g$ with $1 g$. The identity in $G$ is the unit in $\Z G$ and we
shall denote it simply by 1.

Let $JG\subset\Z G$ be the {\em augmentation
ideal},\index{Augmentation ideal} that is, the kernel of the
homomorphism $\Z G\to \Z$ that sends each $g\in G$ to $1$. Elements
of $JG$ are the linear combinations $\sum a_ig_i$ with $\sum a_i=0$.
The powers $J^n G$ of the augmentation ideal form a descending
filtration on $\Z G$.

Let $\Di_k(G)$ 
be the subset of $G$ consisting of all $g\in G$ such that $$g-1\in
J^kG.$$ Obviously, the neutral element of $G$ always belongs to
$\Di_k(G)$. Also, for all $g,h\in \Di_k(G)$ we have
$$gh-1=(g-1)(h-1)+(g-1)+(h-1)$$
and, hence, $\Di_k(G)$ is closed under the product. Finally,
$$g^{-1}-1=-(g-1)g^{-1}$$
which shows that $\Di_k(G)$ is a subgroup of $G$; it is called the
$k$th {\em dimension subgroup} of $G$. Clearly, $\Di_1(G)=G$ and for
each $k$ the subgroup $\Di_{k+1}(G)$ is contained in $\Di_k(G)$.

\begin{xxca} Show that $\Di_k(G)$ is invariant under all
automorphisms of $G$. In particular, it is a normal subgroup of $G$.
\end{xxca}
The descending series of subgroups
$$G=\Di_1(G)\supseteq \Di_2(G) \supseteq \Di_3(G)\supseteq \ldots$$
is called the {\em dimension series of $G$}. \index{Dimension
series}

Consider the {\em group commutator} which can be\footnote{there are
other, equally good, options, such as $[g,h]=ghg^{-1}h^{-1}$.}
defined as
$$[g,h]:=g^{-1}h^{-1}gh.$$ If $g\in\Di_p(G)$ and $h\in\Di_q(G)$, we
have
$$
g^{-1}h^{-1}gh-1=g^{-1}h^{-1}\bigl((g-1)(h-1)-(h-1)(g-1)\bigr),$$
and, hence, $[g,h]\in\Di_{p+q}(G)$. It follows that the group
commutator descends to a bilinear bracket on
$$\Di(G):=\bigoplus_k \Di_k(G)/\Di_{k+1}(G).$$

\noindent{\bf Exercise.} Show that this bracket on $\Di(G)$ is
antisymmetric and satisfies the Jacobi identity. In other words,
show that the commutator endows $\Di(G)$ with the structure of a
{\em Lie ring}.

This exercise implies that $\Di(G)\ot\Q$ is a Lie algebra over the
rationals. The universal enveloping algebra of $\Di(G)\ot\Q$ admits
a very simple description. Denote by $\A_k(G)$ the quotient
$J^kG/J^{k+1}G$. Then the direct sum
$$\A(G):=\bigoplus \A_k(G)$$
is a graded algebra, with the product induced by that of $\Z G$.
Quillen shows in \cite{Quillen-JA} that $\A(G)\ot \Q$ is the
universal enveloping algebra of $\Di(G)\ot\Q$.

The dimension series can be generalized by replacing the integer
coefficients in the definition of the group algebra by coefficients
in an arbitrary ring. The augmentation ideal is defined in the same
fashion as consisting of linear combinations whose coefficients add
up to zero, and the arguments given in this section for integer
coefficients remain unchanged. We denote the $k$th dimension
subgroup of $G$ over a ring $\Ring$ by $\Di_k(G,\Ring)$. Since there
is a canonical homomorphism of the integers to any ring, $\Di_k(G)$
is contained in $\Di_k(G,\Ring)$ for any ring $\Ring$.

\begin{xxca} Show that $\Di_k(G,\Q)/\Di_{k+1}(G,\Q)$ is a torsion-free
abelian group for any $k$.
\end{xxca}

\subsection{Commutators and the lower central series}
For many groups, the dimension subgroups can be described entirely
in terms of  group commutators. For $H,K$ normal subgroups of $G$,
denote by $[H,K]$ the subgroup of $G$ generated by all the
commutators of the form $[h,k]$ with $h\in H$ and $k\in K$. The {\em
lower central series} \index{Lower central series} subgroups
$\gamma_k G$ of a group $G$ are defined inductively by setting
$\gamma_1 G=G$ and
$$\gamma_{k}G=[\gamma_{k-1} G,G].$$ A group $G$ is called {\em
nilpotent} if $\gamma_n G= 1$ for some $n$. The maximal $n$ such
that $\gamma_n G\neq 1$ is called the {\em nilpotency class} of $G$.
If the intersection of all $\gamma_n G$ is trivial, the group $G$ is
called {\em residually nilpotent}.

\begin{xxca} Show that $\gamma_k G$ is invariant under all automorphisms of
$G$. In particular, it is a normal subgroup of $G$.
\end{xxca}
\begin{xxca} Show that for any $p$ and $q$ the subgroup $[\gamma_p G, \gamma_q G]$
is contained in $\gamma_{p+q} G$.
\end{xxca}

We have already seen that the commutator of any two elements of $G$
belongs to $\Di_2 (G)$. Using induction, it is not hard to show that
$\gamma_k G$ is always contained in $\Di_k (G)$. If $\gamma_k G$ is
actually the same thing as the $k$th dimension subgroup of $G$ over
the integers, it is said that $G$ has the {\em dimension subgroup
property}. Many groups have the dimension subgroup property. In
fact, it was conjectured that {\em all} groups have this property
until E.~Rips found a counterexample, published in 1972 \cite{Rips}.
His counterexample was later simplified; we refer to \cite{MPbook}
for the current state of knowledge in this field. In general, if
$x\in \Di_k (G)$, there exists $q$ such that $x^q\in\gamma_k G$, and
the group $\Di_k(G,{\Q} )$ consists of all $x$ with this property,
see Theorem \ref{theorem:truejennings} on page
\pageref{theorem:truejennings} .

The subtlety of the difference between the lower central series and
the dimension subgroups is underlined by the fact that for all
groups $\gamma_kG=\Di_k(G)$ when $k<4$. In order to give the reader
some feeling of the subject let us treat one simple case here:

\begin{xproposition}
For any group $G$ we have $\gamma_2G=\Di_2(G)$.
\end{xproposition}
\begin{proof}
First let us assume that $G$ is abelian, that is, $\gamma_2G=1$ (or,
in additive notation, $\gamma_2G=0$). In this case there is a
homomorphism of abelian groups $s:\Z G\to G$ defined by replacing a
formal linear combination by a linear combination in $G$. The
homomorphism $s$ sends $g-1\in\Z G$ to $g\in G$. On the other hand,
$s(J^2G)=0$. Indeed, it is easy to check that $J^2G$ is additively
spanned by products of the type $(x-1)(y-1)$ with $x,y\in G$; we
have
$$s\bigl((x-1)(y-1)\bigr)=s(xy-x-y+1)=xy-x-y=0$$
since $G$ is abelian. It follows that $g\in \Di_2(G)$ implies that
$g=0$.

Now, let $G$ be an arbitrary group. It can be seen from the
definitions that group homomorphisms respect both the dimension
series and the lower central series. Moreover, it is clear that a
surjective homomorphism of groups induces surjections on the
corresponding terms of the lower central series. This means that if
$\Di_2(G)$ is strictly greater than $\gamma_2G$, the same is true
for $G/\gamma_2G$. On the other hand, $G/\gamma_2G$ is abelian.
\end{proof}

Recall that $H_1(G)$ is the {\em abelianization} of $G$, that is,
its maximal abelian quotient $G/\gamma_2 G$, and that
$H_1(G,\R)=H_1(G)\ot\R$.

\begin{xxca} Show that $H_1(G,\R)$ is canonically isomorphic to $\A_1(G)\ot\R$.
\end{xxca}

\subsection{Filtrations induced by series of subgroups}\label{subsec:canonicalfiltration}

Let $\{G_i\}$  be a descending series of subgroups
$$G=G_1\supseteq G_2\supseteq\ldots$$
of a group $G$ with the property that $[G_p, G_q]\subseteq G_{p+q}$.
For $x\in G$ denote by $\mu(x)$ the maximal $k$ such that $x\in
G_k$. Let $\Q G$  be the group algebra of $G$ with rational
coefficients and $E_nG$ its ideal spanned by the products of the
form $(x_1-1)\cdot\ldots\cdot (x_s-1)$ with
$\sum_{i=1}^s\mu(x_i)\geq n$. We have the filtration of $\Q G$:
$$\Q G\supset J G=E_1G\supseteq E_2G\supseteq\ldots.$$
This filtration is referred to as the {\em canonical filtration
induced by the series $\{G_n\}$}.\index{Filtration!canonical, induced
by a series}\index{Canonical filtration}

\begin{theorem}\label{theorem:jennings}
Let $$G=G_1\supseteq G_2\supseteq\ldots\supseteq G_N=\{1\}$$ be a
finite series of subgroups of a group $G$ with the property that
$[G_p, G_q]\subseteq G_{p+q}$, and such that $G_i/G_{i+1}$ is
torsion-free for all $1\leq i< N$. Then for all $i\geq 1$
$$G_i= G\cap(1+E_iG),$$
where $\{E_iG\}$ is the canonical filtration of $\Q G$ induced by
$\{G_i\}$.
\end{theorem}

As stated above, this theorem can be found in \cite{Passi, Passman}.
The most important case of it has been proved by Jennings
\cite{Jennings}, see also \cite{Hall}. It clarifies the relationship
between the dimension series and the lower central series.

For a subset $H$ of a group $G$ let $\sqrt{H}$ be the set of all
$x\in G$ such that $x^p\in H$ for some $p>0$. If $H$ is a normal
subgroup, and $G/H$ is nilpotent, then $\sqrt{H}$ is again a normal
subgroup of $G$. The set $\sqrt{\{1\}}$ is precisely the set of all
periodic (torsion) elements of $G$; it is a subgroup if $G$ is
nilpotent.

\begin{theorem}\label{lemma:sqrt}
Let $$G=G_1\supseteq G_2\supseteq\ldots \supseteq G_N=\{1\}$$ be a
finite series of subgroups of a group $G$ with the property that
$[G_p, G_q]\subseteq G_{p+q}$. Then $[\sqrt{G_p},
\sqrt{G_q}]\subseteq \sqrt{G_{p+q}}$ and the canonical filtration of
$\Q G$ induced by $\{\sqrt{G_i}\}$ coincides with the filtration
induced by $\{G_i\}$.
\end{theorem}
For the proof see the proofs of Lemmas~1.3 and 1.4 in Chapter IV of
\cite{Passi}.

Now, consider a nilpotent group $G$. We have mentioned that
${\gamma_n G}$ is always contained in $\Di_n(G)$, and, hence, $E_nG$
in this case coincides with $J^nG$. It follows from Theorems
\ref{theorem:jennings} and \ref{lemma:sqrt} that
$$\Di_n(G)=\sqrt{\gamma_n G}$$
for all $n$. The assumption that $G$ is nilpotent can be removed by
considering the group $G/{\gamma_n G}$ instead of $G$, and we get
the following characterization of the dimension series over $\Q$:
\begin{theorem}[Jennings,
\cite{Jennings}]\label{theorem:truejennings} For an arbitrary group
$G$, an element $x$ of $G$ belongs to $\Di_n(G,\Q)$ if and only if
$x^r\in \gamma_n G$ for some $r>0$.
\end{theorem}

\subsection{Semi-direct products}
\label{subsection:semidirect} The augmentation ideals, the dimension
series and the lower central series behave in a predictable way
under taking direct products of groups. When $G=G_1\times G_2$ we
have $$\Z G= \Z G_1\otimes \Z G_2.$$ Moreover,
$$J^kG=\sum_{i+j=k}J^iG_1\otimes J^jG_2,$$
and this implies
$$\Di_k(G)=\Di_k(G_1)\times \Di_k(G_2).$$ It is also easy to see that
$$\gamma_kG=\gamma_kG_1\times \gamma_kG_2.$$

When $G$ is a {\em semi-direct}, rather than direct, product of
$G_1$ and $G_2$ these isomorphisms break down in general. However,
they do extend to one particular case of semi-direct products,
namely, the {\em almost direct product} defined as follows.

Having a semi-direct product $A\ltimes B$ is the same as having an
action of $B$ on $A$ by automorphisms. An action of $B$ on $A$ gives
rise to an action of $B$ on the abelianization of $A$; we say that a
semi-direct product $A\ltimes B$ is {almost direct}\index{Almost
direct product}\label{almostdirect} if this latter action is
trivial.

\begin{xproposition}
For an almost direct product  $G=G_1\ltimes G_2$
$$\gamma_kG=\gamma_kG_1\ltimes
\gamma_kG_2$$ for all $k$. Moreover,
$$J^k(G)=\sum_{i+j=k}J^i(G_1)\otimes J^j(G_2),$$
inside $\Z G$ and, hence, $$\Di_kG=\Di_kG_1\ltimes \Di_kG_2$$ and
$$\A(G)=\A(G_1)\ot\A(G_2)$$
as a graded $\Z$-module.
\end{xproposition}

The proof is not difficult and we leave it as an exercise. The case
of the lower central series is proved in \cite{FR}; for the
dimension subgroups see \cite{Papadima} (or \cite{MostWill} for the
case when $G$ is a pure braid group).

\subsection{The free group} 

Let $x_1,\ldots,x_m$ be a set of free generators of the free group
$F_m$ and set $X_i=x_i-1\in\Z F_m$. Then, for any $k>0$ each element
$w\in F_m$ can be uniquely expressed inside $\Z F_m$ as
$$w=1+\sum_{1\leq i\leq m} a_i X_i+
\ldots+ \sum_{1\leq i_1,\ldots, i_k\leq m}  a_{i_1,\ldots,i_k}
X_{i_1}\ldots X_{i_k}+r(w),$$ where $a_{i_1,\ldots,i_j}$ are
integers and $r(w)\in J^{k+1}F_m$. This formula can be considered as
a Taylor formula for the free group. In fact, the coefficients
$a_{i_1,\ldots,i_j}$ can be interpreted as some kind of derivatives,
see \cite{Fox}.

To show that such formula exists, it is enough to have it for the
generators of $F_m$ and their inverses:
$$x_i=1+X_i$$
and
$$x_i^{-1}=1-X_i+X_i^2-\ldots+(-1)^{k}X_i^{k}+(-1)^{k+1}X_i^{k+1}x_i^{-1}.$$
The uniqueness of the coefficients $a_{i_1,\ldots,i_j}$ will be
clear from the construction below.

Having defined the Taylor formula we can go further and define
something like the Taylor series.

Let $\Z\la X_1,\ldots,X_m \ra$ be the algebra of formal power series
in $m$ non-commuting variables $X_i$. Consider the homomorphism of
$F_m$ into the group of units of this algebra
$$\M: F_m\to \Z \la X_1,\ldots,X_m \ra,$$
which sends the $i$th generator $x_i$ of $F_m$ to $1+X_i$. In
particular,
$$\M(x_i^{-1})=1-X_i+X_i^2-X_i^3+\ldots$$
This homomorphism is called the {\em Magnus expansion}.
\index{Magnus expansion} It is injective: the Magnus expansion of a
reduced word $x_{\alpha_1}^{\varepsilon_1}\ldots
x_{\alpha_k}^{\varepsilon_k}$ contains the monomial
$X_{\alpha_1}\ldots X_{\alpha_k}$ with the coefficient
${\varepsilon_1}\ldots {\varepsilon_k}$, and, hence, the kernel of
$\M$ is trivial.

The Magnus expansion is very useful since it gives a simple test for
an element of the free group to belong to a given dimension
subgroup.
\begin{xlemma}
For $w\in F_m$ the power series $\M(w)-1$ starts with terms of
degree $k$ if and only if $w\in \Di_k (F_m)$ and $w\notin \Di_{k+1}
(F_m)$.
\end{xlemma}
\begin{proof}
Extend the Magnus expansion by linearity to the group algebra $\Z
F_m$. The augmentation ideal is sent by $\M$ to the set of power
series with no constant term and, hence, the Magnus expansion of
anything in $J^{k+1} F_m$ starts with terms of degree at least
$k+1$. It follows that the first $k$ terms of the Magnus expansion
coincide with the first $k$ terms of the Taylor formula. Notice that
this implies the uniqueness of the coefficients in the Taylor
formula. Now, the term of lowest non-zero degree on right-hand side
of the Taylor formula has degree $k$ if and only if $w-1\in J^k
F_m$.
\end{proof}

One can easily see that the non-commutative monomials of degree $k$
in the $X_i$ give a basis for $J^{k} F_m/J^{k+1} F_m$. The Magnus
expansion gives a map
$$\M_k: J^{k} F_m\to \A_k(F_m)= J^{k} F_m/J^{k+1} F_m$$
which sends $x\in J^{k} F_m$ to the degree $k$ term of $\M(x)$. The
following is straightforward:
\begin{lemma}\label{lemma:magnus_uni}
The map $\M_k$ is the quotient map $J^{k} F_m\to \A_{k}(F_m)$.
\end{lemma}

The dimension subgroups for the free group coincide with the
corresponding terms of the lower central series. In other words, the
free group $F_m$ has the dimension subgroup property. A proof can be
found, for example, in Section~5.7 of \cite{MKS}. As a consequence,
we see that the free groups are residually nilpotent since the
kernel of the Magnus expansion is trivial.

\subsection{Chen's iterated integrals}

The Magnus expansion for the free group does not generalize readily
to arbitrary groups. However, there is a general geometric
construction which works in the same way for all finitely generated
groups and detects the terms of the dimension series for any group
just as the Magnus expansion detects them for the free group. This
construction is given by {\em Chen's iterated
integrals},\index{Chen's iterated integral} \cite{Chen1, Chen2}. We
shall only describe it very briefly here since we shall not need it
in the sequel. An accessible introduction to Chen's integrals can be
found in \cite{Ha}.

Let us assume that the group $G$ is the fundamental group of a
smooth manifold $M$. Let $X_1,\ldots, X_m$ be a basis for
$H_1(M,\R)$ and $w_1,\ldots, w_m$ be a set of real closed 1-forms on
$M$ representing the basis of $H^1(M,\R)$ dual to the basis
$\{X_i\}$.

Consider an expression
$$\alpha=\sum_i \alpha_i X_i+ \sum_{i,j} \alpha_{ij} X_iX_j+\sum_{i,j,k} \alpha_{ijk} X_iX_jX_k+\ldots$$
where all the coefficients $\alpha_{*}$ are 1-forms on $M$. We shall
say that $\alpha$ is a {\em $\R\la X_1,\ldots,X_m\ra$-valued 1-form}
on $M$. We refer to  $\sum_i \alpha_i X_i$ as the {\em linear part}
of $\alpha$. Denote by ${\mathfrak{x}}$ the ideal in $\R\la
X_1,\ldots,X_m\ra$ consisting of the power series with no constant
term.

In \cite{Chen1} K.T.Chen proves the following fact:
\begin{xtheorem}
There exists a $\R\la X_1,\ldots,X_m\ra$-valued 1-form  $w$ on $M$
whose linear part is  $\sum_i w_i X_i$ and an ideal ${\mathfrak{j}}$
of $\R\la X_1,\ldots,X_m\ra$ such that there is a ring homomorphism
$$Z:\R\pi_1M\to\R\la X_1,\ldots,X_m\ra/\mathfrak{j}$$
given by
$$
  Z(g) = \sum\limits_{\substack{{0\leq k}\\
            {1\le i_1,\ldots,i_k\le m}}}\
         \int\limits_{0<t_k<\dots<t_1<1}
         w(t_1)\wedge\ldots\wedge w(t_k)\
         ,
$$
where $w(t)$ is the pull-back to the interval $[0,1]$ of the 1-form
$w$ under a map $\gamma: [0,1]\to M$ representing $g$, with the
property that the kernel of the composite map is
$$Z:\R\pi_1M\to\R\la X_1,\ldots,X_m\ra/\mathfrak{j}\to\R\la X_1,\ldots,X_m\ra/(\mathfrak{j}+\mathfrak{x}^n)$$
is precisely $J^n(\pi_1M)\ot\R$.
\end{xtheorem}
We shall call the map $Z$ the {\em Chen expansion}.

In certain important situations the algebra $\R\la
X_1,\ldots,X_m\ra/\mathfrak{j}$ can be replaced by the algebra
$\A(\pi_1M)\ot \R$. Suppose that the algebra $\Lambda M$ of the
differential forms on $M$ has a differential graded subalgebra $A$
with the following properties:
\begin{itemize}
\item the inclusion $A\to \Lambda M$ induces isomorphisms in
cohomology in all dimensions;
\item each element in $H^*(M, \R)$ can be represented by a closed form
in $A$ so that there is a direct sum decomposition
$$A=H^*(M, \R)\oplus A'$$ where $A'$ is an ideal.
\end{itemize}
Chen shows in \cite{Chen2} (Lemma 3.4.2) that in this situation the
ideal $\mathfrak{j}$ is actually homogeneous. As a consequence, the
algebra $\R\la X_1,\ldots,X_m\ra/\mathfrak{j}$ is graded. Since the
Chen expansion sends $J^n(\pi_1M)\ot\R$ to the terms of degree $n$
and higher, it induces an injective map
$$\A(\pi_1M)\ot \R \to \R\la X_1,\ldots,X_m\ra/\mathfrak{j}$$
and the image of the Chen expansion is contained in the graded
completion of the image of this map. This means that the 1-form
$$w=\sum_i w_i X_i+ \sum_{i,j} w_{ij} X_iX_j+\ldots$$
is, actually, $\A(\pi_1M)\ot \R$-valued and we can think of $X_i$ as
the generators of $\A_1(\pi_1M)$.

Examples of manifolds with a subalgebra $A$ satisfying the above
conditions include all compact Kahler manifolds. Another example
which will be of importance for us is the configuration space of $k$
distinct ordered particles $z_1,\ldots, z_k$ in $\C$: its
fundamental group is the pure braid group on $n$ strands. If we
allow complex, rather than real coefficients in the Chen expansion,
we obtain a particularly simple form $w$ which only contains linear
terms:
$$w=\frac{1}{2\pi i} \sum d\log{(z_i-z_j)} \cdot X_{ij},$$
where $X_{ij}$ can be thought of as a chord diagram with one
horizontal chord connecting the $i$th and the $j$th strands.
Comparing the definitions, we see that the Chen expansion of a pure
braid coincides exactly with its Kontsevich integral.

\section{Vassiliev invariants for free groups} \label{section:vasfree}
The main subject of this chapter are the Vassiliev braid invariants,
and, more specifically, the invariants of {\em pure braids}, that
is, the braids whose associated permutation is trivial. Pure braids
are a particular case of tangles and thus we have a general recipe
for constructing their Vassiliev invariants. The only special
feature of braids is the requirement that the tangent vector to a
strand is nowhere horizontal. This leads to the fact that the chord
diagrams for braids have only horizontal chords on a skeleton
consisting of vertical lines; the relations they satisfy are the
horizontal 4T-relations.

We shall start by treating what may seem to be a very particular
case: braids on $m+1$ strands whose all strands, apart from the last
(the rightmost) one, are straight. Such a braid can be thought of as
the graph of a path of a particle in a plane with $m$ punctures.
(The punctures correspond to the vertical strands.) The set of
equivalence classes of such braids can be identified with the
fundamental group of the punctured plane, that is, with the free
group $F_m$ on $m$ generators $x_i$, where $1\leq i\leq m$.
\begin{figure}[htb]
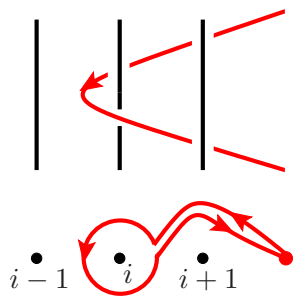

$$\risS{-15}{braids-1}{
\put(-8,3){$i-1$}\put(35,5){$i$}\put(56,3){$i+1$}}{100}{80}{10}$$
\caption{The generator $x_i$ of $F_m$ as a braid and as a path in a
plane with $m$ punctures.}
\end{figure}

A {\em singular path} in the $m$-punctured plane is represented by a
braid with a finite number of transversal double points, whose first
$m$ strands are vertical. Resolving the double points of a singular
path with the help of the Vassiliev skein relation we obtain an
element of the group algebra $\Z F_m$. Singular paths with $k$
double points span the $k$th term of a descending filtration on $\Z
F_m$ which is analogous to the singular knot filtration on $\Z\K$,
defined in Section~\ref{sing_knot_filtr}. A Vassiliev invariant of
order $k$ for the free group $F_m$ is, of course, just a linear map
from $\Z F_m$ to some abelian group that vanishes on singular paths
with more than $k$ double points.

Tangle chord diagrams which correspond to singular paths have a very
specific form: these are horizontal chord diagrams (see
page~\pageref{A^h(n)}) on $m+1$ strands whose all chords have one
endpoint on the last strand. Such diagrams form an algebra, which we
denote temporarily by $\A'(F_m)$, freely generated by $m$ diagrams
of degree 1. We shall see in Section~\ref{vinv_and_magexp} that this
algebra is a subalgebra of $\A(m+1)$, or, equivalently, that the
horizontal 4T relations do not imply any relations in $\A'(F_m)$.

The radical difference between the singular knots and singular paths
(and, for that matter, arbitrary singular braids) lies in the
following
\begin{lemma}\label{lemma:singprod}
A singular path in the $m$-punctured plane with $k$ double points is
a product of $k$ singular paths with one double point each.
\end{lemma}
This is clear from the picture:
$$\ig[width=6cm]{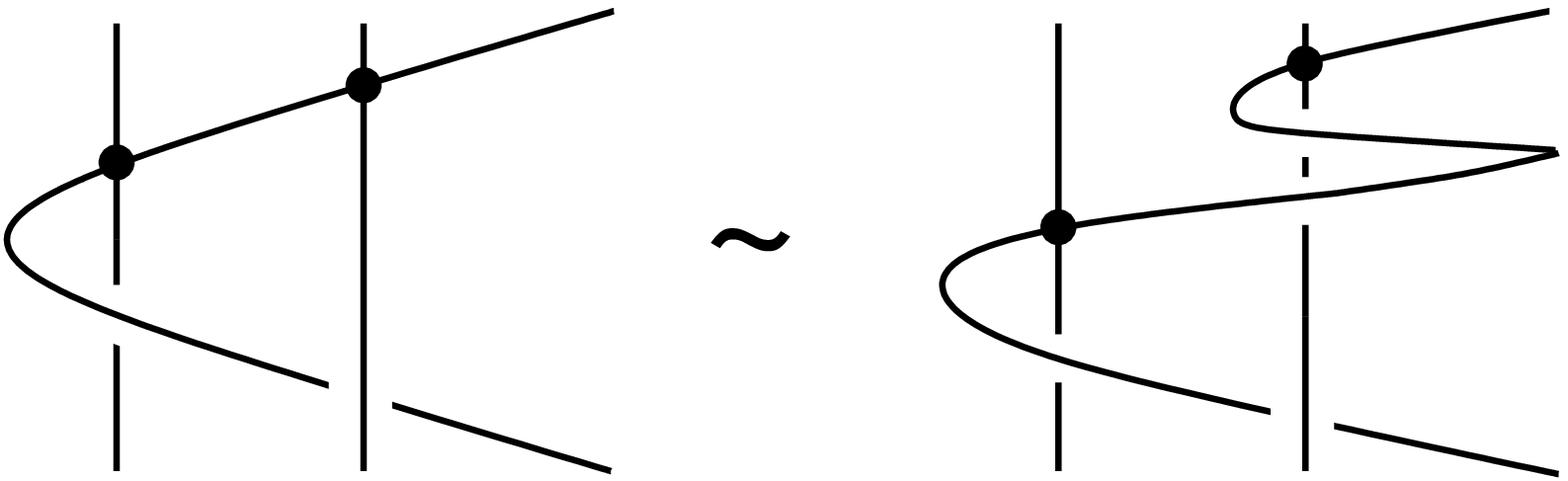}$$
Lemma~\ref{lemma:singprod} allows to describe the singular path
filtration in purely algebraic terms. Namely, singular paths span
the augmentation ideal $J F_m$ in $\Z F_m$ and singular paths with
$k$ double points span the $k$th power of this ideal.

Indeed, each singular path is an alternating sum of non-singular
paths, and, hence, it defines an element of the augmentation ideal
of $F_m$. On the other hand, the augmentation ideal of $F_m$ is
spanned by differences of the form $g-1$ where $g$ is some path. By
successive crossing changes on its braid diagram, the path $g$ can
be made trivial. Let $g_1\, ,\, \ldots\, ,\, g_s$ be the sequence of
paths obtained in the process of changing the crossings from $g$ to
$1$. Then
$$g-1=(g-g_1)+(g_1-g_2)+...+(g_s-1),$$
where the difference enclosed by each pair of brackets is a singular
path with one double point.

We see that the Vassiliev invariants are those that vanish on some
power of the augmentation ideal of $F_m$. The dimension subgroups of
$F_m$ are the counterpart of the Goussarov filtration: $\Di_k F_m$
consists of elements that cannot be distinguished from the unit by
Vassiliev invariants of order less than $k$. We shall refer to these
as to being {\em $k-1$-trivial}.

The algebra $\A'(F_m)$ of chord diagrams for paths is the same thing
as the algebra $$\A(F_m)=\bigoplus_k J^kF_m/J^{k+1}F_m.$$ Indeed,
the set of chord diagrams of degree $k$ is the space of paths with
$k$ double points modulo those with $k+1$ double points. The
generator of $\A(F_m)$ which is the class of the element $x_i-1$,
where $x_i$ is the $i$th generator of $F_m$, is represented by a
chord joining the $i$th and the $m+1$st strands:
$$\rb{-10pt}{$x_i-1\quad =\quad$}\risS{-25}{braids-2a}{
\put(-2,-9){$1$}\put(22,-9){$i$}\put(35,-9){$m+1$}}{50}{10}{30}$$

In fact, the Magnus expansion identifies the algebra $\Z \la
X_1,\ldots,X_m\ra$ with the completion $\widehat{\A}(F_m)$ of the
algebra of the chord diagrams ${\A(F_m)}$. The following statement
is a reformulation of Lemma \ref{lemma:magnus_uni}:

\begin{xtheorem}\index{Vassiliev!invariant!universal for free group}
The Magnus expansion is a universal Vassiliev invariant.
\end{xtheorem}
Since the Magnus expansion is injective, we have
\begin{xcorollary}
The Vassiliev invariants distinguish elements of the free group.
\end{xcorollary}

\subsection{Observation}
If a word $w\in F_m$ contains only positive powers of the generators
$x_i$, the Magnus expansion of $w$ has a transparent combinatorial
meaning: $\M(w)$ is simply the sum of all subwords of $w$, with the
letters capitalized. This is also the logic behind the construction
of the universal invariant for virtual knots discussed in
Chapter~\ref{chapGD}: it associates to a diagram the sum of all its
subdiagrams.

\subsection{The Kontsevich integral}
The Magnus expansion is not the only universal Vassiliev invariant.
(See Exercise~\ref{ex:univas} on page~\pageref{ex:univas}.) Another
important universal invariant is, of course, the Kontsevich
integral. Note that the Kontsevich integral in this case is nothing
but the Chen expansion of $F_m$ where the manifold $M$ is taken to
be the plane $\C$ with $m$ punctures $z_1,\ldots,z_m$ and
$$w=\frac{1}{2\pi i}\cdot\frac{dz}{z-z_j}\cdot X_j.$$
Note that the Kontsevich integral depends on the positions of the
punctures $z_j$ (Exercise~\ref{ex:kinotunique} on
page~\pageref{ex:kinotunique}).

In contrast to the Kontsevich integral, the Magnus expansion has
integer coefficients. We shall see that it also gives rise to a
universal Vassiliev invariant of pure braids with integer
coefficients; however, unlike the Kontsevich integral, this
invariant fails to be multiplicative.

\section{Vassiliev invariants of pure braids}

The interpretation of the Vassiliev invariants for the free group
$F_m$ in terms of the powers of the augmentation ideal in $\Z F_m$
remains valid if the free groups are replaced by the pure braid
groups. One new difficulty is that instead of the free algebra $\A(F_m)$ we
have to study the algebra $\A(P_m)=\A^h(m)$ of horizontal chord
diagrams (see page~\pageref{A^h(n)}). The multiplicative structure
of $\A^h(m)$ is rather complex, but an explicit additive basis for
this algebra can be easily described. This is due to the very
particular structure of the pure braid groups.

\subsection{Pure braids and free groups}\label{subsection:pbfg}
Pure braid groups are, in some sense, very close to being direct
products of free groups.

Erasing one (say, the rightmost) strand of a pure braid on $m$
strands produces a pure braid on $m-1$ strands. This procedure
respects braid multiplication, so, in fact, it gives a homomorphism
$P_m\to P_{m-1}$. Note that this homomorphism has a section
$P_{m-1}\to P_{m}$ defined by adding a vertical non-interacting
strand on the right.

The kernel of erasing the rightmost strand consists of braids on $m$
strands whose first $m-1$ strands are vertical. Such braids are
graphs of paths in a plane with $m-1$ punctures, and they form a
group isomorphic to the free group on $m-1$ letters $F_{m-1}$.

All the above can be re-stated as follows: there is a split
extension
$$ 1\to F_{m-1}\to P_{m}\leftrightarrows P_{m-1}\to 1.$$
It follows that $P_{m}$ is a semi-direct product $F_{m-1}\ltimes
P_{m-1}$, and, proceeding inductively, we see that
\[P_{m}\cong F_{m-1}\ltimes\ldots F_{2}\ltimes F_1.\]
Here $F_{k-1}$ can be identified with the free subgroup of $P_{m}$
formed by pure braids which can be made to be totally straight apart
from the $k$th strand which is allowed to braid around the strands
to the left. As a consequence, every braid in $P_n$ can be written
uniquely as  a product $\beta_{m-1}\beta_{m-2}\ldots\beta_{1}$,
where $\beta_{k}\in F_{k}$. This decomposition is called the {\em
combing} of a pure braid.

\begin{figure}
\ig[width=3cm]{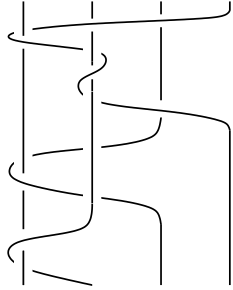} \caption{An example of a combed braid.}
\end{figure}

One can show that the above semi-direct products are not direct (see
Exercise~\ref{ex:pthree} at the end of the chapter). However, they
are almost direct (see the definition on page
\pageref{almostdirect}).

\begin{xlemma}
The semi-direct product $P_{m}=F_{m-1}\ltimes P_{m-1}$ is almost
direct.
\end{xlemma}
\begin{proof}
The abelianization $F^{ab}_{m-1}$ of $F_{m-1}$ is a direct sum of
$m-1$ copies of $\Z$. Given a path $x\in F_{m-1}$, its image in
$F^{ab}_{m-1}$ is given by the $m-1$ linking numbers with each
puncture. The action of a braid $b\in P_{m-1}$ on a generator
$x_i\in F_{m-1}$ consists in ``pushing'' the $x_i$ through the
braid:

$$\ig[width=9cm]{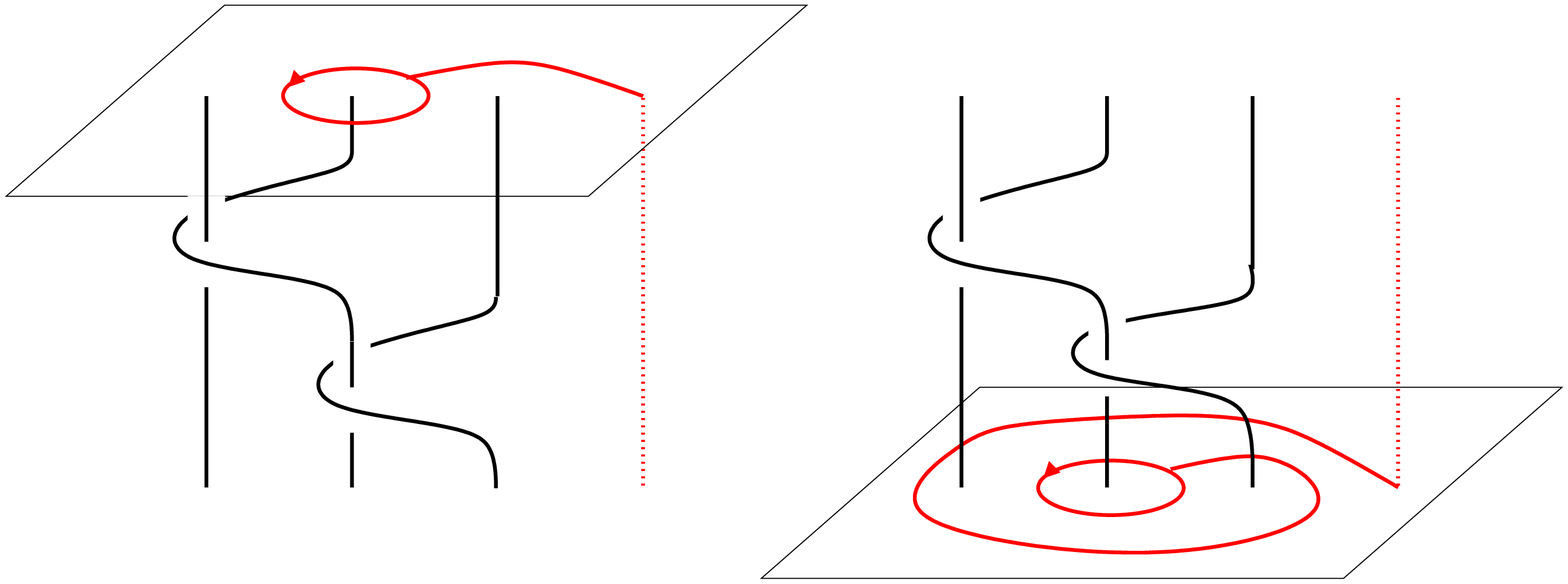}$$

It is clear the linking numbers of the path $b^{-1}x_ib$ with the
punctures in the plane are the same as those of $x_i$, therefore the
action of $P_{m-1}$ on $F^{ab}_{m-1}$ is trivial.
\end{proof}

\begin{xremark}
Strictly speaking, in Section~\ref{link_num} we have only defined
the linking number for two curves in space, while in the above proof
we use the linking number of a point and a loop in a plane. This
number can be defined as the intersection (or incidence) number of
the point with an immersed disk whose boundary is the loop.

Generally, the linking number is defined for two disjoint cycles
(for instance, oriented submanifolds) in $\R^n$ when the sum of the
dimensions of the cycles is one less than $n$, see, for instance,
\cite{Dold}. This linking number is crucial for the definition of
the Alexander duality which we shall use in Chapter~\ref{chap_VSS}.
\end{xremark}

\subsection{Vassiliev invariants and the Magnus
expansion}\label{vinv_and_magexp}

The Vassiliev filtration on the group algebra $\Z P_m$ can be
described in the same algebraic terms as in the
Section~\ref{section:vasfree}. Indeed, singular braids can be
identified with the augmentation ideal $JP_m\subset\Z P_m$. It is
still true that each singular braid with $k$ double points can be
written as a product of $k$ singular braids with one double point
each; therefore, such singular braids span the $k$th power of
$JP_m$. The (linear combinations of) chord diagrams with $k$ chords
are identified with $J^kP_m/J^{k+1}P_m=\A_k(P_m)$ and the Goussarov
filtration on $P_m$ is given by the dimension subgroups
$\Di_k(P_m)$.

Now, since $P_m$ is an almost direct product of $F_{m-1}$ and
$P_{m-1}$ we have that
$$J^k(P_m)=\sum_{i+j=k}J^i(F_{m-1})\otimes J^j(P_{m-1}),$$
$$\A_k(P_m)=\bigoplus_{i+j=k}\A_i(F_{m-1})\otimes \A_j(P_{m-1}),$$
and
$$\Di_k(P_m)=\Di_k(F_{m-1})\ltimes \Di_k(P_{m-1}),$$
see Section \ref{subsection:semidirect}.

These algebraic facts can be re-stated in the language of Vassiliev
invariants as follows.

Firstly, each singular braid with $k$ double points is a linear
combination of {\em combed singular braids} with the same number of
double points. A combed singular braid with $k$ double points is a
product $b_{m-1}b_{m-2}\ldots b_1$ where $b_i$ is a singular path in
$\Z F_i$ with $k_i$ double points, and $k_{m-1}+\ldots+k_1=k$.

Secondly, {\em combed diagrams} form a basis in the space of all
horizontal chord diagrams. A combed diagram $D$\index{Diagram!
combed} is a product $D_{m-1}D_{m-2}\ldots D_1$ where $D_i$ is a
diagram whose all chords have their rightmost end on the $i$th
strand.

Thirdly, a pure braid is $n$-trivial if and only if, when combed, it
becomes a product of $n$-trivial elements of free groups. In
particular, the only braid that is $n$-trivial for all $n$ is the
trivial braid.

Let $\beta\in P_m$ be a combed braid:
$\beta=\beta_{m-1}\beta_{m-2}\ldots\beta_{1}$, where $\beta_{k}\in
F_{k}$. The Magnus expansions of the elements $\beta_i$ can be
``glued together''. Let  $i_k:\A(F_k)\hookrightarrow\A^h(m)$ be the
map that adds $m-k-1$ vertical strands, with no chords on them, to
the right:

$$\ig[width=6cm]{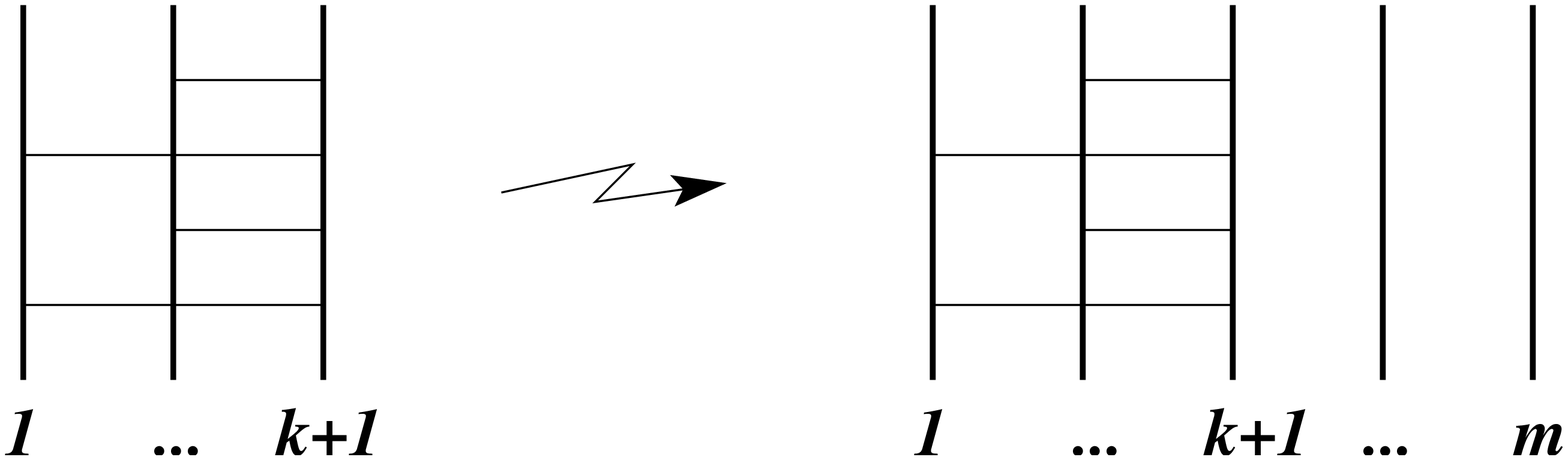}$$

The maps $i_k$ extend to the completions of the algebras $\A(F_k)$
and $\A^h(m)$. Define the Magnus expansion $$\M:P_m\to\Ab^h(m)$$ as
the map sending $\beta$ to $i_{m-1}\M(\beta_{m-1})\ldots
i_{1}\M(\beta_{1})$. For example:

\def\vpicstdheight#1{\rb{-5mm}{\ig[height=13mm]{#1}}}

\begin{align*}
\M\Biggl(\,\vpicstdheight{braids-braid}\,\Biggr) &=\Biggl(
1+\vpicstdheight{braids-diag1} \,\Biggr)
  \Biggl( 1-\vpicstdheight{braids-diag2}+ \vpicstdheight{braids-diag4}
              -\vpicstdheight{braids-diag6}+\ldots  \Biggr)\\
&=1+\vpicstdheight{braids-diag1}-\vpicstdheight{braids-diag2}-
\vpicstdheight{braids-diag3}+\vpicstdheight{braids-diag4}+
\vpicstdheight{braids-diag5}-\vpicstdheight{braids-diag6}+ \ldots.
\end{align*}

\begin{xtheorem}\index{Vassiliev!invariant!universal for pure braids}
The Magnus expansion is a universal Vassiliev invariant of pure
braids.
\end{xtheorem}

As in the case of free groups, the Magnus expansion is injective,
and, therefore, Vassiliev invariants distinguish pure braids. Note
that combing is not multiplicative so the Magnus expansion is not
multiplicative either.

\subsection{A dictionary}
The theory of finite type invariants for the pure braids suggests
the following dictionary between the nilpotency theory for groups
and the theory of Vassiliev invariants:

\begin{center}
\def\lis{\rb{-6pt}{\makebox(0,20){ }}}
\begin{tabular}{|c|c|}
\hline
{\bf Nilpotency theory for groups} & {\bf Vassiliev theory}  \\
\hline\hline

a group $G$ &  a class of tangles \lis \\
&with a fixed skeleton $\boldX$ \\
\hline $A(G)=\oplus J^kG/J^{k-1}G$ & diagram space $\F(\boldX)$\lis \\
\hline functions $\Z G\to\Ring$  & $\Ring$-valued Vassiliev
invariants \lis \\
that vanish on $J^{n+1}G$&  of order $n$\\
\hline Chen expansion & Kontsevich integral\lis \\
\hline
dimension series $\Di_n G$ & filtration by $n$-trivial tangles\lis \\
\hline lower central series $\gamma_nG$ & filtration by
$\gamma_n$-trivial tangles \lis \\
\hline
\end{tabular}
\end{center}
The notion of $\gamma_n$-triviality (that is, $\gamma_n$-equivalence
to the trivial tangle) that appears in the last line will be
discussed later in this chapter, for string links rather than for
general tangles. Note that we do not have a general definition for
the trivial tangle with a given skeleton $\boldX$, so in the last
two lines we should restrict our attention to knots or (string)
links.

The above dictionary must be used with certain care, as illustrated
in the following paragraph.

\subsection{Invariants for the full braid group}

The finite type invariants for braids, considered as tangles, are
defined separately for each permutation. The set of braids on $m$
strands corresponding to the same permutation is in one-to-one
(non-canonical) correspondence with the pure braid group $P_m$:
given a braid $b$ the subset $bP_m\subset B_m$ consists of all the
braids with the same permutation as $b$. This correspondence also
identifies the Vassiliev invariants for $P_m$ with those of $bP_m$.
In particular, the Vassiliev invariants separate braids.

On the other hand, the dimension series for the full braid group
contains very little information. Indeed, it is known from \cite{GL}
that for $m\geq 5$ the lower central series of $B_m$ stabilizes at
$k=2$:
$$\gamma_kB_m=\gamma_2B_m$$
for $k\geq 2$.
\begin{xxca}
Show that for all $m$ the quotient $B_m/\gamma_2B_m$ is an infinite
cyclic group.
\end{xxca}

\section{String links as closures of pure braids}

The Vassiliev invariants for pure braids can be used to prove some
facts about the invariants of knots, and, more generally, string
links.

\subsection{The short-circuit closure.}

String links can be obtained from pure braids by a procedure called
{\em short-circuit closure.}\index{Short-circuit closure}
Essentially, it is a modification of the {\em plat closure}
construction described in \cite{Bir2}.

In the simplest case when string links have one component, the
short-circuit closure produces a long knot out of a pure braid on an
odd number of strands by joining the endpoints of the strands in
turn at the bottom and at the top:
$$\ig[width=8cm]{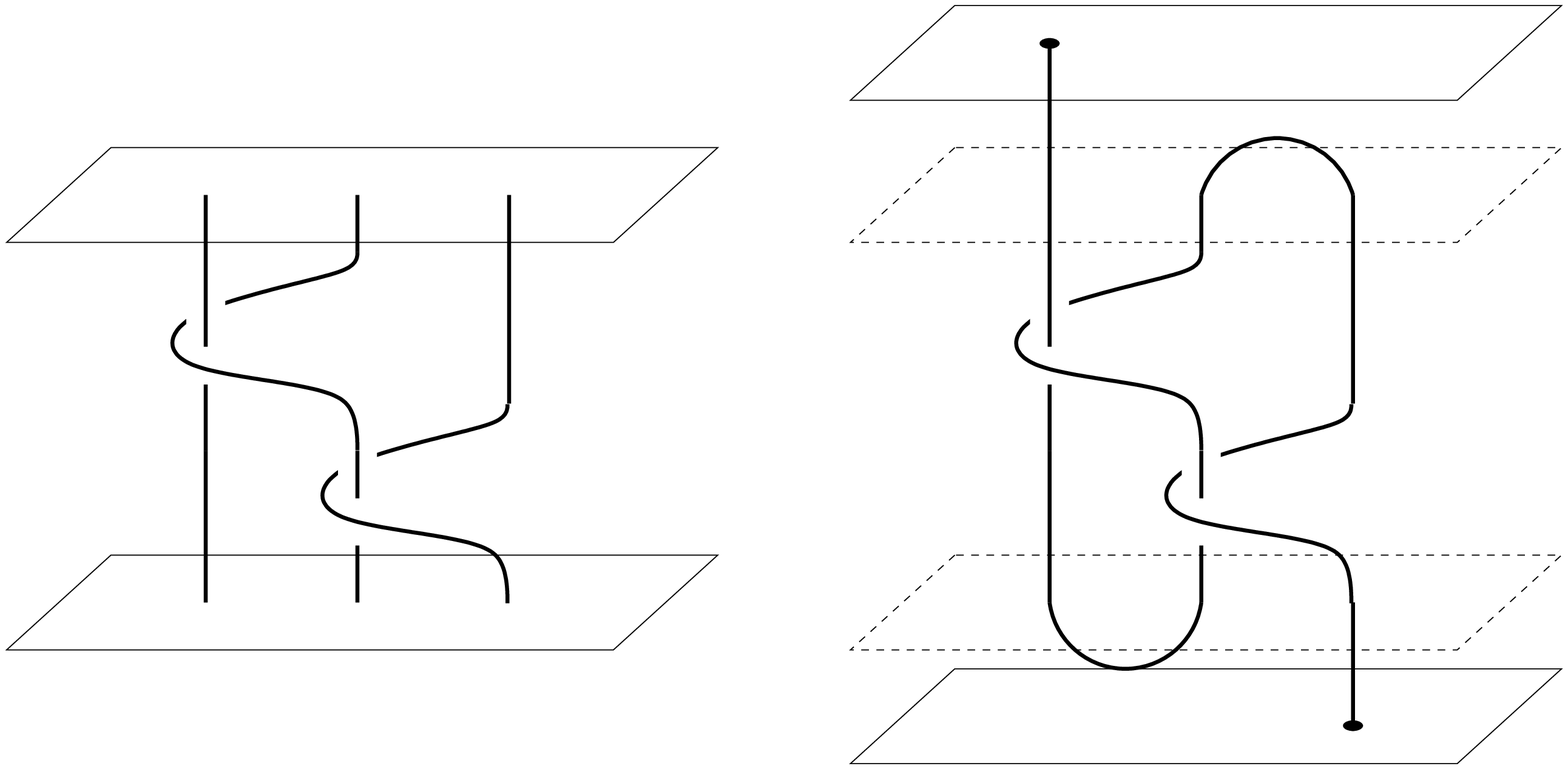}$$
In order to get a string link with $m$ components we have to start
with a pure braid on $(2k+1)m$ strands and proceed as follows.

Draw a braid in such a way that its top and bottom consist of the
integer points of the rectangle $[1,m]\times[0,2k]$ in the plane. A
string link on $m$ strands can be obtained from such a braid by
joining the points $(i,2j-1)$ and $(i,2j)$ (with $0<j\leq k$) in the
top plane and $(i,2j)$ and $(i,2j+1)$ (with $0\leq j< k$) in the
bottom plane by little arcs, and extending the strands at the points
$(i,0)$ in the top plane and $(i,2k)$ in the bottom plane. Here is
an example with $m=2$ and $k=1$:

$$\ig[width=8cm]{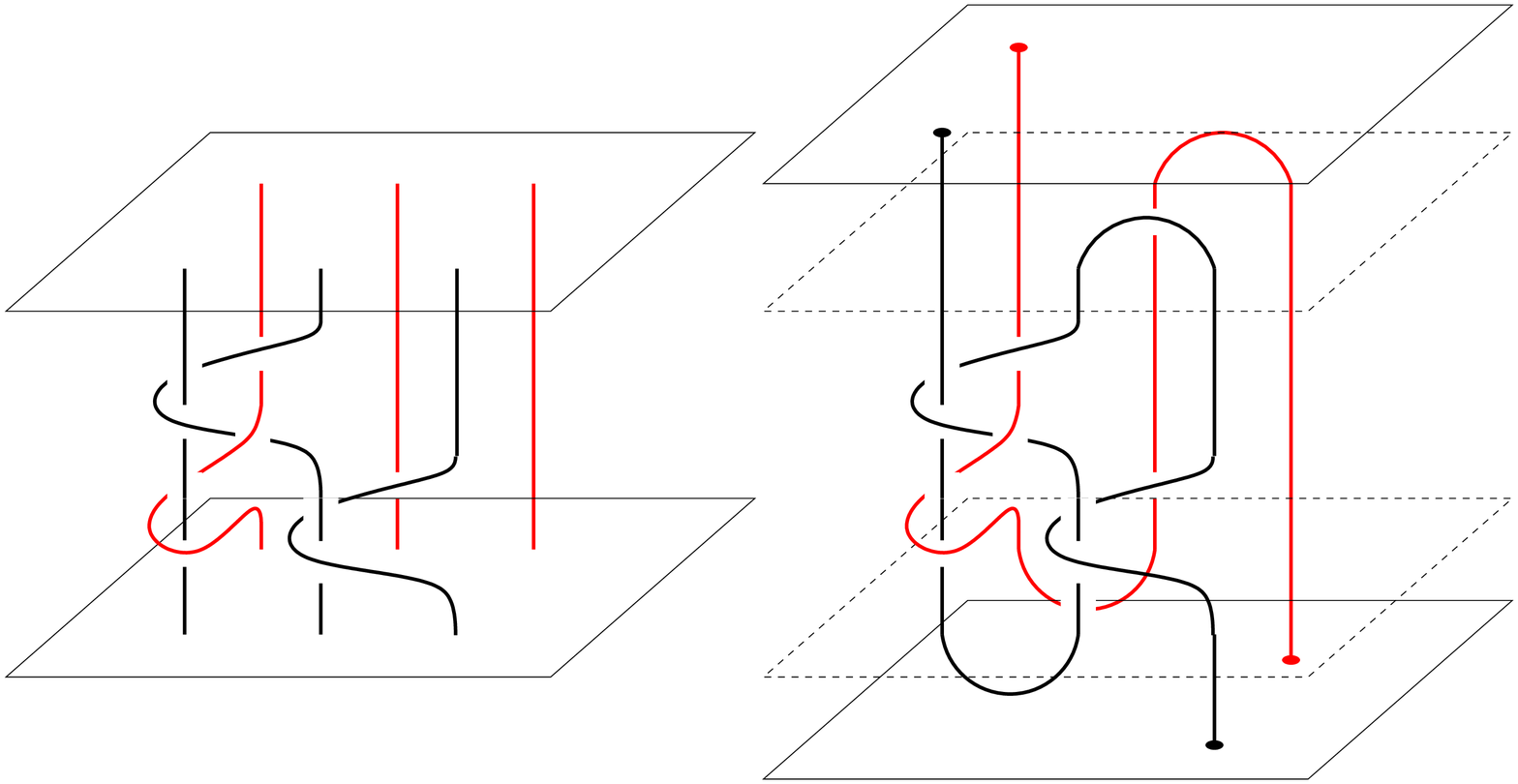}$$

\begin{figure}[htb]
$$\ig[width=10cm]{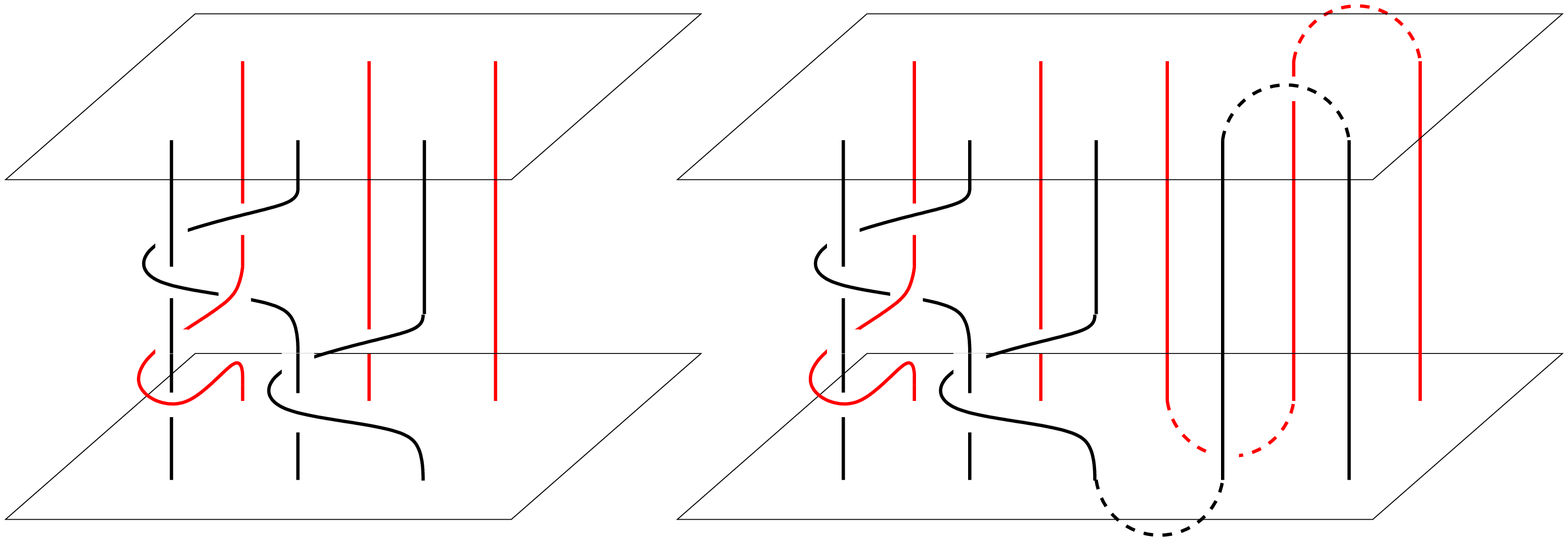}$$
\caption{The stabilization map.} \label{figure:stab}
\end{figure}

The short-circuit closure can be thought of as a map $\mathcal
S_{k}$ from the pure braid group $P_{(2k+1)m}$ to the monoid
$\Links_{m}$ of string links on $m$ strands. This map is compatible
with the {\em stabilization}, which consists of adding $2m$
unbraided strands to the braid on the right, as in
Figure~\ref{figure:stab}.

Therefore, if $P_{\infty}$ denotes the union of the groups
$P_{(2k+1)m}$ with respect to the inclusions $P_{(2k+1)m}\to
P_{(2k+3)m}$, there is a map
$$\mathcal S:P_{\infty}\to \Links_{m}. $$
The map $\mathcal S$ is onto, while $\mathcal S_{k}$, for any finite
$k$, is not\footnote{To show this one has to use the {\em bridge
number} (see page \pageref{ex:bridge-n}) of knots.}.

One can say when two braids in $P_{\infty}$ give the same string
link after the short-circuit closure:
\begin{theorem}\label{theorem:shortcirc}
There exist two subgroups $H^T$ and $H^B$ of $P_{\infty}$ such that
the map $\mathcal S_{n}$ is constant on the double cosets of the
form $H^T x H^B$. The preimage of every string link is a coset of
this form.
\end{theorem}

In other words, $\Links_{m}=H^T\backslash P_{\infty}/ H^B$.

Theorem~\ref{theorem:shortcirc} generalizes a similar statement for
knots (the case $m=1$), which was proved for the first time by J.\
Birman in \cite{Bir2} in the setting of the plat closure. Below we
sketch a proof which closely follows the argument given for knots in
\cite{MSt}.

First, notice that the short-circuit closure of a braid in
$P_{(2k+1)m}$ is not just a string link, but a {\em Morse} string
link: the height in the 3-space is a function on the link with a
finite number of isolated critical points, none of which is on the
boundary. We shall say that two Morse string links are {\em Morse
equivalent} if one of them can be deformed into the other through
Morse string links.
\begin{xlemma}
Assume that the short-circuit closures of  $b_1,b_2\in P_{(2k+1)m}$
are isotopic. There exist $k'\geq k$ such that the short-circuit
closures of the images of $b_1$ and $b_2$ in $P_{(2k'+1)m}$ under
the (iterated) stabilization map are Morse equivalent.
\end{xlemma}
The proof of this Lemma is not difficult; it is identical to the
proof of Lemma~4 in \cite{MSt} and we omit it.

Let us now describe the groups $H^T$ and $H^B$. The group $H^T$ is
generated by elements of two kinds. For each pair of strands joined
on top by the short-circuit map take (a) the full twist of this pair
of strands (b) the braid obtained by taking this pair of strands
around some strand, as in Figure~\ref{figure:ht}:
\begin{figure}[ht]
$$\ig[width=7cm]{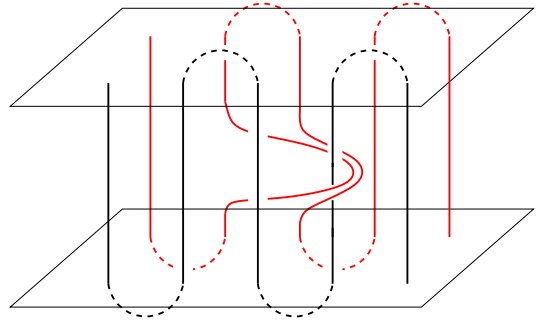}$$
\caption{A generator of $H^T$.}\label{figure:ht}
\end{figure}

The group $H^B$ is defined similarly, but instead of pairs of
strands joined on top we consider those joined at the bottom.
Clearly, multiplying a braid $x$ on the left by an element of $H^T$
and on the right by an element of $H^B$ does not change the string
link $\mathcal{S}(x)$.

Now, given a Morse string link with the same numbers of maxima of
the height function on each component (say, $k$), we can reconstruct
a braid whose short-circuit closure it is, as follows.

\begin{figure}[ht]
$$\ig[width=11cm]{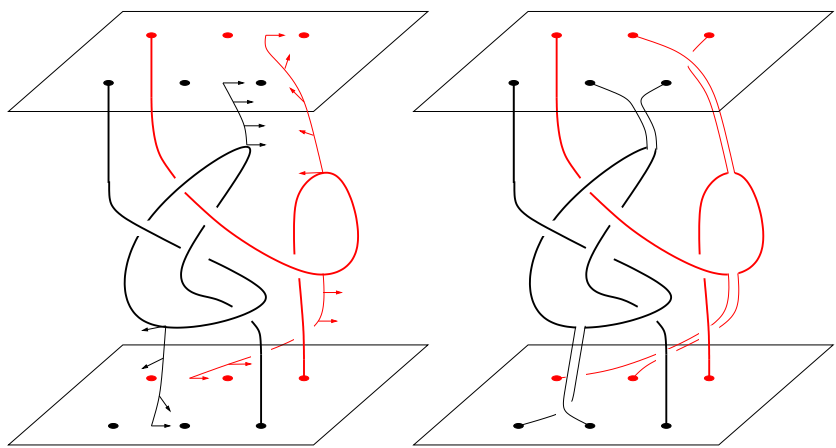}$$
\caption{Obtaining a braid from a string link.}
\label{figure:suspend}
\end{figure}

Suppose that the string link is situated between the top and the
bottom planes of the braid. Without loss of generality we can also
assume that the top point of $i$th strand is the point $(i,0)$ in
the top plane and the bottom point of the same strand is $(i,2k)$ in
the bottom plane. For the $j$th maximum on the $i$th strand, choose
an ascending curve that joins it with the point $(i,2j-1/2)$ in the
top plane, and for the $j$th minimum choose a descending curve
joining it to the point $(i,2k-3/2)$ in the bottom plane. We choose
the curves in such a way that they are all disjoint from each other
and only have common points with the string link at the
corresponding maxima and minima. On each of these curves choose a
framing that is tangent to the knot at one end and is equal to
$(1,0,0)$ at the other end. Then, doubling each of this curves in
the direction of its framing, we obtain a braid as in
Figure~\ref{figure:suspend}.

Each braid representing a given string link can be obtained in this
way. Given two Morse equivalent string links decorated with systems
of framed curves, there exists a deformation of one string link into
the other through Morse links. It extends to a deformation of the
systems of framed curves if we allow a finite number of transversal
intersections of curves with each other or with the string link, all
at distinct values of the parameter of the deformation, and changes
of framing. When a system of framed curves passes such a
singularity, the braid that it represents changes. A change of
framing on a curve ascending from a maximum produces the
multiplication on the left by some power of the twist on the pair of
strands corresponding to the curve. An intersection of the curve
ascending from a maximum with the link or with another curve gives
the multiplication on the left by a braid in $H^T$ obtained by
taking the pair of strands corresponding to the curve around some
other strands. Similarly, singularities involving a curve descending
from a minimum produce multiplications on the right by elements of
$H^B$.

\begin{remark}\label{remark:ht}
The subgroups $H^T$ and $H^B$ can be described in the following
terms. The short-circuit map $\mathcal{S}$  can be thought of as
consisting of two independent steps: joining the top ends of the
strands and joining the bottom ends. A braid belongs to $H^T$ if and
only if the tangle obtained from it after joining the top strands
only is ``trivial'', that is, equivalent to the tangle obtained in
this way from the trivial braid. The subgroup $H^B$ is described in
the same way.
\end{remark}

\subsection{Vassiliev knot invariants as pure braid invariants}
A knot invariant $v$ gives rise to a pure braid invariant $v\circ
\SS$ which is just the pull-back of $v$ with respect to the
short-circuit map. It is clear that if $v$ is of order $n$ the same
is true for $v\circ \SS$ since the short-circuit map sends braids
with double points to singular knots with the same number of double
points.

An example is provided by the Conway polynomial. Each of its
coefficients gives rise to an invariants of pure braids; these
invariants factor through the Magnus expansion since the latter is
the universal Vassiliev invariant. As a result, we get a function on
the chord diagram algebra $\A(2m+1)$ which can be explicitly
described, at least for $m=1$.

Recall that the algebra $\A(3)$ has a basis consisting of diagrams
of the form $w(u_{13},u_{23})\cdot u_{12}^m$ where the $u_{ij}$ are
the generators (horizontal chords connecting the strands $i$ and
$j$) and $w$ is some non-commutative monomial in two variables. Let
$$\chi:\A(3)\to\Z[t]$$ be the map such that for all $x\in P_3$ the
Conway polynomial of $\SS(x)$ coincides with $\chi(\M(x))$. The
following description of $\chi$ is given in \cite{Du4}.

First, it can be shown that $\chi$ vanishes on all the basis
elements of the form $uu_{12}$ and $u_{23} u$ for any $u$, and on
all $u u_{23}^2 u'$ for any $u$ and $u'$. This leaves us with just
two kinds of basis elements: $$[c_1,\ldots, c_k]:=u_{13}^{c_1}
u_{23}\ldots u_{13}^{c_{k-1}}u_{23} u_{13}^{c_k}$$ and
$$[c_1,\ldots, c_k]':=u_{13}^{c_1}u_{23}\cdot\ldots\cdot u_{23}
u_{13}^{c_{k-1}}u_{23}u_{13}^{c_k}u_{23}.$$ The values of $\chi$ on
the elements of the second kind are expressed via those on the
elements of the first kind:
$$\chi([c_1,\ldots, c_k]')=t^{-2}\cdot\chi([c_1,\ldots, c_k, 1]).$$
As for the elements of the first kind, we have
$$\chi([c_1,\ldots , c_k]) = (-1)^{k-1}\left(\prod_{i=1}^{k-1} p_1
p_{c_{i-1}}\right)\cdot p_{c_k},$$ where $p_s = \chi([s])$ is a
sequence of polynomials in $t$ that are defined recursively by $p_0
= 1,$ $p_1 = t^2$ and $p_{s+2} = t^2(p_s + p_{s+1})$ for $s\geq 0$.

\section{Goussarov groups of knots}
There are several facts about the Vassiliev string link invariants
that can be proved by studying the interaction between the
short-circuit closure and the dimension/lower central series for the
pure braid groups. (In view of the results in
Section~\ref{subsection:semidirect} these two series on $P_m$ always
coincide.) In this section we shall consider the case of knots which
is slightly simpler than the general case of string links.

\begin{xdefinition}
Two knots $K_1$ and $K_2$ are {\em $\gamma_n$-equivalent} if there
are $x_1, x_2\in P_{\infty}$ such that $K_i=\mathcal S(x_i)$ and
$x_1x_2^{-1}\in \gamma_{n} P_{\infty}$.
\end{xdefinition}
\begin{xxca}
Show that the connected sum of knots descends to their
$\gamma_n$-equivalence classes.
\end{xxca}
\begin{theorem}[\cite{G1, Ha2}]\label{theorem:groupknots}
For each $n$, the $\gamma_n$-equivalence classes of knots form an
abelian group under the connected sum.
\end{theorem}

The group of knots modulo $\gamma_{n+1}$-equivalence is called the
{\em $n$th Goussarov group} and is denoted by $\K(n)$.
\begin{theorem}\label{theorem:Vassiliev}
Two knots cannot be distinguished by Vassiliev invariants (with
values in any abelian group) of degree at most $n$ if and only if
they define the same element in $\K(n)$.
\end{theorem}
In other words, two knots are $\gamma_{n+1}$-equivalent if and only
if they are $n$-equivalent (see Section~\ref{sing_knot_filtr}).

The rest of this section is dedicated to the proof of Theorems
\ref{theorem:groupknots} and \ref{theorem:Vassiliev}. The main idea
behind the proof of Theorem \ref{theorem:Vassiliev}, which is due to
T.~Stanford \cite{Sta3}, is to interpret knot invariants as pure
braid invariants.

\subsection{The shifting endomorphisms}

For $k>0$, define $\tau_k$ to be the endomorphism of $P_{\infty}$
which replaces the $k$th strand by three parallel copies of itself
as in Figure~\ref{figure:tripleknots}:

\begin{figure}[ht]
$$\ig[width=9cm]{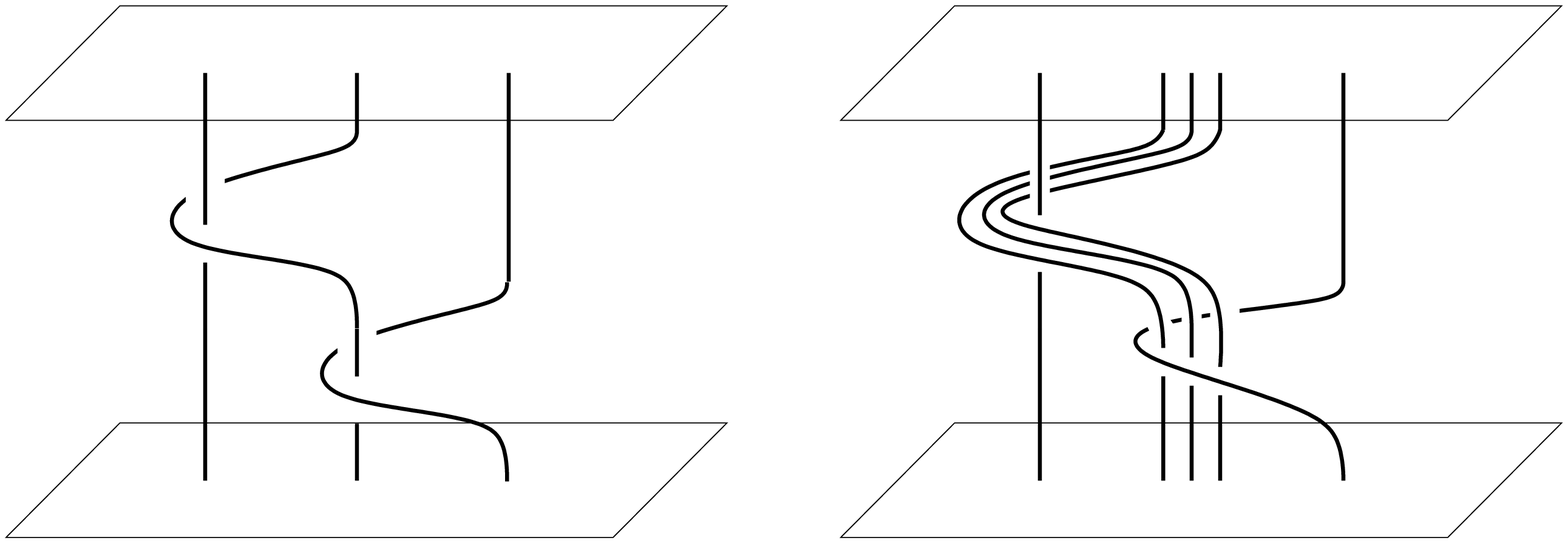}$$
\caption{} \label{figure:tripleknots}
\end{figure}
Denote by $\tau_0$ the endomorphism of $P_{\infty}$ which adds $2$
non-interacting strands to the left of the braid (this is in
contrast to the stabilization map, which adds $2$ strands to the
right and is defined only for $P_{2k+1}$ with finite $k$).

Strand-tripling preserves both $H^T$ and $H^B$. Also, since $\tau_k$
is an endomorphism, it respects the lower central series of
$P_{\infty}$.

\begin{xlemma}\cite{CMS}\label{lemma:shift}
For any $n$ and any $x\in\gamma_n P_{2N-1}$ there exist $t\in
H^T\cap \gamma_n P_{2N+1}$ and $b\in H^B\cap \gamma_n P_{2N+1}$ such
that $\tau_0{(x)}=txb$.
\end{xlemma}

\begin{proof}
Let $t_{2k-1}=\tau_{2k-1}(x)(\tau_{2k}(x))^{-1}$, and let
$b_{2k}=(\tau_{2k+1}(x))^{-1}\tau_{2k}(x)$. Notice that
$t_{2k-1},b_{2k}\in \gamma_n P_{\infty}$. Moreover, $t_{2k-1}$ looks
as in Figure~\ref{figure:t} and, by the Remark~\ref{remark:ht}, lies
in $H^T$. Similarly, $b_{2k}\in H^B$. We have
$$ \tau_{2k-1}(x)=t_{2k-1}\tau_{2k}(x),$$
$$\tau_{2k}(x)=\tau_{2k+1}(x)b_{2k}.$$
There exists $N$ such that $\tau_{2N+1}(x)=x$. Thus the following
equality holds:
$$\tau_0(x)=t_{1}\cdots t_{2N-1}xb_{2N}\cdots b_{0},$$
and this completes the proof.
\begin{figure}[hb]
$$\ig[width=6cm]{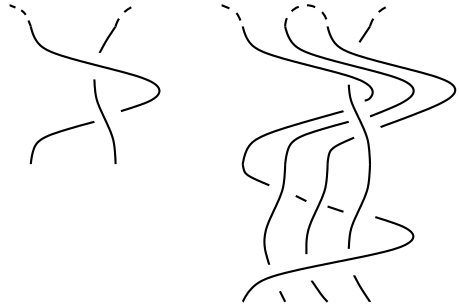}$$
\caption{Braids $x$ and $t_{2k-1}$.} \label{figure:t}
\end{figure}
\end{proof}

\subsection{Existence of inverses}

Theorem~\ref{theorem:groupknots} is a consequence of the following,
stronger, statement:
\begin{xproposition}\label{prop:main}
For any $x\in\gamma_k P_{2N-1}$ and any $n$ there exists $y\in
\gamma_k P_{\infty}$ such that:
\begin{itemize}
\item $y$ is contained in the image of $\tau_0^N$;
\item $xy=thb$ with $h\in\gamma_n P_{\infty}$ and $t,b\in\gamma_k P_{\infty}$.
\end{itemize}
\end{xproposition}
The first condition implies that $\SS(xy)=\SS(x)\#\SS(y)$. It
follows from the second condition that the class of $\SS(y)$ is the
inverse for $\SS(x)$. The fact that $t$ and $b$ lie in $\gamma_k
P_{\infty}$ is not needed for the proof of
Theorem~\ref{theorem:group}, but will be useful for
Theorem~\ref{theorem:nilpotency}.

\begin{proof}
Fix $n$. For $k\geq n$ there is nothing to prove.

Assume there exist braids for which the statement of the proposition
fails; among such braids choose $x$ with the maximal possible value
of $k$. By Lemma~\ref{lemma:shift} we have $\tau_0^N
(x^{-1})=t_1x^{-1}b_1$ with $t_1\in H^T\cap\gamma_k P_{4N-1}$ and
$b_1\in H^B\cap\gamma_k P_{4N-1}$. Then
$$x\tau_0^N (x^{-1})=xt_1x^{-1}b_1=t_1\cdot t_1^{-1}xt_1x^{-1}\cdot b_1.$$
Since $t_1^{-1}xt_1x^{-1}\in\gamma_{k+1}P_{4N-1}$, there exists
$y'\in\gamma_{k+1}P_{\infty}\cap\text{Im}\tau_0^{2N}$ such that
$t_1^{-1}xt_1x^{-1}\cdot y'=t_2hb_2$ where $h\in\gamma_{n}
P_{\infty}$, $t_2\in H^T\cap\gamma_{k+1} P_{\infty}$ and $b_2\in
H^B\cap\gamma_{k+1} P_{\infty}$. Note that $y'$ commutes with $b_1$,
and, hence,
$$x\cdot \tau_0^N(x^{-1})y' = t_1t_2\cdot h\cdot b_2b_1.$$
Setting $y=\tau_0^N(x^{-1})y'$, $t=t_1t_2$  and $b=b_2b_1$ we see
that for $x$ the statement of the proposition is satisfied. We get a
contradiction, and the proposition is proved.
\end{proof}

\subsection{Vassiliev invariants and $\gamma_n$-equivalence}
If we consider knot invariants as braid invariants via the
short-circuit closure, it becomes clear that the value of a knot
invariant of order $n$ or less only depends on the
$\gamma_{n+1}$-equivalence class of the knot. Indeed, multiplying a
pure braid by an element of $\gamma_{n+1}P_m$ amounts to adding an
element of $J^{n+1}P_m$, and this does not affect the invariants of
degree $n$ or less.

\begin{xlemma}
The map $$\K\to \K(n)$$ that sends a knot to its
$\gamma_{n+1}$-equivalence class is an invariant of degree $n$.
\end{xlemma}

This lemma establishes Theorem \ref{theorem:Vassiliev} since it
tautologically implies that Vassiliev invariants of degree at most
$n$ distinguish $\gamma_{n+1}$-equivalence classes of knots.

\begin{proof}
Extend the map $\K\to \K(n)$ by linearity to a homomorphism of
abelian groups $\Z\K\to \K(n)$. The kernel of this map is spanned by
two types of elements:
\begin{itemize}
\item elements of the form $x-y$ where $x$ and $y$ are
$\gamma_{n+1}$-equivalent;
\item elements of the form $x_1\#
x_2-x_1-x_2$.
\end{itemize}
Note that the trivial knot is in the kernel since, up to sign, it is
an element of the second type. The subspace of elements of the
second type in $\Z\K$ coincides with $\K_1\#\K_1$, where $\K_1$ is
the ideal of singular knots.

We need to show that the composite map $$\Z
P_{\infty}\stackrel{\SS}{\longrightarrow} \Z\K\to \K(n)$$ sends
$J^{n+1} P_{\infty}$ to zero.

Define a {\em relator of order $d$ and length $s$} as an element of
$\Z\K$ of the form
$$\SS((x_1-1)(x_2-1)\ldots(x_s-1)y) \label{eq:product}$$
with $y\in P_{\infty}$, $x_i\in\gamma_{d_i}P_{\infty}$ and $\sum d_i
= d$. The greatest $d$ such that a relator is of order $d$ will be
called the {\em exact order}  of a relator. A {\em composite
relator} is an element of $\K_1\#\K_1\subset\Z\K$.

As we noted, the kernel of the map $\Z\K\to \K(n)$ contains all the
relators of length 1 and order $n+1$ and all the composite relators.
On the other hand, an element of $\SS{(J^{n+1}P_{\infty})}$ is a
linear combination of relators of length $n+1$ and, hence, of order
$n+1$. Thus we need to show that any relator of order $n+1$ is a
linear combination of relators of order $n+1$ and length $1$ and
composite relators.

Suppose that there exist relators of order $n+1$ which cannot be
represented as linear combinations of the above form. Among such
relators choose the relator $R$ of minimal length and, given the
length, of maximal exact order.

Assume that $R$ is of the form $(\ref{eq:product})$ as above, with
$y, x_i\in P_{2N-1}$. Choose $t\in H^T$ and $b\in H^B$ such that the
braid $tx_{1}b$ coincides with the braid obtained from $x_{1}$ by
shifting it  by $2N$ strands to the right, that is, with
$\tau_0^N(x_1)$. By Lemma~\ref{lemma:shift} the braids $t$ and $b$
can be taken to belong to the same term of the lower central series
of $P_{\infty}$ as the braid $x$. The relator
\[ R'=\SS((tx_{1}b-1)(x_{2}-1)\ldots(x_{s}-1)y) \]
is a connected sum of two relators and, hence, is a combination of
composite relators. On the other hand,
\[\begin{array}{lcl}
R'-R & = & \SS((tx_{1}b-x_{1})(x_{2}-1)\ldots(x_{s}-1)y)\\
     & = & \SS(x_{1}(b-1)(x_{2}-1)\ldots(x_s-1)y)
\end{array}\]
Notice now that $(b-1)$ can be exchanged with $(x_i-1)$ and $y$
modulo relators of shorter length or higher order. Indeed,
\[(b-1)y=y(b-1)+([b,y]-1)yb\]
and
\[(b-1)(x_i-1)=(x_i-1)(b-1)+ ([b,x_i]-1)(x_{i}b-1)+([b,x_{i}]-1).\]
Thus, modulo  relators of shorter length or higher order
\[\SS(x_{1}(b-1)(x_{2}-1)\ldots(x_{m}-1)y)=
\SS(x_{1}(x_{2}-1)\ldots(x_{m}-1)y(b-1))=0.\] and this means that
$R$ is a linear combination of composite relators and relators of
length 1 and order $n$.
\end{proof}

\section{Goussarov groups of string links}

Much of what was said about the Goussarov groups of knots can be
extended to string links without change. Just as in the case of
knots, two string links $L_1$ and $L_2$ are said to be {\em
$\gamma_n$-equivalent} if there are $x_1, x_2\in P_{\infty}$ such
that $L_i=\mathcal S(x_i)$ and $x_1x_2^{-1}\in \gamma_{n}
P_{\infty}$. A string link is {\em $\gamma_n$-trivial} if it is
$\gamma_n$-equivalent to the trivial link. The product of string
links descends to their $\gamma_n$-equivalence classes.

\begin{theorem}[\cite{G1, Ha2}]\label{theorem:group}
For each $m$ and $n$, the $\gamma_n$-equivalence classes of string
links on $m$ strands form a group under the string link product.
\end{theorem}

These groups are also referred to as {\em Goussarov groups}. We
shall denote the group of string links on $m$ strands modulo
$\gamma_{n+1}$-equivalence by $\Links_{m}(n)$, or simply by
$\Links(n)$, dropping the reference to the number of strands. Let
$\Links(n)_k$ be the subgroup of $\Links(n)$ consisting of the
classes of $k$-trivial links. Note that $\Links(n)_k=1$ for $k>n$.

For string links with more than one component the Goussarov groups
need not be abelian. The most we can say is the following
\begin{theorem}[\cite{G1,Ha2}]\label{theorem:nilpotency}
For all $p,q$ we have
$$[\Links(n)_{p},\Links(n)_q]\subset \Links(n)_{p+q}.$$
In particular, $\Links(n)$ is nilpotent of nilpotency class at most
$n$.
\end{theorem}

As for the relation between $\gamma_{n+1}$-equivalence and
$n$-equivalence for string links, it is not known whether Theorem
\ref{theorem:Vassiliev} is valid for string links in the same form
as for knots. We shall prove a weaker statement:

\begin{theorem}[\cite{Mas}]\label{theorem:VassilievSL}
Two string links cannot be distinguished by $\Q$-valued Vassiliev
invariants of degree $n$ and smaller if and only if the elements
they define in $\Links(n)$ differ by an element of finite order.
\end{theorem}

We refer the reader to \cite{Mas} for further results.

The proof of Theorem~\ref{theorem:group} coincides with the proof of
Theorem~\ref{theorem:groupknots} word for word. The only
modification necessary is in the definition of the shifting
endomorphisms: rather than tripling the $k$th strand, $\tau_k$
triples the $k$th {\em row of strands}. In other words, $\tau_k$
replaces each strand with ends at the points $(i,k-1)$ in the top
and bottom planes, with $1\leq i\leq n$, by three parallel copies of
itself as in Figure~\ref{figure:triple}:
\begin{figure}[ht]
$$\ig[width=9cm]{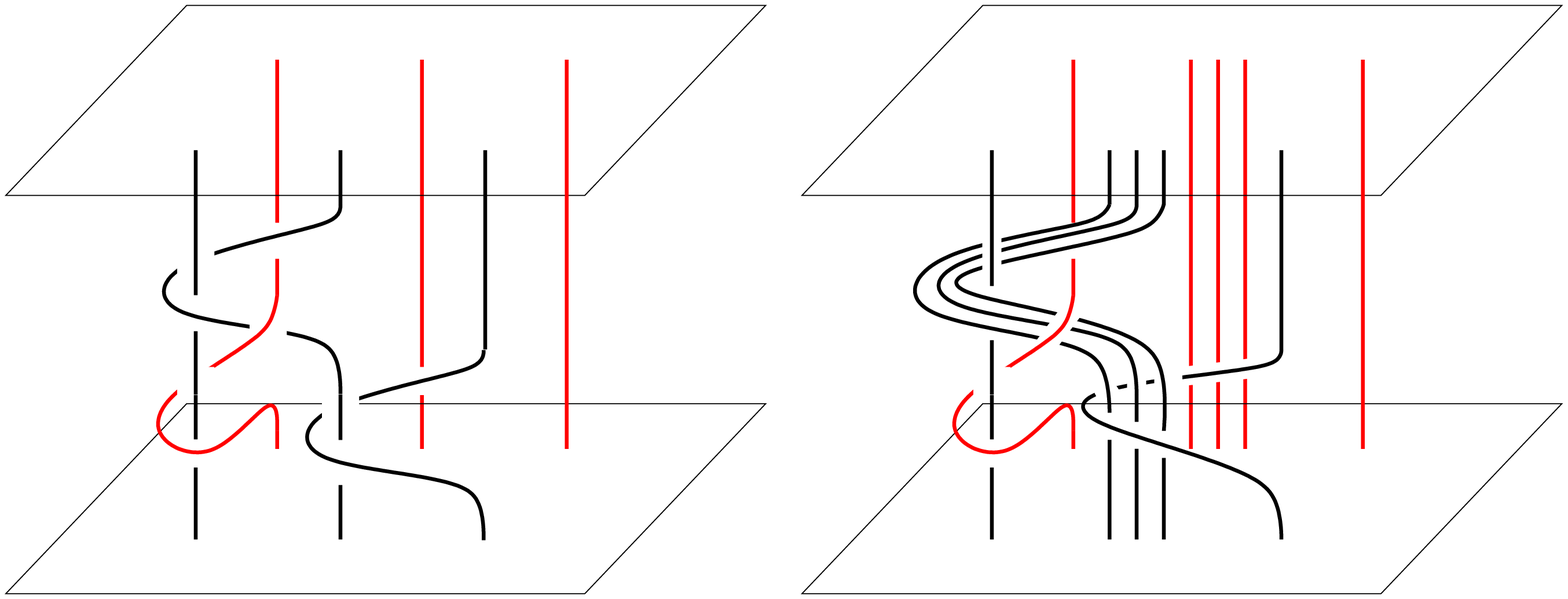}$$
\caption{} \label{figure:triple}
\end{figure}
Similarly, $\tau_0$ adds $2m$ non-interacting strands, arranged in 2
rows, to the left of the braid.

\subsection{The nilpotency of $\Links(n)$}

Let $x\in\gamma_p P_{\infty}$ and $x'\in\gamma_q P_{\infty}$. Choose
the braids $y$ and $y'$ representing the inverses in $\Links(n)$ of
$x$ and $x'$, respectively, such that the conditions of
Proposition~\ref{prop:main} are satisfied, with $n$ replaced by
$n+1$: $xy=t_1h_1b_1$ and $x'y'=t_2h_2b_2$ with $h_i\in\gamma_{n+1}
P_{\infty}$, $t_1,b_1\in\gamma_p P_{\infty}$ and $t_2,b_2\in\gamma_q
P_{\infty}$. Replacing the braids by their iterated shifts to the
right, if necessary, we can achieve that the braids $x$, $x'$, $y$
and $y'$ all involve different blocks of strands, and, therefore,
commute with each other. Then
\begin{multline*}
\SS(x)\cdot\SS(x')\cdot\SS(y)\cdot\SS(y')=\SS(xx'yy')=\SS(xyx'y')\\
=\SS(t_1h_1b_1t_2h_2b_2)=\SS(h_1b_1t_2h_2).
\end{multline*}
The latter link is $n$-equivalent to $\SS(t_2^{-1}b_1t_2b_1^{-1})$
which lives in $\Links(n)_{p+q}$.

It follows that each $n$-fold (that is, involving $n+1$ terms)
commutator in $\Links(n)$ is trivial, which means that $\Links(n)$
is nilpotent of nilpotency class at most $n$.
Theorem~\ref{theorem:nilpotency} is proved.

\subsection{Vassiliev invariants and $\gamma_n$-equivalence}

As in the case of knots, 
the value of any order $n$ Vassiliev invariant on a string link
depends only on the $\gamma_{n+1}$-equivalence class of the link. 
The following proposition is the key to determining when two
different $\gamma_{n+1}$-equivalence classes of string links cannot
be distinguished by Vassiliev invariants of order $n$:

\begin{proposition}\label{prop:relators}
The filtration by the powers of the augmentation ideal $J
P_{\infty}\subset \Q P_{\infty}$ is carried by short-circuit map to
the canonical filtration $\{E_i\Links(n)\}$ of the group algebra $\Q
\Links(n)$, induced by $\{\Links(n)_i\}$.
\end{proposition}
We remind 
that the canonical filtration was defined in Section
\ref{subsec:canonicalfiltration}.
\begin{proof}
We use induction on the power $k$ of  $J P_{\infty}$. For $k=1$
there is nothing to prove.

Any product of the form
$$(x_1-1)(x_2-1)\ldots(x_s-1)y\leqno{(\ast)}$$
with $y\in P_{\infty}$, $x_i\in\gamma_{d_i} P_{\infty}$ and $\sum
d_i = d$ belongs to $J^d P_{\infty}$ since for any $d_i$ we have
$\gamma_{d_i} P_{\infty}-1\subset J^{d_i}P_{\infty}$. We shall refer
to $s$ as the {\em length} of such product, and to $d$ as its {\em
degree}. The maximal $d$ such that a product of the form $(\ast)$ is
of degree $d$, will be referred to as the {\em exact degree} of the
product.

The short-circuit closure of a product of length 1 and degree $k$ is
in $E_k \Links(n)$. Assume there exists a product of the form
$(\ast)$ of degree $k$ whose image $R$ is not in $E_k \Links(n)$;
among such products choose one of minimal length, say $r$, and,
given the length, of maximal exact degree.

There exists $N$ such that
$$R':=\SS((\tau_0^N(x_1)-1)(x_2-1)\ldots(x_r-1)y)=\SS((x_2-1)\ldots(x_r-1)y)\cdot\SS(x_1-1).$$
The length of both factors on the right-hand side is smaller that
$k$, so, by the induction assumption, $R'\in E_k \Links(n)$.
If $\tau_0^N(x_1)=tx_1b$ we have

\[\begin{array}{lcl}
R'-R & = & \SS((tx_{1}b-x_{1})(x_{2}-1)\ldots(x_{m+1}-1)y)\\
     & = & \SS(x_{1}(b-1)(x_{2}-1)\ldots(x_{m+1}-1)y)
\end{array}\]
Notice now that $(b-1)$ can be exchanged with $(x_i-1)$ and $y$
modulo closures of products having shorter length or higher degree.
Indeed,
\[(b-1)y=y(b-1)+([b,y]-1)yb\]
and
\[(b-1)(x_i-1)=(x_i-1)(b-1)+ ([b,x_i]-1)(x_{i}b-1)+([b,x_{i}]-1).\]
Thus, modulo  elements of $E_k \Links(n)$
\[\SS(x_{1}(b-1)(x_{2}-1)\ldots(x_{m+1}-1)y)=
\SS(x_{1}(x_{2}-1)\ldots(x_{m+1}-1)y(b-1))=0.\]
\end{proof}

By Proposition~\ref{prop:relators} the elements of $\Links(n)$ that
cannot be distinguished from the trivial link by the Vassiliev
invariants of degree $n$ form the subgroup
$$\Links(n)\cap (1+E_{n+1}\Links(n)).$$
Since $\Links(n)_{n+1}=1$, by Theorems \ref{theorem:jennings} and
\ref{lemma:sqrt} this subgroup consists of all the elements of
finite order in $\Links(n)$. Finally, if the classes of two links
$L_1$ and $L_2$ cannot be distinguished by invariants of order $n$,
then $L_1-L_2\in E_{n+1}\Links(n)$, and, hence, $L_1 L_2^{-1}-1\in
E_{n+1}\Links(n)$ and $L_1 L_2^{-1}$ is  of finite order in
$\Links(n)$. Theorem~\ref{theorem:VassilievSL} is proved.

\subsection{Some comments}

\begin{xremark}
Rational-valued Vassiliev invariants separate pure braids, and the
Goussarov group of $\gamma_{n+1}$-equivalence classes of pure braids
on $k$ strands is nothing but $P_k/\gamma_{n+1} P_k$, which is
nilpotent of class $n$ for $k>2$. Since this group is a subgroup of
$\Links(n)$, we see that $\Links(n)$ is nilpotent of class $n$ for
links on at least 3 strands. String links on 1 strand are knots, in
this case $\Links(n)$ is abelian. The nilpotency class of
$\Links(n)$ for links on 2 strands is unknown. Note that it follows
from the results of \cite{DK} that $\Links(n)$ for links on 2
strands is, in general, non-abelian.
\end{xremark}

\begin{xremark}
The relation of the Goussarov groups of string links on more than
one strand to integer-valued invariants seems to be a much more
difficult problem. While in Proposition~\ref{prop:relators} the
field $\Q$ can be replaced by the integers with no changes in the
proof, Theorem~\ref{theorem:jennings} fails over $\Z$.
\end{xremark}

\begin{xremark}
Proposition~\ref{prop:relators} shows that the map
$$\Links_{m}\to \Links(n)\to\Q \Links(n)/E_{n+1}\Links(n)$$
is the {\em universal degree $n$ Vassiliev invariant} in the
following sense: each Vassiliev invariant of links in $\Links_{m}$
of degree $n$ can be extended uniquely to a linear function on $\Q
\Links(n)/E_{n+1}\Links(n)$.
\end{xremark}

\section{Braid invariants as string link invariants}

A pure braid is a string link so every finite type string link
invariant is also a braid invariant of the same order (at least). It
turns out that the converse is true:
\begin{theorem}\label{theorem:stringlinx}
A finite type integer-valued pure braid invariant extends to a
string link invariant of the same order.
\end{theorem}
\begin{xcorollary}
The natural map $\A^h(m)\to\A(m)$, where $\A^h(m)$ is the algebra of
the horizontal chord diagrams and $\A(m)$ is the algebra of all
string link chord diagrams, is injective.
\end{xcorollary}
This was first proved in \cite{BN8} by Bar-Natan. He considered
quantum invariants of pure braids, which all extend to string link
invariants, and showed that they span the space of all Vassiliev
braid invariants.

Our approach will be somewhat different. We shall define a map
$$\Links_m(n)\to P_{m}/\gamma_{n+1}P_{m}$$
from the Goussarov group of $\gamma_{n+1}$-equivalence classes of
string links to the group of $\gamma_{n+1}$-equivalence classes of
pure braids on $m$ strands, together with a section
$P_{m}/\gamma_{n+1}P_{m}\to\Links_{m}(n)$. A Vassiliev invariant $v$
of order $n$ for pure braids is just a function on
$P_{m}/\gamma_{n+1}P_{m}$, its pullback to $\Links_{m}(n)$ gives the
extension of $v$ to string links.

\begin{xremark}
Erasing one strand of a string link gives a homomorphism
$\Links_{m}\to\Links_{m-1}$, which has a section. If $\Links_{m}$
were a group, this would imply that string links can be combed, that
is, that $\Links_{m}$ splits as a semi-direct product of
$\Links_{m-1}$ with the kernel of the strand-erasing map. Of course,
$\Links_{m}$ is only a monoid, but it has many quotients that are
groups, and these all split as iterated semi-direct products. For
instance, string links form groups modulo concordance or link
homotopy \cite{HL}; here we are interested in the Goussarov groups.
\end{xremark}

Denote by $\FLinks_{m-1}(n)$ the kernel of the homomorphism
$\Links_{m}({n})\to\Links_{m-1}({n})$ induced by erasing the last
strand. We have semi-direct product decompositions
\[\Links_{m}({n})\cong \FLinks_{m-1}({n})\ltimes\ldots \FLinks_{2}({n})\ltimes \FLinks_{1}({n}).\]

We shall see that any element of $\FLinks_{k}({n})$ can be
represented by a string link on $k+1$ strands whose first $k$
strands are vertical. Moreover, taking the homotopy class of the
last strand in the complement of the first $k$ strands gives a
well-defined map
$$\pi_{k}:\FLinks_{k}({n})\to F_{k}/\gamma_{n+1}F_{k}.$$
Modulo the $n+1$st term of the lower central series, the pure braid
group has a semi-direct product decomposition
\[P_{\infty}/\gamma_{n+1}P_{\infty}\cong
F_{m-1}/\gamma_{n+1} F_{m-1}\ltimes\ldots \ltimes F_{1}/\gamma_{n+1}
F_{1}.\]
The homomorphisms $\pi_i$ with $i<m$ can now be assembled into one
surjective map
$$\Links_{m}({n})\to P_{m}/\gamma_{n+1}P_{m}.$$
Considering a braid as a string link gives a section of this map;
this will establish the theorem stated above as soon as we justify
the our claims about the groups $\FLinks_{k}({n})$. 

\subsection{String links with one non-trivial component.}
\label{subsection:likefree}

The fundamental group of the complement of a string link certainly
depends on the link. However, it turns out that all this dependence
is hidden in the intersection of all the lower central series
subgroups.

Let $X$ be a string link on $m$ strands and $\widetilde X$ be its
complement. The inclusion of the top plane of $X$, punctured at the
endpoints, into  $\widetilde X$ gives a homomorphism $i_t$ of $F_m$
into $\pi_1\widetilde X$.
\begin{xlemma}[\cite{HL}]
For any $n$ the homomorphism
$$F_m/\gamma_n F_m\to\pi_1\widetilde X/\gamma_{n}\pi_1\widetilde X$$
induced by $i_t$ is an isomorphism.
\end{xlemma}

A corollary of this lemma is that for any $n$ there is a
well-defined map
$$\Links_m(n)\to F_{m-1}/\gamma_{n+1} F_{m-1}$$
given by taking the homotopy class of the last strand of a string
link in the complement of the first $m-1$ strands. We must prove
that if two string links represent the same element of
$\FLinks_{m-1}(n)$, their images under this map coincide.

In terms of braid closures, erasing the last strand of a string link
corresponds to erasing all strands of $P_{\infty}$ with ends at the
points $(m,i)$ for all $i\geq 0$. Erasing these strands, we obtain
the group which we denote by $P_{\infty}^{m-1}$; write $\Phi$ for
the kernel of the erasing map. We have a semi-direct product
decomposition
$$P_{\infty}=\Phi\ltimes P_{\infty}^{m-1},$$
and the product is almost direct. In particular, this means that
$$\gamma_k P_{\infty}=\gamma_k \Phi\ltimes \gamma_k
P_{\infty}^{m-1}$$ for all $k$.

\begin{xlemma}
Let $x\in\Phi$, and $h\in \gamma_{n+1}\Phi$. The string links
$\SS(x)$ and $\SS(xh)$ define the same element of
$F_{m-1}/\gamma_{n+1} F_{m-1}$.
\end{xlemma}

\begin{proof}
Each braid in $\Phi$ can be combed: $\Phi$ is an almost direct
product of the free groups $G_i$ which consist of braids whose all
strands, apart from the one with the endpoints at $(m,i)$, are
straight, and whose strands with endpoints at $(m,j)$ with $j<i$ do
not interact. Each element $a$ of $G_i$ gives a path in the
complement of the first $m-1$ strands of the string link, and,
hence, an element $[a]$ of $F_{m-1}$. Notice that this
correspondence is a homomorphism of $G_i$ to $F_{m-1}$. (Strictly
speaking, these copies of $F_{m-1}$ for different $i$ are only
isomorphic, since these are fundamental groups of the same space
with different basepoints. To identify these groups we need a choice
of paths connecting the base points. Here we shall choose intervals
of straight lines.)

Given $x\in\Phi$ we can write it as $x_1x_2\ldots x_r$ with $x_i\in
G_i$. Then the homotopy class of the last strand of $\SS_n(x)$
produces the element $$[x_1][x_2]^{-1}\ldots [x_{r}]^{(-1)^{r-1}}\in
F_{m-1}.$$ Let $x'=xh$ with $h\in\gamma_{n+1}\Phi$. Then the fact
that $\Phi$ is an almost direct product of the $G_i$ implies that if
$x'=x'_1x'_2\ldots x'_r$ with $x_i\in G_i$, then $x_i\equiv x'_i\mod
\gamma_{n+1} G_i$. It follows that the elements of $F_{m-1}$ defined
by $\SS(x)$ and $\SS_n(x')$ differ by multiplication by an element
of $\gamma_{n+1} F_{m-1}$.
\end{proof}

\begin{xlemma}
Let $x\in\Phi$, and $y\in\gamma_{n+1} P_{\infty}^{m-1}$. The string
links $\SS(x)$ and $\SS(xy)$ define the same element of
$F_{m-1}/\gamma_{n+1} F_{m-1}$.
\end{xlemma}

\begin{proof}
Denote by $\widetilde X$ the complement of $\SS(y)$. We shall write
a presentation for the fundamental group of $\widetilde X$. It will
be clear from this presentation that the element of
$$F_{m-1}/\gamma_{n+1}
F_{m-1}=\pi_1\widetilde X /\gamma_{n+1}\pi_1\widetilde X$$ given by
the homotopy class of the last strand of $\SS(xy)$ does not depend
on $y$.

Let us assume that both $x$ and $y$ lie in the braid group
$P_{m(2N+1)}$. Let $H$ be the horizontal plane coinciding with the
top plane of the braid $y$. The plane $H$ cuts the space $\widetilde
X$ into the upper part $H_+$ and the lower part $H_-$. The
fundamental groups of $H_+$, $H_-$ and $H_+\cap H_-$ are free. Let
us denote by $\{\alpha_{i,j}\}$, $\{\beta_{i,j}\}$ y
$\{\gamma_{i,k}\}$ the corresponding free sets of generators (here
$1\leq i < m$, $1\leq j \leq N+1$ and $1\leq k\leq 2N+1$) as in
Figure~\ref{figure:freegen}.
\begin{figure}[htb]
$$\ig[width=9cm]{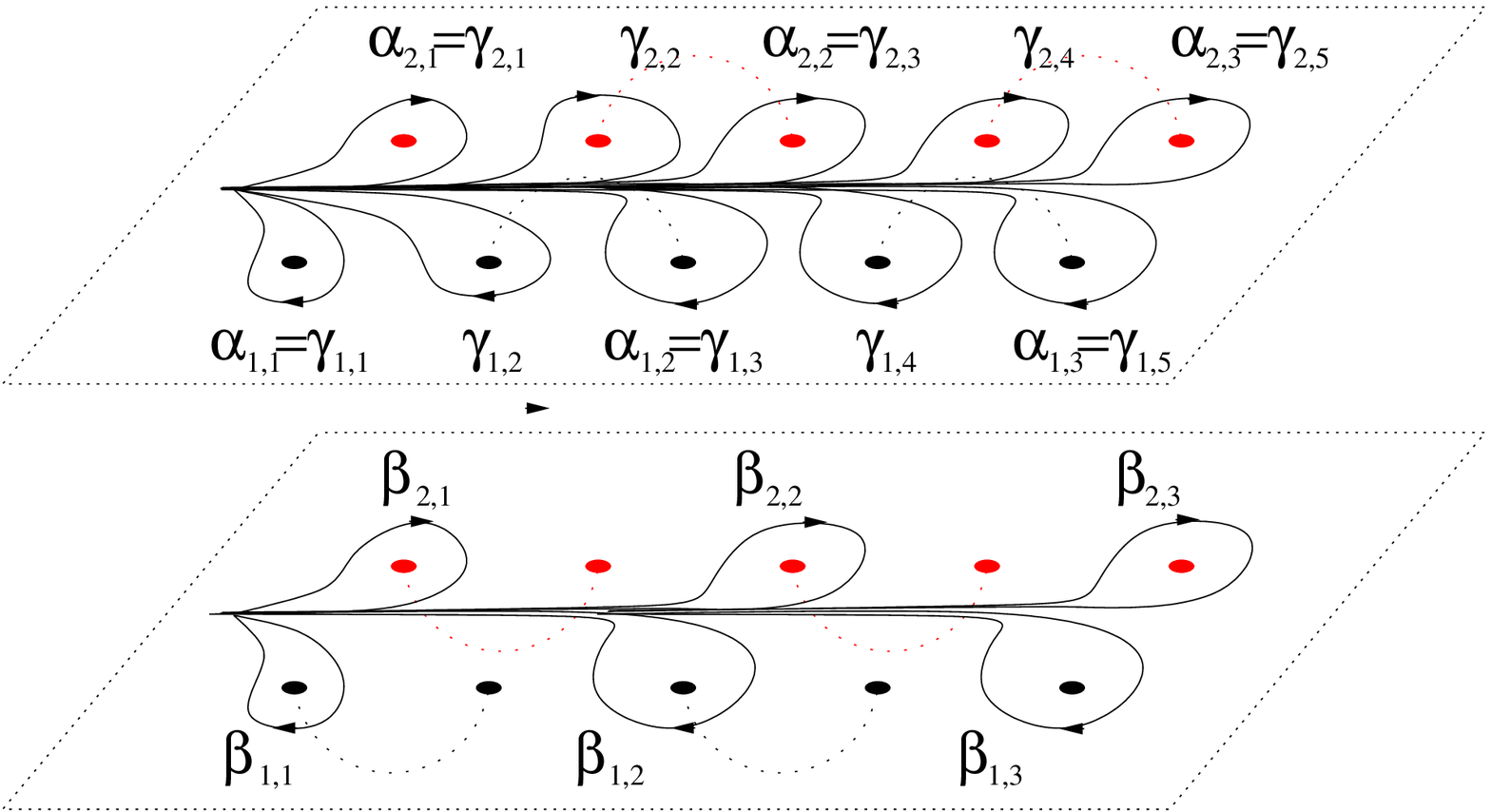}$$
\caption{} \label{figure:freegen}
\end{figure}
By the Van Kampen Theorem, $\pi_1 \widetilde X$ has a presentation
$$
\begin{array}{rrrl}
\langle \alpha_{i,j}, \beta_{i,j}, \gamma_{i,k}\ |\ &
\theta_{y}^{-1}(\gamma_{i,2q-1})=\beta_{i,q},\ &
\theta_{y}^{-1}(\gamma_{i,2q})=\beta_{i,q}^{-1}& \\
&\gamma_{i,2q-1}=\alpha_{i,q}, &\gamma_{i,2q}=\alpha_{i,q+1}^{-1}&
\rangle,
\end{array}
$$
where $1\leq q\leq N+1$ and $\theta_y$ is the automorphism of
$F_{(m-1)(2N+1)}$ given by the braid $y$. Since $y\in\gamma_{n+1}
P_{(m-1)(2N+1)}$, it is easy to see that
$$\theta_{y}^{-1}(\gamma_{i,j})\equiv \gamma_{i,j}\mod{\gamma_{m+1}\pi_1\widetilde X}.$$
Replacing $\theta_{y}^{-1}(\gamma_{i,j})$ by $\gamma_{i,j}$ in the
presentation of $\pi_1\widetilde X$ we obtain a presentation of the
free group $F_{m-1}$.
\end{proof}

Now, a string link that gives rise to an element of
$\FLinks_{m-1}(n)$ can be written as $\SS(xy)$ where $x\in\Phi$ and
$y\in\gamma_{n+1} P_{\infty}^{m-1}$. Any link $n$-equivalent to it
is of the form $\SS(txyb\cdot h)$ where $t\in H^T$, $b\in H^B$ and
$h\in\gamma_{n+1} P_{\infty}$. We have
$$\SS(txyb\cdot h)=\SS(xy\cdot bhb^{-1})=\SS(xh'yh''),$$
where $h'\in \gamma_{n+1}\Phi$ and
$h''\in\gamma_{n+1}P_{\infty}^{m-1}$. It follows from the two
foregoing lemmas that $\SS(xh'yh'')$ and $\SS(xy)$ define the same
element of $F_{m-1}/\gamma_{n+1} F_{m-1}$.

\begin{xcb}{Exercises}

\begin{enumerate}
\item
Show that reducing the coefficients of the Magnus expansion of an
element of $F_n$ modulo $m$, we obtain the universal $\Z_m$-valued
Vassiliev invariant for $F_n$. Therefore, all mod $m$ Vassiliev
invariants for $F_n$ are mod $m$ reductions of integer-valued
invariants.
\item \label{ex:univas}
Let $\M':F_n\to \Z\la X_1,\ldots,X_n\ra$ be any multiplicative map
such that for all $x_i$ we have $\M'(x_i)=1+\alpha_iX_i+\ldots$ with
$\alpha_i\neq 0$. Show that $\M'$ is a universal Vassiliev invariant
for $F_n$.
\item \label{ex:kinotunique}
Show that the Kontsevich integral of an element of a free group
$F_m$ thought of as a path in a plane with $m$ punctures depends on
the positions of the punctures.
\item \label{ex:pthree}
(a) Show that the semi-direct product in the decomposition
$P_3=F_2\ltimes\Z$ given by combing is not direct.\\
(b) Find an isomorphism between $P_3$ and $F_2\times\Z$.
\end{enumerate}

\end{xcb}
 %10 Braids
\chapter{Gauss diagrams} % 13
\label{chapGD}

In this chapter we shall show how the finite type invariants of a
knot can be read off its Gauss diagram. It is not surprising that
this is possible in principle, since the Gauss diagram encodes the
knot completely. However, the particular method we describe,
invented by Polyak and Viro and whose efficiency was proved by
Goussarov, turns out to be conceptually very simple. For a given
Gauss diagram, it involves only counting its subdiagrams of some
particular types.

We shall prove that each finite type invariant arises in this way
and describe several examples of such formulae.

\section{The Goussarov theorem}\label{sec:Gous-th}

Recall that in Chapter~\ref{chapBr} we have constructed a universal
Vassiliev invariant for the free group by sending a word to the sum
of all of its subwords. A similar construction can be performed for
knots if we think of a knot as being ``generated by its crossings''.

Let $\GD$ be the set of all Gauss diagrams (we shall take them to be
based, or long, even though for the moment it is of little
importance).
Denote by $\Z\GD$ the set of all finite linear
combinations of the elements of $\GD$ with integer coefficients. We define the map $I:\Z\GD\to\Z\GD$
by simply sending a diagram to the sum of all its subdiagrams:
$$I(D):=\sum_{D'\subseteq D} D'\label{map-I-gd-ad}
$$
and continuing this definition to the whole of $\Z\GD$ by linearity.
In other terms, the effect of this map can be described as
\smallskip
$$I: \gaussAA\quad\boldsymbol{\longmapsto} \quad\gaussAA\ \ \boldsymbol{+}\ \
\gaussAC$$
\smallskip
For example, we have
\begin{multline*}
I\Bigl(\gaussBA\Bigr)\ =\gaussBA + \gaussBC \\ + \gaussBD + \gaussBE
+ \gaussBF + 2 \gaussBG + \gaussBH
\end{multline*}
\smallskip
Here all signs on the arrows are assumed to be, say, positive.

Define the pairing $\langle A,D\rangle$ of two Gauss diagrams $A$ and $D$ as the coefficient of $A$ in $I(D)$:
$$I(D)=\sum_{A\in\GD} \langle A,D\rangle A.$$
In principle, the integer $\langle A,D\rangle$ may change if a Reidemeister move is
performed on $D$. However, one can find invariant linear
combinations of these integers. For example, in
Section~\ref{subsection:Casson} we have proved that the Casson
invariant $c_2$ of a knot can be expressed as
\def\gdCi#1#2{\bigl\langle \,\risS{-5}{gauss-6A}{\put(0,10){$\scriptstyle #1$}
     \put(35,10){$\scriptstyle #2$}}{45}{0}{0}\,, D\bigr\rangle }
$$\gdCi{+}{+} - \gdCi{-}{+} -\gdCi{+}{-} +\gdCi{-}{-}\ .
$$
More examples of such invariant expressions can be found in Section~\ref{section:examples}.  In fact, as we shall now see, for each Vassiliev knot invariant there exists a formula of this type.

\subsection{The Goussarov Theorem}

Each linear combination of the form
$$\sum_{A\in\GD}c_A\langle A,D\rangle$$  with integer coefficients, considered as a function of $D$, is just the composition $c\circ I$, where $c:\Z\GD\to\Z$ is the linear map with $c(A)=c_A$.

In what follows, usual knots will be referred to as {\em classical} knots, in order to distinguish them from virtual knots.\index{Knot!classical} Gauss diagrams that encode long classical knots, or, in other words, {\em realizable},\index{Gauss diagram!realizable} diagrams, form a subset $\GD^{re}\subset\GD$. Any integer-valued knot invariant $v$ gives rise to a function $\GD^{re}\to \Z$ which extends by linearity to a function $\Z\GD\to\Z$. We also denote this extension by $v$. Here   $\Z\GD^{re}$ is the free abelian group generated by the set  $\GD^{re}$.

\begin{xtheorem}[Goussarov]
For each integer-valued Vassiliev invariant $v$ of classical knots of
order $\le n$ there exists a linear map $c:\Z\GD\to\Z$ such that
$v=c\circ I\!\mid_{\Z\GD^{re}}$ and such that $c$ is zero on each Gauss
diagram with more than $n$ arrows.
\end{xtheorem}

The proof of the Goussarov Theorem is the main goal of this section.

\subsection{Construction of the map $c$}\label{ss:constr-c}

Consider a given Vassiliev knot invariant $v$ of order $\leq n$ as a
linear function $v:\Z\GD^{re}\to\Z$. We are going to
define a map $c:\Z\GD\to\Z$ such that 
the equality $c\circ I = v$ holds on the submodule $\Z\GD^{re}$.

Since the map $I$ is an isomorphism with the inverse being
\begin{equation*}
\label{i-to-mone}
I^{-1}(D)=\sum_{D'\subseteq D} (-1)^{|D-D'|}D',
\end{equation*}
where $|D-D'|$ is the number of arrows of $D$ not contained in $D'$,
the definition appears obvious:
\begin{equation}
c=v\circ I^{-1}\ .
\end{equation}
However, for this equation to make sense we need to extend $v$ from
$\Z\GD^{re}$ to the whole of $\Z\GD$ since the image of $I^{-1}$ contains 
all the subdiagrams of $D$ and a subdiagram of a realizable diagram does not 
have to be realizable.

To make such an extension consistent we need it to satisfy the Vassiliev skein 
relation. Thus we first express this relation in terms of Gauss diagrams 
in Section \ref{subsec:GDwithchords} introducing Gauss diagrams with 
undirected signed chords. It turns out (Section \ref{ss:extens-v}) 
that via Vassiliev skein relation an arbitrary Gauss diagram can be 
presented as a linear combination of some realizable Gauss diagrams 
(which we called {\em descending} and define them in Section 
\ref{section:desc}) plus Gauss diagrams with more than $n$ chords. 
Since $v$ is of order $\le n$, it is natural to extend it by zero 
on Gauss diagrams with more than $n$ chords. Also, we know the values 
of $v$ on descending Gauss diagrams since they are realizable, and,
thus, such a presentation gives us the desired extension.
We complete the proof of the Theorem in Section~\ref{ss:proff-gt}.

\subsection{Gauss diagrams with chords}\label{subsec:GDwithchords}
Gauss diagrams can also be naturally defined for knots with double points.
Apart from the arrows, these diagrams have solid
undirected chords on them, each chord labelled with a sign.
The {\em sign of a chord} is positive if in the positive
resolution of the double
point the overcrossing is passed
first.
(Recall that we are dealing with long Gauss diagrams, and that the
points on a long knot are ordered.)

Gauss diagrams with at most $n$ chords span the space $\Z\GD_n$,
which is mapped to $\Z\GD$ by a version of the Vassiliev skein
relation:
\begin{equation}\label{eq:star}
  \rb{-4mm}{\ig[height=9mm]{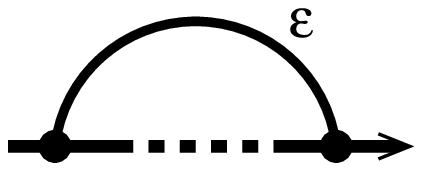}}\quad\boldsymbol{=}\quad
  \eps\ \rb{-4mm}{\ig[height=9mm]{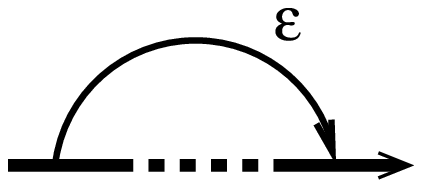}}\quad\boldsymbol{-}\quad
  \eps\ \rb{-4mm}{\ig[height=9mm]{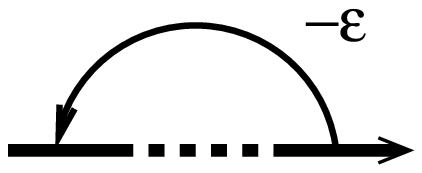}}
\end{equation}
Using this relation, any knot invariant, or, indeed, any
function on Gauss diagrams can be extended to diagrams with chords.
Note that the map $\Z\GD_n\to\Z\GD$ is not injective; in particular,
changing the sign of a chord in a diagram from $\GD_n$ multiplies
its image in $\Z\GD$ by $-1$.
We have a commutative diagram
$$\xymatrix{
\Z\GD_n \ar[rrr]^{\mbox{\scriptsize skein (\ref{eq:star})}}
   \ar[d]_(.5){I}
 &&& \Z\GD \ar[d]^(.5){I} \\
\Z\GD_n \ar[rrr]^{\mbox{\scriptsize skein (\ref{eq:star})}}
    &&& \Z\GD }
$$
where $I:\Z\GD_n\to\Z\GD_n$ is the isomorphism that sends a diagram
to the sum of all its subdiagrams that contain the same 
chords.

\subsection{Descending Gauss diagrams}
\label{section:desc}
We shall draw the diagrams of the long knots in the plane $(x,y)$,
assuming that the knot coincides with the $x$-axis outside some ball.

A diagram of a (classical) long knot is {\em descending}
\index{Gauss diagram!descending}
if for each crossing the overcrossing comes first.
A knot whose diagram is descending is necessarily
trivial. The Gauss diagram corresponding to a descending knot
diagram has all its arrows pointed in the direction of the increase
of the coordinate $x$ (that is, to the right).

The notion of a descending diagram can be generalized to diagrams of
knots with double points.

\begin{xdefinition} A Gauss diagram of a long knot with double
points is called {\em descending}
\index{Gauss diagram!descending!for singular knots}  if
\begin{enumerate}
\item {all the arrows are directed to the right};
\item {no endpoint of an arrow can
be followed by the left endpoint of a chord.}
\end{enumerate}
\end{xdefinition}
In other words, the following situations are forbidden:
$$\ig[width=2cm]{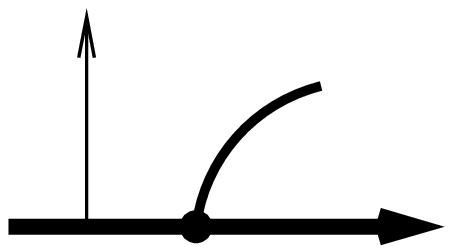}\qquad\ig[width=2cm]{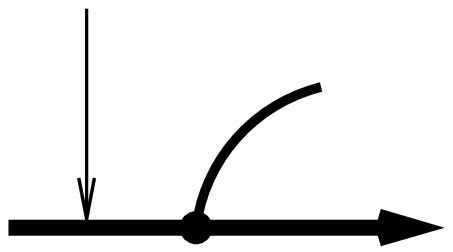}$$
For these two conditions to make sense the Gauss diagram with double
points need not be realizable; we shall speak of descending diagrams
irrespective of whether they can be realized by classical knots with
double points.

Descending diagrams are useful because of the following fact.
\begin{lemma}\label{le:can-repres}
Each long chord diagram with signed chords underlies a unique (up to isotopy)
singular classical long knot that has a descending Gauss diagram.
\begin{proof}
The endpoints of the chords divide the line of the parameter into intervals,
two of which are semi-infinite. Let
us say that such an interval is {\em prohibited} if it is bounded from the
right by a left end of a chord. Clearly, of the two semi-infinite intervals
the left one
is prohibited while the right one is not.  If a chord diagram $D$ underlies a
descending Gauss diagram $G_D$, then $G_D$ has no arrow endpoints
on the prohibited intervals.
We shall refer to the union of all prohibited intervals with some small
neighbourhoods of the chord endpoints (which do not contain endpoints of
other chords or arrows) as the {\em prohibited set}.
$$\risS{0}{proh-set}{}{350}{80}{15}
$$
$$\risS{0}{desc-sing-kn}{}{280}{80}{10}
$$

The prohibited set of a diagram can be immersed into the plane with double points
corresponding to the chords, in such a way that the signs of the chords are respected.
Such an immersion is uniquely defined up to isotopy. 
The image of the prohibited set will be an embedded tree $T$. 

The leaves of $T$ are numbered in the order given by the parameter
along the knot. Note that given $T$, the rest of the plane diagram can be
reconstructed as follows: the leaves of $T$ are joined, in order, by arcs
lying outside of $T$;
these  arcs only touch $T$ at their endpoints and each arc lies
below all the preceding arcs; the last arc
extends to infinity. Such reconstruction is unique since the complement of
$T$ is homeomorphic to a 2-disk, so all possible choices of arcs are homotopic.
\end{proof}
\end{lemma}

\subsection{Extension of $v$}\label{ss:extens-v}

Here, using the Vassiliev skein relation, we extend $v$ not only to  singular long knots (realizable Gauss gauss diagrams with chords) but also to arbitrary Gauss diagrams with signed chords.

If $D$ is a descending Gauss diagram with signed chords, by Lemma~\ref{section:desc} there
exists precisely one singular classical knot $K$ which has a descending diagram with the same signed chords. We set $v(D):=v(K)$. 

If $D$ is an arbitrary diagram, we apply the algorithm which is described below to represent $D$ as a linear combination $\sum a_i D_i$ of descending diagrams modulo diagrams with the number of chords $\ge n$. The algorithm uses the Vassiliev skein relation and has the property that it transforms a realizable Gauss diagram into the linear combination of realizable diagrams. 

Now, if we set $$v(D):=\sum a_i v(D_i);$$ this expression naturally vanishes on diagrams with more than $n$ chords since $v$ is of order $\le n$. Therefore, indeed, we get an extension of $v$.

\medskip

\noindent{\bf The algorithm}\label{pp:desc-alg} consists in the iteration of a certain transformation $P$ of Gauss diagrams which makes a diagram, in a sense, ``more descending''.
The map $P$ works as follows.

Take a diagram $D$. Replace all the arrows of $D$ that point to the
left by the arrows that point to the right (possibly creating new 
chords in the process), using relation (\ref{eq:star}).

Denote by $\sum a_i D'_i$ the resulting linear combination.
Now, each of the $D'_i$ may contain ``prohibited pairs'': these are
the arrow endpoints which are followed by the left endpoint of a chord.
Using the Reidemeister moves a prohibited pair can be transformed as
follows:
$$\risS{-10}{gauss-10A-2}{}{45}{12}{0}\ \risS{0}{totor}{}{25}{0}{0}\
     \risS{-10}{gauss-10A-1}{}{45}{0}{0}\qquad\qquad
\risS{-10}{gauss-10B-2}{}{45}{0}{0}\ \risS{0}{totor}{}{25}{0}{0}\
     \risS{-10}{gauss-10B-1}{}{45}{0}{12}
$$
On a Gauss diagram this transformation can take one of the forms
shown in Figure~\ref{figure:P} where the arrows corresponding to the
new crossings are thinner.

\begin{figure}[h!]
\begin{longtable}{c@{\quad\risS{0}{totor}{}{25}{0}{0}\quad}c}
\risS{-10}{gauss-11-01}{}{100}{20}{0} &
        \risS{-10}{gauss-11-02}{}{100}{0}{0}\\
\risS{-10}{gauss-11-03}{}{100}{40}{0} &
        \risS{-10}{gauss-11-04}{}{100}{0}{0}\\
\risS{-10}{gauss-11-05}{}{100}{40}{0} &
        \risS{-10}{gauss-11-06}{}{100}{0}{0}\\
\risS{-10}{gauss-11-07}{}{100}{40}{0} &
        \risS{-10}{gauss-11-08}{}{100}{0}{0}\\
\risS{-10}{gauss-11-09}{}{100}{40}{0} &
        \risS{-10}{gauss-11-10}{}{100}{0}{0}\\
\risS{-10}{gauss-11-11}{}{100}{50}{0} &
        \risS{-10}{gauss-11-12}{}{100}{0}{0}
\end{longtable}
\caption{}\label{figure:P}\vspace{10pt}
\end{figure}

For each $D'_i$
consider the leftmost prohibited pair, and replace it with the
corresponding configuration of arrows and chords as in
Figure~\ref{figure:P}; denote the resulting diagram by $D''_i$. Set
$P(D):=\sum a_i D''_i$ and extend $P$ linearly to the whole
$\Z\GD_{\infty}=\bigcup_n \Z\GD_n$.

If $D$ is descending, then $P(D)=D$. We claim that applying $P$
repeatedly to any diagram we shall eventually arrive to a linear
combination of descending diagrams, modulo the diagrams with more
than $n$ chords.

Let us order the chords in a diagram by their left endpoints. We say
that a diagram is {\em descending up to the $k$th chord} if the
closed interval from $-\infty$ up to the left end of the $k$th chord
contains neither endpoints of leftwards-pointing arrows, nor
prohibited pairs.

If $D$ is descending up to the $k$th chord, each diagram in $P(D)$
also is. Moreover, applying $P$ either decreases the number of arrow
heads to the left of the left end of the $(k+1)$st chord, or preserves
it. In the latter case, it decreases the number of arrow tails in
the same interval.
It follows that for some finite $m$ each diagram in $P^m(D)$ will be
decreasing up to the $(k+1)$st chord. Therefore, repeating the
process, we obtain after a finite number of steps a combination of
diagrams descending up to the $(n+1)$st chord. Those of them that have
at most $n$ chords are descending, and the rest can be disregarded.

\begin{xremark} By construction, $P$ respects the realizability of
the diagrams. In particular, the above algorithm expresses a long
classical knot as a linear combination of singular classical knots with
descending diagrams.
\end{xremark}

\subsection{Proof of the Goussarov Theorem}\label{ss:proff-gt}

To prove the Goussarov
Theorem we now need to show that $c=v\circ I^{-1}$ vanishes on Gauss diagrams with
more than $n$ arrows.

Let us evaluate $c$ on a descending Gauss diagram $A$ whose total
number of chords and arrows is greater than $n$. We have
$$c(A)=v(I^{-1}(A))=\sum_{A'\subseteq A} (-1)^{|A-A'|}v(A').$$
All the subdiagrams $A'$ of $A$ have the
same chords as $A$ and therefore are descending. Hence, by the construction
of the extension of $v$ to $\Z\GD$,  the values of $v$ on all the $A'$
are equal to $v(A)$. If $A$ has more than $n$ chords, then $v(A)=0$. If
$A$ has at most $n$ chords, it has at least one arrow.  It is easy
to see that in this case $\sum_{A'\subseteq A} (-1)^{|A-A'|}=0$, and
it follows that $c(A)=0$. In particular, $c$ vanishes on all
descending Gauss diagrams with more than $n$ arrows.

In order to treat non-descending Gauss diagrams, we shall introduce
an algorithm, very similar to that of Section~\ref{ss:extens-v} (see page \pageref{pp:desc-alg})
that converts any long Gauss diagram with chords into a combination
of descending diagrams with at least the same total number of chords
and arrows. The algorithm consists in the iteration of a certain map
$Q$, similar to $P$,  which also makes a diagram ``more descending''.
We shall prove that the map $Q$ preserves $c$ in the sense that
$c\circ Q=c$ and does not decrease the total number of chords and arrows.
Then applying $Q$ to a Gauss diagram $A$ enough number
of times we get a linear combination of descending diagrams without
altering the value of $c$. Then the arguments of the previous paragraph
show that $c(A)=0$ which will conclude the proof of the Goussarov Theorem.

Take a Gauss diagram $A$. As in Section~\ref{ss:extens-v},
we replace
all the arrows of $A$ that point leftwards by the arrows that point
to the right, using relation (\ref{eq:star}).

Denote by $\sum a_i A'_i$ the resulting linear combination and check
if the summands $A'_i$ contain prohibited pairs. Here is where our new
construction differs from the previous one. For each $A'_i$
consider the leftmost prohibited pair, and replace it with the sum
of the seven non-empty subdiagrams of the corresponding
diagram from the right column of Figure~\ref{figure:P} containing at least
one of the three arrows. Denote the
sum of these seven diagrams by $A''_i$. For example, if $A'_i$ is the first
diagram from the left column of Figure~\ref{figure:P},
$$A'_i=\risS{-10}{ggauss-11-01}{}{100}{30}{0}, \mbox{\quad then\quad}
A''_i=\risS{-10}{ggauss-11-02-1}{}{100}{0}{0} + $$
$$+ \risS{-10}{ggauss-11-02-2}{}{100}{45}{0} +
\risS{-10}{ggauss-11-02-3}{}{100}{0}{0} +
\risS{-10}{ggauss-11-02-4}{}{100}{0}{0} + $$
$$+\risS{-10}{ggauss-11-02-5}{}{100}{45}{20} +
\risS{-10}{ggauss-11-02-6}{}{100}{0}{0} +
\risS{-10}{ggauss-11-02-7}{}{100}{0}{0}\ .$$
Now, set $Q(A)=\sum a_i A''_i$ and
extend $Q$ linearly to the whole $\Z\GD_{\infty}$.

As before, applying $Q$ repeatedly to any diagram we shall
eventually arrive to a linear combination of descending diagrams,
modulo the diagrams with more than $n$ chords. Note that $Q$ does
not decrease the total number of chords and arrows.

It remains to prove that $Q$ preserves $c$. Since $I:\Z\GD\to\Z\GD$ is
epimorphic, it is sufficient to check this on diagrams of the form $I(D)$.
Assume that we have established
that $c(Q(I(D))))=c(I(D))$ for all Gauss diagrams $D$ with some
chords and at most $k$ arrows. If there are no arrows at all then $D$
is descending and $Q(I(D))=I(D)$. Let now $D$ have $k+1$ arrows. If
$D$ is descending, than again $Q(I(D))=I(D)$ and there is nothing
to prove. If $D$ is not descending, then let us first assume for
simplicity that all the arrows of $D$ point to the right. Denote by
$l$ the arrow involved in the leftmost prohibited pair, and let
$D_l$ be the diagram $D$ with $l$ removed. We have
$$I(P(D))=Q\bigl(I(D)-I(D_l)\bigr)+I(D_l).$$
Indeed, $P(D)$ is a diagram from the right column of Figure~\ref{figure:P}.
Its subdiagrams fall into two categories depending on whether they contain
at least one of the three arrows
indicated on Figure~\ref{figure:P} or none of them. The latter are subdiagrams
of $D_l$ and they are included in $I(D_l)$. The former
can be represented as $Q\bigl(I(D)-I(D_l)\bigr)$.

By the induction assumption, $c(Q(I(D_l)))=c(I(D_l))$.
Therefore, $$c(Q(I(D)))=c(I(P(D)))=v(P(D))\ .$$
But applying $P$ does not change the value of $v$ because of our definition
of the extension of $v$ from Section \ref{ss:constr-c}. Therefore,
$$c(Q(I(D)))=v(P(D))=v(D)=c(I(D))\ ,$$
and, hence $c(Q(A))=c(A)$ for any Gauss diagram $A$.

If some arrows of $D$ point to the left, the argument remains
essentially the same and we leave it to the reader.
\hspace{\fill}$\square$

\subsection{Example. The Casson invariant}
\index{Casson invariant}
We exemplify the proof of the Goussarov theorem by deriving the Gauss diagram
formula for the Casson invariant, that is, the second coefficient of the Conway
polynomial $c_2$. At the beginning of this chapter we already mentioned
a formula for it, first given in Section~\ref{subsection:Casson}.
However, the expression that we are going to derive following the proof of the
Goussarov theorem will be different.

Let $v=c_2$. We shall use the definition $c=v\circ I^{-1}$ to find
the function $c:\Z\GD\to\Z$.

If a Gauss diagram $A$ has at most one arrow then obviously $c(A)=0$.
Also, if $A$ consists of two non-intersecting arrows then $c(A)=0$.
So we need to consider the only situation when $A$ consists of two
intersecting arrows. There are 16 such diagrams differing by
the direction of arrows and signs on them.
The following table shows the values of $c$ on all of them.
\def\gt#1#2#3{\risS{-5}{#3}{\put(0,10){$\scriptstyle #1$}
     \put(35,10){$\scriptstyle #2$}}{45}{13}{6}}
\def\gA#1#2{\gt{#1}{#2}{gauss-6A}}
\def\gB#1#2{\gt{#1}{#2}{gauss-6B}}
\def\gC#1#2{\gt{#1}{#2}{gauss-6C}}
\def\gD#1#2{\gt{#1}{#2}{gauss-6D}}
$$\begin{array}{|c|c|c|c|} \hline
c(\gB{+}{+})=0&c(\gB{-}{+})=0&c(\gB{+}{-})=0&c(\gB{-}{-})=0\\ \hline
c(\gC{+}{+})=0&c(\gC{-}{+})=0&c(\gC{+}{-})=0&c(\gC{-}{-})=0\\ \hline
c(\gA{+}{+})=0&c(\gA{-}{+})=0&c(\gA{+}{-})=0&c(\gA{-}{-})=0\\\hline
c(\gD{+}{+})=1&\!c(\gD{-}{+})\!=\!-1\!\!&\!c(\gD{+}{-})\!=\!-1\!\!
              &c(\gD{-}{-})=1\\ \hline
\end{array}
$$

Let us do the calculation of some of these values in detail.

Take $A=\gA{+}{-}$. According to the definition of $I^{-1}$ from page
\pageref{i-to-mone} we have
$$c(A)=v(\risS{-1}{gauss-c2-0}{}{45}{13}{6})
-v(\risS{-5}{gauss-c2-1l}{\put(30,10){$\scriptstyle +$}}{45}{13}{6})
-v(\risS{-5}{gauss-c2-1r}{\put(5,10){$\scriptstyle -$}}{45}{13}{6})
+v(\gA{+}{-})
$$
The first three values vanish. Indeed, the first and third Gauss diagrams
are descending, so they represent the trivial long knot, and
the value of $c_2$ on it is equal to zero.
For the second value one should use the Vassiliev skein relation
(\ref{eq:star})
$$v(\risS{-5}{gauss-c2-1l}{\put(30,10){$\scriptstyle +$}}{45}{13}{6})
= v(\risS{-5}{gauss-c2-1r}{\put(30,10){$\scriptstyle -$}}{45}{13}{6})
+ v(\risS{-5}{gauss-c2-1c}{\put(30,10){$\scriptstyle -$}}{45}{13}{6})
$$
and then notice that the both diagrams are descending. Moreover, for the second diagram with a single
chord both resolutions of the corresponding double point lead to the trivial knot.

Thus we have
$$c(A)=v(\gA{+}{-})=v(\gB{-}{-})
+ v(\gt{-}{-}{gauss-c2-cl})\ .
$$
The last two Gauss diagrams are descending. Therefore, $c(A)=0$.

Now let us take $A=\gD{+}{+}$. Applying $I^{-1}$ to $A$ we get that the value of $c$ on the first three diagrams is equal to zero as before, and
$$c(\gD{+}{+}) = v(\gt{+}{+}{gauss-6D}).$$ To express the last Gauss diagram
as a combination of descending diagrams first
we should reverse its right arrow using the relation (\ref{eq:star}):
$$
\gt{+}{+}{gauss-6D}\ =\ \gB{+}{-}
\ +\ \gt{+}{-}{gauss-c2-cr}\ .
$$
The first Gauss diagram here is descending. But the second one is not, it has a prohibited pair. So we have to apply the map $P$ from
Section~\ref{ss:extens-v} to it. According to the first case of Figure~\ref{figure:P} we have
$$\gt{+}{-}{gauss-c2-cr}\ =\ \risS{-10}{gauss-c2-4crll}{}{75}{20}{8}
\ =\ \risS{-10}{gauss-c2-4crrl}{}{75}{20}{6} +
     \risS{-10}{gauss-c2-4crcl}{}{75}{20}{6}\ .
$$
In the first diagram we have to reverse one more arrow, and to the second diagram we
need to apply the map $P$ again. After that, the reversion of arrows in it
would not create any problem since the additional terms would have 3
chords, and we can ignore them if we are interested in the second order
invariant $v=c_2$ only. Modulo diagrams with three chords, we have
$$\gt{+}{-}{gauss-c2-cr}\ =\ \risS{-10}{gauss-c2-4crrr}{}{75}{15}{8}
- \risS{-10}{gauss-c2-4crrc}{}{75}{0}{0} +
     \risS{-10}{gauss-c2-6}{}{105}{0}{0}.
$$
The first and third diagrams here are
descending. But with the second one we have a  problem because it has
a prohibited interval with many (three) arrow ends on it. We need to apply
$P$ five times in order to make it descending modulo diagrams with three
chords. The result will be a descending diagram $B$ with two non-intersecting chords, one inside the other. So the value of $v$ on it would be
zero and we may ignore this part of the calculation (see Exercise~\ref{pr:c2-on-B} on page \pageref{pr:c2-on-B}). Nevertheless,  we give the answer here so that
the interested readers can check their understanding of the
procedure:
$$B=\risS{-10}{gauss-c2-14}{}{255}{50}{10}\ .\label{eq:gauss-c2-B}
$$
Combining all these results we have
$$\gt{+}{+}{gauss-6D}\ =\ \gt{+}{-}{gauss-c2-2rr}\ \ - B\ +\
\risS{-10}{gauss-c2-4crrr}{}{75}{15}{8}\ +\
     \risS{-10}{gauss-c2-6}{}{105}{0}{0}
$$
modulo diagrams with at least three chords.
The value of $v$ on the last Gauss diagram is equal to its value on the
descending knot with the same chord diagram,
$\risS{-3}{cat-cd12}{\put(0,8){$\scriptstyle -$}
                   \put(20,8){$\scriptstyle -$}}{30}{0}{0}$\ ,
namely, the knot $$\risS{-7}{cat-dk15}{}{34}{15}{10}.$$%\ .
It is easy to see that the only resolution that gives a non-trivial knot is the
positive resolution of the right double point together with the negative resolution of the left double point; the resulting knot is $4_1$. The value of $v=c_2$ on it is $-1$ according to the Table~\ref{conw_table} on page~\pageref{conw_table}. Thus the value of $v$ on this Gauss diagram is equal to 1. The values of $v$ on the other three descending Gauss diagrams are zero. Therefore, we have
$$c(\gD{+}{+})=1.$$

As an exercise, the reader may wish to check all the other values of $c$ from the table.

This table implies that the value of $c_2$ on a knot $K$ with the Gauss
diagram $D$ is
$$c_2(K)=\bigl\langle \,\gD{+}{+}-\gD{-}{+}-\gD{+}{-}+\gD{-}{-}\,,
D\bigr\rangle \ .
$$\label{otherc2}
This formula differs from the one at the beginning of the chapter by
the orientation of all its arrows.

\section{Canonical actuality tables}
\label{section:can_act_tables}
As a byproduct of the proof of the Goussarov Theorem, namely the Lemma \ref{le:can-repres} on page \pageref{le:can-repres},
we have the following refinement of the notion of an actuality table from Section~\ref{section:actuality_tables}.

In that section we have described a procedure
of calculating a Vassiliev invariant given by an actuality table. This
procedure involves some choices. Firstly, in order to build the
table, we have to choose for each chord diagram a singular knot
representing it. Secondly, when calculating the knot invariant we
have to choose repeatedly sequences of crossing changes that will
express our knot as a linear combination of singular knots from the
table.

It turns out that for long knots these choices can be eliminated. We
shall now define something that can be described as a canonical
actuality table and describe a calculation procedure for Vassiliev
invariants that only depends on the initial Gauss diagram
representing a knot. Strictly speaking, our ``canonical actuality
tables'' are not actuality tables, since they contain one singular
knot for each long chord diagram {\em with signed chords}.

A {\it canonical actuality table}\index{Actuality table!canonical}
for an invariant of order $n$ is the set of its values on all singular
long knots with descending diagrams and at most $n$ double points.

For example, here is the canonical actuality table for the second coefficient
$c_2$ of the Conway polynomial.
$$\begin{array}{c|c|c|c|c|c|c} \hline\hline
\risS{0}{cat-cd01}{}{30}{16}{0} &
\risS{0}{cat-cd02}{\put(2,8){$\scriptstyle +$}}{30}{0}{0} &
\risS{0}{cat-cd02}{\put(2,8){$\scriptstyle -$}}{30}{0}{0} &
\risS{0}{cat-cd04}{\put(-2,6){$\scriptstyle +$}
                   \put(20,6){$\scriptstyle +$}}{30}{0}{0} &
\risS{0}{cat-cd04}{\put(-2,5){$\scriptstyle -$}
                   \put(20,6){$\scriptstyle +$}}{30}{0}{0} &
\risS{0}{cat-cd04}{\put(-2,6){$\scriptstyle +$}
                   \put(22,5){$\scriptstyle -$}}{30}{0}{0} &
\risS{0}{cat-cd04}{\put(-2,5){$\scriptstyle -$}
                   \put(22,5){$\scriptstyle -$}}{30}{0}{0} \\ \hline
\risS{-10}{cat-cd01}{}{30}{14}{15} &
\risS{-10}{cat-dk02}{}{34}{0}{0} &
\risS{-10}{cat-dk03}{}{34}{0}{0} &
\risS{-10}{cat-dk04}{}{34}{0}{0} &
\risS{-10}{cat-dk05}{}{34}{0}{0} &
\risS{-10}{cat-dk06}{}{34}{0}{0} &
\risS{-10}{cat-dk07}{}{34}{0}{0} \\ \hline
0 & 0 & 0 & 0 & 0 & 0 & 0 \\ \hline \hline
\end{array}
$$
$$\begin{array}{c|c|c|c|c|c|c|c} \hline\hline
\risS{0}{cat-cd08}{\put(10,7.5){$\scriptstyle +$}
                   \put(10,17){$\scriptstyle +$}}{30}{25}{5} &
\risS{0}{cat-cd08}{\put(10,6){$\scriptstyle -$}
                   \put(10,17){$\scriptstyle +$}}{30}{0}{0} &
\risS{0}{cat-cd08}{\put(10,7.5){$\scriptstyle +$}
                   \put(10,15){$\scriptstyle -$}}{30}{0}{0} &
\risS{0}{cat-cd08}{\put(10,6){$\scriptstyle -$}
                   \put(10,15){$\scriptstyle -$}}{30}{0}{0} &
\risS{0}{cat-cd12}{\put(0,8){$\scriptstyle +$}
                   \put(20,8){$\scriptstyle +$}}{30}{0}{0} &
\risS{0}{cat-cd12}{\put(0,8){$\scriptstyle -$}
                   \put(20,8){$\scriptstyle +$}}{30}{0}{0} &
\risS{0}{cat-cd12}{\put(0,8){$\scriptstyle +$}
                   \put(20,8){$\scriptstyle -$}}{30}{0}{0} &
\risS{0}{cat-cd12}{\put(0,8){$\scriptstyle -$}
                   \put(20,8){$\scriptstyle -$}}{30}{0}{0} \\ \hline
\risS{-10}{cat-dk08}{}{34}{14}{15} &
\risS{-10}{cat-dk09}{}{34}{0}{0} &
\risS{-10}{cat-dk10}{}{34}{0}{0} &
\risS{-10}{cat-dk11}{}{34}{0}{0} &
\risS{-10}{cat-dk12}{}{34}{0}{0} &
\risS{-10}{cat-dk13}{}{34}{0}{0} &
\risS{-10}{cat-dk14}{}{34}{0}{0} &
\risS{-10}{cat-dk15}{}{34}{0}{0} \\ \hline
0 & 0 & 0 & 0 & 1 & 1 & 1 & 1 \\ \hline \hline
\end{array}
$$

To remove the second ambiguity in the procedure of calculating a Vassiliev
invariant we use the algorithm from Section~\ref{ss:extens-v}.
It expresses an arbitrary Gauss diagram with chords as a linear
combination of descending diagrams, modulo diagrams with more than $n$ chords.
The algorithm reduces a calculation of the value of a Vassiliev invariant to a linear combination of values from the canonical actuality table.

\section{The Polyak algebra for virtual knots}

The fact that all Vassiliev invariants can be expressed with the help of Gauss diagrams suggests that finite type invariants can be actually defined in the setup of Gauss diagrams. This is true and, moreover, there are two (inequivalent) ways to define Vassiliev invariants for virtual
knots: that of \cite{GPV} and that of \cite{Ka5}. Here we review the construction of  \cite{GPV}.
The reader should be warned that it is not known whether this definition coincides with the usual definition on classical knots. However, the logic behind it is very transparent and simple: the universal finite type invariant should send a knot to the sum of all of its  ``subknots''. We have already seen this approach in action in Chapter~\ref{chapBr} where the Magnus expansion of the free group was defined in this precise way.

\subsection{The universal invariant of virtual knots}

The map $I:\Z\GD\to\Z\GD$ from Section~\ref{sec:Gous-th},
sending
a diagram to the sum of all its subdiagrams
$I(D)=\sum_{D'\subseteq D} D'$,
is clearly not invariant under the Reidemeister moves.
However, we can make it invariant by simply taking the quotient of
the image of $I$ by the images of the Reidemeister moves, or their
{\em linearizations}. These linearizations have the following form:
{\flushleft{$\quad \gaussCA\ =\ 0,$}\\
\flushleft{$\quad \gaussCB\ \boldsymbol{+}\ \gaussCC\ \boldsymbol{+}\
\gaussCD\ =\ 0,$}\\
\flushleft{$\quad \gaussCE\ \boldsymbol{+}\ \gaussCF\ \boldsymbol{+}\
\gaussCG\ \boldsymbol{+}\ \gaussCH\ \boldsymbol{=}$}\\
\smallskip

{$\qquad\ \boldsymbol{=}\ \gaussCI\ \boldsymbol{+}\ \gaussCJ\
\boldsymbol{+}\ \gaussCK\ \boldsymbol{+}\ \gaussCL$.}}\\
\smallskip

\noindent The space $\Z\GD$ modulo the linearized Reidemeister moves
is called the {\em Polyak algebra}.
\index{Polyak algebra}\label{PolyakA}
The structure of an algebra
comes from the connected sum of long Gauss diagrams; we shall not
use it here. The Polyak algebra, which we denote by $\Polyak$, looks
rather different from the quotient of $\Z\GD$ by the usual
Reidemeister moves, the latter being isomorphic to the free Abelian group
spanned by the set of all virtual knots $V\K$. Note, however, that
by construction, the resulting invariant $I^*:\Z V\K\to \Polyak$ is
an isomorphism, and, therefore, contains the complete information
about the virtual knot.

It is not clear how to do any calculations in $\Polyak$, since the 
relations are not homogeneous. It may be more feasible to consider the 
(finite-dimensional) quotient $\Polyak_n$ of $\Polyak$ which is obtained by setting all the
diagrams with more than $n$ arrows equal to zero.  In fact, the space $\Polyak_n$ plays an important role in the theory  of Vassiliev invariants for virtual knots. Namely,
the map $I_n: \Z V\K\to \Polyak_n$ obtained by composing $I^*$ with
the quotient map is an order $n$ Vassiliev invariant for virtual knots, universal in 
the sense that any other order $n$ invariant is obtained by composing $I_n$ with 
some linear function on $\Polyak_n$.

Let us now make this statement precise and define the Vassiliev invariants .

While the simplest operation on plane knot diagrams is the crossing
change, for Gauss diagrams there is a similar, but even simpler
manipulation: deleting/inserting of an arrow. An analogue of a knot
with a double point for this operation is a diagram with a dashed
arrow. A dashed arrow can be resolved by means of the following
``virtual Vassiliev skein relation'':

$$\gaussDA\quad\boldsymbol{=} \quad\gaussDB\ \ \boldsymbol{-}\ \
\gaussDC.$$

An invariant of virtual knots is said to be of finite type (or
Vassiliev) of order $n$ if it vanishes on all Gauss diagrams with
more than $n$ dashed arrows.

Observe that the effect of $I$ on a diagram all of whose arrows are
dashed consists in just making all the arrows solid.
More generally, the image under $I$ of a Gauss diagram with some
dashed arrows is a sum of Gauss diagrams all of which contain these
arrows. It follows that $I_n$ is of order $n$: indeed, if a Gauss
diagram has more than $n$ dashed arrows it is sent by $I$ to a
Gauss diagram with at least $n$ arrows, which is zero in
$\Polyak_n$.

\subsection{Dimensions of $\Polyak_n$}
The universal invariants $I_n$, in marked contrast with the Kontsevich integral, 
are defined in a simple combinatorial fashion. However, nothing comes for free: $I_n$
takes its values in the space $\Polyak_n$ which is hard to describe. For small $n$, the 
dimensions of $\Polyak_n$ (over the real numbers) were calculated in \cite{BHLR}:

\begin{center}
\begin{tabular}{c||c|c|c|c|c|c|c|c|c|c|c|c}
$n$         & 1 & 2 & 3 & 4\\
\hline
$\dim \Polyak_n-\dim \Polyak_{n-1}$ & 2 & 7 & 42 & 246
\end{tabular}
\end{center}

\subsection{Open problems}\label{subsection:open}

A finite type invariant of order $n$ for virtual knots gives rise to
a finite type invariant of classical knots of at least the same order.
Indeed, a crossing change can be thought of as deleting an arrow
followed by inserting the same arrow with the direction reversed.

\begin{xxca} Define the Vassiliev invariants for closed (unbased) 
virtual knots and show that the analogue of the space $\Polyak_2$ is
0-dimensional. Deduce that the Casson knot invariant cannot be extended
to a Vassiliev invariant of order 2 for closed virtual knots.
\end{xxca}

It is not clear, however, whether a finite type invariant of classical
knots can be extended to an invariant of virtual {\em long} knots of
the same order. The calculation of \cite{GPV} show that this is true
in orders 2 and 3. 

Given that $I^*$ is a complete invariant for virtual knots, one may
hope that each virtual knot is detected by $I_n$ for some $n$. It is
not known whether this is the case. A positive solution to this
problem would also mean that Vassiliev invariants distinguish
classical knots.

It would be interesting to describe the kernel of the natural
projection $\Polyak_{n}\to\Polyak_{n-1}$ which kills the diagrams
with $n$ arrows. First of all, notice that using the linearization
of the second Reidemeister move, we can get rid of all signs in the
diagrams in $\Polyak_n$ that have exactly $n$ arrows: changing the
sign of an arrow just multiplies the diagram by $-1$. Now, the
diagrams that have exactly $n$ arrows satisfy the following
6T-relation in $\Polyak_n$:
\begin{multline*}
\gaussEA\ \boldsymbol{+}\ \gaussEB\ \boldsymbol{+}\ \gaussEC\ \boldsymbol{=}\\
\boldsymbol{=}\ \gaussED\ \boldsymbol{+}\ \gaussEE\ \boldsymbol{+}\ \gaussEF\ .
\end{multline*}

Consider the space $\vec\A_n$ of chord diagrams with $n$ {\em
oriented} chords, or arrows, modulo the 6T-relation.\label{arrown} There is a map
$i_n:\vec\A_n\to\Polyak_{n}$, whose image is the kernel of the
projection to $\Polyak_{n-1}$. It is not clear, however if $i_n$ is
an inclusion. The spaces $\vec\A_n$ were introduced in \cite{Po1}
where their relation with usual chord diagrams is discussed. A further discussion of 
these spaces and their generalizations can be found in \cite{BN9}.

One more open problem is as follows. Among the linear combinations of
Gauss diagrams of the order no greater than $n$ there are some that produce
a well defined invariant of degree $n$. Obviously, such combinations form
a vector space, call it $L_n$. The combinations that lead to the
identically zero invariant form a subspace $L'_n$. The quotient space
$L_n/L'_n$ is isomorphic to the space of Vassiliev invariants $\V_n$.
The problem is to obtain a description of (or some information about)
the spaces $L_n$ and $L'_n$ and in these terms learn something new
about $\V_n$. For example, we have seen that the Casson invariant
$c_2$ can be given  by two different linear combinations $k_1$, $k_2$
of Gauss diagrams of order 2. It is not difficult to verify that these
two combinations, together with the empty Gauss diagram $k_0$ that
corresponds to the constant 1, span the space $L_2$. The subspace $L'_2$
is spanned by the difference $k_1-k_2$. We see that
$\dim L_2/L'_2=2=\dim\V_2$. For degree 3 the problem is already open. We
know, for instance, three linearly independent combinations of Gauss
diagrams that produce the invariant $j_3$ (see Sections \ref{deg_three} and
\ref{GD_HOMFLY} below), but we do not know if their differences generate 
the space $L'_3$. Neither do we have any description of the space $L_3$.

\section{Examples of Gauss diagram formulae}\label{section:examples}

\subsection{Highest part of the invariant}
Let us start with one observation that will significantly simplify
our formulae.
\begin{xlemma}
Let $c:\Z\GD\to\Z$ be a linear map representing an invariant of
order $n$. If $A_1,A_2\in \GD$ are diagrams with $n$ arrows obtained
from each other by changing the sign of one arrow, then
$c(A_1)=-c(A_2)$.
\end{xlemma}

\begin{proof}
As we noted before, a knot invariant $c$ vanishes on
all linearized Reidemeister moves of the form $I(R)$, where $R=0$ is
a usual Reidemeister move on realizable diagrams.
Consider a linearized second Reidemeister move involving one diagram
$A_0$ with $n+1$ arrows and two diagrams $A_1$ and $A_2$ with $n$
arrows. Clearly, $c$ vanishes on $A_0$, and therefore
$c(A_1)=-c(A_2)$.
\end{proof}

This observation gives rise to the following notation.
Let $A$ be a Gauss diagram with $n$ arrows {\em without signs},
{\em an unsigned Gauss diagram}.\index{Gauss diagram!unsigned}
Given a Gauss diagram $D$, we denote by $\langle A,D\rangle$ the alternating
sum $$\sum_i(-1)^{{\sign}{A_i}}\langle A_i,D\rangle, $$ where the $A_i$ are all
possible Gauss diagrams obtained from $A$ by putting signs on its arrows, and
${\sign}{A_i}$ is the number of chords of $A_i$ whose sign is
negative. Since the value of $c$ on all the $A_i$ coincides, up to
sign, we can speak of the value of $c$ on $A$.

For example, the formula for the Casson invariant of a knot $K$ with
the Gauss diagram $D$ can be written as
$$c_2(K)=\bigl\langle\gD{}{}\,, D\bigr\rangle.$$

\subsection{Invariants of degree 3}
\label{deg_three}
Apart from the Casson invariant, the simplest Vassiliev knot
invariant is the coefficient $j_3(K)$ in the power series expansion
of the Jones polynomial (see Section~\ref{knpol_as_vi}). Many formulae 
for $j_3(K)$ are known; the first such formula was found by
M.~Polyak and O.~Viro in terms of unbased diagrams, see \cite{PV1}.
From the results of \cite{GPV} the following Polyak--Viro
expression for $j_3$ is easily derived:

$$-3\ \Bigl\langle \gaussJA+\gaussJB+\gaussJC+\gaussJD$$
$$+\gaussJE+\gaussJF+\gaussJG+\gaussJH$$
$$+\gaussJI+\gaussJJ+\gaussJK-\gaussJL$$
$$+\ \gD{+}{-}\ -\ 2\ \gD{-}{+}\ -\ \gD{-}{-}\quad ,\quad D\Bigr\rangle .$$
(In this formula a typo of \cite{GPV} is corrected.)
Here the bracket $\langle\cdot,\cdot\rangle$ is assumed to be linear in its
first argument.

S.~Willerton in his thesis \cite{Wil3} found the following formula
for $j_3$:

$$-3\ \Bigl\langle\ \gaussKC+\gaussKD+\gaussKE+\gaussKF$$
$$+\gaussKG-\gaussKH+\gaussKI-\gaussKJ$$
$$\ +2\gaussKK +2\gaussKA+2\gaussKB, D\bigr\rangle.$$

A third Gauss diagram formula for $j_3$ will be given in
Section \ref{GD_HOMFLY}.

Other combinatorial formulae
for $c_2(K)$ and $j_3(K)$ were found earlier by J.~Lannes \cite{Lan}:
they are not Gauss diagram formulae.

\subsection{Coefficients of the Conway polynomial}

Apart from the Gauss diagram formulae for the low degree invariants,
two infinite series of such formulae are currently known: those for the
coefficients of the Conway and the HOMFLY polynomials. The former can be, of
course, derived from the latter, but we start from the discussion of the
Conway polynomial, as it is easier. We shall follow the original
exposition of \cite{CKR}.

\begin{xdefinition}\label{def:one-comp}\rm
A chord diagram $D$ is said to be {\em $k$-component} if after the
parallel doubling of each chord as in the picture
$$\risS{-6}{chord}{}{30}{12}{8}\
 \risS{-2}{totor}{}{25}{0}{0}\ \risS{-6}{dchord}{}{30}{0}{0}\ ,
$$ the resulting curve will have $k$ components. We use the notation
$|D|=k$. (See also Section~\ref{symb_j_n}).
\end{xdefinition}

\begin{xexample}\rm
For chord diagrams with two chords we have:
$$\Bigl|\chd{cd22ch4}\Bigr| =1\ \Longleftarrow\ \chd{cd22-ppar}\ ,\hspace{2cm}
  \Bigl|\chd{cd21ch4}\Bigr| =3\ \Longleftarrow\ \chd{symst4}\ .
$$
We shall be interested in one-component diagrams only.
With four chords, there are four one-component diagrams (the notation is
borrowed from Table~\ref{cd4tab}):
$$d^4_1=\chd{cd4-01}\ ,\quad d^4_5=\chd{cd4-05}\ ,\quad
  d^4_6=\chd{cd4-06}\ ,\quad\mbox{and}\quad d^4_7=\chd{cd4-07}\ .
$$
\end{xexample}

\begin{xdefinition}\rm
We can turn a one-component chord diagram with a base point into an arrow
diagram according to the following rule. {Starting from the base point
we travel along the diagram with doubled chords. During this journey we pass
both copies of each chord in opposite directions. Choose an arrow on each
chord which corresponds to the direction of the first passage
of the chord.} Here is an example.\vspace{-5pt}
$$\chd{cd22ch4}\quad  \risS{2}{totor}{}{25}{15}{15}\quad
  \chd{cd22par}\quad  \risS{2}{totor}{}{25}{0}{0}\quad \chd{cd22arw}\ .
$$
\end{xdefinition}
We call the Gauss diagram obtained in this way {\em ascending}.

\begin{xdefinition}\label{def:cc}\rm
The {\it Conway combination}\index{Conway combination}
$\kr_{2n}$ is the sum of all based
one-component ascending Gauss diagrams with $2n$ arrows. For example,
$$\begin{array}{rcl}
\kr_2 &:=& \risS{-12}{cd22arw}{}{25}{15}{20}\ ,\\
\kr_4 &:=& \ard{cd4-01arw}\
  +\ \ard{cd4-07arw1} + \ard{cd4-07arw2} + \ard{cd4-07arw3} + \ard{cd4-07arw4}+ \\
&&\hspace{-8pt}
   + \ard{cd4-05arw1} + \ard{cd4-05arw2} + \ard{cd4-05arw3} + \ard{cd4-05arw4}
   + \ard{cd4-05arw5} + \ard{cd4-05arw6} + \ard{cd4-05arw7} + \ard{cd4-05arw8} + \\
&&\hspace{-8pt}
   + \ard{cd4-06arw1} + \ard{cd4-06arw2} + \ard{cd4-06arw3} + \ard{cd4-06arw4}
   + \ard{cd4-06arw5} + \ard{cd4-06arw6} + \ard{cd4-06arw7} + \ard{cd4-06arw8}\ .
\end{array}$$
Note that for a given one-component chord diagram we have to consider all possible choices for the base point. However, some choices may lead to the same Gauss diagram. In
$\kr_{2n}$ we list them without repetitions. For instance, all choices of a base point for the diagram $d^4_1$ give the same Gauss diagram. So $d^4_1$ contributes only one Gauss diagram to $\kr_4$. The diagram $d^4_7$ contributes four Gauss diagrams because of its symmetry, while $d^4_5$ and $d^4_6$ contribute eight Gauss diagrams each.
\end{xdefinition}

\begin{xtheorem}
For $n\geqslant 1$, the coefficient $c_{2n}$ of $z^{2n}$ in the Conway polynomial \index{Conway polynomial}
of a knot $K$ with the Gauss diagram $G$ is equal to
$$c_{2n} = \langle \kr_{2n},G \rangle\ .$$
\end{xtheorem}

\begin{xexample}\label{ex:g-2-eva}\rm
Consider the knot $K:=6_2$ and its Gauss diagram $G:=G(6_2)$:
$$\risS{-32}{k6-2or}{\put(-50,30){\mbox{\tt knot $6_2$}}}{70}{35}{40}
\qquad\qquad
G\ =\ \risS{-25}{Gd6-2}{\put(65,30){\mbox{\tt Gauss diagram}}
        \put(90,20){\mbox{\tt of $6_2$}}
             }{60}{20}{25}\hspace{20pt}$$
In order to compute the pairing $\langle \kr_4,G\rangle$ we must match
the arrows of each diagram of $\kr_4$ with the arrows of $G$.
One common property of all terms in $\kr_{2n}$ is that in each term
both endpoints of the arrows that are adjacent to the base point are 
{\em arrowtails}.
This follows from our construction of $\kr_{2n}$.
Hence, the arrow $\{1\}$ of $G$ can not participate in the matching
with any diagram of $\kr_4$. The only candidates to match with the
first arrow of a diagram of $\kr_4$ are the arrows $\{2\}$ and $\{4\}$
of $G$. If it is $\{4\}$, then
$\{1, 2, 3\}$ cannot participate in the matching, and there remain
only 3 arrows to match with the four arrows of $\kr_4$. Therefore, the
arrow of $G$ which matches with the first arrow of a diagram of $\kr_4$
must be $\{2\}$. In a similar way, we can find that the arrow of $G$ which
matches with the last arrow of a diagram of $\kr_4$ must be $\{6\}$. This
leaves three possibilities to match with the
four arrows of $\kr_4$: $\{2,3,4,6\}$, $\{2,3,5,6\}$, and $\{2,4,5,6\}$.
Checking them all we find only one quadruple, $\{2,3,5,6\}$, which
constitute a diagram equal to the second diagram
in the second row of $\kr_4$. The product of the
local writhes of the arrows $\{2,3,5,6\}$ is equal to
$(-1)(-1)(+1)(-1)=-1$. Thus,
$$\langle \kr_4,G \rangle =
  \langle\ \ard{cd4-05arw2}\ ,G \rangle = -1\ ,$$
which coincides with the coefficient $c_4$ of the Conway polynomial
$\nabla(K)= 1-z^2-z^4$.
\end{xexample}

\subsection{Coefficients of the HOMFLY polynomial}
\label{GD_HOMFLY}\index{HOMFLY polynomial}

Let $P(K)$ be the HOMFLY polynomial of the knot $K$.
Substitute $a=e^h$ and take the Taylor expansion in $h$. The
result will be a Laurent polynomial in $z$ and a power series in $h$. Let
$p_{k,l}(K)$ be the coefficient of $h^kz^l$ in that expression.
The numbers $p_{0,l}$ coincide with the coefficients of the Conway
polynomial, since the latter is obtained from HOMFLY by fixing $a=1$.
\medskip

\begin{xremark} It follows from Exercise~\ref{pr:vas-homfly} on page~\pageref{pr:vas-homfly} that

(1) for all nonzero terms the sum $k+l$ is non-negative;

(2) $p_{k,l}$ is a Vassiliev invariant of degree no greater than
$k+l$;

(3) if $l$ is odd, then $p_{k,l}=0$.
\end{xremark}

We shall describe a Gauss diagram formula for 
$p_{k,l}$ following \cite{CP}.

Let $A$ be a (based, or long) Gauss diagram, $S$ a subset of its arrows
(referred to as a {\it state})
and $\alpha$ an arrow of $A$. Doubling all the arrows in $A$ that belong to $S$, in the same fashion as in the preceding section, we obtain a diagram consisting of one or several circles with some signed arrows attached to them.
Denote by $\langle\alpha|A|S\rangle$
the expression in two variables $h$ and $z$ that depends on the sign of the
chord $\alpha$ and the type of the first passage of $\alpha$ (starting from
the basepoint) according to the following table:
$$\begin{array}{c||c|c|c|c}
\rb{7pt}{First passage:}&
 \risS{0}{fp-b}{}{35}{25}{3}&\risS{0}{fp-t}{}{35}{0}{0}&
 \risS{0}{fp-l}{}{35}{0}{0}&\risS{0}{fp-r}{}{35}{0}{0} \\ \hline\hline
\risS{-8}{cr_p-gd}{}{40}{22}{15} & e^{-h}z & 0 & e^{-2h}-1 & 0 \\ \hline
\risS{-8}{cr_m-gd}{}{40}{22}{15} & -e^hz & 0 & e^{2h}-1 & 0
\end{array}
$$
To the Gauss diagram $A$ we then assign a power series $W(A)$ in $h$ and $z$
defined by
$$
W(A)=\sum_S \plwa \left(\frac{e^h-e^{-h}}{z}\right)^{c(S)-1},
$$
where $\plwa=\prod_{\a\in A}\langle\alpha|A|S\rangle$ and
$c(S)$ is the number of components obtained after doubling all the
chords in $S$.
Denote by $w_{k,l}(A)$ the coefficient of $h^kz^l$ in this power series and
consider the following linear combination of Gauss diagrams:
$A_{k,l} := \sum\ w_{k,l}(A)\cdot A$. Note that the number $w_{k,l}(A)$
is non-zero only for a finite number of diagrams $A$.

\begin{xtheorem}
Let $G$ be a Gauss diagram of a knot $L$. Then
$$p_{k,l}(K)=\langle{A_{k,l}},{G}\rangle\ .
$$
\end{xtheorem}

For a proof of the theorem, we refer the reader to the original paper
\cite{CP}. Here we only give one example. To facilitate the practical
application of the theorem, we start with some general remarks.

A state $S$ of a Gauss diagram $A$ is called {\it ascending}, if in
traversing the diagram with doubled arrows we approach the neighbourhood of
every arrow (not only the ones in $S$) first at the arrow head. As follows
directly from the construction, only ascending states contribute to $W(A)$.

Note that since $e^{\pm 2h}-1=\pm 2h + \mbox{(higher degree terms)}$ and
$\pm e^{\mp h}z= \pm z + \mbox{(higher degree terms)}$, the power series
$W(A)$
starts with terms of degree at least $|A|$, the number of arrows of $A$.
Moreover, the $z$-power of $\plwa
\left(\frac{e^h-e^{-h}}{z}\right)^{c(S)-1}$
is equal to $|S|-c(S)+1$.
Therefore, for fixed $k$ and $l$, the weight $w_{k,l}(A)$ of a Gauss
diagram may be non-zero only if $A$ satisfies the following conditions:
\begin{itemize}
\item[(i)]  $|A|$ is at most $k+l$;
\item[(ii)] there is an ascending state $S$ such that $c(S)=|S|+1-l$.
\label{con-iii}
\end{itemize}

For diagrams of the highest degree $|A|=k+l$, the contribution of an
ascending state $S$ to $w_{k,l}(A)$ is equal to $(-1)^{|A|-|S|}2^k\e(A)$,
where $\e(A)$ is the product of signs of all arrows in $A$. If two such
diagrams $A$ and $A'$ with $|A|=k+l$ differ only by signs of the arrows,
their contributions to $A_{k,l}$ differ by the sign $\e(A)\e(A')$. Thus all
such diagrams may be combined to the unsigned diagram $A$, appearing in
$A_{k,l}$ with the coefficient $\sum_S(-1)^{|A|-|S|}2^k$ (where the
summation is over all ascending states of $A$ with $c(S)=|S|+1-l$).
\medskip

\begin{xxca}
Prove that Gauss diagrams with isolated arrows do not contribute
to $A_{k,l}$. ({\sl Hint:}
all ascending states cancel out in pairs.)
\end{xxca}

Now, by way of example, let us find an explicit formula for
$A_{1,2}$. The maximal number of arrows is equal to 3. To get $z^2$ in
$W(A)$ we need ascending states with either $|S|=2$ and $c(S)=1$, or $|S|=3$
and $c(S)=2$. In the first case the equation $c(S)=1$ means that the two
arrows of $S$ must intersect. In the second case the equation $c(S)=2$ does
not add any restrictions on the relative position of the arrows. In the cases
$|S|=|A|=2$ or $|S|=|A|=3$, since $S$ is ascending, $A$ itself must be
ascending as well.

For diagrams of the highest degree $|A|=1+2=3$, we must count ascending
states of unsigned Gauss diagrams with the coefficient $(-1)^{3-|S|}2$, that is,
$-2$ for $|S|=2$ and $+2$ for $|S|=3$. There are only four types of
(unsigned) 3-arrow Gauss diagrams with no isolated arrows:
\def\lhd#1{\ \ \risS{-13}{#1}{}{25}{18}{10}\ \ }
\def\fhd#1{\risS{-20}{#1}{}{25}{12}{25}}
\def\shd#1{\ \risS{-12}{#1}{}{25}{15}{15}}
$$\lhd{bcd35-1}; \qquad \risS{-10}{bcd34-3}{}{27.5}{0}{0}\ \ ,
\qquad \lhd{bcd34-1},\qquad \risS{-10}{bcd34-2}{}{27.5}{0}{0}\ \ .
$$
Diagrams of the same type differ by the directions of arrows.

For the first type, recall that the first arrow should be oriented towards
the base point; this leaves 4 possibilities for the directions of the remaining
two arrows. One of them, namely $\lhd{aiv-n}$, does not have ascending states
with $|S|=2,3$. The remaining possibilities, together with their ascending
states, are shown in the table:

$$\begin{array}[t]{||c|c|c|c||}\hline\hline
 \fhd{aiv-3} & \fhd{aiv-2} & \fhd{aiv-1} & \fhd{aiv-1} \\
 \shd{civ-3} & \shd{civ-2} & \shd{civ-1} & \shd{civ-4} \\
\hline\hline\end{array}
$$
The final contribution of this type of diagrams to $A_{1,2}$ is
equal to\\
$$-2\ \risS{-12}{aiv-3}{}{25}{0}{15}\
  -2\ \risS{-12}{aiv-2}{}{25}{0}{0}\ .
$$

The other three types of degree 3 diagrams differ by the location of the
base point. A similar consideration shows that 5 out of the total of 12
Gauss diagrams of these types, namely
$$\risS{-10}{aiii-n1}{}{27.5}{15}{15}\ ,\quad
  \risS{-10}{aiii-n2}{}{27.5}{0}{0}\ ,\qquad
  \risS{-13}{ai-n}{}{25}{0}{0}\ ,\qquad
  \risS{-10}{aii-n1}{}{27.5}{0}{0}\ ,\quad
  \risS{-10}{aii-n2}{}{27.5}{0}{0}
$$
do not have ascending states with $|S|=2,3$. The remaining possibilities,
together with their ascending states, are shown in the table:
$$\begin{array}[t]{||c|c|c||c|c|c||c|c|c||}\hline\hline
 \risS{-17}{aiii-1}{}{27.5}{10}{25}& \risS{-17}{aiii-2}{}{27.5}{0}{0}
      & \risS{-17}{aiii-2}{}{27.5}{0}{0}
& \fhd{ai-1} & \fhd{ai-2} & \fhd{ai-3} &
 \risS{-17}{aii-2}{}{27.5}{0}{0}& \risS{-17}{aii-1}{}{27.5}{0}{0}
      & \risS{-17}{aii-3}{}{27.5}{0}{0} \\
 \risS{-9}{ciii-1}{}{27.5}{15}{15}& \risS{-9}{ciii-2}{}{27.5}{0}{0}
      & \risS{-9}{ciii-3}{}{27.5}{0}{0}
& \shd{ci-1} & \shd{ci-2} & \shd{ci-3} &
 \risS{-9}{cii-2}{}{27.5}{0}{0}& \risS{-9}{cii-1}{}{27.5}{0}{0}
      & \risS{-9}{cii-3}{}{27.5}{0}{0} \\
\hline\hline\end{array}
$$

The contribution of this type of diagrams to $A_{1,2}$ is thus
equal to
$$-2\ \risS{-10}{aiii-1}{}{27.5}{13}{15}
  -2\ \risS{-13}{ai-1}{}{25}{0}{0}
  -2\ \risS{-13}{ai-2}{}{25}{0}{0}
  +2\ \risS{-13}{ai-3}{}{25}{0}{0}
  -2\ \risS{-10}{aii-2}{}{27.5}{0}{0}\ .
$$

Apart form diagrams of degree 3, some degree 2 diagrams contribute to $A_{1,2}$ as
well. Since $|A|=2<k+l=3$, contributions of 2-diagrams depend also on
their signs. Such diagrams must be ascending (since $|S|=|A|=2$) and should
not have isolated arrows. There are four such diagrams:
$\lhd{ad2}$\!\!,\vspace{4pt} with all choices of the signs $\e_1$, $\e_2$
for the arrows. For each choice we have $\plwa=\e_1\e_2 e^{-(\e_1+\e_2)h}z^2$. If
$\e_1=-\e_2$, then $\plwa=-z^2$, so the coefficient of $hz^2$ vanishes and
such diagrams do not occur in $A_{1,2}$. For the two remaining diagrams with
$\e_1=\e_2=\pm$, the coefficients of $hz^2$ in $\plwa$ are equal to $\mp2$
respectively.

Combining all the above contributions, we finally get
$$A_{1,2} = -2\Bigl( \risS{-12}{aiv-3}{}{25}{0}{15}+
   \risS{-12}{aiv-2}{}{25}{0}{0}+\risS{-10}{aiii-1}{}{27.5}{13}{15}
  +\risS{-13}{ai-1}{}{25}{0}{0}+\risS{-13}{ai-2}{}{25}{0}{0}
  -\risS{-13}{ai-3}{}{25}{0}{0}+\risS{-10}{aii-2}{}{27.5}{0}{0}
    + \shd{ad2pp} - \shd{ad2mm}\Bigr)\ .
$$

At this point we can see the difference between virtual and classical long
knots.  For classical knots the invariant
$\I_{A_{1,2}}=\scp{A_{1,2}}{\cdot}$ can be
simplified further. Note that any classical Gauss diagram $G$ satisfies
$$\scp{ \risS{-13}{ai-2}{}{25}{17}{13}}{G}=
 \scp{ \risS{-13}{ai-3}{}{25}{0}{13} }{G}.$$
This follows from the symmetry of the linking number. Indeed, suppose we
have matched two vertical arrows (which are the same in both diagrams) with
two arrows of $G$. Let us consider the orientation preserving smoothings of
the corresponding two crossings of the link diagram $D$ associated with $G$.
The smoothed diagram $\widetilde{D}$ will have three components. Matchings
of the horizontal arrow of our Gauss diagrams with an arrow of $G$ both
measure the linking number between the first and the third components of
$\widetilde{D}$, using crossings when the first component passes over
(respectively, under) the third one.\vspace{4pt} Thus, as functions on
classical Gauss diagrams, $\langle,\risS{-13}{ai-2}{}{25}{17}{5}\, ,\, \cdot\rangle $ is equal to
$\langle\risS{-13}{ai-3}{}{25}{17}{5}\, ,\, \cdot\rangle$\and we have
$$p_{1,2}(G) = -2\langle \risS{-12}{aiv-3}{}{25}{0}{15}+
   \risS{-12}{aiv-2}{}{25}{0}{0}+\risS{-10}{aiii-1}{}{27.5}{13}{15}
  +\risS{-13}{ai-1}{}{25}{0}{0}+\risS{-10}{aii-2}{}{27.5}{0}{0}
    + \shd{ad2pp} - \shd{ad2mm}\ , G\rangle\ .
$$
For virtual Gauss diagrams this is no longer true.

In a similar manner one may check that $A_{3,0}=-4 A_{1,2}$.

The obtained result implies one more formula for the invariant $j_3$
(compare it with the two other formulae given in Section \ref{deg_three}).
Indeed, $j_3=-p_{3,0}-p_{1,2}=3p_{1,2}$, therefore
\def\shd#1{\ \risS{-12}{#1}{}{25}{15}{15}}
$$j_3(K) = -6\langle \risS{-12}{aiv-3}{}{25}{0}{15}+
   \risS{-12}{aiv-2}{}{25}{0}{0}+\risS{-10}{aiii-1}{}{27.5}{13}{15}
  +\risS{-13}{ai-1}{}{25}{0}{0}+\risS{-10}{aii-2}{}{27.5}{0}{0}
    + \shd{ad2pp} - \shd{ad2mm}\ , G\rangle\ .
$$

\section{The Jones polynomial via Gauss diagrams}

Apart from the Gauss diagram formulae as understood in this chapter, there are many other ways to extract Vassiliev (and other) knot invariants from Gauss diagrams. Here is just one example:
a description of the Jones polynomial (which is essentially a reformulation of the construction from a paper \cite{Zul} by L.~Zulli.) The reader should compare it to the definition of $\so_N$-weight system in Section~\ref{ws_so_n_St}.

Let $G$ be a Gauss diagram representing a knot $K$. Denote by $[G]$
the set of arrows of $G$. The sign  of an arrow $c\in [G]$ can be
considered as a value of the function $\sign: [G]\to \{-1,+1\}$.
A {\em state} $s$ for $G$ is an arbitrary function $s: [G]\to \{-1,+1\}$;
in particular, for a Gauss diagram with $n$ arrows there are $2^n$ states. The function
$\sign(\cdot)$ is one of them. With each state $s$ we associate an immersed
plane curve in the following way. Double every chord $c$ according to
the rule:
$$\risS{-70}{doublech}{
      \put(14,40){\mbox{$\scriptstyle c$}}
      \put(110,65){\mbox{,\qquad if $s(c)=1$}}
      \put(110,15){\mbox{,\qquad if $s(c)=-1$}}
             }{100}{20}{70}$$
Let $|s|$ denote the number of connected components of the curve obtained
by doubling all the chords of $G$. Also, for a state $s$ we
define an integer
$$ p(s) := \sum_{c\in [G]} s(c)\cdot\sign(c)\,.$$

The defining relations for the Kauffman bracket from Section~\ref{kauf_br}
lead to the following expression for the Jones polynomial.

\begin{xtheorem}
$$J(K)=(-1)^{w(K)}t^{3w(K)/4} \sum_s t^{-p(s)/4}
          \bigl(-t^{-1/2}-t^{1/2}\bigr)^{|s|-1}\ ,$$
where the sum is taken over all $2^n$ states for $G$ and
$\displaystyle w(K) = \sum_{c\in [G]} \sign(c)$ is the writhe of $K$.
\end{xtheorem}
This formula can be used to extend the Jones polynomial to virtual knots.

\begin{xexample}
For the left trefoil knot $3_1$ we have the following Gauss diagram.
$$\risS{-18}{31}{\put(19,-5){\mbox{$\scriptstyle 1$}}
                 \put(2,25){\mbox{$\scriptstyle 3$}}
                 \put(35,23){\mbox{$\scriptstyle 2$}}
             }{40}{20}{25}\qquad\qquad
G=\quad\risS{-18}{Gd31}{\put(8,-4){\mbox{$\scriptstyle 1$}}
                        \put(32,38){\mbox{$\scriptstyle 1$}}
                        \put(30,-4){\mbox{$\scriptstyle 2$}}
                        \put(4,38){\mbox{$\scriptstyle 2$}}
                        \put(-7,19){\mbox{$\scriptstyle 3$}}
                        \put(42,19){\mbox{$\scriptstyle 3$}}
                        \put(15,7){\mbox{$-$}}
                        \put(5,22){\mbox{$-$}}
                        \put(27,22){\mbox{$-$}}
             }{40}{20}{25}\qquad\qquad
w(G)=-3$$
There are eight states for such a diagram. Here are the corresponding
curves and numbers $|s|$, $p(s)$.
\def\csps#1#2#3{\begin{array}{c}
                 \ig[width=1.5cm]{curGDtr#1.eps}\\
                 \scriptstyle |s|=#2\\
                 \scriptstyle p(s)=#3
                \end{array}}
$$\csps{1}{2}{-3}\qquad \csps{2}{1}{-1}\qquad \csps{3}{1}{-1}
   \qquad \csps{4}{1}{-1}
$$
$$\csps{5}{2}{1}\qquad \csps{6}{2}{1}\qquad \csps{7}{2}{1}
   \qquad \csps{8}{3}{3}
$$

Therefore,
$$\begin{array}{ccl}
J(3_1)&=&
-t^{-9/4}\Bigl( t^{3/4}\bigl(-t^{-1/2}-t^{1/2}\bigr) + 3t^{1/4}
    +3t^{-1/4}\bigl(-t^{-1/2}-t^{1/2}\bigr) \\
&&\hspace{40pt}
    + t^{-3/4}\bigl(-t^{-1/2}-t^{1/2}\bigr)^2 \Bigr) \\
&=&
-t^{-9/4}\Bigl( -t^{1/4} -t^{5/4} 
    -3t^{-3/4} 
    +t^{-3/4}\bigl(t^{-1}+2+t\bigr) \Bigr) \\
&=& t^{-1}+t^{-3}-t^{-4}\ ,
\end{array}
$$
as we had before in Chapter 2.
\end{xexample}

\begin{xcb}{Exercises}

\begin{enumerate}

\item
Gauss diagrams\index{Gauss diagram!for links}
and Gauss diagram formulae may be defined for links in the same way as for knots.
Prove that for a link $L$ with two components $K_1$ and $K_2$
$$lk(K_1,K_2) = \langle \risS{-8}{adln}{}{45}{10}{20},G(L)\rangle\,.
\index{Linking number}$$

\item\label{pr:c2-on-B}
Find a sequence of Reidemeister moves that transforms the Gauss diagram $B$ from page \pageref{eq:gauss-c2-B} to the diagram
$$\risS{-10}{gauss-c2-4ccrr}{}{75}{20}{6}\ .
$$
Show that this diagram is not realizable. Calculate the value of the
extension, according to \ref{ss:constr-c}, of the invariant $c_2$ on it.

\item
Let $\vec\A$ be the space of arrow (oriented chord) diagrams modulo the 6T relations, see page~\pageref{arrown}. Show that the map $\A\to\vec\A$ which sends a chord diagram to the sum of all the arrow diagrams obtained by putting the orientations on the chords is well-defined. In other words, show that the 6T relation implies the 4T relation.

\item\neresh
Construct analogues of the algebras of closed and open Jacobi diagrams $\F$ and $\B$
consisting of diagrams with oriented edges. (It is known how to do it in the case of closed diagrams with acyclic internal graph, see  \cite{Po1}. )

\end{enumerate}
\end{xcb}
 %13 Gauss
\chapter{Miscellany}%14
\label{chapMisc}

\section{The Melvin--Morton Conjecture}\label{melmor}

\subsection{Formulation}
Roughly speaking, the Melvin--Morton Conjecture says that the
Alexander-Conway polynomial can be read off the highest order part
of the coloured Jones polynomial.

According to Exercise~\ref{ex_sl_2_col_Jones} to Chapter~\ref{LAWS}
(see also \cite{MeMo,BNG}) the coefficient $J^k_n$ of the unframed
coloured Jones polynomial $J^k$ (Section \ref{colJones}) is a
polynomial in $k$, of degree at most $n+1$ and without constant
term. So we may write
$$\frac{J^k_n}{k} = \sum_{0\leq j\leq n} b_{n,j}k^j
       \qquad\mbox{and}\qquad  \frac{J^k}{k} =
 \sum_{n=0}^\infty\quad \sum_{0\leq j\leq n} b_{n,j}k^j h^n\,,$$
where $b_{n,j}$ are Vassiliev invariants of order $\leqslant n$. The
highest order part of the coloured Jones polynomial is a Vassiliev
power series invariant \index{Jones polynomial!coloured!highest part}
$$\MM := \sum_{n=0}^\infty b_{n,n}h^n\,.\label{mel-mor}$$

\medskip

\noindent {\bf The Melvin--Morton Conjecture.} (\cite{MeMo})
\index{Conjecture!Melvin--Morton} {\em The highest order part of the
coloured Jones polynomial $\MM$ is inverse to the Alexander--Conway
power series $A$ defined by equations
(\ref{eq_skein_AC}-\ref{eq_init_AC}). In other words,}
$$\MM(K)\cdot A(K) = 1$$
{\em for any knot $K$.}

\subsection{Historical remarks}

In \cite{Mo} H.~Morton proved the conjecture for torus knots. After
this L.~Ro\-zansky \cite{Roz1} proved the Melvin--Morton Conjecture
on the level of rigour of Witten's path integral interpretation for
the Jones polynomial. The first complete proof was carried out by
D.~Bar-Natan and S.~Garoufalidis \cite{BNG}. They invented a
remarkable reduction of the conjecture to a certain identity on the
corresponding weight systems via canonical invariants; we shall
review this reduction in Section~\ref{MM_red_to_ws}. This identity
was then verified by evaluating the weight systems on chord
diagrams. In fact, Bar-Natan and S.~Garoufalidis  proved a more
general theorem in \cite{BNG} that relates the highest order part of
an arbitrary quantum invariant to the Alexander-Conway polynomial.
Following \cite{Ch2} we shall present another proof of this
generalized Melvin--Morton Conjecture in Section~\ref{MM_for_LA}.
A.~Kricker, B.~Spence and I.~Aitchison \cite{KSA} proved the
Melvin--Morton Conjecture using the cabling operations. Their work
was further generalized in \cite{Kri1}  by A.~Kricker.  Yet another
proof of the Melvin--Morton conjectures appeared in the paper
\cite{Vai1} by A.~Vaintrob. He used calculations on chord diagrams
and the Lie superalgebra $\gl(1|1)$ which gives rise to the
Alexander--Conway polynomial. The idea to use the restriction of the
aforementioned identity on weight systems to the primitive space was
explored in \cite{Ch1,Vai2}. We shall follow \cite{Ch1} in the
direct calculation of the Alexander--Conway weight system in
Section~\ref{AC_ws}.

B.~I.~Kurpita and K.~Murasugi found a different proof of the
Melvin--Morton Conjecture which does not use Vassiliev invariants
and weight systems \cite{KuM}.

Among other things, the works on the Melvin--Morton Conjecture
inspired L.~Rozansky to state his Rationality Conjecture that
describes the fine structure of the Kontsevich integral. This
conjecture was proved by A.~Kricker, and is the subject of
Section~\ref{Roz_r_conj}.)

\subsection{Reduction to weight systems}\label{MM_red_to_ws}

Since both power series Vassiliev invariants $\MM$ and $A$ are
canonical, so is their product (see Exercise~\ref{ex_prod_can} to
Chapter~\ref{advKI}). The constant invariant which is identically
equal to 1 on all knots is also a canonical invariant. We see that
the Melvin--Morton Conjecture states that two canonical invariants
are equal, and it is enough to prove that their symbols coincide.

Introduce the notation\vspace{-5pt}
$$\begin{array}{rcl}
  S_{\MM} &:=&\displaystyle \symb(\MM) =
     \sum_{n=0}^\infty \symb(b_{n,n})\,; \\
  S_A &:=&\displaystyle \symb(A) =
     \symb(\CP) =\sum_{n=0}^\infty \symb(c_n)\,.
\end{array}\label{S_MM}\label{S_A}
$$

The Melvin--Morton Conjecture is equivalent to the relation
$$ S_{\MM}\cdot S_A = \bo_0\,.$$
This is obvious in degrees 0 and 1. So, basically, we must prove
that in degree $\geqslant 2$ the product $S_{\MM}\cdot S_A$ equals
zero. In order to show this we have to establish that $S_{\MM}\cdot
S_A$ vanishes on any product $p_1\cdot\ldots\cdot p_n$ of primitive
elements of degree $>1$.

The weight system $S_{\MM}$ is the highest part of the weight system
$\f'^{V_k}_{\sL_2}/k$ from Exercise~\ref{ex_sl_2_col_Jones} to
Chapter~\ref{LAWS}. The latter is multiplicative as we explained in
Section~\ref{ws_repres}; hence, $S_{\MM}$ is multiplicative too.
Exercise~\ref{ex_consymb} to Chapter~\ref{FT_inv} implies then that
the weight system $S_A$ is also multiplicative. In other words, both
weight systems $S_{\MM}$ and $S_A$ are group-like elements of the
Hopf algebra of weight systems $\ol{\W}$. A product of two
group-like elements is group-like which shows that the weight system
$S_{\MM}\cdot S_A$ is multiplicative. Therefore, it is sufficient to
prove that
$$ S_{\MM}\cdot S_A \bigr|_{\PR_{>1}} = 0\,. $$

By the definitions of the weight system product and of a primitive
element
$$S_{\MM}\cdot S_A(p) = (S_{\MM}\ot S_A)(\d(p))=
  S_{\MM}(p) + S_A(p)\,.$$
Therefore, we have reduced the Melvin--Morton Conjecture to the
equality
$$ S_{\MM}\bigr|_{\PR_{>1}} + S_A\bigr|_{\PR_{>1}} = 0\,.$$

Now we shall exploit the filtration
$$0=\PR_n^1\subseteq  \PR_n^2\subseteq \PR_n^3\subseteq \dots \subseteq
    \PR_n^n=\PR_n\ .
$$
from Section \ref{filtr_pr}. Recall that the wheel $\ol{w_n}$ spans
$\PR_n^n / \PR_n^{n-1}$ for even $n$ and belongs to $\PR_n^{n-1}$
for odd $n$.

The Melvin--Morton Conjecture is a consequence of the following
theorem.

\begin{theorem}\label{MM_comb_th}
The weight systems $S_{\MM}$ and $S_A$ satisfy
\begin{enumerate}
\item $S_{\mathrm{MM}}\bigr|_{\PR^{n-1}_n} = S_A \bigr|_{\PR^{n-1}_n} = 0$;
\item $S_{\mathrm{MM}}(\ol{w_{2m}}) = 2,\quad  S_A(\ol{w_{2m}}) = -2$.
\end{enumerate}
\end{theorem}

The proof is based on several exercises to Chapter~\ref{LAWS}.

First, let us consider the weight system $S_{\mathrm{MM}}$.
Exercise~\ref{ex_sl_2_on_Pkn} implies that for any $D\in
\PR^{n-1}_n$ the weight system $\f_{\sL_2}(D)$  is a polynomial in
$c$ of degree less than or equal to $[(n-1)/2]$. The weight system
of the coloured Jones polynomial is obtained from $\f_{\sL_2}$ by
fixing the representation $V_k$ of $\sL_2$ and deframing. Choosing
the representation $V_k$ means that we have substitute
$c=\frac{k^2-1}{2}$; the degree of the polynomial
$\f^{V_k}_{\sL_2}(D)/k$ in $k$ will be at most $n-1$. Therefore, its
$n$th degree term vanishes and $S_{\MM}\bigr|_{\PR^{n-1}_n} = 0$.
According to Exercise~\ref{ex_sl_2_on_C_wheel}, the highest degree
term of the polynomial $\f_{\sL_2}(\ol{w_{2m}})$ is $2^{m+1}c^m$.
Again, the substitution $c=\frac{k^2-1}{2}$ (taking the trace of the
corresponding operator and dividing the result by $k$) gives that
the highest degree term of $\f^{V_k}_{\sL_2}(\ol{w_{2m}})/k$ is
$\frac{2^{m+1}k^{2m}}{2^m}=2k^{2m}$, and, hence
$S_{\MM}(\ol{w_{2m}}) = 2$.

In order to treat the weight system $S_A$ we use
Exercise~\ref{ex_conw_on_C}, which contains the equality
$S_A(\ol{w_{2m}}) = -2$ as a particular case. It remains to prove
that $S_A \bigr|_{\PR^{n-1}_n} = 0$.

\subsection{Alexander--Conway weight system}\label{AC_ws}

Using the state sum formula for $S_A$ from
Exercise~\ref{ex_conw_on_C} to Chapter \ref{LAWS} we shall prove
that $S_A(D)=0$ for any closed diagram $D\in\PR^{n-1}_n$.

First of all note that any such $D\in\PR^{n-1}_n$ has an internal
vertex which is not connected to any leg by an edge. Indeed, each
leg is connected with only one internal vertex. The diagram $p$ has
at most $n-1$ legs and $2n$ vertices in total, so there must be at
least $n+1$ internal vertices, and only $n-1$ of them can be
connected with legs.

Pick such a vertex connected only with other internal vertices.
There are two possible cases: either all these other vertices are
different or two of them coincide.

Let us start with the second, easier, case. Here we have a
``bubble''
\def\lpb#1{\raisebox{0pt}[15pt][10pt]{ \begin{picture}(40,0)(1,3)
      \put(0,0){\epsfxsize=40pt \epsfbox{#1.eps}} \end{picture}} }
\def\pb#1{\raisebox{-6pt}{ \begin{picture}(60,0)(0,0)
      \put(0,0){\epsfxsize=50pt \epsfbox{#1.eps}} \end{picture}} }
$$\lpb{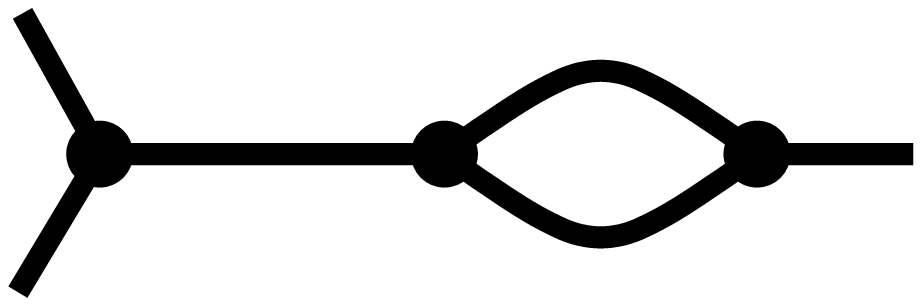}.$$
After resolving the vertices of this fragment according to the state
sum formula and erasing the curves with more than
one component we are left with the linear combination of curves
$$-2\pb{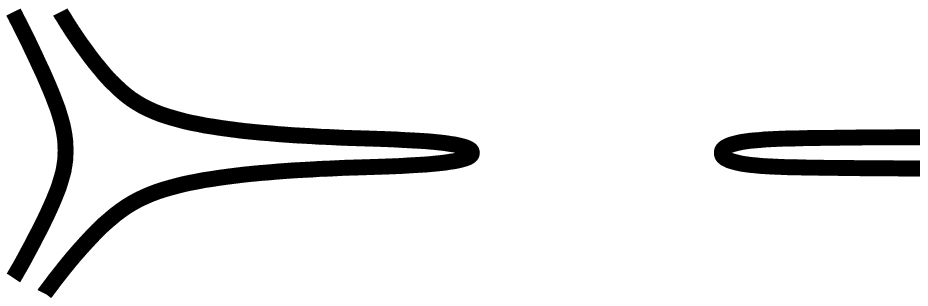} \qquad+\qquad\ 2\pb{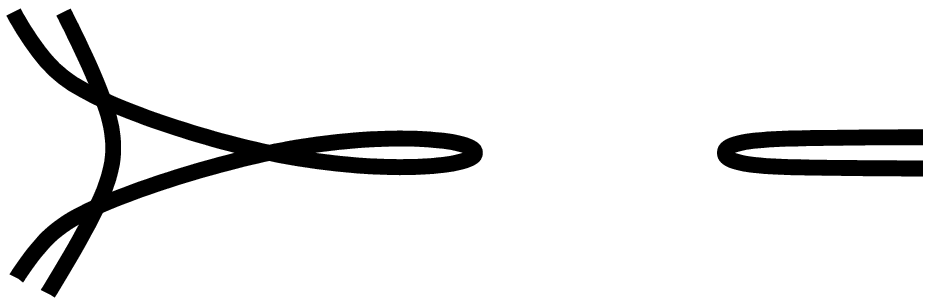}$$
which cancel each other, so $S_A(D)=0$.

\medskip
For the first case we formulate our claim as a lemma.
\def\tr#1{\raisebox{-15pt}{ \begin{picture}(40,0)(-3,0)
      \put(0,0){\epsfxsize=30pt \epsfbox{#1.eps}} \end{picture}} }
\begin{xlemma} \quad
$S_A \left( \tr{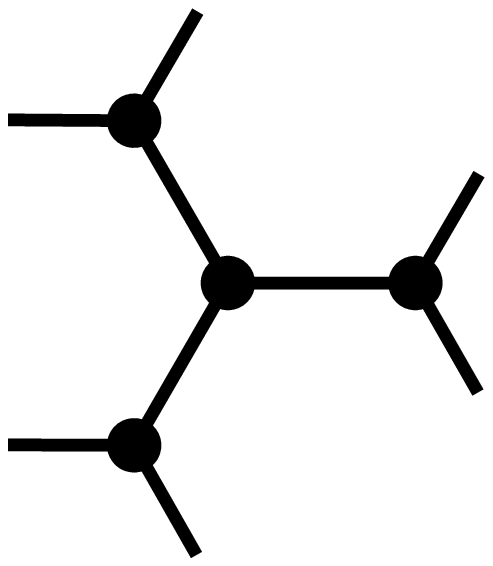} \right) = 0$.
\end{xlemma}

We shall utilize the state surfaces $\Sigma_s(D)$ from
Exercise~\ref{ex_gl_N_on_C_surface} to Chapter \ref{LAWS}. For a
given state, the neighbourhoods of ``$+$''- and ``$-$''-vertices of
look on the surface like three meeting bands:\vspace{8pt}
\def\nev#1#2{\raisebox{0pt}{ \begin{picture}(40,0)(2,19)
      \put(0,0){\epsfxsize=40pt \epsfbox{#1.eps}}
      \put(18,25){\mbox{$\scriptstyle #2$}}  \end{picture}} }
\def\bnev#1{\raisebox{0pt}[15pt][10pt]{ \begin{picture}(40,0)(2,19)
      \put(0,0){\epsfxsize=40pt \epsfbox{#1.eps}} \end{picture}} }
\def\drdr{\raisebox{0pt}{ \begin{picture}(25,0)(2,0)
      \put(0,0){\epsfxsize=25pt \epsfbox{totor.eps}} \end{picture}} }
\begin{equation}
\nev{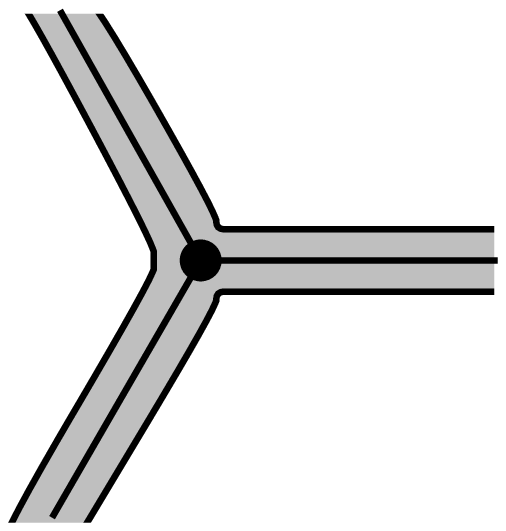}{+};\qquad \nev{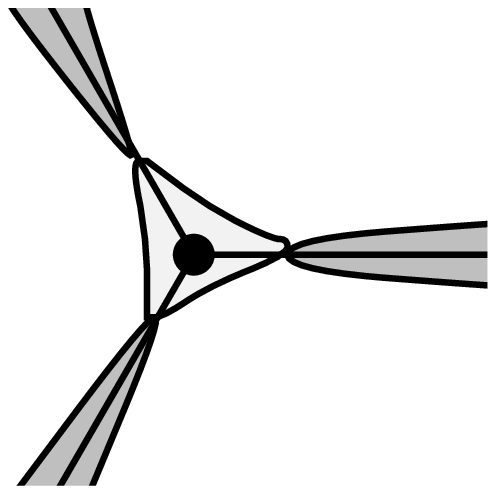}{-} =
 \bnev{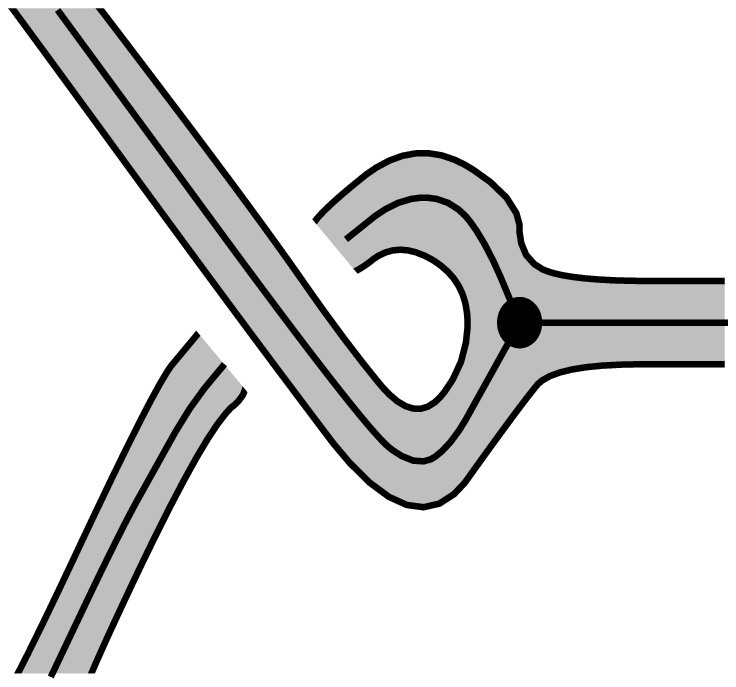} = \bnev{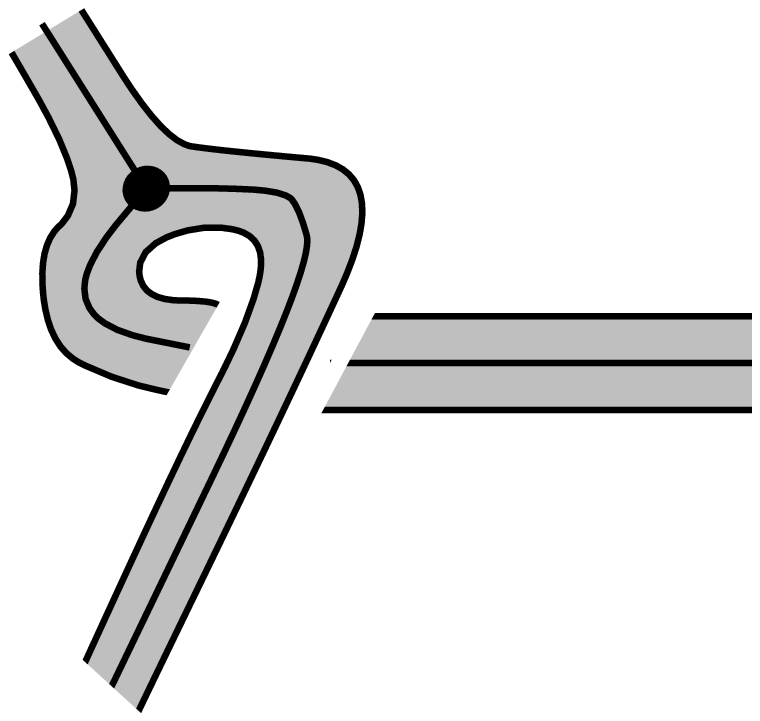} = \bnev{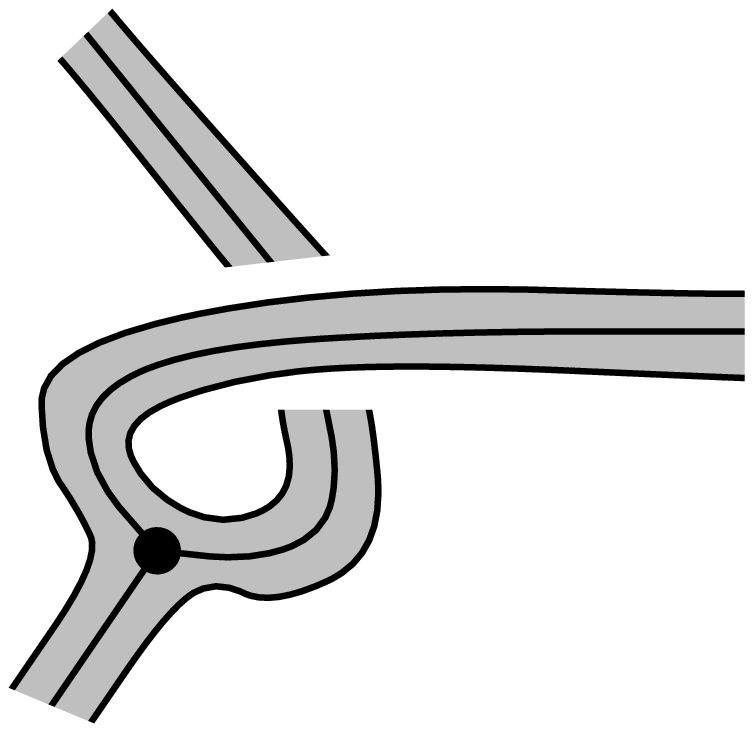}.
\label{eq:melmor_nev}
\end{equation}

Switching a marking (value of the state) at a vertex means reglueing
of the three bands along two chords on the surface: \vspace{8pt}
\def\sdr#1{\ \shortstack{\drdr\\ {\scriptsize #1}}\ }
\def\pcd#1{\raisebox{-13pt}[20pt]{ \begin{picture}(40,0)(2,0)
      \put(0,0){\epsfxsize=40pt \epsfbox{#1.eps}} \end{picture}} }
$$\pcd{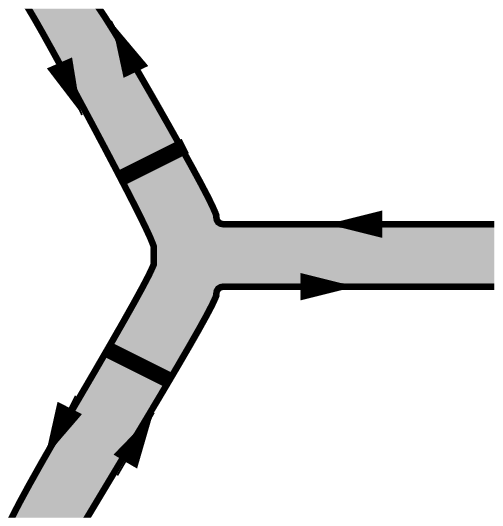} \sdr{cut along chords}
  \pcd{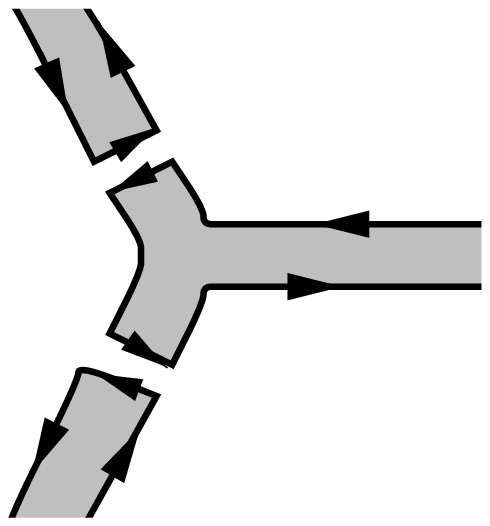} \sdr{interchange}
  \pcd{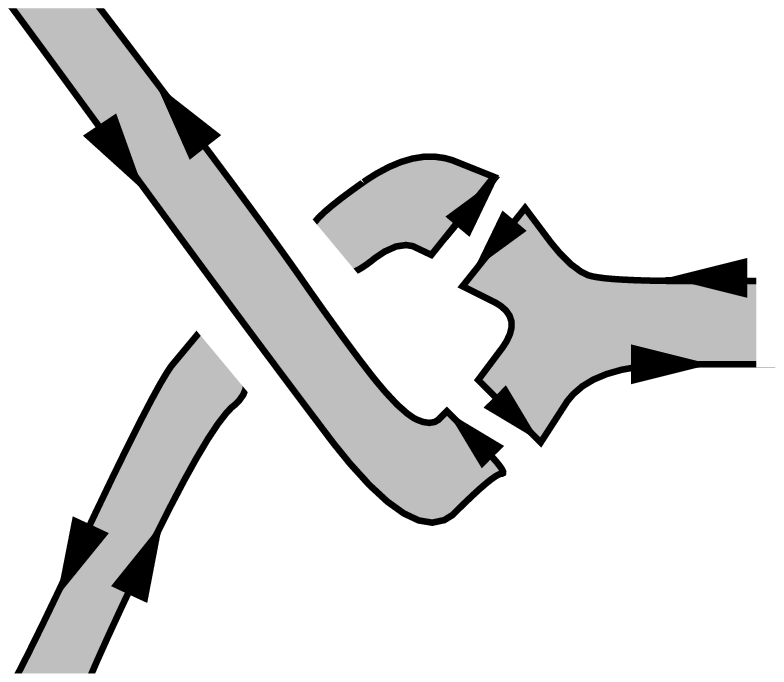} \sdr{glue} \pcd{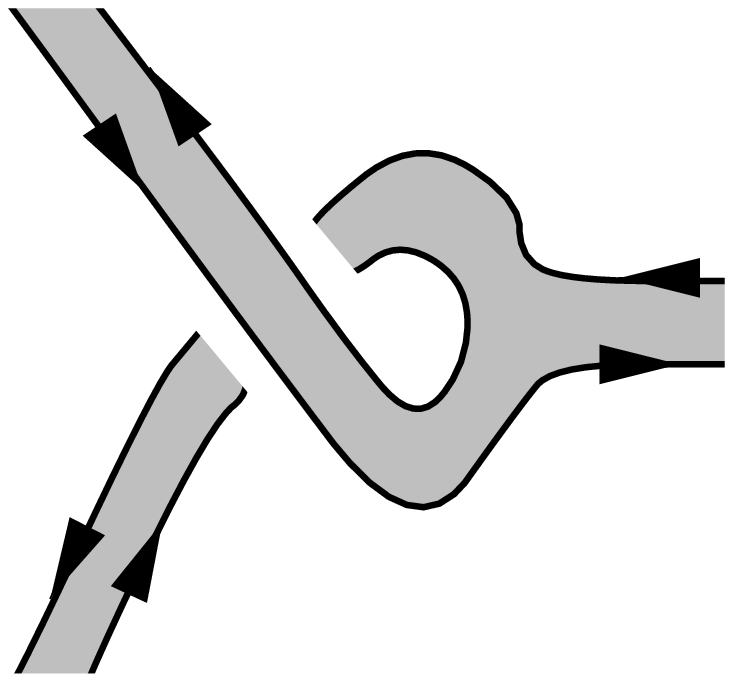}.$$

\begin{proof}
According to Exercise~\ref{ex_conw_on_C}, the symbol of the Conway
polynomial is the coefficient at $N$ in the polynomial
$\f^{St}_{\gl_N}$. In terms of the state surfaces this means that we
only have to consider the surfaces with one boundary component. We
are going to divide the set of all those states $s$ for which the
state surface $\Sigma_s(D)$ has one boundary component into pairs in
such a way that the states $s$ and $s'$ of the same pair differ by
an odd number of markings. The terms of the pairs will cancel each
other and will contribute zero to $S_A(D)$.

In fact, in order to do this we shall adjust only the markings of
the four vertices of the fragment pictured in the statement of the
Lemma. The markings $\e_1$, \dots, $\e_l$ and $\e'_1$, \dots,
$\e'_l$ in the states $s$ and $s'$ will be the same except for some
markings of the four vertices of the fragment. Denote the vertices
by $v$, $v_a$, $v_b$, $v_c$ and their markings in the state $s$ by
$\e$, $\e_a$, $\e_b$, $\e_c$, respectively.

Assume that $\Sigma_s(D)$ has one boundary component. Modifying the
surface as in (\ref{eq:melmor_nev}) we can suppose that the
neighbourhood of the fragment has the form\vspace{-10pt}
\begin{center}
    \begin{picture}(350,75)(0,0)
       \put(20,10){\epsfxsize=50pt \epsfbox{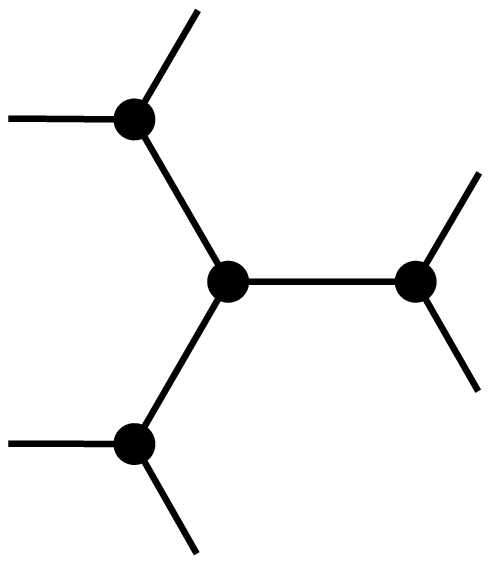}}
       \put(44,42){\mbox{$\scriptstyle v$}}
       \put(25,50){\mbox{$\scriptstyle v_b$}}
       \put(24,25){\mbox{$\scriptstyle v_c$}}
       \put(66,37){\mbox{$\scriptstyle v_a$}}
       \put(180,-8){\epsfxsize=80pt \epsfbox{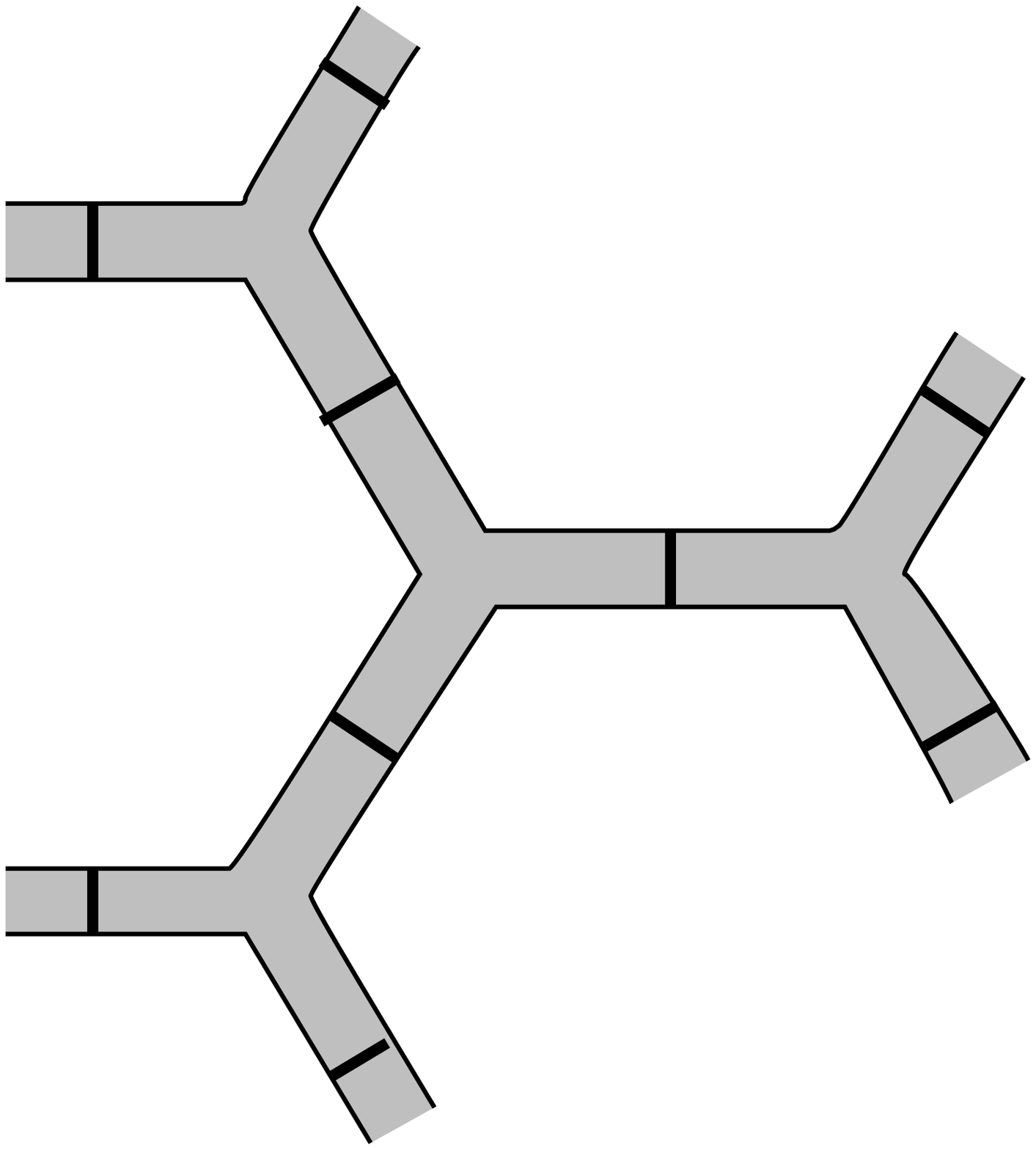}}
       \put(230,42){\mbox{$\scriptstyle a$}}
       \put(214,50){\mbox{$\scriptstyle b$}}
       \put(214,20){\mbox{$\scriptstyle c$}}
       \put(213,68){\mbox{$\scriptstyle b_1$}}
       \put(213,2){\mbox{$\scriptstyle c_2$}}
       \put(241,50){\mbox{$\scriptstyle a_1$}}
       \put(241,22){\mbox{$\scriptstyle a_2$}}
       \put(185,52){\mbox{$\scriptstyle b_2$}}
       \put(185,18){\mbox{$\scriptstyle c_1$}}
    \end{picture}
\end{center}
Draw nine chords $a$, $a_1$, $a_2$, $b$, $b_1$, $b_2$, $c$, $c_1$, $c_2$
on our surface as shown on the picture. The chords $a$, $b$, $c$ are
located near the vertex $v$; $a$, $a_1$, $a_2$ near the vertex $v_a$;
$b$, $b_1$, $b_2$ near $v_b$ and $c$, $c_1$, $c_2$ near $v_c$.

Since the surface has only one boundary component, we can draw this
boundary as a plane circle and $a$, $a_1$, $a_2$, $b$, $b_1$, $b_2$,
$c$, $c_1$, $c_2$ as chords inside it. Let us consider the possible
chord diagrams obtained in this way.

If two, say $b$ and $c$, of three chords located near a vertex, say
$v$, do not intersect, then the surface
$\Sigma_{\dots,-\e,\e_a,\e_b,\e_c, \dots}(D)$ obtained by switching
the marking $\e$ to $-\e$ also has only one boundary component.
Indeed, the reglueing effect along two non-intersecting chords can be
seen on chord diagrams as follows:
\def\pwh#1{\raisebox{-18pt}[30pt]{ \begin{picture}(40,0)(5,0)
      \put(0,0){\epsfxsize=50pt \epsfbox{#1.eps}} \end{picture}} }
\def\pcdt#1{\raisebox{-13pt}[20pt]{ \begin{picture}(30,0)(2,0)
      \put(0,0){\epsfxsize=40pt \epsfbox{#1.eps}} \end{picture}} }
$$\pcdt{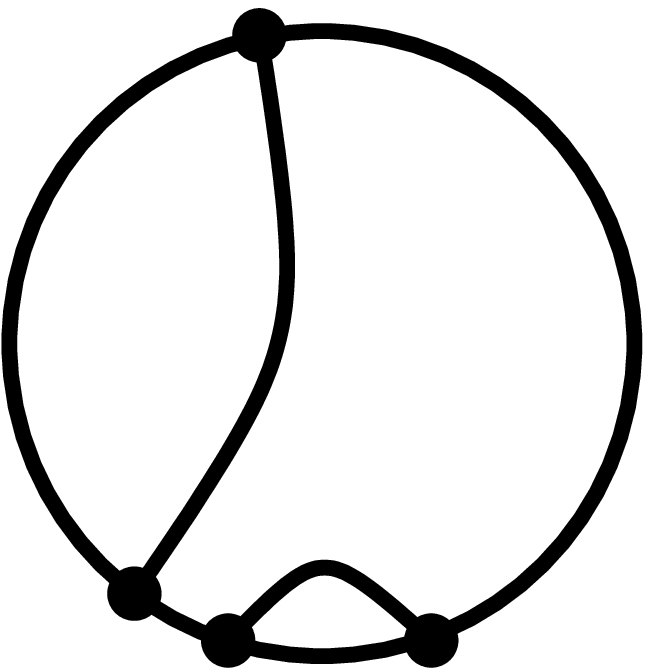}\quad \sdr{cut along chords}
  \pcdt{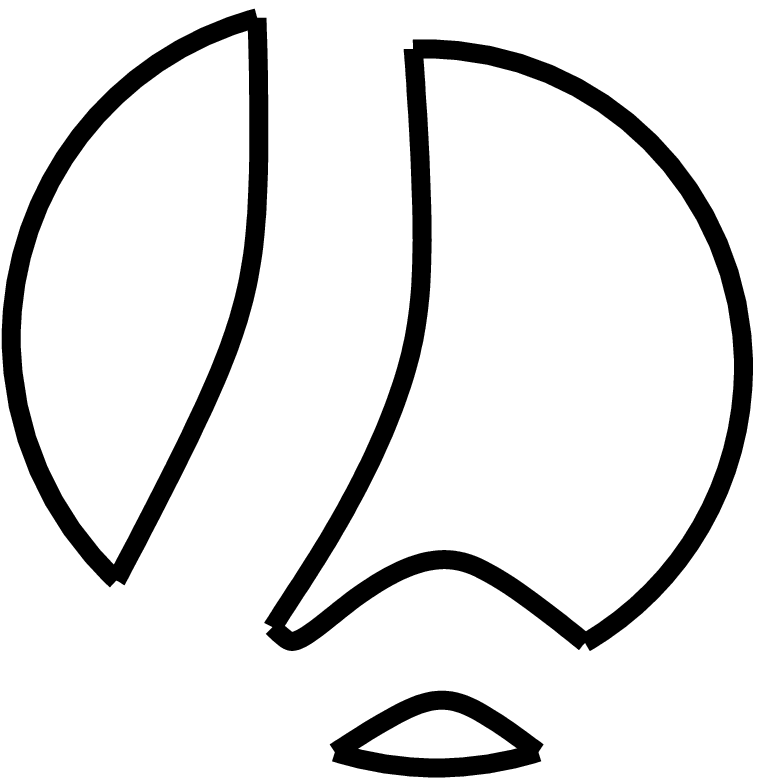}\quad \sdr{interchange}
  \pwh{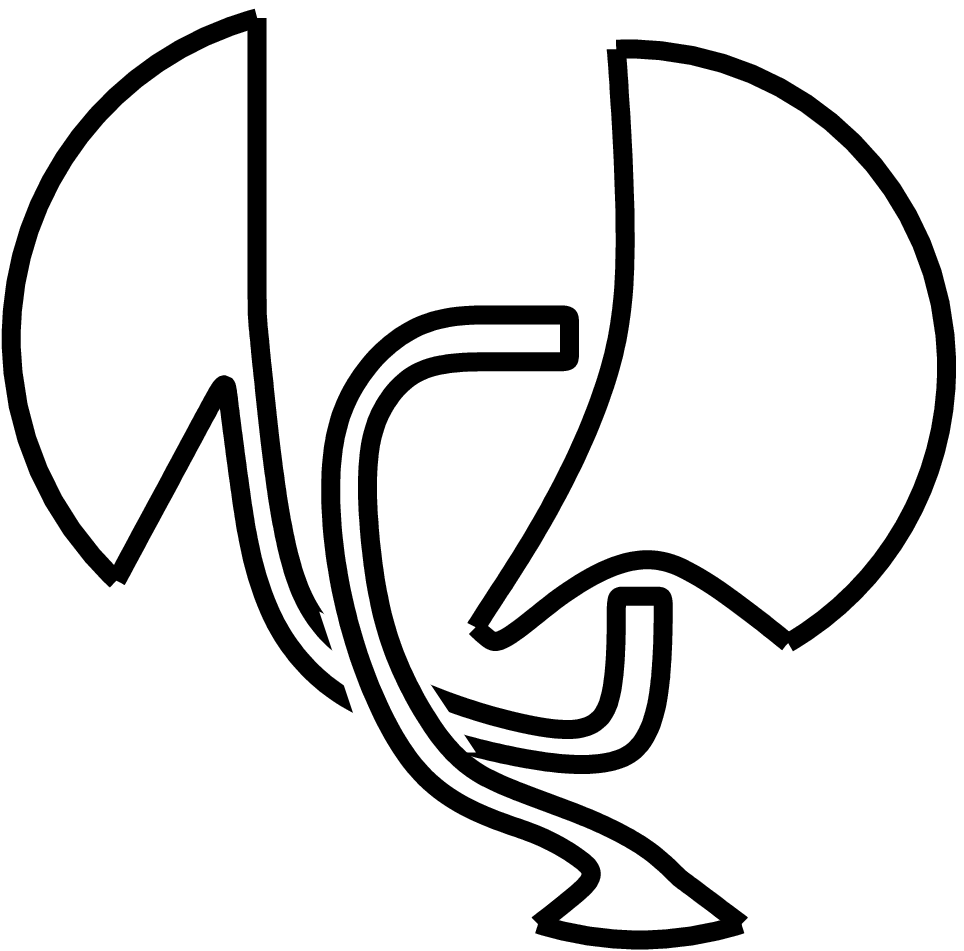}\quad \sdr{glue}\ \pwh{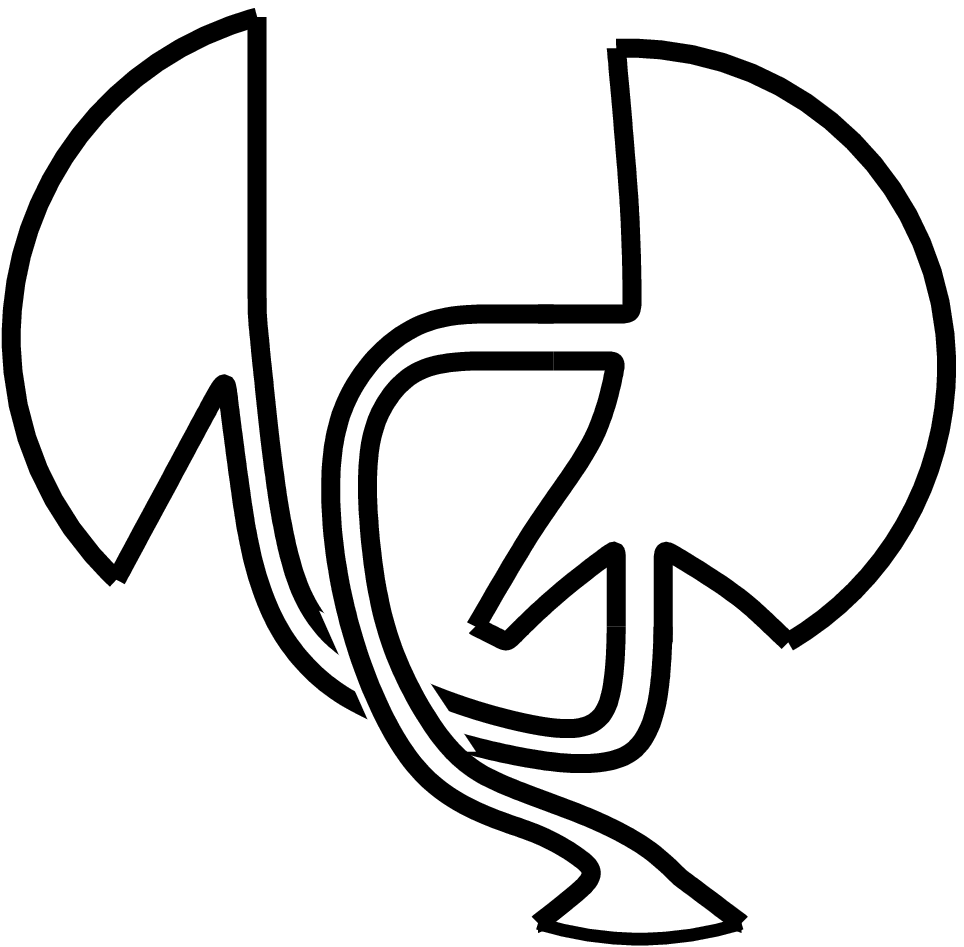}\ \ .$$
Therefore, in this case, the state $s=\{\dots,-\e,\e_a,\e_b,\e_c,
\dots\}$ should be paired with $s'=\{\dots,\e,\e_a,\e_b,\e_c,
\dots\}$.

We see that switching a marking at a vertex we increase the number
of boundary components (so that such a marked diagram may give a
non-zero contribution to $S_A(D)$) if and only if the three chords
located near the vertex intersect pairwise.

Now we can suppose that any two of the three chords
in each triple $(a, b, c)$, $(a, a_1, a_2)$, $(b, b_1, b_2)$,
$(c, c_1, c_2)$ intersect. This leaves us with only one possible chord
diagram:
\begin{center}
    \begin{picture}(80,80)(0,0)
       \put(0,5){\epsfxsize=80pt \epsfbox{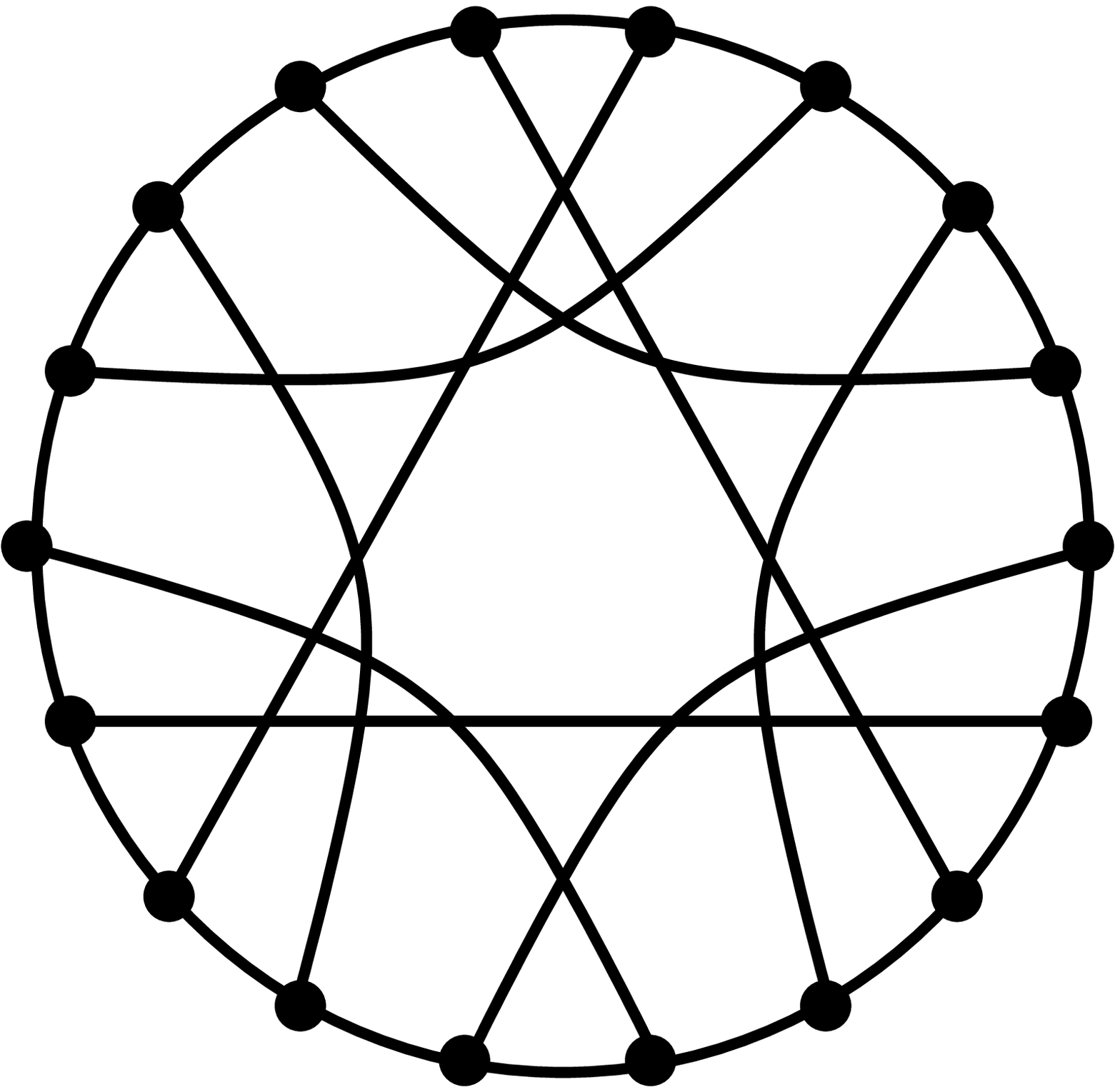}}
       \put(60,79){\mbox{$\scriptstyle a_1$}}
       \put(45,84){\mbox{$\scriptstyle a$}}
       \put(30,83){\mbox{$\scriptstyle b$}}
       \put(13,80){\mbox{$\scriptstyle b_1$}}
       \put(0,69){\mbox{$\scriptstyle a_2$}}
       \put(-6,55){\mbox{$\scriptstyle a_1$}}
       \put(-9,42){\mbox{$\scriptstyle c_2$}}
       \put(-2,29){\mbox{$\scriptstyle c$}}
       \put(5,15){\mbox{$\scriptstyle a$}}
       \put(14,5){\mbox{$\scriptstyle a_2$}}
       \put(30,0){\mbox{$\scriptstyle c_1$}}
       \put(45,0){\mbox{$\scriptstyle c_2$}}
       \put(61,5){\mbox{$\scriptstyle b_2$}}
       \put(71,15){\mbox{$\scriptstyle b$}}
       \put(79,29){\mbox{$\scriptstyle c$}}
       \put(81,42){\mbox{$\scriptstyle c_1$}}
       \put(79,55){\mbox{$\scriptstyle b_1$}}
       \put(72,69){\mbox{$\scriptstyle b_2$}}
   \end{picture}
\end{center}
The boundary curve of the surface connects the ends of our fragment
as in the left picture below.
\begin{center}
\begin{picture}(350,135)(0,0)
  \put(35,5){ \begin{picture}(100,120)(0,0)
       \put(-35,122){\mbox{$\Sigma_{\dots,\e,\e_a,\e_b,\e_c, \dots}(p)$}}
       \put(-20,-10){\epsfxsize=120pt \epsfbox{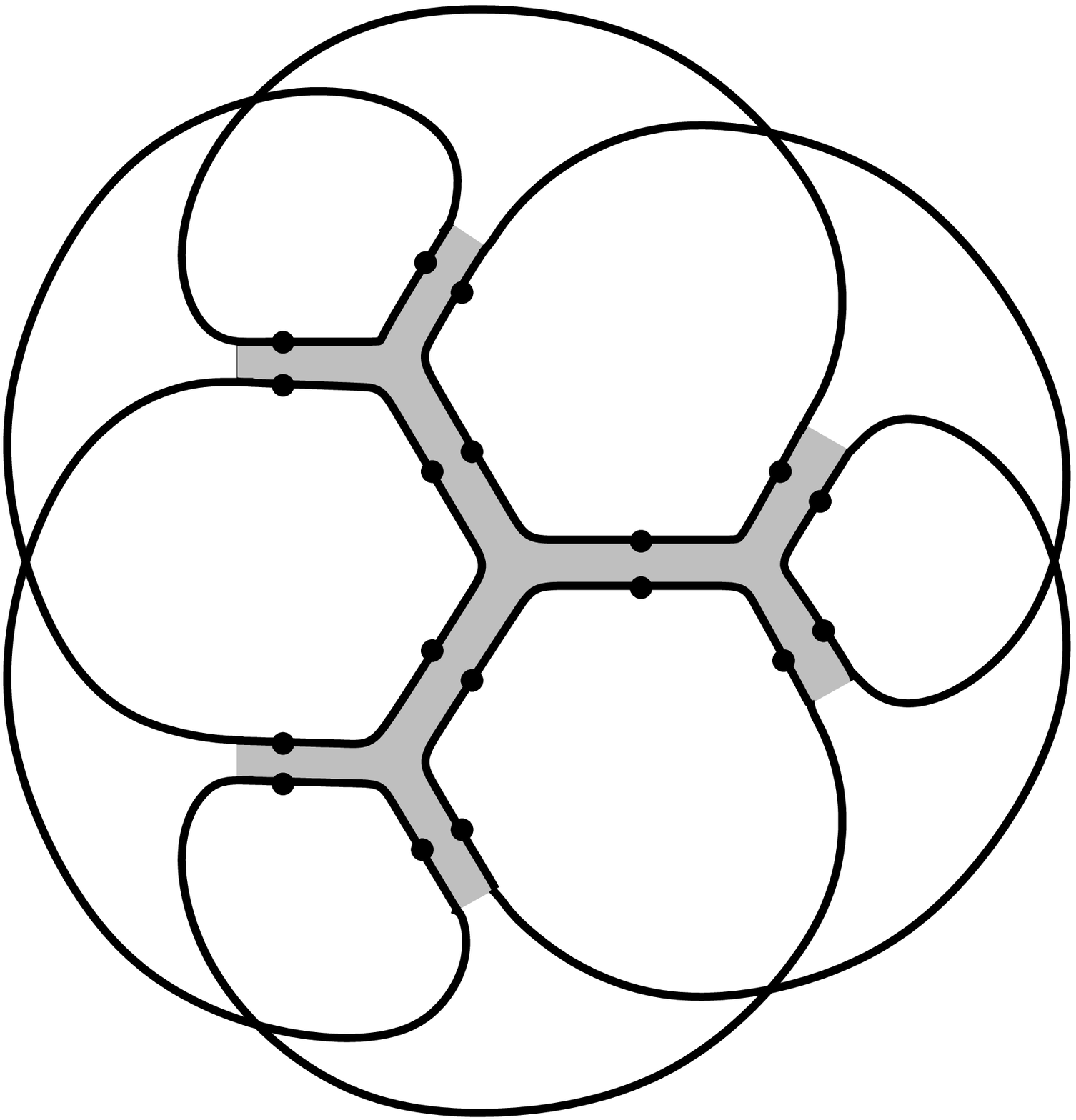}}
       \put(50,58){\mbox{$\scriptstyle a$}}
       \put(50,43){\mbox{$\scriptstyle a$}}
       \put(35,64){\mbox{$\scriptstyle b$}}
       \put(22,60){\mbox{$\scriptstyle b$}}
       \put(35,37){\mbox{$\scriptstyle c$}}
       \put(22,42){\mbox{$\scriptstyle c$}}
       \put(33,76){\mbox{$\scriptstyle b_1$}}
       \put(18,85){\mbox{$\scriptstyle b_1$}}
       \put(34,21){\mbox{$\scriptstyle c_2$}}
       \put(19,16){\mbox{$\scriptstyle c_2$}}
       \put(58,65){\mbox{$\scriptstyle a_1$}}
       \put(72,54){\mbox{$\scriptstyle a_1$}}
       \put(59,36){\mbox{$\scriptstyle a_2$}}
       \put(74,44){\mbox{$\scriptstyle a_2$}}
       \put(8,64){\mbox{$\scriptstyle b_2$}}
       \put(8,80){\mbox{$\scriptstyle b_2$}}
       \put(7,21){\mbox{$\scriptstyle c_1$}}
       \put(8,35){\mbox{$\scriptstyle c_1$}}
              \end{picture}}
  \put(235,5){ \begin{picture}(100,120)(0,0)
       \put(-35,122){\mbox{$\Sigma_{\dots,\e,-\e_a,-\e_b,-\e_c, \dots}(p)$}}
       \put(-20,-10){\epsfxsize=120pt \epsfbox{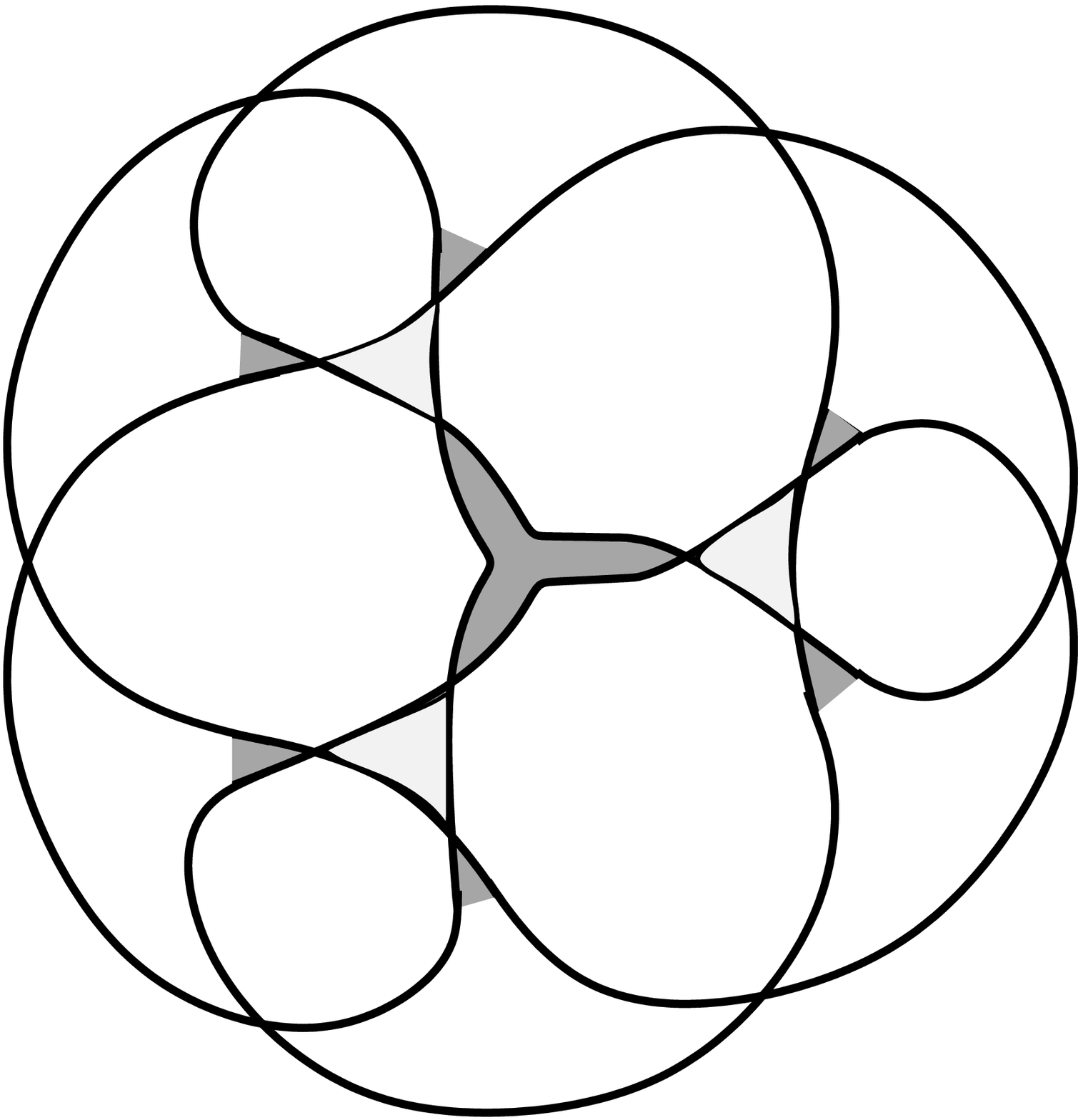}}
               \end{picture}}
\end{picture}
\end{center}
Switching markings at $v_a$, $v_b$, $v_c$ gives a surface which also
has one boundary component as in the right picture above. Pairing
the state $$s=\{\dots,\e,\e_a,\e_b,\e_c, \dots\}$$ up with
$s'=\{\dots,\e,-\e_a,-\e_b,-\e_c, \dots\}$ we get the desired
result.

The Lemma, and thus the Melvin--Morton Conjecture, is proved.
\end{proof}

\subsection{Generalization of the Melvin--Morton Conjecture to other
quantum invariants}\label{MM_for_LA}
\newcommand{\h}{\mathfrak{h}}
\renewcommand{\l}{\lambda}
\def\kf#1#2{{\langle #1,#2\rangle}}
Let $\g$ be a semi-simple Lie algebra and let $V_\l$ be an irreducible
representation of $\g$ of the highest weight $\l$.
Denote by $\h$ a Cartan subalgebra of $\g$, by $R$ the set of all
roots and by $R^+$ the set of positive roots.
Let $\kf{\cdot}{\cdot}$ be the scalar product on $\h^*$
induced by the Killing form.  These data define the unframed quantum invariant $\t^{V_\l}_\g$
which after the substitution $q=e^h$ and the expansion into a power
series in $h$ can be written as
$$ \t^{V_\l}_{\g} = \sum_{n=0}^\infty \t^{\l}_{\g,n}h^n,$$
see Section \ref{can_qi}.
\begin{xtheorem}[\cite{BNG}] \ \\

{\rm (1)} The invariant $\t^{\l}_{\g,n}/\dim(V_\l)$ is a polynomial
in
$\l$ of degree at most $n$.\\

{\rm (2)} Define the Bar-Natan--Garoufalidis function $BNG$ as a
power series in $h$ whose coefficient at $h^n$ is the degree $n$
part of the polynomial $\t^{\l}_{\g,n}/\dim(V_\l)$.
\index{Bar-Natan--Garoufalidis function}\label{BNG} Then for any
knot $K$,
$$BNG(K)\cdot \prod_{\a\in R^+}A_\a(K) = 1\,,$$
where $A_\a$ is the following normalization of the Alexander--Conway
polynomial:
$$\begin{array}{l}\displaystyle
  A_\a\Bigl(\lrints\Bigr)\ -\  A_\a\Bigl(\rlints\Bigr)\quad = \quad
  (e^\frac{\kf\l\a h}2-e^{-\frac{\kf\l\a h}2})A_\a\Bigl(\twoup\Bigr)\,;
     \vspace{10pt}\\ \displaystyle  A_\a\Bigl(\unkn\Bigr)\quad = \quad
   \frac{\kf\l\a h}{e^\frac{\kf\l\a h}2-e^{-\frac{\kf\l\a h}2}}\,.
\end{array}
$$

\end{xtheorem}

\begin{proof}
The symbol $S_{BNG}$ is the highest part (as a function of $\l$) of
the Lie algebra weight system $\f'^{V_\l}_\g$ associated with the
representation $V_\l$. According to Exercise~\ref{ex_ws_subst} on
page~\pageref{ex_ws_subst}, the symbol of $A_\a$ in degree $n$
equals $\kf\l\a^n \symb(c_n)$.

The relation between the invariants can be reduced to the following
relation between their symbols:
$$ S_{BNG}\bigr|_{\PR_n} +
   \sum_{\a\in R^+} \kf\l\a^n \symb(c_n)\bigr|_{\PR_n} = 0\,, \vspace{-10pt}
$$
for $n>1$.

As above, $S_{BNG}\bigr|_{\PR^{n-1}_n} = \symb(c_n)
\bigr|_{\PR^{n-1}_n} = 0$, and $\symb(c_n)(\ol{w_{2m}}) = -2$. Thus
it remains to prove that
$$S_{BNG}(\ol{w_{2m}}) =  2\sum_{\a\in R^+} \kf\l\a^{2m}\,.$$

To prove this equality we shall use the method of
Section~\ref{LAWS_C}. First, we  take the Weyl basis of $\g$ and
write the Lie bracket tensor $J$ in this basis.

Fix the root space decomposition
$\displaystyle \g = \h\op\bigl(\mathop{\op}\limits_{\a\in R}\g_\a\bigr)$.
The Cartan subalgebra $\h$ is orthogonal to all the
${\g_\a}{\mbox{'s}}$ and $\g_\a$ is orthogonal
to $\g_\b$ for $\b \not= -\a$.
Choose the elements $e_\a\in\g_\a$ and $h_\a=[e_\a,e_{-\a}]\in\h$
for each $\a\in R$ in such a way that
$\kf{e_\a}{e_{-\a}} = 2/\kf\a\a$, and for any $\l\in\h^*$,
$\l(h_\a) = 2\kf\l\a / \kf\a\a$.

The elements $\{h_\b,e_\a\}$, where $\b$ belongs to a basis $B(R)$
of $R$ and $\a\in R$, form the Weyl basis of $\g$. The Lie bracket
$[\cdot,\cdot]$
as an element of $\g^*\ot\g^*\ot\g$ can be written as follows:
$$\begin{array}{rcl}
[\cdot,\cdot] &=& \displaystyle{\sum_{\stackrel{\b\in B(R)}{\a\in R}}
\Bigl( h_\b^*\ot e_\a^*\ot\a(h_\b)e_\a - e_\a^*\ot h_\b^*\ot\a(h_\b)e_\a
\Bigr)} \vspace{10pt} \\
&& \!\!\!\!\displaystyle{+ \sum_{\a\in R} e_\a^*\ot e_{-\a}^*\ot h_\a
   + \sum_{\stackrel{\a,\gamma\in R}{\a+\gamma\in R}}
           e_\a^*\ot e_\gamma^*\ot N_{\a,\gamma}e_{\a+\gamma}\,,}
\end{array}$$
where the stars indicate elements of the dual basis. The second sum
is most important because the first and third sums give no
contribution to the Bar-Natan--Garoufalidis weight system $S_{BNG}$.

After identification of $\g^*$ and $\g$ via $\kf\cdot\cdot$
we get
$e_\a^* = (\langle\a,\a\rangle/2)e_{-\a}$. In particular, the second sum
of the tensor $J$ is
$$\sum_{\a\in R} \Bigl( \kf\a\a/2 \Bigr)^2 e_{-\a}\ot e_\a\ot h_\a\,.
$$

According to Section~\ref{LAWS_C}, in order to calculate
$S_{BNG}(\ol{w_{2m}})$ we must assign a copy of the tensor $-J$ to
each internal vertex, perform all the contractions corresponding to
internal edges and, after that, take the product
$\f^{V_\l}_\g(\ol{w_{2m}})$ of the all operators in $V_\l$
corresponding to the external vertices. We have that
$\f^{V_\l}_\g(\ol{w_{2m}})$ is a scalar operator of multiplication
by some constant. This constant is a polynomial in $\l$ of degree at
most $2m$; its part of degree $2m$ is $S_{BNG}(\ol{w_{2m}})$.

We associate the tensor $-J$ with an internal vertex in such a way
that the third tensor factor of $-J$ corresponds to the edge
connecting the vertex with a leg. After that we take the product of
operators corresponding to these external vertices. This means that
we take the product of operators corresponding to the third tensor
factor of $-J$. Of course, we are interested only in those operators
which are linear in $\l$. One can show (see, for example, \cite[
Lemma 5.1]{BNG}) that it is possible to choose a basis in the space
of the representation $V_\l$ in such a way that the Cartan operators
$h_\a$ and raising operators $e_\a$ ($\a\in R^+$) will be linear in
$\l$ while the lowering operators $e_{-\a}$ ($\a\in R^+$) will not
depend on $\l$. So we have to take into account only those summands
of $-J$ that have $h_\a$ or $e_\a$ ($\a\in R^+$) as the third tensor
factor. Further, to calculate the multiplication constant of our
product it is sufficient to act by the operator on any vector. Let
us choose the highest weight vector $v_0$ for this. The Cartan
operators $h_\a$ multiply $v_0$ by $\l(h_\a) =
2\langle\l,\a\rangle/\langle\a,\a\rangle$. So indeed they are linear
in $\l$. But the raising operators $e_\a$ ($\a\in R^+$) send $v_0$
to zero. This means that we have to take into account only those
summands of $-J$ whose third tensor factor is one of the
$h_\a\mbox{'s}$. This is exactly the second sum of $J$ with the
opposite sign:
$$
\sum_{\a\in R} \Bigl( \langle\a,\a\rangle/2 \Bigr)^2 e_\a\ot e_{-\a}\ot h_\a.
$$
Now performing all the contractions corresponding to the edges
connecting the internal vertices of $\ol{w_{2m}}$ we get the tensor
$$
\sum_{\a\in R} \Bigl( \kf\a\a/2 \Bigr)^{2m}
  \underbrace{h_\a\ot \dots \ot h_\a}_{2m \mbox{ {\scriptsize times}}}.
$$
The corresponding element of $U(\g)$ acts on the highest weight
vector $v_0$ as multiplication by
$$S_{BNG}(\ol{w_{2m}}) = \sum_{\a\in R} \kf\l\a^{2m} =
  2\sum_{\a\in R^+} \kf\l\a^{2m}\,.
$$
The theorem is proved.
\end{proof}

\section{The Goussarov--Habiro theory revisited} 
\label{GH}

The term {\em Goussarov-Habiro theory} refers to the study of
$n$-equivalence classes of knots (or, more generally, knotted
graphs), as defined in Section~\ref{sing_knot_filtr}, in terms of
local moves on knot diagrams. It was first developed by M.~Goussarov
who announced the main results in September 1995 at a conference in
Oberwolfach, and, independently, by K.~Habiro \cite{Ha1, Ha2}. (As
often happened with Goussarov's results, his publication on the
subject \cite{G4} appeared several years later.)

There are several different approaches to Goussarov-Habiro theory,
which produce roughly the same results. In Chapter~\ref{chapBr} we
have developed the group-theoretic approach pioneered by T.~Stanford
\cite{Sta4, Sta3} who described $n$-equivalence in terms of the
lower central series of the pure braid groups. Habiro in \cite{Ha1,
Ha2} uses {\em claspers} to define local moves on knots and string
links. Here we shall briefly sketch Goussarov's approach, neither
giving complete proofs, nor striving for maximal generality. Other
versions of theorems of the same type can be found in \cite{CT, TY}.
A proof that the definitions of  Goussarov, Habiro and Stanford are
equivalent can be found in \cite{Ha2}.

\subsection{Statement of the Goussarov--Habiro Theorem}
In what follows we shall use the term {\em tangle} in the sense that
is somewhat different from the rest of this book. Here, by a tangle
\index{Tangle} we shall mean an oriented 1-dimensional submanifold
of a ball in $\R^3$. The isotopy of tangles is understood to fix the
boundary.

\begin{xtheorem}[Goussarov--Habiro]
\label{GHthm} \index{Theorem!Goussarov--Habiro} Let $K_1$ and $K_2$
be two knots. They are $n$-equivalent, that is, $v(K_1)=v(K_2)$ for
any $\Z$-valued Vassiliev invariant $v$ of order $\leq n$ if and
only if $K_1$ and $K_2$ are related by a finite sequence of moves
$\M_n$:
$$\index{Moves!Goussarov--Habiro!$\M_n$}
\index{Goussarov!--Habiro move!$\M_n$}
\underbrace{\ig[width=160pt]{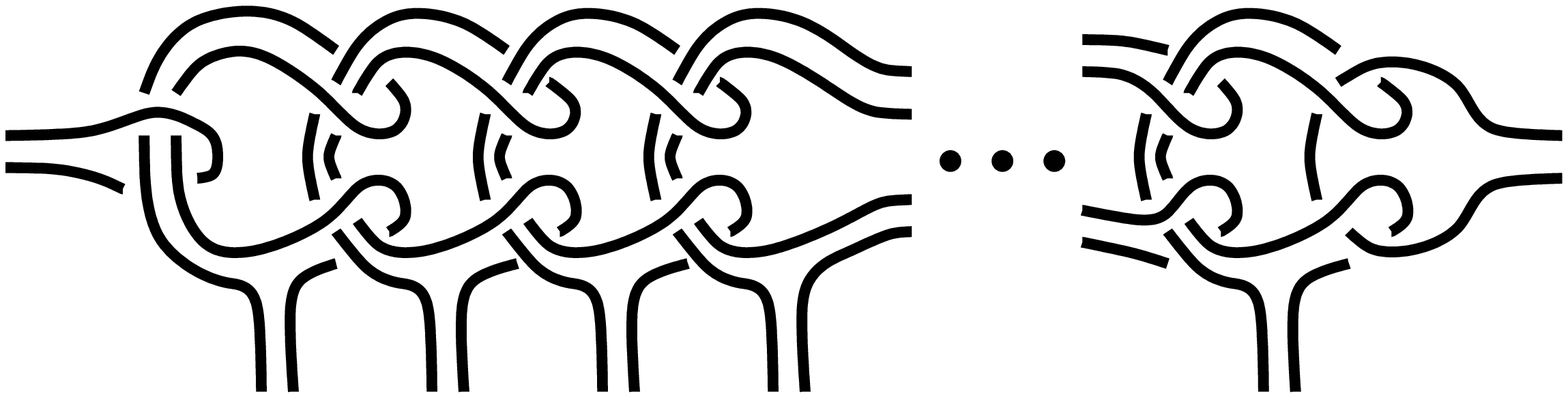}}_{\mbox{\scriptsize $n+2$ components}}
\quad\rb{15pt}{\ig[width=25pt]{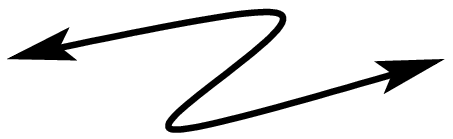}}\quad
\underbrace{\ig[width=120pt]{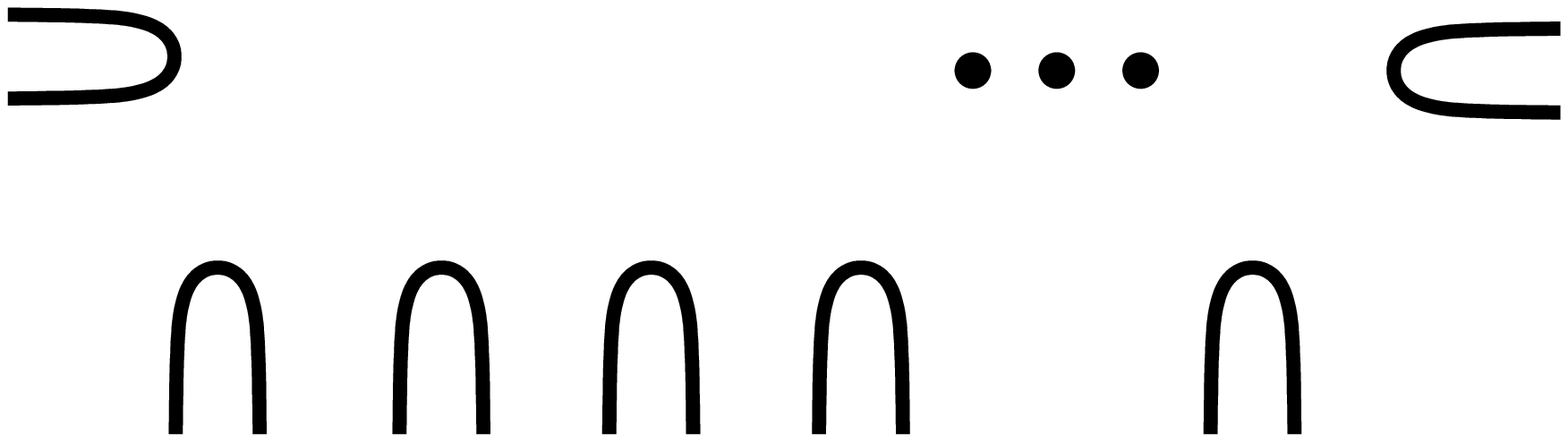}}_{\mbox{\scriptsize $n+2$ components}}
$$
\end{xtheorem}

\def\totonew{\rb{-2pt}{\quad\ig[width=20pt]{totonew.eps}\quad}}
\def\rs2#1#2#3{\rb{#3pt}{\ \ig[width=#1pt]{#2.eps}\ }}

Denote by $\B_n$ and $\T_n$ the tangles on the left and,
respectively, on the right-hand side of the move $\M_n$. The tangle
$\B_n$ is an example of a {\em Brunnian tangle\/} characterized by
the property that removing any of its components makes the remaining
tangle to be isotopic to the trivial tangle $\T_{n-1}$ with $n+1$
components.

The sequence of moves $\M_n$ starts with $n=0$:
$$\M_0:\quad \B_0=\rs2{40}{b0}{-5}\totonew \rs2{30}{t0}{-2}=\T_0.$$
In terms of knot diagrams $\M_0$ consists of a crossing change
followed by a second Reidemeister move.

The move $\M_1$ looks like
$$\M_1:\quad \B_1=\rs2{60}{t3d2}{-22}\totonew \rs2{50}{t3t2}{-10}=\T_1$$
It is also known as the {\em Borromean move}%\quad
\index{Moves!Borromean}\index{Borromean move}
$$\rb{-25pt}{$\rs2{50}{borp}{-20}\rb{5pt}{\totonew}
\rs2{40}{bort}{-20}$} %\vspace{10pt}
$$ Since there are no invariants
of order $\leq 1$ except constants (Proposition \ref{V1}), the
Goussarov--Habiro theorem implies that any knot can be transformed
to the unknot by a finite sequence of Borromean moves $\M_1$; in
other words, $\M_1$ is an\index{Unknotting operation} {\em
unknotting operation}.

\begin{xremark}
The coincidence of all Vassiliev invariants of order $\leq n$
implies the coincidence of all Vassiliev invariants of order $\leq
n-1$. This means that one can accomplish a move $\M_n$ by a sequence
of moves $\M_{n-1}$. Indeed, let us draw the tangle $\B_n$ as shown
below on the left:
$$\ig[width=150pt]{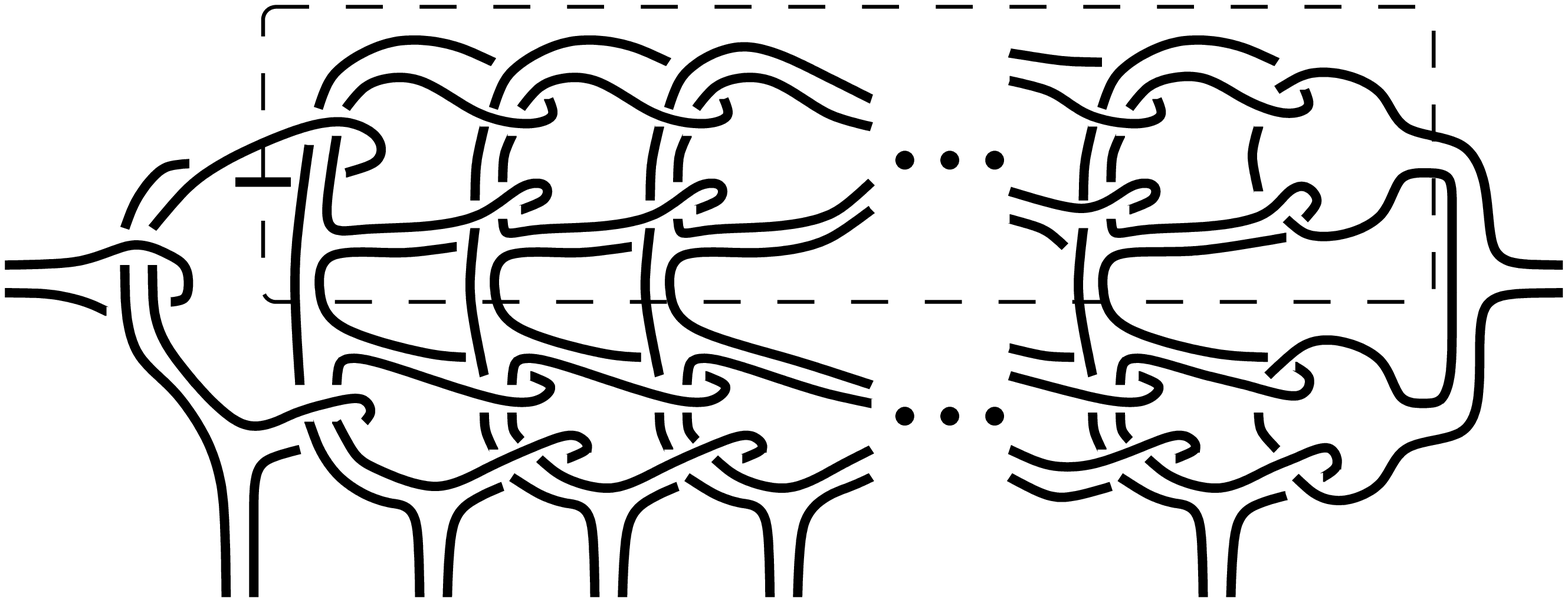}
  \rb{35pt}{$\begin{array}{c}\M_{n-1}\\ \totonew\end{array}$}
  \ig[width=150pt]{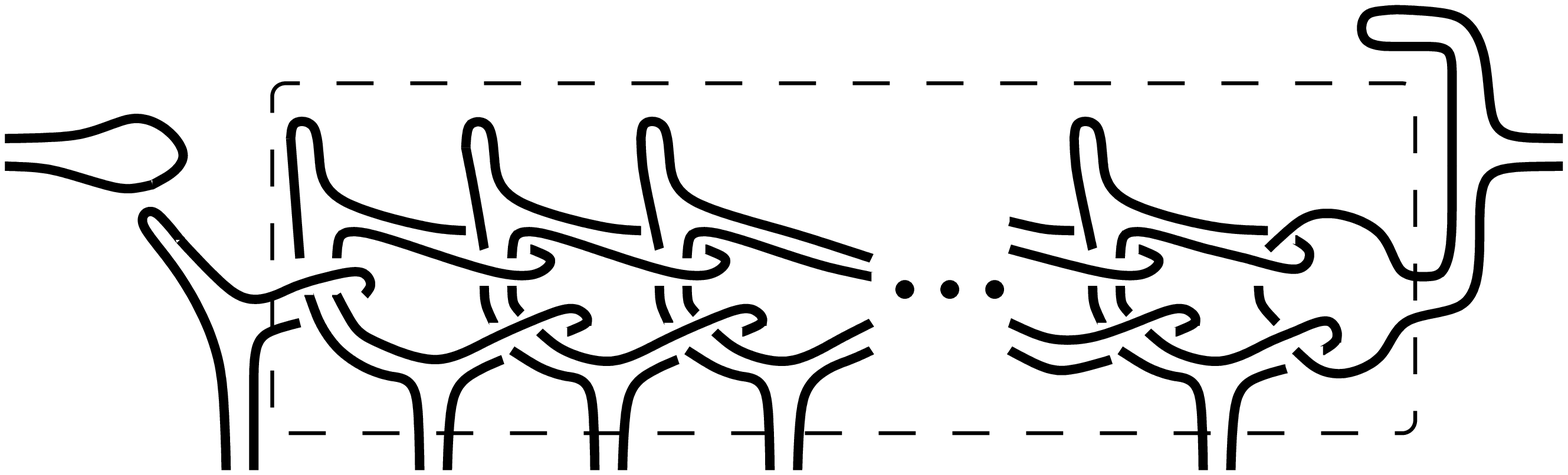}
$$
(In order to see that the tangle on the left is indeed $\B_n$,
untangle the components one by one, working from right to left). The
tangle in the dashed rectangle is $\B_{n-1}$. To perform the move
$\M_{n-1}$ we must replace it with $\T_{n-1}$. This gives us the
tangle on the right also containing $\B_{n-1}$. Now performing once
more the move $\M_{n-1}$ we obtain the trivial tangle $\T_n$.
\end{xremark}

\subsection{Reformulation of the Goussarov--Habiro Theorem}

Recall some notation from \ref{alg_knots} and \ref{sing_knot_filtr}.
We denote by $\K$ the set of all isotopy classes of knots, $\Z\K$ is
the free $\Z$-module (even an algebra) consisting of all finite
formal $\Z$-linear combinations of knots and $\K_n$ stands for the
$n$th term of the singular knot filtration in $\Z\K$. Using the
moves $\M_n$, we can define another filtration in the module $\Z\K$.

Let $\Ha_n$ be the $\Z$-submodule of $\Z\K$ spanned by the
differences of two knots obtained one from another by a single move
$\M_n$ . For example, the difference $3_1-6_3$ belongs to $\Ha_2$.
The Goussarov--Habiro Theorem can be restated as follows.

\begin{xtheorem} \label{GHthm_var}
For all $n$ the submodules $\K_{n+1}$ and $\Ha_n$ coincide.
\end{xtheorem}

\begin{proof}[Proof of the equivalence of the two statements.]
As we have seen in Section~\ref{sing_knot_filtr}, the values of any
Vassiliev invariant of order $\leq n$ are the same on the knots $K$
and $K'$ if and only if the difference $K-K'$ belongs to $\K_{n+1}$.
On the other hand, $K$ can be obtained from $K'$ by a sequence of
$\M_n$-moves if and only if $K-K'$ belongs to $\Ha_n$. Indeed, if
$K-K'\in\Ha_n$, we can write
$$K-K'=\sum_{i}a_i(K_i-K'_{i})$$ 
where $a_i$ are positive integers and  each $K_i$ differs from
$K'_{i}$ by a single $\M_n$ move. Since all the knots in this sum
apart from $K$ and $K'$ cancel each other, we can write it as
$\sum_{i=1}^{n}K_i-K_{i-1}$ with $K_n=K$, $K_0=K'$ and where each
$K_i$ differs from $K_{i-1}$ by a single $\M_n$ move. Then $K$ is
obtained by from $K'$ by a sequence of $\M_n$ moves that change
consecutively $K_{i-1}$ into $K_i$.
\end{proof}

\noindent{\bf Exercise.} \label{GHisS} Prove that two knots are
related by a sequence of $\M_n$-moves if and only if they are
$\gamma_{n+1}$-equivalent (see page~\pageref{theorem:group}).
\medskip

The Goussarov--Habiro Theorem is a corollary of this exercise and
Theorem~\ref{theorem:nilpotency} on
page~\pageref{theorem:nilpotency}. Nevertheless, we shall verify one
part of the Goussarov--Habiro Theorem directly, in order to give the
reader some feel of the Goussarov--Habiro theory. Namely, let us
show that $\Ha_n\subseteq\K_{n+1}$. (The inclusion
$\K_{n+1}\subseteq \Ha_n$ is rather more difficult to prove.)

\subsection{Proof that $\Ha_n$ is contained in $\K_{n+1}$}
In order to prove that $\Ha_n\subseteq \K_{n+1}$ it is sufficient to
represent the difference $\B_n-\T_n$ as a linear combination of
singular tangles with $n+1$ double points each. Let us choose the
orientations of the components of our tangles as shown. Using the
Vassiliev skein relation we shall gradually transform the difference
$\B_n-\T_n$ into the required form.
\vspace{10pt}\\
\def\lr#1#2{\rb{-15pt}{\ig[width=#2pt]{#1.eps}}}
$\begin{array}{lcl}
\B_n-\T_n &=& \lr{bns1}{80}\quad -\quad \lr{tns1}{60}\vspace{20pt}\\
\multicolumn{3}{l}{\qquad
=\quad \lr{bns2}{80}\quad + \quad\lr{bns3}{80}\quad -\quad
\lr{tns1}{60}}
\end{array}$\vspace{15pt}

\noindent
But the difference of the last two
tangles can be expressed as a singular tangle:
$$=\quad \lr{bns2}{80}\quad - \quad\lr{bns4}{80}$$

We got a presentation of $\B_n-\T_n$ as a linear combination of two
tangles with one double point on the first two components. Now we
add and subtract isotopic singular tangles with one double point:
$$\begin{array}{lcl}
\B_n-\T_n &=&
\hspace{-5pt}\left(\lr{bns2}{80}\quad -\quad\lr{bns5}{80}\right)\vspace{15pt}\\
&& -\! \left(\lr{bns4}{80}\quad-\quad\lr{bns6}{80}\right)
\end{array}$$
Then using the Vassiliev skein relation we can see that the
difference in the first pair of parentheses is equal to
$$\begin{array}{cl}
&-\quad\lr{bns7}{80}\quad+\quad\lr{bns8}{80}\quad-\quad\lr{bns5}{80}
    \vspace{10pt}\\
=&-\quad\lr{bns7}{80}\quad+\quad\lr{bns9}{80}
\end{array}$$
Similarly the
difference in the second pair of parentheses would be equal to
$$\begin{array}{cl}
&-\quad\lr{bns10}{80}\quad+\quad\lr{bns11}{80}\quad-\quad\lr{bns6}{80}
    \vspace{10pt}\\
=&-\quad\lr{bns10}{80}\quad+\quad\lr{bns12}{80}
\end{array}$$
Now we have represented $\B_n-\T_n$ as a linear combination of four
singular tangles with two double points each; in each tangle one
double point lies on the first and on the second components and the
other double point --- on the second and on the third components:
$$\begin{array}{lcl}
\B_n-\T_n &=&\quad -\quad\lr{bns7}{80}\quad+\quad\lr{bns9}{80}
     \vspace{10pt}\\
&&\quad+\quad\lr{bns10}{80}\quad-\quad\lr{bns12}{80}
\end{array}$$

Continuing in the same way we arrive to a linear combination of
$2^n$ tangles with $n+1$ double points each; one double point for
every pair of consecutive components. It is easy to see that if we
change the orientations of arbitrary $k$ components of our tangles
$\B_n$ and $\T_n$, then the whole linear combination will be
multiplied by $(-1)^k$.

\begin{xexample}
$$\B_2-\T_2 = \lr{bts1}{70}-\lr{bts2}{70}-
 \lr{bts3}{70}+\lr{bts4}{70}\vspace{15pt}$$
\end{xexample}

\begin{example}
There is only one (up to multiplication by a scalar and adding a
constant) nontrivial Vassiliev invariant of order $\leq 2$, namely
It is the coefficient $c_2$ of the Conway polynomial.

Consider two knots\vspace{-15pt}
$$3_1\quad=\quad
\rb{-20pt}{\ig[width=40pt]{31.eps}}\quad ,\hspace{2cm}
  6_3\quad=\quad
\rb{-20pt}{\ig[width=40pt]{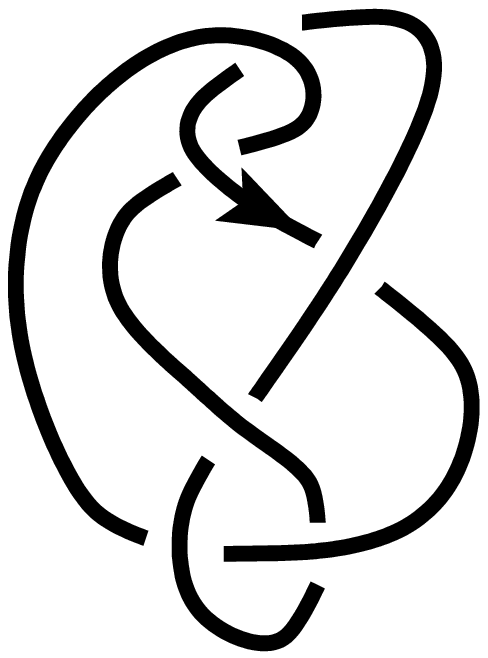}}\ .
$$
We choose the orientations as indicated. Their Conway polynomials
$$
  \CP(3_1)=1+t^2\ ,\qquad\qquad \CP(6_3)=1+t^2+t^4.
$$
have equal coefficients at $t^2$. Therefore for any Vassiliev
invariant $v$ of order $\leq 2$ we have $v(3_1)=v(6_3)$. In this
case the Goussarov--Habiro Theorem states that it is possible to
obtain the knot $6_3$ from the knot $3_1$ by moves $\M_2: \B_2
\totonew \T_2$
$$\M_2:\qquad\rb{-8pt}{\lr{bp2}{80}}\quad\totonew\quad
\rb{-8pt}{\lr{tp2}{60}}$$

Let us show this. We start with the standard diagram of $3_1$, and
then transform it in order to obtain $\B_2$ as a subtangle.
\vspace{10pt}
$$\begin{array}{ccl}
3_1&=&\rs2{32}{31}{-13}\cong \rs2{34}{31-1}{-18}
\cong\rs2{34}{31-2}{-20}\cong\rs2{43}{31-3}{-20}
\cong\rs2{40}{31-4}{-18}
\vspace{15pt}\\
&\cong& \rs2{45}{31-5}{-15}\cong \rs2{55}{31-7}{-15}
\cong\rs2{60}{31-8}{-15}\cong\rs2{60}{31-9}{-15}
\vspace{15pt}\\
&\cong& \rs2{62}{31-10}{-15}
\cong \rs2{70}{31-11}{-20}\cong\rs2{75}{31-12}{-20}
\vspace{15pt}\\
&\cong&\rs2{77}{31-13}{-25}\cong\rs2{88}{31-14}{-25}
\end{array}$$

\medskip
Now we have the tangle $\B_2$ in the dashed oval. Perform the move
$\M_2$ replacing $\B_2$ with the trivial tangle $\T_2$:\vspace{5pt}
$$\rs2{85}{63-1}{-25}\cong\!\rs2{45}{63-2}{-20}\!\!\cong
\!\rs2{50}{63-3}{-25}\!\cong\!\!\rs2{40}{63-4}{-20}\!\cong
\rs2{35}{63-5}{-20}\! = 6_3$$
\end{example}

\subsection{Vassiliev invariants and local moves}

The $mod\ 2$ reduction of $c_2$ is called the {\em Arf invariant}
\index{Arf invariant} of a knot. A description of the Arf invariant
similar to the Goussarov--Habiro description of $c_2$ was obtained
by L.~Kauffman.

\begin{xtheorem} [L.Kauffman \cite{Ka1, Ka2}] $K_1$ and $K_2$
have the same Arf invariant if and only if $K_1$ can be obtained
from $K_2$ by a finite number of so called pass moves:
$$\index{Move!pass}\index{Pass move}
  \lr{pass1}{40}\quad \totonew\quad \lr{pass2}{40}
$$
\end{xtheorem}

The orientations are important. Allowing pass moves with arbitrary
orientations we obtain an unknotting operation (see \cite{Kaw2}).

Actually, one can develop the whole theory of Vassiliev invariants
using the pass move instead of the crossing change in the Vassiliev
skein relation. It turns out, however, that all primitive finite
type invariants with respect to the pass move of order $n$ coincide
with primitive Vassiliev invariants of order $n$ for all $n\geq 1$.
The Arf invariant is the unique finite type invariant of order 0
with respect to the pass move \cite{CMS}.

More generally, in the definition of the finite type invariants one
can replace a crossing change with an arbitrary local move, that is,
a modification of a knot that replaces a subtangle of some fixed
type with another subtangle. For a wide class of moves one obtains
theories of finite type invariants for which the Goussarov--Habiro
Theorem holds, see \cite{TY, CMS}.

One such move is the {\em doubled-delta move}:
$$ \index{Moves!doubled-delta}\index{Doubled-delta move}
 \lr{ddelta1}{80}\quad \rb{20pt}{\totonew}\quad \lr{ddelta2}{80}
$$
S.~Naik and T.~Stanford \cite{NaS} have shown that two knots can be
transformed into each other by doubled-delta moves if and only if
they are {\em S-equivalent}, that is, if they have a common Seifert
matrix, see \cite{Kaw2}. The theory of finite type invariants based
on the doubled-delta move appears to be rather rich. In particular,
for each $n$ there is an infinite number of independent invariants
of order $2n$ which are not of order $2n-1$. We refer the reader to
\cite{CMS} for more details.

\subsection{The Goussarov groups of knots}
There are two main results in the Goussarov-Habiro theory. One is
what we called the Goussarov-Habiro Theorem. The other result says
that classes of knots (more generally, string links) related by
$\M_n$-moves form groups under the connected sum operation.

Modulo the exercise on page~\pageref{GHisS}, we have proved this in
Chapter~\ref{chapBr}, see Theorem~\ref{theorem:Vassiliev}. There we
were mostly interested in applying the technique of braid closures
and the theory of nilpotent groups. Here let us give some concrete
examples.

We shall denote by ${\mathcal G}_n$ the $n$th {\em Goussarov
group,} that is, the set $\K/\Gamma_{n+1}\K$ of $n$-equivalence
classes of knots with the connected sum operation.
\label{gg}\index{Goussarov!group} A {\em $j$-inverse} for a knot $K$
is a knot $K'$ such that $K\#K'$ is $j$-trivial. An $n$-inverse for
$K$ provides an inverse for the class of $K$ in ${\mathcal G}_n$.

Since there are no Vassiliev invariants of order $\leq 1$ except
constants, the zeroth and the first Goussarov groups are trivial.

\subsection{The second Goussarov group ${\mathcal G}_2$.} Consider
the coefficient $c_2$ of the Conway polynomial $\CP(K)$. According
to Exercise~\ref{sum_conway} at the end of Chapter~\ref{kn_inv},
$\CP(K)$ is a multiplicative invariant of the form
$\CP(K)=1+c_2(K)t^2+\dots$. This implies that
$c_2(K_1\#K_2)=c_2(K_1)+c_2(K_2)$, and, hence, $c_2$ is a
homomorphism of ${\mathcal G}_2$ into $\Z$. Since $c_2$ is the only
nontrivial invariant of order $\leq 2$ and there are knots on which
it takes value 1, the homomorphism $c_2: {\mathcal G}_2\to \Z$ is,
in fact, an isomorphism and ${\mathcal G}_2\cong\Z$. From the table
in Section~\ref{conway_tabl} we can see that $c_2(3_1)=1$ and
$c_2(4_1)=-1$. This means that the knot $3_1$ represents a generator
of ${\mathcal G}_2$, and $4_1$ is 2-inverse of $3_1$. The prime
knots with up to 8 crossings are distributed in the second Goussarov
group ${\mathcal G}_2$ as follows:
$$\ig[width=\textwidth]{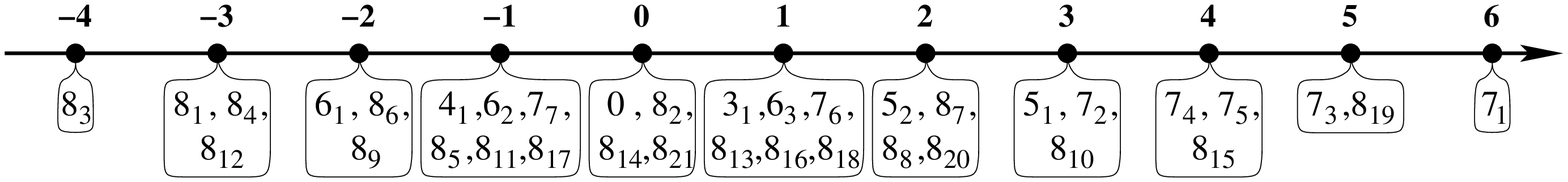} \hspace{-10pt}\rb{35pt}{$c_2$}$$

\subsection{The third Goussarov group ${\mathcal G}_3$.}
\label{third_gg}

In order 3 we have one more Vassiliev invariant; namely, $j_3$, the
coefficient at $h^3$ in the power series expansion of the Jones
polynomial with the substitution $t=e^h$. The Jones polynomial is
multiplicative, $J(K_1\#K_2)=J(K_1)\cdot J(K_2)$ (see
Exercise~\ref{sum_jones} at the end of Chapter~\ref{kn_inv}) and its
expansion has the form $J(K)=1+j_2(K)h^2+j_3(K)h^3+\dots$ (see
Section~\ref{knpol_as_vi}). Thus we can write
$$J(K_1\#K_2)=1+\left(j_2(K_1)+j_2(K_2)\right)h^2+
                \left(j_3(K_1)+j_3(K_2)\right)h^3+\dots$$
In particular, $j_3(K_1\#K_2)=j_3(K_1)+j_3(K_2)$. According to
Exercise~\ref{div_by_6} at the end of Chapter~\ref{FT_inv}, $j_3$ is
divisible by 6. Then $j_3/6$ is a homomorphism from ${\mathcal G}_3$
to $\Z$. The direct sum with $c_2$ gives the isomorphism
$${\mathcal G}_3\cong\Z\oplus\Z=\Z^2;\qquad\qquad
K\mapsto\left(c_2(K),j_3(K)/6\right)$$ Let us identify ${\mathcal
G}_3$ with the integral lattice on a plane. The distribution of
prime knots on this lattice is shown in Figure \ref{distr_2_3};
recall that $\ol{K}$ is the mirror reflection of $K$.

In particular, the 3-inverse of the trefoil $3_1$ can be represented
by $6_2$, or by $\overline{7_7}$. Also, we can see that $3_1\#4_1$
is 3-equivalent to $\overline{8_2}$. Therefore $3_1\#4_1\#8_2$ is
3-equivalent to the unknot, and $4_1\#8_2$ also represents the
3-inverse to $3_1$. The knots $6_3$ and $\overline{8_2}$ represent
the standard generators of ${\mathcal G}_3$.

\subsection*{Open problem.} Is there any torsion in the group
${\mathcal G}_n$?

\section{Willerton's fish and bounds for $c_2$ and $j_3$}
\label{fish}

\begin{figure}[ht]
\rb{-10pt}[480pt][10pt]{\begin{picture}(210,15)(0,0)
  \put(-50,0){\ig[width=300pt,height=470pt]{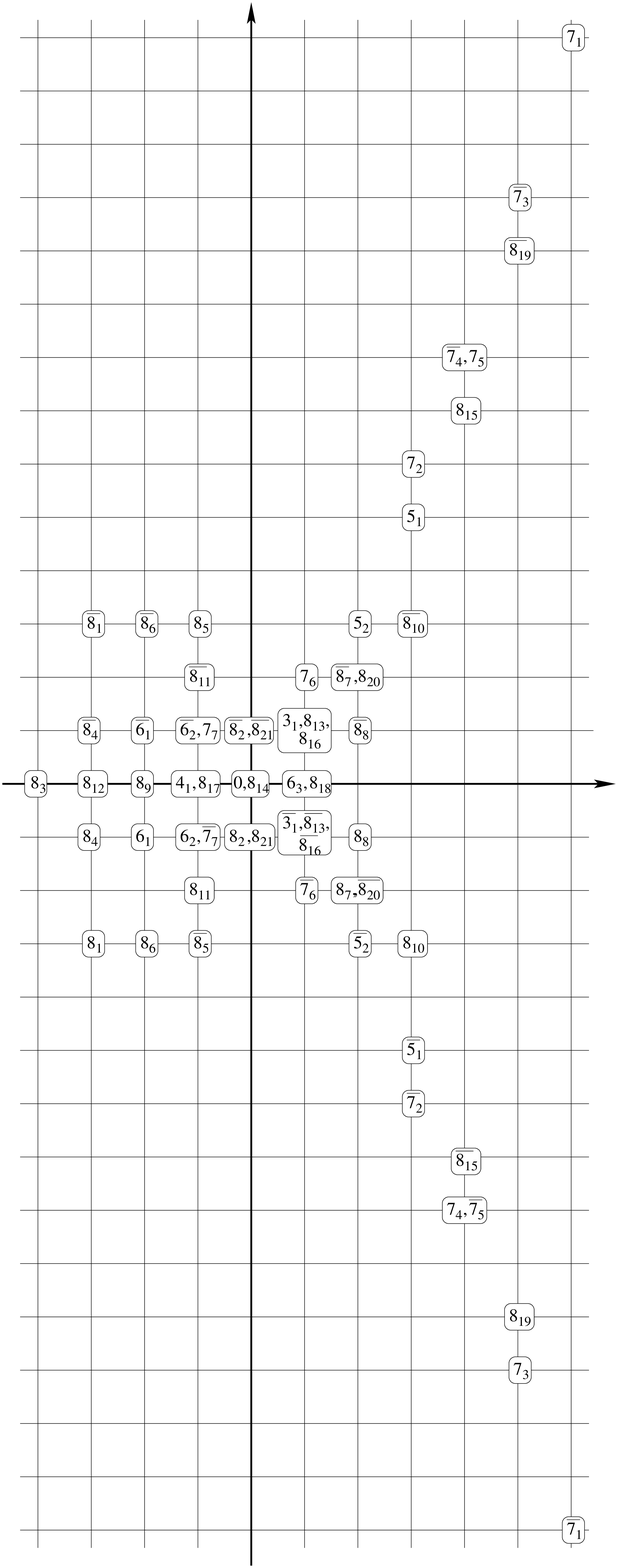}}
  \put(260,245){$c_2$} \put(76,470){$j_3/6$}
     \end{picture}}
\caption{Values of Vassiliev invariants $c_2$ and $j_3$ on prime
knots with up to 8 crossings. The mirror images of $8_{14}$ and of
$8_{17}$, not shown, have the same invariants as the original
knots.} \label{distr_2_3}
\end{figure}

{\em Willerton's fish}\index{Fish ! Willerton's} is a graph where
the Vassiliev invariant $c_2$ is plotted against the invariant $j_3$
for all knots of a given crossing number \cite{Wil5}. The shape of
this graph, at least for the small values of the crossing number
($\leq 14$) where there is enough data to construct it, is
reminiscent of a fish, hence the name. (This shape is already
discernible on Figure~\ref{distr_2_3} which shows all prime knots up
to 8 crossings. The composite knots add no new points to the graph.)

A plausible explanation for the strange shape of these graphs could
involve some inequality on $c_2$, $j_3$ and the crossing number $c$.
At the moment, no such inequality is known. However, there are
several results relating the above knot invariants.
\begin{theorem}[\cite{PV2}] For any knot $K$
$|\, c_2(K)\,|\leq\left[\frac{c(K)^2}{8}\right].$
\end{theorem}
\def\gt#1#2#3{\risS{-5}{#3}{\put(0,10){$\scriptstyle #1$}
     \put(35,10){$\scriptstyle #2$}}{45}{13}{6}}
\def\gdCi#1#2{\bigl\langle \,\risS{-5}{gauss-6A}{\put(0,10){$\scriptstyle #1$}
     \put(35,10){$\scriptstyle #2$}}{45}{0}{0}\,, D\bigr\rangle }
\def\gD#1#2{\gt{#1}{#2}{gauss-6D}}
\begin{proof} Recall the Gauss diagram formulae for $c_2$ on
pages~\pageref{fig:PV} and \pageref{otherc2} which we can write as follows:
$$c_2(K)=\bigl\langle{\gt{}{}{gauss-6A}}\,, D\bigr\rangle = \bigl\langle{\gt{}{}{gauss-6D}}\,, D\bigr\rangle,$$
where $D$ is a based Gauss diagram with $n$ arrows representing the
knot $K$. Let $C^+$ be the set of arrows of $D$ that point forward
(this makes sense since $D$ is based) and let $C^-$ be the set of
backwards-pointing arrows. If $C^+$ consists of $k$ elements, then
$C^-$ has $n-k$ elements.

Now, assume that the diagram $${\gt{}{}{gauss-6A}}$$ appears $n_1$
times as a subdiagram of $D$ and the diagram $${\gt{}{}{gauss-6D}}$$
appears $n_2$ times. Each of these diagrams contains one arrow from
$C^+$ and one from $C^-$. Therefore, we have
$$|\, c_2(K)\,|\leq\min{(n_1,n_2)}\leq\frac{k(n-k)}{2}\leq
\left[\frac{n^2}{8}\right].$$ Now, the smallest possible $n$ in this
formula, that is, the minimal number of arrows in a Gauss diagram
representing $K$, is, by definition, nothing else but the crossing
number $c(K)$.
\end{proof}

\subsection{Invariants of higher degrees}
Similar inequalities exist for all Vassiliev invariants. Indeed,
each invariant of order $n$ can be represented by a Gauss diagram
formula (see Chapter~\ref{chapGD}). This means that its value on a
knot can be calculated by representing this knot by a Gauss diagram
$D$ and counting subdiagrams of $D$ certain types, all with at most
$n$ arrows. The number of such subdiagrams grows as $(\deg{D})^n$,
so for each invariant of degree $n$ there is a bound by a polynomial
of degree $n$ in the crossing number. In particular, S.~Willerton
found the following bound (unpublished):
$$|\,j_3(K)\,|\leq\frac{3}{2}\cdot c(K)(c(K)-1)(c(K)-2).$$

\subsection{Inequalities for torus knots} One particular family of
knots for which $c_2$ and $j_3$ are related by explicit inequalities
are the torus knots \cite{Wil5}. We have
$$24 c_2(K)^3+ 12 c_2(K)^2\leq j_3(K)^2\leq 32 c_2(K)^3+ 4 c_2(K)^2$$
for any torus knot $K$. These bounds are obtained from the explicit
expressions for $c_2$ and $j_3$ for torus knots obtained in
\cite{AL}.

\section{Bialgebra of graphs}
\label{algLando}

It turns out that the natural mapping that assigns to every chord
diagram its intersection graph, can be converted into a homomorphism
of bialgebras $\gamma:\A\to\La$, where $\A$ is the algebra of chord
diagrams and $\La$ is an algebra generated by graphs modulo certain
relations, introduced by S.~Lando \cite{Lnd2}. Here is his
construction.

Let $\Gr$ \label{Gr} be the graded vector space 
spanned by all simple graphs (without
loops or multiple edges) as free generators:
$$
\Gr=\Gr_0\oplus\Gr_1\oplus\Gr_2\oplus\dots,
$$
It is graded by the order (the number of vertices) of a graph.
This space is easily turned into a bialgebra:
\index{Bialgebra!of graphs}

(1) The product is defined as the disjoint union of graphs, then
extended by linearity. The empty graph plays the role of the unit in
this algebra.

(2) The coproduct is defined similarly to the coproduct in the
bialgebra of chord diagrams. If $G$ is a graph, let $V=V(G)$ be the
set of its vertices. For any subset  $U \subset V$ denote by $G(U)$
the graph with the set of vertices $U$ and those vertices of the
graph $G$ whose both endpoints belong to $G$. We set
\begin{equation}
\label{coprlg}
\delta(G)=\sum_{U\subseteq V(G)}G(U)\otimes G(V\setminus U),
\end{equation}
and extend $\delta$ by linearity to the whole of $\Gr$.

The sum in (\ref{coprlg}) is taken over all subsets $U\subset V$
and contains as many as $2^{\#(V)}$ summands.
\medskip

\begin{xexample}
\begin{center}
\begin{picture}(300,60)(10,60)
\put(12,97){$\delta($}
\put(27,100){\circle*{4}}
\put(27,100){\line(1,0){25}}
\put(52,100){\circle*{4}}
\put(52,100){\line(1,0){25}}
\put(77,100){\circle*{4}}
\put(81,97){$)$}
\put(92,97){$=$}
\put(106,97){$1$}
\put(116,97){$\otimes$}
\put(133,100){\circle*{4}}
\put(134,100){\line(1,0){25}}
\put(158,100){\circle*{4}}
\put(158,100){\line(1,0){25}}
\put(183,100){\circle*{4}}
\put(189,97){$+\ 2$}
\put(214,100){\circle*{4}}
\put(221,97){$\otimes$}
\put(237,100){\circle*{4}}
\put(237,100){\line(1,0){25}}
\put(262,100){\circle*{4}}
\put(269,97){$+$}
\put(284,100){\circle*{4}}
\put(291,97){$\otimes$}
\put(307,100){\circle*{4}}
\put(317,100){\circle*{4}}
\put(92,67){$+$}
\put(110,70){\circle*{4}}
\put(120,70){\circle*{4}}
\put(127,67){$\otimes$}
\put(142,70){\circle*{4}}
\put(150,67){$+\ 2$}
\put(175,70){\circle*{4}}
\put(175,70){\line(1,0){25}}
\put(200,70){\circle*{4}}
\put(207,67){$\otimes$}
\put(223,70){\circle*{4}}
\put(230,67){$+$}
\put(245,70){\circle*{4}}
\put(245,70){\line(1,0){25}}
\put(270,70){\circle*{4}}
\put(270,70){\line(1,0){25}}
\put(295,70){\circle*{4}}
\put(302,67){$\otimes\ 1$}
\end{picture}
\end{center}
\end{xexample}

\noindent{\bf Exercise.} Check the axioms of a Hopf algebra for
$\Gr$.
\medskip

The mapping from chord diagrams to intersection graphs does not
extend to a linear operator $\A\to\Gr$ since the combinations of
graphs that correspond to 4-term relations for chord diagrams do not
vanish in $\Gr$. To obtain a linear map, it is necessary to mod out
the space $\Gr$ by the images of the 4 term relations. Here is the
appropriate definition.
\medskip

Let $G$ be an arbitrary graph and $u,v$ an ordered pair of its vertices.
The pair $u,v$ defines two transformations of the graph $G$:
$G\mapsto G_{uv}'$ and $G\mapsto \wt{G}_{uv}$.
Both graphs $G_{uv}'$ and  $\wt{G}_{uv}$ have the same set of
vertices as $G$. They are obtained as follows.

If $uv$ is an edge in $G$, then the graph $G_{uv}'$ is obtained from
$G$ by deleting the edge $uv$; otherwise this edge should be added
(thus, $G\mapsto G_{uv}'$ toggles the adjacency of $u$ and $v$).

The graph $\wt{G}_{uv}$ is obtained from $G$ in a more tricky way.
Consider all vertices $w \in V(G) \setminus \{u,v\}$ which are
adjacent in $G$ with $v$. Then in the graph $\wt{G}_{uv}$ vertices
$u$ and $w$ are joined by an edge if and only if they are not joined
in $G$. For all other pairs of vertices their adjacency in $G$ and
in $\wt{G}_{uv}$ is the same. Note that the two operations applied
at the same pair of vertices, commute and, hence, the graph
$G_{uv}'$ is well-defined.

\begin{xdefinition}\index{Four-term relation!for graphs}
A four-term relation for graphs is
\begin{equation}
  \label{fourtlg}
  G - G_{uv}' = \wt{G}_{uv} - \wt{G}_{uv}'
\end{equation}
\end{xdefinition}
\begin{xexample}
 $$ \ig[width=8cm]{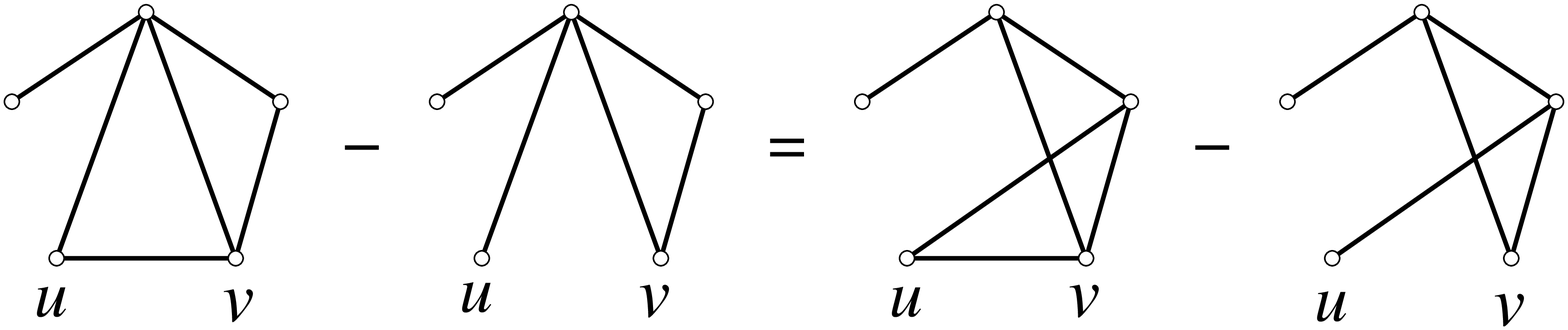}$$
\end{xexample}

\noindent{\bf Exercises.}

(1) Check that, passing to intersection graphs, the four-term
relation for chord diagram carries over exactly into this four-term
relation for graphs.

(2) Find the four-term relation of chord diagrams which is the
preimage of the relation shown in the example above.

\begin{xdefinition}\index{Lando!graph algebra}
The {\em graph bialgebra of Lando} $\La$ \label{La} is the quotient
of the graph algebra $\Gr$ by the ideal generated by all 4-term
relations (\ref{fourtlg}).
\end{xdefinition}

\begin{theorem}
The product and the coproduct defined above induce a bialgebra
structure in the quotient space $\La$.
\end{theorem}

\begin{proof}
The only thing that needs checking is that both the product and the
coproduct respect the 4-term relation (\ref{fourtlg}). For the
product, which is the disjoint union of graphs, this statement is
obvious. In order to verify it for the coproduct it is sufficient to
consider two cases. Namely, let $u,v\in V(G)$ be two distinct
vertices of a graph $G$. The right-hand side summands in the formula
(\ref{coprlg}) for the coproduct split into two groups: the summands
where both vertices $u$ and $v$ belong either to the subset
$U\subset V(G)$ or to its complement $V(G)\setminus U$, and those
where $u$ and $v$ belong to different subsets. By cleverly grouping
the terms of the first kind for the coproduct $\delta(G - G_{uv}' -
\wt{G}_{uv} + \wt{G}_{uv}')$ we can see that they all cancel out in
pairs. The terms of the second kind cancel out in pairs already
within each of the two summands $\delta(G - G_{uv}')$ and
$\delta(\wt{G}_{uv} - \wt{G}_{uv}')$.
\end{proof}

Relations~(\ref{fourtlg}) are homogeneous with respect to the number
of vertices, therefore $\La$ is a graded algebra. By
Theorem~\ref{MMthm} (page~\pageref{MMthm}), the algebra $\La$ is
polynomial with respect to its space of primitive elements.

Now we have a well-defined bialgebra homomorphism
$$
  \gamma:\A\to\La
$$
which extends the assignment of the intersection graph to a chord diagram.
It is defined by the linear mapping between the corresponding primitive
spaces $P(\A)\to P(\La)$.

According to S.~Lando \cite{Lnd2}, the dimensions of the homogeneous
components of $P(\La)$ are known up to degree 7. It turns out that
the homomorphism $\gamma$ is an isomorphism in degrees up to 6,
while the map $\gamma:P_7(\A)\to P_7(\La)$ has a 1-dimensional
kernel. See \cite{Lnd2} for further details and open problems
related to the algebra $\La$.

\section{Estimates for the number of Vassiliev knot invariants} %17
\label{bounds}

Knowing the dimensions of  $\PR_i$ for $i\le d$ is equivalent to
knowing $\dim\A_i$ or $\dim\V_i$ for $i\le d$. These numbers have
been calculated only for small values of $d$ and, at present, their
exact asymptotic behaviour as $d$ tends to infinity is not known.
Below we give a summary of all available results on these
dimensions.

\subsection{Historical remarks: exact results}
The precise dimensions of the spaces related to Vassiliev invariants
are known up to $n=12$, and are listed in the table below. They were
obtained by V.~Vassiliev for $n\le 4$ in 1990 \cite{Va2}, then by
D.~Bar-Natan for $n\le9$ in 1993 \cite{BN1} and by J.~Kneissler, for
$n=10,11,12$, in 1997 \cite{Kn0}. Vassiliev used manual
calculations. Bar-Natan wrote a computer program that implemented
the direct algorithm that solved the system of linear equations
coming from one-term and four-term relations. Kneissler obtained a
lower bound using the marked surfaces \cite{BN1} and an upper bound
using the action of Vogel's algebra $\Lambda$ on the primitive space
$\PR$: miraculously, these two bounds coincided for $n\le 12$.

\begin{center}
\index{Table of!dimensions of!the spaces of Vassiliev invariants}
\begin{tabular}{|*{14}{r|}}
\hline
 $n$       & 0 & 1 & 2 & 3 & 4 & 5  &  6 &  7 &  8 &   9 &  10 &  11 & 12 \\
\hline
$\dim\PR_n$ & 0 & 0 & 1 & 1 & 2 & 3  &  5 &  8 & 12 &  18 &  27 &  39 & 55 \\
\hline
$\dim\A_n$ & 1 & 0 & 1 & 1 & 3 & 4  &  9 & 14 & 27 &  44 &  80 & 132 & 232\\
\hline
$\dim\V_n$ & 1 & 1 & 2 & 3 & 6 & 10 & 19 & 33 & 60 & 104 & 184 & 316 & 548\\
\hline
\end{tabular}
\end{center}

The splitting of the numbers $\dim\PR_n$ for $n\le12$ according to
the second grading in $\PR$ is given in the table on
page~\pageref{dim-of-prim-sp}.
\smallskip

\noindent{\bf Exercise.} Prove that $\dim\A_n^{fr}=\dim\V_n$ for all
$n$.
\smallskip

\subsection{Historical remarks: upper bounds}

A priori it was obvious that $\dim\A_n <
(2n-1)!!=1\cdot3\cdots(2n-1)$, since this is the total number of
linear chord diagrams.

Then, there appeared five papers where this estimate was successively
improved:
\begin{enumerate}
\item
(1993) Chmutov and Duzhin \cite{CD1} proved that $\dim\A_n < (n-1)!$
\item
(1995) K.~Ng in \cite{Ng} replaced $(n-1)!$ by $(n-2)!/2$.
\item
(1996) A.~Stoimenow \cite{Sto1} proved that $\dim \A_n$
grows slower than $n!/a^n$, where $a=1.1$.
\item
(2000) B.~Bollob\'as and O.~Riordan \cite{BR1} obtained the asymptotical
bound $n!/(2\ln(2)+o(1))^n$ (approximately $n!/1.38^n$).
\item
(2001) D.~Zagier \cite{Zag1} improved the last result to
$\frac{6^n\sqrt{n}\cdot n!}{\pi^{2n}}$, which is asymptotically smaller
than $n!/a^n$ for any constant $a<\pi^2/6=1.644...$
\end{enumerate}

For the proofs of these results, we refer the interested reader to
the original papers, and only mention here the methods used to get
these estimates. Chmutov and Duzhin proved that the space $\A_n$ is
spanned by the {\em spine chord diagrams}, that is, diagrams
containing a chord that intersects all other chords, and estimating
the number of such diagrams. Stoimenow did the same with {\em
regular linearized diagrams}; Zagier gave a better estimate for the
number of such diagrams.

\subsection{Historical remarks: lower bounds}

In the story of lower bounds for the number of Vassiliev knot
invariants there is an amusing episode. The first paper by
Kontsevich about Vassiliev invariants (\cite{Kon1}, section 3)
contains the following passage:

``Using this construction\footnote{Of Lie algebra weight systems.},
one can obtain the estimate
$$\dim(\V_n) > e^{c\sqrt{n}}, \quad n \to +\infty$$
for any positive constant $c<\pi\sqrt{2/3}$ (see \cite{BN1a}, Exercise
6.14).''

Here $\V_n$ is a slip of the pen, instead of $\PR_n$, because of the
reference to Exercise 6.14 where primitive elements are considered.
Exercise 6.14 was present, however, only in the first edition of
Bar-Natan's preprint and eliminated in the following editions as
well as in the final published version of his text \cite{BN1}. In
\cite{BN1a} it reads as follows (page 43):

``{\sl Exercise} 6.14. (Kontsevich, [24]) Let $P_{\geq 2}(m)$ denote
the number of partitions of an integer $m$ into a sum of integers
bigger than or equal to 2. Show that $\dim \PR_m \geq P_{\geq
2}(m+1)$.

{\sl Hint} 6.15. Use a correspondence like
\begin{center}
\begin{picture}(300,20)(20,6)
  \bezier{60}(0,25)(80,25)(160,25)
  \bezier{60}(0,0)(80,0)(160,0)
  \bezier{10}(0,0)(0,12.5)(0,25)
  \bezier{10}(40,0)(40,12.5)(40,25)
  \bezier{10}(80,0)(80,12.5)(80,25)
  \bezier{10}(120,0)(120,12.5)(120,25)
  \bezier{10}(160,0)(160,12.5)(160,25)
  \bezier{5}(10,0)(10,5)(10,10)
  \bezier{5}(20,0)(20,5)(20,10)
  \bezier{5}(30,0)(30,5)(30,10)
  \bezier{5}(53.33,0)(53.33,5)(53.33,10)
  \bezier{5}(66.66,0)(66.66,5)(66.66,10)
  \bezier{5}(100,0)(100,5)(100,10)
  \bezier{5}(140,0)(140,5)(140,10)
  \put(19,14){\mbox{\scriptsize 4}}
  \put(59,14){\mbox{\scriptsize 3}}
  \put(99,14){\mbox{\scriptsize 2}}
  \put(139,14){\mbox{\scriptsize 2}}
  \put(190,15){\vector(1,0){20}}
  \put(210,15){\vector(-1,0){20}}
  \put(240,12){\mbox{$10+1 = 4+3+2+2,$}}
\end{picture}
\end{center}
and \dots''

The reference [24] was to ``M.~Kontsevich. Private communication.''! Thus,
both authors referred to each other, and none of them gave any proof. Later,
however, Kontsevich explained what he had in mind (see item 5 below).

Arranged by the date, the history of world records in asymptotic
lower bounds for the dimension of the primitive space $\PR_n$ looks
as follows.
\begin{enumerate}
\item (1994)
$\dim\PR_n \geq 1$ (``forest elements'' found by Chmutov, Duzhin and Lando
\cite{CDL3}).
\item (1995)
$\dim\PR_n \geq [n/2]$ (given by coloured Jones function --- see
Melvin--Morton \cite{MeMo} and Chmutov--Varchenko \cite{CV}).
\item (1996)
$\dim\PR_n \gtrsim n^2/96$ (see Duzhin \cite{Du1}).
\item (1997)
$\dim\PR_n \gtrsim n^{\log n}$, i.~e. the growth is faster than any
polynomial (Chmutov--Duzhin \cite{CD2}).
\item (1997)
$\dim\PR_n > e^{\pi\sqrt{n/3}}$ (Kontsevich \cite{Kon2}).
\item (1997)
$\dim\PR_n > e^{c\sqrt{n}}$ for any constant $c<\pi\sqrt{2/3}$
(Dasbach \cite{Da3}).
\end{enumerate}

Each lower bound for the dimensions of the primitive space $p_n=\dim\PR_n$
implies a certain lower bound for the dimensions of the whole algebra
$a_n=\dim\A_n$.
\begin{xproposition}
We have:
$a_n\gtrsim e^{n/\log_{b}n}$ for any constant $b<\pi^2/6$.
\end{xproposition}
\begin{proof}[Sketch of the proof] 
Fix a basis in each $\PR_k$, suppose that $n=km$ and consider the
elements of $\A_n$ which are products of $m$ basis elements of
$\PR_k$. Finding the maximum of this number over $k$ with fixed $n$,
we get the desired lower bound.
\end{proof}

\smallskip

Note that the best known upper and lower bounds on the dimensions of $\A_n$
are very far apart. Indeed, using the
relation between the generating functions
$$
\sum_{n=0}^\infty a_n t^n=
\prod_{k=1}^\infty(1-t^k)^{-p_k}=
\exp\sum_{n=1}^\infty\big(\sum_{k|n}p_k\big)t^n\,,
$$
one can easily prove (see \cite{Sto3}) that any subexponential lower bound on $p_n$ can only
lead to a subexponential lower bound on $a_n$, while the existing
upper bound is essentially factorial, that is, much greater than exponential.

\subsection{Proof of the lower bound}
\label{lower_bound}

We will sketch the proof of the lower bound for the number of
Vassiliev knot invariants, following \cite{CD2} and then explain how
O.~Dasbach \cite{Da3}, using the same method, managed to improve the
estimate and establish the bound which is still (2011) the best.

The idea of the proof is simple: we construct a large family of open
diagrams whose linear independence in the algebra $\B$ follows from
the linear independence of the values on these diagrams of a certain
polynomial invariant $P$, which is obtained by simplifying the
universal $\gl_N$ invariant.

As we know from Chapter \ref{LAWS}, the $\gl_N$ invariant
$\rho_{\gl_N}$, evaluated on an open diagram, is a polynomial in the
generalized Casimir elements $x_0$, $x_1$, ..., $x_N$. This
polynomial is homogeneous in the sense of the grading defined by
setting $\deg{x_m}=m$. However, in general, it is not homogeneous if
the $x_m$ are considered as variables of degree 1.

\begin{xdefinition}
The polynomial invariant $P:\B\to\Z[x_0,...,x_N]$ is the highest
degree part of $\rho_{\gl_N}$ if all the variables are taken with
degree 1.
\end{xdefinition}

For example, if we had $\rho_{\gl_N}(C)=x_0^2x_2-x_1^2$, then we
would have $P(C)=x_0^2x_2$.

Now we introduce the family of primitive open diagrams whose linear
independence we shall prove.

\begin{xdefinition}\index{Baguette diagram}
The {\em baguette diagram} $B_{n_1,\dots,n_k}$ is
\begin{center}                    
     \begin{picture}(300,90)(0,-20)
      \put(-40,35){\mbox{$B_{n_1,\dots,n_k}\quad =$}}
      \put(50,-10){\ig[width=280pt]{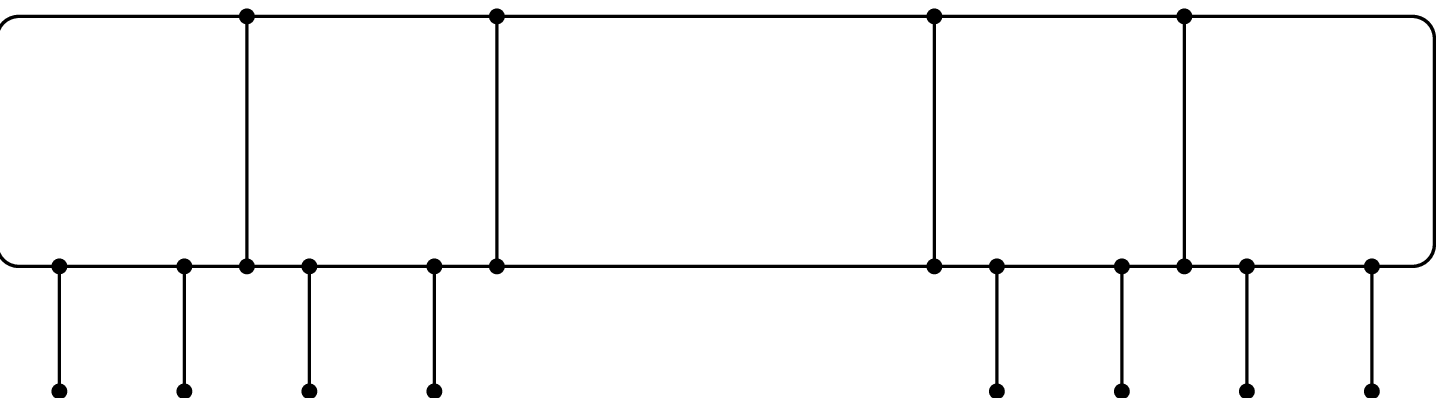}}
      \put(53,-10){\mbox{$\underbrace{\makebox[35pt]{}}_{
                           {\scriptstyle n_1{\rm \ vertices}}}$}}
      \put(103,-10){\mbox{$\underbrace{\makebox[35pt]{}}_{
                           {\scriptstyle n_2{\rm \ vertices}}}$}}
      \put(234,-10){\mbox{$\underbrace{\makebox[35pt]{}}_{\!\!\!\!\!
                           {\scriptstyle n_{k-1}{\rm \ vertices}}}$}}
      \put(282,-10){\mbox{$\underbrace{\makebox[34pt]{}}_{\ \
                           {\scriptstyle n_k{\rm \ vertices}}}$}}
      \put(67,2){\mbox{$\dots$}}
      \put(116,2){\mbox{$\dots$}}
      \put(250,2){\mbox{$\dots$}}
      \put(298,2){\mbox{$\dots$}}
      \put(185,38){\mbox{$\dots$}}
     \end{picture}
  \end{center}
It has a total of $2(n_1+\dots+n_k + k-1)$ vertices,
out of which $n_1+\dots+n_k$ are univalent.
\end{xdefinition}

To write down the formula for the value $P(B_{n_1,\dots,n_k})$, we
shall need the following definitions.

\begin{xdefinition}
Consider $k$ pairs of points arranged in two rows like
\rb{-1mm}{\ig[height=4mm]{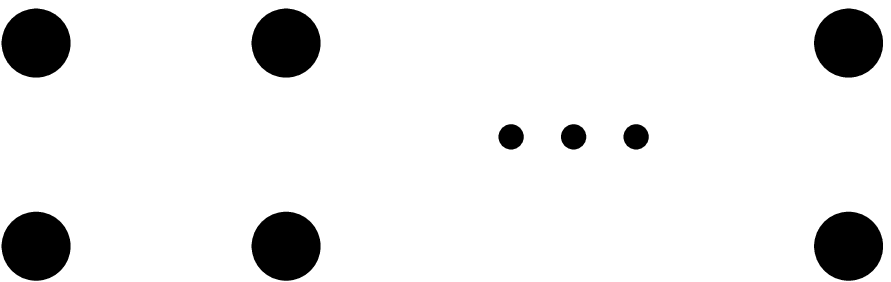}}. Choose one of the $2^{k-1}$
subsets of the set $\{1,\dots,k-1\}$. If a number $s$ belongs to the
chosen subset, then we connect the lower points of $s$th and
$(s+1)$th pairs, otherwise we connect the upper points. The
resulting combinatorial object is called a {\em two-line
scheme}\index{Scheme} of order $k$.
\end{xdefinition}

\begin{xexample} Here is the scheme corresponding to $k=5$ and the
subset $\{2,3\}$:
$$\mbox{\ \begin{picture}(60,15)(0,0)
          \put(0,0){\ig[width=60pt]{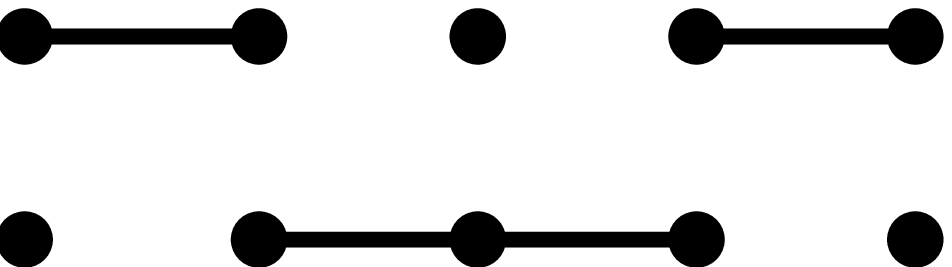}}
          \end{picture}\ }\quad.
$$
\end{xexample}

The number of connected components in a scheme of order $k$ is $k+1$.

\begin{xdefinition}\label{schememon}
\def\s{\sigma}
Let $\s$ be a scheme; $i_1,\dots,i_k$ be non negative integers:
$0\leq i_1\leq n_1, \dots, 0\leq i_k\leq n_k$. We assign $i_s$ to
the lower vertex of the $s$th pair of $\s$ and $j_s=n_s-i_s$ --- to
the upper vertex. For example:\vspace{5pt}
$$\begin{picture}(60,20)(0,0)
          \put(0,0){\epsfxsize=60pt \epsfbox{schth.eps}}
          \put(-1,-7){\mbox{$\scriptstyle i_1$}}
          \put(-1,20){\mbox{$\scriptstyle j_1$}}
          \put(14,-7){\mbox{$\scriptstyle i_2$}}
          \put(14,20){\mbox{$\scriptstyle j_2$}}
          \put(28,-7){\mbox{$\scriptstyle i_3$}}
          \put(28,20){\mbox{$\scriptstyle j_3$}}
          \put(43,-7){\mbox{$\scriptstyle i_4$}}
          \put(43,20){\mbox{$\scriptstyle j_4$}}
          \put(57,-7){\mbox{$\scriptstyle i_5$}}
          \put(57,20){\mbox{$\scriptstyle j_5$}}
          \end{picture}\quad .\vspace{5pt}$$
Then the {\em monomial corresponding to $\s$} is
$x_{\sigma_0}x_{\sigma_1}\dots x_{\sigma_k}$ where $\sigma_t$ is the
sum of integers assigned to the vertices of $t$th connected
component of $\sigma$.
\end{xdefinition}

\begin{xexample} For the above weighted scheme we get the monomial
$$x_{i_1}x_{j_1+j_2}x_{i_2+i_3+i_4}x_{j_3}x_{j_4+j_5}x_{i_5}\ .$$
\end{xexample}

Now the formula for $P$ can be stated as follows.

\begin{proposition}\label{formula_P}

If $N>n_1+\dots+n_k$ then
$$
   P_{gl_N}(B_{n_1,\dots,n_k}) = \sum_{i_1,\dots,i_k}
   (-1)^{j_1+\dots+j_k} \binom{n_1}{i_1} \dots \binom{n_k}{i_k}
   \sum_\sigma x_{\sigma_0}x_{\sigma_1}\dots x_{\sigma_k},
$$
where the external sum ranges over all integers $i_1,\dots,i_k$ such
that $0\leq i_1\leq n_1, \dots, 0\leq i_k\leq n_k$; the internal sum
ranges over all the $2^{k-1}$ schemes, $j_s=n_s-i_s$, and
$x_{\sigma_0}x_{\sigma_1}\dots x_{\sigma_k}$ is the monomial
associated with the scheme $\sigma$ and integers $i_1,\dots,i_k$.
\end{proposition}

\noindent{\bf Examples.}

(1) For the baguette diagram $B_2$ we have $k=1$, $n_1=2$. There is
only one scheme: $\mbox{\ \begin{picture}(10,10)(0,0)
          \put(5,0){\circle*{2}}
          \put(5,10){\circle*{2}}
          \end{picture} }$\ .
The corresponding monomial is $x_{i_1}x_{j_1}$, and
$$\begin{array}{rcl}
P_{gl_N}(B_2) &=& {\displaystyle
   \sum_{i_1=0}^2  (-1)^{j_1}\binom{2}{i_1} x_{i_1}x_{j_1}}\\
    &=&  x_0x_2 - 2x_1x_1 + x_2x_0 \quad
         = \quad 2(x_0x_2 - x_1^2)
\end{array}$$
which agrees with the example given in Chapter \ref{LAWS} on page
\pageref{B2_ex}.

(2) For the diagram $B_{1,1}$ we have $k=2$, $n_1=n_2=1$. There are
two schemes: $\mbox{\
\begin{picture}(10,10)(0,0)
          \put(0,10){\line(1,0){10}}
          \put(0,0){\circle*{2}}
          \put(10,0){\circle*{2}}
          \put(0,10){\circle*{2}}
          \put(10,10){\circle*{2}}
          \end{picture}\ }$ and
$\mbox{\ \begin{picture}(10,10)(0,0)
          \put(0,0){\line(1,0){10}}
          \put(0,0){\circle*{2}}
          \put(10,0){\circle*{2}}
          \put(0,10){\circle*{2}}
          \put(10,10){\circle*{2}}
          \end{picture} }$\ .
The corresponding monomial are
$x_{i_1}x_{i_2}x_{j_1+j_2}$ and $x_{i_1+i_2}x_{j_1}x_{j_2}$.
We have
$$\begin{array}{l}
P_{gl_N}(B_{1,1}) = {\displaystyle
      \sum_{i_1=0}^1 \sum_{i_2=0}^1 (-1)^{j_1+j_2}
       x_{i_1}x_{i_2}x_{j_1+j_2} +
       \sum_{i_1=0}^1 \sum_{i_2=0}^1 (-1)^{j_1+j_2}
       x_{i_1+i_2}x_{j_1}x_{j_2}} \\
=  x_0x_0x_2\! -\! x_0x_1x_1\! -\! x_1x_0x_1\! +\! x_1x_1x_0\! +\!
      x_0x_1x_1\! -\! x_1x_0x_1\! -\! x_1x_1x_0\! +\! x_2x_0x_0 \\
= 2(x_0^2x_2 - x_0x_1^2)
\end{array}$$

\begin{proof}[Sketch of the proof of Proposition \ref{formula_P}] The
diagram $B_{n_1,\dots,n_k}$ has $k$ parts separated by $k-1$ walls.
Each wall is an edge connecting trivalent vertices to which we shall
refer as {\em wall vertices}. The $s$th part has $n_s$ outgoing
legs. We shall refer to the corresponding trivalent vertices as {\em
leg vertices}.

The proof consists of three steps.
\medskip

Recall that in order to evaluate the universal $\gl_N$ weight system on a 
diagram we can use the graphical procedure of ``resolving'' the trivalent vertices of a diagram and associating a tensor to each of these resolutions, see Sections~\ref{univ_ws_glN} and \ref{glN_B}.
At the first step we study the effect of resolutions of the wall vertices. We
prove that the monomial obtained by certain resolutions of these vertices
has the maximal possible degree if and only if for each wall both
resolutions of its vertices have the same sign. These signs are related to
the above defined schemes in the following way. If we take the positive
resolutions at both endpoints of the wall number $s$, then we connect the
lower vertices of the $s$th and the $(s+1)$st pairs in the scheme.
If we take the negative resolutions, then we connect the
upper vertices.

At the second step we study the effect of resolutions of leg vertices. We
show that the result depends only on the numbers of positive resolutions of
leg vertices in each part and does not depend on which vertices in a part
were resolved positively and which negatively. We denote by $i_s$ the number
of positive resolutions in part $s$. This yields the binomial coefficients
$\binom{n_s}{i_s}$ in the formula of Proposition \ref{formula_P}.
The total number
$j_1+\dots+j_k$ of negative resolutions of leg vertices gives the sign
$(-1)^{j_1+\dots+j_k}$.

The first two steps allow us to consider only those cases where the
resolutions
of the left $i_s$ leg vertices in the part $s$ are positive, the rest
$j_s$ resolutions are negative and both resolutions at the ends of each wall
have the same sign. At the third step we prove that such resolutions of wall
vertices lead to monomials associated with corresponding schemes
according to definition \ref{schememon}.
\medskip

We will make some comments only about the first step, because it is
exactly at this step where Dasbach found an improvement of the original
argument of \cite{CD2}.

Let us fix certain resolutions of all trivalent vertices of
$B_{n_1,\dots,n_k}$. We denote the obtained diagram of $n=n_1+\dots+n_k$ pairs of points and $n$ arrows (see page~\pageref{T-diagram})
 by $T$. After a suitable permutation of the
pairs $T$ will look like a disjoint union of certain $x_m$'s. Hence it
defines a monomial in $x_m$'s which we denote by $m(T)$.

Let us close up all arrows in the diagram by connecting the two points in every
pair with an additional short line. We obtain a number of closed curves,
and we can draw them in such a way that they have 3 intersection points in
the vicinity of each negative resolution and do not have other intersections.
Each variable $x_m$ gives precisely one closed curve. Thus the degree of
$m(T)$ is equal to the number of these
closed curves.

Consider an oriented surface $S$ which
has our family of curves as its boundary (the Seifert surface):
$$\mbox{\begin{picture}(60,0)(0,0)
        \put(0,-10){\ig[width=50pt]{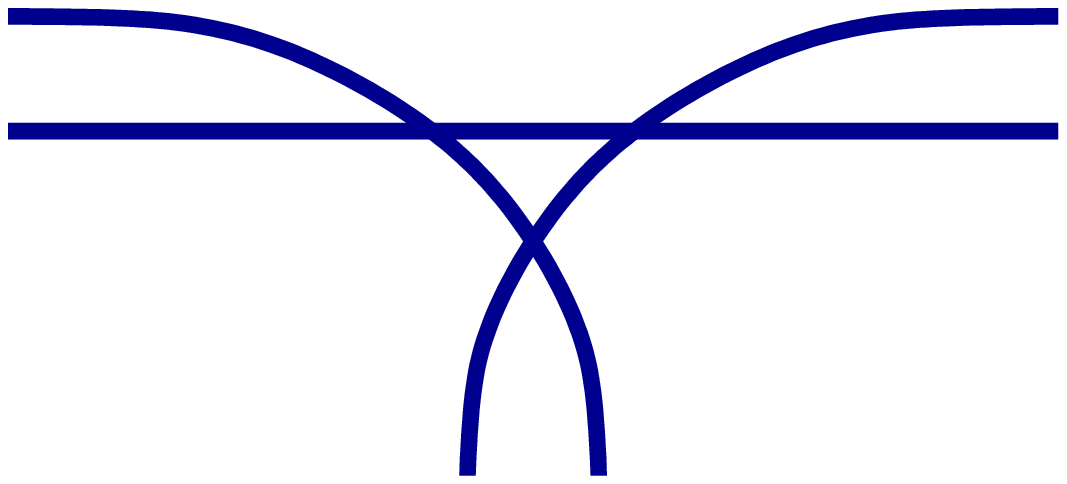}}
        \end{picture}}
   =\quad
  \mbox{\begin{picture}(60,0)(0,0)
        \put(0,-10){\ig[width=50pt]{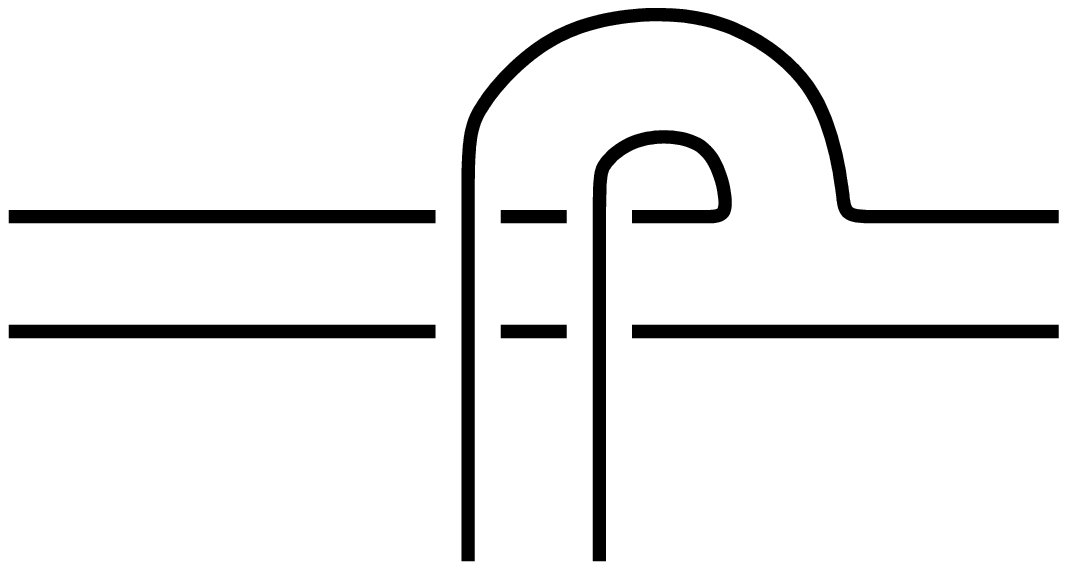}}
        \end{picture}}\quad
  \figsl{toto}\quad\ \
  \rb{-5pt}[0pt][10pt]{\begin{picture}(60,0)(0,0)
        \put(0,-10){\ig[width=50pt]{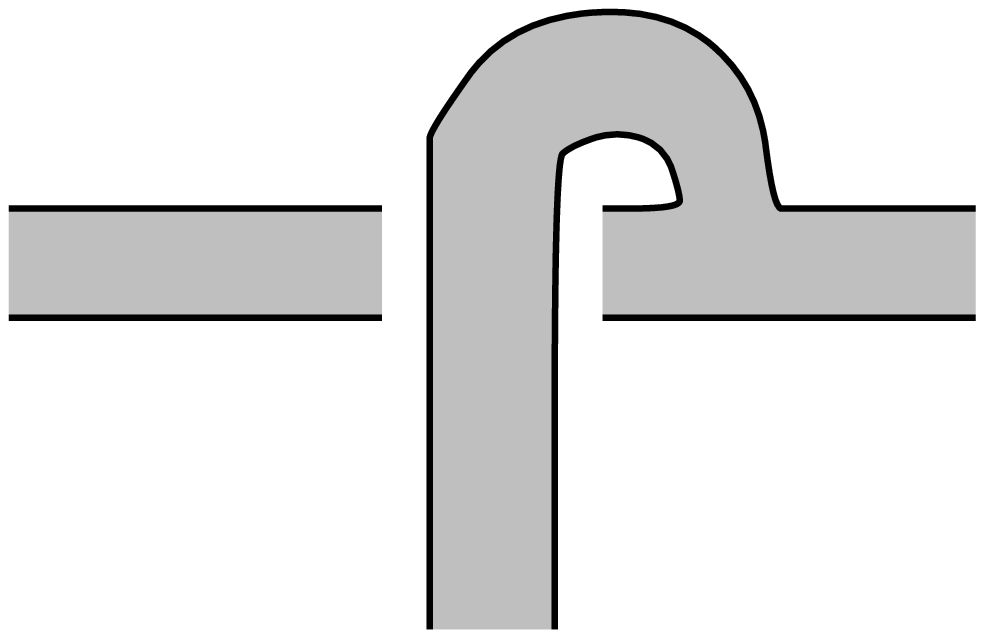}}
        \end{picture}}\ .$$
The degree of $m(T)$ is equal to the number of
boundary components $b$ of $S$.
The whole surface $S$ consists of an annulus
corresponding to the big circle in $B_{n_1,\dots,n_k}$ and $k-1$ bands
corresponding to the walls. Here is an example:
\begin{center}
\begin{picture}(300,68)(0,-10)
     \put(50,-10){\ig[width=200pt]{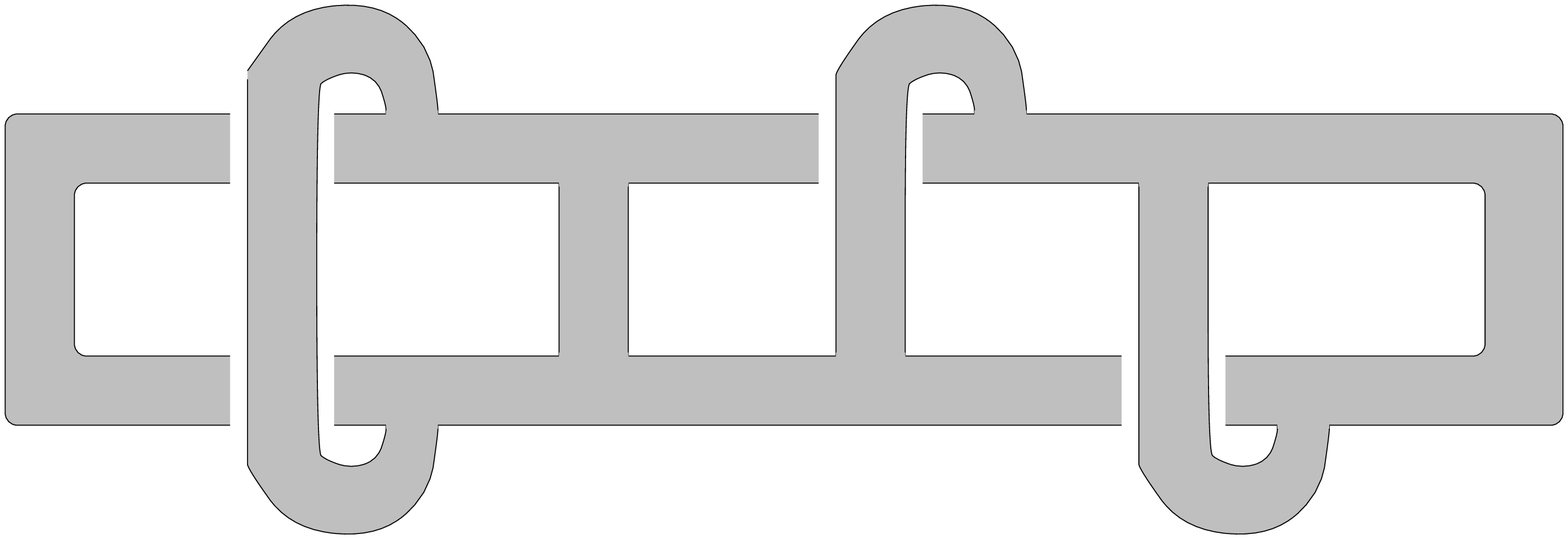}}
\end{picture}
\end{center}
where each of the two walls on the left has the same resolutions at its
endpoints, while the two walls on the right have different
resolutions at their endpoints. The resolutions of the leg vertices
do not influence the surface $S$.

The Euler characteristic $\chi$ of $S$ can be
easily computed. The surface $S$ is contractible
to a circle with $k-1$ chords, thus $\chi = -k+1$. On the
other hand $\chi=2-2g-b$, where $g$ and $b$ are the genus and the number of
boundary components of $S$. Hence $b=k+1-2g$.
Therefore, the degree of $m(T)$, equal to $b$,
attains its maximal value $k+1$ if and only if the surface $S$ has genus 0.

We claim that if there exists a wall whose ends are resolved with the
opposite signs then the genus of $S$ is not zero.
Indeed, in this case we can draw a closed curve in
$S$ which does not separate the surface
(independently of the remaining resolutions):
\begin{center}
\begin{picture}(300,58)(0,0)
     \put(50,0){\ig[width=200pt]{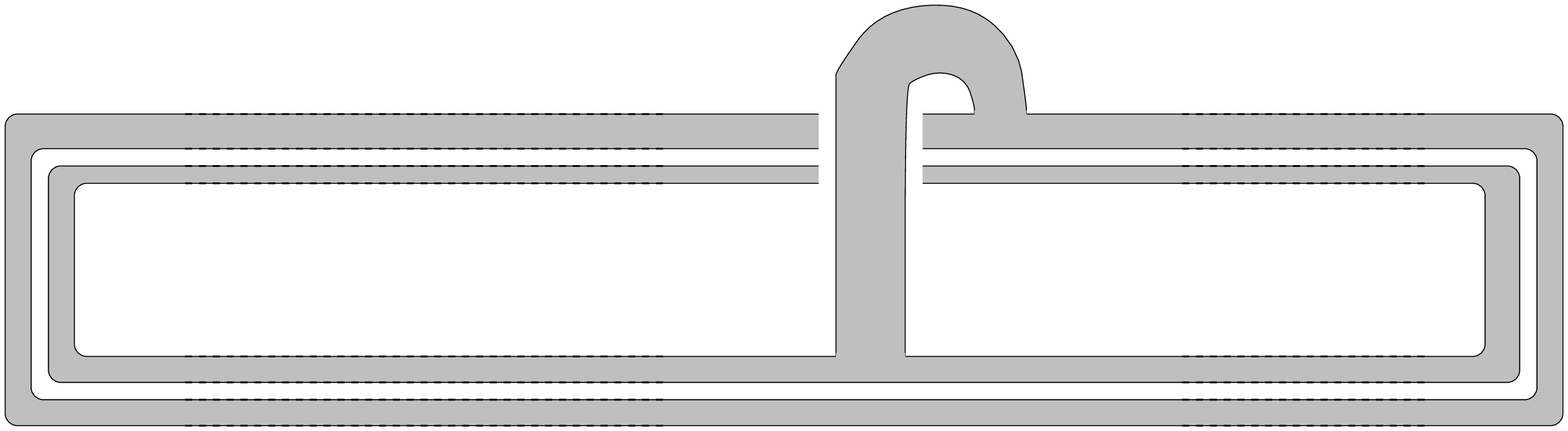}}
\end{picture}
\end{center}

Hence the contribution to $P(B_{n_1,\dots,n_k})$ is given by only
those monomials which come from equal resolutions at the ends of
each wall.
\end{proof}

Now, with Proposition \ref{formula_P} in hand, we can prove the following
result.

\begin{theorem} Let $n=n_1+\dots+n_k$ and $d=n+k-1$.
Baguette diagrams $B_{n_1,\dots,n_k}$ are linearly independent
in $\B$ if $n_1,\dots,n_k$ are all even
and satisfy the following conditions:
\begin{align*}
&    n_1<n_2 \\
&    n_1+n_2<n_3 \\
&    n_1+n_2+n_3<n_4 \\
&    \cdots \cdots \cdots \cdots \cdots \cdots \cdots \\
&    n_1+n_2+\dots+n_{k-2}<n_{k-1} \\
&    n_1+n_2+\dots+n_{k-2}+n_{k-1}<n/3.
\end{align*}
\end{theorem}

The proof is based on the study of the {\it supports\/} of polynomials
$P(B_{n_1,\dots,n_k})$ --- the subsets of $\Z^k$ corresponding to
non-zero terms of the polynomial.

Counting the number of elements described by the theorem, one
arrives at the lower bound $n^{\log(n)}$ for the dimension of the primitive
subspace $\PR_n$ of $\B$.
\medskip

The main difficulty in the above proof is the necessity to consider
the $2^k$ resolutions for the wall vertices of a baguette diagram
that correspond to the zero genus Seifert surface.
O.~Dasbach in \cite{Da3} avoided this difficulty by considering a different
family of open diagrams for which there are only two ways of resolution
of the wall vertices leading to the surface of minimum genus.
These are the {\it Pont-Neuf diagrams\/}:
$$\index{Diagram!Pont-Neuf}\index{Pont-Neuf diagram}\label{pont_neuf}
PN_{a_1,...,a_k,b} \;=\;
\rb{-19mm}{\ig[height=40mm]{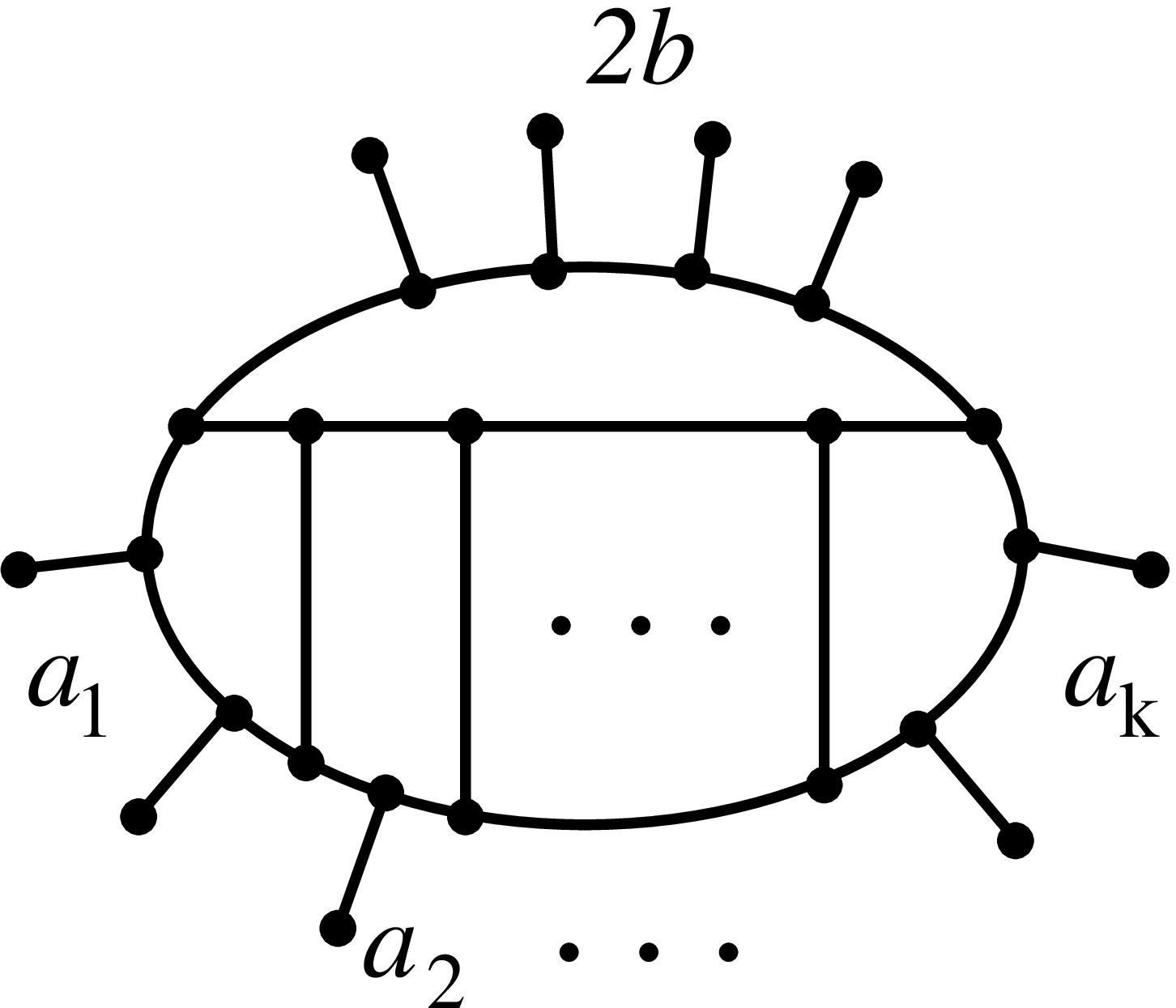}}
$$
(the numbers $a_1$, ..., $a_k$, $2b$ refer to the number of legs
attached to the corresponding edge of the inner diagram).

The reader may wish to check the above property of Pont-Neuf diagrams by way
of exercise. It is remarkable that Pont Neuf diagrams not only lead to
simpler considerations, but they are more numerous, too, and thus lead
to a much better asymptotic estimate for $\dim\PR_n$. The exact statement of
Dasbach's theorem is as follows.

\begin{theorem}
For fixed $n$ and $k$, the diagrams $PN_{a_1,...,a_k,b}$ with
$0\le a_1\le...\le a_k\le b$, $a_1+...+a_k+2b=2n$
are linearly independent.
\end{theorem}

Counting the number of such partitions of $2n$, we obtain precisely the estimate announced by Kontsevich in \cite{Kon1}.
\medskip

\begin{xcorollary}
$\dim\PR_n$ is asymptotically greater than $e^{c\sqrt{n}}$
for any constant $c<\pi\sqrt{2/3}$.
\end{xcorollary}

\begin{xcb}{Exercises}

\begin{enumerate}

\item
\parbox[t]{2.5in}{Show that $\M_1$ is equivalent to the
{\it $\Delta$ move}\index{Moves!$\Delta$}
in the sense that, modulo Reidemeister moves, the $\M_1$ move
can be}
\quad
\rb{-5pt}{$\rs2{50}{delt1}{-20}\!\!\!\!\totonew\!\!\!\!
\rs2{50}{delt2}{-20}\vspace{8pt}$}\vspace{3pt}\\
accomplished by $\Delta$ moves and vise versa.
The fact that $\Delta$ is an unknotting operation was proved in
\cite{Ma,MN}.

\item
Prove that $\M_1$ is equivalent to the move
$$\rs2{75}{vm11}{-20}\totonew\rs2{75}{vm12}{-20}$$

\item
Prove that $\M_2$ is equivalent to the so called {\it clasp-pass}
move\index{Moves!clasp-pass}\index{Clasp-pass move}
$$\rs2{75}{vm21}{-40}\totonew\rs2{75}{vm22}{-40}$$

\item
Prove that $\M_n$ is equivalent to the move ${\mathcal C}_n$:
$$\index{Moves!Goussarov--Habiro!${\mathcal C}_n$}
\index{Goussarov!--Habiro move!${\mathcal C}_n$}
\label{gg-moveCn}
\underbrace{\ig[width=120pt]{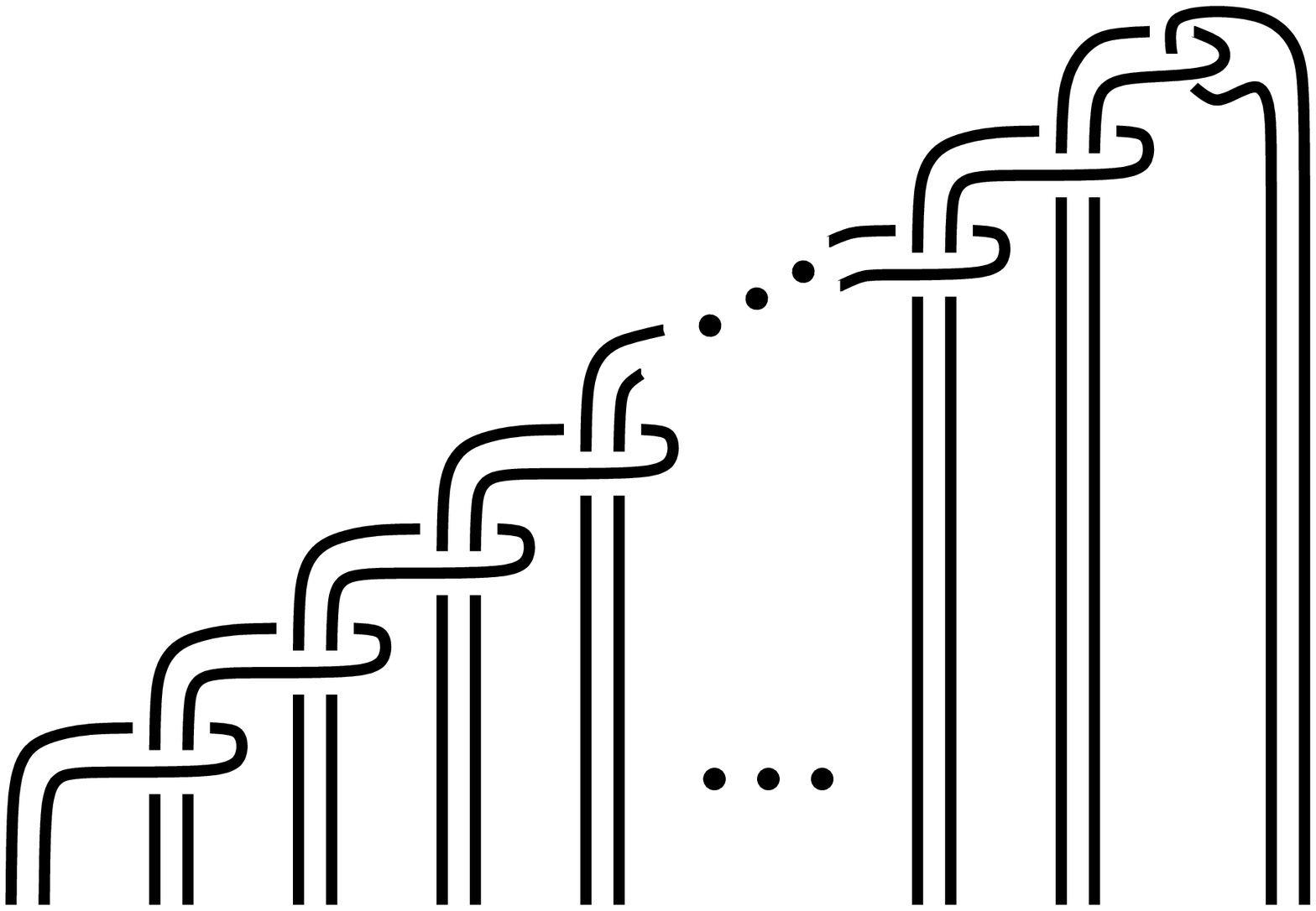}}_{\mbox{\scriptsize $n+2$ components}}
\quad\rb{35pt}{\ig[width=25pt]{totonew.eps}}\quad
\underbrace{\ig[width=120pt]{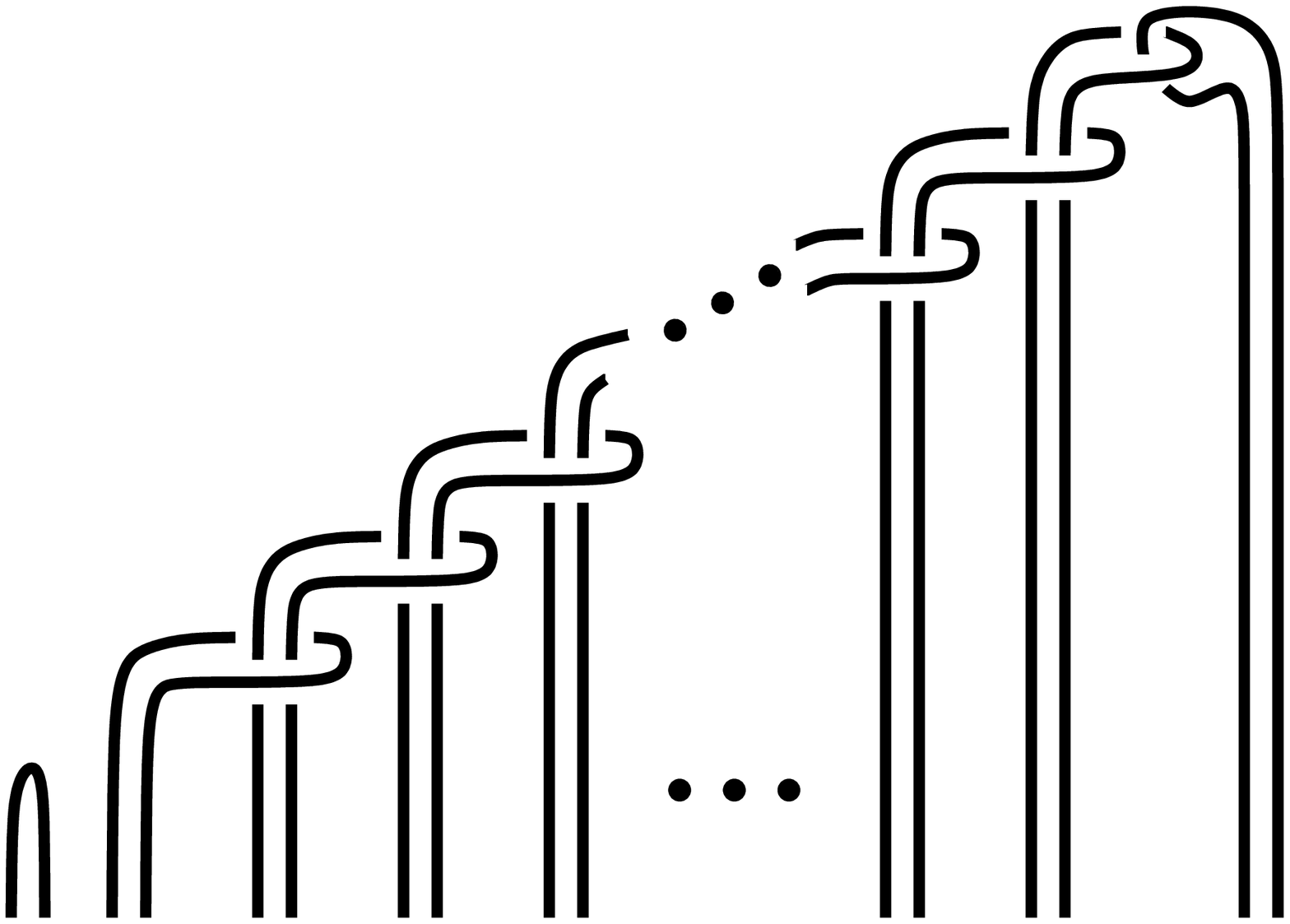}}_{\mbox{\scriptsize $n+2$ components}}
$$

\item Find the inverse element of the knot $3_1$ in the group
$\mathcal{G}_4$.

\item
\begin{enumerate}
\item\neresh (L.~Kauffman) Find a set of moves relating the knots with the
same $c_2$ modulo $n$, for $n=3$, 4,\dots .
\item\neresh  Find a set of moves relating any two knots with the same Vassiliev
invariants modulo 2 (3, 4, ...) up to the order $n$.
\item\neresh Find a set of
moves relating any two knots with the same Conway polynomial.
\end{enumerate}

\item
({\em S. Lando}).
   Let $N$ be a formal variable.
   Prove that $N^{\mathrm{corank}A(G)}$ defines an
algebra homomorphism $\La\to\Z[N]$, where $\La$ is the graph algebra
of Lando and $A(G)$ stands for the adjacency matrix of the graph $G$
considered over the field $\Fi_2$ of two elements.

\item\neresh
Let $\lambda: A \to \La$ be the natural homomorphism from
the algebra of chord diagrams into the graph algebra of Lando.
\begin{itemize}
\item
   Find $\ker\lambda$ (unknown in degrees greater than 7).
\item
   Find $\im\lambda$ (unknown in degrees greater than 7).
\item
   Describe the primitive space $P(\La)$.
\item
   $\La$ is the analog of the algebra of chord diagrams in the case
of intersection graphs. Are there any counterparts of the algebras $\F$
and $\B$?
\end{itemize}

\end{enumerate}
\end{xcb}
 %14 Misc
\chapter{The space of all knots} %15
\label{chap_VSS}

\def\Kspace{{\mathscr K}}

Throughout this book we used the definition of finite type
invariants based on the Vassiliev skein relation. This definition is
justified by the richness of the theory based on it, but it may
appear to be somewhat {\em ad hoc}. In fact, in Vassiliev's original
approach the skein relation is a consequence of a rather
sophisticated construction, which we are going to review briefly in
this chapter.

One basic idea behind Vassiliev's work is that knots, considered as
smooth embeddings $\R^1\to\R^3$, form a topological space $\Kspace$.
An isotopy of a knot can be thought of a continuous path in this
space. Knot invariants are the locally constant functions on
$\Kspace$; therefore, the vector space of $\Ring$-valued invariants,
where $\Ring$ is a ring, is the cohomology group $H^{0}(\Kspace,
\Ring)$. We see that the problem of describing all knot invariants
can be generalized to the following

\noindent{\bf Problem.} {\em Find the cohomology ring
$H^{*}(\Kspace, \Ring)$.}

There are several approaches to this problem. Vassiliev replaces the
study of knots by the study of singular knots with the help of
Alexander duality and then uses simplicial resolutions for the
spaces of singular knots. This method produces a spectral sequence
which can be explicitly described. It is not clear how much
information about the cohomology of the space of knots is contained
in it, but the zero-dimensional classes coming from this spectral
sequence are precisely the Vassiliev invariants.

The second approach is an attempt to build the space of knots out of
the configuration spaces of points in $\R^3$. We are not going to
discuss it here; an instructive explanation of this construction is
given in \cite{Sinha}. Both points of view lead to a new description
of the chord diagram algebra $\A$. It turns out that $\A$ is a part
of an algebraic object which, in a way, is more fundamental than
knots, namely, the {\em Hochshild homology of the Poisson operad}.

It is inevitable that the pre-requisites for this chapter include
rather advanced material such as spectral sequences; at the same
time we delve into less detail. Our goal here is to give a brief
introduction into the subject after which the reader is encouraged
to consult the original sources.

\section{The space of all knots}

First of all, let us give precise definitions.
\begin{xdefinition}
A {\em long curve} is a smooth curve $f:\R\to\R^3$ which at infinity
tends to the diagonal embedding of $\R$ into $\R^3$ :
$$\lim_{t\to\pm\infty} |\,f(t)-(t,t,t)\,| =0.$$
\end{xdefinition}
Here, of course, we could have chosen any fixed linear embedding of
$\R$ into $\R^3$ instead of the diagonal.

There are many ways to organize long curves into a topological
space. For example, one can introduce the $C^0$-metric on the set
$V$ of all long curves with the distance between $f$ and $g$ defined
as
$$d(f,g)=\max_{t\in \R} |\, f(t)-g(t)\,|.$$
Alternatively, let $V_n$ be the set of long curves of the form
$$f(t)=\frac{(P_x(t), P_y(t), P_z(t))}{ (1+t^2)^{-n}},$$
where $P_x$, $P_y$ and $P_z$ are polynomials of the form
$$t^{2n+1}+a_{2n-1} t^{2n-1}+a_{2n-2} t^{2n-2}+\ldots+a_1 t + a_0$$
and $n>0$. (Note the absence of the term of degree $2n$.) We can
consider $V_n$ as a Euclidean space with the coefficients of the
polynomials as coordinates. The space $V_n$ can be identified with
the subspace of $V_{n+1}$ corresponding to the triples of
polynomials divisible by $1+t^2$. Write
$$V_{\infty}=\bigcup_{n=1}^{\infty} V_n$$
\label{V_n}
with the topology of the union (weak topology). We can think of
$V_{\infty}$ as of the space of all {\em polynomial curves}.

\noindent{\bf Exercise.}
\begin{enumerate}
\item
Show that any long curve can be uniformly approximated by polynomial
curves.
\item Is the weak topology on $V_{\infty}$ equivalent to the
topology given by the $C^0$-metric?
\end{enumerate}
\medskip

\begin{xdefinition}
The {\em space of knots} $\Kspace$ is the subset of $V$ consisting
of non-singular curves (smooth embeddings). Similarly, the space
$\Kspace_{n}$, for $n\leq\infty$ is the subspace of smooth
embeddings in $V_{n}$.
\end{xdefinition}
Clearly, $\Kspace$ and $\Kspace_{\infty}$ are just two of the
possible definitions of the space of all knots.

\noindent{\bf Exercise.} Show that the natural map from
$\Kspace_{\infty}$  to $\Kspace$ is a weak homotopy equivalence. In
other words, prove that this map induces a bijection on the set of
connected components and an isomorphism in homotopy groups for each
component. Note that this implies that the cohomology rings of
$\Kspace_{\infty}$  and $\Kspace$ are the same for all coefficients.
\medskip

We shall refer to $\Kspace$  as {\em the space of long
knots}.\index{Knot!long} In
the first chapter we defined long knots as string links on one
string. It is not hard to see that any reasonable definition of a
topology on the space of such string links produces a space that is
weakly homotopy equivalent to $\Kspace$.

\subsection{A remark on the definition of the knot space}
One essential choice that we have made in the definition of the knot
spaces is to consider long curves. As we know, the invariants of
knots in $\R^3$ are the same as those of knots in $S^3$, or those of
long knots. This, however, is no longer true for the higher
invariants of the knot spaces. For example, the component of the
trivial knot in $\Kspace$ is contractible, while in the space of usual 
knots $S^1\to\R^3$ it is not simply-connected, see \cite{Hat1}.

The space of long knots $\Kspace$ has many advantages over the other
types of knot spaces. An important feature of $\Kspace$ is a natural
product
$$\Kspace\times\Kspace\to\Kspace$$
given by the connected sum of long knots. Indeed, the sum of long
knots is defined simply as concatenation, and is well-defined not
just for isotopy classes but for knots as geometric objects\footnote{Well, 
almost. To make this
precise, apply the mapping $(x,y,z)\mapsto(-e^(-x), -e^(-y), -e^(-z))$
to the first knot, the mapping $(x,y,z)\mapsto(e^x, e^y, e^z)$ to the
second, and then glue them together.}. The
connected sum of usual knots, on the contrary, depends on many
choices and is only well-defined as an isotopy class.

\noindent{\bf Exercise.} Show that the product on $\Kspace$ just described
is commutative up to homotopy. Show
that the trivial knot is a unit for this product, up to homotopy.
\medskip

The existence of a homotopy commutative product on $\Kspace$ has
deep consequences for its topology. In fact, it can be shown that
$\Kspace$ is a {\em two-fold loop space}, see \cite{Bud}.

\section{Complements of discriminants}

In this section we shall describe the technical tools necessary for
the construction of the Vassiliev spectral sequence for the space of
knots. This machinery is very general and can be applied in many
situations that are not related to knots in any way; we refer the
reader to Vassiliev's book \cite{Va3} (or its more complete Russian
version \cite{Va7}) for details.

The space $\Kspace$, whose cohomology we are after, is the complement
in the space of all long curves of the closed set whose points
correspond to long curves that fail to be embeddings. In other
words, $\Kspace$ is a complement of a {\em discriminant} in the
space of curves.

The term ``discriminant'' usually denotes the subspace of singular
maps  in the space of all maps between two geometric objects, say,
manifolds. For the discussion that follows the word ``discriminant''
will simply mean ``a closed subvariety in an affine space''.

Vassiliev's construction involves three general technical tools:
Alexander duality, simplicial resolutions and stabilization. Let us
describe them in this order.

\subsection{Alexander duality and the spectral sequence}

If $\Sigma$ is a discriminant, that is, a closed subvariety of an
$N$-dimensional real vector space $V$, the Alexander Duality Theorem
states that
$$\widetilde{H}^q(V-\Sigma,\Z)\simeq \widetilde{H}_{N-q-1}(\onept{\Sigma}),$$
where $0\leq q< N$, the tilde indicates reduced (co)homology and
$\onept{\Sigma}$ is the one-point compactification of $\Sigma$. The
geometric meaning of this isomorphism is as follows: given a cycle
$c$ of dimension $N-q-1$ in $\onept{\Sigma}$ we assign to each
$q$-dimensional cycle in $V-\Sigma$ its linking number with $c$ in
the sphere $S^N:=V\cup \{\infty\}$. This is a $q$-dimensional cocyle
representing the cohomology class dual to the class of $c$. (Here
the integer coefficients can be replaced by coefficients in any
abelian group.)

Now, suppose that the discriminant $\Sigma$ is filtered by closed
subspaces
$$\Sigma_1\subseteq\Sigma_2\subseteq\ldots\subseteq\Sigma_k=\Sigma.$$
Taking one-point compactifications of all terms we get the following
filtration:
$$\onept{\Sigma}_{-1}\subseteq\onept{\Sigma}_0\subseteq\onept{\Sigma}_1
\subseteq\onept{\Sigma}_2\subseteq\ldots\subseteq\onept{\Sigma}_k=\onept{\Sigma},$$
with $\onept{\Sigma}_{-1}=\emptyset$ and
$\onept{\Sigma}_0=\{\infty\}$ being the added point. Then the
homology of $\onept{\Sigma}$ can be studied using the {\em spectral
sequence} arising from this filtration. (We refer the reader to
\cite{Hat2} or \cite{Wei} for basics on spectral sequences.) The
term $E^1_{p,q}$ of this spectral sequence is isomorphic to
${H}_{p+q}(\onept{\Sigma}_p, \onept{\Sigma}_{p-1})$ and the
$E^{\infty}$ term
$$E^{\infty}_m=\bigoplus_{p+q=m} E^{\infty}_{p,q}$$
is associated with $\widetilde{H}_m(\onept{\Sigma})$ in the
following sense: let $_{(i)}\widetilde{H}_m(\onept{\Sigma})$ be the
image of $\widetilde{H}_m(\onept{\Sigma}_i)$ in
$\widetilde{H}_m(\onept{\Sigma})$. Then
$$E^{\infty}_{i,m-i}= {_{(i)}}\widetilde{H}_m(\onept{\Sigma})/ _{(i-1)}\widetilde{H}_m(\onept{\Sigma}).$$

Let us define the cohomological spectral sequence\label{vss} 
$E_{r}^{p,q}$ by setting
$$E_{r}^{p,q}=E^{r}_{-p,N-q-1}$$
and defining the differentials correspondingly, by renaming the
differentials in the homological spectral sequence. According to
Alexander duality, the term $E_{\infty}$ of this sequence is
associated with the cohomology of $V-\Sigma$. All non-zero entries
of this sequence lie in the region $p<0$, $p+q\geq 0$.

The functions on the connected components of $V-\Sigma$ can be
identified with the elements of $H^0(V-\Sigma,\Z)$. The information
about this group is contained in the anti-diagonal entries
$E_{\infty}^{-i,i}$ with positive $i$. Namely, let
$_{(i)}H^0(V-\Sigma,\Z)$ be the subgroup of $H^0$ consisting of the
functions that are obtained as linking numbers with cycles of maximal 
dimension contained in the one-point compactification of $\Sigma_i$; as we shall soon
see, these classes can be thought of as Vassiliev invariants of
order $i$.

\begin{xremark}
The spectral sequence that we just described was first defined by V.~Arnold  \cite{Ar1a} who studied with its help the cohomology of the braid groups. A similar method was later used by G.~Segal in \cite{Seg} to describe the topology of the spaces of rational functions.
\end{xremark}

\subsection{Simplicial resolutions}\label{sec:simpres}

Assume that $f:X\to Y$ is a finite-to-one proper surjective map of
topological spaces. Then $Y$ is obtained from $X$ by identifying the
preimages of each point $y\in Y$. Assume (for simplicity) that there 
exists a constant $R$ such that for any point $y\in Y$ the 
preimage $f^{-1}(y)$ consists of at most $R$ points, and that
$X$ is embedded in some Euclidean space $V$ in such a way that any $k+1$
distinct points of $X$ span a non-degenerate $k$-simplex for all
$k<R$. Then we can form the space $\widetilde{Y}$ as the union, over
all points $y\in Y$, of the convex hulls of the sets $\{f^{-1}(y)\}$
in $V$.

We have a map $\tilde{f}:\widetilde{Y}\to Y$ which assigns to a
point in the convex hull of the set $\{f^{-1}(y)\}$ the point $y\in
Y$. This map is proper and its fibres are simplices, possibly, of
different dimensions. It can be deduced that under mild assumptions
on $Y$ the map $\tilde{f}$ is a homotopy equivalence; it is called a
{\em simplicial resolution of $Y$ associated with $f$}. We shall
refer to the space $\widetilde{Y}$ as the {\em space of the
simplicial resolution $\tilde{f}$}, or, abusing the terminology, as
the simplicial resolution of $Y$.

\begin{xexample} The map of a circle onto the figure eight which
identifies two points has the following simplicial resolution:
$$\ig[width=7cm]{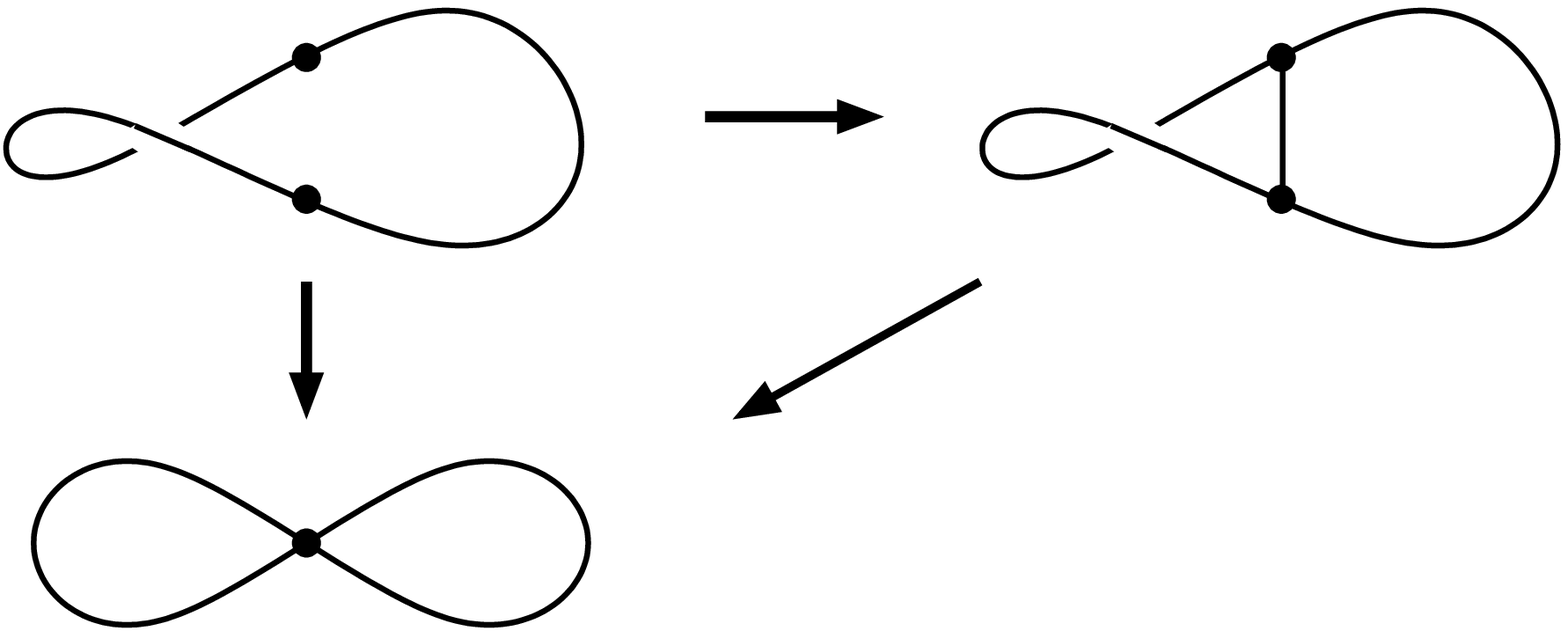}$$
Here the space of the resolution is shown on the right.
\end{xexample}
\noindent{\bf Exercise.} Describe the simplicial resolution
associated with the double cover of a circle by itself. \smallskip

Since we are interested only in calculating the homology groups, we
lose nothing by replacing a space by the space of its simplicial
resolution. On the other hand, simplicial resolutions often have
interesting filtrations on them. For instance, since the space
$\widetilde{Y}$ is a union of simplices, it is natural to consider
the subspaces $\widetilde{Y}_i$ of $\widetilde{Y}$ which are the
unions of the $i-1$-skeleta of these simplices. In the case of the
discriminant in the space of long curves we shall consider another
geometrically natural filtration, see
Section~\ref{sec:vss-explicit}.\smallskip

\noindent{\bf Exercise.} Adapt the results of the preceding section
so that instead of a filtration on a discriminant one can use a
filtration on its simplicial resolution. \smallskip

Simplicial resolutions are especially useful for studying spaces of
functions with singularities of some kind. In such a situation $Y$
is taken to be the space of functions with singularities and $X$ is
the space of all pairs $(\phi,x)$ where $\phi\in Y$ is a function
and $x$ is a point in the domain of $\phi$ where $\phi$ is singular;
the map $X\to Y$ simply forgets the second element of the pair (that
is, the singular point). Various examples of this kind are described
in Vassiliev's book \cite{Va3, Va7}. While many of the ingredients 
of Vassiliev's approach to knot spaces were well-known before 
Vassiliev, the simplicial resolutions are the main innovation of his work.

\subsection{Stabilization}\label{sec:stab} Strictly
speaking, Alexander duality and, as a consequence, all the
foregoing constructions, only makes sense in finite-dimensional
spaces. However, in the case of knots the space $V$ is
infinite-dimensional. This problem can be circumvented by using
finite-dimensional approximations to the space of long curves. For
this we have to understand first how complements of discriminants
behave with respect to inclusions.

Consider two discriminants $\Sigma_1$ and $\Sigma_2$ inside the
Euclidean spaces $V_1$ and $V_2$ respectively. If $V_1\subset V_2$
and $\Sigma_1=V_1\cap\Sigma_2$ we see that $V_1-\Sigma_1$ is a
subspace of $V_2-\Sigma_2$. We would like to describe the induced
map in cohomology
\begin{equation}\label{eq:stab}
H^i(V_2-\Sigma_2,\Z)\to H^i(V_1-\Sigma_1,\Z).
\end{equation}

Assume that $V_1$ intersects $\Sigma_2$ transversally, so that there
exists an $\e$-neighbourhood $V_{\e}$ of $V_1$ such that
$$V_{\e}\cap\Sigma_2=\Sigma_1\times\R^{s},$$
where $s=\dim{V_2}-\dim{V_1}$. There is a homomorphism of reduced
homology groups
$$\widetilde{H}_i(\onept{\Sigma_2})\to \widetilde{H}_{i-s}(\onept{\Sigma_1})$$
where $\onept{X}$, as before, denotes the one-point compactification
of $X$. This homomorphism is known as the {\em Pontrjagin-Thom
homomorphism} and is constructed in two steps. First, we collapse
the part of $\onept{\Sigma_2}$ which lies outside of $\Sigma_2\cap
V_{\e}$ to one point and take the induced homomorphism in homology.
Then notice that the quotient space with respect to this collapsing
map is precisely the $s$-fold suspension of $\onept{\Sigma_1}$, so
we can apply the suspension isomorphism which decreases the degree
by $s$ and lands in the homology of $\onept{\Sigma_1}$.

Since Alexander duality is defined by taking linking numbers, it
follows from this construction that the cohomology map
(\ref{eq:stab}) is dual to the Pontrjagin-Thom homomorphism.

Now, let us consider the situation when both discriminants
$\Sigma_1$ and $\Sigma_2$ are filtered and
$V_{\e}\cap(\Sigma_2)_j=(\Sigma_1)_j\times\R^s$ for all $j$. Then we
have relative Pontrjagin-Thom maps
$H_i((\onept{\Sigma_2})_j,(\onept{\Sigma_2})_{j-1})\to
H_{i-s}((\onept{\Sigma_1})_j,(\onept{\Sigma_1})_{j-1})$.
\begin{xproposition}
If, in the above notation and under the above assumptions, the
relative Pontrjagin-Thom maps are isomorphisms for all $j\leq P$ and
$i >\dim{V}-Q+j$, for some positive $P$ and $Q$, the terms
$E_1^{p,q}$ of the Vassiliev spectral sequences for the cohomology
of $V_1-\Sigma_1$ and $V_2-\Sigma_2$ coincide in the region $-P\leq
p$ and $q\leq Q$.
\end{xproposition}
The proof consists in combining Alexander duality with the
definition of the spectral sequence.

The above proposition will allow us to work in infinite-dimensional
Euclidean spaces as if they had finite dimension, see
Section~\ref{sec:vss-explicit}.

\subsection{Vassiliev invariants}

Suppose that we want to enumerate the connected components of the
complement of a discriminant $\Sigma$ in a vector space $V$; in
other words, we would like to calculate $H^0(V-\Sigma,\Ring)$, the
space of $\Ring$-valued functions on the set of connected components
of $V-\Sigma$. If $\Sigma$ is filtered by closed subspaces
$\Sigma_i$, we can define {\em Vassiliev invariants}
\index{Vassiliev invariant} for the connected components of
$V-\Sigma$ as follows.
\begin{xdefinition}
A Vassiliev invariant of degree $i$ is an element of
$H^0(V-\Sigma,\Ring)$ defined as the linking number with a cycle in
$H_{\dim{V}-1}(\onept{\Sigma_i},\Ring)$.
\end{xdefinition}

This definition also makes sense when we only have the filtration on
the homology of $\Sigma$, rather than on the space $\Sigma$ itself.
Such a situation arises when we consider a filtration on a
simplicial resolution of $\Sigma$. Let us consider the following
rather special situation where the Vassiliev invariants have a
transparent geometric interpretation.

Let $\Sigma$ be the image of a smooth manifold $X$ immersed in a
finite-dimensional vector space $V$, and assume that each point in
$V$, where $\Sigma$ has a singularity, is a point of transversal
$k$-fold self-intersection\footnote{that is, of $k+1$ sheets of $X$} for some $k$. 
Without loss of generality we
can suppose that $\Sigma$ is of codimension 1, since its complement
would be connected otherwise. Locally, $\Sigma$ looks like
$T^{k}\times \R^{\dim{V}-k}$ where $T^{k}$ the union of all
coordinate hyperplanes in $\R^k$. We shall also assume that $\Sigma$
is {\em co-oriented}, that is, that there is a continuous field of
unit normal vectors ({\em co-orientation})\index{Co-orientation} at
the smooth points of $\Sigma$ which extends to the self-intersection
points as a multivalued vector field.

Consider the simplicial resolution ${\sigma}\to\Sigma$ associated
with the map $X\to\Sigma$, and the filtration $\sigma_i$ on
${\sigma}$ by the $i-1$-skeleta of the inverse images of points of
$\Sigma$.  Then we have the following criterion for an element $f\in
H^0(V-\Sigma,\Ring)$ to be a Vassiliev invariant of order $n$.

\begin{figure}[ht]\label{fig:x}
$${\begin{picture}(150,150)(0,0)
  \put(0,0){\ig[width=5cm]{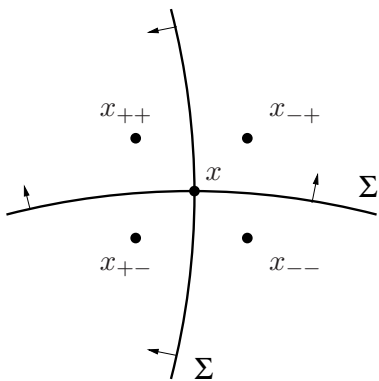}}
  \put(76,76){$x$}
  \put(36,100){$x_{++}$}
  \put(100,100){$x_{-+}$}
  \put(36,42){$x_{+-}$}
  \put(100,42){$x_{--}$}
\end{picture}}
$$
\caption{The neighbourhood of a generic self-intersection point of
the discriminant $\Sigma$.}
\bigskip

\end{figure}

Grouping together the points of the discriminant $\Sigma$ according
to the multiplicity of self-intersection at each point, we get a
decomposition of $\Sigma$ into a union of open strata
$\Sigma_{(i)}$, with $\Sigma_{(i)}$ consisting of points of $i$-fold
intersection and having codimension $i$ in $\Sigma$.
 The function $f$ can be extended from $V-\Sigma$ to a locally
constant function on each stratum of $\Sigma$. If $x$ is a point on
the maximal stratum $\Sigma_{(0)}$ which consists of the points
where $\Sigma$ is smooth, let $x_+$ and $x_-$ be two points in
$V-\Sigma$ obtained by shifting $x$ by ${\pm \e}$ in the direction
of the co-orientation, where $\e$ is small. Then, we set
$f(x)=f(x_+)-f(x_-)$. For $x\in\Sigma_{(1)}$ we take four points
$x_{++}$, $x_{+-}$, $x_{-+}$ and $x_{--}$ obtained from $x$ by
shifting it to each of the four quadrants in $V-\Sigma$, see
Figure~\ref{fig:x}.

For such $x$ we set $f(x)=f(x_{++})-f(x_{+-})-f(x_{-+})+f(x_{--})$.
It is clear how to continue: at a point $x\in\Sigma_{(k)}$ the value
of $f$ is the alternating sum of its values at the $2^{k+1}$ points
obtained by shifting $x$ to the $2^{k+1}$ adjacent quadrants of
$V-\Sigma$.

\begin{xproposition}
An element $f\in H^0(V-\Sigma,\Ring)$ is a Vassiliev invariant of
order $n$ if and only if its extension to the stratum $\Sigma_{(n)}$
of $n$-fold self-intersection
points of $\Sigma$ is identically equal to zero.
\end{xproposition}

This second characterization of the Vassiliev invariants in terms of
their extensions to the strata of the discriminant is the definition
that we used throughout the book. Indeed, a generic point of the
discriminant in the space of long curves is a singular knot with one
simple double point. Knots with two simple double points correspond
to transversal self-intersections of the discriminant, et cetera.
Note, however, that this proposition, as stated, does not apply directly
to the space of knots, since the discriminant in this case has
singularities more complicated than transversal self-intersections.
It turns out that these singularities have no influence on the
homology of the discriminant in the relevant dimensions, and can be
omitted from consideration.

\begin{xxca}
Let $\Sigma$ be the union of the coordinate hyperplanes $x_i=0$ in
$V=\R^n$, $X$ be the disjoint union of these hyperplanes and $X\to
\Sigma$ be the natural projection. Describe the cycles that
represent classes in $H_{n-1}(\onept{\sigma_i})$ and the space
of Vassiliev invariants of degree $i$. Show that the Vassiliev
invariants distinguish the connected components of $V-\Sigma$.
\end{xxca}

\begin{proof}[Sketch of the proof]
We are interested in the cycles of top dimension on
$\onept{\Sigma}$, and these are linear combinations of the (closures
of) connected components of $\Sigma_{(1)}$, that is, the top-dimensional
stratum consisting of smooth points.

By construction of the simplicial resolution, any cycle that locally
is diffeomorphic to the boundary of the ``$k-1$-corner'' $$\{
(x_1,\ldots, x_{n})\, |\, x_1>0,\ldots, x_k>0\}$$ in
$\R^{n}$ defines a homology class in
$H_{n-1}(\onept{\sigma_i},\Ring)$, where $i\geq k$. Conversely,
any cycle in $H_{n-1}(\onept{\sigma_i},\Ring)$ after projection 
to $\onept{\Sigma}$ locally looks
like a linear combination of $k-1$-corners with $k\leq i$.

Now it remains to observe that, locally, the linking number with a
$(k-1)$-corner vanishes on the strata of dimension $k+1$ and,
moreover, this property defines linear combinations of cycles that
locally look like $j$-corners with $j<k$.
\end{proof}

\section{The space of singular knots and Vassiliev invariants}

We want to relate the topology of the space of knots to the
structure of the discriminant in the space of long curves. In order
to use one of our main tools, namely,  Alexander duality, we need
to use finite-dimensional approximations to the space of long
curves. Spaces $V_n$ of polynomial curves of bounded degree 
(see p. \pageref{V_n}) provide
such an approximation; however, it cannot be used directly for
the following reason:
\smallskip

\noindent {\bf Exercise.} Denote by $\Sigma(V_n)$ the discriminant
in the space $V_n$ consisting of non-embeddings. Then the
intersection of $V_{n-1}$ with $\Sigma(V_n)$ inside $V_n$ is not
transversal.
\smallskip

As a consequence, we cannot apply the stabilization procedure
described in Section~\ref{sec:stab}, since it requires
transversality in an essential way. Nevertheless, this is only a
minor technical problem.

\subsection{Good approximations to the space of long curves}
Let $$U_1\subset U_2\subset\ldots \subset U_n\subset\ldots\subset
V_{\infty}$$ be a sequence of finite-dimensional affine subspaces in the space of all
polynomial long curves. Note that each $U_j$ is contained in a
subspace $V_{k}$ for some finite $k$ that depends on $j$. We say
that the sequence $(U_j)$ is a {\em good approximation} to the space
of long curves if
\begin{itemize}
\item each finite-dimensional continuous family of long curves can be uniformly
approximated by a continuous family of curves from $(U_j)$ for some $j$;
\item  for each $U_j$ and each  $V_n$ that contains $U_j$
the intersection of $U_{j}$ with $\Sigma(V_n)$ is transversal inside
$V_n$, where $\Sigma(V_n)$ is the subspace in $V_n$ consisting of
non-embeddings.
\end{itemize}
Write $\Kspace'_j$ for the (topological) subspace of $U_j$ that consists of knots.
The first of the two listed conditions guarantees that the union of
all the spaces $\Kspace'_j$ has the same homotopy, and, hence,
cohomology groups as the space of all knots $\Kspace$. Indeed, it
allows us to approximate homotopy classes of maps $S^n\to\Kspace$ and
homotopies among them by maps and homotopies whose images are
contained in the $\Kspace'_j$. The second condition is to ensure
that we can use the stability criterion from Section~\ref{sec:stab}.

A general position argument gives the following
\begin{xproposition}[Vassiliev \cite{Va3, Va7}]
Good approximations to the space of long curves exist.
\end{xproposition}

The precise form of a good approximation will be unimportant for us.
One crucial property of good approximations is the following:
\smallskip

\noindent {\bf Exercise.} Show that good approximations only contain
long curves with a finite number of singular points (that is, points
where the tangent vector to the curve vanishes) and
self-intersections.
\smallskip

\noindent{\sl Hint:} Show that long curves with an infinite number
of self-intersections and singular points form a subset of infinite
codimension in $V_{\infty}$.

In what follows by the ``space of long curves'' we shall mean the
union $U_{\infty}$ of all the $U_j$ from a good approximation to the
space of long curves and by the ``space of knots'' we shall
understand the space $\Kspace'_{\infty}=\cup_j \Kspace'_j$ 
constructed with the help
of this approximation.

\subsection{Degenerate chord diagrams}

The discriminant in the space of long curves consists of various
parts ({\em strata}) that correspond to various types of {\em
singular knots},\index{Knot!singular} that is, long curves with
self-intersections and singular points.

In our definition of Vassiliev invariants in Chapter~\ref{FT_inv}, we
associated a chord diagram with $n$ chords to a knot with $n$ double
points. In fact, we saw that chord diagrams are precisely the
equivalence classes of knots with double points modulo isotopies and
crossing changes. If we consider knots with more complicated
self-intersections and with singular points we must generalize the
notion of chord diagram.

A {\em degenerate chord diagram}\index{Chord diagram! degenerate} is
a set of distinct pairs $(x_k,y_k)$ of real numbers (called {\em
vertices}) with $x_k\leq y_k$. These pairs can be thought of as
chords on $\R$, with $x_k$ and $y_k$ being the left and right
endpoints of the chords respectively. If all the $x_k$ and $y_k$ are
distinct, we have a usual linear chord diagram.

The ``degeneracy'' of a degenerate chord diagram can be of two
kinds: one chord can degenerate into a singular point ($x_k=y_k$) or
two chords can glue together and share an endpoint. Two degenerate
chord diagrams are {\em combinatorially equivalent} if there is an
self-homeomorphism of $\R$ that preserves the orientation and sends
one diagram to the other.

The vertices of a degenerate chord diagram are of two types: the
{\em singularity} vertices which participate in chords with
$x_k=y_k$, and {\em self-intersection} vertices which participate in
chords with $x_k<y_k$. The same vertex can be a singularity vertex
and a self-intersection vertex at the same time; in this case we
shall count it twice and say that a singularity vertex coincides
with a self-intersection vertex. As with usual chord diagrams, one
can speak of the {\em internal graph}\index{Graph!internal}
of a degenerate diagram: this
is the abstract graph formed by the chords whose ends are distinct.
The self-intersection vertices are divided into {\em groups}: two
vertices belong to the same group if and only if they belong to the
same connected component of the internal graph. Here is an example
of a degenerate chord diagram with two groups of self-intersection
vertices; the singularity vertices are indicated by hollow dots:
$$\ig[width=5cm]{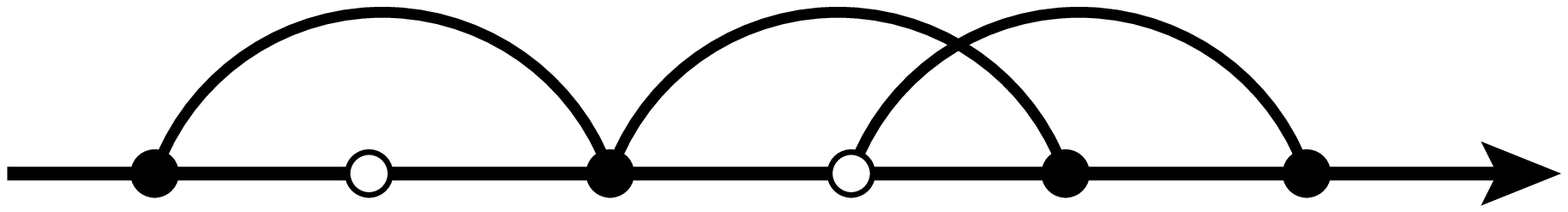}$$

Let us say that two degenerate chord diagrams $D_1$ and $D_2$ are
{\em equivalent} if $D_1$ is combinatorially equivalent to a diagram
that has the same set of singularity vertices and the same groups of
self-intersection vertices as $D_2$. For instance, the following
diagrams are equivalent:
$$\ig[width=2.5cm]{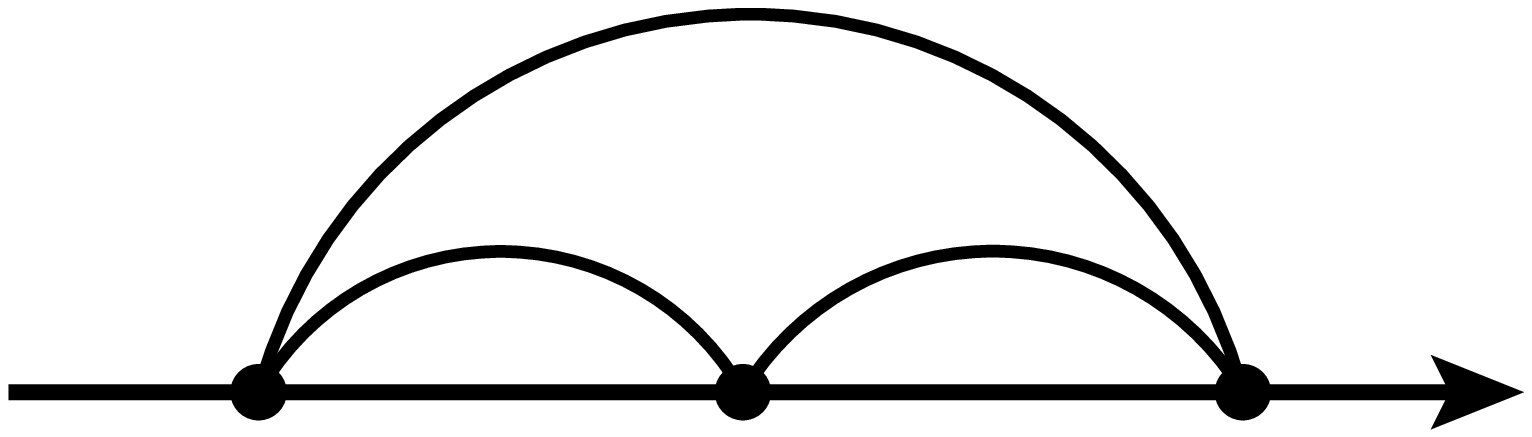}\qquad\ig[width=2.5cm]{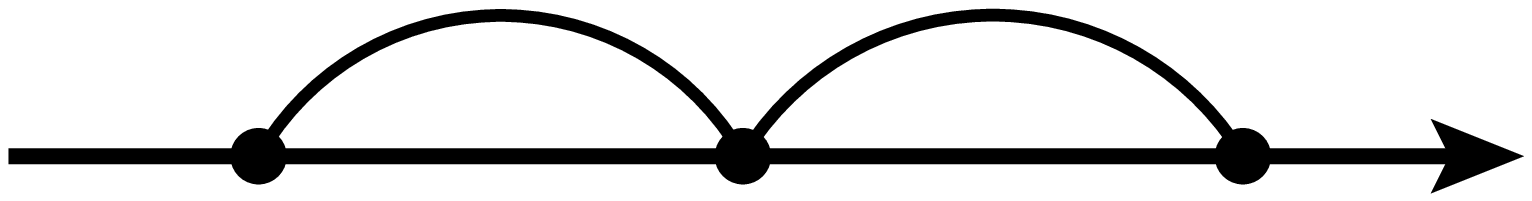}
\qquad\ig[width=2.5cm]{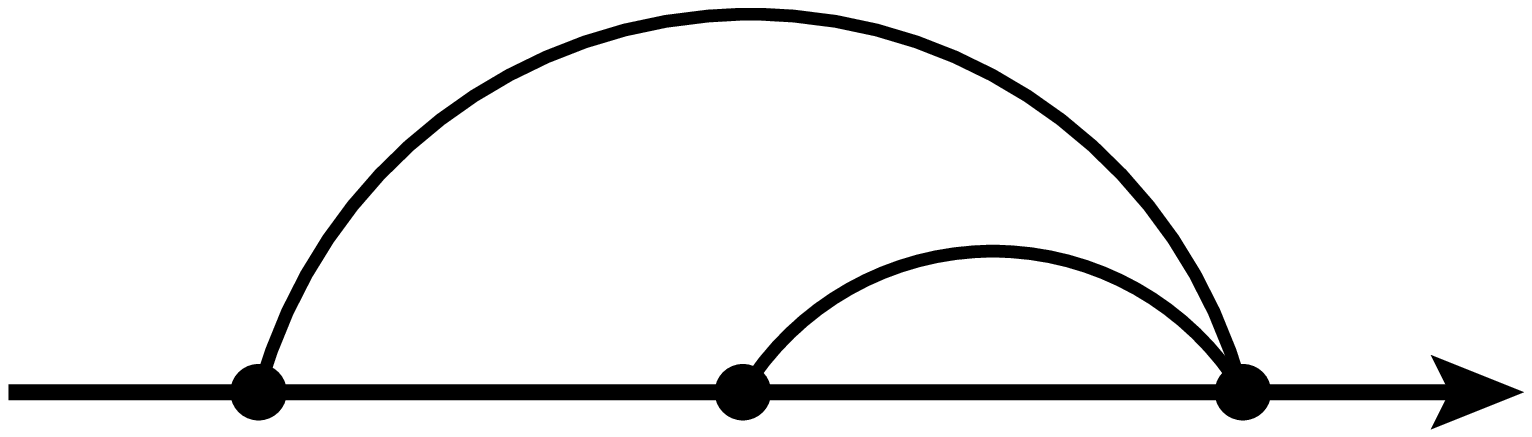}\qquad\ig[width=2.5cm]{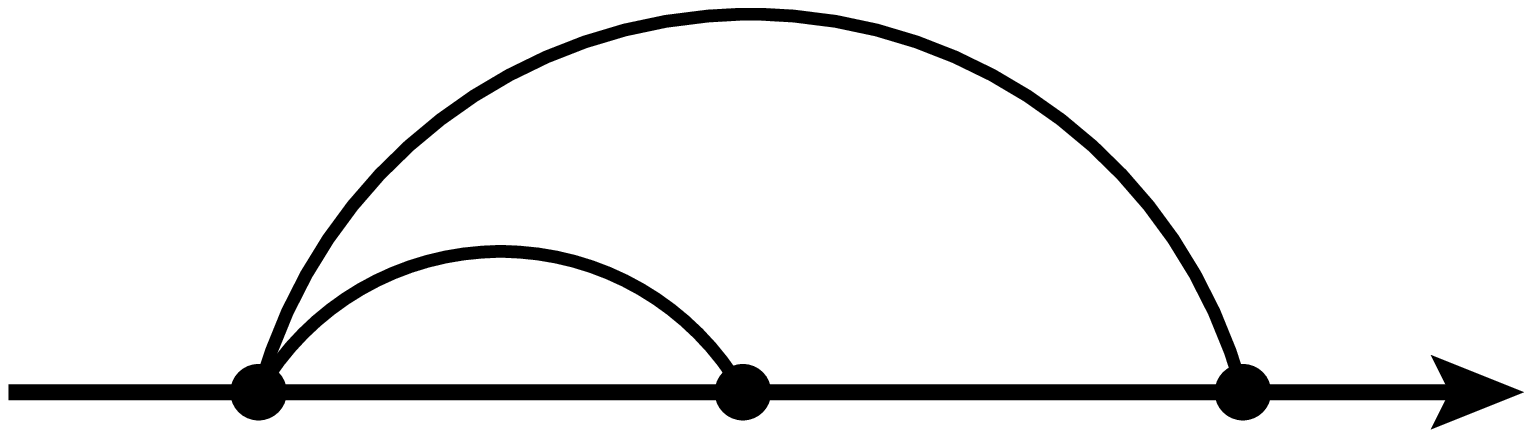}$$
In Vassiliev's terminology, equivalence classes of degenerate chord
diagrams are called {\em
$(A,b)$-configurations}.\index{$(A,b)$-configuration}

Define the {\em complexity}\index{Complexity!of a degenerate chord
diagram} of a degenerate chord diagram as the total number of its
vertices minus the number of groups of its self-intersection
vertices. This number only depends on the equivalence class of the
diagram. The complexity of a usual chord diagram is, clearly, equal
to its degree.

\begin{xremark}
Note that the number of chords of a degenerate chord diagram is not
mentioned in the definitions of equivalence and complexity. This is
only natural, of course, since equivalent diagrams can have
different numbers of chords.
\end{xremark}

An arbitrary singular knot with a finite number of singular points
and self-intersections defines an equivalence class of degenerate
chord diagrams: each group of self-intersection vertices is a
preimage of a self-intersection and the singularity vertices are the
preimages of the singularities. We shall say that a singular knot
$f:\R\to\R^3$ {\em respects} a degenerate chord diagram $D$ if it
glues together all points within each group of self-intersection
vertices of $D$ and its tangent vector is zero at each singularity
vertex of $D$.

\noindent{\bf Exercise.} Show that equivalence classes of degenerate
chord diagrams coincide with the equivalence classes of singular
knots with a finite number of singular points and self-intersections
under isotopies and crossing changes.

\subsection{The discriminant}\label{sec:vss-explicit}
The discriminant in the space of long curves $U_{\infty}$
is a complicated set. Its strata
can be enumerated: they correspond to equivalence classes of
degenerate chord diagrams. However, the structure of these strata is
not easy to describe, since they can (and do) have
self-intersections. The most convenient tool for studying the
discriminant are the simplicial resolutions described in
Section~\ref{sec:simpres}.

In order to tame the multitude of indices, let us write simply
$U$ for the approximating space $U_j$, $N=N_j$ for its  dimension
and $\Sigma=\Sigma(U_j)$ for the discriminant. Write $Sym^2(\R)$ for
the space of all unordered pairs of points in $\R$; this space can
be thought of as the subset of $\R^2$ defined by the inequality
$x\leq y$.

In the product space $\Sigma\times Sym^2(\R)$ consider the subspace
$\widetilde{\Sigma}$ consisting of pairs $(f,(x,y))$ such that
either $x\neq y$ and $f(x)=f(y)$ or $x=y$ and $f'(x)=0$. Forgetting
the pair $(x,y)$ gives a map $\widetilde{\Sigma}\to\Sigma$ which is
finite-to-one and proper, so we can associate a simplicial
resolution with it. (Strictly speaking, in order to define a
simplicial resolution, we must embed $\widetilde{\Sigma}$ in a
Euclidean space in a particular way, but let us sweep this issue
under the carpet and refer to \cite{Va3}.)

Denote the space of this simplicial resolution by $\sigma$. A point
in $\sigma$ is uniquely described by a collection
$$\bigl(f,(x_0, y_0),\ldots,(x_k, y_k), \tau\bigr),$$
where $f$ is a singular knot with $f(x_j)=f(y_j)$ whenever $x_j\neq
y_j$ and $f'(x_j)=0$ when $x_j=y_j$, all pairs $(x_j, y_j)$ are
distinct and $\tau$ is a point in the interior of a $k$-simplex with
vertices labelled by the points $(x_j, y_j)$. In other words, a
point in $\sigma$ is a triple consisting of a singular knot $f$, a
degenerate chord diagram $D$ which $f$ respects, and a point $\tau$
in a simplex whose vertices are labelled by the chords of $D$. Here
$k$ can be arbitrary.

Let $\sigma_i$ be the closed subspace of $\sigma$ consisting of
triples $(f,D,\tau)$ with $D$ of complexity at most $i$. The {\em
cohomological Vassiliev spectral sequence for the space of
knots}\index{Vassiliev!spectral sequence!for the space of knots} is
the spectral sequence that comes from the filtration of
$\sigma$ by the $\sigma_i$ (see page~\pageref{vss}). We have
$$E_{1}^{p,q}=E^{1}_{-p,N-q-1}={H}_{N-(p+q+1)}(\onept{\sigma}_p,
\onept{\sigma}_{p-1})=\widetilde{H}_{N-(p+q+1)}(\onept{\sigma}_p/
\onept{\sigma}_{p-1}).$$ Note that
$\onept{\sigma}_p/\onept{\sigma}_{p-1}$ is homeomorphic to the
one-point compactification of $\sigma_p-\sigma_{p-1}$, and this
space can be described rather explicitly, at least when the
dimension of $U$ is sufficiently large.

Indeed, the condition that a singular knot $f\in U\subset V_d$
respects a degenerate chord diagram $D$ produces several linear
constraints on the coefficients of the polynomials $P_1, P_2$ and
$P_3$ which determine $f$. Namely, if $(x_j,y_j)$ is a chord of $D$,
the polynomials satisfy the conditions
$$\frac{P_{\alpha}(x_j)}{(1+x_j^2)^d}=\frac{P_{\alpha}(y_j)}{(1+y_j^2)^d}$$
when $x_j<y_j$, and
$$\left(\frac{P_{\alpha}(x_j)}{(1+x_j^2)^d}\right)'=0$$
if $x_j=y_j$. Each of these conditions with $x_j$ and $y_j$ fixed
gives one linear equation on the coefficients of each of the
polynomials $P_1, P_2$ and $P_3$.

In general, these linear equations may be linearly dependent.
However, the rank of this system of equations can be explicitly
calculated when the dimension of $U$ is large.
\smallskip

\noindent{\bf Exercise.} Show that for a given degenerate chord
diagram of complexity $p$ there exists $N_0$ such that for $N>N_0$
the number of linearly independent conditions on the coefficients of
the $P_{\alpha}$ is equal to exactly $3p$.

\noindent{\sl Hint.} For any set of distinct real numbers
$x_1,\ldots, x_k$ there exists $d$ such that the vectors
$(1,x_i,x_i^2,\ldots,x_i^d)$ are linearly independent.\smallskip

This exercise shows that, for $N$ sufficiently big, the forgetful
map that sends a triple $(f, D,\tau)\in \sigma_p-\sigma_{p-1}$ to
the pair $(D, \tau)$ is an affine bundle with the fibre of dimension
$N-3p$ over a base $W_p$ which only depends on $p$. In particular,
we have the Thom isomorphism
$$
\widetilde{H}_{N-s}(\onept{(\sigma_p-\sigma_{p-1})})
={H}_{3p-s}(\onept{W_p}).
$$ 
As a consequence, when $p>0$
and $q\geq p$, the first term of the Vassiliev spectral sequence has
the entries
$$
E_{1}^{-p,q}=E^1_{p,N-q-1}=\widetilde{H}_{N+p-q-1}(\onept{(\sigma_p-\sigma_{p-1})})
={H}_{4p-q-1}(\onept{W_p}),
$$ 
and for all other values of
$p$ and $q$ the corresponding entry is zero. The space $\onept{W_p}$
whose homology is, therefore, so important for the theory of
Vassiliev invariants, will be called here {\em the diagram
complex}\index{Diagram complex}.

\section{Topology of the diagram complex}

The diagram complexes $\onept{W_p}$ are constructed out of simplices
and spaces of degenerate chord diagrams of complexity $p$.

\subsection{A cell decomposition for the diagram complex}
The space $\onept{W_p}$ has a cell decomposition with cells indexed
by degenerate chord diagrams (more precisely, combinatorial
equivalence classes of such diagrams) of complexity $p$. The cell
$[D]$ is a product of an open simplex $\Delta_D$ whose vertices are
indexed by the chords of $D$, and the space $E_D$ of all diagrams
combinatorially equivalent to $D$. This latter space is also an open
simplex, of dimension $k$, where $k$ is the number of geometrically distinct vertices of $D$.
Indeed, it is homeomorphic to the configuration space of $k$
distinct points in an open interval.

The boundaries in these cell complexes can also be explicitly
described. Since $$[D]=\Delta_D\times E_D$$ is a product, its
boundary consists of two parts. The first part consists of the cells
that come from $\partial \Delta_D\times E_D$. These are of the form
$[D']$, where $D'$ is obtained form $D$ by removing a number of
chords. The second part comes from $\Delta_D\times \partial E_D$. 
The diagram $D'$ of a
cell $[D']$ of this kind is obtained from $D$ by
collapsing to zero the distance between two adjacent vertices. Note
that by removing chords or glueing together two adjacent vertices we
can decrease the complexity of a diagram; in this case the
corresponding part of $\partial [D]$ is glued to the base point in
$\onept{W_p}$.

If we are interested in the homology of $\onept{W_p}$, we need to
describe the boundaries in the corresponding cellular chain complex
and this involves only those cells $[D']\subset
\partial [D]$ with $\dim{[D']}=\dim{[D]}-1.$ For the dimension of
$[D]$ we have the formula $$\dim{[D]}=\text{no.\ of geometrically
distinct vertices}+\text{no.\ of chords} -1$$ and the complexity
$c(D)$ is given by the expression
$$c(D)=\text{total no.\ of vertices} - \text{no.\ of groups of self-intersection vertices}.$$
These formulae show that if by removing chords of $D$ we obtain a diagram $D'$ with
$c(D')=c(D)$ and $\dim{[D']}=\dim{[D]}-1$, then
$D'$ is obtained from $D$ by removing one chord, and
the endpoints of the removed chord belong to the same group of
self-intersection vertices in $D'$.

Now, suppose that by collapsing two adjacent vertices of $D$ we
obtain a diagram $D'$ with $c(D')=c(D)$ and
$\dim{[D']}=\dim{[D]}-1$. There are several possibilities for this:
\begin{enumerate}
\item both vertices are self-intersection vertices and belong to different
groups, and at most one of the two vertices is also a singularity
vertex:
$$\rb{-6pt}{\ig[height=1cm]{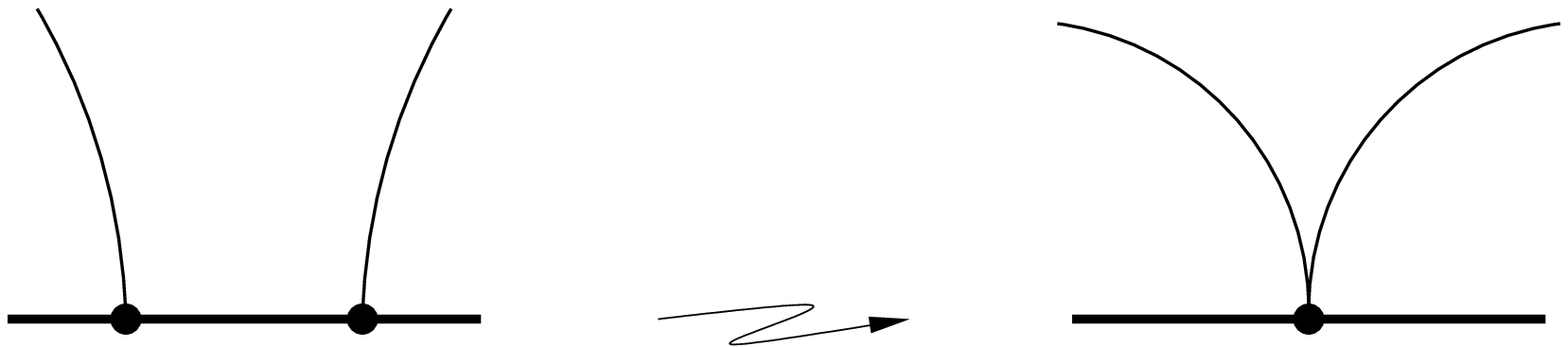}}
\quad\text{or}\quad\rb{-6pt}{\ig[height=1cm]{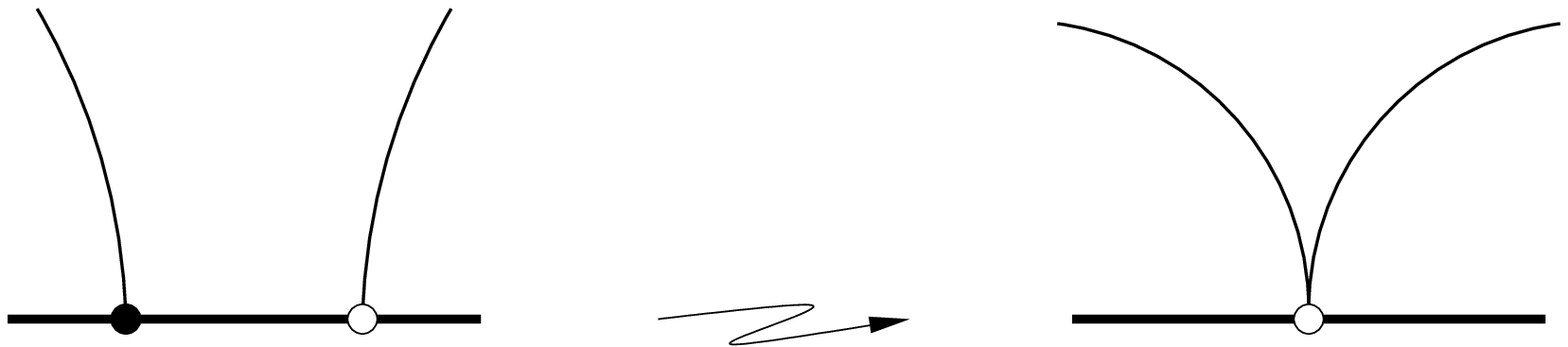}}$$
\item both vertices are endpoints of the same chord and of no other
chord:
$$\ig[height=0.6cm]{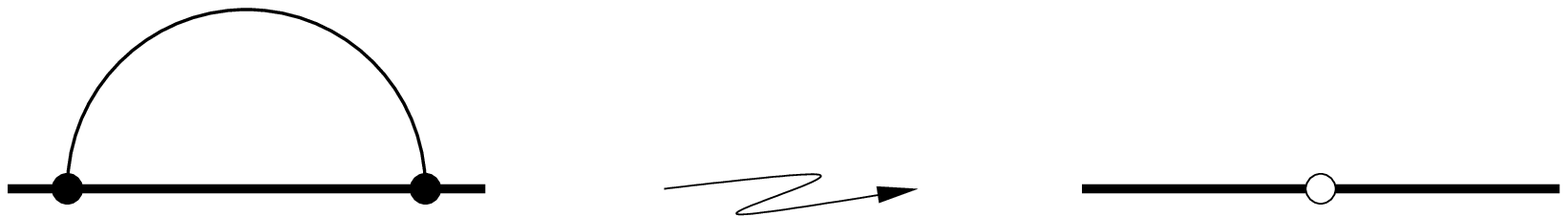}$$
\item one vertex is a self-intersection vertex only and the other is a singularity
vertex only:
$$\ig[height=1cm]{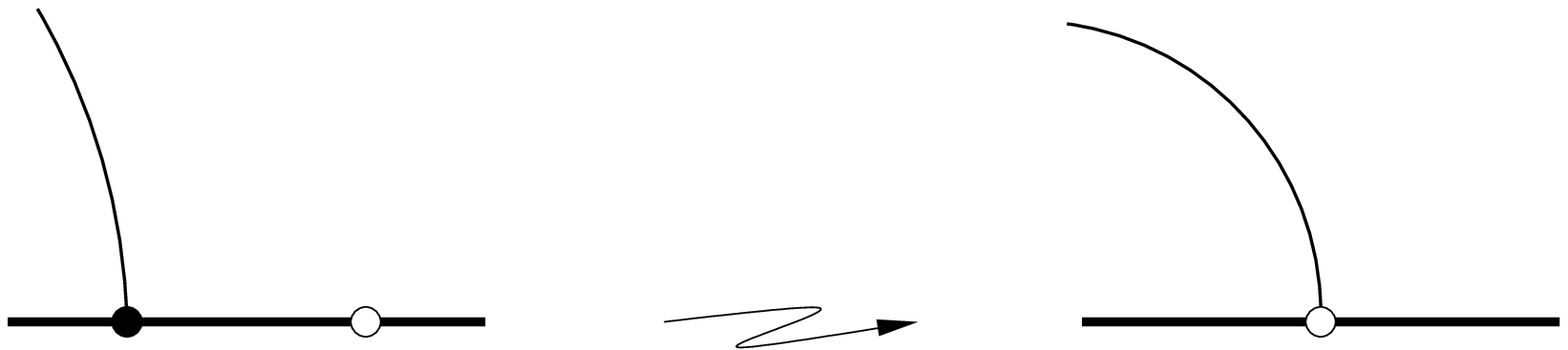}$$
\end{enumerate}
Each of the cells $[D']$ in the above list appears in the boundary
of $[D]$ exactly once. This describes the cellular chain complex for
$\onept{W_p}$ up to signs. This is sufficient if we work modulo
2. In order to calculate the integral homology, we need to fix the
orientations for each $[D]$ and work out the signs.

Recall that the cell $[D]$ is a product of two simplices; therefore,
its orientation can be specified by ordering the vertices of the
factors. $E_D$ is the configuration space of $k$ points in an
interval and its vertices are naturally ordered: the $i$th vertex is
a configuration with $k-i$ points in the left end of the interval
and $i$ points in the right end. (Note that the vertices belong to
closure of $E_D$, but not to $E_D$ itself, and the corresponding
configurations may have coinciding points.)
The vertices of $\Delta_D$ are the chords of $D$. In order to order
them, we first order the chords within each group of
self-intersection vertices: a chord is smaller than another chord if
its left endpoint if smaller; if both chords have the same left
endpoint, the one with the smaller right endpoint is smaller. Next,
we order the groups: a group is smaller if its leftmost vertex is
smaller. It is convenient to consider in this context each
singularity vertex which does not coincide with any
self-intersection vertex as a separate group consisting of one
``degenerate'' chord; the leftmost vertex of such a group is of
course, the singularity vertex itself. Finally, we list the chords
lexicographically: first, all the chords from the first group, then
the chords from the second group, and so on.

Now, it is clear how to assign the signs in the boundaries.
\begin{xexample}
In the cellular chain complex for $\onept{W_2}$ we have
$$d \left(\rb{-5pt}{\ig[width=2cm]{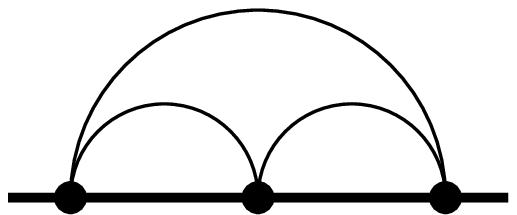}}\right)=
\rb{-5pt}{\ig[width=2cm]{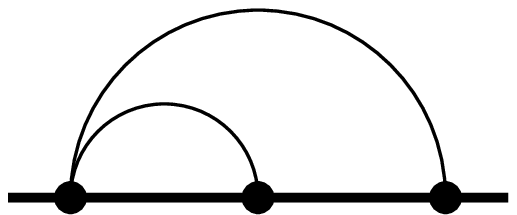}}-\rb{-5pt}{\ig[width=2cm]{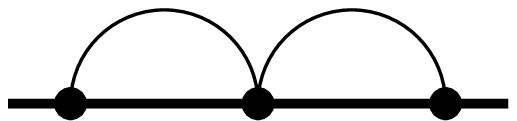}}+
\rb{-5pt}{\ig[width=2cm]{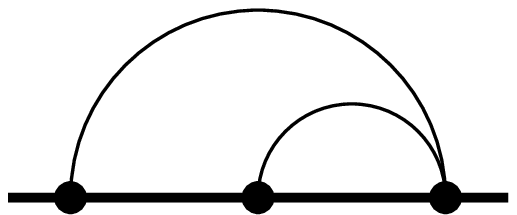}},$$
$$d\left(\rb{-5pt}{\ig[width=2cm]{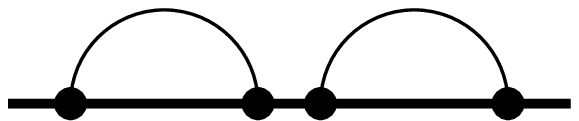}}\right)=
-\rb{-5pt}{\ig[width=2cm]{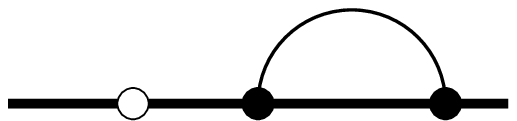}}+\rb{-5pt}{\ig[width=2cm]{abc3}}-
\rb{-5pt}{\ig[width=2cm]{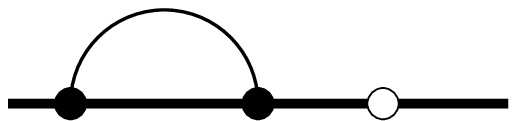}},$$ and
$$d\left(\rb{-5pt}{\ig[width=2cm]{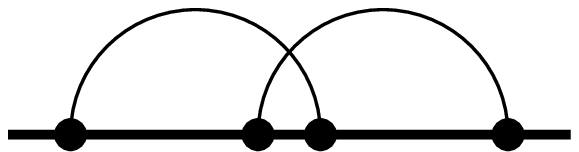}}\right)=
-\rb{-5pt}{\ig[width=2cm]{abc2}}+\rb{-5pt}{\ig[width=2cm]{abc3}}-
\rb{-5pt}{\ig[width=2cm]{abc4}}.$$ Note that our convention for the
orientations of the cells is different from that of \cite{Va3, Va7}.
\end{xexample}
\begin{xxca}
Formulate the general rule for the signs in the boundary of a cell
$[D]$ for an arbitrary $D$.
\end{xxca}

\subsection{The filtration on $\onept{W_p}$ by the number of vertices}

In principle, the homology of $\onept{W_p}$ can be calculated
directly from the cellular chain complex that we have just
described. There is, however, a better way to calculate this
homology.

Let $\onept{W_p}(k)$ be the subspace of $\onept{W_p}$ consisting of
all the cells $[D]$ where $D$ has at most $k$ geometrically distinct
vertices, together with the added basepoint. The smallest $k$ for
which $\onept{W_p}(k)$ is non-trivial is equal to $[p/2]+1$ where
$[\cdot]$ denotes the integer part; this corresponds to the diagrams
all of whose singularity vertices are combined with self-intersection
vertices and the latter are joined into only one group. The maximal
number of distinct vertices is $k=2p$; it is achieved for chord
diagrams of degree $p$. We get the increasing filtration
$$*=\onept{W_p}([p/2])\subset\onept{W_p}([p/2]+1)\subset\ldots\subset\onept{W_p}(2p)=\onept{W_p}$$
by the number of vertices; in \cite{Va3, Va7} it is called the {\em
auxiliary filtration}.\index{Auxiliary filtration}\label{aux-filtr}

The successive quotients in this filtration are bouquets of certain
spaces indexed by equivalence classes of degenerate diagrams:
$$\onept{W_p}(k)/\onept{W_p}(k-1)=\bigvee_{\mathfrak{D}} [\mathfrak{D}],$$
where $\mathfrak{D}$ runs over all equivalence classes of diagrams
with $k$ distinct vertices and $[\mathfrak{D}]$ is the union of all
the cells $[D]$ such that the equivalence class of $D$ is
$\mathfrak{D}$ (the basepoint counted as one of such cells). In
turn, each $[\mathfrak{D}]$ can be constructed out of several
standard pieces.

For a positive integer $a$ define the {\em complex of connected
graphs $\Delta^1(a)$} as follows. Given a set $A$ of $a$ points,
consider the simplex of dimension  $a(a-1)/2-1$ whose vertices are
indexed by chords connecting pairs of points of $A$. Each face of
this simplex corresponds to a graph whose set of vertices is $A$.
Collapsing the union of all those faces that correspond to
non-connected graphs to a point, we get the complex of connected
graphs $\Delta^1(a)$.

The proof of the following statement can be found in \cite{Va3,
Va7}:
\begin{xlemma}
$H_i(\Delta^1(a))=0$ unless $i=0,a-2$, and
$H_{a-2}(\Delta^1(a))=\Z^{(a-1)!}$.
\end{xlemma}

The spaces $[\mathfrak{D}]$ can now be described as follows:
\begin{xlemma}
Let $\mathfrak{D}$ be an equivalence class of diagrams with $m$
groups of self-intersection vertices consisting of $a_1,\ldots, a_m$
vertices respectively, $b$ singularity vertices, and $k$ geometrically
distinct vertices in total. Then
$$[\mathfrak{D}]=\Delta^1(a_1)\wedge\ldots\wedge\Delta^1(a_m)\wedge
S^{k+m+b-1},
$$
where the wedge stands for the reduced product of topological spaces.
\end{xlemma}

\begin{proof}
It is instructive to verify first the case when the diagrams in $\mathfrak{D}$ have
$k$ self-intersection vertices only, all in one group. The space of all possible
sets of vertices for such diagrams is an open $k$-dimensional simplex. Over each
point of this space we have the $k(k+1)/2-1$-dimensional simplex without the
faces corresponding to non-connected graphs, that is, $\Delta^1(k)$ minus the basepoint.
$[\mathfrak{D}]$ is the one-point compactification of this product:
$$[\mathfrak{D}]=S^k\wedge \Delta^1(k)$$
and in this case $m=1$ and $k=a_1$.

In the general case the space of all possible sets of vertices of a diagram is still
a $k$-dimensional simplex. Now, over each set of vertices we have the interior of the join of
all $\Delta^1(a_i)$, taken without their basepoints, and $b$ singularity points. This is nothing
but the product
$$(\Delta^1(a_1)-*)\times\ldots\times(\Delta^1(a_m)-*)\times \R^{m+b-1}.$$
Taking one-point compactification we get the statement of the lemma.
\end{proof}

The two above lemmas imply that the homology of $[\mathfrak{D}]$ vanishes
in all dimensions apart from 0 and $p+k-1$. Now, using the homology exact
sequences, or, which is the same, the spectral sequence associated with the
filtration $\onept{W_p}(k)$, we arrive to the following
\begin{lemma}\label{lemma:threepee}
$$H_{3p-1}(\onept{W_p})=H_{3p-1}(\onept{W_p}/\onept{W_p}(2p-2)).$$
\end{lemma}

\subsection{Chord diagrams and 4T relations}
Cohomology classes of dimension zero, that is, knot invariants,
produced by the Vassiliev spectral sequence correspond to elements
of the groups $E^{-p,p}_{\infty}$ obtained from $E^{-p,p}_{1}$ as
quotients.

A consequence of Lemma~\ref{lemma:threepee} is the following description 
of the group $E^{-p,p}_{1}=H_{3p-1}(\onept{W_p})$:

\begin{xproposition}[\cite{Va3, Va7}]
For any ring $\Ring$ of coefficients, $E^{-p,p}_{1}$ is isomorphic
to the group $\W_p$ of $\Ring$-valued weight systems, that is,
functions on chord diagrams that vanish on the 1T and 4T relations.
\end{xproposition}

Note that, according to the Fundamental Theorem \ref{fund_thm}, over
the rational numbers the group $\W_p$ is isomorphic to the space
$E^{-p,p}_{\infty}$ of Vassiliev invariants of order $\leq p$, modulo
those of order $\leq p-1$.

\begin{proof}
The proof uses the cell decomposition of
$H_{3p-1}(\onept{W_p})$. The $3p-1$-dimensional cells that are not
contained in $\onept{W_p}(2p-2)$ are of the form $[D]$ where $D$ is
either
\begin{itemize}
\item a non-degenerate chord diagram of order $p$;
\item a degenerate chord diagram with $2p-1$ self-intersection vertices
and $p-1$ groups, of which $p-2$ are pairs and one is a group with 3
vertices connected by 3 chords.
\end{itemize}
None of these cells is contained in the boundary of a
$3p$-dimensional cell.

The $3p-2$-dimensional cells that are not contained in
$\onept{W_p}(2p-2)$ are of the form $[D]$ where $D$ is either
\begin{itemize}
\item[(1T)] a diagram with $2p-1$ distinct vertices one of which is a
singularity vertex and the rest
are self-intersection vertices grouped into pairs;
\item[(4T)] a degenerate chord diagram with $2p-1$ self-intersection vertices
and $p-1$ groups, of which $p-2$ are pairs and one is a group with 3
vertices connected by 2 chords.
\end{itemize}

The cellular chain complex consists of free modules, so the kernel
of the boundary on the $3p-1$-cells, is isomorphic to the dual of
the cokernel for the coboundary on the $3p-2$-cells. Unlike the
boundary, the coboundary is easy to calculate.

To be precise, let $d_{3p-1}:C_{3p-1}\to C_{3p-2}$ be the boundary
in the chain complex for $\onept{W_p}/\onept{W_p}(2p-2).$ The dual
modules $\Hom (C_{i}, \Ring)$ can be identified with $C_{i}$; in
particular, they are generated by the same degenerate chord
diagrams. The dual homomorphism $d_{3p-1}^*:C_{3p-2}\to C_{3p-1}$
sends a diagram of the type (1T) to a diagram which has a chord with
adjacent vertices; moreover, every diagram with such a chord is in
the image of $d_{3p-1}^*$.

At this point it will be convenient to modify our convention on the
orientation of the cells. Let us change the orientation of the cells
that correspond to non-degenerate chord diagrams by multiplying it
by $(-1)^r$, where $r$ is the number of intersections among the
chords of the diagram, or, in a slightly fancier language, the
number of edges of its intersection graph. Diagrams of the type (4T)
are sent by $d_{3p-1}^*$ to linear combinations of 3 diagrams:
$$d_{3p-1}^*\left(\rb{-5pt}{\ig[width=2cm]{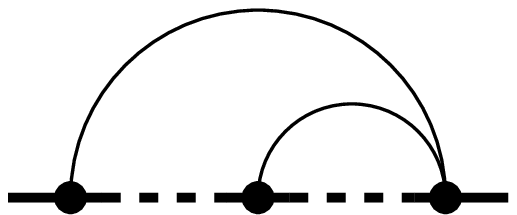}}\right)=
\
\rb{-5pt}{\ig[width=2cm]{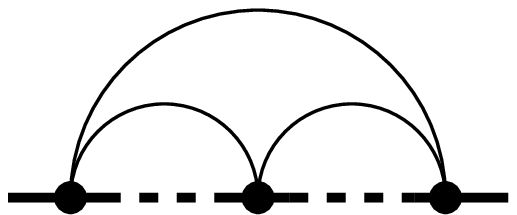}}+(-1)^s\rb{-5pt}{\ig[width=2cm]{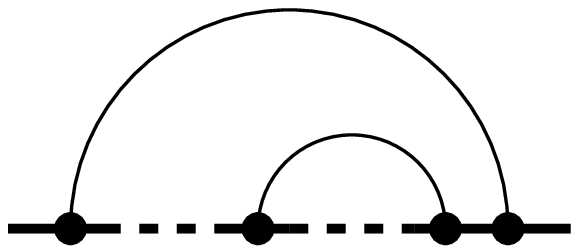}}-
(-1)^s\rb{-5pt}{\ig[width=2cm]{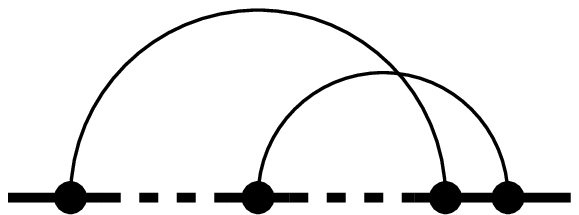}},$$
$$d_{3p-1}^*\left(\rb{-5pt}{\ig[width=2cm]{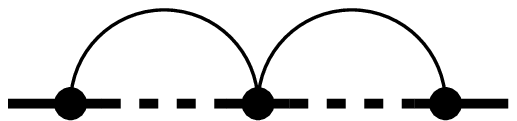}}\right)=
-\rb{-5pt}{\ig[width=2cm]{abc9}}-(-1)^s\rb{-5pt}{\ig[width=2cm]{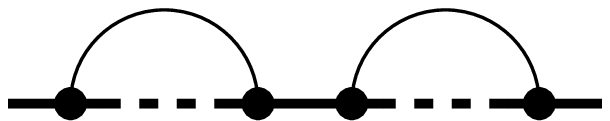}}+
(-1)^s\rb{-5pt}{\ig[width=2cm]{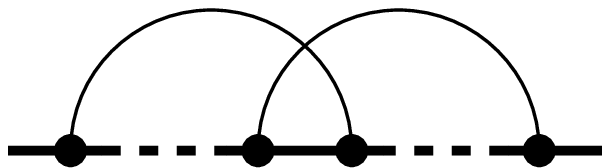}},$$ and
$$d_{3p-1}^*\left(\rb{-5pt}{\ig[width=2cm]{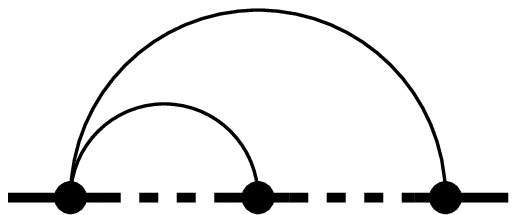}}\right)=
\
\rb{-5pt}{\ig[width=2cm]{abc9}}+(-1)^s\rb{-5pt}{\ig[width=2cm]{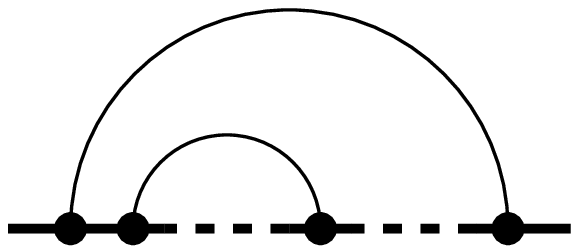}}-(-1)^s
\rb{-5pt}{\ig[width=2cm]{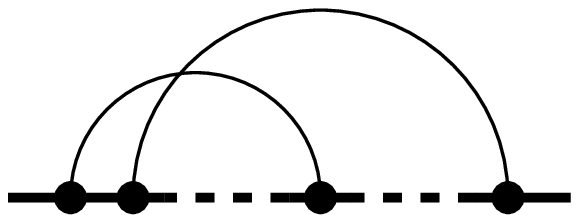}},$$ where $s$ is the same number in
all cases.
\begin{xxca}
Find an expression for $s$ and verify the above formulae.
\end{xxca}
It follows that each 4T relation is in the image of $d_{3p-1}^*$.

Now, let $C_{3p-1}'$ be the subspace of $C_{3p-1}$ spanned by
non-degenerate chord diagrams. Each functional $f$ on
$C_{3p-1}/{d_{3p-1}^*C_{3p-2}}$ can be uniquely reconstructed from
its value on $C_{3p-1}'/\langle {\text{1T, 4T}}\rangle$ since
$$f\left(\rb{-5pt}{\ig[width=2cm]{abc9}}\right)=
-(-1)^s f\left(\rb{-5pt}{\ig[width=2cm]{abc10}}\right)+ (-1)^s
f\left(\rb{-5pt}{\ig[width=2cm]{abc11}}\right)$$ where $s$ is as
above, and, hence, we see that
$H_{3p-1}(\onept{W_p}/\onept{W_p}(2p-2))$ consists precisely of the
$\Ring$-valued weight systems.
\end{proof}

The reader who has survived to this point may note that in
Vassiliev's original approach the road to the combinatorial
description of the weight systems has been long and winding,
especially if compared to the method presented in first chapters 
of this book. We stress, however, that while Vassiliev's approach to the
0-dimensional cohomology classes can be dramatically simplified,
there are no low-tech solutions for classes of higher dimensions.

\section{Homology of the space of knots and Poisson algebras}

The same methods that we have used in this chapter to study the
cohomology of the space of knots can be employed in order to attempt
to describe its homology. In particular, we get a homological
spectral sequence whose first term consists of the cohomology groups
of the diagram complexes $\onept{W}_p$ and can be described
completely in terms of degenerate chord diagrams. The bialgebra of
chord diagrams $\A$ forms a part of this spectral sequence; namely
$\A_p$ is isomorphic over a ring $\Ring$ to the diagonal entry
$E^{1}_{-p,p}=\widetilde{H}^{3p-1}(\onept{W_p},\Ring)$.

It is very interesting to note that the first term of this spectral
sequence has another interpretation, which, at first, seems to be
completely unrelated to knots. Namely, as discovered by V.~Turchin
\cite{Tu1} it is closely related to the {\em Hochshild homology of
the Poisson algebras operad}.

The details of Turchin's work are outside the scope of this book.
Let us just give a rough explanation of how Poisson algebras appear
in the homological Vassiliev spectral sequence.

Recall that a {\em Poisson algebra}\index{Poisson algebra} has two
bilinear operations: a commutative and associative product, and an
antisymmetric bracket satisfying the Jacobi identity. The two
operations are related by the Leibniz rule
$$[ab,c]=a[b,c]+b[a,c].$$
Using the Leibniz rule, one can re-write any composition of products
and brackets as a linear combination of products of iterated
brackets, which we call {\em Poisson monomials}.

In order to describe the cohomology of the diagram complex
$\onept{W_p}$ one can use the auxiliary filtration on it by the
number of vertices of a diagram, see page~\pageref{aux-filtr}. The
successive quotients in this filtration are built out of certain
standard pieces indexed by equivalence classes of degenerate
diagrams. As pointed out in in \cite{Tu1}, these equivalence classes
give rise to Poisson monomials in the following fashion.

Let us restrict our attention to degenerate chord diagrams without
singularity vertices. Label the vertices of such a diagram by
numbers from 1 to $n$ according to their natural order on the real
line. The equivalence class of the diagram is then determined by a
partition of the set $1,\ldots,n$ into several subsets with at least
two elements each. For every such subset $i_1,\ldots, i_k$ form an
iterated bracket $[\ldots[[x_{i_1},x_{i_2}]\ldots, x_{i_k}]$ and
take the product of these brackets over all the subsets. For
instance, the equivalence class of the diagram
$$\ig[height=1.5cm]{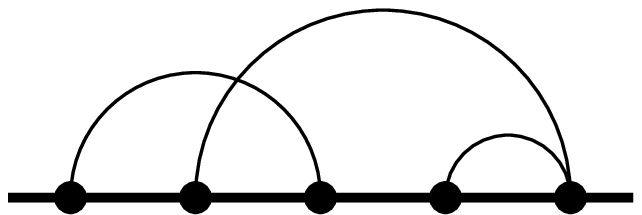}$$
gives rise
to the monomial $[x_1,x_3][[x_2,x_4],x_5]$. A chord diagram of
degree $p$ gives a product of $p$ simple brackets.

Poisson monomials of this type appear in the Hochschild complex for
the Poisson algebra operad. V.~Turchin proves the following result:
\begin{theorem}
The first term of the Vassiliev spectral sequence for the homology
of the space of knots coincides with the Hochschild homology
bialgebra for the operad of Poisson algebras with the Poisson
bracket of degree $-3$, taken modulo two explicit relations, and
with the inverse grading.
\end{theorem}
We refer to \cite{Tu1} for the basics on operads and their
Hochschild homology, and the precise form of this statement. A further reference is the 
paper \cite{Sin1} by D.~Sinha where the relationship between the Vassiliev spectral 
sequence and the Hochschild complex is explained ``on the level of spaces''. 
 %15 All knots

\renewcommand{\thechapter}{A}
\setcounter{section}{0}

\chapter*{Appendix}

\section{Lie algebras and their representations}

\subsection{Lie algebras}
\label{metrized}

A Lie algebra $\g$ over a field $\Fi$ of characteristic zero is a vector space equipped with a bilinear operation ({\em Lie bracket}) $(x,y) \mapsto [x,y]$ subject to the
identities
\begin{gather*}
   [x,y]=-[y,x]\,,\\
   [x,[y,z]]+[y,[z,x]]+[z,[x,y]]=0.
\end{gather*}
In this section we shall only consider finite-dimensional Lie algebras over $\C$.

An {\em abelian} Lie algebra is a vector space with the bracket which is identically 0: $[x,y]=0$ for all $x,y\in \g$. Any vector space with the zero bracket is an abelian Lie algebra.

Considering the Lie bracket as a product, one may speak about homomorphisms of Lie algebras, Lie subalgebras, and so on. In particular, an {\em ideal} in a Lie algebra is a vector subspace stable under taking the bracket with an arbitrary element of the whole algebra. A
Lie algebra is called {\em simple} if it is not abelian and does not
contain any proper ideal. Simple Lie algebras are classified (see,
for example, \cite{FH,Hum}). Over the field of complex numbers $\C$
there are four families of {\em classical algebras}:\index{Lie
algebra!classical}
\begin{center}
\begin{tabular}{c|c|c|l}
Type & $\g$ & $\dim\g$ & description \\ \hline\hline
$A_n$ & $\sL_{n+1}$ &\scriptsize $n^2+2n$ & \scriptsize
      $(n+1)\times (n+1)$ matrices with zero trace, $(n\ge 1)$ \\ \hline
$B_n$ & $\so_{2n+1}$ &\scriptsize $2n^2+n$ & \scriptsize
      skew-symmetric $(2n+1)\times (2n+1)$ matrices, $(n\ge 2)$ \\ \hline
$C_n$ & $\sP_{2n}$ &\scriptsize $2n^2+n$ & \parbox{3in}{\scriptsize
      $2n\times 2n$ matrices $X$ satisfying the relation
      $X^t\cdot M + M\cdot X=0$, where $M$ is the standard
      $2n\times 2n$ skew-symmetric matrix
      $M=\left(\begin{array}{cc}O&Id_n\\-Id_n&0\end{array}\right)$ ,
      $(n\ge 3)$} \\ \hline
$D_n$ & $\so_{2n}$ &\scriptsize $2n^2-n$ & \scriptsize
      skew-symmetric $2n\times 2n$ matrices, $(n\ge 4)$
\end{tabular}
\end{center}
\index{$\sL_N$}\index{$\so_N$}\index{$\sP_N$} and five {\em
exceptional algebras}:\index{Lie algebra!exceptional}
$$\begin{array}{c||@{\quad}c@{\quad}|@{\quad}c@{\quad}|@{\quad}c@{\quad}|
                   @{\quad}c@{\quad}|@{\quad}c@{\quad}}
\mbox{Type}& E_6 & E_7 & E_8 & F_4 & G_2 \\ \hline
\dim\g & 78 & 133 & 248 & 52 & 14
\end{array}
$$\index{$E_6$}\index{$E_7$}\index{$E_8$}\index{$F_4$}\index{$G_2$}
The Lie bracket in the matrix Lie algebras above is the commutator
of matrices.

Apart from the low-dimensional isomorphisms
$$\sP_2\cong\so_3\cong\sL_2;\qquad
  \sP_4\cong\so_5;\qquad
  \so_4\cong\sL_2\op \sL_2;\qquad
  \so_6\cong\sL_4,
$$
all the Lie algebras in the list above are different. The Lie
algebra $\gl_N$ of all $N\times N$ matrices is isomorphic to the
direct sum of $\sL_N$ and the abelian one-dimensional Lie algebra
$\C$.

\subsection{Metrized Lie algebras}
For $x\in\g$ write $\ad_x$ for the linear map $\g\to\g$ given by
$\ad_x(y) = [x,y]$.

The {\em Killing form} \index{Killing form} on a Lie algebra $\g$ is
defined by the equality
$$
  \langle x,y \rangle^K = \mbox{Tr}( \mbox{ad}_x \mbox{ad}_y ).
$$
{\em Cartan's criterion} says that this bilinear form is
non-degenerate if and only if the algebra is {\em semi-simple}, that
is, isomorphic to a direct sum of simple Lie algebras.
\begin{xxca}
Prove that the Killing form is $\ad$-invariant in the sense of the
following definition.
\end{xxca}

\begin{xdefinition}
A bilinear form $\langle\cdot,\cdot\rangle:\g\ot\g\to\C$ is
said to be {\em $\ad$-invariant}\index{ad-invariant!bilinear form}
if it satisfies the identity
$$
  \langle\ad_z(x),y\rangle + \langle x,\ad_z(y)\rangle =0,
$$
or, equivalently,
\begin{equation}
  \langle[x,z],y\rangle = \langle x,[z,y]\rangle. \label{ad-inv}
\end{equation}
for all $x,y,z\in\g$.
\end{xdefinition}

This definition is justified by the fact described in the following
exercise.

\begin{xxca}
Let $\mathfrak{G}$ be the connected Lie group corresponding to the
Lie algebra $\g$ and let $\Ad_g: \g\to\g$ be its adjoint
representation (see, for instance, \cite{AdJ}). Then the
$\ad$-invariance of a bilinear form is equivalent to its
$\Ad$-invariance defined by the natural rule
$$
  \langle\Ad_g(x),\Ad_g(y)\rangle = \langle x,y\rangle
$$
for all $x,y\in\g$ and $g\in\mathfrak{G}$.
\end{xxca}

A Lie algebra is said to be {\em metrized}, \index{Lie
algebra!metrized} if it is equipped with an ad-invariant symmetric
non-degene\-rate bilinear form
$\langle\cdot,\cdot\rangle:\g\ot\g\to\C$. The class of metrized
algebras contains simple Lie algebras with (a multiple of) the
Killing form, abelian Lie algebras with an arbitrary non-degenerate
bilinear form and their direct sums. If a Lie algebra is simple, all
ad-invariant symmetric non-degene\-rate bilinear forms are multiples
of each other.

For the classical simple Lie algebras which consist of matrices it
is often more convenient to use, instead of the Killing form, a
different bilinear form $\langle x,y \rangle = \mbox{Tr}(xy)$, which
is proportional to the Killing form with the coefficient
$\frac{1}{2N}$ for $\sL_N$, $\frac{1}{N-2}$ for $\so_N$, and
$\frac{1}{N+2}$ for $\sP_N$.

\begin{xxca}
Prove that for the Lie algebra $\gl_N$ the Killing form
$\langle{x,y}\rangle=\Tr(\ad_x\cdot\ad_y)$ is degenerate with defect
1 and can be expressed as follows:
$$
  \langle x,y \rangle^K = 2N\mbox{Tr}(xy) - 2\mbox{Tr}(x)\mbox{Tr}(y)\ .
$$
\end{xxca}

\begin{xxca}
Prove that the form $\mbox{Tr}(xy)$ on $\gl_N$ is non-degenerate and
$\ad$-invariant.
\end{xxca}

The bilinear form $\langle\cdot,\cdot\rangle$ is an element of
$\g^*\ot\g^*$. It identifies $\g$ with $\g^*$ and, hence, can be
considered as an element $c\in\g\ot\g$, called the {\em quadratic
Casimir tensor}.\index{Casimir!tensor}\label{casimir_ten} If
$\{e_i\}$ is a basis for $\g$ and $\{e_i^*\}$ is the dual basis, the
Casimir tensor can be written as $$c=\sum_i e_i\ot e_i^*.$$ The
quadratic Casimir tensor is $\ad$-invariant in the sense that
$$\ad_x c:=\sum_i \ad_x e_i\ot e_i^*+ \sum_i  e_i\ot \ad_x e_i^*=0.$$

\subsection{Structure constants}
Given a basis $\{e_i\}$ for the Lie algebra $\g$ the Lie brackets of
the basis elements can be written as $$[e_i,e_j]=\sum_{k=1}^d
c_{ijk}e_k.$$ The numbers $c_{ijk}$ are called the {\em structure
constants of $\g$ with respect to $\{e_i\}$.}
\begin{xlemma}\label{skewsymm}
Let $c_{ijk}$ be the structure constants of a metrized Lie algebra
in a basis $\{e_i\}$, orthonormal with respect to an ad-invariant
bilinear form. Then the constants $c_{ijk}$ are antisymmetric with respect to the
permutations of the indices $i$, $j$ and $k$.
\end{xlemma}

\begin{proof}
The equality $c_{ijk}=-c_{jik}$ is the coordinate expression of the fact that the
commutator is antisymmetric: $[x,y]=-[y,x]$. It remains to prove that
$c_{ijk}=c_{jki}$. This follows immediately from equation (\ref{ad-inv}),
simply by setting $x=e_i$, $y=e_k$, $z=e_j$.
\end{proof}

The Lie bracket, being a bilinear map $\g\ot \g\to\g$, can be
considered as an element of $\g^*\ot\g^*\ot\g$. The metric defines
an isomorphism $\g\simeq\g^*$ and, hence, the Lie bracket of a
metrized Lie algebra produces an element in $J\in\g^{\ot3}$ called
the {\em structure tensor} of $\g$.

\begin{xcorollary}
The structure tensor $J$ of a metrized Lie algebra $\g$ is totally
antisymmetric: $J\in\wedge^3\g$.
\end{xcorollary}

\subsection{Representations of Lie algebras}

A {\em representation} of a Lie algebra $\g$ in a vector space $V$
is a Lie algebra homomorphism of $\g$ into the Lie algebra $\gl(V)$
of linear operators in $V$, that is, a map $\rho:\g\to\gl(V)$ such
that
$$
  \rho([x,y]) = \rho(x)\rho(y)-\rho(y)\rho(x).
$$
It is also said that $V$ is a {\em $\g$-module} and that $\g$ {\em
acts on $V$} by $\rho$. When $\rho$ is understood from the context,
the element $\rho(x)(v)$ can be written as $x(v)$. The {\em
invariants} of an action of $\g$ on $V$ are the elements of $V$ that
lie in the kernel of $\rho(x)$ for all $x\in\g$. The space of all
invariants in $V$ is denoted by $V^{\g}$.

The {\em standard representation} \label{stand_repr} $St$ of a
matrix Lie algebra, such as $\gl_N$ or $\sL_N$, is the
representation in $\C^N$ given by the identity map.

The {\em adjoint representation} \index{Adjoint representation} is
the action $\ad$ of $\g$ on itself according to the rule
$$
  x \mapsto \ad_x \in \Hom(\g,\g)\,,\quad
  \ad_x(y) = [x,y]\,.
$$
It is indeed a representation, since
$\ad_{[x,y]}=\ad_x\cdot\ad_y-\ad_y\cdot\ad_x=[\ad_x,\ad_y]$.

A representation $\rho:\g\to\gl(V)$ is {\em reducible} if there
exist  $\rho_1:\g\to\gl(V_1)$ and  $\rho_2:\g\to\gl(V_2)$ with
$V_i\neq 0$ and $V=V_1\oplus V_2$, such that $\rho=\rho_1\oplus
\rho_2$. A representation that is not reducible is {\em
irreducible}.

\begin{xexample}\label{irr_sl2}
The algebra $\sL_2$ of $2\times 2$-matrices with zero trace has
precisely one irreducible representation of dimension $n+1$ for each
positive $n$. Denote this representation by $V_{n}$. There exist a
basis $e_{0},\ldots,e_n$ for $V_n$ in which the matrices
$$H=\left(\begin{array}{cc} 1&0\\ 0&-1\end{array}\right)\ ,\qquad
  E=\left(\begin{array}{cc} 0&1\\ 0&0\end{array}\right)\ ,\qquad
  F=\left(\begin{array}{cc} 0&0\\ 1&0\end{array}\right)
$$
which span $\sL_2$ act as follows: $$H(e_i)=(n-2i) e_i, \qquad
E(e_i)=(n-i+1)e_{i-1},\qquad  F(e_i)=(i+1)e_{i+1},$$ where it is
assumed that $e_{-1}=e_{n+1}=0$.
\end{xexample}

The {\em Casimir element}\index{Casimir!element}\label{casimir_el}
of a representation $\rho$ of a metrized Lie algebra is the matrix
$$c(\rho)=\sum_i\rho{e_i}\rho{e_i^*}$$
for some basis $\{e_i\}$ of $\g$. If $\rho$ is finite-dimensional,
the trace of the Casimir element of $\rho$ is well-defined; it is
called the {\em Casimir
number}\label{casimir_num}\index{Casimir!number} of $\rho$.
\begin{xxca}
Show that $c(\rho)$ is well-defined and commutes with the image of
$\rho$.
\end{xxca}

\subsection{Tensor algebras}

Let $\g$ be a vector space over a field of characteristic zero. A {\em tensor}\index{Tensor} is
an element of a tensor product of several copies of $\g$ and its
dual space $\g^*$. The number of factors in this product is called
the {\em rank} of the tensor. The canonical map $\g\ot\g^{*}\to\C$
induces maps
$$\g^{\ot p}\ot(\g^{*})^{\ot q}\to \g^{\ot p-1}\ot(\g^{*})^{\ot q-1}$$
called {\em contractions}\index{Contraction of a tensor}, defined
for any pair of factors $\g$ and $\g^{*}$ in the tensor product. %A

Denote by
$$T(\g)=\bigoplus_{n\ge 0} \g^{\ot n}$$ the
{\em tensor algebra}\index{Tensor algebra}\index{Algebra!tensor} of
the vector space $\g$, whose multiplication is given by the tensor
product. In particular, $\g\subset T(\g)$ is the subspace spanned by
the generators of $T(\g)$.

The {\em symmetric algebra}\index{Algebra!symmetric}
\index{Symmetric algebra} of $\g$, denoted by $S(\g)$, is the
quotient of $T(\g)$ by the two-sided ideal generated by all the
elements $x\otimes y-y\otimes x$. The symmetric algebra decomposes
as $$S(\g)=\bigoplus_{n\ge 0} S^n(\g),$$ where the {\em $n$th
symmetric power} $S^n(\g)$ is the images of $\g^{\ot n}$.

 Let $\{e_i\}$ be a
basis of $\g$. Then  $T(\g)$ can be identified with the free algebra
on the generators $e_i$, and $S(\g)$ with the free commutative
algebra on the $e_i$. In particular, elements of $S(e_i)$ can be
thought of as polynomials in the $e_i$ and products
$e_{j_1}e_{j_2}\ldots e_{j_m}$ such that $m\ge 0$ and $j_1\le
{j_2}\le\ldots\le{j_m}$ form an additive basis for $S(\g)$.

The symmetric algebra is a quotient, rather than a subalgebra, of
the tensor algebra. However, it can be identified with the subspace
of {\em symmetric tensors} in $T(\g)$.\index{Tensor!symmetric}\index{Symmetric tensor}\label{symtens} 
Namely, the image of the linear map $S^n(\g)\to \g^{\ot n}$ given by
$$e_{j_1}e_{j_2}\ldots e_{j_m}\to
\frac{1}{m!}\sum_{\sigma\in S_m} e_{\sigma(j_1)}\ot
e_{\sigma(j_2)}\ot \ldots \ot e_{\sigma(j_m)}$$ consists of the
tensors invariant under all permutations of the factors in $\g^{\ot
m}$.

\subsection{Universal enveloping and symmetric algebras}
\label{uea}

Any associative algebra can be considered as a Lie algebra whose Lie
bracket is the commutator
$$[a,b]=ab-ba.$$
While not every Lie algebra is of this form, each Lie algebra is
contained in an associative algebra as a subspace closed under the
commutator.

The {\em universal enveloping algebra}\index{Algebra!universal
enveloping} \index{Universal enveloping algebra} of $\g$, denoted by
$U(\g)$, is the quotient of $T(\g)$ by the two-sided ideal generated
by all the elements $$x\otimes y-y\otimes x-[x,y],$$ $x,y\in\g$. In
other words, we force the commutator of two elements of $\g\subset
T(\g)$ to be equal to their Lie bracket in $\g$. An example of a
universal enveloping algebra is the symmetric algebra $S(\g)$: one
can think of it as the universal enveloping algebra of the abelian
Lie algebra obtained from $\g$ by endowing it with the zero bracket.

The universal enveloping algebra of $\g$ is always
infinite-dimensional. A basis of $\g$ gives rise to an explicit
additive basis of $U(\g)$:
\begin{xtheorem}[Poincar\'e-Birkhoff-Witt]
Let $\{e_i\}$ be a basis of the Lie algebra $\g$. Then all the
products $e_{j_1}e_{j_2}\ldots e_{j_m}$ such that  $m\ge 0$ and
$j_1\le {j_2}\le\ldots\le{j_m}$ form an additive basis for $U(\g)$.
\end{xtheorem}\index{Poincar\'e-Birkhoff-Witt!
theorem}
\begin{xcorollary}
The map $\g\to T(\g)\to U(\g)$ is an inclusion; the restriction to
$\g$ of the commutator on $U(\g)$ coincides with the Lie bracket.
\end{xcorollary}
\begin{xxca}
Show that $U(\g)$ has the following universal property: for each
homomorphism $f$ of $\g$ into a commutator algebra of an associative
algebra $A$ there exists the unique homomorphism of associative
algebras $U(\g)\to A$ whose restriction to $\g$ is $f$.
\end{xxca}

The basis given by the Poincar\'e-Birkhoff-Witt Theorem does not
depend on the Lie bracket of $\g$. In particular, we see that
$$U(\g)\simeq S(\g)$$ as vector spaces.

Further, both $S(\g)$ and $U(\g)$ are $\g$-modules: the adjoint
representation of $\g$ can be extended to $S(\g)$ or $U(\g)$ by the
condition that $\g$ acts by {\em derivations}:
$$\ad_x(yz)=\ad_x(y)z+y\ad_x(z)$$
for all $y,z$ in $S(\g)$ or $U(\g)$. Note that in the case of
$U(\g)$ the element $x$ simply acts by taking the commutator with
$x$. In particular, the Casimir element for this action is simply
the image of the Casimir tensor under the map $T(\g)\to U(\g)$.

The {\em Poincar\'e-Birkhoff-Witt isomorphism}\label{PBW} is the map
$S(\g)\to U(\g)$ defined by
$$e_{j_1}e_{j_2}\ldots e_{j_m}\to
\frac{1}{m!}\sum_{\sigma\in S_m}
e_{\sigma(j_1)}e_{\sigma(j_2)}\ldots e_{\sigma(j_m)}.$$
\begin{xxca}
Show that this is indeed an isomorphism.
\end{xxca}
It follows from the definition that the Poincar\'e-Birkhoff-Witt
isomorphism is an {\em isomorphism of $\g$-modules}, that is, it
commutes with the action of $\g$. In fact, it is also an isomorphism
of coalgebras, see page~\pageref{pbw_coalgebras}. (Clearly, it is
not an algebra isomorphism, since $S(\g)$ is commutative and $U(\g)$
is not, unless $\g$ is abelian.)

\subsection{Duflo isomorphism}\label{duflok}

Since the universal enveloping algebra and the symmetric algebra of
a Lie algebra $\g$ are isomorphic as $\g$-modules, we have an
isomorphism of vector spaces
$$S(\g)^{\g}\simeq U(\g)^{\g}=Z(U(\g))$$
between the subalgebra of invariants in the symmetric algebra and
the centre of the universal enveloping algebra. This map does not
respect the product, but it turns out that $S(\g)^{\g}$ and
$Z(U(\g))$ are actually isomorphic as commutative algebras. The
isomorphism between them is given by the {\em Duflo-Kirillov map},
defined as follows.

A {\em differential operator} $S(\g)\to S(\g)$ is just an element of
the symmetric algebra $S(\g^*)$. The action of $S(\g^*)$ on $S(\g)$
is obtained by extending the pairing of $\g^*$ and $\g$: we
postulate that $$x(ab)=x(a)\cdot b+a \cdot x(b)$$ for any $x\in\g^*$
and $a,b\in S(\g)$, and that $$(xy)(a)= x(y(a))$$ for $x,y\in
S(\g^*)$ and $a\in S(\g)$. An element of $S^k(\g^*)$ is a
differential operator of order $k$: it sends $S^m(\g)$ to
$S^{m-k}(\g)$. We can also speak of differential operators of
infinite order; these are elements of the graded completion of
$S(\g^*)$.

If $\g$ is a metrized Lie algebra, its bilinear form gives an
isomorphism between $S(\g)$ and $S(\g^*)$, which sends elements of
$S(\g)$ to differential operators. Explicitly, if we think of
elements of $S(\g)$ as symmetric tensors, for $a\in S(\g)$ the
operator $\partial_a\in S(\g^*)$ is obtained by taking the sum of
all possible contractions with $a$.

Let $j(x)$ be a formal power series with $x\in\g$ given by
$$j(x)=\det \left( \frac{\sinh{\frac{1}{2}{\ad{\,x}}}}{{\frac{1}{2}\ad{\,x}}}\right).$$
This power series starts with the identity, so we can take the
square root $\sqrt{j}$. It is an element of the graded completion of
$S(\g^*)$ so it can be considered as a differential operator of
infinite order, called the Duflo-Kirillov map:
$${\sqrt{j}}:S(\g)\to S(\g).$$
\begin{xtheorem}[\cite{Duf}, see also  \cite{AT,BLT,Kon3}]
The composition of the Duflo-Kirillov map with the
Poincar\'e-Birkhoff-Witt isomorphism is an isomorphism of
commutative algebras $S(\g)^{\g}\to Z(U(\g))$.
\end{xtheorem}
This isomorphism is known as the\index{Duflo isomorphism}
\index{Isomorphism!Duflo} {\em Duflo isomorphism}.

\subsection{Lie superalgebras}\label{lie_superalgebras}
A {\em super vector space}, or a {\em $\Z_2$-graded vector space} is
a vector space decomposed as a direct sum $$V=V_0\oplus V_1.$$ The
indices (or {\em degrees}) 0 and 1 are thought of as elements of
$\Z_2$; $V_0$ is called the {\em even} part of $V$ and $V_1$ is the
{\em odd} part of $V$. An element $x\in V$ is {\em homogeneous} if
it belongs to either $V_0$ or $V_1$. For $x$ homogeneous we write
$|x|$ for the degree of $x$. The (super) dimension of $V$ is the
pair $(\dim V_0\, |\, \dim V_1)$ also written as $\dim V_0+ \dim
V_1$.

The space $\gl(V)$ of all endomorphisms of a super vector space $V$
is a super vector space itself: $\gl(V)_i$ consists of maps $f:V\to
V$ such that $f(V_j)\subseteq V_{j+i}$; each $f\in \gl(V)$ can be
written as a sum $f_0+f_1$ with $f_i\in \gl(V)_i$. If $V$ is
finite-dimensional the {\em supertrace} of $f$ is defined as
$${\rm sTr} f= \Tr f_0 - \Tr f_1.$$

A {\em superalgebra} is a super vector space $A$ together with a
bilinear product which respects the degree:
$$|xy|=|x|+|y|$$ for all homogeneous $x$ and $y$ in $A$. The {\em
supercommutator} in a superalgebra $A$ is a bilinear operation
defined on homogeneous $x,y\in A$ by
$$[x,y]=xy-(-1)^{|x|\,|y|}yx.$$ The elements of $A$ whose
supercommutator with the whole of $A$ is zero form the {\em super
center} of $A$.

The supercommutator satisfies the following identities:
$$|[x,y]|=|x|+|y|,$$
$$[x,y]=-(-1)^{|x|\,|y|}[y,x]$$
and
$$(-1)^{|z|\,|x|}[x,[y,z]]+(-1)^{|y|\,|z|}[z,[x,y]]+(-1)^{|x|\,|y|}[y,[z,x]]=0,$$
where $x,y,z$ are homogeneous. A super vector space with a bilinear
bracket satisfying these identities is called a {\em Lie
superalgebra}.\index{Lie superalgebra}

Each Lie superalgebra $\g$ can be thought of as a subspace of its
{\em universal enveloping superalgebra} $U(\g)$ defined as the
quotient of the tensor algebra on $\g$ by the ideal generated by
$$x\otimes y- (-1)^{|x|\, |y|} y\otimes x-[x,y],$$
where $x$ and $y$ are arbitrary homogeneous elements of $\g$; the
supercommutator in $U(\g)$ induces the bracket of $\g$.

The theory of Lie superalgebras was developed by V.~Kac \cite{Kac1,
Kac2}; it closely parallels the usual Lie theory.

\begin{xexample}[\cite{FKV}] The Lie superalgebra $\gL_1_1$ consists of the
endomorphisms of the super vector space of dimension 1+1 with the
bracket being the supercommutator of endomorphisms. The supertrace
gives a bilinear form on $\gL_1_1$
$$\langle x, y \rangle= {\rm sTr}(xy),$$
which is non-degenerate and ad-invariant in the same sense as for
Lie algebras:
$$
  \langle[x,z],y\rangle = \langle x,[z,y]\rangle.
$$
Take a basis in the 1+1 - dimensional space whose first vector is
even and the second vector is odd. Then the even part of $\gL_1_1$
is spanned by the matrices
$$H=\left(\begin{array}{cc} 1&0\\ 0&1\end{array}\right)\ ,\qquad
  G=\left(\begin{array}{cc} 0&0\\ 0&1\end{array}\right)\ ,
$$
and the odd part by
$$Q_+=\left(\begin{array}{cc} 0&0\\ 1&0\end{array}\right)\ ,\qquad
  Q_-=\left(\begin{array}{cc} 0&1\\ 0&0\end{array}\right)\ .
$$
The Lie bracket of $H$ with any element vanishes and we have
$$[G,Q_{\pm}]=\pm Q_{\pm}\qquad \text{and}\qquad [Q_+,Q_-]=H.$$
The quadratic Casimir tensor for $\gL_1_1$ is
$$H\ot G+ G\ot H - Q_+Q_-+Q_-Q_+.$$
Its image $c$ in the universal enveloping algebra $U(\gL_1_1)$
together with the image of $H$ under the inclusion of $\gL_1_1$,
which we denote by $h$, generate a polynomial subalgebra of
$U(\gL_1_1)$ which coincides with the super center of $U(\gL_1_1)$.
\end{xexample}

\section{Bialgebras and Hopf algebras}
\label{bialg}

Here we give a brief summary of necessary information about
bialgebras and Hopf algebras. More details can be found in
\cite{Abe, Car3, MiMo}.

\subsection{Coalgebras and bialgebras}

In what follows, all vector spaces and algebras will be considered
over a field $\Fi$ of characteristic zero. First, let us recall the
definition of an algebra in the language of commutative diagrams.

\begin{xdefinition}\index{Algebra}\index{Product}\index{$\mu$}\label{Algebra}
A {\em product}, or a {\em multiplication}, on a vector space $A$ is
a linear map $\mu:A\ot A\to A$. The product $\mu$ on $A$ is {\em
associative} if the diagram
\begin{equation*}\begin{CD}
A\ot A\ot A    @>\mu\ot\id    >>  A\ot A \\
@V{\id\ot\mu}VV             @VV{\mu}V \\
A\ot A    @>>\mu>           A
\end{CD}\end{equation*}
commutes. A {\em unit} for $\mu$ is a linear map $\io:\Fi\to A$
\index{$\io$} (uniquely defined by the element $\io(1)\in A$) that
makes commutative the diagram\index{Unit}
\begin{equation*}\begin{CD}
\Fi\ot A  @>\io\ot\id>>  A\ot A   \\
@AAA                  @VV{\mu}V \\
A        @=           A
\end{CD}\end{equation*}
where the upward arrow is the natural isomorphism. A vector space
with an associative product is called an (associative) {\em
algebra.}
\end{xdefinition}
The unit in an algebra, if exists, is always unique. We shall only
consider associative algebras with a unit.

Reversing the arrows in the above definition we arrive to the notion
of a {\em coalgebra}.

\begin{xdefinition}\index{Coalgebra}\index{Coproduct}\index{$\delta$}\index{$\e$}\label{Coalgebra}
A {coalgebra} is a vector space $A$ equipped with a linear map
$\delta:A\to A\ot A$, referred to as {\em comultiplication}, or {\em
coproduct}, and a linear map $\e:A\to\Fi$, called the {\em counit},
such that the following two diagrams commute:\index{Counit}
\begin{equation*}
\begin{CD}
A\ot A\ot A  @<\delta\ot\id<<  A\ot A \\
@A{\id\ot\delta}AA             @AA{\delta}A \\
A\ot A    @<\delta<<           A
\end{CD}
\hspace{2cm}
\begin{CD}
\Fi\ot A  @<\e\ot\id<<  A\ot A   \\
@VVV                  @AA{\delta}A \\
A        @=           A
\end{CD}
\end{equation*}
\end{xdefinition}

Algebras  (coalgebras) may possess an additional property of
commutativity \index{Commutativity} (respectively,
cocommutativity\index{Cocommutativity}), defined via the following
commutative diagrams:
\begin{equation*}
\begin{CD}
A\ot A  @>\mu>>   A \\
@A{\tau}AA          @| \\
A\ot A  @>\mu>>   A
\end{CD}
\hspace{3cm}
\begin{CD}
A\ot A  @<\delta<<   A \\
@V{\tau}VV         @| \\
A\ot A  @<\delta<<   A
\end{CD}
\end{equation*}
where $\tau:A\ot A \to A\ot A$ is the permutation of the tensor
factors: $$\tau(a\ot b)= b\ot a.$$

\begin{xdefinition}\index{Bialgebra}
A {\em bialgebra} is a vector space $A$ with the structure of an
algebra given by $\mu$, $\io$ and the structure of a coalgebra given
by $\delta$, $\e$ which agree in the sense that the following
identities hold:
\begin{enumerate}
\item $\e(1)=1$;
\item $\delta(1)=1\ot 1$;
\item $\e(ab)=\e(a)\e(b)$;
\item $\delta(ab)=\delta(a)\delta(b)$.
\end{enumerate}
Here $\mu$ is written as a usual product and in the last equation
$\delta(a)\delta(b)$ denotes the component-wise product in $A\ot A$
induced by the product $\mu$ in $A$.
\end{xdefinition}

Note that these conditions, taken in pairs, have the following meaning:

\begin{itemize}
\item (1,3) $\Leftrightarrow$ $\e$ is a  homomorphism of unital algebras.
\item (2,4) $\Leftrightarrow$ $\delta$ is a homomorphism of unital
algebras.
\item (1,2) $\Leftrightarrow$ $\io$ is a homomorphism of
coalgebras.
\item (3,4) $\Leftrightarrow$ $\mu$ is a homomorphism
of coalgebras.
\end{itemize}

The coherence of the two structures in the definition of a bialgebra
can thus be stated in either of the two equivalent ways:
\begin{itemize}
\item $\e$ and $\delta$ are algebra homomorphisms,
\item $\mu$ and $i$ are coalgebra homomorphisms.
\end{itemize}

\begin{xexample}The group algebra $\Fi G $ \index{Group algebra}
of a group $G$ over the field $\Fi$ consists of formal
linear combinations $\sum_{x\in G}\lambda_x x$ where
$\lambda_x\in\Fi$ with the product defined on the basis elements by
the group multiplication in $G$. The coproduct is defined as
$\delta(x)=x\ot x$ for $x\in G$ and then extended by linearity.
Instead of a group $G$, in this example one can, actually, take a
monoid, that is, a semigroup with a unit.
\end{xexample}
\begin{xexample}\label{bialg_ex}
The algebra $\Fi^G$ of $\Fi$-valued functions on a finite group $G$ with
pointwise multiplication
$$
  (fg)(x)=f(x)g(x)
$$
and the comultiplication defined by
$$
\delta(f)(x,y)=f(xy)
$$
where the element $\delta(f)\in\Fi^G\ot\Fi^G$ is understood as a
function on $G\times G$ via the natural isomorphism
$\Fi^G\ot\Fi^G\isom\Fi^{G\times G}$.
\end{xexample}
\begin{xexample} The symmetric algebra $S(V)$ of a vector space
$V$ is a bialgebra with the coproduct defined on the elements $x\in
V=S^1(V)$ by setting $\delta(x)=1\ot x+x\ot 1$ and then extended as
an algebra homomorphism to the entire $S(V)$.
\end{xexample}

\begin{xexample} The {\em completed} symmetric algebra
$$\widehat{S}(\g)=\prod_{n\ge 0} S^n(\g),$$
of a vector space $V$, whose elements are formal power series in the
coordinates of $V$, is a bialgebra whose coproduct extends that of
the symmetric algebra.
\end{xexample}

\begin{xexample}\index{Algebra!universal enveloping} Let $U(\g)$ be the
universal enveloping algebra of a  Lie algebra $\g$ (see
page~\pageref{uea}). Define $\delta(g)=1\ot g+g\ot 1$ for $g\in\g$
and extend it to all of $A$ by the axioms of bialgebra. If $\g$ is
abelian, this example reduces to that of the symmetric algebra.
\end{xexample}
\begin{xxca} Define the appropriate unit and counit in each of the above
examples.
\end{xxca}
\begin{xxca}\label{pbw_coalgebras}
Show that the Poincar\'e-Birkhoff-Witt isomorphism is an
isomorphism of coalgebras (that is, commutes with the counit and the
comultiplication).
\end{xxca}

\subsection{Primitive and group-like elements}\label{ap:pr-gr-el}

In bialgebras there are two remarkable classes of elements:
primitive elements and group-like elements.

\begin{xdefinition}\label{prim}\index{Primitive element}
An element $a\in A$ of a bialgebra $A$ is said to be {\em primitive}
if
$$
  \delta(a) = 1\ot a + a\ot 1.
$$
\end{xdefinition}

The set of all primitive elements forms a vector subspace $\PR(A)$
called the {\em primitive subspace} of the bialgebra $A$. The
primitive subspace is closed under the commutator $[a,b]=ab-ba$,
and, hence, forms a Lie algebra (which is abelian, if $A$ is
commutative). Indeed, since $\delta$ is a homomorphism, the fact
that $a$ and $b$ are primitive implies
\begin{gather*}
  \delta(ab) = \delta(a)\delta(b)=1\ot ab + a\ot b + b\ot a +ab\ot 1,\\
  \delta(ba) = \delta(b)\delta(a)=1\ot ba + b\ot a + a\ot b +ba\ot 1
\end{gather*}
and, therefore,
$$
  \delta([a,b]) = 1\ot [a,b] + [a,b]\ot 1.
$$

\begin{xdefinition}\label{gr_like}\index{Group-like element}
\index{Semigroup-like element} An element $a\in A$ is said to be
{\em semigroup-like} if
$$
   \d(a)=a\otimes a.
$$
If, in addition, $a$ is invertible, then it is called {\em
group-like}.
\end{xdefinition}

The set of all semigroup-like elements in a bialgebra is closed
under multiplication. It follows that the set of all group-like
elements $\GR(A)$ of a bialgebra $A$ forms a multiplicative
group.

Among the examples of bialgebras given above, the notions of the
primitive and group-like elements are especially transparent in the
case $A=\Fi^G$. As follows from the definitions, primitive elements
are the {\em additive} functions ($f(xy)=f(x)+f(y)$) while
group-like elements are the {\em multiplicative} functions
($f(xy)=f(x)f(y)$).

In the example of the symmetric algebra, there is an isomorphism
$$S(V)\ot S(V)\isom S(V\op V)$$ which allows to rewrite the definition of the
coproduct as $\d(x)=(x,x)\in V\op V$ for $x\in V$. It can be even
more suggestive to view the elements of the symmetric algebra $S(V)$
as polynomial functions on the dual space $V^*$ (where homogeneous
subspaces $S^0(V)$, $S^1(V)$, $S^2(V)$ and so on correspond to
constants, linear functions, quadratic functions et cetera on
$V^*$). In these terms, the product in $S(V)$ corresponds to the
usual (pointwise) multiplication of functions, while the coproduct
$\d:S(V)\to S(V\op V)$ acts according to the rule
$$
  \d(f)(\xi,\eta)=f(\xi+\eta),\quad \xi,\eta\in X^*.
$$
Under the same identifications,
$$
  (f\ot g)(\xi,\eta) = f(\xi)g(\eta),
$$
in particular,
\begin{gather*}
  (f\ot 1)(\xi,\eta) = f(\xi),\\
  (1\ot f)(\xi,\eta) = f(\eta).
\end{gather*}
We see that an element of $S(V)$, considered as a function on $V^*$,
is primitive (group-like) if and only if this function is additive
(multiplicative):
\begin{gather*}
  f(\xi,\eta) = f(\xi)+f(\eta),\\
  f(\xi,\eta) = f(\xi)f(\eta).
\end{gather*}
The first condition means that $f$ is a linear function on $V^*$,
that is, it corresponds to an element of $V$ itself; therefore,
$$
  \PR(S(V))=V.
$$
Over a field of characteristic zero, the second condition cannot hold
for polynomial functions except for the constant function equal to 1;
thus
$$
  \GR(S(V))=\{1\}.
$$
The {completed} symmetric algebra $\widehat{S}(V)$, in contrast with
$S(V)$, has a lot of group-like elements. Namely,
$$
  \GR(\widehat{S}(V))=\{\exp(x)\mid x\in V\},
$$
where $\exp(x)$ is defined as a formal power series
$1+x+x^2/2!+\ldots$, see page~\pageref{quillen}.

\begin{xxca}
Describe the primitive and group-like elements in $\Fi G $ and in
$U(\g)$.
\end{xxca}

{\sl Answer}: In $\Fi G$ $\PR=0$, $\GR=G$; in $U(\g)$ $\PR=\g$,
$\GR=\{1\}$.

\def\gra{{\text{gr}}}

\subsection{Filtrations and gradings}
\label{filtration} A decreasing filtration on a vector space $A$ is
a sequence of subspaces $A_i$, $i=0,1,2,...$ such that
$$
A=A_0 \supseteq A_1 \supseteq A_2 \supseteq \dots
$$
The {\em factors} of a decreasing filtration are the quotient spaces
$\gra_i A=A_i/A_{i+1}$.

An increasing filtration on a vector space $A$ is a sequence of
subspaces $A_i$, $i=0,1,2,...$ such that
$$
A_0 \subseteq A_1 \subseteq A_2 \subseteq \dots \subseteq A.
$$
The {factors} of an increasing filtration are the quotient spaces
$\gra_iA=A_i/A_{i-1}$, where by definition $A_{-1}=0$.

A filtration (either decreasing or increasing) is said to be of {\em
finite type} if all its factors are finite-dimensional. Note that in
each case the whole space has a (possibly infinite-dimensional)
``part'' not covered by the factors, namely $\cap_{i=1}^\infty A_i$
for a decreasing filtration and $A/\cup_{i=1}^\infty A_i$ for an
increasing filtration.

A vector space is said to be {\em graded} if it is represented as a
direct sum of its subspaces\index{Graded space}
$$
A=\bigoplus_{i=0}^\infty A_i.
$$
A graded space $A$ has a canonical increasing filtration by the subspaces 
$\oplus_{i=0}^k A_i$ and a canonical decreasing filtration by 
$\oplus_{i=k}^{\infty} A_i$.
 
With a filtered vector space $A$ one can associate a graded vector
space $G(A)$ setting\index{Graded space!associated with a
filtration}
$$
\gra A =\bigoplus_{i=0}^\infty G_iA = \bigoplus_{i=0}^\infty
A_i/A_{i+1}
$$
in case of a decreasing filtration and
$$
\gra A=\bigoplus_{i=0}^\infty G_iA = \bigoplus_{i=0}^\infty
A_i/A_{i-1}
$$
in case of an increasing filtration.

If $A$ is a filtered space of finite type, then the homogeneous
components $G_iA$ are also finite-dimensional; their dimensions have
a compact description in terms of the {\em Poincar\'e series}
$$\sum_{k=0}^\infty \dim(\gra_i A)\, t^k,$$ where $t$ is an auxiliary formal
variable.

\begin{xexample}
The Poincar\'e series of the algebra of polynomials in one variable
is
$$
1+t+t^2+... = \frac{1}{1-t}.
$$
\end{xexample}

\begin{xxca}
Find the Poincar\'e series of the polynomial algebra with $n$
independent variables.
\end{xxca}

One can also speak of filtered and graded algebras, coalgebras and bialgebras: these are filtered (graded) vector spaces with operations that respect the corresponding filtrations (gradings).

\begin{xdefinition}\label{d-filt-bialg}\label{i-filt-bialg}
\index{Bialgebra!filtered} \index{Algebra!filtered} \index{Coalgebra!filtered} 
We say that an algebra $A$ is {\em filtered} if its underlying vector space has a
filtration by subspaces $A_i$ compatible with the product in the sense that
$$A_p A_q \subset A_{p+q} \quad \text{for} \quad p,q\ge0 \quad\text{and}\quad 1\in A_0.$$

A coalgebra $A$ is filtered if it is filtered as a vector space and
$$ \delta(A_n) \subset \sum_{p+q=n} A_p \ot A_q \quad\text{for}\quad n\ge0. 
$$
Finally, a bialgebra is filtered if it is filtered both as an algebra and as a coalgebra, with respect to the same filtration. 
\end{xdefinition}

\begin{xdefinition}
\index{Bialgebra!graded} \index{Algebra!graded} \index{Coalgebra!graded} 
A {\em graded algebra} $A$ is a graded vector space with a product satisfying
$$A_p A_q \subset A_{p+q} \quad \text{for} \quad p,q\ge0.$$
A {\em graded coalgebra} $A$ is a graded vector space with a coproduct satisfying
$$ \delta(A_n) \subset \sum_{p+q=n} A_p \ot A_q \quad\text{for}\quad n\ge0 
\quad\text{and}\quad \e\vert_{A_k}=0 \quad\text{for}\quad k >0.
$$
A {\em graded bialgebra} is a graded vector space which is graded both as an algebra and as a coalgebra. 
\end{xdefinition}

The operations on filtered vector spaces descend to the associated graded spaces.
\begin{xproposition}
The graded vector space associated to a filtered algebra (coalgebra, bialgebra) has a natural structure of a graded algebra (respectively, coalgebra, bialgebra).
\end{xproposition}
The proof is by direct inspection.

\begin{xdefinition}\label{gr-compl}
The {\em graded completion}\index{Graded completion} of a graded
vector space $A=\oplus_{i=0}^{\infty} A_i$ is the vector space
$\widehat{A}=\prod_{i=0}^{\infty} A_i.$
\end{xdefinition}
For instance, the graded completion of the vector space of
polynomials in $n$ variables is the space of formal power series in
the same variables. Note that a priori there is no non-trivial
grading on the graded completion of a graded space.

Note that the product in a graded algebra extends uniquely to its graded completion; the same is true for the coproduct in a graded coalgebra.

\subsection{Dual filtered bialgebra}\index{Filtered bialgebra! dual}
\label{dual_fba} 

Let $A$ be a filtered bialgebra with a decreasing filtration $A_k$ of finite type. For each $k\geq 0$ define  $W_k$ to be the the subspace of $A^*$ consisting of all the linear functions on $A$ that vanish on $A_{k+1}$. Then $W_k$ is contained in $W_{k+1}$ and the union 
$$W=\bigcup_{k\geq 0} W_k$$
 is a filtered vector space with the increasing filtration by the $W_k$. 
\begin{proposition}\label{dual-bialgebra}
$W$ is a bialgebra with an increasing filtration, with the operations induced by duality by those of $A$.
\end{proposition}
We say that $W$ is a bialgebra {\em dual} to $A$. Note that $W\subseteq A^*$ and the equality holds if and only if $A$ is finite-dimensional.
\begin{proof}
If $\mu$ and $\delta$ are the product and the coproduct in $A$, respectively, with $\io$ the unit and $\e$ the counit, the operations in $W$ are as follows:
$$\begin{array}{ll}
\delta^*:   \sum_{k+l=n} W_k\ot W_l  \to W_{n} &
     \mbox{is the product in $W$},\vspace{8pt}\\
  \mu^*: W_n\to 
            \sum_{k+l=n} W_k\ot W_l 
   & \mbox{is the coproduct in $W$},\vspace{10pt}\\
\io^*: W\to \Fi &  \mbox{is the counit in $W$},\vspace{8pt}\\
\e^*: \Fi\to W &  \mbox{is the unit in $W$}.\\
\end{array}
$$

First, let us see that $\delta^*$ is indeed a product which agrees with the filtration. 
$W\ot W$ is a subspace of $A^*\ot A^*$ which, in turn, is a subspace of $(A\ot A)^*$. (The three spaces coincide if and only of $A$ is finite-dimensional.) We need to show that the image of the composition
$$W_k\ot W_l\hookrightarrow (A\ot A)^*\stackrel{\delta^*}{\longrightarrow} A^*$$
lies in $W_{k+l}$.

Take $w_1\in W_k$ and $w_2\in W_l$. The product of these elements is the composition
$$A\stackrel{\delta}{\longrightarrow}A\ot A\stackrel{w_1\ot w_2}{\longrightarrow}\Fi.$$
If $a\in A_{k+l+1}$ then 
$$\delta(a)=\sum_i b_i\ot c_i,$$ 
where for each $i$ we have $b_i\in A_{p}$ and $c_i\in A_{q}$ with $p+q= k+l+1$. As a consequence, either $p>k$ or $q>l$ which implies that $(w_1\ot w_2)(b_i\ot c_i)=0$ for all $i$ and, hence, $\delta^*(w_1\ot w_2)\in W_{k+l}$.

In order to see that $\mu^*$ gives a coproduct on $W$ which respects the filtration, we have to verify that the image of the map
$$W_k\hookrightarrow A^*\stackrel{\mu^*}{\longrightarrow} (A\ot A)^*$$
lies in $\sum_{p+q=k}W_p\ot W_q$. 

Take $w\in W_k$ and consider the composition
$$A\ot A\stackrel{\mu}{\longrightarrow} A\stackrel{w}{\longrightarrow}\Fi.$$
Since $w$ vanishes on $A_{k+1}$, the composition $w\circ\mu$ is equal to zero on the subspace 
$\sum_{p+q=k+1} A_p\ot A_q$ and thus may be considered as a linear function on the quotient vector space
$$A\ot A/\sum_{p+q=k+1} A_p\ot A_q.$$ Since the filtration $A_i$ is of finite type, this quotient 
does not change if we replace $A$ with the finite-dimensional vector space $A/A_{k+1}$. 
Now, for any finite-dimensional vector space $A$ with a descending filtration and for all $k$ the subspaces $$\left(A\ot A/\sum_{p+q=k+1} A_p\ot A_q\right)^*$$ and $$\sum_{p+q=k+1} (A/A_{p+1})^*\ot (A/A_{q+1})^*$$ of $A^*\ot A^*$ coincide. This implies that $\mu^*(w)\in \sum_{p+q=k}W_p\ot W_q$.

We leave checking the bialgebra axioms to the reader. 
\end{proof}

The fact the dual is defined only for filtered bialgebras of finite type and not for bialgebras in general is explained by the following observation. If the vector space $A$ is infinite-dimensional, the inclusion $$A^*\ot A^*\subset (A\ot A)^*$$
is strict. The dual to a coproduct $A\to A\ot A$ is a map $(A\ot A)^*\to A^*$ which restricts to a product $A^*\ot A^*\to A^*$, and, hence, the dual of a coalgebra is an algebra. However, the dual to a product on $A$ is a map $A^*\to A^*\ot A^*$, whose image does not necessarily lie in $A^*\ot A^*$. As a consequence, the dual of an algebra may fail to be a coalgebra.
\begin{xxca}
Give an example of a bialgebra whose product does not induce a coproduct on the dual space.
\end{xxca}

\subsection{Group-like and primitive elements in the dual bialgebra}

Primitive and group-like elements in the dual bialgebra have a
very transparent meaning.

\begin{xproposition}\label{prim_grlike_in_dual}
Primitive (respectively, group-like) elements in the dual of a filtered bialgebra
$A$ are those linear functions which are additive (respectively,
multiplicative), that is, satisfy the respective identities
\begin{align*}
  & f(ab)=f(a)+f(b),\\
  & f(ab)=f(a)f(b)
\end{align*}
for all $a,b\in A$.
\end{xproposition}

\begin{proof}
An element $f$ is primitive if $\d(f)=1\ot f+f\ot 1$.
Evaluating this on an arbitrary tensor product $a\ot b$ with $a,b\in
A$, we obtain
$$
  f(ab) = f(a)+f(b).
$$

An element $f$ is group-like if $\d(f)=f\ot f$. Evaluating
this on an arbitrary tensor product $a\ot b$, we obtain
$$
  f(ab) = f(a)f(b).
$$

In the same way, the additivity (multiplicativity) of a linear function implies that it defines a primitive (respectively, group-like) element. 
\end{proof}

\subsection{Hopf algebras}
\label{hopf_alg}

\begin{xdefinition}\index{Hopf algebra}
A {\em Hopf algebra} is a graded bialgebra.\index{Bialgebra!graded}
This means that $A$ is a graded vector space $A$, with the grading
by non-negative integers
$$
  A = \bigoplus_{k\ge 0} A_k
$$
and the grading is compatible with the operations $\mu$, $\io$, $\delta$, $\e$
in the following sense:
\begin{align*}
&\mu: A_m \ot A_n \to A_{m+n},\\
&\delta:  A_n \to \mathop{\bigoplus}\limits_{k+l=n}A_k \ot A_l.
\end{align*}
A Hopf algebra $A$ is said to be of {\em finite type},
\index{Hopf algebra!of finite type}%
if all its homogeneous components \index{Homogeneous components}
$A_n$ are finite-dimensional. A Hopf algebra is said to be {\em
connected}, \label{connected_alg}
\index{Bialgebra!connected}%
\index{Hopf algebra!connected}%
if 
$\io:\Fi\to A$ is an isomorphism of $\Fi$ onto $A_0\subset A$.
\end{xdefinition}

\begin{xremark} The above definition follows the classical
paper \cite{MiMo}. Nowadays a Hopf algebra is usually defined as a
not necessarily graded bialgebra with an additional operation,
called {\em antipode},\index{Antipode} which\label{antipode} is a
linear map $S:A\to A$ such that $$\mu\circ(S\ot 1)\circ\delta =
\mu\circ(1\ot S)\circ\delta = \io\circ\e.$$ The bialgebras of
interest for us (those that satisfy the premises of
Theorem~\ref{MMthm} below) always have an antipode.
\end{xremark}

\begin{xexample}
Any bialgebra can be considered as a Hopf algebra concentrated in
degree 0.
\end{xexample}
\begin{xexample}
An associated graded bialgebra of a filtered bialgebra is a Hopf algebra.
\end{xexample}
\begin{xexample}
Recall that, given a basis of a vector space $V$, the symmetric
algebra $S(V)$ is spanned by commutative monomials in the elements
of this basis. If $V$ is a graded vector space, and the basis is
chosen to consist of homogeneous elements of $V$, we define the
degree of a monomial to be the sum of the degrees of its factors.
With this grading $S(V)$ is a Hopf algebra.
\end{xexample}

\subsection{Dual Hopf algebra}\index{Hopf algebra! dual}
\label{dual_Ha} 

If $A$ be a Hopf algebra of finite type let $W_k=A_k^*$ and
$$W=\oplus_{k\geq 0} W_k.$$
The space $W$ is also a Hopf algebra; its operations are dual to
those of $A$:
$$\begin{array}{ll}
  \mu^*: W_n\to \mathop{\oplus}\limits_{k+l=n}\!\!\!
                      (A_k \ot A_l)^*\cong\!\!
            \mathop{\oplus}\limits_{k+l=n}\!\!\! W_k\ot W_l 
   & \mbox{is the coproduct in $W$}\vspace{10pt}\\
\delta^*: W_n\ot W_m \to W_{m+n} &
     \mbox{is the product in $W$}\vspace{8pt}\\
\io^*: W\to \Fi &  \mbox{is the counit in $W$}\vspace{8pt}\\
\e^*: \Fi\to W &  \mbox{is the unit in $W$}\\
\end{array}
$$
The Hopf algebra $W$ is called the {\em dual} of $A$. 
\begin{xxca}
Check that this definition is a particular case of the dual of a filtered bialgebra.
\end{xxca}
\begin{xxca}
Show that the dual of the dual of a Hopf algebra $A$ is canonically isomorphic to $A$.
\end{xxca}

\begin{xexample}
The bialgebra of $\Fi$-valued functions on a finite group $G$, considered
as a Hopf algebra concentrated in degree 0, is dual to the bialgebra
$\Fi G$.
\end{xexample}

\subsection{Structure theorem for Hopf algebras}

Is it easy to see that in a Hopf algebra the primitive subspace
$\PR=\PR(A)\subset A$ is the direct sum of its homogeneous
components: $\PR=\Op_{n\ge0}\PR\cap A_n$.

\begin{xtheorem}[Milnor--Moore \cite{MiMo}]
\label{MMthm}\index{Theorem!Milnor--Moore} Any commutative
cocommutative connected Hopf algebra of finite type is canonically
isomorphic to the symmetric algebra on its primitive subspace:
$$
  A = S(\PR(A)).
$$
This isomorphism sends a polynomial in the primitive elements of $A$
into its value in $A$.
\end{xtheorem}

In other words, if a linear basis is chosen in every homogeneous
component $\PR_n=\PR\cap A_n$, then each element of $A$ can be
written uniquely as a polynomial in these variables.

\begin{xremark}
In \cite{MiMo} Milnor and Moore actually consider {\em graded
commutative} Hopf algebras. We do not need this level of generality.
\end{xremark}

\begin{proof}
There are two assertions to prove:
\begin{itemize}
\item[(1)] every element of $A$ can expressed a a polynomial, that is,
as a sum of products, of primitive
elements;
\item[(2)] the value of a nonzero polynomial on a set of linearly independent
homogeneous primitive elements cannot vanish in $A$.
\end{itemize}
First, let us prove assertion (1) for the homogeneous elements of
$A$ by induction on their degree.

Note that under our assumptions the coproduct of a homogeneous
element $x\in A_n$ has the form
\begin{equation}\label{delta}
  \delta(x)=1\ot x +\dots + x\ot 1,
\end{equation}
where the dots stand for an element of
$A_1\ot A_{n-1}+\dots + A_{n-1}\ot A_1$.
Indeed, we can always write $\delta(x) = 1\ot y+\dots+ z\ot 1$.
By cocommutativity $y=z$.
Then, $x=(\e\ot\id)(\delta(x))=y+0+\dots+0=y$.

In particular, for any element $x\in A_1$ equation (\ref{delta})
ensures that $\delta(x)=1\ot x + x\ot 1$, so that $A_1=\PR_1$. (It
may happen that $A_1=0$, but this does not interfere the subsequent
argument!)

Take an element $x\in A_2$. We have
$$
  \delta(x) = 1\ot x +\sum\lambda_{ij}p^1_i\ot p^1_j + x\ot 1,
$$
where $p^1_i$ constitute a basis of $A_1=\PR_1$ and $\lambda_{ij}$
is a symmetric matrix over the ground field. Let
$$
  x'=\frac12\sum\lambda_{ij}p^1_i p^1_j.
$$
Then
$$
  \delta(x') = 1\ot x' +\sum\lambda_{ij}p^1_i\ot p^1_j + x'\ot 1.
$$
It follows that
$$
  \delta(x-x') = 1\ot x' + x'\ot 1,
$$
that is, $x-x'$ is primitive, and $x$ is expressed via primitive
elements as $(x-x')+x'$, which is a polynomial, linear in $\PR_2$
and quadratic in $\PR_1$.

Proceeding in this way, assertion (1) can be proved in degrees 3, 4,
and so on. We omit the formal inductive argument.

Now, assume that there exists a polynomial in the basis elements of
$\PR(A)$ which is equal to zero in $A$. Among all such expressions
there exists one, which we denote by $w$, of the smallest degree; we
can assume that it is homogeneous (lies in $A_n$ for some $n$). In
particular, all monomials of degree smaller than $n$ are linearly
independent. (We remind the reader that we are working in a graded
algebra, so the degree of a polynomial is calculated taking into the
account the degrees of the variables. In particular, a linear
monomial has the degree equal to the degree of the corresponding
variable.)

Let $a$ be a basis primitive element which appears in $w$ as a
factor in at least one of the summands; we can write
$$w=a^k f_k + a^{k-1} f_{k-1}+\ldots f_0,$$
where the $f_i$ for $i>0$ are polynomials in the basis primitive
elements of degree smaller than $n$. Now, $\delta(w)-(1\ot w+ w\ot
1)$ lies in the sum of the terms $A_p\ot A_q$ with $p+q=n$ and
$p,q>0$; the sum of these terms has a basis consisting of
expressions $m_p\ot m_q$ where $m_p,m_q$ are monomials in the basis
elements of $\PR(A)$ of degrees $p,q$ respectively. Inspection shows
that the terms in $\delta(w)-(1\ot w+ w\ot 1)$ with $m_p=a^k$ add up
to $a^k\ot f_k$. Since $\delta(w)$ must be zero in $A\ot A$ this
implies that $f_k=0$ in $A$, which gives a contradiction since the
degree of $f_k$ is smaller than $n$.

This completes the proof.
\end{proof}

\begin{xcorollary}

An algebra $A$ satisfying the assumptions of the theorem

\begin{itemize}
\item[(1)] has no zero divisors,
\item[(2)] has the antipode $S$ defined on primitive elements by
$$
  S(p)=-p\ .
$$
\item[(3)] splits as a direct sum of vector spaces  $A_k=P_k\op R_k$, where $R_k$ is spanned by products of elements of non-zero degrees. 
\end{itemize}
\end{xcorollary}
The space $R=\op R_k$ is called the space of {\em decomposable} elements.

\begin{xremark}
In fact, there is a more general version of the Milnor-Moore theorem
which describes the structure of a general cocommutative connected
Hopf algebra. The primitive subspace of a Hopf algebra is a {\em
graded} Lie algebra; a cocommutative connected Hopf algebra $A$ is
canonically isomorphic to the universal enveloping algebra of
$\PR(A)$.
\end{xremark}

\subsection{Primitive and group-like elements in Hopf algebras}
\label{prim_grl}

As the Milnor-Moore theorem shows, a non-trivial connected
cocommutative Hopf algebra always has a non-empty primitive
subspace. However, the only group-like element in such a Hopf
algebra is the identity. (In the case of commutative algebras, which
are all isomorphic to symmetric algebras, this was noted in
Section~\ref{ap:pr-gr-el}.) As we shall now see, all these Hopf
algebras acquire a wealth of group-like elements {\em after
completion.}

Let $\widehat{A}$ be the graded completion of a Hopf algebra $A$. We
remind that while any element of $A$ can be written as a finite sum
$\sum_{i<N} x_i$ with $x_i\in A_i$, elements of $\widehat{A}$ are
represented by infinite sums $\sum_{i} x_i$ with $x_i\in A_i$. The
operations on $A$ extend to $\widehat{A}$ uniquely; note, however,
that a priori $\widehat{A}$ comes with no non-trivial grading.

\begin{xlemma}\label{quillen}
For the graded completion $\widehat{A}$ of a connected Hopf algebra
$A$ the functions $\exp$ and $\log$, defined by the usual power
series, establish a one-to-one correspondence between the set of
primitive elements $\PR(\widehat{A})$ and the set of group-like
elements $\GR(\widehat{A})$.
\end{xlemma}

\begin{proof}
Let $p\in \PR(\widehat{A})$. Then
$$
  \d(p^n)=(1\ot p+p\ot 1)^n
=\sum_{k+l=n}\frac{n!}{k!l!} p^k\ot p^l
$$
and therefore
$$
  \d(e^p)=\d\left(\sum_{n=0}^\infty \frac{p^n}{n!}\right)
         =\sum_{k=0}^\infty\sum_{l=0}^\infty\frac{1}{k!l!}p^k\ot p^l
=\sum_{k=0}^\infty\frac{1}{k!}p^k\ot\sum_{l=0}^\infty\frac{1}{l!}p^l
=e^p\ot\ e^p
$$
which means that $e^p\in \GR(\widehat{A})$.

Vice versa, assuming that $g\in \GR(\widehat{A})$ we want to prove
that $\log(g)\in \PR(\widehat{A})$. By assumption, our Hopf algebra
$A$ is connected which implies that the graded component $g_0\in
A_0\cong\Fi$ is equal to 1. Therefore we can write $g=1+h$ where
$h\in\prod_{k>0}A_k$. The condition that $g$ is group-like
transcribes as
\begin{equation}
  \d(h)=1\ot h + h\ot 1+ h \ot h. \label{prim_h}
\end{equation}
Now,
$$
  p=\log(g) = \log(1+h) = \sum_{k=1}^\infty \frac{(-1)^{k-1}}{k} h^k
$$
and an exercise in power series combinatorics
shows that equation \ref{prim_h} implies the required property
$$
  \d(p)=1\ot p + p\ot 1.
$$
\end{proof}

\begin{xca}[\cite{Sch, Lnd1}]\label{primitivization}
Define the convolution product of two vector space endomorphisms of a graded connected commutative and co-commutative
Hopf algebra $A$ by
$$
 (f*g)(a) = \sum_{\d(a)=\sum a'_i\ot a''_i} f(a'_i)g(a''_i).
$$
Let $I:A\to A$ be the operator defined as zero on  $A_0$ and as the identity on each $A_i$ with $i>0$.
Show that the map
$$I-\frac{1}{2} I*I + \frac{1}{3} I*I*I -\frac{1}{4} I*I*I*I +\ldots$$ 
is the projector of $A$ onto the subspace
of primitives $P$ parallel to the subspace $R$ of decomposable elements.
\end{xca}

\section{Free algebras and free Lie algebras}

Here we briefly mention the definitions and basic properties of the
the free associative and free Lie algebras. For a detailed treatment
see, for example, \cite{Reu}.

\subsection{Free algebras}

The {\em free algebra} $\Ring\langle x_1, \ldots, x_n \rangle$ over
a commutative unital ring $\Ring$ is the associative algebra of
non-commutative polynomials in the $x_i$ with coefficients in
$\Ring$. If $\Ring=\Fi$ and $V$ is the vector space spanned by the
symbols $x_1,\ldots, x_n$ the free algebra on the $x_i$ is
isomorphic to the tensor algebra $T(V)$.

\begin{xexample}
The algebra $\Ring\langle x_1,x_2\rangle$ consists of finite linear
combinations of the form
$c+c_1^{}x_1^{}+c_2^{}x_2^{}+c_{11}^{}x_1^2+c_{12}^{}x_1^{}x_2^{}
+c_{21}^{}x_2^{}x_1^{}+c_{22}^{}x_2^2+\ldots$, $c_{\a}\in\Ring$,
with natural addition and multiplication.
\end{xexample}

The free algebra $\Ring\langle x_1, \ldots, x_n \rangle$ is
characterized by the following universal property: given an
$\Ring$-algebra $A$ and a set of elements $a_1,\ldots, a_n$ in $A$
there exists a unique map
$$\Ring\langle x_1, \ldots, x_n \rangle \to A$$ which sends $x_i$ to $a_i$
for all $i$. As a consequence, every $\Ring$-algebra generated by
$n$ elements is a quotient of the free algebra $\Ring\langle x_1,
\ldots, x_n \rangle$.

The word {\em free} refers to the above universal property, which is
analogous to the universal property of free groups or Lie algebras
(see below). This property amounts to the fact that the only
identities that hold in $\Ring\langle x_1, \ldots, x_n \rangle$ are
those that follow from the axioms of an algebra, such as
$(x_1^{}+x_2^{})^2=x_1^2+x_1^{}x_2^{}+x_2^{}x_1^{}+x_2^2$.

The algebra $\Ring\langle x_1, \ldots, x_n \rangle$ is graded by the
degree of the monomials; its homogeneous component of degree $k$ has
dimension $n^k$, and its Poincar\'e series is $1/(1-nt)$. The graded
completion of $\Ring\langle x_1, \ldots, x_n \rangle$ is denoted by
$\Ring\la x_1, \ldots, x_n \ra$.

The free algebra $\Ring\langle x_1, \ldots, x_n \rangle$ has a
coproduct $\delta$ defined by the condition that the generators
$x_i$ are primitive:
$$\delta(x_i)=x_i\ot 1+1\ot x_i.$$
This condition determines $\delta$ completely since the coproduct is
an algebra homomorphism. There also exists a counit: it sends a
non-commutative polynomial to its constant term.
\begin{xproposition}
The free algebra $\Fi\langle x_1, \ldots, x_n \rangle$ is a connected
cocommutative Hopf algebra.
\end{xproposition}
The proof is immediate.

\subsection{Free Lie algebras}

Recall that the space of primitive elements in a bialgebra is a Lie
algebra whose Lie bracket is the algebra commutator $[a,b]=ab-ba$.
Let $L(x_1,\ldots, x_n)$ be the Lie algebra of primitive elements in
$\Fi\langle x_1, \ldots, x_n \rangle$. Note that the $x_i$ belong to
$L(x_1,\ldots, x_n)$.

\begin{xproposition}
The Lie algebra $L(x_1,\ldots, x_n)$  has the following universal
property: given a Lie algebra $\g$ and a set of elements $a_1,
\ldots, a_n\in \g$ there exists the unique Lie algebra homomorphism
$L(x_1, \ldots, x_n)\to \g$ sending each $x_i$ to $a_i$.
\end{xproposition}
Indeed, since $\Fi\langle x_1, \ldots, x_n \rangle$ is free, there
exists a unique algebra homomorphism
$$\Fi\langle x_1, \ldots, x_n \rangle\to U(\g)$$
sending the $x_i$ to the $a_i$. Passing to the primitive spaces we
recover the Proposition.

\begin{xdefinition}
The Lie algebra $L(x_1,\ldots, x_n)$  is called {\em the free Lie
algebra on $x_1,\ldots, x_n$.}
\end{xdefinition}

The explicit construction of $L(x_1,\ldots, x_n)$ uses {\em Lie
monomials}, which are defined inductively as follows. A Lie monomial
of degree $1$ in $x_1,\ldots, x_n$ is simply one of these symbols. A
Lie monomial of degree $d$ is an expression of the form $[a,b]$
where $a$ and $b$ are Lie monomials the sum of whose degrees is $d$.

The Lie algebra $L(x_1,\ldots, x_n)$ as a vector space is spanned by
all Lie monomials in $x_1,\ldots, x_n$, modulo the subspace spanned
by all expressions of the form
$$[a,b]-[b,a]$$
and $$[[a,b],c]+[[b,c],a]+[[c,a],b]$$ where $a,b,c$ are Lie
monomials. The Lie bracket is the linear extension of the operation
$[\ ,\ ]$ on Lie monomials. Note that as a vector space a free Lie
algebra is graded by the degree of Lie monomials. Understanding the
bracket as the commutator we get an embedding of $L(x_1,\ldots,
x_n)$ constructed in this way into $\Fi\langle x_1, \ldots, x_n
\rangle$ as the primitive subspace.

Finding a good basis for a free Lie algebra is a non-trivial
problem; it is discussed in detail in \cite{Reu}. 
One explicit basis, the so-called {\em Lyndon basis}, is constructed with the help of
{\em Lyndon words}. The Lyndon words can be defined as follows. 
Take an aperiodic necklace (see page \pageref{necklace}) and choose the lexicographically smallest among all its cyclic shifts. Replacing each bead with the label $i$ by $x_i$ we get a non-commutative monomial  (Lyndon word) in the $x_i$. A Lyndon word $w$ gives rise to an iterated commutator by means of the following recurrent procedure. First,  $w=x_ix_j$ is declared to produce the commutator $[x_i,x_j]$. If $w$ is of degree more than two, among all decompositions of $w$ into a nontrivial product $w=uv$ choose the decomposition with  lexicographically the smallest possible $v$, and take the commutator of the (possibly iterated) commutators that correspond to $u$ and $v$.

Shown below is the Lyndon basis for the free Lie algebra  $L(x,y)$ in small degrees.

$$\begin{array}{ccl}
m & \dim L(x,y)_m & \text{basis} \\
1 &   2     &      x,y \\
2 &   1     &      [x,y] \\
3 &   2     &      [x,[x,y]]\quad [y,[x,y]] \\
4 &   3     &      [x,[x,[x,[x,y]]]]\quad [y,[x,[x,[x,y]]]]\quad [y,[y,[x,[x,y]]]] \\
5 &   6     &      [x,[x,[x,[x,y]]]]\quad [y,[x,[x,[x,y]]]]\quad [y,[y,[x,[x,y]]]] \\
  &         &      [y,[y,[y,[x,y]]]]\quad [[x,y],[x,[x,y]]]\quad [[x,y],[y,[x,y]]]
\end{array}$$

 % appendix

\bibliographystyle{amsalpha}

\chapter*{Notations}

$\Z$, $\Q$, $\R$, $\C$ --- rings of integer, rational, real and
  complex numbers.

$\A$ --- algebra of unframed chord diagrams on the circle,
 p.\pageref{A}.

$\A^{fr}$ --- algebra of framed chord diagrams on the circle,
 p.\pageref{A^{fr}}.

$\A_n$ --- space of unframed chord diagrams of degree $n$,
 p.\pageref{A_n}.

$\A^{fr}_n$ --- space of framed chord diagrams of degree $n$,
 p.\pageref{A_n^{fr}}.

$\A(n)$ --- algebra of chord diagrams on $n$ lines, p.\pageref{A-ot-p}.

%$\A(p,q)$ --- algebra of chord diagrams on $p$ lines and $q$ circles,
%  p.\pageref{A-ot-pq}.

$\A^h(n)$ --- algebra of horizontal chord diagrams, p.\pageref{A^h(n)}.

$\widehat{\A}$ --- graded completion of the algebra of chord diagrams,
  p.\pageref{hat-A}.

$\ChD_n$ --- set of chord diagrams of degree $n$, p.\pageref{ChD_n}.

$A$ --- Alexander-Conway power series invariant,
  p.\pageref{Alex-Conv-pol}.

$\alpha_n$ --- map from $\V_n$ to $\Ring\ChD_n$, symbol of an
  invariant, p.\pageref{alpha_n}.

$\B$ --- algebra of open Jacobi diagrams, p.\pageref{alg-B}.

$\B(m)$ --- space of $m$-coloured open Jacobi diagrams,
  p.\pageref{B-of-m}.

$\Bbig$ --- enlarged algebra $\B$, p.\pageref{B-big}.

$\OD_n$ --- set of open Jacobi diagrams of degree $n$,
  p.\pageref{OD_n}.

$BNG$ --- the Bar-Natan--Garoufalidis function, p.\pageref{BNG}.

$\F$ --- space of closed Jacobi diagrams, p.\pageref{F_n}.

$\F_n$ --- space of closed Jacobi diagrams of degree $n$, p.\pageref{F_n}.

${\mathcal C}_n$ ---  Goussarov--Habiro moves, p.\pageref{gg-moveCn}.

$\F(\xx_1,\ldots,\xx_n\, |\,\yy_1,\ldots,\yy_m)$ --- space of mixed
   Jacobi diagrams, p.\pageref{sp-mix-J-d}.

$\CP$ --- Conway polynomial, p.\pageref{axiconw}.

$\kr_{2n}$ --- Conway combination of Gauss diagrams, p.\pageref{def:cc}.

$\FD_n$ --- set of closed diagrams of degree $n$, p.\pageref{FD_n}.

$c_n$ --- $n$-th coefficient of the Conway polynomial,
 p.\pageref{c_n-Conw-coef}.

%$\ldc_n$ --- $n$-th disconnected cabling of a knot, p.\pageref{cabl-ldc}.

%$\lcc_n$ --- $n$-th connected cabling of a knot, p.\pageref{cabl-lcc}.

$\partial_C$ --- diagrammatic differential operator on $\B$,
  p.\pageref{d-d-o}.

$\partbig_C$ --- diagrammatic differential operator on $\Bbig$,
  p.\pageref{d-d-o-O}.

$\partial_{\Omega}$ --- wheeling map, p.\pageref{wheel-map}.

$\delta$ --- coproduct in a coalgebra, p.\pageref{Coalgebra};
in particular, for the bialgebra $\A^{fr}$ see p.\pageref{cmult-in-A}.

$\Delta^{n_1,\ldots,n_k}$ --- operation
$\A(k)\to\A(n_1+\ldots+n_k)$, p.\pageref{Delta_n1nk}.

$\e$ --- counit in a coalgebra, p.\pageref{Coalgebra};
in particular, for the bialgebra $\A^{fr}$ see p.\pageref{u-cu-in-A}.

$F(L)$ --- unframed two-variable Kauffman polynomial, p.\pageref{t-var-Kauf}.

$\Phi$ --- map $\B({\yy})\to \F(\xx)$, p.\pageref{Phi:B-to-F}.

$\Phi_0$ --- map $\B\to\F$, p.\pageref{Phi-0}.

$\Phi_2$ --- map $\B({\yy}_1,{\yy}_2)\to \F(\xx)$, p.\pageref{Phi-2}.

$\f_\g$ --- universal Lie algebra weight system, p.\pageref{def:phiA}.

$\f_\g^T$ --- Lie algebra weight system associated with the representation,
  p.\pageref{def:pr-phiA}.

${\mathcal G}_n$ --- Goussarov group, p.\pageref{gg}.

$\Gr$ --- bialgebra of graphs, p.\pageref{Gr}.

$\G$ --- algebra of 3-graphs, p.\pageref{alg-3-graphs}.

$\G(D)$ --- intersection graph of a chord diagram $D$,
 p.\pageref{G(D)}.

$H=\left(\begin{array}{cc} 1&0\\ 0&-1\end{array}\right)$ ---
 element of the Lie algebra $\sL_2$, p.\pageref{ws_sl_2_on_A}.

$H$ --- hump unknot, p.\pageref{hump}.

$\bo_n$ --- constant 1 weight system on $\A_n$, p.\pageref{bo_n}.

$I$ --- a map of Gauss diagrams to arrow diagrams,
 p.\pageref{map-I-gd-ad}.

$I(K)$ --- final Kontsevich integral, p.\pageref{ki-I(K)}.

$\I$ --- algebra of knot invariants, p.\pageref{alg_inv}.

$\io$ --- unit in an algebra, p.\pageref{Algebra};
in particular, for the bialgebra $\A^{fr}$ see p.\pageref{u-cu-in-A}.

%$\iota_{pq}$ --- map from $\A(p,q)$ to $\A(p-1,q+1)$,
% p.\pageref{iota-pq}.

$j_n$ --- $n$-th coefficient of the modified Jones polynomial,
 p.\pageref{jones_vi}.

$\K$ --- set of (equivalence classes of) knots, p.\pageref{set-of-knots}.

$\Li$ --- Euler dilogarithm, p.\pageref{E-Li}.

$\La$ --- bialgebra of Lando, p.\pageref{La}.

$\Lambda$ --- Vogel's algebra, p.\pageref{Vogel-algebra}.

$\Lambda(L)$ --- framed two-variable Kauffman polynomial,
 p.\pageref{t-var-Kauf}.

$\nabla$ --- difference operator for Vassiliev invariants,
 p.\pageref{nabla}.

$\M_n$ --- Goussarov--Habiro moves, p.\pageref{GHthm}.

$\MM$ --- highest order part of the coloured Jones polynomial,
 p.\pageref{mel-mor}.

$M_T$ --- mutation of a knot with respect to a tangle $T$,
 p.\pageref{mutkn}.

$\mu$ --- product in an algebra, p.\pageref{Algebra};
in particular, for the bialgebra $\A^{fr}$ see p.\pageref{mult-in-A}.

$\Polyak$ --- Polyak algebra, p.\pageref{PolyakA}.

$\PR_n$ --- primitive subspace of the algebra of chord diagrams,
 p.\pageref{PR_n}.

$P$ --- HOMFLY polynomial, p.\pageref{HOMFLY-def}.

$P^{fr}$ --- framed HOMFLY polynomial, p.\pageref{fr-HOMFLY}.

$p_{k,l}(L)$ --- $k,l$-th coefficient of the modified HOMFLY polynomial, , p.\pageref{pr:vas-homfly}.

$\rho_\g$ ---  universal Lie algebra weight system on $\B$, p.\pageref{def:psiB}.

$\psi_n$ --- $n$-th cabling of a chord diagram, p.\pageref{psiop}.

$\Ring$ --- ground ring (usually $\Q$ or $\C$), p.\pageref{Ring}.

$\Ring(\ChD_n)$ --- $\Ring$-valued functions on chord diagrams,
 p.\pageref{Ring(ChD_n)}.

$R$ --- $R$-matrix, p.\pageref{R-matrix}.

$R$, $R^{-1}$ --- Kontsevich integrals of two braided strings,
  p.\pageref{ki-R}.

$S_A$ --- symbol of the Alexander-Conway invariant $A$,
 p.\pageref{S_A}.

$S_{MM}$ --- symbol of the Melvin-Morton invariant $MM$,
 p.\pageref{S_MM}.

$S_i$ --- operation on tangle (chord) diagrams, p.\pageref{op-S_i}.

$\symb(v)$ --- symbol of the Vasiliev invariant $v$,
 p.\pageref{Symbol}.

$\sigma$ --- mirror reflection of knots, p.\pageref{inv-s-t}.

$\tau$ --- changing the orientation of a knot, p.\pageref{inv-s-t}.

$\tau$ --- reversing the orientation of the Wilson loop,
   p.\pageref{def:tau-rever}.

$\tau$ --- inverse of $\chi:\B\to\F$, p.\pageref{inv-chi}.

$\Theta$ --- the chord diagram with one chord,
    $\risS{-10}{cd1ch4}{}{25}{0}{0}$, p.\pageref{Theta}.

$\t^{fr}$ ---- quantum invariant, p.\pageref{qis4}.

$\t^{fr}$ ---- $\sL_2$-quantum invariant, p.\pageref{qis6}.

$\t^{fr,St}_{\sL_N}$---- $\sL_N$-quantum invariant, p.\pageref{q-sl_N}.

$\V$ --- space of Vassiliev (finite type) invariants,
 p.\pageref{d:Vspace}

$\V_n$ --- space of unframed Vassiliev knot invariants of degree
    $\le n$, p.\pageref{d:Vn}.

$\V^{fr}_n$ --- space of framed Vassiliev knot invariants of degree
   $\le n$, p.\pageref{fr-vi}.

$\V_{\bullet}$ --- space of polynomial Vassiliev invariants,
 p.\pageref{pol-Vas-inv}.

$\widehat{\V}_{\bullet}$ --- space of power series invariants,
   graded completion of $\V_{\bullet}$, p.\pageref{V-grad_compl}.

$\W_n$ --- space of unframed weight systems of degree $n$,
 p.\pageref{W_n}.

$\W^{fr}_n$ --- space of framed weight systems of degree $n$,
 p.\pageref{W_n^{fr}}.

$\widehat{\W}^{fr}$ --- graded completion of the algebra of weight systems,
  p.\pageref{grad_compl}.

$\Zed(\hl)$ --- Kontsevich integral of $\hl$ in algebra $\B(\yy)$,
 p.\pageref{Zed-of-hl}.

$\Zed(\dhl)$ --- Kontsevich integral of $\dhl$ in algebra
 $\B({\yy}_1,{\yy}_2)$, p.\pageref{Zed-of-dhl}.

$\Zed_i(\hl)$ --- $i$-th part of the Kontsevich integral $\Zed(\hl)$,
 p.\pageref{Zed-i-of-hl}.

$Z(K)$ --- Kontsevich integral, p.\pageref{ki-Z(K)}.

$\Z\K$ --- algebra of knots, p.\pageref{alg_knots}.

%\bigskip

%$\Delta_n$

$\chi$ --- symmetrization map $\B\to\F$, p.\pageref{BFisom}.

$\chi_{\yy_m}$ --- map
  $\F(\boldX\, |\, \yy_1,\ldots,\yy_m)\to
  \F(\boldX, \yy_m\, |\, \yy_1,\ldots,\yy_{m-1})$,
  p.\pageref{sym-map-C}.

$\Omega'$ --- part of $\Zed_0(\hl)$ containing wheels,
  p.\pageref{omega-pr-def}.

$\langle\ ,\, \rangle_{\yy}$ --- pairing
  $\Fxy\otimes\B(\yy)\to\F(\xx)$, p.\pageref{pair-C-B}.

$\hl$ --- open Hopf link, p.\pageref{Hopf-hl}.

$\dhl$ --- doubled open Hopf link, p.\pageref{Hopf-dhl}.

$\clhl$ --- closed Hopf link, p.\pageref{Hopf-clhl}.

$\#$ --- connected sum of two knots, p.\pageref{connectedsum} or two diagrams, p.\pageref{connectedsumdiags}; also the action of $\F$ on tangle diagrams, p.\pageref{subsec:modules-tangle}. 

\bigskip

-----------------------

%Add index of open problems
 % notations

\printindex

\end{document}